\pdfoutput=1  

\documentclass[11pt, openany]{book}
\usepackage[margin=.8in,left=.8in]{geometry}
\usepackage{amsthm}
\usepackage{amssymb}
\usepackage{amsmath}
\usepackage{mathtools}
\usepackage{xcolor}
\usepackage{colortbl}    
\usepackage{stmaryrd}    
\usepackage[safe]{tipa} 
\usepackage{adjustbox} 

\usepackage{hyperref}
\hypersetup{
  colorlinks = true,
  allcolors=.           
}

\usepackage{tikz}
\usetikzlibrary{cd}

\DeclareMathAlphabet{\mathpzc}{OT1}{pzc}{m}{it} 
\newcommand\mathscr[1]{\scalebox{1.1}{$\mathpzc{#1}$}}

\newcommand{\underoverset}[3]{\underset{#1}{\overset{#2}{#3}}}

\definecolor{darkblue}{rgb}{0.05,0.25,0.65}
\definecolor{greenii}{RGB}{20,140,10}
\definecolor{lightgray}{rgb}{0.9,0.9,0.9}
\definecolor{orangeii}{RGB}{200,100,5}

\usepackage{multirow}

\usepackage{enumerate} 

\makeatletter
\newcommand\makebig[2]{%
  \@xp\newcommand\@xp*\csname#1\endcsname{\bBigg@{#2}}%
  \@xp\newcommand\@xp*\csname#1l\endcsname{\@xp\mathopen\csname#1\endcsname}%
  \@xp\newcommand\@xp*\csname#1r\endcsname{\@xp\mathclose\csname#1\endcsname}%
}
\makeatother

\makebig{biggg} {3.0}
\makebig{Biggg} {3.5}
\makebig{bigggg}{4.0}
\makebig{Bigggg}{4.5}

\usepackage{mathptmx}
\usepackage{graphicx}
\DeclareRobustCommand{\coprod}{\mathop{\text{\fakecoprod}}}
\newcommand{\fakecoprod}{%
  \sbox0{$\prod$}%
  \smash{\raisebox{\dimexpr.9625\depth-\dp0}{\scalebox{1}[-1]{$\prod$}}}%
  \vphantom{$\prod$}%
}

\def\conical{\rotatebox[origin=c]{70}{$<$}}
\def\smooth{\rotatebox[origin=c]{70}{$\subset$}}
\def\orbisingular{\rotatebox[origin=c]{70}{$\prec$}}

\def\orbisingularG{\raisebox{-3pt}{$\orbisingular^{\hspace{-5.7pt}\raisebox{2pt}{\scalebox{.83}{$G$}}}$}}

\def\orbisingularE{\raisebox{-3pt}{$\orbisingular^{\hspace{-5.5pt}\raisebox{1pt}{\scalebox{.83}{$1$}}}$}}

\def\orbisingularH{\raisebox{-3pt}{$\orbisingular^{\hspace{-5.7pt}\raisebox{2pt}{\scalebox{.83}{\hspace{1pt}$H$}}}$}}
\def\orbisingularK{\raisebox{-3pt}{$\orbisingular^{\hspace{-5.7pt}\raisebox{2pt}{\scalebox{.83}{\hspace{1pt}$K$}}}$}}

\def\orbisingulari{\raisebox{-3pt}{$\orbisingular^{\hspace{-4.2pt}\raisebox{1.3pt}{\scalebox{.83}{$i$}}}$}}

\newcommand{\orbisingularAny}[1]{
  \raisebox{-3pt}{$\orbisingular^{\hspace{-5.7pt}\raisebox{2pt}{\scalebox{.83}{$#1$}}}$}
}

\newcommand{\conicalrelativeG}{
  \conical_{\scalebox{.5}{$\!\!\!\!\orbisingularG$}}
}

\newcommand{\smoothrelativeG}{
  \smooth_{\scalebox{.5}{$\!\!\!\orbisingularG$}}
}

\newcommand{\orbisingularrelativeG}{
  \orbisingular_{\scalebox{.5}{$\!\!\!\!\!\!\!\orbisingularG$}}
}

\newcommand{\Singularities}{
  \mathrm{Snglrt}
}

\newcommand{\Resolvable}{
  \mathrm{rslvbl}
}

\newcommand{\resolvable}{
  \Resolvable
}

\newcommand{\CoverResolvable}{
  {\mathrm{cov}\Resolvable}
}

\newcommand{\ResolvableSingularities}{
  \Resolvable\Singularities
}

\newcommand{\CoverResolvableSingularities}{
  \CoverResolvable\Singularities
}

\newcommand{\stable}{
  \mathrm{stbl}
}

\newcommand{\conjugation}{
  \mathrm{cjg}
}

\newcommand{\Charts}{
  \mathrm{Chrt}
}

\newcommand{\Singular}{
  \mathrm{Snglr}
}

\newcommand{\SingularInfinityGroupoids}{
  \Singular\InfinityGroupoids
}

\newcommand{\SingularSmoothInfinityGroupoids}{
  \Singular\SmoothInfinityGroupoids
}

\newcommand{\Sets}{
  \mathrm{Set}
}

\newcommand{\FiniteSets}{
  \mathrm{Fin}\Sets
}

\newcommand{\underlying}{
  \mathrm{undrlng}
}

\newcommand{\Spaces}{
  \mathrm{Spc}
}

\newcommand{\SimplicialSets}{
  \simplicial\Sets
}

\newcommand{\SimplicialNerve}{
  N
}

\newcommand{\SimplicialGroups}{
  \Groups(\SimplicialSets)
}

\newcommand{\Groups}{
  \mathrm{Grp}
}

\newcommand{\EquivariantClassifyingShape}[2]{
  B_{\scalebox{.6}{$#1$}} #2
}

\newcommand{\NeutralElement}{
  \mathrm{e}
}

\newcommand{\good}{
  \mathrm{good}
}

\newcommand{\wellpointed}{
  \!\mathrm{wellpt}
}

\newcommand{\AbelianGroups}{
  \mathrm{Ab}\Groups
}

\newcommand{\Extensions}[1]{
  {#1}\,\mathrm{Extnss}
}

\newcommand{\SplitExtensions}[1]{
  \mathrm{Splt}\,\Extensions{#1}
}

\newcommand{\TwoGroups}{
  2\Groups
}

\newcommand{\UnitaryGroup}{
  \mathrm{U}
}

\newcommand{\ProjectiveUnitaryGroup}{
  \mathrm{P}\UnitaryGroup
}

\newcommand{\OrthogonalGroup}{
  \mathrm{O}
}

\newcommand{\HilbertSpace}{
  \mathcal{H}
}

\newcommand{\HilbertSpaces}{
  \mathrm{Hilb}
}

\newcommand{\BoundedOperators}{
  \mathcal{B}
}

\newcommand{\CompactOperators}{
  \mathcal{K}
}

\newcommand{\FredholmOperator}{
  F
}

\newcommand{\FredholmOperators}{
  \mathrm{Fred}^{(0)}
}

\newcommand{\UnitaryOperator}{
  U
}

\newcommand{\UH}{
  \UnitaryGroup_{\!\omega}
}

\newcommand{\UHPrime}{
  \UnitaryGroup_{\!\!\omega'}
}

\newcommand{\UHPrimeOptional}{
  \UnitaryGroup_{\!\!\omega^{(\prime)}}
}

\newcommand{\GradedUH}{
  \UH^{\mathrm{gr}}
}

\newcommand{\GradedUHPrime}{
  \UHPrime^{\mathrm{gr}}
}

\newcommand{\OH}{
  \OrthogonalGroup_{\!\omega}
}

\newcommand{\PUH}{
  \mathrm{P}\UH
}

\newcommand{\PUHPrime}{
  \mathrm{P}\UHPrime
}

\newcommand{\PUHPrimeOptional}{
  \mathrm{P}\UHPrimeOptional
}

\newcommand{\POH}{
  \mathrm{P}\OH
}

\newcommand{\GradedOH}{
  \OH^{\mathrm{gr}}
}

\newcommand{\GradedPUH}{
  \PUH^{\mathrm{gr}}
}

\newcommand{\GradedPUHPrime}{
  \PUHPrime^{\mathrm{gr}}
}

\newcommand{\GradedPOH}{
  \POH^{\mathrm{gr}}
}

\newcommand{\KU}{
  \mathrm{KU}
}

\newcommand{\KO}{
  \mathrm{KO}
}

\newcommand{\KR}{
  \mathrm{KR}
}

\newcommand{\CircleGroup}{
  {\UnitaryGroup_{\!1}}
}

\newcommand{\PhaseGroup}{
  \UnitaryGroup_{\!1'}
}

\newcommand{\SOThree}{
  \mathrm{SO}(3)
}

\newcommand{\SUTwo}{
  \mathrm{SU}(2)
}

\newcommand{\SpOne}{
  \mathrm{Sp}(1)
}

\newcommand{\SpinGroup}{
  \mathrm{Spin}
}

\newcommand{\CyclicGroup}[1]{
  \Integers_{
    \scalebox{.7}{
      \hspace{-4pt}
      \scalebox{.9}{/}
      \hspace{-5pt}
      {$#1$}
      \hspace{-5pt}
    }
  }
}

\newcommand{\ZTwo}{
  \CyclicGroup{2}
}

\newcommand{\Homomorphisms}{
  \mathrm{Hom}
}

\newcommand{\CrossedHomomorphisms}{
  \mathrm{Crs}\Homomorphisms
}

\newcommand{\Functors}{
  \mathrm{Fnctr}
}

\newcommand{\ProObjects}{
  \mathrm{Pro}
}

\newcommand{\Hausdorff}{
  \mathrm{Haus}
}

\newcommand{\topological}{
  \mathrm{Top}
}

\newcommand{\PushoutProduct}{
  \,\scalebox{.8}{$\Box$}\,
}

\newcommand{\simplicial}{
  \Delta
}

\newcommand{\Groupoids}{
  \mathrm{Grpd}
}

\newcommand{\OneGroupoids}{
  \Groupoids_1
}

\newcommand{\TwoGroupoids}{
  \mathrm{Grpd}_2
}

\newcommand{\EffectiveEpimorphisms}{
  \mathrm{Atl}
}

\newcommand{\TopologicalGroupoids}{
  \Groupoids(\kTopologicalSpaces)
}

\newcommand{\ConnectedGroupoidComponents}{
  \tau_0
}

\newcommand{\ConstantGroupoid}{
  \mathrm{Cnst}
}

\newcommand{\SpaceOfObjects}{
  (-)_0
}

\newcommand{\CodiscreteGroupoid}{
  \mathrm{Cht}
}

\newcommand{\InfinityGroupoids}{
  \Groupoids_\infty
}
\newcommand{\infinityGroupoids}{
  \InfinityGroupoids
}

\newcommand{\Spectra}{
  \mathrm{Spectra}
}

\newcommand{\TopologicalSpaces}{
  \topological\Spaces
}

\newcommand{\CWComplexes}{
  \mathrm{CWCplx}
}

\newcommand{\GCWComplexes}{
  G\CWComplexes
}

\newcommand{\TopologicalSpace}{
  \mathrm{X}
}

\newcommand{\TopologicalPatch}{
  \mathrm{U}
}

\newcommand{\TopologicalCoefficients}{
  \mathrm{A}
}

\newcommand{\TopologicalPrincipalBundle}{
  \mathrm{P}
}

\newcommand{\TopologicalFiberBundle}{
  \mathrm{E}
}

\newcommand{\kTopologicalSpaces}{
  \mathrm{k}\TopologicalSpaces
}

\newcommand{\DTopologicalSpaces}{
  \mathrm{D}\TopologicalSpaces
}

\newcommand{\DHausdorffSpaces}{
  \mathrm{D}\HausdorffSpaces
}

\newcommand{\ShapeOfSphere}[1]{
  S^{#1}
}

\newcommand{\TopologicalSphere}[1]{
  \ShapeOfSphere{#1}_{\mathrm{top}}
}

\newcommand{\SmoothSphere}[1]{
  \ShapeOfSphere{#1}_{\mathrm{sm}}
}

\newcommand{\RiemannianSphere}[1]{
  \ShapeOfSphere{#1}_{\mathrm{rm}}
}

\newcommand{\TorusGroup}[1]{
  \mathbb{T}^{#1}
}

\newcommand{\Orbi}{
  \mathrm{Orb}
}

\newcommand{\GOrbi}{
  G\Orbi
}

\newcommand{\OrbiSpaces}{
  \Orbi\Spaces
}

\newcommand{\FixedLoci}{
  \mathrm{FxdLoc}
}

\newcommand{\Cartesian}{
  \mathrm{Cart}
}

\newcommand{\CartesianSpaces}{
  \Cartesian\Spaces
}

\newcommand{\NaturalNumbers}{
  \mathbb{N}
}

\newcommand{\Integers}{
  \mathbb{Z}
}

\newcommand{\RealNumbers}{
  \mathbb{R}
}

\newcommand{\ComplexNumbers}{
  \mathbb{C}
}

\newcommand{\ImaginaryUnit}{
  \mathrm{i}
}

\newcommand{\Quaternions}{
  \mathbb{H}
}

\newcommand{\Manifolds}{
  \mathrm{Mfd}
}

\newcommand{\SmoothManifolds}{
  \mathrm{Smth}\Manifolds
}

\newcommand{\SmoothManifold}{
  \TopologicalSpace
}

\newcommand{\RiemannianManifold}{
  \SmoothManifold
}

\newcommand{\RiemannianManifolds}{
  \mathrm{RiemMfd}
}

\newcommand{\VolumeElement}{
  \mathrm{dvol}
}

\newcommand{\Compact}{
 \mathrm{Cpt}
}

\newcommand{\Simplex}[1]{
  \Delta^{#1}
}

\newcommand{\TopologicalSimplex}[1]{
  \Simplex{#1}_{\mathrm{top}}
}

\newcommand{\SmoothSimplex}[1]{
  \Simplex{#1}_{\mathrm{smth}}
}

\newcommand{\CompactSmoothManifolds}{
  \mathrm{Cpt}\SmoothManifolds
}

\newcommand{\Diffeological}{
  \mathrm{Dfflg}
}

\newcommand{\DiffeologicalSpaces}{
  \Diffeological\Spaces
}

\newcommand{\ContinuousDiffeology}{
  \mathrm{Cdfflg}
}

\newcommand{\DTopology}{
  \mathrm{Dtplg}
}

\newcommand{\kification}{
  k
}

\newcommand{\Dification}{
  \ContinuousDiffeology
}

\newcommand{\Categories}{
  \mathrm{Cat}
}

\newcommand{\op}{
  \mathrm{op}
}

\newcommand{\categorybox}[1]{
  \fcolorbox{lightgray}{lightgray}{$#1$}
}

\newcommand{\proofstep}[1]{
  \mbox{\small #1}
}

\newcommand{\Topos}{
  \mathbf{H}
}

\newcommand{\BaseTopos}{
  \mathbf{B}
}

\newcommand{\Cech}{
  C
}

\newcommand{\LocallyConstant}{
  \mathrm{LCnst}
}

\newcommand{\GlobalSections}{
  \mathrm{\Gamma}
}

\newcommand{\SimplicialSite}{
  \mathcal{S}
}

\newcommand{\InfinitySite}{
  \mathbf{S}
}

\newcommand{\GEquivariant}[1]{
  G#1
}

\newcommand{\GloballyEquivariant}[1]{
  \Singular#1
}

\newcommand{\Truncation}[1]{
  \tau_{#1}
}

\newcommand{\Twist}{
  \mathrm{tw}
}

\newcommand{\Modal}[2]{
  {#1}_{\scalebox{.7}{$#2$}}
}

\newcommand{\ModalTopos}[1]{
  \Modal{\Topos}{#1}
}

\newcommand{\Slice}[2]{
  {#1}_{\scalebox{.7}{$/#2$}}
}

\newcommand{\SliceTopos}[1]{
  \Slice{\Topos}{#1}
}

\newcommand{\Discrete}{
  \mathrm{Dsc}
}

\newcommand{\Chaotic}{
  \mathrm{Cht}
}

\newcommand{\Points}{
  \mathrm{Pnt}
}

\newcommand{\Shape}{
  \mathrm{Shp}
}

\newcommand{\shape}{
  \raisebox{1pt}{\rm\normalfont\textesh}
}

\newcommand{\IsomorphismClasses}[1]{
  \left(
    {#1}
  \right)_{/\sim_{\mathrm{iso}}}
}

\newcommand{\ConcordanceClasses}[1]{
  \left(
    {#1}
  \right)_{/\sim_{\mathrm{conc}}}
}

\newcommand{\SingularSimplicialComplex}{
  \mathrm{Pth}
}

\newcommand{\SmoothSingularSimplicialComplex}{
  \mathbf{Pth}
}

\newcommand{\Conical}{
  \mathrm{Cncl}
}

\newcommand{\Space}{
  \mathrm{Spc}
}

\newcommand{\Smooth}{
  \mathrm{Smth}
}

\newcommand{\Singularity}{
  \mathrm{Snglt}
}

\newcommand{\OrbiSingular}{
  \Singularity
}

\newcommand{\Orbisingular}{
  \OrbiSingular
}

\newcommand{\Maps}[3]{
  \mathrm{Map}
  #1(
    #2
    ,\,
    #3
  #1)
}

\newcommand{\SliceMaps}[4]{
  \Maps{#1}{#3}{#4}_{\!\scalebox{.7}{$#2$}}
}

\def\AmbientCategory{\Topos}

\newcommand{\PointsMaps}[3]{
  \AmbientCategory
  #1(
    #2
    ,\,
    #3
  #1)
}

\newcommand{\SlicePointsMaps}[4]{
  \PointsMaps{#1}{#3}{#4}_{\!\scalebox{.7}{$#2$}}
}

\newcommand{\Homs}[3]{
  \mathrm{Hom}
  #1(
    #2
    ,\,
    #3
  #1)
}

\newcommand{\Quotient}[2]{
  #1 / #2
}

\newcommand{\HomotopyQuotient}[2]{
  #1 \!\sslash\! #2
}

\newcommand{\ActionGroupoid}[2]{
  \left(
    {#1} \times #2
      \rightrightarrows
    {#1}
  \right)
}

\newcommand{\DeloopingGroupoid}[1]{
  \left(
    #1
      \rightrightarrows
    \ast
  \right)
}

\newcommand{\ModelCategories}{
  \mathrm{Mdl}\Categories
}

\newcommand{\Quillen}{
  \mathrm{Qu}
}

\newcommand{\QuillenEquivalences}{
  \Quillen\mathrm{Eq}
}

\newcommand{\Enriched}{
  \mathrm{Enr}
}

\newcommand{\EnrichedCategories}[1]{
  {#1}\Enriched\Categories
}

\newcommand{\SimplicialCategories}{
  \simplicial\Categories
}

\newcommand{\SimplicialFunctors}{
  \simplicial\Functors
}

\newcommand{\Localization}[1]{
  \mathrm{Loc}^{\scalebox{.6}{$#1$}}
}

\newcommand{\SimplicialLocalization}[1]{
  \Localization{#1}_{\simplicial}
}

\newcommand{\TwoLocalization}[1]{
  \Localization{#1}_2
}

\newcommand{\TwoLocalizationProjection}[1]{
  \SimplicialLocalization{\WeakEquivalences}
}

\newcommand{\HomotopyCategory}{
  \mathrm{Ho}
}

\newcommand{\TwoHomotopyCategory}{
  \HomotopyCategory
}

\newcommand{\CombinatorialModelCategories}{
  \mathrm{Comb}\ModelCategories
}

\newcommand{\SimplicialModelCategories}{
  \ModelCategories_{\simplicial}
}

\newcommand{\SimplicialCombinatorialModelCategories}{
  \CombinatorialModelCategories_{\simplicial}
}

\newcommand{\CombinatorialSimplicialModelCategories}{
  \SimplicialCombinatorialModelCategories
}

\newcommand{\ProperCombinatorialSimplicialModelCategories}{
  \CombinatorialModelCategories_{\simplicial, \mathrm{prop} }
}

\newcommand{\Toposes}{
  \mathrm{Topos}
}

\newcommand{\infinityToposes}{
  \Toposes_\infty
}

\newcommand{\InfinityToposes}{
  \infinityToposes
}

\newcommand{\PresentableCategories}{
  \mathrm{Pres}\Categories
}

\newcommand{\InfinityCategories}{
  \Categories_\infty
}

\newcommand{\PresentableInfinityCategories}{
  \PresentableCategories_\infty
}

\newcommand{\HomotopyTwoCategory}{
  \HomotopyCategory_{2}
}

\newcommand{\HomotopyTwoCategoryOfPresentableInfinityCategories}{
  \HomotopyTwoCategory
  (
    \PresentableInfinityCategories
  )
}

\newcommand{\WeakEquivalences}{
  \mathrm{WEqs}
}

\newcommand{\Fibrations}{
  \mathrm{Fib}
}

\newcommand{\Cofibrations}{
  \mathrm{Cof}
}

\newcommand{\HurewiczCofibrations}{
  \mathrm{h}\Cofibrations
}

\newcommand{\KanFibrations}{
  \mathrm{Kan}\Fibrations
}

\newcommand{\HomotopyKanFibrations}{
  \mathrm{h}\KanFibrations
}

\newcommand{\Serre}{
  \mathrm{Ser}
}

\newcommand{\SerreFibrations}{
  \Serre\Fibrations
}

\newcommand{\ProperEquivariantSerreFibrations}[1]{
  {#1}\SerreFibrations
}

\newcommand{\Projective}{
  \mathrm{Prj}
}

\newcommand{\Local}{
  \mathrm{Lcl}
}

\newcommand{\ProjectiveWeakEquivalences}{
  \Projective\WeakEquivalences
}

\newcommand{\LocalWeakEquivalences}{
  \Local\WeakEquivalences
}

\newcommand{\ProjectiveFibrations}{
  \Projective\Fibrations
}

\newcommand{\ProjectiveCofibrations}{
  \Projective\Cofibrations
}

\newcommand{\LocalProjectiveCofibrations}{
  \Local\ProjectiveCofibrations
}

\newcommand{\SerreCofibrations}{
  \Serre\Cofibrations
}

\newcommand{\ProperEquivariantSerreCofibrations}[1]{
  {#1}\SerreCofibrations
}

\newcommand{\WeakHomotopyEquivalences}{
  \mathrm{WHmtpEq}
}

\newcommand{\ProperEquivariantWeakHomotopyEquivalences}[1]{
  {#1}\WeakHomotopyEquivalences
}

\newcommand{\HomotopyFiber}{
  \mathrm{hofib}
}

\newcommand{\TopologicalRealization}[2]{
  #1\vert
  \,
    \!#2\!
  \,
  #1\vert
}

\newcommand{\HausdorffSpaces}{
  \Hausdorff\Spaces
}

\newcommand{\kHausdorffSpaces}{
  \mathrm{k}\HausdorffSpaces
}

\newcommand{\LocallyCompactHausdorffSpaces}{
  \mathrm{LC}\HausdorffSpaces
}

\newcommand{\SimplicialTopologicalSpaces}{
  \simplicial\kTopologicalSpaces
}

\newcommand{\TopologicalGroups}{
  \Groups(\kTopologicalSpaces)
}

\newcommand{\SimplicialTopologicalGroups}{
  \Groups(\SimplicialTopologicalSpaces)
}

\newcommand{\Sheaves}{
  \mathrm{Sh}
}

\newcommand{\Presheaves}{
  \mathrm{P}\Sheaves
}

\newcommand{\YonedaEmbedding}{
  y
}

\newcommand{\limit}[1]{
  \underset{
    \underset{#1}{\longleftarrow}
  }{\lim}
  \,
}

\newcommand{\colimit}[1]{
  \underset{
    \underset{#1}{\longrightarrow}
  }{\lim}
  \,
}

\newcommand{\SimplicialPresheaves}{
  \simplicial\Presheaves
}

\newcommand{\proj}{
  \mathrm{proj}
}

\newcommand{\projloc}{
  { \proj }
  \atop
  { \mathrm{loc} }
}

\newcommand{\InfinityPresheaves}{
  \Presheaves_{\infty}
}

\newcommand{\InfinitySheaves}{
  \Sheaves_{\infty}
}

\newcommand{\infinitySheaves}{
  \InfinitySheaves
}

\newcommand{\Orbits}{
  \mathrm{Orb}
}

\newcommand{\OrbitCategory}[1]{
  \Orbits
  \left(
    {#1}
  \right)
}

\newcommand{\Actions}[1]{
  {#1}\,\mathrm{Act}
}

\newcommand{\Principal}{
  \mathrm{Prn}
}

\newcommand{\Fiber}{
  \mathrm{Fib}
}

\newcommand{\Bundles}{
  \mathrm{Bdl}
}

\newcommand{\PrincipalBundles}[1]{
  {#1}\,\Principal\Bundles
}

\newcommand{\PrincipalFiberBundles}[1]{
  {#1}\,\Principal\Fiber\Bundles
}

\newcommand{\Equivariant}{
  \mathrm{Equv}
}

\newcommand{\EquivariantPrincipalBundles}[2]{
  {#1}\,\Equivariant\,\PrincipalBundles{#2}
}

\newcommand{\EquivariantPrincipalFiberBundles}[2]{
  {#1}\,\Equivariant\,\PrincipalFiberBundles{#2}
}

\newcommand{\Formally}{
  \mathrm{Frm}
}

\newcommand{\FormallyPrincipalBundles}[1]{
  \Formally\PrincipalBundles{#1}
}

\newcommand{\FormallyPrincipalFiberBundles}[1]{
  \Formally\PrincipalFiberBundles{#1}
}

\newcommand{\GActionsOnTopologicalSpaces}{
  \Actions{G}(\kTopologicalSpaces)
}

\newcommand{\GEquivariantTopologicalGroups}{
  G\mathrm{EquivTopGrp}
}

\newcommand{\SmoothSets}{
  \Smooth\Groupoids_0
}

\newcommand{\SmoothGroupoids}{
  \Smooth\Groupoids_{1}
}

\newcommand{\SmoothInfinityGroupoids}{
  \Smooth\Groupoids_{\infty}
}

\newcommand{\EquivariantGroups}[1]{
    #1\Equivariant\Groups
}

\newcommand{\Representations}{
  \mathrm{Rep}
}

\newcommand{\Character}{
  \rchi
}

\newcommand{\WeylGroup}{
    \mathrm{W}
      _{\scalebox{.6}{$\!\!G$}}
}

\newcommand{\mapsdown}{\rotatebox[origin=c]{-90}{$\mapsto$}}

\def\acts{\raisebox{1.4pt}{\;\rotatebox[origin=c]{90}{$\curvearrowright$}}\hspace{.5pt}}

\def\rightacts{\rotatebox[origin=c]{180}{$\acts$}}

\DeclareRobustCommand{\rchi}{{\mathpalette\irchi\relax}}
\newcommand{\irchi}[2]{\raisebox{\depth}{$#1\chi$}} 

\makeatletter
\newif\if@sup
\newtoks\@sups
\def\append@sup#1{\edef\act{\noexpand\@sups={\the\@sups #1}}\act}%
\def\reset@sup{\@supfalse\@sups={}}%
\def\mk@scripts#1#2{\if #2/ \if@sup ^{\the\@sups}\fi \else%
  \ifx #1_ \if@sup ^{\the\@sups}\reset@sup \fi {}_{#2}%
  \else \append@sup#2 \@suptrue \fi%
  \expandafter\mk@scripts\fi}
\def\tensor#1#2{\reset@sup#1\mk@scripts#2_/}
\def\multiscripts#1#2#3{\reset@sup{}\mk@scripts#1_/#2%
  \reset@sup\mk@scripts#3_/}
\makeatother

\makeatletter
\newbox\slashbox \setbox\slashbox=\hbox{$/$}
\def\itex@pslash#1{\setbox\@tempboxa=\hbox{$#1$}
  \@tempdima=0.5\wd\slashbox \advance\@tempdima 0.5\wd\@tempboxa
  \copy\slashbox \kern-\@tempdima \box\@tempboxa}
\def\slash{\protect\itex@pslash}
\makeatother

\def\clap#1{\hbox to 0pt{\hss#1\hss}}
\def\mathllap{\mathpalette\mathllapinternal}
\def\mathrlap{\mathpalette\mathrlapinternal}
\def\mathclap{\mathpalette\mathclapinternal}
\def\mathllapinternal#1#2{\llap{$\mathsurround=0pt#1{#2}$}}
\def\mathrlapinternal#1#2{\rlap{$\mathsurround=0pt#1{#2}$}}
\def\mathclapinternal#1#2{\clap{$\mathsurround=0pt#1{#2}$}}

\let\oldroot\root
\def\root#1#2{\oldroot #1 \of{#2}}
\renewcommand{\sqrt}[2][]{\oldroot #1 \of{#2}}

\DeclareSymbolFont{symbolsC}{U}{txsyc}{m}{n}
\SetSymbolFont{symbolsC}{bold}{U}{txsyc}{bx}{n}
\DeclareFontSubstitution{U}{txsyc}{m}{n}

\DeclareSymbolFont{stmry}{U}{stmry}{m}{n}
\SetSymbolFont{stmry}{bold}{U}{stmry}{b}{n}

\DeclareFontFamily{OMX}{MnSymbolE}{}
\DeclareSymbolFont{mnomx}{OMX}{MnSymbolE}{m}{n}
\SetSymbolFont{mnomx}{bold}{OMX}{MnSymbolE}{b}{n}
\DeclareFontShape{OMX}{MnSymbolE}{m}{n}{
    <-6>  MnSymbolE5
   <6-7>  MnSymbolE6
   <7-8>  MnSymbolE7
   <8-9>  MnSymbolE8
   <9-10> MnSymbolE9
  <10-12> MnSymbolE10
  <12->   MnSymbolE12}{}

\usepackage{cleveref}

\crefformat{section}{\S#2#1#3} 
\crefformat{subsection}{\S#2#1#3}
\crefformat{subsubsection}{\S#2#1#3}

\theoremstyle{italics}
\newtheorem{theorem}{Theorem}[section]
\newtheorem{lemma}[theorem]{Lemma}
\newtheorem{proposition}[theorem]{Proposition}
\newtheorem{corollary}[theorem]{Corollary}
\theoremstyle{definition}
\newtheorem{definition}[theorem]{Definition}
\newtheorem{notation}[theorem]{Notation}
\newtheorem{example}[theorem]{Example}
\newtheorem{assumption}[theorem]{Assumption}

\newtheorem{remark}[theorem]{Remark}

\renewcommand{\emph}{\textit}

\begin{document}

\title{
  Equivariant Principal $\infty$-Bundles
}
\author{Hisham Sati, \quad Urs Schreiber}

\maketitle

\thispagestyle{empty}

$\,$

\vspace{3cm}

\begin{center}
\end{center}

\vfill

\noindent

\vspace{2cm}

\noindent  Hisham Sati, {\it Mathematics, Division of Science,
\\
and Center for Quantum and Topological Systems (CQTS), NYUAD Research Institute,
\\
New York University Abu Dhabi, UAE,
}
\\
{\tt hsati@nyu.edu}
\\
\\
\noindent  Urs Schreiber, {\it Mathematics, Division of Science,
\\
and Center for Quantum and Topological Systems (CQTS), NYUAD Research Institute,
\\
New York University Abu Dhabi, UAE;
\\
on leave from Czech Academy of Science, Prague.}
\\
{\tt us13@nyu.edu}

\newpage
\thispagestyle{empty}
$\,$

\vspace{-1.8cm}
\begin{center}
  \bf What this book is about.
\end{center}

\medskip

In this book we prove (Thm. \ref{ProperClassificationOfEquivariantBundlesForResolvableSingularitiesAndEquivariantStructure})
unified classification results for stable equivariant $\Gamma$-principal
bundles when the underlying homotopy type $\shape \, \Gamma$ of the
topological structure group $\Gamma$ is truncated, meaning that its homotopy groups vanish
in and above some degree $n$.
We discuss
(Thm. \ref{BorelClassificationOfEquivariantBundlesForResolvableSingularitiesAndEquivariantStructure}.ii)
how this coincides with the classification of
equivariant higher non-abelian gerbes
and generally of
equivariant {\it principal $\infty$-bundles} \cite{NSS12a}
with structure $n$-group $\shape \, \Gamma$;
and we show how the equivariant homotopy groups
of the respective classifying $G$-spaces are given by the non-abelian
group cohomology
of the equivariance group
with coefficients in $\shape \, \Gamma$
(Thm. \ref{EquivariantHomotopyGroupsOfEquivariantClassifyingSpaces}).

\medskip
The result is proven in a conceptually transparent manner as
a consequence of a
{\it smooth Oka principle}
(Thm. \ref{SmoothOkaPrinciple}, \ref{OrbiSmoothOkaPrinciple},
based on \cite{BBP19}\cite{Pavlov14}),
which is available after faithfully embedding
traditional equivariant topology
into the singular-cohesive homotopy theory
of globally equivariant higher stacks \cite{SS20OrbifoldCohomology}.
This works for
discrete equivariance groups $G$
acting properly on smooth manifolds
(``proper equivariance'' \cite{DHLPS19})
with resolvable singularities (Ntn. \ref{ResolvableOrbiSingularities}),
whence we are equivalently describing principal bundles on good orbifolds.

\medskip

In setting up this proof, we
re-develop the theory of equivariant principal bundles
from scratch by systematic use of Grothendieck's {\it internalization}
(Ntn. \ref{Internalization}).
In particular we prove
(Thm.
\ref{BorelClassificationOfEquivariantBundlesForResolvableSingularitiesAndEquivariantStructure}.i)
that all the intricate
equivariant local triviality conditions
considered in the literature
\cite{tomDieck69}\cite{Bierstone73}\cite{LashofMay86}
are automatically {\it implied}
by regarding $G$-equivariant principal bundles
as principal bundles internal to the $BG$-slice of the
ambient cohesive $\infty$-topos.
We also show that these
conditions are all equivalent
(Thm. \ref{EquivalentNotionsOfEquivariantLocalTriviality}).
Generally we find (see \cref{EquivariantModuliStacks})
that the characteristic subtle phenomena
of equivariant classifying theory all reflect basic {\it modal}
properties of singular-cohesive homotopy theory
(hence of cohesive ``global'' equivariant homotopy theory
\cite{Rezk14}, see e.g.
Prop. \ref{GloballyProperEquivariantClassifyingShape}).

\medskip

A key example is the projective unitary structure group
$\Gamma = \PUH$, in which case
(Ex. \ref{EquivariantBundlesServingAsGeoemtricTwistsOfEquivariantKTheory})
we are classifying
3-twists of equivariant $\KU$-theory
(related to equivariant bundle gerbes, Ex. \ref{EquivariantBundleGerbes}),
recovering the results
\cite[Cor. 2.41]{TXLG04}\cite[Thm. 6.3]{AtiyahSegal04}\cite[Thm. 1.10]{BEJU12}\cite[Thm. 15.17]{LueckUribe14}.
Our general theorem immediately enhances this to
the conjugation-equivariant graded projective unitary structure group
$\ZTwo \acts \, \GradedPUH$ with fixed locus $(\GradedPUH)^{\ZTwo} \,=\, \GradedPOH$
(Ex. \ref{ProjectiveUnitarGroupOnAHilbertSpace}),
where we are classifying the twists of equivariant
$\KR$-theory in degrees 1 and 3 combined, restricting
on ``O-planes'' to the
twists of $\KO$-theory in degrees 1 and 2 combined.
This is the generality in which equivariant K-theory twists
are conjectured to model
quantum symmetries of topological phases of matter
(Rem. \ref{RealEquivariantProjectiveGradedBundlesAsQuantumSymmetries})
and
the B-field in string theory on orbi-orientifolds
(Rem. \ref{GeometricTwistsOfEquivariantKTheoryAsBFieldsOnOrbiOrientifolds}).

\medskip
Our focus on the class of truncated structure groups is largely complementary
to the classical literature on equivariant classifying spaces
\cite{tomDieck69}\cite{Bierstone73}\cite{Lashof82}\cite{LashofMay86}\cite{May90},
which is mostly concerned with
the class of compact Lie structure groups $\Gamma$
(see \hyperlink{Table1}{Tab. 1}).
The two classes intersect when
$\Gamma = \mathbb{T}^d \rtimes K$ is the extension of a finite group $K$
by a compact abelian Lie group,
i.e., the case of equivariant
sheeted torus bundles (``bundle 0-gerbes'').
In this case  our general classification recovers
(Ex. \ref{TruncatedStructureGroupsAndTheirEquivariantClassificationResults})
the results of
\cite[Lem. 1, Thm. 2]{LMS83}\cite[Thm. 3, Thm. 10]{May90}
and \cite[Thm. 1.2]{Rezk16}
and generalizes them to the case of non-trivial $G$-action on $\Gamma$.
\medskip

A classification result for general topological structure groups
had previously been claimed in \cite{MurayamaShimakawa95} but,
as highlighted in \cite[\S 3.12]{GuillouMayMerling17},
the proof had remained open.
Besides characterizing its classification property for truncated
structure groups, we show
that the {\it Murayama-Shimakawa construction}
of \cite{MurayamaShimakawa95} generally produces
the underlying equivariant homotopy type
of the correct equivariant {\it moduli stacks}
(Thm. \ref{MurayamaShimakawaGroupoidIsEquivariantModuliStack}).

\medskip

While compact Lie structure groups $\Gamma$ have received much
attention due to the role of equivariant vector bundles
as {\it cocycles} of equivariant K-theory, the complementary case of
truncated structure groups $\Gamma$ that we discuss
will generally be relevant for equivariant bundles in their role as
{\it twists} of equivariant generalized cohomology theories,
since such twists will typically be considered in a finite number of degrees,
with $\shape \, \Gamma$
a truncation of the $\infty$-group of units of a ring spectrum.
This combined {\it twisted \& equivariant } enhancement of generalized
cohomology theory has previously received little attention
(except for a brief note in \cite{Lind14})
beyond the example of $\KU$-theory. We mean this book to be laying
previously missing foundations (\cref{TwistedEquivariantCohomology}).

\medskip

In particular, via a twisted
generalization of Elmendorf's Theorem (Thm. \ref{TwistedEDKtheorem}),
the equivariant classifying spaces produced
here provide the correct domains for twisted equivariant
enhancements of Chern-Dold character maps
(further generalizing \cite{FSS20CharacterMap}\cite{SS20EquivariantTwistorial})
and hence allow the systematic definition and construction of
twisted equivariant {\it differential} generalized cohomology theories
in equivariant generalization of the construction in
\cite[\S 8]{FSS20CharacterMap}.
This will be discussed in
\cite{TwistedEquivariantChernCharacter}\cite{TwistedEquivariantDifferentialCohomology}.

\newpage
\thispagestyle{empty}

\setcounter{tocdepth}{1}
\tableofcontents

\newpage

\part{Introduction}
\label{IntroductionAndOverview}

We assume that the reader already appreciates principal bundles
and their classifying theory
(e.g., as in \cite{Husemoller94} \cite{Cohen17}\cite[\S 1.1]{RudolphSchmidt17}\cite{Tamaki21})
and is already motivated towards their equivariant generalization
(due to \cite{Stewart61} \cite{tomDieck69}\cite{Bierstone73}\cite{Lashof82}\cite{LashofMay86}),
though we will effectively review and re-develop much  of the classical theory.

\medskip

The content of \cref{InTopologicalSpaces}
is a self-contained discussion
of equivariant principal bundles  entirely in point-set topology,
which may be read independently of \cref{InCohesiveInfinityStacks}
and may be of interest in itself,
as a comprehensive,
streamlined,
and
completed
review.

\medskip
The purpose of \cref{InCohesiveInfinityStacks} is to embed this
classical topic into the relatively new context of higher geometry,
where the nature of equivariant classifying spaces becomes more transparent
and where {\it cohesive homotopy theory}
(\cite[\S 3.1]{SSS09}\cite{dcct}\cite{SchreiberShulman14}\cite{Shulman15}\cite[\S 3]{SS20OrbifoldCohomology})
provides new tools even for the classical theory,
which we use to prove a new classification theorem.
This higher geometric perspective
on principal bundles
(following \cite[\S 3]{SSS09}\cite{NSS12a}\cite{NSS12b})
and on equivariance
(following \cite{Rezk14}\cite{SS20OrbifoldCohomology})
may conceptually be understood as follows.

\medskip

\section{Overview and Summary}
\label{OverviewAndSummary}

\medskip
\noindent
{\bf Higher cohesive geometry as convenient topology.}
For the classical algebraic topologist, one way to appreciate the passage to
higher geometric and to cohesive homotopy theory is as a completion
of the pursuit of ever more {\it convenient categories} of topological spaces
in the famous sense of \cite{Steenrod67}.
This starts out with compactly generated spaces
(used in \cref{InTopologicalSpaces}, see Prop. \ref{CompactlyGeneratedTopologicalSpaces})
but proceeds, along the general lines of \cite{Vogt71}\cite{Wyler73}\cite{FajstrupRosicky08},
to the {\it really convenient category}
\cite{Smith}\cite{Dugger03} of Delta-generated topological spaces
(used in \cref{InCohesiveInfinityStacks}, see Ntn. \ref{ContinuousDiffeologyAndDTopology}).

\medskip
At this point, we may observe \cite[Ex. 3.18]{SS20OrbifoldCohomology}
with \cite[\S 3]{SYH10} (see Prop. \ref{CategoryOfDeltaGeneratedTopologicalSpaces} below)
that
D-topology faithfully embeds first into {\it diffeology}
(\cite{IglesiasZemmour13}, see Ntn. \ref{CartesianSpacesAndDiffeologicalSpaces} below)
and
from here further into the cohesive topos
$\SmoothSets$
over the site of smooth manifolds
and then the cohesive $\infty$-topos
$\SmoothInfinityGroupoids$ (Ntn. \ref{SmoothInfinityGroupoids} below):
$$
  \hspace{-.4cm}
  \begin{tikzcd}[column sep=28pt, row sep=22.5pt]
    \mathclap{
    \mbox{
      \tiny
      \bf
      \def\arraystretch{.9}
      \begin{tabular}{c}
        \phantom{a}
        \\
        \color{gray}
        category of
        \\
        \color{darkblue}
        topological
        \\
        \color{darkblue}
        spaces
      \end{tabular}
    }
    }
    \ar[
      rrrrr,
      shift left=20pt,
      dashed,
      "{
        \mbox{
          \tiny
          \bf
          convenience
        }
        \;\;\;
      }"{description}
    ]
  &
  \mathclap{
    \mbox{
      \tiny
      \bf
      \def\arraystretch{.9}
      \begin{tabular}{c}
        \color{gray}
        cartesian closed
        \\
        \color{gray}
        category of
        \\
        \color{darkblue}
        compactly generated
        \\
        \color{darkblue}
        topological spaces
      \end{tabular}
    }
  }
  &
  \mathclap{
    \mbox{
      \tiny
      \bf
      \def\arraystretch{.9}
      \begin{tabular}{c}
        \color{gray}
        cart. clsd. \& locl. pres.
        \\
        \color{gray}
        category of
        \\
        \color{darkblue}
        $\Delta$-generated
        \\
        \color{darkblue}
        topological spaces
      \end{tabular}
    }
  }
  &
  \mathclap{
    \mbox{
      \tiny
      \bf
      \def\arraystretch{.9}
      \begin{tabular}{c}
        \phantom{a}
        \\
        \color{gray}
        quasi-topos of
        \\
        \color{darkblue}
        diffeological
        \\
        \color{darkblue}
        spaces
      \end{tabular}
    }
  }
  &
  \mathclap{
    \mbox{
      \tiny
      \bf
      \def\arraystretch{.9}
      \begin{tabular}{c}
        \phantom{a}
        \\
        \color{gray}
        cohesive topos of
        \\
        \color{darkblue}
        smooth
        \\
        \color{darkblue}
        sets
      \end{tabular}
    }
  }
  &
  \mathclap{
    \mbox{
      \tiny
      \bf
      \def\arraystretch{.9}
      \begin{tabular}{c}
        {\phantom{a}}
        \\
        \color{gray}
        cohesive $\infty$-topos of
        \\
        \color{darkblue}
        smooth
        \\
        \color{darkblue}
        $\infty$-groupoids
      \end{tabular}
    }
  }
  \\[-25pt]
   \categorybox{\TopologicalSpaces}
    \ar[
      r,
      phantom,
      "{
        \scalebox{.7}{$\bot$}
      }"
    ]
    \ar[
      r,
      shift right=7pt,
      "{
          \kification
      }"{below}
    ]
    \ar[
      rrrrrd,
      rounded corners,
      to path={
        -- ([yshift=-35pt]\tikztostart.south)
        -- node[below]{\scalebox{.7}{$
             \raisebox{1pt}{
               \scalebox{.9}{
               \def\arraystretch{.9}
               \begin{tabular}{c}
                 \color{greenii}
                 \bf
                 form path $\infty$-groupoids
                 \\
                 \scalebox{.85}{(sing. simpl. compl.)}
               \end{tabular}
               }
             }
             \;
             \SingularSimplicialComplex
           $}}
           node[above, yshift=8pt]{\scalebox{1}{$
              \mathrlap{
                \mbox{
                  \hspace{2.5cm}
                  \tiny
                  \color{orangeii}
                  \bf
                  underlying homotopy types are preserved
                }
              }
           $}}
           (\tikztotarget.west)
      }
    ]
    &
    \categorybox{\kTopologicalSpaces}
    \ar[
      l,
      hook',
      shift right=7pt
    ]
    \ar[
      r,
      phantom,
      "{
        \scalebox{.7}{$\bot$}
      }"
    ]
    \ar[
      r,
      shift right=7pt,
      "{
        \scalebox{.7}{$\Dification $}
      }"{below}
    ]
    \ar[
      rrrrd,
      rounded corners,
      to path={
        -- ([yshift=-35pt]\tikztostart.south)
        -- (\tikztotarget.west)
      }
    ]
    &
    \categorybox{\DTopologicalSpaces}
    \ar[
      l,
      shift right=7pt,
      hook'
    ]
    \ar[
      r,
      phantom,
      "{
        \scalebox{.7}{$\bot$}
      }"
    ]
    \ar[
      r,
      hook,
      shift right=7pt
    ]
    \ar[
      rrrd,
      rounded corners,
      to path={
        -- ([yshift=-35pt]\tikztostart.south)
        -- (\tikztotarget.west)
      }
    ]
    &
    \categorybox{\DiffeologicalSpaces}
    \ar[
      l,
      shift right=7pt,
      "{
       \scalebox{.7}{$ \DTopology $}
      }"{above}
    ]
    \ar[
      r,
      hook,
      shift right=7pt,
      "{
        i_{\sharp_1}
      }"{below}
    ]
    \ar[
      rrd,
      rounded corners,
      to path={
        -- ([yshift=-35pt]\tikztostart.south)
        -- (\tikztotarget.west)
      }
    ]
    &
    \categorybox{\SmoothSets}
    \ar[
      r,
      hook,
      shift right=7pt,
      "{
        i_0
      }"{below}
    ]
    &
    \categorybox{\SmoothInfinityGroupoids}
    \ar[
      to path={
        -- ([xshift=-10pt]\tikztostart.south)
        --
            node[below, sloped]{
              \scalebox{.7}{$
                \underset{
                  \raisebox{+1pt}{
                    \color{greenii}
                    \bf
                    shape
                  }
                }{
                  \Shape
                }
              $}
            }
           ([xshift=-10pt, yshift=-24pt]\tikztostart.south)
        --
            node[above, yshift=5pt]{ $\shape$ }
            ([xshift=+10pt, yshift=-24pt]\tikztostart.south)
        --
            node[below, sloped]{ \scalebox{.7}{$\Discrete$}  }
           ([xshift=+10pt]\tikztotarget.south)
      }
    ]
    \\
    &&&&&
    \categorybox{\InfinityGroupoids}
    \\[-25pt]
    &&&&&
    \mathclap{
    \mbox{
      \tiny
      \bf
      \def\arraystretch{.9}
      \begin{tabular}{c}
        \color{gray}
        base $\infty$-topos of
        \\
        \color{darkblue}
        bare $\infty$-groupoids
      \end{tabular}
    }
    }
  \end{tikzcd}
$$
This diagram is
homotopy-commutative, which says
that the usual underlying homotopy type\footnote{
The ``esh''-symbol ``$\shape$''
\eqref{TheModalitiesOnACohesiveInfinityTopos}
stands for \emph{shape}
\cite[3.4.5]{dcct}\cite[9.7]{Shulman15},
following \cite{Borsuk75}, which for the well-behaved topological spaces
of algebraic topology is a convenient synonym for
their \emph{underlying weak homotopy type} \cite[7.1.6]{Lurie09HTT}\cite[4.6]{Wang17}.
It is in this sense that we use the term ``shape''
generalized to objects of cohesive $\infty$-toposes (Def. \ref{CohesiveInfinityTopos}).

This is consistent
(by Prop. \ref{ShapeOfSliceOfCohesiveInfinityToposIsCohesiveShape} below)
with the established generalization
of shape theory to $\infty$-topos theory
(Def. \ref{ShapeOfAnInfinityTopos}),
in that the weak homotopy type of any topological space $\TopologicalSpace$
is equivalently the shape of the slice of
the $\infty$-topos $\SmoothInfinityGroupoids$
over its continuous-diffeological incarnation
$\ContinuousDiffeology(X)$
(by Prop. \ref{SmoothShapeOfSmoothInfinityToposOverContinuousDiffeologyIfWeakHomotopyType}
below).}
$\shape \, \TopologicalSpace$
of topological spaces $\TopologicalSpace$ is preserved
under all these embeddings
(Props.
\ref{DiffeologicalShapeOfContinuousDiffeologyIsWeakHomotopyType},
\ref{SmoothShapeOfTopologicalSpacesIsTheirWeakHomotopyType}
below).
Moreover, the underlying homotopy type of {\it mapping spaces}
coincides with that of the corresponding {\it mapping stacks}
(Prop. \ref{DiffeologicalMappingSpacesHaveCorrectUnderlyingHomotopyType} below):
\vspace{-2mm}
$$
  X,\, Y
  \;\in\;
  \kTopologicalSpaces
  \;\;\;\;\;\;\;\;
    \vdash
  \;\;\;\;\;\;\;\;
  \overset{
    \mathclap{
    \raisebox{5pt}{
      \tiny
      \color{darkblue}
      \bf
      \begin{tabular}{c}
        shape of mapping space
      \end{tabular}
    }
    }
  }{
    \shape
    \,
    \Maps{}
      { X }
      { Y }
  }
  \;\;\;
    \simeq
  \;\;\;
  \overset{
    \mathclap{
    \raisebox{3pt}{
      \tiny
      \color{darkblue}
      \bf
      \begin{tabular}{c}
        shape of topological mapping stack
      \end{tabular}
    }
    }
  }{
    \shape
    \,
    \Maps{\big}
      { \ContinuousDiffeology(X) }
      { \ContinuousDiffeology(Y) }
  }.
$$

\noindent
{\bf The smooth Oka principle.}
After this embedding of differential topology into
cohesive homotopy theory,
a
wonderful theorem becomes available, which we call
the {\it smooth Oka principle}\footnote{
  Our terminology
  {\it smooth Oka principle}
  is meant to rhyme on the
  established
  {\it Oka-Grauert principle} in complex analysis,
  which says
  (reviewed in \cite[Cor. 3.5, 3.2]{ForsternicLarusson11}\cite[(1.1)]{ForsternicLarusson13})
  that over Stein manifolds:
  (a) the homotopy type of the space of holomorphic
  maps into an Oka manifold
  coincides with that of continuous maps
  between the underlying topological spaces;
  (b) the classification of
  holomorphic vector bundles
  coincides with that of topological vector bundles.
  In fact, also the $\infty$-topos over
  the site of complex-analytic manifolds is cohesive
  (this is implicit in \cite[\S 2.1]{HopkinsQuick15}),
  suggesting that
  the classical Oka-Grauert principle and our
  smooth Oka principle are just two special cases
  of one unifying geometric homotopy principle
  in cohesive homotopy theory.
}
(Thm. \ref{SmoothOkaPrinciple} below, based on \cite{BBP19}, see also \cite{Pavlov14}\cite[Thm. B]{Clough21}),
saying that, for smooth maps out of a smooth manifold
into any smooth $\infty$-groupoid, the
underlying homotopy type of the mapping stack is
that of the mapping space of underlying homotopy types:
\vspace{-1mm}
$$
  \left.
  \arraycolsep=2pt
  \begin{array}{cl}
    & \mbox{\small $\SmoothManifold$ a smooth manifold}
    \\
    \mbox{\small \&}
    & \mbox{\small $A$ any smooth $\infty$-groupoid}
  \end{array}
  \right\}
  \;\;\;\;\;\;\;\;\;\;\;\;\;
  \vdash
  \;\;\;\;\;\;\;\;\;\;\;\;\;
  \overset{
    \mathclap{
    \raisebox{8pt}{
      \tiny
      \color{darkblue}
      \bf
      \begin{tabular}{c}
        shape of
        mapping stack
      \end{tabular}
    }
    }
  }{
    \shape
    \,
    \Maps{}
      { \SmoothManifold }
      { A }
    }
  \;\;
  \underset{
    \mathclap{
    \raisebox{-6pt}{
      \tiny
      \color{greenii}
      \bf
       \begin{tabular}{c}
      smooth
      \\
      Oka principle
      \end{tabular}
    }
    }
  }{
    \simeq
  }
  \;\;
  \overset{
    \mathclap{
    \raisebox{4pt}{
      \tiny
      \color{darkblue}
      \bf
      \begin{tabular}{c}
        mapping space
        of shapes
      \end{tabular}
    }
    }
  }{
    \Maps{\big}
      { \shape \, \SmoothManifold }
      { \shape \, A }
  }
  \,.
$$

When $A \,=\, \mathbf{B}\mathcal{G}$ is the moduli $\infty$-stack
for $\mathcal{G}$-principal $\infty$-bundles
(see
Prop. \ref{GroupsActionsAndFiberBundles},
Thm. \ref{DeloopingGroupoidsAreModuliInfinityStacksForPrincipalInfinityBundles})
of a smooth $\infty$-group $\mathcal{G}$,
then the smooth Oka principle immediately implies that its shape
$$
  \begin{tikzcd}[column sep=large]
  \overset{
    \mathclap{
    \raisebox{6pt}{
      \tiny
      \color{darkblue}
      \bf
      \begin{tabular}{c}
        moduli $\infty$-stack for
        \\
        $\mathcal{G}$-principal $\infty$-bundles
      \end{tabular}
    }
    }
  }{
    \mathbf{B}\mathcal{G}
  }
  \quad\coloneqq\quad
  \HomotopyQuotient{\ast}{\mathcal{G}}
  \ar[
    rr,
    "{
      \eta^{\scalebox{.6}{\shape}}_{\mathbf{B}\mathcal{G}}
    }"
    ,
    "{
      \raisebox{-3pt}{
        \tiny
        \color{greenii}
        \bf
        \begin{tabular}{c}
          send bundles to their
          \\
          characteristic classes
        \end{tabular}
      }
    }"{swap, yshift=-2pt}
  ]
  &&
  \quad
  \overset{
    \mathclap{
    \raisebox{6pt}{
      \tiny
      \color{darkblue}
      \bf
      \begin{tabular}{c}
      classifying shape for concordance of
      \\
      smooth $\mathcal{G}$-principal $\infty$-bundles
      \end{tabular}
    }
    }
  }{
  B \mathcal{G}
    \;\,\coloneqq\;\,
  \shape \, \mathbf{B} \mathcal{G}
    \,\;\simeq\;\,
  \mathbf{B} \,\shape\, \mathcal{G}
  }
  \end{tikzcd}
  \;\;\;\;\;
  \in
  \;\;\;
  \InfinityGroupoids
  \xhookrightarrow{ \Discrete }
  \SmoothInfinityGroupoids
$$
is the homotopy type of a classifying space
(the {\it classifying shape})
for concordance classes
(Def. \ref{ConcordanceOfSmoothPrincipalBundles})
of
smooth $\mathcal{G}$-principal $\infty$-bundles (Thm. \ref{ClassificationOfPrincipalInfinityBundles}).
When  $\mathcal{G} \,=\, \Gamma$ is an ordinary Hausdorff-topological group,
this $B \Gamma$ is equivalently its traditional
Milgram classifying space (Prop. \ref{MilgramClassifyingSpaceModelClassifyingShape}).
But, since ordinary topological principal bundles are
isomorphic if they are concordant
(Thm. \ref{ConcordanceClassesOfTopologicalPrincipalBundles}),
this means that the
smooth Oka principle specializes to
an elegant cohesive re-proof of the classical classification theorem
for principal bundles
(Thm. \ref{ClassificationOfPrincipalBundlesAmongPrincipalInfinityBundles}).

\medskip

\noindent
{\bf An orbi-smooth Oka principle.}
Our strategy is to generalize this cohesive bundle theory to the equivariant
context, using that $G$-equivariant principal $\infty$-bundles on
$G \acts \, \SmoothManifold$
are equivalently principal $\infty$-bundles on the
{\it cohesive orbispace} (e.g. {\it orbifold})
$\HomotopyQuotient{\SmoothManifold}{G}$ (\cref{SmoothEquivariantPrincipalInfinityBundles},
see \cite{SS20OrbifoldCohomology} for background on orbifolds in
cohesive homotopy theory).
While the smooth Oka principle cannot hold
for general domains $X$ beyond smooth manifolds,
it will hold with given codomain $A$
for those $\infty$-stacks $X$ that may be
{\it approximated} by smooth manifolds to a degree of accuracy
which the coefficients $A$ do not resolve.
Specifically, if $\Gamma$ is a topological group whose
shape $\shape \, \Gamma$ is $n$-truncated, then $\mathbf{B}\Gamma$
does not distinguish
(Lem. \ref{PrincipalBundlesOnBlowUps})
the
abstract $G$-orbifold singularity $\HomotopyQuotient{\ast}{G}$
from its ``blowup'' $\Quotient{\SmoothSphere{\,n+2}}{G}$
to a smooth spherical space form (Def. \ref{SmoothSphericalSpaceForm})
of dimension above the truncation degree.
Accordingly, an orbi-smooth Oka principle holds
for truncated $\Gamma$-principal bundles
on orbifolds with resolvable singularities
which are ``stable'' (Ntn. \ref{StableEquivariantBundles})
against the blowup
(Thm. \ref{OrbiSmoothOkaPrinciple}):
\begin{equation}
  \label{OrbiSmoothOkaPrincipleInIntroduction}
  \hspace{-6mm}
  \left.
  \arraycolsep=2pt
  \begin{array}{clll}
    &
    \mbox{\small $G$ discrete with resolvable singularities}
    \\
    \mbox{\small \&}
    &
    \mbox{\small $\Gamma$  Hausdorff of truncated shape}
    \\
    \mbox{\small \&}
    &
    \mbox{\small $G \acts \, \SmoothManifold$ a proper smooth $G$-manifold }
  \end{array}
  \right\}
  \;\;\;\;
  \vdash
  \;\;\;\;
  \overset{
    \mathclap{
    \raisebox{4pt}{
      \tiny
      \color{darkblue}
      \bf
      \begin{tabular}{c}
        shape of stable part of slice mapping stack
      \end{tabular}
    }
    }
  }{
    \shape
    \,
    \SliceMaps{\big}{\mathbf{B}G}
      { \HomotopyQuotient{\SmoothManifold}{G} }
      { \HomotopyQuotient{\mathbf{B}\Gamma}{G} }
    ^{\stable}
  }
  \quad
  \underset{
    \mathclap{
    \raisebox{-9pt}{
      \tiny
      \color{greenii}
      \bf
       \begin{tabular}{c}
       \\
        orbi-smooth
        \\
        Oka principle
      \end{tabular}
    }
    }
  }{
    \simeq
  }
  \quad
  \overset{
    \mathclap{
    \raisebox{4pt}{
      \tiny
      \color{darkblue}
      \bf
      \begin{tabular}{c}
        slice mapping space of shapes
      \end{tabular}
    }
    }
  }{
    \SliceMaps{\big}{B G}
      { \HomotopyQuotient{ \shape \, \SmoothManifold }{G} }
      { \HomotopyQuotient{B \Gamma}{G} }
  }.
\end{equation}
As before, this immediately specializes to the classification theorem
for equivariant bundles up to concordance,
and this
--
since also equivariant principal bundles are isomorphic if they are
concordant (Lem. \ref{ConcordantEquivariantPrincipalBundlesConstructedFromCechCocyclesAreIsomorphic})
--
proves the full classification theorem
(Thm. \ref{BorelClassificationOfEquivariantBundlesForResolvableSingularitiesAndEquivariantStructure})
under the given assumptions.

\medskip
Notice that {\it every}
$\infty$-group is the shape
$\shape \, \Gamma$ of some Hausdorff topological group
$\Gamma$ (Prop. \ref{AllBareInfinityGroupsAreShapesOfTopologicalGroups}),
and that the assumption of
discrete equivariance group $G$ and
truncated shape $\shape \, \Gamma$
is generically verified in applications
of equivariant principal bundles
as twists for generalized cohomology of good orbifolds
(Ex. \ref{TruncatedStructureGroupsAndTheirEquivariantClassificationResults},
Ex. \ref{EquivariantBundlesServingAsGeoemtricTwistsOfEquivariantKTheory}):

\vspace{.1cm}
\begin{center}
\hypertarget{Table1}{}
\def\arraystretch{1.7}
\begin{tabular}{|c|c||c|c|}
\hline
  \multicolumn{2}{|c||}{
    \color{darkblue}
    \bf
       $G$-equivariant
       $\Gamma$-principal
       bundles:
  }
  &
 \multicolumn{2}{c|}{ \bf Structure group $\Gamma$ }
 \\
 \cline{3-4}
  \multicolumn{2}{|c||}{
    \color{darkblue}
    \bf
    examples and applications
  }
  &
  {\bf Compact Lie}
  &
  {\bf Truncated}
  \\
  \hline
  \hline
  \multirow{3}{*}{
    \rotatebox{90}{
      \bf
        Equivariance
        group $G$
      \hspace{-10pt}
    }
  }
  &
  \multirow{2}{*}{
    \begin{tabular}{c}
      {\bf Discrete}
    \end{tabular}
  }
  &
  \def\arraystretch{1}
  \begin{tabular}{l}
    \phantom{$\mathclap{\vert^{\vert^{\vert^{\vert^{\vert}}}}}$}
    -- torus-principal bundles
    \\
    \phantom{}
   -- finite covering spaces
    \\
    \phantom{}
    -- K-theory cocycles
    \\
    \phantom{}
    -- tangential twists
  \end{tabular}
  &
  \def\arraystretch{1}
  \begin{tabular}{l}
    \phantom{$\mathclap{\vert^{\vert^{\vert^{\vert^{\vert}}}}}$}
    -- torus-principal bundles
    \\
    \phantom{}
    \color{orangeii}
    -- ordinary cohomology twists
    \\
    \phantom{}
    \color{orangeii}
    -- K-theory twists
    \\
    \phantom{}
    \color{orangeii}
    -- differential twists
  \end{tabular}
  \\
  &&
  \multicolumn{2}{c|}{
    \raisebox{24pt}{
    $
      \underset{
        \mbox{
          on good orbifolds
        }
      }{
      \underbrace{
        \phantom{--------------------------}
      }
      }
    $
    }
  }
  \\
  \cline{2-4}
  &
  {\bf Compact Lie}
  &
  \multicolumn{2}{c|}
  {
    \def\arraystretch{1}
    \begin{tabular}{c}
      \phantom{$\mathclap{\vert^{\vert^{\vert}}}$}
      as above,
      \\
      on bad effective orbifolds \& in FHT theory
    \end{tabular}
  }
  \\
  \hline
\end{tabular}

\vspace{.1cm}

{\bf Table 1.}
\end{center}
\vspace{.1cm}

\medskip

\noindent
{\bf Proper equivariant cohesive homotopy theory.}
A remarkable aspect of
the classification result Thm. \ref{BorelClassificationOfEquivariantBundlesForResolvableSingularitiesAndEquivariantStructure}
is (see \cref{EquivariantModuliStacks})
that
 {\bf (1)} the
moduli $\infty$-stacks of equivariant principal bundles
always exist in {\it Borel-equivariant} cohesive homotopy theory,
and that
 {\bf (2)}
under the truncation condition on $\shape \, \Gamma$
also the corresponding classifying spaces exist in Borel equivariant homotopy theory.
While the first of these is a general fact,
the second is a happy coincidence
(often used without further amplification, e.g. in \cite{FreedHopkinsTeleman02ComplexCoefficients})
conditioned on stability and the truncation of $\shape \, \Gamma$.

General equivariant classifying spaces
require more fine-grained proper-equivariant homotopy theory
of spaces parameterized over the category
of orbits (Ntn. \ref{GOrbitCategory})
of the equivariance group
(e.g. \cite[\S 8]{tomDieck79}\cite{May96}\cite{Blu17}).
Following \cite{Rezk14}\cite{SS20OrbifoldCohomology}, we may understand this
(see \cref{GeneralSingularCohesion})
through an embedding of the category of topological $G$-spaces
(Ntn. \ref{GActionOnTopologicalSpaces})
into,
once again, a yet more convenient ambient category of $G$-actions, namely
the slice of a {\it singular-cohesive $\infty$-topos}
(Def. \ref{SingularCohesiveInfinityTopos})
over the generic
$G$-orbi-singularity $\orbisingularG$ (Ntn. \ref{Singularities}):
$$
  \hspace{-.6cm}
  \begin{tikzcd}
    \mathclap{
    \mbox{
      \tiny
      \bf
      \begin{tabular}{c}
        \color{gray}
        category of
        \\
        \color{darkblue}
        topological $G$-spaces
      \end{tabular}
    }
    }
    \ar[
      rrr,
      shift left=13pt,
      dashed,
      "{
        \mbox{
          \tiny
          \bf
          convenience
        }
        \;\;\;
      }"{description}
    ]
    &[-6pt]
    \mathclap{
    \mbox{
      \tiny
      \bf
      \begin{tabular}{c}
        \color{gray}
        $\infty$-topos of
        \\
        \color{darkblue}
        smooth $G$-$\infty$-actions
      \end{tabular}
    }
    }
    &[-6pt]
    \mathclap{
    \mbox{
      \tiny
      \bf
      \begin{tabular}{c}
        \color{gray}
        slice of cohesive $\infty$-topos of
        \\
        \color{darkblue}
        smooth $\infty$-groupoids over $\mathbf{B}G$
      \end{tabular}
    }
    }
    &
    \mathclap{
    \mbox{
      \tiny
      \bf
      \begin{tabular}{c}
        \color{gray}
        slice of singular-cohesive $\infty$-topos of
        \\
        \color{darkblue}
        $G$-orbi-singular $\infty$-groupoids
      \end{tabular}
    }
    }
    \\[-20pt]
    \categorybox{
      \!\!
      \Actions{G}(
        \kTopologicalSpaces
      )
      \!\!
    }
    \ar[r,  "\ContinuousDiffeology"]
    \ar[
      drrr,
      rounded corners,
      to path={
        -- ([yshift=-39pt]\tikztostart.south)
        --  node[swap, yshift=-7pt]{ \scalebox{.7}{
             \color{greenii}
             \bf
             underlying proper $G$-equivariant homotopy type
           } }
           ([xshift=-00pt]\tikztotarget.west)
      }
    ]
    &
    \categorybox{
      \!\!
      \Actions{G}(\SmoothInfinityGroupoids)
      \!\!
    }
    \ar[r, "\HomotopyQuotient{(-)}{G}", "\sim"{swap}]
    \ar[
      drr,
      rounded corners,
      to path={
        -- ([yshift=-39pt]\tikztostart.south)
        -- ([xshift=-00pt]\tikztotarget.west)
      }
    ]
    &
    \categorybox{
      \!\!
      \Slice{\SmoothInfinityGroupoids}{\mathbf{B}G}
      \!\!
    }
    \ar[
      from=r,
      to path={
        -- ([yshift=10pt]\tikztostart.west)
        -- node[swap, yshift=+5pt]{\scalebox{.7}{$
                \Smooth
              $}}
           ([xshift=+00pt, yshift=+10pt]\tikztotarget.east)
        -- node[xshift=15pt]{ $\orbisingular$ }
            ([xshift=+00pt, yshift=-10pt]\tikztotarget.east)
        -- node[yshift=-10pt]{ \scalebox{.7}{$
               \underset{
                 \raisebox{-3pt}{
                   \color{greenii}
                   \bf
                   orbi-singular
                 }
               }{
                 \Singularity
               }
              $}}
           ([yshift=-10pt]\tikztostart.west)
      }
    ]
    &
    \categorybox{
      \!\!
      \Slice{\GloballyEquivariant\SmoothInfinityGroupoids}{\orbisingularG}
      \!\!
    }
    \ar[
      r,
      to path={
        -- ([yshift=-10pt]\tikztostart.east)
        -- node[swap, yshift=-12pt]{\scalebox{.7}{$
               \underset{
                 \raisebox{-5pt}{
                   \color{greenii}
                   \bf
                   $G$-orbi-spatial
                 }
               }{
                G\Orbi\Smooth
               }
              $}}
           ([xshift=-10pt, yshift=-10pt]\tikztotarget.west)
        -- node[swap, xshift=-13pt]{ $\smoothrelativeG$ }
            ([xshift=-10pt, yshift=+10pt]\tikztotarget.west)
        -- node[swap, yshift=5pt]{ \scalebox{.7}{$G\Orbi\Space$}  }
           ([yshift=+10pt]\tikztostart.east)
      }
    ]
    \ar[
      to path={
        -- ([xshift=-10pt]\tikztostart.south)
        --
            node[below, sloped]{
              \scalebox{.7}{$
                \underset{
                  \raisebox{+0pt}{
                    \color{greenii}
                    \bf
                    shape
                  }
                }{\Shape}
              $}
            }
           ([xshift=-10pt, yshift=-24pt]\tikztostart.south)
        --
            node[above, yshift=5pt]{ $\shape$ }
            ([xshift=+10pt, yshift=-24pt]\tikztostart.south)
        --
            node[below, sloped]{ \scalebox{.7}{$\Discrete$}  }
           ([xshift=+10pt]\tikztotarget.south)
      }
    ]
    &[9pt]
    \!\!\!\!\!\!\!\!\!\!
    \categorybox{
      \!\!
      \GEquivariant\SmoothInfinityGroupoids
      \!\!
    }
    \!\!\!\!\!\!\!\!\!\!\!
    \ar[d, "\Shape"{sloped}]
    \\[+3pt]
    &
    &
    &
    \categorybox{
      \!\!
      \Slice{\GloballyEquivariant\InfinityGroupoids}{\orbisingularG}
      \!\!
    }
    \ar[r, "G\Orbi\Smooth"{swap}]
    &
    \!\!\!
    \categorybox{
      \!\!
      \GEquivariant\InfinityGroupoids
      \!\!
    }
    \!\!\!
    \\[-24pt]
    &&&
    \mathclap{
    \mbox{
      \tiny
      \bf
      \begin{tabular}{c}
        \color{gray}
        slice of base $\infty$-topos of
        \\
        \color{darkblue}
        global equivariant
        homotopy theory
      \end{tabular}
    }
    }
    &
    \mathclap{
    \mbox{
      \tiny
      \bf
      \begin{tabular}{c}
        \color{gray}
        sub-$\infty$-topos of
        \\
        \color{darkblue}
        proper $G$-equivariant
        \\
        \color{darkblue}
        homotopy theory
      \end{tabular}
    }
    }
  \end{tikzcd}
$$

\medskip

\noindent
{\bf Elmendorf's theorem as a generalized equivariant Oka principle.}
On the right, $\GEquivariant\InfinityGroupoids$
(Ex. \ref{ClassicalEquivariantHomotopyTheory})
denotes the $\infty$-presheaves over the {\it category of $G$-orbits} (Def.
\ref{GOrbitCategory}). This is the modern context of $G$-equivariant homotopy theory, traditionally motivated by {\it Elmendorf's theorem} (\cite{Elmendorf}\cite{DwyerKan84}, recalled as Prop. \ref{ElmendorfDwyerKanTheorem} below). With the hindsight of cohesive homotopy theory, this classical theorem may conceptually be understood as {\it enforcing} an ``equivariant Oka principle'': For $G \acts \TopologicalSpace$ a $G$-CW-complex (Ex. \ref{GCWComplexesAreCofibrantObjectsInProperEquivariantModelcategory}) and $G \acts Y$ any topological $G$-space, the simplicial enhancement of Elmendorf's theorem (\cite[Thm. 3.1]{DwyerKan84}, see also \cite{CordierPorter96}\cite[Thm. 1.3.8]{Blu17}) says, after the above embedding into singular-cohesive homotopy theory (by Prop. \ref{OrbiSpaceIncarnationOfGSpaceIsOrbisingularizationOfHomotopyQuotient} below), that one may take the shape operation inside the equivariant mapping stack {\it if}
in the process one enhances $G$-orbifolds $\HomotopyQuotient{\mathrm{X}}{G}$ to
their {\it orbi-singularization}
$\orbisingular(\HomotopyQuotient{\mathrm{X}}{G})$
(Ex. \ref{CohesiveFormulationOfEDKTheoremForDiscreteEquivarianceGroups} below):
\begin{equation}
  \label{ElmendorfTheorem}
  \overset{
    \mathclap{
    \raisebox{3pt}{
      \tiny
      \color{darkblue}
      \bf
      \def\arraystretch{.9}
      \begin{tabular}{c}
        \hspace{-.2cm}
        shape of equivariant
        mapping stack
      \end{tabular}
    }
    }
  }{
  \shape
  \,
  \SliceMaps{\big}{\mathbf{B}G}
    {
      \HomotopyQuotient
        { \mathrm{X} }
        { G }
    }
    {
      \HomotopyQuotient
        { \mathrm{Y} }
        { G }
    }
  }
  \;\;\;
  \underset{
    \mathclap{
    \raisebox{-12pt}{
      \tiny
      \color{greenii}
      \bf
      \def\arraystretch{.9}
      \begin{tabular}{c}
        Elmendorf-Dwyer-Kan
        \\
        theorem
      \end{tabular}
     }
    }
  }{
    \simeq
  }
  \;\;\;
  \smooth
  \overset{
    \raisebox{4pt}{
      \tiny
      \color{darkblue}
      \bf
      \begin{tabular}{c}
        mapping space of
        equivariant shapes
      \end{tabular}
    }
  }{
  \SliceMaps{\big}{\orbisingularG}
    {
      \shape
      \orbisingular
      (
      \HomotopyQuotient
        { \mathrm{X} }
        { G }
      )
    }
    {
      \shape
      \orbisingular
      (
      \HomotopyQuotient
        { \mathrm{Y} }
        { G }
      )
    }
    \mathrlap{\,.}
  }
\end{equation}

\medskip

\noindent
{\bf Proper equivariant classifying spaces.}
Therefore, chasing a topological $G$-space
through the above diagram produces its
usual proper-equivariant homotopy type as a
$G$-orbi-space (Prop. \ref{OrbiSpaceIncarnationOfGSpaceIsOrbisingularizationOfHomotopyQuotient}).
But we may now also feed the above moduli stack
$\HomotopyQuotient{\mathbf{B}\Gamma}{G}$ of $G$-equivariant $\Gamma$-principal
bundles through this machine, and we find
(Thm. \ref{MurayamaShimakawaGroupoidIsEquivariantModuliStack})
that the resulting
equivariant homotopy type is that of the Murayama-Shimakawa construction
(\cite{MurayamaShimakawa95}\cite{GuillouMayMerling17}, \cref{ConstructionOfUniversalEquivariantPrincipalBundles}):
$$
  \overset{
    \mathclap{
    \raisebox{6pt}{
      \tiny
      \color{darkblue}
      \bf
      \begin{tabular}{c}
        proper equivariant classifying space
        \\
        for equivariant principal bundles
      \end{tabular}
    }
    }
  }{
    \EquivariantClassifyingShape{G}{\Gamma}
    \;\;\coloneqq\;\;
    \smoothrelativeG
    \;\,
    \shape
    \,\,
    \orbisingular
    \,\,
    (\HomotopyQuotient{\mathbf{B}\Gamma}{G})
  }
  \;\;\;
    :
  \;\;\;
  G/H
  \;\;\;
  \longmapsto
  \;\;\;\;
  \overset{
    \mathclap{
    \raisebox{6pt}{
      \tiny
      \color{darkblue}
      \bf
      \begin{tabular}{c}
        concordances of equivariant
        \\
        bundles over the point
      \end{tabular}
    }
    }
  }{
    \shape
    \,
    \SliceMaps{}{\mathbf{B}G}
      { \mathbf{B}H }
      { \HomotopyQuotient{\mathbf{B}\Gamma}{G} }
  }
  \;\;\;
    \simeq
  \;\;\;
  \overset{
    \mathclap{
    \raisebox{4pt}{
      \tiny
      \color{darkblue}
      \bf
      \begin{tabular}{c}
        Murayama-Shimakawa construction
      \end{tabular}
    }
    }
  }{
    \SingularSimplicialComplex
    \,
    \TopologicalRealization{}
    {
      \Maps{}
        { \mathbf{E}G }
        { \mathbf{B}\Gamma }
    }^H
  }
$$

In the special case when
($G$-singularities are resolvable and)
$\shape \, \Gamma$ is truncated, the
above orbi-smooth Oka principle
\eqref{OrbiSmoothOkaPrincipleInIntroduction}
applies
to these values of the equivariant classifying space and gives
(Prop. \ref{EquivariantClassifyingShapeOfTruncatedTopologicalGroupsCoincidesWithThatOftheirShape}):

$$
\hspace{-2mm}
  \left.
  \arraycolsep=2pt
  \begin{array}{clll}
    &
    \mbox{\small $G$ discrete with resolvable singularities}
    \\
    \mbox{\small \&}
    &
    \mbox{\small $\Gamma$  topological of truncated shape}
    \\
    \mbox{\small \&}
    &
    \mbox{\small $G \acts \, \SmoothManifold$ a proper smooth $G$-manifold}
  \end{array}
  \right\}
  \;
  \vdash
  \quad
  \EquivariantClassifyingShape{G}{\Gamma}
  \;\simeq\;
  \EquivariantClassifyingShape{G}{(\shape \, \Gamma)}
  \;
  :
  \;
  G/H
  \;
  \longmapsto
  \quad
  \def\arraystretch{0}
  \begin{aligned}
    &   \SliceMaps{}{B G}
      { B H }
      { \HomotopyQuotient{B \Gamma}{G} }
   \\
    &
    \phantom{AAA}
    \simeq
    \;
   \SingularSimplicialComplex
    \,
    \TopologicalRealization{\big}
    {
      \Maps{}
        { \TopologicalRealization{}{\mathbf{E}G} }
        { \TopologicalRealization{}{\mathbf{B}\Gamma} }
    }
    ^H
     .
  \end{aligned}
$$

\medskip

\noindent
{\bf Classification statement in proper-equivariant cohomology.}
These proper-equivariant homotopy types are the coefficient
systems that represent Borel-equivariant cohomology inside proper-equivariant cohomology
(Prop. \ref{ProperEquivariantCohomologySubsumesBorelEquivariantCohomology}),
so that our classification result may be re-stated in the following
equivalent forms (Thm. \ref{ProperClassificationOfEquivariantBundlesForResolvableSingularitiesAndEquivariantStructure}):
$$
  \hspace{-1mm}
  \left.
  \arraycolsep=2pt
  \begin{array}{clll}
    &
    \small  \mbox{$G$ discrete with resolvable singularities}
    \\
    \mbox{\small \&}
    &
    \small    \mbox{$\Gamma$ Hausdorff of truncated shape}
    \\
    \mbox{\small \&}
    &
    \mbox{\small $G \acts \, \SmoothManifold$ a proper smooth $G$-manifold}
  \end{array}
  \right\}
  \;
  \vdash
  \;
  \left\{
  \hspace{-.6cm}
  \def\arraystretch{1.6}
  \begin{array}{lll}
    &
    \IsomorphismClasses{
      \EquivariantPrincipalFiberBundles{G}{\Gamma}
      (\DTopologicalSpaces)_{\SmoothManifold}^{\stable}
    }
    &
      \hspace{-.4cm}
      \raisebox{1pt}{
        \tiny
        \color{darkblue}
        \bf
        \def\arraystretch{.9}
        \begin{tabular}{l}
          isom. classes of
          equivariantly locally trivial
          \\
          stable equivariant principal topol. bundles
        \end{tabular}
      }
    \\
    &
    \;\;
    \simeq
    \;\;
      H^1_G
      (X;\, \shape\,\Gamma)
      \;
      =
      \;
      H^0_G
      (X;\, B \Gamma)
      &
      \hspace{-.35cm}
      \def\arraystretch{.9}
      \raisebox{4pt}{
        \tiny
        \color{darkblue}
        \bf
        \begin{tabular}{l}
          Borel-equivariant cohomology with
          \\
          coefficients in classifying space
        \end{tabular}
      }
    \\
    &
    \;\;
    \simeq
    \;\;
      H^0_{\scalebox{.7}{$\orbisingularG$}}
      \left(X;\, \EquivariantClassifyingShape{G}{(\shape\, \Gamma)}\right)
      \;
      \simeq
      \;
      H^0_{\scalebox{.7}{$\orbisingularG$}}
      \left(
        X;\,
        (\EquivariantClassifyingShape{G}{\Gamma})^{\stable}
      \right)
      &
      \hspace{-.4cm}
      \raisebox{1pt}{
        \tiny
        \color{darkblue}
        \bf
        \def\arraystretch{.9}
        \begin{tabular}{l}
          proper-equivariant cohomology with
          \\
          coefficients in equivariant classifying space.
        \end{tabular}
      }
  \end{array}
  \right.
$$
Specialized to trivial $G$ action on $\Gamma$ and the
cases where
$\Gamma$ is either a 1-truncated compact Lie group or
the projective unitary group $\PUH$,
this theorem reproduces
(when $\HomotopyQuotient{\SmoothManifold}{G}$ is a good orbifold with resolvable
singularities)
a series of statements found in the literature
(Ex. \ref{TruncatedStructureGroupsAndTheirEquivariantClassificationResults},
Ex. \ref{EquivariantBundlesServingAsGeoemtricTwistsOfEquivariantKTheory}).

\medskip
\medskip

\noindent
{\bf Outline.}
In order to make the presentation of these results reasonably self-contained
also for the non-expert reader, we lay out a fair bit of the required
background in \cref{EquivariantTopology} (equivariant differential topology)
and \cref{Generalities} (equivariant cohesive homotopy theory).
The new constructions and results are the content of
\cref{EquivariantPrincipalTopologicalBundles} (equivariant topological bundles)
and, particularly, of \cref{EquivariantInfinityBundles}
(equivariant $\infty$-bundles).

\medskip
\medskip

\noindent
{\bf Conclusion and outlook.}
In conclusion, we find that the unifying mechanism
behind a large class of equivariant classification results and their
generalization to higher truncated structure groups is a smooth Oka principle
in cohesive homotopy theory, which seamlessly embeds the theory of equivariant bundles
into a transparent modal homotopy theory of higher geometry,
in particular into the context of higher principal bundles over orbifolds
and more general cohesive orbispaces.

\medskip
By way of outlook, notice that
there is a large supply of equivariance groups with resolvable singularities
(Prop. \ref{ExistenceOfSmoothSphericaSpaceForms})
including key examples of interest in applications
(Ex. \ref{ADEGroupsHaveSphericalSpaceForms});
and
there is a large supply of truncated structure groups,
as every $\infty$-group is the shape
of some Hausdorff topological group (Prop. \ref{AllBareInfinityGroupsAreShapesOfTopologicalGroups}).
In particular, for $R$ a ring spectrum and $\mathrm{GL}(1,R)$ its
$\infty$-group of units (see \cite[Ex. 2.37]{FSS20CharacterMap} for pointers),
there is for every $n$ a
Hausdorff group $\Gamma$ with shape its $n$-truncation $\shape \, \Gamma \;\simeq\;
\Truncation{n} \mathrm{GL}(1,R)$. The corresponding $G$-equivariant $\Gamma$-principal
bundles are candidate geometric twists for equivariant $R$-cohomology theory,
generalizing the archetypical case of twisted complex K-theory
(where $R = \KU$, $\Gamma = \GradedPUH$
with $\shape \, \GradedPUH \,\simeq\, \Truncation{2} \, \mathrm{GL}(1,\KU)$).

\medskip

In all these cases, the classification theorem of \cref{EquivariantModuliStacks}
shows, in particular, that the equivariant classifying spaces of twists are, generically,
equivariantly non-simply connected,
which means that many traditional tools, notably
of rational homotopy theory, do not apply without extra care.
For example, in the case of twisted complex K-theory,
this fact
(see Ex. \ref{OrbiSmoothOkaPrincipleForPUHCoefficientsOverThePoint}
and \eqref{EquivariantHomotopyGroupsOfClassifyingShapeOfStableEquivariantPUHBundles}),
explains, we claim,
the otherwise somewhat unexpected
(cf. \cite[(3.22)]{FreedHopkinsTeleman02ComplexCoefficients})
appearance of local systems in the twisted equivariant Chern character
(\cite[Def. 3.10]{TuXu06}). Generally, one may use the equivariant
classifying theory developed here to give a general construction of
twisted equivariant Chern-Dold character maps and hence of
twisted equivariant differential cohomology theories, in
equivariant generalization (along the lines of \cite[\S 3]{SS20EquivariantTwistorial})
of the construction in \cite[Def. 5.4]{FSS20CharacterMap}:
$$
  \begin{tikzcd}[column sep=huge]
    \overset{
      \mathclap{
      \raisebox{3pt}{
        \tiny
        \color{darkblue}
        \bf
        \def\arraystretch{.9}
        \begin{tabular}{c}
          universal local coefficient bundle
          \\
          for $\Gamma$-twisted $G$-equivariant
          \\
          $A$-cohomology theory
        \end{tabular}
      }
      }
    }{
      \shape
      \,
      \orbisingular
      \,
      \left(
        \HomotopyQuotient{A}{(\HomotopyQuotient{\Gamma}{G})}
      \right)
    }
    \ar[
      rr,
      "{
        \mbox{\hspace{-7mm}
          \tiny
          \begin{tabular}{c}
            \color{greenii}
            \bf
            proper equivariant sliced rationalization
            \\
            representing the
            \\
            \color{greenii}
            \bf
            twisted equivariant $A$-character map
          \end{tabular}
        }
      }"{swap, yshift=-7.5pt, xshift=12pt}
    ]
    \ar[d]
    &&
    L_{\mathbb{Q}}
    \left(
      \shape
      \,
      \orbisingular
      \,
      \left(
        \HomotopyQuotient{A}{(\HomotopyQuotient{\Gamma}{G})}
      \right)
    \right)
    \ar[d]
    \\
    \underset{
      \mathclap{
      \def\arraystretch{.9}
      \raisebox{-5pt}{
        \tiny
        \color{darkblue}
        \begin{tabular}{c}
          \bf
          equivariant classifying
          \\
          {\bf
          space of twists }
          {\color{black}
          (\cref{EquivariantModuliStacks})
          }
        \end{tabular}
      }
      }
    }{
      \EquivariantClassifyingShape{G}{\Gamma}
    }
    \ar[rr]
    &&
    \underset{
      \raisebox{-3pt}{
        \tiny
        \color{darkblue}
        \bf
        \begin{tabular}{c}
        \end{tabular}
      }
    }{
    L_{\mathbb{Q}}
    \left(
      \EquivariantClassifyingShape{G}{\Gamma}
    \right)
    }
  \end{tikzcd}
$$
The homotopy pullback
(formed in $\SingularSmoothInfinityGroupoids$)
of stacks of flat equivariant $L_\infty$-algebra valued differential forms
(\cite[Def. 3.57]{SS20EquivariantTwistorial})
along this rationalization operation
(in direct equivariant generalization of \cite[Def. 4.38]{FSS20CharacterMap})
solves the open problem of providing a general construction of
$G$-equivariant $\Gamma$-twisted {\it differential} $A$-cohomology.
This application will be discussed in detail in
\cite{TwistedEquivariantChernCharacter}\cite{TwistedEquivariantDifferentialCohomology}.

\newpage

\section{Tools and Techniques}
\label{ToolsAndTechniques}

\noindent
{\bf Higher cohesive geometry as intrinsically equivariant geometry.}
The point of
{\it higher homotopical geometry}
is (see \cite[p. 4-5]{SS20OrbifoldCohomology})
that the notion and presence of {\it gauge transformations} (homotopies) and
{\it higher gauge-of-gauge transformations}
is natively built into
the theory, so that absolutely every concept formulated in higher
geometry is {\it intrinsically} equivariant with respect to all relevant
symmetries.
This makes higher geometry the natural context for laying foundations
for equivariant algebraic topology.

\vspace{-3mm}
\begin{equation}
  \label{GaugeTransformations}
  \overset{
    \mathclap{
    \raisebox{4pt}{
      \tiny
      \color{darkblue}
      \bf
      \def\arraystretch{.9}
      \begin{tabular}{c}
        $\Sigma$-shaped probes of
        \\
        higher geometric space $\mathcal{X}$
        \\
        \phantom{a}
      \end{tabular}
    }
    }
  }{
    \underset{
      \mathclap{
      \raisebox{-4pt}{
        \tiny
        \color{darkblue}
        \bf
        \def\arraystretch{.9}
        \begin{tabular}{c}
        \\
          value of $\infty$-stack $\mathcal{X}$
          \\
          on site-object $\Sigma$
        \end{tabular}
      }
      }
    }{
      \mathcal{X}(\Sigma)
    }
  }
  \qquad
    =
  \;\;
  \left\{\quad
  {\phantom{\mbox{\tiny\bf domain}}}
  \begin{tikzcd}[column sep=large]
    \mathllap{
      \mbox{
        \tiny
        \color{darkblue}
        \bf
        \begin{tabular}{c}
          probe
          \\
          space
          \\
         (`` brane")
        \end{tabular}
      }
      \!\!\!\!
    }
    \Sigma
    \ar[
      rr,
      bend left=60,
      "\mbox{\tiny\color{greenii}\bf configuration}"{description},
      "\ "{name=sl, xshift=-12pt,below},
      "\ "{name=sr, xshift=+12pt,below}
    ]
    \ar[
      rr,
      bend right=60,
      "\mbox{\tiny\color{greenii}\bf configuration}"{description},
      "\ "{name=tl, xshift=-12pt, above},
      "\ "{name=tr, xshift=+12pt, above}
    ]
    \ar[
      from=sr,
      to=tr,
      Rightarrow,
      bend left=50,
      "\mbox{\tiny\color{orangeii}\bf gauge transf.}"{sloped, description},
      "\ "{name=t2, left}
    ]
    \ar[
      from=sl,
      to=tl,
      Rightarrow,
      bend right=50,
      "\rotatebox{0}{\tiny\color{orangeii}\bf gauge transf.}"{sloped, description},
      "\ "{name=s2, right}
    ]
    \ar[
      from=s2,
      to=t2,
      dashed,
      Rightarrow,
      "\mbox{\tiny\color{purple}\bf gauge-of-gauge}"{above},
      "\mbox{\tiny\color{purple}\bf transformations}"{below}
    ]
    \ar[
      from=s2,
      to=t2,
      -,
      dashed
    ]
    &{\phantom{A}}&
    \mathcal{X}
    \mathrlap{
      \!\!\!\!
      \mbox{
        \tiny
        \color{darkblue}
        \bf
        \begin{tabular}{c}
          field
          \\
          space
        \end{tabular}
      }
    }
  \end{tikzcd}
  {\phantom{\mbox{\tiny\bf codomain}}}
  \!\!\!\!\!\!\right\}
  \;\;
  \;\in\; \InfinityGroupoids
  \,.
\end{equation}

\medskip
\noindent
{\bf The machinery of $\infty$-category theory.}
While higher stacks and their higher geometry are often perceived as an
esoteric and convoluted mathematical subject, this is rather a property
of their {presentation} by models in simplicial homotopy theory
(reviewed in \cref{AbstractHomotopyTheory}, \cite[\S A]{FSS20CharacterMap})
and must be understood as a reflection of
their precious richness, not as of their intractability.
Indeed, with
{\it $\infty$-category theory}
\cite{Joyal08Notes}\cite{Joyal08Theory}\cite{Lurie09HTT}\cite{Cisinski19}\cite{RiehlVerity21}
(see \cref{AbstractHomotopyTheory}),
and specifically with
{\it $\infty$-topos theory}
\cite{Simpson99}\cite{Lurie03}\cite{ToenVezzosi05}\cite{Joyal08Logoi}\cite{Lurie09HTT}\cite{Rezk10}
(see \cref{ToposTheory}),
there is a high-level language,
abstracting away from the zoo of models (e.g. \cite{Bergner07Survey}),
to admit efficient reasoning about higher stacks via elementary categorical logic.
This point is made fully manifest
by the existence of an elementary internal logic of $\infty$-toposes \cite{Shulman19},
now known as Homotopy Type Theory \cite{UFP13}
(in our context see \cite{SchreiberShulman14}\cite[p. 5]{SS20OrbifoldCohomology})
which condenses all such reasoning to coding in a kind of programming language.

\medskip

For example, once Cartesian (pullback) squares are understood as
homotopy Cartesian squares, namely filled with a homotopy which
exhibits the expected unique factorization property up to homotopy and in the sense of
a contractible space of homotopy-factorizations:
\begin{equation}
  \label{HomotopyCartesianSquare}
  \begin{tikzcd}
    X \times_B Y
    \ar[
      r,
      "\ "{swap, pos=.8, name=s}
    ]
    \ar[d]
    &
    Y
    \ar[d]
    \\
    X
    \ar[
      r,
      "\ "{pos=.2, name=t},
      "{
        \mbox{
          \tiny
          \color{darkblue}
          \bf
          \begin{tabular}{c}
            homotopy Cartesian square
            {\color{black}/}
            \\
            homotopy pullback square
          \end{tabular}
        }
      }"{swap, yshift=-6pt}
    ]
    &
    B
    \ar[
      from=s,
      to=t,
      Rightarrow,
      "{ \mbox{\tiny\rm(pb)} }"
    ]
  \end{tikzcd}
  \;\;\;
  \Rightarrow
  \qquad
  \left(
  \begin{tikzcd}
    Q
    \ar[
      drr,
      bend left=30,
      "\ "{swap, name=s2, pos=.7}
    ]
    \ar[
      ddr,
      bend right=30,
      "\ "{name=t2, pos=.7}
    ]
    \\[-15pt]
    &[-15pt]
    &
    Y
    \ar[d]
    \\
    &
    X
    \ar[
      r,
      "\ "{pos=.2, name=t}
    ]
    &
    B
    \ar[
      from=s2,
      to=t2,
      Rightarrow
    ]
  \end{tikzcd}
  \;\;\;
  \Leftrightarrow
  \;\;\;
  \begin{tikzcd}
    Q
    \ar[
      drr,
      bend left=30,
      "\ "{swap, name=s2, pos=.4}
    ]
    \ar[
      ddr,
      bend right=30,
      "\ "{name=t2, pos=.35}
    ]
    \ar[
      dr,
      dashed,
      "\ "{name=t3},
      "\ "{swap, name=s3},
    ]
    \\[-10pt]
    &[-15pt]
    X \times_B Y
    \ar[
      r,
      "\ "{swap, pos=.4, name=s}
    ]
    \ar[d]
    &
    Y
    \ar[d]
    \\
    &
    X
    \ar[
      r,
      "\ "{pos=.2, name=t}
    ]
    &
    B
    \ar[
      from=s2,
      to=t3,
      Rightarrow
    ]
    \ar[
      from=s3,
      to=t2,
      Rightarrow
    ]
    \ar[
      from=s,
      to=t,
      Rightarrow
    ]
  \end{tikzcd}
  \right)
  \,,
\end{equation}
then they follow patterns familiar from 1-category theory:
For instance, the {\it pasting law} in 1-categories (recalled as Prop. \ref{PastingLaw} below) continues to hold verbatim in any $\infty$-category \cite[Lem. 4.4.2.1]{Lurie09HTT}:
\begin{equation}
  \label{HomotopyPastingLaw}
  \begin{tikzcd}
    X
    \ar[
      r,
      "\ "{swap, name=s1, pos=.9}
    ]
    \ar[d]
    &
    Y
    \ar[d]
    \ar[
      r,
      "\ "{swap, name=s2, pos=.9}
    ]
    &
    Z
    \ar[d]
    \\
    A
    \ar[
      r,
      "\ "{name=t1, pos=.1}
    ]
    &
    B
    \ar[
      r,
      "\ "{name=t2, pos=.1}
    ]
    &
    C
    \ar[
      from=s1,
      to=t1,
      Rightarrow,
      "\mbox{\tiny\rm(pb)}"{description}
    ]
    \ar[
      from=s2,
      to=t2,
      Rightarrow,
    ]
  \end{tikzcd}
  \quad \Rightarrow\quad
  \left(
  \begin{tikzcd}
    X
    \ar[
      rr,
      "\ "{swap, name=s, pos=.9}
    ]
    \ar[d]
    &
    &
    Z
    \ar[d]
    \\
    A
    \ar[
      rr,
      "\ "{name=t, pos=.1}
    ]
    &
    &
    C
    \ar[
      from=s,
      to=t,
      Rightarrow,
      "\mbox{\tiny\rm(pb)}"{description}
    ]
  \end{tikzcd}
  \;\;\Leftrightarrow\;\;
  \begin{tikzcd}
    X
    \ar[
      r,
      "\ "{swap, name=s1, pos=.9}
    ]
    \ar[d]
    &
    Y
    \ar[d]
    \ar[
      r,
      "\ "{swap, name=s2, pos=.9}
    ]
    &
    Z
    \ar[d]
    \\
    A
    \ar[
      r,
      "\ "{name=t1, pos=.1}
    ]
    &
    B
    \ar[
      r,
      "\ "{name=t2, pos=.1}
    ]
    &
    C
    \ar[
      from=s1,
      to=t1,
      Rightarrow,
      "\mbox{\tiny\rm(pb)}"{description}
    ]
    \ar[
      from=s2,
      to=t2,
      Rightarrow,
      "\mbox{\tiny\rm(pb)}"{description}
    ]
  \end{tikzcd}
  \right)
  \,;
\end{equation}
and it still makes sense, for instance, to say that a morphism $f$ is a {\it monomorphism}
if and only if its homotopy fiber product with itself is equivalently its domain
(\cite[p. 575]{Lurie09HTT}\cite[p. 10]{Rezk19}, see also Ex. \ref{MonomorphismsOfInfinityGroupoids}):

\vspace{-.3cm}

\begin{equation}
  \label{InfinityMonomorphism}
  \begin{tikzcd}
    X
    \ar[
      r,
      hook,
      "f",
      "{
        \mbox{
          \tiny
          \rm
          \color{greenii}
          monomorphism
        }
      }"{swap}
    ]
    &[+12pt]
    Y
  \end{tikzcd}
  \;\;\;\;\;\;\;\;
  \Leftrightarrow
  \;\;\;\;\;\;\;\;
  \begin{tikzcd}[column sep=large]
    X
    \ar[r, "\mathrm{id}", "\ "{swap, pos=.9, name=s}]
    \ar[d, "\mathrm{id}"{swap}]
    &
    X
    \ar[d, "f"]
    \\
    X
    \ar[r, "f"{swap}, "\ "{pos=.1, name=t}]
    &
    Y
    \mathrlap{\,.}
    \ar[
      from=s,
      to=t,
      Rightarrow,
      "{ \mbox{\tiny\rm(pb)} }"{description}
    ]
  \end{tikzcd}
\end{equation}

\noindent
(In the following we will leave the homotopies filling these squares notationally implicit.)

\medskip
As another example: For any $\infty$-category $\mathbf{C}$ and any pair
$X, A$ of its objects, we have their {\it hom(omorphism) $\infty$-groupoids}
(see \cite{DuggerSpivak11})
\vspace{-2mm}
\begin{equation}
  \label{HomSpace}
  \overset{
    \mathclap{
    \raisebox{3pt}{
      \tiny
      \color{darkblue}
      \bf
      \begin{tabular}{c}
        hom-$\infty$-groupoid
        \\
        from $X$ to $A$
        \\
        in $\infty$-category $\mathbf{C}$
        \\
        \phantom{a}
      \end{tabular}
    }
    }
  }{
    \mathbf{C}(X,A)
  }
  \quad
  =
  \;\;
  \left\{
  {\phantom{\mbox{\tiny\bf domain}}}
  \begin{tikzcd}[column sep=large]
    \mathllap{
      \mbox{
        \tiny
        \color{darkblue}
        \bf
        \begin{tabular}{c}
          domain
          \\
          object
        \end{tabular}
      }
      \!\!\!\!
    }
    X
    \ar[
      rr,
      bend left=60,
      "\mbox{\tiny\color{greenii}\bf homomorphism}"{description},
      "\ "{name=sl, xshift=-8pt,below},
      "\ "{name=sr, xshift=+8pt,below}
    ]
    \ar[
      rr,
      bend right=60,
      "\mbox{\tiny\color{greenii}\bf homomorphism}"{description},
      "\ "{name=tl, xshift=-8pt, above},
      "\ "{name=tr, xshift=+8pt, above}
    ]
    \ar[
      from=sr,
      to=tr,
      Rightarrow,
      bend left=50,
      "\rotatebox{180}{\tiny\color{orangeii}\bf homotopy}"{sloped, description},
      "\ "{name=t2, left}
    ]
    \ar[
      from=sl,
      to=tl,
      Rightarrow,
      bend right=50,
      "\rotatebox{0}{\tiny\color{orangeii}\bf homotopy}"{sloped, description},
      "\ "{name=s2, right}
    ]
    \ar[
      from=s2,
      to=t2,
      dashed,
      Rightarrow,
      "\mbox{\tiny\color{purple}\bf higher}"{above},
      "\mbox{\tiny\color{purple}\bf homotopies}"{below}
    ]
    \ar[
      from=s2,
      to=t2,
      -,
      dashed
    ]
    &{\phantom{A}}&
    A
    \mathrlap{
      \!\!\!\!
      \mbox{
        \tiny
        \color{darkblue}
        \bf
        \begin{tabular}{c}
          codomain
          \\
          object
        \end{tabular}
      }
    }
  \end{tikzcd}
  {\phantom{\mbox{\tiny\bf codomain}}}
  \right\}
  \;\;
  \;\in\; \InfinityGroupoids
  \,,
\end{equation}
which are well-defined (model-independent) up to
(weak homotopy-)equivalence,
and homotopy functorial in the two arguments in the expected way.

\medskip
For instance, an $\infty$-functor is {\it fully faithful} if it
induces a natural equivalence on values of hom-$\infty$-functors:
\vspace{-2mm}
\begin{equation}
  \label{FullyFaithfulInfinityFunctor}
  \begin{tikzcd}
    \mathbf{C}
    \ar[
      rr,
      hook,
      "{F}",
      "{
        \mbox{
          \tiny
          \color{greenii}
          \bf
          fully faithful
        }
      }"{swap}
    ]
    &&
    \mathbf{D}
  \end{tikzcd}
  \hspace{1cm}
  \Leftrightarrow
  \hspace{1cm}
  \begin{tikzcd}
    \mathbf{C}(-,\,-)
    \ar[
      rr,
      "F_{(-,-)}",
      "\sim"{swap}
    ]
    &&
    \mathbf{D}
    \left(
      F(-)
      ,\,
      F(-)
    \right)\,.
  \end{tikzcd}
\end{equation}

\vspace{-1mm}
\noindent
Moreover, the hom $\infty$-functor satisfies the expected category-theoretic properties
in homotopy-theoretic form:

\vspace{1mm}
\noindent
{\bf (i)}
hom-$\infty$-groupoids respect homotopy (co)-limits via natural equivalences of the form
\vspace{-1mm}
\begin{equation}
  \label{HomFunctorRespectsLimits}
  \mathbf{C}
  \Big(
    \underset{\underset{i \in \mathcal{I}}{\longrightarrow}}{\lim}
    \,
    X_i,
    \,
    \underset{\underset{j \in \mathcal{J}}{\longleftarrow}}{\lim}
    \,
    A_j
  \Big)
  \;\;
  \simeq
  \;\;
  \underset{\underset{i \in \mathcal{I}}{\longleftarrow}}{\lim}
  \;
  \underset{\underset{j \in \mathcal{J}}{\longleftarrow}}{\lim}
  \;
  \mathbf{C}
  \left(
    X_i,
    \,
    A_j
  \right)
  \,;
\end{equation}

\vspace{-2mm}
\noindent
{\bf (ii)}
for a pair of {\it adjoint functors} between $\infty$-categories,
there is a natural equivalence between these hom-$\infty$-groupoids \eqref{HomSpace},
of the usual form
\cite[p. 159]{Joyal08Theory}\cite[Def. 5.2.2.7]{Lurie09HTT}\cite[Prop. F.5.6]{RiehlVerity21}:
\vspace{-2mm}
\begin{equation}
  \label{AdjunctionAndHomEquivalence}
  \begin{tikzcd}[column sep=large]
    \mathbf{D}
    \ar[
      rr,
      shift right=7pt,
      "R"{description},
      "\mbox{\tiny\color{greenii}\bf right adjoint}"{below,yshift=-1pt}
    ]
    \ar[
      rr,
      phantom,
      "\scalebox{.7}{$\bot$}"
    ]
    &&
    \mathbf{C}
    \ar[
      ll,
      shift right=7pt,
      "L"{description},
      "\mbox{\tiny\color{greenii}\bf left adjoint}"{above,yshift=1pt}
    ]
  \end{tikzcd}
  \;\;\;\;\;\;\;\;
  \Leftrightarrow
  \;\;\;\;\;\;\;\;
  \underset{
    \raisebox{-2pt}{
      \tiny
      \color{darkblue}
      \bf
      natural equivalence of hom $\infty$-groupoids
    }
  }{
    \mathbf{C}
    \left(
      C, R(D)
    \right)
    \;\;
    \simeq
    \;\;
    \mathbf{D}
    \left(
      L(C), D
    \right)
  }
  \,;
\end{equation}

\noindent
{\bf (iii)} in usual consequence, right (left) $\infty$-adjoint
$\infty$-functors
preserve
$\infty$-limits ($\infty$-colimits) via natural equivalences:
\begin{equation}
  \label{InfinityAdjointPreservesInfinityLimits}
  R
  \big(\,
    \underset{
      \underset{i \in \mathcal{I}}{\longleftarrow}
    }{\mathrm{lim}}
    \,
    X_i
  \big)
  \;\;
  \simeq
  \;\;
  \underset{
    \underset{i \in \mathcal{I}}{\longleftarrow}
  }{\mathrm{lim}}
  \;
  R
  (
    X_i
  )
  \,,
  \phantom{AAAAA}
  L
  \big(\,
    \underset{
      \underset{i \in \mathcal{I}}{\longrightarrow}
    }{\mathrm{lim}}
    \,
    X_i
  \big)
  \;\;
  \simeq
  \;\;
  \underset{
    \underset{i \in \mathcal{I}}{\longrightarrow}
  }{\mathrm{lim}}
  \;
  L
  (
    X_i
  )
  \,.
\end{equation}

Therefore, with a good supply of systems of adjoint
$\infty$-functors \eqref{AdjunctionAndHomEquivalence}
--
which we gain by invoking {\it modal} and specifically
{\it cohesive} homotopy theory
(\cite[\S 3.1]{SSS09}\cite{dcct}\cite{SS20OrbifoldCohomology},
see \cref{CohesiveHomotopyTheory} below)
--
many proofs in higher geometry,
which may be formidable when done in simplicial components,
are reduced to formal manipulations
yielding strings of such natural equivalences.
This is how we prove the main theorems in
\cref{EquivariantInfinityBundles}.

\medskip

\noindent
{\bf Homotopy fibers and $\infty$-groups.}
While $\infty$-category theory thus parallels category theory in the abstract,
a key difference is that fiber sequences in $\infty$-categories
(whenever they exist)
are, generically, {\it long}:
\vspace{-2mm}
\begin{equation}
  \label{GenericHomotopyFiberSequence}
  \begin{tikzcd}[row sep=large, column sep=huge]
    &
    \cdots
    \ar[r, "\mathrm{fib}^6(f)"]
    \ar[
      d,
      phantom,
      ""{coordinate, name=t1}
    ]
    &
    \Omega^2 C
    \ar[
      dll,
      rounded corners,
      "\mathrm{fib}^5(f)"{description},
      to path={
        -- ([xshift=8pt]\tikztostart.east)
        |- (t1) [pos=1]\tikztonodes
        -| ([xshift=-8pt]\tikztotarget.west)
        -- (\tikztotarget)
      }
    ]
    &
    \\
    \Omega A
    \ar[r, "\mathrm{fib}^4(f)"]
    &
    \Omega B
    \ar[r, "\mathrm{fib}^3(f)"]
    \ar[
      d,
      phantom,
      ""{coordinate, name=t}
    ]
    &
    \Omega C
    \ar[
      dll,
      rounded corners,
      "\mathrm{fib}^2(f)"{description},
      to path={
        -- ([xshift=8pt]\tikztostart.east)
        |- (t) [pos=1]\tikztonodes
        -| ([xshift=-8pt]\tikztotarget.west)
        -- (\tikztotarget)
      }
    ]
    \\
    A
    \ar[r, "\mathrm{fib}(f)"]
    &
    B
    \ar[r, "f"]
    &
    C
    \,.
  \end{tikzcd}
\end{equation}
In particular, the (homotopy-)fiber of a point inclusion is not
in general trivial (as it necessarily is in 1-category theory),
but is the {\it looping}
\vspace{-2mm}
\begin{equation}
  \label{LoopingInIntroduction}
  \begin{tikzcd}[column sep=large]
    \Omega_x X
    \ar[r]
    \ar[d]
    \ar[
      dr,
      phantom,
      "\mbox{\tiny\rm(pb)}"
    ]
    &
    \ast
    \ar[
      d,
      "x"{right}
    ]
    \\
    \ast
    \ar[
      r,
      "x"{below}
    ]
    &
    X
    \mathrlap{\,,}
  \end{tikzcd}
\end{equation}
whence the long homotopy fiber sequences
\eqref{GenericHomotopyFiberSequence}
follow by the pasting law \eqref{HomotopyPastingLaw}:
\vspace{-1mm}
$$
  \begin{tikzcd}[row sep=25pt, column sep=60pt]
    \Omega^2 C
    \ar[r]
    \ar[
      d,
      "\scalebox{.85}{$\mathrm{fib}^5(f)$}"{description}
    ]
    \ar[
      dr,
      phantom,
      "\mbox{\tiny\rm(pb)}"
    ]
    &
    \ast
    \ar[d]
    \\
    \Omega A
    \ar[d]
    \ar[
      r,
      "\scalebox{.85}{$\mathrm{fib}^4(f)$}"{description}
    ]
    \ar[
      dr,
      phantom,
      "\mbox{\tiny\rm(pb)}"
    ]
    &
    \Omega B
    \ar[r]
    \ar[
      d,
      "\scalebox{.85}{$\mathrm{fib}^3(f)$}"{description}
    ]
    \ar[
      dr,
      phantom,
      "\mbox{\tiny\rm(pb)}"
    ]
    &
    \ast
    \ar[d]
    \\
    \ast
    \ar[r]
    &
    \Omega C
    \ar[
      r,
      "\scalebox{.85}{$\mathrm{fib}^2(f)$}"{description}
    ]
    \ar[d]
    \ar[
      dr,
      phantom,
      "\mbox{\tiny\rm(pb)}"
    ]
    &
    A
    \ar[r]
    \ar[
      d,
      "\scalebox{.85}{$\mathrm{fib}(f)$}"{description}
    ]
    \ar[
      dr,
      phantom,
      "\mbox{\tiny\rm(pb)}"
    ]
    &
    \ast
    \ar[d]
    \\
    &
    \ast
    \ar[
      r
    ]
    &
    B
    \ar[
      r,
      "\scalebox{.85}{$f$}"{description}
    ]
    &
    C
    \mathrlap{\,.}
  \end{tikzcd}
$$

\medskip
\noindent
{\bf $\infty$-Toposes.}
In $\infty$-categories $\Topos$ of higher stacks \eqref{GaugeTransformations},
namely in {\it $\infty$-toposes}
\cite{Lurie03}\cite{ToenVezzosi05}\cite{Lurie09HTT}\cite{Rezk10},
every group object
($\infty$-group) $\mathcal{G} \,\in\, \Groups(\Topos)$ arises,
uniquely up to equivalence, as the looping \eqref{LoopingInIntroduction}
of its {\it delooping stack} $\mathbf{B}\mathcal{G} \,\in\, \Topos$
(recalled in Prop. \ref{GroupsActionsAndFiberBundles}):
\vspace{-1mm}
\begin{equation}
  \label{LoopingDeloopingInIntroduction}
  \mathcal{G} \;\simeq\; \Omega_\ast \mathbf{B} \mathcal{G}
  \,.
\end{equation}
Accordingly,
the connected components $\Truncation{0}(-)$ (see Prop. \ref{nTruncation} below)
of hom-$\infty$-groupoids
\eqref{HomSpace}
in $\infty$-toposes $\Topos$
may be understood
as (non-abelian, generalized) {\it cohomology theories}
(\cite[\S 2]{FSS20CharacterMap}\cite[p. 6]{SS20OrbifoldCohomology}):
\vspace{0mm}
$$
  \underset{
    \raisebox{-5pt}{
      \tiny
      \color{darkblue}
      \bf
      \def\arraystretch{.9}
      \begin{tabular}{c}
        domain
        \\
        space
      \end{tabular}
    }
  }{
    X
  },
  \underset{
    \raisebox{-5pt}{
      \tiny
      \color{darkblue}
      \bf
      \def\arraystretch{.9}
      \begin{tabular}{c}
        coefficient stack{\color{black}/}
        \\
        classifying stack
      \end{tabular}
    }
  }{
    A
  }\,\in\, \Topos
  \;\;\;\;\;\;\;\;\;\;\;\;\;\;\;\;\;
  \vdash
  \;\;\;\;\;\;\;\;\;\;\;\;\;\;\;\;\;
  \underset{
    \mathclap{
    \raisebox{-4pt}{
      \tiny
      \color{orangeii}
      \bf
      \begin{tabular}{c}
        $A$-cohomology of $X$
        \\
        in degree $\bullet$
      \end{tabular}
    }
    }
  }{
    H^{\bullet}
    \scalebox{1.1}{$($}
      X,
      \,
      A
    \scalebox{1.1}{$)$}
  }
  \;\;\coloneqq\;\;
  \tau_0
  \,
  \Topos
  \scalebox{1.1}{$($}
    X,
    \,
    \Omega^{-\bullet} A
  \scalebox{1.1}{$)$}
  \,.
$$
Specifically, for $\mathcal{G} \,\in\, \Groups(\Topos)$
we have {\it first non-abelian cohomology sets}:
$$
  H^1(X;\, \mathcal{G})
  \;=\;
  \Truncation{0}
  \Topos(X; \mathbf{B}\mathcal{G})
  \,.
$$

Moreover, the {\it fundamental theorem of $\infty$-topos theory}
(\cite[Prop. 6.5.3.1]{Lurie09HTT}, see around Prop. \ref{SliceInfinityTopos} below) says that for every object
$B \,\in\, \Topos$ in an $\infty$-topos, the slice $\infty$-category
$\SliceTopos{B}$, with
\begin{equation}
  \label{SliceHomInIntroduction}
  \overset{
    \mathclap{
    \raisebox{6pt}{
      \tiny
      \color{darkblue}
      \bf
      \begin{tabular}{c}
        slice hom $\infty$-groupoid
        \\
        in $\Topos$ over $B$
      \end{tabular}
    }
    }
  }{
    \SlicePointsMaps{\big}{B}
      { (X_1, p_1) }
      { (X_2, p_2) }
  }
  \;\;
    =
  \;\;
  \left\{
  \begin{tikzcd}
    X_1
    \ar[
      dr,
      " "{name=t2, pos=.8},
      " "{name=t2prime, pos=.6},
      "p_1"{swap},
      shorten=-4pt
    ]
    \ar[
      rr,
      bend left=30pt,
      "\ "{name=s, swap},
      "\ "{name=sprime, pos=.7, swap},
      shorten=-4pt
    ]
    \ar[orangeii,
      from=sprime,
      to=t2prime,
      dashed,
      Rightarrow,
      bend left=15pt,
      shorten=-3pt
    ]
    \ar[
      rr,
      bend right=30pt,
      "\ "{name=t},
      "\ "{name=s2,
      pos=.65, swap},
      crossing over,
      shorten=-4pt
    ]
    &&
    X_2
    \ar[
      dl,
      "p_2",
      shorten=-4pt
    ]
    \\[+16pt]
    &
    B
    \ar[
      from=ul,
      shorten=-4pt
    ]
    \ar[orangeii,
      from=s,
      to=t,
      Rightarrow,
      bend left=50pt,
      shift left=7pt,
      crossing over,
      shorten=-1pt,
      "\ "{name=s3, swap}
    ]
    \ar[orangeii,
      from=s,
      to=t,
      Rightarrow,
      bend right=50pt,
      shift right=7pt,
      crossing over,
      shorten=-1pt,
      "\ "{name=t3}
    ]
    \ar[
      from=t3,
      to=s3,
      Rightarrow,
      shorten=-3pt,
      crossing over
    ]
    \ar[
      from=t3,
      to=s3,
      Rightarrow,
      shorten=-3pt,
      crossing over
    ]
    \ar[
      from=t3,
      to=s3,
      -,
      shorten=-3pt
    ]
 \ar[orangeii,
      from=s2,
      to=t2,
      Rightarrow,
      bend left=20pt,
      crossing over,
      shorten=-3pt
    ]
  \end{tikzcd}
  \right\}
  \;\;\;\;
  \in
  \;
  \InfinityGroupoids
  \,,
\end{equation}
is itself an $\infty$-topos.
When $B \,=\, \mathbf{B}\mathcal{G}$,  the cohomology in the
slice $\SliceTopos{\mathbf{B}\mathcal{G}}$ is
{\it Borel-$\mathcal{G}$-equivariant cohomology} of $\Topos$
(see Def. \ref{BorelEquivariantAndProperEquivariantCohomologyInCohesiveInfinityTopos} below):
\vspace{-1mm}
$$
  \underset{
    \mathclap{
    \raisebox{-3pt}{
      \tiny
      \color{darkblue}
      \bf
      \begin{tabular}{c}
        Borel-equivariant
        \\
        $A$-cohomology of $X$
      \end{tabular}
    }
    }
  }{
    H^\bullet_{\mathcal{G}}(X;\, A)
  }
  \;\;
  \coloneqq
  \;\;
  \Truncation{0}
  \,
  \SlicePointsMaps{\big}{\mathbf{B}\mathcal{G}}
    { \HomotopyQuotient{X}{\mathcal{G}} }
    {
      \Omega_{\mathbf{B}\mathcal{G}}^{-\bullet}
      (
        \HomotopyQuotient{A}{\mathcal{G}}
      )
    }
  \,.
$$
This is a consequence of the following fact:

\medskip
\noindent
{\bf Transformation groups in an $\infty$-topos.}
The key fact which propels our general theory in
\cref{EquivariantInfinityBundles}
is
(following \cite{NSS12a} and \cite[\S 2.2]{SS20OrbifoldCohomology}):
that the notions of {\it groups},
of {\it group actions},
of {\it principal bundles}
and their {\it moduli stacks}
and associated {\it fiber bundles}
are all {\it native to $\infty$-topos theory}, in that these concepts and
their pertinent properties are available internal to any $\infty$-topos
without needing further axiomatization:

\begin{proposition}[Groups, actions and principal bundles in any $\infty$-topos
({Props. \ref{LoopingAndDeloopingEquivalence},
  \ref{HomotopyQuotientsAndPrincipaInfinityBundles},
Thm. \ref{DeloopingGroupoidsAreModuliInfinityStacksForPrincipalInfinityBundles}
})]
\label{GroupsActionsAndFiberBundles}
$\,$

\noindent
Let $\Topos$ be an $\infty$-topos.

\noindent
{\bf (i)} The operation
of forming loop space objects \eqref{LoopingInIntroduction}
constitutes an equivalence\footnote{This is the \emph{May recognition theorem} \cite{May72}
generalized from $\InfinityGroupoids$ to $\infty$-toposes by
\cite[7.2.2.11]{Lurie09HTT}\cite[6.2.6.15]{Lurie17}.}
of
group objects with pointed connected objects in $\Topos$ (Ntn. \ref{ConnectedObject}):

\vspace{-6mm}
\begin{equation}
  \label{LoopingAndDelooping}
  \begin{tikzcd}[row sep=0pt]
    \Groups(\Topos)
    \ar[
      rr,
      shift right=5pt,
      "\mathbf{B}"{below}
    ]
    \ar[
      rr,
      phantom,
      "{\scalebox{.8}{$\simeq$}}"
    ]
    &&
    \Topos^{\ast/}_{\geq 1}
    \ar[
      ll,
      shift right=5pt,
      "\Omega"{above}
    ]
    \\
  \scalebox{0.7}{$  \mathcal{G} $}
    \ar[
      rr,
      |->
    ]
    &&
  \scalebox{0.7}{$  \ast \!\sslash\! \mathcal{G} $}
  \end{tikzcd}
  {\phantom{AAAA}}
  \mbox{\rm i.e.}
  {\phantom{AAAAA}}
  \begin{tikzcd}[column sep=huge]
    \overset{
    \mathclap{
      \raisebox{3pt}{
        \tiny
        \color{darkblue}
        \bf
        group stack
      }
      }
    }{
      \mathcal{G}
    }
    \ar[r]
    \ar[d]
    \ar[
      dr,
      phantom,
      "\mbox{\tiny\rm(pb)}"
    ]
    &
    \overset{
      \mathclap{
      \raisebox{5pt}{
        \tiny
        \color{darkblue}
        \bf
        \begin{tabular}{c}
          base
          \\
          point
        \end{tabular}
      }
      }
    }{
      \ast
    }
    \ar[
      d,
      "\mathrm{pt}_{\mathbf{B}G}"
    ]
    \\
    \ast
    \ar[
      r,
      "\mathrm{pt}_{\mathbf{B}G}"
    ]
    &
    \mathbf{B}\mathcal{G}
  \end{tikzcd}
 \end{equation}
hence
\begin{equation}
  \label{PointedConnectedObjectEquivalentToDeloopingOfItsLoopSpaceObject}
  \mbox{
    $
    X
    \,\in\,
    \Topos^{\ast/}
    $
    is connected
  }
  {\phantom{AAAA}}
  \Leftrightarrow
  {\phantom{AAAA}}
  X
  \;\simeq\;
  \mathbf{B} \Omega X
  \,.
\end{equation}

\vspace{-1mm}
\noindent
{\bf (ii)}
For $\mathcal{G} \,\in\, \Groups(\Topos)$,
the $\mathcal{G}$-actions (Def. \ref{ActionObjectsInAnInfinityTopos})
and $\mathcal{G}$-principal bundles (Def. \ref{PrincipalInfinityBundles})
are both identified with the slice objects \eqref{SliceHomInIntroduction}
over the delooping $\mathbf{B} \mathcal{G}$ \eqref{GroupsActionsAndFiberBundles},
as follows:
\vspace{-2mm}
\begin{equation}
  \label{EquivalenceBetweenActionsAndPrincipalBundlesAndSlices}
  \hspace{0cm}
  \begin{tikzcd}[row sep=-1pt, column sep=8pt]
    \Actions{\mathcal{G}}(\Topos)
    \ar[
      rr,
      "\sim"{above, yshift=1pt},
      "(-) \!\sslash\! G"{below}
    ]
    &&
    \PrincipalBundles{\mathcal{G}}(\Topos)
    &&
    \Topos_{/\mathbf{B}\mathcal{G}}
    \ar[
      ll,
      "\sim"{above, yshift=1pt},
      "\mathrm{fib}"{below}
    ]
    \\
\scalebox{0.7}{$    G \acts \, P $}
    &\longmapsto&
  \scalebox{.7}{$  \left(\!\!\!\!
    \def\arraystretch{.9}
    \begin{array}{c}
      P
      \\
      \downarrow
      \\
      P \!\sslash\! \mathcal{G}
    \end{array}
    \!\!\!\! \right)
    $}
    &\longmapsfrom&
    \scalebox{.7}{$ \left(\!\!\!\!
    \def\arraystretch{.9}
    \begin{array}{c}
      P \!\sslash\! \mathcal{G}
      \\
      \downarrow
      \\
      \mathbf{B}\mathcal{G}
    \end{array}
    \!\!\!\! \right)
    $}
  \end{tikzcd}
  \end{equation}

\vspace{-2mm}
\noindent  i.e.,
\vspace{-2mm}
\begin{equation}
  \begin{tikzcd}[column sep=80pt]
    \mathllap{
      \mbox{
        \tiny
        \color{greenii}
        \bf
        action
      }
      \;\;\;\;\;\;
    }
    \overset{
      \mathclap{
      \raisebox{6pt}{
        \tiny
        \color{darkblue}
        \bf
        \begin{tabular}{c}
          $\mathcal{G}$-principal
          \\
          bundle
          \\
          \phantom{.}
        \end{tabular}
      }
      }
    }{
      P
    }
      \ar[out=180-60+90, in=60+90, looseness=3.8, "\scalebox{.77}{$\mathclap{
        \mathcal{G}
      }$}"{pos=.41, description},shift right=1]
    \ar[r]
    \ar[
      d
    ]
    \ar[
      dr,
      phantom,
      "\mbox{\tiny\rm(pb)}"
    ]
    &
    \overset{
      \mathclap{
      \raisebox{6pt}{
        \tiny
        \color{darkblue}
        \bf
        \begin{tabular}{c}
          universal
          \\
          $\mathcal{G}$-principal bundle
        \end{tabular}
      }
      }
    }{
      \ast
      \mathrlap{
        \;
        \simeq \mathcal{G} \!\sslash\! \mathcal{G}
      }
    }
    \ar[
      d,
      "\mathrm{pt}_{\mathbf{B}\mathcal{G}}"{right}
    ]
    \\
    \mathllap{
      \mathllap{
      \mbox{
        \tiny
        \color{darkblue}
        \bf
        \begin{tabular}{c}
          action quotient-
          \\
          / base-stack
        \end{tabular}
      }
      \!\!\!\!\!
      }
    }
    P \!\sslash\! \mathcal{G}
    \ar[
      r,
      "\vdash P"{above},
      "\mbox{\tiny \color{greenii}\bf cocycle}"{below}
    ]
    &
    \mathbf{B}\mathcal{G}
    \mathrlap{
      \!\!\!\!\!\!
      \mbox{
        \tiny
        \color{darkblue}
        \bf
        \begin{tabular}{c}
          universal
          \\
          moduli stack
        \end{tabular}
      }
    }
  \end{tikzcd}
\end{equation}
\end{proposition}

\medskip

\noindent
{\bf Equivariant principal $\infty$-bundles.} With the conceptualization
of Prop. \ref{GroupsActionsAndFiberBundles} in hand, there is an
evident general-abstract definition of $G$-equivariant principal $\infty$-bundles
(Def. \ref{GEquivariantGammaPrincipalBundles}):
These must simply be the principal $\infty$-bundles internal to
a slice $\infty$-topos $\SliceTopos{\mathbf{B}G}$ over the delooping $\mathbf{B}G$:

\vspace{-.2cm}
$$
  \overset{
    \mathclap{
    \raisebox{7pt}{
      \tiny
      \color{darkblue}
      \bf
      \begin{tabular}{c}
        $G$-equivariant $\Gamma$-principal
        $\infty$-bundles in $\Topos$
      \end{tabular}
    }
    }
  }{
    \EquivariantPrincipalBundles{G}{\Gamma}(\Topos)_X
  }
  \quad
  \coloneqq
  \;\;
  \overset{
    \mathclap{
    \raisebox{3pt}{
      \tiny
      \color{darkblue}
      \bf
      \begin{tabular}{c}
        $\HomotopyQuotient{\Gamma}{G}$-principal
        $\infty$-bundles in $\SliceTopos{\mathbf{B}G}$
        \\
        \phantom{.}
      \end{tabular}
    }
    }
  }{
    \PrincipalBundles{ (\HomotopyQuotient{\Gamma}{G}) }
    (
      \SliceTopos{\mathbf{B}G}
    )_{ \HomotopyQuotient{X}{G} }
  }
  \,.
$$

We prove
(in \eqref{EquivalenceOfGroupoidOfEquivariantCechCocyclesIntoGroupoidOfEquivariantPrincipalBundles} of Thm. \ref{BorelClassificationOfEquivariantBundlesForResolvableSingularitiesAndEquivariantStructure})
that, for the case $\Topos \,\coloneqq\,\SmoothInfinityGroupoids$
and restricted to topological $G$-spaces $\TopologicalSpace$
and topological structure groups $\Gamma$,
this canonical $\infty$-topos-theoretic definition
recovers the traditional definition
of topological equivariant principal bundles,
including their equivariant local triviality property
(which are all reviewed and developed in \Cref{EquivariantPrincipalTopologicalBundles}).

\medskip
By extension, this means that for more general $\Topos$ and/or more general
$X, \Gamma \,\in\, \SliceTopos{\mathbf{B}G}$, we obtain sensible generalizations
of these classical definitions.
For example, by taking $\Topos$ to be the cohesive $\infty$-topos of super-geometric $\infty$-groupoids
(\cite[\S 3.1.3]{SS20OrbifoldCohomology}) the general theory
developed here immediately produces a good theory of equivariant
higher super-gerbes (as needed, e.g., in super-string theory
on super-orbifolds \cite{FSS13}\cite{HSS18}).

\medskip
Finally, this abstract definition in combination with the orbi-smooth Oka principle
\eqref{OrbiSmoothOkaPrincipleInIntroduction} implies the classification
of stable $G$-equivariant $\Gamma$-principal bundles,
at least in the case that $G$-singularities are cover-resolvable
(Ntn. \ref{ResolvableOrbiSingularities})
and $\Gamma$ is a truncated Hausdorff group.
This is the content of our main
Theorems \ref{ClassificationOfPrincipalBundlesAmongPrincipalInfinityBundles} ands \ref{ProperClassificationOfEquivariantBundlesForResolvableSingularitiesAndEquivariantStructure} below.

\medskip

Even though this classification result concerns only equivariant ``1-bundles''
instead of more general equivariant $\infty$-bundles,
the cohesive $\infty$-topos theory drives the proof:
For example, our generalization of the
existing classification results to non-trivial discrete $G$-action on the structure group
$\Gamma$
is a direct consequence
(see the end of the proof of Thm. \ref{ShapeOfMappingStackOutOfOrbiSingularityIsMappingStackIntoShape})
of the general fact (Prop. \ref{ShapeFunctorPreservesHomotopyFibersOverDiscreteObjects})
that the shape operation in
cohesive $\infty$-toposes preserves homotopy fiber products over geometrically
discrete groupoids (such as $\mathbf{B}G$ for discrete $G$).

\newpage

\section{Notation and Terminology}

\phantom{AA} {\bf Categories and functors.}

\vspace{1mm}
\def\arraystretch{1.3}
\begin{tabular}{lll}
\hline
  {\bf Category} & {\bf of}
  \\
  \hline
  \hline
  \rowcolor{lightgray}
  $\kTopologicalSpaces$ & cg topological spaces
  &
  Ntn. \ref{CompactlyGeneratedTopologicalSpaces}
  \\
  $\kHausdorffSpaces$ & cg Hausdorff spaces & Ntn. \ref{CompactlyGeneratedTopologicalSpaces}
  \\
  \rowcolor{lightgray}
  $\Actions{G}(\kTopologicalSpaces)$ & topological $G$-spaces
  &
  Ntn. \ref{GActionOnTopologicalSpaces}
  \\
  $\Groups(\Sets)$ & discrete groups
  \\
  \rowcolor{lightgray}
  $\Groups(\kHausdorffSpaces)$ & Hausdorff groups &
  \\
  $\FormallyPrincipalBundles{\Gamma}(\mathcal{C})$
  &
  formally principal internal bundles
  & Ntn. \ref{InternalizationOfPrincipalBundleTheory}
  \\
  \rowcolor{lightgray}
  $\EquivariantPrincipalBundles{G}{\Gamma}$ & equivariant principal bundles
  &
  Def. \ref{EquivariantPrincipalBundle}
  \\
  $\EquivariantPrincipalFiberBundles{G}{\Gamma}$ & ... equivariantly locally trivial
  &
  Def. \ref{TerminologyForPrincipalBundles}
  \\
  \rowcolor{lightgray}
  \hline
  $\CartesianSpaces$ & Cartesian spaces &
  \\
  $\DiffeologicalSpaces$ & diffeological spaces
  & Ntn. \ref{CartesianSpacesAndDiffeologicalSpaces}
  \\
  \rowcolor{lightgray}
  $\DTopologicalSpaces$ & D-topological spaces
  &
  Ntn. \ref{DeltaGeneratedTopologicalSpaces}
  \\
  \hline
  \hline
  $\SimplicialSets$ & simplicial sets
  &
  Ntn. \ref{SimplicialSets}
  \\
  \rowcolor{lightgray}
    $\SimplicialCategories$
  & simplicial categories
  & Ntn. \ref{SimplicialCategories}
  \\
    $\SimplicialPresheaves$
  & simplicial presheaves
  & Ntn. \ref{ModelCategoriesOfSimplicialPresheaves}
  \\
  \hline
  \hline
  {\bf 2-category} & {\bf of}
  \\
  \rowcolor{lightgray}
  \hline
  $\Groupoids$ & groupoids &
  \\
    $\TopologicalGroupoids$
  & topological groupoids
  & Ntn. \ref{TopologicalGroupoids}
  \\
  \rowcolor{lightgray}
  $\HomotopyTwoCategory(\PresentableInfinityCategories)$
  &
  presentable $\infty$-categories
  &
  Prop. \ref{HomotopyCategoryOfPresInfinityCategoriesIsThatOfCombinatorialModelCategories}
  \\
  \hline
  \hline
  {\bf $\infty$-category} & {\bf of}
  \\
  \hline
  \hline
    $\InfinityGroupoids$
  & $\infty$-groupoids
  &
  Ntn. \ref{SimplicialSetsAndInfinityGroupoids}
  \\
  \rowcolor{lightgray}
    $\SmoothInfinityGroupoids$
  & smooth $\infty$-groupoids
  & Ntn. \ref{SmoothInfinityGroupoids}
  \\
    $\SingularSmoothInfinityGroupoids$
  & singular smooth $\infty$-groupoids
  & Ntn. \ref{SingularSmoothInfinityGroupoids}
  \\
  \hline
  \hline
  {\bf Functor} & {\bf producing}
  \\
  \rowcolor{lightgray}
  \hline
  $N$ & simplicial nerve &
  Ntn. \ref{NerveOfTopologicalGroupoids}
  \\
  $\TopologicalRealization{}{-}$
  &
  topological realization
  &
  Ntn. \ref{TopologicalRealizationFunctors}
  \\
  \rowcolor{lightgray}
    $\ContinuousDiffeology$
  & continuous diffeology
  & Ex. \ref{ContinuousDiffeologyAndDTopology}
  \\
    $\DTopology$
  & D-topology
  & Ex. \ref{ContinuousDiffeologyAndDTopology}
  \\
    \hline
\end{tabular}

\vspace{5mm}

\phantom{.} {\bf Types of groups.}

\vspace{1mm}
\def\arraystretch{1.4}
\begin{tabular}{llll}
  \hline
  \hline
  \rowcolor{lightgray}
  $G$
  & $\in \Groups(\Sets) \xhookrightarrow{\;} \Groups(\Topos)$
  & equivariance group
  & Ntn. \ref{GActionOnTopologicalSpaces}
  \\
  \hline
  $G \acts \, \Gamma$
  &
  $\in \Groups(\kTopologicalSpaces)$
  &
  equivariant structure group
  &
  Def. \ref{EquivariantTopologicalGroup}, Lem. \ref{EquivariantTopologicalGroupsAreSemidirectProductsWithG}
  \\
  \hline
  \rowcolor{lightgray}
    $\mathcal{G}$
  & $\in \Groups(\Topos)$
  & structure $\infty$-group
  & Prop. \ref{GroupsActionsAndFiberBundles}
  \\
  \hline
  $G \acts \, \Gamma$
  & $\in \Groups \left( \Actions{G}(\Topos)\right)$
  &
  \hspace{-2pt}
  \multirow{2}{*}{
    equivariant structure $\infty$-group
  }
    &
  \multirow{2}{*}{
    Def. \ref{GGroups}
  }
  \\
  \cline{1-2}
  $\Gamma \!\sslash\! G$
  & $\in \Groups(\Topos_{/\mathbf{B}G})$
  &
  &
  \\
  \hline
\end{tabular}

\newpage

\noindent
{\bf Ambient $\infty$-toposes.}

\vspace{.2cm}

$$
  \begin{tikzcd}[row sep=50pt, column sep=large]
    &
    &[+20pt]
    \mathllap{
      \raisebox{4pt}{
        \tiny
        \color{darkblue}
        \bf
        \def\arraystretch{.9}
        \begin{tabular}{c}
          slice over
          \\
          $G$-orbi-singularity
        \end{tabular}
      }
      \hspace{+1pt}
    }
    \categorybox{
      \SliceTopos{\orbisingularG}
    }
    \ar[
      rr,
      <->
    ]
    \ar[
      dl,
      <->,
      end anchor={[xshift=+2pt, yshift=+2pt]}
    ]
    \ar[
      dd,
      <->,
      crossing over,
      "{
        \def\arraystretch{.6}
        \begin{array}{c}
          G\Conical
          \\
          \bot
          \\
          G\Space
          \\
          \bot
          \\
          G\Smooth
          \\
          \bot
          \\
          G\Orbisingular
        \end{array}
      }"{swap, pos=.6, xshift=4pt}
    ]
    &&
    \categorybox{
      \GloballyEquivariant\InfinityGroupoids{}_{/\scalebox{.7}{$\orbisingularG$}}
    }
    \ar[
      dd,
      <->
    ]
    \ar[
      dl,
      <->
    ]
    \\[-45pt]
    &
    \mathllap{
      \mbox{
        \tiny
        \cref{GeneralSingularCohesion}
        \color{darkblue}
        \bf
        \begin{tabular}{c}
          singular cohesive
          \\
          $\infty$-topos
        \end{tabular}
      }
      \hspace{+0pt}
    }
    \categorybox{
      \Topos
    }
    \ar[
      rr,
      <->,
      crossing over
    ]
    \ar[
      dd,
      <->,
      "{
        \def\arraystretch{.6}
        \begin{array}{c}
          \Conical
          \\
          \bot
          \\
          \Space
          \\
          \bot
          \\
          \Smooth
          \\
          \bot
          \\
          \Orbisingular
        \end{array}
      }"{swap, xshift=5pt}
    ]
    &&
    \categorybox{
      \GloballyEquivariant\InfinityGroupoids
    }
    \mathrlap{
      \hspace{-3pt}
      \mbox{
        \tiny
        \color{darkblue}
        \bf
        \def\arraystretch{.9}
        \begin{tabular}{c}
          global
          \\
          homotopy theory
        \end{tabular}
      }
    }
    \\[+20pt]
    &
    &
    \categorybox{
      \GEquivariant\ModalTopos{\smooth}
    }
    \ar[
      dl,
      <->,
      end anchor={[xshift=+2pt, yshift=+2pt]}
    ]
    \ar[
      rr,
      <->,
      crossing over
    ]
    &&
    \categorybox{
      \GEquivariant\InfinityGroupoids
    }
    \mathrlap{
      \hspace{-3pt}
      \mbox{
        \tiny
        \color{darkblue}
        \bf
        \def\arraystretch{.9}
        \begin{tabular}{c}
          $G$-equivariant
          \\
          homotopy theory
        \end{tabular}
      }
    }
    \ar[
      dl,
      <->
    ]
    \\[-45pt]
    &
    \mathllap{
      \mbox{
        \tiny
        \cref{GeneralCohesion}
        \color{darkblue}
        \bf
        \def\arraystretch{.9}
        \begin{tabular}{c}
          smooth cohesive
          \\
          $\infty$-topos
        \end{tabular}
      }
      \hspace{+1pt}
    }
    \categorybox{
      \ModalTopos{\smooth}
    }
    \ar[
      rr,
      <->,
      "{
        \Shape
        \;\dashv\;
        \Discrete
        \;\dashv\;
        \Points
        \;\dashv\;
        \Chaotic
      }"{swap}
    ]
    &&
    \categorybox{
      \InfinityGroupoids
    }
    \mathrlap{
      \hspace{-3pt}
      \raisebox{-3pt}{
        \tiny
        \color{darkblue}
        \bf
        \def\arraystretch{.9}
        \begin{tabular}{c}
          base
          \\
          $\infty$-topos
        \end{tabular}
      }
    }
    \ar[
      from=uu,
      <->,
      crossing over
    ]
  \end{tikzcd}
$$

\vspace{2mm}
Here:
\vspace{-2mm}
$$
  \def\arraystretch{1.5}
  \begin{array}{lcrcll}
  \ModalTopos{\smooth}
  &:=&
  \SmoothInfinityGroupoids
  &:=&
  \InfinitySheaves(\CartesianSpaces)
  &
  \proofstep{(Ntn. \ref{SmoothInfinityGroupoids},
    {\cite[Ex. 3.18]{SS20OrbifoldCohomology}}
 \!)},
  \\
  \Topos
  &:=&
  \SingularSmoothInfinityGroupoids
  &:=&
  \InfinitySheaves(\CartesianSpaces \times \Singularities)
  &
  \proofstep{(Ntn. \ref{SingularSmoothInfinityGroupoids},
    {\cite[Ex. 3.56]{SS20OrbifoldCohomology}}
  \!)}.
  \end{array}
$$

\vspace{5mm}

\noindent {\bf Modalities.}

\vspace{2mm}

\begin{center}
\def\arraystretch{3}
\begin{tabular}{|c||l|l|l|}
  \hline
  \def\arraystretch{1}
  \begin{tabular}{c}
    {\bf Cohesion}
    \\
    Def. \ref{CohesiveInfinityTopos}
  \end{tabular}
  &
  $
  \overset{
    \mathrlap{
      \raisebox{4pt}{
        \tiny
        \color{darkblue}
        \bf
        shape
      }
    }
  }{
  \shape
    \,\coloneqq\,
  }
  \Discrete \circ \Shape
  $
  &
  $
  \overset{
    \mathrlap{
      \raisebox{4pt}{
        \tiny
        \color{darkblue}
        \bf
        discrete
      }
    }
  }{
  \flat
    \,\coloneqq\,
  }
  \Discrete \circ \Points
  $
  &
  $
  \overset{
    \mathrlap{
      \raisebox{4pt}{
        \tiny
        \color{darkblue}
        \bf
        sharp
      }
    }
  }{
    {\rm chaotic}
      \,\coloneqq\,
  }
  \Chaotic \circ \Points
  $
  \\
  \hline
  \def\arraystretch{1}
  \begin{tabular}{c}
    {\bf Singularities}
    \\
    Def. \ref{GEquivariantAndGloballyEquivariantHomotopyTheories}
  \end{tabular}
  &
  $
  \overset{
    \mathrlap{
      \raisebox{4pt}{
        \tiny
        \color{darkblue}
        \bf
        conical
      }
    }
  }{
  \conical
    \,\coloneqq\,
  }
  \Space \circ \Conical
  $
  &
  $
  \overset{
    \mathrlap{
      \raisebox{4pt}{
        \tiny
        \color{darkblue}
        \bf
        smooth
      }
    }
  }{
  \smooth
    \,\coloneqq\,
  }
  \Space \circ \Smooth
  $
  &
  $
  \overset{
    \mathrlap{
      \raisebox{4pt}{
        \tiny
        \color{darkblue}
        \bf
        orbisingular
      }
    }
  }{
    \orbisingular
      \,\coloneqq\,
  }
  \Singularity \circ \Smooth
  $
  \\
  \hline
  \def\arraystretch{1}
  \begin{tabular}{c}
 {\bf    $G$-singularities }
    \\
    Def. \ref{ModalitieswithrespecttoGOrbiSingularities}
  \end{tabular}
  &
  $
  \overset{
    \mathrlap{
      \raisebox{4pt}{
        \tiny
        \color{darkblue}
        \bf
        $G$-conical
      }
    }
  }{
  \conicalrelativeG
    \,\coloneqq\,
  }
  G\Space \circ G\Conical
  $
  &
  $
  \overset{
    \mathrlap{
      \raisebox{4pt}{
        \tiny
        \color{darkblue}
        \bf
        $G$-smooth
      }
    }
  }{
  \smoothrelativeG
    \,\coloneqq\,
  }
  G\Space \circ G\Smooth
  $
  &
  $
  \overset{
    \mathrlap{
      \raisebox{4pt}{
        \tiny
        \color{darkblue}
        \bf
        $G$-orbisingular
      }
    }
  }{
    \orbisingularrelativeG
      \,\coloneqq\,
  }
  G\Singularity \circ G\Smooth
  $
  \\
  \hline
\end{tabular}
\end{center}

\def\arraystretch{1}

\vspace{5mm}

\noindent {\bf Notions of space.}

\vspace{-.1cm}

$$
  \begin{tikzcd}[row sep=17pt, column sep=23pt]
    \overset{
      \mathclap{
      \raisebox{4pt}{
        \tiny
        Ntn. \ref{CompactlyGeneratedTopologicalSpaces}
      }
      }
    }{
      \categorybox{\kTopologicalSpaces}
    }
    \ar[
      rr,
      hook
    ]
    \ar[
      ddr,
      "\SingularSimplicialComplex"{swap, pos=.65},
      "{
        \mbox{
          \tiny
          \begin{tabular}{c}
            Ntn. \ref{DiffeologicalSingularComplex}
            \\
            \color{greenii}
            \bf
            sing. simp. complex
          \end{tabular}
        }
      }"{sloped, pos=.54    }
    ]
    &&
    \overset{
      \raisebox{4pt}{
        \tiny
        Ntn. \ref{GActionOnTopologicalSpaces}
      }
    }{
      \categorybox{
        \Actions{G}(\kTopologicalSpaces)
      }
    }
    \ar[
      ddr,
      "{
        \shape \orbisingular (\HomotopyQuotient{-}{G})
      }",
      "\mbox{
        \tiny
        \begin{tabular}{c}
          \color{greenii}
          \bf
          equivariant shape
          \\
          \eqref{EquivariantShape}
        \end{tabular}
      }"{sloped, swap, pos=.4}
    ]
    \ar[
      rr,
      "{
        \orbisingular (\HomotopyQuotient{-}{G})
      }",
      "{
        \mbox{
          \tiny
          \color{greenii}
          \bf
          orbi-singularized homotopy quotient
        }
      }"{swap}
    ]
    &&
    \overset{
      \raisebox{4pt}{
        \tiny
        Ntn. \ref{SingularSmoothInfinityGroupoids}
      }
    }{
      \categorybox{
        (
          \SingularSmoothInfinityGroupoids
        )_{/\scalebox{.7}{$\orbisingularG$}}
      }
    }
  \ar[
    ddl,
    "{
      \mbox{
        \tiny
        \begin{tabular}{c}
          \color{greenii}
          \bf
          shape
        \end{tabular}
      }
    }"{swap, sloped},
    "{
      \shape
    }"{swap}
  ]
  \\[-22pt]
 \hspace{-5mm}
 \underset{
    \mathclap{
    \raisebox{-6pt}{
      \tiny
      \color{darkblue}
      \bf
      \begin{tabular}{c}
        topological
        \\
        space
      \end{tabular}
      }
    }
  }{
    \scalebox{.93}{$\TopologicalSpace$}
  }
  &&
  \hspace{-1cm}
  \underset{
    \mathclap{
    \raisebox{-6pt}{
      \tiny
      \color{darkblue}
      \bf
      \begin{tabular}{c}
        topological
        \\
        $G$-space
      \end{tabular}
      }
    }
  }{
    \scalebox{.93}{$G \acts \, \TopologicalSpace$}
  }
  &&
  \underset{
    \mathclap{
    \raisebox{-6pt}{
      \tiny
      \color{darkblue}
      \bf
      \begin{tabular}{c}
        orbi-stack {\color{black}/}
        \\
        cohesive orbi-space
      \end{tabular}
      }
    }
  }{
    \mathcal{X}
  }
    \\
    &
    \mathllap{
    \mbox{ \tiny
      Ntn. \ref{SimplicialSetsAndInfinityGroupoids}
      }
    \;\;\;}
    \categorybox{\InfinityGroupoids}
    \ar[
      rr,
      hook
    ]
    &&
    \categorybox{
      (
        \SingularInfinityGroupoids
      )_{/\scalebox{.7}{$\orbisingularG$}}
    }
    \\[-20pt]
  &
  \underset{
    \mathclap{
    \raisebox{-6pt}{
      \tiny
      \color{darkblue}
      \bf
      \begin{tabular}{c}
        (shape of)
        \\
        space
      \end{tabular}
      }
    }
  }{
    X
  }
  &&
  \underset{
    \mathclap{
    \raisebox{-6pt}{
      \tiny
      \color{darkblue}
      \bf
      \begin{tabular}{c}
        orbi-
        \\
        space
      \end{tabular}
      }
    }
  }{
    \scalebox{1.06}{$\mathscr{X}$}
  }
\end{tikzcd}
$$

\newpage

\part{In topological spaces}
 \label{InTopologicalSpaces}

\chapter{Equivariant topology}
\label{EquivariantTopology}

We recall and develop basics of equivariant algebraic topology
(see \cite{Illman72}\cite{Bredon72}\cite{tomDieck87}\cite{May96}\cite{Blu17})
that we invoke below in \cref{EquivariantPrincipalTopologicalBundles}.
The expert reader may want to skip this chapter and refer back to it just as need be.

\begin{itemize}

\vspace{-.2cm}
\item[{\cref{TopologicalGActions}}] recalls basics of equivariant point-set topology,
highlighting how the change-of-group adjoint triple governs the theory,
and establishing lemmas needed in \cref{EquivariantPrincipalTopologicalBundles}.

\vspace{-.2cm}
\item[{\cref{GActionsOnTopologicalGroupoids}}] recalls basics of topological groupoids,
generalizing to equivariant topological groupoids,
and establishing lemmas needed in \cref{ConstructionOfUniversalEquivariantPrincipalBundles}.

\vspace{-.2cm}
\item[{\cref{GEquivariantHomotopyTypes}}] recalls basics of proper equivariant homotopy theory,
recording lemmas needed in \cref{ConstructionOfUniversalEquivariantPrincipalBundles}
and \cref{EquivariantLocalTrivializationIsImplies}.

\end{itemize}

Throughout, we make extensive use of
a hierarchy of {\it internalizations}
of mathematical structures (Ntn. \ref{Internalization})
into {\it categories with pullbacks} (Ntn. \ref{CartesianSquares}),
starting in a {\it convenient category of topological spaces}
(Ntn. \ref{CompactlyGeneratedTopologicalSpaces}).
Here is the basic terminology and notation that we are using:

\medskip
\noindent
{\bf Categories.}
In \cref{TopologicalGActions} and \cref{EquivariantPrincipalTopologicalBundles}
we need just basic notions of (co)limits and adjoint functors
in plain category theory (e.g., \cite{AHS90}\cite{Borceux94I}\cite{Borceux94II}),
while in \cref{GActionsOnTopologicalGroupoids}
and \cref{EquivariantLocalTrivializationIsImplies}
we need these notions in their $\Groupoids$-enriched enhancement
(e.g., \cite[\S 6]{Borceux94II}\cite[\S 3]{Riehl14}\cite[\S 1.3]{JohnsonYau21})
--
but we only need the most basic concepts:
enriched adjunctions and conical (i.e., non-weighted) enriched (co-)limits,
and here mostly just finite ones.

\begin{notation}[Basic categories]
  \label{CategoryOfGroupoids}
  We write

  \noindent
  {\bf (i)} $\Sets$ for the category of sets with functions between them;

  \noindent
  {\bf (ii)}
  $\Groupoids \;\coloneqq\; \Groupoids(\Sets)$
  for the category of small groupoids with functors between them.
\end{notation}

\begin{notation}[Morphisms]
  \label{BasicNotationForCategories}
  Let $\mathcal{C}$ be a category and
  $X, Y \,\in\, \mathcal{C}$
  a pair of its objects. We write

  \noindent
  {\bf (i)}
  $\mathcal{C}(X,Y) \,\in\, \Sets$
  for the set of morphisms $X \to Y$ in $\mathcal{C}$ (the {\it hom-set});

  \noindent
  {\bf (ii)}
  $X \xrightarrow{\;\; \sim \;\;} Y$
  to indicate {\it isomorphisms}.

\end{notation}

\begin{notation}[$\Groupoids$-enriched categories]
  \label{Strict2Categories}
With $\Groupoids$ (Ntn. \ref{CategoryOfGroupoids}) regarded as a
cartesian monoidal category,
a
{\it strict (2,1)-category},
namely a
{\it $\Groupoids$-enriched category}
\cite{FanthamMoore83},
has for each pair of objects $X , Y \,\in\, \mathcal{C}$
a {\it hom-groupoid}
$
  \mathcal{C}(X,Y)
  \;\in\;
  \Groupoids
$.
Equivalently, this is a
{\it strict 2-category} or
{\it $\mathrm{Cat}$-enriched category}
(e.g. \cite[\S 2.3]{JohnsonYau21}\cite[\S 9.5]{Richter20})
whose hom-categories happen to be groupoids.
\end{notation}

\begin{notation}[Adjoint functors]
  \label{AdjointFunctors}
  We denote pairs of adjoint functors as shown on the left here:
  \vspace{-2mm}
  \begin{equation}
    \label{FormingAdjuncts}
    \begin{tikzcd}
      \mathcal{D}
      \ar[
        rr,
        "R"{below},
        shift right=5.5pt
      ]
      \ar[
        rr,
        phantom,
        "\scalebox{.7}{$\bot$}"{description}
      ]
      &&
      \mathcal{C}
      \ar[
        ll,
        "L"{above},
        shift right=5pt
      ]
    \end{tikzcd}
    {\phantom{AAAA}}
    \Leftrightarrow
    {\phantom{AAAA}}
    \begin{tikzcd}[row sep=-4pt, column sep=1pt]
    \mathcal{C}
    \left(
      c, R(d)
    \right)
    &
    \simeq
    &
    \mathcal{D}
    \left(
      L(c),\, d
    \right)
    \\
    c \xrightarrow{f} R(d)
    &\leftrightarrow&
    L(c) \xrightarrow{\tilde f} d
    \end{tikzcd}
  \end{equation}

  \vspace{-3mm}
  \noindent
  meaning that for all objects $c \in \mathcal{C}$ and $d \in \mathcal{C}$ there is
  a natural isomorphism (``forming adjuncts'')
  between the hom-objects
  (Ntn. \ref{BasicNotationForCategories}, \ref{Strict2Categories})
  out of the image of left adjoint functor $L$
  and into the image of the right adjoint functor $E$,
  as shown on the right.
\end{notation}

\begin{notation}[Cartesian/pullback squares]
  \label{CartesianSquares}
  We indicate that a commuting square of morpisms
  in some category $\mathcal{C}$ (Ntn. \ref{BasicNotationForCategories})
  -- here typically in the category of
  $\GActionsOnTopologicalSpaces$ \eqref{GActionsOnTopologicalSpaces} --
  is a {\it pullback square} (also: {\it Cartesian square} or {\it fiber product})
  by putting the symbols ``{\color{orangeii} (pb)}'' at its center.
  This means that each pair of morphisms forming another commuting square with its
  right and bottom morphism (a ``cone'') factors uniquely through the its
  top left object
  such that the resulting triangles commute:
  
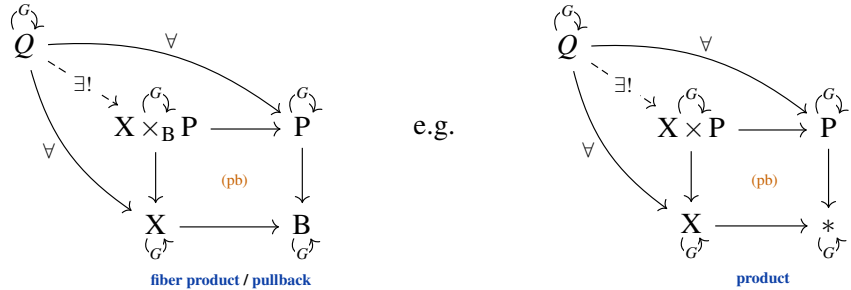
\vspace{-3mm}
  $$
    \begin{tikzcd}
      Q
      \ar[out=180-66, in=66, looseness=3.5, "\scalebox{.77}{$\mathclap{
        G
      }$}"{description},shift right=1]
      \ar[
        drr,
        bend left=20,
        "{ \forall }"{above}
      ]
      \ar[
        ddr,
        bend right=20,
        "{ \forall }"{left, xshift=-2pt}
      ]
      \ar[
       dr,
       dashed,
       "\exists !"{description}
      ]
      \\[-6pt]
      &[-6pt]
      \TopologicalSpace \times_{\mathrm{B}} \mathrm{P}
      \ar[r]
      \ar[d]
      \ar[
        dr,
        phantom,
        "\mbox{\color{orangeii} \tiny\rm(pb)}"{description}
      ]
      \ar[out=180-66, in=66, looseness=3.5, "\scalebox{.77}{$\mathclap{
        G
      }$}"{description},shift right=1]
      &
      \mathrm{P}
     \ar[out=180-66, in=66, looseness=3.5, "\scalebox{.77}{$\mathclap{
        G
      }$}"{description},shift right=1]
      \ar[d]
      \\
      &
      \TopologicalSpace
      \ar[r]
      \ar[out=-180+66, in=-66, looseness=3.5, "\scalebox{.77}{$\mathclap{
        G
      }$}"{description},shift left=1]
      \ar[
        r,
        phantom,
        "{
          \mbox{
            \tiny
            \color{darkblue}
            \bf
              fiber product {\color{black}/}
              pullback
          }
        }"{below, yshift=-15pt}
      ]
      &
      \mathrm{B}
      \ar[out=-180+66, in=-66, looseness=3.5, "\scalebox{.77}{$\mathclap{
        G
      }$}"{description},shift left=1]
    \end{tikzcd}
    {\phantom{AAAAA}}
    \mbox{e.g.}
    {\phantom{AAAAA}}
    \begin{tikzcd}
      Q
      \ar[out=180-66, in=66, looseness=3.5, "\scalebox{.77}{$\mathclap{
        G
      }$}"{description},shift right=1]
      \ar[
        drr,
        bend left=20,
        "{ \forall }"{above}
      ]
      \ar[
        ddr,
        bend right=20,
        "{ \forall }"{left, xshift=-2pt}
      ]
      \ar[
       dr,
       dashed,
       "\exists !"{description}
      ]
      \\[-6pt]
      &[-6pt]
      \TopologicalSpace \times \mathrm{P}
      \ar[r]
      \ar[d]
      \ar[
        dr,
        phantom,
        "\mbox{\color{orangeii} \tiny\rm(pb)}"{description}
      ]
      \ar[out=180-66, in=66, looseness=3.5, "\scalebox{.77}{$\mathclap{
        G
      }$}"{description},shift right=1]
      &
      \mathrm{P}
     \ar[out=180-66, in=66, looseness=3.5, "\scalebox{.77}{$\mathclap{
        G
      }$}"{description},shift right=1]
      \ar[d]
      \\
      &
      \TopologicalSpace
      \ar[r]
      \ar[out=-180+66, in=-66, looseness=3.5, "\scalebox{.77}{$\mathclap{
        G
      }$}"{description},shift left=1]
      \ar[
        r,
        phantom,
        "{
          \mbox{
            \tiny
            \color{darkblue}
            \bf
            product
          }
        }"{below, yshift=-15pt}
      ]
      &
      \ast
      \ar[out=-180+66, in=-66, looseness=3.5, "\scalebox{.77}{$\mathclap{
        G
      }$}"{description},shift left=1]
    \end{tikzcd}
  $$

  \vspace{-2mm}
  \noindent
  Over the {\it terminal object}, denoted by a point:
  $
    \begin{tikzcd}
      Q
      \ar[out=180-66, in=66, looseness=3.5, "\scalebox{.77}{$\mathclap{
        G
      }$}"{description},shift right=1]
      \ar[
        r,
        dashed,
        "{ \exists ! }"
      ]
      &
      \ast
      \ar[out=180-66, in=66, looseness=3.5, "\scalebox{.77}{$\mathclap{
        G
      }$}"{description},shift right=1]
    \end{tikzcd}
  $,
  a fiber product is a plain {\it product}, as shown on the right.
\end{notation}
\begin{proposition}[Finite limits, {e.g. \cite[Prop. 2.8.2]{Borceux94I}}]
  \label{FiniteLimits}
    In the presence of a terminal object,
  every {\it finite limit} is an iteration of
  pullbacks (Ntn. \ref{CartesianSquares}).
\end{proposition}

\begin{example}[Pullback preserves isomorphisms]
  \label{PullbackPreservesIsomorphisms}
  A commuting square with a bottom isomorphism
  (Ntn. \ref{BasicNotationForCategories})
  is a pullback square
  (Ntn. \ref{CartesianSquares}) if and only if also the top morphism is an
  isomorphism:
  \vspace{-2mm}
  $$
    \begin{tikzcd}[row sep=small, column sep=large]
      \mathrm{A}
      \ar[out=180-66, in=66, looseness=3.5, "\scalebox{.77}{$\mathclap{
        G
      }$}"{description},shift right=1]
      \ar[
        r,
        "f"
      ]
      \ar[d]
      \ar[
        dr,
        phantom,
        "\mbox{\tiny\rm(pb)}"{description}
      ]
      &
      \mathrm{P}
      \ar[out=180-66, in=66, looseness=3.5, "\scalebox{.77}{$\mathclap{
        G
      }$}"{description},shift right=1]
      \ar[d]
      \\
      \TopologicalSpace
      \ar[out=-180+66, in=-66, looseness=4.2, "\scalebox{.77}{$\mathclap{
        G
      }$}"{description},shift left=1]
      \ar[
        r,
        "\sim"{below}
      ]
      &
      \mathrm{B}
      \ar[out=-180+66, in=-66, looseness=4.2, "\scalebox{.77}{$\mathclap{
        G
      }$}"{description},shift left=1]
    \end{tikzcd}
    {\phantom{AAA}}
      \Leftrightarrow
    {\phantom{AAA}}
    \begin{tikzcd}
      \mathrm{A}
      \ar[
        r,
        "f"{above},
        "\sim"{below}
      ]
      &
      \mathrm{P}\;.
    \end{tikzcd}
  $$
\end{example}

\begin{proposition}[Right adjoint functors preserve limits]
  \label{RightAdjointFunctorsPreserveFiberProducts}
  A right adjoint functor $R$ (Ntn. \ref{AdjointFunctors}) preserves all limits,
  in particular it preserves all finite limits (Prop. \ref{FiniteLimits})
  and hence terminal objects and
  pullbacks (Ntn. \ref{CartesianSquares}):
  \vspace{-2mm}
  $$
    R
    (
      \TopologicalSpace \times_{\mathrm{B}} \mathrm{P}
    )
    \;\simeq\;
    R(\TopologicalSpace)
      \times_{R(\mathrm{B})}
    R(\mathrm{P})
    \,.
  $$
\end{proposition}

\begin{proposition}[Pasting law (e.g. {\cite[Prop. 11.10]{AHS90}})]
  \label{PastingLaw}
  Given two adjacent commuting squares where the right one is a pullback
  (Ntn. \ref{CartesianSquares})
  then the left square is a pullback if and
  only if the total rectangle is:
  \vspace{-3mm}
  $$
    \begin{tikzcd}[row sep=small, column sep=large]
      \mathrm{P}_1
      \ar[out=180-66, in=66, looseness=3.5, "\scalebox{.77}{$\mathclap{
        G
      }$}"{description},shift right=1]
      \ar[r]
      \ar[d]
      &
      \mathrm{P}_2
      \ar[out=180-66, in=66, looseness=3.5, "\scalebox{.77}{$\mathclap{
        G
      }$}"{description},shift right=1]
      \ar[r]
      \ar[d]
      \ar[
        dr,
        phantom,
        "\mbox{\tiny\rm(pb)}"
      ]
      &
      \mathrm{P}_3
      \ar[out=180-66, in=66, looseness=3.5, "\scalebox{.77}{$\mathclap{
        G
      }$}"{description},shift right=1]
      \ar[d]
      \\
      \TopologicalSpace_1
      \ar[out=-180+66, in=-66, looseness=4.2, "\scalebox{.77}{$\mathclap{
        G
      }$}"{description},shift left=1]
      \ar[r]
      &
      \TopologicalSpace_2
      \ar[out=-180+66, in=-66, looseness=4.2, "\scalebox{.77}{$\mathclap{
        G
      }$}"{description},shift left=1]
      \ar[r]
      &
      \TopologicalSpace_3
      \ar[out=-180+66, in=-66, looseness=4.2, "\scalebox{.77}{$\mathclap{
        G
      }$}"{description},shift left=1]
    \end{tikzcd}
  $$
\end{proposition}

\begin{notation}[Effective epimorphism {\cite[p. 101]{Grothendieck61}\cite[Def. 2.5.3]{Borceux94II}}]
  \label{EffectiveEpimorphism}
  A morphism $p$ is called a {\it regular epimorphism}
  if it is the coequalizer of \emph{some} parallel pair of morphisms,
  and an {\it effective epimorphism},
  to be denoted by double-headed arrows,
  if it is the
  coequalizer specifically of the two projections out of
  its pullback (Ntn. \ref{CartesianSquares}) along itself:
  \vspace{-2mm}
  $$
    \begin{tikzcd}
      \mathrm{P}
        \times_{\TopologicalSpace}
      \mathrm{P}
      \ar[
        rr,
        shift left=3pt,
        "\mathrm{pr}_1"{above}
      ]
      \ar[
        rr,
        shift right=3pt,
        "\mathrm{pr}_2"{below}
      ]
      &&
      \mathrm{P}
      \ar[
        rr,
        ->>,
        "p"{above},
        "\mathrm{coeq}"{below}
      ]
      &&
      \TopologicalSpace \;.
    \end{tikzcd}
  $$
\end{notation}

\begin{definition}[Regular categories ({\cite{BarrGrilletyOsdol}\cite[\S 2]{Borceux94II}\cite{Gran21}})]
  \label{RegularCategory}
  A category is called {\it regular} if

\noindent {\bf (i)}
for every morphism $X \to Y$,

\vspace{-3mm}
\begin{itemize}
\setlength\itemsep{-2pt}
  \item[{\bf (a)}]
  the fiber product
  $X \times_Y X$ exists (the ``kernel pair'');

  \item[{\bf (b)}]
  the coequalizer
  $
    X \times_Y X \rightrightarrows X \xrightarrow{\;\mathrm{coeq}\;} X/(X\times_Y X)
  $
  exists (the {\it image});
\end{itemize}

\vspace{-2mm}
\noindent {\bf (ii)}
pullbacks of regular epimorphisms (Ntn.  \ref{EffectiveEpimorphism})
exist and are again regular epimorphisms.
\end{definition}
In this Part I we are interested in the exceptional example of the category of compactly-generated topological spaces (Prop. \ref{CompactyGeneratedTopologicalSpacesFormARegularCategory} below). A generic class of examples of regular categories are toposes (such as the category of simplicial sets, which is of interest in Part II, see Ntn. \ref{SimplicialSets} below):
\begin{example}[Toposes are regular {\cite[p. 17]{BarrGrilletyOsdol}\cite[Prop. 3.4.14]{Borceux94III}\cite[p. 92]{Johnstone02a}}]
  \label{ToposesAreRegular}
  Every topos (hence in particular every category of presheaves) is a regular category (Def. \ref{RegularCategory}). Moreover, in toposes the classes of (i) epimorphisms, (ii) regular epimorphisms and (iii) effective epimorphisms
  (Ntn. \ref{EffectiveEpimorphism})
  all coincide (e.g. \cite[\S IV.7, Thm. 8 (p. 197)]{MacLaneMoerdijk92}\cite[Prop. 3.4.13, 3.4.15]{Borceux94III}).
\end{example}

\begin{lemma}[Effective epimorphisms in regular categories
(e.g. {\cite[Prop. 2.5.7]{Borceux94I}\cite[Prop. 2.3.3]{Borceux94I}})]
  \label{EffectiveEpimorphismsArePreservedByPullbackInRegularCategories}
  In a regular category (Def. \ref{RegularCategory}),
  the notions of regular and of effective epimorphisms (Ntn. \ref{EffectiveEpimorphism})
  coincide, and pullback along any morphism $f$ preserves effective epimorphisms $p$
  together with their coequalizer diagrams:
  $$
    \begin{tikzcd}[column sep=large]
      (f^\ast\TopologicalPrincipalBundle)
        \times_{\TopologicalSpace'}
      (f^\ast\TopologicalPrincipalBundle)
      \ar[d, shift right=4pt, "\mathrm{pr}_1"{swap}]
      \ar[d, shift left=4pt, "\mathrm{pr}_2"]
      \ar[r]
      \ar[dr, phantom, "\mbox{\tiny\rm(pb)}"{pos=.4}]
      &
      \TopologicalPrincipalBundle
        \times_{\TopologicalSpace}
      \TopologicalPrincipalBundle
      \ar[d, shift right=4pt, "\mathrm{pr}_1"{swap}]
      \ar[d, shift left=4pt, "\mathrm{pr}_2"]
      \\
      f^\ast\TopologicalPrincipalBundle
      \ar[d, ->>, "f^\ast p", "\mathrm{coeq}"{swap}]
      \ar[r]
      \ar[
        dr,
        phantom,
        "\mbox{\tiny\rm (pb)}"
      ]
      &
      \TopologicalPrincipalBundle
      \ar[d, ->>, "p", "\mathrm{coeq}"{ swap}]
      \\[+10pt]
      \TopologicalSpace'
      \ar[r, "f"{swap}]
      &
      \TopologicalSpace
    \end{tikzcd}
  $$
\end{lemma}
In regular categories, there are partial reverses
to the implications of
Ex. \ref{PullbackPreservesIsomorphisms}
and Prop. \ref{PastingLaw} (see also the $\infty$-category theoretic version in Lem. \ref{ReversePastingLawForInfinityPullbacks} below):
\begin{lemma}[Reverse pasting law in regular categories (e.g. {\cite[Lem. 1.15]{Gran21}})]
  \label{ReversePastingLawInRegularCategories}
  Given a commuting diagram
  in a regular category (Def. \ref{RegularCategory})
  of the form
  \vspace{-1mm}
  $$
    \begin{tikzcd}
      {}
      \ar[r]
      \ar[d]
      \ar[dr, phantom, "\mbox{\tiny \rm (pb)}"]
      &
      {}
      \ar[r]
      \ar[d]
      &
      {}
      \ar[d]
      \\
      {}
      \ar[r, ->>]
      &
      {}
      \ar[r]
      &
      {}
    \end{tikzcd}
  $$

  \vspace{-1mm}
  \noindent
  where the left square is Cartesian (Ntn. \ref{CartesianSquares})
  and the bottom left morphism is an effective epimorphism
  (Ntn. \ref{EffectiveEpimorphism}), then
  the right square is Cartesian if and only if the total rectangle is
  Cartesian.
\end{lemma}

\begin{lemma}[Local recognition of isomorphisms in regular categories]
  \label{IsomorphismsOnRegularCategoriesAreDetectedOnEffectiveCovers}
  In a regular category (Def. \ref{RegularCategory}),
  if the pullback $p^\ast f$ of a morphism $f$
  along an effective epimorphism
  $p$ (Ntn.  \ref{EffectiveEpimorphism})
  is an isomorphism, then $f$ was already an isomoprhism itself.
\end{lemma}
\begin{proof}
By assumption, we have a pullback square as on the bottom of the following diagram:
   \vspace{-1mm}
  $$
    \begin{tikzcd}[column sep=large]
      \widehat X
        \times_X
      \widehat X
      \ar[r, "\sim"]
      \ar[d, shift left=3pt]
      \ar[d, shift right=3pt]
      \ar[dr, phantom, "\mbox{\tiny\rm (pb)}"{pos=.4}]
      &
      \widehat Y
        \times_Y
      \widehat Y
      \ar[d, shift left=3pt]
      \ar[d, shift right=3pt]
      \\
      \widehat{X}
      \ar[r, "\sim"{swap}, "p^\ast f"]
      \ar[d, ->>, "f^\ast p"{swap}]
      \ar[dr, phantom, "\mbox{\tiny\rm (pb)}"]
      &
      \widehat{Y}
      \ar[d, ->>, "p"]
      \\
      X
      \ar[r, "f"{below}]
      &
      Y
    \end{tikzcd}
  $$

  \vspace{-1mm}
  \noindent
  Here the bottom left morphism is an effective epimorphism by
  Lem. \ref{EffectiveEpimorphismsArePreservedByPullbackInRegularCategories},
  since $p$ is so by assumption.
  Since limits commute with each other,
  we get the top pullback squares, where the topmost morphism
  is an isomorphism as the pullback of the isomorphism $p^\ast f$
  (by Ex. \ref{PullbackPreservesIsomorphisms}).
  By the nature of effective epimorphisms (Ntn.  \ref{EffectiveEpimorphism}),
  this now exhibits $f$ as the
  image under passage to coequalizers of an isomorphism of coequalizer
  diagrams, hence as an isomorphism.
\end{proof}

\medskip

\noindent
{\bf Topological spaces.}
We use the following {\it convenient category of topological spaces}
\cite{Steenrod67} which has become the standard foundation for algebraic topology:

\begin{notation}[Category of compactly-generated topological spaces]
  \label{CompactlyGeneratedTopologicalSpaces}
  We write
  \vspace{-3mm}
  \begin{equation}
    \label{CategoryOfTopologicalSpaces}
    \begin{tikzcd}
    \overset{
      \mathclap{
      \raisebox{4pt}{
        \tiny
        \color{darkblue}
        \bf
        \begin{tabular}{c}
          topological spaces
        \end{tabular}
      }
      }
    }{
    \categorybox{\TopologicalSpaces}
    }
    \;\;
    \ar[
      r,
      phantom,
      "{ \scalebox{.7}{$\bot$} }"
    ]
    \ar[
      r,
      shift right=5.5pt,
      "{ k }"{below}
    ]
    &
    \;\;
    \overset{
      \mathclap{
      \raisebox{4pt}{
        \tiny
        \color{orangeii}
        \bf
        \begin{tabular}{c}
          topological k-spaces
        \end{tabular}
      }
      }
    }{
    \categorybox{\kTopologicalSpaces}
    }
    \ar[
      l,
      hook',
      shift right=5.5pt
    ]
    &
    \;\;
    \overset{
      \mathclap{
      \raisebox{4pt}{
        \tiny
        \color{darkblue}
        \bf
        \begin{tabular}{c}
          Hausdorff k-spaces
        \end{tabular}
      }
      }
    }{
      \categorybox{\kHausdorffSpaces}
    }
    \;\;
    \ar[
      l,
      hook'
    ]
    &
    \;\;
    \overset{
      \mathclap{
      \raisebox{3pt}{
        \tiny
        \color{darkblue}
        \bf
        \begin{tabular}{c}
          locally compact
          \\
          Hausdorff spaces
        \end{tabular}
      }
      }
    }{
      \categorybox{\LocallyCompactHausdorffSpaces}
    }
    \ar[
      l,
      hook'
    ]
    \end{tikzcd}
  \end{equation}

  \vspace{-1mm}
  \noindent
  for
  the coreflective subcategory of
  those topological spaces
  which are the colimits of all images of compact spaces inside them
  ({\it k-spaces} \cite[\S 1]{Gale50}),
  and its full subcategory of
  Hausdorff spaces among these
  (subsuming all locally compact Hausdorff spaces),
  hence of {\it compactly generated} topological spaces
  (e.g. \cite[\S XI.9]{Dugundji66}\cite{Steenrod67}\cite{Lewis78} \cite[\S 3.4]{HerrlichStrecker97},
  concise practical review is in \cite[\S 0]{FHT00},
  and specifically for the equivariant context in \cite[\S 16]{LueckUribe14}).
\end{notation}

\begin{remark}[Mapping spaces]
  In particular, the category of k-spaces is Cartesian closed,
  which means that, for $\TopologicalSpace, \mathrm{Y} \,\in\, \kTopologicalSpaces$
  \eqref{CategoryOfTopologicalSpaces},
  the {\it mapping space}
  \vspace{-2mm}
  \begin{equation}
    \label{MappingSpace}
    \mathrm{Maps}(\TopologicalSpace,\mathrm{Y})
    \;\;\;
    \in
    \;
    \kTopologicalSpaces
  \end{equation}

  \vspace{-1mm}
  \noindent
  (namely the set of continuous functions $\TopologicalSpace \xrightarrow{\;} \mathrm{Y}$
  equipped with the $k$-ified compact-open topology) serves as an
  {\it exponential object}
  or
  {\it cartesian internal hom},
  (e.g., \cite[\S 9]{Niefield78}\cite[Thm. B.10]{Piccinini92}\cite[\S 7.1-7.2]{Borceux94II}).
  That is,
  $\mathrm{Maps}(\TopologicalSpace,-)$ is right adjoint (Ntn. \ref{AdjointFunctors})
  to forming the categorical product (the k-ified product topological space) with $\TopologicalSpace$:

  \vspace{-5mm}
  \begin{equation}
    \label{MappingSpaceAdjunction}
    \begin{tikzcd}[column sep=large]
      \kTopologicalSpaces
      \ar[
        rr,
        shift right=5pt,
        "{
          \mathrm{Maps}(\TopologicalSpace,\,-)
        }"{below}
      ]
      \ar[
        rr,
        phantom,
        "{
          \scalebox{.7}{$\bot$}
        }"
      ]
      &&
      \kTopologicalSpaces
      \mathrlap{\,.}
      \ar[
        ll,
        shift right=5pt,
        "{
          \TopologicalSpace \times (-)
        }"{above}
      ]
    \end{tikzcd}
  \end{equation}
  \vspace{-.4cm}
\end{remark}
Besides making the adjunction \eqref{MappingSpaceAdjunction} work,
compactly generated topological spaces
behave essentially like plain topological spaces:

\begin{remark}[Colimits of compactly generated topological spaces]
Colimits of k-spaces (Ntn. \ref{CompactlyGeneratedTopologicalSpaces})
are computed as
usual colimits of topological spaces.
For instance:

\noindent
{\bf (i)}
Orbit spaces --
i.e., the usual quotient topological spaces of continuous group actions
(see Lem. \ref{HausdorffQuotientSpaces} and Cor. \ref{QuotientCoprojectionOfFreeProperActionIsLocallyTrivial} below)
--
are (split) coequalizers in $\kTopologicalSpaces$ \eqref{CategoryOfTopologicalSpaces}.

\noindent
{\bf (ii)}
It is still true that
the functor
  $\pi_0 \,\colon\, \kTopologicalSpaces \xrightarrow{\;} \Sets$
  assigning sets of path-connected components
  preserves coequalizers
  (quotients by relations) and
  it preserves finite products:
  \vspace{-1mm}
  \begin{equation}
    \label{PathConnectedComponentsPreserveQuotientsAndFiniteProducts}
    \TopologicalSpace, \mathrm{Y}
    \;\;
    \in
    \;
    \kTopologicalSpaces
    \;\;\;\;
    \Rightarrow
    \;\;\;\;
    \left\{\!\!\!\!
    \begin{array}{rcl}
    \pi_0
    \left(
      \mathrm{coeq}
      (\mathrm{R} \rightrightarrows \TopologicalSpace)
    \right)
    &
    \;\simeq\;
    &
    \mathrm{coeq}
    \left(
      \pi_0(\mathrm{R}) \rightrightarrows \pi_0(\TopologicalSpace)
    \right)
    \\
    \pi_0
    (
      \TopologicalSpace
        \times
      \mathrm{Y}
    )
    &
    \;
    \simeq
    \;
    &
    \pi_0(\TopologicalSpace)
      \times
    \pi_0(\mathrm{Y})
    \end{array}
    \right.
    \;\;\;\;
    \in
    \;\;
    \Sets
    \,.
  \end{equation}
\end{remark}

\begin{remark}[Relation to locally compact Hausdorff spaces]
$\,$

\noindent
{\bf (i)}
If $\TopologicalSpace \,\in\,
 \LocallyCompactHausdorffSpaces \xhookrightarrow{\;} \kTopologicalSpaces$
 \eqref{CategoryOfTopologicalSpaces}
then, for any $\mathrm{Y} \,\in\, \kTopologicalSpaces$,
the usual product topological space is the category-theoretic product
$\TopologicalSpace \!\times\! \mathrm{Y} \,\in\, \kTopologicalSpaces$
in k-spaces
(\cite[Lem. 2.4]{Lewis78}\cite[Thm. B.6]{Piccinini92}).

\noindent
{\bf (ii)}
Any topological space is a k-space if and only if it is a
quotient topological space of a locally compact Hausdorff space.
(\cite[Lem. 3.2 (v)]{EscardoLawsonSimpson04}, strengthening \cite[\S XI, Thm. 9.4]{Dugundji66}\cite[Thm. B.4]{Piccinini92}).

\end{remark}

Moreover:
\begin{proposition}[Compactly-generated topological spaces form a regular category
{\cite[p. 3]{CagliariMantovaniVitale95}}]
 \label{CompactyGeneratedTopologicalSpacesFormARegularCategory}
 $\,$

 \noindent
 The categories $\kTopologicalSpaces$
 and $\kHausdorffSpaces$ (Ntn. \ref{CompactlyGeneratedTopologicalSpaces})
 are regular (Def. \ref{RegularCategory}).
\end{proposition}

  \begin{example}[Open covers are effective epimorphisms]
  \label{OpenCoversAreEffectiveEpimorphisms}
  For $\TopologicalSpace \,\in\, \kTopologicalSpaces$
  and $\{ \TopologicalPatch_i \xhookrightarrow{\;} \TopologicalSpace  \}_{i \in I}$
  an open cover, then the canonical map
  $\!\!
    \begin{tikzcd}
      \underset{i \in I}{\sqcup}
      \TopologicalPatch
      \ar[r, ->>]
      &
      \TopologicalSpace
    \end{tikzcd}
  \!\!$
  is an effective epimorphism (Def. \ref{EffectiveEpimorphisms}), in that
  the canonical maps
    \vspace{-3mm}
  $$
    \begin{tikzcd}
      \underset{j_1, j_2 \in I}{\sqcup}
      \TopologicalPatch_{j_1} \cap \TopologicalPatch_{j_2}
      \ar[r, shift left=3pt]
      \ar[r, shift right=3pt]
      &
      \underset{i \in I}{\coprod} \TopologicalPatch_i
      \ar[r]
      &
      \TopologicalSpace
    \end{tikzcd}
  $$

    \vspace{-2mm}
\noindent
  make a coequalizer diagram.
  (Namely, this  is the case on underlying sets, by the fact that the
  covering is surjective;
  hence the remaining condition
  is that a subset of $\TopologicalSpace$ is open precisely if its intersection with
  each $\TopologicalPatch_i$ is open, which is the case by the fact that the covering is
  by open subsets.)
\end{example}

\begin{example}[Isomorphism of bundles is detected on covers]
\label{IsomorphismOfBundlesDetectedOnOpenCovers}
Given a commuting diagram in $\kTopologicalSpaces$ (Ntn. \ref{CompactlyGeneratedTopologicalSpaces})
of the form
\vspace{-2mm}
$$
  \begin{tikzcd}[row sep=small]
    \TopologicalPrincipalBundle
    \ar[rr, "f"]
    \ar[dr]
    &&
    \TopologicalPrincipalBundle'
    \ar[dl]
    \\
    &
    \TopologicalSpace
  \end{tikzcd}
$$

\vspace{-2mm}
\noindent
then $f$ is an isomorphism as soon as it is locally so, hence
if its pullback to any open cover
$\widehat X \,\coloneqq\, \underset{i \in I}{\coprod} \TopologicalPatch_i$
is an isomorphism:
\vspace{-1mm}
$$
  \begin{tikzcd}[column sep={between origins, 60pt}]
    \TopologicalPrincipalBundle|_{\widehat X}
    \ar[rr]
    \ar[dr, "\sim"{sloped}]
    \ar[ddr, bend right=24]
    \ar[drrr, phantom, "\mbox{\tiny\rm(pb)}"{pos=.4}]
    &&
    \TopologicalPrincipalBundle
    \ar[dr, "f"]
    \ar[ddr, bend right=24]
    \\
    &
    \TopologicalPrincipalBundle'|_{ \widehat \TopologicalSpace}
    \ar[rr, ->>, crossing over]
    \ar[d]
    \ar[drr, phantom, "\mbox{\tiny\rm(pb)}"{pos=.3}]
    &{}&
    \TopologicalPrincipalBundle'
    \ar[d]
    \\
    &
    \widehat X
    \ar[rr, ->>]
    &&
    \TopologicalSpace
    \,.
  \end{tikzcd}
$$

\vspace{0mm}
\noindent
Namely,
here the bottom morphism is an effective epimorphism by
Ex. \ref{OpenCoversAreEffectiveEpimorphisms},
hence the middle morphism is an effective epimorphism by
regularity of $\kTopologicalSpaces$ (Prop. \ref{CompactyGeneratedTopologicalSpacesFormARegularCategory}).
Now with the rear square also the top square is Cartesian, by the
pasting law (Prop. \ref{PastingLaw}),
whence the top square exhibits the pullback of $f$ along
an effective epimorphism as an isomorphism, so that
Lem. \ref{IsomorphismsOnRegularCategoriesAreDetectedOnEffectiveCovers},
implies that it is itself an isomorphism.
\end{example}

\medskip

\noindent
{\bf Internal mathematical structures.} We find below that equivariant topology
is both a beautiful example of and is itself beautified by a systematic use
of {\it internalization} of mathematical structures into ambient categories
other than $\Sets$; a basic point that seems not to have received due attention.

\begin{notation}[Internal mathematical structures]
  \label{Internalization}
  For $S$ a mathematical structure
  expressible in terms of finite limits
  (a ``finite limit sketch''
  \cite{BastianiEhresmann72}\cite[\S 4]{BarrWells83}\cite[\S 1.49]{AdamekRosicky94}),
  hence by operations on fiber products (Prop. \ref{FiniteLimits}),
  and for $\mathcal{C}$ any category,

  \noindent {\bf (i)}   we write
  $S(\mathcal{C})$ for the category of $S$-models in $\mathcal{C}$,
  hence the category of $S$-structures
  {\it internal} to $\mathcal{C}$ in the
  original sense of \cite[p. 370]{Grothendieck60II}.

  \noindent {\bf (ii)}  For $F$ a functor that preserves finite limits (denoted {\it lex}),
  there is the evident induced functor on $S$-structures, which we denote as follows:
  \vspace{-2mm}
  \begin{equation}
    \label{FunctorOnStructuresInducedFromLexFunctor}
    F \,:\, \mathcal{C} \xrightarrow{\mathrm{lex}} \mathcal{D}
    \qquad
    \vdash
    \qquad
    S(F)
    \;:\;
    S(\mathcal{C})
    \longrightarrow
    S(\mathcal{D})
    \,.
  \end{equation}

  \end{notation}

  \begin{example}[Archetypical examples of internal structures]
  We have the categories:

  \vspace{-.1cm}
  \begin{itemize}

  \vspace{-.2cm}
  \item
  $\Groups(\mathcal{C})$
  of internal groups
  (e.g., \cite{EckmannHilton61}\cite{EckmannHilton62}\cite{EckmannHilton63}\cite[\S 4.1]{BarrWells83},
  see Def. \ref{EquivariantTopologicalGroup} below);

  \vspace{-.2cm}
  \item
  $\Actions{G}(\mathcal{C})$
  of internal group actions,
  (e.g., \cite[\S 7]{Boardman95}\cite[p. 8]{BorceuxJanelidzeKelly05}, see Def. \ref{EquivariantPrincipalBundle} below);

  \vspace{-.2cm}
  \item
  $\Groupoids(\mathcal{C})$
  of internal groupoids
  (e.g., \cite[\S 8]{Borceux94I}\cite[\S 1]{NiefieldPronk19}, see Ntn. \ref{TopologicalGroupoids} below);
  \end{itemize}
  \vspace{-.2cm}

  \noindent
  all originally due to \cite[\S 4]{Grothendieck61}.
  We consider these notions mainly internal to the category
  $\mathcal{C} = \GActionsOnTopologicalSpaces$  \eqref{GActionsOnTopologicalSpaces}
  of, in turn, group actions internal to
  topological spaces \eqref{CategoryOfTopologicalSpaces}.

\end{example}

\begin{notation}[Internalization of principal bundle theory]
  \label{InternalizationOfPrincipalBundleTheory}
  The key example of internal structures for our purposes is the category

  \vspace{-.1cm}
  \begin{itemize}

    \vspace{-.2cm}
    \item
    $\FormallyPrincipalBundles{\Gamma}(\mathcal{C})$
    of
    {\it formally principal bundles}
    (\cite[p. 312 (15 of 30)]{Grothendieck60}\cite[p. 9 (293)]{Grothendieck71},
    \\
    \phantom{$\FormallyPrincipalBundles{\Gamma}(\mathcal{C})$}
    also: {\it pseudo-torsors} \cite[\S 16.5.15]{Grothendieck67}),
  \end{itemize}
  \vspace{-.3cm}

 \noindent
 which is the subcategory of $\Actions{\Gamma}(\mathcal{C})$
 (Ntn. \ref{Internalization})
 on those actions that are
 either fiberwise principal or empty \eqref{PrincipalityConditionAsShearMapBeingAnIsomorphism},
 see Rem. \ref{PseudoTorsorCondition} below.
\end{notation}
We observe in Cor. \ref{InternalDefinitionOfGPrincipalBundlesCoicidesWithtomDieckDefinition}
 that these (formally) principal bundles, when internalized in
 $\mathcal{C} \,=\, \Actions{G}(\kTopologicalSpaces)$,
 are equivalently equivariant principal bundles
 in the original and general sense of \cite{tomDieck69},
 see Rem. \ref{LiteratureOnEquivariantPrincipalBundles} below.
 While this might not be a surprising observation for experts with the relevant background,
 it is, we find, absolutely foundational to the subject of equivariant bundle theory,
 and seems not to have been made before in existing literature.

\newpage

\section{$G$-Actions on topological spaces}
\label{TopologicalGActions}

\begin{notation}[Equivariant topology (``transformation groups'', ``$G$-spaces'', e.g.{\cite{Bredon72}\cite{tomDieck79}\cite{tomDieck87}})]
  \label{GActionOnTopologicalSpaces}
  $\,$

\noindent Throughout:

\vspace{-3mm}
\begin{itemize}
\setlength\itemsep{-2pt}
\item
$
  G \;\in\;
  \mathrm{Grps}
  (
    \kHausdorffSpaces
  )
$
denotes a Hausdorff topological group, with group operation denoted $(-)\cdot (-)$;

\item
$H \subset G$ denotes a topological subgroup
(necessarily Hausdorff, since $G$ is);

\item
$N\!(H) \subset G$ denotes its {\it normalizer subgroup}, and

\item $W\!(H) \coloneqq N\!(H)\!/H$ its {\it Weyl group} (e.g. \cite[p. 13]{May96}):

\vspace{-5mm}
\begin{equation}
 \label{NormalizerAndWeylGroup}
    \begin{tikzcd}
        \mathclap{
        \raisebox{3pt}{
          \tiny
          \color{darkblue}
          \bf
          \begin{tabular}{c}
            normalizer
            \\
            subgroup
          \end{tabular}
        }
        }
        \qquad
      {
        N(H)
      }
      \ar[
        r,
        ->>
      ]
      &
        N(H)/H
        \,=:\,
        W(H)
        \qquad
        \mathclap{
        \raisebox{1pt}{
          \tiny
          \color{darkblue}
          \bf
          \begin{tabular}{c}
            Weyl group
          \end{tabular}
        }
        }
    \end{tikzcd}
\end{equation}

  \vspace{-3mm}
  \item
  We write
  \vspace{-1mm}
  \begin{equation}
    \label{GActionsOnTopologicalSpaces}
    \hspace{-1cm}
    \GActionsOnTopologicalSpaces
    \;\coloneqq\;
    \left\{\!\!
      \left(
        \arraycolsep=1pt
        \begin{array}{ccl} \small
          \TopologicalSpace &\in& \kTopologicalSpaces \,,
          \\
          \rho &:& G \xrightarrow{\rho} \mathrm{Aut}(\TopologicalSpace)
        \end{array}
      \right)
    \!\!\right\}
    \;=\;
    \Big\{
      G \acts \, \TopologicalSpace
      \;\coloneqq\;
     \big(\!\!\!
     \begin{tikzcd}[column sep=40pt]
        G \times \TopologicalSpace
        \ar[
          r,
          "
         \scalebox{.7}{$   (-)\cdot(-) $}
          "{
            above
          },
          "
           \mbox{
             \tiny
             \color{greenii}
             \bf
             \begin{tabular}{c}
               continuous
               \\
               actions
             \end{tabular}
           }
          "{
            below
          }
        ]
        &\TopologicalSpace
      \end{tikzcd}
      \!\!\!
      \big)
      \!
    \Big\}
  \end{equation}

  \vspace{-3mm}
  \noindent
  for the category whose objects
  $G \acts \, \TopologicalSpace$ are topological spaces $\TopologicalSpace$
  equipped with continuous left $G$-actions
  and
  whose morphisms are $G$-equivariant contionuous functions between these
  (often: ``$G$-spaces'', for short).

  \item
  For $G \acts \, \TopologicalSpace \in \GActionsOnTopologicalSpaces$,
  we write

  \vspace{-.1cm}
  \begin{itemize}

    \vspace{-.2cm}
    \item[$\circ$]
    $
      \TopologicalSpace^G
        \,\coloneqq\,
      \big\{
        x \in \TopologicalSpace
          \,\vert\,
        \underset{g \in G}{\forall}
        \;
        g \cdot x \,=\, x
      \big\}
        \xhookrightarrow{\;}
      \TopologicalSpace
      \;\;\;\;\;\;\;\;\;\;\;\;\;\;
      \in
      \;
      \kTopologicalSpaces
    $
    for the {\it $G$-fixed subspace};

    \vspace{0mm}
    \item[$\circ$]
    $
      \TopologicalSpace
      \relbar\joinrel\twoheadrightarrow
      X_{G}
      \,\coloneqq\,
      X/G
      \,\coloneqq\,
      \left\{
        [x] \coloneqq G \cdot x
        \;\vert\;
        \TopologicalSpace \in \TopologicalSpace
      \right\}
      \;\;
      \in
      \;
      \kTopologicalSpaces
    $
    for the {\it $G$-quotient space} ({\it $G$-orbit space}).
  \end{itemize}
  \vspace{-.1cm}

  \item
For $G\acts \, \TopologicalSpace \in \GActionsOnTopologicalSpaces$ and $x \in \TopologicalSpace$,
its {\it isotropy subgroup} is denoted
\vspace{-2mm}
\begin{equation}
  \label{StabilizerSubgroupInEquivarianceGroup}
  G_x
  \;\coloneqq\;
  \mathrm{Stab}_{{}_{G}}(x)
  \;\coloneqq\;
  \{
    g \in G \,\vert\, g \cdot x = x
  \}
  \;\subset\;
  G\;.
\end{equation}

\end{itemize}
\end{notation}

From \cref{NotionsOfEquivariantLocalTrivialization} on, we make the
following further assumptions on the equivariance group:

\begin{assumption}[Proper equivariant topology (following {\cite{DHLPS19}\cite{SS20OrbifoldCohomology}})]
\label{ProperEquivariantTopology}
We speak of {\it proper equivariant topology} if:

\vspace{-4mm}
\begin{itemize}
\setlength\itemsep{-2pt}

\item
equivariance groups $G$ are Lie groups with compact connected components;

\item
subgroups $H \subset G$ are compact;

\item
domain spaces $\TopologicalSpace$ are locally compact and Hausdorff;

\item
equivariance actions $G \acts \, \TopologicalSpace$ are proper. \footnote{
  Under the previous assumption that domain spaces are locally compact and
  Hausdorff, all notions of proper actions agree \cite[Thm. 1.2.9]{Palais61};
  see also \cite[Rem. 5.2.4]{Karppinen16}.
}
\end{itemize}
\end{assumption}

\begin{lemma}[Equivariance subgroups in proper equivariant topology]
  \label{EquivarianceSubgroupsInProperEquivariantTopology}
  Under Assumption \ref{ProperEquivariantTopology},

  \vspace{-4mm}
  \begin{enumerate}[{\bf (i)}]
  \setlength\itemsep{-2pt}
  \item
  every
  $H \subset G$
  (namely every compact subgroup of a Lie group with compact connected components)
  is:

  \vspace{-4mm}
  \begin{enumerate}[{\bf (a)}]
\setlength\itemsep{-2pt}
   \item
   a closed subgroup;

   \item
   a compact Lie group;

  \end{enumerate}
  \vspace{-3mm}

  \item
  every $G_x \subset G$
  (namely the isotropy subgroup \eqref{StabilizerSubgroupInEquivarianceGroup}
  of a proper action at any point $x$)
  is of this form.

  \end{enumerate}
  \vspace{-.2cm}

\end{lemma}
\begin{proof}
  For the first statement,
  it is sufficient to consider the connected components of the neutral element.
  Here statement (a) follows since Lie groups are Hausdorff spaces
  and compact subspaces of compact Hausdorff spaces are equivalently closed subspaces.
  With this, statement (b) follows from Cartan's closed-subgroup theorem
  (e.g. \cite[Thm. 10.12]{Lee12}).
  The assumption that $\TopologicalSpace$ is locally compact and Hausdorff
  ensures that all notions of proper action agree, and it follows that
  all stabilizer subgroups of points are compact.
  With this, the second statement follows from the first.
\end{proof}

\medskip

\noindent
{\bf Basic examples of $G$-actions.} To fix notation and conventions,
we make explicit the following basic $G$-actions.

\begin{example}[Left and right-inverse multiplication action]
  \label{LeftAndInverseRightMultiplicationAction}
  Each $G \,\in\, \Groups(\kTopologicalSpaces)$
  carries canonical left $G$-actions \eqref{GActionsOnTopologicalSpaces},
  by left multiplication
  and by inverse right multiplication, respectively,

  \vspace{-.3cm}
  \begin{equation}
    \label{LeftMultiplicationAndInverseRightMultiplicationActionsOnATopologicalGroup}
    G \acts \,  G^L \,\in\, \GActionsOnTopologicalSpaces
    {\phantom{AAAAAAA}}
    G \acts  \, G^R \,\in\, \GActionsOnTopologicalSpaces
  \end{equation}
  \vspace{-.5cm}
  $$
    \begin{tikzcd}[row sep=-5pt, column sep=4pt]
\scalebox{0.8}{$      G \times G^L $}
      \ar[rr]
      &&
   \scalebox{0.8}{$     G^L $}
      \\
   \scalebox{0.8}{$     (g,h) $} &\longmapsto&  \scalebox{0.8}{$   g \cdot h $}
    \end{tikzcd}
    {\phantom{AAAAAAA}}
    \begin{tikzcd}[row sep=-5pt, column sep=4pt]
    \scalebox{0.8}{$    G \times G^R $}
      \ar[rr]
      &&
   \scalebox{0.8}{$     G^R $}
      \\
 \scalebox{0.8}{$       (g,h) $} &\longmapsto& \scalebox{0.8}{$  h \cdot g^{-1}$}
      \mathrlap{\,.}
    \end{tikzcd}
  $$

  \vspace{-1mm}
\noindent  Under inversion, these two actions are isomorphic:

  \vspace{-.4cm}
  \begin{equation}
    \label{IsomorphismBetweenLeftMultiplicationAndInverseRightMultiplicationAction}
    \begin{tikzcd}
      G^L
      \ar[
        r,
        "{ (-)^{-1} }"{above}
        "{\sim}"{below, yshift=+1pt}
      ]
      &
      G^R
      \;\;\;
      \in
      \;
      \GActionsOnTopologicalSpaces \;.
    \end{tikzcd}
  \end{equation}
  \vspace{-.3cm}
\end{example}

\begin{example}[Diagonal action]
  \label{DiagonalActionOnProductGSpaces}
  For
  $
    G \acts \, \TopologicalSpace_1, \, G \acts \, \TopologicalSpace_2
      \,\in\,
    \GActionsOnTopologicalSpaces
  $,
  one can consider the {\it diagonal action}
  $G \acts \,  (\TopologicalSpace_1 \times \TopologicalSpace_2)$
  on the product space of the underlying spaces:
    \vspace{-3mm}
  $$
    \begin{tikzcd}[row sep=-5pt]
      G \times (\TopologicalSpace_1 \times \TopologicalSpace_2)
      \ar[
        rr
      ]
      &&
\scalebox{0.8}{$        \TopologicalSpace_1 \times \TopologicalSpace_2$}
      \\
\scalebox{0.8}{$        \left(g, (x_1, x_2)\right) $}
      &\longmapsto&
\scalebox{0.8}{$        \left( g\cdot x_1, \, g \cdot x_2\right)$}.
    \end{tikzcd}
  $$
\end{example}

\begin{example}[Conjugation action on mapping space (e.g. {\cite[p. 5]{GuillouMayRubin13}})]
  \label{ConjugationActionOnMappingSpaces}
  Let $G \acts \, \TopologicalSpace_1, \, G \acts \, \TopologicalSpace_2 \,\in \, \GActionsOnTopologicalSpaces$.

  \noindent
  {\bf (i)}
  The mapping space \eqref{MappingSpace} of the underlying topological spaces
  carries a $G$-action
  given by {\it conjugation}:

  \vspace{-.3cm}
  \begin{equation}
    \label{ConjugationActionOnMappingSpace}
    G
      \acts \;
    \mathrm{Maps}(\TopologicalSpace_1, \TopologicalSpace_2)
    \;\;\;
    \in
    \;
    \GActionsOnTopologicalSpaces
  \end{equation}
  \vspace{-.6cm}

  \vspace{-.5cm}
  \begin{equation}
    \label{ConjugationActionOnMapsBetweenGSpaces}
    \hspace{-2mm}
    f \,\in\, \mathrm{Maps}(\TopologicalSpace_1, \TopologicalSpace_2)
    \;\;\;\;\;\;
    \vdash
    \;\;\;\;\;\;
    \underset{g \in G}{\forall}
   \quad
    \begin{tikzcd}
      \TopologicalSpace_1
      \ar[
        rr,
        "{
          g \cdot f
        }"
      ]
      \ar[
        d,
        "{
          g^{-1}\cdot(-)
        }"{left}
      ]
      &&
      \TopologicalSpace_2
      \\
      \TopologicalSpace_1
      \ar[
        rr,
        "{
          f
        }"{below}
      ]
      &&
      \TopologicalSpace_2
      \ar[
        u,
        "{
          g \cdot (-)
        }"
      ]
    \end{tikzcd}
    {\phantom{AA}}
    \mbox{i.e.,}
    \;\;\;\;
    \underset{x \in \TopologicalSpace_1}{\forall}
    \;
    (g \cdot f)(x)
    \;=\;
    g
      \cdot
    (
    f
    (
      g^{-1} \cdot x
    )
    )\;.
  \end{equation}
  \vspace{-.4cm}

  \noindent
  {\bf (ii)} The fixed locus of the conjugation action
  \eqref{ConjugationActionOnMapsBetweenGSpaces} is the subspace of
  {\it $G$-equivariant functions}

  \vspace{-3mm}
  \begin{equation}
    \label{EquivariantFunctions}
    \overset{
      \mathclap{
      \raisebox{3pt}{
        \tiny
        \color{darkblue}
        \bf
        \def\arraystretch{.9}
        \begin{tabular}{c}
          subspace of $G$-equivariant maps
        \end{tabular}
      }
      }
    }{
    \big\{
      f \,\in\,
      \mathrm{Maps}(\TopologicalSpace_1, \TopologicalSpace_2)
      \,\big\vert\,
      f(-) = g^{-1}\cdot f(g\cdot -)
    \big\}
    }
    \;\;
      =
    \;\;
    \overset{
      \raisebox{3pt}{
        \tiny
        \color{darkblue}
        \bf
        \begin{tabular}{c}
          $G$-fixed subspace
          \\
          of conjugation action
        \end{tabular}
      }
    }{
      \mathrm{Maps}
      (
        \TopologicalSpace_1,
        \,
        \TopologicalSpace_2
      )^G
    }
   \;\;  \xhookrightarrow{\quad}
   \;\;
    \overset{
      \mathclap{
      \raisebox{3pt}{
        \tiny
        \color{darkblue}
        \bf
        \def\arraystretch{.9}
        \begin{tabular}{c}
          space of all
          \\
          continuous maps
        \end{tabular}
      }
      }
    }{
      \mathrm{Maps}
      (
        \TopologicalSpace_1,
        \,
        \TopologicalSpace_2
      )\;.
    }
  \end{equation}
  \vspace{-.4cm}

\noindent
{\bf (iii)}
This construction \eqref{ConjugationActionOnMappingSpace} is functorial in both
arguments, contravariantly in the first.
With \eqref{IsomorphismBetweenLeftMultiplicationAndInverseRightMultiplicationAction}
with means, in particular, for $G \acts  \, X \in \GActionsOnTopologicalSpaces$ that

\vspace{-3mm}
\begin{equation}
  \mathrm{Maps}
  (
    G^L,
    \TopologicalSpace
  )
  \;\simeq\;
  \mathrm{Maps}
  (
    G^R,
    \TopologicalSpace
  )
  \;\;\;
  \in
  \;
  \GActionsOnTopologicalSpaces \,.
\end{equation}

\vspace{-1mm}
\noindent
{\bf (iv)}
With the first argument fixed, this construction
\eqref{ConjugationActionOnMappingSpace} is a right adjoint to the
product operation from Ex. \ref{DiagonalActionOnProductGSpaces}:
\vspace{-3mm}
\begin{equation}
  \label{InternalHomInGSpaces}
  \begin{tikzcd}[column sep=35pt]
    \Actions{G}(\kTopologicalSpaces)
    \ar[
      rr,
      shift right=6pt,
      "{
        \scalebox{.7}{$
          G \acts \, \mathrm{Maps}(\TopologicalSpace,\, -)
        $}
      }"{below}
    ]
    \ar[
      rr,
      phantom,
      "\scalebox{.7}{$\bot$}"
    ]
    &&
    \Actions{G}(\kTopologicalSpaces) \;.
    \ar[
      ll,
      shift right=6pt,
      "{
        \scalebox{.7}{$
          G \acts \,  \TopologicalSpace \,\times\, (-)
        $}
      }"{above}
    ]
  \end{tikzcd}
\end{equation}

\vspace{-.4cm}
\end{example}

\medskip

\noindent
{\bf Change of equivariance group.}
Much of our formulation of equivariant topology proceeds by applying the
{\it change of equivariance group}
adjoint triple from the following Lem. \ref{InducedAndCoinducedActions}, in numerous ways.

\begin{lemma}[Change of equivariance group (e.g. {\cite[\S I.1]{May96}\cite[p. 9]{DHLPS19}})]
  \label{InducedAndCoinducedActions}
  Given a continuous homomorphism of topological groups
  \vspace{-2mm}
  $$
    \begin{tikzcd}
      G_1
        \ar[
          r,
          "\phi"
        ]
      &
      G_2\;,
    \end{tikzcd}
  $$

  \vspace{-2mm}
  \noindent
  we have a triple of adjoint functors (Ntn. \ref{AdjointFunctors})
  between their categories of
  continuous actions (Ntn. \ref{GActionOnTopologicalSpaces}):

  \vspace{-3mm}
  \begin{equation}
    \label{AdjointTripleOfChangeOfEquivarianceGroup}
    \begin{tikzcd}[column sep=55pt]
      G_1\mathrm{Act}
      (
        \mathrm{TopSp}
      )
      \ar[
        rrr,
        shift right=13pt,
        "{
          \mathrm{Maps}
          (
            G_2,
            \,
        -
          )^{G_1}
          \;\coloneqq\;
          \mathrm{Maps}
          \left(
            \phi^\ast
            (
              G^L_2
            )
            ,\,
            -
          \right)^{G_1}
        }"{
          below
        },
        "\scalebox{.7}{$\bot$}"{above}
      ]
      \ar[
        rrr,
        shift left=13pt,
        "{
          G_2 \times_{{}_{G_1}} (-)
          \;\coloneqq\;
          \left(
            \phi^\ast
            (
              G^R_2
            )
              \times
            (-)
          \right)_{G_1}
        }"{
          above
        },
        "\scalebox{.7}{$\bot$}"{below}
      ]
      &&&
      G_2\mathrm{Act}
      (
        \mathrm{TopSp}
      )
      \mathrlap{\,,}
      \ar[
        lll,
        "\;\phi^\ast"{description}
      ]
    \end{tikzcd}
\end{equation}
\vspace{-.4cm}

\noindent
  where:

   -- the $G_1$-{\rm pullback action} on $\phi^\ast Y$ is through $\phi$,
    on the same underlying topological space;

  -- the {\rm induced $G_2$-action} on $G_2 \times_{G_1} X$ is
   that given by left multiplication of $G_2$ on the $G_2$-factor;

  -- the {\rm co-induced $G_2$-action} on $\mathrm{Maps}(G_2,X)^{G_1}$
   is given by right multiplication on the $G_2$-argument.

\end{lemma}

\begin{example}[Quotient spaces, fixed loci and trivial action]
  \label{QuotientAndFixedLociFromChangeOfGroupAdjunction}
  For $G \xrightarrow{\;} 1$ the unique group homomorphism to the trivial group, the corresponding pullback action (Lemma \ref{InducedAndCoinducedActions}) is the trivial $G$-action, whose adjoints \eqref{AdjointTripleOfChangeOfEquivarianceGroup}
  form the quotient space $(-)_G$ and the $G$-fixed space
  $(-)^G$ (Ntn. \ref{GActionOnTopologicalSpaces}), respectively:
  \begin{equation}
    \hspace{-.7cm}
    \begin{tikzcd}[column sep=50pt]
      \Actions{G}(\kTopologicalSpaces)
      \ar[
        rr,
        shift left=16pt,
        "{
          \mbox{
            \tiny
            \color{greenii}
            \bf
            quotient space
          }
        }"{description},
        "{
          (-)/G
        }"{yshift=1pt}
      ]
      \ar[
        from=rr,
        "{
          \mbox{
            \tiny
            \color{greenii}
            \bf
            trivial $G$-action
          }
        }"{description}
      ]
      \ar[
        rr,
        shift right=16pt,
        "{
          \mbox{
            \tiny
            \color{greenii}
            \bf
            fixed locus
          }
        }"{description},
        "{
          (-)^G
        }"{swap,yshift=-1pt}
      ]
      \ar[
        rr,
        phantom,
        "{\scalebox{.75}{$\bot$}}",
        shift left=8pt
      ]
      \ar[
        rr,
        phantom,
        "{\scalebox{.75}{$\bot$}}",
        shift right=8pt
      ]
      &&
      \kTopologicalSpaces
      \mathrlap{\,.}
    \end{tikzcd}
  \end{equation}
\end{example}

\begin{example}[Underlying topological spaces and (co-)free actions]
  \label{ForgettingGActionsAsPullbackAction}
  For $1 \xhookrightarrow{\;} G$ the unique inclusion of the trivial group,
  the corresponding pullback action (Lemma \ref{InducedAndCoinducedActions})
  is the forgetful functor from continuous $G$-actions to their underlying
  topological spaces, whose adjoints \eqref{AdjointTripleOfChangeOfEquivarianceGroup}
  form the free action and cofree action, respectively:
  \vspace{-2mm}
  \begin{equation}
    \label{FreeForgetfulAdjunctionForGAction}
    \hspace{.7cm}
    \begin{tikzcd}[column sep=50pt]
      \kTopologicalSpaces
      \ar[
        rr,
        shift left=2*8pt,
        "{
          \mbox{
            \tiny
            \color{greenii}
            \bf
            free action
          }
        }"{description},
        "{
          G \times (-)
        }"{above, yshift=+1pt}
      ]
      \ar[
        rr,
        shift right=2*8pt,
        "{
          \mbox{
            \tiny
            \color{greenii}
            \bf
            co-free action
          }
        }"{description},
        "{
          \mathrm{Maps}(G,-)
        }"{below, yshift=-1pt}
      ]
      \ar[
        rr,
        shift left=1*8pt,
        phantom,
        "\scalebox{.7}{$\bot$}"{description}
      ]
      \ar[
        rr,
        shift right=1*8pt,
        phantom,
        "\scalebox{.7}{$\bot$}"{description}
      ]
      \ar[
        rr,
        shift right=1*8pt,
        phantom,
        "\scalebox{.7}{$\bot$}"{description}
      ]
      &&
      \GActionsOnTopologicalSpaces
      \mathrlap{\;.}
      \ar[
        ll,
        "\mbox{\tiny\color{greenii} \bf forget $G$-action}"{description}
      ]
    \end{tikzcd}
  \end{equation}
\end{example}

\begin{lemma}[Forgetting $G$-action creates limits and colimits]
  \label{ForgetfulFunctorFromTopologicalGSpacesToGSpaces}
  The forgetful functor from topological $G$-actions
  to underlying topological spaces
  (Example \ref{ForgettingGActionsAsPullbackAction})
  {\it creates limits and colimits}, in that a diagram of
  topological $G$-actions is a (co)limiting (co)cone diagram precisely if
  its underlying diagram of topological spaces is:
  \vspace{-2mm}
  $$
    \begin{tikzcd}[row sep=-2pt]
      \GActionsOnTopologicalSpaces
      \ar[
        rr,
        "
         \mbox{
           \tiny
           \rm
           \color{greenii}
          \bf  forget $G$-action
         }
        "
      ]
      &&
      \kTopologicalSpaces
      \\
      \scalebox{0.8}{$
        G \acts \, \TopologicalSpace
        \;\simeq\;
        \limit{i}
        (
          G \acts \, \TopologicalSpace_i
        )
      $}
      &\Leftrightarrow&
      \scalebox{0.8}{$
        \TopologicalSpace
        \;\simeq\;
        \limit{i}
        (
          \TopologicalSpace_i
        )
      $}
      \\
      \scalebox{0.8}{$
        G \acts \, \TopologicalSpace
        \;\simeq\;
        \colimit{i}
        (
          G \acts \, \TopologicalSpace_i
        )
      $}
      &\Leftrightarrow&
      \scalebox{0.8}{$
        \TopologicalSpace
        \;\simeq\;
        \colimit{i}
        (
          \TopologicalSpace_i
        )
      $}
      \;.
    \end{tikzcd}
  $$
\end{lemma}
\begin{proof}
  That the forgetful functor {\it preserves} all limits
  and colimits
  follows
  (Prop. \ref{RightAdjointFunctorsPreserveFiberProducts})
  from it being a right and a left adjoint \eqref{FreeForgetfulAdjunctionForGAction}.
  A general abstract way to see that it also reflects
  and hence creates all limits and colimits is to notice that $G$-actions
  are algebras for the monad $G \times (-)$, and that monadic functors
  create all limits which exist in their codomain, and create all  colimits which
  exist and are preserved by the monad
  (e.g., \cite[pp. 137-138]{MacLane70}). But the monad here is the composite
  of the two left adjoints in the change of group adjoint triple
  \eqref{AdjointTripleOfChangeOfEquivarianceGroup}
  along $1 \to G$
  and hence preserves all colimits.
  The resulting claim also appears in \cite[\S B]{Schwede18}.
  For the record, we spell out the reflection of pullbacks/fiber products
  (Ntn. \ref{CartesianSquares}, the general proof is directly analogous),
   which is the main case of interest in \cref{EquivariantPrincipalTopologicalBundles}:

  Consider a commuting square of topological $G$-actions whose
  underlying square of topological spaces is a pullback, and consider
  a cone with tip $G \acts \, Q$ over this square, as shown on the left here:
  \vspace{-2mm}
  \begin{equation}
    \label{ReflectedPullbackOfGActions}
    \hspace{-8mm}
    \begin{tikzcd}
      Q
      \ar[out=180-66, in=66, looseness=3.5, "\scalebox{.77}{$\mathclap{
        G
      }$}"{description},shift right=1]
      \ar[
        drr,
        bend left=20
      ]
      \ar[
        ddr,
        bend right=20
      ]
      \ar[
       dr,
       dashed
      ]
      \\
      &
      \TopologicalSpace \times_{\mathrm{B}} \mathrm{P}
      \ar[r]
      \ar[d]
      \ar[out=180-66, in=66, looseness=3.5, "\scalebox{.77}{$\mathclap{
        G
      }$}"{description},shift right=1]
      &
      \mathrm{P}
     \ar[out=180-66, in=66, looseness=3.5, "\scalebox{.77}{$\mathclap{
        G
      }$}"{description},shift right=1]
      \ar[d]
      \\
      &
      \TopologicalSpace
      \ar[r]
      \ar[out=-180+66, in=-66, looseness=3.5, "\scalebox{.77}{$\mathclap{
        G
      }$}"{description},shift left=1]
      &
      \mathrm{B}
      \ar[out=-180+66, in=-66, looseness=3.5, "\scalebox{.77}{$\mathclap{
        G
      }$}"{description},shift left=1]
    \end{tikzcd}
    {\phantom{AAAAAAAA}}
    \begin{tikzcd}[column sep=tiny, row sep=-4pt]
      & &[10pt]
      & &
      G \times \mathrm{P}
      \ar[ddr]
      \ar[dddd]
      \\
      G
        \times
      \mathrm{Q}
      \ar[
        rr,
        dashed,
        "\mathrm{id} \times f"
      ]
      \ar[
        dddd
      ]
      &&
      G \times
      (
        \TopologicalSpace \times_{\mathrm{B}} \mathrm{P}
      )
      \ar[
        dddd
      ]
      \ar[urr]
      \ar[ddr]
      \\
      & & & & &
      G \times \mathrm{B}
      \ar[dddd, "\mbox{\tiny\color{greenii} group action}"{sloped}]
      \\
      & & &
      G \times \TopologicalSpace
      \ar[
        urr,
        crossing over
      ]
      \\
      & & & &
      \mathrm{P}
      \ar[ddr]
      \\
      \mathrm{Q}
      \ar[
        rr,
        dashed,
        "f"
      ]
      &&
      \TopologicalSpace \times_{\mathrm{B}} \mathrm{P}
      \ar[
        urr
      ]
      \ar[ddr]
      \\
      && &&&
      \mathrm{B}
      \\
      & & &
      \TopologicalSpace
      \ar[urr]
      \ar[
        from=uuuu,
        crossing over
      ]
    \end{tikzcd}
  \end{equation}

  \vspace{-2mm}
  \noindent
  We need to show that there exists a unique dashed morphism making the
  full diagram on the left commute. But since the underlying square of topological spaces
  is a pullback, there exists a unique such continuous function,
  in the bottom part of the diagram on the right of
  \eqref{ReflectedPullbackOfGActions}.
  Hence it remains to show that this unique function is necessarily $G$-equivariant.
  But the functor $G \times (-)$, preserves limits,
  being itself a limit,
  so that also the top square on the right is a pullback.
  Therefore also the top dashed morphism exists uniquely, and makes
  the square commute as shown, by functoriality of limits.

  The argument for colimits is analogous, now using that
  $G \times (-)$ also preserves colimit diagrams, by \eqref{MappingSpaceAdjunction}.
\end{proof}

\begin{proposition}[Compactly generated topological $G$-actions form a regular category]
  \label{CompactlyGeneratedTopologicalGActionsFormARegularCategory}
  For $G \,\in\, \Groups(\kTopologicalSpaces)$,
  the category $\Actions{G}(\kTopologicalSpaces)$
  \eqref{GActionOnTopologicalSpaces}
  is regular (Def. \ref{RegularCategory}).
\end{proposition}
\begin{proof}
  Since regularity is entirely a condition on limits and colimits of a category,
  it transfers through any forgetful functor which creates all limits and colimits.
  Therefore the statement follows by the combination of
  Lem. \ref{ForgetfulFunctorFromTopologicalGSpacesToGSpaces}
  with Prop. \ref{CompactyGeneratedTopologicalSpacesFormARegularCategory}.
\end{proof}

In generalization of Example \ref{ForgettingGActionsAsPullbackAction}, we have:

\begin{example}[Restricted actions]
  \label{RestrictedActions}
  Let $H \xhookrightarrow{\;} G$ any subgroup inclusion. Then
  the corresponding pullback $H$-action  (Lemma \ref{InducedAndCoinducedActions})
  of a $G \acts \, \mathrm{Y} \in \GActionsOnTopologicalSpaces$
  is just its restriction to the action of $H \subset G$
      \vspace{-2mm}
  \begin{equation}
  \label{InducedRestrictedActionAdjunction}
   \begin{tikzcd}[row sep=0pt]
     H\mathrm{Act}
     (
       \mathrm{TopSp}
     )
     \ar[
       rr,
       shift left=5pt,
       "G \times_H (-)"{above}
     ]
     \ar[
       rr,
       phantom,
       "\scalebox{.7}{$\bot$}"
     ]
     &&
     G\mathrm{Act}
     (
       \mathrm{TopSp}
     )
     \ar[
       ll,
       shift left=5pt,
       "\mbox{\tiny\color{greenii} \bf restricted action}"
     ]
     \\
     \scalebox{0.7}{$       H \acts \, \mathrm{Y} $}
     &\longmapsfrom&
     \scalebox{0.7}{$        G \acts \, \mathrm{Y} $}
     \mathrlap{\,,}
   \end{tikzcd}
  \end{equation}

  \vspace{-2mm}
\noindent
  and the adjunction \eqref{AdjointTripleOfChangeOfEquivarianceGroup}
  corresponds to a natural bijection \eqref{FormingAdjuncts} of hom-sets of the following form:
      \vspace{-2mm}
  \begin{equation}
    \label{HomIsomorphismForRestrictedActionAndInducedAction}
    \bigg\{
    \begin{tikzcd}[row sep=-3pt]
      \TopologicalSpace
      \ar[out=180-66, in=66, looseness=3.5, "\scalebox{.77}{$\mathclap{
        H
      }$}"{description},shift right=1]
      \ar[r]
      &
      \mathrm{Y}
      \ar[out=180-66, in=66, looseness=3.5, "\scalebox{.77}{$\mathclap{
        H
      }$}"{description},shift right=1]
      \\
  \scalebox{.7}{$    x $}
      \ar[
        r,
        phantom,
        "{\overset{}{\longmapsto}}"{description}
      ]
      &
   \scalebox{.7}{$   f(x)  $}
    \end{tikzcd}
    \bigg\}
    \qquad
    \mathrel{\mathop{
    \xleftrightarrow{
      \qquad \quad
      \widetilde{(-)}
      \qquad \quad
    }}_{\scalebox{0.6}{\bf \color{greenii} induction/restriction}}
    }
    \qquad
    \bigg\{\!\!\!
    \begin{tikzcd}[row sep=-3pt]
      G \times_H
      \TopologicalSpace
      \ar[out=180-60, in=60, looseness=3.0, "\scalebox{.77}{$\mathclap{
        G
      }$}"{description},shift right=1]
      \ar[r]
      &
      \mathrm{Y}
      \ar[out=180-66, in=66, looseness=3.5, "\scalebox{.77}{$\mathclap{
        G
      }$}"{description},shift right=1]
      \\
  \scalebox{.7}{$    {[g,x]} $}
      \ar[
        r,
        phantom,
        "{\overset{}{\longmapsto}}"{description}
      ]
      &
  \scalebox{.7}{$    g \cdot f(x)  $}
    \end{tikzcd}
   \!\!\! \bigg\}.
  \end{equation}
\end{example}

\begin{example}[Fixed loci with residual Weyl-group action]
  \label{FixedLociWithResidualWeylGroupAction}
  Let $H \subset G$ be a subgroup inclusion. Consider the functor which
  sends a $G$-action to its $H$-fixed subspace equipped with its
  residual Weyl group action (Ntn. \ref{GActionOnTopologicalSpaces})
  \vspace{-2mm}
  \begin{equation}
    \label{HFixedLocusFunctor}
    \hspace{-.6cm}
    \begin{tikzcd}[row sep=-4pt, column sep=large]
     N(H)
     \mathrm{Act}
     (
       \mathrm{TopSp}
     )
      \ar[
        rr,
        "{
         \scalebox{0.75}{$    (-)^H
          \;\coloneqq\;
          \mathrm{Maps}
          \left(
            N(H)/H, \,
            X
          \right)^{N(H)}
          $}
        }"{above}
      ]
      &&
      W(H)
      \mathrm{Act}
      (\mathrm{TopSp})
      \\
 \scalebox{0.7}{$        N(H) \acts \, \TopologicalSpace $}
      &\longmapsto&
\scalebox{0.7}{$         W(H)
      \; \acts\;
      \TopologicalSpace^H
        \,=:\,
        \big\{
          x \in \TopologicalSpace
          \; \big\vert \;\;
            \underset{
              h \in H \subset G
            }{\forall} \;
            h \cdot x = x
        \big\}.
        $}
    \end{tikzcd}
  \end{equation}

\vspace{-3mm}
\noindent
This fixed locus functor $(-)^H$ is, equivalently, the pull-push of change-of-equivariance-groups
(Lemma \ref{InducedAndCoinducedActions}) through the normalizer correspondence
    \vspace{-4mm}
$$
  \begin{tikzcd}[row sep=-4pt]
    &
    N\!(H)
    \ar[dl, hook]
    \ar[dr,->>]
    \\
    G
    &&
    W\!(H)\;,
  \end{tikzcd}
$$

\vspace{-2mm}
\noindent
in that it is the composite right adjoint
in the following composite of
change-of-equivariance group adjunction \eqref{AdjointTripleOfChangeOfEquivarianceGroup}:

\vspace{-.5cm}
\begin{equation}
\label{RightAdjointWeylGroupValuedFixedLocusFunctor}
\hspace{-5mm}
\begin{tikzcd}[column sep=7em]
  G\mathrm{Act}
  \ar[
    rr,
    shift right=7pt,
    "(N(H) \hookrightarrow G)^\ast"{below}
  ]
  \ar[
    rr,
    phantom,
    "\scalebox{.7}{$\bot$}"
  ]
  \ar[
    rrrr,
    rounded corners,
    to path={
         -- ([yshift=-15pt]\tikztostart.south)
         --node[below]{\scalebox{.7}{$(-)^H$}} ([yshift=-12pt]\tikztotarget.south)
         -- (\tikztotarget.south)}
  ]
  &&
  N\!(H)\mathrm{Act}
  \ar[
    ll,
    shift right=7pt,
    "G \times_{N\!(H)} (-)"{above}
  ]
  \ar[
    rr,
    shift right=7pt,
    "{
      \mathrm{Maps}
      (
        N\!(H)\!/H,
        -
      )^{N\!(H)}
    }"{below}
  ]
  \ar[
    rr,
    phantom,
    "\scalebox{.7}{$\bot$}"
  ]
  &&
  \big(
    N\!(H)\!/H
  \big)\mathrm{Act}
  \mathrlap{\,.}
  \ar[
    ll,
    shift right=7pt,
    "{
      \left( N\!(H)\, \twoheadrightarrow N\!(H)\!/H \right)^\ast
    }"{above}
  ]
  \ar[
    llll,
    rounded corners,
    to path={
         -- ([yshift=+13pt]\tikztostart.north)
         --node[above]{\scalebox{.7}{$G/H \times_{N\!(H)\!/H}(-)$}} ([yshift=+15pt]\tikztotarget.north)
         -- (\tikztotarget.north)}
  ]
\end{tikzcd}
\end{equation}
\vspace{-.4cm}
\end{example}

\begin{example}[Coset spaces (e.g. {\cite[p. 34]{Bredon72}})]
  \label{CosetSpacesAsActions}
  For $H \hookrightarrow G$ a subgroup inclusion,
  the induced $G$-action \eqref{InducedRestrictedActionAdjunction}
  of the unique trivial $H$-action on the point
  is the coset space
      \vspace{-2mm}
  \begin{equation}
    \label{CosetSpace}
    G \times_H \ast
      \;=\;
    G/H
      \;\coloneqq\;
    \{ g H \subset G  \,\vert\, g \in G \}
      \;\in\;
    \GActionsOnTopologicalSpaces
  \end{equation}

  \vspace{-2mm}
  \noindent
  equipped with its $G$-action by left multiplication of representatives in $G$.
\end{example}

\medskip

\noindent
{\bf Quotient spaces.} In generalization of Example \ref{QuotientAndFixedLociFromChangeOfGroupAdjunction} we have:
\begin{example}[Partial quotient spaces]
  \label{QuotientSpaces}
  For
  $G, G' \in \mathrm{Grps}(\mathrm{HausSp})$
  and $G \times G' \xrightarrow{\mathrm{pr}_2} G'$ the projection
  homomorphism out of their direct product,
  the corresponding pullback action in Lemma \ref{InducedAndCoinducedActions}
  assigns trivial $G$-actions
  \vspace{-2mm}
  $$
    \begin{tikzcd}[column sep=large]
      (
        G \times G'
      )
      \mathrm{Act}
      (
        \mathrm{TopSp}
      )
      \ar[
        rr,
        shift left=5pt,
        "(-)/G"
      ]
      \ar[
        rr,
        phantom,
      "\scalebox{.7}{$\bot$}"
      ]
      &&
      G'\mathrm{Act}
      (
        \mathrm{TopSp}
      )
      \ar[
        ll,
        shift left=5pt,
        "\mbox{\tiny\color{greenii}\bf trivial $G$-action}"
      ]
    \end{tikzcd}
  $$

   \vspace{0mm}
\noindent
  and its left adjoint \eqref{AdjointTripleOfChangeOfEquivarianceGroup}
  forms
  $G$-quotients
  $
    G' \times_{G \times G'} (-)
    \;=\;
    (-)/G
    \,.
  $
  The unit of this adjunction is the natural transformation which sends
  any $G \times G'$-action to the coprojection
  $q_{\TopologicalSpace} : \TopologicalSpace \xrightarrow{\;} \TopologicalSpace/G$
  onto its $G$-quotient space, and
  any $G$-equivariant continuous function $f : \TopologicalSpace \xrightarrow{\;} \mathrm{Y}$
  to a commuting square of $G'$-$\mathrm{actions}$:
   \vspace{-2mm}
  \begin{equation}
    \label{QuotientSpaceNaturalitySquare}
    \begin{tikzcd}[column sep=large, row sep=small]
      \TopologicalSpace
      \ar[out=180-66, in=66, looseness=3.5, "\scalebox{.77}{$\mathclap{
        G'
      }$}"{description},shift right=1]
      \ar[
        r,
        "f"{above}
      ]
      \ar[
        d,
        "q_{\TopologicalSpace}"{left}
      ]
      &
      \mathrm{Y}
      \ar[out=180-66, in=66, looseness=3.5, "\scalebox{.77}{$\mathclap{
        G'
      }$}"{description},shift right=1]
      \ar[
        d,
        "q_{\mathrm{Y}}"
      ]
      \\
      \TopologicalSpace/G
      \ar[out=-180+66, in=-66, looseness=3.6, "\scalebox{.77}{$\mathclap{
        G'
      }$}"{description},shift left=1]
      \ar[
        r,
        "f/G"{below}
      ]
      &
      \mathrm{Y}/G
      \ar[out=-180+66, in=-66, looseness=3.6, "\scalebox{.77}{$\mathclap{
        G'
      }$}"{description},shift left=1]
    \end{tikzcd}
  \end{equation}
\end{example}
\begin{lemma}[Hausdorff quotient spaces (e.g. {\cite[Thm. 3.1]{Bredon72}})]
  \label{HausdorffQuotientSpaces}
  If the equivariance group $G$ is compact
  and the topological space underlying
  $\TopologicalSpace \,\in\, \GActionsOnTopologicalSpaces$ (Ntn. \ref{GActionOnTopologicalSpaces})
  is Hausdorff, then the quotient space $X/G$ is Hausdorff.
\end{lemma}

\begin{remark}[Recognition of pullbacks of quotient coprojections]
{\bf (i)} The quotient coprojection squares \eqref{QuotientSpaceNaturalitySquare}
are not in general pullbacks (Ntn. \ref{CartesianSquares}); and it is important
to recognize those situations
in which they are.

\vspace{-1mm}
\item {\bf (ii)}
A general recognition principle applies to
{\it compact} quotient groups (Lemma \ref{RecognitionOfCartesianQuotientProjections} below)
which is however of little value in the applications
to twisted cohomology theory, where the quotient groups
generically are topological representatives of general $\infty$-groups
(topological realizations of general simplicial groups) and thus rarely compact.

\vspace{-1mm}
\item {\bf (iii)}
Without assuming compactness of the quotient group we
may still recognize pullbacks of {\it free} (principal) quotients
(Lemma \ref{HomomorphismsOfLocallyTrivialPrincipalBundlesArePullbackSquares})
from just the fact that the left one is locally trivial below
(compare Ntn. \ref{TerminologyForPrincipalBundles}).
This is a basic fact of principal bundle theory, but rarely, if
ever, stated in the general form of Lemma
 \ref{HomomorphismsOfLocallyTrivialPrincipalBundlesArePullbackSquares}
in which it drives much of the proofs in \cref{NotionsOfEquivariantLocalTrivialization}.
\end{remark}

\begin{lemma}[Recognition of pullbacks of compact group quotients ({\cite[Prop. 4.1]{BykovFlores15}})]
  \label{RecognitionOfCartesianQuotientProjections}
  If $G$ is compact and the underlying topological spaces of
  $\TopologicalSpace, \mathrm{Y} \,\in\, \GActionsOnTopologicalSpaces$
  are Hausdorff
  then, for any morphism $f : \TopologicalSpace \xrightarrow{\;} \mathrm{Y}$,
  its quotient naturality square
  \eqref{QuotientSpaceNaturalitySquare}
  is a pullback square (Ntn. \ref{CartesianSquares})
  if and only if $f$ preserves isotropy groups (as subgroups of $G$):
  \vspace{-1mm}
  $$
    \begin{tikzcd}[column sep=large]
      \TopologicalSpace
      \ar[
        r,
        "f"{above}
      ]
      \ar[
        d,
        "q_{\TopologicalSpace}"{left}
      ]
      \ar[
        dr,
        phantom,
        "\mbox{\tiny\rm(pb)}"
      ]
      &
      \mathrm{Y}
      \ar[
        d,
        "q_{\mathrm{Y}}"
      ]
      \\
      \TopologicalSpace/G
      \ar[
        r,
        "f/G"{below}
      ]
      &
      \mathrm{Y}/G
    \end{tikzcd}
    {\phantom{AAA}}
    \Leftrightarrow
    {\phantom{AAA}}
    \underset{x \in \TopologicalSpace}{\forall}
    \left(
      G_x \simeq G_{f(x)}
    \right)
    \,.
  $$
\end{lemma}

\begin{lemma}[Recognition of pullbacks of principal quotients]
  \label{HomomorphismsOfLocallyTrivialPrincipalBundlesArePullbackSquares}
  A homomorphism of $G$-principal fibrations
  (covering any morphism of base spaces)
  is a pullback square
  (Ntn. \ref{CartesianSquares})
  as soon as the domain is a locally trivializable fiber bundle:

  \vspace{-7mm}
  \begin{equation}
    \label{MorphismOfPrincipalBundlesIsPullback}
    {
    \begin{tikzcd}[column sep=large, row sep=20pt]
      \\
      \mathrm{P}_1
      \ar[
        rr,
        "f"
      ]
      \ar[
        d,
        "p_1"{left}
      ]
      \ar[out=180-66, in=66, looseness=3.5, "\scalebox{.77}{$\mathclap{
        G
      }$}"{description}, shift right=1]
      &&
      \mathrm{P}_2
      \ar[out=180-66, in=66, looseness=3.5, "\scalebox{.77}{$\mathclap{
        G
      }$}"{description}, shift right=1]
      \ar[
        d,
        "\,p_2"
      ]
      \\
      \underset{
        \mathclap{
        \raisebox{-4pt}{
          \tiny
          \color{darkblue}
          \bf
          \def\arraystretch{1}
          \begin{tabular}{c}
            $G$-principal \&
            \\
            locally trivial
          \end{tabular}
        }
        }
      }{
        \TopologicalSpace_1
      }
      \ar[
        rr,
        "f/G"
      ]
      &&
      \underset{
        \mathclap{
        \raisebox{-4pt}{
          \tiny
          \color{darkblue}
          \bf
          \def\arraystretch{1}
          \begin{tabular}{c}
            $G$-principal
          \end{tabular}
        }
        }
      }{
        \TopologicalSpace_2
      }
    \end{tikzcd}
    }
    {\phantom{AAAA}}
    \Rightarrow
    {\phantom{AAAA}}
    \begin{tikzcd}[column sep=large]
      \mathrm{P}_1
      \ar[
        rr
      ]
      \ar[
        d,
        "p_1"{left}
      ]
      \ar[
        drr,
        phantom,
        "\mbox{\tiny\rm(pb)}"{description}
      ]
      &&
      \mathrm{P}_2
      \ar[
        d,
        "\,p_2"
      ]
      \\
      \TopologicalSpace_1
      \ar[rr]
      &&
      \TopologicalSpace_2
    \end{tikzcd}
  \end{equation}
\end{lemma}
\begin{proof}
  First, consider the special case when the domain bundle is
  actually trivial and that the base morphism is the identity
  \vspace{-2mm}
  $$
    \begin{tikzcd}[column sep=large, row sep=small]
      G \times \TopologicalSpace
      \ar[
        d,
        "\; \mathrm{pr}_2"{left}
      ]
      \ar[
        r,
        "\sigma"
      ]
      &
      \mathrm{P}
      \ar[
        d,
        "\; p"
      ]
      \\
      \TopologicalSpace
      \ar[r,-,shift left=1pt]
      \ar[r,-,shift right=1pt]
      &
      \TopologicalSpace
    \end{tikzcd}
  $$
  (equivalently a global section $\sigma(\NeutralElement)(-)$).
  In that case, a continuous inverse of $\sigma$ is given by the composite
      \vspace{-2mm}
  $$
    \begin{tikzcd}
      G \times \TopologicalSpace
      &
      G \times \mathrm{P}
      \ar[
        l,
        "{\mathrm{id} \times p}"{above}
      ]
      &
      \mathrm{P} \times_{\TopologicalSpace} \mathrm{P}
      \ar[
        l,
        "{\sim}"{above}
      ]
      &&&
      \mathrm{P}
      \mathrm{\,,}
      \ar[
       lll,
        "{
          \left(
            \sigma(\NeutralElement,\,p(-))
            ,\,
            \mathrm{id}
          \right)
        }"{above}
      ]
    \end{tikzcd}
  $$

  \vspace{-1mm}
  \noindent
  where the second map is the inverse of the shear map (see
  \eqref{PrincipalityConditionAsShearMapBeingAnIsomorphism}) of the codomain bundle,
  and where $\NeutralElement \in G$ denotes the neutral element.
    From this it follows
  that morphisms out of any locally trivial principal bundle
  over the identity base morphism are isomorphisms,
  by local recognition of homeomorphisms
  (Ex. \ref{IsomorphismOfBundlesDetectedOnOpenCovers}).
  Finally, in the general case the universal comparison morphism
  \vspace{-3mm}
  $$
    \begin{tikzcd}[row sep=small]
      \mathrm{P}_1
      \ar[
        r,
        dashed,
        "\sim"
      ]
      \ar[dr]
      &
      \mathrm{P}_2 \times_{{}_{\TopologicalSpace_1}} \mathrm{X_2}
      \ar[rr]
      \ar[d]
      \ar[
        drr,
        phantom,
        "\mbox{\tiny\rm(pb)}"{description}
      ]
      &&
      \mathrm{P}_2
      \ar[d]
      \\
      &
      \TopologicalSpace_1
      \ar[rr]
      &&
      \TopologicalSpace_2
    \end{tikzcd}
  $$

  \vspace{-2mm}
\noindent  from $\mathrm{P}_1$ to the pullback of $\mathrm{P}_2$ in
\eqref{MorphismOfPrincipalBundlesIsPullback}
  is such a homomorphism over a common
  base space $\TopologicalSpace_1$, hence is an isomorphism, thus exhibiting
  $\mathrm{P}_1$ as a pullback.
\end{proof}

\medskip

\noindent
{\bf Slices of $G$-orbits.} The existence of local
{\it slices through families of orbits} of group actions
is guaranteed by the Slice Theorem (Prop. \ref{SliceTheorem})
below and serves to ensure or detect local triviality of
plain principal bundles (Cor. \ref{QuotientCoprojectionOfFreeProperActionIsLocallyTrivial}) below
and of equivariant principal bundles (around Ntn. \ref{SlicesInsideBierstonePatches} below).

\begin{definition}[Slices of $G$-orbits]
  \label{SliceOfTopologicalGSpace}
  For $G \acts \; \mathrm{U} \in \GActionsOnTopologicalSpaces$
  and $H \subset G$ a subgroup, an $H$-subspace
  \vspace{-2mm}
  \begin{equation}
    \label{SliceSubspace}
    \begin{tikzcd}
      \mathrm{S}
      \ar[out=180-66, in=66, looseness=3.5, "\scalebox{.77}{$\phantom{}\mathclap{
        H
      }\phantom{}$}"{description}, shift right=1]
      \ar[
        r,
        hook,
        "\iota"
      ]
      &
      \mathrm{U}
      \ar[out=180-66, in=66, looseness=3.5, "\scalebox{.77}{$\phantom{}\mathclap{
        H
      }\phantom{}$}"{description}, shift right=1]
    \end{tikzcd}
  \end{equation}

  \vspace{-2mm}
  \noindent
  is called a {\it slice} through its $G$-orbit modulo $H$
  if its
  induction/restriction-adjunct \eqref{HomIsomorphismForRestrictedActionAndInducedAction}
  is an isomorphism
  \vspace{-3mm}
  \begin{equation}
    \label{SliceIsomorphism}
    \begin{tikzcd}[row sep=-2pt]
      G \times_H \mathrm{S}
      \ar[out=180-66, in=66, looseness=3.5, "\scalebox{.77}{$\phantom{}\mathclap{
        G
      }\phantom{}$}"{description}, shift right=1]
      \ar[
        r,
        "\sim"{below},
        "\tilde \iota"{above}
      ]
      &
      \mathrm{U}
      \ar[out=180-66, in=66, looseness=3.5, "\scalebox{.77}{$\phantom{}\mathclap{
        G
      }\phantom{}$}"{description}, shift right=1]
      \\
   \scalebox{.7}{$   {[g , s]} $}
   \;\;   \ar[
        r,
        phantom,
        "\longmapsto"{description}
      ]
      &
  \scalebox{.7}{$    g \cdot s $}
      \,.
    \end{tikzcd}
  \end{equation}

  \vspace{-2mm}
  \noindent
  Specifically, for
  $G \acts \,  \TopologicalSpace \in \GActionsOnTopologicalSpaces$
  and $x \in \TopologicalSpace$ a point,
  by a {\it slice through $x$} one means
  (e.g. {\cite[\S II, Def. 4.1]{Bredon72}})
  a slice \eqref{SliceSubspace}
  relative to the isotropy group $G_x$ \eqref{StabilizerSubgroupInEquivarianceGroup}
  through $x$ of an open $G$-neighborhood $G \acts \; \mathrm{U}_x$ of $x$
  \vspace{-2mm}
  \begin{equation}
    \label{ASliceThroughAPoint}
    \begin{tikzcd}[column sep=-1]
      x
      &\in
      &
      \mathrm{S}_x
      \ar[out=180-66, in=66, looseness=3.5, "\scalebox{.77}{$\phantom{}\mathclap{
        G_x
      }\phantom{\cdot}$}"{description}, shift right=1]
      \ar[
        rr,
        hook
      ]
      &{\phantom{AA}}&
      \mathrm{U}_x
      \mathrlap{\,.}
      \ar[out=180-66, in=66, looseness=3.5, "\scalebox{.77}{$\phantom{}\mathclap{
        G_x
      }\phantom{\cdot}$}"{description}, shift right=1]
    \end{tikzcd}
  \end{equation}
\end{definition}
\begin{proposition}[Slice Theorem ({\cite[Prop. 2.3.1]{Palais61}\cite[Thm. 6.2.7]{Karppinen16}})]
  \label{SliceTheorem}
  Under the assumption
  \ref{ProperEquivariantTopology}
  of proper equivariance,
  given $G \acts \, \TopologicalSpace \in \GActionsOnTopologicalSpaces$
  then for every $x \in \TopologicalSpace$
  there exists a slice through $x$ (Def. \ref{SliceOfTopologicalGSpace}).
\end{proposition}
\begin{remark}[Technical conditions in the slice theorem]
\label{TechnicalConditionsInTheSliceTheorem}
The slice theorem
for {\it compact} Lie group actions
is due to \cite[Thm. 2.1]{Mostow57}\cite[Cor. 1.7.19]{Palais60},
and for proper actions of general Lie groups it is due to \cite{Palais61},
reviewed in \cite{Karppinen16}.
Beware that \cite{Palais61} goes to some length
to further generalize beyond proper actions, which leads to a wealth of technical
conditions that, it seems, have been of rare use in practice.
But, under the assumption \ref{ProperEquivariantTopology}
that all $G$-spaces are locally compact, all these conditions reduce to
properness \cite[Thm 1.2.9]{Palais61}\cite[Rem. 5.2.4]{Karppinen16},
and thus the theorem reduces to the statement in Prop. \ref{SliceTheorem}.
\end{remark}
\begin{corollary}[Quotient coprojection of free proper action is locally trivial
{\cite[\S 4.1]{Palais61}}]
  \label{QuotientCoprojectionOfFreeProperActionIsLocallyTrivial}
  Under the assumption \ref{ProperEquivariantTopology} of proper equivariance,
  the quotient space coprojection
  $P \xrightarrow{q} P/G$
  of a \emph{free} action $G \acts \,  P$ admits local sections.
\end{corollary}

\medskip

\noindent
{\bf Equivariant open covers.}

\begin{definition}[Properly equivariant open cover]
  \label{ProperEquivariantOpenCover}
  Given $G \acts \, \TopologicalSpace \,\in\, \Actions{G}(\kTopologicalSpaces)$
  \eqref{GActionOnTopologicalSpaces},
  we say that an open cover
  (Ex. \ref{OpenCoversAreEffectiveEpimorphisms}) of the underlying space
  \vspace{-1mm}
  \begin{equation}
    \label{OpenCoverToBeEquivariant}
    \widehat{\TopologicalSpace}
    \,=\,
    \underset{i \in I}{\sqcup}
    \TopologicalPatch_i
    \twoheadrightarrow
    \TopologicalSpace
  \end{equation}

  \vspace{-1mm}
  \noindent
  is

  \noindent
  {\bf (i)}
  {\it equivariant} if the $G$-action on $\TopologicalSpace$ pulls back to
  $\widehat{\TopologicalSpace}$
  \vspace{-2mm}
  \begin{equation}
    \label{AProperEquivariantOpenCover}
    \begin{tikzcd}
      G \acts \; \widehat{\TopologicalSpace}
      \ar[r, ->>, "p"]
      &
      G \acts \, \TopologicalSpace
    \end{tikzcd}
    \;\;\;
    \in
    \;
    \Actions{G}(\kTopologicalSpaces) \;;
  \end{equation}

  \vspace{-2mm}
  \noindent
  {\bf (ii)}
  {\it regular}
  if there is a $G$-action on the index set such that
  \begin{equation}
    \label{RegularityConditionsOnEquivariantOpenCover}
    \begin{aligned}
   {\bf  (a)}
    \quad
    &
    \underset{i,j \in I}{\forall}
    \Bigg(
      \begin{tikzcd}[row sep=7pt]
        \TopologicalPatch_i
        \ar[r, "\sim"{swap}]
        \ar[d, hook]
        &
        \TopologicalPatch_{g \cdot j}
        \ar[d, hook]
        \\
        \TopologicalSpace
        \ar[r, -, shift left=1pt]
        \ar[r, -, shift right=1pt]
        &
        \TopologicalSpace
      \end{tikzcd}
      \;\;\;\;
      \Rightarrow
      \;\;\;\;
      i \,=\, j
    \Bigg);
    \\
    {\bf (b)}
    \quad
    &
       \underset{i \in I}{\forall}
       \,\,
       \underset{g \in G}{\forall}
       \;
       \bigg(
         U_i
          \cap
         g \cdot U_j
         \,\neq\,
         \varnothing
         \;\;\Rightarrow\;\;
        \begin{tikzcd}[row sep=7pt]
          \TopologicalPatch_i
          \ar[r, "\sim"{swap}]
          \ar[d, hook]
          &
          \TopologicalPatch_{i}
          \ar[d, hook]
          \\
          \TopologicalSpace
          \ar[r, "g", "\sim"{swap}]
          &
          \TopologicalSpace
        \end{tikzcd}
       \bigg);
       \\
      {\bf  (c)}
       \quad
       &
       \underset{n \in \mathbb{N}}{\forall}
       \;
       \underset{
          {i_0, \cdots, i_n \in I}
          \atop
          {g_0, \cdots, g_n \,\in\, G }
       }{\forall}
       \left(
         \begin{array}{rcl}
           U_{i_0} \cap \cdots \cap U_{i_n}
           & \neq
           & \varnothing,
           \\
           g_0 \cdot U_{i_0} \cap \cdots g_n \cdot U_{i_n}
           & \neq
           & \varnothing
         \end{array}
         \;\;\Rightarrow\;\;
         \underset{g \in G}{\exists}
         \;
         \underset{0 \leq k \leq n}{\forall}
         \;
         g \cdot U_{i_k} \,=\, g_k \cdot U_{i_k}
       \right);
    \end{aligned}
  \end{equation}

  \noindent
  {\bf (iii)}
  {\it properly equivariant}
  if, in addition, each $H$-fixed locus of $\widehat{\TopologicalSpace}$
  is an open cover of that of $\TopologicalSpace$:
  \vspace{-2mm}
  \begin{equation}
    \label{EquivariantOpenCoverRestrictedToFixedLoci}
    \underset{
      H
        \underset{\mathclap{\mathrm{clsd}}}{\subset}
      G
    }{\forall}
    \;\;\;
    \begin{tikzcd}
      \widehat{\TopologicalSpace}^H
      \ar[rr, ->>, "p^H", "\mbox{\tiny \rm open cover}"{swap}]
      &&
      \TopologicalSpace^H
      ;
    \end{tikzcd}
  \end{equation}

  \vspace{-2mm}
  \noindent
  {\bf (iv)}
  properly equivariantly {\it good}
  if all the restrictions
  \eqref{EquivariantOpenCoverRestrictedToFixedLoci}
  are good open covers
  (Def. \ref{GoodOpenCovers}).
\end{definition}
\begin{proposition}[Smooth $G$-manifolds admit properly equivariant regular good open covers]
  \label{SmoothGManifoldsAdmitProperlyEquivariantGoodOpenCovers}
  At least for

 \noindent  $G \,\in\, \Groups(\FiniteSets) \xhookrightarrow{\Groups(\Discrete)}
  \Groups(\kTopologicalSpaces)$,
  every smooth $G$-manifold
  $
    G \acts \, \TopologicalSpace
      \,\in\,
    \Groups(\SmoothManifolds)
    \xhookrightarrow{\;}
    \Groups(\kTopologicalSpaces)
  $
  admits a regular properly equivariant and good open cover
  (Def. \ref{ProperEquivariantOpenCover}).
\end{proposition}
\begin{proof}
  This follows with the
  equivariant triangulation theorem
  \cite[Thm. 3.1]{Illman72}\cite{Illman83}; see \cite[Thm. 2.11]{Yang14}.
\end{proof}

\begin{remark}
  \label{StabilizerSubgroupOfPointAlsoFixedIndexOfPatchInRegularEquivariantOpenCover}
  If an equivariant open cover \eqref{EquivariantOpenCover}
  is regular \eqref{RegularityConditionsOnEquivariantOpenCover}
  then for $x \,\in\, \TopologicalPatch_i \xhookrightarrow{\;} \TopologicalSpace$
  the
  stabilizer group $G_x \,\coloneqq\, \mathrm{Stab}_G(x)$
  of $x \,\in\, X$ also fixes the index $i$:
  \vspace{-2mm}
  $$
    \begin{tikzcd}[row sep=small, column sep=large]
      \mathllap{x \in \;}
      \TopologicalPatch_i
      \ar[r, "\sim"{swap}]
      \ar[d, hook]
      &
      \TopologicalPatch_i
      \ar[d, hook]
      \\
      \TopologicalSpace
      \ar[r, "g \,\in\, G_x"]
      &
      \TopologicalSpace
      \,.
    \end{tikzcd}
  $$
\end{remark}

\section{$G$-Actions on topological groupoids}
\label{GActionsOnTopologicalGroupoids}

The theory of {\it universal} equivariant bundles turns out to
be most naturally formulated
in the language not just of topological spaces equipped with
$G$-actions, but of topological {\it groupoids} equipped with $G$-actions.
This observation, which is due to \cite{MurayamaShimakawa95}
and was only
more recently amplified again in \cite{GuillouMayMerling17},
we turn  in \cref{ConstructionOfUniversalEquivariantPrincipalBundles} and then especially in \cref{InCohesiveInfinityStacks} below.
Here we recall and develop some basics of equivariant topological groupoids
that will make the theory of universal equivariant bundles
in \cref{ConstructionOfUniversalEquivariantPrincipalBundles}
be transparent and run smoothly.
(The material here is not needed for \cref{EquivariantPrincipalTopologicalBundles}.)

\medskip
\noindent
{\bf Topological groupoids.}

\begin{notation}[Topological groupoids]
  \label{TopologicalGroupoids}
  We write
  $\TopologicalGroupoids$ for the
  strict (2,1)-category (Ntn. \ref{Strict2Categories})
  of groupoid objects internal (Ntn. \ref{Internalization})
  to
  $\kTopologicalSpaces$ \eqref{CategoryOfTopologicalSpaces},
  hence of {\it topological groupoids}
  (\cite{Ehresmann59}, survey in \cite[\S II]{Mackenzie87},
  exposition in \cite[p. 6]{Weinstein96}):

  \noindent
  {\bf (i)} Its objects are diagrams of topological spaces

  \vspace{-4mm}
  \begin{equation}
    \label{DiagramForTopologicalGroupoid}
    \begin{tikzcd}[row sep=1pt, column sep=30pt]
      \mathclap{
      \mbox{
        \tiny
        \color{darkblue}
        \bf
        \begin{tabular}{c}
          space of composable
          \\
          pairs of morphisms
        \end{tabular}
      }
      }
      &[-16pt]
      \mathclap{
      \mbox{
        \tiny
        \color{greenii}
        \bf
        \begin{tabular}{c}
          composition
          \\
          map
        \end{tabular}
      }
      }
      &
      \quad
          \mathclap{
      \mbox{
        \tiny
        \color{darkblue}
        \bf
        \begin{tabular}{c}
          space of
          \\
          morphisms
        \end{tabular}
      }
      }
           \quad
     &
      \mathclap{
      \mbox{
        \tiny
        \color{greenii}
        \bf
        \begin{tabular}{c}
          source, target \& unit
          \\
          maps
        \end{tabular}
      }
      }
      &
      \quad
      \mathclap{
      \mbox{
        \tiny
        \color{darkblue}
        \bf
        \begin{tabular}{c}
          space of
          \\
          objects
        \end{tabular}
      }
      }
      \\
      (\TopologicalSpace_1)
      \,
      {}_{t}\!\!\underset{\TopologicalSpace_0}{\times_s}
      (\TopologicalSpace_1)
      \ar[
        rr,
        "\circ"{pos=.4}
      ]
      &&
      \quad
      \TopologicalSpace_1
      \ar[out=180-60+180, in=60+180, looseness=3.0, "\scalebox{.77}{$\mathclap{
        \mathllap{
          \mbox{
            \tiny
            \color{greenii}
            \bf
            \begin{tabular}{c}
              inversion
              \\
              map
            \end{tabular}
          }
          \!\!\!
        }
        (-)^{-1}
      }$}"{below}]
      \ar[
        rr,
        shift left=5pt,
        "{s}"{above}
      ]
      \ar[
        rr,
        shift right=5pt,
        "{t}"{below}
      ]
      \quad
      &&
      \quad
        \TopologicalSpace_0
      \ar[
        ll,
        "{\mathrm{e}}"{description}
      ]
    \end{tikzcd}
  \end{equation}

  \vspace{-2mm}
  \noindent
  such that the composition operation $\circ$ is associative, unital
  with respect to $\mathrm{e}$ and with inverses given by $(-)^{-1}$.
  We will mostly denote such an object by

  \vspace{-.4cm}
  $$
    \TopologicalSpace_1
    \rightrightarrows
    \TopologicalSpace_0
    \;\;\;
    \in
    \;
    \TopologicalGroupoids
    \,,
  $$

  \vspace{-1mm}
  \noindent
  with the rest of the structure understood from the given context.

  \vspace{0mm}
  \noindent
  {\bf (ii)} Its morphisms are
  {\it continuous functors}
  hence continuous functions $F_0$, $F_1$
  compatible with all this structure:
      \vspace{-2mm}
  \begin{equation}
    \label{ContinuousFunctors}
    \begin{tikzcd}[row sep=14pt, column sep=30pt]
      (\TopologicalSpace_1)
      \,
      {}_{t}\!\!\underset{\TopologicalSpace_0}{\times_s}
      (\TopologicalSpace_1)
      \;\;
      \ar[
        r,
        "\circ"
      ]
      \ar[
        d,
        "\scalebox{0.8}{$
          F_1 \, \underset{F_0}{{}_t\!\!\times_s} \,  F_1
        $}"{left}
      ]
      &
      \quad
      \TopologicalSpace_1
      \ar[out=180-60, in=60, looseness=3.0, "\scalebox{.77}{$\mathclap{
        (-)^{-1}
      }$}"{above}]
      \quad
      \ar[
        rr,
        shift left=5pt,
        "{s}"{above}
      ]
      \ar[
        rr,
        shift right=5pt,
        "{t}"{below}
      ]
      \ar[
        d,
        "{F_1}"
      ]
      &&
      \;\;
      \TopologicalSpace_0
      \ar[
        ll,
        "{\mathrm{e}}"{description}
      ]
      \ar[
        d,
        "F_0"
      ]
      \\
      (\mathrm{Y}_1)
      \,
      {}_{t}\!\!\underset{\TopologicalSpace_0}{\times_s}
      (\mathrm{Y}_1)
      \;\;
      \ar[
        r,
        "\circ"
      ]
      &
      \quad
      \mathrm{Y}_1
      \ar[in=180-60+180, out=60+180, looseness=3.0, "\scalebox{.77}{$\mathclap{
        (-)^{-1}
      }$}"{below}]
      \quad
      \ar[
        rr,
        shift left=5pt,
        "{s}"{above}
      ]
      \ar[
        rr,
        shift right=5pt,
        "{t}"{below}
      ]
           &&
           \;\;
      \mathrm{Y}_0
      \ar[
        ll,
        "{\mathrm{e}}"{description}
      ]
    \end{tikzcd}
  \end{equation}

  \vspace{-3mm}
  \noindent
  {\bf (iii)}
  Its {\it 2-morphism} $\eta \,\colon\, F \Rightarrow F'$
  are {\it continuous natural transformations}, hence
  continuous functions
  $\eta(-) \,\colon\, \TopologicalSpace_0 \xrightarrow{\;} \mathrm{Y}_1$
  making all the naturality squares commute:
      \vspace{-2mm}
  \begin{equation}
    \label{NaturalitySquareForTopologicalGroupoids}
    \begin{tikzcd}
      (
        \TopologicalSpace_1
        \rightrightarrows
        \TopologicalSpace_0
      )
      \ar[
        rr,
        bend left=30,
        "{F}"{above},
        " "{below,name=s}
      ]
      \ar[
        rr,
        bend right=30,
        "{F'}"{below},
        " "{above,name=t}
      ]
      &&
      (
        \mathrm{Y}_1
        \rightrightarrows
        \mathrm{Y}_0
      )
      \ar[
        from=s,
        to=t,
        Rightarrow,
        "\eta"{xshift=1pt}
      ]
    \end{tikzcd}
    {\phantom{AA}}
    \colon
    {\phantom{AA}}
    \begin{tikzcd}[row sep=16pt]
      x
      \ar[
        dd,
        "{\gamma}"
      ]
      &[-15pt]&[-15pt]
      F(x)
      \ar[
        rr,
        "{ \eta(x) }"
      ]
      \ar[
        dd,
        "F(\gamma)"
      ]
      &&
      F'(x)
      \ar[
        dd,
        "F'(\gamma)"
      ]
      \\[-16pt]
      &
        \overset{
          \mbox{
            \tiny
            \rm
            cts.
          }
        }{\longmapsto}
      &
      \\[-16pt]
      x'
      &&
      F(x')
      \ar[
        rr,
        "\eta(x')"
      ]
      &&
      F'(x')
      \mathrlap{\,.}
    \end{tikzcd}
  \end{equation}

  \vspace{-2mm}
  \noindent
  {\bf (iv)}
  An {\it isomophism of topological groupoids}
  $(\TopologicalSpace_1 \rightrightarrows \TopologicalSpace_0) \;\simeq\; (\mathrm{Y}_1 \rightrightarrows \mathrm{Y}_0)$
  is an isomorphsim in the underlying 1-category (ignoring the 2-morphisms).

  \noindent
  {\bf (v)}
  An {\it equivalence of topological groupoids} is a pair of
  morphisms going back and forth between them, together with
  2-morphisms \eqref{NaturalitySquareForTopologicalGroupoids}
  relating their composites to the identity morphism:
      \vspace{-2mm}
  \begin{equation}
    \label{EquivalenceOfTopologicalGroupoids}
    \begin{tikzcd}
    (\TopologicalSpace_1 \rightrightarrows \TopologicalSpace_0)
    \,\underset{\mathrm{hmtpy}}{\simeq}\,
    (\mathrm{Y}_1 \rightrightarrows \mathrm{Y}_0)
    \qquad
      \Leftrightarrow
    \qquad
    (
      \TopologicalSpace_1
        \rightrightarrows
      \TopologicalSpace_0
    )
    \ar[
      r,
      shift right=2pt,
      "R"{below}
    ]
    &
    (
      \mathrm{Y}_1
        \rightrightarrows
      \mathrm{Y}_0
    )
    \mathrlap{\,,}
    \ar[
      l,
      shift right=2pt,
      "L"{above}
    ]
    \end{tikzcd}
    \;\;\;
    L \circ R \Rightarrow \mathrm{id}\,,
    \;\;\;\;\;
    \mathrm{id} \Rightarrow R \circ L
    \,.
  \end{equation}
\end{notation}

\begin{example}[Topological spaces as constant topological groupoids]
  \label{TopologicalSpacesAsTopologicalGroupoids}
  Each $\TopologicalSpace \,\in\, \kTopologicalSpaces$ \eqref{CategoryOfTopologicalSpaces}
  becomes a topological groupoid (Ntn. \ref{TopologicalGroupoids})
  \vspace{-3mm}
  $$
    \ConstantGroupoid(\mathrm{C})
    \;\coloneqq\;
    \big(
    \TopologicalSpace
      \underoverset
        {\mathrm{id}}
        {\mathrm{id}}
        {\rightrightarrows}
    \TopologicalSpace
    \big)
    \;\;\;
    \in
    \;
    \TopologicalGroupoids
  $$

        \vspace{-2mm}
\noindent  by taking all structure maps \eqref{DiagramForTopologicalGroupoid}
  to be the identity on $\TopologicalSpace$. This construction
  constitutes a full subcategory inclusion
      \vspace{-2mm}
  \begin{equation}
    \label{FullInclusionOfTopologicalSpacesIntoTopologicalGroupoids}
    \kTopologicalSpaces \;
    \xhookrightarrow{\;\; \ConstantGroupoid \;\;}
    \;
    \Groupoids(\kTopologicalSpaces)
    \mathrlap{\,.}
  \end{equation}
\end{example}

\begin{example}[Topological pair groupoid]
  \label{TopologicalPairGroupoid}
  For $\TopologicalSpace \,\in\, \kTopologicalSpaces$, its
  {\it chaotic groupoid} or
  {\it pair groupoid}
  is the topological groupoid (Ntn. \ref{TopologicalGroupoids})
  whose space of morphisms is the product of $\TopologicalSpace$ with itself
  (the space of pairs of elements of $\TopologicalSpace$),
  with source and target given by the two canonical projection maps
      \vspace{-2mm}
  $$
    \CodiscreteGroupoid(\TopologicalSpace)
    \;\coloneqq\;
    \big(
    \TopologicalSpace \times \TopologicalSpace
    \underoverset
      {\mathrm{pr}_2}
      {\mathrm{pr}_1}
      {\rightrightarrows}
    \TopologicalSpace
    \big)
    \;\;\;
    \in
    \;
    \TopologicalGroupoids
  $$

      \vspace{-2mm}
\noindent
  and equipped with the unique admissible composition operation:
     \vspace{-2mm}
  $$
    (\TopologicalSpace \times \TopologicalSpace)
      \underset{
        \TopologicalSpace
      }{
        {}_{t}\!\times_s
      }
    (\TopologicalSpace \times \TopologicalSpace)
    \;=\;
    \TopologicalSpace \times \TopologicalSpace \times \TopologicalSpace
    \xrightarrow{ ( \mathrm{pr}_1 , \, \mathrm{pr}_3  ) }
    \TopologicalSpace \times \TopologicalSpace
    \,.
  $$

      \vspace{-2mm}
\noindent
  This construction constitutes {\it another} full subcategory inclusion
     \vspace{-2mm}
  $$
    \kTopologicalSpaces \;
    \xhookrightarrow{\; \CodiscreteGroupoid \;} \;
    \Groupoids(\kTopologicalSpaces)
    \,.
  $$
\end{example}

\begin{definition}[Space of components of a topological groupoid]
\label{SpaceOfConnectedComponentsOfTopologicalGroupoid}
For $(\TopologicalSpace_1 \rightrightarrows \TopologicalSpace_0)$
a topological groupoid (Ntn. \ref{TopologicalGroupoids}),
its {\it space of connected components}
(or: {\it 0-truncation})
\vspace{-2mm}
$$
  \tau_0
  (
    \TopologicalSpace_1
    \rightrightarrows
    \TopologicalSpace_0
  )
  \;\;\;
  \in
  \kTopologicalSpaces
$$

\vspace{-2mm}
\noindent is the quotient space by the source/target relation, hence
the coequalizer of its source and target  maps:
 \vspace{-2mm}
\begin{equation}
  \begin{tikzcd}
    \TopologicalSpace_1
    \ar[
      r,
      shift left=3pt,
      "{s}"{above}
    ]
    \ar[
      r,
      shift right=3pt,
      "{t}"{below}
    ]
    &
    \TopologicalSpace_0
    \ar[
      rr,
      "{\mathrm{coeq}(s,t)}"{above}
    ]
    &&
    \tau_0(\TopologicalSpace_1 \rightrightarrows \TopologicalSpace_0)
    \,.
  \end{tikzcd}
\end{equation}
\end{definition}

All these basic notions are unified as follows:
\begin{proposition}[Adjunctions between topological groupoids and topological spaces]
  \label{AdjunctionBetweenTopologicalGroupoidsAndTopologicalSpaces}
  The 1-category of topological groupoids (Ntn. \ref{TopologicalGroupoids})
  is related to that of topological spaces (Ntn. \ref{CompactlyGeneratedTopologicalSpaces})
  by a quadruple of adjoint functors (Ntn. \ref{AdjointFunctors})
  \vspace{-2mm}
  $$
    \begin{tikzcd}
    \Groupoids(\kTopologicalSpaces)
    \ar[
      rr,
      shift left=4*8pt,
      "{
        \ConnectedGroupoidComponents
      }"{description}
    ]
    \ar[
      rr,
      "{
        \SpaceOfObjects
      }"{description}
    ]
    &&
    \kTopologicalSpaces
    \mathrlap{\,,}
    \ar[
      ll,
      hook',
      shift right=2*8pt,
      "{
        \ConstantGroupoid
      }"{description}
    ]
    \ar[
      ll,
      hook',
      shift left=2*8pt,
      "{
        \CodiscreteGroupoid
      }"{description}
    ]
    \ar[
      ll,
      phantom,
      shift right = 3*8pt,
      "{\scalebox{.6}{$\bot$}}"
    ]
    \ar[
      ll,
      phantom,
      shift right = 1*8pt,
      "{\scalebox{.6}{$\bot$}}"
    ]
    \ar[
      ll,
      phantom,
      shift left = 1*8pt,
      "{\scalebox{.6}{$\bot$}}"
    ]
    \end{tikzcd}
  $$

  \vspace{-2mm}
  \noindent
  where

  \vspace{-4mm}
  \begin{itemize}
  \setlength\itemsep{-4pt}

    \item
    $\ConstantGroupoid$ assigns constant groupoids
    in the sense of
    Ex. \ref{TopologicalSpacesAsTopologicalGroupoids};

    \item
    $\CodiscreteGroupoid$ assigns pair groupoids
    in the sense of
    Ex. \ref{TopologicalPairGroupoid};

    \item
    $\SpaceOfObjects$ assigns spaces of objects \eqref{DiagramForTopologicalGroupoid};

    \item
    $\ConnectedGroupoidComponents$ assigns spaces of
    {\it connected components} (Def. \ref{SpaceOfConnectedComponentsOfTopologicalGroupoid}).
  \end{itemize}
  \vspace{-.2cm}

\end{proposition}
\begin{proof}
  The hom-isomorphisms \eqref{FormingAdjuncts}
  are readily seen by unwinding the definitions:

  \noindent
  {\bf (1)}
  For $\ConnectedGroupoidComponents \dashv \ConstantGroupoid$,
  the natural bijection
      \vspace{-2mm}
  $$
    \kTopologicalSpaces
    \big(
      \ConnectedGroupoidComponents(\TopologicalSpace_1 \rightrightarrows \TopologicalSpace_0)
      ,\,
      \mathrm{Y}
    \big)
    \;\simeq\;
    \Groupoids(\kTopologicalSpaces)
    \big(
      (\TopologicalSpace_1 \rightrightarrows \TopologicalSpace_0)
      ,\,
      \ConstantGroupoid(\mathrm{Y})
    \big)
  $$

  \vspace{-2mm}
  \noindent
  exhibits the universal property of the coequalizer:
      \vspace{-2mm}
  $$
    \begin{tikzcd}[row sep=14pt]
      \TopologicalSpace_1
      \ar[rr]
      \ar[
        d,
        shift left=3pt,
        "{s}"{right}
      ]
      \ar[
        d,
        shift right=3pt,
        "{t}"{left}
      ]
      &&
      \mathrm{Y}
      \ar[
        d,-,
        shift left=1pt
      ]
      \ar[
        d,-,
        shift right=1pt
      ]
      \\
      \TopologicalSpace_0
      \ar[rr]
      \ar[
        d,
        "{\mathrm{coequ}(s,t)}"{left}
      ]
      &&
      \mathrm{Y}
      \\
      \ConnectedGroupoidComponents(\TopologicalSpace_1 \rightrightarrows \TopologicalSpace_0)
      \ar[
        urr,
        dashed,
        "{\exists !}"{below}
      ]
    \end{tikzcd}
  $$

    \vspace{-2mm}
  \noindent
  {\bf (2)} For $\ConstantGroupoid \dashv \SpaceOfObjects$,
  the natural bijection
      \vspace{-2mm}
  $$
    \Groupoids(\kTopologicalSpaces)
    \big(
      \ConstantGroupoid(\TopologicalSpace)
      ,\,
      (\mathrm{Y}_1 \rightrightarrows \mathrm{Y}_0)
    \big)
    \;\simeq\;
    \kTopologicalSpaces
    (
      \TopologicalSpace
      ,\,
      \mathrm{Y}_0
    )
  $$

    \vspace{-2mm}
\noindent
  reflects the unitality of the groupoid composition:
      \vspace{-2mm}
  $$
    \begin{tikzcd}[row sep=small]
      \TopologicalSpace
      \ar[
        rr,
        dashed,
        "{\exists !}"{above}
      ]
      &&
      \mathrm{Y}_1
      \\
      \TopologicalSpace
      \ar[
        rr
      ]
      \ar[
        u,-,
        shift left=1pt,
        "{\mathrm{e}}"{left}
      ]
      \ar[
        u,-,
        shift right=1pt
      ]
      &&
      \mathrm{Y}_0
      \ar[
        u,
        "{e}"
      ]
    \end{tikzcd}
  $$

    \vspace{-2mm}
  \noindent
  {\bf (3)}
  For $\SpaceOfObjects \dashv \CodiscreteGroupoid$, the natural bijection
      \vspace{-2mm}
  $$
    \kTopologicalSpaces
    (
      \TopologicalSpace_0
      ,\,
      \mathrm{Y}
    )
    \;\simeq\;
    \Groupoids(\kTopologicalSpaces)
    \big(
      (\TopologicalSpace_1 \rightrightarrows \TopologicalSpace_0)
      ,\,
      \CodiscreteGroupoid(\mathrm{Y})
    \big)
  $$

      \vspace{-2mm}
\noindent
  reflects the universal property of the Cartesian product:
      \vspace{-2mm}
  $$
    \begin{tikzcd}[row sep=small]
      \TopologicalSpace_1
      \ar[
        rr,
        dashed,
        "{\exists !}"
      ]
      \ar[
        d,
        shift right=3pt,
        "{s}"{left}
      ]
      \ar[
        d,
        shift left=3pt,
        "{\,t}"{right}
      ]
      &&
      \mathrm{Y} \times \mathrm{Y}
      \ar[
        d,
        shift right=3pt,
        "{\mathrm{pr}_1}"{left}
      ]
      \ar[
        d,
        shift left=3pt,
        "{\, \mathrm{pr}_2}"{right}
      ]
      \\
      \TopologicalSpace_0
      \ar[
        rr
      ]
      &&
      \mathrm{Y}
      \mathrlap{\,.}
    \end{tikzcd}
  $$

  $\,$
  \vspace{-1cm}

\end{proof}

We continue to list some classes of examples of topological groupoids
that we need later on.
\begin{example}[Topological action groupoid]
  \label{TopologicalActionGroupoid}
  For $G \acts \, \TopologicalSpace \,\in\, \GActionsOnTopologicalSpaces$
  \eqref{GActionsOnTopologicalSpaces},
  the corresponding {\it action groupoid} is the topological groupoid
  (Ntn. \ref{TopologicalGroupoids}) given by
  \vspace{-3mm}
  \begin{equation}
    \label{LeftActionGroupoid}
    \hspace{-5mm}
    \begin{tikzcd}[row sep=0pt]
    \TopologicalSpace \times G^{\mathrm{op}} \times G^{\mathrm{op}}
    \ar[
      rr,
      "{
        \mathrm{id}
          \times
        \scalebox{.8}{$($}(-) \cdot (-) \scalebox{1.1}{$)$}
      }"
    ]
    &&
    \;\;
    \TopologicalSpace \times G^{\mathrm{op}}
    \;\;
      \ar[
        rr,
        shift left=13pt,
        "{
          \mathrm{pr}_1
        }"
      ]
      \ar[
        rr,
        shift right=5pt,
        "{
          \scalebox{0.7}{$(-) \cdot (-) $}
        }"{description}
      ]
    &&
  \quad   \TopologicalSpace
    \ar[
      ll,
      shift right=5pt,
      "{
        \mathrm{id} \times e
      }"{description}
    ]
    \\
\scalebox{0.8}{$    (x, g_1, g_2) $}
    &\longmapsto&
 \scalebox{0.8}{$   (x, g_2 \cdot g_1) $}
   &\;\;\;\;\; \longmapsto&
 \quad  \scalebox{0.8}{$   (g_2 \cdot g_1 \cdot x) $}
    \end{tikzcd}
    \quad
    \in
    \;
    \TopologicalGroupoids
  \end{equation}

  \vspace{-1mm}
  \noindent
  with composition given by the {\it reverse} of the group operation in $G$.
\end{example}

\begin{example}[Topological delooping groupoid]
  \label{TopologicalDeloopingGroupoid}
  For $\Gamma \,\in\, \Groups(\kTopologicalSpaces)$,
  its {\it delooping groupoid} is the
  topological left action groupoid (Ex. \ref{TopologicalActionGroupoid})
  of the unique $\Gamma^{\mathrm{op}}$-action
   $\Gamma^{\mathrm{op}} \acts \, \ast \,\in\, \Actions{\Gamma}(\kTopologicalSpaces)$
  on the point space:
      \vspace{-1mm}
  $$
    \mathbf{B}\Gamma
    \;\coloneqq\;
    (
      \ast \times \Gamma
      \rightrightarrows
      \ast
    )
    \;=\;
    (
    \Gamma
    \rightrightarrows
    \ast
    )
    \;\;\;
    \in
    \;
    \TopologicalGroupoids
    \,,
  $$

  \vspace{-1mm}
  \noindent
  with composition given by the group operation in $\Gamma$:
  \vspace{-2mm}
  $$
    \begin{tikzcd}
      \bullet
      \ar[
        r,
        "{g_1}"
      ]
      \ar[
        rr,
        rounded corners,
        to path={
          -- ([yshift=+5pt]\tikztostart.north)
          --node[above]{
              \scalebox{.7}{$g_1 \cdot g_2$}
            }
            ([yshift=+5pt]\tikztotarget.north)
          -- (\tikztotarget.north)}
      ]
      &
      \bullet
      \ar[
        r,
        "{g_2}"
      ]
      &
      \bullet
    \end{tikzcd}
    \;\;\;
    \in
    \;
    \mathbf{B}\Gamma
    \,.
  $$
\end{example}

\begin{example}[Action groupoid of group multiplication is pair groupoid]
  \label{ActionGroupoidOfLeftGroupMultiplicationIsPairGroupoid}
  For $\Gamma \,\in\, \Groups(\kTopologicalSpaces)$,
  the topological pair groupoid (Ex. \ref{TopologicalPairGroupoid})
  on its underlying topological space
  is
  isomorphic to the action groupoid (Ex. \ref{TopologicalActionGroupoid})
  of both the left and inverse-right action of $G^{\mathrm{op}}$
  on itself \eqref{LeftMultiplicationAndInverseRightMultiplicationActionsOnATopologicalGroup}
      \vspace{-2mm}
  $$
    \begin{tikzcd}[row sep=-5pt]
      \overset{
        \mathclap{
        \raisebox{3pt}{
          \tiny
          \color{darkblue}
          \bf
          \def\arraystretch{.9}
          \begin{tabular}{c}
            left multiplication
            \\
            action groupoid
          \end{tabular}
        }
        }
      }{
        (
          G \times G
          \rightrightarrows G
        )
      }
      &&
      \overset{
        \mathclap{
        \raisebox{3pt}{
          \tiny
          \color{darkblue}
          \bf
          \def\arraystretch{.9}
          \begin{tabular}{c}
            pair groupoid
          \end{tabular}
        }
        }
      }{
        (
          G \times G
          \rightrightarrows G
        )
      }
      \ar[
        ll,
        "\sim"{above, yshift=-1pt}
      ]
      \ar[
        rr,
        "\sim"{above, yshift=-1pt}
      ]
      &&
      \overset{
        \mathclap{
        \raisebox{3pt}{
          \tiny
          \color{darkblue}
          \bf
          \def\arraystretch{.9}
          \begin{tabular}{c}
            right multiplication
            \\
            action groupoid
          \end{tabular}
        }
        }
      }{
        (
          G \times G
          \rightrightarrows G
        )
      }
      \\
 \scalebox{0.7}{$     \big(
        g_1
        \xrightarrow{ ( g_1, \, g_1 \cdot g_2^{-1} ) }
        g_2
      \big)
      $}
      &\longmapsfrom&
  \scalebox{0.7}{$         \big(
        g_1
        \xrightarrow{ (g_1, g_2) }
        g_2
      \big)
      $}
      &\longmapsto&
     \scalebox{0.7}{$      \big(
        g_1
        \xrightarrow{ ( g_1, \, g_1^{-1} \cdot g_2 ) }
        g_2
      \big)
      $}
    \end{tikzcd}
  $$
\end{example}
\begin{example}[Topological mapping groupoid]
  \label{TopologicalFunctorGroupoid}
  Given
  a pair of topological groupoids (Ntn. \ref{TopologicalGroupoids}),
  their {\it mapping groupoid} or {\it functor groupoid}
  (e.g. \cite[\S 2]{NiefieldPronk19})
  \vspace{-2mm}
  \begin{equation}
    \label{MappingGroupoidOfTopologicalGroupoids}
    \mathrm{Maps}
    (
      \TopologicalSpace_1 \rightrightarrows \mathrm{X_0},
      \,
      \mathrm{Y}_1 \rightrightarrows \mathrm{Y_0}
    )
    \;\;\;
    \in
    \;
    \TopologicalGroupoids
  \end{equation}

  \vspace{-2mm}
  \noindent
  has as object space the subspace of the
  product of mapping spaces
  $\mathrm{Maps}(\TopologicalSpace_0,\mathrm{Y}_0) \times \mathrm{Maps}(\TopologicalSpace_1, \mathrm{Y}_1)$
  \eqref{MappingSpace}
  on the elements that satisfy the functoriality condition \eqref{ContinuousFunctors},
  and
  as morphism space the subspace on the product of that space with the
  mapping space $\mathrm{Maps}(\TopologicalSpace_0, \TopologicalSpace_1)$
  on those elements which satisfy the naturality condition \eqref{NaturalitySquareForTopologicalGroupoids}.

  This construction is a $\Groupoids$-enriched functor in both arguments
  contravariantly so in the first:
  \vspace{-2mm}
  $$
    \mathrm{Maps}
    (- ,\, - )
    \;:\;
    \begin{tikzcd}
    \Groupoids(\kTopologicalSpaces)^{\mathrm{op}}
      \times
    \Groupoids(\kTopologicalSpaces)
    \ar[r]
    &
    \Groupoids(\kTopologicalSpaces)
    \,.
    \end{tikzcd}
  $$

  \vspace{-2mm}
  \noindent
  With the first argument fixed, it constitutes a
  $\Groupoids$-enriched
  right adjoint (Ntn. \ref{AdjointFunctors})
  to the product functor (e.g. \cite[Prop. 3.1]{NiefieldPronk19}):
  \vspace{-2mm}
  \begin{equation}
    \label{InternalHomAdjunctionForTopologicalGroupoids}
    \begin{tikzcd}[column sep=40pt]
    \Groupoids(\kTopologicalSpaces)
    \ar[
      rr,
      shift right=6pt,
      "{
       \scalebox{.7}{$ \mathrm{Maps}
        \left(
          (\TopologicalSpace_1 \rightrightarrows \, \TopologicalSpace_0)
          ,\,
          -
        \right)
        $}
      }"{below}
    ]
    \ar[
      rr,
      phantom,
      "{\scalebox{.7}{$\bot$}}"
    ]
    &&
    \Groupoids(\kTopologicalSpaces)
    \mathrlap{\,.}
    \ar[
      ll,
      shift right=6pt,
      "{
         \scalebox{.7}{$ (\TopologicalSpace_1 \rightrightarrows \, \TopologicalSpace_0)
          \times
        (-)$}
      }"{above}
    ]
    \end{tikzcd}
  \end{equation}
\end{example}

\begin{example}[Mapping groupoid between delooping groupoids]
  \label{MappingGroupoidBetweenDeloopingGroupoids}
  For $G, \Gamma \,\in\, \Groups(\kTopologicalSpaces)$,
  the mapping groupoid (Ex. \ref{TopologicalFunctorGroupoid})
  between their topological delooping groupoids
  (Ex. \ref{TopologicalDeloopingGroupoid})
  is isomorphic to the topological action groupoid
  (Ex. \ref{TopologicalActionGroupoid})
  of the adjoint action of $\Gamma$
  on the hom-set of group homomorphisms $G \xrightarrow{\;} \Gamma$
  (topologized as a subspace of $\mathrm{Maps}(G,\Gamma)$):
      \vspace{-5mm}
  \begin{equation}
    \label{MappingGroupoidOfDeloopingGroupoidsIsAdjointActionGroupoid}
    \hspace{2cm}
    \mathrm{Maps}
    (
      G \rightrightarrows \ast
      ,\,
      \Gamma \rightrightarrows \ast
    )
    \;\;
      \simeq
    \;\;
    \big(
      \Groups(G,\, \Gamma)
        \times
      \Gamma^{\mathrm{op}}
        \; \rightrightarrows \;
      \Groups(G,\, \Gamma)
    \big)
    \,.
  \end{equation}
  \vspace{-6mm}
  $$
  \hspace{-2cm}
    \begin{tikzcd}[column sep=large]
      \bullet
      \ar[
        d,
        "g_1"
      ]
      \ar[
        dd,
        rounded corners,
        to path={
           -- ([xshift=-10pt]\tikztostart.west)
           --node[below, sloped]{\rotatebox{0}{\scalebox{.7}{$
               g_1 \cdot g_2
             $}}} ([xshift=-10pt]\tikztotarget.west)
           -- (\tikztotarget.west)}
      ]
      &[-10pt]&[-10pt]
      \bullet
      \ar[
        d,
        "\phi(g_1)"
      ]
      \ar[
        dd,
        rounded corners,
        to path={
           -- ([xshift=-10pt]\tikztostart.west)
           --node[below, sloped]{\rotatebox{0}{\scalebox{.7}{$
               \phi(g_1 \cdot g_2)
             $}}} ([xshift=-10pt]\tikztotarget.west)
           -- (\tikztotarget.west)}
      ]
      \ar[
        rr,
        "\gamma"
      ]
      &&
      \bullet
      \ar[
        d,
        "\phi'(g_1)"{left}
      ]
      \\
      \bullet
      \ar[
        d,
        "g_2"
      ]
      &&
      \bullet
      \ar[
        d,
        "\phi(g_2)"
      ]
      \ar[
        rr,
        "\gamma"
      ]
      &&
      \bullet
      \ar[
        d,
        "\phi'(g_2)"{left}
      ]
      \\
      \bullet
      &&
      \bullet
      \ar[
        rr,
        "\gamma"
      ]
      &&
      \bullet
    \end{tikzcd}
    \;\;\;\;\;\;\longmapsto\;\;\;\;\;\;
 \scalebox{0.7}{$    \big(
      \phi
      \xrightarrow{ (\gamma,\phi) }
      \phi'
      =
      \mathrm{ad}_\gamma \circ \phi
    \big)
    {\phantom{AAAAAAAAAAAAA}}
    $}
  $$
\end{example}

\medskip

\noindent
{\bf Crossed homomorphisms and first non-abelian group cohomology.}
The classical notion of {\it crossed homomorphisms}, recalled as
Def. \ref{CrossedHomomorphismsAndFirstNonAbelianGroupCohomology} below,
turns out to play a pivotal role in equivariant bundle theory
(e.g., Lem. \ref{EquivariantPrincipalTwistedProductBundles}
and Prop. \ref{FixedLociOfBaseOfUniversalEquivariantPrincipalGroupoid}),
often secretly so (Rem. \ref{GraphsOfCrossedHomomorphismsInTheLiterature} below).
Here we highlight a transparent groupoidal understanding
of crossed homomorphisms with crossed conjugations
between them (Prop. \ref{ConjugationGroupoidOfCrossedHomomorphismsIsSectionsOfDeloopedSemidirectProductProjection} below).

\begin{definition}[Crossed homomorphisms and first non-abelian group cohomology]
  \label{CrossedHomomorphismsAndFirstNonAbelianGroupCohomology}
  Let $G \,\in\, \Groups(\kTopologicalSpaces)$
  and $G \acts \, \Gamma \,\in\, \Actions{G}(\kTopologicalSpaces)$,
  with
  $\alpha : G \xrightarrow{\;} \mathrm{Aut}_{\mathrm{Grp}}(\Gamma)$
  the underlying automorphism action (Lem. \ref{EquivariantTopologicalGroupsAreSemidirectProductsWithG}).

  \noindent
  {\bf (i)} A continuous {\it crossed homomorphism} from $G$ to $\Gamma$
  is a continuous map $\phi \colon  G \xrightarrow{\;} \Gamma$
  which satisfies the following {\it $G$-crossed} homomorphism property:
  \vspace{-2mm}
  \begin{equation}
    \label{GCrossedHomomorphismProperty}
    \underset{
      g_1, g_2 \in G
    }{\forall}
    \;\;\;
    \phi(g_1 \cdot g_2)
    \;=\;
    \phi(g_1)
      \cdot
    \alpha(g_1)
    \left(
      \phi(g_2)
    \right)
    \,.
  \end{equation}

  \vspace{-2mm}
  \noindent
  We write
  \begin{equation}
    \label{SpaceOfCrossedHomomorphisms}
    \CrossedHomomorphisms(G, \, G \acts \, \Gamma)
      \;\subset\;
    \mathrm{Maps}(G,\Gamma)
    \;\;\;
    \in
    \;
    \kTopologicalSpaces
  \end{equation}
  for the subspace of the mapping space \eqref{MappingSpace}
  on the crossed homomorphisms.

  \noindent
  {\bf (ii)} A {\it crossed conjugation} between two crossed homomorphisms
  $\phi \xrightarrow{\;} \phi'$ is an element $\gamma \,\in\, \Gamma$ such that
  \vspace{-2mm}
  \begin{equation}
    \label{CrossedConjugationAction}
    \underset{g \in G}{\forall}
    \;\;\;
    \phi'(g)
      \;=\;
    \gamma^{-1}
      \cdot
    \phi(g)
      \cdot
    \alpha(g)(\gamma)
    \,.
  \end{equation}

  \vspace{-2mm}
\noindent  We denote the continuous $\Gamma$-action by crossed conjugation  by
 \vspace{-2mm}
  \begin{equation}
    \label{ActionOnSpaceOfCrossedHomomorphisms}
    \Gamma
      \acts \;
    \CrossedHomomorphisms(G,\, G \acts \, \Gamma)
    \;\;\;
    \in
    \;
    \Actions{\Gamma}(\kTopologicalSpaces)
    \,.
  \end{equation}

  \vspace{-2mm}
  \noindent  We write $\phi \sim_{\mathrm{ad}} \phi'$
  for the corresponding equivalence relation.

  \noindent
  {\bf (iii)} The {\it non-abelian group cohomology} of $G$
  in degree 1
  with coefficients
  in $G \acts \, \Gamma$ is
  -- at least when $G$ is discrete\footnote{
  For non-discrete domain groups
  the notion of crossed homomorphisms need no longer capture all
  1-cocycles in non-abelian group cohomology, when the latter is formulated in
  proper stacky generality; see \cite{WagemannWockel11} for pointers.} --
  the set of connected components\footnote{
    The passage to connected components in \eqref{NonAbelianGroup1Cohomology}
    seems not to be considered in existing literature.
    It makes no difference when the coefficient group is discrete,
    (which tends to be tacitly understood in this context,
    but is not the most general case of interest),
    as well as under other sufficient conditions discussed in
    Prop. \ref{DiscreteSpacesOfCrossedConjugacyClassesOfCrossedHomomorphisms} below.
    But in general
    the correct homotopy-meaningful definition of non-abelian group 1-cohomology
    (see Rem. \ref{GeneralAbstractPerspectiveOnNonAbelianGroup1Cohomology} below)
    is
    only obtained with passage to connected components included.
    (The statement in \cite[\S 4.3]{GuillouMayMerling17},
    that any groupoid is equivalent to the coproduct of its automorphism sub-groupoids,
    is patently false for topological groupoids, in general.
    It does hold under suitable extra conditions, such
    as in Prop. \ref{DiscreteSpacesOfCrossedConjugacyClassesOfCrossedHomomorphisms} below).
    For further discussion of this point see the companion article (cite).
  }
  of the quotient space by crossed conjugation classes \eqref{CrossedConjugationAction}
  of the space \eqref{SpaceOfCrossedHomomorphisms}
  of crossed homomorphisms \eqref{GCrossedHomomorphismProperty}:
  \vspace{-2mm}
  \begin{equation}
    \label{NonAbelianGroup1Cohomology}
    H^1_{\mathrm{Grp}}
    (
      G
      ,\,
      G \acts \, \Gamma
    )
    \;\;
    \coloneqq
    \;\;
    \pi_0
    \big(
      \CrossedHomomorphisms(G, \, G \acts \, \Gamma)
      /\!\sim_{\mathrm{ad}}
    \!\!\big)
    \;\;\;\;
    \in
    \;
    \Sets^{\ast/}
    \,.
  \end{equation}
\end{definition}
\begin{remark}[Crossed homomorphisms in the literature]
\label{CrossedHomomorphismsInTheliterature}
  Since the notion of crossed homomorphisms,
  in the generality that we need them here
  (Def. \ref{CrossedHomomorphismsAndFirstNonAbelianGroupCohomology}),
  tends to be neglected in the literature, we record some pointers:
Crossed homomorphisms \eqref{GCrossedHomomorphismProperty}
appear first, already in full non-abelian generality,
in \cite[(3.1)]{Whitehead49}.
Much later, following \cite[\S IV.2]{MacLane75},
they became widely appreciated only in the special case when $\Gamma$ is an abelian group,
as a tool in ordinary group cohomology (e.g. \cite[p. 45]{Brown82}).
Crossed homomorphisms in their non-abelian generality
appear again
in \cite[\S 2.1]{tomDieck69} (not using the ``crossed'' terminology, though)
and in \cite[p. 2]{MurayamaShimakawa95}\cite[Def. 4.1]{GuillouMayMerling17},
all in the context of equivariant bundle theory
(in which we consider them in \cref{ConstructionOfUniversalEquivariantPrincipalBundles}).
Textbook accounts in this generality are in \cite[p. 16]{NSW08}\cite[\S 15.a-b]{Milne17}.
The corresponding definition
\eqref{NonAbelianGroup1Cohomology}
of non-abelian group 1-cohomology is
rarely made explicit;
exceptions are
\cite[Def. 2.3.2]{GilleSzamuely06}\cite[Def. 4.17]{GuillouMayMerling17}\cite[\S 3.k]{Milne17}\footnote{
These are sections 16.a-b \& 27.a in the expanded version of Milne's book at
\href{https://www.jmilne.org/math/CourseNotes/iAG200.pdf}{\tt www.jmilne.org/math/CourseNotes/iAG200.pdf}}.
\end{remark}
We also need to recall the following standard fact (e.g. \cite[Ex. 15.1]{Milne17}):

\begin{lemma}[Crossed homomorphisms are sections of the semidirect product projection]
  \label{CrossedHomomorphismsAreSectionsOftheSemidirectProductProjection}
  {\bf (i)}
  Crossed homomorphisms $\phi \colon G \to \Gamma$ (Def. \ref{MappingGroupoidBetweenDeloopingGroupoids})
  are in bijective correspondence to homomorphic sections
  of the semidirect group projection \eqref{SplitGroupExtensionOfGByGamma}:
  \vspace{-2mm}
  $$
    \begin{tikzcd}[column sep=60pt, row sep=small]
      &
      \Gamma \rtimes G
      \ar[
        d,
        "\; \mathrm{pr}_2"
      ]
      \\
      G
      \ar[
        r,-,
        shift left=1pt
      ]
      \ar[
        r,-,
        shift right=1pt
      ]
      \ar[
        ur,
        dashed,
        "{
          g \,\mapsto\,
          (
            \phi(g)
            \,,
            g
          )
        }"{above, yshift=1pt, sloped}
      ]
      &
      G
    \end{tikzcd}
  $$

  \vspace{-2mm}
\noindent  {\bf (ii)} Under this identification, crossed conjugations \eqref{CrossedConjugationAction}
  are equivalently plain conjugations with elements in
  $\Gamma \xrightarrow{i} \Gamma \rtimes G$.
\end{lemma}

\begin{proof}
  Having a section means that
  $g \,\mapsto\, ( \phi(g),\, g ) \,\in\, \Gamma \rtimes G$,
  and this being a homomorphism means that
    \vspace{-1mm}
  \begin{equation}
    \label{CrossedHomomorphismAsPlainHomomorphismsToSemidirectProductGroup}
    (
      \phi(g_1 \cdot g_2)
      ,\,
      g_1 \cdot g_2
    )
    \;=\;
    (
      \phi(g_1)
      ,\,
      g_1
    )
    \cdot
    (
      \phi(g_2)
      ,\,
      g_2
    )
    \;=\;
    \scalebox{1.3}{$($}
      \phi(g_1)
        \cdot
      \alpha(g_2)
      \scalebox{1.15}{$($}
        \phi(g_2)
      \scalebox{1.15}{$)$}
      ,\,
      g_1 \cdot g_2
    \scalebox{1.3}{$)$}
    \,,
  \end{equation}

    \vspace{-1mm}
\noindent
  where the second equality on the right is the definition of the semidirect product group operation,
  evidently reproducing the defining condition \eqref{GCrossedHomomorphismProperty}
  in the first argument.

  Analogously, plain conjugation in the semidirect product group with elements of the form
  $(\gamma, \mathrm{e}) \,\in\, \Gamma \rtimes G$
  gives
    \vspace{-2mm}
  \begin{equation}
    \label{CrossedConjugationAsPlainConjugationInSemidirectProductGroup}
    (\gamma, \mathrm{e})^{-1}
    \cdot
    \left(
      \phi(g) ,\, g
    \right)
    \cdot
    (\gamma, \mathrm{e})
    \;\;
    =
    \;\;
    \big(
      \gamma^{-1}
      \cdot
      \phi(g)
      \cdot
      \alpha(g)(\gamma)
      ,\,
      g
    \big)
    \mathrlap{\,,}
  \end{equation}
   reproducing the formula
  \eqref{CrossedConjugationAction} in the first argument.
\end{proof}
Similarly elementary, but maybe less widely appreciated
(see Rem. \ref{GraphsOfCrossedHomomorphismsInTheLiterature}),
is the following:

\begin{lemma}[Graphs of crossed homomorphisms {\cite[\S 2.1]{tomDieck69}\cite[Lem. 4.5]{GuillouMayMerling17}}]
  \label{CrossedHomomorphismsAreEquivalentlyMaySubgroupsOfSemidirectProducts}
  $\,$

  \noindent
  {\bf (i)}
  The graph of a crossed homomorphism $\phi \;:\; G \xrightarrow{\;} \Gamma$
  (Def. \ref{CrossedHomomorphismsAndFirstNonAbelianGroupCohomology})
  is a subgroup
  \vspace{-1mm}
  \begin{equation}
    \label{MaySubgroupsOfSemidirectProducts}
    \widehat{G}
    \;\in\;
    \Gamma \rtimes G
    \,,
    \;\;\;\;\;\;\;
    \mbox{\rm such that}
    \;\;\;\;\;
    \mathrm{pr}_2
    (
      \widehat{G}
    )
    \,\simeq\,
    G
    \;\;\;\;
    \mbox{\rm and}
    \;\;\;\;
    \widehat{G}
      \cap
    i(\Gamma)
    \,\simeq\,
    \{
      (\mathrm{e},\,\mathrm{e})
    \}
    \,,
  \end{equation}

\vspace{-2mm}
\noindent
  where $i : \Gamma \xhookrightarrow{\;} \Gamma \rtimes G$ is the canonical
subgroup inclusion \eqref{SplitGroupExtensionOfGByGamma}.

  \noindent
  {\bf (ii)}
  Every such subgroup \eqref{MaySubgroupsOfSemidirectProducts}
  is the graph of a unique crossed homomorphism.
\end{lemma}

\begin{proof}
  The first statement is immediate from the definitions.
  For the converse statement (ii),
  consider a subgroup $\widehat{G}$ as in \eqref{MaySubgroupsOfSemidirectProducts}.
  Then the subgroup property implies that
  \vspace{-2mm}
  $$
    (\gamma,\,g)
    ,\,
    (\gamma{\;}',\,g)
    \;\in\;
    \widehat{G}
    \;\;\;\;\;\;\;
    \Rightarrow
    \;\;\;\;\;\;\;
    (\gamma{\;}',\, g)
    \cdot
    (\gamma,\, g)^{-1}
    \;=\;
    (\gamma{\;}',\, g)
    \cdot
    \left(
      \alpha(g^{-1})(\gamma^{-1})
      ,\,
      g^{-1}
    \right)
    \;=\;
    \big(
      \gamma{\;}' \cdot \gamma^{-1}
      ,\,
      \mathrm{e}
    \big)
    \;\;\;
    \in
    \;
    \widehat{G}
    \,.
  $$

\vspace{-2mm}
\noindent
  From this, the second condition in \eqref{MaySubgroupsOfSemidirectProducts}
  implies that
  \vspace{-2mm}
  $$
    (\gamma,\,g)
    ,\,
    (\gamma{\;}',\,g)
    \;\in\;
    \widehat{G}
    \;\;\;\;\;\;\;
    \Rightarrow
    \;\;\;\;\;\;\;
    \gamma \,=\, \gamma{\;}'
    .
  $$

\vspace{-2mm}
\noindent
  Together with the first condition in \eqref{MaySubgroupsOfSemidirectProducts},
  this implies that
  $\widehat G$ is the graph of a function $\phi \,:\, G \times \Gamma$.
  From this the claim follows by
  Lem. \ref{CrossedHomomorphismsAreSectionsOftheSemidirectProductProjection}.
\end{proof}

\begin{remark}[Graphs of crossed homomorphisms in the literature]
  \label{GraphsOfCrossedHomomorphismsInTheLiterature}
  Subgroups of the form \eqref{MaySubgroupsOfSemidirectProducts}
  were used in early articles on equivariant bundle
  theory (e.g. \cite[Thm. 10]{LashofMay86}\cite[Thm. 7]{May90}).
  That these are equivalently (graphs of) crossed homomorphisms
  (Lem. \ref{CrossedHomomorphismsAreEquivalentlyMaySubgroupsOfSemidirectProducts})
  and hence homomorphic sections of the semidirect product group
  (Lem. \ref{CrossedHomomorphismsAreSectionsOftheSemidirectProductProjection})
  may have been
  (in view of
  Prop. \ref{ConjugationGroupoidOfCrossedHomomorphismsIsSectionsOfDeloopedSemidirectProductProjection}
  below)
  one of the key observations that led to the construction in
  \cite{MurayamaShimakawa95}
  (discussed in \cref{ConstructionOfUniversalEquivariantPrincipalBundles} below);
 however, Lem. \ref{CrossedHomomorphismsAreEquivalentlyMaySubgroupsOfSemidirectProducts}
  is still not made explicit there.
\end{remark}

\begin{notation}[Conjugation groupoid of crossed homomorphisms]
  \label{ConjugationGroupoidOfCrossedHomomorphisms}
  We write
  \vspace{-2mm}
  $$
    \CrossedHomomorphisms(G,\,G \acts \, \Gamma) \sslash_{\!\mathrm{ad}} \Gamma
    \;\;\coloneqq\;\;
    \big(
      \CrossedHomomorphisms(G,\, G \acts \,  \Gamma)
        \times
      \Gamma
     \;\; \rightrightarrows \;\;
      \CrossedHomomorphisms(G ,\, G \acts \, \Gamma)
    \big)
  $$

  \vspace{-2mm}
\noindent
  for the topological action groupoid (Ex. \ref{TopologicalActionGroupoid})
  of crossed conjugations \eqref{CrossedConjugationAction}
  acting on the space
  \eqref{SpaceOfCrossedHomomorphisms} of
  crossed homomorphisms \eqref{GCrossedHomomorphismProperty}.
\end{notation}

\begin{definition}[Topological groupoid of sections of delooped semidirect product projection]
  \label{GroupoidOfSectionsOfDeloopedSemidirectProductProjection}
  For $G \,\in\, \Groups(\kTopologicalSpaces)$
  and $G \acts \, \Gamma \,\in\, \Actions{G}( \kTopologicalSpaces)$,
  consider the topological mapping groupoid (Ex. \ref{TopologicalFunctorGroupoid})
  from the delooping groupoid (Ex. \ref{TopologicalDeloopingGroupoid})
  of $G$ into that of the semidirect product group $\Gamma \rtimes G$
  (Lem. \ref{EquivariantTopologicalGroupsAreSemidirectProductsWithG})
  restricted to those functors and natural transformations
  which are sections, in that their
  composition with the projection \eqref{SplitGroupExtensionOfGByGamma}
  back to $\mathbf{B}G$ is the identity:
  \vspace{-3mm}
  \begin{equation}
    \label{GroupoidOfCrossedHomomorphisms}
    \hspace{-2mm}
    \SliceMaps{\big}{\mathbf{B}G}
      { \mathbf{B}G }
      { \mathbf{B}(\Gamma \rtimes G) }
    \;\coloneqq\;
    \mathrm{Maps}
    \left(
      \mathbf{B}G
      ,\,
      \mathbf{B}(\Gamma \rtimes G)
    \right)
    \!\!\!
    \underset{
      \raisebox{-5pt}{
        \scalebox{.7}{$
          \mathrm{Maps}
          (
            \mathbf{B}G
            ,\,
            \mathbf{B}G
          )
        $}
      }
    }{\times}
    \!\!\!
    \{\mathrm{id}\}
    \;=\;
    \Biggg\{\!\!
       \begin{tikzcd}
      &&
      \mathbf{B}(\Gamma \rtimes G)
      \ar[
        d,
        "\mathbf{B}\mathrm{pr}_2"
      ]
      \\
      \mathbf{B}G
      \ar[
        rr,-,
        shift right=1pt
      ]
      \ar[
        rr,-,
        shift left=1pt
      ]
      \ar[
        urr,
        dashed,
        bend left=30,
        "\ "{right, name=s}
      ]
      \ar[
        urr,
        dashed,
        bend right=10,
        "\ "{left, name=t}
      ]
      &&
      \mathbf{B}G
      \ar[
        from=s,
        to=t,
        dashed,
        Rightarrow
      ]
    \end{tikzcd}
        \!\!\!\! \Biggg\}
    \,.
  \end{equation}
\end{definition}

\begin{proposition}[Conjugation groupoid of crossed homomorphisms is sections of delooped semidirect product]
  \label{ConjugationGroupoidOfCrossedHomomorphismsIsSectionsOfDeloopedSemidirectProductProjection}
  The groupoid of sections
  of the delooped semidirect product projection (Def. \ref{GroupoidOfSectionsOfDeloopedSemidirectProductProjection})
  is isomorphic to
  the conjugation groupoid of crossed homomorphisms
  (Ntn. \ref{ConjugationGroupoidOfCrossedHomomorphisms}):
  \vspace{-2mm}
  \begin{equation}
    \label{ConjugationGroupoidOfCrossedHomomorphismsIsSlicedMappingGroupoidFromBGToBGammaRtimesG}
    \SliceMaps{}{\mathbf{B}G}
      { \mathbf{B}G }
      { \mathbf{B}(\Gamma \rtimes G) }
    \;\simeq\;
    \CrossedHomomorphisms(G,\, G \acts \, \Gamma) \sslash_{\!\!\mathrm{ad}} \Gamma
    \,.
  \end{equation}

  \vspace{-1mm}
  \noindent
  Hence its connected components (Def. \ref{SpaceOfConnectedComponentsOfTopologicalGroupoid})
  are bijective to the
  non-abelian group 1-cohomology \eqref{NonAbelianGroup1Cohomology}:
  \vspace{-1mm}
  $$
    \tau_0
      \SliceMaps{\big}{\mathbf{B}G}
        { \mathbf{B}G }
        { \mathbf{B}(\Gamma \rtimes G) }
    \;\simeq\;
    H^1_{\mathrm{Grp}}
    (
      G
      ,\,
      G \acts \, \Gamma
    )
    \,.
  $$
  \vspace{-.5cm}

\end{proposition}
\begin{proof}
By definition, a morphism in the groupoid
\eqref{GroupoidOfCrossedHomomorphisms} is a commuting diagram
of functors and natural transformations as shown on the left of the following:
\vspace{-4mm}
  $$
    \begin{tikzcd}[row sep=52pt]
      \\[-30pt]
      &&
      \mathbf{B}(\Gamma \!\!\rtimes_\alpha\! G)
      \ar[
        d,
        "\mathbf{B}\mathrm{pr}_2"
      ]
      \\
      \mathbf{B}G
      \ar[
        rr,-,
        shift right=1pt
      ]
      \ar[
        rr,-,
        shift left=1pt
      ]
      \ar[
        urr,
        dashed,
        bend left=40,
        "{F}"{description, pos=.525},
        "{
          \mbox{
            \tiny
            \color{darkblue}
            \bf
            sliced
          }
        }"{above, yshift=+9pt, sloped},
        "{
          \mbox{
            \tiny
            \color{darkblue}
            \bf
            functor
          }
        }"{above, yshift=+3pt, sloped},
        "\ "{right, name=s}
      ]
      \ar[
        urr,
        dashed,
        bend right=20,
        "{F'}"{description},
        "\ "{left, name=t}
      ]
      &&
      \mathbf{B}G
      \ar[
        from=s,
        to=t,
        dashed,
        Rightarrow,
        "{
          \mbox{
            \tiny
            \color{greenii}
            \bf
            sliced
          }
        }"{above, yshift=7pt, sloped},
        "{
          \mbox{
            \tiny
            \color{greenii}
            \bf
            transf.
          }
        }"{above, yshift=2pt, sloped}
      ]
    \end{tikzcd}
    \raisebox{-14pt}{
     $:$
    }
    \;\;
    \begin{tikzcd}[row sep=32pt, column sep=50pt]
      \bullet
      \ar[
        d,
        "{
          g_1
        }"{right}
      ]
      \ar[
        dd,
        rounded corners,
        to path={
          -- ([xshift=-10pt]\tikztostart.west)
          --node[above, sloped]{\scalebox{.7}{
             $g_1 \cdot g_2$
          }} ([xshift=-10pt]\tikztotarget.west)
          -- (\tikztotarget.west)}
      ]
      &[-15pt]
      &[-5pt]
      \overset{
        \hspace{-20pt}
        \mathclap{
        \raisebox{6pt}{
          \tiny
          \color{darkblue}
          \bf
          \def\arraystretch{.9}
          \begin{tabular}{c}
            crossed
            \\
            homo-
            \\
            morphism
          \end{tabular}
        }
        }
      }{
        \bullet
      }
      \ar[
        rr,
        "{
          \scalebox{1.3}{$($}
            \gamma, \, \mathrm{e}
          \scalebox{1.3}{$)$}
        }"{description},
        "{
          \mbox{
            \tiny
            \color{greenii}
            \bf
            \def\arraystretch{.9}
            \begin{tabular}{c}
              crossed
              \\
              conjugation
            \end{tabular}
          }
        }"{above,yshift=6pt}
      ]
      \ar[
        d,
        "{
          \scalebox{1.3}{$($}
            \phi(g_1), \, g_1
          \scalebox{1.3}{$)$}
        }"
      ]
      \ar[
        dd,
        rounded corners,
        to path={
          -- ([xshift=-12pt]\tikztostart.west)
          --node[above, sloped]{\scalebox{.6}{
             $
               \scalebox{1.3}{$($}
                 \phi(g_1)
                 \!\cdot\!
                 \alpha(g_1)
                 \scalebox{1.15}{$($}
                   \phi(g_2)
                 \scalebox{1.15}{$)$}
                 ,\,
                 g_1 \!\cdot\! g_2
               \scalebox{1.3}{$)$}
             $
          }} ([xshift=-12pt]\tikztotarget.west)
          -- (\tikztotarget.west)}
      ]
      &&
      \bullet
      \ar[
        d,
        "{
          \scalebox{1.3}{$($}
            \phi'(g_1), \, g_1
          \scalebox{1.3}{$)$}
          \,=\,
          \scalebox{1.3}{$($}
            \gamma^{-1}
            \cdot
            \phi(g_1)
            \cdot
            \alpha(g_1)
            (\gamma)
          \scalebox{1.3}{$)$}
        }"
      ]
      \\
      \bullet
      \ar[
        d,
        "{g_2}"{right}
      ]
      & \hspace{-1cm} \mapsto &
      \bullet
      \ar[
        rr,
        "{
          \scalebox{1.3}{$($}
            \gamma, \, \mathrm{e}
          \scalebox{1.3}{$)$}
        }"{description}
      ]
      \ar[
        d,
        "{
          \scalebox{1.3}{$($}
            \phi(g_2), \, g_2
          \scalebox{1.3}{$)$}
        }"
      ]
      &&
      \bullet
      \ar[
        d,
        "{
          \scalebox{1.3}{$($}
            \phi'(g_2), \, g_2
          \scalebox{1.3}{$)$}
          \,=\,
          \scalebox{1.3}{$($}
            \gamma^{-1}
            \cdot
            \phi(g_2)
            \cdot
            \alpha(g_2)
            (\gamma)
          \scalebox{1.3}{$)$}
        }"
      ]
      \\
      \bullet
      &
      &
      \bullet
      \ar[
        rr,
        "{
          \scalebox{1.3}{$($}
            \gamma, \, \mathrm{e}
          \scalebox{1.3}{$)$}
        }"{description}
      ]
      &&
      \bullet
    \end{tikzcd}
  $$
  By unwinding the definition
  and using Lem. \ref{CrossedHomomorphismsAreSectionsOftheSemidirectProductProjection}
  the claim follows,
  as indicated on the right.
\end{proof}

\begin{remark}[General abstract perspective on non-abelian group 1-cohomology]
  \label{GeneralAbstractPerspectiveOnNonAbelianGroup1Cohomology}
  Prop. \ref{ConjugationGroupoidOfCrossedHomomorphismsIsSectionsOfDeloopedSemidirectProductProjection}
  says that non-abelian group cohomology
  (here in degree 1) is an example of the general notion of
  twisted non-abelian cohomology in \cite[\S 2.2]{SS20OrbifoldCohomology},
  for local coefficient bundle being the
  fibration $B (\Gamma \rtimes G) \xrightarrow{\;} B G$.
\end{remark}

\begin{proposition}[Discrete spaces of crossed-conjugacy classes of crossed homomorphisms]
\label{DiscreteSpacesOfCrossedConjugacyClassesOfCrossedHomomorphisms}
let

\noindent  {(a)} $G$ be a {\it compact} Lie group (e.g. a finite group),

\noindent  {(b)} $\Gamma$ be any Lie group, and

\noindent  \hypertarget{AssumptionThatAlphaIsTrivialOnTheCenter}{(c)}
  $\alpha \,\colon\, G \to \mathrm{Aut}_{\mathrm{Grp}}(\Gamma)$
    be trivial on the center of $G$.

\noindent
Then:

\noindent
{\bf (i)}
the crossed-conjugation quotient space in \eqref{NonAbelianGroup1Cohomology}
is discrete, hence already is the group 1-cohomology as a set:
\vspace{-2mm}
\begin{equation}
  \label{ConjugationQuotientOfCrossedHomomorphismSpaceHappeningToBeDiscrete}
  H^1_{\mathrm{Grp}}(G,\, G \acts \, \Gamma)
  \;\;
  \simeq
  \;\;
  \CrossedHomomorphisms(G ,\, G \acts \, \Gamma)/\sim_{\!\mathrm{ad}}
  \quad
  \in
  \;
  \Sets
    \xhookrightarrow{\; \mathrm{Disc} \;}
  \kTopologicalSpaces
  \,.
\end{equation}

\vspace{-2mm}
\noindent
{\bf (ii)}
The crossed homomorphism space equipped with its
crossed conjugation action
\eqref{ActionOnSpaceOfCrossedHomomorphisms},
decomposes as a disjoint union,
indexed by this group 1-cohomology set
\eqref{ConjugationQuotientOfCrossedHomomorphismSpaceHappeningToBeDiscrete},
of the coset spaces (Ex. \ref{CosetSpacesAsActions})
of $\Gamma$
by the  crossed conjugation stabilizer subgroups $C_\Gamma(\phi) \subset \Gamma$:
\vspace{-2mm}
\begin{equation}
  \label{CrossedHomomorphismSpaceDecomposed}
  \CrossedHomomorphisms(G,\, G \acts \, \Gamma)
  \;\;\;
  \simeq
  \quad
  \underset{
    {[\phi] \in}
    \atop
    \mathclap{
      \scalebox{.7}{$
        H^1_{\mathrm{Grp}}(G,\, G \acts \, \Gamma)
      $}
    }
  }{
    \coprod
  }
  \;\;\;
  \Gamma/C_\Gamma(\phi)
  \;\;\;\;\;
  \in
  \Actions{\Gamma}(\kTopologicalSpaces)
  \,.
\end{equation}
\end{proposition}
\begin{proof}
  For the special case when the action $\alpha$ is entirely trivial,
  the statement is observed in \cite[Rem. 2.2.1]{Rezk14}.
  To deduce from this the more general statement,
  identify crossed homomorphisms $\phi : G \xrightarrow{\;} \Gamma$
  with constrained plain homomorphisms $G \xrightarrow{\;} \Gamma \rtimes G$
  according to Lem. \ref{CrossedHomomorphismsAreSectionsOftheSemidirectProductProjection}.
  Applied to these, the previous version of the statement guarantees
  for nearby $\phi$, $\phi'$ an element $(\gamma, h) \,\in\, \Gamma \rtimes G$,
  such that:
  \vspace{-2mm}
  $$
    \underset{
      g \in G
    }{\forall}
    \;\;\;\;\;\;
    \left(
      \phi'(g),\,(g)
    \right)
    \;\;
    =
    \;\;
    (\gamma,\, h)^{-1}
    \cdot
    \left(
      \phi(g),\, g
    \right)
    \cdot
    (\gamma,\, h)
    \,.
  $$

  \vspace{-2mm}
  \noindent
  In the second component this implies that $h \,\in\, C(G)$
  is in the center of $G$.
  Therefore, further conjugation of this equation with
  $(\mathrm{e}, h) \,\in\, \Gamma \rtimes G$,
  and using assumption \hyperlink{AssumptionThatAlphaIsTrivialOnTheCenter}{(iii)},
  yields
  \vspace{-2mm}
  $$
    \underset{
      g \in G
    }{\forall}
    \;\;\;\;\;\;
    \left(
      \phi'(g),\,(g)
    \right)
    \;\;
    =
    \;\;
    (\gamma,\, \mathrm{e})^{-1}
    \cdot
    \left(
      \phi(g),\, g
    \right)
    \cdot
    (\gamma,\, \mathrm{e})
    \,,
  $$

  \vspace{-2mm}
\noindent
  which implies
  the claim as in \eqref{CrossedConjugationAsPlainConjugationInSemidirectProductGroup}.
\end{proof}

\begin{proposition}[Weyl group action on group 1-cohomology of subgroup]
\label{WeylGroupActionOnGroup1CohomologyOfSubgroup}
  For $G \,\in\, \Groups(\kTopologicalSpaces)$
  and $G \acts \, \Gamma \,\in\, \Groups ( \Actions{G}(\kTopologicalSpaces))$,
  let $H \,\subset\, G$ be a subgroup, with
  $H \acts \,  \Gamma \,\in\, \Groups (\Actions{H}(\kTopologicalSpaces))$
  denoting the restricted action (Ex. \ref{RestrictedActions}).
  Then:

  \noindent
  {\bf (i)}
  The assignment
  \vspace{-3mm}
  \begin{equation}
    \label{NormalizerGroupActionOnCrossedHomomorphismsOutOfSubgroup}
    \begin{tikzcd}[column sep=4pt, row sep=-4pt]
      N(H)
        \,\times\,
      \CrossedHomomorphisms(H,\, H \acts \, \Gamma)
      \ar[rr]
      &&
      \CrossedHomomorphisms(H,\, H \acts \, \Gamma)
      \\
   \scalebox{0.8}{$    ( n,\,\phi ) $}
      &\longmapsto&
      \hspace{-2cm}
  \scalebox{0.8}{$     \phi_n
      \mathrlap{
        \;\colon\;
        h
        \,\mapsto\,
        \alpha(n)
        (
          n^{-1} \cdot h \cdot n
        )
      }
      $}
    \end{tikzcd}
  \end{equation}

  \vspace{-2mm}
\noindent  is a continuous action of the normalizer subgroup $N(H)$
  (Ntn. \ref{GActionOnTopologicalSpaces})
  on the space \eqref{SpaceOfCrossedHomomorphisms} of crossed homomorphisms
  out of $H$.

  \noindent
  {\bf (ii)}
  This descends to the quotient
  by crossed conjugations, hence to the
  non-abelian group 1-cohomology set
  \eqref{NonAbelianGroup1Cohomology},

  \noindent
  {\bf (iii)} where it passes through an action of the Weyl group
  $W(H) \,\coloneqq\, N(H)/H$
  (Ntn. \ref{GActionOnTopologicalSpaces}):
  \begin{equation}
    \label{WeylGroupActionOnSetOfGroup1CohomologyOfSubgroup}
    G
      \acts \,
    \Gamma
      \,\in\,
    \Groups\left( \Actions{H}(\kTopologicalSpaces) \right)
    \;\;\;\;\;\;\;
    \vdash
    \;\;\;\;\;\;\;
    \underset{H \subset G}{\forall}
    \;\;\;\;
    W(H)
      \;\acts\;
    H^1_{\mathrm{Grp}}
    (
      H,\, H \acts \, \Gamma
    )
    \;\;\;
    \in
    \;
    \Actions{ W(H) }(\Sets)
    \,.
  \end{equation}
\end{proposition}
\begin{proof}
  The first two statements are immediate in the equivalent
  incarnation of crossed homomorphisms $\phi \,\colon\, H \xrightarrow{\;} \Gamma$
  as plain homomorphisms of the form
  $( \phi(-),\, (-)) \,\colon\, G \xrightarrow{\;} \Gamma \rtimes G$
  (Lem. \ref{CrossedHomomorphismsAreSectionsOftheSemidirectProductProjection}),
  observing that, under this identification,
  the assignment \eqref{NormalizerGroupActionOnCrossedHomomorphismsOutOfSubgroup}
  is just the ``conjugation action by the adjoint action'':
    \vspace{-2mm}
  $$
    \left(
      \phi_n(-),\, (-)
    \right)
    \;\;
    =
    \;\;
    \left(
      n,\, \mathrm{e}
    \right)
    \cdot
    \left(
      \phi_n
      \left(
        n^{-1} \cdot (-) \cdot n
      \right)
      ,\,
      \left(
        n^{-1} \cdot (-) \cdot n
      \right)
    \right)
    \cdot
    \left(
      n,\, \mathrm{e}
    \right)
    \,.
  $$

    \vspace{-2mm}
\noindent
  For the last statement we need to exhibit a
  crossed conjugation between $\phi_n$ and $\phi$ whenever
  $n \,\in\, H \,\subset\, N(H)$. The
  following direct computation shows that this is given by
  $\phi(n)$ (which is defined for such $n$):
    \vspace{-2mm}
  $$
    \label{CrossedConjugationToImageOfHActionViaActionOfNormalizerSubgroup}
    \begin{array}{lll}
      \phi_n(h)
      & \;=\;
      \alpha(n)
      \big(
        \phi ( n^{-1} \cdot h \cdot n )
      \big)
      &
      \mbox{\small by definition \eqref{NormalizerGroupActionOnCrossedHomomorphismsOutOfSubgroup}}
      \\
      & \;=\;
      \alpha(n)
      \Big(
        \phi( n^{-1})
        \cdot
        \alpha( n^{-1})
        \big(
          \phi(h)
          \cdot
          \alpha(h)
          (
            \phi(n)
          )
        \big)
      \Big)
      &
      \mbox{\small by crossed hom. property \eqref{GCrossedHomomorphismProperty}}
      \\
      & \;=\;
      \phi(n^{-1})
       \cdot
      \phi(h)
       \cdot
      \alpha(h)
      \left(
        \phi(n)
      \right)
      &
      \mbox{\small by action property of $\alpha$}
      \mathrlap{\,.}
    \end{array}
  $$

    \vspace{-6mm}
\end{proof}


\noindent
{\bf Simplicial topological spaces.}

\begin{notation}[Simplicial topological spaces]
  \label{CategoryOfSimplicialTopologicalSpaces}
  We denote by
  \vspace{-2mm}
  \begin{equation}
    \label{SimplicialTopologicalSpaces}
    \SimplicialTopologicalSpaces
    \;\coloneqq\;
    \Functors(\Delta^{\mathrm{op}}, \kTopologicalSpaces)
  \end{equation}

  \vspace{-2mm}
\noindent  the category of simplicial topological spaces,
  hence of simplicial objects internal to (Ntn. \ref{Internalization})
  the category
  of compactly generated weak Hausdorff spaces \eqref{CategoryOfTopologicalSpaces}.
\end{notation}

\begin{definition}[Homotopy Kan fibrations {\cite[Def. 3]{Lurie11}\cite[p. 2]{MazelGee14}}]
  \label{HomotopyKanFibrations}
  We say that a morphism
  of simplicial topological
  spaces \eqref{SimplicialTopologicalSpaces}
  $f_\bullet \;\colon\; \TopologicalSpace_\bullet \xrightarrow{\;} \mathrm{Y}_\bullet$
  is a {\it homotopy Kan fibrations}, to be denoted
  \vspace{-2mm}
    $$
      f_\bullet \;\in\; \HomotopyKanFibrations
      \;\subset\;
      \mathrm{Mor}
      (
        \SimplicialTopologicalSpaces
      )
      \,,
    $$

     \vspace{-2mm}
    \noindent
    if for all positive  $n \in \NaturalNumbers_+$
    every solid homotopy-commutative diagram as
    follows has a dashed lift up to homotopy, as shown:
    \vspace{-2mm}
    \begin{equation}
      \label{HomotopyKanFibrationProperty}
      \begin{tikzcd}
        \Lambda^n_k
        \ar[
          rr,
          hook,
          "{\ }"{below, name=s1, pos=.7}
        ]
        \ar[d]
        &&
        \TopologicalSpace_\bullet
        \ar[
          d,
          "{f_\bullet}"{right}
        ]
        \\
        \Delta^n
        \ar[
          rr,
          "{\ }"{above, name=t2, pos=.3}
        ]
        \ar[
          urr,
          dashed,
          "{\ }"{above, name=t1},
          "{\ }"{below, name=s2}
        ]
        &&
        \mathrm{Y}_{\bullet}
        \ar[
          from=s1,
          to=t1,
          Rightarrow,
          dashed
        ]
        \ar[
          from=s2,
          to=t2,
          Rightarrow,
          dashed
        ]
      \end{tikzcd}
      {\phantom{AAAAA}}
      \Leftrightarrow
      {\phantom{AAAAA}}
      \mbox{$
      \pi_0
      \Big(\!
      X(\Delta^n)
      \xrightarrow{\;\;}
      X(\Lambda^n_k)
        \underset{
          Y(\Lambda^n_k)
        }{\times^h}
      Y(\Delta^n)
      \!\Big)$
      is surjective,
      }
    \end{equation}

  \vspace{-3mm}
  \noindent
  where the morphism on the far left is the $(n,k)$-horn inclusion
  of simplicial sets regarded as degreewise discrete simplicial topological spaces.
  The equivalent condition on the right of \eqref{HomotopyKanFibrationProperty}
  simplifies when horn filling in $\mathrm{Y}_\bullet$ is a Serre fibration,
  in which case the homotopy fiber product in \eqref{HomotopyKanFibrationProperty}
  is given by the plain fiber product, so that:
  \vspace{-2mm}
  \begin{equation}
    \label{HomotopyKanFibrationPropertyInTheCaseThatHornFillingInCodomainIsSerreFibration}
    \hspace{-4mm}
    \bigg(
    \underset{
      { n \,\in\, \mathbb{N}_+ }
      \atop
      { 0 \leq k \leq n }
    }{\forall}
    \TopologicalSpace(\Delta^n)
    \xrightarrow{\!\in \, \SerreFibrations \!}
    \mathrm{Y}(\Lambda^n_k)
    \bigg)
    \;
       \Rightarrow
    \;
    \left(\!\!
      \big(
      f_\bullet
        \;\in\;
      \HomotopyKanFibrations
      \big)
    \Leftrightarrow
    \bigg(
    \underset{
      { n \,\in\, \mathbb{N}_+ }
      \atop
      { 0 \leq k \leq n }
    }{\forall}
    \pi_0
    \big(
    \TopologicalSpace(\Delta^n)
    \xrightarrow{\;}
    \TopologicalSpace(\Lambda^n_k)
    \underset{
      \mathrm{Y}(\Lambda^n_k)
    }{\times}
    \mathrm{Y}(\Delta^n)
    \big)
    \,
    \mbox{surj.}
    \!\bigg)
    \!\!\!\right)
    .
  \end{equation}
\end{definition}

\begin{notation}[Nerve of topological groupoids]
  \label{NerveOfTopologicalGroupoids}
  We denote the {\it nerve} functor from
  topological groupoids (Ntn. \ref{TopologicalGroupoids})
  to simplicial topological spaces \eqref{CategoryOfSimplicialTopologicalSpaces}
  by
  \vspace{-2mm}
  \begin{equation}
    \label{SimplicialTopologicalNerveOfTopologicalGroupoids}
    \begin{tikzcd}[row sep=-5pt]
      \mathllap{
        \SimplicialNerve
        \;:\;
      }
      \Groupoids(\kTopologicalSpaces)
      \ar[rr]
      &&
      \SimplicialTopologicalSpaces
      \\
  \scalebox{0.8}{$    \big(
        \TopologicalSpace_1
        \underoverset
          {t}
          {s}
          {\rightrightarrows}
        \TopologicalSpace_1
      \big)
      $}
      &\longmapsto&
     \scalebox{0.8}{$   \big(
      [n]
        \,\mapsto\,
      \underset{
        \mbox{\tiny \rm $n$ factors}
      }{
      \underbrace{
      (\TopologicalSpace_1)
        \, {}_{t}\!\! \underset{\TopologicalSpace_0}{\times_{s}}
        \cdots
        \, {}_{t}\!\! \underset{\TopologicalSpace_0}{\times_{s}}
      \TopologicalSpace_1
      }
      }
      \big)
      $}
    \end{tikzcd}
  \end{equation}
\end{notation}

\vspace{-3mm}
\begin{example}[Nerves of action groupoids have horn filling by Serre fibrations]
  \label{NervesOfActionGroupoids}
  For $\Gamma \,\in\, \Groups(\kTopologicalSpaces)$
  and $\TopologicalSpace \, \rightacts \Gamma
  \,\in\, \Actions{\Gamma^{\mathrm{op}}}(\kTopologicalSpaces)$,
  the nerve \eqref{SimplicialTopologicalNerveOfTopologicalGroupoids}
  of the corresponding action groupoid (Ex. \ref{TopologicalActionGroupoid})
  is of the form
  \vspace{-1mm}
  $$
    N
    (
      \TopologicalSpace \times \Gamma \rightrightarrows \Gamma
    )_n
    \;\;\;
      \simeq
    \;\;\;
    \TopologicalSpace \times \Gamma^{\times_n}
    \;\;\;
    \in
    \;
    \kTopologicalSpaces
    \,.
  $$

  \vspace{-2mm}
  \noindent
  from which it is readily seen that the
  horn filler maps are all Serre fibrations
  (for $n = 1$ by Lem.  \ref{LocallyTrivialBundlesAreSerreFibrations}):
  \vspace{-2mm}
  \begin{equation}
    \label{HornFillerMapsInNerveOfActionGroupoidAreSerreFibrations}
     \hspace{-1cm}
    \begin{tikzcd}[row sep=2pt]
          &
      N
      (
        \TopologicalSpace \times \Gamma \rightrightarrows \Gamma
      )(\Delta^1)
      \;\simeq\;
      \TopologicalSpace \times \Gamma
      \;\;
      \ar[
        r,
        "{ \mathrm{pr}_1 }"{above},
        "{ \in \SerreFibrations }"{below}
      ]
      &
      \;\;
      \TopologicalSpace
      \phantom{ \times \Gamma }
      \;\simeq\;
      N
      (
        \TopologicalSpace \times \Gamma \rightrightarrows \Gamma
      )(\Lambda^1_1)
      \\
      \underset{
        { n \geq 2 }
        \atop
        { 0 \leq k \leq n }
      }{\forall}
      &[-24pt]
      \hspace{-1cm}
      N
      (
        \TopologicalSpace \times \Gamma \rightrightarrows \Gamma
      )(\Delta^n)
      \;\simeq\;
      \TopologicalSpace \times \Gamma^{\times_n}
     \;\;
      \ar[
        r,
        "{ \sim }"{above},
        "{ \in \SerreFibrations }"{below}
      ]
      &
      \;\;
      \TopologicalSpace \times \Gamma^{\times_{n}}
      \;\simeq\;
      N
      (
        \TopologicalSpace \times \Gamma \rightrightarrows \Gamma
      )(\Lambda^n_k)
      \mathrlap{\,.}
    \end{tikzcd}
  \end{equation}
\end{example}

\vspace{-3mm}
\begin{example}[Nerves of morphisms of delooping groupoids are homomotopy Kan fibrations]
  \label{NervesOfMorphismsOfOfDeloopingGroupoids}
  Let $\Gamma_1 \xrightarrow{\phi} \Gamma_2 \;\; \in \Groups(\kTopologicalSpaces)$
  be a homomorphism of topological groups
  which is surjective on connected components of underlying topological spaces:
  \vspace{-2mm}
  $$
    \begin{tikzcd}[column sep=large]
      \pi_0(\Gamma_1)
      \ar[
        r,
        ->>,
        "{ \pi_0(\phi) }"
      ]
      &
      \pi_0(\Gamma_2)
      \mathrlap{\,.}
    \end{tikzcd}
  $$

    \vspace{-1mm}
\noindent
  Then its
  image under taking nerves \eqref{SimplicialTopologicalNerveOfTopologicalGroupoids}
  of delooping groupoids (Ex. \ref{TopologicalDeloopingGroupoid})
  is a homotopy Kan fibration (Def. \ref{HomotopyKanFibrations}):
  \vspace{-2mm}
  \begin{equation}
    \label{AssumptionThatTopologicalGroupHomomorphismIsSurjectiveOnConnectedComponents}
    \begin{tikzcd}
      N
      (
        \Gamma_1 \rightrightarrows \ast
      )_\bullet
      \ar[
        rr,
        "{ N( \phi \rightrightarrows \ast )_\bullet }"{above},
        "{ \in \; \HomotopyKanFibrations } "{below}
      ]
      &&
      N
      (
        \Gamma_2 \rightrightarrows \ast
      )_\bullet
      \mathrlap{\,.}
    \end{tikzcd}
  \end{equation}
\end{example}
\begin{proof}
  By \eqref{HornFillerMapsInNerveOfActionGroupoidAreSerreFibrations}
  in Ex. \ref{NervesOfActionGroupoids}
  the condition to be checked is
  that on the right of
  \eqref{HomotopyKanFibrationPropertyInTheCaseThatHornFillingInCodomainIsSerreFibration},
  which in degree 1 is satisfied by assumption, since here the comparison map is
  $\phi$ itself:
    \vspace{-2mm}
  $$
    \begin{tikzcd}
      N
      (
        \Gamma_1 \rightrightarrows \ast
      )(\Delta^1)
      \;\simeq\;
      \Gamma_1
      \;\;
      \ar[
        r,
        "\phi"
      ]
      &
      \;\;
      \Gamma_2
      \;\simeq\;
      N
      (
        \Gamma_1 \rightrightarrows \ast
      )(\Delta^0)
      \quad \;\;
      \underset{
        \mathclap{
          \raisebox{-3pt}{
          \scalebox{.7}{$
            N
            (
              \Gamma_2 \rightrightarrows \ast
            )(\Delta^0)
          $}}
        }
      }{\times}
      \quad \;\;\;
      (
        \Gamma_2 \rightrightarrows \ast
      )(\Delta^1)
      \mathrlap{\,,}
    \end{tikzcd}
  $$

  \vspace{-3mm}
  \noindent
  while in higher degrees the condition is trivial, as here
  the comparison map is an isomorphism
  (since the pullbacks of the isomorphisms in \eqref{HornFillerMapsInNerveOfActionGroupoidAreSerreFibrations}
  are isomorphisms):
    \vspace{-2mm}
  $$
    \underset{
      { n \geq 2 }
      \atop
      { 0 \leq k \leq n }
    }{\forall}
    \begin{tikzcd}
      N
      (
        \Gamma_1 \rightrightarrows \ast
      )(\Delta^n)
      \;\simeq\;
      \Gamma_1^{\times_n}
    \;\;  \ar[
        r,
        "\sim"
      ]
      &
      \;\;
      \Gamma_1^{\times_n}
      \;\simeq\;
      N
      (
        \Gamma_1 \rightrightarrows \ast
      )(\Lambda^n_k)
      \quad \;\;
      \underset{
        \mathclap{
          \raisebox{-3pt}{
          \scalebox{.7}{$
          N
          (
            \Gamma_2 \rightrightarrows \ast
          )(\Lambda^n_k)
          $}}
        }
      }{\times}
      \quad \;\;\;
      (
        \Gamma_2 \rightrightarrows \ast
      )(\Delta^n)
      \mathrlap{\,.}
    \end{tikzcd}
  $$

  \vspace{-8mm}
\end{proof}

\begin{example}[Base morphisms out of nerves of action groupoids are homotopy Kan fibrations]
  \label{BaseMorphismsOutOfNervesOfActionGroupoidsAreHomotopyKanFibrations}
  $\,$

  \noindent
  For $\Gamma \,\in\, \Groups(\kTopologicalSpaces)$
  and $\TopologicalSpace \, \rightacts \Gamma \,\in\, \Actions{\Gamma^{\mathrm{op}}}(\kTopologicalSpaces)$,
  the image under taking nerves \eqref{SimplicialTopologicalNerveOfTopologicalGroupoids}
  of the canonical morphism
  from the action groupoid (Ex. \ref{TopologicalActionGroupoid})
  to the delooping groupoid (Ex. \ref{TopologicalDeloopingGroupoid})
  is a homotopy Kan fibration (Def. \ref{HomotopyKanFibrations}):
  \vspace{-3mm}
  $$
    \begin{tikzcd}
      N
      (
        \TopologicalSpace \times \Gamma
          \rightrightarrows
        \TopologicalSpace
      )
      \ar[
        rr,
        "{
        }"{above},
        "{
          \in \, \HomotopyKanFibrations
        }"{below}
      ]
      &&
      N
      (
        \Gamma
        \rightrightarrows
        \ast
      )\;.
    \end{tikzcd}
  $$
\end{example}
\begin{proof}
  Again, by \eqref{HornFillerMapsInNerveOfActionGroupoidAreSerreFibrations}
  in Ex. \ref{NervesOfActionGroupoids}
  the condition to be checked is
  that on the right of
  \eqref{HomotopyKanFibrationPropertyInTheCaseThatHornFillingInCodomainIsSerreFibration}.
  One readily sees that the comparison maps are in fact isomorphisms,
  hence in particular surjective on connected components.
\end{proof}


\noindent
{\bf Topological realization of topological groupoids.}

\begin{notation}[Topological realization functors]
  \label{TopologicalRealizationFunctors}
  We denote

  \noindent
  {{\bf (i)}}
  the cosimplicial topological space of
  {\it standard topological simplices} by
  \vspace{-2mm}
  $$
    \begin{tikzcd}[row sep=-5pt, column sep=17pt]
      \mathllap{
        \Delta^\bullet_{\mathrm{top}}
        \;:\;
      }
      \Delta
      \ar[rr]
      &&
      \kTopologicalSpaces
      \\
      \scalebox{0.8}{${[n]}$}
      &\longmapsto&
      \scalebox{0.8}{$
      \big\{
        \vec x \,\in\, \mathbb{R}^n
        \,\big\vert\,
        \sum_i x_i \,=\, 1
        \;\mbox{and}\;
        \forall_i (0 \leq x_i \leq 1)
      \big\}
      $}
    \end{tikzcd}
    \hspace{-2.3cm}
  $$

  \vspace{-2mm}
  \noindent
  {\bf (ii)} the {\it topological realization of simplicial topological spaces}
  (\cite[\S 6]{MacLane70}\cite[\S 11]{May72})
  by
  \vspace{-2mm}
  \begin{equation}
    \label{TopologicalRealizationOfSimplicialTopologicalSpaces}
    \begin{tikzcd}[row sep=-5pt]
      {\vert
       -
      \vert}
      \;:\;
      \SimplicialTopologicalSpaces
      \ar[rr]
      &&
      \kTopologicalSpaces
      \\
  \scalebox{0.8}{$      \TopologicalSpace_\bullet $}
      &\longmapsto&
    \scalebox{0.8}{$    \overset{[n] \in \Delta}{\int}
      \TopologicalSpace_n \times \Delta^n_{\mathrm{top}}
      $}
    \end{tikzcd}
  \end{equation}

  \vspace{-2mm}
  \noindent
  {\bf (iii)}
  the {\it topological realization of topological groupoids}
  by the same symbol, leaving
  the nerve \eqref{SimplicialTopologicalNerveOfTopologicalGroupoids}
  notationally implicit:
  \vspace{-2mm}
  \begin{equation}
    \label{TopologicalRealizationOfTopologicalGroupoids}
    \vert
    -
    \vert
    \;\;:\;
    \begin{tikzcd}
      \Groupoids(\kTopologicalSpaces) \;.
      \ar[
        r,
        "{N}"
      ]
      &
      \SimplicialTopologicalSpaces
      \ar[
        r,
        "{{\vert - \vert}}"
      ]
      &
      \kTopologicalSpaces \;.
    \end{tikzcd}
  \end{equation}
\end{notation}
\begin{remark}[Classifying spaces and realization]
  The realization \eqref{TopologicalRealizationOfTopologicalGroupoids}
  of topological groupoids
  has traditionally been denoted $B(-)$.
  This may seem like a good idea if one already notationally conflates
  topological groups $\Gamma$
  with their delooping groupoids $\Gamma \rightrightarrows \ast$
  (Ex. \ref{TopologicalDeloopingGroupoid}), because then $B \Gamma$
  denotes
  the intended classifying space.
  But, since a group is crucially {not} the same as its delooping groupoid,
  we give the latter its own symbol
  (following \cite{NSS12a}\cite{SS20OrbifoldCohomology}):
  $\mathbf{B}\Gamma \,\coloneqq\, ( \Gamma \rightrightarrows \ast)$,
  and denote topological realization of topological groupoids
  by the same symbol $\left\vert-\right\vert$ as that of topological spaces
  (see corresponding discussion in \cite[p. 4]{GuillouMayMerling17}).
  Therefore, in our notation,
  the classifying space \eqref{QuotientCoprojectionOfUniversalPrincipalBundle}
  of a topological group $\Gamma$
  (rather: of $\Gamma$-principal bundles, when $\Gamma$ is well-behaved)
  arises as
  \vspace{-1mm}
  $$
    B \Gamma
      \;=\;
    \vert \mathbf{B} \Gamma \vert
    \,.
  $$

  \vspace{-2mm}
  \noindent
  This typesetting serves to express that the classifying space $B G$
  is but a shadow (namely: the {\it shape})
  of the richer {\it moduli stack} $\mathbf{B}G$
  (see \cite{SS20OrbifoldCohomology}).
\end{remark}

\begin{example}[Topological realization of constant groupoids]
  \label{TopologicalRealizationOfConstantGroupoids}
  The topological realization \eqref{TopologicalRealizationOfTopologicalGroupoids}
  of constant topological groupoids (Ex. \ref{TopologicalSpacesAsTopologicalGroupoids})
  is the topological realization \eqref{TopologicalRealizationOfSimplicialTopologicalSpaces}
  of the corresponding constant simplicial space, and hence is
  homeomorphic to the underlying topological space:
  \vspace{-2mm}
  \begin{equation}
    \label{TopologicalRealizationOfConstantGroupoidIsomorphicToSPaceOfObjects}
    \vert
      \ConstantGroupoid(\TopologicalSpace)
    \vert
    \;\simeq\;
    \TopologicalSpace
    \;\;\;
    \in
    \;
    \kTopologicalSpaces
    \,.
  \end{equation}
\end{example}

\begin{lemma}[Topological realization preserves finite limits]
  \label{TopologicalRealizationPreservesFiniteLimits}
  The topological realization functors (Ntn. \ref{TopologicalRealizationFunctors}),
  both
  {(i)} of simplicial topological spaces \eqref{TopologicalRealizationOfSimplicialTopologicalSpaces}
  and {(ii)} of topological groupoids \eqref{TopologicalRealizationOfTopologicalGroupoids},
  preserve finite limits:
      \vspace{-2mm}
  $$
    \left\vert-\right\vert
    \;:\;
    \begin{tikzcd}
    \Groupoids(\kTopologicalSpaces)
    \ar[
      r,
      "{\mathrm{lex}}"{above}
    ]
    &
    \kTopologicalSpaces
    \,.
    \end{tikzcd}
  $$
\end{lemma}
\begin{proof}
  The second statement follows from the first by the fact that the
  nerve \eqref{SimplicialTopologicalNerveOfTopologicalGroupoids}
  is a right adjoint (\cite[\S 3]{Kan58}) and hence preserves all limits.

  For the first statement it is sufficient
  (e.g. \cite[Prop. 2.8.2]{Borceux94I}) to see that topological realization
  {(a)} preserves the terminal object (here: the point) and
  {(b)} preserves all pullbacks.
  Here (a) follows by Ex. \ref{TopologicalRealizationOfConstantGroupoids}
  and (b) is the statement of \cite[Cor. 11.6]{May72}
  (making crucial use of compact generation \eqref{CategoryOfTopologicalSpaces},
  see \cite[p. 1]{May72}).
\end{proof}

\begin{lemma}[Topological realization of equivalence is homotopy equivalence]
  \label{TopologicalRealizationOfEquivalenceOfGroupoidsIsHomotopyEquivalence}
  An equivalence \eqref{EquivalenceOfTopologicalGroupoids} of topological groipoids
  induces a homotopy equivalence between their topological realizations
  \eqref{TopologicalRealizationOfTopologicalGroupoids}:
      \vspace{-1mm}
  $$
    (\TopologicalSpace_1 \rightrightarrows \TopologicalSpace_0)
    \;\underset{\mathrm{htpy}}{\simeq}\;
    (\mathrm{Y}_1 \rightrightarrows \mathrm{Y}_0)
    \;\;\;\;\;\;\;\;\;\;\;
    \Rightarrow
    \;\;\;\;\;\;\;\;\;\;\;
    \left\vert(\TopologicalSpace_1 \rightrightarrows \TopologicalSpace_0)\right\vert
    \;\underset{\mathrm{htpy}}{\simeq}\;
    \left\vert(\mathrm{Y}_1 \rightrightarrows \mathrm{Y}_0)\right\vert
    \,.
  $$
\end{lemma}
\begin{proof}
  Observe that a 2-morphism
  $F \xRightarrow{\eta} F'$
  \eqref{NaturalitySquareForTopologicalGroupoids}
  of topological groipoids is equivalently
  a morphism $\widetilde \eta$ out of the product of the domain groupoid with the
  codiscrete groupoid (Ex. \ref{TopologicalPairGroupoid}) on the 2-element set:
      \vspace{-2mm}
  $$
    \begin{tikzcd}[row sep=-5pt, column sep=2pt]
      \mathllap{
        \widetilde \eta
        \;:\;
      }
      (\TopologicalSpace_1 \rightrightarrows \TopologicalSpace_0)
        \times
      \CodiscreteGroupoid(\{\ast, \ast' \})
      \ar[
        rr
      ]
      &&
      (\mathrm{Y}_1 \rightrightarrows \mathrm{Y}_0)
      \\
 \scalebox{0.8}{$
      (x \xrightarrow{\gamma} x') \times \ast
      $}
      &\longmapsto& \quad
    \scalebox{0.8}{$   F{\phantom{'}}(x{\phantom{'}}) \xrightarrow{F{\phantom{'}}(\gamma)} F{\phantom{'}}(x')
    $}
      \\
   \scalebox{0.8}{$    (x \xrightarrow{\gamma} x') \times \ast'
   $}
      &\longmapsto& \quad
    \scalebox{0.8}{$   F'(x{\phantom{'}}) \xrightarrow{F'(\gamma)} F'(x')
    $}
      \\
   \scalebox{0.8}{$    x \times (\ast \xrightarrow{\;} \ast')
   $}
      &\longmapsto& \quad
 \scalebox{0.8}{$      F{\phantom{'}}(x) \xrightarrow{ \;\eta(x)\; } F'(x{\phantom{'}})
 $}
      \mathrlap{\,.}
    \end{tikzcd}
  $$

  \vspace{-2mm}
  \noindent
  The topological realization of this morphism is,
  by Lem. \ref{TopologicalRealizationPreservesFiniteLimits},
  of the following form
  \vspace{-2mm}
  $$
    \begin{tikzcd}[column sep=large]
    \vert
      (\TopologicalSpace_1 \rightrightarrows \TopologicalSpace_0)
    \vert
    \ar[
      d,
      "{
        \mathrm{id}(-) \times \ast
      }"{left}
    ]
    \ar[
      drr,
      "{
        \left\vert
          F
        \right\vert
      }"{above}
    ]
    \\
    \vert
      (\TopologicalSpace_1 \rightrightarrows \TopologicalSpace_0)
    \vert
    \times
    \left\vert \CodiscreteGroupoid(\{\ast,\ast'\}) \right\vert
    \ar[
      rr,
      "{
        \left\vert
          \widetilde \eta
        \right\vert
      }"{description}
    ]
    &&
    \vert
      (\mathrm{Y}_1 \rightrightarrows \mathrm{Y}_0)
    \vert \;.
    \\
    \vert
      (\TopologicalSpace_1 \rightrightarrows \TopologicalSpace_0)
    \vert
    \ar[
      u,
      "{
        \mathrm{id}(-) \times \ast'
      }"{left}
    ]
    \ar[
      urr,
      "{
        \left\vert
          F'
        \right\vert
      }"{below}
    ]
    \end{tikzcd}
  $$

    \vspace{-2mm}
\noindent
  Observing that
      \vspace{-2mm}
  $$
   \begin{tikzcd}
      \{\ast, \ast'\}
      \ar[
        r,
        hook
      ]
      \ar[
        rr,
        rounded corners,
        to path={
             -- ([yshift=+6pt]\tikztostart.north)
             --node[above]{\scalebox{.7}{$
               \nabla_\ast
             $}} ([yshift=+10pt]\tikztotarget.north)
             -- (\tikztotarget.north)}
      ]
      &
      \left\vert \CodiscreteGroupoid (\{\ast,\ast'\}) \right\vert\
      \ar[r]
      &
      \ast
    \end{tikzcd}
  $$
  is a cylinder object for the point in the classical model structure
  on topological spaces, this exhibits a homotopy of continuous functions
  \vspace{-2mm}
  $$
    \left\vert F \right\vert
    \Rightarrow
    \vert F' \vert
    \,.
  $$
  Applying this to the two homotopies involved in the data \eqref{EquivalenceOfTopologicalGroupoids}
  constituting an equivalence of topological groupoids,
  this yields a homotopy equivalence of topological spaces as claimed.
\end{proof}

\medskip

\noindent
{\bf Topological 2-groups.}

\begin{notation}[Topological strict 2-groups]
  \label{StrictTwoGroups}
  An internal group (Ntn. \ref{Internalization})
  in topological groupoids (Ntn. \ref{TopologicalGroupoids})
  is also called a topological {\it strict 2-group}:
    \vspace{-1mm}
  \begin{equation}
    \label{CategoryOfStrict2GroupsInTopologicalSpaces}
    \TwoGroups(\kTopologicalSpaces)
    \;\coloneqq\;
    \Groups
    \left(
      \Groupoids(\kTopologicalSpaces)
    \right)
    \,.
  \end{equation}

  \vspace{-2mm}
\noindent  By Lem. \ref{TopologicalRealizationPreservesFiniteLimits},
  topological realization \eqref{TopologicalRealizationOfTopologicalGroupoids}
  induces \eqref{FunctorOnStructuresInducedFromLexFunctor}
  a functor from topological 2-groups to topological groups
\vspace{-2mm}
  \begin{equation}
    \label{TopologicalRealizationOfTopological2Groups}
    \TwoGroups(\kTopologicalSpaces)
    \;=\;
    \Groups\left( \Groupoids(\kTopologicalSpaces) \right)
      \xrightarrow{ \;\; \vert-\vert \;\; }
    \Groups(\kTopologicalSpaces)
    \mathrlap{\,.}
  \end{equation}
\end{notation}

\begin{example}[Topological groups as topological 2-groups]
  \label{TopologicalGroupsAsTopologicalTwoGroups}
  For $\Gamma \,\in\, \Groups(\kTopologicalSpaces)$,
  its constant groupoid (Ex. \ref{TopologicalSpacesAsTopologicalGroupoids})
  carries the structure of a topological strict 2-group (Ntn. \ref{StrictTwoGroups})
\vspace{-1mm}
  $$
    (
      \Gamma
      \rightrightarrows
      \Gamma
    )
    \;\;\;
    \in
    \;
    \TwoGroups(\kTopologicalSpaces)
  $$

  \vspace{-1mm}
\noindent  with group structure given degreewise by that of $\widehat \Gamma$.
\end{example}

\begin{example}[Strict 2-groups delooping abelian groups]
  \label{Strict2GroupsDeloopingAbelianGroups}
  For $A \,\in\, \AbelianGroups(\kTopologicalSpaces)$
  an abelian topological group, its delooping groupoid (Ex. \ref{TopologicalDeloopingGroupoid})
  carries the structure of a topological strict 2-group
  (Ntn. \ref{StrictTwoGroups})
  \vspace{-1mm}
  $$
    (
      A
      \rightrightarrows
      1
    )
    \;\;\;
    \in
    \;
    \TwoGroups(\kTopologicalSpaces)
  $$

  \vspace{-2mm}
\noindent
  with group structure given degreewise by that of $A$.
  The topological group which is the topological realization
  \eqref{TopologicalRealizationOfTopological2Groups} of this 2-group
  \vspace{-1mm}
  $$
    B A
    \;\coloneqq\;
    \vert
      A \rightrightarrows \ast
    \vert
    \;\;\;
    \in
    \;
    \Groups(\kTopologicalSpaces)
  $$
  has as underlying space the Milgram classifying space
  \eqref{QuotientCoprojectionOfUniversalPrincipalBundle} of $A$.
\end{example}

\begin{example}[Long exact sequence of 2-groups from central extension of groups]
  \label{Strict2GroupsFromCentralExtensions}
  Given a central extension of topological groups
  \vspace{-2mm}
  \begin{equation}
    \label{CentralExtensionOfTopologicalGroups}
    \begin{tikzcd}[row sep=6pt, column sep={between origins, 30pt}]
      1
      \ar[r]
      &
      A
      \ar[
        rr,
        hook,
        "{i}"
      ]
      \ar[
        dr,
        hook
      ]
      &&
      \widehat{\Gamma}
      \ar[
        rr,
        ->>,
        "{p}"
      ]
      &&
      \Gamma
      \ar[r]
      &
      1
      \\
      & &
      Z(\widehat \Gamma)
      \ar[
        ur,
        hook
      ]
    \end{tikzcd}
  \end{equation}

  \vspace{-2mm}
\noindent  the action groupoid (Ex. \ref{TopologicalActionGroupoid})
  of the canonical $A$ action on $\widehat \Gamma$ on
  becomes a strict 2-group (Ntn. \ref{StrictTwoGroups})
  under the direct product group structure:
  \vspace{-1mm}
  $$
   (\,
     \widehat \Gamma \times A
     \rightrightarrows
     \widehat \Gamma
   )
   \;\;\;
   \in
   \;
   \TwoGroups(\kTopologicalSpaces)
   \,.
  $$
  This sits in the following diagram of 2-groups
    \vspace{-5mm}
  \begin{equation}
  \label{LongDuagramOf2GroupFromCentralExtensionOfGroups}
  \hspace{-1cm}
  \begin{tikzcd}
    &&
    &&
    (
      \widehat \Gamma \times A
        \rightrightarrows
      \widehat \Gamma
    )
    \ar[
      rr,
      "{ (\mathrm{pr}_2 \,\rightrightarrows\, 1) }"
    ]
    \ar[
      d,
      "{ ( ( p\circ \mathrm{pr}_1) \,\rightrightarrows\, p) }"
    ]
    &&
    (
      A
      \rightrightarrows
      1
    )
    \\
    (
      A \rightrightarrows A
    )
    \ar[
      rr,
      "{ (i \,\rightrightarrows\, i) }"
    ]
    &&
    (
      \widehat \Gamma \rightrightarrows \widehat \Gamma
    )
    \ar[
      rr,
      "{ (p \,\rightrightarrows\, p) }"
    ]
    \ar[
      urr,
      "{
        (
          (\mathrm{id}, \mathrm{e}) \,\rightrightarrows\, \mathrm{id}
        )
      }"
    ]
    &&
    (
      \Gamma \rightrightarrows \Gamma
    )
  \end{tikzcd}
  \end{equation}

  \vspace{-2mm}
  \noindent
  with the constant 2-groups
  from Ex. \ref{TopologicalGroupsAsTopologicalTwoGroups} on the bottom,
  and
  the delooping 2-group from Ex. \ref{Strict2GroupsDeloopingAbelianGroups} on the top right.
\end{example}

\medskip

\noindent
{\bf Equivariant topological groupoids.}

\begin{definition}[Equivariant topological groupoids]
  \label{EquivariantTopologicalGroupoids}
  We write
  $\Groupoids(\GActionsOnTopologicalSpaces)$
  for the category of internal groupoids (Ntn. \ref{Internalization})
  internal to equivariant topological spaces \eqref{GActionsOnTopologicalSpaces},
  hence of diagrams of the form \eqref{DiagramForTopologicalGroupoid}
  but where all spaces involved are equipped with continuous $G$-actions
  and all maps involved are $G$-equivariant.
With the equivariance group
  regarded as a group object in topological groupoids
  via the inclusion \eqref{FullInclusionOfTopologicalSpacesIntoTopologicalGroupoids}
  of Ex. \ref{TopologicalSpacesAsTopologicalGroupoids}
      \vspace{-1mm}
  $$
    G
      \,\in\,
      \Groups (\kTopologicalSpaces)
        \xhookrightarrow{\;\;
          \Groups( \ConstantGroupoid )
       \;\; }
      \Groups\left(\TopologicalGroupoids\right),
  $$

      \vspace{-1mm}
\noindent
  this is equivalent to the category of
  internal $G$-actions (Ntn. \ref{Internalization})
  in
  $\TopologicalGroupoids$ (Ntn. \ref{TopologicalGroupoids}):
      \vspace{-2mm}
  \begin{equation}
    \label{EquivalentIncarnationsOfGActionsOnTopologicalGroupoids}
    \Groupoids
    \left(
      \GActionsOnTopologicalSpaces
    \right)
    \;\;
    \simeq
    \;\;
    \Actions{G}
    \left(
      \TopologicalGroupoids
    \right)
    \,.
  \end{equation}
\end{definition}

\begin{example}[Right action groupoid inherits left group actions]
  \label{RightActionGroupoidInheritsLeftGroupActions}
  Let $G_L, G_R \,\in\, \Groups(\kTopologicalSpaces)$
  and
      \vspace{-2mm}
  $$
    G_L \acts \, \TopologicalSpace \; \rightacts G_R
    \;\;\;
    \in
    \;
    \Actions{
    \left(
      G_L \times G_R^{\mathrm{op}}
    \right)
    }
    (\kTopologicalSpaces)
    \,.
  $$

      \vspace{-2mm}
\noindent
  Then the $G^{\mathrm{op}}_R$-action groupoid of $\TopologicalSpace$ (Ex. \ref{TopologicalActionGroupoid})
  inherits the left $G_L$-action to become a
  $G_L$-eqivariant topological groupoid (Def. \ref{EquivariantTopologicalGroupoids}):
        \vspace{-2mm}
  $$
    G_L
      \acts \;
    (
      \TopologicalSpace
        \times
      G_R
        \rightrightarrows
      \TopologicalSpace
    )
    \;\;\;
    \in
    \;
    \Actions{G_L}
    \left(
      \Groupoids(\kTopologicalSpaces)
    \right)
    \,.
  $$
\end{example}

In generalization of Ex. \ref{ConjugationActionOnMappingSpaces}, we have:

\begin{example}[Conjugation action on mapping groupoid]
  \label{ConjugationActionOnMappingGroupoids}
  Given a pair of
  $G$-equivariant topological groupoids
  (Def. \ref{EquivariantTopologicalGroupoids})
  $$
    G \acts \,
    (
      \TopologicalSpace_1 \rightrightarrows \TopologicalSpace_0
    )
    ,\;
    G \acts \,
    (
      \mathrm{Y}_1 \rightrightarrows \mathrm{Y}_0
    )
    \;\;\;
    \in
    \;
    \Actions{G}(\TopologicalGroupoids)
  $$
  the mapping groupoid \eqref{MappingGroupoidOfTopologicalGroupoids}
  of their underlying topological groupoids
  (Ntn. \ref{TopologicalGroupoids})
  inherits the conjugation action
  \vspace{-2mm}
  $$
    G
    \acts \,
    \mathrm{Maps}
    \big(
      (\TopologicalSpace_1
      \rightrightarrows
      \TopologicalSpace_0)
      ,\,
      (\mathrm{Y}_1
      \rightrightarrows
      \mathrm{Y}_0)
    \big)
    \;\;\;
    \in
    \Actions{G}(\TopologicalGroupoids)
  $$

  \vspace{-2mm}
  \noindent
  given on the spaces of morphisms and of objects by the
  restriction of the ordinary conjugation action \eqref{ConjugationActionOnMapsBetweenGSpaces}.
With the first argument fixed, this construction constitutes
a right adjoint to the product operation,
in joint generalization of
\eqref{InternalHomInGSpaces} and
\eqref{InternalHomAdjunctionForTopologicalGroupoids}:
    \vspace{-5mm}
\begin{equation}
  \label{InternalHomInGEquivariantTopologicalGroupoids}
  \begin{tikzcd}[column sep=50pt]
    \Actions{G}
    \left(
      \Groupoids(\kTopologicalSpaces)
    \right)
    \ar[
      rr,
      phantom,
      "{\scalebox{.7}{$\bot$}}"
    ]
    \ar[
      rr,
      shift right=5pt,
      "{
        G \acts \, \mathrm{Maps}
        \left(
          (\TopologicalSpace_1 \rightrightarrows \TopologicalSpace_0 )
          , \,
          -
        \right)
      }"{below}
    ]
    &&
    \Actions{G}
    \left(
      \Groupoids(\kTopologicalSpaces)
    \right)
    \mathrlap{\,.}
    \ar[
      ll,
      shift right=5pt,
      "{
        G \acts \, (\TopologicalSpace_1 \rightrightarrows \TopologicalSpace_0)
        \times (-)
      }"{above}
    ]
  \end{tikzcd}
\end{equation}
\end{example}

\noindent
{\bf Equivariant topological realization of equivariant topological groupoids.}
The simplicial nerve functor \eqref{SimplicialTopologicalNerveOfTopologicalGroupoids}
is a right adjoint and as such preserves all limits; but
it preserves only a small class of colimits, among them
the class of the following Lem. \ref{NervePreservesLeftQuotientsOfRightActionGroupoids},
where this preservation secretly underlies classical constructions
(Prop. \ref{TopologicalRealizationOfUniversalPrincipalGroupoidIsUniversalPrincipalBundle} below)
of universal principal bundles via realization of topological groupoids
(as amplified in \cite[p. 11]{GuillouMayMerling17}):

\begin{lemma}[Nerve preserves left quotients of right action groupoids]
  \label{NervePreservesLeftQuotientsOfRightActionGroupoids}
  Let $\Gamma_L, \Gamma_R \,\in\, \Groups(\kTopologicalSpaces)$
  and
      \vspace{-2mm}
  $$
    \Gamma_L \acts \, \TopologicalSpace \; \rightacts \Gamma_R
    \;\;\;
    \in
    \;
    \Actions{
    \left(
      \Gamma_L \times \Gamma_R^{\mathrm{op}}
    \right)
    }
    (\kTopologicalSpaces)
    \,.
  $$

      \vspace{-2mm}
\noindent
  Then the canonical comparison morphism
  from
  the nerve \eqref{SimplicialTopologicalNerveOfTopologicalGroupoids}
  of
  the left $\Gamma_L$-quotient
  of the equivariant right action groupoid (Ex. \ref{RightActionGroupoidInheritsLeftGroupActions})
  to the $\Gamma_L$-quotient of the nerve
  is an isomorphism:
      \vspace{-2mm}
  $$
    N
    \left(
      \Gamma_L
        \backslash
      (
        \TopologicalSpace
          \times
        \Gamma_R
          \rightrightarrows
        \TopologicalSpace
      )
    \right)
    \;\;
    \simeq
    \;\;
    \Gamma_L
      \backslash
    \left(
      N
      (
        \TopologicalSpace
          \times
        \Gamma_R
          \rightrightarrows
        \TopologicalSpace
      )
    \right)
    \,.
  $$
\end{lemma}
\begin{proof}
  The point is that
  in each degree $n \in \mathbb{N}$,
  the $\Gamma_L$-action is non-trvial
  only on the leftmost factor of the $n$-th component space of the nerve:

  \vspace{-7mm}
  $$
    \begin{tikzcd}[column sep=-2pt]
    N
    (
      \TopologicalSpace
        \times
      \Gamma
      \rightrightarrows
      \TopologicalSpace
    )_n
    \ar[
      d,
      "{
        \mbox{
          \tiny
          \color{greenii}
          \bf
          \begin{tabular}{c}
            quotient
            \\
            coprojection
          \end{tabular}
        }
      }"{left}
    ]
    &
    \simeq
    &
    \TopologicalSpace
    \times
    \overset{
      \mbox{\tiny \rm $n$ factors}
    }{
      \overbrace{
        \Gamma \times \cdots \times \Gamma
      }
    }
    \ar[d]
    \\
    \Gamma
      \backslash
    \left(
      N
      (
        \TopologicalSpace \times \Gamma
        \rightrightarrows
        \TopologicalSpace
      )_n
    \right)
    &\simeq&
    (
      \Gamma
        \backslash
      \TopologicalSpace
    )
    \times
    \underset{
      \mbox{\tiny \rm $n$ factors}
    }{
      \underbrace{
        \Gamma \times \cdots \times \Gamma
      }
    }
    &\simeq&
      N
      \left(
        \Gamma
          \backslash
        (
          \TopologicalSpace
            \times
          \Gamma
            \rightrightarrows
          \TopologicalSpace
        )
      \right)_n
      \mathrlap{\,.}
    \end{tikzcd}
  $$

  \vspace{-7mm}
\end{proof}

\begin{proposition}[Topological realization of topological groupoids respects equivariance]
  \label{TopologicalRealizationOfTopologicalGroupoidsRespectsEquivariance}
  $\,$

  \noindent
  For any $G \,\in\, \Groups(\kTopologicalSpaces)$,
  topological realization \eqref{TopologicalRealizationOfTopologicalGroupoids}
  of topological groupoids
  extends to a functor
  from $G$-actions on topological groupoids (Def. \ref{EquivariantTopologicalGroup})
  to $G$-action on topological spaces (Ntn. \ref{GActionOnTopologicalSpaces}):
      \vspace{-1mm}
  $$
    \left\vert
      -
    \right\vert
    \;:\;
    \Actions{G}
    \left(
      \Groupoids(\kTopologicalSpaces)
    \right)
    \xrightarrow{\;\;\;\;\;}
    \Actions{G}
    (\kTopologicalSpaces)
    \,.
  $$
\end{proposition}
\begin{proof}
  Under the equivalence \eqref{EquivalentIncarnationsOfGActionsOnTopologicalGroupoids}
  the functor in question is of the form
    \vspace{-2mm}
  $$
    \begin{tikzcd}
      \Actions{G}
      \left(
        \Groupoids(\kTopologicalSpaces)
      \right)
      \ar[
        rr,
        " \scalebox{0.8}{$
          \Actions{(\vert- \vert)}
          (
            \vert- \vert
          )
        $}"
      ]
      &&
      \Actions{G}
      \left(
        \kTopologicalSpaces
      \right)
      \,,
    \end{tikzcd}
  $$

      \vspace{-2mm}
\noindent
  in the notation \eqref{FunctorOnStructuresInducedFromLexFunctor},
  and hence exists since $\left\vert-\right\vert$ preserves finite limits
  (here in particular: finite products),
  by Lemma \ref{TopologicalRealizationFunctors}.
  (Here we are using that form of the group object is indeed preserved, by
  Ex. \ref{TopologicalRealizationOfConstantGroupoids}).
\end{proof}

\medskip

\noindent
{\bf Equivariantly equivariant topological groupoids.}
In equivariant generalization of Def. \ref{EquivariantTopologicalGroup}, we have:

\begin{definition}[Equivariant topological groupoid]
Let $G \,\in\, \Groups(\kTopologicalSpaces)$
and
$
  G \acts \, \Gamma_L
  \,\in\,
  \Groups
  \left(
    \Actions{G}(\kTopologicalSpaces)
  \right)
$
be a $G$-equivariant group (Def. \ref{EquivariantTopologicalGroup}).
With $G$-equivariant topological groups regarded as
constant $G$-equivariant topological groupoids
    \vspace{-2mm}
$$
  (G \acts \, \Gamma)
  \;\in\;
  \Groups
  \left(
    \Actions{G}(\kTopologicalSpaces)
  \right)
  \xhookrightarrow{\;\;
    \Groups
    \left(
      \Actions{
        (\ConstantGroupoid(G))
      }
      (
        \ConstantGroupoid(-)
      )
    \right)
  \;\;}
\;  \Groups
  \left(\!
    \Actions{G}
    \left(
      \Groupoids(\kTopologicalSpaces)
    \right)
  \!\right),
$$

    \vspace{-2mm}
\noindent
we obtain the category
of $(G \acts \, \Gamma)$-equivariant $G$-equivariant topological groupoids:
    \vspace{-1mm}
$$
  \Actions{(G \acts \, \Gamma)}
  \left(\!
    \Groupoids
    \left(
      \Actions{G}(\kTopologicalSpaces)
    \right)
  \!\right)
  \;\;
    \simeq
  \;\;
  \Actions{(G \acts \, \Gamma)}
  \left(\!
    \Actions{G}
    \left(
      \Groupoids
        (\kTopologicalSpaces)
    \right)
  \!\right)
  \,.
$$
\end{definition}

In equivariant generalization of Ex. \ref{RightActionGroupoidInheritsLeftGroupActions}, we get:

\begin{example}[Right equivariant action groupoid inherits left equivariant group action]
\label{RightEquivariantActionGroupoidInheritsLeftEquivariantGroupAction}
For $G \,\in\, \Groups(\kTopologicalSpaces)$, let
    \vspace{-1mm}
$$
  G \acts \, \Gamma_L
  ,\,
  G \acts \, \Gamma_R
  \quad \in\,
  \Groups
  \left(
    \Actions{G}(\kTopologicalSpaces)
  \right)
$$

    \vspace{-1mm}
\noindent
be two $G$-equivariant groups (Def. \ref{EquivariantTopologicalGroup}) and
consider a topological $G$-space $\TopologicalSpace$ \eqref{GActionsOnTopologicalSpaces}
equipped with commuting equivariant left and
right actions by these, respectively:
    \vspace{-2mm}
$$
  (G \acts \, \Gamma_L)
  \acts
  \;
  (
    G \acts \, \TopologicalSpace
  )
  \;
  \rightacts
  (G \acts \, \Gamma_R)
  \;\;\;
  \in\;
  \Actions{
    \big(\!
      (G \acts \, \Gamma_1)
      \times
      (G \acts \, \Gamma_2)^{\mathrm{op}}
    \big)
  }
  \left(
    \Actions{G}(\kTopologicalSpaces)
  \right)
  \,.
$$

    \vspace{-1mm}
\noindent
Then the evident $(G \acts  \, \Gamma_R)$-action $G$-groupoid
(formed just as in Ex. \ref{TopologicalActionGroupoid}
with Ex. \ref{RightActionGroupoidInheritsLeftGroupActions},
only that
now all component spaces carry $G$-action and all structure maps are $G$-equivariant)
inherits
a (further, $G$-equivariant) $(G \acts \, \Gamma_L)$-action:
    \vspace{-2mm}
$$
  \underset{
    \mbox{
      \tiny
      \color{orangeii}
      \bf
      \begin{tabular}{c}
        $G$-equivariant
        \\
        left $\Gamma$-action
      \end{tabular}
    }
  }{
  \underbrace{
  (G \acts \, \Gamma_L)
  \,\acts
  }}
  \;
  \underset{
    \mathclap{
    \raisebox{0pt}{
      \tiny
      \color{darkblue}
      \bf
      \begin{tabular}{c}
        $G$-equivariant
        \\
        right action groupoid
      \end{tabular}
    }
    }
  }{
  \underbrace{
  \big(
    (G \acts \, \TopologicalSpace)
    \times
    (G \acts \, \Gamma_R)
    \rightrightarrows
    (G \acts \, \TopologicalSpace)
  \big)
  }
  }
  \;\;\;
  \in
  \;
  \Actions{
    (G \acts \, \Gamma_L)
  }
  \left(\!
    \Groupoids
    \left(
      \Actions{G}
      (\kTopologicalSpaces)
    \right)
  \!\right)
  \,.
$$
\end{example}

In equivariant generalization of Ex. \ref{TopologicalDeloopingGroupoid},
we get:

\begin{example}[Equivariant topological delooping groupoid]
  \label{EquivariantTopologicalDeloopingGroupoid}
  For $(G \acts \, \Gamma) \,\in\, \Groups ( \Actions{G}(\kTopologicalSpaces))$,
  its {\it equivariant delooping groupoid} is the
  $G$-equivariant topological left action groupoid
  (Ex. \ref{RightEquivariantActionGroupoidInheritsLeftEquivariantGroupAction})
  of the unique $(G \acts \, \Gamma)$-action on $(G \acts \, \ast)$:
      \vspace{-1mm}
  $$
    G \acts \, \mathbf{B}\Gamma
    \;\coloneqq\;
    G \acts \,
    (
      \ast \times \Gamma^{\mathrm{op}}
        \rightrightarrows
      \ast
    )
    \;=\;
    G \acts \,
    (
      \Gamma^{\mathrm{op}}
      \rightrightarrows
      \ast
    )
    \;\;\;
    \in
    \;
    \Groupoids
    \left(
      \Actions{G}(\kTopologicalSpaces)
    \right)
    \,.
  $$
\end{example}

\section{$G$-Equivariant homotopy types}
\label{GEquivariantHomotopyTypes}

We recall and develop some of equivariant homotopy theory
needed
in \cref{ConstructionOfUniversalEquivariantPrincipalBundles}.
for identifying equivariant homotopy groups
of equivariant classifying spaces.

\begin{notation}[Classical model structure on topological spaces]
  \label{ClassicalModelStructureOnTopologicalSpaces}
  {\bf (i)}
  We write
 \vspace{-2mm}
  \begin{equation}
    \label{TheClassicalModelCategoryOnTopologicalSpaces}
    \kTopologicalSpaces_{\mathrm{Qu}}
    \;\;\;
    \in
    \;
    \ModelCategories
  \end{equation}

   \vspace{-2mm}
\noindent  for the classical model category on
  (compactly generated, Ntn. \ref{CompactlyGeneratedTopologicalSpaces})
  topological spaces, whose weak equivalences are the weak homotopy equivalences
  $\WeakHomotopyEquivalences$
  and whose fibrations are the Serre fibrations $\SerreFibrations$
  (e.g., \cite[Thm. 2.4.23]{Hovey99}).

  \noindent {\bf (ii)}  We write
    \vspace{-1mm}
  \begin{equation}
    \label{TheClassicalHomotopyCategory}
    \HomotopyCategory
    (
      \kTopologicalSpaces_{\mathrm{Qu}}
    )
    \;\;\;
    \in
    \;
    \Categories
  \end{equation}

  \vspace{-1mm}
  \noindent
  for the classical homotopy category; the localization of
  \eqref{TheClassicalModelCategoryOnTopologicalSpaces} at the weak equivalences.
\end{notation}
For example:
\begin{lemma}[Locally trivial bundles are Serre fibrations {\cite[p. 130, Thm. 6.3.3]{tomDieck08}}]
  \label{LocallyTrivialBundlesAreSerreFibrations}
  A map $\mathrm{E} \xrightarrow{p} \TopologicalSpace \;\; \in \; \kTopologicalSpaces$,
  which over some open cover
  $\widehat {\TopologicalSpace}
    \, = \, \underset{i \in I}{\sqcup} \TopologicalPatch_i
    \twoheadrightarrow \TopologicalSpace
  $  (Ex. \ref{OpenCoversAreEffectiveEpimorphisms})
  restricts to
  a Cartesian product,
  is a Serre fibration:
  \vspace{-2mm}
  $$
    \begin{tikzcd}
      \mathrm{U} \times \mathrm{F}
      \ar[rr]
      \ar[
        d,
        "{
          \mathrm{pr}_1
        }"
      ]
      \ar[
        drr,
        phantom,
        "\mbox{\tiny\rm (pb)}"
      ]
      &&
      \mathrm{E}
      \ar[
        d,
        "p"
      ]
      \\
      \mathrm{U}
      \ar[
        rr,
        ->>,
        "\mathrm{open}"{below}
      ]
      &&
      \TopologicalSpace
    \end{tikzcd}
    \;\;\;\;\;\;\;\;\;\;
    \Rightarrow
    \;\;\;\;\;\;\;\;\;\;
    p \,\in\, \SerreFibrations
    \,.
  $$
\end{lemma}

\begin{notation}[Homotopy fiber sequence]
  \label{HomtopyFiberSequenceOfTopologicalSpaces}
  A sequence of morphisms
    \vspace{-2mm}
  $$
    \begin{tikzcd}[column sep=small]
    F
    \ar[
      r,
      "i"
    ]
    &
    \mathrm{E}
    \ar[
      r,
      "{p}"
    ]
    &
    \TopologicalSpace
    \end{tikzcd}
    \;\;\;\;
    \in
    \;
    \mathcal{C}
  $$

    \vspace{-3mm}
\noindent
  in a model category $\mathcal{C}$
  is a {\it homotopy fiber sequence}
  for a given point $\ast \xrightarrow{x} \TopologicalSpace$, to be denoted
    \vspace{-2mm}
  $$
    i \,\simeq\, \HomotopyFiber{x}(p)
    \,,
  $$

  \vspace{-2mm}
\noindent
  if, for any choice of fibration replacement
  $\hat p \,\colon\, \widehat{\TopologicalSpace} \xrightarrow{\;} \TopologicalSpace$ of $p$,
  the sequence is isomorphic, in the homotopy category,
  to the ordinary $x$-fiber mapping into the fibration:
    \vspace{-2mm}
  \begin{equation}
    \label{CharacterizationOfHomotopyKanFibration}
    \begin{tikzcd}[row sep=20pt, column sep=huge]
      F
      \ar[
        r,
        "i"
      ]
      \ar[
        d,
        dashed,
        "\in \, \mathrm{W}"{right},
        "\exists"{left}
      ]
      &[14pt]
      \mathrm{E}
      \ar[
        r,
        "{p}"
      ]
      \ar[
        d,
        "\in \, \mathrm{W}"{right}
      ]
      &
      \mathrm{Y}
      \ar[
        d,-,
        shift left=1pt
      ]
      \ar[
        d,-,
        shift right=1pt
      ]
      \\
      \widehat{E}_{x_0}
      \ar[
        r,
        "\mathrm{fib}_{x}(\hat p)"
      ]
      &
      \widehat {\mathrm{E}}
      \ar[
        r,
        " \hat p "{above},
        "\in \, \mathrm{Fib}"{below}
      ]
      &
      \TopologicalSpace
    \end{tikzcd}
    \;\;\;\;\;
    \in
    \;
    \HomotopyCategory
    (
      \mathcal{C}
    )
    \,.
  \end{equation}
\end{notation}

\medskip

\noindent
{\bf Proper equivariant homotopy theory of $G$-spaces.}

\begin{proposition}[Proper equivariant model category of $G$-spaces {\cite[\S 1.2]{DwyerKan84}}]
  \label{ProperEquivariantModelCategoryOfGSpaces}
  For $G \,\in\, \Groups(\CompactSmoothManifolds)
   \xhookrightarrow{\;} \Groups(\kTopologicalSpaces)$
  the topological group underlying a compact Lie group,
  there exists a model category structure
  on the category of topological $G$-spaces (Ntn. \ref{GActionOnTopologicalSpaces}),
  to be denoted
  \vspace{-2mm}
  \begin{equation}
    \label{FineModelStructureOnTopologicalGSpaces}
    \Actions{G} (\kTopologicalSpaces_{\mathrm{Qu}})_{\mathrm{prop}}
    \;\;\;
    \in
    \;
    \ModelCategories
  \end{equation}

  \vspace{-2mm}
  \noindent
  whose weak equivalences
  $\ProperEquivariantWeakHomotopyEquivalences{G}$
  and
  fibrations
  $\ProperEquivariantSerreFibrations{G}$
  are those morphisms whose
  underlying continuous maps between $H$-fixed loci,
  for all closed \footnote{This model category structure itself exists for any
  topological group and any non-empty
  subset of subgroups,
  but for it to be equivalent to $G$-CW-complexes
  localized at the $G$-equivariant homotopy equivalences
  one needs to make sufficient restrictions, such as to
  compact Lie groups and their closed subgroups.}
  subgroups $H \,\underset{\mathclap{\mathrm{clsd}}}{\subset}\, G$, are
  weak equivalences or fibrations in the
  Serre-Quillen model structure on topological spaces,
  (Ntn. \ref{ClassicalModelStructureOnTopologicalSpaces})
  hence are weak homoropy equivalences or Serre fibrations, respectively.
\end{proposition}

\begin{proposition}[Cofibrant generation of the proper equivariant model structure]
\label{CofibrantGenerationOfProperEquivariantModelStructure}
$\,$

\noindent
  The model category
  $\Actions{G}(\kTopologicalSpaces_{\mathrm{Qu}})_{\mathrm{prop}}$
  from Prop. \ref{ProperEquivariantModelCategoryOfGSpaces} is

\vspace{-4mm}
  \begin{enumerate}[{\bf (i)}]
  \setlength\itemsep{-4pt}

    \item
    a proper model category;

    \item
    a cofibrantly generated model category (see e.g. \cite[\S 11]{Hirschhorn02})
    whose class of
    generating (acyclic) cofibrations is the product of
    coset spaces (Ex. \ref{CosetSpacesAsActions})
    with the generating (acyclic) cofibrations of the
    Serre-Quillen nodel strutucre
    \vspace{-3mm}
    \begin{align}
      \label{GeneratingCofibrationsOfProperEquivariantModelCategory}
      I_{\mathrm{prop}}
      &
      \;\coloneqq\;
      \big\{
        G/H \,\times\, S^{n-1}
        \xhookrightarrow{\;\; \mathrm{id}_{G/H} \times  i_{n} \;\;}
        G/H \,\times\, D^n
      \big\}_{
        n \in \mathbb{N},\, H \underset{\mathrm{clsd}}{\subset} G
      }
      \\
      \label{GeneratingAcyclicCofibrationsOfProperEquivariantModelCategory}
      I_{\mathrm{prop}}
      &
      \;\coloneqq\;
      \big\{
        G/H \,\times\, D^n \times \{0\}
        \xhookrightarrow{\;\; \mathrm{id}_{G/H} \times  j_{n} \;\;}
        G/H \,\times\, D^n \times [0,1]
      \big\}_{
        n \in \mathbb{N},\, H \underset{\mathrm{clsd}}{\subset} G
        \,;
      }
    \end{align}

  \vspace{-.2cm}
  \item
  an enriched model category
  (see e.g. \cite[\S 4.3]{GuillouMay11})
  over $\kTopologicalSpaces_{\mathrm{Qu}}$ (Ntn. \ref{ClassicalModelStructureOnTopologicalSpaces}),
  in that

  \vspace{-2mm}
  \begin{enumerate}[{\bf (a)}]
\item   the functor assigning spaces \eqref{EquivariantFunctions}
  of equivariant maps
    \vspace{-2mm}
  $$
    \mathrm{Maps}(-,-)^G
    \;\colon\;
    \Actions{G}(\kTopologicalSpaces_{\mathrm{Qu}})^{\mathrm{op}}
      \times
    \Actions{G}(\kTopologicalSpaces_{\mathrm{Qu}})
      \longrightarrow
    \kTopologicalSpaces_{\mathrm{Qu}}
  $$

  \vspace{-3mm}
  \noindent
  is a right Quillen bifunctor;
 \item for $G \acts \, \TopologicalSpace \xrightarrow{ c \in \mathrm{Cof} } G \acts  \, \mathrm{Y}$
  a cofibration,
  and
  $G \acts \, \mathrm{A} \xrightarrow{f \in \mathrm{Fib}} G \acts  \, \mathrm{B}$ a fibration, the
  morphism
  \vspace{-2mm}
  \begin{equation}
    \label{QuillenBifunctorPropertyForEquivariantMappingSpace}
    \mathrm{Maps}
      (\mathrm{Y}
      ,\,
      \mathrm{A} )
    \xrightarrow{
      \;\;
      (c^\ast, f_\ast)
    }
    \mathrm{Maps}(\TopologicalSpace,\, \mathrm{A})
    \underset{
      \mathrm{Maps}( \TopologicalSpace ,\, \mathrm{B} )
    }{\times}
    \mathrm{Maps}( \mathrm{Y} ,\, \mathrm{B} )
  \end{equation}

  \vspace{-1mm}
    \noindent
  is a fibration,
  \item
  and is in addition a weak equivalence (in $\kTopologicalSpaces_{\mathrm{Qu}}$)
  if
  $c$ or $f$ is a weak equivalence.
\end{enumerate}
\end{enumerate}
\end{proposition}
\begin{proof}
  The first two statement may be found as \cite[Prop. 2.11]{Fausk08}\cite[Prop. 2.6]{Stephan13}.
  With the third statement included this appears in
  \cite[Thm. 3.7]{GuillouMayRubin13}\cite[Prop. B.7]{Schwede18}.
\end{proof}
\begin{example}[$G$-CW complexes are cofibrant objects in the proper equivariant model category]
\label{GCWComplexesAreCofibrantObjectsInProperEquivariantModelcategory}
$\,$

\noindent
{\bf (i)} A {\it $G$-CW complex} \cite{Matumoto71a}
is, by definition,
a $G$-space obtained by a
sequence, monotone in the dimension $n$, of cell attachments
with the generating cofibrations \eqref{GeneratingCofibrationsOfProperEquivariantModelCategory}.
Hence Prop. \ref{CofibrantGenerationOfProperEquivariantModelStructure} implies that
$G$-CW complexes are cofibrant objects in the proper equivariant model structure
(Prop. \ref{ProperEquivariantModelCategoryOfGSpaces}).

\noindent
{\bf (ii)} By the equivariant triangulation theorem \cite[Thm. 3.1]{Illman72}\cite{Illman83},
every smooth manifold with a smooth action by a compact Lie group $G$ admits
an equivariant triangulation, hence is a $G$-CW complex and hence a cofibrant
object in the proper equivariant model category:

    \vspace{-2mm}
    $$
      \Actions{G}(\SmoothManifolds)
      \xhookrightarrow{\quad}
      \GCWComplexes
      \;\xhookrightarrow{\quad}\;
      \left(
        \Actions{G}
        (
          \kTopologicalSpaces_{\mathrm{Qu}}
        )_{\mathrm{prop}}
      \right)^{\mathrlap{\mathrm{cof}}}
      \;\;\;\;.
    $$
\end{example}

\begin{proposition}[Monoidalness of the proper equivariant model structure]
  \label{MonoidalnessOfProperEquivariantModelCategory}
  The model category
  $\Actions{G}(\kTopologicalSpaces_{\mathrm{Qu}})_{\mathrm{prop}}$
  from Prop. \ref{ProperEquivariantModelCategoryOfGSpaces} is
  a monoidal model category, in that its Cartesian product is
  a left Quillen bifunctor, meaning that for
  for any pair of cofibrations
  $G \acts \, \TopologicalSpace_1 \xrightarrow{f_1 \in \mathrm{Cof}} \mathrm{Y}_1$,
  $G \acts \, \TopologicalSpace_2 \xrightarrow{f_1 \in \mathrm{Cof}} G \acts \, \mathrm{Y}_2$
  the {\it pushout-product} morphism
  \vspace{-1mm}
  \begin{equation}
    \label{PushoutProductMorphism}
    \begin{tikzcd}
      \TopologicalSpace_1 \times \mathrm{Y}_2
      \underset{
        \mathrm{Y}_1 \times \mathrm{Y}_2
      }{\coprod}
      \mathrm{Y}_1 \times \TopologicalSpace_2
      \ar[
        rr,
        "{
          \scalebox{.7}{$
            f_1 \PushoutProduct f_2
          $}
        }"{above},
        "{\in \, \mathrm{Cof}}"{below}
      ]
      &&
      \TopologicalSpace_1 \times \TopologicalSpace_2
    \end{tikzcd}
  \end{equation}

  \vspace{-2mm}
  \noindent
  is a cofibration, and is in addition a weak equivalence if
  $f_1$ or $f_2$ is such.
\end{proposition}
\begin{proof}
  It is folklore that this follows
  from the fact that products of coset spaces
  $G/H_1 \times G/H_2$
  admit a $G$-CW complex structure,
  by Ex. \ref{GCWComplexesAreCofibrantObjectsInProperEquivariantModelcategory}
  (making crucial use of the assumption that $G$ is a Lie group,
  so that its coset spaces are smooth manifolds with smooth $G$-action).
  The conclusion has more recently been
  made explicit in \cite[Prop. 1.1.3 (iii)]{DHLPS19}.
\end{proof}

\begin{proposition}[Equivariant mapping space Quillen adjunction]
  \label{EquivariantMappingSpaceQuillenAdjunction}
  For $\TopologicalSpace \,\in\,
  \left(\Actions{G}(\kTopologicalSpaces_{\mathrm{Qu}})_{\mathrm{prop}}\right)^{\mathrm{cof}}$
  a cofibrant $G$-space,
  the equivariant mapping space adjunction \eqref{InternalHomInGSpaces}
  for $\TopologicalSpace$
  is a Quillen adjunction from the proper equivariant model structure
  (Prop. \ref{ProperEquivariantModelCategoryOfGSpaces}) to itself:
  \vspace{-3mm}
  $$
  \begin{tikzcd}[column sep=35pt]
    \Actions{G}(\kTopologicalSpaces_{\mathrm{Qu}})_{\mathrm{prop}}
    \ar[
      rr,
      shift right=6pt,
      "{
        \scalebox{.7}{$
          G \acts \, \mathrm{Maps}(\TopologicalSpace,\, -)
        $}
      }"{below}
    ]
    \ar[
      rr,
      phantom,
      "\scalebox{.7}{$\bot_{\mathrlap{\mathrm{Qu}}}$}"
    ]
    &&
    \Actions{G}(\kTopologicalSpaces_{\mathrm{Qu}})_{\mathrm{prop}}\;.
    \ar[
      ll,
      shift right=6pt,
      "{
        \scalebox{.7}{$
          G \acts \,  \TopologicalSpace \,\times\, (-)
        $}
      }"{above}
    ]
  \end{tikzcd}
  $$
\end{proposition}
\begin{proof}
 This is a standard consequence of the monoidal
 model category structure due to Prop. \ref{MonoidalnessOfProperEquivariantModelCategory},
 which immediately implies that the operation of forming the cartesian products with a cofibrant
 object is a left Quillen functor.
\end{proof}

\begin{example}[Equivariant classifying shapes for pure shape structure]
\label{TopologicalGSpaceModellingHomotopyGFixedLociInAClassifyingSpace}
$\,$

\noindent For $G \,\in\, \Groups(\FiniteSets)$
and $\mathcal{G} \,\in\, \Groups(\SimplicialTopologicalSpaces)_{\wellpointed}$
(Ntn. \ref{WellPointedTopologicalGroup}),
the topological $G$-space
\vspace{-2mm}
$$
  \Maps{}
    { \TopologicalRealization{}{ \mathbf{E}G  }  }
    { \TopologicalRealization{}{ \mathbf{B} \mathcal{G} } }
  \;\;
  \in
  \;\;
  \Actions{G}(\kTopologicalSpaces)
$$

\vspace{-2mm}
\noindent
has fixed loci at $H \subset G$ naturally weakly homotopy equivalent to
\vspace{-2mm}
\begin{align*}
  \SingularSimplicialComplex
  \big(
    \Maps{}
      { \TopologicalRealization{}{ \mathbf{E}G }  }
      { \TopologicalRealization{}{ \mathbf{B} \mathcal{G} } }
    ^{H}
  \big)
  \;\;
  &\simeq
  \;\;
  \SingularSimplicialComplex
  \big(
    \Maps{}
      { \TopologicalRealization{}{ \mathbf{E}G } / H  }
      { \TopologicalRealization{}{ \mathbf{B} \mathcal{G} } }
  \big)
\\
  \;\;
 & \simeq
  \;\;
  \SingularSimplicialComplex
  \big(
    \Maps{}
      { \TopologicalRealization{}{ \mathbf{B}H }  }
      { \TopologicalRealization{}{ \mathbf{B} \mathcal{G} } }
  \big)
  \\
  \;
  &\simeq\;
  \Shape
  \,
  \Maps{}
    { B H }
    { B \mathcal{G} }
  \,.
\end{align*}

\vspace{-2mm}
\noindent
Hence
\vspace{-2mm}
$$
  \begin{tikzcd}[row sep=-1pt, column sep=large]
    \Actions{G}(\DTopologicalSpaces)
    \ar[r, "\FixedLoci"{above}]
    &
    \GEquivariant\SmoothInfinityGroupoids
    \ar[r, "\Shape"{above}]
    &
    \GEquivariant\InfinityGroupoids \;.
    \\
  \scalebox{0.7}{$
    G \acts \;
    \Maps{}
      { \TopologicalRealization{}{ \mathbf{E}G } }
      { \TopologicalRealization{}{ \mathbf{B}\mathcal{G} } }
      $}
          \ar[r, phantom, "\longmapsto"]
    &
    \scalebox{0.7}{$
    \Maps{}
      { \TopologicalRealization{}{ \mathbf{B}(-) } }
      { \TopologicalRealization{}{ \mathbf{B}\mathcal{G} } }
      $}
    \ar[r, phantom, "\longmapsto"]
    &
    \scalebox{0.7}{$
    \Maps{}
      { B(-) }
      { B \mathcal{G} }
      $}
  \end{tikzcd}
$$
\end{example}

\medskip

\noindent
{\bf Borel equivariant homotopy theory of $G$-spaces.} In fact, the
above statements and proofs about the proper equivariant model structure
apply more generally to any choice of subset of closed subgroups of $G$
containing the trivial subgroup, the above being the case of the
maximal such subset. For the mininal subset containing only the
trivial group one obtains the {\it coarse} or {\it Borel}
equivariant model structure:

\begin{proposition}[Borel model structure on $G$-spaces]
  \label{BorelModelStructureOnGSpaces}
  For $G \,\in\, \Groups(\kTopologicalSpaces)$,
  there exists a model category structure
  on the category of topological $G$-spaces (Ntn. \ref{GActionOnTopologicalSpaces}),
  to be denoted
  \vspace{-1.5mm}
  \begin{equation}
    \label{CoarseModelStructureOnTopologicalGSpaces}
    \Actions{G} (\kTopologicalSpaces_{\mathrm{Qu}})_{\mathrm{coarse}}
    \;\;\;
    \in
    \;
    \ModelCategories
    \,,
  \end{equation}

   \vspace{-1.5mm}
\noindent  whose weak equivalences and fibrations are those
  $G$-equivariant functions whose underlying maps
  are weak equivalences or fibrations, respectively,
  in the classical model structure on topological spaces
  (Ntn. \ref{ClassicalModelStructureOnTopologicalSpaces}).
\end{proposition}
\begin{proof}
  This is also the special case of the
  general projective model structure on
  $\kTopologicalSpaces_{\mathrm{Qu}}$-enriched presheaf categories
  (\cite[Thm, 5.4]{Piacenza91}) for the site being the
  topological delooping groupoid $(G \rightrightarrows \ast)$; see Ex. \ref{TopologicalDeloopingGroupoid}.
\end{proof}

\begin{proposition}[Cofibrant generation of the Borel equivariant model structure]
$\,$

  \noindent
  The Borel model structure
  $\Actions{G}(\kTopologicalSpaces)_{\mathrm{coarse}}$
  from Prop. \ref{BorelModelStructureOnGSpaces} is

  \noindent
  {\bf (i)} Cofibrantly generated with generating cofibrations
    \vspace{-2mm}
    \begin{equation}
      \label{GeneratingCofibrationsOfCoarseEquivariantModelCategory}
      I_{\mathrm{coarse}}
      \;\coloneqq\;
      \big\{
        G \,\times\, S^{n-1}
        \xhookrightarrow{ \mathrm{id}_{G} \times  i_{n} }
        G \,\times\, D^n
      \big\}_{
        n \in \mathbb{N}
             }.
    \end{equation}

   \vspace{-1mm}
  \noindent
  {\bf (ii)} Cartesian monoidal, in that the pushout-product
   $f_1 \PushoutProduct f_2$
  \eqref{PushoutProductMorphism}
  of a pair of cofibrations is itself a cofibration,
  and in addition is a weak if $f_1$ or $f_2$ is so.
\end{proposition}

\medskip

\noindent
{\bf Good simplicial topological spaces.}

\begin{notation}[H-cofibrations]
  \label{HurewiczClosedCofibrations}
  We write
    \vspace{-2mm}
  $$
    \HurewiczCofibrations
    \;\subset\;
    \mathrm{Mor}(\kTopologicalSpaces)
  $$

    \vspace{-2mm}
    \noindent
  for the class of {\it Hurewicz cofibrations}
  -- to be called {\it h-cofibrations}, following \cite[p. 16]{MandellMaySchwedeShipley01} --
  in the
  category \eqref{CategoryOfTopologicalSpaces} of compactly generated weak Hausdorff spaces,
  hence for those maps $i \,\colon\, \mathrm{A} \xrightarrow{\;} \TopologicalSpace$,
  such that every solid commuting diagram as follows
  admits a dashed lift, as shown (e.g. \cite[p. 43]{May99}):
    \vspace{-2mm}
     $$
    \begin{tikzcd}
      \mathrm{A}
      \ar[
        d,
        "{i}"{left}
      ]
      \ar[
        r
      ]
      &
      \mathrm{Maps}([0,1] ,\, \mathrm{Y})
      \ar[
        d,
        "{ \mathrm{Maps}(i_0 ,\, \mathrm{Y}) }"
      ]
      \\
      \TopologicalSpace
      \ar[
        r
      ]
      \ar[
        ur,
        dashed,
        "{\exists}"{description}
      ]
      &
      \mathrm{Y}
      \mathrlap{\,.}
    \end{tikzcd}
  $$
  Notice that, due to the weak Hausdorff property assumed in
  \eqref{CategoryOfTopologicalSpaces}, such maps are necessarily
  injections with closed images,
  hence these are equivalently {\it closed cofibrations}
  (e.g. \cite[p. 44]{May99}).
\end{notation}

\begin{proposition}[Equivalent characterizations of h-cofibrations]
  \label{EquivalentCharacterizationsOfClosedCofibrations}
  $\,$

  \noindent
  A closed subspace inclusion
  $i \,\colon\, \mathrm{A} \xhookrightarrow{\;} \TopologicalSpace$
  is
  a Hurewicz cofibration, $i \,\in\,\HurewiczCofibrations$ (Ntn. \ref{HurewiczClosedCofibrations})
  if and only if the following equivalent conditions hold:

\noindent
{\bf (i)}  (\cite[Thm. 2]{Strom68}, review in \cite[\S VII, Thm. 1.3]{Bredon93}\cite[p. 45]{May99}):
The pushout-product $i \PushoutProduct 0$ of $i$ with the endpoint inclusion
    $0 \,\colon\, \ast \xrightarrow{\;} [0,1]$  admits a retraction map $r$:
      \vspace{-2mm}
    $$
      \begin{tikzcd}
        \TopologicalSpace \!\times\! \{0\}
        \;\cup\;
        \mathrm{A} \!\times\! [0,1]
        \ar[
          r,
          "{ i \PushoutProduct 0  }"
        ]
        \ar[
          rr,
          rounded corners,
          to path={
            -- ([yshift=-9pt]\tikztostart.south)
            --node[above]{
                \scalebox{.7}{$\mathrm{id}$}
              }
              ([yshift=-9pt]\tikztotarget.south)
            -- (\tikztotarget.south)}
        ]
        &
        \TopologicalSpace \times [0,1]\;.
        \ar[
          r,
          dashed,
          "{r}"
        ]
        &
        \TopologicalSpace \!\times\! \{0\}
        \;\cup\;
        \mathrm{A} \!\times\! [0,1] \;.
      \end{tikzcd}
    $$

\noindent
{\bf (ii)}  (\cite[Thm. 2]{Strom66}, review in \cite[\S VII, Thm. 1.5]{Bredon93}):
There exists

  \vspace{-2mm}
\begin{itemize}
\setlength\itemsep{-2pt}
 \item[{\bf (a)}]
 a neighborhood
$\mathrm{A} \xhookrightarrow{\;} \mathrm{U} \xhookrightarrow{\;} \TopologicalSpace$
of $\mathrm{A}$ in $\TopologicalSpace$;

  \item[{\bf (b)}]
a homotopy $\eta$ deforming this neighborhood into $A$ relative to $A$,
in that its makes the following diagram commute:
  \vspace{-2mm}
$$
  \begin{tikzcd}[column sep=20pt]
    &
    A \times [0,1]
    \ar[
      dr,
      ->>
    ]
    \\[-25pt]
    U
    \ar[
      d,
      hook,
      "{
        (\mathrm{id}, 0)
      }"{left}
    ]
    \ar[
      drr,
      hook,
      "{
      }"
    ]
    &
    &
    A
    \ar[
      d,
      hook
    ]
    \\
    U \times [0,1]
    \ar[
      rr,
      dashed,
      "{\eta}"{description}
    ]
    \ar[
      from=uur,
      hook,
      crossing over
    ]
    &&
    X
    \\
    U
    \ar[
      u,
      hook,
      "{
        (\mathrm{id}, 1)
      }"{left}
    ]
    \ar[
      rr,
      dashed
    ]
    &&
    A
    \ar[
     u,
     hook
    ]
  \end{tikzcd}
$$

  \vspace{-.2cm}
  \item[{\bf (c)}]
a continuous function $\phi \,\colon\, \TopologicalSpace \xrightarrow{\;} [0,1]$
making the following diagram commute:
\vspace{-2mm}
$$
  \begin{tikzcd}[row sep=small]
    A
    \ar[
      d,
      hook
    ]
    \ar[rr]
    &&
    \{0\}
    \ar[
      d,
      hook
    ]
    \\
    X
    \ar[
      rr,
      dashed,
      "{\phi}"{description}
    ]
    &&
    {[0,1]}
    \\
    X \setminus U
    \ar[
      u,
      hook
    ]
    \ar[rr]
    &&
    \{1\}\;.
    \ar[
      u,
      hook
    ]
  \end{tikzcd}
$$
\end{itemize}
\end{proposition}

\begin{lemma}[Composition preserves h-cofibrations (e.g. {\cite[Ex. 4.2.17]{AguilarGitlerPrieto02}})]
  \label{CompositionPreservesHCofibrations}
  If $\mathrm{A} \xrightarrow{i_1} \mathrm{B}$
  and
  $\mathrm{B} \xrightarrow{i_2} \TopologicalSpace$ are h-cofibrations
  (Ntn. \ref{HurewiczClosedCofibrations})
  then so is their composite $i_2 \circ i_1$.
\end{lemma}

\begin{lemma}[Products preserve h-cofibrations]
  \label{ProductsPreserveClosedCofibrations}
  For $\mathrm{A} \xrightarrow{i \,\in\, \HurewiczCofibrations} \TopologicalSpace$
  an h-cofibration (Ntn. \ref{HurewiczClosedCofibrations}) and any
  $\mathrm{Y} \,\in\, \kTopologicalSpaces$, also their product is a closed cofibration:
\vspace{-2mm}
  $$
    \mathrm{Y} \times \mathrm{A}
    \xhookrightarrow{ \;\; \mathrm{id} \times i \,\in\, \HurewiczCofibrations \;\; }
    \mathrm{Y} \times \TopologicalSpace
    \,.
  $$
\end{lemma}
\begin{proof}
  By Prop. \ref{EquivalentCharacterizationsOfClosedCofibrations}
  it is sufficient to show that from a retraction $r$ of $i \PushoutProduct 0$
  we obtain a retraction of $(\mathrm{id}_{\mathrm{Y}} \times i) \PushoutProduct 0$.
  But since the product operation $\mathrm{Y} \times (-)$
  is a left adjoint \eqref{MappingSpaceAdjunction},
  it preserves the pushout-product:
  \vspace{-1mm}
  $$
    (\mathrm{id}_{\mathrm{Y}} \times i) \PushoutProduct 0
    \;=\;
    \mathrm{id}_{\mathrm{Y}} \times ( i \PushoutProduct 0 )
    \,.
  $$

  \vspace{-1mm}
\noindent  Therefore, the required retraction is
  given by $\mathrm{id}_{\mathrm{Y}} \times r$.
\end{proof}

\begin{definition}[Good simplicial topological space {\cite[Def. B.4]{Segal74}}]
  \label{GoodSimplicialTopologicalSpace}
  {\bf (i)}
  A simplicial topological space
  $\TopologicalSpace_\bullet \,\in\, \SimplicialTopologicalSpaces$
  \eqref{SimplicialTopologicalSpaces} is
  called {\it good} if all its degeneracy maps
  are h-cofibrations (Ntn. \ref{HurewiczClosedCofibrations}, Prop. \ref{EquivalentCharacterizationsOfClosedCofibrations}).
  \vspace{-2mm}
  $$
    \underset{
      {n \in \mathbb{N}}
      \atop
      { 0 \leq i \leq n }
    }{\forall}
    \TopologicalSpace_n \xrightarrow{\;\; \sigma_i \,\in\, \HurewiczCofibrations \;\;} \TopologicalSpace_{n+1}
    \,.
  $$

  \vspace{-1mm}
  \noindent
  {\bf (ii)} We denote the full subcategory of good simplicial topological spaces by:
  \vspace{-2mm}
  $$
    \SimplicialTopologicalSpaces_{\good}
    \xhookrightarrow{\quad}
    \SimplicialTopologicalSpaces \;.
  $$
\end{definition}

\begin{notation}[Well-pointed simplicial topological group]
  \label{WellPointedTopologicalGroup}
  {\bf (i)} A topological group $\Gamma \,\in\, \Groups(\kTopologicalSpaces)$
  is called {\it well-pointed}
  (e.g., \cite[\S VII, Def. 1.8]{Bredon93})
  if the inclusion of its neutral element
  is an h-cofibration (Ntn. \ref{HurewiczClosedCofibrations}):
  \vspace{-2mm}
  \begin{equation}
    \label{TheGoodPointOfAWellPointedTopologicalGroup}
    \{\mathrm{e}\}
    \xhookrightarrow{\;\;\; \in\, \HurewiczCofibrations \;\;\;}
    \Gamma
    \,.
  \end{equation}
  \noindent {\bf (ii)}
  We denote the full subcategory of well-pointed topological groups by:
  $$
    \Groups(\kTopologicalSpaces)_{\wellpointed} \;
    \xhookrightarrow{\quad}
    \Groups(\kTopologicalSpaces)
    \,.
  $$
 \noindent {\bf (iii)}   More generally (see Rem. \ref{WellPointedGroupsAreThoseWithGoodNerves}),
  a simplicial topological group
  $\mathcal{G}_\bullet \,\in\, \Groups( \SimplicialTopologicalSpaces)$
  is called {\it well-pointed} if
  all component maps of the
  the unique homomorphism from
  $1 \,\coloneqq\, \mathrm{const}(\{\mathrm{e}\})_\bullet \,\in\,\Groups( \SimplicialTopologicalSpaces)$
  are h-cofibrations:
  \vspace{-2mm}
  \begin{equation}
    \label{TheGoodPointsOfAWellPointedSimplicialTopologicalGroup}
    \underset{n \in \mathbb{N}}{\forall}
    \;\;\;
    \{e\}
    \;
   \xhookrightarrow{\;\;\; \in\, \HurewiczCofibrations \;\;\;}
      \mathcal{G}_n
      \mathrlap{\,.}
  \end{equation}

  \vspace{-1mm}
  \noindent {\bf (iv)} We denote the full subcategory of well-pointed simplicial topological groups by
  $$
    \Groups(\SimplicialTopologicalSpaces)_{\wellpointed} \;
    \xhookrightarrow{\quad}
    \Groups(\SimplicialTopologicalSpaces)
    \,.
  $$
\end{notation}

\begin{proposition}[Banach Lie groups are well-pointed]
  \label{BanachLieGroupsAreWellPointed}
  Every (paracompact Hausdorff) Banach Lie group,
  hence in particular every ordinary Lie group, is well-pointed
  (Ntn. \ref{WellPointedTopologicalGroup}).
\end{proposition}
\begin{proof}
  Paracompact Banach manifolds are ``absolute neighborhood retracts'' (ANRs),
  by \cite[Cor. to Thm. 5, p. 3]{Palais66};
  and closed inclusions of ANRs are h-cofibrations,
  by \cite[Thm. 4.2.15]{AguilarGitlerPrieto02}.
  The claim follows since
  points are ANRs, trivially, and are closed by assumption of Hausdorffness.
\end{proof}

\medskip

\noindent
{\bf Infinite projective groups}
The following infinite-projective groups are examples of a well-pointed groups
which are {\it not} special cases of the class in Prop. \ref{BanachLieGroupsAreWellPointed}:
\begin{example}[Projective groups on countably-dimensional Hilbert spaces]
  \label{ProjectiveUnitarGroupOnAHilbertSpace}
  We write
  \vspace{-1mm}
  \begin{equation}
    \label{TheGroupUH}
    \UH
    \;:=\;
    \UnitaryGroup(\HilbertSpace)
    \;\;
    \,\in\, \Groups(\kTopologicalSpaces)
  \end{equation}

  \vspace{-1mm}
  \noindent
  for the group of unitary operators on any
  countably infinite-dimensional complex Hilbert space
  equipped with its strong operator topology
  (which is equal to the weak operator topology
  \cite[Cor. 9.4]{HilgertNeeb92}
  as well as to the
  compact-open topology \cite{EspinozaUribe14}\cite{Schottenloher18},
  but strictly coarser
  than the norm topology that would make it a Banach Lie group).
  Notice that this group
  (which is famously contractible in the norm topology by Kuiper's theorem)
  is still contractible in the strong operator topology
  (\cite[Lem. 3 on p. 251]{DixmierDouady63}\cite[p. 4]{Schottenloher18}):
  \vspace{-1mm}
  \begin{equation}
    \label{UHIsContractible}
    \UH \;\underset{\mathrm{htpy}}{\simeq}\; \ast
    \,.
  \end{equation}

  \vspace{-2mm}
  \noindent
  {\bf (i)} We write
  \begin{equation}
    \label{TheGroupPUH}
    \PUH
    \;\coloneqq\;
    \UH/\CircleGroup
    \;\;\;
    \in
    \;
    \Groups(\kTopologicalSpaces)_{\wellpointed}
  \end{equation}
  for the topological quotient group of \eqref{TheGroupUH}
  by its subgroup of operators acting by multiplication with a complex number.
  This {\it projective unitary group} \eqref{TheGroupPUH}
  is well-pointed (Ntn. \ref{WellPointedTopologicalGroup}),
  see \cite[p. 23]{HebestreitSagave20}.
  The quotient coprojection is a locally trivial $\mathrm{U}(1)$-principal bundle
  \cite[Thm. 1]{Simms70}:
  \vspace{-2mm}
  \begin{equation}
    \label{PUHFiberSequence}
    \begin{tikzcd}[row sep=small]
      \CircleGroup
      \ar[r, hook]
      &
      \UH
      \ar[d]
      \\
      &
      \PUH
      \mathrlap{\,.}
    \end{tikzcd}
  \end{equation}

\vspace{-2mm}
\noindent
{\bf (ii)} Under the canonical $\ZTwo$-action by complex conjugation,
  these
  are compatibly
  $\ZTwo$-equivariant topological groups (Def. \ref{EquivariantTopologicalGroup})
\begin{equation}
  \label{ComplexConjugationActionOnProjectiveUnitaryGroup}
  \begin{tikzcd}
    \CircleGroup
    \ar[r, hook]
    \ar[out=180-66, in=66, looseness=3.5, "\scalebox{.77}{$\;\mathclap{
      \ZTwo
    }\;$}"{description},shift right=1]
    &
    \UH
    \ar[r, ->>]
    \ar[out=180-66, in=66, looseness=3.5, "\scalebox{.77}{$\;\mathclap{
      \ZTwo
    }\;$}"{description},shift right=1]
    &
    \PUH
    \ar[out=180-66, in=66, looseness=3.5, "\scalebox{.77}{$\;\mathclap{
      \ZTwo
    }\;$}"{description},shift right=1]
    \end{tikzcd}
    \;\in\;
    \Actions{\ZTwo}
    \left(
      \Groups(\kTopologicalSpaces)
    \right)
    \,.
\end{equation}
  Their fixed locus is the group of orthogonal operators on any
  countably infinity-dimensional {\it real} Hilbert space
  equipped with its operator topology,
  \begin{equation}
    \label{TheGroupOH}
    \OH
    \;\simeq\;
    (\UH)^{\ZTwo}
    \,\in\,
    \Groups(\kTopologicalSpaces)
    \,.
  \end{equation}
  In turn, the quotient of the latter
  by the subgroup
  of operators acting by multiplication with real units,
  is the infinite
  {\it projective orthogonal group}
  \cite[\S 3]{Ros}\cite{MMS}:
  \begin{equation}
    \label{TheGroupPOH}
    \POH
    \;=\;
    (\PUH)^{\ZTwo}
    \;=\;
    \OH/\ZTwo
    \;\;\;
    \in
    \;
    \Groups(\kTopologicalSpaces)_{\wellpointed}
    \,.
  \end{equation}
  This is again well-pointed, see \cite[p. 23]{HebestreitSagave20}.

\noindent
{\bf (iii)} These statements generalize to the $\ZTwo$-graded projective group
  $\GradedPUH$
  \cite[Prop. 2.2]{Parker88} (see also \cite[p. 5]{CareyWang08}),
  which is obtained from the graded unitary group
$$
  \GradedUH
  \;:=\;
  \Bigg\{
    \left(
    \!\!\!
    \arraycolsep=2pt
    \begin{array}{ll}
      \UnitaryOperator_{{}_{++}}
      &
      0
      \\
      0
      &
      \UnitaryOperator_{{}_{--}}
    \end{array}
    \!\!\!
    \right)
    ,\,
    \left(
    \!\!\!
    \arraycolsep=2pt
    \begin{array}{ll}
      0
      &
      \UnitaryOperator_{{}_{+-}}
      \\
      \UnitaryOperator_{{}_{-+}}
      &
      0
    \end{array}
    \!\!\!
    \right)
    \,\Bigg\vert\,
    U_{{}_{\bullet, \bullet}}
    \,\in\,
    \UnitaryGroup(\HilbertSpace)
  \Bigg\}
  \;\;
  \subset
  \;\;
  \UnitaryGroup
  \big(
    \HilbertSpace \otimes \ComplexNumbers^2
  \big)
  \,,
$$
that itself is a $\UH \times \UH$-extension of $\ZTwo$
\begin{equation}
  \label{GradedUnitaryGroupAsExtension}
  \begin{tikzcd}[row sep=-2pt]
    \UH^2
    \ar[rr, hook]
    &&
    \GradedUH
    \ar[rr, ->>, "{ c }"]
    &&
    \ZTwo
    \\
 \scalebox{0.7}{$   (U_{{}_{++}}, U_{{}_{--}}) $}
    &\longmapsto&
  \scalebox{0.7}{$    \left(
      \!\!\!\!\!
      \def\arraystretch{1.2}
      \begin{array}{l}
        U_{{}_{++}}
        \\
        0
      \end{array}
      \!\!
      \def\arraystretch{1.2}
      \begin{array}{l}
        0
        \\
        U_{{}_{--}}
      \end{array}
      \!\!\!\!\!\!
    \right)
    $}
    &\longmapsto&
 \scalebox{0.7}{$     \NeutralElement $}
    \\
    &&
 \scalebox{0.7}{$     \left(
      \!\!\!\!\!
      \def\arraystretch{1.2}
      \begin{array}{l}
        0
        \\
        U_{{}_{+-}}
      \end{array}
      \!\!
      \def\arraystretch{1.2}
      \begin{array}{l}
        U_{{}_{-+}}
        \\
        0
      \end{array}
      \!\!\!\!\!\!
    \right)
    $}
    &\longmapsto&
    \scalebox{0.7}{$  \mathrm{odd} $}
  \end{tikzcd}
\end{equation}
by quotienting out the
diagonal subgroup $\CircleGroup \xhookrightarrow{\;} \UH \xhookrightarrow{\;} \GradedUH$:
\vspace{-2mm}
\begin{equation}
  \label{TheGroupGradedPUH}
  \begin{tikzcd}
    \CircleGroup
    \ar[r, hook]
      &
    \GradedUH
    \ar[r, ->>]
      &
    \GradedPUH
    \mathrlap{\,.}
  \end{tikzcd}
\end{equation}

\vspace{-2mm}
\noindent
{\bf (iv)}
The $\ZTwo$-action \eqref{ComplexConjugationActionOnProjectiveUnitaryGroup}
evidently extends to these graded groups \eqref{TheGroupGradedPUH}
(acting trivially on the $\{\pm 1\}$-grading)
  $$
    \ZTwo \acts \, \GradedPUH
    \;\;\;
    \in
    \;
    \Actions{\ZTwo}\big(\Groups(\kTopologicalSpaces)_{\wellpointed}\big)
  $$
  with fixed locus being the analogous graded projective orthogonal group:
  $$
    \GradedPOH
    \;=\;
    (\GradedPUH)^{\ZTwo}
    \;\simeq\;
    \GradedOH/\ZTwo
    \;\;\;
    \in
    \;
    \Groups(\kTopologicalSpaces)
    \,.
  $$

\noindent
{\bf (v)}
This allows to form semidirect products of all these (projective, graded)
unitary groups with their complex conjugation action, yielding the following
system of short exact sequences of well-pointed topological groups:
\begin{equation}
  \label{ProjectiveGradedExtensionOfZTwoTimesZTwo}
  \begin{tikzcd}
    \CircleGroup
    \ar[rr,-,shift left=1pt]
    \ar[rr,-,shift right=1pt]
    \ar[d, hook]
    &&
    \CircleGroup
    \ar[rr, ->>]
    \ar[d, hook]
    &&
    1
    \ar[d]
    \\
    \UH^2
    \ar[rr, hook]
    \ar[d, ->>]
    &&
    \GradedUH \rtimes \ZTwo
    \ar[rr, ->>, "{ (c, \mathrm{pr}_2) }"]
    \ar[d, ->>]
    &&
    \ZTwo \times \ZTwo
    \ar[d,-,shift left=1pt]
    \ar[d,-,shift right=1pt]
    \\
    \UH^2/\CircleGroup
    \ar[rr, hook]
    &&
    \GradedPUH \rtimes \ZTwo
    \ar[rr, ->>, "{ (c,\, \mathrm{id}) }"]
    &&
    \ZTwo \times \ZTwo
    \,.
  \end{tikzcd}
\end{equation}

\noindent
{\bf (vi)}
The projective graded unitary group \eqref{TheGroupGradedPUH}
(and thus its semidirect product \eqref{ProjectiveGradedExtensionOfZTwoTimesZTwo})
has a canonical continuous action
on the following space of Fredholm operators\footnote{The invertibility
up to compact operators, in \eqref{TheSpaceOfFredholmOperators}, is equivalent
to the more traditional definition of Fredholm operators, by Atkinson's theorem,
e.g., \cite[Thm. 2.1]{Murphy94}. }
(\cite[Def. 3.2]{AtiyahSegal04}, see also \cite[Def. A.39]{FreedHopkinsTelement07TwistedKTheoryI}):
\begin{equation}
  \label{TheSpaceOfFredholmOperators}
  \hspace{-2mm}
  \FredholmOperators
  \coloneqq
  \left\{
    \FredholmOperator
    \,\in\,
    \mathcal{B}\big(\HilbertSpace_+ \oplus \HilbertSpace_-\big)
    \,\middle\vert\,
    \def\arraystretch{1.2}
    \begin{array}{l}
      \FredholmOperator^\dagger \,=\, \FredholmOperator,
      \\
      \mathrm{deg}(\FredholmOperator) \;=\; \mathrm{odd} \in \ZTwo,
      \\
      \FredholmOperator^2 - 1
        \,\in\,
      \CompactOperators
      \big(\HilbertSpace_+ \oplus \HilbertSpace_-\big)
    \end{array}
 \!\! \right\}
  \xhookrightarrow{\;
    \FredholmOperator
      \,\mapsto\,
    (
      \FredholmOperator
      , \,
      \FredholmOperator^2  - 1
    )
 \; }
  \BoundedOperators_{\mathrm{co}} \times \CompactOperators_{\mathrm{nrm}}
    \big(\HilbertSpace_+ \oplus \HilbertSpace_-\big)
\end{equation}
(where ``$\BoundedOperators_{\mathrm{co}}$'' denotes the space of
bounded operators with the compact-open topology and
``$\CompactOperators_{\mathrm{nrm}}$'' denotes
space of compact operators with the normal topology), namely by conjugation:
\begin{equation}
  \label{ConjugationActionOfProjectiveGradedUnitaryGroup}
  \begin{tikzcd}[row sep=-1pt, column sep=small]
    \big(
      \GradedPUH
      \,\rtimes\,
      \ZTwo
    \big)
      \times
    \FredholmOperators
    \ar[rr]
    &&
    \FredholmOperators
    \\
   \scalebox{0.7}{$   \Bigg(
    \bigg(
    \left[
    \!\!\!\!\!
    \def\arraystretch{1.3}
    \begin{array}{c}
      U_{{}_{++}}
      \\
      U_{{}_{-+}}
    \end{array}
    \!\!
    \def\arraystretch{1.3}
    \begin{array}{c}
      U_{{}_{+-}}
      \\
      U_{{}_{--}}
    \end{array}
    \!\!\!\!\!
    \right]
    ,
    \sigma
    \bigg)
    ,\,
    \left(
    \!\!\!\!\!
    \def\arraystretch{1.3}
    \begin{array}{l}
      0
      \\
      \FredholmOperator^\dagger_{{}_{+-}}
    \end{array}
    \!\!
    \def\arraystretch{1.3}
    \begin{array}{l}
      \FredholmOperator_{{}_{+-}}
      \\
      0
    \end{array}
    \!\!\!\!\!
    \right)
    \Bigg)
    $}
    &\longmapsto&
   \scalebox{0.7}{$   \left(
    \!\!\!\!\!
    \def\arraystretch{1.3}
    \begin{array}{c}
      U_{{}_{++}}
      \\
      U_{{}_{-+}}
    \end{array}
    \!\!
    \def\arraystretch{1.3}
    \begin{array}{c}
      U_{{}_{+-}}
      \\
      U_{{}_{--}}
    \end{array}
    \!\!\!\!\!
    \right)
    \circ
    \left(
    \!\!\!\!\!
    \def\arraystretch{1.3}
    \begin{array}{l}
      0
      \\
      \FredholmOperator^{\dagger \sigma}_{{}_{+-}}
    \end{array}
    \!\!
    \def\arraystretch{1.3}
    \begin{array}{l}
      \FredholmOperator^\sigma_{{}_{+-}}
      \\
      0
    \end{array}
    \!\!\!\!\!
    \right)
    \circ
    \left(
    \!\!\!\!\!
    \def\arraystretch{1.3}
    \begin{array}{c}
      U_{{}_{++}}
      \\
      U_{{}_{-+}}
    \end{array}
    \!\!
    \def\arraystretch{1.3}
    \begin{array}{c}
      U_{{}_{+-}}
      \\
      U_{{}_{--}}
    \end{array}
    \!\!\!\!\!
    \right)^{-1}
    $}
    \!.
  \end{tikzcd}
\end{equation}
(Here either
$\UnitaryOperator_{{}_{+-}}, \UnitaryOperator_{{}_{-+}} = 0$
or
$\UnitaryOperator_{{}_{++}}, \UnitaryOperator_{{}_{--}} = 0$;
the square bracket denotes the $\mathrm{diag}(\CircleGroup)$-equivalence class
of the matrix;
and $\FredholmOperator^{\sigma}$ equals $\FredholmOperator$ when $\sigma = \NeutralElement$
and equals its complex conjugate operator otherwise, see also  \cite[\S 5.B]{Matumoto71}.)
\end{example}

\medskip

\begin{lemma}[Good connected simplicial groups are well-pointed]
  \label{GoodConnectedSimplicialGroupsAreWellPointed}
  If $\mathcal{G}_\bullet \;\in\; \Groups(\SimplicialTopologicalSpaces)$
  is trivial in degree 0, $\mathcal{G}_0 \,=\, 1$, and such that
  its underlying simplicial space
  is good (Def. \ref{GoodSimplicialTopologicalSpace}),
  then it is well-pointed (Ntn. \ref{WellPointedTopologicalGroup}).
\end{lemma}
\begin{proof}
  For all $n \in \mathbb{N}$, the inclusion of the neutral element in that degree
  is the composition of a sequence of degeneracy maps with the inclusion of the
  neutral element in degree 0:
  \vspace{-2mm}
  $$
    \begin{tikzcd}[column sep=large]
      \{\mathrm{e}\}
      \ar[
        r,-,
        shift left=1pt
      ]
      \ar[
        r,-,
        shift right=1pt
      ]
      &[-15pt]
      \mathcal{G}_0
      \ar[
        r,
        "{
          \sigma
        }"{above},
        "{
          \in \, \HurewiczCofibrations
        }"{below}
      ]
      &
      \mathcal{G}_1
      \ar[
        r,
        "{
          \sigma
        }"{above},
        "{
          \in \, \HurewiczCofibrations
        }"{below}
      ]
      &
      \cdots
      \ar[
        r,
        "{
          \sigma
        }"{above},
        "{
          \in \, \HurewiczCofibrations
        }"{below}
      ]
      &
      \mathcal{G}_n
      \,.
    \end{tikzcd}
  $$

  \vspace{-2mm}
  \noindent  By assumption, all these morphisms are h-cofibrations, and hence
  so is their composite, by Lem. \ref{CompositionPreservesHCofibrations}.
\end{proof}
\begin{proposition}[Nerves of action groupoids of well-pointed topological group actions are good]
  \label{NervesOfActionGroupoidsOfWellPointedTopologicalGroupActionsAreGood}
  $\,$

  \noindent
  Let $\Gamma \,\in\, \Groups(\kTopologicalSpaces)$
  be well-pointed (Ntn. \ref{WellPointedTopologicalGroup}).
  Then for every $\TopologicalSpace \,\in\, \Actions{\Gamma}(\kTopologicalSpaces)$
  the nerve \eqref{SimplicialTopologicalNerveOfTopologicalGroupoids}
  of the topological action groupoid
  (Ex. \ref{TopologicalActionGroupoid}) is a good simplicial space
  (Def. \ref{GoodSimplicialTopologicalSpace}).
\end{proposition}
\begin{proof}
  This follows by Lemma \ref{ProductsPreserveClosedCofibrations}
  with the observation that all degeneracy maps in the action groupoid are
  of the form
  \vspace{-1mm}
  $$
    \mathrm{id}_{
      \TopologicalSpace \times \Gamma^{\times_n}
    }
    \times
    \big(\!
      \{\mathrm{e}\}
      \xhookrightarrow{ \in \, \HurewiczCofibrations}
      \Gamma
    \big)
    \;\;\;
    \in
    \;
    \in \HurewiczCofibrations
    \,.
  $$

  \vspace{-8mm}
\end{proof}
\begin{remark}[Well-pointed groups are those with good delooping]
  \label{WellPointedGroupsAreThoseWithGoodNerves}
  Specialized to the trivial action on the point,
  Prop. \ref{NervesOfActionGroupoidsOfWellPointedTopologicalGroupActionsAreGood}
  says that a topological group is well-pointed (Ntn. \ref{WellPointedTopologicalGroup})
  if and only if the nerve of its delooping groupoid (Ex. \ref{TopologicalDeloopingGroupoid})
  is good (Def. \ref{GoodSimplicialTopologicalSpace}),
  and if and only if it is well-pointed
  in the sense \eqref{TheGoodPointsOfAWellPointedSimplicialTopologicalGroup}
  of a (constant) simplicial topological group.
\end{remark}
Moreover:

\begin{lemma}[Well-pointed simplicial groups are good simplicial spaces {\cite[Prop. 3 (1)]{RobertsStevenson12}}]
  \label{WellPointedSimplicialGroupsAreGoodSimplicialSpaces}
  Underlying any well-pointed simplicial topological group
  (Ntn. \ref{WellPointedTopologicalGroup})
  is a good simplicial topological space (Ntn. \ref{GoodSimplicialTopologicalSpace}).
  \vspace{-2mm}
  $$
    \begin{tikzcd}[row sep=10pt]
      \SimplicialTopologicalGroups_{\wellpointed}
      \ar[d, hook]
      \ar[rr, dashed]
      &&
      \SimplicialTopologicalSpaces_{\good}
      \ar[d, hook]
      \\
      \SimplicialTopologicalGroups
      \ar[
        rr,
        "{
          \mbox{\tiny \color{greenii} \bf underlying}
        }"
      ]
      &&
      \SimplicialTopologicalSpaces
    \end{tikzcd}
  $$
\end{lemma}

\begin{proposition}[Good resolutions of simplicial topological spaces {\cite[p. 308-309]{Segal74}}]
  For every $\TopologicalSpace_\bullet \;\in\; \SimplicialTopologicalSpaces$
  there exists a good simplicial space
  $\TopologicalSpace^{\good}_\bullet \;\in\; \SimplicialTopologicalSpaces_{\good}$
  (Def. \ref{GoodSimplicialTopologicalSpace})
  and a morphism
    \vspace{-2mm}
  $$
    \TopologicalSpace_\bullet^{\good}
    \xrightarrow{\;\; \in \WeakEquivalences  \;\;}
    \TopologicalSpace_\bullet
  $$

  \vspace{-2mm}
 \noindent which is degreewise a weak homotopy equivalence.
\end{proposition}

\newpage
\noindent
{\bf Topological realization of good simplicial spaces.}

\begin{proposition}[Realization of well-pointed simplicial group is well-pointed
{\cite[Prop. 3]{RobertsStevenson12}\cite[Lem. 1]{BaezStevenson09}}]
  \label{RealizationOfWellPointedSimplicialGroupIsWellPointed}
  If  a simplicial topological group
   is well-pointed (Ntn. \ref{WellPointedTopologicalGroup}),
   then its topological realization \eqref{TopologicalRealizationOfSimplicialTopologicalSpaces}
   --
   with its induced structure \eqref{FunctorOnStructuresInducedFromLexFunctor}
   of a topological group, by Lem. \ref{TopologicalRealizationPreservesFiniteLimits}
   --
   is again well-pointed \eqref{TheGoodPointOfAWellPointedTopologicalGroup}:
     \vspace{-2mm}
   $$
     \begin{tikzcd}[row sep=10pt]
       \Groups(\SimplicialTopologicalSpaces)_{\wellpointed}
       \ar[
         rr,
         dashed
       ]
       \ar[
         d,
         hook
       ]
       &&
       \Groups(\kTopologicalSpaces)_{\wellpointed}
       \ar[
         d,
         hook
       ]
       \\
       \Groups(\SimplicialTopologicalSpaces)
       \ar[
         rr,
         "{
           \Groups( \vert - \vert)
         }"
       ]
       &&
       \Groups(\kTopologicalSpaces)
       \mathrlap{\,.}
     \end{tikzcd}
   $$
\end{proposition}

\begin{remark}[Topological realization of good simplicial spaces models their homotopy colimit]
\label{TopologicalRealizationOfGoodSimplicialSpacesModelsTheirHomotopyColimit}
The point of good simplicial spaces (Def. \ref{GoodSimplicialTopologicalSpace})
is that their topological realization (Ntn. \ref{TopologicalRealizationFunctors})
is weakly homotopy equivalent to the ``fat'' realization \cite[Prop. A.1 (iv)]{Segal74}
(also \cite[Prop. 1]{tomDieck74} with \cite[\S A]{RobertsStevenson12})
which, in turn, is a standard model for their homotopy colimit
(e.g., \cite[Ex. 6.4]{ArkhipovOrsted18}):
\vspace{-2mm}
\begin{equation}
  \label{GoodImpliesThatRealizationIsHomotopyColimit}
  \mathrm{X_\bullet}
  \;\;
  \in
  \;
  \SimplicialTopologicalSpaces_{\good}
  \;\;\;\;\;\;\;\;
  \Rightarrow
  \;\;\;\;\;\;\;\;
  \SingularSimplicialComplex
  \,
  \TopologicalRealization{}
  {
    \TopologicalSpace_\bullet
  }
  \;\simeq\;
  \underset{ [n] \in \Delta^{\mathrm{op}}  }{\mathrm{hocolim}}
  \left(
    \SingularSimplicialComplex
    (\TopologicalSpace_n)
  \right)
  \;\;\;
  \in
  \;
  \HomotopyCategory
  (\SimplicialSets_{\mathrm{Qu}})
  \,.
\end{equation}
\end{remark}
In view of this fact \eqref{GoodImpliesThatRealizationIsHomotopyColimit},
we quote the following statements about homotopy colimits of simplicial spaces
(Prop. \ref{SufficientConditionsForTopologicalRealizationOfSimplicialSpacesToPreserveHomotopyPullbacks})
in terms of topological realizations of good simplicial spaces:

\begin{proposition}[Sufficient conditions for topological realization to preserve homotopy fibers]
  \label{SufficientConditionsForTopologicalRealizationOfSimplicialSpacesToPreserveHomotopyPullbacks}
  Sufficient conditions for a homotopy fiber sequence
  (Ntn. \ref{HomtopyFiberSequenceOfTopologicalSpaces})
  of {\it good} simplicial topological spaces (Def. \ref{GoodSimplicialTopologicalSpace}),
  hence a degreewise homotopy fiber product of topological spaces,
  to remain a homotopy fiber sequence under topological realization
  (Ntn. \ref{TopologicalRealizationFunctors})
    \vspace{-2mm}
  $$
    \left(\!\!\!
    \begin{tikzcd}
      \mathrm{F}_\bullet
        \ar[
          rr,
          "{ i_\bullet }"{above},
          "{ \simeq \, \HomotopyFiber(p_\bullet) }"{below}
        ]
      &&
      \mathrm{E}_\bullet
        \ar[
          r,
          "p_\bullet"
        ]
      &
      \TopologicalSpace_\bullet
    \end{tikzcd}
    \;\;\;
    \in
    \;
    \SimplicialTopologicalSpaces_{\good}
    \right)
    \;\;
    \Rightarrow
    \;\;
    \left(\!\!\!
    \begin{tikzcd}
      \vert \mathrm{F}_\bullet \vert
        \ar[
          rr,
          "{ \vert i_\bullet \vert }"{above},
          "{ \simeq \, \HomotopyFiber(\vert p_\bullet\vert ) }"{below}
        ]
      &&
      \vert \mathrm{E}_\bullet \vert
        \ar[
          r,
          "\vert p_\bullet \vert"
        ]
      &
      \vert \TopologicalSpace_\bullet \vert
    \end{tikzcd}
    \;\;\;
    \in
    \;
    \kTopologicalSpaces
    \right)
  $$

    \vspace{-2mm}
\noindent
  include any of the following conditions on $p_\bullet$:

  \vspace{-2mm}
\begin{enumerate}[{\bf (i)}]
\setlength\itemsep{-4pt}

  \vspace{-.2cm}
  \item
  \cite[p. 2]{Anderson78}:

  \vspace{-2mm}
  \begin{itemize}

    \vspace{-.2cm}
    \item
    on simplicial sets of connected components, $\pi_0(p)_\bullet$ is a Kan fibration,

    \vspace{-.2cm}
    \item
    the component spaces $\mathrm{E}_n$, $\mathrm{Y}_n$ ($n \in \mathbb{N}$) are connected or discrete.

  \end{itemize}

  \vspace{-.2cm}
  \item
  \cite[Prop. 5.4]{Rezk14HomotopyColimits}\cite[Lem. 5.5.6.17]{Lurie17}:

  \vspace{-.1cm}
  \begin{itemize}

    \vspace{-.3cm}
    \item
    the simplicial set of connected components of the base is constant:
     $\pi_0(Y)_\bullet \,\simeq\, \mathrm{const}\big( \pi_0(Y_0) \big)_\bullet$.

  \end{itemize}

  \vspace{-.2cm}
  \item
  \cite[Cor. 6.7]{MazelGee14}\cite[Prop. 10]{Lurie11}:

  \vspace{-.1cm}
  \begin{itemize}

    \vspace{-.3cm}
    \item
    $p_\bullet$ is a homotopy Kan fibration:
    $p_\bullet \;\in\; \HomotopyKanFibrations$
    (Def. \ref{HomotopyKanFibrations}).
  \end{itemize}

  \end{enumerate}
\end{proposition}

\newpage

\chapter{Equivariant principal bundles}
\label{EquivariantPrincipalTopologicalBundles}

We give a streamlined account of
basic notions of topological equivariant principal bundles:

-- \cref{PrincipalBundlesInternalToTopologicalGActions}
  discusses the basic definitions.

-- \cref{NotionsOfEquivariantLocalTrivialization}
  discusses equivariant local triviality.

-- \cref{ConstructionOfUniversalEquivariantPrincipalBundles}
  constructs universal equivariant principal bundles.

\medskip

\section{As bundles internal to $G$-actions}
\label{PrincipalBundlesInternalToTopologicalGActions}

\noindent
{\bf Definition of equivariant principal topological bundles.}
We consider the definition of equivariant principal topological bundles
as principal bundles {\it internal} (Ntn. \ref{Internalization})
to
topological $G$-actions (Def. \ref{EquivariantPrincipalBundle} below),
extract what this means externally
(Cor.  \ref{InternalDefinitionOfGPrincipalBundlesCoicidesWithtomDieckDefinition} below),
and explain how this compares to previous definitions
found in the literature (Remark \ref{LiteratureOnEquivariantPrincipalBundles} below).

\begin{remark}[Assumption of local trivalizability]
\label{AssumptionOfLocalTrivializability}
Here we consider (formally) principal bundles which are not required to be
locally trivial, and we shall say ``principal {\it fiber} bundles'' for those that are,
discussed further below in \cref{NotionsOfEquivariantLocalTrivialization},
see Def. \ref{TerminologyForPrincipalBundles}.
Notice that this terminology in line with:

\noindent
{\bf (i)}  Cartan's original definition of
``principal bundle'', which did {\it not} include the
local triviality clause
(whence  \cite[Def. 1.1.2]{Palais61} speaks of ``Cartan principal bundles''
in this more general case),

\noindent
{\bf (ii)} tradition in equivariant topology,
where the equivariant local triviality condition is
more subtle (see \cref{NotionsOfEquivariantLocalTrivialization})
and not included
(see \cite{tomDieck69}\cite{Bierstone73}\cite{Lashof82})
in the bare definition of equivariant principal bundles.

In \cref{InCohesiveInfinityStacks} we show that this issue
is an artefact of internalization into
a 1-category of topological spaces,
instead of into an $\infty$-topos of higher geometric spaces:
When ($G$-equivariant) ordinary principal bundles are regarded internal
to (the $\mathbf{B}G$-slice of) the $\infty$-topos $\SmoothInfinityGroupoids$,
then their (equivariant)
local triviality is automatically implied
(Thm. \ref{OrdinaryPrincipalBundlesAmongPrincipalInfinityBundles},
Thm. \ref{BorelClassificationOfEquivariantBundlesForResolvableSingularitiesAndEquivariantStructure} below).
\end{remark}

\begin{definition}[Equivariant topological group]
  \label{EquivariantTopologicalGroup}
  We say that an {\it equivariant topological group}
  is a group object
  internal
  (Ntn. \ref{Internalization})
  to $\GActionsOnTopologicalSpaces$
  (Ntn. \ref{GActionOnTopologicalSpaces}), hence:

  \vspace{-3mm}
  \begin{enumerate}[{\bf (i)}]
  \setlength\itemsep{-1pt}
  \item
  an object
  $
    \begin{tikzcd}
      \Gamma
      \ar[out=180-66, in=66, looseness=3.5, "\scalebox{.77}{$\mathclap{
        G
      }$}"{description},shift right=1]
    \end{tikzcd}
      \in
    \;
    \GActionsOnTopologicalSpaces
  $;

  \vspace{-.2cm}
  \item
  morphisms
  $
    \begin{tikzcd}[column sep=large]
      \ast
      \ar[
        r,
        "\NeutralElement"{above},
        "\mbox{\tiny\color{greenii} \bf neutral element}"{below}
      ]
      &
      \Gamma
      \mathrlap{\,,}
      \ar[out=180-66, in=66, looseness=3.5, "\scalebox{.77}{$\mathclap{
        G
      }$}"{description},shift right=1]
    \end{tikzcd}
    {\phantom{AAAA}}
    \begin{tikzcd}[column sep=-5pt]
      \Gamma
      \ar[out=180-66, in=66, looseness=3.5, "\scalebox{.77}{$\mathclap{
        G
      }$}"{description},shift right=1]
      &
      \times
      &
      \Gamma
      \ar[out=180-66, in=66, looseness=3.5, "\scalebox{.77}{$\mathclap{
        G
      }$}"{description},shift right=1]
      \ar[
        rr,
        "m"{above},
        "\mbox{\tiny\color{greenii} \bf multiplication}"{below}
      ]
      &{\phantom{AAAAAA}}&
      \Gamma
      \mathrlap{\,,}
      \ar[out=180-66, in=66, looseness=3.5, "\scalebox{.77}{$\mathclap{
        G
      }$}"{description},shift right=1]
    \end{tikzcd}
    {\phantom{AAAA}}
    \begin{tikzcd}[column sep=large]
      \Gamma
      \ar[out=180-66, in=66, looseness=3.5, "\scalebox{.77}{$\mathclap{
        G
      }$}"{description},shift right=1]
      \ar[
        r,
        "(-)^{-1}"{above},
        "\mbox{\tiny\color{greenii} \bf inverses}"{below}
      ]
      &
      \Gamma
      \mathrlap{\,,}
      \ar[out=180-66, in=66, looseness=3.5, "\scalebox{.77}{$\mathclap{
        G
      }$}"{description},shift right=1]
    \end{tikzcd}
  $
  \end{enumerate}

   \vspace{-2mm}
\noindent
  such that the following diagrams commute:

\vspace{-2mm}
  \begin{enumerate}[{\bf (a)}]
       \setlength\itemsep{2pt}
    \item
    {\bf (Unitality)}
   $\quad$
    $
      \begin{tikzcd}[column sep=-5pt, row sep=small]
        \Gamma
        \ar[out=180-66, in=66, looseness=3.5, "\scalebox{.77}{$\mathclap{
          G
        }$}"{description},shift right=1]
        &
        \times
        \ar[
          d,
          "\simeq"
        ]
        &
        \ast
        \ar[
          rr,
          "\mathrm{id} \times \NeutralElement"
        ]
        &
        {\phantom{AAAA}}
        &
        \Gamma
        \ar[out=180-66, in=66, looseness=3.5, "\scalebox{.77}{$\mathclap{
          G
        }$}"{description},shift right=1]
        &
        \times
        \ar[
          d,
          "\, m"
        ]
        &
        \Gamma
        \ar[out=180-66, in=66, looseness=3.5, "\scalebox{.77}{$\mathclap{
          G
        }$}"{description},shift right=1]
        \\
        &
        \Gamma
        \ar[out=-180+66, in=-66, looseness=4.2, "\scalebox{.77}{$\mathclap{
          G
        }$}"{description},shift left=1]
        \ar[
          rrrr,
          -,
          shift left=1pt
        ]
        \ar[
          rrrr,
          -,
          shift right=1pt
        ]
        &&
        &&
        \Gamma
        \ar[out=-180+66, in=-66, looseness=4.2, "\scalebox{.77}{$\mathclap{
          G
        }$}"{description},shift left=1]
      \end{tikzcd}
    $
 $\quad$
    and
     $\quad$
    $
      \begin{tikzcd}[column sep=-5pt, row sep=small]
        \ast
        &
        \times
        \ar[
          d,
          "\simeq"
        ]
        &
        \Gamma
        \ar[out=180-66, in=66, looseness=3.5, "\scalebox{.77}{$\mathclap{
          G
        }$}"{description},shift right=1]
        \ar[
          rr,
          "\mathrm{id} \times \NeutralElement"
        ]
        &
        {\phantom{AAAA}}
        &
        \Gamma
        \ar[out=180-66, in=66, looseness=3.5, "\scalebox{.77}{$\mathclap{
          G
        }$}"{description},shift right=1]
        &
        \times
        \ar[
          d,
          "\, m"
        ]
        &
        \Gamma
        \ar[out=180-66, in=66, looseness=3.5, "\scalebox{.77}{$\mathclap{
          G
        }$}"{description},shift right=1]
        \\
        &
        \Gamma
        \ar[out=-180+66, in=-66, looseness=4.2, "\scalebox{.77}{$\mathclap{
          G
        }$}"{description},shift left=1]
        \ar[
          rrrr,
          -,
          shift left=1pt
        ]
        \ar[
          rrrr,
          -,
          shift right=1pt
        ]
        &&
        &&
        \Gamma
        \ar[out=-180+66, in=-66, looseness=4.2, "\scalebox{.77}{$\mathclap{
          G
        }$}"{description},shift left=1]
      \end{tikzcd}
    $
    \item
    {\bf (Associativity)}   $\quad$
    $
      \begin{tikzcd}[column sep =-5pt]
        \Gamma
        \ar[out=180-66, in=66, looseness=3.5, "\scalebox{.77}{$\mathclap{
          G
        }$}"{description},shift right=1]
        &\times&
        \Gamma
        \ar[out=180-66, in=66, looseness=3.5, "\scalebox{.77}{$\mathclap{
          G
        }$}"{description},shift right=1]
        \ar[
          d,
          "m \times \mathrm{id}"
        ]
        &\times&
        \Gamma
        \ar[out=180-66, in=66, looseness=3.5, "\scalebox{.77}{$\mathclap{
          G
        }$}"{description},shift right=1]
        \ar[
          rr,
          "\mathrm{id} \times m"
        ]
        &{\phantom{AAAAA}}&
        \Gamma
        \ar[out=180-66, in=66, looseness=3.5, "\scalebox{.77}{$\mathclap{
          G
        }$}"{description},shift right=1]
        &
        \times
        \ar[
          d,
          "m"
        ]
        &
        \Gamma
        \ar[out=180-66, in=66, looseness=3.5, "\scalebox{.77}{$\mathclap{
          G
        }$}"{description},shift right=1]
        \\
        &
        \Gamma
        \ar[out=-180+66, in=-66, looseness=4.2, "\scalebox{.77}{$\mathclap{
          G
        }$}"{description},shift left=1]
        &\times&
        \Gamma
        \ar[out=-180+66, in=-66, looseness=4.2, "\scalebox{.77}{$\mathclap{
          G
        }$}"{description},shift left=1]
        \ar[
          rrrr,
          "m"
        ]
        &
        &&
        &
        \Gamma
        \ar[out=-180+66, in=-66, looseness=4.2, "\scalebox{.77}{$\mathclap{
          G
        }$}"{description},shift left=1]
      \end{tikzcd}
    $
    \item
    {\bf (Invertibility)}   $\quad$
    $
      \begin{tikzcd}[column sep=-5pt, row sep=3pt]
        &{\phantom{AAA}}&
        \Gamma
        \ar[out=180-66, in=66, looseness=3.5, "\scalebox{.77}{$\mathclap{
          G
        }$}"{description},shift right=1]
        &\times&
        \Gamma
        \ar[out=180-66, in=66, looseness=3.5, "\scalebox{.77}{$\mathclap{
          G
        }$}"{description},shift right=1]
        \ar[
          rr,
          "(-)^{-1} \times \mathrm{id}"
        ]
        &{\phantom{AAAAAAA}}&
        \Gamma
        \ar[out=180-66, in=66, looseness=3.5, "\scalebox{.77}{$\mathclap{
          G
        }$}"{description},shift right=1]
        &
        \times
        \ar[
          dd,
          "m"
        ]
        &
        \Gamma
        \ar[out=180-66, in=66, looseness=3.5, "\scalebox{.77}{$\mathclap{
          G
        }$}"{description},shift right=1]
        \\
        \Gamma
        \ar[out=180-66, in=66, looseness=3.5, "\scalebox{.77}{$\mathclap{
          G
        }$}"{description},shift right=1]
        \ar[
          urr,
          "\mathrm{diag}"{near end}
        ]
        \ar[
          drrr
        ]
        \\
        &&&
        \ast
        \ar[
          rrrr,
          "\NeutralElement"
        ]
        &
        &&
        &
        \Gamma
        \ar[out=-180+66, in=-66, looseness=4.2, "\scalebox{.77}{$\mathclap{
          G
        }$}"{description},shift left=1]
      \end{tikzcd}
    $
      $\quad$
    and
      $\quad$
    $
      \begin{tikzcd}[column sep=-5pt, row sep=3pt]
        &{\phantom{AAA}}&
        \Gamma
        \ar[out=180-66, in=66, looseness=3.5, "\scalebox{.77}{$\mathclap{
          G
        }$}"{description},shift right=1]
        &\times&
        \Gamma
        \ar[out=180-66, in=66, looseness=3.5, "\scalebox{.77}{$\mathclap{
          G
        }$}"{description},shift right=1]
        \ar[
          rr,
          "\mathrm{id} \times(-)^{-1}"
        ]
        &{\phantom{AAAAAAA}}&
        \Gamma
        \ar[out=180-66, in=66, looseness=3.5, "\scalebox{.77}{$\mathclap{
          G
        }$}"{description},shift right=1]
        &
        \times
        \ar[
          dd,
          "m"
        ]
        &
        \Gamma
        \ar[out=180-66, in=66, looseness=3.5, "\scalebox{.77}{$\mathclap{
          G
        }$}"{description},shift right=1]
        \\
        \Gamma
        \ar[out=180-66, in=66, looseness=3.5, "\scalebox{.77}{$\mathclap{
          G
        }$}"{description},shift right=1]
        \ar[
          urr,
          "\mathrm{diag}"{near end}
        ]
        \ar[
          drrr
        ]
        \\
        &&&
        \ast
        \ar[
          rrrr,
          "\NeutralElement"
        ]
        &
        &&
        &
        \Gamma
        \ar[out=-180+66, in=-66, looseness=4.2, "\scalebox{.77}{$\mathclap{
          G
        }$}"{description},shift left=1]
      \end{tikzcd}
    $
  \end{enumerate}

  \vspace{-3mm}
  \noindent
  A homomorphism of equivariant groups is a morphism
  $
    \begin{tikzcd}
        \Gamma_1
         \ar[out=180-66, in=66, looseness=3.5, "\scalebox{.77}{$\mathclap{
          G
        }$}"{description},shift right=1]
           \ar[r,"f"]
        &
        \Gamma_1
 \ar[out=180-66, in=66, looseness=3.5, "\scalebox{.77}{$\mathclap{
          G
        }$}"{description},shift right=1]
    \end{tikzcd}
        \in
        \GActionsOnTopologicalSpaces
  $
  which makes commuting diagrams with the structure morphisms,
  in the evident way. This yields the category
  \vspace{-2mm}
  \begin{equation}
    \label{CategoryOfGEquivariantTopologicalGroups}
    \GEquivariantTopologicalGroups
    \;:=\;
    \mathrm{Grps}
    \left(
      \GActionsOnTopologicalSpaces
    \right)
    \,.
  \end{equation}
\end{definition}

\begin{definition}[Equivariant principal bundle]
  \label{EquivariantPrincipalBundle}
  We say that
  a (topological) {\it $G$-equivariant princial bundle}
  with {\it equivariant structure group}
  $G \acts \, \Gamma \,\in\, \Groups\big( G\mathrm{Actions}(\kTopologicalSpaces)\big)$
  (Def. \ref{EquivariantTopologicalGroup})
  is a {\it formally principal} bundle
  (\cite[p. 312 (15 of 30)]{Grothendieck60}\cite[p. 9 (293)]{Grothendieck71},
  also: {\it pseudo-torsor} \cite[\S 16.5.15]{Grothendieck67},
  see Rem. \ref{PseudoTorsorCondition})
  internal
  (Ntn. \ref{Internalization})
  to $\GActionsOnTopologicalSpaces$
  (Ntn. \ref{GActionOnTopologicalSpaces}),
  hence:

  \vspace{-3mm}
 \begin{enumerate}[{\bf (i)}]
\setlength\itemsep{-2pt}
  \item
  {\bf (Bundle)}
  a morphism
  \vspace{-2mm}
  \begin{equation}
    \label{InternalBundleInGActions}
    \begin{tikzcd}
      \mathrm{P}
      \ar[out=180-66, in=66, looseness=3.5, "\scalebox{.77}{$\mathclap{
        G
      }$}"{description},shift right=1]
      \ar[rr,"p"]
      &&
      \TopologicalSpace
      \ar[out=180-66, in=66, looseness=3.5, "\scalebox{.77}{$\mathclap{
        G
      }$}"{description},shift right=1]
    \end{tikzcd}
    \phantom{AAAA}
    \in \, \GActionsOnTopologicalSpaces \;.
  \end{equation}

  \item
  {\bf (Action)}
  an internal left $(G \acts \, \Gamma)$-action
  (e.g. \cite[\S V.6]{MacLaneMoerdijk92}\cite[p. 8]{BorceuxJanelidzeKelly05}),
  on its total space,
  namely a morphism
  \begin{equation}
    \label{InternalActionInGActions}
    \begin{tikzcd}[column sep=-5pt]
      \Gamma
      \ar[out=180-66, in=66, looseness=3.5, "\scalebox{.77}{$\mathclap{
        G
      }$}"{description},shift right=1]
      &\times&
      \mathrm{P}
      \ar[out=180-66, in=66, looseness=3.5, "\scalebox{.77}{$\mathclap{
        G
      }$}"{description},shift right=1]
      \ar[
        rr,
        "\rho"
      ]
      &{\phantom{AAAA}}&
      \mathrm{P}
      \ar[out=180-66, in=66, looseness=3.5, "\scalebox{.77}{$\mathclap{
        G
      }$}"{description},shift right=1]
    \end{tikzcd}
    \phantom{AAAA}
    \in \, \GActionsOnTopologicalSpaces
  \end{equation}
  making the following diagrams commute:
  \vspace{-2mm}
  $$
    \mbox{\bf(unitality)}
    \;\;
      \begin{tikzcd}[column sep=-5pt]
        \ast
        &
        \times
        \ar[
          d,
          "\simeq"
        ]
        &
        P
        \ar[out=180-66, in=66, looseness=3.5, "\scalebox{.77}{$\mathclap{
          G
        }$}"{description},shift right=1]
        \ar[
          rr,
          "\NeutralElement \times \mathrm{id}"
        ]
        &
        {\phantom{AAAA}}
        &
        \Gamma
        \ar[out=180-66, in=66, looseness=3.5, "\scalebox{.77}{$\mathclap{
          G
        }$}"{description},shift right=1]
        &
        \times
        \ar[
          d,
          "\rho"
        ]
        &
        P
        \ar[out=180-66, in=66, looseness=3.5, "\scalebox{.77}{$\mathclap{
          G
        }$}"{description},shift right=1]
        \\
        &
        P
        \ar[out=-180+66, in=-66, looseness=4.2, "\scalebox{.77}{$\mathclap{
          G
        }$}"{description},shift left=1]
        \ar[
          rrrr,
          -,
          shift left=1pt
        ]
        \ar[
          rrrr,
          -,
          shift right=1pt
        ]
        &&
        &&
        P
        \ar[out=-180+66, in=-66, looseness=4.2, "\scalebox{.77}{$\mathclap{
          G
        }$}"{description},shift left=1]
      \end{tikzcd}
      \phantom{AAAAAA}
    \mbox{\bf(action property)}
    \;\;
      \begin{tikzcd}[column sep =-5pt]
        \Gamma
        \ar[out=180-66, in=66, looseness=3.5, "\scalebox{.77}{$\mathclap{
          G
        }$}"{description},shift right=1]
        &\times&
        \Gamma
        \ar[out=180-66, in=66, looseness=3.5, "\scalebox{.77}{$\mathclap{
          G
        }$}"{description},shift right=1]
        \ar[
          d,
          "m \times \mathrm{id}"
        ]
        &\times&
        P
        \ar[out=180-66, in=66, looseness=3.5, "\scalebox{.77}{$\mathclap{
          G
        }$}"{description},shift right=1]
        \ar[
          rr,
          "\mathrm{id} \times \mathrm{\rho}"
        ]
        &{\phantom{AAAAA}}&
        \Gamma
        \ar[out=180-66, in=66, looseness=3.5, "\scalebox{.77}{$\mathclap{
          G
        }$}"{description},shift right=1]
        &
        \times
        \ar[
          d,
          "\rho"
        ]
        &
        P
        \ar[out=180-66, in=66, looseness=3.5, "\scalebox{.77}{$\mathclap{
          G
        }$}"{description},shift right=1]
        \\
        &
        \Gamma
        \ar[out=-180+66, in=-66, looseness=4.2, "\scalebox{.77}{$\mathclap{
          G
        }$}"{description},shift left=1]
        &\times&
        P
        \ar[out=-180+66, in=-66, looseness=4.2, "\scalebox{.77}{$\mathclap{
          G
        }$}"{description},shift left=1]
        \ar[
          rrrr,
          "\rho"
        ]
        &
        &&
        &
        P
        \ar[out=-180+66, in=-66, looseness=4.2, "\scalebox{.77}{$\mathclap{
          G
        }$}"{description},shift left=1]
      \end{tikzcd}
    $$
  \end{enumerate}

  \vspace{-2mm}
  \noindent
  such that, moreover:

\vspace{-3mm}
  \begin{enumerate}[{\bf (a)}]
  \setlength\itemsep{-8pt}
  \item
  {\bf (Fiberwise action)}
  the action $\rho$ \eqref{InternalActionInGActions}
  is fiberwise relative to $p$ \eqref{InternalBundleInGActions},
  in that the following diagram commutes:
  \vspace{-2mm}
  $$
    \begin{tikzcd}[column sep=-1pt, row sep=small]
      \Gamma
      \ar[out=180-66, in=66, looseness=3.5, "\scalebox{.77}{$\mathclap{
        G
      }$}"{description},shift right=1]
      &
      \times
      \ar[
        drr,
        "p \,\circ\, \mathrm{pr}_2\;\;\;"{below, xshift=-4pt}
      ]
      &
      \mathrm{P}
      \ar[out=180-66, in=66, looseness=3.5, "\scalebox{.77}{$\mathclap{
        G
      }$}"{description},shift right=1]
      \ar[
        rr,
        "\rho"
      ]
      &{\phantom{AAAAAAA}}&
      P
      \ar[out=180-66, in=66, looseness=3.5, "\scalebox{.77}{$\mathclap{
        G
      }$}"{description},shift right=1]
      \ar[
        dl,
        "p"{below, xshift=3pt}
      ]
      \\
      &&&
      \TopologicalSpace
      \ar[out=-180+66, in=-66, looseness=4.2, "\scalebox{.77}{$\mathclap{
        G
      }$}"{description},shift left=1]
    \end{tikzcd}
  $$
  \item
  {\bf (Principality)}
  The resulting commuting square is
  a pullback square (shown on the left), meaning equivalently that the
  fiberwise {\it shear map}
  (the universal morphism factoring through the fiber product, shown on the right) is an isomorphism:
  \vspace{-5mm}
  \begin{equation}
    \label{PrincipalityConditionAsShearMapBeingAnIsomorphism}
    \begin{tikzcd}[column sep=-5pt]
      \Gamma
      \ar[out=180-66, in=66, looseness=3.5, "\scalebox{.77}{$\mathclap{
        G
      }$}"{description},shift right=1]
      &
      \times
      \ar[
        d,
        "\mathrm{pr}_2"{left}
      ]
      \ar[
        rrrd,
        phantom,
        "\mbox{\tiny\rm(pb)}"
      ]
      &
      \mathrm{P}
      \ar[out=180-66, in=66, looseness=3.5, "\scalebox{.77}{$\mathclap{
        G
      }$}"{description},shift right=1]
      \ar[
        rr,
        "\rho"
      ]
      &
      {\phantom{AAAAA}}
      &
      \mathrm{P}
      \ar[out=180-66, in=66, looseness=3.5, "\scalebox{.77}{$\mathclap{
        G
      }$}"{description},shift right=1]
      \ar[
        d,
        "p"
      ]
      \\
      &
      \mathrm{P}
      \ar[out=-180+66, in=-66, looseness=4.2, "\scalebox{.77}{$\mathclap{
        G
      }$}"{description},shift left=1]
      \ar[
        rrr,
        "p"{below}
      ]
      &
      &&
      \TopologicalSpace
      \ar[out=-180+66, in=-66, looseness=4.2, "\scalebox{.77}{$\mathclap{
        G
      }$}"{description},shift left=1]
    \end{tikzcd}
    {\phantom{AAAA}}
    \Leftrightarrow
    {\phantom{AAAA}}
    \begin{tikzcd}[column sep=-5pt]
      \Gamma
      \ar[out=180-66, in=66, looseness=3.5, "\scalebox{.77}{$\mathclap{
        G
      }$}"{description},shift right=1]
      \ar[
        ddrrrrr,
        bend right=20,
        "\mathrm{pr}_2"{left}
      ]
      &
      \times
      \ar[
        drrr,
        dashed,
        "\sim"{above, sloped},
        "
          \mbox{
            \tiny
            \color{greenii}
            \bf
            shear map
          }
        "{below, sloped}
      ]
      &
      \mathrm{P}
      \ar[out=180-66, in=66, looseness=3.5, "\scalebox{.77}{$\mathclap{
        G
      }$}"{description},shift right=1]
      \ar[
        rrrrrrd,
        bend left=17,
        "\rho"{above}
      ]
      &&
      \\
      &&
      &{\phantom{AAAA}}&
      P
      \ar[out=180-66, in=66, looseness=3.5, "\scalebox{.77}{$\mathclap{
        G
      }$}"{description},shift right=1]
      &
      \times_{\TopologicalSpace}
      \ar[
        d,
        "\mathrm{pr}_2"{left}
      ]
      \ar[
        rrrd,
        phantom,
        "\mbox{\tiny\rm(pb)}"
      ]
      &
      \mathrm{P}
      \ar[out=180-66, in=66, looseness=3.5, "\scalebox{.77}{$\mathclap{
        G
      }$}"{description},shift right=1]
      \ar[
        rr,
        "\mathrm{pr}_1"
      ]
      &
      {\phantom{AAAAA}}
      &
      \mathrm{P}
      \ar[out=180-66, in=66, looseness=3.5, "\scalebox{.77}{$\mathclap{
        G
      }$}"{description},shift right=1]
      \ar[
        d,
        "p"
      ]
      \\
      &&
      &&
      &
      \mathrm{P}
      \ar[out=-180+66, in=-66, looseness=4.2, "\scalebox{.77}{$\mathclap{
        G
      }$}"{description},shift left=1]
      \ar[
        rrr,
        "p"{below}
      ]
      &
      &&
      \TopologicalSpace
      \ar[out=-180+66, in=-66, looseness=4.2, "\scalebox{.77}{$\mathclap{
        G
      }$}"{description},shift left=1]
    \end{tikzcd}
  \end{equation}
\end{enumerate}
We denote the category of these objects by:
\begin{equation}
  \label{CategoryOfTopologicalEquivariantPrincipalBundles}
  \EquivariantPrincipalBundles{G}{\Gamma}(\kTopologicalSpaces)
  \;\;
  \coloneqq
  \;\;
  \FormallyPrincipalBundles{(G\acts \, \Gamma)}
  \big(
    \Actions{G}(\kTopologicalSpaces)
  \big)
  \,.
\end{equation}
\end{definition}

We extract what the internal definition (Def. \ref{EquivariantPrincipalBundle}) of
equivariant principal bundles means externally:
\begin{lemma}[$G$-Equivariant groups are semidirect products with $G$]
  \label{EquivariantTopologicalGroupsAreSemidirectProductsWithG}
  There is a fully faithful functor \eqref{FullyFaithfulInfinityFunctor}
   \vspace{-2mm}
  \begin{equation}
    \label{IdentifyingEquivariantGroupsWithSemidirectProductGroups}
    \begin{tikzcd}[row sep=-5pt]
      \GEquivariantTopologicalGroups
      \ar[
        rr,
        hook
      ]
      &&
      \TopologicalGroups^{G/}_{/G}
      \\
 \scalebox{0.8}{$      \Gamma
      \ar[out=-180+50, in=-50, shift left=1, looseness=3.7,
        "\scalebox{.77}{$\mathclap{
          \;\; G \;\;
        }$}"{description},
        "\alpha\;\;\;"{below,very near start}
      ]
      $}
     \ar[
        rr,
        phantom,
        "\qquad \longmapsto"
      ]
      &&
    \scalebox{0.8}{$   \Gamma \rtimes_{\alpha} G $}
    \end{tikzcd}
  \end{equation}

   \vspace{-2mm}
\noindent   which identifies
  $G$-equivariant topological groups (Def. \ref{EquivariantTopologicalGroup})
  with the
  topological semidirect product groups with $G$,
  regarded as pointed objects
  in the slice over $G$, via the canonical group homomorphisms
  \vspace{-2mm}
  \begin{equation}
    \label{SplitGroupExtensionOfGByGamma}
    \begin{tikzcd}[row sep=-1pt, column sep=large]
   \scalebox{0.8}{$ g $}
      \ar[
        r,
        phantom,
        "\longmapsto"{description}
      ]
      &
      \scalebox{0.8}{$
        (\NeutralElement, g)
      $}
      &&
      \scalebox{0.8}{$
        g
      $}
      \ar[
        ll,
        phantom,
        "\longmapsfrom"
      ]
      \\
      G
      \ar[
        r,
        hook,
        "{i}"
      ]
      &
      \Gamma \rtimes_\alpha G
      \ar[
        rr,
        ->>,
        shift right=5pt,
        "\mathrm{pr}_2"{description}
      ]
      &&
      G
      \ar[
        ll,
        hook,
        shift right=5pt,
        "s"{description}
      ]
      \\
      &
      \scalebox{0.8}{$
        (\gamma,g)
      $}
      \ar[
        rr,
        phantom,
        "\longmapsto"{
          description
        }
      ]
      &&
      \scalebox{0.8}{$
        g
      $}
    \end{tikzcd}
  \end{equation}

     \vspace{-2mm}
\noindent
  that jointly exhibit the semidirect product as a split
  group extension of
  $G$ by $\Gamma$.
\end{lemma}
\begin{proof}
  This is a matter of straightforward unwinding of the definitions:

  First, observe that an equivariant topological group in the internal
  sense of Def. \ref{EquivariantTopologicalGroup} is equivalently a
  topological group $\Gamma$ equipped with a continuous left $G$-action
  by group automorphisms $\alpha : G \xrightarrow{\;} \mathrm{Aut}_{\mathrm{Grp}}(\Gamma)$.

  In terms of this identification,
  the functor is defined on objects by sending the pair $(\Gamma,\alpha)$ to the semidirect product $\Gamma\rtimes_{\alpha}G$.
  The functor sends a morphism $f:(\Gamma,\alpha)\to (\Gamma',\alpha')$, determined by an equivariant continuous homomorphism $f:\Gamma\to \Gamma'$, to the map
  $f_* : \Gamma\rtimes_{\alpha}G\to \Gamma'\rtimes_{\alpha'}G$, defined by $f_*(\gamma,g)= (f(\gamma),g)$, which is a morphism in the pointed slice over $G$.
  This association is clearly faithful. Any morphism $h:\Gamma\rtimes_{\alpha}G\to\Gamma'\rtimes_{\alpha'}G$ in the slice over $G$ determines a map  $h':\Gamma\to \Gamma'$. The group structure on the semidirect product forces $h'$ to be an  equivariant homomorphism. Hence the functor is also full.
\end{proof}

\begin{lemma}[Equivariant group actions seen externally]
  Under the identification from Lemma \ref{EquivariantTopologicalGroupsAreSemidirectProductsWithG},
  actions of $G$-equivariant topological groups
  $(\Gamma,\alpha)$ (Def. \ref{EquivariantTopologicalGroup})
  on $G$-equivariant topological spaces
  $(\TopologicalSpace,\rho)$ (Ntn. \ref{GActionOnTopologicalSpaces})
  are equivalently actions of the corresponding semidirect product groups
  \eqref{IdentifyingEquivariantGroupsWithSemidirectProductGroups}
  on the underlying plain topological space:
     \vspace{-2mm}
  \begin{equation}\label{FunctorGammaActionstoSemidirectProduct}
  \begin{tikzcd}[row sep=-5pt, column sep=small]
    (
      G \acts \, \Gamma
    )
    \mathrm{Act}
    \left(
      \GActionsOnTopologicalSpaces
    \right)
    \ar[
      rr,
      "\sim"
    ]
    &&
    (
      \Gamma\rtimes_{\alpha} G
    )
    \mathrm{Act}
    (
      \kTopologicalSpaces
    )
    \\
    \scalebox{0.8}{$    (
      \Gamma,
      \alpha
    )
    \xrightarrow{R}
    \mathrm{Aut}
    (
      \TopologicalSpace, \rho
    )
    $}
    &\longmapsto&
    \scalebox{0.8}{$  \Gamma
     \rtimes_\alpha
    G
    \xrightarrow{(R,\rho)}
    \mathrm{Aut}(X)
    $}
    \,,
  \end{tikzcd}
  \end{equation}

     \vspace{-1mm}
\noindent
  where the semidirect product group action on the right is
  $
    (R,\rho)(\gamma,g)(x)
    \;\coloneqq\;
    R(\gamma)
    \left(
      \rho(g)(x)
    \right)
    \,.
  $
\end{lemma}
\begin{proof}
  This is again a matter of straightforward unwinding of the definitions. The morphism $R$ can be
  identified with a continuous homomorphism $R:\Gamma\to {\rm Aut}(X)$ that is equivariant
  with respect to $\alpha$ and the conjugacy action on automorphisms. Given an action by the
  semidirect product $\Gamma \rtimes_{\alpha}G$ on $X$, we can recover $R$ by restricting the
   action to elements $(\gamma,1)\in \Gamma\rtimes_{\alpha}G$. This map is indeed equivariant
   since for all $x\in X$,
       \vspace{-2mm}
  $$
  \left(\alpha(g)\gamma\right) x=\left(\alpha(g)\gamma,1\right)x=(1,g)\left(\gamma,g^{-1}\right)x
  =g\gamma \left(g^{-1}x\right)\;.
  $$

    \vspace{-2mm}
\noindent
With these identifications, it is trivial to verify that the
functor \eqref{FunctorGammaActionstoSemidirectProduct} is fully faithful.
\end{proof}

\begin{corollary}[External description of equivariant principal bundles]
  \label{InternalDefinitionOfGPrincipalBundlesCoicidesWithtomDieckDefinition}
  For $(G \acts \, \Gamma)$ a $G$-equivariant group, hence
  a topological group $\Gamma$ equipped with an action
  $\alpha : G \xrightarrow{\;} \mathrm{Aut}(\Gamma)$, the
  category of equivariant $(G \acts  \,\Gamma)$-principal bundles
  in the internal sense of Def. \ref{EquivariantPrincipalBundle}
  is equivalently described as follows:

  \vspace{-3mm}
  \begin{enumerate}[{\bf (i)}]
  \setlength\itemsep{-4pt}

    \item {\bf (Bundle)}
    Objects are topological bundles, i.e., continuous functions
    $\!\!
      \begin{tikzcd}[column sep=small]
        \mathrm{P}
         \ar[r,"p"]
         &
         \TopologicalSpace\;.
        \end{tikzcd}
      $

    \item {\bf (Action)}
    equipped with continuous actions
        \vspace{-3mm}
    $$
      \begin{tikzcd}[row sep=-2pt]
        (
          \Gamma \rtimes_\alpha G
        )
        \times \mathrm{P}
        \ar[
          rr,
          "(-)\cdot(-)"
        ]
        &&
        \mathrm{P}\;,
        \\
        G
        \times \TopologicalSpace
        \ar[
          rr,
          "(-)\cdot(-)"
        ]
        &&
        \TopologicalSpace\;,
      \end{tikzcd}
    $$

        \vspace{-4mm}
\noindent
    such that

\vspace{-3mm}
    \begin{enumerate}[{\bf (a)}]
   \setlength\itemsep{-6pt}

  \item

     {\bf (Fiberwise action)}
       the bundle projection is $G$-equivariant and
       the $\Gamma$-action is $G$-equivariant and fiberwise,
       in that this square commutes:
       \vspace{-4mm}
       \begin{equation}
        \label{ExternalEquivariancePropertyOfEquivariantPrincipalBundle}
       \begin{tikzcd}
         (
           \Gamma \rtimes_\alpha G
         )
           \times
         \mathrm{P}
         \ar[
           rr,
           "(-)\cdot (-)"
         ]
         \ar[
           d,
           "\mathrm{pr}_2 \times p"{left}
         ]
         &&
         \mathrm{P}
         \ar[
           d,
           "p"
         ]
         \\
         G \times \TopologicalSpace
         \ar[
           rr,
           "(-)\cdot (-)"
         ]
         &&
         \TopologicalSpace
       \end{tikzcd}
      \end{equation}

      \item
      {\bf (Principality)} the fiberwise shear map is an isomorphism
      (i.e., a homeomorphism)
         \vspace{-2mm}
          \begin{equation}
            \label{FiberwiseShearMapIsomorphism}
            \begin{tikzcd}[row sep=-2pt]
              \Gamma \times \mathrm{P}
              \ar[
                rr,
                "\sim"
              ]
              &&
              \mathrm{P} \times_{{}_{\TopologicalSpace}} \mathrm{P}
              \\
       \scalebox{0.8}{$         (\gamma, p) $}
              \ar[
                rr,
                phantom,
                "\longmapsto"{
                  description
                }
              ]
              &&
        \scalebox{0.8}{$        (\gamma \cdot p, p) $}
              \,.
            \end{tikzcd}
          \end{equation}
      \end{enumerate}
  \end{enumerate}

\vspace{-3mm}
\noindent     Morphisms
     between these objects
     are $\Gamma \rtimes_\alpha G$-equivariant continuous
     functions:
     \vspace{-2mm}
\begin{equation}
\label{MorphismsOfEquivariantBundlesExternally}
\begin{tikzcd}[row sep=small]
  \mathrm{P}_1
  \ar[out=180-66, in=66, looseness=3.5, "\scalebox{.77}{$\phantom{\cdot}\mathclap{
    \Gamma \rtimes_\alpha G
  }\phantom{\cdot}$}"{description}, shift right=1]
  \ar[
    rr,
    "f_{\mathrm{P}}"
  ]
  \ar[d]
  \ar[
    drr,
    phantom
  ]
  &&
  \mathrm{P}_2
  \ar[out=180-66, in=66, looseness=3.5, "\scalebox{.77}{$\phantom{\cdot}\mathclap{
    \Gamma \rtimes_\alpha G
  }\phantom{\cdot}$}"{description}, shift right=1]
  \ar[d]
  \\
  \TopologicalSpace_1
  \ar[out=-180+66, in=-66, looseness=3.5, "\scalebox{.77}{$\mathclap{
    G
  }$}"{description}, shift left=1]
  \ar[
    rr,
    "f_{\TopologicalSpace}"{below}
  ]
  &&
  \TopologicalSpace_2
  \ar[out=-180+66, in=-66, looseness=3.5, "\scalebox{.77}{$
    \mathclap{G}
  $}"{description}, shift left=1]
\end{tikzcd}
\end{equation}
\vspace{-.7cm}

\end{corollary}
For archetypical classes of examples
of equivariant principal bundles,
see Lemma \ref{EquivariantPrincipalTwistedProductBundles}
and Prop. \ref{TopologicalRealizationOfEquivariantPrincipalGroupoid}
below.

\begin{remark}[Comparison to the literature on equivariant principal bundles]
  \label{LiteratureOnEquivariantPrincipalBundles}
 The external description of equivariant principal bundles
 obtained
 in Cor.  \ref{InternalDefinitionOfGPrincipalBundlesCoicidesWithtomDieckDefinition}
 from the canonical internal definition in Def. \ref{EquivariantPrincipalBundle}:

\vspace{-1mm}
 \item {\bf (i)} recovers the definition
 of equivariant bundles due to
 \cite[\S 1.1]{tomDieck69}\cite[\S I.8]{tomDieck87},
 followed in e.g. \cite{Nishida78}\cite{MurayamaShimakawa95};

\vspace{-1mm}
\item {\bf (ii)} while much of the literature
(e.g. \cite{Bierstone73}\cite{Lashof82}\cite{LueckUribe14})
considers the definition recovered in
Cor.  \ref{InternalDefinitionOfGPrincipalBundlesCoicidesWithtomDieckDefinition}
only in the special case when $\alpha$
\eqref{IdentifyingEquivariantGroupsWithSemidirectProductGroups} is trivial,
hence when the semidirect product group $\Gamma \rtimes_\alpha G$
reduces to the direct product group
$\Gamma \times G$, hence when the action of
the structure group $\Gamma$  and the equivariance group $G$
commute with each other.

\vspace{-1mm}
\item {\bf (iii)} In the other direction,
\cite{LashofMay86}\cite[IV.1]{LewisMaySteinberger86}
proposed to consider group extensions of $G$ by $\Gamma$
more general than semidirect products (which are the split extensions);
but followups \cite{May90}\cite{GuillouMayMerling17}
fall back to the semidirect product group actions
originally considered in \cite{tomDieck69} and recovered in
Cor.  \ref{InternalDefinitionOfGPrincipalBundlesCoicidesWithtomDieckDefinition}
from our Def. \ref{EquivariantPrincipalBundle}.
\end{remark}

In fact, the internal definition subsumes one degenerate boundary case which is
traditionally disregarded, but which is important to include for the
equivariant theory to work well:

\begin{remark}[The role of formal principality]
  \label{PseudoTorsorCondition}
  Beware that
  the formal principality condition
  \cite[p. 312 (15 of 30)]{Grothendieck60}\cite[p. 9 (293)]{Grothendieck71}
  in
  Def. \ref{EquivariantPrincipalBundle}
  does not
  explicitly demand that the base space $\TopologicalSpace$
  be the quotient $\mathrm{P}/\Gamma$ of the total space by the structure group
  action, as would be demanded for actual principal bundles/torsors.
  We claim that this is the right definition for
  equivariant principal bundles, for two reasons:

  \vspace{-2mm}
  \begin{enumerate}[{\bf (i)}]
  \setlength\itemsep{-1pt}
  \item
  The internal Def. \ref{EquivariantPrincipalBundle} involves only finite limits
  (products and fiber product, but no colimits such as quotients),
  which guarantees (Prop. \ref{RightAdjointFunctorsPreserveFiberProducts})
  that every right adjoint functor
  (Ntn. \ref{AdjointFunctors}) on the ambient category of
  $\GActionsOnTopologicalSpaces$,
  hence notably the fixed locus functor (Ex. \ref{FixedLociWithResidualWeylGroupAction}),
  preserves these principal bundle objects
  (crucial in the form of Cor. \ref{FixedLociOfEquivariantPrincipalBundles} below).

  \item
  Once equivariant local triviality is imposed
  (in \cref{NotionsOfEquivariantLocalTrivialization} below)
  the quotient condition
  is both implied as well as circumvented, as need be.
  Namely, in that case the underlying bundle $\mathrm{P} \xrightarrow{\;} \TopologicalSpace$
  (Cor. \ref{UnderlyingPrincipalBundles}) is a locally trivial fiber bundle
  (Prop. \ref{tomDieckLocalTrivialityImpliesOrdinaryLocalTriviality}),
  whence there are two cases:

    \vspace{-3mm}
  \begin{enumerate}[{\bf (a)}]
    \setlength\itemsep{-1pt}
    \item
    If fibers are inhabited (i.e., not empty),
    it follows that $\mathrm{P} \xrightarrow{\;} \TopologicalSpace$ is an
    effective epimorphism (Ntn.  \ref{EffectiveEpimorphism}),
    in which case the shear map isomorphism
    \eqref{FiberwiseShearMapIsomorphism} {\it implies} the quotient condition:
    $\TopologicalSpace \simeq \mathrm{P}/\TopologicalSpace$ (Lem. \ref{EffectiveEquivariantPrincipalBundles}).

    \item
    If the fibers are empty
    then, while the quotient condition is of course not met
    (unless $\TopologicalSpace$ itself is empty), the bundle is still
    formally principal
    in the sense of Def. \ref{EquivariantPrincipalBundle}, Cor. \ref{InternalDefinitionOfGPrincipalBundlesCoicidesWithtomDieckDefinition},
    as the shear map is indeed an isomorphism

    \vspace{-.7cm}
    $$
      \begin{tikzcd}
        \Gamma \times \varnothing
        \ar[
          rr,
          "\exists !"{above},
          "\sim"{below}
        ]
        &&
        \varnothing \times_{\TopologicalSpace} \varnothing
      \end{tikzcd}
    $$
    \vspace{-.6cm}

    (since both its domain and codomain are the empty space). Hence:
    {\it Empty bundles are formally principal.}

    While this degenerate case is irrelevant in ordinary topological bundle theory
    and hence traditionally ignored
    (not so in algebraic geometry \cite[\S 16.5.15]{Grothendieck67}),
    its inclusion is crucial for equivariant principal bundle theory to work well,
    since the fixed loci (Cor. \ref{FixedLociOfEquivariantPrincipalBundles})
    of equivariant principal bundles are frequently empty
    (see Example \ref{FixedLociInBierstoneLocalModelBundles}).
  \end{enumerate}

  \end{enumerate}
  We see this effect clearly brought out in the
  universal equivariant principal bundles,
  Rem. \ref{PseudoPrincipalNatureOfFixedLociInUniversalEquivariantPrincipalBundles} below.
\end{remark}

\begin{corollary}[Fixed loci bundles of equivariant principal bundles]
  \label{FixedLociOfEquivariantPrincipalBundles}
  For $H \subset G$ a subgroup inclusion,
  passage to $H$-fixed loci constitutes a functor
  from $G$-equivariant $\Gamma$-principal bundles
  to $W\!(H)$-equivariant (Ntn. \ref{GActionOnTopologicalSpaces})
  $\Gamma^H$-principal bundles:
  \vspace{-5mm}
  $$
    \begin{tikzcd}
      \FormallyPrincipalBundles{(G\acts \,\Gamma)}
      (
        G \mathrm{Actions}
      )
      \ar[
        rr,
        "(-)^H"
      ]
      &&
      \FormallyPrincipalBundles{
      \left(
        W\!(H)
        \acts \, \Gamma^H
      \right)}
      \left(
        W\!(H) \mathrm{Act}
      \right)
      \;.
    \end{tikzcd}
  $$
\end{corollary}
\begin{proof}
  By Example \ref{FixedLociWithResidualWeylGroupAction}, the functor
  $(-)^H : G\mathrm{Actions} \xrightarrow{\;} W\!(H)\mathrm{Actions}$
  is a right adjoint, hence preserves finite limits
  (Prop. \ref{RightAdjointFunctorsPreserveFiberProducts}),
  hence induces \eqref{FunctorOnStructuresInducedFromLexFunctor}
  a functor of
  internal groups (Def. \ref{EquivariantTopologicalGroup}),
  and internal formally principal bundles (Def. \ref{EquivariantPrincipalBundle}).
\end{proof}

For completeness, we also record:

\begin{corollary}[Underlying principal bundles of equivariant principal bundles]
  \label{UnderlyingPrincipalBundles}
  Passage to underlying topological spaces (forgetting the equivariance group action)
  (Example \ref{ForgettingGActionsAsPullbackAction})
  constitutes a functor from equivariant to ordinary principal bundles:
    \vspace{-2mm}
  \begin{equation}
    \label{ForgetfulFunctorFromEquivariantToUnderlyingPrincipalBundles}
    \begin{tikzcd}
      \FormallyPrincipalBundles{(G \acts \,  \Gamma)}\left(\GActionsOnTopologicalSpaces\right)
      \ar[
        rr,
        "
          \mbox{
            \tiny
            \color{greenii}
            \rm \bf
            forget $G$-action
          }
        "{below}
      ]
      &&
      \FormallyPrincipalBundles{\Gamma}(\kTopologicalSpaces)
      \,.
    \end{tikzcd}
  \end{equation}
\end{corollary}
\begin{proof}
  This is immediate from inspection of the definitions, but it
  also follows abstractly, as in the proof of Cor. \ref{FixedLociOfEquivariantPrincipalBundles},
  from the fact that the forgetful functor \eqref{FreeForgetfulAdjunctionForGAction}
  is a right adjoint, hence preserves
  limits
  (Prop. \ref{RightAdjointFunctorsPreserveFiberProducts}),
  and therefore \eqref{FunctorOnStructuresInducedFromLexFunctor}
  preserves
  internal groups and internal principal bundles.
\end{proof}

\medskip

\noindent
{\bf Basic properties of equivariant principal topological bundles.}
We discuss some basic aspects of equivariant bundles related to
effective epimorphy.

The following Lem. \ref{EffectiveEquivariantPrincipalBundles}
is a slight strengthening of the fact that
effective epimorphisms in topological spaces are equivalently quotient maps,
which is worth recording
(and holds more generally for formally principal bundles internal
to any ambient category):
\begin{lemma}[Effective equivariant principal bundles]
  \label{EffectiveEquivariantPrincipalBundles}
  Given an equivariant principal bundle
  $G \acts \, \mathrm{P} \xrightarrow{p} G \acts \, \TopologicalSpace
    \,\in\,
    \FormallyPrincipalBundles{\Gamma}(\kTopologicalSpaces)\,$
  in the sense of
  Def. \ref{EquivariantPrincipalBundle}
  (i.e. not requiring any local trivializability)
  the following two conditions are equivalent:

\noindent  {\bf (i)} $\mathrm{P} \xrightarrow{p} \TopologicalSpace$
     is the $\Gamma$-quotient coprojection.

\noindent  {\bf (ii)} $\mathrm{P} \xrightarrow{p} \TopologicalSpace$
    is an effective epimorphism (Ntn.  \ref{EffectiveEpimorphism}).
\end{lemma}
\begin{proof}
  By definition, the two conditions mean equivalently that
  the following vertical diagrams on the left or right, respectively,
  are coequalizer diagrams

  \vspace{-4mm}
  $$
    \begin{tikzcd}[row sep=small]
      \Gamma
        \times
      \mathrm{P}
      \ar[
        rr,
        "\sim"{above,yshift=-1pt},
        "{
          \mbox{
            \tiny
            \color{greenii}
            \bf
            shear map
          }
        }"{below}
      ]
      \ar[
        d,
        shift right=5pt,
        "{\rho}"{left}
      ]
      \ar[
        d,
        shift left=5pt,
        "\, {\mathrm{pr}_2}"{right}
      ]
      &&
      \mathrm{P} \times_{\TopologicalSpace} \mathrm{P}
      \ar[
        d,
        shift right=5pt,
        "{ \mathrm{pr}_1 }"{left}
      ]
      \ar[
        d,
        shift left=5pt,
        "\, { \mathrm{pr}_2 }"{right}
      ]
      \\
      \mathrm{P}
      \ar[
        rr,-,
        shift left=1pt
      ]
      \ar[
        rr,-,
        shift right=1pt
      ]
      \ar[
        d,
        "\, p"
      ]
        &&
      \mathrm{P}
      \ar[
        d,
        "\, p"
      ]
      \\
      \TopologicalSpace
      \ar[
        rr,-,
        shift left=1pt
      ]
      \ar[
        rr,-,
        shift right=1pt
      ]
        &&
      \TopologicalSpace
    \end{tikzcd}
  $$

  \vspace{-2mm}
\noindent  But the shear map \eqref{PrincipalityConditionAsShearMapBeingAnIsomorphism}
  is readily seen to provide a morphism between the top diagrams
  (makes the two parallel squares commute) and is an isomorphism
  by principality. Therefore the two vertical diagrams are equivalent,
  and one is a coequalizer diagram if and only if the other is.
\end{proof}

The following may seem obvious, but does need an argument given that
we define principality via shear isomorphisms:
\begin{lemma}[Locally trivial $\Gamma$-actions are principal bundles]
  \label{LocallyTrivialActionsArePrincipal}
   Let $\Gamma \,\in\, \Groups(\kTopologicalSpaces)$
   and $\Gamma \acts \, \TopologicalPrincipalBundle \,\in\, \Actions{\Gamma}(\kTopologicalSpaces)$
   such that
   $\TopologicalPrincipalBundle \xrightarrow{q} \TopologicalPrincipalBundle/\Gamma
   =: \TopologicalSpace$
   is locally, over an open cover, $\Gamma$-equivariantly isomorphic to a
   trivial $\Gamma$-principal bundle. Then $\TopologicalPrincipalBundle$
   is itself a (locally trivial) $\Gamma$-principal bundle, in that its
   shear map is an isomorphism.
\end{lemma}
\begin{proof}
Let $\{ \TopologicalPatch_i \xhookrightarrow{\;} \TopologicalSpace\}_{i \in I}$
be a trivializing open cover and abbreviate $\widehat{\TopologicalSpace}
\coloneqq \underset{i \in I}{\sqcup} \TopologicalPatch_i$.
Consider the following pasting diagram of pullbacks,
whose left side follows by the assumption of local triviality:
\vspace{-2mm}
$$
  \begin{tikzcd}[row sep=small, column sep={between origins, 60pt}]
    \Gamma
      \times
    (
      \widehat{\TopologicalSpace} \times \Gamma
    )
    \ar[rr]
    \ar[dr, "\sim"{sloped}]
    \ar[ddr, bend right=24]
    \ar[drrr, phantom, "\mbox{\tiny\rm(pb)}"{pos=.4}]
    &&
    \Gamma \times \TopologicalPrincipalBundle
    \ar[dr, "\mathrm{shear}"{sloped}]
    \ar[ddr, bend right=24]
    \\
    &
    (\widehat \TopologicalSpace \times \Gamma)
    \underset{\widehat{\TopologicalSpace}}{\times}
    (\widehat \TopologicalSpace \times \Gamma)
    \ar[rr, ->>, crossing over]
    \ar[d]
    \ar[drr, phantom, "\mbox{\tiny\rm(pb)}"{pos=.3}]
    &{}&
    \TopologicalPrincipalBundle
      \times_{\TopologicalSpace}
    \TopologicalPrincipalBundle
    \ar[d]
    \\
    &
    \widehat X
    \ar[rr, ->>]
    &&
    \TopologicalSpace
    \,.
  \end{tikzcd}
$$

\vspace{-2mm}
\noindent By local recognition of homeomorphisms (Ex. \ref{IsomorphismOfBundlesDetectedOnOpenCovers})
it follows that the shear map itself is an isomorphism.
\end{proof}

\section{Equivariant local triviality}
\label{NotionsOfEquivariantLocalTrivialization}

We consider
equivariant principal bundles
which are equivariantly locally trivial
(in Thm. \ref{EquivalentNotionsOfEquivariantLocalTriviality}, Def. \ref{TerminologyForPrincipalBundles} below,
recall Rem. \ref{AssumptionOfLocalTrivializability})
and prove that
isomorphism classes
these equivariant principal {\it fiber bundles}
coincide with concordance classes
(Thm. \ref{ConcordanceClassesOfTopologicalPrincipalBundles} below).
This is a key ingredient in the proof of the classification theorem
for equivariant principal bundles fiber bundles in
Thm. \ref{BorelClassificationOfEquivariantBundlesForResolvableSingularitiesAndEquivariantStructure}
below.

\medskip

\noindent
{\bf Notions of equivariant local triviality.}
There is a curious subtlety in the notion of
equivariant local trivializability
(whose resolution is provided by the perspective of higher geometry
in  \cref{AsInfinityBundlesInternalToSliceOverBG} below):
While internalization (Ntn. \ref{Internalization})
of the notion of principal bundles into equivariant topology
is the method of choice
for understanding the plain definition of equivariant principal bundles
(as discussed in \cref{PrincipalBundlesInternalToTopologicalGActions}),
it {\it fails} when it comes to understanding their local triviality.

\medskip
That is, using that $\Actions{G}(\kTopologicalSpaces)$ is a regular category
(Prop. \ref{CompactlyGeneratedTopologicalGActionsFormARegularCategory}),
internal local triviality of
an equivariant principal bundle $G \acts \, \TopologicalPrincipalBundle$
would be the
existence of an effective epimorphism
$G \acts \, \widehat{\mathrm{X}} \twoheadrightarrow G \acts \, \TopologicalSpace$
and of a pullback diagram of the form
\vspace{-2mm}
\begin{equation}
  \label{NaiveEquivariantLocalTrivialization}
  \begin{tikzcd}[column sep=-9pt, row sep=small]
    \Gamma
      \ar[out=180-66, in=66, looseness=3.5, "\scalebox{.77}{$\mathclap{
        G
      }$}"{description},
      "\scalebox{.6}{$\alpha$}"{pos=.3, yshift=-5pt}
      shift right=1]
    &\times&
    \widehat{\TopologicalSpace}
      \ar[out=180-66, in=66, looseness=2.3, "\scalebox{.77}{$\mathclap{
        G
      }$}"{description},shift right=1]
    \ar[d, "\mathrm{pr}_2"{swap}]
    \ar[dr, phantom, "\mbox{\tiny\rm(pb)}"]
    \ar[r]
    &[50pt]
    \TopologicalPrincipalBundle
      \ar[out=180-66, in=66, looseness=3.5, "\scalebox{.77}{$\mathclap{
        G
      }$}"{description},shift right=1]
    \ar[d]
    \\
    &&
    \widehat{\TopologicalSpace}
      \ar[shift left=4pt, out=-180+66, in=-66, looseness=3.5, "\scalebox{.77}{$\mathclap{
        G
      }$}"{description},shift right=1]
    \ar[r, ->>]
    &
    \TopologicalSpace
      \ar[shift left=4pt, out=-180+66, in=-66, looseness=3.5, "\scalebox{.77}{$\mathclap{
        G
      }$}"{description},shift right=1]
  \end{tikzcd}
  \;\;\;\;
  \in
  \;
  \Actions{G}(\kTopologicalSpaces)
  \,,
\end{equation}
(such that the top morphism is also $\Gamma$-equivariant).
That this notion of local triviality would be too restrictive
in practice is readily seen in the special case when the action
$\alpha$ of $G$ on $\Gamma$ is trivial: In this case
\eqref{NaiveEquivariantLocalTrivialization}
essentially says that the $G$-action is
trivial on total spaces of equivariant principal bundles,
relative to the action on the base.

\medskip
Putting the systematics of internalization aside for the moment,
one may guess what the ``correct'' equivariant local trivialization condition
should be, and several authors have done so. Below we review the existing
definitions in streamlined form, generalize them where necessary to the case when the
structure group carries a non-trivial $G$-action, and then show that the
resulting definitions are all equivalent to each other (Thm. \ref{EquivalentNotionsOfEquivariantLocalTriviality}).
Further below in \cref{AsInfinityBundlesInternalToSliceOverBG}
we show that these equivariant local triviality conditions
do follow from systematic internalization after all:
though not in the  1-category of topological spaces but in a
higher category of, in particular, D-topological stacks.

\medskip
One way to motivate the correct definition of equivariant local
trivialization (without yet passing to $\infty$-topos theory)
is to observe that equivariantly there
is no longer a single local model for principal bundles:
A trivial underlying
$\Gamma$-principal bundle $\Gamma \times \mathrm{U} \to \mathrm{U}$
in general carries several inequivalent $G$-actions that make it
an equivariant $(G \acts \, \Gamma)$-principal bundle.
We consider three different ways of parametrizing the possible
equivariant local model bundles, which we name after the original
authors who put these models into focus:

\begin{center}
\begin{tabular}{|c|c|c|}
  \hline
  \multicolumn{3}{|c|}{
    \bf
    Local model equivariant principal bundles
  }
  \\
  \hline
  \hline
  \begin{tabular}{c}
    {tom Dieck}
    \\
    \small
    \cref{tomDieckLocalTrivializations}
  \end{tabular}
  &
  \begin{tabular}{c}
    {Lashof-May}
    \\
    \small
    \cref{LashofMyLocalTrivialization}
  \end{tabular}
  &
  \begin{tabular}{c}
    {Bierstone}
    \\
    \small
    \cref{BierstoneLocalTrivializations}
  \end{tabular}
  \\
  \hline
  \hline
  \begin{tikzcd}[row sep=small]
    \widehat G / \widehat H
    \ar[d]
    \\
    G/H
  \end{tikzcd}
  &
  \begin{tikzcd}[row sep=small]
    \widehat G \times_{\widehat H} \mathrm{S}
    \ar[d]
    \\
    G \times_H \mathrm{S}
  \end{tikzcd}
  &
  \begin{tikzcd}[row sep=small]
   \widehat G_x \times_{\widehat {G_x}} \mathrm{U}_x
   \ar[d]
   \\
   G_x \times_{G_x}
   \mathrm{U}_x
  \end{tikzcd}
  \\
  \hline
\end{tabular}
\end{center}

\medskip

Arguably, tom Dieck's definition is the most conceptual
(here one derives the general form of the local model bundles, instead of
prescribing them, see Prop. \ref{CharacterizationOfEquivariantBundlesOverCosetSpaces} below)
and also the most general (its underlying local triviality
follows for $G$ any topological group, by using the first item of Lemma
\ref{CosetSpaceCoprojectionsAdmittingLocalSections} in
Lemma \ref{SemidirectProductCosetBundles}),
while Bierstone's definition makes manifest
(Lem. \ref{EquivariantDirectProductBundlesFromHatGxTwistedProductBundles} below)
that the underlying bundles are indeed locally trivial,
only that the $G$-action on the local models is not, in general.
In particular, this makes manifest that
equivariant local triviality is preserved by passage to fixed loci
(Prop. \ref{BierstoneLocalTrivialityIsPreservedByPassageToHFixedLoci} below).
Finally, Bierstone's definition (in the special case when $\alpha$ is trivial)
is the one traditionally used (albeit without attribution to Bierstone)
in the seminal example of the equivariant degree-3 twisting bundle of complex
K-theory, following \cite[p. 28]{AtiyahSegal04}.

\medskip
Notice that \cite{Bierstone73} and \cite{Lashof82} considered only the special case
when $\alpha$ in \eqref{IdentifyingEquivariantGroupsWithSemidirectProductGroups} is trivial;
we generalize their definitions and proofs;
see Remark \ref{BierstoneLocalModelBundlesGeneralized} and
Remark \ref{ComparingLocalTrivializationProofOfLashofMayBundlesToTheLiterature} below.
Moreover, the relation to tom Dieck's original notion of equivariant local trivialization
seems to not have been discussed before.

\begin{theorem}[Equivalent notions of equivariant local triviality]
\label{EquivalentNotionsOfEquivariantLocalTriviality}
Under the assumption of proper equivariance (Assump. \ref{ProperEquivariantTopology}),
the three notions of local trivializability
of equivariant principal bundles (Def. \ref{EquivariantPrincipalBundle}) are all equivalent:
\vspace{0mm}
$$
  \begin{tikzcd}[column sep=-10pt, row sep=12pt]
    \mbox{
      \def\arraystretch{1}
      \begin{tabular}{c}
      \fbox{
        \rm
        tom Dieck
        }
        \\
        \tiny
        \rm
        (Def. \ref{tomDieckLocalTrivializability})
      \end{tabular}
    }
    \ar[
      rr,
      Rightarrow,
      shift left=7pt,
      "
        \mbox{
          \tiny
          \rm
          Prop. \ref{tomDieckLocalTrivializationImpliesBierstone}
        }
      "
    ]
    &&
    \mbox{
      \def\arraystretch{1}
      \begin{tabular}{c}
      \fbox{
        \rm
        Bierstone}
        \\
        \tiny
        \rm
        (Def. \ref{BierstoneEquivariantLocalTrivializability})
      \end{tabular}
    }
    \ar[
      dl,
      Rightarrow,
      start anchor={[xshift=+10pt, yshift=8pt]},
      "\mbox{
        \tiny
        \rm
        Prop. \ref{BierstoneLocalTrivialityImpliesLashofMay}
      }"{swap, yshift=-1pt, sloped, pos=.45}
    ]
    \\
    &
    \mbox{
      \def\arraystretch{1}
      \begin{tabular}{c}
       \fbox{ \rm
        Lashof-May }
        \\
        \tiny
        \rm
        (Def. \ref{LashofMayEquivariantLocalTrivializability})
      \end{tabular}
    }
    \ar[
      ul,
      Rightarrow,
      end anchor={[xshift=-10pt, yshift=8pt]},
      "
        \mbox{
          \tiny
          \rm
          Prop. \ref{LashoMayLocalTrivialityImpliestomDieck}
        }
      "{swap, yshift=-1pt, sloped}
    ]
  \end{tikzcd}
$$
\end{theorem}

\begin{definition}[Equivariant principal fiber bundles]
\label{TerminologyForPrincipalBundles}
We write
\vspace{-2mm}
\begin{equation*}
  \hspace{-5mm}
  \begin{tikzcd}
    \underset{
      \raisebox{-3pt}{
        \tiny
        \color{darkblue}
        \bf
        \begin{tabular}{c}
          equivariantly locally trivial equivariant principal bundles
        \end{tabular}
      }
    }{
      \EquivariantPrincipalFiberBundles{G}{\Gamma}(\kTopologicalSpaces)
    }
    \hspace{-5mm}
    \ar[r, hook]
    &
    \underset{
      \raisebox{-3pt}{
        \tiny
        \color{darkblue}
        \bf
        equivariant principal bundles
      }
    }{
      \EquivariantPrincipalBundles{G}{\Gamma}(\kTopologicalSpaces)
    }
    \;\;
    :=
    \;\;
    \underset{
      \raisebox{-3pt}{
        \tiny
        \color{darkblue}
        \bf
        \begin{tabular}{c}
          formally principal internal actions
        \end{tabular}
      }
    }{
      \FormallyPrincipalBundles{(G \acts \, \Gamma)}
      \left(
        \GActionsOnTopologicalSpaces
      \right)
    }
  \end{tikzcd}
\end{equation*}
for the full subcategory of \eqref{CategoryOfTopologicalEquivariantPrincipalBundles}
on those equivariant principal bundles
(Def. \ref{EquivariantPrincipalBundle}, see Rem. \ref{AssumptionOfLocalTrivializability})
which satisfy the equivalent equivariant local triviality conditions of
Thm. \ref{EquivalentNotionsOfEquivariantLocalTriviality}.
\end{definition}

\medskip

\noindent
{\bf Immediate consequences of equivariant local trivialization.}

\begin{proposition}[Passage to fixed loci preserves equivariant local triviality]
  \label{PassageToFixedLociPreservesEquivariantLocalTriviality}
  For $H \subset G$ any closed subgroup, the
  functor (Cor. \ref{FixedLociOfEquivariantPrincipalBundles})
  sending $G$-equivariant principal bundles to their
  $\WeylGroup{H}$-equivariant $H$-fixed loci bundles
  respects
  the subcategories in Def. \ref{TerminologyForPrincipalBundles}
  of equivariant principal fiber bundles
  (Def. \ref{TerminologyForPrincipalBundles}):
\begin{equation}
  \label{PassageToFixedLociPreservesFullSubcategoryOftomDieckLocallyTrivialEquivariantBundles}
  \hspace{-4.5mm}
  \begin{tikzcd}[column sep=large]
    \overset{
      \raisebox{4pt}{
        \tiny
        \color{darkblue}
        \bf
        \begin{tabular}{c}
          equivariantly locally trivial
        \end{tabular}
      }
    }{
      \EquivariantPrincipalFiberBundles{G}{\Gamma}(\kTopologicalSpaces)
    }\!\!
    \ar[
      r,
      hook
    ]
    \ar[
      d,
      "(-)^H"{right}
    ]
    &
\!\!    \overset{
      \raisebox{4pt}{
        \tiny
        \color{darkblue}
        \bf
        \begin{tabular}{c}
          equivariant principal bundles
        \end{tabular}
      }
    }{
      \EquivariantPrincipalBundles{G}{\Gamma}(\kTopologicalSpaces)
    }
    \ar[
      d,
      "(-)^H"{right},
      "{
        \mbox{
          \tiny
          \color{greenii}
          \bf
          \begin{tabular}{c}
            passage to
            \\
            $H$-fixed loci
          \end{tabular}
        }
      }"{left}
    ]
    \\
    \EquivariantPrincipalFiberBundles{\WeylGroup(H)}{\Gamma^H}(\kTopologicalSpaces)
    \ar[
      r,
      hook
    ]
    &
    \EquivariantPrincipalBundles{\WeylGroup(H)}{\Gamma^H}(\kTopologicalSpaces)
  \end{tikzcd}
\end{equation}
\end{proposition}
\begin{proof}
  This is manifest in Bierstone's formulation
  (Prop. \ref{BierstoneLocalTrivialityIsPreservedByPassageToHFixedLoci})
  and
  hence holds generally, by Theorem \ref{EquivalentNotionsOfEquivariantLocalTriviality}.
\end{proof}

\begin{lemma}[Equivariant principal fiber bundles are proper equivariant Serre fibrations]
  \label{EquivariantLocallyTrivialPrincipalBundlesAreProperEquivariantFibrations}
  An equivariant principal bundle
  $G \acts \, \TopologicalPrincipalBundle
   \xrightarrow{\;p\;} G \acts  \, \TopologicalSpace$
  which is equivariantly locally trivial
  (Def. \ref{TerminologyForPrincipalBundles})
  is a fibration in the proper equivariant model structure
  (Prop. \ref{ProperEquivariantModelCategoryOfGSpaces}):
  \vspace{-2mm}
  $$
    \EquivariantPrincipalFiberBundles{G}{\Gamma}(\kTopologicalSpaces)
    \;\subset\;
    \ProperEquivariantSerreFibrations{G}
    \,.
  $$
\end{lemma}
\begin{proof}
  By Thm. \ref{EquivalentNotionsOfEquivariantLocalTriviality}
  and Prop. \ref{PassageToFixedLociPreservesEquivariantLocalTriviality},
  we have for all $H \,\subset G\,$ that
  $
    \WeylGroup{H} \acts  \, \TopologicalPrincipalBundle^H
      \xrightarrow{\;}
    \WeylGroup{H} \acts \, \TopologicalSpace
  $
  satisfies Bierstone's condition (Def. \ref{BierstoneEquivariantLocalTrivializability}).
By Lem. \ref{EquivariantDirectProductBundlesFromHatGxTwistedProductBundles}, this
  implies that the underlying principal bundle
  $\TopologicalPrincipalBundle^H \xrightarrow{\;} \TopologicalSpace^H$
  \eqref{ForgetfulFunctorFromEquivariantToUnderlyingPrincipalBundles}
  is locally trivial, hence a Serre fibration by
  Lem. \ref{LocallyTrivialBundlesAreSerreFibrations}.
\end{proof}

\begin{lemma}[Quotiented fiber product of equivariant principal fiber bundles is
  proper equivariant fibration]
  \label{QuotientedFiberProductOfEquivariantPrincipalFiberBundlesIsProperEquivariantFibration}
  Given a pair of equivariant principal fiber bundles
  $G\acts \, \TopologicalPrincipalBundle_1
   ,\,
   G\acts \, \TopologicalPrincipalBundle_2
   \;
   \in
   \;
   \EquivariantPrincipalFiberBundles{G}{\Gamma}(\kTopologicalSpaces)
  $
  (Def. \ref{TerminologyForPrincipalBundles}),
  the $G\acts \Gamma$-quotient of their fiber product
  is a proper equivariant fibration
  (according to Prop. \ref{ProperEquivariantModelCategoryOfGSpaces}):
  \vspace{-2mm}
  \begin{equation}
    \label{QuotientedFiberProductOfEquivariantPrincipalFiberBundlesEquivalentlyProperEquivariantFibration}
    \begin{tikzcd}
      \Quotient
      {
        (
          \TopologicalPrincipalBundle_1
            \times_{\TopologicalSpace}
          \TopologicalPrincipalBundle_2
        )
      }{
        \Gamma
      }
       \ar[out=180-66, in=66, looseness=3.5, "\scalebox{.77}{$\mathclap{
        G
      }$}"{description},shift right=1]
     \ar[d]
      \\
      \TopologicalSpace
      \ar[shift left=4pt, out=-180+66, in=-66, looseness=3.5, "\scalebox{.77}{$\mathclap{
        G
      }$}"{description},shift right=1]
    \end{tikzcd}
    \;\;\;
    \in
    \;
    \ProperEquivariantSerreFibrations{G}\;.
  \end{equation}
\end{lemma}
\begin{proof}
By Prop. \ref{PassageToFixedLociPreservesEquivariantLocalTriviality}
and
Lem. \ref{EquivariantDirectProductBundlesFromHatGxTwistedProductBundles},
we have for all $H \underset{\mathrm{cpt}}{\subset} G$
that the underlying $\Gamma^H$-principal bundle
of the fixed locus bundle $\TopologicalPrincipalBundle_i^H$
is locally trivial, for $i \,\in\, \{0,1\}$.
Let $\widehat{\TopologicalSpace^H} \twoheadrightarrow \TopologicalSpace^H$
be an open cover over which both trivialize and consider the following commuting
diagram:
\vspace{-2mm}
$$
  \begin{tikzcd}[column sep=7pt]
    (\widehat{\TopologicalSpace^H} \times \Gamma^H)
      \times_{\TopologicalSpace^H}
    (\widehat{\TopologicalSpace^H} \times \Gamma^H)
      \times
    \Gamma^H
    \ar[rr]
    \ar[d, shift left=4pt]
    \ar[d, shift right=4pt]
    &
    {}
    \ar[ddd, phantom, "\mbox{\tiny\rm(pb)}"{pos=.33-.165}]
    \ar[ddd, phantom, "\mbox{\tiny\rm(pb)}"{pos=.66-.165}]
    \ar[ddd, phantom, "\mbox{\tiny\rm(pb)}"{pos=.99-.165}]
    &
    (
    \TopologicalPrincipalBundle^H_1
      \times_{\TopologicalSpace^H}
    \TopologicalPrincipalBundle^H_2
    )
      \times
    \Gamma
    \ar[d, shift left=4pt]
    \ar[d, shift right=4pt]
    \\
    (\widehat{\TopologicalSpace^H} \times \Gamma^H)
      \times_{\TopologicalSpace^H}
    (\widehat{\TopologicalSpace^H} \times \Gamma^H)
    \ar[rr]
    \ar[d, ->>]
    &&
    \TopologicalPrincipalBundle^H_1
      \times_{\TopologicalSpace^H}
    \TopologicalPrincipalBundle^H_2
    \ar[d, ->>]
    \\
    \widehat{\TopologicalSpace^H} \times \Gamma^H
    \ar[rr]
    \ar[d]
    &&
    \Quotient{
      (
      \TopologicalPrincipalBundle^H_1
        \times_{\TopologicalSpace^H}
      \TopologicalPrincipalBundle^H_2
      )
    }
    { \Gamma^H }
    \ar[d]
    \\
    \widehat{\TopologicalSpace^H}
    \ar[rr, ->>]
    &{}&
    \TopologicalSpace^H
  \end{tikzcd}
$$

\vspace{-1mm}
\noindent Here, the right vertical column is given and we are pulling back to the cover.
The form of the top square follows by the pasting law (Prop. \ref{PastingLaw})
and
the left middle morphism is an effective epimorphism by
regularity of $\Actions{G}(\kTopologicalSpaces)$ (Prop. \ref{CompactlyGeneratedTopologicalGActionsFormARegularCategory}).
This identifies is codomain with
the quotient $\widehat{\TopologicalSpace^H} \times \Gamma^H$,
as shown. In conclusion, the bottom square
shows that each $H$-fixed locus of
\eqref{QuotientedFiberProductOfEquivariantPrincipalFiberBundlesEquivalentlyProperEquivariantFibration}
is a locally trivial fiber bundle, hence a Serre fibration by
Lem. \ref{LocallyTrivialBundlesAreSerreFibrations}.
\end{proof}

\noindent
{\bf Concordance of equivariant principal bundles.}
Below, in \cref{EquivariantInfinityBundles}, we find that
the classification of equivariant principal bundles
{\it up to concordance}
follows
on general abstract grounds,
after embedding them into cohesive $\infty$-topos theory,
from the {\it orbi-smooth Oka principle}
(Thm. \ref{SmoothOkaPrinciple}, Thm. \ref{OrbiSmoothOkaPrinciple}).
Here we show that
for topological equivariant principal fiber bundles,
their concordance classes actually coincide with their isomorphism classes
(Thm. \ref{ConcordanceClassesOfTopologicalPrincipalBundles} below).
Together with the orbi-smooth Oka principle,
this hence implies the full classification theorem
(Thm. \ref{BorelClassificationOfEquivariantBundlesForResolvableSingularitiesAndEquivariantStructure} below).

\begin{definition}[Concordance of equivariant principal bundles]
  \label{ConcordanceOfEquivariantPrincipalBundles}
  For $G \acts \, \TopologicalSpace \,\in\, \Actions{G}(\kTopologicalSpaces)$,

 \noindent {\bf (i)} we say that a pair
  of equivariant principal fiber bundles
  $G \acts \, \TopologicalPrincipalBundle_0,
   \,
   G \acts \,  \TopologicalPrincipalBundle_1
    \in\,
  \EquivariantPrincipalFiberBundles{G}{\Gamma}_{\TopologicalSpace}$
  (Def. \ref{TerminologyForPrincipalBundles})
  are {\it concordant} if there exists
  an equivariant principal fiber bundle
  \vspace{-2mm}
  $$
    G \acts \; \widehat{ \TopologicalPrincipalBundle }
    \;\in\;
    \EquivariantPrincipalBundles{G}{\Gamma}_{\TopologicalSpace \times [0,1]}
  $$

  \vspace{-2mm}
  \noindent
  on the cylinder $G \acts \, \TopologicalSpace \times [0,1]$
  with trivial $G$-action on the topological interval,
  whose restriction to the endpoints is isomorphic to these bundles:
    \vspace{-2mm}
  $$
    G \acts \; \TopologicalPrincipalBundle_0
    \;
    \simeq
    \;
    G \acts \; {\widehat{\TopologicalPrincipalBundle}}\vert_{\TopologicalSpace \times \{0\}}
        \;\;\;\;\;
    \mbox{and}
    \;\;\;\;\;
    G \acts \; \TopologicalPrincipalBundle_1
    \;
    \simeq
    \;
    G \acts \; {\widehat{\TopologicalPrincipalBundle}}|_{\TopologicalSpace \times \{1\}}
    \,.
  $$

 \vspace{-2mm}
 \noindent {\bf (ii)}
  We denote the sets of equivalences classes of
  equivariant principal bundles under isomorphisms and under
  concordance, respectively, by:
  \vspace{-2mm}
  \begin{equation}
    \label{GeneralProjectionFromIsomorphismToConcordanceClassesOfEquivariantPrincipalBundles}
    \begin{tikzcd}
      \underset{
        \raisebox{-3pt}{
          \tiny
          \color{darkblue}
          \bf
          isomorphism classes
        }
      }{
      \IsomorphismClasses
      {
        \EquivariantPrincipalFiberBundles{G}{\Gamma}_{\TopologicalSpace}
      }
      }
      \ar[r, ->>]
      &
      \underset{
        \raisebox{-3pt}{
          \tiny
          \color{darkblue}
          \bf
          concordance classes
        }
      }{
      \ConcordanceClasses
      {
        \EquivariantPrincipalFiberBundles{G}{\Gamma}_{\TopologicalSpace}
      }
      }
      \,.
    \end{tikzcd}
  \end{equation}
\end{definition}
We are going to show (Thm. \ref{ConcordanceClassesOfTopologicalPrincipalBundles})
that
\eqref{ProjectionFromIsomorphismToConcordanceClassesOfEquivariantPrincipalBundles}
is
in fact a bijection. For this purpose, we need the following
Lemma \ref{IsomorphismsOfPrincipalBundlesAsSectionsOfTheirQuotientedFiberProduct},
whose statement is not surprising, but whose proof requires some care:

\begin{lemma}[Isomorphisms of equivariant principal bundles as sections of their quotiented fiber product]
  \label{IsomorphismsOfPrincipalBundlesAsSectionsOfTheirQuotientedFiberProduct}
  Given $G \acts \TopologicalSpace \,\in\, \Actions{G}(\kTopologicalSpaces)$,
  and
 $
    G\acts \,  \TopologicalPrincipalBundle_1,
    G\acts \, \TopologicalPrincipalBundle_2
      \,\in\,
    \EquivariantPrincipalBundles{G}{\Gamma}(\kTopologicalSpaces)_{\TopologicalSpace}
  $
  a pair of $G$-equivariant $\Gamma$-principal fiber bundles
  \footnote{
     In fact, the proof of
     Lem. \ref{IsomorphismsOfPrincipalBundlesAsSectionsOfTheirQuotientedFiberProduct}
     only uses that the
     underlying principal bundles are locally trivial,
     not that this also the case for all $H$-fixed loci.
  }
  (Def. \ref{TerminologyForPrincipalBundles}),
  then
  there is a natural bijection
  \vspace{-2mm}
  \begin{equation}
    \label{BijectionBetweenPrincipalBundleMorphismsAndSectionsOfTheirQuotientedFiberProduct}
    \left\{
    \begin{tikzcd}[row sep=small]
      \TopologicalPrincipalBundle_1
      \ar[rr, dashed, "f", "\sim"{swap}]
      \ar[out=180-66, in=66, looseness=3.5, "\scalebox{.77}{$\;\mathclap{
        \Gamma \rtimes G
      }\;$}"{description}, shift right=1]
      \ar[dr]
      &&
      \TopologicalPrincipalBundle_2
      \ar[out=180-66, in=66, looseness=3.5, "\scalebox{.77}{$\;\mathclap{
        \Gamma \rtimes G
      }\;$}"{description}, shift right=1]
      \ar[dl]
      \\
      &
      \TopologicalSpace
      \ar[out=-180+66, in=-66, looseness=3.5, "\scalebox{.77}{$\mathclap{
        G
      }$}"{description},shift left=1]
    \end{tikzcd}
    \right\}
    \;\;\;\;\;\;
    \leftrightarrow
    \;\;\;\;\;\;
    \Bigggg\{
    \begin{tikzcd}[column sep=large]
      &
      (
      \TopologicalPrincipalBundle_1
        \times_{\TopologicalSpace}
      \TopologicalPrincipalBundle_2
      ) / \Gamma
      \ar[out=180-66, in=66, looseness=3.5, "\scalebox{.77}{$\;\mathclap{
        G
      }\;$}"{description}, shift right=1]
      \ar[d]
      \\
      \TopologicalSpace
      \ar[r, -, shift left=1pt]
      \ar[r, -, shift right=1pt]
      \ar[ur, dashed, "\sigma"]
      \ar[out=-180+66, in=-66, looseness=3.5, "\scalebox{.77}{$\mathclap{
        G
      }$}"{description},shift left=1]
      &
      \TopologicalSpace
      \ar[out=-180+66, in=-66, looseness=3.5, "\scalebox{.77}{$\mathclap{
        G
      }$}"{description},shift left=1]
    \end{tikzcd}
    \Bigggg\}
  \end{equation}

  \vspace{-1mm}
  \noindent
  between
  their homomorphisms
  (necessarily isomorphisms, by Lem. \ref{HomomorphismsOfLocallyTrivialPrincipalBundlesArePullbackSquares})
  and
  the continuous sections of the diagonal $\Gamma$-quotient of their
  fiber product.
\end{lemma}

\begin{proof}
We spell out a proof in the non-equivariant case, i.e.
for the case $G = 1$,
in a way that uses nothing but the fact (Prop. \ref{CompactyGeneratedTopologicalSpacesFormARegularCategory})
that $\kTopologicalSpaces$ is a regular category.
The general case
then follows by observing that these arguments lift through the
forgetful functor
(since this creates limits and colimits, Lem. \ref{ForgetfulFunctorFromTopologicalGSpacesToGSpaces})
to the regular category $\Actions{G}(\kTopologicalSpaces)$
(Prop. \ref{CompactlyGeneratedTopologicalGActionsFormARegularCategory}).
In other words, the general proof is verbatim the following proof with
$G \acts \, (-)$ adjoined to all objects appearing.

Noticing that the underlying bundles of
$\TopologicalPrincipalBundle_1, \TopologicalPrincipalBundle_2$
are locally trivial (Lem. \ref{EquivariantDirectProductBundlesFromHatGxTwistedProductBundles}),
first observe that the quotient coprojection
\vspace{-1mm}
\begin{equation}
  \label{QuotientCoprojectionOfFiberProductOfPrincipalBundlesIsLocallyTrivial}
  \begin{tikzcd}[row sep=10pt]
    \TopologicalPrincipalBundle_1
      \times_{\TopologicalSpace}
    \TopologicalPrincipalBundle_2
    \ar[d, ->>]
    \\
    (
      \TopologicalPrincipalBundle_1
        \times_{\TopologicalSpace}
      \TopologicalPrincipalBundle_2
    )
      /
    \Gamma
  \end{tikzcd}
  \;\;\;\;\;\;\;\;
  \mbox{\small is locally trivial and an effective epimorphism.}
\end{equation}
Namely, restriction of the fiber product bundle
to any patch $\TopologicalPatch \xhookrightarrow{\;} \TopologicalSpace$
over which both factor bundles
$\TopologicalPrincipalBundle_1$
and
$\TopologicalPrincipalBundle_2$
trivialize gives the pullback diagram shown on the left here:
$$
  \begin{tikzcd}
    \TopologicalPatch
      \times
    \Gamma^{\times^3}
    \ar[d, shift left=4pt]
    \ar[d, shift right=4pt]
    \ar[rr]
    &
    {}
    \ar[ddd, phantom, "\mbox{\tiny \rm (pb)}"{pos=.15}]
    \ar[ddd, phantom, "\mbox{\tiny \rm (pb)}"{pos=.66}]
    &
    (
    \TopologicalPrincipalBundle_1
      \times_{\TopologicalSpace}
    \TopologicalPrincipalBundle_2
    )
    \times \Gamma
    \ar[d, shift left=4pt]
    \ar[d, shift right=4pt]
    \\
    \TopologicalPatch
      \times
    \Gamma^{\times^2}
    \ar[rr]
    \ar[dd]
    &{}&
    \TopologicalPrincipalBundle_1
      \times_{\TopologicalSpace}
    \TopologicalPrincipalBundle_2
    \ar[dd]
    \\
    \phantom{
    \TopologicalPatch
      \times
    \Gamma
    }
    &{}&
    \phantom{
    (
      \TopologicalPrincipalBundle_1
        \times_{\TopologicalSpace}
      \TopologicalPrincipalBundle_2
    )
      /
    \Gamma
    }
    \ar[d]
    \\
    \TopologicalPatch
    \ar[rr, hook]
    &{}&
    \TopologicalSpace
  \end{tikzcd}
  \;\;\;\;\;\;\;\;\;\;\;
  \simeq
  \;\;\;\;\;\;\;\;\;\;\;
  \begin{tikzcd}
    \TopologicalPatch
      \times
    \Gamma^{\times^3}
    \ar[d, shift left=4pt]
    \ar[d, shift right=4pt]
    \ar[rr]
    &
    {}
    \ar[ddd, phantom, "\mbox{\tiny \rm (pb)}"{pos=.15}]
    \ar[ddd, phantom, "\mbox{\tiny \rm (pb)}"{pos=.5}]
    &
    (
    \TopologicalPrincipalBundle_1
      \times_{\TopologicalSpace}
    \TopologicalPrincipalBundle_2
    )
    \times \Gamma
    \ar[d, shift left=4pt]
    \ar[d, shift right=4pt]
    \\
    \TopologicalPatch
      \times
    \Gamma^{\times^2}
    \ar[rr]
    \ar[d, ->>, "\mathrm{coeq}"]
    &{}&
    \TopologicalPrincipalBundle_1
      \times_{\TopologicalSpace}
    \TopologicalPrincipalBundle_2
    \ar[d, ->>, "\mathrm{coeq}"]
    \\
    \TopologicalPatch
      \times
    \Gamma
    \ar[d]
    \ar[rr, dashed]
    &{}&
    (
      \TopologicalPrincipalBundle_1
        \times_{\TopologicalSpace}
      \TopologicalPrincipalBundle_2
    )
      /
    \Gamma
    \ar[d]
    \\
    \TopologicalPatch
    \ar[rr, hook]
    &{}&
    \TopologicalSpace
    \,.
  \end{tikzcd}
$$
Passage to coequalizers of the top morphism pairs
factors
the left diagram as shown on the right, which exhibits the resulting
middle square as a homomorphism of $\Gamma$-torsors over
$\TopologicalSpace$,
whose domain is a trivial principal bundle.
This implies that the middle square on the right is a
pullback,
by Lem. \ref{HomomorphismsOfLocallyTrivialPrincipalBundlesArePullbackSquares}.
Therefore
regularity of $\kTopologicalSpaces$
(Prop. \ref{CompactyGeneratedTopologicalSpacesFormARegularCategory})
implies with Lem. \ref{ReversePastingLawInRegularCategories}
that also the bottom square on the right is a pullback.
Since this holds for all
patches $\TopologicalPatch$ in any joint trivializing cover
of the two bundles, this proves
the local triviality in \eqref{QuotientCoprojectionOfFiberProductOfPrincipalBundlesIsLocallyTrivial}.
From this, the effective epimorphy
follows with Lem. \ref{LocallyTrivialActionsArePrincipal}
and Lem. \ref{EffectiveEquivariantPrincipalBundles}.

Now,  given a section $\sigma$ as on the right of
\eqref{BijectionBetweenPrincipalBundleMorphismsAndSectionsOfTheirQuotientedFiberProduct},
complete it to the following commuting diagram,
for $i \,\in\, \{1,2\}$:
\begin{equation}
  \label{TurningSectionOfQuotientedFiberProductBundleIntoBundleIsomorphisms}
  \begin{tikzcd}[column sep=70pt]
    \TopologicalPrincipalBundle
    \times
    \Gamma
    \ar[d, shift right=4pt, "\mathrm{pr}_1"{swap}]
    \ar[d, shift left=4pt, "\rho"]
    \ar[r]
    \ar[dr, phantom, "\mbox{\tiny\rm(pb)}"]
    &
    (
    \TopologicalPrincipalBundle_1
      \times_{\TopologicalSpace}
    \TopologicalPrincipalBundle_2
    )
    \times
    \Gamma
    \ar[d, shift right=4pt, "\mathrm{pr}_1"{swap}]
    \ar[d, shift left=4pt, "\rho"]
    \ar[r, " \mathrm{pr}_i \times \mathrm{id} "]
    &
    \TopologicalPrincipalBundle_i
    \times
    \Gamma
    \ar[d, shift right=4pt, "\mathrm{pr}_1"{swap}]
    \ar[d, shift left=4pt, "\rho"]
    \\
    \TopologicalPrincipalBundle
    \ar[r, "{ (\phi, f \circ \phi^{-1} )} "{description} ]
    \ar[d, ->>]
    \ar[dr, phantom, "\mbox{\tiny\rm(pb)}"]
    &
    \TopologicalPrincipalBundle_1
      \times_{\TopologicalSpace}
    \TopologicalPrincipalBundle_2
    \ar[r, "\mathrm{pr}_i"]
    \ar[d, ->>]
    &
    \TopologicalPrincipalBundle_i
    \ar[d, ->>]
    \\
    \TopologicalSpace
    \ar[r, dashed, "\sigma"{description}]
    \ar[
      rr,
      rounded corners,
      to path={
        -- ([yshift=-10pt]\tikztostart.south)
        -- node[below]{\scalebox{.7}{$\mathrm{id}$}}
           ([yshift=-10pt]\tikztotarget.south)
        -- ([yshift=-00pt]\tikztotarget.south)
      }
    ]
    &
    (
    \TopologicalPrincipalBundle_1
      \times_{\TopologicalSpace}
    \TopologicalPrincipalBundle_2
    ) / \Gamma
    \ar[r, "\mathrm{pr}_i/\Gamma"{above}]
    &
    \TopologicalSpace
    \mathrlap{\,.}
  \end{tikzcd}
\end{equation}
Here the squares on the left are pullbacks by definition,
the top right square commutes by nature of the diagonal action,
and the bottom left square is induced by passage to quotients.

But the right part of this diagram exhibits the projection $\mathrm{pr}_i$
as a homomorphism of $\Gamma$-torsors over $\TopologicalSpace$
out of a locally trivial principal bundle
\eqref{QuotientCoprojectionOfFiberProductOfPrincipalBundlesIsLocallyTrivial},
whereby Lem. \ref{HomomorphismsOfLocallyTrivialPrincipalBundlesArePullbackSquares}
implies that also the bottom right square is a pullback.
From this the pasting law (Prop. \ref{PastingLaw}) implies that
the total bottom rectangle is a pullback, hence
(by Ex. \ref{PullbackPreservesIsomorphisms}), that the
total horizontal middle morphisms are isomorphism
of $\Gamma$-torsors over $\TopologicalSpace$, for $i \,\in\, \{1,2\}$,
which we may suggestively denote by:
\vspace{-2mm}
$$
  \begin{tikzcd}
    \TopologicalPrincipalBundle
    \ar[r, "\phi", "\sim"{swap}]
    &
    \TopologicalPrincipalBundle_1
  \end{tikzcd}
  \,,
  \;\;\;\;
  \begin{tikzcd}[column sep=large]
    \TopologicalPrincipalBundle
    \ar[r, "f\circ \phi^{-1}", "\sim"{swap}]
    &
    \TopologicalPrincipalBundle_2
  \end{tikzcd}
  \,,
  \;\;\;\;
  \mbox{\small hence:}
  \;\;\;\;
  \begin{tikzcd}
    \TopologicalPrincipalBundle_1
    \ar[r, "f", "\sim"{swap}]
    &
    \TopologicalPrincipalBundle_2
    \,.
  \end{tikzcd}
$$

Conversely, given such $f$, form the
top part of the diagram
\eqref{TurningSectionOfQuotientedFiberProductBundleIntoBundleIsomorphisms}
with $\TopologicalPrincipalBundle \,\coloneqq\, \TopologicalPrincipalBundle_1$
and $\phi \,\coloneqq\, \mathrm{id}$. Then passage to coequalizers
as in the bottom part of
\eqref{TurningSectionOfQuotientedFiberProductBundleIntoBundleIsomorphisms}
yields a section $\sigma$. By effective epimorphy of the
bottom vertical morphisms, these two constructions
are inverse to each other, as claimed.
\end{proof}

\begin{theorem}[Concordant equivariant principal fiber bundles bundles are isomorphic]
  \label{ConcordanceClassesOfTopologicalPrincipalBundles}
  For

  -
  $G \acts \, \Gamma
      \,\in\,
    \Groups\left(\Actions{G}(\kTopologicalSpaces)\right)$,

  - $G \acts \, \TopologicalSpace \,\in\,
     \GCWComplexes \xhookrightarrow{\quad}
     \Actions{G}(\kTopologicalSpaces)$
    (e.g., a smooth $G$-manifold),

  \noindent
  concordance classes
  (Def. \ref{ConcordanceOfEquivariantPrincipalBundles})
  of
  topological $G$-equivariant $\Gamma$-principal fiber bundles
  over $\TopologicalSpace$
  (Def. \ref{TerminologyForPrincipalBundles})
  coincide with their isomorphism classes in that the
  quotient projection
  \eqref{GeneralProjectionFromIsomorphismToConcordanceClassesOfEquivariantPrincipalBundles}
  is a bijection:
  \begin{equation}
    \label{ProjectionFromIsomorphismToConcordanceClassesOfEquivariantPrincipalBundles}
    \begin{tikzcd}
      \IsomorphismClasses
      {
        \EquivariantPrincipalFiberBundles{G}{\Gamma}(\kTopologicalSpaces)_{\TopologicalSpace}
      }
      \ar[
        r,
        ->>,
        "{\sim}"
      ]
      &
      \ConcordanceClasses
      {
        \EquivariantPrincipalFiberBundles{G}{\Gamma}(\kTopologicalSpaces)_{\TopologicalSpace}
      }
      \,.
    \end{tikzcd}
  \end{equation}
\end{theorem}
The following proof adapts the idea of the proof of
\cite[Cor. 15]{RobertsStevenson12} to equivariant bundles.
\begin{proof}
Given a concordance, as shown on the left,
we need to produce an isomorphism, as shown on the right here:
\vspace{-3mm}
\begin{equation}
  \label{IsomorphismOfPrincipalBundlesFromConcordance}
  \begin{tikzcd}
    \mathrm{P}_0
      \ar[out=180-66, in=66, looseness=3.5, "\scalebox{.77}{$\mathclap{
        G
      }$}"{description},shift right=1]
    \ar[r]
    \ar[d, "p_0"]
    \ar[dr, phantom, "\mbox{\tiny\rm(pb)}"]
    &
    \mathrm{P}
      \ar[out=180-66, in=66, looseness=3.5, "\scalebox{.77}{$\mathclap{
        G
      }$}"{description},shift right=1]
    \ar[d, "p"]
    \ar[dr, phantom, "\mbox{\tiny\rm(pb)}"]
    &
    \mathrm{P}_1
      \ar[out=180-66, in=66, looseness=3.5, "\scalebox{.77}{$\mathclap{
        G
      }$}"{description},shift right=1]
    \ar[d, "p_1"]
    \ar[l]
    \\
    \TopologicalSpace \times \{0\}
      \ar[shift left=4pt, out=-180+66, in=-66, looseness=3.5, "\scalebox{.77}{$\mathclap{
        G
      }$}"{description},shift right=1]
    \ar[r, hook]
    &
    \TopologicalSpace \times [0,1]
      \ar[shift left=4pt, out=-180+66, in=-66, looseness=3.5, "\scalebox{.77}{$\mathclap{
        G
      }$}"{description},shift right=1]
    &
    \TopologicalSpace \times \{1\}
      \ar[shift left=4pt, out=-180+66, in=-66, looseness=3.5, "\scalebox{.77}{$\mathclap{
        G
      }$}"{description},shift right=1]
    \ar[l, hook']
  \end{tikzcd}
  \qquad
  \rightsquigarrow
  \qquad
  \begin{tikzcd}[row sep=small]
    \mathrm{P}_0
      \ar[out=180-66, in=66, looseness=3.5, "\scalebox{.77}{$\;\mathclap{
        \Gamma \rtimes G
      }\;$}"{description},shift right=1]
    \ar[dr, "p_0"{below, xshift=-2pt}]
    \ar[rr, "\sim"]
    &&
    \mathrm{P}_1
    \ar[dl, "p_1"{below, xshift=2pt}]
      \ar[out=180-66, in=66, looseness=3.5, "\scalebox{.77}{$\;\mathclap{
        \Gamma \rtimes G
      }\;$}"{description},shift right=1]
    \\
    &
    \TopologicalSpace
      \ar[shift left=4pt, out=-180+66, in=-66, looseness=3.5, "\scalebox{.77}{$\mathclap{
        G
      }$}"{description},shift right=1]
      \end{tikzcd}
\end{equation}
Consider the following solid commuting square in
$\Actions{G}(\kTopologicalSpaces)$,
showing a local section $\sigma_0$
which exhibits,
under Lem. \ref{IsomorphismsOfPrincipalBundlesAsSectionsOfTheirQuotientedFiberProduct},
that the restriction of $P_0 \times [0,1]$ to $\{0\} \subset [0,1]$
is isomorphic to $P_0$, by construction:
\vspace{-2mm}
\begin{equation}
  \label{LiftingProblemInShowingThatConcordanceOfTopPrinBundImpliesIsomorphism}
  \begin{tikzcd}[column sep=large]
    \TopologicalSpace \times \{0\}
      \ar[out=180-66, in=66, looseness=3.5, "\scalebox{.77}{$\mathclap{
        G
      }$}"{description},shift right=1]
    \ar[rr, "\sigma_0"]
    \ar[
      d,
      "{
        \def\arraystretch{.5}
        \begin{array}{r}
        \ProperEquivariantSerreCofibrations{G}
        \\
        \,\cap\,
        \ProperEquivariantWeakHomotopyEquivalences{G}
        \end{array}
        \,\ni
      }"{left}]
    &&
    \Quotient
    {
      \big(
        \mathrm{P}
          \underset
            {\TopologicalSpace \times [0,1]}
            {\times}
        (\mathrm{P}_0 \times [0,1])
      \big)
    }
    { \Gamma }
      \ar[out=180-66, in=66, looseness=3.5, "\scalebox{.77}{$\mathclap{
        G
      }$}"{description},shift right=1]
    \ar[d, "\in \, \ProperEquivariantSerreFibrations{G}"]
    \\
    \TopologicalSpace \times [0,1]
      \ar[shift left=4pt, out=-180+66, in=-66, looseness=3.5, "\scalebox{.77}{$\mathclap{
        G
      }$}"{description},shift right=1]
    \ar[rr, -, shift right=1pt]
    \ar[rr, -, shift left=1pt]
    \ar[urr, dashed, "\exists"]
    &&
    \TopologicalSpace \times [0,1]
      \ar[shift left=4pt, out=-180+66, in=-66, looseness=3.5, "\scalebox{.77}{$\mathclap{
        G
      }$}"{description},shift right=1]
  \end{tikzcd}
\end{equation}
Now observe that, with respect to
the proper equivariant model structure
(from Prop. \ref{ProperEquivariantModelCategoryOfGSpaces}):

\vspace{-2.5mm}
\begin{enumerate}[{\bf (i)}]
\setlength\itemsep{-3pt}

\item
the
left morphism is
an acyclic cofibration, since
it is the image of the generating acyclic cofibration
$G/G \times (D^0 \xhookrightarrow{} D^0 \times [0,1])$
\eqref{GeneratingAcyclicCofibrationsOfProperEquivariantModelCategory}
under the
the product functor $\TopologicalSpace \times (-)$,
which is left Quillen (by Prop. \ref{EquivariantMappingSpaceQuillenAdjunction})
since
$G$-CW-complexes are cofibrant (Ex. \ref{GCWComplexesAreCofibrantObjectsInProperEquivariantModelcategory});

\item
the
right vertical morphism
is a fibration, by
Lem. \ref{QuotientedFiberProductOfEquivariantPrincipalFiberBundlesIsProperEquivariantFibration}.
\end{enumerate}
\vspace{-2.5mm}

\noindent Therefore, the proper equivariant model structure implies that
a dashed lift in
\eqref{LiftingProblemInShowingThatConcordanceOfTopPrinBundImpliesIsomorphism}
exists. The resulting commutativity of the
bottom right triangle in
\eqref{LiftingProblemInShowingThatConcordanceOfTopPrinBundImpliesIsomorphism}
means that this lift is a section
which,
under Lem. \ref{IsomorphismsOfPrincipalBundlesAsSectionsOfTheirQuotientedFiberProduct},
exhibits an isomorphism
$G \acts \, \TopologicalPrincipalBundle
  \,\simeq\,
  G \acts \, \TopologicalPrincipalBundle_0 \times [0,1]$
of equivariant principal bundles over
$G \acts \, \TopologicalSpace \times [0,1]$.
The restriction of this isomorphism to $\{1\} \subset [0,1]$ is
of the required form \eqref{IsomorphismOfPrincipalBundlesFromConcordance}.
\end{proof}

\begin{remark}[Generalized conditions such that concordance classes coincide with isomorphism classes]
  The proof of Thm. \ref{ConcordanceClassesOfTopologicalPrincipalBundles}
  clearly works internal to any (regular) model category in which fiber bundles
  are fibrations and $X \times \{0\} \to X \times [0,1]$ exists and is
  an acyclic cofibration. In particular, there ought to be
  a suitable such model structure on diffeological spaces\footnote{
    This is claimed to be the case by Dmitri Pavlov, in private conversation.
  }
  which would imply that concordance classes coincide with isomorphism classes
  also for the case of diffeological principal bundles --
  in particular of smooth principal bundles with structure Lie groups.
  Moreover, if a proper equivariant version of such a model structure
  on diffeological spaces existed, it would imply the
  identification of concordance classes with isomorphism classes
  also for smooth equivariant principal bundles.
  When used further below in the proof of
  Thm. \ref{BorelClassificationOfEquivariantBundlesForResolvableSingularitiesAndEquivariantStructure},
  this would
  imply, under the assumptions used there,
  the classification if smooth equivariant principal bundles
  not just for D-topological but for more general diffeological structure.
\end{remark}

\medskip

\noindent
{\bf Proof of equivalence of notions of equivariant local triviality.}
We now turn to the proof of equivalence of the
three notions of equivariant local triviality (Thm. \ref{EquivalentNotionsOfEquivariantLocalTriviality}).
The following notion plays a pivotal role in all three perspectives:

\begin{notation}[Lifts of equivariance subgroups to semidirect product with structure group]
\label{LiftsOfEquivarianceSubgroupsToSemidirectProductWithStructureGroup}
Given a $G$-equivariant topological group $(\Gamma,\alpha)$
(Def. \ref{EquivariantTopologicalGroup}),
we will now
abbreviate the semidirect product group
\eqref{IdentifyingEquivariantGroupsWithSemidirectProductGroups}
as $\widehat G$. Moreover, given a closed subgroup $H \subset G$,
we will write $\widehat H \subset \widehat G$ for any {\it choice} of
lift of $H$ to a subgroup of $\widehat G$:
\vspace{-3mm}
\begin{equation}
  \label{LiftOfSubgroupsHToSemidirectProductGroup}
  \begin{tikzcd}[column sep=tiny]
    \widehat H
    \ar[
      d,
      "\simeq"{
        right
      },
      shift left=1pt,
      bend left=10
    ]
    &\subset&
    \widehat G
    \mathrlap{
      \;
      \coloneqq
      \Gamma \rtimes_\alpha G
    }
    \ar[
      d,
      "\mathrm{pr}_2"
    ]
    &{\phantom{AAAAAAA}}&
    \widehat G / \widehat H
    \ar[
      d,
      "\mathrm{pr}_2/\widehat H"{
        near start
      }
    ]
    \\
    H
    \ar[
      u,
      shift left=1pt,
      bend left=10,
      "
        \scalebox{.7}{$
          \widehat{(-)}
        $}
      "{
        left
      }
    ]
    &\subset&
    G
    &&
    G/H
  \end{tikzcd}
\end{equation}

  \vspace{-3mm}
\noindent
Equivalently,
such $\widehat H$ is
(by Lem. \ref{CrossedHomomorphismsAreEquivalentlyMaySubgroupsOfSemidirectProducts})
the graph of a crossed homomorphism $H \xrightarrow{\;} \Gamma$
(Def. \ref{CrossedHomomorphismsAndFirstNonAbelianGroupCohomology}).
\end{notation}
These choices of lifts $\widehat H$
are at the heart of the local characterization
of equivariant principal bundles:

\begin{lemma}[Equivariant principal $\widehat H$-twisted product bundles]
  \label{EquivariantPrincipalTwistedProductBundles}
%
  Consider $H \subset G$ a subgroup,
  with a lift
  $\widehat H \subset \widehat G \coloneqq \Gamma \rtimes_\alpha G$
  (Ntn. \ref{LiftsOfEquivarianceSubgroupsToSemidirectProductWithStructureGroup})
  and
  $\mathrm{S} \in H\mathrm{Act}(\mathrm{TopSp})$.
  Then sufficient conditions for the canonical projection of twisted product spaces
    \vspace{-3mm}
  \begin{equation}
  \hspace{-3cm}
    \begin{tikzcd}[column sep=small]
      \widehat G
        \times_{\widehat H}
      \mathrm{S}
      \ar[out=180-66, in=66, looseness=3.5, "\scalebox{.77}{$\phantom{}\mathclap{
        \widehat G
      }\phantom{}$}"{description}, shift right=1]
      \ar[
        dd,
        "\mathrm{pr}_2 \times_{\widehat H} \mathrm{id}_{\mathrm{S}}"{left}
      ]
      &
      \coloneqq
      &
      \{
        (\gamma,g), s
      \}_{
        \mathrlap{\!\!
        \big/
        \left(
          (
            (\gamma,g),\, s
          )
        \;\sim\;
          \left(
            (\gamma, g) \cdot \hat h, \, h^{-1} \cdot s
          \right)
        \right)
        }
      }
      \ar[
        dd,
        |->,
        "{
          \def\arraystretch{.8}
          \begin{array}{c}
            {[(\gamma,g),s]}
            \\
            \mapsdown
            \\
            {[g,s]}
          \end{array}
        }"{left}
      ]
      \\
      \\
      G \times_H \mathrm{S}
      \ar[out=-180+66, in=-66, looseness=3.5, "\scalebox{.77}{$\phantom{}\mathclap{
        G
      }\phantom{}$}"{description}, shift left=1]
      &\coloneqq&
      \{
        (g,s)
      \}_{
        \mathrlap{\!\!
        \big/
        \left(
          (g,s)
            \;\sim\;
          (g \cdot h,\, h^{-1}\cdot s)
        \right)
        }
      }
    \end{tikzcd}
  \end{equation}

  \vspace{-2mm}
\noindent   (equipped with their canonical
  left multiplication actions)
  to be a $G$-equivariant $(\Gamma,\alpha)$-principal bundle,
  according to Def. \ref{EquivariantPrincipalBundle} and via Cor. \ref{InternalDefinitionOfGPrincipalBundlesCoicidesWithtomDieckDefinition}, are:

   \vspace{-3mm}
  \begin{itemize}
\setlength\itemsep{-2pt}
      \item
    the underlying $\Gamma$-principal bundle is locally trivial;
    \item
    or $\Gamma$ is compact.
  \end{itemize}
\end{lemma}

\begin{proof}
  The equivariance conditions \eqref{ExternalEquivariancePropertyOfEquivariantPrincipalBundle}
  are immediate by construction.
  Principality \eqref{FiberwiseShearMapIsomorphism},
  follows from consideration of the following diagram:
    \vspace{-.2cm}
  \begin{equation}
    \label{RecognizingTwistedProductBundlesAsEquivariantPrincipal}
    \hspace{-6mm}
    \begin{tikzcd}[column sep=large]
      \Gamma
        \times
      \left(
        (
          \Gamma \rtimes_\alpha G
        )
        \times_{\widehat H}
      \mathrm{S}
      \right)
      \ar[
        rr,
        shift left=5pt,
        "{
          \left(
          \gamma',
          [
            (\gamma,\, g),\, s
          ]
          \right)
         \, \mapsto \,
          [
            ( \gamma' \cdot \gamma,\, \gamma,\, g),\, s
          ]
        }"
      ]
      \ar[
        rr,
        phantom,
        "\scalebox{.7}{$\sim$}"{description}
      ]
      \ar[
        drr,
        bend right=15,
        "\mathrm{pr}_2\;\;\;\;\;\;\;\;\;\;\;"{left, }
      ]
      &{\phantom{AAAAAAA}}&
      \left(
        (\Gamma \times \Gamma) \rtimes_\alpha G
      \right)
      \times_{\widehat H}
      \mathrm{S}
      \ar[
        ll,
        shift left=4pt,
        "{
          \left(
            \gamma_1 \cdot \gamma_2^{-1},\,
            [
              (\gamma_2,\, g)
              ,\,
              s
            ]
          \right)
         \, \mapsfrom \,
          [
            (\gamma_1,\, \gamma_2,\, g)
            ,\,
            s
          ]
        }"
      ]
      \ar[r]
      \ar[
        d,
        "p"
      ]
      \ar[
        dr,
        phantom,
        "\mbox{\tiny\rm(pb)}"
      ]
      &
      (
        \Gamma \rtimes_\alpha G
      )
      \times_{\widehat H}
      \mathrm{S}
      \ar[
        d,
        "\mathrm{pr}_2 \times_{\widehat H} \mathrm{id}_{\mathrm{S}}"{right}
      ]
      \\
      &&
      (
        \Gamma \rtimes_\alpha G
      )
      \times_{\widehat H}
      \mathrm{S}
      \ar[
        r,
        "\mathrm{pr}_2 \times_{\widehat H} \mathrm{id}_{\mathrm{S}}"{below}
      ]
      &
      G \times_H \mathrm{S}
      \mathrlap{\,.}
    \end{tikzcd}
  \end{equation}

  \vspace{-2mm}
\noindent  Here the square on the right,
  regarded in $\Gamma\mathrm{Act}(\mathrm{TopSp})$
  with respect to the left multplication action on the left of the two $\Gamma$-factors,
  is recognized as a pullback:
  \vspace{-.25cm}
  \begin{itemize}
\setlength\itemsep{-1pt}
  \item
  if the underlying $\Gamma$-principal bundle is locally trivial:
  by Lemma \ref{HomomorphismsOfLocallyTrivialPrincipalBundlesArePullbackSquares},
  observing that any
  local trivialization of
  $\mathrm{pr}_2 \times_{\widehat H} \mathrm{id}_{\mathrm{S}}$
  induces one of $p$

  \item
  if $\Gamma$ is compact: by
  Lemma \ref{RecognitionOfCartesianQuotientProjections},
  observing that all spaces are Hausdorff by
  assumption on $G$, $H$, and $\Gamma$ and by Lemma \ref{HausdorffQuotientSpaces}.

  \end{itemize}

    \vspace{-3mm}
\noindent
  With this, the shear map \eqref{FiberwiseShearMapIsomorphism}
  has an evident continuous inverse, as shown on the top left of \eqref{RecognizingTwistedProductBundlesAsEquivariantPrincipal};
  where the point is to observe that these formulas are indeed compatible with the
  given quotient space structures.
\end{proof}

Accordingly, in order to apply Lemma \ref{EquivariantPrincipalTwistedProductBundles}
we consider now classes of choices of $H$ and $\mathrm{S}$
such that the twisted product projections \eqref{RecognizingTwistedProductBundlesAsEquivariantPrincipal}
are locally trivial as $\Gamma$-principal bundles.

\subsection{tom Dieck local trivialization}
\label{tomDieckLocalTrivializations}

\begin{definition}[Equivariant open cover]
  \label{EquivariantOpenCover}
  For $G \acts \, \TopologicalSpace \,\in\, \GActionsOnTopologicalSpaces$,
  a {\it $G$-equivariant open cover} is an index set $I \,\in\, \Sets$,
  and an $I$-indexed set of open $G$-subspaces
  \vspace{-2mm}
  $$
    \raisebox{3pt}{\big\{}
      \begin{tikzcd}[column sep=large]
        \mathrm{U}_i
        \ar[out=180-66, in=66, looseness=3.5, "\scalebox{.77}{$\mathclap{
          G
        }$}"{description},shift right=1]
        \ar[
          r,
          hook,
          "\iota_i"{above},
          "\mathrm{open}"{below}
        ]
        &
        \TopologicalSpace
        \ar[out=180-66, in=66, looseness=3.5, "\scalebox{.77}{$\mathclap{
          G
        }$}"{description},shift right=1]
      \end{tikzcd}
        \raisebox{3pt}{\big\}}_{i \in I}
    \,,
    {\phantom{AAAAAA}}
\begin{tikzcd}[column sep=large]
  \mathllap{
    \underset{i \in I}{\sqcup}
  }
  \mathrm{U}_i
  \ar[
    out=180-66,
    in=66,
    looseness=3.5,
    "
    \scalebox{.77}{$
      \mathclap{
        G
      }
    $}
    "{
      description
    },
    shift right=1
  ]
  \ar[
    r,
    ->>,
    "\mathrm{open}"{
     below
    }
  ]
  &
  \TopologicalSpace
  \ar[
    out=180-66,
    in=66,
    looseness=3.5,
    "
    \scalebox{.77}{$
      \mathclap{
        G
      }
    $}
    "{
      description
    },
    shift right=1
  ]
\end{tikzcd}
  $$

    \vspace{-3mm}
\noindent
  such that, forgetting the $G$-action, the underlying open subsets cover $\TopologicalSpace$
  (cf to ``proper equivariant cover'' below in Def. \ref{ProperEquivariantOpenCover})
\end{definition}

\begin{definition}[tom Dieck's equivariant local trivializability]
  \label{tomDieckLocalTrivializability}
{\bf (i)}   A $G$-equivariant principal bundle $\mathrm{P} \xrightarrow{p} \TopologicalSpace$
  (Def. \ref{EquivariantPrincipalBundle}, Cor. \ref{InternalDefinitionOfGPrincipalBundlesCoicidesWithtomDieckDefinition})
  is {\it locally trivial} in the
  sense of \cite[\S 2.1]{tomDieck69}\footnote{
In fact, \cite{tomDieck69} does not require the square in \eqref{EquivariantLocalTrivialityCondition} to be a pullback, but instead
adds to the definition of equivariant bundles in the sense of
Cor. \ref{InternalDefinitionOfGPrincipalBundlesCoicidesWithtomDieckDefinition}
the condition that the underlying ordinary principal bundles are locally trivial.
These conditions are immediately equivalent to the ones we use, as shown by
Prop. \ref{tomDieckLocalTrivialityImpliesOrdinaryLocalTriviality} and
Prop. \ref{HomomorphismsOftomDieckLocallyTrivialEquivariantBundlesArePullbacks}.
We find the order of the conditions as used here a bit more systematic.}
  if
  there exists a $G$-equivariant open cover
  (Def. \ref{EquivariantOpenCover})
  \vspace{-3mm}
  $$
\begin{tikzcd}
  \mathllap{
    \underset{i \in I}{\sqcup}
  }
  \mathrm{U}_i
  \ar[out=180-66, in=66, looseness=3.5, "\scalebox{.77}{$\mathclap{
    G
  }$}"{description}, shift right=1]
  \ar[
    rr,
    ->>,
    "\mathrm{open}"{
     below
    }
  ]
  &&
  \TopologicalSpace
  \ar[
    out=180-66,
    in=66,
    looseness=3.5,
    "
    \scalebox{.77}{$
      \mathclap{
        G
      }
    $}
    "{
      description
    },
    shift right=1
  ]
\end{tikzcd}
  $$

  \vspace{-2mm}
  \noindent
  of its base space
  such that for each $i \in I$
the restriction of $\mathrm{P}$ to $\mathrm{U}_i$ is the pullback of
an equivariant principal bundle
$\mathrm{P}_i$
over a coset space $G/H_i$
for some closed subgroup $H_i \subset G$:
\vspace{-2mm}
\begin{equation}
  \label{EquivariantLocalTrivialityCondition}
  \begin{tikzcd}[row sep=1.3em]
  \mathrm{P}_{\vert \mathrm{U}_i}
  \ar[d]
  \ar[
    drr,
    phantom,
    "\mbox{\tiny\rm(pb)}"
  ]
  \ar[out=180-66, in=66, looseness=3.5, "\scalebox{.77}{$\phantom{\cdot}\mathclap{
        \Gamma \rtimes_\alpha G
  }\phantom{\cdot}$}"{description}, shift right=1]
  \ar[
    rr
  ]
  &&
  \mathrm{P}_i
  \ar[d]
  \ar[
    out=180-66,
    in=66,
    looseness=3.5,
    "
    \scalebox{.77}{$
      \phantom{\cdot}
      \mathclap{
        \Gamma \rtimes_\alpha G
      }
      \phantom{\cdot}
    $}
    "{
      description
    },
    shift right=1
  ]
  \ar[
    d,
    "\exists"{left},
    "\mathrm{principal}"{right}
  ]
  \\
  \mathrm{U}_i
  \ar[out=-180+66, in=-66, looseness=3.5, "\scalebox{.77}{$
        \mathclap{G}
  $}"{description}, shift left=1]
  \ar[
    rr,
    "\exists"{
      below
    }
  ]
  &&
  G/H_i
  \mathrlap{\,.}
  \ar[out=-180+66, in=-66, looseness=3.5, "\scalebox{.77}{$
        \mathclap{G}
  $}"{description}, shift left=1]
  \end{tikzcd}
\end{equation}
\end{definition}

\begin{lemma}[Coset space coprojections admitting local sections]
  \label{CosetSpaceCoprojectionsAdmittingLocalSections}
  Let $G$ be a topological group and $H$
  a topological subgroup. Then the following are
  sufficient conditions\footnote{
    Regarding necessity
    of these conditions, see counter-examples given in \cite[\S 3]{Karube58}.
  }
  for the coset space coprojection
  $G \xrightarrow{q} G/H$ to admit local sections:

\begin{enumerate}[{\bf (i)}]
  \vspace{-.2cm}
  \item
  $G$ is arbitrary, and
  \\
  $H$ is a compact Lie group;
  \vspace{-.2cm}
  \item
  $G$ is a locally compact separable metric space of finite dimension, and
  \\
  $H$ is a closed subgroup;
  \vspace{-.2cm}
  \item
  $G$ is a Lie group;
  \\
  $H$ is a closed subgroup.
\end{enumerate}
\vspace{-.3cm}

\end{lemma}

\begin{proof}
  The first statement is \cite[Thm. 4.1]{Gleason50},
  the second is \cite[Thm. 3]{Mostert53}, see also \cite[Thm. 2]{Karube58}.
  The third statement may be found in \cite[Thm. 4.3]{tomDieckBrocker85}.
\end{proof}

\begin{lemma}[Semidirect product coset bundles {\cite[Lem. 2.1]{tomDieck69}} are locally trivial]
  \label{SemidirectProductCosetBundles}
  Let
  $H \subset G$ a compact Lie group and
  $\widehat H$ any lift of $H$ to $\widehat G \coloneqq \Gamma \rtimes_\alpha G$
  (Ntn. \ref{LiftsOfEquivarianceSubgroupsToSemidirectProductWithStructureGroup}).
  Then the
  ordinary $\Gamma$-principal bundle,
  underlying
  (by Cor. \ref{UnderlyingPrincipalBundles})
  the equivariant principal $\widehat H$-twisted coset space bundle from
  Lemma \ref{EquivariantPrincipalTwistedProductBundles} for
  $\mathrm{S} = \ast$ (Ex. \ref{QuotientSpaces}):
    \vspace{-3mm}
  $$
   \begin{tikzcd}
     \widehat G
       \times_{\widehat H}
     \ast
     \ar[r,-, shift right=1pt]
     \ar[r,-, shift left=1pt]
     \ar[
       d,
       "\mathrm{pr}_2 \times_{\widehat H}\mathrm{id}_\ast"{left}
     ]
     &
     \widehat G
     /
     \widehat H
      \ar[out=180-66, in=66, looseness=3.5, "\scalebox{.77}{$\phantom{\cdot}\mathclap{
        \Gamma \rtimes_\alpha G
      }\phantom{\cdot}$}"{description}, shift right=1]
     \ar[
       d,
       "
         \mathrm{pr}_2/\widehat H
         \mathrlap{
           \mbox{
             \tiny
             \color{darkblue}
             \bf
             \def\arraystretch{1}
             \begin{tabular}{c}
               local model
               \\
               $G$-equivariant
               \\
               $\Gamma$-principal
               \\
               bundle
             \end{tabular}
           }
         }
       "
     ]
     \\
     G \times_H \ast
     \ar[r,-, shift right=1pt]
     \ar[r,-, shift left=1pt]
     &
     G/H
      \ar[out=-180+66, in=-66, looseness=3.5, "\scalebox{.77}{$
        \mathclap{G}
      $}"{description}, shift left=1]
   \end{tikzcd}
  $$

\vspace{-3mm}
\noindent
is locally trivial.
\end{lemma}
\noindent
(Compare the following proof with its re-derivation in geometric homotopy theory, below in Ex. \ref{GEquivariantPrincipalBundlesOverGCosetSpacesViaGeometricHomotopy}.)
\begin{proof}
  Under the given assumption on $H \subset G$,
  Lemma \ref{CosetSpaceCoprojectionsAdmittingLocalSections} says that
  there exists an open cover of $G/H$ such that over each of its charts
  $U \subset G/H$ we have a continuous section $\sigma$
\vspace{-2mm}
  $$
    \begin{tikzcd}[row sep=small]
      &
      G_{\vert U}
      \ar[d]
      \ar[r]
      \ar[
        dr,
        phantom,
        "\mbox{\tiny\rm(pb)}"
      ]
      &
      G
      \ar[d]
      \\
      U
      \ar[
        r,
        -,
        shift left=1pt
      ]
      \ar[
        r,
        -,
        shift right=1pt
      ]
      \ar[
        ur,
        dashed,
        "\sigma\;"{
          left
        }
      ]
      &
      U
      \ar[
        r,
        hook
      ]
      &
      G/H
    \end{tikzcd}
  $$

  \vspace{-2mm}
  \noindent
  With this, observe that
  \vspace{-4mm}
  \begin{equation}
    \label{LocalGammaTrivializationOfEquivariantCosetSpaceBundle}
    \begin{tikzcd}[row sep=1.5em, column sep=huge]
      \Gamma \times U
      \ar[
        rrr,
        "{
          (\gamma,u)
         \; \mapsto \;
          [\gamma,\sigma(u)]
        }"{
          above
        },
        shift left=3pt
      ]
      \ar[
        d,
        "\mathrm{pr}_1\;\;"{left}
      ]
      &&&
      \left(
        (
          \Gamma \rtimes_\alpha G
        )/\widehat H
      \right)\big\vert_{U}
      \ar[
        d
      ]
      \ar[
        lll,
        shift left=3pt,
        "{
          \big(
          \mathrm{pr}_1
          \big(
            (\gamma, g) \cdot \widehat{g^{-1} \sigma([g])} )
          \big),
          [g]
          \big)
   \;         \longmapsfrom \;
          [\gamma,g]
        }"{below}
      ]
      \\
      U
      \ar[rrr,-,shift left=1pt]
      \ar[rrr,-,shift right=1pt]
      &&&
      U
    \end{tikzcd}
  \end{equation}

  \vspace{-2mm}
  \noindent
  is a $\Gamma$-equivariant homeomorphism over $U$:
  The top map is clearly $\Gamma$-equivariant,
  the bottom map is clearly its inverse and both are continuous,
  by the continuity of all the maps from which they are composed, as shown.
\end{proof}

\begin{lemma}[Isomorphisms of local model equivariant principal bundles
  {\cite[\S 2.2]{tomDieck69}} ]
  \label{IsomorphismsOfLocalModelEquivariantPrincipalBundles}
  Let $H \subset G$ be any subgroup and
  $\widehat H_1, \widehat H_2 \subset \widehat G \coloneqq \Gamma \rtimes_\alpha G$
 be two lifts of $H$ as in \eqref{LiftOfSubgroupsHToSemidirectProductGroup}.
  There exists an isomorphism of equivariant principal bundles
  between the corresponding
  semidirect product coset space bundles from Lemma \ref{SemidirectProductCosetBundles},
  covering an isomorphism of their base spaces,
  precisely if $\widehat H_1$ and $\widehat H_2$ are conjugate in $\widehat G$
  by an element $\hat g \in \widehat G$:
  \vspace{-2mm}
  $$
    \begin{tikzcd}
      \widehat G / \widehat H_1
      \ar[d]
      \ar[
        rr,
        "{
          [\hat g]
          \,\mapsto\,
          [\hat g \hat g_0]
        }"{below},
        "\sim"{above}
      ]
     \ar[out=180-66, in=66, looseness=3.5, "\scalebox{.77}{$\mathclap{
        \widehat G
      }$}"{description},shift right=1]
       &&
      \widehat G / \widehat H_2
     \ar[out=180-66, in=66, looseness=3.5, "\scalebox{.77}{$\mathclap{
        \widehat G
      }$}"{description},shift right=1]
      \ar[d]
      \\
      G/H_1
      \ar[out=-180+66, in=-66, looseness=3.5, "\scalebox{.77}{$
        \mathclap{G}
      $}"{description}, shift left=1]
      \ar[
        rr,
        "{
          [g]
          \,\mapsto\,
          [g g_0]
        }"{below},
        "\sim"{above}
      ]
      &&
      G/H_2
      \ar[out=-180+66, in=-66, looseness=3.5, "\scalebox{.77}{$
        \mathclap{G}
      $}"{description}, shift left=1]
      \mathrlap{\,.}
    \end{tikzcd}
    {\phantom{AAAA}}
    \Leftrightarrow
    {\phantom{AAAA}}
    \begin{tikzcd}[column sep=0pt]
      \widehat H_1
      \ar[
        d,
        "\mathrm{pr}_2"
      ]
      &\subset&
      \hat g_0 \widehat H_2 \hat g_0^{-1}
      &\subset&
      \widehat G
      \\
      H_1
      &=&
      g_0 H_2 g_0^{-1}
      &\subset&
      G
    \end{tikzcd}
  $$
\end{lemma}
\begin{proof}
  Any $\widehat G$-equivariant map of $\widehat G$-coset spaces is
  determined, via equivariance, by its image $[\hat g_0]$ of the element $[e_{\widehat G}]$;
  and the condition that this assignment is consistent,
  in that it descends to the quotient space
  \vspace{-3mm}
  $$
    \begin{tikzcd}[row sep=small]
      \widehat G
      \ar[out=180-66, in=66, looseness=4.2, "\scalebox{.77}{$\mathclap{
        \widehat G
      }$}"{description},shift right=1]
      \ar[
        rr,
        "{
          \hat g \mapsto [\hat g \hat g_0]
        }"
      ]
      \ar[d]
      &&
      \widehat G / \widehat H_2
      \mathrlap{\,,}
      \ar[out=180-66, in=66, looseness=4.2, "\scalebox{.77}{$\mathclap{
        \widehat G
      }$}"{description},shift right=1]
      \\
      \widehat G / \widehat H_1
      \ar[
        urr,
        dashed
      ]
    \end{tikzcd}
  $$

  \vspace{-3mm}
  \noindent  is equivalent to the condition $\widehat H_1 \hat g_0 \subset \hat g_0 \widehat H_2$.
  If this inclusion is an isomorphism of subgroups, as assumed and as shown on the left here:
  $$
    \widehat H_2
      \;=\;
    \hat g_0^{-1}
      \cdot
    \widehat H_1
      \cdot
    \hat g_0
    \;\;\;\;\;\;\;\;
    \overset{\mathrm{pr}_2}{\Longrightarrow}
    \;\;\;\;\;\;\;\;
    H_2
      \;=\;
    g_0^{-1}
      \cdot
    H_1
      \cdot
    g_0\;,
  $$
  then the analogous argument
  (with $\widehat g$ replaced by $\widehat g^{-1}$) shows that
  we have an inverse map. Under the equivariant projection
  map, $\mathrm{pr}_2$, the analogous statement holds for the image subgroups in $G$.
\end{proof}

\begin{proposition}[Characterization of equivariant principal bundles over orbits {\cite[\S 2.1]{tomDieck69}}]
  \label{CharacterizationOfEquivariantBundlesOverCosetSpaces}
  Given a subgroup $H \subset G$, every $G$-
  equivariant $\Gamma$-principal bundle
  $\mathrm{P} \xrightarrow{\;} G/H$ (Def. \ref{EquivariantPrincipalBundle},
  Cor. \ref{InternalDefinitionOfGPrincipalBundlesCoicidesWithtomDieckDefinition})
  whose base space is the coset space $G/H$
  is either empty (Rem. \ref{PseudoTorsorCondition})
  or isomorphic to a semidirect product coset bundle from Lemma \ref{SemidirectProductCosetBundles},
  for some lift $\hat H$ \eqref{LiftOfSubgroupsHToSemidirectProductGroup}:
  \vspace{-2mm}
  \begin{equation}
    \label{EquivariantPrincipalBundleOverCosetAsSemidirectProductCosetProjection}
    \begin{tikzcd}[row sep=small]
      \mathllap{
        \exists
        \;\;\;\;\;
      }
      (
        \Gamma \rtimes_\alpha G
      )/\widehat H
      \ar[out=180-66, in=66, looseness=3.5, "\scalebox{.77}{$\phantom{\cdot}\mathclap{
        \Gamma \rtimes_\alpha G
      }\phantom{\cdot}$}"{description}, shift right=1]
      \ar[
        rr,
        "
          \sim
        "
      ]
      \ar[
        d
      ]
      &&
      P
      \ar[
        d,
        "
          \;\; p
        "{
          right
        }
      ]
      \ar[out=180-66, in=66, looseness=3.5, "\scalebox{.77}{$\phantom{\cdot}\mathclap{
        \Gamma \rtimes_\alpha G
      }\phantom{\cdot}$}"{description}, shift right=1]
      \\
      G/H
      \ar[rr,-,shift left=1pt]
      \ar[rr,-,shift right=1pt]
      &&
      G/H
    \end{tikzcd}
  \end{equation}
\end{proposition}
\begin{proof}
  In the case when the bundle is not empty,
  there exists an element in some fiber; and hence, by $G$-equivariance and
  transitivity of the $G$-action on $G/H$, we may find a point
  $\widehat{[e_G]} \in P_{[H]}$ in the fiber over $[e_G] \in G/H$.
  Let then
  \vspace{-2mm}
  \begin{equation}
    \label{StabilizerOfPointInFiber}
    \widehat H
    \;\coloneqq\;
    \mathrm{Stab}_{\Gamma \rtimes_\alpha G}
    \left(
      \widehat{[e_G]}
    \right)
    \;\subset\;
    \widehat G
    \;\coloneqq\;
    \Gamma \rtimes_\alpha G
  \end{equation}

      \vspace{-2mm}
\noindent
  be its isotropy group under the action of the full semidirect product group.
  We observe that this $\widehat H$ is a lift of $H$ as required in \eqref{LiftOfSubgroupsHToSemidirectProductGroup}:

   \vspace{-3mm}
  \begin{enumerate}[{\bf (i)}]
 \setlength\itemsep{-3pt}
  \item
    $\widehat H \subset \widehat G$ is a closed subgroup, since it is the preimage
    of the closed (by Hausdorffness, Assump. \ref{ProperEquivariantTopology})
    singleton subset $\{ \widehat{[e_G]}\} \subset \mathrm{P}$
    under the continuous function
    \begin{tikzcd}
      \widehat G
      \ar[
        rr,
        "
          \hat g \,\mapsto\, \hat g \cdot \widehat{[e_G]}
        "
      ]
      &&
      \mathrm{P}
      \mathrlap{\,.}
    \end{tikzcd}

  \item The restriction of $\widehat G \xrightarrow{\mathrm{pr}_2} G$
  to $\widehat H$ factors through $H \subset G$, since
   for all $\hat h \in H$ we have
      \vspace{-1mm}
   $$
     \begin{aligned}
       [e_G]
         & =
       p
       \big(
         \hat h \cdot \widehat{[e_G]}
       \big)
       \\
       & =
       \mathrm{pr}_2
       (
         \hat h
       )
       \cdot
       p
       \big(
         \widehat{[e_G]}
       \big)
       \\
       & =
       \mathrm{pr}_2
       (
         \hat h
       )
       \cdot [e_G]
       \,.
     \end{aligned}
   $$

\vspace{-2mm}
   \item
   This map $\mathrm{pr}_2 : \widehat H \xrightarrow{\;} H$
   is an isomorphism since (Example \ref{PullbackPreservesIsomorphisms})
   it is a pullback of
   (a restriction of) the $\Gamma$-principality isomorphism \eqref{FiberwiseShearMapIsomorphism}:
   \vspace{-4mm}
   $$
     \begin{tikzcd}[column sep=large]
       \widehat H
       \ar[
         rr,
         "{
           (\gamma,g) \,
           \mapsto \,
           g
         }"
       ]
       \ar[
         dd,
         "{
           \def\arraystretch{.7}
           \begin{array}{c}
             (\gamma,g)
             \\
             \mapsdown
             \\
             \gamma^{-1}
           \end{array}
         }"{
           left
         }
       ]
       \ar[
         ddrr,
         phantom,
         "\mbox{\tiny\rm(pb)}"
       ]
       &&
       \{
         \widehat{[e_G]}
       \}
         \times
       H
       \ar[
         dd,
         "{
           \def\arraystretch{.7}
           \begin{array}{c}
             g
             \\
             \mapsdown
             \\
             (e_\Gamma , g)\cdot \widehat {[e]}
           \end{array}
         }
         "
       ]
       \\
       \\
       \big\{
         \widehat{[e_G]}
       \big\}
         \times
       \Gamma
       \ar[
         rr,
         "{
           \gamma
             \;\mapsto\,
           (\gamma, e_G) \cdot \widehat{[e_G]}
         }"{
           above
         },
         "\simeq"{
           below
         }
       ]
       &&
       P_{[e]}
       \mathrlap{\,.}
     \end{tikzcd}
   $$

  \end{enumerate}

\vspace{-1cm}
\end{proof}

\begin{example}[Adjusting the classifying maps in a tom Dieck local trivialization]
\label{AdjustingTheClassifyingMapsInAtomDieckLocalTrivialization}
Given a tom Dieck local trivialization (Def. \ref{tomDieckLocalTrivializability})
over a $G$-equivariant patch $\mathrm{U}_i$,
consider any point $x \in \mathrm{U}_i$
and write $[g_0] \in G/H_i$ for its image under
the given classifying map \eqref{EquivariantLocalTrivialityCondition}.
Now choosing any identification of the local model bundle with a semidirect product
coset projection $\widehat G / \widehat H_i$ (Lemma \ref{CharacterizationOfEquivariantBundlesOverCosetSpaces})
and choosing any lift $[\hat g_0] \in \widehat G / \widehat {H_i}$
(Ntn. \ref{LiftOfIsotropyGroupsToSemidirectProductWithStructureGroup}),
we may then form the pasting composite of this
pullback square with the corresponding isomorphism
of local model bundles from Lemma \ref{IsomorphismsOfLocalModelEquivariantPrincipalBundles}:
 \vspace{-2mm}
\begin{equation}
  \label{PastingCompositeAdjustingTheClassifyingMapOfAtomDieckLocalTrivialization}
\hspace{-1cm}
  \begin{tikzcd}
    \mathrm{P}_{\vert \mathrm{U}_i}
    \ar[
      rr,
    ]
    \ar[
      d
    ]
    \ar[
      drr,
      phantom,
      "\mbox{\tiny\rm(pb)}"{description}
    ]
    &&
    \widehat G / \widehat{H_i}
    \ar[d]
    \ar[
      rr,
      "\sim"{below},
      "{
        [\hat g]
        \,\mapsto\,
        [\hat g \hat g^{-1}_0]
      }"
    ]
    \ar[
      drr,
      phantom,
      "\mbox{\tiny\rm(pb)}"
    ]
    &&
    \widehat G / \widehat{H'_i}
    \ar[d]
    \\
    \mathrm{U}_i
    \ar[
      rr,
      "{
        x \,\mapsto\, [g_0]
      }"{below}
    ]
    &&
    G/H_i
    \ar[
      rr,
      "\sim"{above},
      "{
        [g]
          \,\mapsto\,
        [g g^{-1}_0]
      }"{below}
    ]
    &&
    G/H'_i
  \end{tikzcd}
  {\phantom{AAAA}}
  \mbox{with}
  {\phantom{AAAA}}
  \begin{tikzcd}[column sep=0pt]
    \widehat {H'_i}
    &\coloneqq&
    \hat g_0 \widehat{H_i} \hat g^{-1}_0
    \ar[d]
    &\subset&
    \widehat G
    \\
    H'_i
    &\coloneqq&
    g_0 H_i g^{-1}_0
    &\subset&
    G
  \end{tikzcd}
\end{equation}

 \vspace{-1mm}
\noindent
Since the square on the right is a pullback (by Example \ref{PullbackPreservesIsomorphisms}),
the total rectangle is a pullback (by the pasting law, Prop. \ref{PastingLaw})
and hence exhibits another, equivalent, tom Dieck local trivialization,
whose classifying map takes $x$ to $[e_G]$.
\end{example}

\begin{proposition}[tom Dieck's local triviality implies underlying ordinary local triviality]
  \label{tomDieckLocalTrivialityImpliesOrdinaryLocalTriviality}
  An  equivariant principal bundle (Def. \ref{EquivariantPrincipalBundle},
  Cor. \ref{InternalDefinitionOfGPrincipalBundlesCoicidesWithtomDieckDefinition}),
  which is locally trivial in the sense of tom Dieck
  (Def. \ref{tomDieckLocalTrivializability}),
  has an underlying principal bundle
  (Cor. \ref{UnderlyingPrincipalBundles})
  which is locally trivial in the ordinary sense:
   \vspace{-1mm}
  $$
    \begin{tikzcd}
      \FormallyPrincipalFiberBundles{(G\acts \, \Gamma)}
      \left(
        \GActionsOnTopologicalSpaces
      \right)
      \ar[
        rr,
        "
          \mbox{
            \tiny
            \color{greenii}
            \bf
            forget $G$-action
          }
        "
      ]
      &&
      \FormallyPrincipalFiberBundles{\Gamma}
      (
        \kTopologicalSpaces
      )\;.
    \end{tikzcd}
  $$
\end{proposition}
\begin{proof}
  By Prop. \ref{CharacterizationOfEquivariantBundlesOverCosetSpaces}
  and Lemma \ref{SemidirectProductCosetBundles},
  the statement is true for the local model bundles over coset spaces.
  By Lemma \ref{ForgetfulFunctorFromTopologicalGSpacesToGSpaces}
  the pullback of this local model bundle
  is given by the pullback of the underlying topological spaces.
  With this,
  the condition \eqref{EquivariantLocalTrivialityCondition}
  implies the claim,
  since any pullback bundle of an ordinary locally trivial bundle is again locally trivial.
\end{proof}

\begin{proposition}[Homomorphisms of tom Dieck-locally trivial equivariant bundles are pullbacks]
  \label{HomomorphismsOftomDieckLocallyTrivialEquivariantBundlesArePullbacks}
  Let $G$ be a Lie group.
  Then every homomorphism $f$ \eqref{MorphismsOfEquivariantBundlesExternally}
  between
  tom Dieck-locally trivial
  $G$-equivariant $\Gamma$-principal bundles (Def. \ref{tomDieckLocalTrivializability})
  is a pullback square:
   \vspace{-2mm}
  $$
\begin{tikzcd}[row sep=small]
  \mathrm{P}_1
  \ar[out=180-66, in=66, looseness=3.5, "\scalebox{.77}{$\phantom{\cdot}\mathclap{
    \Gamma \rtimes_\alpha G
  }\phantom{\cdot}$}"{description}, shift right=1]
  \ar[
    rr,
    "f"
  ]
  \ar[d]
  \ar[
    drr,
    phantom,
    ""{description}
  ]
  &&
  \mathrm{P}_2
  \ar[out=180-66, in=66, looseness=3.5, "\scalebox{.77}{$\phantom{\cdot}\mathclap{
    \Gamma \rtimes_\alpha G
  }\phantom{\cdot}$}"{description}, shift right=1]
  \ar[d]
  \\
  \TopologicalSpace_1
  \ar[out=-180+66, in=-66, looseness=3.5, "\scalebox{.77}{$\mathclap{
    G
  }$}"{description}, shift left=1]
  \ar[
    rr,
    "f/\Gamma"{
      below
    }
  ]
  &&
  \TopologicalSpace_2
  \ar[out=-180+66, in=-66, looseness=3.5, "\scalebox{.77}{$
    \mathclap{G}
  $}"{description}, shift left=1]
\end{tikzcd}
  {\phantom{AAAA}}
  \Rightarrow
  {\phantom{AAAA}}
\begin{tikzcd}[row sep=small]
  \mathrm{P}_1
  \ar[out=180-66, in=66, looseness=3.5, "\scalebox{.77}{$\mathclap{
    G
  }$}"{description}, shift right=1]
  \ar[
    rr,
    "f"
  ]
  \ar[d]
  \ar[
    drr,
    phantom,
    "\mbox{\tiny\rm(pb)}"{description}
  ]
  &&
  \mathrm{P}_2
  \ar[out=180-66, in=66, looseness=3.5, "\scalebox{.77}{$\mathclap{
    G
  }$}"{description}, shift right=1]
  \ar[d]
  \\
  \TopologicalSpace_1
  \ar[out=-180+66, in=-66, looseness=3.5, "\scalebox{.77}{$\mathclap{
    G
  }$}"{description}, shift left=1]
  \ar[
    rr,
    "f/\Gamma"{
      below
    }
  ]
  &&
  \TopologicalSpace_2
  \ar[out=-180+66, in=-66, looseness=3.5, "\scalebox{.77}{$
    \mathclap{G}
  $}"{description}, shift left=1]
\end{tikzcd}
$$
\end{proposition}
\begin{proof}
  Prop. \ref{tomDieckLocalTrivialityImpliesOrdinaryLocalTriviality}
  says that, under the given assumptions,
  the ordinary principal bundles underlying the given equivariant principal
  bundles are locally trivial. With this, the claim
  follows from the fact that morphisms of ordinary locally trivial
  bundles are pullbacks (Lemma \ref{HomomorphismsOfLocallyTrivialPrincipalBundlesArePullbackSquares})
  and that underlying pullbacks reflect equivariant pullbacks
  (Lemma \ref{ForgetfulFunctorFromTopologicalGSpacesToGSpaces}).
\end{proof}

\subsection{Lashof-May local trivialization}
\label{LashofMyLocalTrivialization}

\begin{definition}[Lashof-May's equivariant local trivialization]
  \label{LashofMayEquivariantLocalTrivializability}
  An equivariant principal bundle $\mathrm{P} \xrightarrow{p} \TopologicalSpace$
  (Def. \ref{EquivariantPrincipalBundle},
  Cor. \ref{InternalDefinitionOfGPrincipalBundlesCoicidesWithtomDieckDefinition})
  is {\it locally trivial} in the sense of
  \cite[p. 258]{Lashof82}\cite[p. 267]{LashofMay86}
  if
  there exists a $G$-equivariant open cover
  (Def. \ref{EquivariantOpenCover})
  $
  {\phantom{\underset{i \in I}{\sqcup}}}
\begin{tikzcd}
  \mathllap{
    \underset{i \in I}{\sqcup}
  }
  \mathrm{U}_i
  \ar[out=180-66, in=66, looseness=3.5, "\scalebox{.77}{$\mathclap{
    G
  }$}"{description}, shift right=1]
  \ar[
    r,
    ->>,
    "\mathrm{open}"{
     below
    }
  ]
  &
  \TopologicalSpace
  \ar[out=180-66, in=66, looseness=3.5, "\scalebox{.77}{$\mathclap{
    G
  }$}"{description}, shift right=1]
\end{tikzcd}
  $
  by orbits of $H_i$-slices (Def. \ref{SliceOfTopologicalGSpace})
     \vspace{-7mm}
  \begin{equation}
    \label{OpenCoverByOrbitsOfSlices}
    \begin{tikzcd}
    \mathrm{U}_i
    \;=\;
    G \cdot \mathrm{S}_i
    &
    G \times_{H_i} \mathrm{S}_i
    \mathrlap{\,,}
    \ar[
      l,
      "\sim"{above},
      "\scalebox{0.7}{$(-)\cdot (-)$}"{below}
    ]
    \end{tikzcd}
    \qquad
    \begin{tikzcd}
      \mathrm{S}_i
      \ar[out=180-66, in=66, looseness=3.5, "\scalebox{.77}{$\phantom{}\mathclap{
        H_i
      }\phantom{}$}"{description}, shift right=1]
      \ar[
        r,
        hook
      ]
      &
      \TopologicalSpace
      \ar[out=180-66, in=66, looseness=3.5, "\scalebox{.77}{$\phantom{}\mathclap{
        H
      }\phantom{}$}"{description}, shift right=1]
    \end{tikzcd}
  \end{equation}

    \vspace{-2mm}
\noindent
  for closed subgroups
    \vspace{-2mm}
  \begin{equation}
    \label{ClosedSubgroupsInLashofMayLocalTrivialityCondition}
      \begin{tikzcd}
        H_i
        \ar[
          r,
          hook,
          "\,\mathrm{closed}\,"{above}
        ]
        &
        G
      \end{tikzcd}
    ,
  \end{equation}

    \vspace{-2mm}
\noindent
  such that, for each $i \in I$,

  \vspace{-2mm}
  \begin{itemize}
  \setlength\itemsep{-3pt}
\item   either the restriction $\mathrm{P}_{\vert \mathrm{U}_i}$ is empty
  (Rem. \ref{PseudoTorsorCondition}),

  \item or there is a lift $\widehat{H_i}$ \eqref{LiftOfSubgroupsHToSemidirectProductGroup}
  of $H_i$
  \eqref{ClosedSubgroupsInLashofMayLocalTrivialityCondition}
  to $\widehat G \coloneqq \Gamma \rtimes_\alpha G$
  (Ntn. \ref{LiftsOfEquivarianceSubgroupsToSemidirectProductWithStructureGroup}),
  such that the
  the restriction of $\mathrm{P}$ to $\mathrm{U}_i$ is isomorphic,
  as a $G$-equivariant $\Gamma$-principal bundle,
  to the equivariant principal
  $\widehat{H_i}$-twisted product
  bundle of the slice,
  from Lemma \ref{EquivariantPrincipalTwistedProductBundles}:
  \end{itemize}

     \vspace{-6mm}
  \begin{equation}
    \label{LashofMayLocalModelBundle}
    \begin{tikzcd}[column sep=large]
      \mathrm{P}_{\vert \mathrm{U}_i}
      \ar[out=180-66, in=66, looseness=4.2, "\scalebox{.77}{$\mathclap{
        \widehat G
      }$}"{description},shift right=1]
      \ar[
        r,
        "\sim"{above}
      ]
      \ar[
        d,
        "
         p_{\vert \mathrm{U}_i}
        "
      ]
      &
      \widehat G
        \times_{\widehat {H_i}}
      \mathrm{S}_i
      \ar[out=180-66, in=66, looseness=3.5, "\scalebox{.77}{$\mathclap{
        \widehat G
      }$}"{description},shift right=1]
      \ar[
        d,
        "\mathrm{pr}_2 \times_{\widehat H} \mathrm{id}_{\mathrm{S}_i}"{right}
      ]
      \\
      \mathrm{U}_i
  \ar[
    out=-180+66,
    in=-66,
    looseness=3.5,
    "
      \scalebox{.77}{$
        \mathclap{G}
      $}
    "{
      description
    },
    shift left=1
  ]
      \ar[
        r,
        "\sim"{above}
      ]
      &
      G \times_H \mathrm{S}_i
  \ar[
    out=-180+66,
    in=-66,
    looseness=3.5,
    "
      \scalebox{.77}{$
        \mathclap{G}
      $}
    "{
      description
    },
    shift left=1
  ]
    \end{tikzcd}
  \end{equation}
\end{definition}

\begin{lemma}[Local triviality of $\widehat H$-twisted product bundles]
  \label{LashofBundlesOverSliceOrbitsAreLocallyTrivialAsOrdinaryPrincipalBundles}
  Let $H \subset G$ be a compact subgroup (Assump. \ref{ProperEquivariantTopology})
  and  $H \acts \,  S \in H \mathrm{Act}( \mathrm{TopSp})$.

\vspace{-3mm}
 \begin{enumerate}[{\bf (i)}]
    \setlength\itemsep{-1pt}

  \item There exists a tubular neighborhood $D \times \mathrm{S}$ of
  $\mathrm{S}$ in $G \times_H S$ ($D$ an open ball), and hence
  an open cover of $G \times_H \mathrm{S}$ by its $G$-translates.

\item
  For every
  $\widehat H \subset \widehat G \coloneqq \Gamma \times_\alpha G$
  (Ntn. \ref{LiftsOfEquivarianceSubgroupsToSemidirectProductWithStructureGroup}),
  the induced $\widehat H$-twisted product bundle (Lemma \ref{EquivariantPrincipalTwistedProductBundles})
  has underlying $\Gamma$-principal bundle (Cor. \ref{UnderlyingPrincipalBundles})
  which is locally trivial over this open cover:
 \vspace{-3mm}
  $$
   \begin{tikzcd}[column sep=large]
      \Gamma
        \times
      D
        \times
      \mathrm{S}
      \ar[out=180-66, in=66, looseness=4.0, "\scalebox{.77}{$\mathclap{
        \Gamma
      }$}"{description},shift right=1]
     \ar[
       r
     ]
     \ar[
       d,
       "\mathrm{pr}_2"{left}
     ]
     \ar[
       dr,
       phantom,
       "\mbox{\tiny\rm(pb)}"{description}
     ]
     &
     \widehat G \times_{\widehat H } \mathrm{S}
     \ar[out=180-66, in=66, looseness=3.5, "\scalebox{.77}{$\mathclap{
       \Gamma
     }$}"{description},shift right=1]
     \ar[
       d
     ]
     \\
     D
       \times
     \mathrm{S}
     \ar[
       r,
       hook,
       "\mathrm{open}"{below}
     ]
     &
     G \times_H \mathrm{S}
   \end{tikzcd}
  $$
  \end{enumerate}
\end{lemma}
\begin{proof}
  The core of this argument is indicated inside the proof of \cite[Lem. 1.1]{Lashof82}
  for the special case when $\alpha$ \eqref{IdentifyingEquivariantGroupsWithSemidirectProductGroups}
  is trivial. We spell it out and generalize it.

  Since the connected component of $G$ is a compact Lie group (Assump. \ref{ProperEquivariantTopology}),
  we may find a bi-invariant Riemannian metric on $G$
  (by \cite[Cor. 1.4]{Milnor76}). With respect to this metric, let
    \vspace{-2mm}
  \begin{equation}
    \label{OpenBallAtNeutralElementNormalToClosedSubgroupH}
    D_e^\epsilon
    \;\coloneqq\;
    \begin{tikzcd}
      D^\epsilon N_e H
    \end{tikzcd}
  \end{equation}

    \vspace{-2mm}
\noindent
  be an open $\epsilon$-ball for any $\epsilon \in \mathbb{R}_+$
  in the normal bundle at the neutral element,
  and hence in a tubular neighborhood of $H$.
  By right invariance of the metric, the multiplication action
  yields a homomorphism
    \vspace{-2mm}
  \begin{equation}
    \label{IsomorphismOfDtimesHToDH}
    \begin{tikzcd}[row sep=-3pt, column sep=tiny]
      D \times H
      \ar[
        rr,
        "\sim"{below},
        "\kappa"{above}
      ]
      &&
      D \cdot H
      \\
    \scalebox{0.8}{$   (d,h) $}
      &\longmapsto&
   \scalebox{0.8}{$    d \cdot h $}
    \end{tikzcd}
  \end{equation}

    \vspace{-2mm}
\noindent
  which is surjective (since for each $h$ it restricts to an isometry
  $D \simeq D_h$) and continuously invertible
  (with inverse given by right multiplication with the inverse of the
  normal projection into $H$)
  so that by
  dimensional reasons we have that
  $$
    D \times \mathrm{S}
    \;\simeq\;
    (
      D \cdot H
   ) \times_H S
    \;\subset\;
    G \times_H S
  $$
  is an open neighborhood of $S \subset G \times_H S$.
  This proves the first statement.

  For the second statement, consider the following diagram,
  whose top morphisms are $\Gamma$-equivariant:
  \vspace{-3mm}
  \begin{equation}
   \label{PullbackRecognizingLashofModelBundlesAsLocallyTrivial}
   \hspace{-1cm}
   \begin{tikzcd}[column sep=large]
     \Gamma
       \times
     D
       \times
     \mathrm{S}
     \ar[
       r,
       "\sim"
     ]
     \ar[d]
     \ar[
       dr,
       phantom,
       "\mbox{\tiny\rm(pb)}"
     ]
     &
     \left(
       \Gamma
       \rtimes_\alpha
       (D \times H)
     \right)
       \times_{\widehat H}
     \mathrm{S}
     \ar[
       r,
       "\sim"
     ]
     \ar[d]
     \ar[
       dr,
       phantom,
       "\mbox{\tiny\rm(pb)}"
     ]
     &
     \left(
       \Gamma
         \rtimes_\alpha
         D \cdot H
     \right)
       \times_{\widehat H}
     \mathrm{S}
     \ar[r]
     \ar[d]
     \ar[
       dr,
       phantom,
       "\mbox{\tiny\rm(pb)}"
     ]
     &
     \widehat G \times_{\widehat H} \mathrm{S}
     \ar[d]
     \\
     D \times \mathrm{S}
     \ar[
       r,
       "\sim"
     ]
     &
     D \times H \times_H \mathrm{S}
     \ar[
       r,
       "\sim"{above},
       "\kappa \times_{H} \mathrm{id}_{\mathrm{S}}"{below}
     ]
     &
     (
       D \cdot H
     )
     \times_H \mathrm{S}
     \ar[
       r,
       hook,
       "\mathrm{open}"{below}
     ]
     &
     G \times_H \mathrm{S}
     \mathrlap{\,.}
   \end{tikzcd}
  \end{equation}

  \vspace{-2mm}
  \noindent
  Here the two squares on the right
  are pullback squares since they have parallel isomorphisms
  (by Example \ref{PullbackPreservesIsomorphisms}),
  induced from \eqref{IsomorphismOfDtimesHToDH}.
  This witnesses the second vertical morphism from the right
  as a trivialized $\Gamma$-principal bundle.
  With this, the right square is a pullback by
  Lemma \ref{HomomorphismsOfLocallyTrivialPrincipalBundlesArePullbackSquares}.
  In conclusion, the total rectangle in \eqref{PullbackRecognizingLashofModelBundlesAsLocallyTrivial}
  is a pullback by the pasting law (Prop. \ref{PastingLaw}),
  and hence exhibits the restriction of
  the twisted product bundle to the open subset $D \times \mathrm{S}$ of
  its base space as a Cartesian product projection.
\end{proof}

\begin{proposition}[Lashof-May's local trivializability implies tom Dieck's]
  \label{LashoMayLocalTrivialityImpliestomDieck}
  If a $G$-equivariant $(G\acts \, \Gamma)$-principal bundle
  (Def. \ref{EquivariantPrincipalBundle},
  Cor. \ref{InternalDefinitionOfGPrincipalBundlesCoicidesWithtomDieckDefinition},
  Assump. \ref{ProperEquivariantTopology})
  is locally trivial
  in the sense of Lashof-May (Def. \ref{LashofMayEquivariantLocalTrivializability})
  then it is also locally trivial in the sense of tom Dieck
  (Def. \ref{tomDieckLocalTrivializability}).
\end{proposition}
\begin{proof}
  Observe the evident commuting square
  of equivariant maps shown on the left in the following \eqref{ConstructingThePullbackOftomDieckLocalModelsToLashofMayLocalModels},
  where $\widehat G$ acts by left multiplication on itself
  (and on $G$, through $\mathrm{pr}_2$),
  while $\widehat {H_i}$ acts by right inverse multiplication on $\widehat G$
  (and on $G$ through $\mathrm{pr}_2$),
  by its given action on $\mathrm{S}_i$
  and diagonally on the Cartesian products:
   \vspace{-2mm}
  \begin{equation}
   \label{ConstructingThePullbackOftomDieckLocalModelsToLashofMayLocalModels}
    \begin{tikzcd}[column sep=large]
      \widehat G
      \times
      \mathrm{S}_i
      \ar[out=180-60, in=60, looseness=3.5, "\scalebox{.77}{$\phantom{-}\mathclap{
        \widehat G
        \times
        \widehat {H_i}
      }\phantom{-}$}"{description}, shift right=1]
      \ar[
        rr,
        "\mathrm{id}_{\widehat G} \times \mathrm{pr} "
      ]
      \ar[
        d,
        "\mathrm{pr}_2 \times \mathrm{id}"
      ]
      &&
      \widehat G
      \times
      \ast
      \ar[out=180-60, in=60, looseness=3.5, "\scalebox{.77}{$\phantom{-}\mathclap{
        \widehat G
        \times
        \widehat {H_i}
      }\phantom{-}$}"{description}, shift right=1]
      \ar[
        d,
        "\mathrm{pr}_2"
      ]
      \\
      G
        \times
      \mathrm{S}_i
      \ar[out=-180+56, in=-56, looseness=3.8, "\scalebox{.77}{$\phantom{\cdot}\mathclap{
          G
          \times
          \widehat {H_i}
      }\phantom{-}$}"{description}, shift left=1]
      \ar[
        rr,
        " \mathrm{id}_G \times \mathrm{pr} "
      ]
      &&
      G
      \ar[out=-180+56, in=-56, looseness=5, "\scalebox{.77}{$\phantom{\cdot}\mathclap{
        G
        \times
        \widehat {H_i}
      }\phantom{-}$}"{description}, shift left=1]
    \end{tikzcd}
    {\phantom{AAA}}
   \xmapsto{\;\;
      (-)/\widehat {H_i} \;\;}
    {\phantom{AAA}}
    \begin{tikzcd}
      \widehat G \times_{\widehat {H_i}} \mathrm{S}_i
      \ar[out=180-66, in=66, looseness=3.5, "\scalebox{.77}{$\phantom{\cdot}\mathclap{
        \widehat G
      }\phantom{\cdot}$}"{description}, shift right=1]
      \ar[rr]
      \ar[d]
      &&
      \widehat G / \widehat {H_i}
      \ar[out=180-66, in=66, looseness=3.5, "\scalebox{.77}{$\phantom{\cdot}\mathclap{
        \widehat G
      }\phantom{\cdot}$}"{description}, shift right=1]
      \ar[d]
      \\
      G \times_{H_i} \mathrm{S}_i
  \ar[
    out=-180+66,
    in=-66,
    looseness=3.5,
    "
      \scalebox{.77}{$
        \mathclap{
          G
        }
      $}
    "{
      description
    },
    shift left=1
  ]
      \ar[rr]
      &&
      G/H_i
  \ar[
    out=-180+66,
    in=-66,
    looseness=3.5,
    "
      \scalebox{.77}{$
        \mathclap{
          G
        }
      $}
    "{
      description
    },
    shift left=1
  ]
    \end{tikzcd}
  \end{equation}

   \vspace{-2mm}
\noindent  Passage to $\widehat {H_i}$-quotients (Example \ref{QuotientSpaces})
  yields the commuting square shown on the right, where
  the vertical map on the left is the
  Lashof-May local model bundle \eqref{LashofMayLocalModelBundle},
  while the vertical map on the right is a
  tom Dieck local model bundle \eqref{EquivariantLocalTrivialityCondition},
  by Lemma \ref{SemidirectProductCosetBundles}.

  Hence it is now sufficient to see that the square on the right
  of \eqref{ConstructingThePullbackOftomDieckLocalModelsToLashofMayLocalModels}
  is in fact a
  pullback square, as required in
  tom Dieck's local trivialization condition \eqref{EquivariantLocalTrivialityCondition}.
  By Lemma \ref{ForgetfulFunctorFromTopologicalGSpacesToGSpaces}, this
  is equivalent to showing that the underlying square of topological spaces
  is a pullback.
  But the square is a homomorphism of
  ordinary $\Gamma$-principal bundles
  (since $\Gamma \subset \widehat G$)
  which are both locally trivial,
  by
  Prop. \ref{tomDieckLocalTrivialityImpliesOrdinaryLocalTriviality}
  and
  by Lemma \ref{LashofBundlesOverSliceOrbitsAreLocallyTrivialAsOrdinaryPrincipalBundles}.
  Therefore its pullback property follows by
  Lemma \ref{HomomorphismsOfLocallyTrivialPrincipalBundlesArePullbackSquares}.
\end{proof}

\subsection{Bierstone local trivialization}
\label{BierstoneLocalTrivializations}

\begin{notation}[Lift of isotropy groups to semidirect product with structure group]
\label{LiftOfIsotropyGroupsToSemidirectProductWithStructureGroup}
Given an action $G \acts \, \TopologicalSpace$ and a point $x \in \TopologicalSpace$,
recall from \eqref{StabilizerSubgroupInEquivarianceGroup} the isotropy group $G_x \subset X$.
Thinking of the semidirect product group
$\widehat G \coloneqq \Gamma \rtimes_\alpha G$
(Ntn. \ref{LiftsOfEquivarianceSubgroupsToSemidirectProductWithStructureGroup})
as acting on $\TopologicalSpace$ through the projection $\mathrm{pr}_2$,
we write
\vspace{-2mm}
\begin{equation}
  \label{StabilizerSubgroupInSemidirectProductGroup}
  \widehat G_x
  \;\coloneqq\;
  \mathrm{Stab}_{{}_{\widehat G}}(x)
  \;\coloneqq\;
  \Gamma
    \rtimes_\alpha
  G_x
  \,.
\end{equation}

\vspace{-3mm}
\noindent
Hence in further following Ntn. \ref{LiftsOfEquivarianceSubgroupsToSemidirectProductWithStructureGroup},
a choice of lift \eqref{LiftOfSubgroupsHToSemidirectProductGroup} of
$G_x \subset G$
to $\widehat G_x$ is to be denoted $\widehat {G_x}$
(note the scope of the hat):
  \vspace{-2mm}
\begin{equation}
  \label{LiftOfIsotorpyGroupToSemidirectProductGroup}
  \begin{tikzcd}[column sep=tiny]
    \widehat {G_x}
    \ar[
      d,
      "\simeq"{
        right
      },
      shift left=1pt,
      bend left=10
    ]
    &\subset&
    \widehat G_x
    \mathrlap{
      \;
      \coloneqq
      \Gamma \rtimes_\alpha G_x
    }
    \ar[
      d,
      "\mathrm{pr}_2"
    ]
    \\
    G_x
    \ar[
      u,
      shift left=1pt,
      bend left=10,
      "
        \scalebox{.7}{$
          \widehat{(-)}
        $}
      "{
        left
      }
    ]
    &=&
    G_x
  \end{tikzcd}
\end{equation}
\end{notation}

\begin{lemma}[Beck-Chevalley for lift of isotropy groups]
 \label{BeckChevalleyConditionForLiftOfIsotropyGroup}
Given an isotropy subgroup $G_x \xhookrightarrow{\;} G$,
its two inclusions into $\widehat G$ (via Ntn. \ref{LiftOfIsotropyGroupsToSemidirectProductWithStructureGroup})
form a commuting square as shown on the left here:
  \vspace{-2mm}
\begin{equation}
  \label{BeckChevalleyTransformForIsotropySubgroup}
  \begin{tikzcd}
    &
    G_x
    \ar[
      dr,
      hook,
      "i_x"
    ]
    \ar[
      dl,
      hook,
      "s\;\;\;"{left, near start}
    ]
    \ar[
      dd,
      phantom,
      "\mbox{\tiny\rm(pb)}"{description}
    ]
    \\
    \widehat G_x
    \ar[
      dr,
      hook,
      "\hat i_x\;\;"{below}
    ]
    &&
    G
    \ar[
      dl,
      hook,
      "s"
    ]
    \\
    &
    \widehat G
  \end{tikzcd}
  {\phantom{AAAA}}
  \begin{tikzcd}[column sep=-10pt]
    &
    G_x \mathrm{Actions}
    \ar[
      dr,
      "\;G \times_{G_x}(-)"{right},
      ""{name=s, left}
    ]
    \\
    \widehat G_x
    \mathrm{Actions}
    \ar[
      ur,
      "s^\ast\;\;\;"{left},
    ]
    \ar[
      dr,
      "\widehat G \times_{\widehat G_x}(-)\;\;\;\;"{left, near end},
      ""{name=t, right}
    ]
    \ar[
      from=s,
      to=t,
      Rightarrow,
      "\raisebox{-11pt}{$\sim$}"{right, yshift=-1pt, sloped}
    ]
    &&
    G
    \mathrm{Actions}
    \\
    &
    \widehat G
    \mathrm{Actions}
    \ar[
      ur,
      "\;\;s^\ast"{right, near start}
    ]
  \end{tikzcd}
\end{equation}

  \vspace{-2mm}
\noindent
and the pull-push change-of-group operations
(Lemma \ref{InducedAndCoinducedActions})
via inductions/restriction (Example \ref{RestrictedActions})
from left to right through this square coincide to yield the
commuting diagram of functors shown on the right of \eqref{BeckChevalleyTransformForIsotropySubgroup}.
\end{lemma}
\begin{proof}
This is a direct consequence
of
the basic fact that
$
  \widehat G_x / \Gamma
  \;=\;
  (
    \Gamma \rtimes_\alpha G_x
  )
  /
  \Gamma
  \;\simeq\;
  G_x
  \,.
$
Explicitly,
we claim that the natural transformation filling the square
on the right of \eqref{BeckChevalleyTransformForIsotropySubgroup} is
\vspace{-3mm}
\begin{equation}
  \label{BeckChevalleyNaturalTransformation}
 \begin{tikzcd}[column sep=10pt, row sep=-3pt]
   \mathrm{P}
   \ar[out=180-60, in=60, looseness=3.2, "\scalebox{.77}{$\,\mathclap{
     \widehat G_x
   }\phantom{\cdot}$}"{description},shift right=2]
   \ar[rr]
   &&
   G \times_{G_x} \mathrm{P}
   \ar[out=180-60, in=60, looseness=3.2, "\scalebox{.77}{$\,\mathclap{
     G
   }$}"{description},shift right=1]
   \ar[
     rr,
     "\sim"{below}
   ]
   &&
   \widehat G \times_{\widehat G_x} \mathrm{P}
   \ar[out=180-60, in=60, looseness=2.2, "\scalebox{.77}{$\,\mathclap{
     G
   }$}"{description},shift right=1]
   \\
   &&
   \scalebox{0.7}{${[g,p]}$}
   &\longmapsto&
    \scalebox{0.7}{$ {[(\NeutralElement_\Gamma, g), p]} $}
   \mathrlap{\,.}
 \end{tikzcd}
\end{equation}

\vspace{-3mm}
\noindent This is manifestly equivariant under the left multiplication action by
$G \xhookrightarrow{s} \Gamma \rtimes_\alpha G$ and it is manifestly
natural under $\widehat G_x$-equivariant maps $\mathrm{P} \xrightarrow{\;} \mathrm{P}'$
on the right.
Moreover, \eqref{BeckChevalleyNaturalTransformation} is surjective because the general element
$[ (\gamma,g), p ]$
on the right is equal to
\vspace{-3mm}
$$
  [ (\gamma,g),\, p ]
  \,=\,
  \big[
    (\NeutralElement_\Gamma,\, g)
      \cdot
    ( \alpha(g^{-1})(\gamma),\, \NeutralElement_G )
    ,\,
    p
  \big]
  \,=\,
  \big[
    (\NeutralElement_\Gamma,\, g)
    ,\,
    ( \alpha(g^{-1})(\gamma),\, \NeutralElement_G )
      \cdot
    p'
  \big]
  \mathrlap{\,,}
$$

\vspace{-2mm}
\noindent
and it is injective because
\vspace{-2mm}
$$
  \begin{aligned}
  [
    (\NeutralElement_\Gamma, g_1), p_1
  ]
  \,=\,
  [
    (\NeutralElement_\Gamma, g_2), p_2
  ]
  \;\;\;
  \Leftrightarrow
  \;\;
  \underset{
    {g' \in G_x}
    \atop
    {\gamma \in \Gamma}
  }{\exists}
  \left\{\!\!\!
    \begin{array}{l}
      (\NeutralElement_\Gamma, g_1) \cdot (\gamma, g') = (\NeutralElement_\Gamma, g_1)\,,
      \\
      (\gamma, g')^{-1} \cdot p_1 = p_2
    \end{array}
  \right.
\;\;  \Rightarrow
  &
  \quad
  \underset{
    {g' \in G_x}
  }{\exists}
  \left\{
    \begin{array}{l}
      (\NeutralElement_\Gamma, g_1) \cdot (\NeutralElement_\Gamma, g')
        =
      (\NeutralElement_\Gamma, g_1)\;.
      \\
      (\NeutralElement_\Gamma, g')^{-1} \cdot p_1 = p_2
    \end{array}
  \right.
  \\
  \Leftrightarrow
  &
  \;\;
  [g_1 ,\, g_1] \,=\, [g_2 ,\, p_2]
  \end{aligned}
  $$

  \vspace{-2mm}
\noindent
In conclusion, we have a $G$-equivariant natural bijection, as claimed.
\end{proof}

\begin{definition}[Bierstone's equivariant local trivializability]
  \label{BierstoneEquivariantLocalTrivializability}
  A $G$-equivariant $(G \acts \, \Gamma)$-principal bundle $\mathrm{P} \xrightarrow{p} \TopologicalSpace$
  (Def. \ref{EquivariantPrincipalBundle},
  Cor. \ref{InternalDefinitionOfGPrincipalBundlesCoicidesWithtomDieckDefinition})
  is {\it locally trivial} in the
  sense of \cite[\S 4]{Bierstone73}
  generalized to non-trivial $\alpha$ \eqref{IdentifyingEquivariantGroupsWithSemidirectProductGroups},
  if for every point $x \in \TopologicalSpace$ there exists a
  $G_x$-equivariant \eqref{StabilizerSubgroupInEquivarianceGroup}
  open neighborhood
  $
    \begin{tikzcd}
      x
      \in
      \mathrm{U}_x
      \ar[out=180-66, in=66, looseness=4.2, "\scalebox{.77}{$\mathclap{
        G_x
      }$}"{description},shift right=1]
      \ar[
        r,
        hook,
        "\mathrm{open}"{below}
      ]
      &
      \TopologicalSpace
      \mathrlap{\,,}
      \ar[out=180-66, in=66, looseness=4.2, "\scalebox{.77}{$\mathclap{
        G_x
      }$}"{description},shift right=1]
    \end{tikzcd}
  $
  and a lift $\widehat {G_x} \subset \widehat G_x \coloneqq \Gamma \rtimes_\alpha G_x$
  \eqref{LiftOfIsotorpyGroupToSemidirectProductGroup}
  such that the restriction of $\mathrm{P}$ to $\mathrm{U}_x$ is

  \vspace{-3mm}
    \begin{itemize}
  \setlength\itemsep{-3pt}
\item
  either empty (Rem. \ref{PseudoTorsorCondition}),
\item
  or
  isomorphic,
  as a $G_x$-equivariant $(\Gamma,\alpha)$-principal bundle
  under restriction along $G_x \hookrightarrow G$ (Ex. \ref{RestrictedActions}),
  to the equivariant direct product bundle from
  Lemma \ref{EquivariantDirectProductBundlesFromHatGxTwistedProductBundles}:
  \end{itemize}
  \vspace{-5mm}
  \begin{equation}
    \label{BierstoneLocalIdentification}
    \begin{tikzcd}[column sep=large, row sep=15pt]
      \mathrm{P}_{\vert \mathrm{U}_x}
      \ar[out=180-66, in=66, looseness=4.2, "\scalebox{.77}{$\mathclap{
        \widehat G_x
      }$}"{description},shift right=1]
      \ar[d]
      \ar[
        r,
        "\exists"{above},
        "\sim"{below}
      ]
      &
      \widehat G_x \times_{\widehat {G_x}} \mathrm{U}_x
      \ar[out=180-66, in=66, looseness=4.2, "\scalebox{.77}{$\mathclap{
        \widehat G_x
      }$}"{description},shift right=1]
      \ar[d]
      \\
      \mathrm{U}_x
      \ar[out=-180+66, in=-66, looseness=4.2, "\scalebox{.77}{$\mathclap{
        G_x
      }$}"{description},shift left=1]
      \ar[r,-,shift left=1pt]
      \ar[r,-,shift right=1pt]
      &
      \mathrm{U}_x
      \mathrlap{\,.}
      \ar[out=-180+66, in=-66, looseness=4.2, "\scalebox{.77}{$\mathclap{
        G_x
      }$}"{description},shift left=1]
    \end{tikzcd}
  \end{equation}
\end{definition}

\begin{lemma}[Underlying principal bundles of Bierstone local model bundles are trivial]
  \label{EquivariantDirectProductBundlesFromHatGxTwistedProductBundles}
  Any Bierstone local model bundle \eqref{BierstoneLocalIdentification}
  is isomorphic, as a bundle of topological $G_x$-actions,
  to the projection out of the
  Cartesian product of $G_x \acts \, \mathrm{U}_x$ with a
  $G_x$-action $\rho_x$ on $\Gamma$:
  \vspace{-3mm}
  \begin{equation}
    \label{EquivariantPrincipalHatGxTwistedProductBundleIsCartesianProductProjectionOfActions}
    \begin{tikzcd}[column sep=-5pt, row sep=9pt]
      \widehat G_x
        \times_{\widehat {G_x}}
      \mathrm{U}_x
      \ar[out=180-60, in=60, looseness=3.2, "\scalebox{.77}{$\,\mathclap{
        G_x
      }\phantom{\cdot}$}"{description},shift right=2]
      \ar[
        dd,
        " \mathrm{pr}_2 \times_{\widehat {G_x}} \mathrm{id}_{\mathrm{U}_x}"{left}
      ]
      \ar[
        rr,
        "\sim"
      ]
      &{\phantom{AAAAAAAA}}&
      \Gamma
      \ar[out=180-66, in=66, looseness=4.5, "\scalebox{.77}{$\,\mathclap{
        G_x
      }\phantom{\cdot}$}"{description}, "\scalebox{.7}{$\!\!\rho_x$}"{very near end},shift right=1]
      &\times&
      \mathrm{U}_x
      \ar[out=180-66, in=66, looseness=4.5, "\scalebox{.77}{$\,\mathclap{
        G_x
      }\phantom{\cdot}$}"{description},shift right=1]
      \ar[
        dd,
        "\mathrm{pr}_2"
      ]
      \\
      \\
      G_x \times_{G_x} \mathrm{U}_x
      \ar[out=-180+60, in=-60, looseness=4.2, "\scalebox{.77}{$\,\,\mathclap{
        G_x
      }$}"{description},shift left=1]
      \ar[
        rrrr,
        "\sim"
      ]
      &&
      &&
      \mathrm{U}_x
      \ar[out=-180+60, in=-60, looseness=4.2, "\scalebox{.77}{$\,\,\mathclap{
        G_x
      }$}"{description},shift left=1]
    \end{tikzcd}
  \end{equation}
\end{lemma}
\begin{proof}
The point is that the $\widehat {G_x}$-action is transitive on the
$G_x$-factor in $\widehat G_x$, which implies that there are
canonical representatives of the elements in the twisted product quotient
(namely those with the neutral element in the $G_x$-factor)
passing to which yields,
first of all, an isomorphism of underlying topological spaces:
\vspace{-3mm}
\begin{equation}
  \label{IsomorphismRealizingBierstoneLocalModelsAsCartesianProductProjections}
  \begin{tikzcd}[row sep=-5pt]
    \widehat G_x
      \times_{\widehat {G_x}}
    \mathrm{U}_x
    \ar[
      rr,
      "\sim"
    ]
    &&
    \Gamma \times S
    \\
  \scalebox{0.8}{$   \left[
      (\gamma, \NeutralElement_{G}), s
    \right]
    $}
    &\longmapsto&
    \scalebox{0.8}{$
      (\gamma, s)
    $}
    \mathrlap{\,.}
  \end{tikzcd}
\end{equation}

\vspace{-2mm}
\noindent
One readily checks that, under this identification, the $G_x$-action on the left
turns into the claimed Cartesian product action:
\vspace{-3mm}
\begin{equation}
  \label{ComputingTheCartesianProductActionOnBierstoneLocalModelBundles}
  g
  \cdot
  \left[
    (\gamma,\, \NeutralElement_{{}_G}),\, u
  \right]
  \;=\;
  \left[
    (\alpha(g)(\gamma),\, g),\, u
  \right]
  \;=\;
  \left[
    (\alpha(g)(\gamma),\, g) \cdot \hat g^{-1},\, g \cdot u
  \right]
  \;=\;
  \left[
    (\rho_x(g)(\gamma),\, \NeutralElement_{{}_{G}}),\, g \cdot u
  \right]
  \,,
\end{equation}

\vspace{-1mm}
\noindent
with $(\gamma,g) \mapsto \rho_x(g)(\gamma) \in \Gamma$ defined by the
rightmost equation.
\end{proof}

\begin{remark}[Bierstone's local model bundles and generalization]
  \label{BierstoneLocalModelBundlesGeneralized}
  In the special case when $\alpha$ \eqref{IdentifyingEquivariantGroupsWithSemidirectProductGroups}
  is trivial,
  the equivariant bundles appearing on the right of
  \eqref{EquivariantPrincipalHatGxTwistedProductBundleIsCartesianProductProjectionOfActions}
  are the local model principal bundles considered in \cite{Bierstone73}.
  We may regard Lemma \ref{EquivariantDirectProductBundlesFromHatGxTwistedProductBundles}
  as providing the consistent generalization of Bierstone's notion
  to equivariant principal bundles with nontrivial action $\alpha$
  of the equivariance group on the structure group.
\end{remark}

\begin{example}[Fixed loci in Bierstone local model bundle for trivial $\alpha$]
  \label{FixedLociInBierstoneLocalModelBundles}
  Consider the special case when
  $\alpha$ \eqref{IdentifyingEquivariantGroupsWithSemidirectProductGroups}
  is trivial.
  Then any lift $\widehat {G_x}$ in \eqref{StabilizerSubgroupInSemidirectProductGroup}
  is equivalently the graph of an injective group homomorphism $\phi$
  \vspace{-2mm}
  $$
    \widehat {G_x}
    \;=\;
    \left\{
      (\phi(g),\, g)
      \,\vert\,
      g \in G_x
    \right\}
    \;\subset\;
    \Gamma \times G_x
    \,,
    {\phantom{AA}}
    G_x \xhookrightarrow{\phi} \Gamma
  $$

  \vspace{-2mm}
 \noindent
and the action $\rho_x$ in
\eqref{ComputingTheCartesianProductActionOnBierstoneLocalModelBundles}
  is just the right inverse multiplication action through this homomorphism
  \vspace{-2mm}
  $$
    \begin{tikzcd}[row sep=-5pt]
      G \times
      (
         \Gamma \times \mathrm{U}_x
      )
      \ar[
        rr,
        "\rho_x"
      ]
      &&
      \Gamma \times \mathrm{U}_x
      \\
   \scalebox{0.7}{$     (
        g, (\gamma, s)
      )
      $}
      &\longmapsto&
  \scalebox{0.7}{$      \left(
        \gamma \cdot \phi(g)^{-1},
        \,
        g \cdot s
      \right)
      $}
      \mathrlap{\,.}
    \end{tikzcd}
  $$

  \vspace{-2mm}
  \noindent
  Since this action is free, it has no $H$-fixed points as soon as
  $H \subset G_x$ is nontrivial.
  Therefore
  the $H$-fixed locus bundle
  (Cor. \ref{FixedLociOfEquivariantPrincipalBundles})
  of a Bierstone local model bundle
  for trivial $\alpha$ is necessarily
  of the following form:
  \vspace{-2mm}
  $$
    \big(
      \widehat G_x
        \times_{\widehat{G_x}}
      \mathrm{S}
    \big)^H
    \;\simeq\;
    \left\{
    \begin{array}{ccc}
      {
      \begin{tikzcd}[column sep=-4pt]
        \Gamma
        \ar[out=180-66, in=+66, looseness=3.8, "\scalebox{.77}{$\mathclap{
          \,G_x\,
        }\phantom{\cdot}$}"{description},"\scalebox{.7}{$\rho_x$}"{very near end},shift right=1]
        &\times&
        \mathrm{U}_x
        \ar[out=180-66, in=+66, looseness=3.8, "\scalebox{.77}{$\mathclap{
          G_x
        }\phantom{\cdot}$}"{description},shift right=1]
      \end{tikzcd}
      }
      &\vert&
      H = 1
      \\
      {
      \begin{tikzcd}[column sep=-4pt]
        \varnothing
        &\times&
        \mathrm{U}^H_x
        \ar[out=180-66, in=+66, looseness=3.8, "\scalebox{.77}{$\;\mathclap{
          W\!(H)
        }\phantom{-}$}"{description},shift right=1]
      \end{tikzcd}
      }
      &\vert&
      \mbox{otherwise}
    \end{array}
    \right.
  $$
\end{example}

\begin{proposition}[Bierstone local triviality is preserved by passage to $H$-fixed loci]
  \label{BierstoneLocalTrivialityIsPreservedByPassageToHFixedLoci}
  If a $G$-equivariant principal bundle is locally trivial in the sense of
  Bierstone (Def. \ref{BierstoneEquivariantLocalTrivializability}),
  then for all $H \subset G$ so is its $W\!(H)$-equivariant $H$-fixed locus
  bundle (Cor. \ref{FixedLociOfEquivariantPrincipalBundles}).
\end{proposition}
\begin{proof}
  Let $x \in \TopologicalSpace^H \subset \TopologicalSpace$ be an $H$-fixed point
  and consider a Bierstone local trivialization \eqref{BierstoneLocalIdentification}
  over an open neighborhood
  $\mathrm{U}_x \subset \TopologicalSpace$ with given lift $\widehat {G_x}$.
  It is sufficient to show that the
  $H$-fixed locus of this $G$-equivariant Bierstone local model bundle is $W\!(H)$-equivariantly
  isomorphic to a $W\!(H)$-equivariant Bierstone local model bundle:

  Since the fixed locus functor
  $(-)^H$ is a right adjoint \eqref{RightAdjointWeylGroupValuedFixedLocusFunctor},
  it preserves
  (Prop. \ref{RightAdjointFunctorsPreserveFiberProducts})
  the Cartesian product on the right hand side of
  the identification
  \eqref{EquivariantPrincipalHatGxTwistedProductBundleIsCartesianProductProjectionOfActions}
  of the Bierstone local model model with a Cartesian product
  (Lemma \ref{EquivariantDirectProductBundlesFromHatGxTwistedProductBundles}).
  Hence we just have to observe that this isomorphism
  \eqref{IsomorphismRealizingBierstoneLocalModelsAsCartesianProductProjections}
  passes to the fixed locus:
  \vspace{-3mm}
  \begin{equation}
    \label{FixedLocusOfBierstoneLocalModelBundle}
    \hspace{-2mm}
    \begin{tikzcd}[column sep=-5pt, row sep=20pt]
      \widehat G_x
        \times_{\widehat {G_x}}
      \mathrm{U}_x
      \ar[out=180-60, in=60, looseness=3.2, "\scalebox{.77}{$\,\mathclap{
        G_x
      }\phantom{\cdot}$}"{description},shift right=2]
      \ar[
        dd,
        " \mathrm{pr}_2 \times_{\widehat {G_x}} \mathrm{id}_{\mathrm{U}_x}"{left}
      ]
      \ar[
        rr,
        "\sim"
      ]
      &{\phantom{AAAAAAAA}}&
      \Gamma
      \ar[out=180-66, in=66, looseness=4.5, "\scalebox{.77}{$\,\mathclap{
        G_x
      }\phantom{\cdot}$}"{description}, "\scalebox{.7}{$\!\!\rho_x$}"{very near end},shift right=1]
      &\times&
      \mathrm{U}_x
      \ar[out=180-66, in=66, looseness=4.5, "\scalebox{.77}{$\,\mathclap{
        G_x
      }\phantom{\cdot}$}"{description},shift right=1]
      \ar[
        dd,
        "\mathrm{pr}_2"
      ]
      \\
      \\
      G_x \times_{G_x} \mathrm{S}
      \ar[out=-180+60, in=-60, looseness=4.2, "\scalebox{.77}{$\,\,\mathclap{
        G_x
      }$}"{description},shift left=1]
      \ar[
        rrrr,
        "\sim"
      ]
      &&
      &&
      \mathrm{U}_x
      \ar[out=-180+60, in=-60, looseness=4.2, "\scalebox{.77}{$\,\,\mathclap{
        G_x
      }$}"{description},shift left=1]
      \\
   \scalebox{0.8}{$   [
        ( \gamma,\, e_{{}_{G_x}} ), u
      ]
      $}
      &
      \longmapsto
      &
      &
 \scalebox{0.8}{$     \mathclap{
      (
        \gamma, u
      )
      }
      $}
    \end{tikzcd}
    {\phantom{AAAA}}
    \Rightarrow
    {\phantom{AAAA}}
    \begin{tikzcd}[column sep=-5pt, row sep=20pt]
      (\widehat {W\!(H)})_x
        \times_{\widehat {W(H)_x}}
      \mathrm{U}^H_x
      \ar[out=180-60, in=60, looseness=3.2, "\scalebox{.77}{$\phantom{-}\mathclap{
        W\!(H)_x
      }\phantom{-}$}"{description},shift right=2]
      \ar[
        dd,
        " \mathrm{pr}_2 \times_{\widehat {W(H)_x}} \mathrm{id}_{\mathrm{U}^H_x}"{left}
      ]
      \ar[
        rr,
        "\sim"
      ]
      &{\phantom{AAAAAAAA}}&
      \Gamma^H
      \ar[out=180-66, in=66, looseness=4.5, "\scalebox{.77}{$\,\phantom{-}\mathclap{
        W\!(H)_x
      }\phantom{-}$}"{description},shift right=1]
      &\times&
      \mathrm{U}^H_x
      \ar[out=180-66, in=66, looseness=4.5, "\scalebox{.77}{$\phantom{-}\mathclap{
        W\!(H)_x
      }\phantom{-}$}"{description},shift right=1]
      \ar[
        dd,
        "\mathrm{pr}^H_2"
      ]
      \\
      \\
      W\!(H)
        \times_{W(H)_x}
      \mathrm{U}^H_x
      \ar[out=-180+60, in=-60, looseness=3.8, "\scalebox{.77}{$\phantom{-}\mathclap{
        W(H)_x
      }\phantom{-}$}"{description},shift left=1]
      \ar[
        rrrr,
        "\sim"
      ]
      &&
      &&
      \mathrm{U}^H_x
      \ar[out=-180+60, in=-60, looseness=4.2, "\scalebox{.77}{$\phantom{-}\mathclap{
        W\!(H)_x
      }\phantom{-}$}"{description},shift left=1]
      \\
\scalebox{0.8}{$      \left[
        ( \gamma,\, e_{{}_{W(H)_x}}), u
      \right]
      $}
      &
      \longmapsto
      &
      &
   \scalebox{0.8}{$   \mathclap{
      (
        \gamma, u
      )
      }
      $}
    \end{tikzcd}
  \end{equation}

    \vspace{-2mm}
\noindent
Here in the left column of the right square we systematically follow our notation conventions
(Ntn.
\ref{GActionOnTopologicalSpaces},
\ref{LiftsOfEquivarianceSubgroupsToSemidirectProductWithStructureGroup},
\ref{LiftOfIsotropyGroupsToSemidirectProductWithStructureGroup})
to the new equivariance group $W\!(H)$ at the fixed locus
  \vspace{-2mm}
$$
  W\!(H)_x
  \;\coloneqq\;
  \mathrm{Stab}_{W(H)}(x)
  \;=\;
  \left(
    G_x \cap N\!(H)
  \right)/ (G_x \cap H)
  \,,
  \;\;\;\;\;\;
  (\widehat {W(H)})_x
  \;\coloneqq\;
  \Gamma^H \rtimes_\alpha W\!(H)_x
$$

  \vspace{-2mm}
\noindent
and take the lift to be the subquotient of the given lift
  \vspace{-2mm}
$$
  \widehat {W(H)_x}
  \;\coloneqq\;
  \widehat{
    \scalebox{.8}{$G_x \cap N\!(H)$}
  }
  /
  H
  \,.
$$

  \vspace{-2mm}
\noindent
The computation \eqref{ComputingTheCartesianProductActionOnBierstoneLocalModelBundles}
now applies verbatim to check that
the homeomorphism \eqref{FixedLocusOfBierstoneLocalModelBundle}
is indeed $W\!(H)_x$-equivariant
  \vspace{-1mm}
\begin{equation}
  w
  \cdot
  \left[
    (\gamma,\, e_{{}_{W(H)}}),\, u
  \right]
  \;=\;
  \big[
    ( \alpha(w)(\gamma),\, w),\, u
  \big]
  \;=\;
  \left[
    ( \alpha(w)(\gamma),\, w) \cdot \hat w^{-1},\, w \cdot u
  \right]
  \;=\;
  \left[
    \big(\rho_x(w)(\gamma),\, e_{{}_{W(H)}}\big),\, w \cdot u
  \right]
  \,.
\end{equation}

  \vspace{0mm}
\noindent In conclusion, this shows that if a $G$-equivariant principal bundle
is isomorphic in a neighborhood $\mathrm{U}_x$ of an $H$-fixed point $x$
to a Bierstone local model bundle of the form \eqref{EquivariantPrincipalHatGxTwistedProductBundleIsCartesianProductProjectionOfActions},
then its $H$-fixed locus is $W\!(H)$-equivariantly isomorphic over $\mathrm{U}^H_x$
to the Bierstone local model \eqref{FixedLocusOfBierstoneLocalModelBundle}.
\end{proof}

\medskip

We proceed to discuss the relation of Bierstone's local trivalizability condition to
that of tom Dieck and of Lashof-May.

\vspace{-3mm}
\begin{lemma}[Intersection with open neighborhoods preserves slices through points]
  \label{IntersectionWithOpenNeighborhoodsPreservesSlicesThroughPoints}
Let
$
\!\!
  \begin{tikzcd}
      x
      \in
      \mathrm{U}_x
      \ar[out=180-66, in=66, looseness=4.2, "\scalebox{.77}{$\mathclap{
        G_x
      }$}"{description},shift right=1]
      \ar[
        r,
        hook,
        "\mathrm{open}"{above}
      ]
      &
      \TopologicalSpace
      \ar[out=180-66, in=66, looseness=4.2, "\scalebox{.77}{$\mathclap{
        G_x
      }$}"{description},shift right=1]
    \end{tikzcd}
  $
  be a $G_x$-equivariant open neighborhood
  in a given $G \acts \, \TopologicalSpace \in \GActionsOnTopologicalSpaces$,
  of some $x \in \TopologicalSpace$.
Then, with every
  slice
  \newline
  $
    \begin{tikzcd}[column sep=-5]
      x
      &\in
      &
      \mathrm{S}'_x
      \ar[out=180-66, in=66, looseness=3.5, "\scalebox{.77}{$\phantom{}\mathclap{
        G_x
      }\phantom{\cdot}$}"{description}, shift right=1]
      \ar[
        rr,
        hook,
        "\iota"
      ]
      &{\phantom{AA}}&
      G \cdot \mathrm{S}'_x
      \ar[out=180-66, in=66, looseness=3.5, "\scalebox{.77}{$\phantom{}\mathclap{
        G_x
      }\phantom{\cdot}$}"{description}, shift right=1]
    \end{tikzcd}
   $
   \eqref{ASliceThroughAPoint}
   through $x$,  the intersection
  $
    \begin{tikzcd}[column sep=-5pt]
      \mathrm{S}_x
      \ar[out=180-66, in=66, looseness=3.5, "\scalebox{.77}{$\phantom{}\mathclap{
        G_x
      }\phantom{\cdot}$}"{description}, shift right=1]
      &
       \;
        \coloneqq
       \;
      &
      \mathrm{S}'_x
      \ar[out=180-66, in=66, looseness=3.5, "\scalebox{.77}{$\phantom{}\mathclap{
        G_x
      }\phantom{\cdot}$}"{description}, shift right=1]
      &\cap&
      \mathrm{U}_x
      \ar[out=180-66, in=66, looseness=3.5, "\scalebox{.77}{$\phantom{}\mathclap{
        G_x
      }\phantom{\cdot}$}"{description}, shift right=1]
    \end{tikzcd}
  $
  is also a slice through $x$ (Def. \ref{SliceOfTopologicalGSpace}).
\end{lemma}
\begin{proof}
  We need to see that the image
    \vspace{-4mm}
  $$
    \begin{tikzcd}
      G \times_{G_x}
      (
        \mathrm{S}'_x \cap \mathrm{U}_x
      )
      \ar[
        r,
        "{\widetilde {\iota_{\vert \mathrm{U}_x}}}"{above},
        "\sim"{below}
      ]
      &
      G \cdot ( \mathrm{S}'_x \cap \mathrm{U}_x )
      \ar[
        r,
        hook
      ]
      &
      G \cdot \mathrm{S}'_x
      \ar[
        r,
        hook,
        "\mathrm{open}"{above}
      ]
      &
      \TopologicalSpace
    \end{tikzcd}
  $$

  \vspace{-3mm}
\noindent
  is still an open subset of $\TopologicalSpace$, hence an open neighborhood of $x$.
  Observe that for any slice $\mathrm{S}$, under Assump. \ref{ProperEquivariantTopology},
  the map
  $
  \!\!
    \begin{tikzcd}[column sep =40pt]
      G \times \mathrm{S}
      \ar[
        r,
        "{(g,s) \mapsto g\cdot s}"
      ]
      &
      G \cdot \mathrm{S}
    \end{tikzcd}
    \!\!
  $
  is an open map (\cite[Thm. 2.1 (1)]{Antonyan17}),
  hence sends open subsets to open subsets.
  With this, and since
  $
    \begin{tikzcd}
      G \times
      (
        \mathrm{S}'_x \cap \mathrm{U}_x
      )
      \ar[r,hook,"\mathrm{open}"{above}]
      &
      G \times \mathrm{S'}_x
    \end{tikzcd}
  $
  is open by assumption on $\mathrm{U}_x$, it follows that its image
  is open in $G \cdot \mathrm{S}'_x$, which in turn is open
  in $\TopologicalSpace$ by assumption on $\mathrm{S}'_x$.
\end{proof}

\begin{notation}[Slices inside Bierstone patches]
  \label{SlicesInsideBierstonePatches}
  Given an equivariant local trivialization
  \eqref{BierstoneLocalIdentification}
  in the sense of Bierstone
  (Def. \ref{BierstoneEquivariantLocalTrivializability}),
  we
  may choose
  (by the slice theorem, Prop. \ref{SliceTheorem}, using Assump. \ref{ProperEquivariantTopology})
  for each point $x \in \TopologicalSpace$
  a $G_x$-slice $\mathrm{S}'_x$ through that point. Furthermore,  by
  Lemma \ref{IntersectionWithOpenNeighborhoodsPreservesSlicesThroughPoints},
  we obtain from this a slice $\mathrm{S}_x \coloneqq \mathrm{S}'_x \cap \mathrm{U}_x$
  which is still through $x$ but also contained in $\mathrm{U}_x$:
  \vspace{-5mm}
  \begin{equation}
    \label{SliceThroughBierstonePatch}
    \begin{tikzcd}[column sep=small, row sep=1pt]
      &[-18pt] &[-18pt] &&
      \mathrm{U}_x
      \ar[out=180-66, in=66, looseness=3.5, "\scalebox{.77}{$\phantom{}\mathclap{
        G_x
      }\phantom{\cdot}$}"{description}, shift right=1]
      \ar[
        drr,
        hook
      ]
      \\
      x
      &
      \in
      &
      \mathrm{S}_x
      \ar[out=180-66, in=66, looseness=3.6, "\scalebox{.77}{$\phantom{}\mathclap{
        G_x
      }\phantom{\cdot}$}"{description}, shift right=1]
      \ar[
        drr,
        hook
      ]
      \ar[
        urr,
        hook
      ]
      &{\phantom{AA}}&
      &{\phantom{AA}}&
      \TopologicalSpace
      \ar[out=-180+66, in=-66, looseness=4.2, "\scalebox{.77}{$\phantom{}\mathclap{
        G
      }$}"{description}, shift left=1]
      \\
      && &&
      G \cdot \mathrm{S}'_x
      \ar[out=-180+66, in=-66, looseness=3.6, "\scalebox{.77}{$\phantom{}\mathclap{
        G
      }$}"{description}, shift left=1]
      \ar[
        urr,
        hook
      ]
    \end{tikzcd}.
  \end{equation}
\end{notation}

\begin{lemma}[Bierstone local trivializations restricted to slices of the base are slices in the total space]
  \label{BierstoneLocalTrivializationsRestrictedToSlicesAreSlicesInTotalSpace}
  Given an equivariant local trivialization
  \eqref{BierstoneLocalIdentification}
  in the sense of Bierstone
  (Def. \ref{BierstoneEquivariantLocalTrivializability}),
  then the further restrictions of the equivariant bundle
  to slices $\mathrm{S}_x \hookrightarrow \mathrm{U}_x$ inside the Bierstone patches
  (according to Nota \ref{SlicesInsideBierstonePatches})
  are slices (Def. \ref{SliceOfTopologicalGSpace})
  in the total space of the bundle.
  That is, the induction/restriction-adjunct
  \eqref{HomIsomorphismForRestrictedActionAndInducedAction}
  of their inclusions are isomorphisms:
    \vspace{-3mm}
  \begin{equation}
    \label{AdjunctOfRestrictionOfPrincipalBundleToSlice}
    \begin{tikzcd}[row sep=small]
      \mathrm{P}\vert_{\mathrm{S}_x}
      \ar[out=180-66, in=66, looseness=3.5, "\scalebox{.77}{$\mathclap{
        G_x
      }\phantom{\cdot}$}"{description}, shift right=1]
      \ar[d]
      \ar[
        rr,
        hook,
        "P_{\vert \iota_x}"
      ]
      \ar[
        rrd,
        phantom,
        "\mbox{\tiny\rm(pb)}"{description}
      ]
      &&
      \mathrm{P}\vert_{ G \cdot \mathrm{S}_x}
      \ar[out=180-66, in=66, looseness=3.5, "\scalebox{.77}{$\mathclap{
        G_x
      }\phantom{\cdot}$}"{description}, shift right=1]
      \ar[
        d
      ]
      \\
      \mathrm{S}_x
      \ar[out=-180+66, in=-66, looseness=4.2, "\scalebox{.77}{$\phantom{\cdot}\mathclap{
        G_x
      }$}"{description}, shift left=1]
      \ar[
        rr,
        hook,
        "\iota_x"{below}
      ]
      &&
      G \cdot \mathrm{S}_x
      \ar[out=-180+66, in=-66, looseness=4.2, "\scalebox{.77}{$\phantom{\cdot}\mathclap{
        G_x
      }$}"{description}, shift left=1]
    \end{tikzcd}
    {\phantom{AAAAA}}
    \xleftrightarrow{
      \;\;\;\;
      \widetilde{(-)}
      \;\;\;\;
    }
    {\phantom{AAAAA}}
    \begin{tikzcd}[row sep=small]
      G
        \times_{G_x}
      \mathrm{P}\vert_{ \mathrm{S}_x}
      \ar[out=180-66, in=66, looseness=3.5, "\scalebox{.77}{$\mathclap{
        G
      }$}"{description}, shift right=2]
      \ar[d]
      \ar[
        rr,
        "\sim"{below},
        "{
          \widetilde{
            P\vert_{ \iota_x}
          }
        }"{above}
      ]
      &&
      \mathrm{P}\vert_{ G\cdot \mathrm{S}_x}
      \ar[out=180-66, in=66, looseness=3.5, "\scalebox{.77}{$\mathclap{
        G
      }$}"{description}, shift right=1]
      \ar[d]
      \\
      G \times_{G_x} \mathrm{S}_x
      \ar[out=-180+66, in=-66, looseness=4.2, "\scalebox{.77}{$\mathclap{
        G
      }$}"{description}, shift left=1]
      \ar[
        rr,
        "\sim"{above},
        "\widetilde {\iota_x}"{below}
      ]
      &&
      G \cdot \mathrm{S}_x
      \ar[out=-180+66, in=-66, looseness=4.2, "\scalebox{.77}{$\mathclap{
        G
      }$}"{description}, shift left=1]
    \end{tikzcd}
  \end{equation}
\end{lemma}
\begin{proof}
  First we claim that the underlying $\Gamma$-principal bundle
  (Cor. \ref{UnderlyingPrincipalBundles}) of
  $G \times _{G_x} \mathrm{P}\vert_{ \mathrm{S}_x}$
  is locally trivializable:

  By Lemma \ref{LashofBundlesOverSliceOrbitsAreLocallyTrivialAsOrdinaryPrincipalBundles},
  there exists a tubular neighborhood
  $D \times \mathrm{S}_x \xhookrightarrow{\;} G \cdot \mathrm{S}_x $
  (using Assump. \ref{ProperEquivariantTopology}).
  Since this is contractible to $\mathrm{S}_x$,
  and since
  the underlying $\Gamma$-principal bundle
  of $P\vert_{ \mathrm{S}_x}$ is trivializable -- being the
  restriction, along \eqref{SliceThroughBierstonePatch},
  of that of $\mathrm{P}\vert_{ \mathrm{U}_x}$, which is trivializable by
  assumption \eqref{BierstoneLocalIdentification} and
  by Lemma \ref{EquivariantDirectProductBundlesFromHatGxTwistedProductBundles} --
  it follows that
  the underlying bundle of $\mathrm{P}\vert_{ D \times \mathrm{S}_x}$
  is trivializable. By $G$-equivariance the same is then true for all
  its $G$-translates; and since the $G$-translates of its base space
  $D \times \mathrm{S}_x \xhookrightarrow{\;} G \cdot \mathrm{S}_x$
  form a cover, this implies the claim.

  Second, we claim that the adjunct square on the right of
  \eqref{AdjunctOfRestrictionOfPrincipalBundleToSlice} is a morphism of
  $\Gamma$-principal bundles, hence that its top map is $\Gamma$-equivariant.
  A formal way to see this is to notice,
  with Lemma \ref{BeckChevalleyConditionForLiftOfIsotropyGroup},
  that the top morphism
  is equivalently the induction/restriction-adjunct \eqref{HomIsomorphismForRestrictedActionAndInducedAction}
  along the lifted inclusion
  $\widehat G_x \xhookrightarrow{\;} \widehat G$ (Ntn. \ref{LiftOfIsotropyGroupsToSemidirectProductWithStructureGroup}).

  These two claims together imply that the square on the right of
  \eqref{AdjunctOfRestrictionOfPrincipalBundleToSlice} is a pullback square,
  by Lemma \ref{HomomorphismsOfLocallyTrivialPrincipalBundlesArePullbackSquares}.
  Since its bottom morphism is an isomorphism
  \eqref{SliceIsomorphism}
  by the assumption that $\mathrm{S}_x$ is a slice, the top morphism
  is exhibited as the pullback of an iosmorphism and thus is itself an
  isomorphism (Example \ref{PullbackPreservesIsomorphisms}),
  which is the statement to be shown.
\end{proof}

\begin{remark}[Comparison to the literature]
  \label{ComparingLocalTrivializationProofOfLashofMayBundlesToTheLiterature}
  In \cite[Proof of Lem. 1.3]{Lashof81} the adjunct \eqref{AdjunctOfRestrictionOfPrincipalBundleToSlice}
  is written down in the special case when $\alpha$ is trivial
  (Rem. \ref{BierstoneLocalModelBundlesGeneralized})
  but no reason is offered for it being an isomorphism.
  While it is clear by point-set analysis that the underlying function is a
  bijection, its homeomorphy is a little subtle.
  The above category-theoretic proof of Lem.
  \ref{BierstoneLocalTrivializationsRestrictedToSlicesAreSlicesInTotalSpace}
  makes this transparent and immediately generalizes the argument to
  equivariant bundles internal to other ambient categories,
  such as to differentiable equivariant bundles.
\end{remark}

\begin{proposition}[Bierstone's local triviality implies Lashof-May's]
  \label{BierstoneLocalTrivialityImpliesLashofMay}
  An equivariant principal bundle over $\TopologicalSpace$
  (Def. \ref{EquivariantPrincipalBundle},
  Cor. \ref{InternalDefinitionOfGPrincipalBundlesCoicidesWithtomDieckDefinition})
  that is locally trivial
  in the sense of Bierstone (Def. \ref{BierstoneEquivariantLocalTrivializability})
  is also locally trivial in the sense of Lashof \&  May
  (Def. \ref{LashofMayEquivariantLocalTrivializability}).
\end{proposition}
\begin{proof}
  For $x \in \TopologicalSpace$, the given
  Bierstone trivialization \eqref{BierstoneLocalIdentification}
  over a Bierstone patch $\mathrm{U}_x$ further restricted to
  a slice $\mathrm{S}_x$ through $x$ inside $\mathrm{U}_x$
  (Ntn. \ref{SlicesInsideBierstonePatches}) yields the following diagram,
    by Lemma \ref{RecognitionOfCartesianQuotientProjections}
   \vspace{-2mm}
  \begin{equation}
    \label{BierstoneLocalPatchRestrictedToSlice}
    \begin{tikzcd}[row sep=small]
      \mathrm{P}\vert_{\mathrm{S}_x}
      \ar[
        r,
        "\sim"
      ]
      \ar[d]
      &
      \widehat G_x
        \times_{\widehat{G_x}}
      \mathrm{S}_x
      \ar[
        d
      ]
      \\
      \mathrm{S}_x
      \ar[r,-,shift left=1pt]
      \ar[r,-,shift right=1pt]
      &
      \mathrm{S}_x
    \end{tikzcd}
  \end{equation}

\vspace{-2mm}
\noindent
using that
  the underlying $\Gamma$-principal bundle (Cor. \ref{UnderlyingPrincipalBundles}) of
  $\mathrm{P}\vert_{\mathrm{S}_x}$ is trivializable, by Lemma \ref{EquivariantDirectProductBundlesFromHatGxTwistedProductBundles}.
Consider then the following pasting composite
  of the
  image of this diagram
  \eqref{BierstoneLocalPatchRestrictedToSlice}
  under the induction functor $G \times_{G_x} (-)$ \eqref{InducedRestrictedActionAdjunction}
  with, on the left,
  the inverse of the identification \eqref{AdjunctOfRestrictionOfPrincipalBundleToSlice}
  from Lemma \ref{BierstoneLocalTrivializationsRestrictedToSlicesAreSlicesInTotalSpace},
  and, on the right,
  the equivalence induced via Lemma \ref{BeckChevalleyConditionForLiftOfIsotropyGroup}:
  \vspace{-2mm}
  $$
    \begin{tikzcd}[column sep=large]
      \mathrm{P}\vert_{G \cdot \mathrm{S}_x}
      \ar[
        r,
        "\sim"
      ]
      \ar[
        d
      ]
      &
      G
        \times_{G_x}
      \mathrm{P}\vert_{\mathrm{S}_x}
      \ar[
        r,
        "\sim"
      ]
      \ar[d]
      &
      G
        \times_{G_x}
      \widehat G_x
        \times_{\widehat{G_x}}
      \mathrm{S}_x
      \ar[
        d
      ]
      \ar[
        r,
        "\sim"
      ]
      &
      \widehat G
        \times_{\widehat {G_x}}
      \mathrm{S}_X
      \ar[
        d,
        " \mathrm{pr}_2 \times_{\widehat{G_x}} \mathrm{id}"
      ]
      \\
      G \cdot \mathrm{S}_x
      \ar[
        r,
        "\sim"
      ]
      &
      G
        \times_{G_x}
      \mathrm{S}_x
      \ar[r,-,shift left=1pt]
      \ar[r,-,shift right=1pt]
      &
      G
        \times_{G_x}
      \mathrm{S}_x
      \ar[r,-, shift left=1pt]
      \ar[r,-, shift right=1pt]
      &
      G
        \times_{G_x}
      \mathrm{S}_x
      \mathrlap{\,.}
    \end{tikzcd}
  $$

  \vspace{-1mm}
\noindent
  The resulting composite rectangle manifestly exhibits a
  Lashof-May local trivialization \eqref{LashofMayLocalModelBundle},
  with open subsets indexed by the points of $\TopologicalSpace$.
\end{proof}

\begin{proposition}[tom Dieck's local trivialization implies Bierstone's]
  \label{tomDieckLocalTrivializationImpliesBierstone}
  An equivariant principal bundle over $\TopologicalSpace$ (Def. \ref{EquivariantPrincipalBundle},
  Cor. \ref{InternalDefinitionOfGPrincipalBundlesCoicidesWithtomDieckDefinition})
  that is locally trivial in the sense of tom Dieck (Def. \ref{tomDieckLocalTrivializability})
  is also locally trivial in the sense of Bierstone
  (Def. \ref{BierstoneEquivariantLocalTrivializability}).
\end{proposition}
\begin{proof}
Consider a tom Dieck local trivialization (Def. \ref{tomDieckLocalTrivializability}),
over patches $\mathrm{U}_i$ \eqref{EquivariantLocalTrivialityCondition}
given via Prop. \ref{CharacterizationOfEquivariantBundlesOverCosetSpaces}
by $\widehat {H_i}$-coset space coprojections over $H_i$-cosets
(Ntn. \ref{LiftsOfEquivarianceSubgroupsToSemidirectProductWithStructureGroup})
for closed subgroups $H_i \subset G$.

In view of
Ntn. \ref{LiftsOfEquivarianceSubgroupsToSemidirectProductWithStructureGroup}
for the given tom Dieck local models
and
Ntn.
\ref{LiftOfIsotropyGroupsToSemidirectProductWithStructureGroup}
for the Bierstone local models that we have to reproduce,  it is suggestive to
write $\widehat G_i$ for the corresponding semidirect product subgroups:
\vspace{-3mm}
\begin{equation}
\label{tomDieckToBierstoneSemidirectProductSubgroups}
\begin{tikzcd}[column sep=0pt, row sep=small]
  \widehat G_i
  \ar[
    d,
    hook
  ]
  &\coloneqq&
  \Gamma \rtimes_\alpha H_i
  \ar[
    d,
    hook
  ]
  \\
  \widehat G
  &\coloneqq&
  \Gamma \rtimes_\alpha H_i
\end{tikzcd}
\end{equation}

\vspace{-1mm}
\noindent
By Lemma \ref{LashofBundlesOverSliceOrbitsAreLocallyTrivialAsOrdinaryPrincipalBundles}
(for $\mathrm{S} = H_i$),
we may find an open ball
\eqref{OpenBallAtNeutralElementNormalToClosedSubgroupH}
around $\NeutralElement \in G$ normal to $H_i \subset G$
\vspace{-2mm}
$$
  D_i \,=\, D^\epsilon N_{\NeutralElement} H_i \;\subset\; G
$$

\vspace{-2mm}
\noindent
which is an $H$-equivariant subspace under the conjugation action
of $H$ on itself
\vspace{-2mm}
\begin{equation}
  \begin{tikzcd}
    \label{ConjugationActionOnNormalBall}
    D_i
    \ar[out=180-60, in=60, looseness=3.2, "\scalebox{.77}{$\,\mathclap{
      H_i
    }\phantom{\cdot}$}"{description},"\scalebox{.6}{$\mathrm{Ad}$}"{very near end},shift right=1]
    \ar[
      r,
      hook
    ]
    &
    H_i
    \ar[out=180-60, in=60, looseness=3.2, "\scalebox{.77}{$\,\mathclap{
      H_i
    }\phantom{\cdot}$}"{description},"\scalebox{.6}{$\mathrm{Ad}$}"{very near end},shift right=1]
  \end{tikzcd}
\end{equation}

\vspace{-3mm}
\noindent and such that we have an isomorphism
\vspace{-4mm}
\begin{equation}
  \label{SplittingOfDcdotHi}
  \begin{tikzcd}
    D_i \times H_i
    \ar[
      rr,
      "\sim"{below},
      "{
        (d,h) \,\mapsto\, d \cdot h
      }"
    ]
    &&
    D_i \cdot H_i
  \end{tikzcd}
\end{equation}

\vspace{-2mm}
\noindent
and hence an $H_i$-equivariant open neighborhood
\vspace{-4mm}
\begin{equation}
  \label{OpenNeighborhoodInCosetSpace}
  \begin{tikzcd}[column sep=40pt]
    D_i
    \ar[out=180-60, in=60, looseness=3.2, "\scalebox{.77}{$\,\mathclap{
      H_i
    }\phantom{\cdot}$}"{description},"\scalebox{.6}{$\mathrm{Ad}$}"{very near end},shift right=1]
    \ar[
      r,
      hook,
      "\mathrm{open}"{below}
    ]
    &
    G/H_i
    \mathrlap{\,.}
    \ar[out=180-60, in=60, looseness=3.2, "\scalebox{.77}{$\,\mathclap{
      H_i
    }\phantom{\cdot}$}"{description},shift right=1]
  \end{tikzcd}
\end{equation}

\vspace{-3mm}
\noindent
Using this isomorphism \eqref{SplittingOfDcdotHi},
observe the following $G \times \widehat{H_i}$-equivariant homeomorphism
(where $h \mapsto \widehat h$ is the given lift \eqref{LiftOfSubgroupsHToSemidirectProductGroup}):
\vspace{-6mm}
$$
  \begin{tikzcd}[column sep=small, row sep=0pt]
    \Gamma \rtimes_\alpha (H_i \cdot D_i)
    \ar[out=180-60, in=60, looseness=4, "\scalebox{.77}{$\phantom{.}\mathclap{
      G \times \widehat{H_i}
    }\phantom{-}$}"{description},shift right=1]
    \ar[
      rr,
      "\sim"
    ]
    &&
    \left(
      \Gamma \rtimes_\alpha H_i
    \right)
      \times
    D_i
    \ar[out=180-60, in=60, looseness=4, "\scalebox{.77}{$\phantom{.}\mathclap{
      G \times \widehat{H_i}
    }\phantom{-}$}"{description},shift right=1]
    \\
\scalebox{0.7}{$    (
      \gamma,\,
      h' \cdot d
    )
    $}
    \ar[
      dd,
      |->,
      "\widehat h \,\in\, \widehat{H_i}"
    ]
    &
      \longmapsto
    &
  \scalebox{0.7}{$  \big(
      (
        \alpha(d^{-1})(\gamma)
        \,,
        h'
      )
      ,\,
      d
    \big)
    $}
    \ar[
      dd,
      |->,
      "\widehat h \,\in\, \widehat{H_i} "
    ]
    \\
    {\phantom{A}}
    \\
\scalebox{0.7}{$    \big(
      \gamma \cdot \alpha(d \cdot h')(\mathrm{pr}_1(\hat h^{-1}))
      ,\,
      h' h^{-1} \cdot  h d h^{-1}
    \big)
    $}
    &
      \longmapsto
    &
\scalebox{0.7}{$    \Big(
      \big(
        \alpha(d^{-1})(\gamma)
          \cdot
        \alpha(h')(\mathrm{pr}_1(\hat h^{-1}))
        ,\,
        h' h^{-1}
      \big)
      ,\,
      h  d h^{-1}
    \Big)
    $}
    \mathrlap{\,.}
      \end{tikzcd}
$$

\vspace{-3mm}
\noindent
Here the vertical maps define the $\widehat {H_i}$-action,
given by inverse right multiplication on $\Gamma \rtimes_\alpha G$
and by the conjugation action \eqref{ConjugationActionOnNormalBall}
(through $\mathrm{pr}_2$), and where the remaining $G$-action is
by left multiplication (through $s$).
This is such that passage
to $\widehat{H_i}$-quotients (Example \ref{QuotientSpaces})
yields a $G$-equivariant homeomorphism of the following form,
where on the right we use the adapted notation \eqref{tomDieckToBierstoneSemidirectProductSubgroups}:
\vspace{-3mm}
\begin{equation}
  \label{RestrictionOftomDieckLocalModelToInvarinantOpenNeighborhood}
  \begin{tikzcd}
    \left(
      \Gamma
        \rtimes_\alpha
      (H_i \cdot D_i)
    \right)
      /
    \widehat {H_i}
    \ar[
      rr,
      "\sim"
    ]
    &&
    \widehat G_i
      \times_{ \widehat{H_i} }
    D_i \;.
  \end{tikzcd}
\end{equation}

\vspace{-3mm}
\noindent
It follows (with Lemma \ref{RecognitionOfCartesianQuotientProjections})
that the following square is a pullback:
\vspace{-2mm}
\begin{equation}
  \label{RestrictionOftomDieckLocalModelToInvariantNeighborhood}
  \begin{tikzcd}[row sep=small]
    \widehat G_i
      \times_{H_i}
    D_i
    \ar[
      rr
    ]
    \ar[
      d
    ]
    \ar[
      drr,
      phantom,
      "\mbox{\tiny\rm(pb)}"
    ]
    &&
    \widehat G / \widehat{H_i}
    \ar[
      d
    ]
    \\
    D_i
    \ar[
      rr,
      hook
    ]
    &&
    G/H_i
  \end{tikzcd}
\end{equation}

\vspace{-2mm}
\noindent
because the top morphism is manifestly $\Gamma$-equivariant, while the
left vertical morphism is trivial as a $\Gamma$-principal bundle
(by Lemma \ref{EquivariantDirectProductBundlesFromHatGxTwistedProductBundles}).

Now consider any point $x \in \mathrm{U}_i$.
Without restriction of generality, we may assume that this point is
taken by the classifying map of the given
tom Dieck local trivialization to
$[\NeutralElement_G] \in G/H_i$,
for if not
then we may adjust the local trivialization,
by  Ex. \ref{AdjustingTheClassifyingMapsInAtomDieckLocalTrivialization}.

Then consider the open neighborhood of $x$ which is the preimage of
$\mathrm{D}_i$ under the classifying map:
\vspace{-2mm}
\begin{equation}
  \label{OpenNeighourhoodOfPointIntomDieckCoverFromOpenNeighbouthodOfeInCosetSpace}
  \begin{tikzcd}
    \mathrm{U}_x
    \ar[out=180-66, in=66, looseness=3.2, "\scalebox{.77}{$\mathclap{
      H_i
    }\,$}"{description}, shift right=1]
    \ar[rr]
    \ar[
      d,
      hook
    ]
    \ar[
      drr,
      phantom,
      "\mbox{\tiny\rm(pb)}"
    ]
    &&
    D_i
    \ar[
      d,
      hook
    ]
    \ar[out=180-66, in=66, looseness=3.2, "\scalebox{.77}{$\mathclap{
      H_i
    }\,$}"{description}, shift right=1]
    \\
    \mathrm{U}_i
    \ar[
      rr,
      "f_i"{below}
    ]
    \ar[out=-180+66, in=-66, looseness=3.6, "\scalebox{.77}{$\,\mathclap{
      H_i
    }$}"{description}, shift left=1]
    &&
    G/H_i
    \ar[out=-180+66, in=-66, looseness=3.6, "\scalebox{.77}{$\,\mathclap{
      G_x
    }$}"{description}, shift left=1]
  \end{tikzcd}
  {\phantom{AAAAAAA}}
  \begin{tikzcd}
    \{x\}
    \ar[out=180-66, in=66, looseness=3.2, "\scalebox{.77}{$\mathclap{
      G_x
    }\,$}"{description}, shift right=1]
    \ar[
      r,
      hook
    ]
    \ar[
      dr,
      hook
    ]
    &
    \mathrm{U}_x
    \ar[out=180-66, in=66, looseness=3.2, "\scalebox{.77}{$\mathclap{
      G_x
    }\,$}"{description}, shift right=1]
    \ar[rr]
    \ar[
      d,
      hook
    ]
    \ar[
      drr,
      phantom,
      "\mbox{\tiny\rm(pb)}"
    ]
    &&
    D_i
    \ar[
      d,
      hook
    ]
    \ar[out=180-66, in=66, looseness=3.2, "\scalebox{.77}{$\mathclap{
      G_x
    }\,$}"{description}, shift right=1]
    \\
    &
    \mathrm{U}_i
    \ar[
      rr,
      "f_i"{below}
    ]
    \ar[out=-180+66, in=-66, looseness=3.6, "\scalebox{.77}{$\,\mathclap{
      G_x
    }$}"{description}, shift left=1]
    &&
    G/H_i
    \ar[out=-180+66, in=-66, looseness=3.6, "\scalebox{.77}{$\,\mathclap{
      H_i
    }$}"{description}, shift left=1]
  \end{tikzcd}
\end{equation}

\vspace{-2mm}
\noindent
This inherits $H_i$-equivariance from $D_i$, but
we may restrict this (Ex. \ref{RestrictedActions})
to the isotropy group $G_x \subset H_i \subset H$, shown on the right.

Finally, consider the following diagram:
\vspace{-2mm}
\begin{equation}
  \label{DiagramProvingBierstoneLocalTrivialityFromTomDieckLocalTriviality}
  \begin{tikzcd}[row sep=small]
    &
    \mathrm{P}\vert_{ \mathrm{U}_i}
    \ar[
      rrrr
    ]
    \ar[
      dd
    ]
    &&&&
    \widehat G / \widehat {H_i}
    \ar[
      dd
    ]
    \\
    \widehat G_x
      \times_{\widehat G_x}
    \mathrm{U}_x
    \ar[
      dd
    ]
    \ar[
      ur
    ]
    \ar[
      rrrr,
      crossing over
    ]
    &&&&
    \widehat G_i \times_{H_i} D
    \ar[
      ur
    ]
    \\
    &
    \mathrm{U}_i
    \ar[
      rrrr,
      near start, below,
      "f_i"
    ]
    &&&&
    G/H_i
    \\
    \mathrm{U}_x
    \ar[
      ur,
      hook
    ]
    \ar[
      rrrr
    ]
    && &&
    D_i
    \ar[
      ur,
      hook
    ]
    \ar[
      from=uu,
      crossing over
    ]
  \end{tikzcd}
\end{equation}

  \vspace{-2mm}
\noindent
Here the rear face is the given tom Dieck local trivialization and,
as such, is a pullback. The bottom square is the above pullback \eqref{OpenNeighourhoodOfPointIntomDieckCoverFromOpenNeighbouthodOfeInCosetSpace}
which defines the open neighborhood of any given point $x$ in the given cover.
To prove that we have a local trivialization in the sense of Bierstone,
it is sufficient to prove it for this particular open neighborhood $\mathrm{U}_x$,
which means to produce a pullback square as shown on the left of
\eqref{DiagramProvingBierstoneLocalTrivialityFromTomDieckLocalTriviality}.

But the right square of \eqref{DiagramProvingBierstoneLocalTrivialityFromTomDieckLocalTriviality}
is a pullback by \eqref{RestrictionOftomDieckLocalModelToInvariantNeighborhood},
while the front square is similarly recognized as a pullback by Lemma \ref{HomomorphismsOfLocallyTrivialPrincipalBundlesArePullbackSquares},
using that the underlying front left $\Gamma$-principal bundle is trivial,
by Lemma \ref{EquivariantDirectProductBundlesFromHatGxTwistedProductBundles}.
Therefore, the rear, right, and front squares are all pullbacks. Hence
it follows that the left square is indeed a pullback, by the pasting law
(Prop. \ref{PastingLaw}).
\end{proof}

\medskip

This concludes the proofs of the propositions that
make up Thm. \ref{EquivalentNotionsOfEquivariantLocalTriviality}.

\section{Equivariant classifying spaces}
\label{ConstructionOfUniversalEquivariantPrincipalBundles}

\medskip
We give a streamlined account of the
{\it Murayama-Shimakawa construction}
(\cite{MurayamaShimakawa95}\cite{GuillouMayMerling17})
for
would-be equivariant classifying spaces
of
equivariant principal topological bundles
(their actual classifying property is discussed further below in \cref{EquivariantModuliStacks}).
A transparent formulation is obtained by first establishing all relevant properties
for equivariant principal bundles {\it internal} (Ntn. \ref{Internalization})
to equivariant topological {\it groupoids} (discussed in \cref{GActionsOnTopologicalGroupoids}),
where the proofs amount to elementary verifications,
and from where most of the desired statements are functorially induced
\eqref{FunctorOnStructuresInducedFromLexFunctor}
by the fact that topological realization preserves all finite limits
(Lem. \ref{TopologicalRealizationPreservesFiniteLimits}).

\medskip

\noindent
{\bf Ordinary classifying spaces.}
As a transparent template for the following equivariant discussion,
we start by briefly discussing the classical construction of
ordinary (i.e., not equivariant) universal principal bundles
via topological realization of a universal
formally principal bundle
$\mathbf{E}\Gamma \xrightarrow{q_{{}_{\Gamma}}} \mathbf{B}\Gamma$
{\it internal to topological groupoids} (Ex. \ref{UniversalPrincipalGroupoid} below).
The elegance of this construction rests in the fact that
the pertinent properties of the bundle
$\mathbf{E}G \xrightarrow{\;} \mathbf{B}G$
{\it internal to topological groupoids}
are  elementary and immediate to check explicitly
(Lem. \ref{UniversalPrincipalGroupoidIsEquivalentToThePoint},
Lem. \ref{UniversalPrincipalGroupoidIsFormallyPrincipal}),
from where they are then simply inherited
due to the good abstract properties of the topological realization functor
that we have established in \cref{GActionsOnTopologicalGroupoids}.
The resulting topological realization
\vspace{-1mm}
$$
  E \Gamma
  \;\coloneqq\;
  \left\vert
    \mathbf{E}\Gamma
  \right\vert
    \xrightarrow{ \vert q_{{}_{\Gamma}} \vert }
  \left\vert
    \mathbf{B}\Gamma
  \right\vert
  \;=:\;
  B G
$$

\vspace{-1mm}
\noindent is the classical universal $\Gamma$-principal bundle
in its incarnation as the {\it Milgram bar-construction} \cite{Milgram67},
which \cite[\S 6]{MacLane70}
observed to be the realization coend \eqref{TopologicalRealizationOfTopologicalGroupoids}
applied to the simplicial topological space which we
recognize as
the nerve \eqref{SimplicialTopologicalNerveOfTopologicalGroupoids}
of the topological groupoid $\mathbf{E}\Gamma$.

\begin{example}[Universal principal topological groupoid]
  \label{UniversalPrincipalGroupoid}
  Consider any $\Gamma \,\in\, \Groups(\kTopologicalSpaces)$
  with its combined left and right multiplication action:
  $$
    \Gamma
      \acts \;
    \Gamma^{\mathrm{L},\mathrm{R}} \;
    \rightacts
    \Gamma
    \;\;\;
    \in
    \;
    \Actions{
      (
        \Gamma \times \Gamma^{\mathrm{op}}
      )
    }
    (\kTopologicalSpaces)
    \,.
  $$
  We denote the resulting left $\Gamma$-equivariant
  right $\Gamma$-action groupoid (Ex. \ref{RightActionGroupoidInheritsLeftGroupActions})
  by:
  \begin{equation}
    \label{TheUniversalPrincipalTopologicalGroupoid}
    \Gamma
      \acts
      \,
    \mathbf{E}\Gamma
    \;\coloneqq\;
    \Gamma \acts \;
    \Big(
      \Gamma \times \Gamma
    \;\;  \underoverset
        {\scalebox{0.5}{$\mathclap{(-)\cdot(-)}$}}
        {\mathrm{pr}_1}
        {\rightrightarrows} \;\;
      \Gamma
    \Big)
    \;\;\;
    \in
    \;
    \Actions{\Gamma}
    \left(
      \Groupoids(\kTopologicalSpaces)
    \right)
    \,.
  \end{equation}
  Its left $\Gamma$-quotient
  is the delooping groupoid
  of $\Gamma$ (Ex. \ref{TopologicalDeloopingGroupoid}):
  \begin{equation}
    \label{QuotientCoprojectionOfUniversalPrincipalGammaGroupoid}
    \mathbf{E}\Gamma
    \xrightarrow{\;\;\; q \;\;\;}
    \Gamma \,
      \backslash
    (\mathbf{E}\Gamma)
    \;\simeq\;
    \mathbf{B}\Gamma
    \;\;\;
    \in
    \;
    \Groupoids(\kTopologicalSpaces)
    \,.
  \end{equation}

  We discussed below (Rem. \ref{UniversalPrincipalBundleOverTheModuliStack})
  that $\mathbf{E}\Gamma \to \mathbf{B}\Gamma$
  \eqref{TheUniversalPrincipalTopologicalGroupoid} is the
  universal $\Gamma$-principal bundle
  not over the classifying space $B \Gamma$, but over the {\it moduli stack}
  $\mathbf{B}\Gamma$
  (see Prop. \ref{DeloopingGroupoidsOfDTopologicalGroupsAreLocallyFibrant}).
\end{example}

\begin{lemma}[Universal principal groupoid is formally principal]
  \label{UniversalPrincipalGroupoidIsFormallyPrincipal}
  The quotient projection \eqref{QuotientCoprojectionOfUniversalPrincipalGammaGroupoid}
  exhibits the universal $\Gamma$-principal groupoid
  \eqref{TheUniversalPrincipalTopologicalGroupoid}
  as a
  formally principal bundle (Ntn. \ref{InternalizationOfPrincipalBundleTheory})
  internal  to topological groupoids (Ntn. \ref{TopologicalGroupoids}),
  in that the shear map is an isomorphism:
  \vspace{-2mm}
  $$
    \begin{tikzcd}[column sep=large]
      \Gamma \times \mathbf{E}\Gamma
      \ar[
        rr,
        "{
         \scalebox{0.8}{$   \big(
            (-)\cdot(-)
            ,\,
            \mathrm{pr}_2
          \big)
          $}
        }"{above},
        "{\sim}"{below,yshift=+1pt}
      ]
      &&
      (
        \mathbf{E}\Gamma
      )
      \underset{
        \mathbf{B}\Gamma
      }{\times}
      (
        \mathbf{E}\Gamma
      )
    \end{tikzcd}
    \;\;\;\;\;
    \in
    \;
    \Groupoids(\kTopologicalSpaces)
    \mathrlap{\,.}
  $$
\end{lemma}
\begin{proof}
  Since limits of topological groupoids are computed degree-wise,
  this is equivalent to the following two
  horizontal component morphisms being
  isomorphisms of topological spaces (i.e., homeomorphisms):
    \vspace{-2mm}
  \begin{equation}
    \label{ComponentIsomorphismsExhibitingShearMapIsoForUniversalPrincipalGroupoid}
    \begin{tikzcd}[column sep=50pt]
      \Gamma_L \times ( \Gamma \times \Gamma_R )
      \ar[
        rr,
        "{
          (\gamma_\ell, \; \gamma, \; \gamma_r)
          \,\mapsto\,
          ( \gamma_\ell \cdot \gamma, \; \gamma, \; \gamma_r )
        }"
      ]
      \ar[
        d,
        shift left=5pt
      ]
      \ar[
        d,
        shift right=5pt
      ]
      &&
      \Gamma_1 \times \Gamma_2 \times \Gamma_R
      \ar[
        d,
        shift left=5pt
      ]
      \ar[
        d,
        shift right=5pt
      ]
      \\
      \Gamma_L \times \Gamma
      \ar[
        rr,
        "{
          (\gamma_\ell, \; \gamma)
          \,\mapsto\,
          (\gamma_\ell \cdot \gamma, \; \gamma)
        }"
      ]
      &&
      \Gamma_1 \times \Gamma_2
    \end{tikzcd}
    \;\;\;
    \in
    \;
    \kTopologicalSpaces
    \,,
    \;\;\;\;\;\;\;\;\;\;\;\;\;\;
  \end{equation}

    \vspace{-2mm}
\noindent   which is clearly the case.
  (We show subscripts only for readability, all these objects
  $\Gamma_{(-)}$
  are copies of $\Gamma$.)
\end{proof}

\begin{lemma}[Universal principal groupoid is equivalent to the point]
  \label{UniversalPrincipalGroupoidIsEquivalentToThePoint}
  There is an equivalence \eqref{EquivalenceOfTopologicalGroupoids} of topological groupoids
  between the universal $\Gamma$-principal groupoid (Def. \ref{UniversalPrincipalGroupoid})
  and the point $\ast \,:=\, \ConstantGroupoid(\ast)$.
    \vspace{-1mm}
  $$
    \mathbf{E}\Gamma
    \;\underset{\mathrm{htpy}}{\simeq}\;
    \ast
    \;\;\;
    \in
    \;
    \Groupoids(\kTopologicalSpaces)
    \,.
  $$
\end{lemma}
\begin{proof}
  The canonical choice of equivalence
  \eqref{EquivalenceOfTopologicalGroupoids} is
    \vspace{-2mm}
  \begin{equation}
    \label{ContractionOfUniversalPrincipalGroupoid}
    \begin{tikzcd}[column sep=large]
    \big(
      \Gamma
        \times
      \Gamma
   \;\;   \underoverset
        {\mathclap{\scalebox{0.5}{$(-)\cdot (-)$}}}
        {\mathrm{pr}_1}
        {\rightrightarrows} \;\;
      \Gamma
    \big)
    \ar[
      r,
      shift right=4pt,
      "\exists! \; R"{below}
    ]
    &
    (
      \ast
      \rightrightarrows
      \ast
    )
    \,,
    \ar[
      l,
      shift right=4pt,
      "{
        L \,\coloneqq\, \mathrm{e}
      }"{above}
    ]
    \end{tikzcd}
    \;\;\;\;\;
    L \circ R
    \xRightarrow{ (-)^{-1} }
    \mathrm{id}
    ,\,
    \;\;\;\;\;
    \mathrm{id}
    \;=\;
    R \circ L
    \,.
  \end{equation}
  \vspace{-8mm}

\end{proof}

These lemmas readily imply the following classical fact:

\begin{proposition}[Topological realization of universal principal groupoid is universal principal bundle]
  \label{TopologicalRealizationOfUniversalPrincipalGroupoidIsUniversalPrincipalBundle}
  $\,$

  \noindent
  For $\Gamma \,\in\, \Groups(\kTopologicalSpaces)$,
   the topological realization \eqref{TopologicalRealizationOfTopologicalGroupoids}
  \vspace{-2mm}
  \begin{equation}
    \label{UniversalPrincipalBundleAsTopologicalRealizationOfUniversalPrincipalGroupoid}
    E \Gamma
    \;\coloneqq\;
    \left\vert
      \mathbf{E}\Gamma
    \right\vert
    \;\;\;
    \in
    \;
    \kTopologicalSpaces
  \end{equation}

  \vspace{-2mm}
  \noindent
  of the universal $\Gamma$-principal groupoid $\mathbf{E}\Gamma$
  (Ex. \ref{UniversalPrincipalGroupoid})

\noindent {\bf (i)}    is contractible
  \vspace{-1mm}
  \begin{equation}
    \label{UniversalPrincipalBundleIsContractible}
    E \Gamma
    \;\underset{\mathrm{htpy}}{\simeq}\;
    \ast\;;
  \end{equation}

  \vspace{-2mm}
 \noindent {\bf (ii)}   inherits a $\Gamma$-action
  $
    \Gamma \acts \, E \Gamma
    \;\;\;
    \in
    \;
    \Actions{\Gamma}(\kTopologicalSpaces)
  $
  whose quotient coprojection
  coincides
  with the topological realization of
  the groupoid quotient \eqref{QuotientCoprojectionOfUniversalPrincipalGammaGroupoid}

  \vspace{-.5cm}
  \begin{equation}
    \label{QuotientCoprojectionOfUniversalPrincipalBundle}
    E \Gamma
    \xrightarrow{ \;\;\;\;\;\; }
    \Gamma \backslash (E \Gamma)
    \;\simeq\;
    \overset{
      \raisebox{3pt}{
        \tiny
        \color{darkblue}
        \bf
        Milgram classifying space
      }
    }{
      \vert \mathbf{B}\Gamma \vert
      \;=:\;
      B \Gamma\;;
    }
  \end{equation}

  \vspace{-2mm}
\noindent {\bf (iii)}    and is principal, in that the shear map is a homeomorphism:
  \vspace{-2mm}
  $$
    \begin{tikzcd}
      \Gamma \times E \Gamma
      \ar[
        rr,
        "{
         \scalebox{0.8}{$   ( (-)\cdot (-), \; \mathrm{pr}_2 ) $}
        }"{above},
        "{\sim}"{below, yshift=1pt}
      ]
      &&
      E \Gamma
        \underset{B \Gamma}{\times}
      E \Gamma
    \end{tikzcd}
    \;\;\;\;\;\;\;
    \in
    \;
    \kTopologicalSpaces
    \,.
  $$
\end{proposition}
\begin{proof}

  \noindent
  {(i)}
  That $E \Gamma$ is contractible follows
  from Lem. \ref{TopologicalRealizationOfEquivalenceOfGroupoidsIsHomotopyEquivalence}
  applied to Lem. \ref{UniversalPrincipalGroupoidIsEquivalentToThePoint}.

  \noindent
  {(ii)}
  That $E \Gamma \xrightarrow{\;} B \Gamma$ is the $\Gamma$-quotient projection
  follows
  from Lem. \ref{NervePreservesLeftQuotientsOfRightActionGroupoids}
  applied to Ex. \ref{UniversalPrincipalGroupoid},
  which yields that the nerve preserves the quotient projection,
  followed by the fact that the remaining topological realization of simplicial spaces
  \eqref{TopologicalRealizationOfSimplicialTopologicalSpaces} preserves all colimits.

  \noindent
  {(iii)}
  That the $\Gamma$-action on $E \Gamma$ is principal over $B \Gamma$
  follows from the principal bundle structure on $\mathbf{E} \Gamma$
  (from Lem. \ref{UniversalPrincipalGroupoidIsFormallyPrincipal})
  since topological realization preserves finite limits,
  by Lem. \ref{TopologicalRealizationPreservesFiniteLimits},
  thus inducing a functor \eqref{FunctorOnStructuresInducedFromLexFunctor}
  of internal formally principal bundles:
  \vspace{-2mm}
  $$
    \begin{tikzcd}[column sep=large]
      \FormallyPrincipalBundles{\Gamma}
      \left(
        \Groupoids(\kTopologicalSpaces)
      \right)
      \ar[
        rr,
        "{
         \scalebox{0.8}{$   \FormallyPrincipalBundles{(\vert- \vert)}
          (
            \vert- \vert
          )$}
        }"
      ]
      &&
      \FormallyPrincipalBundles{\Gamma}
      (
        \kTopologicalSpaces
      )
      \,.
    \end{tikzcd}
  $$
  \vspace{-11mm}

\end{proof}
\begin{example}[Topological groupoid refinement of the Borel construction]
\label{GroupoidLevelBorelConstruction}
For
$\TopologicalSpace  \, \rightacts  \Gamma
  \;\in\; \Actions{\Gamma^{\mathrm{op}}}(\kTopologicalSpaces)
)$,
its right $\Gamma$-action groupoid
(Ex. \ref{RightActionGroupoidInheritsLeftGroupActions})
is isomorphic to the quotient, in the 1-category $\Groupoids(\kTopologicalSpaces)$,
of the product of
the constant groupoid on $\TopologicalSpace$ (Ex. \ref{TopologicalSpacesAsTopologicalGroupoids}),
regarded with its induced inverse left action,
with the universal principal groupoid on $\Gamma$ (Ex. \ref{UniversalPrincipalGroupoid}),
regarded with its canonical left $\Gamma$-action
\eqref{TheUniversalPrincipalTopologicalGroupoid},
by the resulting product action:
\vspace{-2mm}
\begin{equation}
  (
    \TopologicalSpace \times \Gamma
    \rightrightarrows
    \TopologicalSpace
  )
  \;\;
    \simeq
  \;\;
  \frac{
    \ConstantGroupoid(\TopologicalSpace)
      \times
    \mathbf{E}\Gamma
  }{
    \Gamma
  }
  \;\;\;
  \in
  \;
  \Groupoids(\kTopologicalSpaces)
  \,.
\end{equation}
Since topological realization \eqref{TopologicalRealizationOfTopologicalGroupoids}
preserves all quotients, finite products
(by Lem. \ref{TopologicalRealizationPreservesFiniteLimits})
and constant groupoids (by Ex. \ref{TopologicalRealizationOfConstantGroupoids}),
this means that the topological realization of the right action groupoid on
$\TopologicalSpace$ is its {\it Borel construction}:
\begin{equation}
  \label{BorelConstructionAsTopologicalRealizationOfActionGroupoid}
  \left\vert
    \TopologicalSpace \times \Gamma
    \rightrightarrows
    \TopologicalSpace
  \right\vert
  \;\;
    \simeq
  \;\;
  \frac{
    \vert\ConstantGroupoid(\TopologicalSpace)\vert
      \times
    \vert\mathbf{E}\Gamma\vert
  }{
    \left\vert \Gamma \right\vert
  }
  \;\;
    \simeq
  \;\;
  \frac{
    \TopologicalSpace \times E \Gamma
  }{
    \Gamma
  }
  \;\;\;
  \in
  \;
  \kTopologicalSpaces
  \,,
\end{equation}
whose path-connected components are those of the plain quotient,
by \eqref{PathConnectedComponentsPreserveQuotientsAndFiniteProducts} and
\eqref{UniversalPrincipalBundleIsContractible}:
\begin{equation}
  \label{PathConnectedComponentsOfBorelConstruction}
  \pi_0
  \left(
    \frac{
      \TopologicalSpace
        \times
      E \Gamma
    }{
      \Gamma
    }
  \right)
  \;\;
  \simeq
  \;\;
  \pi_0
  \left(
    \frac{
      \TopologicalSpace
    }{
      \Gamma
    }
  \right)
  \;\;\;
  \in
  \;
  \Sets
\end{equation}
\end{example}

\begin{example}[Higher classifying spaces]
  \label{HigherClassifyingSpaces}
  For $A \,\in\, \AbelianGroups(\kTopologicalSpaces)$,
  its classifying space \eqref{QuotientCoprojectionOfUniversalPrincipalBundle}
  itself carries the structure of a topological group, by
  Ex. \ref{Strict2GroupsDeloopingAbelianGroups}.
  Therefore the
  classifying space construction \eqref{QuotientCoprojectionOfUniversalPrincipalBundle}
  may be iterated, and we write, recursively:
  $$
    B^0 A \,\coloneqq\, A
  $$
  and
  \begin{equation}
    \label{HigherClassifyingSpacesOfAbelianTopologicalGroup}
    B^{n+1} A
    \;\coloneqq\;
    B
    (
      B^n A
    )
    \;=\;
    B \,
    \vert
      (B^n A) \rightrightarrows \ast
    \vert
    \;\;\;\;
    \in
    \;
    \kTopologicalSpaces
    \,.
  \end{equation}
\end{example}

\medskip

\noindent
{\bf The Murayama-Shimakawa construction.}
In equivariant generalization of Ex. \ref{UniversalPrincipalGroupoid}, we make
the following Def. \ref{EquivariantUniversalPrincipalGroupoid},
which, upon its topological realization
(Prop. \ref{TopologicalRealizationOfEquivariantPrincipalGroupoid} below),
is a streamlined rendering of
the construction due to \cite[p. 1293 (5 of 7)]{MurayamaShimakawa95},
specialized to discrete equivariance groups $G$ according to the
remark at the bottom of \cite[p. 1294 (6 of 7)]{MurayamaShimakawa95}
(which seems to be the greatest validated generality,
as observed in \cite[Scholium 3.12]{GuillouMayMerling17}):

\begin{definition}[Universal equivariant principal groupoid via Murayama-Shimakawa construction]
  \label{EquivariantUniversalPrincipalGroupoid}
  $\,$

  \noindent
  {\bf (i)}
  Given a discrete group
  $G \,\in\, \Groups(\Sets) \xhookrightarrow{\;} \Groups(\kTopologicalSpaces)$,
  consider any
  $G \acts \, \Gamma \,\in\, \Groups( \Actions{G}(\kTopologicalSpaces) )$
  \eqref{CategoryOfGEquivariantTopologicalGroups}
  with its combined left and right multiplication action
      \vspace{-2mm}
  $$
    (G \acts \,  \Gamma)
    \,\acts\;
    \left(
      G \acts \, \Gamma^{L,R}
    \right)
    \; \raisebox{-2pt}{$\rightacts$} \,
    (G \acts \, \Gamma)
    \;\;\;
    \in
    \;
    \Actions{
    \big(
      (G \acts \, \Gamma)
      \times
      (G \acts \, \Gamma^{\mathrm{op}})
    \big)
    }
    \left(
      \Actions{G}(\kTopologicalSpaces)
    \right)
    .
  $$

      \vspace{-2mm}
 \noindent
  By Ex. \ref{RightEquivariantActionGroupoidInheritsLeftEquivariantGroupAction},
  we obtain the resulting
  left $(G \acts \, \Gamma_L)$-equivariant
  right $(G \acts \,  \Gamma_R)$-action
  $G$-groupoid
      \vspace{-2mm}
  \begin{equation}
    \label{EquivariantPreUniversalPrincipalGroupoid}
    \hspace{-4mm}
    (G \acts \, \Gamma)
    \,\acts\,
    \mathbf{E}(G \acts \, \Gamma)
    \, \coloneqq
    (G \acts \, \Gamma)
    \,\acts\,
    \Big(\!\!
      (G \acts \, \Gamma)
      \times
      (G \acts \, \Gamma)
      \;\underoverset
        {\scalebox{0.5}{$\mathclap{(-)\cdot (-)}$}}
        {\mathrm{pr}_1}
        {\rightrightarrows} \;
      (G \acts \, \Gamma)
\!\!    \Big)
    \;
    \in
    \;
    \Actions{
      (G \acts \, \Gamma)
    }
    \left(\!
      \Groupoids
      \left(
        \Actions{G}(\kTopologicalSpaces)
      \right)
   \! \right)
    .
  \end{equation}

 \vspace{-3mm}
 \noindent
 Notice that after forgetting the $(G \acts \Gamma)$-action on
 \eqref{EquivariantPreUniversalPrincipalGroupoid} there still remains the
 $G$-action inherited from the $G$-action on $\Gamma$:
     \vspace{-2mm}
\begin{equation}
  \label{GEquivariantUniversalGammaPrincipalGroupoid}
  G \acts \, \mathbf{E}\Gamma
  \;=\;
  G \acts \;
  \big(
    \Gamma \times \Gamma
  \;\;  \underoverset
      {\scalebox{0.5}{$\mathclap{(-)\cdot (-)}$}}
      {\mathrm{pr}_1}
      {\rightrightarrows} \;\;
    \Gamma
  \big)
  \;\;\;
  \in
  \;
  \Groupoids
  \left(
   \Actions{G}(\kTopologicalSpaces)
  \right).
\end{equation}

    \vspace{-2mm}
 \noindent
 The left $(G \acts \, \Gamma)$-quotient of \eqref{EquivariantPreUniversalPrincipalGroupoid}
 is the equivariant delooping groupoid
 (Ex. \ref{EquivariantTopologicalDeloopingGroupoid}):
     \vspace{-2mm}
 \begin{equation}
   \label{EquivariantCoprojectionOutOfEGamma}
   G \acts \,
   \mathbf{E}\Gamma
   \xrightarrow{
     \;
     \scalebox{.7}{$G \acts \, q$}
     \;
   }
   (G \acts \, \Gamma)
     \backslash
   \mathbf{E}\Gamma
   \;\simeq\;
   G \acts \,  \mathbf{B}\Gamma
   \;\;\;
   \in
   \;
   \Actions{G}
   \left(
     \Groupoids(\kTopologicalSpaces)
   \right).
 \end{equation}

    \vspace{-2mm}
\noindent
{\bf (ii)}
We say that
the {\it universal equivariant principal groupoid} over $(G \acts \, \Gamma)$
is the
$G$-equivariant mapping groupoid
(Ex. \ref{ConjugationActionOnMappingGroupoids})
of $G \acts \, \mathbf{E}G$ \eqref{TheUniversalPrincipalTopologicalGroupoid}
into $G \acts \, \mathbf{E}\Gamma$ \eqref{GEquivariantUniversalGammaPrincipalGroupoid}
with its $(G \acts \, \Gamma)$-action inherited from \eqref{EquivariantPreUniversalPrincipalGroupoid}:
\vspace{-2mm}
 \begin{equation}
  \label{UniversalEquivariantPrincipalGroupoid}
  (G \acts \, \Gamma)
  \,\acts\;
  \mathrm{Maps}
  (
    \mathbf{E}G
    ,\,
    \mathbf{E}\Gamma
  )
  \;\;\;
  \in
  \;
  \Actions{
    (G \acts \, \Gamma)
  }
  \left(
    \Groupoids
    \left(
      \Actions{G}(\kTopologicalSpaces)
    \right)
  \right)
  .
\end{equation}
\end{definition}

\begin{remark}[Conceptual nature of universal equivariant principal groupoid]
  It may be surprising on first sight --
  and is the key innovation of \cite{MurayamaShimakawa95} --
  that Def. \ref{EquivariantUniversalPrincipalGroupoid}
  involves not just the evident $G$-groupoid $G \acts \, \mathbf{E}\Gamma$
  \eqref{GEquivariantUniversalGammaPrincipalGroupoid},
  but the mapping groupoid $G \acts \, \mathrm{Maps}(\mathbf{E}G, \, \mathbf{E}\Gamma)$
  \eqref{UniversalEquivariantPrincipalGroupoid}.
  Of course, one recognizes the fixed locus
  $\mathrm{Maps}(E G, \, \TopologicalSpace)^G \,\simeq\, \TopologicalSpace^{\mathrm{h}G}$
  as a model for the {\it homotopy fixed locus}
  (going back to \cite{Thomason83}, see also \cite{Virk21}),
  but this is more of a hint than an explanation.
  Instead, the construction had been justified
  (in \cite{MurayamaShimakawa95}\cite{GuillouMayMerling17})
  {\it a posteriori}
  by proof that its topological realization
  does happen to satisfy (under suitable conditions on $G$ and $\Gamma$)
  abstract properties known to characterize
  universal equivariant principal bundles.
  Below in Thm. \ref{MurayamaShimakawaGroupoidIsEquivariantModuliStack}
  we give a more conceptual derivation of
  Def. \ref{EquivariantUniversalPrincipalGroupoid},
  showing that this construction arises naturally
  as the  {\it equivariant shape} \eqref{EquivariantShape}
  of the proper-equivariant {\it moduli stack} of equivariant principal bundles.
\end{remark}

Using our internalization technology, we prove in
Prop. \ref{TopologicalRealizationOfEquivariantPrincipalGroupoid}
that the topological realization of Def. \ref{EquivariantUniversalPrincipalGroupoid}
is an equivariantly contractible equivariant principal bundle,
as a formal consequence of the following sequence of Lemmas,
which establish elementary properties of the underlying
equivariant principal groupoid.

\begin{lemma}[Universal equivariant principal groupoid is right action groupoid]
  \label{MappingGroupoidIntoCodiscreteGroupoidIsRightActionGroupoid}
The universal equivariant principal groupoid \eqref{UniversalEquivariantPrincipalGroupoid}
is a left $(G \acts \, \Gamma)$-equivariant right action $(G \acts \, \Gamma)$-groupoid
of the form of Ex. \ref{RightActionGroupoidInheritsLeftGroupActions}:
    \vspace{-2mm}
\begin{equation}
  (G \acts \, \Gamma)
  \,\acts\;
  \mathrm{Maps}
  (\mathbf{E}G,\,\mathbf{E}\Gamma)
  \;\;\simeq\;\;
  \Big(
    \mathrm{Fnctr}(\mathbf{E}G,\,\mathbf{E}\Gamma)^{L,R}
    \times
    \mathrm{Maps}(G,\,\Gamma)
    \;\;\;
    \underoverset
      {\scalebox{0.5}{$\mathclap{(-)\cdot(-)}$}}
      {\mathrm{pr}_1}
      {\rightrightarrows}
      \;\;\;
    \mathrm{Fnctr}(\mathbf{E}G,\,\mathbf{E}\Gamma)
 \! \Big)
  \,,
\end{equation}

    \vspace{-3mm}
 \noindent
where
    \vspace{-2mm}
$$
  \mathrm{Fnctr}(\mathbf{E}G,\,\mathbf{E}\Gamma)
  \;\coloneqq\;
  \mathrm{Maps}(\mathbf{E}G,\,\mathbf{E}\Gamma)_0
$$
is equipped with the argument-wise induced left and right $(G \acts \, \Gamma)$-actions.
\end{lemma}
\begin{proof}
  This is a direct unwinding of the definitions:
  The point to observe is only that the space of
  natural transformations
  out of any fixed functor $\mathbf{E}G \xrightarrow{\;} \mathbf{E}\Gamma$
  is already isomorphic to that of all possible component functions,
  which is
  \vspace{-2mm}
  $$
    \mathrm{Maps}
    \left(
      (\mathbf{E}G)_0
      ,\,
      (\mathbf{E}\Gamma)_1
    \right)
    \;\simeq\;
    \mathrm{Maps}(G,\Gamma)
    \,.
  $$

  \vspace{-8mm}
\end{proof}

\begin{lemma}[Mapping groupoid into universal principal groupoid preserves canonical quotient]
  \label{MappingGroupoidIntoUniversalPrincipalGroupoidPreservesCanonicalQuotient}
  $\,$

  \noindent
  Let $\Gamma \,\in\, \Groups(\kTopologicalSpaces)$
  and let
  $( \TopologicalSpace_1 \rightrightarrows \TopologicalSpace_0)
    \,\in\,
    \Groupoids(\kTopologicalSpaces)$
    be topologically discrete and connected as a groupoid.
  Then pushforward $q_\ast$ along the quotient coprojection
  $\mathbf{E}\Gamma \xrightarrow{q} \Gamma \backslash \mathbf{E}\Gamma = \mathbf{B}\Gamma$
  \eqref{QuotientCoprojectionOfUniversalPrincipalGammaGroupoid}
  of the mapping groupoid
  \eqref{InternalHomAdjunctionForTopologicalGroupoids}
  into $\Gamma \acts \, \mathbf{E}\Gamma$ \eqref{TheUniversalPrincipalTopologicalGroupoid}
  exhibits the quotient of the latter
  by its induced pointwise left $\Gamma$-action:
      \vspace{-3mm}
  $$
    \begin{tikzcd}[row sep=small]
    \mathrm{Maps}
    \left(
      (\TopologicalSpace_1 \rightrightarrows \TopologicalSpace_0)
      ,\,
      \mathbf{E}\Gamma
    \right)
    \ar[
      d
    ]
    \ar[
      drr,
      "q_\ast"
    ]
    \\
    \Gamma
      \,\backslash\,
    \mathrm{Maps}
    \left(
      (\TopologicalSpace_1 \rightrightarrows \TopologicalSpace_0)
      ,\,
      \mathbf{E}\Gamma
    \right)
    \ar[
      rr,
      "{\sim}"{yshift=-1pt}
    ]
    &&
    \mathrm{Maps}
    \left(
      (\TopologicalSpace_1 \rightrightarrows \TopologicalSpace_0)
      ,\,
      \mathbf{B}\Gamma
    \right).
    \end{tikzcd}
  $$

\end{lemma}
\noindent
This statement corresponds to the last part of \cite[Thm. 2.7]{GuillouMayMerling17}.
We mean to give a transparent proof:
\begin{proof}
  By the assumption that the domain groupoid is connected
  and discrete,
  we may choose a base object $x_0$ and
  a continuous assignment $f_{(-)}$ of
  morphisms connecting all objects to this one:
      \vspace{-3mm}
  $$
    x_0 \,\in\, \TopologicalSpace_0
    \,,
    {\phantom{AAAAAAAA}}
    \begin{tikzcd}[row sep=-6pt]
      \mathllap{f_{(-)} \;:\; }
      \TopologicalSpace_0
      \ar[
        r
      ]
      &
      \TopologicalSpace_1
      \\
  \scalebox{0.8}{$      x $}
      \ar[r, |->]
      &
 \scalebox{0.8}{$       (x_0 \xrightarrow{f_x} x) $}
    \end{tikzcd}
     $$
  \vspace{-.5cm}

  \noindent
  Using this, we obtain the following
  component functions,
  where $q \,:\, \mathbf{E}G \xrightarrow{\;} \mathbf{B}G$
  is the coprojection \eqref{QuotientCoprojectionOfUniversalPrincipalGammaGroupoid}
  from Ex. \ref{UniversalPrincipalGroupoid}:

  \vspace{-.6cm}
  \begin{equation}
    \label{DecomposingAFunctorToEG}
    \begin{tikzcd}[row sep=-5pt]
 \scalebox{0.7}{$      F $}
      &\longmapsto&
      \!\!\!\!\!\!\!\!\!\!\!\!\!\!
\scalebox{0.7}{$      \left(
        F_0(x_0)
        ,\,
        q \circ F
      \right)
      $}
      \\
      \mathrm{Maps}
      \left(
        (\TopologicalSpace_1 \rightrightarrows \TopologicalSpace_0)
        ,\,
        \mathbf{E}\Gamma
      \right)
      \ar[
        rr,
        shift left=6pt
      ]
      &&
      (
        \Gamma \rightrightarrows \Gamma
      )
        \times
      \mathrm{Maps}
      \left(
        (\TopologicalSpace_1 \rightrightarrows \TopologicalSpace_0)
        ,\,
        \mathbf{B} \Gamma
      \right)\,.
      \ar[
        ll
      ]
      \\
  \scalebox{0.7}{$     \Big(
        \big(
          x_1 \xrightarrow{f} x_2
        \big)
        \,\mapsto\,
        \big(
          \gamma
            \cdot
          \phi(f_{x_1})
          \xrightarrow{\;}
          \gamma
            \cdot
          \phi(f_{x_2})
        \big)
      \Big)
      $}
      &\raisebox{6pt}{\rotatebox{180}{$\longmapsto$}}&
     \scalebox{0.7}{$  (
        \gamma
        ,\
        \phi
      )
      $}
      \mathrlap{\,.}
    \end{tikzcd}
  \end{equation}

  \vspace{-2mm}
  \noindent
  Unwinding the definitions readily reveals that
  this is a pair of inverse continuous funcors. Moreover, the
  top morphism is evidently $\Gamma$-equivariant, with respect to the
  canonical left $\Gamma$-action on the $(\Gamma \rightrightarrows \Gamma)$-factor
  and the trivial action on the remaining factor on the right.
  Therefore, passing the quotient operation by $\Gamma$ along the top isomorphism
  implies the claim.
\end{proof}

\begin{lemma}[Quotient coprojection of universal equivariant principal groupoid]
  \label{QuotientCoprojectionOfUniversalEquivariantPrincipalGroupoid}
  Let $G \,\in\, \Groups(\Sets) \xhookrightarrow{\;} \Groups(\kTopologicalSpaces)$
  be a discrete group
  and $(G \acts \, \Gamma) \,\in\, \Groups( \Actions{G}(\kTopologicalSpaces) )$,
  then pushforward along the
  $G$-equivariant quotient coprojection
  $G \acts \, \mathbf{E} \Gamma \xrightarrow{\scalebox{.7}{$G \acts \, q$}} G \acts \, \mathbf{B}\Gamma$
  \eqref{EquivariantCoprojectionOutOfEGamma}
  of the equivariant mapping groupoid \eqref{InternalHomInGEquivariantTopologicalGroupoids}
  out of $G \acts \, \mathbf{E}G$ \eqref{TheUniversalPrincipalTopologicalGroupoid}
  exhibits the equivariant quotient of the latter
  by its induced left $\Gamma$-action:
      \vspace{-2mm}
  \begin{equation}
    \label{TheQuotientCoprojectionOnUniversalEquivariantPrincipalgroupoid}
    \begin{tikzcd}[row sep=small]
    \mbox{
      \tiny
      \color{darkblue}
      \bf
      \begin{tabular}{c}
        universal equivariant
        \\
        principal groupoid
      \end{tabular}
    }
    &
    G
      \acts \;
    \mathrm{Maps}
    (
      \mathbf{E}G
      ,\,
      \mathbf{E}\Gamma
    )
    \ar[
      d
    ]
    \ar[
      drr,
      "\scalebox{.7}{$G \acts \, q_\ast$}"
    ]
    &&
    \\
    \mbox{
      \tiny
      \color{darkblue}
      \bf
      \begin{tabular}{c}
        universal equivariant
        \\
        {\color{orangeii} base }
        groupoid
      \end{tabular}
    }
    &
    (G \acts \, \Gamma)
      \,\backslash\,
    \mathrm{Maps}
    (
      \mathbf{E}G
      ,\,
      \mathbf{E}\Gamma
    )
    \ar[
      rr,
      "{\sim}"{yshift=-1pt}
    ]
    &&
    G \acts \;
    \mathrm{Maps}
    (
      \mathbf{E}G
      ,\,
      \mathbf{B}\Gamma
    )
    \,.
    \end{tikzcd}
  \end{equation}
\end{lemma}
\begin{proof}
  Since the functor which forgets the $G$-action is a left adjoint
  (Ex. \ref{ForgettingGActionsAsPullbackAction}) and hence preserves
  colimits such as the quotient considered here,
  and since $\mathbf{E}G$ is connected (by definition)
  and topologically discrete (by assumption on $G$),
  the statement follows with Lem. \ref{MappingGroupoidIntoUniversalPrincipalGroupoidPreservesCanonicalQuotient}.
\end{proof}

In equivariant generalization of Lem. \ref{UniversalPrincipalGroupoidIsEquivalentToThePoint},
we have:
\begin{lemma}[Universal equivariant principal groupoid is equivalent to the point]
  There is an equivalence of $G$-equivariant topological groupoids
  between the universal equivariant
  principal groupoid \eqref{UniversalEquivariantPrincipalGroupoid}
  from
  (Def. \ref{EquivariantUniversalPrincipalGroupoid})
  and the point:
  \begin{equation}
    \label{EquivariantContractionOfUniversalEquivariantPrincipalGroupoid}
    G \acts \;
    \mathrm{Maps}
    (
      \mathbf{E}G
      ,\,
      \mathbf{E}\Gamma
    )
    \;\;
    \underset{\mathrm{htpy}}{\simeq}
    \;\;
    G \acts \; \ast
    \;\;\;\;
    \in
    \;
    \Actions{G}
    (
      \Groupoids(\kTopologicalSpaces)
    )
    \,.
  \end{equation}
  In particular, for all subgroups $H \subset G$ the $H$-fixed subgroupoid
  of the mapping groupoid is equivalent to the point:
  \vspace{-2mm}
  $$
    \underset{H \subset G}{\forall}
    \;\;\;\;
    \left(
    G
      \acts \;
    \mathrm{Maps}
    (
      \mathbf{E}G
      ,\,
      \mathbf{E}\Gamma
    )
    \right)^H
    \;
     \underset{\mathrm{htpy}}{\simeq}
    \;
    \ast
    \;\;\;
    \in
    \;
    \Actions{G}
    \left(
      \Groupoids
      (\kTopologicalSpaces)
    \right)
    \,.
  $$
\end{lemma}
\begin{proof}
  Since the $G$-action on $\Gamma$ is by group automorphisms
  (Lem. \ref{EquivariantTopologicalGroupsAreSemidirectProductsWithG})
  the equivalence \eqref{ContractionOfUniversalPrincipalGroupoid}
  contracting the universal principal groupoid $\mathbf{E}\Gamma$ is clearly
  $G$-equivariant, so that
  it is also contractible as a
  $G$-groupoid \eqref{GEquivariantUniversalGammaPrincipalGroupoid}
  \vspace{-1mm}
  \begin{equation}
    \label{EquivariantContractionOfUniversalPrincipalGroupoid}
    G \acts  \, \mathbf{E}\Gamma
    \;\underset{\mathrm{htpy}}{\simeq}\;
    \ast
    \;\;\;
    \in
    \;
    \Groupoids
    (
      \Actions{G}(\kTopologicalSpaces)
    )
    \,.
  \end{equation}

  \vspace{-2mm}
\noindent  The required equivalence \eqref{EquivariantContractionOfUniversalEquivariantPrincipalGroupoid}
  is obtained as the
  image of this equivalence \eqref{EquivariantContractionOfUniversalPrincipalGroupoid}
  under the equivariant mapping groupoid 2-functor
  $G \acts \, \mathrm{Maps}(\mathbf{E}G , -)$ \eqref{InternalHomInGEquivariantTopologicalGroupoids}.
\end{proof}

\begin{lemma}[Equivariant mapping groupoid out of universal principal groupoid preserves
 constant groupoids]
  \label{EquivariantMappingGroupoidOutOfUniversalPrincipalGroupoidPreservesConstantGroupoids}
  For $(G \acts \; \Gamma) \,\in\, \Groups( \Actions{G}(\kTopologicalSpaces))$,
  we have an isomorphism of group objects
  \vspace{-2mm}
  $$
    \begin{tikzcd}
      G
        \acts \;
      \mathrm{Maps}
      (
        \mathbf{E}G
        ,\,
        \ConstantGroupoid(\Gamma)
      )
      \ar[
        rr,
        "{\mathrm{ev}_{\mathrm{e}}}"{above},
        "{\sim}"{below,yshift=.5pt}
      ]
      &&
      G \acts \, \Gamma
    \end{tikzcd}
    \;\;\;
    \in
    \;
    \Groups
    \left(
      \Actions{G}
      (
        \Groupoids(\kTopologicalSpaces)
      )
    \right)
  $$

 \vspace{-2mm}
 \noindent
  between $\Gamma$ and the equivariant mapping groupoid
  (Ex. \ref{ConjugationActionOnMappingGroupoids})
  out of $G \acts \, \mathbf{E}G$ \eqref{TheUniversalPrincipalTopologicalGroupoid}
  into the constant groupoid (Ex. \ref{TopologicalSpacesAsTopologicalGroupoids})
  on $\Gamma$.
\end{lemma}
\begin{proof}
$\,$
  \vspace{-2mm}
  $$
    \begin{array}{lll}
      {\rm Maps}
      \left(  { \mathbf{E}G },
        { \ConstantGroupoid(\Gamma) }
        \right)
      & \;\simeq\;
      {\rm Maps}
      \left(
        { \Truncation{0} (\mathbf{E}G) },
        { \Gamma }
      \right)
      &
      \proofstep{ by Prop. \ref{AdjunctionBetweenTopologicalGroupoidsAndTopologicalSpaces} }
      \\
      & \;\simeq\;
      \Maps{}
        { \ast }
        { \Gamma }
      \\
      & \;\simeq\;
      \Gamma
      \mathrlap{\,.}
    \end{array}
  $$

  \vspace{-6mm}
\end{proof}

In equivariant generalization of Lem. \ref{UniversalPrincipalGroupoidIsFormallyPrincipal},
we have:
\begin{lemma}[Universal equivariant principal groupoid is formally principal]
  \label{UniversalEquivariantPrincipalGroupoidIsFormallyPrincipal}
  The quotient coprojection
  \eqref{TheQuotientCoprojectionOnUniversalEquivariantPrincipalgroupoid}
  exhibits the universal equivariant principal groupoid
  \eqref{UniversalEquivariantPrincipalGroupoid}
  as a formally principal bundle (Ntn. \ref{InternalizationOfPrincipalBundleTheory})
  internal to $G$-equivariant topological groupoids
  (Def. \ref{EquivariantTopologicalGroupoids}),
  \vspace{-2mm}
  \begin{equation}
    \label{EquivariantPreUniversalPrincipalGroupoidIsFormallyPrincipal}
    (G \acts  \, \Gamma)
      \acts \;
    \mathrm{Maps}
    (
      \mathbf{E}G
     \,,
     \mathbf{E}\Gamma
    )
    \xrightarrow{ \;q_\ast\; }
    G
      \acts \;
    \mathrm{Maps}
    (
      \mathbf{E}G
     \,,
     \mathbf{B}\Gamma
    )
    \;\;\;
    \in
    \;
    \FormallyPrincipalBundles{(G \acts  \, \Gamma)}
    \Big(
      \Actions{G}
      \big(
        \Groupoids(\kTopologicalSpaces)
      \big)
     \!
   \Big)
    \,,
  \end{equation}

  \vspace{-2mm}
  \noindent
  in that the
  shear map is an isomorphism:
  \vspace{-2mm}
  $$
    \begin{tikzcd}
    (G \acts \, \Gamma)
    \times
    \big(
      G \acts \;
      \mathrm{Maps}(\mathbf{E}G,\,\mathbf{E}\Gamma)
    \big)
    \ar[
      rr,
      "{
    \scalebox{0.7}{$    \big(
          (-)\cdot(-)
          ,\,
          \mathrm{pr}_2
        \big)
        $}
      }"{above},
      "{\sim}"{below, yshift=1pt}
    ]
    &&
    \big(
      G \acts \;
      \mathrm{Maps}(\mathbf{E}G,\,\mathbf{E}\Gamma)
    \big)
      \!\!
      \underset{
        \scalebox{.6}{$
          G \acts \;
          \mathrm{Maps}(\mathbf{E}G,\,\mathbf{B}\Gamma)
        $}
      }{\times}
      \!\!
    \big(
      G \acts \;
      \mathrm{Maps}(\mathbf{E}G,\,\mathbf{E}\Gamma)
    \big)
    \,.
    \end{tikzcd}
  $$
\end{lemma}
\begin{proof}
  First, observe that
  the coprojection
  \eqref{EquivariantCoprojectionOutOfEGamma}
  is formally $(G \acts \, \Gamma)$-principal:
  \vspace{-2mm}
  \begin{equation}
    \label{PrecursorEquivariantPreUniversalPrincipalGroupoidIsFormallyPrincipal}
    (G \acts \, \Gamma) \acts \, \mathbf{E}\Gamma \xrightarrow{\;\;q\;\;} G \acts \, \mathbf{B}\Gamma
    \;\;\;
    \in
    \;
    \FormallyPrincipalBundles{(G \acts \, \Gamma)}
    \Big(
      \Actions{G}
      \big(
        \Groupoids(\kTopologicalSpaces)
      \big)
     \!
   \Big)
    \,.
  \end{equation}

  \vspace{-2mm}
  \noindent
  This follows by Lem. \ref{UniversalPrincipalGroupoidIsFormallyPrincipal}
  combined with Lem. \ref{ForgetfulFunctorFromTopologicalGSpacesToGSpaces},
  observing that the
  isomorphisms \eqref{ComponentIsomorphismsExhibitingShearMapIsoForUniversalPrincipalGroupoid}
  are clearly $G$-equivariant (Lem. \ref{EquivariantTopologicalGroupsAreSemidirectProductsWithG}).

  Now, since the equivariant mapping groupoid functor
  $G\acts  \; \mathrm{Maps}(\mathbf{E}G,\, -)$
  is a
  right adjoint \eqref{InternalHomInGEquivariantTopologicalGroupoids}
  it preserves finite limits (Prop. \ref{RightAdjointFunctorsPreserveFiberProducts})
  and hence induces \eqref{FunctorOnStructuresInducedFromLexFunctor}
  a functor of formally principal bundles,
  which by Lem. \ref{EquivariantMappingGroupoidOutOfUniversalPrincipalGroupoidPreservesConstantGroupoids}
  is an endo-functor:
  $$
    \FormallyPrincipalBundles{(G \acts \, \Gamma)}
    \Big(
      \Groupoids
      \big(
        \Actions{G}(\kTopologicalSpaces)
      \big)
    \!\!
    \Big)
    \xrightarrow{
      \Maps{}{ \mathbf{E}G }{-}
    }
    \FormallyPrincipalBundles{(G \acts \, \Gamma)}
    \Big(
      \Groupoids
      \big(
        \Actions{G}(\kTopologicalSpaces)
      \big)
    \!\!
    \Big)
    \mathrlap{\,.}
  $$
  The image of
  \eqref{PrecursorEquivariantPreUniversalPrincipalGroupoidIsFormallyPrincipal}
  under this functor gives the required structure
  \eqref{EquivariantPreUniversalPrincipalGroupoidIsFormallyPrincipal}.
\end{proof}

After this series of lemmas, we may finally conclude,
in equivariant generalization of Prop. \ref{TopologicalRealizationOfUniversalPrincipalGroupoidIsUniversalPrincipalBundle}:

\begin{notation}[Murayama-Shimakawa construction]
  \label{MurayamaShimakawaConstruction}
  $\,$

  \noindent
  For
  $G \,\in\, \Groups(\Sets) \xhookrightarrow{\;} \Groups(\kTopologicalSpaces)$
  and
  $\Gamma \,\in\, \Groups( \Actions{G}(\kTopologicalSpaces))$,
  we denote the topological realization \eqref{TopologicalRealizationOfTopologicalGroupoids}
  with its induced $G$-action (Prop. \ref{TopologicalRealizationOfTopologicalGroupoidsRespectsEquivariance})
  of the $G$-equivariant universal $(G \acts \, \Gamma)$-principal groupoid
  $\mathrm{Maps}(\mathbf{E}G, \mathbf{E}\Gamma)$
  (Def. \ref{EquivariantUniversalPrincipalGroupoid})
  by
  \vspace{-2mm}
  \begin{equation}
    \label{DefinitionOfUniversalEquivariantPrincipalBundle}
    G \acts \,
    E(G \acts \Gamma)
    \;\coloneqq\;
    G \acts \;
    \left\vert
      \Maps{}
        { \mathbf{E}G }
        { \mathbf{E}\Gamma }
    \right\vert
    \;\;\;
    \in
    \;
    \Actions{G}( \kTopologicalSpaces )\;.
  \end{equation}
\end{notation}
\begin{proposition}[Murayama-Shimakawa construction as equivariant principal bundle]
  \label{TopologicalRealizationOfEquivariantPrincipalGroupoid}
  $\,$

  \noindent
  Let
  $G \,\in\, \Groups(\Sets) \xhookrightarrow{\;} \Groups(\kTopologicalSpaces)$
  and
  $\Gamma \,\in\, \Groups( \Actions{G}(\kTopologicalSpaces))$.

  \vspace{-3mm}
\begin{enumerate}[{\bf (i)}]
\setlength\itemsep{-7pt}
\item
  Then the Murayama-Shimakawa construction
  \eqref{DefinitionOfUniversalEquivariantPrincipalBundle}
  is
  $G$-equivariantly contractible
  \vspace{-2mm}
  \begin{equation}
    \label{EquivariantContractibilityOfUniversalEquivariantBundle}
    G
      \acts \,
    E(G \acts \, \Gamma)
    \;\;
      \underset{\mathrm{htpy}}{\simeq}
    \;\;
    G
      \acts
    \ast
    \mathrlap{\,.}
  \end{equation}

  \item
  Hence, in particular, it has contractible fixed loci
  \vspace{-2mm}
  \begin{equation}
    \label{ContractibleFixedLociOfEquivariantUniversalBundle}
    \underset{H \subset G}{\forall}
    \;\;
    (
      E(G \acts \, \Gamma)
    )^H
    \;\;
      \underset{\mathrm{htpy}}{\simeq}
    \;\;
    \ast
    \mathrlap{\,,}
  \end{equation}

  \vspace{-6mm}
  and inherits a $(G \acts \, \Gamma)$-action
    \vspace{-4mm}
  $$
    (G \acts \, \Gamma)
    \,\acts\;
    E(G \acts \, \Gamma)
    \;\;\;
    \in
    \;
    \Actions{G}(\kTopologicalSpaces) \;.
  $$

\vspace{-2mm}
  \item
  The corresponding quotient coprojection
  \vspace{-2mm}
  \begin{equation}
    \label{QuotientCoprojectionOfUniversalEquivariantPrincipalBundle}
    E(G \acts \, \Gamma)
    \longrightarrow
    \Gamma
      \backslash
    E(G \acts \, \Gamma)
    \;=:\;
    B(G \acts \, \Gamma)
  \end{equation}

  \vspace{-3mm}
  \noindent
  is a $G$-equivariant $\Gamma$-principal bundle (Def. \ref{EquivariantPrincipalBundle})
      \vspace{-4mm}
  $$
    \begin{tikzcd}[row sep=0pt]
      \Maps{\big}
      {
        \mathbf{E}G
       }
       {
        \mathbf{E}\Gamma
       }
      \ar[out=180-82.8+90, in=82.8+90, looseness=2.6,
        shift right=2pt,
        "\scalebox{.77}{$\mathclap{
          G
        }$}"{description, pos=.4},
      ]
      \ar[
        dd,
        "{
          q_\ast
        }"{right}
      ]
      &&
      \overset{
        \;\;\;
        \scalebox{.7}{$
          E(G \acts \,  \Gamma)
        $}
      }{
      \overbrace{
        \TopologicalRealization{\big}
        {
          \Maps{}
          {
            \mathbf{E}G
          }
          {
            \mathbf{E}\Gamma
          }
        }
      }
      }
      \ar[out=180-82.8+90, in=82.8+90, looseness=2.6,
        shift right=2pt,
        "\scalebox{.77}{$\mathclap{
          G
        }$}"{description, pos=.4},
      ]
      \ar[
        dd,
        "\Gamma \backslash (-)"
      ]
      \\[-5pt]
      &
      \xmapsto{\;\;
        \left\vert-\right\vert
      \;\;}
      &
      \\[+5pt]
      \Maps{}
      {
        \mathbf{E}G
      }
      {
        \mathbf{B}\Gamma
      }
      \ar[out=180-82.8+90, in=82.8+90, looseness=2.6,
        shift right=2pt,
        "\scalebox{.77}{$\mathclap{
          G
        }$}"{description, pos=.4},
      ]
      &&
      \underset{
        \;\;\;
        \scalebox{.7}{$
          B(G \acts \,  \Gamma)
        $}
      }{
      \underbrace{
        \TopologicalRealization{\big}
        {
          \Maps{}
          {
            \mathbf{E}G
          }
          {
            \mathbf{B}\Gamma
          }
        }
      }
      }
      \ar[out=180-82.8+90, in=82.8+90, looseness=2.6,
        shift right=2pt,
        "\scalebox{.77}{$\mathclap{
          G
        }$}"{description, pos=.4},
      ]
    \end{tikzcd}
    \;\;\;\;\;
    \in
    \;
    \FormallyPrincipalBundles{(G \acts \, \Gamma)}
    \big(
      \Actions{G}(\kTopologicalSpaces)
    \big)
    \,.
  $$
  \end{enumerate}
\end{proposition}
\begin{proof}
 {\bf  (i)} The equivariant contraction
  is the image
  under
  Lem. \ref{TopologicalRealizationOfEquivalenceOfGroupoidsIsHomotopyEquivalence}
  and
  Prop. \ref{TopologicalRealizationOfTopologicalGroupoidsRespectsEquivariance}
  of the contraction as an equivariant groupoid
  from Lem. \ref{UniversalPrincipalGroupoidIsEquivalentToThePoint}.

  \noindent
  {\bf (ii)} The quotient coprojection is the image of
  the coprojection of equivariant groupoids
  from Lem. \ref{QuotientCoprojectionOfUniversalEquivariantPrincipalGroupoid},
  using that
  this is preserved first by the nerve operation \eqref{SimplicialTopologicalNerveOfTopologicalGroupoids},
  due to
  Lem. \ref{NervePreservesLeftQuotientsOfRightActionGroupoids}
  with
  Lem. \ref{MappingGroupoidIntoCodiscreteGroupoidIsRightActionGroupoid},
  and then by topological realization \eqref{TopologicalRealizationOfSimplicialTopologicalSpaces},
  which preserves all colimits.

  \noindent
{\bf  (iii)} Formal principality follows from that
  that of Lem. \ref{UniversalEquivariantPrincipalGroupoidIsFormallyPrincipal}
  under the functor which is induced \eqref{FunctorOnStructuresInducedFromLexFunctor}
  from the fact that topological realization preserves finite limits
  (Prop. \ref{TopologicalRealizationPreservesFiniteLimits})
  and using that it sends the constant groupoid
  $\ConstantGroupoid(\Gamma)$ to its underlying topological space $\Gamma$
  (Ex. \ref{TopologicalRealizationOfConstantGroupoids})
      \vspace{-2mm}
  $$
    \begin{tikzcd}[column sep=40pt]
      \FormallyPrincipalBundles{\Gamma}
      \left(\!
        \Actions{G}
        \left(
          \Groupoids(\kTopologicalSpaces)
        \right)
      \!\right)
      \ar[
        rr,
        "{
        \scalebox{0.8}{$    \FormallyPrincipalBundles{
            (
              \vert- \vert
            )
          }
          (
            \vert-\vert
          )
          $}
        }"
      ]
      &&
      \FormallyPrincipalBundles{\Gamma}
      \left(
        \Actions{G}
        (
          \kTopologicalSpaces
        )
      \right)
      \mathrlap{\,.}
    \end{tikzcd}
  $$

  \vspace{-8mm}
\end{proof}

\medskip

\noindent
{\bf Fixed loci of equivariant classifying spaces.} For the purposes of equivariant homotopy theory, the key property of the above equivariant universal bundles are (the homotopy types of) their fixed

\begin{lemma}[Fixed loci of base of universal equivariant principal groupoid]
  \label{FixedLociOfBaseOfUniversalEquivariantPrincipalGroupoid}
  For
  $G \in \Groups(\Sets) \xhookrightarrow{\;} \Groups(\kTopologicalSpaces)$
  a discrete group and
  $\Gamma \,\in\, \Groups( \Actions{G}(\kTopologicalSpaces))$,
  we have for each subgroup $H \subset G$
  an equivalence \eqref{EquivalenceOfTopologicalGroupoids}
  of $W(H)$-equivariant (Ntn. \ref{GActionOnTopologicalSpaces})
  topological groupoids (Def. \ref{EquivariantTopologicalGroupoids}):
  \vspace{-2mm}
  \begin{align}
    &
    \Maps{}
    {
      \mathbf{E}G
    }
    {
      \mathbf{B}\Gamma
    }^H
    \\
    \label{FixedLocusOfMappingGroupoidOutOfPairGroupoidOfEquivarianceGroup}
    &
    \underset{\mathclap{\mathrm{htpy}}}{\simeq}
    \;
    \SliceMaps{\big}{\mathbf{B}H}
    {
      \mathbf{B}H
    }
    {
      \mathbf{B}(\Gamma \rtimes H)
    }
    \\
    \label{FixedLociOfBaseOfMurayamaShimakawaConstructionInTermsOfCrossedHomomorphisms}
    &
      \simeq
      \ActionGroupoid
        { \CrossedHomomorphisms(H, \, H \acts \, \Gamma)_{\mathrm{ad}} }
        { \Gamma }
    \qquad
    \in
    \;
    \Actions{W(H)}\left(
      \TopologicalGroupoids
    \right)
  \end{align}

  \vspace{-3mm}
  \noindent
  between

  {\bf (a)}
  the $H$-fixed sub-groupoid of the
  universal equivariant base groupoid \eqref{TheQuotientCoprojectionOnUniversalEquivariantPrincipalgroupoid},
  and

  {\bf (b)} the groupoid of sections of $\mathbf{B}(\Gamma \rtimes H)$
  (Def. \ref{GroupoidOfSectionsOfDeloopedSemidirectProductProjection}),
   hence (by Prop. \ref{ConjugationGroupoidOfCrossedHomomorphismsIsSectionsOfDeloopedSemidirectProductProjection})

   {\bf (c)}
  the conjugation groupoid of crossed homomorphisms (Ntn. \ref{ConjugationGroupoidOfCrossedHomomorphisms}).
\end{lemma}
This statement corresponds to that of \cite[Thm. 4.14, Cor. 4.15]{GuillouMayMerling17}.
We mean to give a detailed proof. Its ingredients are needed below in
the proof of Prop. \ref{WeylGroupActionOnConnectedComponentsOfFixedLociInEquivariantClassifyingSpace}.
\begin{proof}
  First, observe that a $G$-equivariant function \eqref{EquivariantFunctions}
  on $G^L \times G^L$ \eqref{LeftMultiplicationAndInverseRightMultiplicationActionsOnATopologicalGroup}
  is equivalently a general function on $G$:
  \vspace{-3mm}
  $$
    \begin{tikzcd}[row sep=-7pt]
      \mathrm{Maps}
      \left(
        G^L \times G^L
        ,\,
        \Gamma
      \right)^G
      \ar[
        rr,
        <->,
        "{\sim}"
      ]
      &&
      \mathrm{Maps}
      (
        G
        ,\,
        \Gamma
      )
      \ar[
        rr,
        <->,
        "{\sim}"
      ]
      &[-20pt]&[-20pt]
      \mathrm{Maps}_{{}_{/G}}
      (
        G
        ,\,
        \Gamma \times G
      )
      \\
      \scalebox{0.7}{$
      \big(
        (g_1,g_2)
          \,\mapsto\,
        F(g_1,g_2)
      \big)
      $}
      &\longmapsto&
         \scalebox{0.7}{$ \left(
        g
        \,\mapsto\,
        F(e,g)
      \right)
      $}
      &\longmapsto&
         \scalebox{0.7}{$ \left(
        g
          \,\mapsto\,
        \left(
          F(e,g)
          ,\,
          g
        \right)
      \right)
      $}
      \\
          \scalebox{0.7}{$
      \left(
        (g_1,g_2)
        \,\mapsto\,
        g_1
          \cdot
        f
        \left(
          g_1^{-1} \cdot g_2
        \right)
      \right)
      $}
      &\rotatebox{180}{$\longmapsto$}&
         \scalebox{0.7}{$ \left(
        g
          \,\mapsto\,
        f(g)
      \right)
      $}
      &\rotatebox{180}{$\longmapsto$}&
    \scalebox{0.7}{$
            \left(
        g
          \,\mapsto\,
        \left(
          f(g)
          ,\,
          g
        \right)
      \right)
    $}
    \end{tikzcd}
  $$

  \vspace{-2mm}
  \noindent
  Using the isomorphism of $\mathbf{E}G$ with the pair groupoid on $G$
  (Ex. \ref{ActionGroupoidOfLeftGroupMultiplicationIsPairGroupoid}),
  this applies to the component functions of
  functors and natural transformations
  to yield the claim for the special case $H = G$:
  \vspace{-3mm}
  \begin{equation}
    \label{IsomorphismOfGFixedMappingGroupoidOutOfGPairGroupoidWithMappingOutOfGDelooping}
    \hspace{-1cm}
    \begin{tikzcd}[row sep=0pt, column sep=40pt]
      \mathrm{Maps}
      \left(
        (
        G^L \times G^L
        \rightrightarrows
        G^L
        )
        ,\,
        \Gamma^\alpha
          \rightrightarrows
        \ast
      \right)^G
      \ar[
        rr,
        "{\sim}"{above,yshift=-1pt}
      ]
      &&
      \mathrm{Maps}_{{}_{/(G \rightrightarrows \ast)}}
      \left(
        G \rightrightarrows \ast
        ,\,
        \Gamma \rtimes_\alpha G \rightrightarrows \ast
      \right)
    \end{tikzcd}
  \end{equation}
  \vspace{-.4cm}
  $$
    \left(
    \!
    \begin{tikzcd}[column sep=30pt]
      \mathrm{e}
      \ar[
        d
      ]
      &[-20pt]
      &[-5pt]
      \bullet
      \ar[
        rr,
        "{\eta(\mathrm{e})}"
      ]
      \ar[
        d,
        "{
          F(\mathrm{e},g_1)
        }"{left}
      ]
      &&[30pt]
      \bullet
      \ar[
        d,
        "{
          F'(\mathrm{e},g_1)
        }"{left}
      ]
      \\
      g_1
      \ar[
        d
      ]
      &\mapsto&[20pt]
      \bullet
      \ar[
        rr,
        "{
          \eta(g_1)
        }"{description}
      ]
      \ar[
        d,
        "{
          F(g_1,g_2)
        }"{left}
      ]
      &&
      \bullet
      \ar[
        d,
        "{
          F(g_1,g_2)
        }"{left}
      ]
      \\
      g_2
      &
      &
      \bullet
      \ar[
        rr,
        "{
          \eta(g_2)
        }"{below}
      ]
      &&
      \bullet
    \end{tikzcd}
    \!
    \right)
    \;\;\;\;\;
    \longmapsto
    \;\;\;\;\;
    \left(
    \!
    \begin{tikzcd}[column sep=30pt]
      \bullet
      \ar[
        d,
        "{g_1}"{right}
      ]
      &[-15pt]
      &
      \bullet
      \ar[
        rr,
        "{
          \scalebox{1.3}{$($}
            \eta(\mathrm{e})
            ,\
            \mathrm{e}
          \scalebox{1.3}{$)$}
        }"
      ]
      \ar[
        d,
        "{
          \scalebox{1.3}{$($}
            F(\mathrm{e},g_1)
            ,\,
            g_1
          \scalebox{1.3}{$)$}
        }"{left}
      ]
      &&[40pt]
      \bullet
      \ar[
        d,
        "{
          \scalebox{1.3}{$($}
            F'(\mathrm{e},g_1)
            ,\,
            g_1
          \scalebox{1.3}{$)$}
        }"{left}
      ]
      \\
      \bullet
      \ar[
        d,
        "{
          {
            g' :=
          }
          \atop
          {
            g_1^{-1} \cdot g_2
          }
        }"{right}
      ]
      &\mapsto \phantom{AAA}&
      \bullet
      \ar[
        rr,
        "{
          \scalebox{1.3}{$($}
            \eta(\mathrm{e})
            ,\,
            \mathrm{e}
          \scalebox{1.3}{$)$}
        }"{description}
      ]
      \ar[
        d,
        "{
          \scalebox{1.3}{$($}
            F(\mathrm{e},\, g')
            ,\,
            g'
          \scalebox{1.3}{$)$}
        }"{left}
      ]
      &&
      \bullet
      \ar[
        d,
        "{
          \scalebox{1.3}{$($}
            F'(\mathrm{e},\, g')
            ,\,
            g'
          \scalebox{1.3}{$)$}
        }"{left}
      ]
      \\
      \bullet
      &
      &
      \bullet
      \ar[
        rr,
        "{
          \scalebox{1.3}{$($}
            \eta(\mathrm{e})
            ,\,
            \mathrm{e}
          \scalebox{1.3}{$)$}
        }"{below}
      ]
      &&
      \bullet
    \end{tikzcd}
    \!
    \right)
  $$
  Notice that it is the semidirect product group structure
  that does make the assignment on the right
  be functorial, given the data on the left,
  e.g.:
  \vspace{-2mm}
  $$
    (F(\mathrm{e},g_1), g_1) \cdot ( \eta(\mathrm{e}), \mathrm{e} )
    \,=\,
    \left(
      F(\mathrm{e},g_1) \cdot \alpha(g_1)(\eta(\mathrm{e}))
      ,\,
      g_1
    \right)
    \,=\,
    \left(
      F(\mathrm{e},g_1) \cdot \eta(g_1)
      ,\,
      g_1
    \right).
  $$

  \vspace{-2mm}
\noindent  In order to generalize this proof to arbitrary $H \subset G$,
  choose a section of the coset space projection
  \vspace{-2mm}
  \begin{equation}
    \label{SectionOfDiscreteCosetSpaceCoprojection}
    \begin{tikzcd}[row sep=small]
      & G
      \ar[
        d,
        ->>
      ]
      \\
      G/H
      \ar[r,-,shift left=1pt]
      \ar[r,-,shift right=1pt]
      \ar[ur,dashed, "\sigma"]
      &
      G/H
      \mathrlap{\,,}
      &
      \;\;\;\;\;\;\;
      \mbox{such that $\sigma([\mathrm{e}]) = \mathrm{e}$}
      \,,
    \end{tikzcd}
  \end{equation}

      \vspace{-2mm}
\noindent
  which exists and is continuous
  by the assumption that $G$ is discrete.
  This induces a decomposition of the underlying set of $G$ into $H$-orbits
  $$
    G
    \;\;
    \,\simeq\,
    \underset{
      [g] \in G/H
    }{\bigcup}
    H \cdot \sigma([g])
    \;\;\;\;\;\;
    \in
    \;
    \Sets
  $$
  and thus implies that the pair groupoid
  $(G^L \times G^L \rightrightarrows G^L)$ is generated,
  under (i) composition, (ii) taking inverses and (iii) acting with elements of $H$,
  by the following two classes of morphisms:
      \vspace{-2mm}
  \begin{equation}
    \label{GeneratingMorphismsOfPairGroupoidUnderCompositionInversesAndHAction}
    \left\{
    (
      \mathrm{e}
      \xrightarrow{  }
      h
    )
    \;\vert\;
    h \,\in\, H
    \right\}
    ,\;\;\;\;
    \left\{
    \left(
      \mathrm{e}
      \xrightarrow{  }
      \sigma([g])
    \right)
    \; \vert\;
    [g] \,\in\, G/H
    \right\}
    \;\;\;\;\;\;
    \subset
    \;
    G \times G
    \mathrlap{\,.}
  \end{equation}

  \vspace{-2mm}
  \noindent
  Using this, consider the following expression for a pair of
  continuous functors:
  \vspace{-3mm}
  \begin{equation}
    \label{FunctorsBetweeHFixedMappingGroupoidSectionsGroupoidOfDeloopedSemidirectProductWithH}
    \begin{tikzcd}[row sep=-8pt]
\scalebox{0.7}{$      \left(
        (g_1,g_2)
        \,\mapsto\,
        F(g_1, g_1)
      \right)
      $}
      &\longmapsto&
      \;\;\;\;\;\;\;\;\;\;\;\;\;\;
 \scalebox{0.7}{$     \left(
        h \,\mapsto\, \left( F(\mathrm{e},h),\, h \right)
      \right)
      $}
      \\
      \mathrm{Maps}
      \left(
        (
          G \times G
            \rightrightarrows
          G
        )
        ,\,
        (
          \Gamma
            \rightrightarrows
          \ast
        )
      \right)^H
      \ar[
        rr,
        shift left=2pt,
        "{L}"
      ]
      &&
      \mathrm{Maps}_{{}_{/(H \rightrightarrows \ast)}}
      \left(
        (
          H \rightrightarrows \ast
        )
        ,\,
       (
          \Gamma \rtimes H \rightrightarrows \ast
        )
      \right)
      \ar[
        ll,
        shift left=2pt,
        "{R}"{below}
      ]
      \\
    \scalebox{0.7}{$  \left(\!\!\!
        \def\arraystretch{.9}
        \begin{array}{l}
          \,
          (\mathrm{e},h) \;\;\;\mapsto\; \phi(h)
          \\
          \scalebox{1.1}{$($}
            \mathrm{e}, \sigma([g])
          \scalebox{1.1}{$)$}
          \mapsto
          \mathrm{e}
        \end{array}
      \!\!\!\right)
      $}
      &\rotatebox{180}{$\longmapsto$}&
      \;\;\;\;\;\;\;\;\;\;\;\;\;
      \scalebox{0.7}{$
      \left(
        h \,\mapsto\, \left( \phi(h) \,, h \right)
      \right)
      $}
      \mathrlap{\,,}
    \end{tikzcd}
  \end{equation}

  \vspace{-2mm}
\noindent  where $L$ is restriction along
  $
    ( H \times H \rightrightarrows H)
      \xhookrightarrow{\;}
    ( G \times G \rightrightarrows G)
  $
  followed by the isomorphism \eqref{IsomorphismOfGFixedMappingGroupoidOutOfGPairGroupoidWithMappingOutOfGDelooping}
  for $G = H$, while $R$ is given on morphisms
  as follows:
  \vspace{-2mm}
  \begin{equation}
  \label{RightComparisonFunctorFromSectionedFunctorsToFixedLociInEquivariantClassifyingGroupoid}
  \hspace{-3mm}
  \left(\!\!\!\!
  \begin{tikzcd}
    \mathrm{e}
    \ar[
      d
    ]
    &[-20pt]\mapsto&[-20pt]
    \bullet
    \ar[
      rrr,
      "{
        \eta(\mathrm{e})
        \,\coloneqq\,
        \eta(\bullet)
      }"
    ]
    \ar[
      d,
      "{
        \phi(h)
      }"{left}
    ]
    &&&
    \bullet
    \ar[
      d,
      "\phi'(h)"{left}
    ]
    \\
    h
    &\mapsto&
    \bullet
    \ar[
      rrr,
      "{
        \eta(h)
        \,\coloneqq\,
        \alpha(h)\left(\eta(\bullet)\right)
      }"{below}
    ]
    &&&
    \bullet
    \\[-10pt]
    \mathrm{e}
    \ar[d]
    &\mapsto&
    \bullet
    \ar[
      rrr,
      "{
        \eta(\mathrm{e})
        \coloneqq\,
        \eta(\bullet)
      }"{below  }
    ]
    \ar[
      d,
      "{\mathrm{e}}"{left}
    ]
    &&&
    \bullet
    \ar[
      d,
      "{\mathrm{e}}"{left}
    ]
    \\
    \sigma([g])
    &\mapsto&
    \bullet
    \ar[
      rrr,
      "{
        \eta\left(\sigma([g])\right)
        \,\coloneqq\,
        \eta(\bullet)
      }"{below}
    ]
    &&&
    \bullet
  \end{tikzcd}
  \right)
  \;\;\;\;\;\;\;\;
  \overset{R}{\rotatebox{180}{$\longmapsto$}}
  \;\;\;\;\;\;\;\;\;
  \left(
  \begin{tikzcd}
    \bullet
    \ar[
      d,
      "{h}"
    ]
    &[-20pt]\mapsto&[-20pt]
    \bullet
    \ar[
      rrr,
      "{
        \left(
          \eta(\bullet)
          ,\,
          \mathrm{e}
        \right)
      }"
    ]
    \ar[
      d,
      "{
        \left(
          \phi(h)
          ,\,
          h
        \right)
      }"{right}
    ]
    &&&
    \bullet
    \ar[
      d,
      "{
        \left(
          \phi'(h)
          ,\,
          h
        \right)
      }"{right}
    ]
    \\
    \bullet
    &\mapsto&
    \bullet
    \ar[
      rrr,
      "{
        \left(
          \eta(\bullet)
          ,\,
          \mathrm{e}
        \right)
      }"{below}
    ]
    &&&
    \bullet
    \mathrlap{\,.}
  \end{tikzcd}
  \!\! \right)
  \end{equation}

\vspace{-2mm}
  \noindent
  One readily checks that this is well-defined
  and that $L \circ R \,=\, \mathrm{id}$.
  Therefore, it now suffices to give a
  continuous natural transformation
  $\mathrm{id} \,\xRightarrow{\eta}\, R \circ L$.
  This is obtained by choosing for any functor $F$
  in the groupoid on the left
  of \eqref{FunctorsBetweeHFixedMappingGroupoidSectionsGroupoidOfDeloopedSemidirectProductWithH}
  a natural transformation $\eta_F \,\colon\, F \Rightarrow R \circ L(F)$,
  given by the following component function:
  \vspace{-2mm}
  $$
    \begin{tikzcd}
      \mathrm{e}
      \ar[
        d
      ]
      &\mapsto&&
      \bullet
      \ar[
        rrrr,
        "{
          \eta_{{}_{F}}\!(\mathrm{e}) \,\coloneqq\, \mathrm{e}
        }"{above}
      ]
      \ar[
        d,
        "{
          F(\mathrm{e},h)
        }"{left}
      ]
      &&&&
      \bullet
      \ar[
        d,
        "{
          \left(R \circ L (F)\right)(\mathrm{e},h) \,=\, F(\mathrm{e},h)
        }"
      ]
      \\[-6pt]
      h
      &\mapsto&&
      \bullet
      \ar[
       rrrr,
       "{
         \eta_{{}_F}\!(h) \,\coloneqq\, \mathrm{e}
       }"{above}
      ]
      &&&&
      \bullet
      \\[-5pt]
      \mathrm{e}
      \ar[
        d
      ]
      &\mapsto&&
      \bullet
      \ar[
        rrrr,
        "\eta_{{}_{F}}\!(\mathrm{e}) \,\coloneqq\, \mathrm{e} "{below}
      ]
      \ar[
        d,
        "{
          F\scalebox{1.1}{$($} \mathrm{e}, \sigma([h]) \scalebox{1.1}{$)$}
        }"{left}
      ]
      &&&&
      \bullet
      \ar[
        d,
        "{
          \left(R \circ L (F)\right)
          \scalebox{1.1}{$($} \mathrm{e}, \sigma([g])\scalebox{1.1}{$)$}
            \,=\,
          \mathrm{e}
        }"
      ]
      \\
      \sigma([g])
      &\mapsto&&
      \bullet
      \ar[
       rrrr,
       "{
         \eta_{{}_F}\!\scalebox{1.1}{$($} \sigma([g]) \scalebox{1.1}{$)$}
         \,\coloneqq\,
         F\scalebox{1.1}{$($} \mathrm{e}, \sigma([g]) \scalebox{1.1}{$)$}^{-1}
         \mathrlap{\,.}
       }"{below}
      ]
      &&&&
      \bullet
    \end{tikzcd}
  $$

  \vspace{-1mm}
  \noindent
  It only remains to see that this is indeed natural in $F$,
  which amounts to checking that
  for any $H$-equivariant natural transformation $\beta$ from $F$ to $F'$
  the following two squares on the right commute
  (in $\mathbf{B}\Gamma$):
    \vspace{-2mm}
  $$
    \begin{tikzcd}
      F
      \ar[
        d,
        "{
          \beta
        }"
      ]
      \ar[
        rr,
        "{
          \eta_F
        }"
      ]
      &&
      R \circ L(F)
      \ar[
        d,
        "{
          R \circ L(\beta)
        }"
      ]
      \\
      F'
      \ar[
        rr,
        "{
          \eta_{F'}
        }"
      ]
      &&
      R \circ L(F')
    \end{tikzcd}
    \;\;\;\;
    \colon
    \;\;\;\;
    \left\{
    \begin{array}{l}
    \;\;\;\;\;
    h
    \;\;\;\;\;
    \;\;\;\;
      \mapsto
    \;\;\;\;
    \begin{tikzcd}[column sep=huge]
      \bullet
      \ar[
        d,
        "{ \beta(h) }"
      ]
      \ar[
        rr,
        "{
          \eta_{{}_{F}}\!(h) \,=\, \mathrm{e}
        }"
      ]
      &&
      \bullet
      \ar[
        d,
        "{
          \left(R \circ L(\beta)\right)(h)
          \,=\,
          \alpha(h)\left(\beta(e)\right)
        }"
      ]
      \\
      \bullet
      \ar[
        rr,
        "{
          \eta_{{}_{F'}}\!(h) \,=\, \mathrm{e}
        }"{below}
      ]
      &&
      \bullet
    \end{tikzcd}
    \\
    \sigma([g])
    \;\;\;\;
      \mapsto
    \;\;\;\;
    \begin{tikzcd}[column sep=huge]
      \bullet
      \ar[
        d,
        "{ \beta\scalebox{1.1}{$($} \sigma([g]) \scalebox{1.1}{$)$} }"
      ]
      \ar[
        rr,
        "{
          \mathclap{
          \eta_{{}_{F}}\!\scalebox{1.1}{$($}\sigma([g]) \scalebox{1.1}{$)$}
          \,=\,
          F(
            e
            , \,
            \sigma([g])
          )^{-1}
          }
        }"
      ]
      &&
      \bullet
      \ar[
        d,
        "{
          \left(R \circ L(\beta)\right)
          \scalebox{1.1}{$($}
            \sigma([g])
          \scalebox{1.1}{$)$}
          \,=\,
          \beta(\mathrm{e})
        }"
      ]
      \\
      \bullet
      \ar[
        rr,
        "{
          \mathclap{
          \eta_{{}_{F'}}\!\scalebox{1.1}{$($}\sigma([g]) \scalebox{1.1}{$)$}
          \,=\,
          F'(
            e
            , \,
            \sigma([g])
          )^{-1}
          }
        }"{below}
      ]
      &&
      \bullet
    \end{tikzcd}
    \end{array}
    \right.
  $$

\vspace{-2mm}
\noindent
But the top square commutes by the $H$-equivariance of $\beta$ \eqref{EquivariantFunctions},
the bottom one by the naturality of $\beta$ \eqref{NaturalitySquareForTopologicalGroupoids}.
\end{proof}

\begin{proposition}[Connected components of fixed loci in equivariant classifying spaces]
  \label{ConnectedComponentsOfHFixedEquivariantClassifyingSpaceWhenCrossedConjugationQuotientIsDiscrete}
$\,$

\noindent {\bf (i)}   Let

  {\bf (a)}
  $
    G
      \,\in\,
    \Groups(\FiniteSets)
    \xhookrightarrow{\;}
    \Groups(\kTopologicalSpaces)
  $
  be a finite group;

  {\bf (b)}
  $\Gamma \,\in\, \Groups\left( \Actions{G}(\kTopologicalSpaces) \right)$;

  {\hypertarget{AssumptionThatH1GrpIsDiscrete}{}\bf (c)} such that
  $\CrossedHomomorphisms(G, \, G\acts \, \Gamma)/\sim_{\mathrm{ad}}$
  \eqref{NonAbelianGroup1Cohomology} is a discrete space (e.g. via Prop. \ref{DiscreteSpacesOfCrossedConjugacyClassesOfCrossedHomomorphisms}).

  \noindent
  Then for each subgroup $H \,\subset\, G$,
  the connected components of
  the $H$-fixed locus of the base space
  $G \acts \, B(G \acts \, \Gamma)$ \eqref{QuotientCoprojectionOfUniversalEquivariantPrincipalBundle}
  of the universal $G$-equivariant $\Gamma$-principal bundle;
  are the non-abelian group cohomology in degree 1
  \eqref{NonAbelianGroup1Cohomology}
  of $H$ with coefficients in the restricted action $H \acts \, \Gamma$:
    \vspace{-2mm}
  \begin{equation}
    \label{ConnectedComponentsOfFixedLociOfEquivariantClassifyingSpace}
    \pi_0
    \left(
      \left(
        B(G \acts \, \Gamma)
      \right)^H
    \right)
    \;\;
    \simeq
    \;\;
    H^1_{\mathrm{Grp}}
    (
      H
      ,\,
      H \acts \, \Gamma
    )
    \;\;\;
    \in
    \;
    \Sets
    \,.
  \end{equation}

\vspace{-2mm}
\noindent
{\bf (ii)} If, moreover,

  {\bf (d)} $\Gamma$ is a Lie group and
    $\alpha \,:\, G \to \mathrm{Aut}_{\mathrm{Grp}}(\Gamma)$
    restricts to the identity on the center of $G$,

\noindent
then the underlying topological space of the
fixed locus itself is homotopy equivalent to
\vspace{-2mm}
\begin{equation}
  \left(
    B(G \acts \, \Gamma)
  \right)^H
  \;\;
  \underset{
    \mathrm{htpy}
  }{\simeq}
  \qquad
  \underset
    {
      {[\phi] \in }
      \atop
      \mathclap{
        \scalebox{.6}{$
          H^1_{\mathrm{Grp}}(H,\, H \acts \, \Gamma)
        $}
      }
    }
    {\coprod}
    \quad
    \frac{
      \Gamma/C_{{}_{\Gamma}}(\phi) \times E \Gamma
    }
    {\Gamma}
  \;\;\;
  \in
  \;
  \kTopologicalSpaces
  \,,
\end{equation}

\vspace{-2mm}
\noindent
hence (when the $\Gamma \to \Gamma/C_\Gamma(\phi)$ admits local sections...)
to a disjoint union, indexed by classes $[\phi]$ in the group 1-cohomology set,
of ordinary classifying spaces of the stabilizer groups of the cocycles $\phi$:
\vspace{-2mm}
$$
  \left(
    B(G \acts \, \Gamma)
  \right)^H
  \;\;\;
  \simeq
  \qquad
  \underset{
    { [\phi] \in }
    \atop
    \mathclap{
      \scalebox{.6}{$
        H^1_{\mathrm{Grp}}(H ,\, H \acts \,  \Gamma)
      $}
    }
  }{\coprod}
  \quad
  B
  (
    C_{{}_{\Gamma}}(\phi)
  )
  \;\;\;\;
  \in
  \;
  \HomotopyCategory(\kTopologicalSpaces_{\mathrm{Qu}})
  \,.
$$
\end{proposition}
\begin{proof}
First, observe that we have a
homotopy equivalence
of the fixed locus with
the Borel construction \eqref{BorelConstructionAsTopologicalRealizationOfActionGroupoid}
of the space of crossed homomorphisms $H \to \Gamma$ \eqref{SpaceOfCrossedHomomorphisms}
by the $\Gamma$-action of crossed conjugations \eqref{CrossedConjugationAction}:
    \vspace{-2mm}
\begin{equation}
  \label{FixedLociInEquivariantClassifyingSpaceInTermsOfCrossedHomomorphisms}
    \left(
      B(G \acts \, \Gamma)
    \right)^H
    \;\;
    \underset{
      \mathrm{htpy}
    }{\simeq}
    \;\;
    \frac{
      \scalebox{1}{$
        \CrossedHomomorphisms(G,\, G \acts \, \Gamma)
      $}
    }{
      \Gamma
    }
    \;\;\;
    \in
    \;
    \kTopologicalSpaces
    \,.
\end{equation}

    \vspace{-2mm}
\noindent
obtained as the following composite:
\vspace{-3mm}
$$
  \def\arraystretch{1.8}
  \begin{array}{lll}
    \left(
     B(G \acts \, \Gamma)
    \right)^H
    &
    \;\simeq\;
    \big(
      \Gamma
        \backslash
      \left\vert
        \mathrm{Maps}
        (\mathbf{E}G ,\, G \acts \,  \Gamma)
      \right\vert
    \big)^H
    &
    \mbox{\small by \eqref{QuotientCoprojectionOfUniversalEquivariantPrincipalBundle}}
    \\
    &
    \;\simeq\;
    \left\vert
      \Gamma
        \backslash
      \mathrm{Maps}
      (\mathbf{E}G ,\, \mathbf{E}\Gamma)
    \right\vert^H
    &
    \mbox{\small by
     Lem. \ref{NervePreservesLeftQuotientsOfRightActionGroupoids}
     \&
     Lem. \ref{MappingGroupoidIntoCodiscreteGroupoidIsRightActionGroupoid},
    }
    \\
    &
    \;\simeq\;
    \left\vert
      \mathrm{Maps}
      (\mathbf{E}G ,\, \mathbf{B}\Gamma)
    \right\vert^H
    &
    \mbox{\small by Lem. \ref{QuotientCoprojectionOfUniversalEquivariantPrincipalGroupoid}}
    \\
    &
    \;\simeq\;
    \left\vert
      \left(
        \mathrm{Maps}
        (\mathbf{E}G ,\, \mathbf{B}\Gamma)
      \right)^H
    \right\vert
    &
    \mbox{\small by Lem. \ref{TopologicalRealizationPreservesFiniteLimits}}
    \\
    &
    \;
    \underset{
      \mathrm{htpy}
    }{
      \simeq
    }
    \;
    \left\vert
      \CrossedHomomorphisms(H ,\, H \acts \, \Gamma )
        \sslash_{\! \mathrm{ad}}
      \Gamma
    \right\vert
    &
    \mbox{\small by Lem. \ref{FixedLociOfBaseOfUniversalEquivariantPrincipalGroupoid}}
    \\
    & \;\simeq\;
    \frac{
            \CrossedHomomorphisms(H ,\, H \acts \, \Gamma)
        \times
        E \Gamma
          }{
      \Gamma
    }
    &
    \mbox{\small by Ex. \ref{GroupoidLevelBorelConstruction}}
    \mathrlap{\,.}
  \end{array}
$$

  \vspace{-2mm}
\noindent
With this,
the first claim arises by the following sequence of bijections:
  \vspace{-2mm}
$$
  \hspace{-2cm}
  \def\arraystretch{1.8}
  \begin{array}{lll}
    \pi_0
    \left(\!\!
      \left(
         B(G \acts \, \Gamma)
      \right)^H
    \right)
    &
    \;\simeq\;
    \pi_0
    \left(
    \frac{
        \CrossedHomomorphisms(H ,\, H \acts \, \Gamma)
        \times
        E \Gamma
    }{
      \Gamma
    }
    \right)
    &
    \mbox{\small by \eqref{FixedLociInEquivariantClassifyingSpaceInTermsOfCrossedHomomorphisms}}
    \\
    & \;\simeq\;
    \pi_0
    \left(
      \frac{
          \CrossedHomomorphisms(H,\, H \acts \, \Gamma)
      }{
        \Gamma
      }
    \right)
    &
    \mbox{\small by \eqref{PathConnectedComponentsOfBorelConstruction}}
    \\
    & \;\simeq\;
    \frac{
        \CrossedHomomorphisms(H ,\, H \acts \, \Gamma)
    }{
      \Gamma
    }
    &
    \mbox{\small by$\;$ \hyperlink{AssumptionThatH1GrpIsDiscrete}{(iii)} }
    \\
    &
    \;\simeq\;
    H^1_{\mathrm{Grp}}(H,\, H \acts \, \Gamma)
    &
    \mbox{\small by \eqref{NonAbelianGroup1Cohomology}}
    \,,
  \end{array}
$$

  \vspace{-2mm}
\noindent
and the second by the following homeomorphisms:
  \vspace{-3mm}
$$
  \hspace{-2cm}
  \def\arraystretch{1.8}
  \begin{array}{lll}
    \left(
       B(G \acts \, \Gamma)
    \right)^H
    &
    \;\simeq\;
    \frac{
        \CrossedHomomorphisms(H,\, H \acts \, \Gamma)
        \times
        E \Gamma
    }{
      \Gamma
    }
    & \;
    \mbox{\small by \eqref{FixedLociInEquivariantClassifyingSpaceInTermsOfCrossedHomomorphisms}}
    \\
    & \;\simeq\qquad
    \underset
      {
        {[\phi] \in }
        \atop
        \mathclap{
          \scalebox{.7}{$
            H^1_{\mathrm{Grp}}(H,\, H \acts \, \Gamma)
          $}
        }
      }
      {\coprod}
      \quad \frac
        {
          \Gamma/C_{{}_{/\Gamma}}(\phi) \,\times\, E \Gamma
        }
        {
          \Gamma
        }
      &
      \mbox{
        by Prop. \ref{DiscreteSpacesOfCrossedConjugacyClassesOfCrossedHomomorphisms}
      }.
  \end{array}
$$

\vspace{-7mm}
\end{proof}

\begin{proposition}[Weyl group action on connected components of fixed loci in equivariant classifying space]
  \label{WeylGroupActionOnConnectedComponentsOfFixedLociInEquivariantClassifyingSpace}
  Under the equivalence of Prop. \ref{ConnectedComponentsOfHFixedEquivariantClassifyingSpaceWhenCrossedConjugationQuotientIsDiscrete}
  the residual $W(H)$-action on $( B(G \acts \, \Gamma))^H$
  (Ex. \ref{FixedLociWithResidualWeylGroupAction})
  is, on the set of connected components, the
  Weyl group action
  on non-abelian first group cohomology
  from Prop. \ref{WeylGroupActionOnGroup1CohomologyOfSubgroup}:
      \vspace{-2mm}
  $$
    W(H)
      \,\acts\;
    \pi_0\left(\!\! \left( B(G \acts \, \Gamma)\right)^H\right)
    \;\;
    \simeq
    \;\;
    W(H)
      \,\acts\;
    H^1_{\mathrm{Grp}}(H,\, H \acts \, \Gamma)
    \;\;\;
    \in
    \;
    \Actions{W(H)}(\Sets)
    \,.
  $$
\end{proposition}
\begin{proof}
  We may transport the action along the explicit equivalence
  established inside the proof of Lem. \ref{FixedLociOfBaseOfUniversalEquivariantPrincipalGroupoid}:
  For $n \in N(H)$ and $\phi \,\in\, \CrossedHomomorphisms(H,\, H \acts \, \Gamma)$
  we need to check that
    \vspace{-2mm}
  \begin{equation}
    \label{WeylGroupActionOnFixedLocusOfEquClassSpaceToInduceActionOnGroup1Cohomology}
    [
      L( n \cdot (R\phi) )
    ]
    \;=\;
    [\phi_n]
    \;\;\;
    \in
    \;
    H^1_{\mathrm{Grp}}(H,\, H \acts \, \Gamma)\;,
  \end{equation}

    \vspace{-2mm}
\noindent
  where
  $L$ and $R$ are from \eqref{FunctorsBetweeHFixedMappingGroupoidSectionsGroupoidOfDeloopedSemidirectProductWithH},
  $\phi_n$ is from \eqref{NormalizerGroupActionOnCrossedHomomorphismsOutOfSubgroup},
  and $n \cdot (-)$ is the conjugation action by $n \in G$ on
  $\mathrm{Maps}(\mathbf{E}G,\,\mathbf{B}\Gamma)$.

  Notice in the case $n \in H \subset N(H)$
  that the $n$-action on the $H$-fixed locus is trivial by definition,
  while that on $H^1_{\mathrm{Grp}}$ is trivial by
  \eqref{WeylGroupActionOnSetOfGroup1CohomologyOfSubgroup}
  in Prop. \ref{WeylGroupActionOnGroup1CohomologyOfSubgroup},
  so that there is nothing further to be proven.

  Therefore, we may assume now that $n$ is not in $H$.
  Since $L$ does not depend on the choice of section $\sigma$
  in \eqref{SectionOfDiscreteCosetSpaceCoprojection},
  the bijection of connected components induced by $L$ and $R$
  is independent of this choice,
  so that {\it for any $n$}
  we may choose any such section convenient for the analysis.
  We now choose
  a $\sigma$ \eqref{SectionOfDiscreteCosetSpaceCoprojection}
  that takes $n^{-1}$ to be the representative of its
  coset class:
    \vspace{-2mm}
  $$
    \sigma([ n^{-1} ])
    \;
    \coloneqq
    \;
    n^{-1}
    .
  $$

    \vspace{-1mm}
\noindent
  With this choice, the definition of $R$
  \eqref{RightComparisonFunctorFromSectionedFunctorsToFixedLociInEquivariantClassifyingGroupoid}
  says that $R\phi$ sends morphisms (in $\mathbf{E}G$)
  between $n^{-1}$ and the neutral element to
  the neutral element:
    \vspace{-1mm}
  \begin{equation}
    \label{RphiOnMorphismsBetweenInverseOfnAndNeutralElement}
    (R\phi)(\mathrm{e}, n^{-1})
    \;=\;
    \mathrm{e}
    \,,
    \;\;\;\;\;
    (R\phi)(n^{-1}, \mathrm{e})
    \;=\;
    \mathrm{e}
    \,.
  \end{equation}

    \vspace{-1mm}
\noindent
  Now for any $h \in H$, we compute as follows:
    \vspace{-3mm}
  \begin{equation}
    \def\arraystretch{1.4}
    \begin{array}{lll}
      L\left( n \cdot (R\phi) \right)(h)
      & \;=\;
      \left( n \cdot (R\phi) \right)(\mathrm{e}, \, h )
      &
      \mbox{\small by \eqref{FunctorsBetweeHFixedMappingGroupoidSectionsGroupoidOfDeloopedSemidirectProductWithH}}
      \\
      & \;=\;
      \alpha(n)
      \left(
        (R\phi)
        (
          n^{-1}, \, n^{-1} \cdot h
        )
      \right)
      &
      \mbox{\small by \eqref{ConjugationActionOnMapsBetweenGSpaces}}
      \\
      & \;=\;
      \alpha(n)
      \big(
        (R\phi)(\mathrm{e},\, n^{-1} \cdot h)
      \big)
      &
      \mbox{\small by \eqref{RphiOnMorphismsBetweenInverseOfnAndNeutralElement}}
      \\
      &
      \;=\;
      \alpha(n^{-1} \cdot h \cdot n \cdot n)
      \big(
        (R\phi)(n^{-1} \cdot h^{-1} \cdot n,\, n^{-1})
      \big)
      &
      \mbox{\small by $H$-equivariance of $R\phi$}
      \\
      &
      \;=\;
      \alpha(n^{-1} \cdot h \cdot n \cdot n)
      \big(
        (R\phi)(n^{-1} \cdot h^{-1} \cdot n,\, \mathrm{e})
      \big)
      &
      \mbox{\small by \eqref{RphiOnMorphismsBetweenInverseOfnAndNeutralElement}}
      \\
      &
      \;=\;
      \alpha(n)
      \big(
        (R\phi)(\mathrm{e},\, n^{-1} \cdot h \cdot n)
      \big)
      &
      \mbox{\small by $H$-equivariance of $R\phi$}
      \\
      & \;=\;
      \alpha(n)
      \big(
        \phi(n^{-1} \cdot h \cdot n)
      \big)
      &
      \mbox{\small by \eqref{RightComparisonFunctorFromSectionedFunctorsToFixedLociInEquivariantClassifyingGroupoid}}
      \\
      & \;=\;
      \phi_n(h)
      &
      \mbox{\small by \eqref{NormalizerGroupActionOnCrossedHomomorphismsOutOfSubgroup}}
      \,.
    \end{array}
  \end{equation}

    \vspace{-2mm}
\noindent
  This manifestly implies the desired \eqref{WeylGroupActionOnFixedLocusOfEquClassSpaceToInduceActionOnGroup1Cohomology}.
\end{proof}

\begin{remark}[Fixed loci of equivariant classifying spaces for Lie groups]
  Applied to the special case when
  both $G$ and $\Gamma$ are
  compact Lie groups,
  Prop. \ref{ConnectedComponentsOfHFixedEquivariantClassifyingSpaceWhenCrossedConjugationQuotientIsDiscrete}
  with
  Prop. \ref{WeylGroupActionOnConnectedComponentsOfFixedLociInEquivariantClassifyingSpace}
  reproduces
  most of the
  content of \cite[Thm. 10, Thm. 11]{LashofMay86}.
  The statement of
  Prop. \ref{ConnectedComponentsOfHFixedEquivariantClassifyingSpaceWhenCrossedConjugationQuotientIsDiscrete}
  is essentially that of \cite[Thm. 4.23]{GuillouMayMerling17}
  (though it seems some topological conditions are missing there).
\end{remark}

\begin{remark}[Pseudo-principal nature of fixed loci in universal equivariant principal bundle]
  \label{PseudoPrincipalNatureOfFixedLociInUniversalEquivariantPrincipalBundles}
  Since the total space of the equivariant universal bundle
  has contractible fixed loci \eqref{ContractibleFixedLociOfEquivariantUniversalBundle},
  while
  Prop. \ref{ConnectedComponentsOfHFixedEquivariantClassifyingSpaceWhenCrossedConjugationQuotientIsDiscrete}
  says that its base space has, in general, fixed loci with several connected components,
  it follows that all except one of these base components
  carry {\it empty} fibers:
    \vspace{-4mm}
  $$
    \begin{tikzcd}[column sep=2pt, row sep=7pt]
      \left(
        E(G \acts  \, \Gamma)
      \right)^H
      \ar[
        d
      ]
      &\simeq&
      E
      (
        \Gamma^H
      )
      \ar[
        d
      ]
      &[-2pt]
        \sqcup
      &[-6pt]
      \varnothing
      \ar[
        d
      ]
      \\
      \left(
        B(G \acts \, \Gamma)
      \right)^H
      &\simeq&
      B
      (
        \Gamma^H
      )
      &
      \sqcup
      &
      \;\;\;
      \underset{
        { [\phi \neq \mathrm{e}] \in }
        \atop
        \mathclap{\scalebox{.6}{$
            H^1_{\mathrm{Grp}}(H,\,H \acts  \, \Gamma)
          $}
        }
      }{\coprod}
      \;
      B\left(C_{{}_{\Gamma}}(\phi)\right)
    \end{tikzcd}
  $$

  \vspace{-3mm}
  \noindent
  While a bundle with empty fibers cannot be principal,
  these fixed loci of the universal equivariant principal bundle
  are still $W(H)$-equivariant {\it formally principal} $\Gamma^H$-principal
  bundles (Rem. \ref{PseudoTorsorCondition}),
  by Cor. \ref{FixedLociOfEquivariantPrincipalBundles}.
\end{remark}

\medskip

\noindent
{\bf Homotopy types of ordinary classifying spaces.}
We begin by establishing some homotopy theoretic properties
of ordinary (not equivariant) classifying spaces \eqref{QuotientCoprojectionOfUniversalPrincipalBundle}.

\begin{remark}[Local triviality]
Prop. \ref{TopologicalRealizationOfUniversalPrincipalGroupoidIsUniversalPrincipalBundle}
does not establish local trivializability
of the Milgram $\Gamma$-princial bundle (compare Rem. \ref{AssumptionOfLocalTrivializability}),
and this cannot be expected to hold for general $\Gamma$.
Prop. \ref{MilgramBundleIsHomotopyFiberSequence} below is classical
only in the case when $E \Gamma \xrightarrow{\;} B \Gamma$
has been shown to be locally trivial
(using that then it
is guaranteed to be a Serre fibration, by Lem. \ref{LocallyTrivialBundlesAreSerreFibrations}).
While
Prop. \ref{TopologicalRealizationOfUniversalPrincipalGroupoidIsUniversalPrincipalBundle}
follows, by Rem. \ref{TopologicalRealizationOfGoodSimplicialSpacesModelsTheirHomotopyColimit} and
also by general facts of $\infty$-topos theory discussed in \cite{NSS12a},
the following is a more
topological derivation (though it does crucially invoke
Prop. \ref{SufficientConditionsForTopologicalRealizationOfSimplicialSpacesToPreserveHomotopyPullbacks})
which is useful in the present context.
\end{remark}

\begin{lemma}[Homotopy fibration property of Borel construction]
  \label{HomotopyFibrationOfBorelConstruction}
  Let $\Gamma \,\in\, \Groups(\kTopologicalSpaces)$
  be well-pointed (Ntn. \ref{WellPointedTopologicalGroup})
  then for any $\TopologicalSpace \, \rightacts \Gamma \,\in\, \Actions{\Gamma^{\mathrm{op}}}(\kTopologicalSpaces)$
  the Borel construction (Ex. \ref{GroupoidLevelBorelConstruction})
  sits in a homotopy fiber sequence (Ntn. \ref{HomtopyFiberSequenceOfTopologicalSpaces})
  of the form
    \vspace{-2mm}
  \begin{equation}
    \label{HomotopyFiberSequenceOfTopologicalBorelConstruction}
    \begin{tikzcd}[row sep=12pt]
      \TopologicalSpace
        \ar[
          rr,
          hook,
          "{ \HomotopyFiber(b) }"
        ]
      &&
        ({\TopologicalSpace \times E \Gamma})
        /
        {G}
      \ar[
        d,
        "{ \; p }"
      ]
      \\
      &&
      B \Gamma
    \end{tikzcd}
  \end{equation}
\end{lemma}
\begin{proof}
  By Ex. \ref{GroupoidLevelBorelConstruction},
  the sequence \eqref{HomotopyFiberSequenceOfTopologicalBorelConstruction}
  is the
  image under topological realization
  (Ntn. \ref{TopologicalRealizationFunctors})
  of the sequence of simplicial topological spaces which,
  in turn,
  is the image under taking nerves
  \eqref{SimplicialTopologicalNerveOfTopologicalGroupoids}
  of the evident sequence of topological action groupoids:
  \vspace{-2mm}
  \begin{equation}
  \label{SimplicialSpaceAvatarOfTopologicalBorelConstructionFiberSequence}
  N
  \left(\!\!\!\!
  \begin{tikzcd}[column sep=small]
    (
      \TopologicalSpace
      \rightrightarrows
      \TopologicalSpace
    )
    \ar[
      r,
      hook
    ]
    &
    (
      \TopologicalSpace \times \Gamma^{\mathrm{op}}
      \rightrightarrows
      \TopologicalSpace
    )
    \ar[
      r,
      ->>
    ]
    &
    (
      \Gamma^{\mathrm{op}}
      \rightrightarrows
      \ast
    )
  \end{tikzcd}
 \!\!\!\! \right)
  \;\;
  =
  \;\;
  \big( \!\!\!
  \begin{tikzcd}[column sep=small]
    \TopologicalSpace
    \ar[
      r,
      hook
    ]
    &
    \TopologicalSpace \times \Gamma^{\times_{\bullet}}
    \ar[
      r,
      ->>
    ]
    &
    \Gamma^{\times_{\bullet}}
  \end{tikzcd}
\!\!\!  \big).
\end{equation}

 \vspace{-2mm}
\noindent Therefore, the claim follows by
Prop. \ref{SufficientConditionsForTopologicalRealizationOfSimplicialSpacesToPreserveHomotopyPullbacks}
as soon as we verify the following conditions:

{(i)} All simplicial spaces in \eqref{SimplicialSpaceAvatarOfTopologicalBorelConstructionFiberSequence}
 are good (Def. \ref{GoodSimplicialTopologicalSpace}).

{(ii)} For each $n \in \mathbb{N}$ we have that
  $\TopologicalSpace
      \xhookrightarrow{\;}
   \TopologicalSpace \times \Gamma^{\times_n}
      \xrightarrow{\;}
  \Gamma^{\times}_n$
  is a homotopy fiber sequence (Ntn. \ref{HomtopyFiberSequenceOfTopologicalSpaces}).

{(iii)} The morphism on the right of
 \eqref{SimplicialSpaceAvatarOfTopologicalBorelConstructionFiberSequence}
 satisfies the homotopy Kan fibration property
 \eqref{HomotopyKanFibrationProperty}.


\noindent
Condition (i) follows by Prop. \ref{NervesOfActionGroupoidsOfWellPointedTopologicalGroupActionsAreGood}.
Condition (ii) is immediate from the fact that the projections
$\TopologicalSpace \times \Gamma^{\times_n} \xrightarrow{\;} \Gamma^{\times_n}$
are Serre fibrations (by Lem. \ref{LocallyTrivialBundlesAreSerreFibrations})
with ordinary fiber $\TopologicalSpace$.
Condition (iii) follows by Ex. \ref{BaseMorphismsOutOfNervesOfActionGroupoidsAreHomotopyKanFibrations}.
\end{proof}

\begin{proposition}[Milgram bundle is homotopy fiber sequence]
  \label{MilgramBundleIsHomotopyFiberSequence}
  If $\Gamma \,\in\, \Groups(\kTopologicalSpaces)_{\wellpointed}$ is
  a well-pointed topological group (Def. \ref{WellPointedTopologicalGroup}),
  then
  the sequence \eqref{QuotientCoprojectionOfUniversalPrincipalBundle}
    \vspace{-3mm}
  $$
    \begin{tikzcd}
      \Gamma
      \ar[
        rr,
        "{
          \in \, \HomotopyFiber(q)
        }"{below}
      ]
      &&
      E \Gamma
      \ar[
        r,
        "q"
      ]
      &
      B \Gamma
    \end{tikzcd}
  $$

    \vspace{-3mm}
\noindent
  is a homotopy fiber sequence (Ntn. \ref{HomtopyFiberSequenceOfTopologicalSpaces}).
\end{proposition}
\begin{proof}
  This is Lem. \ref{HomotopyFibrationOfBorelConstruction}
  for the case when $\TopologicalSpace \,=\, \Gamma$.
\end{proof}
\begin{lemma}[Borel construction for free locally trivial actions]
  \label{BorelConstructionForFreeAndLocallyTrivialActions}
  If  $\Gamma \in \Groups(\kTopologicalSpaces)$
  and
  $\Gamma \acts \, \TopologicalSpace \,\in \Actions{G}(\kTopologicalSpaces)$
  are such that

 \noindent {\bf (i)} both underlying spaces are connected, $\pi_0 = \ast$,

  \noindent {\bf (ii)} The coprojection $\TopologicalSpace \xrightarrow{q} \TopologicalSpace/G$
   is a locally trivial $G$-fiber bundle.

  \noindent
  Then the Borel construction \eqref{BorelConstructionAsTopologicalRealizationOfActionGroupoid}
  is weakly homotopy equivalent to the plain quotient space:
  \vspace{-2mm}
  $$
    (
      \TopologicalSpace \times E \Gamma
    )/\Gamma
    \xrightarrow{ \;\;\,\in\, \WeakHomotopyEquivalences \;\;}
    \TopologicalSpace/\Gamma
    \;\;\;\;\;\;
    \in
    \;\;
    \kTopologicalSpaces_{\mathrm{Qu}}
    \,.
  $$
\end{lemma}
\begin{proof}
  Observe that any local trivialization of $\TopologicalSpace \xrightarrow{q} \TopologicalSpace/G$
  lifts to a local trivialization of
  $\TopologicalSpace \times E G \xrightarrow{\;} (\TopologicalSpace \times  E G)/G$
  through the morphism that projects out the $E G$-factor.
  This makes the resulting commuting diagram a morphism of Serre fiber sequences
  and hence of homotopy fiber sequences (Ntn. \ref{HomtopyFiberSequenceOfTopologicalSpaces}).
  The statement now follows
  from the five-lemma (in its generality for possibly non-abelian groups)
  applied to the resulting morphism of
  long exact sequences of homotopy groups.
\end{proof}

\begin{lemma}[Classifying spaces of realizations of homotopy fibers of
simplicial topological groups]
 \label{ClassifyingSpacesOfRealizationsOfHomotopyFibersOfSimplicialTopologicalGroups}
  If
\vspace{-2mm}
  \begin{equation}
    \label{AHomotopyFiberSequenceOfSimplicialTopologicalGroups}
    \begin{tikzcd}
      \mathcal{H}_\bullet
      \ar[
        rr,
        " i_\bullet "{above},
        " \simeq \, \HomotopyFiber(p_\bullet) "{below}
      ]
      &&
      \widehat{\mathcal{G}}_\bullet
      \ar[
        rr,
        "p_\bullet"{above},
        "\in \; \HomotopyKanFibrations"{below}
      ]
      &&
      \mathcal{G}_{\bullet}
    \end{tikzcd}
    \;\;\;\;\;\;
    \in
    \;
    \Groups
    (
      \SimplicialTopologicalSpaces
    )_{\wellpointed}
  \end{equation}

  \vspace{-2mm}
  \noindent is a sequence of morphisms
  simplicial topological groups such that, as indicated:

 \noindent { \bf (i)} all three simplicial groups are well-pointed (Ntn. \ref{WellPointedTopologicalGroup});

\noindent { \bf (ii)}  the underlying sequence of simplicial spaces is
  a homotopy fiber sequence (Ntn. \ref{HomtopyFiberSequenceOfTopologicalSpaces});

\noindent { \bf (iii)}  the underlying morphism $p_\bullet$ is a homotopy Kan fibration
   (Def. \ref{HomotopyKanFibrations});

  \noindent
  and in addition

\noindent { \bf (iv)}
  under topological realization \eqref{TopologicalRealizationOfSimplicialTopologicalSpaces},
  $\vert p_\bullet \vert $ is surjective on connected components;

  \noindent
  then the image
  of \eqref{AHomotopyFiberSequenceOfSimplicialTopologicalGroups}
  under passage to
  classifying spaces $B(-)$ \eqref{QuotientCoprojectionOfUniversalPrincipalBundle}
  of topologically realized groups $\vert - \vert$
  \eqref{TopologicalRealizationOfSimplicialTopologicalSpaces}
  is still a homotopy fiber sequence
  \eqref{CharacterizationOfHomotopyKanFibration},
  now of topological classifying spaces:
    \vspace{-2mm}
  \begin{equation}
    \label{HomotopyFiberSequenceOfClassifyingSpacesOfRealizationsOfWellPointedSImplicialTopologicalGroups}
    \begin{tikzcd}
      B
      \vert \mathcal{H}_\bullet \vert
      \ar[
        rr,
        " B \vert i_\bullet \vert "{above},
        "{
          \scalebox{.9}{$
            \simeq
              \,
              \HomotopyFiber
              \scalebox{1.2}{$($}
                B \vert p_\bullet \vert
              \scalebox{1.1}{$)$}
          $}
        }"{below}
      ]
      &&
      B \vert \widehat{\mathcal{G}}_\bullet \vert
      \ar[
        rr,
        "B \vert p_\bullet \vert"
      ]
      &&
      B \vert \mathcal{G}_{\bullet} \vert
    \end{tikzcd}
    \;\;\;\;\;\;\;
    \in
    \;
    \kTopologicalSpaces
    \,.
  \end{equation}
\end{lemma}
\begin{proof}
  First, with assumption (i),
  Lem. \ref{WellPointedSimplicialGroupsAreGoodSimplicialSpaces}
  implies that
  the underlying simplicial topological spaces in
  \eqref{AHomotopyFiberSequenceOfSimplicialTopologicalGroups}
  are good (Def. \ref{GoodSimplicialTopologicalSpace}).
  Together with assumptions (ii) and (iii),
  this implies, by Prop. \ref{SufficientConditionsForTopologicalRealizationOfSimplicialSpacesToPreserveHomotopyPullbacks},
  that the topological realization
  of \eqref{HomotopyFiberSequenceOfClassifyingSpacesOfRealizationsOfWellPointedSImplicialTopologicalGroups}
  is a homotopy fiber sequence (Ntn. \ref{HomtopyFiberSequenceOfTopologicalSpaces}):
    \vspace{-2mm}
  \begin{equation}
    \label{SequenceOfTopologicalGroupsRealizingSequenceOfSimplicialTopologicalGroups}
    \begin{tikzcd}
      \vert
        \mathcal{H}_\bullet
      \vert
      \ar[
        rr,
        " \vert i_\bullet \vert "{above},
        " \simeq \; \HomotopyFiber (\vert p_\bullet\vert) "{below}
      ]
      &&
      \vert
        \widehat{\mathcal{G}}_\bullet
      \vert
      \ar[
        rr,
        "\vert p_\bullet \vert"
      ]
      &&
      \vert \mathcal{G}_{\bullet} \vert
    \end{tikzcd}
    \;\;\;\;\;\;
    \in
    \;
    \Groups
    (
      \kTopologicalSpaces
    )_{\wellpointed}
    \,,
  \end{equation}

    \vspace{-2mm}
  \noindent
  all of whose entries are well-pointed simplicial topological groups,
  by  Prop. \ref{RealizationOfWellPointedSimplicialGroupIsWellPointed}
  with assumption (i).

  Next, by condition (iv) with Ex. \ref{NervesOfMorphismsOfOfDeloopingGroupoids},
  this implies that the
  image of \eqref{SequenceOfTopologicalGroupsRealizingSequenceOfSimplicialTopologicalGroups}
  under forming nerves \eqref{SimplicialTopologicalNerveOfTopologicalGroupoids}
  of delooping groupoids (Ex. \ref{TopologicalDeloopingGroupoid})
  is the homotopy fiber sequence
  of a homotopy Kan fibration (Def. \ref{HomotopyKanFibrations}):
    \vspace{-2mm}
  \begin{equation}
    \label{NerveOfDeloopingOfSequenceOfTopologicalGroupsRealizingSequenceOfSimplicialTopologicalGroups}
    \begin{tikzcd}[column sep=25pt]
      N
      (
      \vert
        \mathcal{H}_\bullet
      \vert
        \rightrightarrows
      1
      )
      \ar[
        rr,
        "{
          \scalebox{.8}{$
            N
            \scalebox{1.2}{$($}
              \vert i_\bullet \vert \,\rightrightarrows\, 1
            \scalebox{1.2}{$)$}
          $}
        }"{above},
        "{
          \scalebox{.8}{$
            \in \; \HomotopyFiber
            \scalebox{1.2}{$($}
              \vert p_\bullet\vert
            \scalebox{1.2}{$)$}
          $}
        }"{below}
      ]
      &&
      N
      (
        \vert
          \widehat{\mathcal{G}}_\bullet
        \vert
        \rightrightarrows
        1
      )
      \ar[
        rr,
        "{
          \scalebox{.8}{$
            N
            \scalebox{1.2}{$($}
              \vert p_\bullet \vert
            \scalebox{1.2}{$)$}
          $}
        }"{above},
        "{
          \in
          \;
          \HomotopyKanFibrations
        }"{below}
      ]
      &&
      N
      (
        \vert \mathcal{G}_{\bullet} \vert
        \rightrightarrows
        1
      )
    \end{tikzcd}
    \;\;\;\;\;\;
    \in
    \;
    \SimplicialTopologicalSpaces_{\good}
    \,,
  \end{equation}

    \vspace{-2mm}
  \noindent
  all of whose entries
  are good simplicial spaces (Def. \ref{GoodSimplicialTopologicalSpace}),
  by Prop. \ref{NervesOfActionGroupoidsOfWellPointedTopologicalGroupActionsAreGood}.

  Therefore Prop.
  \ref{SufficientConditionsForTopologicalRealizationOfSimplicialSpacesToPreserveHomotopyPullbacks}
  applies also to
  \eqref{NerveOfDeloopingOfSequenceOfTopologicalGroupsRealizingSequenceOfSimplicialTopologicalGroups}
  and thus yields the claim
  \eqref{HomotopyFiberSequenceOfClassifyingSpacesOfRealizationsOfWellPointedSImplicialTopologicalGroups}.
\end{proof}

\begin{proposition}[2-Cocycle from topological central extension]
  \label{TopologicalTwoCocycleFromTopologicalCentralCentralExtension}
  Let
  $\!\!\begin{tikzcd}[column sep=small]
    A
    \ar[r,hook,"i"]
    &
    \widehat \Gamma
    \ar[r,->>,"p"]
    &
    \Gamma
  \end{tikzcd}\!\!$
  be a central extension
  \eqref{CentralExtensionOfTopologicalGroups}
  of topological groups
  such that

 \noindent {\bf (i)} $\widehat \Gamma \xrightarrow{\;} \Gamma$ is a locally trivial $A$-fiber bundle;

 \noindent {\bf (ii)} all groups are well-pointed (Ntn. \ref{WellPointedTopologicalGroup});

 \noindent {\bf (iii)} $A$ and $\widehat \Gamma$ are connected topological groups.

  \noindent
  Then there is a homotopy fiber sequence (Ntn. \ref{HomtopyFiberSequenceOfTopologicalSpaces}) of
  classifying spaces \eqref{QuotientCoprojectionOfUniversalPrincipalBundle}
  and higher classifying spaces (Eq. \ref{HigherClassifyingSpaces})
  of this form:
  \vspace{-3mm}
  $$
    \begin{tikzcd}
      B
      \widehat \Gamma
      \ar[
        rr,
        "{
          \simeq \, \HomotopyFiber(c)
        }"{below}
      ]
      &&
      B \Gamma
      \ar[
        r,
        "c"
      ]
      &
      B^2 A
    \end{tikzcd}
    \;\;\;
    \in
    \;
    \HomotopyCategory
    (
      \kTopologicalSpaces_{\mathrm{Qu}}
    )
    \,.
  $$
\end{proposition}
\begin{proof}
Consider the following image under topological realization
\eqref{TopologicalRealizationOfTopologicalGroupoids}
of the diagram \eqref{LongDuagramOf2GroupFromCentralExtensionOfGroups}
of topological 2-groups from Ex. \ref{Strict2GroupsFromCentralExtensions}:
\vspace{-2mm}
$$
  \begin{tikzcd}[column sep=large]
    \vert
      \widehat \Gamma
      \rightrightarrows
      \widehat \Gamma
    \vert
    \ar[
      rr,
      "{
        \vert
          (\mathrm{id},\mathrm{e}) \,\rightrightarrows\, \mathrm{id}
        \vert
      }"{above},
      "{
        \simeq \, \HomotopyFiber(\vert \mathrm{pr}_2 \rightrightarrows 1 \vert)
      }"{below}
    ]
    \ar[
      d,-,
      shift left=1pt
    ]
    \ar[
      d,-,
      shift right=1pt
    ]
    &&[17pt]
    \big\vert
      \widehat \Gamma \times A
      \rightrightarrows
      \widehat \Gamma
    \big\vert
    \ar[
      rr,
      "{
        \vert \mathrm{pr}_2 \rightrightarrows 1 \vert
      }"
    ]
    \ar[
      d,
      "{\vert p\circ \mathrm{pr}_1 \,\rightrightarrows\, p \vert}"{left},
      "\in \, \mathrm{W}"{right}
    ]
    &&
    \vert
      A
      \rightrightarrows
      \ast
    \vert
    \ar[
      d,-,
      shift left=1pt
    ]
    \ar[
      d,-,
      shift right=1pt
    ]
    \\
    \widehat \Gamma
    \ar[rr]
    &&
    \Gamma
    &&
    B A
  \end{tikzcd}
  \;\;\;\;\;
  \in
  \;
  \Groups
  (
    \kTopologicalSpaces
  )
  \,,
$$

\vspace{-1mm}
\noindent where in the bottom row we have used
\eqref{TopologicalRealizationOfConstantGroupoidIsomorphicToSPaceOfObjects}
and
\eqref{QuotientCoprojectionOfUniversalPrincipalBundle}.
Observe that, as indicated:

  {\bf (i)} the top sequence is a homotopy fiber sequence (Ntn. \ref{HomtopyFiberSequenceOfTopologicalSpaces}),
 by Lem. \ref{HomotopyFibrationOfBorelConstruction};

 {\bf (ii)} the middle morphism is a weak homotopy equivalence, by Lem. \ref{BorelConstructionForFreeAndLocallyTrivialActions}.

\noindent
It follows by Lem. \ref{ClassifyingSpacesOfRealizationsOfHomotopyFibersOfSimplicialTopologicalGroups}
that the same two properties are enjoyed by
the further image of this diagram under passage to classifying spaces:
\vspace{-1mm}
$$
  \begin{tikzcd}[column sep=large]
    B
    \vert
      \widehat \Gamma
      \rightrightarrows
      \widehat \Gamma
    \vert
    \ar[
      rr,
      "{
        B
        \vert
          (\mathrm{id},\mathrm{e}) \,\rightrightarrows\, \mathrm{id}
        \vert
      }"{above},
      "{
        \simeq \, \HomotopyFiber(B\vert \mathrm{pr}_2 \rightrightarrows 1 \vert)
      }"{below}
    ]
    \ar[
      d,-,
      shift left=1pt
    ]
    \ar[
      d,-,
      shift right=1pt
    ]
    &&[17pt]
    B
    \left\vert
      \widehat \Gamma \times A
      \rightrightarrows
      \widehat \Gamma
    \right\vert
    \ar[
      rr,
      "{
        B\vert \mathrm{pr}_2 \rightrightarrows 1 \vert
      }"
    ]
    \ar[
      d,
      "{B\vert p\circ \mathrm{pr}_1 \,\rightrightarrows\, p \vert}"{left},
      "\in \, \mathrm{W}"{right}
    ]
    &&
    B
    \vert
      A
      \rightrightarrows
      \ast
    \vert
    \ar[
      d,-,
      shift left=1pt
    ]
    \ar[
      d,-,
      shift right=1pt
    ]
    \\
    B
    \widehat \Gamma
    \ar[rr]
    &&
    B
    \Gamma
    \ar[
      rr,
      dashed,
      "c"{above}
    ]
    &&
    B^2 A
  \end{tikzcd}
  \;\;\;\;
  \in
  \HomotopyCategory
  (
    \kTopologicalSpaces_{\mathrm{Qu}}
  )
  \,,
$$

\vspace{-1mm}
\noindent
where the dashed morphism is the unique morphism
making the right square commute in the homotopy category.
Now the commutativity of the left square implies the claim.
\end{proof}
\begin{example}[Homotopy fiber sequence of the circle group]
Regarding the real line with
its additive group structure and
its Euclidean topology as a topological group
$$
  \mathbb{R}
  \;\;\;
  \in
  \;
  \Groups(\kTopologicalSpaces)
  \,,
$$
consider the defining sequence of the circle group $\CircleGroup$
$$
  \begin{tikzcd}
    \mathbb{Z}
    \ar[r, hook]
    &
    \mathbb{R}
    \ar[r, ->>]
    &
    \CircleGroup
    \,.
  \end{tikzcd}
$$
Under passage to higher classifying spaces (Ex. \ref{HigherClassifyingSpaces})
and using that topological realization preserves finite limits
(Lem. \ref{TopologicalRealizationPreservesFiniteLimits})
and hence fibers, this induces for each $n \,\in\, \mathbb{N}$
a fiber sequence of topological groups
\begin{equation}
  \label{FiberSequenceOfHigherCircleClassifyingSpaces}
  \begin{tikzcd}
    B^n \mathbb{Z}
    \ar[r, hook]
    &
    B^n \mathbb{R}
    \ar[r, ->>]
    &
    B^n \CircleGroup
    \,.
  \end{tikzcd}
\end{equation}
Prop. \ref{TopologicalTwoCocycleFromTopologicalCentralCentralExtension}
says that there is a continuous map $B^n \CircleGroup \xrightarrow{\;c\;} B^{n+1} \mathbb{Z}$
which extends this to a homotopy fiber sequence of the form
\begin{equation}
  \label{HomotopyFiberSequenceOfHigherCircleClassifyingSpaces}
  \begin{tikzcd}
    B^n \mathbb{Z}
    \ar[r]
    &
    B^n \mathbb{R}
    \ar[r, "\mathrm{hofib}(c)"]
    &
    B^n \CircleGroup
    \ar[r, "c", "\sim"{swap}]
    &
    B^{n+1} \mathbb{Z}
  \end{tikzcd}
  \;\;\;
  \in
  \;
  \HomotopyCategory(
    \kTopologicalSpaces_{\mathrm{Qu}}
  )
  \,.
\end{equation}
where on the left we are using \eqref{GenericHomotopyFiberSequence}
in order to identify the homotopy fiber of the homotopy fiber of $c$.
But $\mathbb{R} \underset{\mathrm{htpy}}{\simeq} \ast$
implies $B^n \mathbb{R} \underset{\mathrm{htpy}}{\simeq} B^n 1 \,\simeq\, \ast$
for all $n \in \mathbb{N}$,
from which the long exact sequence of homotopy groups induced by
\eqref{FiberSequenceOfHigherCircleClassifyingSpaces} implies that
the map $c$ is a weak homotopy equivalence, as indicated.
\end{example}

\begin{example}[Homotopy fiber sequence of the projective unitary group]
  \label{HomotopyFiberSequenceOfTheProjectiveUnitaryGroup}
  From the defining fiber sequence \eqref{PUHFiberSequence} of
  projective unitary group from Ex. \ref{ProjectiveUnitarGroupOnAHilbertSpace},
  Prop. \ref{TopologicalTwoCocycleFromTopologicalCentralCentralExtension} implies
  a homotopy fiber sequence of classifying spaces of the form
  $$
    \begin{tikzcd}
      B \UH
      \ar[rr]
      &&
      B \PUH
      \ar[r, "c", "\sim"{swap}]
      &
      B^2 \CircleGroup
      \ar[r, "\sim"{swap}]
      &
      B^3 \mathbb{Z}
    \end{tikzcd}
    \;\;\;
    \in
    \;
    \HomotopyCategory
    (
      \kTopologicalSpaces_{\mathrm{Qu}}
    )
    \,,
  $$
  where on the right we have appended the equivalence
  \eqref{HomotopyFiberSequenceOfHigherCircleClassifyingSpaces}.
  But since $\UH \underset{\mathrm{htpy}}{\simeq} \ast$
  \eqref{UHIsContractible} also $B \UH \underset{\mathrm{htpy}}{\simeq}$,
  whence the corresponding long exact sequence of homotopy groups implies
  that $c$ is an equivalence, as indicated.

  In summary, the combined
  2-cocycle maps from Prop. \ref{TopologicalTwoCocycleFromTopologicalCentralCentralExtension}
  witness that the classifying space of the infinite complex projective group
  is an Eilenberg-MacLane space concentrated in degree 3:
  \begin{equation}
    \label{ClassifyingSpaceOfPUHIsEMSpace}
    B \PUH \,\simeq\, B^3 \mathbb{Z} \;=\; K(\mathbb{Z},3)
    \;\;\;
    \in
    \;
    \HomotopyCategory
    (
      \kTopologicalSpaces_{\mathrm{Qu}}
    )
    \,.
  \end{equation}
  We discuss a stacky refinement of this situation below in Ex. \ref{ProjectiveRepresentationsAndTheirCentralExtensions}.
\end{example}

\newpage

\part{In cohesive $\infty$-stacks}
 \label{InCohesiveInfinityStacks}

\chapter{Equivariant $\infty$-topos theory}
\label{Generalities}

We recall and develop some background in modal
and specifically in cohesive and globally equivariant $\infty$-topos theory
that are needed below in \cref{EquivariantInfinityBundles}.
Some of the following material is recalled from
\cite{Lurie09HTT}\cite{dcct} \cite{SS20OrbifoldCohomology} \cite{FSS20CharacterMap},
but several statements are new or at least not readily quotable from the literature.
Nevertheless, the expert reader may want to skip this chapter and
just refer back to it as need be.

\medskip

-- \cref{AbstractHomotopyTheory}:  Abstract homotopy theory (presentable $\infty$-categories).

-- \cref{ToposTheory}: Geometric homotopy theory (general $\infty$-toposes).

-- \cref{CohesiveHomotopyTheory}: Cohesive homotopy theory (cohesive $\infty$-toposes).

\section{Abstract homotopy theory}
\label{AbstractHomotopyTheory}

We record and develop some facts of the
abstract homotopy theory
of presentable $\infty$-categories
in the guise of
combinatorial simplicial model category theory.
Then we highlight some aspects of
the homotopy theory of simplicial transformation groups.

\medskip

-- \cref{OnModelCategories}: The homotopy 2-category of presentable $\infty$-categories.

-- \cref{TransformationInfinityGroups}: Homotopy theory of simplicial transformation groups.

\subsection{The homotopy 2-category of presentable $\infty$-categories}
 \label{OnModelCategories}

Besides the original
\cite{Quillen67},  standard textbook accounts of model category theory include
\cite{Hovey99}\cite{Hirschhorn02}\cite[A.2]{Lurie09HTT}.
For review streamlined towards our context and applications see \cite[\S A]{FSS20CharacterMap}.

\medskip
We briefly put this in perspective:

\medskip

\noindent
{\bf The idea of abstract homotopy theory.}
An {\it abstract homotopy theory} is
meant to be a kind of category which adheres to
what we may call the
{\it higher gauge principle},
where
instead of asking if any two
parallel morphisms are equal one must ask
whether there is a {\it gauge transformation} between
them -- here called a {\it homotopy} -- and then for any
two such homotopies one is to ask for {\it higher}
gauge transformations, namely higher homotopies
between these;  and so forth. See diagrams \eqref{GaugeTransformations}
and \eqref{HomSpace} in \cref{ToolsAndTechniques}.

\medskip

\noindent
{\bf Model categories.}
This informal idea had famously been formalized by the notion of
{\it model categories} \cite{Quillen67},
among which the ``combinatorial''
ones have come to be understood
(see Rem. \ref{HomotopyCategoryOfPresInfinityCategoriesIsThatOfCombinatorialModelCategories} below)
as presenting exactly the
``presentable'' homotopy theories \cite{Lurie09HTT},
essentially meaning: those
that have any relevance in practice.
Since model categories {\it present}
homotopy theories in terms of ordinary categories equipped with
an ingenious construction principle for higher homotopies,
they have proven to serve as useful toolboxes for
operationally handling abstract homotopy in practice.
In fact, those combinatorial model categories which are
also locally cartesian closed (in an evident model category-theoretic
sense) are  so very well-adapted to this task that
their presentation of homotopy theories can be entirely mechanized
via a kind of programming language now known as {\it homotopy type theory}
(\cite{UFP13}\cite{Shulman19}).

\medskip

\noindent
{\bf Homotopy theory of homotopy theories.}
But in the all-encompassing spirit of the gauge principle,
there ought to be a {\it homotopy theory of homotopy theories},
necessary for making gauge-invariant sense of crucial notions such as that
reflecting upon sub-homotopy theories of objects
with special qualitative properties.
By analogy (in fact: by de-categorification), ordinary categories form a (``very large'')
{\it 2-category} $\Categories$ which provides the
formal context (\cite{Gray74})
for making invariant sense of notions
such as that of reflective subcategories.

\medskip
Accordingly, it eventually came to
light that homotopy theories may be understood as
homotopy-theoretic categories, now technically known as
{\it $(\infty,1)$-categories} (or {\it $\infty$-categories}, for short,
see \cite{Lurie09HTT}\cite{Cisinski19}\cite{Riehl19}\cite{Land21} for introduction),
whose {\it homotopy theory of homotopy theories}
is a (``very large'') homotopical 2-category
$\InfinityCategories$, hence an ``$(\infty,2)$-category''.
While it transpired that the resulting $(\infty,1)$-category theory formally behaves
much along the lines familiar from ordinary category theory,
thus providing the
expert with a powerful way of effectively reasoning inside it,
the sheer volume of technical detail needed to rigorously set up
the $(\infty,2)$-category of $(\infty,1)$-categories, in any of its multitude
of interdependent models, can be prohibitive.

\medskip

\noindent
{\bf The homotopy 2-category of homotopy theories.}
However, as highlighted in \cite{RiehlVerity13},
this is often more than is
necessary in practice, as much of the
relevant formal theory of homotopy theories is reflected already within
just the {\it homotopy 2-category} $\HomotopyTwoCategory(\InfinityCategories)$
that is obtained from the $(\infty,2)$-category by
identifying gauge equivalent 2-morphisms \cite{RiehlVerity16}\cite{RiehlVerity21}.
In particular, adjunctions \eqref{AdjunctionAndHomEquivalence}
between homotopy theories, such as
the above-mentioned reflections onto sub-theories,
may be fully understood already in $\HomotopyTwoCategory(\InfinityCategories)$
\cite[\S 4.4.5]{RiehlVerity13}.

\medskip
\noindent
{\bf The homotopy 2-category of combinatorial model categories.}
We observe here that this opens the door to short-cutting the
rigorous construction of
$\infty$-category theory in a dramatic way, if
the homotopy 2-category of homotopy theories has a usefully direct
construction from model category theory.
This is indeed the case, at least for the practically relevant
{\it presentable} homotopy theories
(Prop. \ref{HomotopyCategoryOfPresInfinityCategoriesIsThatOfCombinatorialModelCategories} below).
The homotopy 2-category of presentable homotopy theories
turns out
to be simply the {\it 2-localization} of the 2-category of
combinatorial model categories at the Quillen equivalences,
i.e. the 2-category obtained from that of combinatorial model categories
by universally forcing its Quillen equivalences to become actual
homotopy equivalences.
This is a technical theorem
(recently proven by \cite{Pavlov21}, inspired by \cite{Renaudin06})
if one starts with existing definitions
of $(\infty,1)$-categories.
But we may thus take it to be the very definition
of what might be called the {\it formal category theory of homotopy theories}:
\vspace{-2mm}
$$
  \overset{
    \mathclap{
    \raisebox{7pt}{
      \tiny
      \color{darkblue}
      \bf
      \begin{tabular}{c}
        homotopy 2-category of
        \\
        presentable homotopy theories
      \end{tabular}
    }
    }
  }{
    \HomotopyTwoCategory
    (\PresentableInfinityCategories)
  }
  \;\;
  \coloneqq
  \;\;
  \overset{
    \mathclap{
    \raisebox{7pt}{
      \tiny
      \color{darkblue}
      \bf
      \begin{tabular}{c}
        2-localization 2-category of
        \\
        combinatorial model categories
      \end{tabular}
    }
    }
  }{
    \TwoLocalization{\Quillen}
    (\CombinatorialModelCategories)
  }
  \,.
$$
This definition turns out to give a powerful but convenient handle
on abstract homotopy theory.

\medskip

\begin{notation}[Category of simplicial sets]
\label{SimplicialSets}
  We write $\Delta \,\in\, \Categories$ for the simplex category,
  and
  $
    \SimplicialSets
    \;\coloneqq\;
    \Presheaves(\Delta,\, \Sets)
  $
  for the category of simplicial sets
  (see e.g. \cite{Friedman08}\cite[\S 1]{GoerssJardine09}).
  For $n \,\in\, \mathbb{N}$ we write $[n] \,\in\, \Delta$
  and then $\Delta[n] \,\coloneqq\, y([n]) \,\in\, \SimplicialSets$
  for the objects of the simplex category, regarded as the standard simplicial
  simplices.
\end{notation}
\begin{remark}[Cartesian closed structure on simplicial sets]
  \label{CartesianClosedStructureOnSimplicialSets}
  $\SimplicialSets$ being a category of presheaves means that the product of simplicial
  sets exists and is given for $X ,\, Y \,\in\, \SimplicialSets$ by
  \vspace{-2mm}
  $$
    (X \times Y)_n \;\simeq\; X_n \times Y_n
    \,.
  $$
\end{remark}

\begin{notation}[Simplicial categories]
  \label{SimplicialCategories}
  $\,$

  \noindent
 {\bf (i)}  All categories in the following are simplicial categories
  (\cite[\S 9.1.1]{Hirschhorn02}),
  i.e., enriched (\cite{Kelly82}\cite[\S 6]{Borceux94I}) in
  the category of simplicial sets (Ntn. \ref{SimplicialSets})
  with respect to its cartesian monoidal structure (Ntn. \ref{CartesianClosedStructureOnSimplicialSets}).

  We write $\SimplicialCategories$ for the large category of simplicial categories.

  \noindent {\bf (ii)}
  This means, in particular, that for
  any pair of objects of some $\mathcal{C} \,\in\, \SimplicialCategories$,
  we have their {\it simplicial hom-complex}, to be denoted
  \begin{equation}
    \label{SimplicialHomComplex}
     X, A \;\in\; \mathcal{C}
     \;\;\;
     \vdash
     \;\;\;
    \mathcal{C}(X,A)
    \;\in\;
    \SimplicialSets
    \,.
  \end{equation}

 \noindent {\bf (iii)}  For $\mathcal{C} \in \SimplicialCategories$, its
{\it homotopy category} $\HomotopyCategory(\mathcal{C})$
is the 1-category with the same objects and with hom-sets the connected components
of the hom-complexes \eqref{SimplicialHomComplex}:
\vspace{-2mm}
\begin{equation}
  \label{HomotopyCategoryOfASimplicialCategory}
  \HomotopyCategory(\mathcal{C})
  (
    X
    ,\,
    A
  )
  \;\coloneqq\;
  \pi_0
  \left(
    \mathcal{C}
    (
      X
      ,\,
      A
    )
  \right)
  \,.
\end{equation}
\end{notation}

\begin{notation}[2-Category of simplicial combinatorial model categories]
\label{CategoryOfCombinatorialModelCategories}
$\,$

\noindent
{\bf (i)} We write $\CombinatorialSimplicialModelCategories$
for the very large strict 2-category whose

\vspace{-3mm}
\begin{itemize}
 \setlength\itemsep{-2pt}

\item[\bf (a)]
objects are are combinatorial model categories
(\cite[Def. 2.1]{Dugger01Combinatorial}\cite[\S A.2.6]{Lurie09HTT}
review in \cite{Raptis14})
that are also
simplicial model categories
(\cite[\S II.2]{Quillen67}\cite[\S 9.1.5]{Hirschhorn02}\cite[Ex. A.3.1.4]{Lurie09HTT});

\item[\bf (b)]
morphisms are left Quillen functors, to be denoted
\vspace{-3mm}
$$
  \begin{tikzcd}
    \mathcal{C}
    \ar[
      rr,
      shift left=7pt,
      "L"{above}
    ]
    &&
    \mathcal{D}\;;
    \ar[
      ll,
      shift left=7pt,
      "\mathllap{\exists\,} R"{below}
    ]
    \ar[
      ll,
      phantom,
      "\scalebox{.7}{$\bot_{\mathrlap{\Quillen}}$}"
    ]
  \end{tikzcd}
$$

\vspace{-3mm}
\item[\bf (iii)]
2-morphisms are natural transformations between these left adjoints
\vspace{-2mm}
$$
  \begin{tikzcd}
    \mathcal{C}
    \ar[
      rr,
      bend left=25,
      "L"{above},
      "{ }"{name=s, below}
    ]
    \ar[
      rr,
      bend right=25,
      "L'"{below},
      "{ }"{name=t, above}
    ]
    &&
    \mathcal{D} \;.
    \ar[
      from=s,
      to=t,
      Rightarrow
    ]
  \end{tikzcd}
$$
\end{itemize}
\end{notation}


\begin{definition}[Homotopy 2-category of presentable homotopy theories]
\label{HomotopyTwoCategoryOfPresentableHomotopyTheories}
We write
\vspace{-2mm}
\begin{equation}
  \label{TheHomotopy2CategoryOfPresentableModelCategories}
  \begin{tikzcd}
    \SimplicialCombinatorialModelCategories
    \ar[
      rr,
      "\SimplicialLocalization{\WeakEquivalences}"
    ]
    &&
    \TwoLocalization{\QuillenEquivalences}
    \left(
      \CombinatorialSimplicialModelCategories
    \right)
  \end{tikzcd}
\end{equation}

\vspace{-2mm}
\noindent
for the
2-localization \cite[\S 2]{Pronk96}\cite[\S 1.2]{Renaudin06}
of the 2-category $\SimplicialCombinatorialModelCategories$
(Ntn. \ref{CategoryOfCombinatorialModelCategories})
at the class of left Quillen equivalences.
\end{definition}
This localization exists by
\cite[\S 2.3]{Renaudin06}.\footnote{Considered in \cite{Renaudin06} is the 2-localization of
the category of all combinatorial model categories,
but by \cite[Thm 1.1]{Dugger01Combinatorial}
every combinatorial model category is Quillen equivalent to
a (left proper) simplicial combinatorial model category,
so that the
above localization \eqref{TheHomotopy2CategoryOfPresentableModelCategories} exists equivalently.}

\begin{proposition}[Recovering the homotopy 2-category of presentable $\infty$-categories {\cite{Pavlov21}}]
  \label{HomotopyCategoryOfPresInfinityCategoriesIsThatOfCombinatorialModelCategories}
  The 2-localization of the 2-category of simplicial combinatorial model categories in
  Def. \ref{HomotopyTwoCategoryOfPresentableHomotopyTheories}
  is equivalent to
  the homotopy 2-category of
  the $(\infty,2)$-category $\PresentableInfinityCategories$
  of presentable $\infty$-categories with left adjoint $\infty$-functors,
  according to \cite[\S 5]{Lurie09HTT}:
  \begin{equation}
    \label{HomotopyTwoCategoryOfPresentableInfinityCategoires}
    \TwoLocalization{\QuillenEquivalences}
    (
      \SimplicialCombinatorialModelCategories
    )
    \;\;\simeq\;\;
    \HomotopyTwoCategory
    (
      \PresentableInfinityCategories
    )\;.
  \end{equation}
\end{proposition}
\begin{remark}
Prop. \ref{HomotopyCategoryOfPresInfinityCategoriesIsThatOfCombinatorialModelCategories}
is the culmination of several well-known results:

\vspace{1mm}
\noindent
{\bf (i)}
That presentable $\infty$-categories
in the sense of \cite[\S 5.5]{Lurie09HTT},
are, up to equivalence,
the simplicial localizations of combinatorial model categories is
\cite[Prop. A.3.7.6]{Lurie09HTT}\cite[Thm. 7.11.16, Rem. 7.11.17]{Cisinski19}.

\vspace{1mm}
\noindent
{\bf (ii)}
That the image under simplicial localization
of left Quillen functors between combinatorial model categories are
left adjoint $\infty$-functors between these $\infty$-categories
in the sense of \cite[\S 5.2]{Lurie09HTT},
is \cite{Mazel-Gee15}. More specifically, left Quillen functors between
left proper combinatorial model categories functorially lift
(by \cite[Prop. A.3]{BlumbergRiehl14}) to
simplicial Quillen adjunctions between Quillen-equivalent simplicial model categories,
and these in turn lift functorially (by \cite[Prop. 5.2.4.6]{Lurie09HTT})
to adjunctions between the corresponding quasi-categories.

\vspace{1mm}
\noindent
{\bf (iii)}
That every equivalence of presentable $\infty$-categories
arises from a composite of (such simplicial localizations of) simplicial Quillen adjunctions
between simplicial combinatorial model categories
is claimed in \cite[Rem. A.3.7.7]{Lurie09HTT}.
\end{remark}

\medskip

\noindent
{\bf $\infty$-Groupoids.}
In higher analogy to how $\Sets$ is the archetypical category,
the archetypical $\infty$-category is that of $\infty$-groupoids
(e.g. \cite[\S 1.1.2]{Lurie09HTT}):

\begin{notation}[Kan-Quillen model category of simplicial sets]
\label{KanQuillenModelCategoryOfSimplicialSets}
\noindent
The classical Kan-Quillen model category of simplicial sets
(\cite[\S II.3]{Quillen67}\cite{GoerssJardine09}),
regarded as a simplicial combinatorial model category
(Ntn. \ref{CategoryOfCombinatorialModelCategories}),
we denote
\vspace{-2mm}
\begin{equation}
  \label{SimplicialSetsAsAnEnrichedCategory}
  \SimplicialSets_{\Quillen}
    \;\in\;
  \SimplicialCombinatorialModelCategories
  \,.
\end{equation}
\end{notation}

\begin{remark}[Simplicial tensoring Quillen adjunction]
For every $\mathcal{C} \,\in\, \SimplicialCombinatorialModelCategories$
(Ntn. \ref{CategoryOfCombinatorialModelCategories})
and $X \,\in\, \mathcal{C}^{\mathrm{co}}$, we have a morphism
$\SimplicialSets_{\Quillen} \xrightarrow{ X } \mathcal{C}$ in
$\SimplicialCombinatorialModelCategories$,
namely the left Quillen functor
given by the tensoring (co-powering) of the simplicial model category
$\mathcal{C}$ over simplicial sets (Ntn. \ref{KanQuillenModelCategoryOfSimplicialSets}):
\vspace{-2mm}
$$
  \begin{tikzcd}
    \mathcal{C}
    \ar[
      rr,
      shift right=7pt,
      "{ \mathcal{C}(X,-) }"{below}
    ]
    &&
    \SimplicialSets_{\Quillen}\;.
    \ar[
      ll,
      shift right=7pt,
      "{ (-) \cdot X }"{above}
    ]
    \ar[
      ll,
      phantom,
      "\scalebox{.7}{$\bot_{\Quillen}$}"
    ]
  \end{tikzcd}
$$
\end{remark}

\begin{notation}[$\infty$-Groupoids]
  \label{SimplicialSetsAndInfinityGroupoids}
  We denote the homotopy theory (Def. \ref{HomotopyTwoCategoryOfPresentableHomotopyTheories})
  presented by the classical Kan-Quillen model category of simplicial sets
  (Ntn. \ref{KanQuillenModelCategoryOfSimplicialSets})
  by
  \vspace{-2mm}
  \begin{equation}
    \label{InfinityGroupoids}
    \InfinityGroupoids
      \;\coloneqq\;
    \TwoLocalizationProjection{\Quillen}
    (
      \SimplicialSets_{\Quillen}
    )
    \;\;
    \in
    \;
    \HomotopyTwoCategoryOfPresentableInfinityCategories
    \,.
  \end{equation}
\end{notation}

\begin{definition}[Hom $\infty$-groupoid]
  \label{HomInfinityGroupoid}
  For
  $
    \mathbf{C}
      \;\simeq\;
    \SimplicialLocalization{\WeakEquivalences}(\mathcal{C})
    \,\in\,
  \HomotopyTwoCategoryOfPresentableInfinityCategories$
  (Def. \ref{HomotopyCategoryOfPresInfinityCategoriesIsThatOfCombinatorialModelCategories})
  and for $X,A \,\in\, \mathcal{C}^{\mathrm{cof}}_{\mathrm{fib}}$
  two bifibrant objects, we write
  \vspace{-3mm}
  $$
    \mathbf{C}
    (X,A)
    \;\coloneqq\;
    \mathcal{C}(X,A)
    \;\;\;
    \in
    \;
    \SimplicialSets_{\mathrm{Qu}}
    \xrightarrow{\;\;
      \scalebox{.7}{$
        \eta^{\SimplicialLocalization{\mathrm{W}}}
      $}
    \;\;}
    \InfinityGroupoids
  $$

  \vspace{-2mm}
 \noindent for the image of their simplicial hom-complex,
 regarded as an $\infty$-groupoid via Ntn. \ref{SimplicialSetsAndInfinityGroupoids}.
\end{definition}
In fact, it suffices
that the domain be cofibrant and the codomain fibrant
to obtain the homotopy type of the hom-$\infty$-groupoid, as shown by the
following standard observation:

\begin{lemma}[Hom $\infty$-groupoid from cofibrant domains into fibrant codomains]
  \label{HomInfinityGroupoidFromCofibrantDomainAndFibrantCodomain}
  For $\mathcal{C} \,\in\, \SimplicialModelCategories$,
  let $X \in \mathcal{C}^{\mathrm{cof}}$ be a cofibrant object
  and $A \in \mathcal{C}_{\mathrm{fib}}$ be a fibrant object,
  and consider any fibrant replacement
  $\!\!
    \begin{tikzcd}
      X
      \ar[
        r,
        "q_X"{above},
        "\in \mathrm{Cof} \cap \mathrm{W}"{below}
      ]
      &
      P X
    \end{tikzcd}
  \!\!$
  and
  cofibrant replacement
  $\!\!
    \begin{tikzcd}
      Q A
      \ar[
        r,
        "p_A"{above},
        "\in \mathrm{Fib} \cap \mathrm{W}"{below}
      ]
      &
      A
    \end{tikzcd}
  \!\!$.
  Then the canonical comparison map
  \vspace{-3mm}
  $$
    \begin{tikzcd}
      \mathcal{C}(X,A)
      \ar[
        rr,
        "{
          \mathcal{C}(q_X, p_A)
        }",
        "\in \mathrm{W}"{below}
      ]
      &&
      \mathcal{C}
      (
        P X,
        \,
        Q A
      )
      \;=\;
      \mathbf{C}
      (
        P X, \, Q A
      )
      \;\in\;
      \InfinityGroupoids
    \end{tikzcd}
  $$

  \vspace{-1mm}
  \noindent is an equivalence to the hom-$\infty$-groupoid
  from Def. \ref{HomInfinityGroupoid}
\end{lemma}
\begin{proof}
  Since $p_X$ is an acylic cofibration and
  $q_A$ is an acylic fibration,
  the axioms)
  (e.g. \cite[Rem. A.3.6.1 (2')]{Lurie09HTT})
  of simplicially enriched model categories
  imply that
  both morphisms in the the following factorization are equivalences:
  \vspace{-3mm}
  $$
    \mathcal{C}(p_X, q_A)
    \;:\;
    \begin{tikzcd}[column sep=large]
      \mathcal{C}(X,A)
      \ar[
        rr,
        "{\mathcal{C}(\mathrm{id}_X, q_A)}"{above},
        "\in \mathrm{W} \cap \mathrm{Fib}"{below}
      ]
      &&
      \mathcal{C}(P X, Q A)
      \ar[
        rr,
        "{\mathcal{C}(q_X, \mathrm{id}_A)}"{above},
        "\in \mathrm{W} \cap \mathrm{Fib}"{below}
      ]
      &&
      \mathcal{C}(P X, QA)\;.
    \end{tikzcd}
  $$

  \vspace{-1mm}
 \noindent Hence the composite is an equivalence.
\end{proof}

\begin{proposition}[$\infty$-Groupoids form a cartesian closed $\infty$-category]
  \label{InfinityGroupoidsFormCartesianClosedInfinityCategory}
  For $S \,\in\, \InfinityGroupoids$ (Ntn. \ref{SimplicialSetsAndInfinityGroupoids}),
  there is a pair of adjoint $\infty$-functors \eqref{AdjunctionAndHomEquivalence}
  of the form
  \vspace{-3mm}
  \begin{equation}
    \label{HomAdjunctionInInfinityGroupoids}
    \begin{tikzcd}
      \InfinityGroupoids
      \ar[
        rr,
        phantom,
        "{\scalebox{.6}{$\bot$}}"
      ]
      \ar[
        rr,
        shift right=5pt,
        "{
         \scalebox{.7}{$ \InfinityGroupoids(S,-) $}
         }"{below}
      ]
      &&
     \InfinityGroupoids \;.
      \ar[
        ll,
        shift right=5pt,
        "{
        \scalebox{.7}{$   S \times (-) $}
        }"{above}
      ]
    \end{tikzcd}
  \end{equation}
\end{proposition}

\begin{notation}[Objects and equivalences in $\infty$-categories]
  For $\mathbf{C} = \SimplicialLocalization{\mathrm{W}}(\mathcal{C})$
  (Def. \ref{HomotopyTwoCategoryOfPresentableHomotopyTheories}),

 \vspace{-3mm}
\begin{enumerate}[{\bf (i)}]
\setlength\itemsep{-4pt}
\item  we say that the homotopy 2-category of $\mathbf{C}$ is
  \vspace{-1mm}
  $$
    \TwoHomotopyCategory(\mathbf{C})
    \;:=\;
    \TwoHomotopyCategory(\PresentableInfinityCategories)
    \left(
      \InfinityGroupoids, \mathbf{C}
    \right);
  $$
\item   we say that objects $X \,\in\, \mathcal{C}$
  are {\it presentations} for their images, up to isomorphism, in the
  homotopy category of the model category $\mathcal{C}$ which we
  may identify with the homotopy category of $\mathcal{C}$
  \vspace{-2mm}
  \begin{equation}
    \label{PresentingObjectInModelCategory}
    \begin{tikzcd}[row sep=0pt]
      X
      \,\in\,
      \mathcal{C}
      \ar[
        rr,
        "\Localization{\WeakEquivalences}"
      ]
      &&
      \HomotopyCategory(\mathcal{C})
      \,=:\,
      \HomotopyCategory(\mathbf{C}) \;.
    \end{tikzcd}
  \end{equation}
\end{enumerate}

\end{notation}

\begin{example}[Points of $\infty$-Groupoids]
\label{PointsOfInfinityGroupoids}
All $S \in \InfinityGroupoids$
(Ntn. \ref{SimplicialSetsAndInfinityGroupoids})
are equivalent to the hom-$\infty$-groupoid (Def. \ref{HomInfinityGroupoid})
into themselves out of the point:
\vspace{-2mm}
\begin{equation}
  \label{GroupoidsAreTheirOwnGroupoidsOfPoints}
  S \;\simeq\;
  \infinityGroupoids(\ast, S)
  \;\;\;
  \in
  \;
  \infinityGroupoids
  \,.
\end{equation}
\end{example}

\begin{example}[Monomorphisms of $\infty$-groupoids]
  \label{MonomorphismsOfInfinityGroupoids}
  A morphism of $\infty$-groupoids is a
  monomorphism \eqref{InfinityMonomorphism}
  precisely if it is fully faithful \eqref{FullyFaithfulInfinityFunctor},
  hence if it is, up to equivalence, the inclusion of a disjoint summand
  of connected components:
  \begin{equation}
    \label{FullyFaithfulMorphismOfInfinityGroupoidsCharacterized}
    S
    \xhookrightarrow{}
    S'
    \;\;\;\;\;\;\;\;\;
    \Leftrightarrow
    \;\;\;\;\;\;\;\;\;
    S' \,\simeq\, S \;\sqcup\; S' \setminus S
    \;\;\;\;\;\;\;\;\;\;
    \in
    \;
    \InfinityGroupoids.
  \end{equation}

  Moreover, by the homotopy-pullback definition
  \eqref{InfinityMonomorphism}
  of monomorphisms, a morphism
  $f$
  in any $\infty$-category
  $\mathbf{C}$ is a monomorphism iff
  the induced map on hom-$\infty$-groupoids
  out of any object
  is a monomorphism
  of $\infty$-groupoids \eqref{FullyFaithfulMorphismOfInfinityGroupoidsCharacterized}:
  $$
      A
      \xhookrightarrow{\;i\;}
      B
      \;
      \in
      \;
      \mathbf{C}
    \;\;\;\;\;\;\;\;\;\;\;\;\;
    \Leftrightarrow
    \;\;\;\;\;\;\;\;\;\;\;\;\;
    \underset{X \in \mathbf{C}}{\forall}
    \;\;
    \mathbf{C}(X,\,A)
    \xhookrightarrow{ \mathbf{C}(X,i) }
    \mathbf{C}(X,\,B)
    \;\;\;
    \in
    \;
    \InfinityGroupoids
    \,.
  $$
\end{example}

\begin{notation}[1-Groupoids]
  \label{GroupoidNotation}
  Let $\mathcal{C}$ be a category with finite products.
  For $G \in \Groups(\mathcal{C})$, we denote

\vspace{-.3cm}
  \begin{itemize}
 \setlength\itemsep{-3pt}
  \item
  its  {\it delooping groupoid} by:
  \begin{equation}
    \label{DeloopingGroupoid}
    ( G \rightrightarrows \ast)
    \;\in\;
    \Groupoids(\mathcal{C})\;;
  \end{equation}

  \item
  its {\it action groupoid} of right multiplication action on itself by:
  \vspace{-2mm}
  \begin{equation}
    \label{ActionGroupoidOfRightMultiplicationAction}
    \left( G \times G \rightrightarrows G \right)
    \;\coloneqq\;
    \Big(
      G \times G
      \underoverset
        {\scalebox{.5}{$(-)\cdot (-)$}}
        {\mathrm{pr}_1}
        {\rightrightarrows}
      G
    \Big)
    \;\in\;
    \Groupoids(\mathcal{C}) \;.
  \end{equation}

  \end{itemize}

  \vspace{-3mm}
  \noindent More generally, for $G \acts \, X \,\in\, \Actions{G}(\mathcal{C})$
  a left action
 \vspace{-2mm}
$$
  \begin{tikzcd}[row sep=-6pt]
    G \times X
    \ar[
      rr,
      "\mbox{
         \tiny
         \color{greenii}
         \bf
         left action
      }"
    ] && X
    \mathrlap{\,,}
    \\
  \scalebox{0.7}{$  (g,x)  $}
  &\longmapsto&
  \scalebox{0.7}{$ g \cdot x $}
  \end{tikzcd}
  {\phantom{AAAAAA}}
  \begin{tikzcd}[row sep=-6pt]
    X \times G
    \ar[
      rr,
      "\mbox{
         \tiny
         \color{greenii}
         \bf
         induced right action
      }"
    ]
    && X
    \\
 \scalebox{0.7}{$   x,g $}
      &\longmapsto&
  \scalebox{0.7}{$  x \cdot g \mathrlap{ \; \;\coloneqq\; g^{-1} \cdot x  } $}
  \end{tikzcd}
$$

\vspace{-2mm}
\noindent we denote

\vspace{-3mm}
\begin{itemize}
  \item
  its {\it right action groupoid} by

  \vspace{-6mm}
  \begin{equation}
    \label{ActionGroupoid}
    \left( X \times G \rightrightarrows X \right)
    \;\coloneqq\;
    \Big(
      X \times G
      \underoverset
        {\scalebox{.5}{$(-)\cdot (-)$}}
        {\mathrm{pr}_1}
        {\rightrightarrows}
      X
    \Big)
    \;\in\;
    \Groupoids(\mathcal{C}) \;,
  \end{equation}
  whose nerve hence has 2-cells of the following form:
 \vspace{-2mm}
  $$
    N
    \left(
      X \times G \rightrightarrows X
    \right)_2
    \;\;\;
    =
    \;\;\;
    \left\{
    \left.
    \begin{tikzcd}[row sep=40pt]
      &
      x \cdot g_1
      \ar[
        dr,
        "{
          (x\cdot g_1,\, g_0)
        }"{above, sloped, pos=.4},
        "\ "{below, pos=.01, name=s}
      ]
      \\
      x
      \ar[
        ur,
        "{
          (x,\, g_1)
        }"{above, sloped}
      ]
      \ar[
        rr,
        "{
          (x,\, g_1 \cdot g_0)
        }"{below},
        "\ "{above, name=t}
      ]
      &&
      x \cdot g_1 \cdot g_0
      \ar[
        from=s,
        to=t,
        Rightarrow,
        "{(x,\,g_1,\, g_0)}"{description}
      ]
    \end{tikzcd}
    \;
    \right\vert
    \;
    \begin{array}{l}
      x \in X,
      \\
      g_1, g_0 \in G
    \end{array}
    \right\}.
  $$

\end{itemize}

  \vspace{-.2cm}
\noindent
For $\mathcal{C} = \Sets$, we denote the
images of these groupoids in $\InfinityGroupoids$ \eqref{InfinityGroupoids}
as follows:
\vspace{-2mm}
$$
  \begin{tikzcd}[row sep=-4pt]
    \Groupoids(\Sets)
    \ar[rr, "N"]
    &&
    \InfinityGroupoids
    \\
  \scalebox{0.7}{$  (X \times G \rightrightarrows G)  $}
  &\longmapsto&
  \scalebox{0.7}{$ X \!\sslash\! G $}
    \\
\scalebox{0.7}{$    (G \rightrightarrows \ast) $}
&\longmapsto&
\scalebox{0.7}{$\ast \!\sslash\! G \,=:\, B G $}
    \,.
  \end{tikzcd}
$$

\end{notation}

\medskip

\noindent
{\bf Limits and colimits.} We record some facts about homotopy (co-)limits
(e.g. \cite{DHKS04}) that we will need.

\begin{example}[Homotopy pullback preserves projections out of Cartesian products]
  \label{HomotopyPullbackPreserbesProductProjections}
  We have the following diagram:
  \vspace{-2mm}
  $$
    \begin{tikzcd}[column sep=large]
      X \times A
      \ar[r, "{ \mathrm{id} \times f }"]
      \ar[d, "{ \mathrm{pr}_2 }"{left}]
      \ar[dr, phantom, "{ \mbox{\tiny\rm(hbp)} }"]
      &
      X \times B
      \ar[d, "\mathrm{pr}_2"]
      \\
      A
      \ar[r, "{ f }"{below}]
      &
      B
    \end{tikzcd}
  $$
\end{example}
\begin{proof}
  Since homotopy limits commute over each other, this
  homotopy Cartesian diagram is the objectwise product of the
  following two homotopy Cartesian diagrams:
  \vspace{-2mm}
  $$
    \begin{tikzcd}[column sep=large]
      X
      \ar[r, "{ \mathrm{id} }"]
      \ar[d, "{ }"]
      \ar[dr, phantom, "{ \mbox{\tiny\rm(hbp)} }"]
      &
      X
      \ar[d, ""]
      \\
      \ast
      \ar[r, "{  }"]
      &
      \ast
      \mathrlap{\,,}
    \end{tikzcd}
    {\phantom{AAAAA}}
    \begin{tikzcd}[column sep=large]
      A
      \ar[r, "{ f }"]
      \ar[d, "{ \mathrm{id} }"{left}]
      \ar[dr, phantom, "{ \mbox{\tiny\rm(hbp)} }"]
      &
      B
      \ar[d, " \mathrm{id} "]
      \\
      A
      \ar[r, "{ f }"{below}]
      &
      B
    \end{tikzcd}
  $$

  \vspace{-8mm}
\end{proof}

\begin{example}[Colimits over simplicial diagrams of 0-truncated objects are coequalizers (e.g., {\cite[Ex. 8.3.8]{Riehl14}})]
  \label{ColimitsOverSimplicialDiagramsOfZeroTruncatedObjectsAreCoequalizers}
  The colimit of a simplicial diagram of 0-truncated objects
  \vspace{-2mm}
  $$
    X_\bullet
      \,:\,
    \Delta^{\mathrm{op}} \longrightarrow \Topos_{\leq 0}
  $$

  \vspace{-2mm}
  \noindent
  is the coequalizer of the first two face maps,
  hence the quotient of the set $X_0$ by the equivalence relation generated
  by $X_1$:
  \vspace{-2mm}
  $$
    \underset{\longrightarrow}{\mathrm{lim}} \, X_\bullet
    \;\simeq\;
    \mathrm{coeq}
    \big(\!\!\!\!
      \begin{tikzcd}
        X_1
        \ar[
          rr,
          shift left=5pt,
          "d_0"{description}
        ]
        \ar[
          rr,
          shift right=5pt,
          "d_1"{description}
        ]
        &&
        X_0
      \end{tikzcd}
    \!\!\!\!\big)
    \;=:\;
    X_0 / X_1
    \,.
  $$
\end{example}

\medskip

\noindent
{\bf Free homotopy sets.}
We will need to rephrase some statements about homotopy groups,
hence about homotopy classes of pointed maps of $n$-spheres,
in terms of {\it free homotopy sets}, namely
homotopy classes of unconstrained (unpointed) maps out of $n$-spheres.

\begin{lemma}[Detecting weak homotopy equivalences via free homotopy sets {\cite[Thm. 2]{MatumotoMinamiSugawara84}}]
  \label{DetectingWeakHomotopyEquivalencesViaFreeHomotopySets}
  $\,$

  \noindent
  A map
  $X \xrightarrow{\;f\;} Y \,\in\, (\InfinityGroupoids)_{\geq 1}$
  of connected $\infty$-groupoids
  is an equivalence, namely a weak homotopy equivalence
  \vspace{-2mm}
  $$
    \underset{
      n \in \mathbb{N}_+
    }{\forall}
    \;
    \;
    \begin{tikzcd}
      \pi_n(X)
      \ar[r, "{f_\ast}"{above}, "\sim"{below}]
      &
      \pi_n(Y)
      \,,
    \end{tikzcd}
  $$

  \vspace{-2mm}
  \noindent if and only if it induces

  \noindent
  {\bf(a)} isomorphisms
  on all free homotopy sets
  \emph{and}

  \noindent
  {\bf (b)}
  a surjection
  on homotopy classes of maps out of the wedge of circles
  indexed by (the underlying set of) $\pi_1(Y)$:
    \vspace{-2mm}
   $$
    \underset{n \in \mathbb{N}_+}{\forall}
    \;
    \begin{tikzcd}
    \tau_0
    \,
    \Maps{}
      { \ShapeOfSphere{n} }
      { X }
    \ar[r, "{f_\ast}"{above}, "\sim"{below}]
    &
    \tau_0
    \,
    \Maps{}
      { \ShapeOfSphere{n} }
      { Y }
    \end{tikzcd}
    \;\;\;\;\;
    \mbox{\rm and}
    \;\;\;\;\;
    \begin{tikzcd}
    \tau_0
    \,
    \Maps{\big}
      {\, \underset{ \pi_1(Y) }{\vee} \ShapeOfSphere{1} }
      { X }
    \ar[r, ->>,  "{f_\ast}"{above}]
    &
    \tau_0
    \,
    \Maps{\big}
       { \, \underset{ \pi_1(Y) }{\vee} \ShapeOfSphere{1} }
      { Y }
    \,.
    \end{tikzcd}
  $$
\end{lemma}

\begin{lemma}[Truncation and $n$-fold free loop spaces]
  \label{TruncationAndnFoldFreeLoopSpaces}
  For $X \,\in\, \InfinityGroupoids$,
  and
  $n \,\in\, \mathbb{N}$, the following are equivalent:

\vspace{-3mm}
  \begin{enumerate}[{\bf (i)}]
  \setlength\itemsep{-4pt}
    \item
      For all $d > n$
      and all $x \,\in\, X$,
      we have
      $\pi_{d}(X, x) \,=\, \ast$, that is,  $X$ is $n$-truncated, $\tau_n(X) \,\simeq\, X$.

    \item
     For all $d > n$ the point evaluation map is a weak equivalence:
     $
       \begin{tikzcd}
         \Maps{}
           { \ShapeOfSphere{d+1} }
           { X }
         \ar[
           rr,
           "(\mathrm{pt}_d)^\ast"{above},
           "\in\,\WeakHomotopyEquivalences"{below}
         ]
         &&
         X
         \,.
       \end{tikzcd}
     $

    \item
     For all $d > n$ the equatorial evaluation map is a weak equivalence:
     $
       \begin{tikzcd}[column sep=17pt]
         \Maps{}
           { \ShapeOfSphere{d+1} }
           { X }
         \ar[
           rr,
           "(i_d)^\ast"{above},
           "\in\,\WeakHomotopyEquivalences"{below}
         ]
         &&
         \Maps{}
           { \ShapeOfSphere{d} }
           { X }
         \,.
       \end{tikzcd}
     $

  \end{enumerate}

\end{lemma}
\noindent
This is an elementary argument using standard ingredients,
but since the statement it is rarely made explicit in this form
we spell out the proof:
\begin{proof}
First, to amplify the standard fact
that the homotopy groups of a based loop space
at {\it each} loop $\ell \,\in\, \Omega_x X$
(not necessarily the constant loop)
are the shifted homotopy
groups of the underlying space at the given basepoint $x \,\in\, X$,
since the standard fiber sequence with the
contractible based path space $P_x X$
yields a long exact sequence of homotopy groups of this form:
\vspace{-2mm}
\begin{equation}
  \label{HomotopyGroupsOfLoopSpaceAreShiftedHomotopyGroupsOfOriginalSpace}
  \begin{tikzcd}[column sep=10pt]
    \ast
    \ar[r,phantom,"\simeq"]
    &
    \pi_{\bullet + 1}(P_x X, \ell)
    \ar[rr]
    &&
    \pi_{\bullet+1}(X,\,x)
    \ar[rr, "\sim"{above}]
    &&
    \pi_{\bullet}
    (
      \Omega_x X
      ,\,
      \ell
    )
    \ar[
      rr
    ]
    &&
    \pi_{\bullet}
    (
      P X
      ,\,
      \ell
    )
    \ar[r,phantom,"\simeq"]
    &
    \ast
    \,.
  \end{tikzcd}
\end{equation}

\vspace{-2mm}
\noindent
By induction, this implies
that for all iterated loops $\sigma \,\in\, \Omega^d_x X$
we have
\vspace{-2mm}
\begin{equation}
  \label{HomotopyGroupsOfIteratedLoopSpaceAreShiftedHomotopyGroupsOfOriginalSpace}
  \pi_\bullet
  \big(
    \Omega^d_x X
    ,\,
    \sigma
  \big)
  \;\;
  \simeq
  \;\;
  \pi_{\bullet + d}
  (
    X
    ,\,
    x
  )
  \,.
\end{equation}

\vspace{-2mm}
\noindent
Similarly, for any choice of
$\sigma \,\in\, \Maps{}{ \ShapeOfSphere{d} }{ X }$
with
basepoint
$x \,\coloneqq\, \mathrm{pt}_d^\ast(\sigma) \,\in\, X$,
the $d$-fold loop space is the homotopy fiber
of the evaluation map
$
  (\mathrm{pt}_d)^\ast
  \,=\,
  \Maps{}
    { \ast \xrightarrow{\mathrm{pt}_d } \ShapeOfSphere{d} }
    { X }
$
out of the mapping space
from $\ShapeOfSphere{d}$ into $X$:
\vspace{-2mm}
$$
  \begin{tikzcd}
    \Omega_x^d X
    \ar[
      rr,
      "{ \mathrm{fib}_x\left( (\mathrm{pt}_d)^\ast \right) }"
    ]
    &&
    \Maps{}
      { \ShapeOfSphere{d} }
      { X }
    \ar[
      r,
      "{
        (\mathrm{pt}_d)^\ast
      }"
    ]
    &
    X
    \;.
  \end{tikzcd}
$$

\vspace{-2mm}
\noindent
With \eqref{HomotopyGroupsOfIteratedLoopSpaceAreShiftedHomotopyGroupsOfOriginalSpace},
this induces a long exact sequence of homotopy groups of the following form:
\vspace{-2mm}
\begin{equation}
  \label{LongExactSequenceOfHomotopyGroupsForHigherFreeLooSpaceFibration}
  \hspace{-5mm}
  \begin{tikzcd}[column sep=13pt]
    \cdots
    \ar[r]
    &
    \pi_{\bullet + 1}(X,x)
    \ar[r]
    &
    \pi_{\bullet + d}(X,x)
    \ar[
      r
    ]
    &
    \pi_{ \bullet }
    \big(
      \Maps{}
        { \ShapeOfSphere{d} }
        { A }
      ,\,
      \sigma
    \big)
    \ar[
      rr,
      "{
        \scalebox{.7}{$
          \scalebox{1.1}{$($}
            (\mathrm{pt}_d)^\ast
          \scalebox{1.1}{$)$}_\ast
        $}~
      }"
    ]
    &&
    \pi_{\bullet}(X,x)
    \ar[
      r
    ]
    &
    \pi_{\bullet + d - 1}(X,x)
    \ar[r]
    &
    \cdots
    \,.
  \end{tikzcd}
\end{equation}

\vspace{-2mm}
\noindent
The implication
$\mbox{(i)} \Rightarrow \mbox{(ii)}$
follows immediately from the exactness of \eqref{LongExactSequenceOfHomotopyGroupsForHigherFreeLooSpaceFibration}.

Next, since the equatorial embeddings $S^d \xhookrightarrow{ \;i_d\; } S^{d+1}$
clearly respects the basepoint inclusion,
we have for each $x \,\in\, X$ a commuting diagram of this form:
\vspace{-2mm}
\begin{equation}
  \label{CommutingDiagramForComparingHigherFreeLoopSpaces}
  \begin{tikzcd}[column sep=large]
    \Omega_x^{d+1} X
    \ar[
      rr,
      "{
        \mathrm{fib}_x
        \left(
          (\mathrm{pt}_d)^\ast
        \right)
      }"{above}
    ]
    \ar[d]
    &&
    \Maps{}
      { \ShapeOfSphere{d+1} }
      { X }
    \ar[
      r,
      "{
        (\mathrm{pt}_d)^\ast
      }"{above}
    ]
    \ar[
      d,
      "{ (i_d)^\ast }"{left}
      "\in \, \WeakHomotopyEquivalences"{right}
    ]
    &
    X
    \ar[d,-, shift left=1pt]
    \ar[d,-, shift right=1pt]
    \\
    \Omega_x^{d} X
    \ar[
      rr,
      "{
        \mathrm{fib}_x
        \left(
          (\mathrm{pt}_d)^\ast
        \right)
      }"{above}
    ]
    &&
    \Maps{}
      { \ShapeOfSphere{d} }
      {X}
    \ar[
      r,
      "{
        (\mathrm{pt}_d)^\ast
      }"{above}
    ]
    &
    X
  \end{tikzcd}
\end{equation}

\vspace{-2mm}
\noindent This gives the implication $\mbox{(ii)} \Rightarrow \mbox{(iii)}$ by
the 2-out-of-3 property for weak equivalences.

Moreover, for each $\sigma \,\in\, \Maps{}{ \ShapeOfSphere{d+1} }{ X }$
the morphism of long exact sequences of homotopy groups
induced from \eqref{CommutingDiagramForComparingHigherFreeLoopSpaces} is
of this form:
\vspace{-2mm}
$$
\hspace{-2mm}
  \begin{tikzcd}
    \pi_{\bullet + 1}
    \big(
      \Maps{}
        { \ShapeOfSphere{d+1} }
        { X }
      ,\,
      \sigma
    \big)
    \ar[r]
    \ar[
      d,
      "{
        \scalebox{.9}{$
          \scalebox{1.1}{$($}
            (i_d)^\ast
          \scalebox{1.1}{$)$}_\ast
        $}
      }"
    ]
    &
    \pi_{\bullet + 1}(X,\, x)
    \ar[r]
    \ar[d,-,shift left=1pt]
    \ar[d,-,shift right=1pt]
    &
    \pi_{\bullet + d +1 }(X ,\, x)
    \ar[
      d,
      "0"
    ]
    \ar[r]
    &
    \pi_{\bullet}
    \big(
      \Maps{}
        { \ShapeOfSphere{d+1} }
        { X }
      ,\,
      \sigma
    \big)
    \ar[r]
    \ar[
      d,
      "{
        \scalebox{.9}{$
          \scalebox{1.1}{$($}
            (i_d)^\ast
          \scalebox{1.1}{$)$}_\ast
        $}
      }"
    ]
    &
    \pi_{\bullet}(X ,\, x)
    \ar[d,-,shift left=1pt]
    \ar[d,-,shift right=1pt]
    \\
    \pi_{\bullet}
    \big(
      \Maps{}
        { \ShapeOfSphere{d} }
        { X }
      ,\,
      \sigma
    \big)
    \ar[r]
    &
    \pi_{\bullet+1}(X ,\, x)
    \ar[r]
    &
    \pi_{\bullet + d}(X ,\, x)
    \ar[r]
    &
    \pi_{\bullet}
    \big(
      \Maps{}
        { \ShapeOfSphere{d} }
        { X }
      ,\,
      \sigma
    \big)
    \ar[r]
    &
    \pi_{\bullet}(X ,\, x)
    \,.
  \end{tikzcd}
$$

\vspace{-2mm}
\noindent
Here the middle vertical morphism is the function constant on 0,
as shown,
since an equatorial extension of a based map from $S^d$ to $S^{d+1}$
is in particular an extension to $D^{d+1}$ and hence a contracting homotopy.

So when $(i_d)^\ast$ is a weak homotopy equivalence,
then the (non-abelian) five lemma for the above diagram
implies that this 0-map is an isomorphism,
which means that $\pi_{\bullet + d}(X, x) \,=\, 0$ for all $x$.
This is the implication $\mbox{(iii)} \Rightarrow \mbox{(i)}$
and hence concludes the proof.
\end{proof}

\subsection{Homotopy theory of simplicial transformation groups}
\label{TransformationInfinityGroups}

\noindent
{\bf Simplicial Groups.} In analogy to how plain simplicial sets model
bare $\infty$-groupoids, so simplicial groups
(e.g. \cite[\S IV.17]{May67}\cite[\S 3]{Curtis71})
model bare $\infty$-groups:

\begin{notation}[Homotopy theory of simplicial groups ({\cite[\S II 3.7]{Quillen67}\cite[\S V]{GoerssJardine09}})]
  \label{HomotopyTheoryOfSimplicialGroups}
  We write
  \vspace{-2mm}
  \begin{equation}
    \label{CategoryOfSimplicialGroups}
    \Groups(\SimplicialSets)_{\mathrm{proj}}
    \;\in\;
    \ModelCategories
  \end{equation}

  \vspace{-2mm}
  \noindent
  for the model category of simplicial groups
  (\cite[\S 17]{May67}\cite[\S 3]{Curtis71}),
  whose
  weak equivalences and fibrations are those of
  the underlying  $\SimplicialSets_{\mathrm{Qu}}$.
\end{notation}
\begin{lemma}[Simplicial groups are Kan complexes ({\cite[Thm. 3]{Moore54}\cite[Thm. 17.1]{May67}\cite[Lem. 3.1]{Curtis71}})]
  \label{UnderlyingSimplicialSetOfASimplicialGroupIsKanComplex}
  The underling simplicial set of any simplicial group
  \eqref{CategoryOfSimplicialGroups}
  is a Kan complex.
\end{lemma}

The following model category theoretic version
(Prop. \ref{QuillenEquivalenceBetweenSimplicialGroupsAndReducedSimplicialSets})
of the general
looping/delooping equivalence \eqref{LoopingAndDelooping}
relates simplicial groups to reduced simplicial sets:

\begin{notation}[Homotopy theory of reduced simplicial sets {(\cite[\S V Prop. 6.12]{GoerssJardine09})}]
  \label{HomotopyTheoryOfReducedSimplicialSets}
  We write
  \vspace{-2mm}
  $$
    \SimplicialSets_{\geq 1, \mathrm{inj}}
    \;\;
    \in
    \;\;
    \ModelCategories
  $$

  \vspace{-2mm}
\noindent
  for the model category of reduced simplicial sets
  (those $S \in \SimplicialSets$ with a single vertex, $S_0 = \ast$ )
  whose weak equivalences and cofibrations are those of the underlying
  $\SimplicialSets_{\mathrm{Qu}}$ (Ntn. \ref{KanQuillenModelCategoryOfSimplicialSets}).
\end{notation}

\begin{lemma}[Fibrant reduced simplicial sets are Kan complexes {\cite[\S V, Lem. 6.6]{GoerssJardine09}}]
  \label{FibrantReducedSimplicialSetsAreKanComplexes}
  While the forgetful functor
  \vspace{-2mm}
  $$
    \begin{tikzcd}
      (\SimplicialSets_{\geq 1})_{\mathrm{inj}}
      \ar[
        rr,
        "\underlying"
      ]
      &&
      \SimplicialSets_{\mathrm{Qu}}
    \end{tikzcd}
  $$

  \vspace{-2mm}
  \noindent
  does not preserve all fibrations, it does preserve fibrant objects.
\end{lemma}

\begin{proposition}[Quillen equivalence between simplicial groups and reduced simplicial sets
 {\cite[\S V, Prop. 6.3]{GoerssJardine09}}]
  \label{QuillenEquivalenceBetweenSimplicialGroupsAndReducedSimplicialSets}
  The simplicial classifying space construction (Def. \ref{StandardModelOfUniversalSimplicialPrincipalComplex})
  is the right adjoint of a Quillen equivalence between
  the
  projective model structure on simplicial groups (Ntn. \ref{HomotopyTheoryOfSimplicialGroups})
  and the injective model structure on reduced simplicial sets
  (Nota \ref{HomotopyTheoryOfReducedSimplicialSets}):
  \vspace{-2mm}
  $$
    \begin{tikzcd}
      \Groups(\SimplicialSets)_{\mathrm{proj}}
      \ar[
        rr,
        shift right=6pt,
        "\overline{W}(-)"{below}
      ]
      \ar[
        rr,
        phantom,
        "\scalebox{.7}{$\simeq_{\mathrlap{\rm Qu}}$}"
      ]
      &&
      (\SimplicialSets_{\geq 1})_{\mathrm{inj}}\;.
      \ar[
        ll,
        shift right=6pt
      ]
    \end{tikzcd}
  $$
\end{proposition}

\noindent
{\bf Universal principal simplicial complex.}
With $\infty$-groups presented by simplicial groups $\mathcal{G}$,
$\mathcal{G}$-principal $\infty$-bundles are presented by the
construction traditionally known denoted ``$W \mathcal{G}$'':

\begin{definition}[Universal principal simplicial complex
{\cite[Def. 10.3]{Kan58}\cite[p. 87-788]{May67}\cite[p. 269]{GoerssJardine09}}]
 \label{StandardModelOfUniversalSimplicialPrincipalComplex}
  Let $\mathcal{G} \in \Groups(\SimplicialSets)$.

  \noindent
  {\bf (i)} Its {\it standard universal principal complex} is the simplicial set
  \vspace{-2mm}
  $$
    W \mathcal{G}
    \;\;
    \in
    \;
    \SimplicialSets
  $$

  \vspace{-3mm}
  \noindent
  whose

  \vspace{-4mm}
  \begin{itemize}
  \setlength\itemsep{-6pt}
  \item
  component sets are
  $$
    (W \mathcal{G})_n
    \;\coloneqq\;
    \mathcal{G}_{n}
      \times
    \mathcal{G}_{n-1}
      \times
    \cdots
      \times
    \mathcal{G}_0
    \,,
  $$

  \item
  face maps are given by
  \vspace{-2mm}
  \begin{equation}
    \label{FaceMapsOfWG}
    \hspace{-5mm}
    d_i
    \left(
      \gamma_n, \, \gamma_{n-1},\, \cdots ,\, \gamma_0
    \right)
    \;\coloneqq\;
    \left\{\!\!\!\!
    \begin{array}{lcl}
      \big(
        d_i(\gamma_n),
        \,
        d_{i-1}(\gamma_{n-1}),
        \,
        \cdots,
        \,
        d_0(\gamma_{n-i}) \cdot \gamma_{n-i-1},
        \,
        \gamma_{n-i-2},
        \,
        \cdots,
        \,
        \gamma_0
      \big)
      &\mbox{for}&
      0 < i < n
      \\[5pt]
      \big(
        d_n(\gamma_n),
        \,
        d_{n-1}(\gamma_{n-1}),
        \,
        \cdots,
        \,
        d_1(\gamma_1)
      \big)
      &\mbox{for}&
      i = n
      \,,
    \end{array}
    \right.
  \end{equation}

  \item
  degeneracy maps are given by
  \vspace{-2mm}
  \begin{equation}
    \label{DegeneracyMapsOfWG}
    s_i
    \left(
      \gamma_n,
      \,
      \gamma_{n-1},
      \,
      \cdots,
      \,
      \gamma_0
    \right)
    \;:=\;
    \big(
      s_i(\gamma_n),
      \,
      s_{i-1}(\gamma_{n - 1}),
      \,
      \cdots,
      \,
      s_0(\gamma_{n - i}),
      \,
      e,
      \,
      \gamma_{n - i - 1},
      \,
      \cdots,
      \,
      \gamma_0
    \big)
    \,,
  \end{equation}

  \item
  and equipped with the left $\mathcal{G}$-action
  (Ex. \ref{UniversalPrincipalSimplicialComplexInGActions})
  given by
  \vspace{-2mm}
  \begin{equation}
    \label{LeftActionOnUniversalSimplicialPrincipalSpace}
    \begin{tikzcd}[row sep=-3pt]
      \mathcal{G} \times W\mathcal{G}
      \ar[rr]
      &&
      W \mathcal{G}
      \\
       \scalebox{0.7}{$
       \big(
        h_n,
        \,
        (\gamma_n, \gamma_{n-1}, \cdots, \gamma_0)
       \big)
      $}
      &\longmapsto&
      \scalebox{0.7}{$   \big(
        h_n \cdot \gamma_n,
        \,
        \gamma_{n-1},
        \,
        \cdots,
        \,
        \gamma_0
      \big)
      $}
      \,.
    \end{tikzcd}
  \end{equation}

  \end{itemize}

  \vspace{-4mm}
  \noindent
  {\bf (ii)}  Its standard {\it simplicial delooping} or
  {\it simplicial classifying complex}
  $\overline{W}\mathcal{G}$ is the quotient
  by that action \eqref{LeftActionOnUniversalSimplicialPrincipalSpace}:
  \vspace{-2mm}
  \begin{equation}
    \label{StandardSimplicialDeloopingAsQuotient}
    W\mathcal{G}
    \xrightarrow{ \;\;q_{{}_{W\mathcal{G}}}\;\; }
    \overline{W}\mathcal{G}
    \;:=\;
    (
      W \mathcal{G}
    )/\mathcal{G}
    \,.
  \end{equation}
\end{definition}

\begin{example}[Low-dimensional cells of universal simplicial principal complex]
Unwinding the definition \eqref{FaceMapsOfWG} of the face maps of $W \mathcal{G}$
(Def. \ref{StandardModelOfUniversalSimplicialPrincipalComplex}) shows that its
1-simplices are of the form
$$
  (W \mathcal{G})_1
  \;\;
  =
  \;\;
  \Big\{\!\!\!\!
  \begin{tikzcd}
    d_1(g_1)
    \ar[
      rr,
      "{(g_1, g_0)}"
    ]
    &&
    d_0(g_1)\cdot g_0
  \end{tikzcd}
  \;\Big\vert\;
  \begin{array}{l}
    g_0 \,\in\, \mathcal{G}_0
    \\
    g_1 \,\in\, \mathcal{G}_1
  \end{array}
\!\!  \Big\}
$$
and its 2-simplices are of this form:
$$
  (W\mathcal{G})_2
  \;\;
  =
  \;\;
  \left\{
  \left.
  \begin{tikzcd}[column sep=-10pt, row sep=large]
    &&
    \scalebox{.9}{$
      {d_0 d_2(g_2) \cdot d_1(g_1)}
      \,=\,
      {d_1 d_0(g_2) \cdot d_1(g_1)}
    $}
    {}
    \ar[
      ddrr,
      "{
        \scalebox{1}{$($}
          d_0(g_2) \cdot g_1,
          \,
          g_0
        \scalebox{1}{$)$}
      }"{above, sloped, pos=.6},
      "\ "{below, name=s, pos=.02}
    ]
    \\
    \\
      \scalebox{1.2}{$    {
      {\phantom{=\,} d_1 d_2(g_2)}
      \atop
      {=\, d_1 d_1(g_2)}
    }
    $}
    \ar[
      rrrr,
      "{
        \scalebox{1}{$($}
        d_1(g_2),
        \,
        d_0(g_1)\cdot g_0
        \scalebox{1}{$)$}
      }"{below},
      "{\ }"{above, name=t}
    ]
    \ar[
      uurr,
      "{
        \scalebox{1}{$($}
          d_2(g_2),
          \,
          d_1(g_1)
        \scalebox{1.3}{$)$}
      }"{above, sloped, pos=.4}
    ]
    &&&&
      \scalebox{1.2}{$  {
      {\phantom{=,} d_0 d_0(g_2) \cdot d_0(g_1) \cdot g_0}
      \atop
      {=\,  d_0 d_1(g_2)\cdot d_0(g_1) \cdot g_0}
    }
    $}
    \ar[
      from=s,
      to=t,
      Rightarrow,
      "{
        (g_2, g_1, g_0)
      }"{description}
    ]
  \end{tikzcd}
  \;
  \right\vert
  \;
  \begin{array}{l}
    g_0 \,\in\, \mathcal{G}_0,
    \\
    g_1 \,\in\, \mathcal{G}_1,
    \\
    g_2 \,\in\, \mathcal{G}_2
  \end{array}
  \right\}.
$$

\end{example}

\begin{example}[Universal principal simplicial complex for ordinary group $G$]
  \label{UniversalSimplicialPrincipalComplexForOrdoinaryGroupG}
  If
  \vspace{-2mm}
  $$
    G
    \;\in\;
    \Groups
    \xhookrightarrow{\;\;\;}
    \Groups
    (
      \SimplicialSets
    )
  $$

  \vspace{-2mm}
  \noindent
  is an ordinary discrete group, regarded as a simplicial group
  (hence the functor constant on $G$ on the opposite simplex category),
  then the  standard model of its universal principal complex (Def. \ref{StandardModelOfUniversalSimplicialPrincipalComplex})
  is isomorphic to the nerve of the
  action groupoid \eqref{ActionGroupoidOfRightMultiplicationAction}
  of the right multiplication action of $G$ on itself:
  \vspace{-2mm}
  \begin{equation}
    \label{WGForOrdinaryGroupIsNerveOfActionGroupoid}
    W G
    \;=\;
    N
    (
      G \times G
        \rightrightarrows
      G
    )
    \,.
  \end{equation}
\vspace{-4mm}
$$
(W G)_2
\;\;
=
\;\;
\left\{
\left.
\begin{tikzcd}
  & g_2 g_1
  \ar[
    ddr,
    "{(g_2 g_1, g_0)}"{sloped},
    ""{name=s, below, pos=0.01}
  ]
  \\
  \\
  g_2
  \ar[
    uur,
    "{(g_2, g_1)}"{sloped}
  ]
  \ar[
    rr,
    "{(g_2, g_1 g_0)}"{below},
    ""{name=t, above}
  ]
  &
  {}
  &
  g_2 g_1 g_0
  \ar[
    from=s, to=t,
    Rightarrow,
    "{(g_2, g_1, g_0)}"{description}
  ]
\end{tikzcd}
\;\right\vert\;
 g_0, g_1, g_2 \,\in\, G
\right\}.
$$
Accordingly, the standard simplicial delooping \eqref{StandardSimplicialDeloopingAsQuotient}
of an ordinary group is isomorphic to the nerve of its delooping groupoid
\eqref{DeloopingGroupoid}:
$$
  \overline{W}G
  \;\simeq\;
  N( G \rightrightarrows \ast )
  \,.
$$
\end{example}
\begin{proposition}[Basic properties of standard simplicial principal complex {\cite[\S V, Lem. 4.1, 4.6, Cor. 6.8]{GoerssJardine09}}]
  \label{BasicPropertiesOfStandardSimplicialPrincipalComplex}
  For $\mathcal{G} \in \Groups(\mathrm{SimplSets})$,
  its standard universal principal complex (Def. \ref{StandardModelOfUniversalSimplicialPrincipalComplex})
  has the following properties:

  {\bf (i)} $W \mathcal{G}$ is contractible;

  {\bf (ii)} $W \mathcal{G}$ and $\overline{W} \mathcal{G}$ are Kan complexes.
\end{proposition}
\begin{proof}
  That $\overline{W}\mathcal{G}$ is Kan fibrant
  follows as the combination of
  Lem. \ref{UnderlyingSimplicialSetOfASimplicialGroupIsKanComplex},
  Prop. \ref{QuillenEquivalenceBetweenSimplicialGroupsAndReducedSimplicialSets},
  and
  Lem. \ref{FibrantReducedSimplicialSetsAreKanComplexes}.
  This implies that $W \mathcal{G}$ is Kan fibrant since
  $W \mathcal{G} \xrightarrow{q} \overline{W}\mathcal{G}$
  is a Kan fibration \eqref{TheUniversalPrincipalSimplicialBundle}
  (by Prop. \ref{QuillenEquivalenceBetweenSimplicialGroupsAndReducedSimplicialSets},
  see Ex. \ref{CoprojectionsOutOfBorelConstructionAreKanFibrations}).
\end{proof}

\begin{proposition}[{\cite[\S V, Cor. 6.9]{GoerssJardine09}}]
  \label{EssentiallySurjectiveKanFibrationsOfSimplicialGroupsInducedKanFibrationsUnderBarW}
  Let $\mathcal{G}_1 \xrightarrow{\phi} \mathcal{G}_2$
  be a homomorphism of simplicial groups which is a Kan fibration
  of underlying simplicial sets. Then the induced morphism of
  simplicial classifying spaces
  $\overline{W}\mathcal{G}_1 \xrightarrow{ \overline{W}\phi } \overline{W}\mathcal{G}_2$
  (Def. \ref{StandardModelOfUniversalSimplicialPrincipalComplex})
  is a Kan fibration if and only if $\phi$ is a surjection on connected components.
\end{proposition}
For plain groups $\mathcal{G} \,\in\, \Groups(\Sets) \hookrightarrow \Groups(\SimplicialSets)$,
the statement of Prop. \ref{EssentiallySurjectiveKanFibrationsOfSimplicialGroupsInducedKanFibrationsUnderBarW}
is also readily seen by direct inspection.
Elementary as it is, it has some important consequences
(see Lem. \ref{InfinityActionBaseChangeComonadAlongDiscreteGroupExtensions},
Prop. \ref{PullbackOfActionsAlongSurjectiveGroupHomomorphismsIsFullyFaithful} below).

\medskip

\noindent
{\bf Simplicial group actions.} With bare $\infty$-groups presented by
simplicial groups, their $\infty$-actions are presented by
the homotopy theory of simplicial group actions (Ntn. \ref{ProjectiveModelStructureOnActionsOfSimplicialGroups} below),
sometimes referred to as the {\bf Borel model structure}.

\begin{notation}[Simplicial group actions]
  \label{SimplicialGroupActions}
  For $\mathcal{G} \in \Groups(\SimplicialSets)$,
  we denote

  \noindent
  {\bf (i)}
  by

  \vspace{-.5cm}
  \begin{equation}
    \label{SimplicialDeloopingOfSimplicialGroup}
    \mathbf{B}\mathcal{G}
    \;\in\;
    \EnrichedCategories{\SimplicialSets}
    \,,
    \;\;\;\;\;
    \mathbf{B}\mathcal{G}(\ast,\ast) \;:=\; \mathcal{G}
    \,,
  \end{equation}
  \vspace{-.4cm}

  \noindent
  the simplicial groupoid with a single object $\ast$, with
  $\mathcal{G}$ as its unique hom-object
  and with composition ``$\circ$''
  given by the {\it reverse} of the group product ``$\cdot$''

  \vspace{-.5cm}
  \begin{equation}
    \label{CompositionInDeloopedSimplicialGroup}
    \begin{tikzcd}[row sep =-4pt]
      \mathcal{G} \times \mathcal{G}
      \ar[
        rr,
        "\circ"
      ]
      &&
      \mathcal{G}
      \\
   \scalebox{0.7}{$      (g_n , h_n) $}
      &\longmapsto&
  \scalebox{0.7}{$       h_n \cdot g_n $\;;}
    \end{tikzcd}
  \end{equation}
  \vspace{-.4cm}

  \noindent
  {\bf (ii)}
  the category of $\mathcal{G}$-actions on simplicial sets by:
  \vspace{-2mm}
  \begin{equation}
    \label{CategoryOfSimplicialGroupActions}
    \Actions{\mathcal{G}}
    (
      \SimplicialSets
    )
    \;:=\;
    \SimplicialFunctors
    (
      \mathbf{B} \mathcal{G},
      \,
      \SimplicialSets
    )
    \,,
  \end{equation}

  \vspace{-2mm}
\noindent
  identified
  with the category of $\SimplicialSets$-enriched functors
  from the delooping \eqref{SimplicialDeloopingOfSimplicialGroup}
  of $\mathcal{G}$ to $\SimplicialSets$
  \eqref{SimplicialSetsAsAnEnrichedCategory}.
\end{notation}

\begin{remark}[Simplicial group actions are from the left]
  \label{SimplicialGroupActionsAreFromTheLeft}
  The convention \eqref{CompositionInDeloopedSimplicialGroup}
  for the delooping $\mathbf{B}\mathcal{G}$ \eqref{SimplicialDeloopingOfSimplicialGroup}
  implies
  that the simplicial $\mathcal{G}$-actions
  \eqref{CategoryOfSimplicialGroupActions} are {\it left} actions:
  \vspace{-2mm}
  $$
    \begin{tikzcd}[row sep=0pt]
      \mathbf{B}\mathcal{G}
      \ar[
        rr,
        "\mathcal{G} \acts \, X"
      ]
      &&
      \SimplicialSets
      \\[2pt]
      \bullet
      \ar[
        d,
        "g_1"{left}
      ]
      \ar[
        dd,
        rounded corners,
        to path={
           -- ([xshift=-8pt]\tikztostart.west)
           --node[below, sloped]{
               \scalebox{.7}{$
                 g_1 \circ g_2
                 \,\coloneqq\,
                 g_2 \cdot g_1
               $}
             } ([xshift=-8pt]\tikztotarget.west)
           -- (\tikztotarget.west)}
      ]
      &&
      X
      \ar[
        d,
        "g_1 \cdot (-)"
      ]
      \ar[
        dd,
        rounded corners,
        to path={
           -- ([xshift=+21pt]\tikztostart.east)
           --node[above, sloped]{
               \scalebox{.7}{$
                 (g_2 \cdot g_1) \cdot (-)
               $}
             } ([xshift=+21pt]\tikztotarget.east)
           -- (\tikztotarget.east)}
      ]
      \\[17pt]
      \bullet
      \ar[
        d,
        "g_2"{left}
      ]
      &&
      X
      \ar[
        d,
        "g_2 \cdot (-)"
      ]
      \\[17pt]
      \bullet
      &&
      X
    \end{tikzcd}
    {\phantom{AAAAAA}}
    \begin{tikzcd}[row sep=-2pt]
      \mathcal{G} \times X
      \ar[
        rr,
        "(-) \cdot (-)"
      ]
      &&
      X
      \\
    \scalebox{0.7}{$     (g_n, x_n) $}
      &\longmapsto&
     \scalebox{0.7}{$    g_n \cdot x_n $}
      \mathrlap{\,.}
    \end{tikzcd}
  $$
  \vspace{-.4cm}
\end{remark}
\begin{example}[Universal principal simplicial complex in $\mathcal{G}$-actions]
  \label{UniversalPrincipalSimplicialComplexInGActions}
  For $\mathcal{G} \in \SimplicialGroups$,
  the universal principal simplicial complex
  $W \mathcal{G}$ (Def. \ref{StandardModelOfUniversalSimplicialPrincipalComplex})
  becomes an object of \eqref{CategoryOfSimplicialGroupActions}
  by the formula \eqref{LeftActionOnUniversalSimplicialPrincipalSpace}.
  \vspace{-2mm}
  \begin{equation}
    \label{UniversalPrincipalSimplicialComplexAsLeftSimplicialGroupAction}
    \mathcal{G} \acts \; W \mathcal{G}
    \,\in\,
    \Actions{\mathcal{G}}(\SimplicialSets)\;.
  \end{equation}
    \end{example}

Making explicit the following elementary Ex. \ref{SimplicialGroupCanonicallyActingOnItself} is,
serves to straighten out a web of conventions about (simplicial) group actions.
\begin{example}[Simplicial group canonically acting on itself]
  \label{SimplicialGroupCanonicallyActingOnItself}
  Any $\mathcal{G} \in \Groups(\SimplicialSets)$ becomes
  an object of the category of simplicial $\mathcal{G}$-actions
  \eqref{CategoryOfSimplicialGroupActions}
  in three canonical ways:

  \vspace{-6mm}
  \begin{equation}
  \hspace{-5mm}
    \begin{tikzcd}[column sep=18pt, row sep=-4pt]
      \mathcal{G} \times \mathcal{G}
      \ar[
        rr,
        "\mathclap{\mbox{
          \tiny
          \color{greenii}
          \bf
          \def\arraystretch{.9}
          \begin{tabular}{c}
            left
            \\
            multiplication
            action
          \end{tabular}
        }}"
      ]
      &&
      \mathcal{G}\;,
      \\
    \scalebox{0.7}{$     (g_n, h_n) $}
      &\longmapsto&
     \scalebox{0.7}{$    g_n \cdot h_n $}
    \end{tikzcd}
    \qquad  \quad
    \begin{tikzcd}[column sep=18pt, row sep=-4pt]
      \mathcal{G} \times \mathcal{G}
      \ar[
        rr,
        "\mathclap{\mbox{
          \tiny
          \color{greenii}
          \bf
          \def\arraystretch{.9}
          \begin{tabular}{c}
            right inverse
            \\
            multiplication
            action
          \end{tabular}
        }}"
      ]
      &&
      \mathcal{G}\;,
      \\
     \scalebox{0.7}{$    (g_n, h_n) $}
      &\longmapsto&
     \scalebox{0.7}{$    h_n \cdot g_n^{-1} $}
    \end{tikzcd}
    \qquad  \quad
    \begin{tikzcd}[column sep=15pt, row sep=-4pt]
      \mathcal{G} \times \mathcal{G}
      \ar[
        rr,
        "\mathclap{\mbox{
          \tiny
          \color{greenii}
          \bf
          \def\arraystretch{.9}
          \begin{tabular}{c}
            adjoint/conjugation
            \\
            action
          \end{tabular}
        }}"
      ]
      &&
      \mathcal{G}\;.
      \\
      \scalebox{0.7}{$   (g_n, h_n) $}
      &\longmapsto&
      \scalebox{0.7}{$   g_n \cdot h_n \cdot g_n^{-1} $}
    \end{tikzcd}
  \end{equation}

  \vspace{-2mm}
\noindent
  The first two are isomorphic in $\Actions{\mathcal{G}}(\SimplicialSets)$
  via the inversion operation:
  \vspace{-2mm}
  \begin{equation}
    \begin{tikzcd}[row sep=1pt, column sep=14pt]
      \scalebox{0.7}{$   (g_n, h_n) $}
      \ar[
       rrrr,
       |->
      ]
      \ar[
        ddd,
        |->
      ]
      &&[40pt] &&
   \scalebox{0.7}{$      g_n \cdot h_n $}
      \ar[
        ddd,
        |->
      ]
      \\
      &
      \mathcal{G} \times \mathcal{G}
      \ar[
        rr,
        "{
          \mbox{
            \tiny
            \color{greenii}
            \bf
            \def\arraystretch{.9}
            \begin{tabular}{c}
              left
              \\
              multiplication
            \end{tabular}
          }
        }"{above}
      ]
      \ar[
        d,
        "{
          \mathrm{id} \times (-)^{-1}
        }"{left},
        "\sim"{above, sloped}
      ]
      &&
      \mathcal{G}
      \ar[
        d,
        "{
          (-)^{-1}
        }"{right},
        "\sim"{below, sloped}
      ]
      \\[20pt]
      &
      \mathcal{G} \times \mathcal{G}
      \ar[
        rr,
        "{
          \mbox{
            \tiny
            \color{greenii}
            \bf
            \def\arraystretch{.9}
            \begin{tabular}{c}
              right inverse
              \\
              multiplication
            \end{tabular}
          }
        }"{below}
      ]
      &&
      \mathcal{G}
      \\[2pt]
       \scalebox{0.7}{$  (g_n, h_n^{-1}) $}
      \ar[
        rrrr,
        |->
      ]
      &&&&
   \scalebox{0.7}{$      h_n^{-1} \cdot g_n^{-1} $}
    \end{tikzcd}
  \end{equation}
  \vspace{-.4cm}

\end{example}

\begin{notation}[Model category of simplicial group actions {(\cite[\S 2]{DDK80}\cite[\S 5]{Guillou06}\cite[\S V Thm. 2.3]{GoerssJardine09})}]
  \label{ProjectiveModelStructureOnActionsOfSimplicialGroups}
  For $\mathcal{G} \in \Groups( \SimplicialSets)$,
  we have on
  the category of $\mathcal{G}$-actions \eqref{CategoryOfSimplicialGroupActions}
  the projective model structure
  (the {\it coarse-} or {\it Borel- equivariant model structure})
  whose fibrations and weak equivalences are those
  of the underlying $\SimplicialSets_{\mathrm{Qu}}$
  (Ntn. \ref{KanQuillenModelCategoryOfSimplicialSets}),
  which we denote as:

  \vspace{-.5cm}
  \begin{equation}
    \label{BorelModelStructure}
    \Actions{\mathcal{G}}
    (
      \SimplicialSets
    )_{\mathrm{proj}}
    \;:=\;
    \SimplicialFunctors
    (
      \mathbf{B} \mathcal{G},
      \,
      \SimplicialSets
    )_{\mathrm{proj}}
    \,.
  \end{equation}
  \vspace{-.4cm}
\end{notation}

\begin{lemma}[Cofibrations of simplicial group actions {\cite[Prop. 2.2 (ii)]{DDK80}\cite[Prop. 5.3]{Guillou06}\cite[\S V Lem. 2.4]{GoerssJardine09}}]
  \label{CofibrationsOfSimplicialGroupActions}
  The cofibrations of
  $
    \Actions{\mathcal{G}}
    (
      \SimplicialSets
    )_{\mathrm{proj}}
  $ \eqref{BorelModelStructure} are the monomorphisms
  such that the $\mathcal{G}$-action on the simplices not in their image is free.
\end{lemma}

\begin{lemma}[Equivariant equivalence of simplicial universal principal complexes]
  \label{EquivariantEquivalenceOfSimplicialUniversalPrincipalComplexes}
  For $\mathcal{H} \xhookrightarrow{\;i\;} \mathcal{G}$
  a simplicial subgroup inclusion, the induced inclusion

  \vspace{-.4cm}
  $$
    \begin{tikzcd}
      W \mathcal{H}
      \ar[
        rr,
        "W(i)"{above},
        "\in \mathrm{W}"{below}
      ]
      &&
      W \mathcal{G}
    \end{tikzcd}
    \;\;\;
    \in
    \;
    \Actions{\mathcal{H}} (\SimplicialSets)_{\mathrm{proj}}
  $$
  \vspace{-.5cm}

  \noindent
  of their standard simplicial principal complexes
  (Def. \ref{StandardModelOfUniversalSimplicialPrincipalComplex})
  equipped with their canonical $\mathcal{H}$-action
  \eqref{LeftActionOnUniversalSimplicialPrincipalSpace}
  is a weak equivalence
  in the Borel-equivariant model structure \eqref{BorelModelStructure}.
\end{lemma}
\begin{proof}
The underlying simplicial sets of both are contractible,
  by Prop. \ref{BasicPropertiesOfStandardSimplicialPrincipalComplex}),
so that underlying any equivariant morphism between them is an
simplicial weak homotopy equivalence.
\end{proof}

\begin{proposition}[Quillen equivalence between Borel model structure and the slice over
classifying complex]
\label{QuillenEquivalenceBetweenBorelModelStructureAndSliceOverClassifyingComplex}
\noindent {\bf (i)} For
any $\mathcal{G} \in \Groups( \SimplicialSets)$,
there is a simplicial adjunction

\vspace{-3mm}
\begin{equation}
  \label{BorelConstructionAdjunction}
  \mathclap{
    \raisebox{3pt}{
      \tiny
      \color{darkblue}
      \bf
      homotopy fiber
    }
    }
    \qquad \;\;
  {
    (-) \times_{\overline{W}\mathcal{G}} W \mathcal{G}
  }
  \quad
  \dashv
  \quad
  {
    \left(
      (-) \times W \mathcal{G}
    \right)/\mathcal{G}
  }
  \qquad \quad
    \mathclap{
    \raisebox{3pt}{
      \tiny
      \color{darkblue}
      \bf
      Borel construction
    }
  }
\end{equation}
\vspace{-5mm}

\noindent
between the Borel model structure
\eqref{BorelModelStructure}
and the slice model structure
of $\SimplicialSets_{\mathrm{Qu}}$ \eqref{KanQuillenModelCategoryOfSimplicialSets}
over the simplicial classifying complex $\overline{W}\mathcal{G}$ \eqref{StandardSimplicialDeloopingAsQuotient}.

\noindent {\bf (ii)} Hence there is a natural isomorphism of hom-complexes
\vspace{-2mm}
\begin{equation}
  \label{BorelSliceAdjunctionAsNaturalIsoOfSimplicialHomComplexes}
  \hspace{-3mm}
  \begin{aligned}
  \Actions{\mathcal{G}}
  (\SimplicialSets)
  \big(
    (-) \times_{\overline{W}\mathcal{G}} W \mathcal{G},
    \,
    (-)
  \big)
  \;\simeq\;
  &
  \SimplicialSets_{/\overline{W}\mathcal{G}}
  \left(\!
    (-),
    \,
    \left(
      (-) \times W \mathcal{G}
    \right)/\mathcal{G}
  \right)
  \end{aligned}
  \!\!\!\!\!\!\!
  \in
  \;
  \SimplicialSets
  \,,
\end{equation}

  \vspace{-2mm}
\noindent
which is a Quillen equivalence:
\vspace{-4mm}
\begin{equation}
  \label{QuillenEquivalenceBetweenBorelModelStructureAndSliceOverSimplicialClassifyingComplex}
  \begin{tikzcd}[column sep=large]
    \Actions{\mathcal{G}}
    (
      \SimplicialSets
    )_{\mathrm{proj}}
    \ar[
      rr,
      shift right=7pt,
      "{
        (
          (-) \times W \mathcal{G}
        )/ \mathcal{G}
      }"{below}
    ]
    \ar[
      rr,
      phantom,
      "{\simeq_{\mathrlap{\mathrm{Qu}}}}"
    ]
    &&
    (
      \SimplicialSets_{\mathrm{Qu}}
    )_{/\overline{W}\mathcal{G}}
    \mathrlap{\,.}
    \ar[
      ll,
      shift right=7pt,
      "{
        (-) \times_{\overline{W}\mathcal{G}} W\mathcal{G}
      }"{above}
    ]
  \end{tikzcd}
\end{equation}
\end{proposition}
\begin{proof}
  The plain adjunction constituting a Quillen equivalence \eqref{QuillenEquivalenceBetweenBorelModelStructureAndSliceOverSimplicialClassifyingComplex}
  is \cite[Prop. 2.3]{DDK80}. The simplicial enrichment \eqref{BorelSliceAdjunctionAsNaturalIsoOfSimplicialHomComplexes},
  hence the natural bijections
  \vspace{-2mm}
  $$
    \mathrm{Hom}
    \big(\!\!
      \left(
        (-) \times_{\overline{W}\mathcal{G}}
      \right)
      \times \Delta[k]
      ,\,
      (-)
    \big)
    \;\simeq\;
    \mathrm{Hom}
    \big(
      (-) \times \Delta[k]
      ,\,
      \left(
        (-) \times W \mathcal{G}
      \right)/\mathcal{G}
    \big)
    \;\;\;
    \in
    \;
    \Sets
    \,,
  $$

  \vspace{-1mm}
  \noindent
  are left somewhat implicit in \cite[Prop. 2.4]{DDK80}, but it follows
  readily from the plain adjunction via the natural isomorphism
  $$
    \left(
      (-) \times_{\overline{W}\mathcal{G}} W G
    \right) \times \Delta[k]
    \;\simeq\;
    \left(
      (-) \times \Delta[k]
    \right) \times_{\overline{W}\mathcal{G}} W \mathcal{G}
    \,,
  $$
  which, in turn, follows from the pasting law \eqref{HomotopyPastingLaw}:
  \vspace{-2mm}
  $$
  \hspace{-1cm}
  \begin{tikzcd}[column sep=small, row sep=small]
      { \left(X \times_{\overline{W}G} W G\right) \times \Delta[k] } \!\!\!
      &
      { \mathllap{\simeq \;} \;\; (X \times \Delta[k]) \times_{\overline{W}G} W G  }
  \ar[rrrrr]
  \ar[d]
  \ar[drrrrr,phantom,"\mbox{\tiny\rm(pb)}"]
  &&&&&
  X \times \Delta[k]
  \ar[
    d,
    "\;\mathrm{pr}_1"
  ]
  \\
  &
  X \times_{\overline{W}G} W G
  \ar[rrrrr]
  \ar[d]
  \ar[drrrrr,phantom,"\mbox{\tiny\rm(pb)}"]
  &&&&&
  X
  \ar[
    d
  ]
  \\
  &
  W G
  \ar[rrrrr]
  &&&&&
  \overline{W}G
  \mathrlap{\,.}
\end{tikzcd}
  $$

  \vspace{-7mm}
\end{proof}

\begin{example}[Coprojections out of Borel construction are Kan fibrations]
  \label{CoprojectionsOutOfBorelConstructionAreKanFibrations}
  For $\mathcal{G} \acts \, X \,\in\, \Actions{\mathcal{G}}(\SimplicialSets)$
  such that the underlying simplicial set $X$ being a Kan complex, hence such that
  \vspace{-2mm}
  $$
    \mathcal{G} \acts \, X
    \xrightarrow{ \;\in\; \mathrm{Fib} }
    \ast
    \;\;\;
    \in
    \;
    \Actions{\mathcal{G}}(\SimplicialSets)_{\mathrm{proj}}
    \,,
  $$

  \vspace{-2mm}
  \noindent
  the projection from the Borel construction
  \eqref{BorelConstructionAdjunction}
  to the simplicial classifying space
  (Def. \ref{StandardModelOfUniversalSimplicialPrincipalComplex})
  is a Kan fibration, due to the right Quillen functor property \eqref{QuillenEquivalenceBetweenBorelModelStructureAndSliceOverSimplicialClassifyingComplex}:
  \vspace{-2mm}
  \begin{equation}
    \label{TheUniversalAssociatedSimplicialBundle}
    \begin{tikzcd}
      X
      \ar[
        r,
        "{\mathrm{fib}(q)}"
      ]
      &
      (W\mathcal{G} \times X)/\mathcal{G}
      \ar[
        d,
        "{\in \mathrm{Fib}}"{right},
        "q"{left}
      ]
      \\
      &
      \overline{W}\mathcal{G}
    \end{tikzcd}
    \;\;\;
    =
    \;\;\;
    \bigg(\!\!
    W\mathcal{G}
    \times
    \left(\!\!\!
     \scalebox{0.7}{$
      \def\arraystretch{.6}
      \begin{array}{c}
        X
        \\
        \downarrow
        \\
        \ast
      \end{array}
      $}
   \!\!\! \right)
   \!\!\! \bigg)
   \big/ \mathcal{G}
   \,.
  \end{equation}

 \vspace{-2mm}
\noindent  The fiber of this fibration, hence the {\it homotopy fiber}
  is clearly $X$.

  In the special case where $\mathcal{G} \acts \, X = \mathcal{G} \acts \, \mathcal{G}_L$
  is the multiplication action
  of the simplicial group on itself, this construction reduces
  to the {\it universal principal simplicial bundle} \eqref{StandardSimplicialDeloopingAsQuotient}
   \vspace{-2mm}
  \begin{equation}
    \label{TheUniversalPrincipalSimplicialBundle}
    \begin{tikzcd}
      \mathcal{G}
      \ar[r]
      &
      W\mathcal{G}
      \ar[
        d,
        "q"{left},
        "\in \mathrm{Fib}"{right}
      ]
      \\
      &
      \overline{W}\mathcal{G}
    \end{tikzcd}
    {\phantom{AAA}}
    =
    {\phantom{AAA}}
    \bigg(\!\!
    W\mathcal{G}
    \times
    \left(\!\!\!
    \scalebox{0.7}{$
      \def\arraystretch{.6}
      \begin{array}{c}
        \mathcal{G}_L
        \\
        \downarrow
        \\
        \ast
      \end{array}
      $}
   \!\!\! \right)
   \!\!\! \bigg) \big/ \mathcal{G}
   \;.
  \end{equation}
  By Prop. \ref{BasicPropertiesOfStandardSimplicialPrincipalComplex} this
  \eqref{TheUniversalPrincipalSimplicialBundle}
  exhibits $\mathcal{G}$ as the looping \eqref{LoopingInIntroduction} of $\overline{W}$,
  whereby \eqref{TheUniversalAssociatedSimplicialBundle} exhibits the simplicial model for the abstractly characterizedhomotopy quotient according to Prop. \ref{HomotopyQuotientsAndPrincipaInfinityBundles}.
\end{example}
Therefore it makes sense to use the following:
\begin{notation}[Homotopy quotient of simplicial group actions]
  \label{HomotopyQuotientOfSimplicialGroupActions}
  For $\mathcal{G} \in \Groups(\SimplicialSets)$,
  we denote the right derived functor of the
  {\it Borel construction} right Quillen functor
  \eqref{QuillenEquivalenceBetweenBorelModelStructureAndSliceOverSimplicialClassifyingComplex}
  by
   \vspace{-2mm}
  $$
    \HomotopyQuotient
      { (-) }
      { \mathcal{G} }
    \;\coloneqq\;
    \mathbb{R}
    \Big(\!\!
      \big(
        (-) \times W\mathcal{G}
      \big)/\mathcal{G}
    \Big)
    \;\;\;
      :
    \;\;\;
    \HomotopyCategory
    \big(
      \Actions{\mathcal{G}}(\SimplicialSets)_{\mathrm{proj}}
    \big)
    \longrightarrow
    \HomotopyCategory
    (
      \SimplicialSets_{\mathrm{Qu}}
    )
    \,.
  $$

\vspace{-2mm}
\noindent  Applied to the point $\mathcal{G} \acts \, \ast$, we also write
\vspace{-1mm}
  \begin{equation}
    \label{DeloopingAsHomotopyQuotientOfThePoint}
    \mathbf{B}\mathcal{G}
    \;\;
    \coloneqq
    \;\;
    \ast \!\sslash\! \mathcal{G}
    \,.
  \end{equation}
\end{notation}

As an example:
\begin{proposition}[$G$-Sets in the homotopy theory over $B G$]
  \label{GSetsInTheHomotopyTheoryOverBG}
  For $G \in \Groups(\Sets) \xhookrightarrow{\;} \Groups(\SimplicialSets)$,
  we have an equivalence between $G$-actions on sets and 0-truncated
  objects in the homotopy theory over $B G$ (Ntn. \ref{GroupoidNotation}):

  \vspace{-.4cm}
  $$
    \begin{tikzcd}[row sep=-6pt]
      \Actions{G}(\Sets)
      \ar[r, phantom, "\simeq"]
      &
      \left(
        (\InfinityGroupoids)_{/B G}
      \right)_{\leq 0}
      \ar[r, hook]
      &
      (\InfinityGroupoids)_{/B G}
      \\
    \scalebox{0.7}{$     G \acts \, X $}
      &\longmapsto&
     \scalebox{0.7}{$    \big( X \!\sslash\! G \xrightarrow{\;} B G  \big) $}
    \end{tikzcd}
  $$
  \vspace{-.4cm}

\end{proposition}
\begin{proof}
  Since $( X \times G \rightrightarrows X ) \xrightarrow{\;} (G \rightrightarrows \ast)$
  is clearly a Kan fibration with fiber $X$, the latter is the
  homotopy fiber of $X \!\sslash\! G \xrightarrow B G$.
  With this, the statement follows by Prop. \ref{QuillenEquivalenceBetweenBorelModelStructureAndSliceOverClassifyingComplex}.
\end{proof}

\medskip

\section{Geometric homotopy theory}
\label{ToposTheory}

We recall and develop basics of higher topos theory
\cite{Simpson99}\cite{Lurie03}\cite{ToenVezzosi05}\cite{Joyal08Logoi}\cite{Lurie09HTT}\cite{Rezk10}
with focus on the discussion of transformation groups in this context.

\medskip

-- \cref{SheavesAndStacks}: Simplicial sheaves and $\infty$-stacks.

-- \cref{GeneralConstructionsInInfinityToposes}: General constructions in $\infty$-toposes.

-- \cref{TransformationGroupsInInfinityToposes}: Transformation groups in $\infty$-toposes.

\medskip

\noindent
{\bf $\infty$-Toposes.} As usual, we say {\it $\infty$-topos} for
{$\infty$-categories of $\infty$-sheaves} (of $\infty$-stacks), in contrast to the broader
concept of {\it elementary $\infty$-toposes}.

\medskip
For our purposes, the reader may take the following characterization
of $\infty$-toposes to be their definition:

\begin{proposition}[$\infty$-Giraud theorem {\cite[Prop. 6.1.0.6]{Lurie09HTT}\cite[\S 2.2]{Lurie03}}]
\label{InfinityGiraudTheorem}
A presentable $\infty$-category
$\Topos \in \HomotopyTwoCategory(\PresentableInfinityCategories)$
(Def. \ref{HomotopyTwoCategoryOfPresentableHomotopyTheories},
Prop. \ref{HomotopyCategoryOfPresInfinityCategoriesIsThatOfCombinatorialModelCategories})
is an \emph{$\infty$-topos} if and only if it satisfies the following three conditions:

{\bf (a)} {\it coproducts are disjoint},

{\bf (b)} {\it colimits are universal} (see Prop. \ref{ColimitsAreUniversalInAnInfinityTopos}),

{\bf (c)}  {\it groupoid objects are effective} (see Prop. \ref{EquivalentPerspectivesOnGroupoidObjectsInAnInfinityTopos}).
\end{proposition}

\begin{example}[The canonical base $\infty$-topos of $\infty$-groupoids]
  \label{InfinityGroupoidsAsAnInfinityTopos}
  The $\infty$-category $\InfinityGroupoids$ of $\infty$-groupoids
  (Ntn. \ref{SimplicialSetsAndInfinityGroupoids}) is an
  $\infty$-topos (Prop. \ref{InfinityGiraudTheorem}), being the
  $\infty$-category $\infty$-sheaves on the point:
  \vspace{-2mm}
  $$
    \InfinityGroupoids
    \;\simeq\;
    \InfinitySheaves(\ast)
    \,.
  $$
\end{example}

\begin{proposition}[Terminal $\infty$-geometric morphism {\cite[Prop. 6.3.4.1 with Def. 6.1.0.4]{Lurie09HTT}}]
\label{TerminAlInfinityGeometricMorphism}
Given an $\infty$-topos $\Topos$ (Prop. \ref{InfinityGiraudTheorem})
there is an essentially unique pair of adjoint $\infty$-functors
\eqref{AdjunctionAndHomEquivalence} between $\Topos$ and
$\InfinityGroupoids$ (Ex. \ref{InfinityGroupoidsAsAnInfinityTopos})
  \vspace{-2mm}
\begin{equation}
  \label{TerminalGeometricMorphism}
  \begin{tikzcd}[column sep=large]
    \Topos
    \ar[
      rr,
      phantom,
      "{\scalebox{.6}{$\bot$}}",
    ]
    \ar[
      rr,
      shift right=5pt,
      "\GlobalSections"{below},
      "\mathclap{\scalebox{.6}{\colorbox{white}{\rm acc}}}"{description, pos=.05}
    ]
    &&
    \InfinityGroupoids \;,
    \ar[
      ll,
      shift right=5pt,
      "{ \LocallyConstant }"{above},
      "\mathclap{\scalebox{.6}{\colorbox{white}{\rm lex}}}"{description, pos=.05}
    ]
  \end{tikzcd}
\end{equation}

  \vspace{-2mm}
  \noindent
such that, in addition to the default (co-)limit preservation \eqref{InfinityAdjointPreservesInfinityLimits}:

{\bf (a)} the left adjoint $\LocallyConstant$ preserves finite $\infty$-limits;

{\bf (b)}  the right adjoint $\GlobalSections$ is accessible (\cite[Def. 5.4.2.5]{Lurie09HTT}).
\end{proposition}

\begin{remark}[Constructing locally constant $\infty$-stacks preserves finite products]
\label{ConstructingLocallyConstantInfinityStacksPreservesFiniteProduct}
That $\LocallyConstant$ preserves finite limits means in particular that
it preserves, via natural equivalences:

{\bf (a)}  the terminal object
\begin{equation}
  \label{LConstPreservesTerminalObject}
  \LocallyConstant(\ast) \;\simeq\; \ast_{\Topos}\;;
\end{equation}

{\bf (b)}  finite products
\begin{equation}
  \LocallyConstant
  (
    \mathcal{S}
      \times
      \cdots
      \times
    \mathcal{F}
  )
  \;\simeq\;
  \LocallyConstant(\mathcal{S})
    \times
    \cdots
    \times
  \LocallyConstant(\mathcal{S})
  \,.
\end{equation}
\end{remark}
\begin{example}[Global sections are co-represented by the terminal object]
\label{GlobalSectionsAreCoRepresentedByTheTerminalObject}
Global sections (Prop. \ref{TerminAlInfinityGeometricMorphism}) are equivalently
maps out of the terminal object:
\vspace{-2mm}
\begin{equation}
  \label{GlobalSectionsInInfinityToposAsMapsOutOfTerminalObject}
  \begin{aligned}
    \GlobalSections(-)
    \;\simeq\;
    \PointsMaps{}{ \ast }{-}
    \,.
  \end{aligned}
\end{equation}

\vspace{-2mm}
\noindent This is because:
$$
 \def\arraystretch{1.2}
  \begin{array}{lll}
    \GlobalSections(X)
    & \;\simeq\;
    \InfinityGroupoids
    \left(
      \ast
      ,\,
      \GlobalSections(X)
    \right)
    &
    \proofstep{ by Ex. \ref{PointsOfInfinityGroupoids} }
    \\
    & \;\simeq\;
   {\rm Grpd}_\infty \left(*, \Gamma(X)\right)
       &
    \proofstep{ by \eqref{TerminalGeometricMorphism} with \eqref{AdjunctionAndHomEquivalence} }
    \\
    & \;\simeq\;
    \PointsMaps{}
      { \ast }
      { X }
    &
    \proofstep{ by \eqref{LConstPreservesTerminalObject}. }
  \end{array}
$$
\end{example}

\begin{proposition}[Colimits are universal in an $\infty$-topos {\cite[p. 532]{Lurie09HTT}}]
\label{ColimitsAreUniversalInAnInfinityTopos}
In an $\infty$-topos $\Topos$ (Prop. \ref{InfinityGiraudTheorem}),
given a morphism $f \,;\, A \xrightarrow{\;} B$
and a diagram $X_\bullet \,:\, \mathcal{I} \xrightarrow{\;} \Slice{\Topos}{B}$,

\noindent {\bf (i)} pullback along $f$ preserves its $\infty$-colimit:
\vspace{-3mm}
\begin{equation}
  \label{UniversalityOfColimits}
  \begin{tikzcd}
    [row sep=small, column sep=large]
    \mathllap{
      f^\ast
      \big(
        \colimit{}
        \,
        X_\bullet
      \big)
      \simeq\;\;
    }
    \colimit{}
    \,
    f^\ast(X_\bullet)
    \ar[r]
    \ar[d, start anchor={[yshift=8pt]}]
    \ar[
      dr,
      phantom,
      "{ \mbox{\tiny\rm(pb)} }"{pos=.3}
    ]
    &
    \colimit{}
    \;
    X_\bullet
    \ar[
      d,
      start anchor={[yshift=4pt]}    ]
    \\
    A
    \ar[
      r,
      "f"
    ]
    &
    B
  \end{tikzcd}
\end{equation}

\vspace{-2mm}
\noindent {\bf (ii)} In particular, products preserve colimits:
\vspace{-3mm}
\begin{equation}
  \begin{tikzcd}
    [row sep=small, column sep=large]
    \mathllap{
      A \times
      \big(
        \colimit{}
        \,
        X_\bullet
      \big)
     \simeq\;\;
    }
    \colimit{}
    \,
    (A \times X_\bullet)
    \ar[r]
    \ar[d, start anchor={[yshift=8pt]}]
    \ar[
      dr,
      phantom,
      "{ \mbox{\tiny\rm(pb)} }"{pos=.3}
    ]
    &
    \colimit{}
    \;
    X_\bullet
    \ar[
      d,
      start anchor={[yshift=4pt]}    ]
    \\
    A
    \ar[
      r,
      "{}"
    ]
    &
    \ast
  \end{tikzcd}
\end{equation}
\end{proposition}

\medskip

\noindent
{\bf Groupoid objects in an $\infty$-topos.}

\begin{definition}[Groupoid objects in an $\infty$-topos {\cite[Def. 6.1.2.7]{Lurie09HTT}}]
  \label{GroupoidObjectInAnInfinityTopos}
  Given an $\infty$-topos $\Topos$ (Prop. \ref{InfinityGiraudTheorem}),

  \noindent
  {\bf (i)}
  a simplicial object
  \vspace{-1mm}
  $$
    \mathcal{X}_\bullet
    \;:\;
    \Delta^{\mathrm{op}}
      \longrightarrow
    \Topos
  $$
  is called a {\it groupoid object} if it satifies the
  groupoidal Segal conditions, hence  if for all $n \,\in\, \mathbb{N}$
  and all partitions (not necessarily order-preserving)
  \vspace{-3mm}
  $$
    \{0, 1, \cdots, n\}
    \;\simeq\;
    \{
      i_0, \cdots,  i_k
    \}
    \underset{ \{ i_k \} }{\sqcup}
    \{
      i_k, \cdots, i_n
    \}
  $$

  \vspace{-2mm}
 \noindent the corresponding diagram of evaluations of $\mathcal{X}_\bullet$
  is homotopy Cartesian
  \vspace{-2mm}
  \begin{equation}
    \label{SegalConditionsForGroupoidObjectInInfinityTopos}
    \begin{tikzcd}
      \mathcal{X}_{[n]}
      \ar[rr, "{ (i_0, \cdots, i_k) }^\ast"]
      \ar[d, "{ (i_k, \cdots, i_{n-k}) }^\ast"{left}]
      \ar[drr, phantom, "{ \mbox{\tiny \rm (pb)} }"]
      &&
      \mathcal{X}_{[k]}
      \ar[d, "{ (i_k)^\ast }"]
      \\
      \mathcal{X}_{[n-k]}
      \ar[rr, "{ (i_k)^\ast }"{below}]
      &&
      \mathcal{X}_{[0]}
      \,.
    \end{tikzcd}
  \end{equation}

   \vspace{-2mm}
 \noindent
  Here $(i_0, \cdots, i_k) \,\colon\, \Delta[k] \xrightarrow{\;} \Delta[n]$
  denotes the unique unique order-preserving map
  of finite sets, whose image contains
  the vertices $i_0, \cdots, i_k$.

\noindent {\bf (ii)}     We write
 \vspace{-2mm}
  $$
    \Groupoids(\Topos)
    \xhookrightarrow{\quad}
    \InfinityPresheaves(\Delta, \, \Topos)
  $$
  for the full sub-$\infty$-category of groupoid objects.
\end{definition}

\begin{example}[2-Horn fillers in groupoid objects]
  For $n = 2$ the condition \eqref{SegalConditionsForGroupoidObjectInInfinityTopos}
  says that
   \vspace{-2mm}
  \begin{align*}
   & \mathcal{X}_{\{0 <  1 <  2\}}
    \;\simeq\;
    \mathcal{X}_{\{0 < 1\}}
    \underset{
      \mathcal{X}_{\{1\}}
    }{\times}
    \mathcal{X}_{\{1 < 2\}}
    \,,
   \\
    \;\;
        \;\;
    &\mathcal{X}_{\{0 <  1 <  2\}}
    \;\simeq\;
    \mathcal{X}_{\{0 < 2\}}
    \underset{
      \mathcal{X}_{\{0\}}
    }{\times}
    \mathcal{X}_{\{0 < 1\}}
    \,,
    \;\;
   \\
    \mbox{and}
       \;\;
&    \mathcal{X}_{\{0 <  1 <  2\}}
    \;\simeq\;
    \mathcal{X}_{\{1 < 2\}}
    \underset{
      \mathcal{X}_{\{2\}}
    }{\times}
    \mathcal{X}_{\{0 < 2\}}\;,
  \end{align*}

  \vspace{-2mm}
 \noindent expressing the fact that morphisms
  of the groupoid objects (namely the elements of $\mathcal{X}_{[1]})$
  have

  (a) essentially unique composites when composable,

  (b) essentially unique left inverses,  and

  (c) essentially unique right inverses.
\end{example}

\begin{example}[Plain objects as constant groupoid objects in an $\infty$-topos]
  \label{PlainObjectsAsConstantGroupoidObjectsInAnInfinityTopos}
  For an object $X \,\in\, \Topos$, the simplicial object
  $\ConstantGroupoid(X)_\bullet$ which is
  constant on this
  value is a groupoid object (Def. \ref{GroupoidObjectInAnInfinityTopos}).
  This construction extends to a full-sub-$\infty$-category inclusion
  $$
    \ConstantGroupoid(-)_\bullet
    \;\colon\;
    \Topos
    \xhookrightarrow{\quad}
    \Groupoids(\Topos)
    \,.
  $$
  In particular, the groupoid constant on the terminal object $\ast \,\in\, \Topos$
  is the terminal groupoid object
  \begin{equation}
    \label{TerminalGroupoidObjectInAnInfinityTopos}
    \ConstantGroupoid(\ast)_\bullet
    \;
    \simeq
    \;
    \ast
    \;\;\;
    \in
    \;
    \Groupoids(\Topos)
    \,.
  \end{equation}
\end{example}

\begin{example}[{\v C}ech groupoids]
  \label{CechGroupoidsInAnInfinityTopos}
  Given a morphism $X \xrightarrow{\;p\;} \mathcal{X}$ in $\Topos$,
  its iterated fiber products with itself
  form a groupoid object (Def. \ref{GroupoidObjectInAnInfinityTopos}):
  $$
    \Cech(p)_\bullet
    \;\coloneqq\;
    X^{\times^\bullet_{\mathcal{X}}}
    \;\in\;
    \Groupoids(\Topos)
    \,.
  $$
\end{example}

\medskip
\noindent
{\bf Effective epimorphisms and atlases.}

\begin{definition}[Effective epimorphism in an $\infty$-category {\cite[\S 6.2.3]{Lurie09HTT}}]
  \label{EffectiveEpimorphisms}
  A morphism in an $\infty$-category with finite limits and
  simplicial colimits
  is called an {\it effective epimorphism},
  to be denoted in the form
  $\!\!\begin{tikzcd} X \ar[r, ->>,  "p"] & \mathcal{X} \end{tikzcd}\!\!$,
   if
  it is equivalently the $\infty$-colimit coprojection out of its
  {\v C}ech groupoid (Ex. \ref{CechGroupoidsInAnInfinityTopos})
  \vspace{-2mm}
  $$
    \begin{tikzcd}[row sep=15pt]
      X
      \ar[rr, -, shift left=1pt]
      \ar[rr, -, shift right=1pt]
      \ar[d, ->>, "{\, p }"]
      &&
      X
      \ar[d, "{\, q_0 }"]
      \\
      \colimit{ [n] \in \Delta^{\mathrm{op}} }
      \,
      X^{\times^n_{\mathcal{X}}}
      \ar[rr, "{ \sim }"]
      &&
      \mathcal{X}
    \end{tikzcd}
  $$

 \vspace{-3mm}
 \noindent
  In view of Prop. \ref{EquivalentPerspectivesOnGroupoidObjectsInAnInfinityTopos},
    we also say that an effective epimorphism $X \twoheadrightarrow \mathcal{X}$
  is an {\it atlas} for $\mathcal{X}$ (see {\cite[p. 27]{SS20OrbifoldCohomology}})
  and we write
   \vspace{-1mm}
  $$
    \EffectiveEpimorphisms(\Topos)
    \xhookrightarrow{\quad}
    \Presheaves(\Delta^1, \, \Topos)
  $$

   \vspace{-1mm}
 \noindent
  for the sub-$\infty$-category of the arrow category on the effective
  epimorphisms.
\end{definition}
\begin{notation}[Inhabited objects]
  \label{InhabitedObject}
  An object $X \,\in\, \Topos$ is called {\it inhabited}
  the morphism to the terminal object is an effective epimorphism:
  $X \twoheadrightarrow \ast$ (Def. \ref{EffectiveEpimorphisms}).
\end{notation}

\begin{proposition}[Detecting effective epimorphisms on 0-truncations {\cite[Prop. 7.2.1.14]{Lurie09HTT}}]
  \label{DetectingEffectiveEpimorphismsOnZeroTruncation}
  A morphism $f \,\colon\, X \xrightarrow{\;} \mathcal{X}$
  in an $\infty$-topos $\Topos$, is an effective epimorphism
  (Def. \ref{EffectiveEpimorphisms}), precisely if its
  0-truncation
  $\Truncation{0}(f) \,\colon\, \Truncation{0}(X)
  \xrightarrow{\;} \Truncation{0}(\mathcal{X})$
  is an effective epimorphism in $\Topos_0$.
\end{proposition}
\begin{example}[Pointed connected objects are those with point atlas]
  \label{PointedConnectedObjectsHaveEffectiveEpimorphicPointInclusion}
  In an $\infty$-topos $\Topos$,
  given a pointed object $\ast \xrightarrow{\mathrm{pt}} \mathcal{X}$,
  it is connected,
  in that $\Truncation{0}(\mathcal{X}) \,\simeq \, \ast$ (Ntn. \ref{ConnectedObject}),
  precisely if the pointing is
  an effective epimorphism (Def. \ref{EffectiveEpimorphisms}).
  Therefore, the pointed connected objects in $\Topos$
  form a full sub-$\infty$-category of the pointed objects in
  $\EffectiveEpimorphisms(\Topos)$:
   \vspace{-2mm}
  \begin{equation}
    \label{RegardingAPointedConnectedObjectAsAPointedAtlas}
    \begin{tikzcd}[row sep=-4pt]
      \Topos^{\ast/}_{\geq 0}
      \ar[rr, hook]
      &&
      \EffectiveEpimorphisms(\Topos)^{(\ast \to \ast)/}
      \\
     \scalebox{0.7}{$ \big(
        \ast \overset{\mathrm{pt}~}{\to} \mathcal{X}
      \big)
      $}
      &\longmapsto&
   \scalebox{0.7}{$   \left(
      \begin{array}{c}
        \ast
        \\
        \downarrow
        \\
        \ast
      \end{array}
      \begin{array}{c}
        \twoheadrightarrow
        \\
        {\phantom{\downarrow}}
        \\
        \overset{\mathrm{pt}}{\twoheadrightarrow}
      \end{array}
      \begin{array}{c}
        \ast
        \\
        \downarrow\mathrlap{\scalebox{.7}{$\mathrm{pt}$}}
        \\
        \mathcal{X}
      \end{array}
      \right)
      $}
    \end{tikzcd}
  \end{equation}
\end{example}
\begin{proof}
  By Prop. \ref{DetectingEffectiveEpimorphismsOnZeroTruncation}
  and using that $\Truncation{0}$ preserves the terminal object,
  it is sufficient to show that $\mathcal{X}$ is connected
  precisely if $\ast \xrightarrow{\; \Truncation{0}(\mathrm{pt})\; } \Truncation{0}(\mathcal{X})$
  is an effective epimorphism in the 1-category $\Topos_0$. But in
  a 1-category any fiber product of the terminal object with itself is
  the terminal object itself. Therefore, the {\v C}ech nerve of
  $\Truncation{0}(\mathrm{pt})$ in $\Topos_0$ is $\ConstantGroupoid(\ast)_\bullet$,
  whose colimit coprojection is $\ast \xrightarrow{\;} \ast$.
\end{proof}

\begin{proposition}[Groupoid objects in an $\infty$-topos are effective {\cite[Cor. 6.2.3.5]{Lurie09HTT}}]
  \label{EquivalentPerspectivesOnGroupoidObjectsInAnInfinityTopos}
  $\,$

  \noindent
  Given an $\infty$-topos $\Topos$,
  the operation of forming {\v C}ech groupoids
  (Ex. \ref{CechGroupoidsInAnInfinityTopos})
  constitutes an equivalence of $\infty$-categories
  \vspace{-2mm}
  \begin{equation}
    \label{EquivalenceBetweenEffectiveEpimorphismsAndGroupoidObjects}
    \begin{tikzcd}[column sep=large]
      \EffectiveEpimorphisms(\Topos)
      \ar[
        rr,
        shift right=5pt,
        "{ \Cech(-)_\bullet }"{below}
      ]
      \ar[rr, phantom, "\sim"]
      &&
      \Groupoids(\Topos)
      \ar[
        ll,
        shift right=5pt,
        "{
          (-)_0 \; \to \; \colimit{}(-)
        }"{above, yshift=-4pt}
      ]
    \end{tikzcd}
  \end{equation}
  between effective epimorphisms
  and groupoid objects (Def. \ref{GroupoidObjectInAnInfinityTopos})
\end{proposition}

\begin{proposition}[The effective-epi/mono factorization system]
  \label{TheEffectiveEpiMonoFactorizationSystem}
  In an $\infty$-topos $\Topos$,
  the classes
  of $(-1)$-truncated morphisms
  (hence of monomorphisms \eqref{InfinityMonomorphism})
  and
  of $(-1)$-connected morphisms
  (hence of
  effective epimorphisms (Def. \ref{EquivalentPerspectivesOnGroupoidObjectsInAnInfinityTopos})),
  form an orthogonal factorization system.
  In particular:

  \noindent
  {\bf (i)}
  Diagrams of the following form have essentially unique lifts
  \vspace{-2mm}
  \begin{equation}
    \label{MinusOneConnectedTruncatedLiftingProblem}
    \begin{tikzcd}[column sep=large]
      X
      \ar[d, ->>]
      \ar[r, "\forall"]
      &
      A
      \ar[d, hook]
      \\
      \mathcal{X}
      \ar[r, "\forall"{below}]
      \ar[ur, dashed, "\exists !"]
      &
      \mathcal{A}
    \end{tikzcd}
  \end{equation}

  \vspace{-2mm}
  \noindent
  {\bf (ii)}
  An $\infty$-limit in the arrow category over $(-1)$-truncated morphisms
  is again $(-1)$-truncated
  \vspace{-2mm}
  \begin{equation}
    \label{LimitOfMonomorphismsInArrowCategoryIsMono}
    \limit{i \in I}
    \big(
      X_i \xhookrightarrow{\; f_i \;} Y_i
    \big)
    \;\;\;\simeq\;\;\;
    \bigg(\!\!
    \big(
      \limit{i \in I}
      X_i
    \big)
    \xhookrightarrow{
      \;\;
       \scalebox{0.7}{$   \limit{i \, \in \, I} f_i $}
      \;\;
    }
    \big(
      \limit{i \in I}
      Y_i
    \big)
    \!\!\bigg)
    \,.
  \end{equation}
\end{proposition}

\medskip

\noindent
{\bf Shape theory of $\infty$-toposes.}

\begin{definition}[Pro-$\infty$-groupoids {\cite[Def. 7.1.6.1]{Lurie09HTT}}]
  We write
\vspace{-2mm}
  \begin{equation}
    \label{ProInfinityGroupoids}
    \ProObjects\InfinityGroupoids
    \xhookrightarrow{\quad}
    \Functors(\InfinityGroupoids,\InfinityGroupoids)^{\mathrm{op}}
  \end{equation}

  \vspace{-2mm}
\noindent  for the opposite of the
  full sub-$\infty$-category on those $\infty$-functors
  which preserve finite $\infty$-limits and $\kappa$-filtered $\infty$-colimits
  for some regular cardinal $\kappa$.
\end{definition}

\begin{lemma}[Factorization through pro-$\infty$-groupoids]
  The $\infty$-Yoneda embedding \eqref{InfinityYonedaEmbedding}
  of $\InfinityGroupoids$ factors through pro-$\infty$-groupoids
  \eqref{ProInfinityGroupoids}:
    \vspace{-2mm}
  \begin{equation}
    \label{InfinityGroupoidsInsideProInfinityGroupoids}
    \begin{tikzcd}[row sep=-3pt, column sep=small]
      \InfinityGroupoids
      \ar[
        rr,
        hook,
        "{\YonedaEmbedding}"
      ]
      &&
      \ProObjects\InfinityGroupoids \;.
      \\
     \scalebox{.7}{$ S  $}
      &\longmapsto&
    \scalebox{.7}{$  \InfinityGroupoids(S,-) $}
    \end{tikzcd}
  \end{equation}
\end{lemma}
\begin{proof}
  By cartesian closure,
  $\infinityGroupoids(S,-) \colon \InfinityGroupoids \xrightarrow{\;} \InfinityGroupoids$
  is a right adjoint
  and as such preserves finite $\infty$-limits
  \eqref{InfinityAdjointPreservesInfinityLimits}
  and $\kappa$-filtered colimits for some regular cardinal $\kappa$
  (by \cite[Prop. 5.4.7.7]{Lurie09HTT}).
\end{proof}

\begin{definition}[Shape of an $\infty$-topos {\cite[Def. 7.1.6.3]{Lurie09HTT}}]
\label{ShapeOfAnInfinityTopos}
  The {\it shape} of an $\infty$-topos
of the form  \eqref{TerminalGeometricMorphism}
  is the pro-$\infty$-groupoid
  \eqref{ProInfinityGroupoids} given by
    \vspace{-3mm}
  \begin{equation}
    \label{FormulaForLurieShapeOfAnInfinityTopos}
    \Shape(\mathbf{H})
    \;\;\coloneqq\;\;
    \GlobalSections(-)
    \circ
    \LocallyConstant
    \;\;\;
    \in
    \;
    \ProObjects\InfinityGroupoids
    \,.
  \end{equation}
\end{definition}
\begin{remark}[Shape of $\infty$-topos is well-defined]
  The formula \eqref{FormulaForLurieShapeOfAnInfinityTopos}
  is well-defined, in that $\GlobalSections \circ \LocallyConstant$
  preserves finite $\infty$-limits and $\kappa$-filtered $\infty$-colimits,
  by \eqref{TerminalGeometricMorphism} and \eqref{InfinityAdjointPreservesInfinityLimits}.
\end{remark}

\begin{proposition}[Shape of $\infty$-topos is {\'e}tale pro-homotopy type of its terminal object
  {\cite[Def. 2.3]{Hoyois15}}]
  \label{ShapeOfInfinityToposIsEtaleProHomotopyTypeOfItsTerminalObject}
  The shape of an $\infty$-topos in the sense of
  Def. \ref{ShapeOfAnInfinityTopos}
  is equivalently the \emph{{\'e}tale pro-homotopy type}
  of its terminal object $\ast_{\Topos} \,\in\, \Topos$,
  namely the image of the terminal object
  under the pro-left adjoint to $\LocallyConstant$:
  \vspace{-2mm}
  \begin{equation}
    \label{ProLeftAdjointOfLConst}
    \begin{tikzcd}[row sep=-4pt, column sep=small]
      \Shape
      \;:
       &[-12pt]
       \Topos
       \ar[rr]
       &&
       \InfinityGroupoids
       \\
       &
  \scalebox{0.7}{$     X $}
       &\longmapsto&
    \scalebox{0.7}{$     \PointsMaps{}
         { X }
         { \LocallyConstant(-) }
         $}
    \end{tikzcd}
  \end{equation}

  \vspace{-3mm}
\noindent
 in that:
 \vspace{-1mm}
  $$
    \overset{
      \mathclap{
      \raisebox{3pt}{
        \tiny
        \color{darkblue}
        \bf
        \begin{tabular}{c}
          shape of
          \\
          $\infty$-topos
        \end{tabular}
      }
      }
    }{
      \Shape(\Topos)
    }
    \;\;\;\;\;
    \simeq
    \;\;\;\;\;
    \overset{
      \mathclap{
      \raisebox{3pt}{
        \tiny
        \color{darkblue}
        \bf
        \begin{tabular}{c}
          {\'e}tale pro-homotopy type
          \\
          of terminal object
        \end{tabular}
      }
      }
    }
    {
      \Shape(\ast_{\Topos})
    }
    \;\;\;
    \in
    \;
    \ProObjects\InfinityGroupoids
    \,.
  $$
\end{proposition}
\begin{proof}
This is the composite of the following sequence of
natural equivalences:
\vspace{-2mm}
$$
  \def\arraystretch{1.2}
  \begin{array}{lll}
    \Shape(\Topos)
    &
    \;\coloneqq\;
    \GlobalSections \circ \LocallyConstant(-)
    &
    \mbox{\small by Def. \ref{ShapeOfAnInfinityTopos} }
    \\
    & \;\simeq\;
    \InfinityGroupoids
    \left(
      \ast
      ,\,
      \GlobalSections \circ \LocallyConstant(-)
    \right)
    &
    \mbox{\small  by Prop. \ref{InfinityGroupoidsFormCartesianClosedInfinityCategory} }
    \\
    & \;\simeq\;
     \mathbf{H}\left(
      { \LocallyConstant(\ast) },
      { \LocallyConstant(-) }
      \right)
    &
    \mbox{\small  by \eqref{TerminalGeometricMorphism} }
    \\
    & \;\simeq\;
   \mathbf{H}\left(
      { \ast_{\Topos} },
      { \LocallyConstant(-) }
      \right)
    &
    \mbox{\small  by \eqref{LConstPreservesTerminalObject} }
    \\
    &
    \;=:\;
    \Shape(\ast_{\Topos})
    &
    \mbox{\small  by \eqref{ProLeftAdjointOfLConst} }
    .
  \end{array}
$$

\vspace{-7mm}
\end{proof}

\subsection{Simplicial Sheaves and $\infty$-stacks}
\label{SheavesAndStacks}

With $\infty$-groupoids presented by
simplicial sets, $\infty$-stacks are presented by simplicial presheaves,
an observation that really goes back to \cite{Brown73}.
A textbook account with focus on the injective model structures
on simplicial presheaves
may be found in \cite{Jardine15}.
Here we need to focus on their projective model structures whose
role was first highlighted in \cite{Dugger01Universal}.
A more detailed review of the specific constructions of interest here may be
found in \cite[pp. 105]{FSS20CharacterMap}.

\medskip

\noindent
{\bf Simplicial sites and $\infty$-Sites.}
  When discussing $\infty$-toposes
  we need some would-be $\infty$-categories which need not and typically are not presentable, and as such are outside the scope of our defining Prop. \ref{HomotopyCategoryOfPresInfinityCategoriesIsThatOfCombinatorialModelCategories}.
  These are the {\it $\infty$-sites}
  (Def. \ref{InfinitySite})
  over with $\infty$-(pre-)sheaves are defined and which hence serve to present $\infty$-toposes (see Prop. \ref{PresentationOfInfinityToposesBySimplicialPresheaves} below).
  But in fact, their presentation of $\infty$-toposes is practically the only purpose of $\infty$-sites, and for this it is sufficient to model the $\infty$-site as a simplicially enriched (preferably Kan complex-enriched) category (Ntn. \ref{SimplicialCategories}). Last not least, the $\infty$-Yoneda lemma (Prop. \ref{InfinityYonedaLemma}) says that every $\infty$-site in this sense is a full  sub-$\infty$-category of its associated $\infty$-topos, hence of a presentable $\infty$-category (where either the $\infty$-site is small or else the ambient presentable $\infty$-category is ``very large'').

Therefore:
\begin{definition}[Simplicial sites and $\infty$-sites {\cite[Def. 3.1.1]{ToenVezzosi05}}]
\label{SimplicialSite}
A {\it simplicial site} $(\mathcal{S}, J)$ is a simplicial category $\mathcal{S} \,\in\, \SimplicialCategories$ (Ntn. \ref{SimplicialCategories}) equipped with a Grothendieck topology on its homotopy category \eqref{HomotopyCategoryOfASimplicialCategory}.
\end{definition}
\begin{notation}[$\infty$-Sites]
\label{InfinitySite}
In order to notationally indicate that a construction depends only on the corresponding $\infty$-category of an $\infty$-site, we write $\InfinitySite$ for it, as in \eqref{SimplicialLocalizationOfLocalModelCategoryOfSimplicialPresheaves} below.
\end{notation}
\begin{example}[$\infty$-Sites from full sub-$\infty$-categories of presentable $\infty$-categories]
  \label{InfinitySitesFromFullSubInfinityCategoriesOfPresentableInfinityCategories}
  For $\mathbf{C} \,\in\, \PresentableInfinityCategories$,  Prop. \ref{HomotopyCategoryOfPresInfinityCategoriesIsThatOfCombinatorialModelCategories} provides a presenting simplicial combinatorial model category $\mathcal{C}$.

  Hence for any class of objects of $\mathcal{C}$ we obtain an $\infty$-site (with trivial Grothendieck topology) given by the corresponding simplicial full subcategory
  $$
    \mathcal{S}
    \;\hookrightarrow\;
    \mathcal{C}_{\mathrm{fib}}^{\mathrm{cof}}
    \;\hookrightarrow\;
    \mathcal{C}
  $$
  of the full simplicial subcategory of fibrant- and cofibrant objects of $\mathcal{C}$
  (i.e. with hom-$\infty$-groupoids as in Def. \ref{HomInfinityGroupoid}).
  The equivalence class of the simplicial categories $\mathcal{S}$ arising this way as the presentation $\mathcal{C}$ varies up to equivalence in
  $
    \TwoLocalization{\QuillenEquivalences}
    (
      \SimplicialCombinatorialModelCategories
    )
  $
  is labeled by the corresponding $\infty$-site $\mathbf{S}$ (Ntn. \ref{InfinitySite}).
\end{example}

The main instance of this class examples (Exp. \ref{InfinitySitesFromFullSubInfinityCategoriesOfPresentableInfinityCategories}) that is of interest below is the {\it 2-category of orbi-singularities} (or {\it global orbit category}), see Ntn. \ref{Singularities}.
A variation of this example is the orbit category of a topological group, which will appear for us
(in Def. \ref{ProperTopologicalOrbitInfinityCategory} below) as a special case of the following variant of the above construction:

\begin{example}[$\infty$-Sites from Cartesian closed presentable $\infty$-categories]
  \label{InfinitySitesFromCartesianClosedInfinityCategories}

  If $\mathbf{C} \,\in\, \PresentableInfinityCategories$ is Cartesian closed (as in Prop. \ref{MappingStacks} below) with mapping object $\infty$-functor
  $$
    \Maps{}{(-)}{(-)}
    \;:\;
    \mathbf{C}^{\mathrm{op}}
    \times
    \mathbf{C}
    \xrightarrow{\;}
    \mathbf{C}
  $$
  and if
  $$
    R
    \;:\;
    \mathbf{C}
    \xrightarrow{\;}
    \InfinityGroupoids
  $$
  is any $\infty$-functor which preserves products, then we obtain an $\infty$-site $\mathbf{S}$ with objects any sub-class of that of $\mathbf{C}$ and
  mapping complex for for $X$, $Y$ a pair of objects given by
  $$
    \mathbf{S}(X,Y)
    \;:=\;
    R\big(\Maps{}{X}{Y}\big)
  $$
  and with composition given by
  $$
    \mathbf{S}(X,Y)
    \times
    \mathbf{S}(Y,Z)
    \;\simeq\;
    R\big(\Maps{}{X}{Y}\big)
    \times
    R\big(
      \Maps{}{Y}{Z}
    \big)
    \;\simeq\;
    R
    \big(
      \Maps{}{X}{Y}
      \times
      \Maps{}{Y}{Z}
    \big)
    \xrightarrow{
      R\big((-)\circ(-)\big)
    }
    \Maps{}{X}{Z}
    \,.
  $$
\end{example}

\medskip

\noindent
{\bf Simplicial presheaves and $\infty$-stacks.}
\begin{notation}[Model categories of simplicial presheaves]
  \label{ModelCategoriesOfSimplicialPresheaves}
  Given a simplicial site $(\SimplicialSite, J)$
  (Def. \ref{SimplicialSite}), we write:

\vspace{1mm}
  \noindent
  {\bf (i)}
  $
    \SimplicialPresheaves(\SimplicialSite)
    \;\in\;
    \SimplicialCategories
  $
  for the simplicial category of simplicial presheaves on $\SimplicialSite$,

  \vspace{1mm}
  \noindent
  {\bf (ii)}
  $
    \SimplicialPresheaves(\SimplicialSite)_{ \proj }
    \;\in\;
    \ProperCombinatorialSimplicialModelCategories
  $
  for the global projective model category structure,

\vspace{1mm}
  \noindent
  {\bf (iii)}
  $
    \SimplicialPresheaves(\SimplicialSite, J)_{ \projloc }
    \;\in\;
    \ProperCombinatorialSimplicialModelCategories
  $
  for the $J$-local projective model category structure,
  whose weak equivalences are the $J$-stalk-wise weak equivalences in
  $\SimplicialSets_{\mathrm{Qu}}$ (the ``hypercomplete'' local model structure).
\end{notation}
For our purposes, the following may be taken to be the definition of $\infty$-toposes:
\begin{proposition}[Presentation of $\infty$-stacks by simplicial presheaves {(\cite{ToenVezzosi05}\cite{Lurie09HTT}, following \cite{Dugger01Universal}\cite{Brown73}})]
  \label{PresentationOfInfinityToposesBySimplicialPresheaves}
  Let $(\mathcal{S}, J)$ be a simplicial site.

  \noindent {\bf (i)} The
  {\v C}ech/stalk-local
  projective model category structure on
  simplicial presheaves over $\SimplicialSite$ presents
  the topological/hypercomplete $\infty$-topos over $\SimplicialSite$:
  \vspace{-2mm}
  \begin{equation}
    \label{SimplicialLocalizationOfLocalModelCategoryOfSimplicialPresheaves}
    \SimplicialPresheaves(\SimplicialSite)_{\projloc}
    \xhookrightarrow{\;\; \Localization{\LocalWeakEquivalences} \;\;}
    \SimplicialLocalization{\LocalWeakEquivalences}
    \big(
      \SimplicialPresheaves(\SimplicialSite, J)_{\projloc}
    \big)
    \;\simeq\;
    \infinitySheaves(\InfinitySite,J)
    \;\;\;
    \in
    \;
    \HomotopyTwoCategory
    (
      \PresentableInfinityCategories
    )
    \,,.
  \end{equation}

  \vspace{-1mm}
  \noindent
  {\bf (ii)}
  In particular, for all cofibrant $X$ and fibrant $A$
  in $\SimplicialPresheaves(\mathcal{C})_{\projloc}$ there is a natural equivalence of
  hom-$\infty$-groupoids
  \vspace{-1mm}
  \begin{equation}
    \SimplicialPresheaves(\SimplicialSite)
    (
      X,\, A
    )
    \;\simeq\;
    \InfinitySheaves
    \left(
      \Localization{\LocalWeakEquivalences}(X)
      ,\,
      \Localization{\LocalWeakEquivalences}(A)
    \right)
    \,.
  \end{equation}
\end{proposition}
\begin{notation}[Presentation of $\infty$-presheaves by simplicial presheaves]
\label{PresentationOfInfinityPresheavesBySimplicialPresheaves}
When the Grothendieck topology $J$ on $\SimplicialSite$
(in Prop. \ref{PresentationOfInfinityToposesBySimplicialPresheaves})
is trivial, we write, as usual,
$$
  \InfinityPresheaves(\InfinitySite)
  \;:=\;
  \InfinitySheaves(\InfinitySite, \mathrm{triv})
  \;\;\;
  \in
  \;
  \HomotopyTwoCategory
  (
    \PresentableInfinityCategories
  )
$$
for the $\infty$-category of $\infty$-presheaves over the $\infty$-category $\InfinitySite$.
\end{notation}
\begin{proposition}[$\infty$-Yoneda lemma {\cite[Prop. 8.2.1.3, Thm. 8.2.5.4]{LurieYonedaLemma}\cite[Prop. 5.1.3.1, 5.5.2.1]{Lurie09HTT}\cite[Thm. 7.2.22]{RiehlVerity09}}]
  \label{InfinityYonedaLemma}
  Let $\SimplicialSite \,\in\, \SimplicialCategories$,
  not necessarily small.

  \noindent
  {\bf (i)}
  The assignment of representable $\infty$-presheaves
  is a fully faithful embedding \eqref{FullyFaithfulInfinityFunctor}
  \vspace{-2mm}
  \begin{equation}
    \label{InfinityYonedaEmbedding}
    \begin{tikzcd}[row sep=-4pt, column sep=small]
      \InfinitySite
      \ar[
        rr,
        hook,
        "\YonedaEmbedding"
      ]
      &&
      \InfinityPresheaves(\InfinitySite)\;.
      \\
      \scalebox{0.7}{$U$}
        &\longmapsto&
      \scalebox{0.7}{$ \InfinitySite(-,\,U) $}
    \end{tikzcd}
  \end{equation}

  \vspace{-3mm}
  \noindent
  {\bf (ii)}
  For $U  \,\in\, \InfinitySite$ and $X \in \InfinityPresheaves(\InfinitySite)$,
  \vspace{-2mm}
  \begin{equation}
    \label{YonedaEquivalence}
    \mathrm{PSh}_\infty
    \left(
      y(U),
      \,
      X
    \right)
    \;\simeq\;
    X(U)
    \,.
  \end{equation}
\end{proposition}

\medskip
\noindent
{\bf Handling $\infty$-stacks via presenting simplicial presheaves.}
We discuss a few techniques for reasoning about $\infty$-stacks
-- on the right of \eqref{SimplicialLocalizationOfLocalModelCategoryOfSimplicialPresheaves} -- by constructions in the 1-category of simplicial presheaves -- on the left of \eqref{SimplicialLocalizationOfLocalModelCategoryOfSimplicialPresheaves}.

\begin{lemma}[Objectwise connected $\infty$-sheaves are connected]
 \label{ObjectwiseConnectedInfinitySheavesAreConnected}
 A sufficient (but not necessary) condition for an
 object $X \,\in\, \InfinitySheaves(\InfinitySite)$
 to be connected is that it is connected as an object of
 $\InfinityPresheaves(\InfinitySite)$, i.e., that
 for all $U \,\in\,\InfinitySite$ we have that
 $X(U) \,\in\, \InfinityGroupoids$ is connected.
\end{lemma}
\begin{proof}
  Consider a presentation of $\InfinitySheaves(\InfinitySite)$
  by a model category of simplicial presheaves.
  Then 0-truncation $\tau_0$ is the derived functor
  of objectwise truncation of simplicial sets followed
  by ordinary sheafification. But under the given assumption,
  the objectwise truncation is the presheaf with constant value the
  singleton set. This is a sheaf, the terminal sheaf
  (by the fact that inclusion into presheaves preserves the empty limit),
  hence is already the derived truncation operation in question.
\end{proof}

\begin{lemma}[Simplicial presheaves model homotopy colimit of their component sheaves]
  \label{SimplicialPresheavesModelHomotopyColimitOfTheirComponentSheaves}
  The localization
  \eqref{SimplicialLocalizationOfLocalModelCategoryOfSimplicialPresheaves}
  of any
  $X_\bullet \,\in\, \SimplicialPresheaves(\SimplicialSite)$
  is equivalently the homotopy colimit of its
  components
  $
    X_n
      \,\in\,
    \Presheaves(\SimplicialSite)
      \xhookrightarrow{\mathrm{const}}
    \SimplicialPresheaves(\SimplicialSite)
  $:
  $$
    \Localization{\Local\WeakEquivalences}(X_\bullet)
    \;\;
    \simeq
    \;\;
    \colimit{[n] \in \Delta^\op}
    \Localization{\LocalWeakEquivalences}(X_n)
    \;\;\;
    \in
    \;
    \SimplicialLocalization{\LocalWeakEquivalences}
    \big(
      \SimplicialPresheaves(\SimplicialSite)_{\projloc}
    \big)
    \,.
  $$
\end{lemma}
\begin{proof}
  The analogous statement for simplicial sets (i.e. for $\SimplicialSite = \ast$)
  is classical. This implies the statement for the global projective
  model structure since its homotopy colimits are computed objectwise.
  The full statement now follows since Bousfield localization to the
  local model structure is a (Quillen/derived/$\infty$-) left adjoint.
\end{proof}

The following Lem \ref{ComputingHomotopyPullbacksOfInfinityStacks}
is clearly a special case of a much more general statement,
but we highlight it in this form for definiteness, as we will have crucial
use of this, for instance in the proof of Prop. \ref{EquivariantPrincipalBundlesFromCechCocycles} below:

\begin{lemma}[Computing homotopy pullbacks of $\infty$-stacks]
  \label{ComputingHomotopyPullbacksOfInfinityStacks}
  Given a simplicial site
  $(\SimplicialSite, J)$
  and a diagram of simplicial presheaves
  (Ntn. \ref{ModelCategoriesOfSimplicialPresheaves})
  of the form
  \vspace{-3mm}
  \begin{equation}
    \label{HomotopyPullbackInSimplicialPresheaves}
    \begin{tikzcd}[row sep=small]
      &&
      P
      \ar[d, "\in \ProjectiveWeakEquivalences"]
      \\
      X \times_A \widehat P
      \ar[rr]
      \ar[d]
      \ar[drr, phantom, "\mbox{\tiny\rm(pb)}"]
      &&
      \widehat{P}
      \ar[
        d,
        "\in \ProjectiveFibrations"
      ]
      \\
      X
      \ar[rr]
      &&
      A
    \end{tikzcd}
    \;\;\;
    \in
    \;
    \SimplicialPresheaves(\SimplicialSite)
  \end{equation}
  then the homotopy pullback
  \eqref{HomotopyCartesianSquare}
  of $P$ with respect to the \emph{local} projective model
  structure $\SimplicialPresheaves(\SimplicialSite,J)_{\projloc}$
  is presented already by the homotopy pullback in the
  \emph{global} projective model
  structure, which in turn is presented by the ordinary pullback
  of any projective fibration resolution $\widehat P$:
  $$
    \Localization{\LocalWeakEquivalences}
    (
      X \times_A \widehat P
    \,)
    \;\;\;
    \simeq
    \;\;\;
    \Localization{\LocalWeakEquivalences}(X)
    \underset
      {
        \scalebox{.7}{$
          \Localization{\LocalWeakEquivalences}(A)
        $}
      }
      { \times }
    \Localization{\LocalWeakEquivalences}(P)
    \;\;\;\;\;
    \in
    \;
    \SimplicialLocalization{\WeakEquivalences}
    \left(
      \SimplicialPresheaves(\SimplicialSite, J)
    \right).
  $$
\end{lemma}
\begin{proof}
  The point is that the left Bousfield localization
  of the model categories has underlying identity functors
  but models the lex localization of presentable $\infty$-categories
  which preserves finite $\infty$-limits, in particular homotopy pullbacks \eqref{HomotopyCartesianSquare}:
  \vspace{-2mm}
  $$
    \SimplicialLocalization{\WeakEquivalences}
    \Big(\!\!\!
    \begin{tikzcd}
      \SimplicialPresheaves(\SimplicialSite, J)_{\projloc}
      \ar[from=r, shift right=5pt, "\mathrm{id}"{swap}]
      \ar[r, shift right=5pt, "\mathrm{id}"{swap}]
      \ar[r, phantom, "\scalebox{.7}{$\bot$}"]
      &
      \SimplicialPresheaves(\SimplicialSite, J)_{\proj}
    \end{tikzcd}
    \!\!\!\Big)
    \;\;
      \simeq
    \;\;
    \Big(\!\!\!
    \begin{tikzcd}
      \InfinitySheaves\left(\SimplicialLocalization{\WeakEquivalences}(\SimplicialSite)\right)
      \ar[from=r, shift right=5pt, "\mathrm{lex}"{swap}]
      \ar[r, shift right=5pt, hook]
      \ar[r, phantom, "\scalebox{.7}{$\bot$}"]
      &
      \InfinityPresheaves\left(\SimplicialLocalization{\WeakEquivalences}(\SimplicialSite)\right)
    \end{tikzcd}
    \!\!\!\!\!\Big)
    \,.
  $$

  \vspace{-2mm}
  \noindent
  This shows that the local homotopy pullback is given,
  by the image of the
  global homotopy pullback under the left derived functor of the identity,
  hence by a global cofibrant resolution of the global homotopy pullback.
  But left Bousfield localization does not change the class of cofibrations
  (just the class of acyclic cofibrations) so that
  global cofibrant resolution is implied by local cofibrant resolution,
  so that, under simplicial localization, the plain global homotopy
  pullback does present the local homotopy pullback.

  Hence it just remains to see that for computing the global homotopy
  pullback of a cospan diagram of simplicial presheaves it is sufficient
  to resolve one of the two morphisms by a global fibration.
  Since the global projective model strcuture is evidently right proper
  (since $\SimplicialSets_{\mathrm{Qu}}$ is so),
  this follows by a classical result (e.g. \cite[Prop. A.2.4.4]{Lurie09HTT}).
\end{proof}
This Lemma \ref{ComputingHomotopyPullbacksOfInfinityStacks} is particularly useful since there is a convenient way to achieve the factorization in \eqref{HomotopyPullbackInSimplicialPresheaves}:
\begin{notation}[Canonical projective path object for simplicial presheaves]
\label{CanonicalProjectivePathObjectForSimplicialPresheaves}
For $A \,\in\, \SimplicialPresheaves(\mathcal{S})$
(Ntn.
ref{(Ntn. \ref{ModelCategoriesOfSimplicialPresheaves})}), we write
\begin{equation}
  \label{CanonicalProjectivePathObjectViaFunctionComplex}
  A^I
  \,:\,
  s \,\mapsto\,
  \big(
    A(s)
  \big)^I
  \,:=\,
  \Maps{\big}{ \Simplex{1} }{A(s)}
  \,:=\,
  \SimplicialSets
  \big(
    (-)\times \Simplex{1}
    ,\,
    A(s)
  \big)
  \;\;\;
  \in
  \;
  \SimplicialSets
\end{equation}
for the simplicial presheaf $A^I \,\in\, \SimplicialPresheaves(\mathcal{S})$ which assigns to any object $s \,\in\, \mathcal{S}$ the simplicial function complex from $\Simplex{1}$ into the value $A(s) \,\in\, \SimplicialSets$ of the presheaf on that object.

This is evidently a {\it path space object} for $A$
(in the sense of \cite[Def. I.4]{Quillen67}, review in \cite[Def. A.0.12]{FSS20CharacterMap})
in that with the evident face and degeneracy maps it sits in a diagram of simplicial presheaves of this form:
\begin{equation}
  \label{CanonicalProjectivePathObjectFactorization}
  \begin{tikzcd}
    A
    \ar[
      rr,
      "\sigma",
      "{\in \ProjectiveWeakEquivalences}"{swap}
    ]
    \ar[
      rrrr,
      rounded corners,
      to path={
           ([yshift=-0pt]\tikztostart.south)
        -- ([yshift=-10pt]\tikztostart.south)
        -- node[yshift=5pt]{\scalebox{.7}{$\mathrm{diag}$}}
           ([yshift=-9pt]\tikztotarget.south)
        -- ([yshift=-0pt]\tikztotarget.south)
      }
    ]
    &&
    A^I
    \ar[
      rr,
      "{
        (\mathrm{ev}_0, \mathrm{ev}_1)
      }",
      "{
        \in \ProjectiveFibrations
      }"{swap}
    ]
    &&
    A \times A
  \end{tikzcd}
\end{equation}
\end{notation}
\begin{lemma}[Factorization lemma for projective simplicial presheaves (specializing {\cite[p. 421]{Brown73}})]
  \label{FactorizationLemmaForProjectiveSimplicialPresheaves}
  Given a morphism
  $$
    P \xrightarrow{\;f\;} A
    \;\;\;
    \in
    \;
    \big(
      \SimplicialPresheaves(\SimplicialSite)_{\proj}
    \big)_{\mathrm{fib}}
  $$
  of projectively fibrant simplicial presheaves
  (Ntn. \ref{ModelCategoriesOfSimplicialPresheaves}), the following pullback
  (formed
  by using the path space object $A^I$ from Ntn. \ref{CanonicalProjectivePathObjectForSimplicialPresheaves})
  exhibits a factorization of this morphism through a weak equivalence followed by a fibration in $\SimplicialPresheaves(\mathcal{S})_{\proj}$:
  \begin{equation}
    \label{FactorizationLemmaDiagram}
    \begin{tikzcd}
      P
      \ar[
        dr,
        dashed,
        shorten >=-5pt,
        "{
          (\mathrm{id}_P, \sigma_A \circ f)
        }"{sloped, pos=.66},
        "{
          \in \ProjectiveWeakEquivalences
        }"{swap, sloped, pos=.6}
      ]
      \ar[
        drr,
        bend left=20,
        "{
          \mathrm{id}_P
        }"
      ]
      \ar[
        dddr,
        bend right=20,
        rounded corners,
        to path={
               ([yshift=-4pt]\tikztostart.south)
            -- node[xshift=-6pt]{
               \scalebox{.7}{$f$}
            }
            ([xshift=-45pt]\tikztotarget.west)
            -- ([xshift=-0pt]\tikztotarget.west)
         },
        "{f}"{description}
      ]
      \\
      &
      P \underset{A}{\times} A^I
      \ar[r]
      \ar[d]
      \ar[
        dr,
        phantom,
        "{\mbox{\tiny\rm(pb)}}"{pos=.33}
      ]
      \ar[
        dd,
        shorten <=-8pt,
        shorten >=-4pt,
        shift right=6pt,
        bend right=20pt,
        "{
          \in \ProjectiveFibrations
        }"{description, sloped}
      ]
      &
      P
      \ar[
        d,
        "{f}"
      ]
      \\
      &
      A^I
      \ar[
        r,
        "{
          \mathrm{ev}_0
        }"{swap}
      ]
      \ar[
        d,
        "{
          \mathrm{ev}_1
        }"
      ]
        &
      A
      \\
      &
      A
      \mathrlap{\,.}
    \end{tikzcd}
  \end{equation}
\end{lemma}
\begin{example}[Quotient stack of coset space is delooping of subgroup]
\label{QuotientStackOfCosetSpacesIsDeloopingOfSubgroup}
For $(\mathcal{S},J)$ a (simplical) site, consider a presheaf $G$ of ordinary groups, presenting a group stack  in the  $\infty$-topos over the site,
and a subgroup-object $H$:
$$
  H
    \xhookrightarrow{i_H}
  G
  \;\;\in\;\;
  \Groups
  \big(
    \Presheaves(\SimplicialSite)
  \big)
  \xhookrightarrow{\;}
  \Groups
  \big(
    \SimplicialPresheaves(\SimplicialSite)
  \big)
  \xrightarrow{
    \Groups\big(
    \Localization
      {\LocalWeakEquivalences}
     \big)
  }
  \Groups
  \big(
    \InfinitySheaves(\InfinitySite, J)
  \big)
  \,.
$$
Then the homotopy fiber of the
$\infty$-stackification of the
induced morphism of presheaves of nerves of delooping groupoids
$$
  N\DeloopingGroupoid{H}
  \xrightarrow{
    N\DeloopingGroupoid{i_H}
  }
  N\DeloopingGroupoid{G}
  \;\;\;
  \in
  \;
  \SimplicialPresheaves(\mathcal{S})
  \xrightarrow{
    \Localization
      {\LocalWeakEquivalences}
  }
  \InfinitySheaves(\mathcal{S}, J)
$$
is computed as follows, where we are applying Lemma \ref{ComputingHomotopyPullbacksOfInfinityStacks} to the projective fibration resolution of the above morphism which is provided by the factorization lemma \ref{FactorizationLemmaForProjectiveSimplicialPresheaves}:
$$
  \begin{tikzcd}[
    column sep=8pt
  ]
    N
    \big(
      (G/H)
      \rightrightarrows
      (G/H)
    \big)
    \ar[
      from=d,
      "{
        \in
        \ProjectiveWeakEquivalences
      }"{swap}
    ]
    &&
    N
    \DeloopingGroupoid{H}
    \ar[
      d,
      "{
        N
        \Big(
        \big(
          \mathrm{const}_{\NeutralElement},
          (i_H,(-)^{-1})
        \big)
        \,\rightrightarrows\,
        \mathrm{const}_{\NeutralElement}
        \Big)
        \,\in\,
        \ProjectiveWeakEquivalences
      }"{description}
    ]
    \\[+20pt]
    N
    \ActionGroupoid
      { G }
      { H^{\mathrm{op}} }
    \ar[rr]
    \ar[d]
    \ar[
      drr,
      phantom,
      "{
        \mbox{\tiny\rm(pb)}
      }"
    ]
    &&
    N
    \big(
      G \times (G \times H^{\mathrm{op}})
      \rightrightarrows
      G
    \big)
    \ar[
      d,
      "{
        N
        \big(
          \mathrm{pr}_2
          \,\rightrightarrows\,
          \mathrm{const}_\ast
        \big)
        \,\in\,
        \ProjectiveFibrations
      }"
    ]
    \\
    N(\ast \,\rightrightarrows\, \ast)
    \ar[rr]
    &&
    N
    \DeloopingGroupoid{G}
  \end{tikzcd}
  \hspace{.7cm}
  \raisebox{-20pt}{$
    \overset{
      \Localization
        {\LocalWeakEquivalences}
    }{\longmapsto}
  $}
  \hspace{.7cm}
  \begin{tikzcd}[column sep=10pt]
    G/H
    \ar[d, "\sim"{sloped}]
    \ar[rr]
    &&
    \mathbf{B}H
    \ar[
      d,
      "{\sim}"{sloped}
    ]
    \\
    \HomotopyQuotient{G}{H}
    \ar[rr]
    \ar[d]
    \ar[
      drr,
      phantom,
      "{
        \mbox{\tiny\rm(pb)}
      }"{pos=.4}
    ]
    &&
      \HomotopyQuotient
        {(G/H)}
        { G }
   \ar[d]
    \\
    \ast
    \ar[rr]
    &&
    \mathbf{B}G
  \end{tikzcd}
$$
Comparison with Prop. \ref{GroupsActionsAndFiberBundles} on the right shows first of all (for the case $H = 1$ the trivial subgroup), that the $\infty$-stackification of the nerve of the delooping groupoid of $G$ is the $\infty$-topos theoretic delooping of the $\infty$-stackification of $G$:
\begin{equation}
  \Localization{\LocalWeakEquivalences}
  \big(
    N\DeloopingGroupoid{G}
  \big)
  \;\;
  \simeq
  \;\;
  \mathbf{B}G
  \,;
\end{equation}
and then, more generally, that
 the homotopy quotient of any $H \subset G$-coset object in the $\infty$-topos by its residual $G$-action is equivalently the delooping stack of $H$:
 \begin{equation}
   \label{HomotopyQuotientOfCosetSpace}
   \HomotopyQuotient
     {(G/H)}
     {G}
  \;\simeq\;
  \mathbf{B}H
   \;\;
   \in
   \;
   \InfinitySheaves(\InfinitySite, J)
   \,.
 \end{equation}
\end{example}
\begin{example}[Double homotopy coset stacks]
  \label{DoubleHomotopyCosetStacks}
  In direct generalization of Exp. \ref{QuotientStackOfCosetSpacesIsDeloopingOfSubgroup}
  the Factorization Lemma \ref{FactorizationLemmaForProjectiveSimplicialPresheaves}
  combined with Lemma \ref{ComputingHomotopyPullbacksOfInfinityStacks}
  shows that for a pair of presheaves of sub-groups
  $$
    H_1, H_2
    \,\in\,
    \Groups
    \big(
      \Presheaves(\SimplicialSite)
    \big)
    \xhookrightarrow{\;\;}
    \Groups
    \big(
      \SimplicialPresheaves(\SimplicialSite)
    \big)
    \xhookrightarrow{
      \Groups
      \big(
        \Localization{\LocalWeakEquivalences}
      \big)
    }
  $$
  the homotopy fiber product of the delooping of their
  inclusion is the homotopy double coset object:
  $$
    \begin{tikzcd}
      \HomotopyQuotient
        { (G/H_1) }
        { H_2 }
      \ar[rr]
      \ar[d]
      \ar[
        drr,
        phantom,
        "{\mbox{\tiny\rm(pb)}}"
      ]
      &&
      \mathbf{B}H_1
      \ar[
        d,
        "{
          \mathbf{B}
          \big(
            i_{H_1}
          \big)
        }"
      ]
      \\
      \mathbf{B}H_2
      \ar[
        rr,
        "{
          \mathbf{B}\big( i_{H_2}\big)
        }"{swap}
      ]
      &&
      \mathbf{B}G
    \end{tikzcd}
    \;\;\;\;\;\;
    \in
    \;\;
    \InfinitySheaves(\InfinitySite, J)
    \,.
  $$
  Here $H_1$ acts by multiplication from one side, $H_2$ by residual multiplication from the other side. Notice that the symmetry of the Cartesian diagram shows as once that
  $$
    \HomotopyQuotient
      { (G/H_1) }
      { H_2 }
    \;\;
    \simeq
    \;\;
    \HomotopyQuotient
      { (G/H_2) }
      { H_1 }
    \,.
  $$
\end{example}

\medskip

This way, while the  Factorization Lemma \ref{FactorizationLemmaForProjectiveSimplicialPresheaves}
serves to handle most our fibrancy needs (together with Prop. \ref{RecognitionOfFibrantObjectsOverCartSp}, Prop. \ref{LocalFibrancyOfSimplicialDeloopingsOfLocallyFibrantSimplicialGroups} below), the following Prop. \ref{DuggerCofibrancyRecognition} is our tool for ensuring cofibrancy (in Exp. \ref{GoodOpenCoversAreProjectivelyCofibrantResolutionsOfSmoothManiolds} and Exp. \ref{CechActionGroupoidOfEquivariantGoodOpenCoverIsLocalCofibrantResolution} below):

\begin{proposition}[Dugger's cofibrancy recognition {\cite[Cor . 9.4]{Dugger01Universal}}]
  \label{DuggerCofibrancyRecognition}
  Let $\mathcal{S}$ be a 1-site.
  A sufficient condition for
  $X \in \mathrm{SimplPSh}(\mathcal{S})_{ \projloc }$
  (Ntn. \ref{ModelCategoriesOfSimplicialPresheaves})
  to be
  projectively cofibrant
  is that in each simplicial degree $k$,
  the component presheaf $X_k \in \mathrm{PSh}(\mathcal{S})$ is

  {\bf (i)} a coproduct $X_k \simeq \underset{i_k}{\coprod} U_{i_k}$
    of representables $U_{i_k} \in \mathcal{S} \xhookrightarrow{\; y \;} \mathrm{PSh}(\mathcal{S})$;

  {\bf (ii)} whose degenerate cells split off as a disjoint summand:
    $X_k \simeq N_k \coprod \mathrm{im}(\sigma)$
    for some $N_k$.
\end{proposition}

\begin{example}[Basic examples of projectively cofibrant simplicial presheaves]
  \label{ExamplesOfProjectiveleCofibrantSimplicialPresheaves}
  Let $\mathcal{S}$ be a 1-site.
  Examples of projectively cofibrant simplicial presheaves
  over $\mathcal{S}$ include:

  \vspace{-3mm}
  \begin{enumerate}[{\bf (i)}]
  \setlength\itemsep{-3pt}
  \item  every representable $U \in \mathcal{Y} \xhookrightarrow{\;y\;} \mathrm{SimplPSh}(\mathcal{S})$;

  \item every constant
    simplicial presheaf
    $S \in \mathrm{SimpSets} \xhookrightarrow{\; \rm const\; } \mathrm{SimplPSh}(\mathcal{S})$;

  \noindent
  and in joint generalization of these two cases:

  \item  every product $U \times S$ of a representable with a
  simplicial set.
\end{enumerate}

\vspace{-3mm}
\noindent
  In all cases the defining lifting property is readily checked.
  Alternatively, these follow with
  Prop. \ref{DuggerCofibrancyRecognition}.
\end{example}

\medskip
\noindent
{\bf Kan extension to presheaves.} A central operation in
\cref{CohesiveHomotopyTheory} below is the induction,
by {\it Kan extension} of
quadruples of adjoint $\infty$-functors between $\infty$-stacks
from pairs of adjoint functors between simplicial presheaves.
For completeness, we make explicit the elementary but crucial
operations.

\medskip

\begin{lemma}[Kan extension of $\infty$-functor to $\infty$-presheaves]
  \label{InfinityKanExtension}
  If $\InfinitySite_1 \xrightarrow{f} \InfinitySite_2$
  is an $\infty$-functor between small $\infty$-sites (Ntn. \ref{InfinitySite}), then
  precomposition
  $\phi^\ast
     \;\coloneqq\;
    \InfinityPresheaves(\InfinitySite_2)
      \xrightarrow{ (-) \circ \phi }
    \InfinityPresheaves(\InfinitySite_1)
  $
  preserves both limits and colimits, and hence has
  left- and right adjoint $\infty$-functors
  \vspace{-4mm}
  $$
    \begin{tikzcd}
      \InfinityPresheaves(\InfinitySite_1)
      \ar[
        rr,
        shift left=14pt,
        "{ f_! }"{above}
      ]
      \ar[
        rr,
        shift right=13.5pt,
        "{ f_\ast }"{below}
      ]
      &&
      \InfinityPresheaves(\InfinitySite_2) \;.
      \ar[
        ll,
        "{ f^\ast }"{ description }
      ]
      \ar[
        ll,
        phantom,
        shift left=8.5pt,
        "{ \scalebox{.6}{$\bot$} }"{description}
      ]
      \ar[
        ll,
        phantom,
        shift right=9pt,
        "{ \scalebox{.6}{$\bot$} }"{description}
      ]
    \end{tikzcd}
  $$
\end{lemma}
\begin{lemma}[Left Kan extension on representables is original functor]
  \label{LeftKanExtensionOnRepresentablesIsOriginalFunctor}
  On representables, the left Kan extension $f_!$
  (Lem. \ref{InfinityKanExtension})
  acts as the original functor $f$:
  \vspace{-2mm}
  $$
    f_!
    \big(
      y(s_1)
    \big)
    \;\simeq\;
    y
    \big(
      f(s_1)
    \big)
    \,.
  $$
\end{lemma}
\begin{proof}
For $s_1\,\in\, \InfinitySite_1$ and $X \,\in\, \InfinityPresheaves(\InfinitySite_2)$
we have the following sequence of natural equivalences:
\vspace{-2mm}
 $$
   \begin{aligned}
         \InfinityPresheaves(\InfinitySite_2)
     \left(
       f_!
       \big(
         y(s_1)
       \big)
       ,\,
       X
     \right)
         & \;\simeq\;
     \InfinityPresheaves(\InfinitySite_1)
     \big(
       y(s_1)
       ,\,
       f^\ast(X)
     \big)
     \\
     &
     \;\simeq\;
     X
     \big(
       f(s_1)
     \big)
     \\
     &
     \;\simeq\;
     \InfinityPresheaves(\InfinitySite_2)
     \Big(
       y\big(f(s_1)\big)
       ,\,
       X
     \Big)
     \,.
   \end{aligned}
 $$

 \vspace{-2mm}
\noindent Since these equivalences are natural in $X$, the
 $\infty$-Yoneda lemma (Prop. \ref{InfinityYonedaLemma})
 for $\InfinityPresheaves(\InfinitySite_2)^{\mathrm{op}}$,
 implies the claim.
\end{proof}

\begin{lemma}[Left Kan extension preserves binary products if original functor does]
  \label{LeftKanExtensionPreservesBinaryProductsIfOriginalFunctorDoes}
  If a pair of small $\infty$-categoris $\InfinitySite_1, \InfinitySite_2$ has binary products
  and an $\infty$-functor $\InfinitySite_1 \xrightarrow{ \;f\; } \InfinitySite_2  $
  preserves these, in that for $s, s' \,\in\, \InfinitySite_1 $
  there is a natural equivalence $f(s \times s') \,\simeq\, f(s) \times f(s')$,
  then its left Kan extension also preserves binary products. That is,
  for $ X, \, X' \,\in\, \InfinityPresheaves(\InfinitySite_1)$,
  there is a natural equivalence
  \vspace{-2mm}
  $$
    f_!
    (
      X \times X'
    )
    \;\;
    \simeq
    \;\;
    f_!(X)
    \times
    f_1(X')
    \,.
  $$
\end{lemma}
\begin{proof}
We have the sequence of natural equivalences
\vspace{-2mm}
$$
  \begin{aligned}
      f_!(X \times X')
       & \;\simeq\;
    f_!
    \Big(\!
      \big(
        \underset{
          \underset{ s \to X}{\longrightarrow}
        }{\lim}
        \,
        y(s)
      \big)
      \times
      \big(
        \underset{
          \underset{ s' \to X}{\longrightarrow}
        }{\lim}
        \,
        y(s')
      \big)
    \!\Big)
    \\
    & \;\simeq\;
    f_!
    \Big(
        \underset{
          \underset{ s \to X}{\longrightarrow}
        }{\lim}
        \;
        \underset{
          \underset{ s' \to X}{\longrightarrow}
        }{\lim}
        \,
      \big(
        y(s)
        \times
        y(s')
      \big)
    \!\!\Big)
    \\
    & \;\simeq\;
        \underset{
          \underset{ s \to X}{\longrightarrow}
        }{\lim}
        \;
        \underset{
          \underset{ s' \to X}{\longrightarrow}
        }{\lim}
    \,
    f_!
      \left(
        y(s)
        \times
        y(s')
      \right)
    \\
    & \;\simeq\;
        \underset{
          \underset{ s \to X}{\longrightarrow}
        }{\lim}
        \;
        \underset{
          \underset{ s' \to X}{\longrightarrow}
        }{\lim}
    \,
    f_!
      \left(
        y(s)
      \right)
   \times
    f_!
      \left(
        y(s')
      \right)
    \\
    & \;\simeq\;
    \Big(\,
        \underset{
          \underset{ s \to X}{\longrightarrow}
        }{\lim}
    \,
    f_!
      \left(
        y(s)
      \right)
   \!\!\Big)
   \times
   \Big(\,
        \underset{
          \underset{ s' \to X}{\longrightarrow}
        }{\lim}
    \,
    f_!
      \left(
        y(s')
      \right)
  \!\!\Big).
  \end{aligned}
$$

\vspace{-7mm}
\end{proof}

\medskip
\begin{lemma}[Left Kan extension of fully faithful functor is fully faithful]
\label{LeftKanExtensionOfFullyFaithfulFunctorIsFullyFaithful}
The left Kan extension $f_!$ of a fully faithful functor $f$ between small $\infty$-categories
is itself fully faithful \eqref{FullyFaithfulInfinityFunctor}:
\vspace{-2mm}
$$
  \begin{tikzcd}
    \InfinitySite_1
    &
    \InfinitySite_2
    \ar[
      l,
      hook',
      "{i}"{above}
    ]
  \end{tikzcd}
  \;\;\;\;\;\;\;\;\;
  \Rightarrow
  \;\;\;\;\;\;\;\;\;
  \begin{tikzcd}
    \InfinityPresheaves(\InfinitySite_1)
    &
    \InfinityPresheaves(\InfinitySite_2)\;.
    \ar[
      l,
      hook',
      "{i_!}"{above}
    ]
  \end{tikzcd}
$$
\end{lemma}
\begin{proof}
For $X ,\, X' \,\in\, \InfinityPresheaves(\InfinitySite_2)$,
we have the following sequence of natural equivalences:
\vspace{-2mm}
$$
  \begin{aligned}
    \InfinityPresheaves(\InfinitySite_1)
    \left(
      i_!(X)
      ,\,
      i_!(X')
    \right)
    &
    \;\simeq\;
    \InfinityPresheaves(\InfinitySite_1)
    \Big(
      i_!
      \big(\,
      \colimit{s \to X}
        y(s)
      \big)
      ,\,
      i_!
      \big(\,
        \colimit{s' \to X'}
        y(s')
      \big)
    \!\Big)
    \\
    &
    \;\simeq\;
    \InfinityPresheaves(\InfinitySite_1)
    \Big(\,
      \colimit{s \to X}
        i_!
        \left(
          y(s)
        \right)
      ,\,
      \colimit{s' \to X'}
        i_!
        \left(
          y(s')
        \right)
   \!\! \Big)
    \\
    &
    \;\simeq\;
    \InfinityPresheaves(\InfinitySite_1)
    \Big(\,
      \colimit{s \to X}
        y
        \left(
          i
          (s)
        \right)
      ,\,
      \colimit{s' \to X'}
        y
        \left(
          i(s')
        \right)
   \!\! \Big)
    \\
    &
    \;\simeq\;
    \limit{s \to X}
    \colimit{s' \to X'}
    \InfinityPresheaves(\InfinitySite_1)
    \left(
        y
        \left(
          i
          (s)
        \right)
      ,\,
        y
        \left(
          i(s')
        \right)
    \right)
    \\
    &
    \;\simeq\;
    \limit{s \to X}
    \colimit{s' \to X'}
    \InfinitySite_1
    \left(
      i(s)
      ,\,
      i(s')
    \right)
    \\
    &
    \;\simeq\;
    \limit{s \to X}
    \colimit{s' \to X'}
    \InfinitySite_2
    (
      s
      ,\,
      s'
    )
    \\
    &
    \;\simeq\;
    \limit{s \to X}
    \colimit{s' \to X'}
    \InfinityPresheaves(\InfinitySite_2)
    \left(
      y(s)
      ,\,
      y(s')
    \right)
    \\
    &
    \;\simeq\;
    \InfinityPresheaves(\InfinitySite_2)
    \Big(\,
      \colimit{s \to X}
      y(s)
      ,\,
      \colimit{s' \to X'}
      y(s')
    \Big)
    \\
    &
    \;\simeq\;
    \InfinityPresheaves(\InfinitySite_2)
    \left(
      X
      ,\,
      X'
    \right).
  \end{aligned}
$$

\vspace{-8mm}
\end{proof}

\begin{lemma}[Kan extensions of an adjoint pair yield an adjoint quadruple]
\label{KanExtensionOfAdjointPairOfInfinityFunctors}
Given a pair of adjoint $\infty$-functors between small $\infty$-categories
\vspace{-2mm}
$$
  \begin{tikzcd}
    \InfinitySite_1
    \ar[
      rr,
      shift left=5pt,
      "{ \ell }"{above}
    ]
    &&
    \InfinitySite_2
    \ar[
      ll,
      shift left=5pt,
      "{ r }"{below}
    ]
    \ar[
      ll,
      phantom,
      "{ \scalebox{.6}{$\bot$} }"
    ]
  \end{tikzcd}
$$

\vspace{-2mm}
\noindent the corresponding pullback functors on presheaves
are adjoint to each other,
$\ell^\ast \,\dashv\, r^\ast$,
and their adjoint triples of
Kan extensions (from Lem. \ref{InfinityKanExtension})
overlap:
\vspace{-4mm}
$$
    \begin{tikzcd}[column sep=large]
      \InfinityPresheaves(\InfinitySite_1)
      \;\;
      \ar[
        rr,
        shift left=-7pt,
        "{ \ell_\ast \,\simeq\, r^\ast    }"{description}
      ]
      \ar[
        rr,
        shift left=32pt-7pt,
        "{ \ell_! }"{description, pos=.4}
      ]
      &&
      \;\;
      \InfinityPresheaves(\InfinitySite_2)\;.
      \ar[
        ll,
        shift right=16pt-7pt,
        "{ \ell^\ast \,\simeq\, r_! }"{description}
      ]
      \ar[
        ll,
        shift right=-16pt-7pt,
        "{ r_\ast }"{description, pos=.4}
      ]
      \ar[
        ll,
        phantom,
        shift right=8pt-7pt,
        "{\scalebox{.6}{$\bot$}}"
      ]
      \ar[
        ll,
        phantom,
        shift right=24pt-7pt,
        "{\scalebox{.6}{$\bot$}}"
      ]
      \ar[
        ll,
        phantom,
        shift right=-8pt-7pt,
        "{\scalebox{.6}{$\bot$}}"
      ]
    \end{tikzcd}
$$
\end{lemma}
\begin{proof}
  The characterizing hom-equivalence of the adjunction $\ell^\ast \dashv r^\ast$
  is the composition of the following sequence of natural equivalences,
  for $X_i \,\in\, \InfinityPresheaves(\InfinitySite_i)$
  is checked along the lines of the proof of Lem. \ref{LeftKanExtensionOnRepresentablesIsOriginalFunctor}.
  From this, the overlap of the adjoint triples of Kan extensions
  follows by essential uniqueness of adjoints.
\end{proof}

\medskip

\noindent
{\bf Systems of local sections.}

\begin{notation}[$\infty$-Yoneda embedding of slice sites]
  \label{InfinityYonedaEmbeddingOfSliceSites}
  For $\InfinitySite$ a small $\infty$-category and $X \,\in\, \InfinitySite$ an object,
  we have both the $\infty$-Yoneda embedding \eqref{InfinityYonedaEmbedding}
  of the slice $\InfinitySite_{/X}$
  \vspace{-2mm}
  \begin{equation}
    \label{InfinityYonedaEmbeddingOfSliceSite}
    \begin{tikzcd}[row sep=-4pt]
      \InfinitySite_{/X}
      \ar[
        rr,
        hook,
        "{\YonedaEmbedding_{(\InfinitySite_{/X})}}"
      ]
      &&
      \InfinityPresheaves(\InfinitySite_{/X})
      \\
  \scalebox{0.7}{$    (U \overset{\phi}{\to} X) $}
      &\longmapsto&
 \scalebox{0.7}{$        \InfinitySite_{/X}
      \big(
        -,
        \,
        (U\overset{\phi}{\to} X)
      \big)
      $}
    \end{tikzcd}
  \end{equation}

  \vspace{-2mm}
\noindent  and also the slicing of the Yoneda embedding on $\InfinitySite$:
  \vspace{-2mm}
  \begin{equation}
    \label{SlicedInfinityYonedaEmbedding}
    \begin{tikzcd}[row sep=-4pt, column sep=small]
      \InfinitySite_{/X}
      \ar[
        rr,
        hook,
        "{ (\YonedaEmbedding_{\InfinitySite})_{/X} }"
      ]
      &&
      \InfinityPresheaves(\InfinitySite)_{/\YonedaEmbedding_{\InfinitySite}(X)}
      \\
\scalebox{0.7}{$         (U \overset{\phi}{\to}X) $}
      &\longmapsto&
 \scalebox{0.7}{$        \big(
      \InfinitySite(-,\, U )
      \xrightarrow{ \InfinitySite(-,\phi) }
      \InfinitySite(-,\, X)
      \Big).
      $}
    \end{tikzcd}
  \end{equation}
\end{notation}

\begin{proposition}[Systems of local sections of bundles internal to $\infty$-presheaves]
  \label{SystemsOfLocalSectionsOfBundlesOfInfinityPresheaves}
  For $\BaseTopos$ an $\infty$-topos, $\InfinitySite$ a small $\infty$-category
  and $X \,\in\, \InfinitySite$ an object, the $\infty$-functor which
  assembles local sections of bundles of $\infty$-presheaves
  over (the image under the $\infty$-Yoneda embedding of) $X$ is
  an equivalence of $\infty$-categories:
  \vspace{-2mm}
  \begin{equation}
    \label{LocalSectionsAsEquivalenceOfInfinityCategories}
    \begin{tikzcd}[column sep=80pt]
      \Presheaves
      (
        \InfinitySite
        ,\,
        \BaseTopos
      )_{/ \YonedaEmbedding_{\InfinitySite}(X)  }
      \ar[
        rr,
        "{
          \Gamma_{(-)}(-)
          \,\coloneqq\,
          \Presheaves(
            \InfinitySite,
            \BaseTopos)_{/\YonedaEmbedding_{\InfinitySite}(X)}
          \left(
            (\YonedaEmbedding_{\InfinitySite})_{/X}(-)
            ,\,
            -
          \right)
        }"{above},
        "{ \sim }"{below}
      ]
      &&
      \Presheaves
      (
        \InfinitySite_{/X}
        ,\,
        \BaseTopos
      )
      \,,
    \end{tikzcd}
  \end{equation}

\vspace{-2mm}
\noindent
where $(\YonedaEmbedding_{\InfinitySite})_{/X}$ is
from \eqref{InfinityYonedaEmbeddingOfSliceSite}.
\end{proposition}
\begin{proof}
  For the analogous equivalence for 1-categories
  (which is classical, e.g. \cite[Lem. 1.4.12]{KashiwaraSchapira06})
  one readily checks\footnote{This is spelled out also at \href{https://ncatlab.org/nlab/show/slice+of+presheaves+is+presheaves+on+slice}{\tt ncatlab.org/nlab/show/slice+of+presheaves+is+presheaves+on+slice}}
  that:

  (1) the analogous functor is a right adjoint, and

  (2) generalizes to a simplicial adjunction of simplicial presheaves over simplicial sites,

  (3) where it is a Quillen equivalence for the
    projective model structure and its slice model structure.

  \noindent
  This implies the claim,
  by Prop. \ref{HomotopyCategoryOfPresInfinityCategoriesIsThatOfCombinatorialModelCategories}.
\end{proof}
   An alternative general argument for some equivalence of $\infty$-categories as in Prop. \ref{SystemsOfLocalSectionsOfBundlesOfInfinityPresheaves}
   is given in \cite[Cor. 5.1.6.12]{Lurie09HTT}, in terms of quasi-categories.
On the other hand, in our formulation the proof of Prop. \ref{SystemsOfLocalSectionsOfBundlesOfInfinityPresheaves} applies verbatim also in the case that the base presheaf is not necessarily representable (not necessarily in the image of the Yoneda embedding), if we use the following more general notion of slice $\infty$-site:

\begin{definition}[Slice $\infty$-site]
  \label{SliceInfinitySite}
  Given an $\infty$-site (Ntn. \ref{InfinitySite}) $\InfinitySite$ and an $\infty$-presheaf $X \,in\, \InfinityPresheaves(\InfinitySite)$ (Ntn. \ref{PresentationOfInfinityPresheavesBySimplicialPresheaves}) the {\it slice $\infty$-site} $\Slice{\InfinitySite}{X}$ is the full sub-$\infty$-category of the slice $\Slice{\InfinityPresheaves(\InfinitySite)}{X}$
  \begin{equation}
    \label{SliceInfinitySiteEmbedding}
    \begin{tikzcd}
      \Slice{\InfinitySite}{X}
      \ar[
        rr,
        hook
      ]
      &&
      \Slice{\InfinityPresheaves{\InfinitySite}}{X}
      \,.
    \end{tikzcd}
  \end{equation}
  on the morphisms
  with representable domain, i.e. those of the form $y_{\InfinitySite}(U) \xrightarrow{\phi} X$ for some $U \,\in\, \InfinitySite$.
\end{definition}
By Prop. \ref{HomSpaceInSliceToposAsFiberProduct} this means that the hom-$\infty$-groupoids in a slice $\infty$-site
\eqref{SliceInfinitySiteEmbedding}
are given by the following homotopy fiber product of $\infty$-groupoids:
\begin{equation}
  \label{HomInfinityGroupoidOfSliceInfinitySite}
  \begin{tikzcd}
    \Slice{\InfinitySite}{X}
    \Big(
      (y_{\InfinitySite}(U) \xrightarrow{\phi} X)
      ,\,
      (y_{\InfinitySite}(U') \xrightarrow{\phi'} X)
    \Big)
    \ar[rr]
    \ar[d]
    \ar[
      drr,
      phantom,
      "{\mbox{\tiny\rm(pb)}}"{pos=.4}
    ]
    &&
    \InfinityPresheaves(\InfinitySite)
    \big(
      y_{\InfinitySite}(U)
      ,\,
      y_{\InfinitySite}(U')
    \big)
    \ar[
      d,
      "{
        \phi' \,\circ\, (-)
      }"
    ]
    \\
    \ast
    \ar[
      rr,
      "{
        \vdash \, \phi
      }"{swap}
    ]
    &&
    \InfinityPresheaves(\InfinitySite)
    \big(
      y_{\InfinitySite}(U)
      ,\,
      X
    \big)
  \end{tikzcd}
\end{equation}

With this, the proof of Prop. \ref{SystemsOfLocalSectionsOfBundlesOfInfinityPresheaves} generalizes to:
\begin{proposition}[Fundamental theorem of $\infty$-presheaf $\infty$-topos theory]
  \label{FundamentalTheoremOfPresheafToposTheory}
  For $\InfinitySite$ an $\infty$-site and any
  $X \,\in\, \InfinityPresheaves(\InfinitySite)$
  there is an equivalence of $\infty$-categories between the slice of all $\infty$-presheaves over $X$ and the $\infty$-presheaves over the slice $\infty$-site (Def. \ref{SliceInfinitySite}):
  \begin{equation}
    \label{EquivalenceForFundamentalTheoremOfSliceInfinityToposes}
    \Slice{\InfinityPresheaves(\InfinitySite)}{X}
    \;\;\simeq\;\;
    \InfinityPresheaves
    \big(
      \Slice{\InfinitySite}{X}
    \big)
    \,.
  \end{equation}
\end{proposition}
A proof for this statement in the special case that $X$ is 1-truncated also essentially appears as \cite[Thm. 6.1(a)]{Hollander08}.

  The equivalence \eqref{EquivalenceForFundamentalTheoremOfSliceInfinityToposes} in
  Prop. \ref{FundamentalTheoremOfPresheafToposTheory}
  shows in particular that every slice of a presheaf $\infty$-topos
  (on the left) is again an $\infty$-topos (manifest on the right).
  This is the archetypical special case of the
  {\it fundamental theorem of $\infty$-topos theory}
  (e.g. \cite[Rem. 2.2.10]{Vergura19}, following the terminology for 1-toposes \cite{McLarty}),
  which says that any slice of any $\infty$-topos is again an $\infty$-topos:

\begin{proposition}[Fundamental theorem of $\infty$-topos theory
{\cite[Prop. 6.3.5.1]{Lurie09HTT}}]
  \label{SliceInfinityTopos}
  $\,$

  \noindent
  For $\Topos$ an $\infty$-topos and $B \,\in\, \Topos$ an object,
   the slice $\infty$-category $\ModalTopos{/B}$ \eqref{SliceHomInIntroduction}
   is also an $\infty$-topos.
\end{proposition}

\begin{example}[Basic structures in slice infinity topos]
  \label{BasicStructuresInSliceInfinityTopos}
  In $\SliceTopos{B}$ the terminal object is $(B, \mathrm{id}_B)$.
\end{example}

\medskip

\noindent
{\bf Base change of $\infty$-toposes.} In higher generalization of
Lem. \ref{InducedAndCoinducedActions}, much of the theory of
(equivariant) $\infty$-toposes is driven by their base change adjoint triples:

\begin{proposition}[Base change]
  \label{BaseChange}
  For $\Topos \,\in\, \InfinityToposes$
  and $B_1 \xrightarrow{f} B_2$ a morphism in $\Topos$,
  there is an adjoint triple of $\infty$-functors between the
  slice $\infty$-toposes (Prop. \ref{SliceInfinityTopos}):
  \vspace{-1mm}
  \begin{equation}
    \label{BaseChangeAdjointTriple}
    \begin{tikzcd}[column sep=large]
      \Topos_{/B_1}
      \ar[
        rr,
        shift left=14pt,
        "f_! \,=\, \sum_f"
      ]
      \ar[
        rr,
        phantom,
        shift left=9pt,
        "\scalebox{.5}{$\bot$}"
      ]
      \ar[
        rr,
        phantom,
        shift right=8pt,
        "\scalebox{.5}{$\bot$}"
      ]
      \ar[
        rr,
        shift right=14pt,
        "f_\ast \,=\, \prod_f"{below}
      ]
      &&
      \Topos_{/B_2}
      \ar[
        ll,
        "f^\ast"{description}
      ]
      \ar[
        rr, shift left=14pt,
        phantom,
        "\mbox{
          \tiny
          \color{greenii}
          \bf
          left base change
        }"
      ]
      \ar[
        rr,
        phantom,
        "\mbox{
          \tiny
          \color{greenii}
          \bf
          \phantom{left} base change
        }"
      ]
      \ar[
        rr, shift right=14pt,
        phantom,
        "\mbox{
          \tiny
          \color{greenii}
          \bf
          right base change
        }"
      ]
      &&
      {}
    \end{tikzcd}
  \end{equation}

  \vspace{-1mm}
  \noindent
  where $f^\ast$ is given by pullback along $f$ and
  $f_!$ is given by postcomposition with $f$.
\end{proposition}
\begin{example}[Base change to absolute context]
  \label{BaseChangeToAbsoluteContext}
 For $B \in \Topos$ any object
 the base change (Prop. \ref{BaseChange})
 along its terminal morphism $B \xrightarrow{\exists !} \ast$
 is the operation of taking the Cartesian product with $B$,
 and the corresponding
 adjoint triple \eqref{BaseChangeAdjointTriple} looks as follows:
 \vspace{-1mm}
  \begin{equation}
    \label{BaseChangeAdjointTripleToAbsoluteContext}
    \begin{tikzcd}[column sep=large]
      \Topos_{/B}
      \;\;
      \ar[
        rr,
        shift left=13pt,
        "\mathrm{dom} \,=\, \sum_B"
      ]
      \ar[
        rr,
        phantom,
        shift left=9pt,
        "\scalebox{.5}{$\bot$}"
      ]
      \ar[
        rr,
        phantom,
        shift right=8pt,
        "\scalebox{.5}{$\bot$}"
      ]
      \ar[
        rr,
        shift right=13pt,
        "\prod_B"{below}
      ]
      &&
\;\;      \Topos_{/\ast}
      \mathrlap{\; \simeq \Topos\;.}
      \ar[
        ll,
        "(-) \times B"{description}
      ]
    \end{tikzcd}
  \end{equation}
\end{example}

\begin{lemma}[Frobenius reciprocity]
  For $B_1 \xrightarrow{f} B_2$ a morphism in an $\infty$-topos $\Topos$,
  we have for $E_1 \in \Topos_{/B_1}$ and $E_2 \in \Topos_{/B_2}$
  a natural equivalence
  \vspace{-3mm}
  \begin{equation}
    \label{FrobeniusReciprocityEquivalence}
    f_!
    \left(
      E_1 \times f^\ast(E_2)
    \right)
    \;\simeq\;
    f_!(E_1) \times_{B_2} E_2
    \;\;\;
    \in
    \;
    \Topos_{/B_2}
  \end{equation}

  \vspace{-2mm}
\noindent
  of left base changes \eqref{BaseChangeAdjointTriple}.
\end{lemma}

\begin{proof}
  Since, by Prop. \ref{BaseChange}, the left base change is
  given by postcomposition with $f$, this follows by the pasting law
  \eqref{HomotopyPastingLaw}:
  \vspace{-2mm}
  $$
    \begin{tikzcd}[row sep=small, column sep=large]
      f_!(E_1) \times_{B_2} E_2
      \ar[r]
      \ar[d]
      \ar[
        dr,
        phantom,
        "\mbox{\tiny\rm(pb)}"
      ]
      \ar[
      rr,
      rounded corners,
      to path={
           -- ([yshift=+10pt]\tikztostart.north)
           --node[above]{
               \scalebox{.7}{$
                 p_{ f_!( E_1 \times_{B_1} f^\ast(E_2) ) }
               $}
             } ([yshift=+12pt]\tikztotarget.north)
           -- (\tikztotarget.north)}
      ]
      &
      f^\ast(E_2)
      \ar[r]
      \ar[d]
      \ar[
        dr,
        phantom,
        "\mbox{\tiny\rm(pb)}"
      ]
      &
      E_2
      \ar[
        d,
        "p_{E_2}"
      ]
      \\
      E_1
      \ar[
        r,
        "p_{E_1}"{below}
      ]
    \ar[
      rr,
      rounded corners,
      to path={
           -- ([yshift=-12pt]\tikztostart.south)
           -- node[below]{
              \scalebox{.7}{$
                p_{f_! E_1}
              $}
             } ([yshift=-12pt]\tikztotarget.south)
           -- (\tikztotarget.south)}
    ]
      &
      B_1
      \ar[
        r,
        "f"{below}
      ]
      &
      B_2
    \end{tikzcd}
  $$

  \vspace{-6mm}
\end{proof}

\subsection{General constructions in $\infty$-toposes}
\label{GeneralConstructionsInInfinityToposes}

\noindent
{\bf Truncation and connectivity.}

\begin{proposition}[$n$-Truncation modality {\cite[5.5.6.18, 6.5.1.2]{Lurie09HTT}}]
\label{nTruncation}
For $n \in \{-2,-1, 0, 1, \cdots\}$ and $\Topos$ an $\infty$-topos,
its full sub-$\infty$-category of $n$-truncated objects is reflective,
and the reflector $\tau_n$ preserves finite products:
\vspace{-2mm}
\begin{equation}
  \label{nTruncationReflection}
  \begin{tikzcd}
    \Topos_{n}
    \ar[
      rr,
      hook,
      shift right=5pt,
      "i_n"{below}
    ]
    \ar[
      rr,
      phantom,
      "\scalebox{.7}{$\bot$}"
    ]
    &&
    \Topos
    \mathrlap{\,.}
    \ar[
      ll,
      shift right=5pt,
      "\mathclap{\times}"{description, pos=0},
      "\tau_n"{above}
    ]
  \end{tikzcd}
\end{equation}
\end{proposition}

\begin{notation}[Connected objects]
  \label{ConnectedObject}
  An object $X \,\in\,\Topos$ is {\it connected} if
  $\Truncation{0} \, X \;\simeq\; \ast$. We write
  \vspace{-1mm}
  $$
    \ModalTopos{\geq }^{\ast/}
    \xhookrightarrow{\quad}
    \ModalTopos{\geq }
    \xhookrightarrow{\quad}
    \Topos
  $$

    \vspace{-1mm}
\noindent
  for the full sub-$\infty$-categories of (pointed and)
  connected objects.
\end{notation}

\begin{lemma}[Truncation of $\infty$-stacks modeled by simplicial coskeleta of simplicial presheaves]
  \label{TruncationModeledBySimplicialCoskeleta}
  For $n \in \NaturalNumbers$, the $n$-(co-)skeleton adjunction on simplicial sets
  (see \cite[Prop. A.0.41]{FSS20CharacterMap} for review and further pointers) extends objectwise to simplicial presheaves

  \vspace{-.2cm}
  \begin{equation}
    \label{CoSkeletaAdjunction}
    \begin{tikzcd}
      \SimplicialSets
      \ar[
        from=rr,
        shift right=7pt,
        "{\mathrm{sk}_k}"{swap}
      ]
      \ar[
        rr,
        shift right=7pt,
        "{\mathrm{cosk}_k}"{swap}
      ]
      \ar[
        rr,
        phantom,
        "{
          \scalebox{.7}{$\bot$}
        }"
      ]
      &&
      \SimplicialSets
      \mathrlap{\,,}
    \end{tikzcd}
    \hspace{1cm}
    \begin{tikzcd}
      \SimplicialPresheaves(\CartesianSpaces)
      \ar[
        from=rr,
        shift right=7pt,
        "{\mathrm{sk}_k}"{swap}
      ]
      \ar[
        rr,
        shift right=7pt,
        "{\mathrm{cosk}_k}"{swap}
      ]
      \ar[
        rr,
        phantom,
        "{
          \scalebox{.7}{$\bot$}
        }"
      ]
      &&
      \SimplicialPresheaves(\CartesianSpaces)
      \;\in\;
      \Topos
      \mathrlap{\,,}
    \end{tikzcd}
  \end{equation}
  \vspace{-.2cm}

  \noindent
  and under the presentation of $\infty$-stacks by simplicial presheaves (Prop. \ref{PresentationOfInfinityToposesBySimplicialPresheaves}), the $n+1$ coskeleton operation models $n$-truncation (Prop. \ref{nTruncation}):

  \vspace{-.2cm}
  \begin{equation}
    \label{TruncationOfInfinitySheafEquivalentToCoskeletonOfRepresentingSimplicialPresheaf}
    \mathcal{X}
    \,\in\,
    \SimplicialPresheaves(\SimplicialSite,J)_{\projloc}
    \xrightarrow{
      \Localization{\LocalWeakEquivalences}
    }
    \Topos
    \;\;\;\;\;\;\;\;
    \vdash
    \;\;\;\;\;\;\;\;
    \Localization{\LocalWeakEquivalences}
    \big(
      \mathrm{cosk}_{n+1}
      \mathcal{X}
    \big)
    \;\simeq\;
    \tau_n
    \big(
      \Localization{\LocalWeakEquivalences}
      \mathcal{X}
    \big)
    \,.
  \end{equation}
\end{lemma}
\begin{proof}
  Since $\infty$-sheafification preserves all colimits and finite homotopy limits, the recursive definition
  of $n$-truncation \cite[Lem. 5.5.6.1]{Lurie09HTT} implies that we may equivalently compute it on $\infty$-presheaves. Here it is objectwise given by $n$-truncation in $\SimplicialLocalization{\WeakHomotopyEquivalences}(\SimplicialSets) \,\simeq\, \InfinityGroupoids$, which in turn is given by the $(n+1)$-coskeleton construction according to classical results (e.g. \cite[Prop. A.0.41]{FSS20CharacterMap}).
\end{proof}

\medskip
\noindent
{\bf Mapping stacks.}

\begin{proposition}[Mapping stacks]
  \label{MappingStacks}
  For $\Topos$ an $\infty$-topos and $X \in \Topos$ any object,
  the Cartesian product functor $X \times (-)$ has a right adjoint
  \vspace{-2mm}
  \begin{equation}
    \label{InternalHomAdjunction}
    \begin{tikzcd}[column sep=large]
      \Topos
      \ar[
        rr,
        shift right=5pt,
        "{
          \Maps{}{X}{-}
        }"{below}
      ]
      \ar[
        rr,
        phantom,
        "{\scalebox{.7}{$\bot$}}"
      ]
      &&
      \Topos
      \ar[
        ll,
        shift right=5pt,
        "{X \times (-)}"{above}
      ] \;.
    \end{tikzcd}
  \end{equation}

  \vspace{-2mm}
  \noindent
  If $\InfinitySite$ is an $\infty$-site for
  $
   \!\!\!
   \begin{tikzcd}
    \Topos
    \ar[
      r,
      hook,
      shift right=5.5pt
    ]
    \ar[
      r,
      phantom,
      "\scalebox{.7}{$\bot$}"{description}
    ]
    &
    \mathrm{PSh}_\infty(\InfinitySite)
    \ar[
      l,
      shift right=5.5pt,
      "\mathrm{lex}"{description}
    ]
   \end{tikzcd}
   \!\!\!,
  $
  then for $U \,\in\, \InfinitySite$ the value of the {\it mapping stack}
  $\Maps{}{X}{Y} \,\in\, \Topos$ is naturally equivalent to
  \vspace{-2mm}
  \begin{equation}
    \label{ValuesOfMappingStackAsHomSpaces}
    \Maps{}{X}{Y}(U)
    \;\simeq\;
    \Topos
    \left(
      U
      ,\,
      \Maps{}{X}{Y}
    \right)
    \;\simeq\;
    \Topos
    (
      U \times X, \, Y
    )
    \;\;\;\;\;
    \in
    \;
    \InfinityGroupoids\;.
  \end{equation}
\end{proposition}

\begin{notation}[Evaluation map on mapping stacks]
  The counit of the mapping stack adjunction \eqref{InternalHomAdjunction}
  is the {\it evaluation map}:
  \vspace{-2mm}
  \begin{equation}
    \label{EvaluationMap}
    X \times
    \Maps{}
      { X }
      {A}
    \xrightarrow{\; \mathrm{ev}_{X,A} \;}
    A
    \,.
  \end{equation}
\end{notation}

\vspace{1mm}
\begin{example}[Mapping space consists of global points of mapping stack]
\label{MappingSpaceConsistsOfPointsOfMappingStack}
\vspace{-2mm}
\begin{equation}
  \label{HomIsGlobalPointsOfMappingStacks}
     \mathbf{H}\left(
    { \ast },
    {
      \Maps{}
        { X }
        { A }
    }
    \right)
  \;\;\simeq\;\;
  \PointsMaps{}
    { X }
    { A }
  \,.
\end{equation}
\end{example}

\begin{lemma}[$\infty$-Stacks consist of their internal points]
  \label{StacksConsistOfTheirInternalPoints}
  For $X \,\in\, \mathbf{H}$ there is a natural equivalence
  \vspace{-2mm}
  $$
    X
    \;\simeq\;
    \Maps{}
      { \ast }
      { X }
    \,.
  $$
\end{lemma}
\begin{proof}
  For $U \,\in\, \InfinitySite$, there is the natural equivalence \eqref{ValuesOfMappingStackAsHomSpaces}
 \vspace{-2mm}
  $$
    \PointsMaps{}
      { U }
      { X }
    \;\;
    \simeq
    \;\;
    \PointsMaps{}
      { U \times \ast }
      { X }
    \;\;
    \simeq
    \;\;
    \mathbf{H}\left(
      { U },
      {
        \Maps{}
          { \ast }
          { X }
      }
      \right)
    .
  $$

  \vspace{-2mm}
\noindent  With this, the claim follows by the $\infty$-Yoneda lemma (Prop. \ref{InfinityYonedaLemma}).
\end{proof}

\begin{lemma}[Looping of mapping stack is mapping stack into looping]
  \label{LoopingOfMappingStackIsMappingStackIntoLooping}
  We have a natural equivalence
  \vspace{-2mm}
  $$
    \Omega
    \Maps{}
      { X }
      { A }
    \;\simeq\;
    \Maps{}
      { X }
      { \Omega A }\;.
  $$
\end{lemma}
\begin{proof}
  Since the looping operation is a limit,
  this follows from \eqref{InternalHomAdjunction}
  by \eqref{InfinityAdjointPreservesInfinityLimits}.
\end{proof}

\begin{lemma}[Internal hom-adjointness]
  \label{InternalHomAdjointness}
  The kind of hom-equivalence that characterizes the right adjoint
  \eqref{InternalHomAdjunction} also holds internally, in that
  for all $X, Y, A \,\in\, \Topos$
  there is a natural equivalence
  \vspace{-2mm}
  $$
    \Maps{}
      { X \times Y }
      { A }
    \;\;
      \simeq
    \;\;
    \Maps{\big}
      { X }
      {
        \Maps{}
          { Y }
          { A }
      }
    \,.
  $$
\end{lemma}
\begin{proof}
  For $U \,\in\, \InfinitySite$, there is the natural equivalence \eqref{ValuesOfMappingStackAsHomSpaces}
    \vspace{-2mm}
  $$
    \PointsMaps{\big}
    { U }
    {
      \Maps{\big}
        { X }
        {
          \Maps{}
            { Y }
            { A }
        }
    }
    \;\;\simeq\;\;
    \PointsMaps{\big}
      { U \times X \times Y }
      { A }
    \;\;\simeq\;\;
    \PointsMaps{\big}
      { U }
      {
        \Maps{}
          { X \times Y }
          { A }
      }
    \,.
  $$

  \vspace{-2mm}
  \noindent
  With this, the claim follows by the $\infty$-Yoneda lemma (Prop. \ref{InfinityYonedaLemma}).
\end{proof}

\begin{proposition}[Mapping stack construction preserves $\infty$-limits]
  \label{MappingStackConstructionPreservesLimits}
  The mapping stack construction (Prop. \ref{MappingStacks})
  preserves not only $\infty$-limits in its second argument,
  but also
  turns $\infty$-colimits in its first argument into $\infty$-limits:
  For $A \,\in\, \Topos$
  and
  $X_{(-)} \,\colon\, I \xrightarrow{\;} \mathbf{H}$ a small diagram,
  we have a natural equivalence
  \vspace{-2mm}
  $$
    \Maps{\big}
      {
        \underset{ \underset{ i \in I }{\longrightarrow} }{\mathrm{lim}}
        \,
        X_i
      }
      { A }
    \;\;
    \simeq
    \;\;
    \underset{ \underset{ i \in I }{\longleftarrow} }{\mathrm{lim}}
    \,
    \Maps{}
      { X_i }
      { A }
    \;\;\;
    \in
    \;
    \Topos
    \,.
  $$
\end{proposition}
\begin{proof}
  For $\InfinitySite$ any $\infty$-site for $\Topos$
  consider the following sequence of natural equivalences
  for $U \,\in\, \InfinitySite$ any object:
    \vspace{-2mm}
  $$
    \begin{array}{lll}
      \PointsMaps{\Big}
        { \! U }
        {
          \Maps{\big}
            {
              \underset{ \underset{ i \in I }{\longrightarrow} }{\mathrm{lim}}
              \,
              X_i
            }
            { A }
        \!\!}
      &
      \;\simeq\;
      \PointsMaps{\Big}
        {
          \big(
            \underset{ \underset{ i \in I }{\longrightarrow} }{\mathrm{lim}}
            \,
            X_i
          \big)
          \times U
        }
        { A }
      &
      \mbox{\small by \eqref{ValuesOfMappingStackAsHomSpaces} }
      \\
      &
      \;\simeq\;
      \PointsMaps{\Big}
        {
          \underset{ \underset{ i \in I }{\longrightarrow} }{\mathrm{lim}}
          \big(
            X_i
            \times U
          \big)
        }
        { A }
      &
      \mbox{\small by \eqref{UniversalityOfColimits} }
      \\
      &
      \;\simeq\;
      \underset{ \underset{ i \in I }{\longleftarrow} }{\mathrm{lim}}
      \,
      \PointsMaps{}
        { X_i \times U }
        { A }
      &
      \mbox{\small by \eqref{HomFunctorRespectsLimits}  }
      \\
      &
      \;\simeq\;
      \underset{ \underset{ i \in I }{\longleftarrow} }{\mathrm{lim}}
      \,
      \PointsMaps{\big}
        { U }
        {
          \Maps{}
            { X_i }
            { A }
        }
      &
      \mbox{\small by \eqref{ValuesOfMappingStackAsHomSpaces} }
      \\
      &
      \;\simeq\;
      \PointsMaps{\big}
        { U }
        {
          \underset{ \underset{ i \in I }{\longleftarrow} }{\mathrm{lim}}
          \,
          \Maps{}
            { X_i }
            { A }
        }
      &
      \mbox{\small by \eqref{HomFunctorRespectsLimits}  }
    \end{array}
  $$

  \vspace{-2mm}
  \noindent
  With this, the claim follows by the $\infty$-Yoneda lemma
  (Prop. \ref{InfinityYonedaLemma}).
\end{proof}

\begin{proposition}[Base change preserves mapping stacks]
\label{BaseChangePreservesMappingStacks}
For $B_1 \xrightarrow{f} B_2$ any morphism in
an $\infty$-topos $\Topos$
we have for $X,Y \in \Topos_{/B_2}$ a natural equivalence
\vspace{-2mm}
\begin{equation}
  \label{EquivalenceExhibitingBaseChangePreservingMappingStacks}
  f^\ast \Maps{}{X}{Y}
  \;\simeq\;
  \Maps{}
    { f^\ast(X) }
    { f^\ast(Y) }
  \;\;\;
  \in
  \;
  \Topos_{/B_!}
\end{equation}

\vspace{-2mm}
\noindent between the base change (Def. \ref{BaseChange})
of the mapping stack \eqref{InternalHomAdjunction}
and the mapping stack of the base changes.
\end{proposition}
\begin{proof}
For any $U \in \Topos_{/B_1}$, consider the following
sequence of natural equivalences of $\InfinityGroupoids$
\vspace{-2mm}
$$
\def\arraystretch{1.5}
  \begin{array}{lll}
    \Topos_{/B_1}
    \big(
      U,
      f^\ast \Maps{}{X}{Y}
    \big)
    & \;\simeq\;
    \Topos_{/B_2}
    \big(
      f_! U
      ,\,
      \Maps{}{X}{Y}
    \big)
    &
    \mbox{\small  by \eqref{BaseChangeAdjointTriple} }
    \\
    & \;\simeq\;
    \Topos_{/B_2}
    \big(
      (f_! U) \times_{B_2} X
      ,
      Y
    \big)
    &
    \mbox{\small  by \eqref{InternalHomAdjunction} }
    \\
    & \;\simeq\;
    \Topos_{/B_2}
    \big(
      f_! ( U \times f^\ast (X) )
      ,
      Y
    \big)
    &
    \mbox{\small  by \eqref{FrobeniusReciprocityEquivalence} }
    \\
    & \;\simeq\;
    \Topos_{/B_1}
    \big(
      U \times f^\ast (X)
      ,
      f^\ast(Y)
    \big)
    &
    \mbox{\small  by \eqref{BaseChangeAdjointTriple} }
    \\
    & \simeq
    \Topos_{/B_1}
    \big(
      U
      ,\,
      \Maps{}
        { f^\ast (X) }
        { f^\ast(Y) }
    \big)
    &
    \mbox{\small  by \eqref{InternalHomAdjunction} }
    .
  \end{array}
$$

  \vspace{-2mm}
\noindent
With this, the claim follows by the $\infty$-Yoneda lemma
(Prop. \ref{InfinityYonedaLemma}).
\end{proof}

\noindent
{\bf Slice mapping stacks.}
\begin{proposition}[Hom-$\infty$-groupoids in slice toposes {\cite[Prop. 5.5.5.12]{Lurie09HTT}}]
  \label{HomSpaceInSliceToposAsFiberProduct}
  Let $\Topos$ an $\infty$-topos over a base $\infty$-topos
  $\BaseTopos$.
  Then for $B \in \Topos$ and $(E_1,p_1), (E_2,p_2) \in \Topos_{/B}$,
  the hom-$\infty$-groupoid \eqref{HomSpace}
  in the slice $\Topos_{/B}$ (Prop. \ref{SliceInfinityTopos}) is
  naturally identified with the fiber of the plain hom-$\infty$-groupoid:
  \vspace{-3mm}
  \begin{equation}
    \label{HomSpaceInSliceAsFiberProduct}
    \hspace{-1cm}
    \Topos_{/B}
      \left(
        E_1, \,  E_2
      \right)
    \;\simeq\;
    \Topos(E_1, E_2)
      \underset
        {\Topos(E_1, B)}
        {\times}
      \{p_1\}
      \qquad \mbox{\rm i.e.,}\qquad
      \begin{tikzcd}
        \Topos_{/B}(E_1, E_2)
        \ar[rr]
        \ar[d]
        \ar[
          drr,
          phantom,
          "\mbox{\tiny\rm(pb)}"
        ]
        &&
        \Topos(E_1, E_2)
        \ar[
          d,
          "{ \Topos(E_1, p_2) }"
        ]
        \\
        \ast
        \ar[
          rr,
          "\vdash p_1"
        ]
        &&
        \Topos(E_1, B)
      \end{tikzcd}
  \end{equation}
\end{proposition}

\begin{proposition}[Slice mapping space is space of sections of pullback bundle]
  \label{SliceMappingStackIsStackOfSectionsOfPullbackBundle}
  There is a natural equivalence
  \vspace{-2mm}
  $$
    \Topos_{/B}
    \left(E_1, E_2 \right)
    \;\simeq\;
    \Topos_{/E_1}
    \left(
      E_1,\, p_1^\ast E_2
    \right)
    \,.
  $$
\end{proposition}
\begin{proof}
Factoring
the bottom morphism $\vdash p_1$ in
\eqref{HomSpaceInSliceAsFiberProduct}
as $p_1 \circ (\vdash \mathrm{id}_{E_1})$
and using the pasting law \eqref{HomotopyPastingLaw}
yields the following
identification:
\vspace{-2mm}
$$
  \begin{aligned}
      \begin{tikzcd}
        \Topos_{/B}(E_1, E_2)
        \ar[rr]
        \ar[d]
        \ar[
          drr,
          phantom,
          "\mbox{\tiny\rm(pb)}"
        ]
        &&
        \Topos(E_1, E_2)
        \ar[
          d,
          "{ \Topos(E_1, p_2) }"{description}
        ]
        \\
        \ast
        \ar[
          rr,
          "\vdash p_1"
        ]
        &&
        \Topos(E_1, B)
      \end{tikzcd}
    &
    \quad \simeq\qquad
      \begin{tikzcd}
        \Topos_{/E_1}(E_1, p_1^\ast E_2)
        \ar[r]
        \ar[d]
        \ar[
          dr,
          phantom,
          "\mbox{\tiny\rm(pb)}"
        ]
        &
        \Topos(E_1, p_1^\ast E_2)
        \ar[r]
        \ar[
          dr,
          phantom,
          "\mbox{\tiny\rm(pb)}"
        ]
        \ar[d]
        &
        \Topos(E_1, E_2)
        \ar[
          d,
          "{ \Topos(E_1, p_2) }"
        ]
        \\
        \ast
        \ar[
          r,
          "\vdash \mathrm{id}_{E_1}"{below}
        ]
        &
        \Topos(E_1,E_1)
        \ar[
          r,
          "{\Topos(E_1,p_1)}"{below}
        ]
        &
        \Topos(E_1, B)
      \end{tikzcd}
  \end{aligned}
$$

\vspace{-2mm}
\noindent
Here the pullback square on the far right is
identified as shown by the fact that
$\Topos(E_1,-)$ preserves all limits, hence in particular the pullback
of $p_2$ along $p_1$.  Hence the claim follows by recognizing the left
square on the right as the claimined hom-space over $E_1$, by
Prop. \ref{HomSpaceInSliceToposAsFiberProduct}.
\end{proof}

\begin{definition}[Slice mapping stack]
\label{SliceMappingStack}
For $B \in \Topos$ and $(E_1,p_1), (E_2,p_2) \in \SliceTopos{B}$
a pair of objects in the slice over $B$ (Prop. \ref{SliceInfinityTopos}),
we say that the
{\it slice mapping stack} between them is the object of $\Topos$
which is the fiber of the $\Topos$-internal hom between their total space objects:
\vspace{-4mm}
\begin{equation}
  \label{TheSliceMappingStack}
  \hspace{-5mm}
  \SliceMaps{}{B}
   {E_1}{E_2}
  \;\coloneqq\;
  \Maps{}
    {E_1}{E_2}
    \underset
      {
        \Maps{}
          {E_1}{B}
      }
      {\times}
   \!\! \{ p_1 \},
  \quad
  \mbox{\rm i.e.,}
  \quad
  \begin{tikzcd}[column sep=small]
    {
      \SliceMaps{}{B}
        {E_1}{E_2}
    }
    \ar[rr]
    \ar[d]
    \ar[
      drr,
      phantom,
      "{\mbox{\tiny\rm(pb)}}"
    ]
    &&
    {
      \Maps{}
        {E_1}{E_2}
    }
    \ar[
      d,
      "{
        \Maps{}
          {E_1}{p_2}
      }"
    ]
    \\
    \ast
    \ar[
      rr,
      "\vdash p_1"
    ]
    &&
    {
      \Maps{}
        {E_1}{B}
    }
  \end{tikzcd}
  \;\;\;\;
  \in
  \;
  \Topos\;.
\end{equation}
\end{definition}
When $E_1 = B$, this may also be called the {\it stack of sections} of $E_2$:
\begin{equation}
  \label{SliceMappingStackOutOfBaseIsSpaceOfSections}
  \Gamma_B(E)
  \;:=\;
  \SliceMaps{}{B}
    {B}{E}
  \,.
\end{equation}
\begin{remark}
  The slice mapping stack (Def. \ref{SliceMappingStack}) only depends on the connected component $\mathbf{B}\Omega_{p_1} \Maps{}{E_1}{B}$
  of $\vdash p_1$ in $\Maps{}{E_1}{B}$ in that it is equivalently the top fiber product in the following diagram (obtained by applying the pasting law \eqref{HomotopyPastingLaw} to the homotopy image factorization as recalled e.g.in \cite[Ex. 2.67]{SS20OrbifoldCohomology})
  \begin{equation}
    \label{SliceMappingStackAsFiberOverConnectedComponentOfCodomainProjection}
    \begin{tikzcd}
      \SliceMaps{}{B}{E_1}{E_2}
      \ar[rr]
      \ar[d,->>]
      \ar[
        drr,
        phantom,
        "{\mbox{\tiny(pb)}}"
      ]
      &&
      \{p_1\}
      \ar[d,->>]
      \\
      \Maps{}{E_1}{E_2}^{p_1}
      \ar[d,hook]
      \ar[rr]
      \ar[
        drr,
        phantom,
        "{
          \mbox{\tiny(pb)}
        }"
      ]
      &&
      \mathbf{B}
      \Omega_{p_1}
      \Maps{}
        {E_1}{B}
      \ar[d,hook]
      \\
      \Maps{}{E_1}{E_2}
      \ar[
        rr,
        "{
          \Maps{}{E_1}{p_2}
        }"{swap}
      ]
      &&
      \Maps{}{E_1}{B}
    \end{tikzcd}
  \end{equation}
\end{remark}

\begin{example}[Slice mapping stack into product projection]
  \label{SliceMappingStackIntoProductProjection}
  In the case when $E_2 \,=\, X \times B$ with $p_2 \,=\, \mathrm{pr}_2$
  the projection on the second factor, we have
  $\Maps{}{E_1}{E_2} \,\simeq\, \Maps{}{E_1}{X} \times \Maps{}{E_1}{B}$
  with $\Maps{}{E_1}{p_2}$
  being again the projection onto the second factor,
  so that the slice mapping stack in Def. \ref{SliceMappingStack}
  reduces to the plain mapping stack
  (by Ex. \ref{HomotopyPullbackPreserbesProductProjections}):
  $$
    \SliceMaps{}{B}{E_1}{X \times B}
    \;\simeq\;
    \Maps{}{E_1}{X}
    \,.
  $$
\end{example}

\begin{lemma}[Plots of slice mapping stack are slice homs]
  \label{PlotsOfSliceMappingStackAreSliceHoms}
The slice mapping stack (Def. \ref{SliceMappingStack})
is equivalently the right base change (Prop. \ref{BaseChange})
of the internal hom in the slice $\infty$-topos down to the base,
in that

\noindent {\bf (i)} for
$B, U \,\in\, \Topos$ and
$(E_1, p_1), (E_2,p_2) \,\in\, \SliceTopos{B} $ (Prop. \ref{SliceInfinityTopos}),
we have a natural equivalence
\begin{equation}
  \label{SliceMappingStackIsRightBaseChangeOfMappingStackInSlice}
  \SliceMaps{}{B}
    {E_1}{E_2}
  \;\simeq\;
  \underset{B}{\prod} \,
  \Maps{}
    {p_1}{p_2}
  \,.
\end{equation}

\vspace{-1mm}
\noindent
{\bf (ii)} This means that there is a natural equivalence
between the value of the slice mapping stack
  at $U$ and the slice hom-space (Def. \ref{HomSpaceInSliceToposAsFiberProduct})
  between $(U \times E_1,  \mathrm{pr}_2 \circ p_1 )$ and $(E_2, p_2)$:
  \vspace{-.5mm}
  \begin{equation}
    \label{PlotsOfSliceMappingStack}
   \PointsMaps{\big}
      { U }
      {
        \SliceMaps{}{B}
          { E_1}{ E_2 }
      }
    \;\simeq\;
    \SlicePointsMaps{}{B}
      { U \times E_1 }
      { E_2 }
    \,.
  \end{equation}

\vspace{-1mm}
\noindent
{\bf (iii)}
In particular, the global points of the slice mapping stack
constitute the slice mapping space \eqref{HomSpaceInSliceAsFiberProduct}:
\vspace{-1mm}
\begin{equation}
  \label{GlobalPointsOfSliceMappingStackIsSliceHom}
  \PointsMaps{\big}
    { \ast }
    {
      \SliceMaps{}{B}
      {E_1}{E_2}
    }
    \;\;
    \simeq
    \;\;
    \SlicePointsMaps{}{B}
      { E_1 }
      { E_2}
    \,.
\end{equation}
\noindent
{\bf (iv)}
  and the space of sections \eqref{SliceMappingStackOutOfBaseIsSpaceOfSections} is equivalently the right base change to the point:
  \begin{equation}
    \label{SpaceOfSectionsIsRightBaseChangeToPoint}
    \Gamma_B(E)
    \;=\;
    \SliceMaps{}{B}
      { B }{ E }
    \;\simeq\;
    \underset{B}{\prod}
    \,
    (E,p)
    \,.
  \end{equation}

\vspace{-.3cm}
\noindent
{\bf (v)}
Generally, thinking of $E_1 = U \xrightarrow{i_U} B$
as a local patch of $B$, the slice mapping stack
gives the
{\it $U$-local sections}:
\begin{equation}
  \label{SpaceOfLocalSectionsViaSliceHom}
  \Gamma_U(E)
  \;:=\;
  \Gamma_U
  \big(
    i_U^\ast E
  \big)
  \;\simeq\;
  \SliceMaps{}{B}
    { U }
    { E }
  \;\simeq\;
  \underset{B}{\prod}
    \;
    i_U^\ast (E,p)
  \,.
\end{equation}
\end{lemma}
\vspace{-.4cm}
\begin{proof}
The first equivalence is the pre-image
under the $\infty$-Yoneda lemma (Prop. \ref{InfinityYonedaLemma})
of the composite of the following sequence of
natural equivalences for $U \,\in\, \Topos$:
\vspace{-2mm}
$$
\hspace{-1mm}
  \def\arraystretch{1.7}
  \begin{array}{lll}
    \PointsMaps{\big}
      { U }
      {  \underset{B}{\prod}
        {\rm Maps}
       \left( {(E_1,p_1)}, {(E_2,p_2)} \right) \!\!
          }
    &
     \!\!\!\! \simeq
    \SlicePointsMaps{\big}{B}
      { (U \times B, \mathrm{pr}_2) }
      {
      {\rm Maps}
       \left( {(E_1,p_1)}, {(E_2,p_2)} \right) \!\!
       }
    &
    \proofstep{ by \eqref{AdjunctionAndHomEquivalence} with Prop. \ref{BaseChange} }
    \\
    &
         \!\!\!\!\simeq
    \SlicePointsMaps{\big}{B}
      { (U \times E_1, p_1\circ \mathrm{pr}_2) }
      { (E_2, p_2) }
    &
    \proofstep{ by \eqref{ValuesOfMappingStackAsHomSpaces} with Prop. \ref{SliceInfinityTopos} }
    \\
    &
         \!\!\!\! \simeq
    \PointsMaps{}
      { U \times E_1 }
      { E_2 }
    \underset{
      \PointsMaps{}
        { U \times E_1 }
        { B }
    }{\times}
    \{ p_1 \circ \mathrm{pr}_2 \}
    &
    \proofstep{ by Prop. \ref{HomSpaceInSliceToposAsFiberProduct} }
    \\
    &
         \!\!\!\!\simeq
    \PointsMaps{\big}
      { U }
      { \Maps{}{E_1}{E_2} }
    \underset{
      \PointsMaps{}
        { U }
        { \Maps{}{E_1}{B} }
    }{\times}
    \PointsMaps{}
      { U }
      {
        \{ p_1 \}
      }
    &
    \proofstep{ by \eqref{ValuesOfMappingStackAsHomSpaces}  }
    \\
    &
    \!\!\!\!\simeq
    \PointsMaps{\Big}
      { U }
      {
        \Maps{}{E_1}{E_2}
          \underset{
            \Maps{}{E_1}{B}
          }{\times}
        \{ p_1  \}
      }
    &
    \proofstep{ by \eqref{InfinityAdjointPreservesInfinityLimits} }
    \\
    &
    \!\!\!\!\simeq
    \PointsMaps{\big}
      { U }
      {
        \SliceMaps{}{B}
          { E_1 }
          { E_2 }
      }
    & \proofstep{ by Def. \ref{SliceMappingStack}. }
  \end{array}
$$

\vspace{-1mm}
\noindent With this, the second statement is equivalently the first step
in this sequence.
\end{proof}

\medskip

\subsection{Transformation groups in $\infty$-toposes}
\label{TransformationGroupsInInfinityToposes}

We discuss basics of {\it higher transformation group theory},
i.e., of higher groups and their group actions, internal to $\infty$-toposes.

\medskip

\noindent
{\bf Group objects in an $\infty$-topos.}
\begin{definition}[Group objects in an $\infty$-topos]
  \label{GroupObjectsInAnInfinityTopos}
  Given an $\infty$-topos $\Topos$,

\noindent {\bf (i)}  a {\it group object} $\mathcal{G}_{\bullet}$
  is a groupoid object (Def. \ref{GroupoidObjectInAnInfinityTopos}),
  to be denoted $(\mathbf{B}\mathcal{G})_\bullet$, which is
  equipped with a morphism
  $\ConstantGroupoid(\ast)_\bullet \xrightarrow{\mathrm{pt}} (\mathbf{B}\mathcal{G})_\bullet$
  from the terminal groupoid object
  \eqref{TerminalGroupoidObjectInAnInfinityTopos},
  such that this is an equivalence in degree 0:
  $(\mathbf{B}\mathcal{G})_0 \,\simeq\, \ast$.

 \noindent {\bf (ii)}  We write
  \begin{equation}
    \label{InfinityCategoryOfGroupObjects}
    \Groups(\Topos)
    \xhookrightarrow{\quad}
    \Groupoids(\Topos)^{\ast/}
  \end{equation}
  for the full sub-$\infty$-category
  of group objects.
  among the pointed groupoid objects.
\end{definition}

\begin{proposition}[Looping and delooping equivalence]
  \label{LoopingAndDeloopingEquivalence}
  Under the equivalence of groupoid objects with atlases
  (Prop. \ref{EquivalentPerspectivesOnGroupoidObjectsInAnInfinityTopos})
   group objects $\mathcal{G}$ (Def. \ref{GroupObjectsInAnInfinityTopos})
  are identified as the loop objects
  $\mathcal{G} \,\simeq\, \Omega \mathbf{B}\mathcal{G}$ of
  pointed connected objects $\ast \twoheadrightarrow{\;} \mathbf{B}\mathcal{G}$,
  \vspace{-2mm}
  \begin{equation}
    \begin{tikzcd}[row sep=30pt, column sep=huge]
      \Topos^{\ast/}_{\geq  0}
      \ar[rr, phantom, "\sim"]
      \ar[rr, shift right=5pt,  "\Omega(-)"{below}]
      \ar[from=rr, shift right=5pt, "\mathbf{B}(-)"{above}]
      \ar[d, hook]
      &&
      \Groups(\Topos)
      \ar[d, hook]
      \\
      \EffectiveEpimorphisms(\Topos)^{(\ast \to \ast)/}
      \ar[from=rr, shift right=5pt, " (-)_0 \,\to \, \colimit{}(-) "{above, yshift=-5pt}]
      \ar[rr, phantom, "\sim"]
      \ar[rr, shift right=5pt, "{ \Cech(-)^{(\ast \to \ast)/}_\bullet }"{below}]
      &&
      \Groupoids(\Topos)^{\ConstantGroupoid(\ast)_{\bullet}/}
      \mathrlap{\,,}
    \end{tikzcd}
  \end{equation}
\vspace{-1mm}
\noindent
where the {\it delooping} $\mathbf{B}\mathcal{G}$ is the homotopy colimit:
\vspace{-1mm}
\begin{equation}
  \label{DeloopingOfInfinityGroupAsColimit}
  \begin{tikzcd}
    \mathbf{B}\mathcal{G}
    \;\coloneqq\;
    \colimit{[n] \in \Delta^{\mathrm{op}}}
    \mathcal{G}^{\times_n}
    &&
    \ast
    \ar[
      ll,
      ->>,
      "{
        q
      }"{above}
    ]
    \ar[
      r
    ]
    &
    \mathcal{G}
    \ar[
      r,
      shift left=4pt
    ]
    \ar[
      r,
      shift right=4pt
    ]
    \ar[
      l,
      shift right=4pt
    ]
    \ar[
      l,
      shift left=4pt
    ]
    &
    \mathcal{G} \times \mathcal{G}
    \ar[
      l,
      shift left=8pt
    ]
    \ar[
      l
    ]
    \ar[
      l,
      shift right=8pt
    ]
    \ar[r, -, dotted]
    &[-15pt]
    {}
    \;.
  \end{tikzcd}
\end{equation}
\end{proposition}
\begin{proof}
By the equivalence at the bottom (Prop. \ref{EquivalentPerspectivesOnGroupoidObjectsInAnInfinityTopos})
and by the
definition of $\Groups(\Topos)$ (Def. \ref{GroupObjectsInAnInfinityTopos}),
the $\infty$-category
in the top left is the full sub-$\infty$-category of
that on the bottom left on the point atlases.
These are equivalently the
pointed connected objects,
by Prop. \ref{PointedConnectedObjectsHaveEffectiveEpimorphicPointInclusion}.
Under this identification
\eqref{RegardingAPointedConnectedObjectAsAPointedAtlas}
the {\v C}ech nerve construction on a point atlas
$\ast \twoheadrightarrow{\;} \mathbf{B}\mathcal{G}$
clearly is the loop space construction:
$$
  \ast
    \underset{\mathbf{B}\mathcal{G}}{\times}
  \ast
  \;\simeq\;
  \Omega \mathbf{B}\mathcal{G}
  \;\simeq\;
  \mathcal{G}
  \,.
$$

\vspace{-7mm}
\end{proof}

\begin{example}[Delooping preserves products]
  The operation of delooping (Prop. \ref{LoopingAndDeloopingEquivalence})
  preserves products, in that the delooping of the direct product of
  $\mathcal{G}_1, \, \mathcal{G}_2 \,\in\, \Groups(\Topos)$ is the
  naturally equivalent to the product of the separate deloopings:
  \begin{equation}
    \label{DeloopingPreservesProducts}
    \mathbf{B}(\mathcal{G}_1 \times \mathcal{G}_2)
    \;\;
    \simeq
    \;\;
    (\mathbf{B}\mathcal{G}_1)
    \times
    (\mathbf{B}\mathcal{G}_2)
    \,.
  \end{equation}
  This follows from the inverse looping equivalence in Prop. \ref{LoopingAndDeloopingEquivalence},
  the fact that $\mathbf{B}\mathcal{G}_1 \times \mathbf{B}\mathcal{G}_2$ is connected,
  and using that looping, being itself a limit operation, commutes with products:
  $$
    \Omega
    \big(
      (\mathbf{B}\mathcal{G}_1)
      \times
      (\mathbf{B}\mathcal{G}_2)
    \big)
    \;\simeq\;
    (\Omega \mathbf{B}\mathcal{G}_1)
    \times
    (\Omega \mathbf{B}\mathcal{G}_2)
    \;\simeq\;
    \mathcal{G}_1 \times \mathcal{G}_1
    \,.
  $$
\end{example}

\begin{lemma}[$\infty$-Groups presented by presheaves of simplicial groups
 {\cite[Prop. 3.35, 3.73, Rem. 3.67]{NSS12b}}]
\label{InfinityGroupsPresentedByPresheavesOfSimplicialGroups}
Let $(\mathcal{S}, J)$ be a 1-site with a terminal object and
$\Topos \coloneqq Sh_{\infty}\left((\mathcal{C},J)\right)$ its $\infty$-topos.

\noindent
{\bf (i)} Presheaves of simplicial groups represent $\infty$-groups
\vspace{-2mm}
$$
  \mathcal{G}
  \;\in\;
  \Groups
  \left(
    \SimplicialPresheaves
    (
      \mathcal{S}
    )
  \right)
  \xrightarrow{\;\Localization{J}\;}
  \Groups
  (
    \Topos
  )
$$

\vspace{-1mm}
\noindent
with delooping equivalent to the objectwise $\overline{W}(-)$
\eqref{StandardSimplicialDeloopingAsQuotient}
\vspace{-2mm}
$$
  \mathbf{B} (\mathcal{G})
  \;\simeq\;
  \Localization{J}
  \left(
    \overline{W}(\mathcal{G})
  \right)
  \,.
$$

\vspace{-2mm}
\noindent
{\bf (ii)} For $\mathcal{G} \acts \, X$, the homotopy quotient is represented by
the right derived
simplicial Borel constuction functor
\eqref{QuillenEquivalenceBetweenBorelModelStructureAndSliceOverSimplicialClassifyingComplex}
given by the diagonal quotient of the Cartesian product with
$W(-)$ (Def. \ref{StandardModelOfUniversalSimplicialPrincipalComplex}):

\vspace{-2mm}
$$
  \Localization{J}(X) \sslash \Localization{J}(\mathcal{G})
  \;\simeq\;
  \Localization{J}
  \left(
    (X \times W\mathcal{G})/\mathcal{G}
  \right)
  \,.
$$
\end{lemma}
We consider examples of this construction in Ex. \ref{ProjectiveRepresentationsAndTheirCentralExtensions}
below.

\begin{proposition}[All bare $\infty$-groups are shapes of topological groups]
  \label{AllBareInfinityGroupsAreShapesOfTopologicalGroups}
  For every bare $\infty$-group
  $\mathcal{G} \,\in\, \Groups(\InfinityGroupoids)$,
  there exists a Hausdorff topological group
  $\Gamma \,\in\, \Groups(\kHausdorffSpaces) \xrightarrow{ \ContinuousDiffeology }
  \Groups(\SmoothInfinityGroupoids)$,
  such that
    \vspace{-2mm}
  $$
    \shape \,\Gamma\,
    \;\simeq\;
    \mathcal{G}
    \;\;\;
    \in
    \;
    \Groups(\InfinityGroupoids)
    \xhookrightarrow{ \Groups(\Discrete) }
    \Groups(\SmoothInfinityGroupoids)
    \,.
  $$
\end{proposition}
\begin{proof}
  By Lem. \ref{InfinityGroupsPresentedByPresheavesOfSimplicialGroups},
  we have
  \vspace{-2mm}
  $$
    \Groups(\InfinityGroupoids)
    \;\;
    \simeq
    \;\;
    \SimplicialLocalization{\WeakHomotopyEquivalences}
    \left(
      \Groups(\SimplicialSets)
    \right)
    \,.
  $$

    \vspace{-2mm}
\noindent
  Under this equivalence, we may think of the given
  $\infty$-group as a simplicial group $\mathcal{G}$.

  Now since topological realization \eqref{TopologicalRealizationOfSimplicialTopologicalSpaces}
  preserves finite products (Lem. \ref{TopologicalRealizationPreservesFiniteLimits}),
  the adjoint pair Quillen equivalence
    \vspace{-2mm}
  \begin{equation}
    \label{QuillenEquivalenceBetweenClassicalModelStructures}
    \begin{tikzcd}
      \TopologicalSpaces_{\mathrm{Qu}}
      \ar[
        from=rr,
        shift right=6pt,
        "{\TopologicalRealization{}{-}}"{swap}
      ]
      \ar[
        rr,
        shift right=6pt,
        "\SingularSimplicialComplex"{swap}
      ]
      \ar[
        rr,
        phantom,
        "\scalebox{.7}{$\simeq_{\mathrlap{\mathrm{Qu}}}$}"
      ]
      &&
      \SimplicialSets_{\mathrm{Qu}}
    \end{tikzcd}
  \end{equation}
  induces \eqref{FunctorOnStructuresInducedFromLexFunctor}
  an adjoint pair between group objects
  \vspace{-2mm}
  $$
    \begin{tikzcd}
      \Groups(\TopologicalSpaces)
      \ar[
        from=rr,
        shift right=5pt,
        "{\Groups(\TopologicalRealization{}{-})}"{swap}
      ]
      \ar[
        rr,
        shift right=5pt,
        "\Groups(\SingularSimplicialComplex)"{swap}
      ]
      \ar[
        rr,
        phantom,
        "\scalebox{.7}{$\bot$}"
      ]
      &&
      \Groups(\SimplicialSets) \;.
    \end{tikzcd}
    $$

  \vspace{-2mm}
\noindent
  Hence
  $$
    \Gamma
      \,\coloneqq\,
    \TopologicalRealization{}{\mathcal{G}}
    \,,
  $$
  is a Hausdorff topological group (since all CW complexes are Hausdorff,
  e.g. \cite[Porp. A.3]{Hatcher02})
  such that
  $$
    \begin{array}{lll}
      \shape \, \Gamma
      &
      \;\simeq\;
      \SingularSimplicialComplex
      \,
      \TopologicalRealization{}{\mathcal{G}}
      &
      \proofstep{ by Prop. \ref{SmoothShapeGivenBySmoothPathInfinityGroupoid} }
      \\
      &
      \;\simeq\;
      \mathcal{G}
      \qquad
      \in
      \;
      \SimplicialLocalization{\WeakHomotopyEquivalences}
      \left(
        \Groups(\SimplicialSets)
      \right)
      &
      \proofstep{ by \eqref{QuillenEquivalenceBetweenClassicalModelStructures} }
      \,.
    \end{array}
  $$

  \vspace{-6mm}
\end{proof}

\medskip

\noindent
{\bf Group actions in an $\infty$-topos.}

\begin{definition}[Group actions internal to an $\infty$-topos]
  \label{ActionObjectsInAnInfinityTopos}
  Given an $\infty$-topos $\Topos$,
  for
  $X \,\in\, \Topos$
  and
  $\mathcal{G} \,\in\, \Groups(\Topos)$,

\noindent {\bf (i)}   we say that a {\it $\mathcal{G}$-action groupoid} on $X$ is
  a groupoid object (Def. \ref{GroupoidObjectInAnInfinityTopos})
  which in degree $k$ is equivalent to the cartesian product of $X$
  with $k$ copies of $\mathcal{G}$:
    \vspace{-2mm}
  $$
    (\HomotopyQuotient{X}{\mathcal{G}})_\bullet
    \;\in\;
    \Groupoids(\Topos)
    \,,
    \;\;\;\;\;\;\;\;\;\;
    (\HomotopyQuotient{X}{\mathcal{G}})_k
    \;\simeq\;
    X \times \mathcal{G}^{\times^k}
    \,,
  $$

    \vspace{-2mm}
\noindent
  compatibly with the group structure on $\mathcal{G}$, hence for which
  the two simplicial objects
  form a diagram of homotopy Cartesian squares:
  \vspace{-3mm}
  \begin{equation}
    \label{ActionGroupoidInInfinityToposAsCartesianMorphismOfGroupoidObjectsToGroupObject}
    \begin{tikzcd}
      {}
      \ar[d, -, dotted]
      &&
      {}
      \ar[d, -, dotted]
      \\[-10pt]
      \mathllap{
        \mathcal{G} \times \mathcal{G}
          \times
        X
        \;\simeq\;
      }
      (\HomotopyQuotient{X}{\mathcal{G}})_2
      \ar[d, shift left=10pt]
      \ar[from=d, shift left=5pt]
      \ar[d]
      \ar[from=d, shift right=5pt]
      \ar[d, shift right=10pt]
      \ar[rr, "{ c_2 }"{description}]
      \ar[drr, phantom, "{ \mbox{\tiny \rm (pb)} }"]
      &&
      \mathcal{G} \times \mathcal{G}
      \ar[d, shift left=10pt]
      \ar[from=d, shift left=5pt]
      \ar[d]
      \ar[from=d, shift right=5pt]
      \ar[d, shift right=10pt]
      \\
      \mathllap{
        \mathcal{G}
          \times
        X
        \;\simeq\;
      }
      (\HomotopyQuotient{X}{\mathcal{G}})_{1}
      \ar[d, shift left=5pt]
      \ar[from=d]
      \ar[d, shift right=5pt]
      \ar[rr, "{ c_1 }"{description}]
      \ar[d]
      \ar[drr, phantom, "{ \mbox{\tiny \rm (pb)} }"]
      &&
      \mathcal{G}
      \ar[d, shift left=5pt]
      \ar[from=d]
      \ar[d, shift right=5pt]
      \\
      \mathllap{
        X
        \;\simeq\;
      }
      (\HomotopyQuotient{X}{\mathcal{G}})_{0}
      \ar[rr, "{ c_0 }"{description}]
      &&
      \ast
    \end{tikzcd}
  \end{equation}

\noindent {\bf (ii)}   We write
\vspace{-1mm}
  \begin{equation}
    \label{InfinityCategoryOfInfinityActionsInInfinityTopos}
    \begin{tikzcd}[row sep=-3pt, column sep=3pt]
      \Actions{\mathcal{G}}(\Topos)
      \ar[rr, hook]
      &
      &
      \Groupoids(\Topos)_{/(\mathbf{B}\mathcal{G})_\bullet}
      \\
  \scalebox{0.7}{$    G \acts \, X $}
      &\longmapsto&
       \scalebox{0.7}{$  (\HomotopyQuotient{X}{\mathcal{G}})_\bullet $}
    \end{tikzcd}
  \end{equation}

  \vspace{-1mm}
\noindent
  for the full sub-$\infty$-category of the slice of
  that of groupoid objects over the groupoid object underlying $\mathcal{G}$.
\end{definition}

\begin{proposition}[Homotopy quotients of $\infty$-actions]
  \label{HomotopyQuotientsAndPrincipaInfinityBundles}
  Under the equivalence of groupoid objects with atlases
  (Prop. \ref{EquivalentPerspectivesOnGroupoidObjectsInAnInfinityTopos}),
  $\mathcal{G}$-actions (Def. \ref{ActionObjectsInAnInfinityTopos})
  are identified with objects in the slice of $\Topos$ over $\mathbf{B}\mathcal{G}$:
  \vspace{-2mm}
  $$
    \begin{tikzcd}[column sep=huge]
      \Slice{\Topos}{\mathbf{B}\mathcal{G}}
      \ar[
        from=rr,
        shift right=5pt,
        "{
          (X \to \HomotopyQuotient{X}{\mathcal{G}})
             \;\;\mapsfrom\;\;
          G \acts \, X
        }"{above}
      ]
      \ar[rr, phantom, "\sim"]
      \ar[rr, shift right=5pt]
      \ar[d, hook]
      &&
      \Actions{\mathcal{G}}(\Topos)
      \ar[d, hook]
      \\
      \EffectiveEpimorphisms(\Topos)_{/(\ast \to \mathbf{B}\mathcal{G})}
      \ar[from=rr, shift right=5pt, "(-)_0 \, \to  \, \colimit{}(-)"{above, yshift=-5pt}]
      \ar[rr, phantom, "\sim"]
      \ar[rr, shift right=5pt, "{ \Cech(-)_\bullet^{/(\ast \to \mathbf{B}\mathcal{G})} }"{below}]
      &&
      \Groupoids(\Topos)_{/(\mathbf{B}\mathcal{G})_\bullet}
      \mathrlap{\,,}
    \end{tikzcd}
  $$

  \vspace{-2mm}
\noindent
where the base is the {\it homotopy quotient}
\vspace{-2mm}
\begin{equation}
  \label{HomotopyQuotientAsHomotopyColimit}
  \begin{tikzcd}
    \HomotopyQuotient{ X }{\mathcal}{G}
    \;\simeq\;
    \colimit{[n] \in \Delta^{\mathrm{op}}}
    \mathcal{G}^{\times_n}
      \times
    X
    &&
    X
    \ar[
      ll,
      ->>,
      "{
        q
      }"{above}
    ]
    \ar[
      r
    ]
    &
    \mathcal{G}
      \times
    X
    \ar[
      r,
      shift left=4pt
    ]
    \ar[
      r,
      shift right=4pt
    ]
    \ar[
      l,
      shift right=4pt
    ]
    \ar[
      l,
      shift left=4pt
    ]
    &
    \mathcal{G} \times \mathcal{G}
      \times
    X
    \ar[
      l,
      shift left=8pt
    ]
    \ar[
      l
    ]
    \ar[
      l,
      shift right=8pt
    ]
    \ar[r, -, dotted]
    &[-15pt]
    {}.
  \end{tikzcd}
\end{equation}
\end{proposition}
\begin{proof}
  By the bottom equivalence (Prop. \ref{EquivalentPerspectivesOnGroupoidObjectsInAnInfinityTopos})
  and the definition of $\Actions{\mathcal{G}}(\Topos)$ (Def. \ref{ActionObjectsInAnInfinityTopos}),
  the $\infty$-category
  that appears in the top left must
  be the full sub-$\infty$-category of
  $\Slice{\mathrm{Atl}(\Topos)}{(\ast \to \mathbf{B}\mathcal{G})}$
  whose objects are those squares
  \vspace{-2mm}
  $$
    \begin{tikzcd}[row sep=small, column sep=large]
      X
      \ar[r, ->>]
      \ar[d]
      &
      \HomotopyQuotient{X}{\mathcal{G}}
      \ar[d]
      \\
      \ast
      \ar[r, ->>]
      &
      \mathbf{B}\mathcal{G}
      \mathrlap{\,,}
    \end{tikzcd}
  $$
  which arise under the colimit operation
  from Cartesian morphisms of
  groupoid objects \eqref{ActionGroupoidInInfinityToposAsCartesianMorphismOfGroupoidObjectsToGroupObject}.
  But these are precisely the Cartesian such squares, by the theory of
  equifibered  transformations (\cite[6.1.3.9(4)]{Lurie09HTT}\cite[6.5]{Rezk10},
  see \cite[Prop. 2.32]{SS20OrbifoldCohomology}).
  (This is consistent, since effective epimorphisms
  are stable under pullback; so that for every such Cartesian square the top horizontal morphism is an
  effective epimorphism, since $\ast \twoheadrightarrow \mathbf{B}\mathcal{G}$ is.)
  But by the universal property of pullbacks,
  morphisms of pullback squares over
  $\ast \twoheadrightarrow \mathbf{B}\mathcal{G}$ are equivalent to morphsims of their right vertical component over $\mathbf{B}\mathcal{G}$, hence to morphisms in $\SliceTopos{\mathbf{B}\mathcal{G}}$.
\end{proof}
\begin{remark}[$\infty$-Topos of $\infty$-Actions]
\label{InfinityToposOfInfinityActions}
A key consequences of Prop. \ref{HomotopyQuotientsAndPrincipaInfinityBundles}
is that $\mathcal{G}$-actions in an
$\infty$-topos form themselves an $\infty$-topos:
$$
  \Actions{\mathcal{G}}(\Topos)
  \;\;
  \simeq
  \;\;
  \SliceTopos{\mathbf{B}\mathcal{G}}
  \;\in\;
  \InfinityToposes
  \,,
$$
by the Fundamental Theorem of $\infty$-topos theory
(Prop. \ref{SliceInfinityTopos}).
\end{remark}

In higher generalization of Exp. \ref{QuotientAndFixedLociFromChangeOfGroupAdjunction} we have:
\begin{example}[Fixed loci of $\infty$-actions]
  \label{FixedLociOfInfinityActions}
  For $\mathcal{G} \,\in\, \Groups(\Topos)$,
  the operations of forming

  (i) the homotopy quotient of $\mathcal{G}$-actions \eqref{HomotopyQuotientAsHomotopyColimit}

  (ii) the trivial $\mathcal{G}$-action

  (iii) the fixed points of a $\mathcal{G}$-action
   \cite[Def. 2.97]{SS20OrbifoldCohomology}

  \noindent
  constitute an adjoint triple of $\infty$-functors
  as shown on the left below, which
  under the equivalence of Prop. \ref{HomotopyQuotientsAndPrincipaInfinityBundles} is identified with the base change adjoint triple
  \eqref{BaseChangeAdjointTriple} along $\mathbf{B}\mathcal{G} \xrightarrow{\;} \ast$:
  \begin{equation}
    \label{InfinityFixedPointsViaBaseChange}
    \begin{tikzcd}[column sep=40pt]
      \Actions{\mathcal{G}}(\Topos)
      \ar[
        rr,
        shift right=14pt,
        "{
          (-)^{\mathcal{G}}
        }"{swap}
      ]
      \ar[
        from=rr,
        "{\mathrm{trivial}}"{description}
      ]
      \ar[
        rr,
        shift left=14pt,
        "{
          \HomotopyQuotient{(-)}{\mathcal{G}}
        }"
      ]
      \ar[
        rr,
        phantom,
        shift left=8pt,
        "{\scalebox{.7}{$\bot$}}"
      ]
      \ar[
        rr,
        phantom,
        shift right=8pt,
        "{\scalebox{.7}{$\bot$}}"
      ]
      &&
      \Topos
    \end{tikzcd}
    \hspace{1.4cm}
    \simeq
    \hspace{1.4cm}
    \begin{tikzcd}[column sep=40pt]
      \SliceTopos{\mathbf{B}\mathcal{G}}
      \ar[
        rr,
        shift right=14pt,
        "{
          \prod_{\mathbf{B}\mathcal{G}}
        }"{swap}
      ]
      \ar[
        from=rr,
        "{\mathbf{B}\mathcal{G}\times (-)}"{description}
      ]
      \ar[
        rr,
        shift left=14pt,
        "{
          \sum_{\mathbf{B}\mathcal{G}}
        }"
      ]
      \ar[
        rr,
        phantom,
        shift left=8pt,
        "{\scalebox{.7}{$\bot$}}"
      ]
      \ar[
        rr,
        phantom,
        shift right=8pt,
        "{\scalebox{.7}{$\bot$}}"
      ]
      &&
      \Topos
      \,.
    \end{tikzcd}
  \end{equation}
  (For the homotopy quotient this is part of the result of Prop. \ref{HomotopyQuotientsAndPrincipaInfinityBundles}, for the trivial action it is evident, and for the fixed locus it is the very definition \cite[Def. 2.97]{SS20OrbifoldCohomology}.)
\end{example}
\begin{example}[Fixed loci as slice mapping stacks into homotopy quotients over deloopings]
  \label{FixedLociAsSliceMappingStackss}
  Via  \eqref{SpaceOfSectionsIsRightBaseChangeToPoint},
  the fixed loci
  \eqref{InfinityFixedPointsViaBaseChange}
  of any $\mathcal{G} \acts X$
  are equivalently given by the slice mapping stack (Def. \ref{SliceMappingStack})
  out of the base $\mathbf{B}\mathcal{G}$
  into the homotopy quotient (cf. \cite[Ex. 2.98]{SS20OrbifoldCohomology}):
  \begin{equation}
    \label{FixedPointsAsSliceMappingStack}
    X^{\mathcal{G}}
    \;\simeq\;
    \SliceMaps{\big}{\mathbf{B}\mathcal{G}}
      { \mathbf{B}\mathcal{G} }
      {
        \HomotopyQuotient{X}{\mathcal{G}}
      }
    \,.
  \end{equation}
More generally, via
\eqref{SpaceOfLocalSectionsViaSliceHom},
for
  $\mathcal{H} \xrightarrow{\;i\;} \mathcal{G}$
  any homomorphism of $\infty$-groups,
hence $\mathbf{B}i : \mathbf{B}\mathcal{H} \xhookrightarrow{\;} \mathbf{B}\mathcal{G}$,
the fixed loci of the induced action $\mathcal{H} \acts X$
are given by the slice mapping stack out of $\mathbf{B}\mathcal{H}$ over $\mathbf{B}\mathcal{G}$:
  \begin{equation}
    \label{FixedPointsForSubgroupAsSliceMappingStack}
    X^{\mathcal{H}}
    \;\simeq\;
    \SliceMaps{\big}{\mathbf{B}\mathcal{G}}
      { \mathbf{B}\mathcal{H} }
      {
        \HomotopyQuotient{X}{\mathcal{G}}
      }
    \,.
  \end{equation}
\end{example}

\begin{example}[$\infty$-Stack of equivariant maps]
 \label{EquivariantMappingStackAsSliceMappingStack}
Given $\mathcal{G} \,\in\, \Groups(\Topos)$ and $\mathcal{G} \acts X, \, \mathcal{G} \acts Y \,\in\, \Actions{\mathcal{G}}(\Topos)$ (Prop. \ref{HomotopyQuotientsAndPrincipaInfinityBundles}), the $\mathcal{G}$-equivariant mapping stack between them is equivalently the slice mapping stack
(Def. \ref{SliceMappingStack}) between their homotopy quotients:
\begin{equation}
  \label{EquivariantMappingStackAsSliceMappingStackEquivalence}
  \def\arraystretch{1.8}
  \begin{array}{ll}
  \Maps{\big}
    { \mathcal{G} \acts X }
    { \mathcal{G} \acts Y }
  ^{ \mathcal{G} }
  \\
  \;\simeq\;
  \underset{\mathbf{B}\mathcal{G}}{\prod}
  \Maps{\big}
    { \HomotopyQuotient{X}{\mathcal{G}} }
    { \HomotopyQuotient{Y}{\mathcal{G}} }
  &
  \proofstep{
    by
    \eqref{InfinityFixedPointsViaBaseChange}
  }
  \\
  \;\simeq\;
  \SliceMaps{\big}{\mathbf{B}\mathcal{G}}
    { \HomotopyQuotient{X}{\mathcal{G}} }
    { \HomotopyQuotient{Y}{\mathcal{G}} }
  &
  \proofstep{
    by
    \eqref{SliceMappingStackIsRightBaseChangeOfMappingStackInSlice}.
  }
  \end{array}
\end{equation}
\end{example}

In higher generalization of Exp. \ref{ForgettingGActionsAsPullbackAction} we have:
\begin{proposition}[Free- and cofree $\mathcal{G}$-action]
\label{UnderlyingObjectsOfGAction}
The forgetful functor from $\mathcal{G}$-actions
(Def. \ref{ActionObjectsInAnInfinityTopos})
to underlying objects
has a left adjoint and a right adjoint of the following form:
\vspace{-2mm}
$$
  \begin{tikzcd}[row sep=0pt]
    \Actions{\mathcal{G}}(\Topos)
    \ar[
      rr,
      "\underlying"
    ]
    &&
    \Topos
    \\
   \scalebox{0.8}{$ \mathcal{G}\acts \, P $}
    \ar[rr, |->]
    &&
   \scalebox{0.8}{$ P $}
  \end{tikzcd}
  {\phantom{AAAA}}
    \begin{tikzcd}[column sep=huge]
      \Topos
      \ar[
        rr,
        shift left=12pt,
        "{\mathcal{G} \acts \, \left( (-) \times \mathcal{G} \right)}"
      ]
      \ar[
        rr,
        phantom,
        shift left=9pt,
        "\scalebox{.5}{$\bot$}"
      ]
      \ar[
        rr,
        phantom,
        shift right=8pt,
        "\scalebox{.5}{$\bot$}"
      ]
      \ar[
        rr,
        shift right=12pt,
        "{
          \mathcal{G} \acts \, \Maps{}{\mathcal{G}}{-}
        }"{below}
      ]
      &&
      \Actions{\mathcal{G}}(\Topos)\;.
      \ar[
        ll,
        "\underlying"{description}
      ]
    \end{tikzcd}
$$
\end{proposition}
\begin{proof}
Observe that, under Prop. \ref{GroupsActionsAndFiberBundles},
the underlying object functor is
equivalently the pullback
along $\ast \xrightarrow{\mathrm{pt}_{\mathbf{B}\mathcal{G}}} \mathbf{B}\mathcal{G}$:
 \vspace{-2mm}
\begin{equation}
  \label{UnderlyingObjectOfActionViaBaseChange}
  P
  \;\simeq\;
  \underlying (G \acts \, P)
  \;\simeq\;
  \mathrm{fib}
  \scalebox{0.7}{$
  \left(\!\!\!\!
    \begin{array}{c}
      P \!\sslash\! \mathcal{G}
      \\
      \downarrow
      \\
      \mathbf{B}G
    \end{array}
  \!\!\!\! \right)
  $}
  \;\simeq\;
  (
    \mathrm{pt}_{\mathbf{B}\mathcal{G}}
  )^\ast
   \scalebox{0.7}{$
  \left(\!\!\!\!
    \begin{array}{c}
      P \!\sslash\! \mathcal{G}
      \\
      \downarrow
      \\
      \mathbf{B}\mathcal{G}
    \end{array}
  \!\!\!\! \right)
  $}.
\end{equation}

\noindent Therefore, the existence of the adjoints is a special case
of the base change adjoint triple (Prop. \ref{BaseChange}).
By the discussion there, the left adjoint, seen on the slice topos, is given
by postcomposition, which implies
--
by the pasting law \eqref{HomotopyPastingLaw}
and using the looping equivalence \eqref{GroupsActionsAndFiberBundles} --
that its underlying object is $(-) \times G$:
 \vspace{-2mm}
$$
  \begin{tikzcd}[row sep=small, column sep=large]
    X \times G
    \ar[r]
    \ar[d]
    \ar[
      dr,
      phantom,
      "\mbox{\tiny\rm(pb)}"
    ]
    \ar[
      dd,
      bend right=20,
      rounded corners,
      to path={
        -- ([xshift=-5pt]\tikztostart.west)
        -- node[left]{
            \scalebox{.7}{$
             {(\mathrm{pt}_{\mathbf{B}\mathcal{G}})^\ast
              (\mathrm{pt}_{\mathbf{B}\mathcal{G}})_!(X)}
             $}
           }
           ([xshift=-16pt]\tikztotarget.west)
        -- (\tikztotarget.west)
      }
    ]
    &
    X
    \ar[d]
    \ar[
      dd,
      rounded corners,
      to path={
        -- ([xshift=10pt]\tikztostart.east)
        -- node[right, yshift=1pt]
           { \scalebox{.7}{$
               (\mathrm{pt}_{\mathbf{B}\mathcal{G}})_!(X)
             $}
             }
               ([xshift=6pt]\tikztotarget.east)
        -- (\tikztotarget.east)
        }
    ]
    \\
    G
    \ar[r]
    \ar[d]
    \ar[
      dr,
      phantom,
      "\mbox{\tiny\rm(pb)}"
    ]
    &
    \ast
    \ar[d]
    \\
    \ast
    \ar[
      r,
      "\mathrm{pt}_{\mathbf{B}\mathcal{G}}"{below}
    ]
    &
    \mathbf{B}\mathcal{G}
  \end{tikzcd}
$$

\vspace{-2mm}
\noindent
Now observing the induced adjoint pair
 \vspace{-3mm}
$$
  \begin{tikzcd}[column sep=huge]
    \Topos
    \ar[
      rr,
      shift right=5pt,
      "{
        ( \mathrm{pt}_{\mathbf{B}\mathcal{G}})^\ast
          \circ
        ( \mathrm{pt}_{\mathbf{B}\mathcal{G}})_\ast
      }"{below}
    ]
    \ar[
      rr,
      phantom,
      "\scalebox{.7}{$\bot$}"
    ]
    &{\phantom{AAAA}}&
    \Topos
    \ar[
      ll,
      shift right=5pt,
      "{
        ( \mathrm{pt}_{\mathbf{B}\mathcal{G}})^\ast
          \circ
        ( \mathrm{pt}_{\mathbf{B}\mathcal{G}})_!
        \,\simeq\,
        (-) \times \mathcal{G}
      }"{above}
    ]
  \end{tikzcd}
$$

 \vspace{-2mm}
\noindent
implies the form of the right adjoint by the defining adjunction
\eqref{InternalHomAdjunction} and the essential uniqueness of adjoints.
\end{proof}

\medskip
\noindent
{\bf Some properties of $\infty$-group actions.} We establish some properties of $\infty$-actions in $\infty$-toposes that we will need later on.

\begin{lemma}[Homotopy quotient preserves connectedness]
  \label{HomotopyQuotientPreservesConnectedness}
  If
  $\mathcal{G} \acts \,  X \,\in \, \Actions{\mathcal{G}}(\Topos)$
  is an action on a connected object
  $X \in \Topos_{\geq 1}$, then also the homotopy quotient
  $X \!\sslash\! \mathcal{G}$ is connected.
\end{lemma}
\begin{proof}
  By \cite[Prop. 6.5.1.12]{Lurie09HTT}, an object is connected
  precisely if its image under 0-truncation
  \vspace{-2mm}
  $$
    \begin{tikzcd}
      \Topos_{\leq 0}
      \ar[
        rr,
        hook,
        shift right=5pt
      ]
      \ar[
        rr,
        phantom,
        "\scalebox{.7}{$\bot$}"
      ]
      &&
      \Topos
      \ar[
        ll,
        shift right=5pt,
        "\tau_0"{above}
      ]
    \end{tikzcd}
  $$

  \vspace{-1mm}
  \noindent
  is equivalent to the
  terminal object. But truncation is a left adjoint and
  hence sends homotopy quotients to homotopy quotients,
  $
    \tau_0
    \big(
      X \!\sslash\! \mathcal{G}
    \big)
    \;\simeq\;
    \tau_0(X) / \tau_0(\mathcal{G})
    \,.
  $
  But homotopy quotient in the subcategory $H_{\leq 0}$ of 0-truncated
  objects are plain quotients, which send the point to the point.
\end{proof}

In higher generalization of Ex. \ref{RestrictedActions}, we have:
\begin{proposition}[Restricted and left induced $\infty$-actions]
  \label{RestrictedAndLeftInducedInfinityActions}
  For $G \,\in\, \Groups(\Sets)$
  and $K \xhookrightarrow{\;i\;} G$ a subgroup inclusion,
  the induced left base change (via Prop. \ref{HomotopyQuotientsAndPrincipaInfinityBundles}
  and Prop. \ref{BaseChange})
  \vspace{-2mm}
$$
  \begin{tikzcd}
    \Actions{G}(\Topos)
    \,\simeq\,
    \SliceTopos{\mathbf{B} G}
    \ar[r, " (B i)_! "]
    &
    \SliceTopos{\mathbf{B} G'}
    \,\simeq\,
    \Actions{G'}(\Topos)
  \end{tikzcd}
$$

\vspace{-2mm}
\noindent satisfies
\begin{equation}
  \label{NatureOfLeftInducedInfinityAction}
  (B i)_!
  (
    K \acts X
  )
  \;\simeq\;
  G \acts \;
  (
    X \times_K G
  )
  \;\;\;\;\;
  \mbox{and}
  \;\;\;\;\;
  \HomotopyQuotient
    { X }
    { K }
  \;\simeq\;
  \HomotopyQuotient
    { (X \times_K G) }
    { G  }
  \,.
\end{equation}
\end{proposition}
\begin{proof}
  Unwinding the definitions, and using Prop. \ref{GroupsActionsAndFiberBundles}
  and the pasting law \eqref{HomotopyPastingLaw}, this follows from the following
  pasting diagram
$$
  \begin{tikzcd}
    X \times_K G
    \ar[r]
    \ar[d]
    \ar[
      dr,
      phantom,
      "\mbox{\tiny\rm(pb)}"
    ]
    &
    \HomotopyQuotient
      { X }
      { K }
    \ar[d]
    \ar[r, phantom, "\simeq"]
    &[-14pt]
    \HomotopyQuotient
      { (X \times_K G) }
      { G }
    \\
    G/K
    \ar[r]
    \ar[d]
    \ar[
      dr,
      phantom,
      "\mbox{\tiny\rm(pb)}"
    ]
    &
    \mathbf{B}K
    \ar[d, "{ \mathbf{B}i }"]
    \\
    \ast
    \ar[r]
    &
    \mathbf{B}G
    \mathrlap{\,.}
  \end{tikzcd}
$$

\vspace{-6mm}
\end{proof}

\begin{lemma}[Presentation of homotopy fixed point spaces of equivariant moduli stacks]
  \label{PresentationOfHomotopyFixedPointSpacesOfEquivariantClassifyingStacks}
  For
  $U \in \CartesianSpaces$,
  $G \in \Groups(\SimplicialSets)
     \xrightarrow{\;}
   \Groups(\SmoothInfinityGroupoids)
  $, and
  $(\Gamma, \rho) \in
   \Actions{G}\left(\Groups(\SimplicialPresheaves)\right)
     \longrightarrow
   \Actions{G}\left(\Groups(\SmoothInfinityGroupoids)\right)
  $,
  we have a natural equivalence of $\InfinityGroupoids$
   \vspace{-2mm}
  $$
    \left(
      \SmoothInfinityGroupoids
    \right)_{\!/\mathbf{B}G}
    \left(
      U \times \mathbf{B}G,
      \,
      (
        \mathbf{B}\Gamma
      )\sslash G
    \right)
    \;\simeq\;
    \SimplicialPresheaves(\CartesianSpaces)_{/\overline{W}G}
    \big(
       U \times \overline{W}G,
       \,
       \left(
       \overline{W}\Gamma
         \times
       W G
       \right)/G
    \big)
    \,.
  $$
\end{lemma}
\begin{proof}
  By Lemma \ref{InfinityGroupsPresentedByPresheavesOfSimplicialGroups},
  we have that $\overline{W}(-)$ represents $\mathbf{B}(-)$
  and that $\left((-) \times W G\right)/G$ represents $(-) \sslash G$.
  It hence remains to see that the simplicial hom-set on the right
  has the correct homotopy type of the derived hom-set.
  For this it is sufficient, by Lemma \ref{HomInfinityGroupoidFromCofibrantDomainAndFibrantCodomain},
  to show that:

\noindent  {\bf (i)} $U \times \overline{W}G$ is projectively cofibrant,

\noindent   {\bf (ii)} $\left(\overline{W}\mathcal{G} \times W G\right) / G$ is projectively fibrant.

  \noindent
  The first statement is Example \ref{ExamplesOfProjectiveleCofibrantSimplicialPresheaves}.
  The second statement follows since
  $\overline{W}\mathcal{G}$ is projectively fibrant as a simplicial presheaf (Prop. \ref{BasicPropertiesOfStandardSimplicialPrincipalComplex}),
  hence is projectively fibrant as a $G$-action \eqref{BorelModelStructure},
  and since $\left((-) \times W G\right)/G$ is a right Quillen functor
  (Prop. \ref{QuillenEquivalenceBetweenBorelModelStructureAndSliceOverClassifyingComplex})
  and thus preserves fibrancy.
\end{proof}

\medskip

\noindent
{\bf Base change of $\infty$-actions along discrete group covers.}

\medskip
\noindent We discuss (Prop. \ref{PullbackOfActionsAlongSurjectiveGroupHomomorphismsIsFullyFaithful} below)
the base change
of $\infty$-actions (Def. \ref{ActionObjectsInAnInfinityTopos})
along coverings of discrete groups, hence
along deloopings of surjective homomorphisms
  \begin{equation}
    \label{ADiscreteGroupEpimorphismInAnInfinityTopos}
    \hspace{-6mm}
    \begin{tikzcd}
      \widehat{G}
      \ar[r, ->>, "p"]
      &
      G
      \;\;\;\;
      \in
      \;
      \Groups(\Sets)
      \xrightarrow{ \Groups(\LocallyConstant) }
      \Groups(\ModalTopos)
      \,,
      \;\;\;
      \mbox{\small delooping to}
      \;\;\;
      B\widehat{G}
      \ar[r,  "B p"]
      &
      B G
      \;\;\;\;
      \in
      \;
      \OneGroupoids
      \xrightarrow{ \LocallyConstant \mbox{\tiny \cref{TerminalGeometricMorphism}} }
      \ModalTopos{1}
      \,.
    \end{tikzcd}
  \end{equation}
The following Prop. \ref{PullbackOfActionsAlongSurjectiveGroupHomomorphismsIsFullyFaithful}
is one key aspect of ``globalizing'' $G$-equivariant homotopy theory
(see also Rem. \ref{GloballyEquivariantNatureOfEquivariantPrincipalBundles}).

\medskip

\begin{lemma}[Terminal epimorphisms of delooping 1-groupoids]
  \label{TerminalEpimorphismsOfOneGroupoids}
  For all $G \,\in\, \Groups(\Sets)$ the terminal map out of its delooping groupoid
  $$
    B G \xrightarrow{ \; } \ast
    \;\;\;
    \in
    \;
    \OneGroupoids
  $$
  is an {\it epimorphism} in $\OneGroupoids$, in that for all
  all $\mathcal{X} \,\in\, \OneGroupoids \xhookrightarrow{\;} \InfinityGroupoids$ the
  induced map
  $$
    \Maps{}
      { B G \to \ast }
      { \mathcal{X} }
    \;\;
    :
    \;\;
    \mathcal{X}
      \xrightarrow{\;\;}
    \Maps{}
      { B G }
      { \mathcal{X} }
  $$
  is a monomorphism \eqref{InfinityMonomorphism},
  namely a fully faithful functor, i.e. the full inclusion
  of a connected component.
\end{lemma}
\begin{remark}[Epimorphisms of $n$-groupoids]
  \label{EpimorphismsOfNGroupoids}
  The proof of Lem. \ref{TerminalEpimorphismsOfOneGroupoids}
  is immediate by inspection, due to the fact that a natural transformation
  between functors out of
  $B G \simeq  \Localization{\WeakHomotopyEquivalences} (G \rightrightarrows \ast)$
  has a single component corresponding to the point $\ast \to B G$.

  The subtlety to notice here is just that this situation crucially relies
  on the ambient category being $\OneGroupoids$, as it fails in
  $\TwoGroupoids$ and thus in
  $\InfinityGroupoids$: If $\mathcal{X}$ is a 2-groupoid (or higher)
  then a (pseudo-)natural transformation between maps $B N \to \mathcal{X}$
  has, in general, not just a component on the point, but also components
  over each element of $G$, which are not faithfully reflected on $\ast \to B G$.
  Indeed, for $B G \to \ast$ to be an epimorphism in $\InfinityGroupoids$
  the group $G$ must be perfect, in that its abelianization is trivial
  (\cite{Raptis19}\cite[Lem. 3]{Hoyois19}).
\end{remark}

\begin{lemma}[$\infty$-Action base change comonad along discrete group extension]
  \label{InfinityActionBaseChangeComonadAlongDiscreteGroupExtensions}
  If $\begin{tikzcd} \widehat{G} \ar[r, ->>, "p"] & G\end{tikzcd}$
  is a surjective homomorphisms in $\Groups(\Sets)$
  \eqref{ADiscreteGroupEpimorphismInAnInfinityTopos}, then
  the induced base change comonad (Prop. \ref{BaseChange},
  via Prop. \ref{GroupsActionsAndFiberBundles}) of $G$-$\infty$-actions in $\Topos$
  along $B p$
  \vspace{-2mm}
  \begin{equation}
    \label{LeftBaseChangeOfInfinityActionsAlongCoverOfEquivarianceGroup}
    \begin{tikzcd}
      \Actions{\widehat{G}}(\Topos)
      \,\simeq\,
      \SliceTopos{\mathbf{B}\widehat{G}}
      \ar[
        rr,
        shift left=7pt,
        "{
          (B p)_!
        }"
      ]
      \ar[
        rr,
        phantom,
        "\scalebox{.7}{$\bot$}"
      ]
      \ar[
        from=rr,
        shift left=7pt,
        "{
          (B p)^\ast
        }"
      ]
      &&
      \SliceTopos{\mathbf{B}G}
      \,\simeq\,
      \Actions{G}(\Topos)
    \end{tikzcd}
  \end{equation}

  \vspace{-2mm}
  \noindent
  is naturally equivalent to the cartesian product with the delooping groupoid
  $B N$ of the kernel $N \,\coloneqq\, \mathrm{ker}(p) \,\subset\, \widehat G $
  equipped with its conjugation action:
  \begin{equation}
    \label{LeftBaseChangeComonadOnInfinityAction}
    (B p)_!
    \,
    (B p)^\ast
    \left(
      G \acts (-)
    \right)
    \;\simeq\;
    G \acts
    \left(
      (-) \times B N
    \right)
    \,.
  \end{equation}
  In particular, on 0-truncated objects the operation is the identity, in that
  for $G \acts W \,\in\, \Actions{G}(\ModalTopos{0})$ we have:
  \begin{equation}
    \label{ZeroTruncationOfPullPushOfActionAlongDeloopedSurjection}
    \Truncation{0}
    \,
    (B p)_!
    \,
    (B p)^\ast
    (
      G\acts \, W
    )
    \;\simeq\;
    G \acts \, W\;.
  \end{equation}
\end{lemma}
\begin{proof}
Consider the following diagram in $\Topos$:
$$
  \begin{tikzcd}
    W \times B N
    \ar[r]
    \ar[d]
    \ar[
      dr,
      phantom,
      "{
        \mbox{
          \tiny
          \rm
          (pb)
        }
      }"
    ]
    &
    \HomotopyQuotient
      {\left( (B p)^\ast W \right)}
      {\widehat G}
    \ar[d]
    \ar[r]
    \ar[
      dr,
      phantom,
      "{
        \mbox{
          \tiny
          \rm
          (pb)
        }
      }"
    ]
    &
    \HomotopyQuotient
      { W }
      { G }
    \ar[d]
    \\
    B N
    \ar[d]
    \ar[r]
    \ar[
      dr,
      phantom,
      "{
        \mbox{
          \tiny
          \rm
          (pb)
        }
      }"
    ]
    &
    B \widehat{G}
    \ar[d, "B p"{left}]
    \ar[r, "B p"{description}]
    &
    B G
    \ar[dl , -, shift left=1pt]
    \ar[dl , -, shift right=1pt]
    \\
    \ast
    \ar[r]
    &
    B G
    \mathrlap{\,.}
  \end{tikzcd}
$$
Here the total middle vertical morphism is
$(B p)_! \, (B p^\ast) \, (W)$, by Prop. \ref{BaseChange},
so that the total left
Cartesian rectangle computes, in its top left corner, its
underlying object. The bottom left rectangle is Cartesian by
Prop. \ref{InfinityGroupsPresentedByPresheavesOfSimplicialGroups}
with
Lemma \ref{EssentiallySurjectiveKanFibrationsOfSimplicialGroupsInducedKanFibrationsUnderBarW}.
Hence, by repeated use of the pasting law \eqref{HomotopyPastingLaw},
the top left object
is equivalently
the pullback exhibited by the top left square and hence that by the
total top rectangle. But by the factorization at the bottom,
this total top pullback is the Cartesian product of the
underlying object of $G \acts \, W$ with $B N$, as claimed.
\end{proof}

\begin{proposition}[Base change of 2-actions along surjective group homomorphisms is faithful]
  \label{PullbackOfActionsAlongSurjectiveGroupHomomorphismsIsFullyFaithful}
  If $p$ is a surjective homomorphisms of discrete groups
  \eqref{ADiscreteGroupEpimorphismInAnInfinityTopos}, then
  the induced base change (Prop. \ref{BaseChange}, via Prop. \ref{GroupsActionsAndFiberBundles})
  along $B p$ is

  \noindent
  {\bf (i)}
  fully faithful
  \eqref{FullyFaithfulInfinityFunctor}
  on 0-truncated objects
  \vspace{-2mm}
  \begin{equation}
    \label{BaseChangeEquivalenceOfOneActionsAlongSurjectiveHomomorphism}
      \begin{tikzcd} \widehat{G} \ar[r, ->>, "p"] & G\end{tikzcd}
      \hspace{.8cm}
      \vdash
      \hspace{.8cm}
      \begin{tikzcd}
      \Actions{\widehat G}(\ModalTopos{0})
      \,\simeq\,
      \left(
        \SliceTopos{B \widehat{G}}
      \right)_0
      \ar[
        from=rr,
        hook',
        "(B p)^\ast"{swap}
      ]
      &&
      \left(
        \SliceTopos{B G}
      \right)_0
      \,\simeq\,
      \Actions{G}(\ModalTopos{0})
      \,.
    \end{tikzcd}
  \end{equation}

  \vspace{-1mm}
  \noindent
  {\bf (ii)}
  faithful on 1-truncated objects, in the sense that $(B p)^\ast$ induces a
  monomorphism, hence (Ex. \ref{MonomorphismsOfInfinityGroupoids}) an
  injection of connected components of hom-$\infty$-groupoids:
  \vspace{-1.5mm}
  \begin{equation}
    \label{BaseChangeOfTwoActionsAlongSurjectiveHomomorphism}
    \left.
    \begin{array}{rl}
      &
      \begin{tikzcd} \widehat{G} \ar[r, ->>, "p"] & G \end{tikzcd}
      \\
      \mbox{and}
      &
      G\acts \, W, G\acts \, V \,\in\,
      \Actions{G}(\ModalTopos{1})
    \end{array}
    \right\}
    \hspace{.3cm}
    \vdash
    \hspace{.3cm}
    \begin{tikzcd}
    \Homs{\big}
      {
        (B p)^\ast
        (W)
      }
      {
        (B p)^\ast
        (V)
      }
    \ar[
      from=rr,
      hook',
      "{
        (B p)^\ast_{W,V}
      }"{swap}
    ]
    &&
    \Homs{}
      { W }
      { V }
    \,.
    \end{tikzcd}
  \end{equation}
\end{proposition}
\begin{proof}
For the first statement, consider the
following sequence of natural equivalences, for
0-truncated
$G \acts \, V, \, G \acts \, W \,\in\, \Actions{G}(\ModalTopos{0})$:
\vspace{-3.5mm}
$$
  \def\arraystretch{1.6}
  \begin{array}{lll}
    \Homs{\big}
      {(B p)^\ast W}
      {(B p)^\ast V}
    &
    \;\simeq\;
    \Homs{\big}
      {(B p)_! (B p)^\ast W}
      {V}
    &
    \proofstep{ by \eqref{AdjunctionAndHomEquivalence} }
    \\
    & \;\simeq\;
    \Homs{\big}
      {\Truncation{0} (B p)_! (B p)^\ast W}
      {V}
    &
    \proofstep{ Prop. \ref{nTruncation} with \eqref{AdjunctionAndHomEquivalence} }
    \\
    & \;\simeq\;
    \Homs{}
      {W}
      {V}
    &
    \proofstep{ by \eqref{ZeroTruncationOfPullPushOfActionAlongDeloopedSurjection}. }
  \end{array}
$$

\vspace{-1.5mm}
\noindent
By the $\infty$-Yoneda lemma (Prop. \ref{InfinityYonedaLemma})
the composite natural
equivalence is a left inverse to
$(B p)^\ast_{(-,-)}$
and hence exhibits fully faithfulness on 0-truncated objects.

The second statement follows by analogous reasoning.
For readability, we display the proof only for the special case when
$\widehat{G} \xrightarrow{p} G$ is a {\it central} extension, so that
the $G$-action on $B N$ is trivial, hence the case when
\vspace{-1mm}
\begin{equation}
  \label{AssuminTrivialGActionOnBN}
  G \acts \, B N
  \;\simeq\;
  (B p_G)^\ast (B N)
  \,,
  \;\;\;\;
  \mbox{
    for $p_G : G \xrightarrow{\;} \ast$.
  }
\end{equation}

\vspace{-1.5mm}
\noindent Then we have the following natural equivalences:
\vspace{-1.5mm}
$$
  \def\arraystretch{1.47}
  \begin{array}{lll}
  &
        \Actions{\widehat{G}}(\ModalTopos{1})
    \big(
      {(B p)^\ast (G \acts \, W)}
      ,\,
      {(B p)^\ast (G \acts \, V)}
    \big)
\\
    &
    \;\simeq\;
    \Actions{G}(\ModalTopos{1})
    \big(
      {(B p)_! (B p)^\ast (G \acts \, W)}
      ,\,
      {G \acts \, V}
    \big)
    &
    \proofstep{ by \eqref{AdjunctionAndHomEquivalence} }
    \\
    & \;\simeq\;
    \Actions{G}(\ModalTopos{1})
    \left(
      { G \acts \, W \times B N }
      ,\,
      { G \acts \, V }
    \right)
    &
    \proofstep{
      by \eqref{LeftBaseChangeComonadOnInfinityAction}
      in Lem. \ref{InfinityActionBaseChangeComonadAlongDiscreteGroupExtensions}
    }
    \\
    & \;\simeq\;
    \Actions{G}(\ModalTopos{1})
    \left(
      G \acts \, \ast
      ,\,
     { \rm Maps} \left(
        { G \acts \, W \times B N },
        { G \acts \, V }
        \right)
        \right)
    &
    \proofstep{ by \eqref{HomIsGlobalPointsOfMappingStacks} }
    \\
    & \;\simeq\;
    \Actions{G}(\ModalTopos{1})
    \big(
      G \acts \, \ast
      ,\,
     {\rm Maps} \left(
        { G \acts \, B N },
        {
          \Maps{}
            { G \acts \, W }
            { G \acts \, V }
        }
        \right)\!\!
    \big)
    &
    \proofstep{
      by Lem. \ref{InternalHomAdjointness}
    }
    \\
    & \;\simeq\;
    \Actions{G}(\ModalTopos{1})
    \big(
      { G \acts \, B N }
      ,\,
      \Maps{}
        { G \acts \, W }
        { G \acts \, V }
    \big)
    &
    \proofstep{ by \eqref{InternalHomAdjunction} with \eqref{AdjunctionAndHomEquivalence} }
    \\
    & \;\simeq\;
    \Actions{G}(\ModalTopos{1})
    \big(
      { (B p_{G})^\ast B N }
      ,\,
      \Maps{}
        { G \acts \, W }
        { G \acts \, V }
    \big)
    &
    \proofstep{ by \eqref{AssuminTrivialGActionOnBN} }
    \\
    & \;\simeq\;
    \ModalTopos{1}
    \big(
      { B N }
      ,\,
      (B p_G)_\ast
      \Maps{}
        { G \acts \, W }
        { G \acts \, V }
    \big)
    &
    \proofstep{ by Prop. \ref{UnderlyingObjectsOfGAction} }
    \\
    & \;\simeq\;
    \ModalTopos{1}
    \big(
      { \LocallyConstant(B N) }
      ,\,
      (B p_G)_\ast
      \Maps{}
        { G \acts \, W }
        { G \acts \, V }
    \big)
    &
    \proofstep{ by \eqref{ADiscreteGroupEpimorphismInAnInfinityTopos} }
    \\
    & \;\simeq\;
    \OneGroupoids
    \big(
      { B N }
      ,\,
      \GlobalSections
      \,
      (B p_G)_\ast
      \,
      \Maps{}
        { G \acts \, W }
        { G \acts \, V }
    \big)
    &
    \proofstep{ by \eqref{TerminalGeometricMorphism} with \eqref{AdjunctionAndHomEquivalence} }
    \\
    & \;\simeq\;
    \OneGroupoids
    \big(
      { B N }
      ,\,
      \Actions{G}(\ModalTopos{1})
      (
        { G \acts \, W }
        ,\,
        { G \acts \, V }
      )
    \big)
    &
    \proofstep{
      by
      \eqref{GlobalSectionsInInfinityToposAsMapsOutOfTerminalObject}
      with
      \eqref{HomIsGlobalPointsOfMappingStacks}
    }
    \\
    &
    \;=\;
    {\rm Maps} \left(
      { B N },
      {
        \Homs{}
          { W }
          { V }
      }
      \right)
    \,.
  \end{array}
$$

\vspace{-.5mm}
\noindent With this, the claim follows by
Lem. \ref{TerminalEpimorphismsOfOneGroupoids}.
\end{proof}

\medskip

\noindent
{\bf Higher shear maps.}
In preparation of the discussion of
principal $\infty$-bundles (Def. \ref{PrincipalInfinityBundles} below),
we discuss the higher analogs of the shear map \eqref{PrincipalityConditionAsShearMapBeingAnIsomorphism}:

\begin{definition}[Higher shear maps]
  \label{HigherShearMaps}
  Given $\mathcal{G} \acts \, X \,\in\, \Actions{G}(\Topos)$
  (Def. \ref{ActionObjectsInAnInfinityTopos}),

  \noindent {\bf (i)} we say that its {\it shear map} is the universal morphism $\mathrm{shear}_1$
  in the following pasting diagram, whose front and rear faces are
  Cartesian, by construction:
  \vspace{-2mm}
  \begin{equation}
    \label{ShearMapOfAnInfinityAction}
    \begin{tikzcd}[row sep=small]
      \mathcal{G} \times X
      \ar[rr, "\rho"]
      \ar[dd, "\mathrm{pr}_1"{left}]
      \ar[dr, dashed, "\mathrm{shear}_1"{sloped, below}]
      &&
      X
      \ar[dr, -, shift left=1pt]
      \ar[dr, -, shift right=1pt]
      \ar[dd, ->>]
      \\
      &
      X \times X
      \ar[rr, crossing over, "\mathrm{pr}_1"{pos=.17}]
      \ar[ddrr, phantom, "{\mbox{\tiny\rm (pb)}}"{pos=.1}]
      &&
      X
      \ar[dd]
      \\
      X
      \ar[dr, -, shift left=1pt]
      \ar[dr, -, shift right=1pt]
      \ar[rr, ->>]
      &&
      X \!\sslash\! \mathcal{G}
      \ar[dr]
      \\
      &
      X \ar[rr]
      \ar[from=uu, crossing over, "\mathrm{pr}_2"{left, pos=.25}]
      &&
      \ast
    \end{tikzcd}
  \end{equation}

\vspace{-2mm}
\noindent
{\bf (ii)}   More generally, we say that the
  {\it higher shear maps}
  \vspace{-2mm}
  \begin{equation}
    \mathcal{G}^{\times_n}
      \times
    X
    \xrightarrow{\;\mathrm{shear}_{n}\;}
    X^{\times n+1}
    \,,
    \;\;\;\;\;
    n \in \mathbb{N}
  \end{equation}

  \vspace{-2mm}
  \noindent
  are the universal morphisms
  given by functorially
  forming the {\v C}ech nerves
  of $\mathrm{id}_X$ canonically mapping into
  $\HomotopyQuotient{X}{G} \xrightarrow{\;}\ast$
  \eqref{EquivalenceBetweenEffectiveEpimorphismsAndGroupoidObjects}:
  \vspace{-3mm}
  \begin{equation}
    \label{DefiningDiagramForHigherShearMaps}
    \begin{tikzcd}[row sep=small, column sep=large]
      {}
      \ar[d, -, shift left=15pt, dotted]
      \ar[from=d, -, shift left=10pt, dotted]
      \ar[d, -, shift left=5pt, dotted]
      \ar[from=d, -, dotted]
      \ar[d, -, shift right=5pt, dotted]
      \ar[from=d, -, shift right=10pt, dotted]
      \ar[d, -, shift right=15pt, dotted]
      &
      {}
      \ar[d, -, shift left=15pt, dotted]
      \ar[from=d, -, shift left=10pt, dotted]
      \ar[d, -, shift left=5pt, dotted]
      \ar[from=d, -, dotted]
      \ar[d, -, shift right=5pt, dotted]
      \ar[from=d, -, shift right=10pt, dotted]
      \ar[d, -, shift right=15pt, dotted]
      \\
      \mathcal{G} \times \mathcal{G}
        \times
      X
      \ar[d, shift left=10pt]
      \ar[from=d, shift left=5pt]
      \ar[d]
      \ar[from=d, shift right=5pt]
      \ar[d, shift right=10pt]
      \ar[r, "\mathrm{shear}_2"]
      &
      X \times X \times X
      \ar[d, shift left=10pt]
      \ar[from=d, shift left=5pt]
      \ar[d]
      \ar[from=d, shift right=5pt]
      \ar[d, shift right=10pt]
      \\
      \mathcal{G}
        \times
      X
      \ar[d, shift left=5pt]
      \ar[from=d]
      \ar[d, shift right=5pt]
      \ar[r, "\mathrm{shear}_1"]
      &
      X \times X
      \ar[d, shift left=5pt]
      \ar[from=d]
      \ar[d, shift right=5pt]
      \\
      X
      \ar[r, -, shift right=1pt]
      \ar[r, -, shift left=1pt, "\mathrm{shear}_0"{above}]
      \ar[d, ->>]
      &
      X
      \ar[d]
      \\
      \HomotopyQuotient{X}{\mathcal{G}}
      \ar[r]
      &
      \ast
    \end{tikzcd}
  \end{equation}
\end{definition}
\begin{example}[Recovering the traditional notion of shear map]
  Let $\mathbf{H} \,=\, (\SmoothInfinityGroupoids)_{/X}$
  be the slice $\infty$-topos over some
  $X \,\in\, \DTopologicalSpaces \xhookrightarrow{\;} \SmoothInfinityGroupoids$,
  and let
  $\Gamma \,\in\, \Groups(\DTopologicalSpaces) \xhookrightarrow{ (-) \times X }
  \Groups(\mathbf{H})$
  and
  $\Gamma \acts \, P \,\in\,
  \Actions{G}(\DTopologicalSpaces_{/X})$ $
  \xhookrightarrow{\;} \Actions{G}(\mathbf{H})$
  be a fiberwise topological group action on a topological bundle $P$ over $X$.

  Then its shear map
  \eqref{ShearMapOfAnInfinityAction}
  according to Def. \ref{HigherShearMaps}
  is the universal dashed morphism in the following diagram:
  \vspace{-2mm}
  $$
    \begin{tikzcd}[row sep=small]
      \Gamma \times P
      \ar[rr, "\rho"]
      \ar[dd, "\mathrm{pr}_1"{left}]
      \ar[dr, dashed, "\mathrm{shear}_1"{sloped, below}]
      &&
      P
      \ar[dr, -, shift left=1pt]
      \ar[dr, -, shift right=1pt]
      \ar[dd, ->>]
      \\
      &
      P \times_X P
      \ar[rr, crossing over, "\mathrm{pr}_1"{pos=.17}]
      \ar[ddrr, phantom, "{\mbox{\tiny\rm (pb)}}"{pos=.1}]
      &&
      P
      \ar[dd]
      \\
      P
      \ar[dr, -, shift left=1pt]
      \ar[dr, -, shift right=1pt]
      \ar[rr, ->>]
      &&
      P \!\sslash_{{}_{\!X}}\! \Gamma
      \ar[dr]
      \\
      &
      P \ar[rr]
      \ar[from=uu, crossing over, "\mathrm{pr}_2"{left, pos=.25}]
      &&
      X
    \end{tikzcd}
  $$

  \vspace{-1mm}
\noindent
  This manifestly coincides with the traditional shear map \eqref{PrincipalityConditionAsShearMapBeingAnIsomorphism}.
\end{example}

\medskip

\noindent
{\bf Free $\infty$-actions.}

\begin{definition}[Free, transitive and regular $\infty$-actions]
  \label{FreeInfinityAction}
  We say that $\mathcal{G} \acts \, X \,\in\, \Actions{\mathcal{G}}(\Topos)$
  (Def. \ref{ActionObjectsInAnInfinityTopos})
  is

  \vspace{-2mm}
  \begin{itemize}
  \setlength\itemsep{-6pt}

  \item
  {\it free} if its shear map
  \eqref{ShearMapOfAnInfinityAction}
  is $(-1)$-truncated:
    \vspace{-2mm}
  \begin{equation}
    \label{FreeActionWitnessedByHigherShearMaps}
    \mbox{
      $\mathcal{G} \acts \, X$
      is free
    }
    \;\;\;\;
    \Leftrightarrow
    \;\;\;\;
    \begin{tikzcd}[column sep=large]
      \mathcal{G}
      \times
      X
      \ar[r, hook, "\mathrm{shear}_{1}"{above}]
      &
      X \times X
      \,;
    \end{tikzcd}
  \end{equation}

  \item
  {\it transitive} if its shear map
  is $(-1)$-connected:
  \vspace{-2mm}
  \begin{equation}
    \label{TransitiveActionWitnessedByHigherShearMaps}
    \mbox{
      $\mathcal{G} \acts \, X$
      is transitive
    }
    \;\;\;\;
    \Leftrightarrow
    \;\;\;\;
    \begin{tikzcd}[column sep=large]
      \mathcal{G}
      \times
      X
      \ar[r, ->>,  "\mathrm{shear}_{1}"{above}]
      &
      X \times X
      \,;
    \end{tikzcd}
  \end{equation}

  \item
  {\it regular} if its shear map
  is an equivalence:
  \vspace{-2mm}
  \begin{equation}
    \label{RegularActionWitnessedByHigherShearMaps}
    \mbox{
      $\mathcal{G} \acts \, X$
      is regular
    }
    \;\;\;\;
    \Leftrightarrow
    \;\;\;\;
    \begin{tikzcd}[column sep=large]
      \mathcal{G}
      \times
      X
      \ar[r, "\sim"{below},  "\mathrm{shear}_{1}"{above}]
      &
      X \times X
      \,.
    \end{tikzcd}
  \end{equation}

\end{itemize}
\vspace{-.1cm}

\end{definition}

\begin{lemma}[Higher shear maps are $(-1)$-truncated if first shear map is]
\label{HigherShearMapsAreMinusOneTruncatedIfFirstShearMapIs}
  Given $\mathcal{G} \acts  \, X \,\in\, \Actions{G}(\Topos)$,
  if its shear map $\mathrm{shear}_1$ is $(-1)$-truncated
  or an equivalence,
  then so are all higher shear
  maps (Def. \ref{HigherShearMaps}), respectively:
  \vspace{-2mm}
  $$
    \begin{tikzcd}[column sep=large]
      \mathcal{G}
      \times
      X
      \ar[r, hook, "\mathrm{shear}_1"{above}]
      &
      \mathcal{X} \times \mathcal{X}
    \end{tikzcd}
    \;\;\;\;\;
    \Rightarrow
    \;\;\;\;\;
    \underset{n \in \mathbb{N}}{\forall}
    \begin{tikzcd}[column sep=large]
      \mathcal{G}
      \times
      X
      \ar[r, hook, "\mathrm{shear}_n"{above}]
      &
      \mathcal{X} \times \mathcal{X}
      \mathrlap{,}
    \end{tikzcd}
  $$
  \vspace{-6mm}
  $$
    \begin{tikzcd}[column sep=large]
      \mathcal{G}
      \times
      X
      \ar[r, "\sim"{below}, "\mathrm{shear}_1"{above}]
      &
      \mathcal{X} \times \mathcal{X}
    \end{tikzcd}
    \;\;\;\;\;
    \Rightarrow
    \;\;\;\;\;
    \underset{n \in \mathbb{N}}{\forall}
    \begin{tikzcd}[column sep=large]
      \mathcal{G}
      \times
      X
      \ar[r, "{\sim}"{below}, "\mathrm{shear}_n"{above}]
      &
      \mathcal{X} \times \mathcal{X}
      \mathrlap{.}
    \end{tikzcd}
  $$
\end{lemma}
\begin{proof}
Observe that for all effective epimorphisms $X \twoheadrightarrow \mathcal{X}$
we have, for all $n \in \mathbb{N}$, homotopy-cartesian squares as follows:
  \vspace{-2mm}
$$
  \begin{tikzcd}[row sep=small, column sep=25pt]
    &&
    X^{\times^{n+2}_{\mathcal{X}}}
    \ar[dl]
    \ar[dd, phantom, "\mbox{\tiny \rm (pb)}"]
    \ar[dr]
    \\
    &
    X^{\times^{n+1}_{\mathcal{X}}}
    \ar[dl]
    \ar[dr]
    \ar[dd, phantom, "\mbox{\tiny \rm (pb)}"]
    &&
    X^{\times^2_{\mathcal{X}}}
    \ar[dl]
    \ar[dr]
    \ar[dd, phantom, "\mbox{\tiny \rm (pb)}"]
    \\
    X^{\times^n_{\mathcal{X}}}
    \ar[dr]
    &&
    X
    \ar[dl]
    \ar[dr]
    &&
    X
    \ar[dl]
    \\
    &
    \mathcal{X}
    &&
    \mathcal{X}
  \end{tikzcd}
$$

  \vspace{-1mm}
\noindent
Namely, the two bottom squares are cartesian by construction of
{\v C}ech nerves,
while the top square is cartesian by the groupoidal
Segal conditions \eqref{SegalConditionsForGroupoidObjectInInfinityTopos}
satisfied by groupoid objects.

Applied to the morphism of effective epimorphisms
\vspace{-2mm}
$$
  \begin{tikzcd}[row sep=small]
    X
    \ar[r, -, shift left=1pt]
    \ar[r, -, shift right=1pt]
    \ar[d, ->>]
    &
    X
    \ar[d]
    \\
    \HomotopyQuotient{X}{\mathcal{G}}
    \ar[r]
    &
    \ast
  \end{tikzcd}
$$

  \vspace{-2mm}
\noindent
this implies that the $n+1$st shear map is the following homotopy
pullback in the arrow category of $\Topos$:
  \vspace{-2mm}
$$
  \begin{tikzcd}[row sep=small]
    X \times G^{\times^{n+1}}
    \ar[rr]
    \ar[dr, dashed, "\mathrm{shear}_{n+1}"{sloped, below}]
    \ar[dd]
    &&
    X \times G
    \ar[dr, "\mathrm{shear}_1"{sloped}]
    \ar[dd, "\mathrm{pr}_1"{left, pos=.3}]
    \\
    &
    X^{\times^{n+2}}
    \ar[rr, crossing over]
    &&
    X \times X
    \ar[dd, "\mathrm{pr}_1"{left, pos=.3}]
    \\
    X \times G^{\times^n}
    \ar[rr]
    \ar[dr, "\mathrm{shear}_n"{sloped, below}]
    &&
    X
    \ar[dr, -, shift left=1pt]
    \ar[dr, -, shift right=1pt]
    \\
    &
    X^{\times^{n+1}}
    \ar[rr]
    \ar[from=uu, crossing over]
    &&
    X
  \end{tikzcd}
$$

  \vspace{-2mm}
\noindent
Since the classes of (-1)-truncated morphisms
and of equivalences
(include all identity morphisms)
and both are stable
\eqref{LimitOfMonomorphismsInArrowCategoryIsMono}
under $\infty$-limits in the arrow category
(by Prop. \ref{TheEffectiveEpiMonoFactorizationSystem}), the claim follows by induction.
\end{proof}

\begin{proposition}[Inhabited homotopy quotients of regular $\infty$-actions are terminal]
  \label{InhabitedHomotopyQuotientsOfRegularInfinityActionsAreTerminal}
  If $X \,\in\, \Topos$ is inhabited (Ntn. \ref{InhabitedObject})
  then an $\infty$-action $\mathcal{G} \acts \, X \,\in\, \Actions{\mathcal{G}}(\Topos)$
  (Def. \ref{ActionObjectsInAnInfinityTopos})
  is regular (Def. \ref{FreeInfinityAction})
  iff
  its homotopy quotient (Prop. \ref{HomotopyQuotientsAndPrincipaInfinityBundles})
  is terminal:
  $$
    (
      X \twoheadrightarrow \ast
    )
    \;\;\;\;\;\;\;\;\;
    \Rightarrow
    \;\;\;\;\;\;\;\;\;
    \big(
      \mbox{\rm $\mathcal{G} \acts \, X$ is regular}
      \;\;\;\;\;
      \Leftrightarrow
      \;\;\;\;\;
      \HomotopyQuotient
        { X }
        { \mathcal{G} }
      \;\simeq\;
      \ast
    \big)
    \,.
  $$
\end{proposition}
\begin{proof}
  This follows by inspection of the defining diagram
  \eqref{DefiningDiagramForHigherShearMaps} of the higher shear maps,
  noticing that,
  under the given assumption on $X$, not only its left but also its
  right bottom morphism is an effective epimorphism, so that
  the bottom rectangle is the image under forming homotopy quotients
  of the morphisms of simplicial objects. Therefore:
  If the action is regular in that
  $\mathrm{shear}_1$ is an equivalence \eqref{RegularActionWitnessedByHigherShearMaps},
  then all top horizontal morphisms in
  \eqref{DefiningDiagramForHigherShearMaps}
  are equivalences, by
  Lem. \ref{HigherShearMapsAreMinusOneTruncatedIfFirstShearMapIs},
  whence all top squares are Cartesian, so that the bottom morphism
  is the image of an equivalence under the $\infty$-colimit functor,
  and hence itself an equivalence.

  Conversely, if the homotopy quotient is terminal,
  then the bottom morphism is an equivalence, and hence so is the shear map,
  by its defining diagram \eqref{ShearMapOfAnInfinityAction}.
\end{proof}

\begin{proposition}[Homotopy quotient of free actions is ordinary quotient of 0-truncation]
  \label{HomotopyQuotientOfFreeActionsIsOrdinaryQuotientOfZeroTruncation}
  For $\Topos$ an $\infty$-topos,
  with
  $\mathcal{G} \,\in\, \Groups(\Topos)$
  if
  $\mathcal{G} \acts \, X \,\in\, \Actions{\mathcal{G}}(\Topos)$
  is free  (Def. \ref{FreeInfinityAction})
  and $X$ is inhabited (Ntn. \ref{InhabitedObject}),
  then the homotopy quotient \eqref{HomotopyQuotientAsHomotopyColimit}
  coincides with the
  ordinary quotient of the 0-truncation:
  \vspace{-2mm}
  $$
    \HomotopyQuotient{ X }{ \mathcal{G} }
    \;\simeq\;
    i_0
    \big(
      \tau_0(X) / \tau_0(\mathcal{G})
    \big)
    \;\coloneqq\;
    i_0
    \Big(
    \mathrm{coeq}
    \big(
      \!\!\!
      \begin{tikzcd}[column sep=large]
        \tau_0(\mathcal{G})
        \times
        \tau_0(X)
        \ar[r, shift left=3pt, "{\mathrm{pr_1}}"{above}]
        \ar[r, shift right=3pt, "{\mathrm{\rho}}"{below}]
        &
        \tau_0(X)
      \end{tikzcd}
      \!\!\!
    \big)
    \Big)
    \;\;\;
    \in
    \Topos_0
    \xhookrightarrow{\; i_0 \;}
    \Topos
    \,.
  $$
\end{proposition}
\begin{proof}
First we show that $\HomotopyQuotient{X}{\mathcal{G}}$ is 0-truncated.
For this it is sufficient to demonstrate that,
for all $U \,\in\, \Topos$
and
$\mathcal{K} \,\in\, \Groups(\InfinityGroupoids)
\xrightarrow{\Groups(\LocallyConstant)}
\Groups{\Topos}$,
every morphism
out of $U \times \mathbf{B}\mathcal{K}$ into it factors through the projection
onto $U$:
\vspace{-2mm}
$$
  \begin{tikzcd}
    U \times \mathbf{B}\mathcal{K}
    \ar[rr, "\forall"{above}]
    \ar[d, "\mathrm{pr}_1"{left}]
    &&
    \HomotopyQuotient{X}{\mathcal{G}}
    \ar[d]
    \\
    U
    \ar[rr]
    \ar[urr, dashed, "\exists"]
    &&
    \ast
    \,.
  \end{tikzcd}
$$

\vspace{-2mm}
\noindent
To see this, consider the following square diagram of augmented simplicial
objects, obtained by choosing atlases as shown in the bottom row and
then forming the four {\v C}ech groupoids
running vertically (Ex. \ref{CechGroupoidsInAnInfinityTopos}):
\vspace{-4mm}
\begin{equation}
  \label{DiagramOfCechNervesInProofThatHomotopyQuotientOfFreeActionIsZeroTruncated}
  \begin{tikzcd}[row sep=small]
    {}
    \ar[dd, -, dotted]
    &&
    {}
    \ar[dd, -, dotted]
    \\
    &
    {}
    &&
    {}
    \ar[dd, -, dotted]
    \\[-22pt]
    \mathcal{K}
    \times
    \mathcal{K}
    \times
    U
    \ar[dd, shift left=12pt, start anchor={[yshift=-4pt]}]
    \ar[from=dd, shift left=6pt]
    \ar[dd]
    \ar[from=dd, shift right=6pt]
    \ar[dd, shift right=12pt]
    \ar[rr]
    \ar[dr, ->>]
    &&
    \mathcal{G}
    \times
    \mathcal{G}
    \times
    X
    \ar[dd, shift left=12pt]
    \ar[from=dd, shift left=6pt]
    \ar[dd]
    \ar[from=dd, shift right=6pt]
    \ar[dd, shift right=12pt]
    \ar[dr, hook]
    \\
    &
    U
    \ar[rr, crossing over]
    \ar[ur, dashed]
    \ar[from=uu, -, dotted, crossing over]
    &&
    X \times X \times X
    \ar[dd, shift left=12pt]
    \ar[from=dd, shift left=6pt]
    \ar[dd]
    \ar[from=dd, shift right=6pt]
    \ar[dd, shift right=12pt]
    \\
    \mathcal{K}
    \times
    U
    \ar[rr]
    \ar[dr, ->>]
    \ar[dd, shift left=6pt]
    \ar[from=dd]
    \ar[dd, shift right=6pt]
    &&
    \mathcal{G}
    \times
    X
    \ar[dr, hook]
    \ar[dd, shift left=6pt]
    \ar[from=dd]
    \ar[dd, shift right=6pt]
    \\
    &
    U
    \ar[rr, crossing over]
    \ar[from=uu, shift left=12pt, crossing over, end anchor={[yshift=7pt]}]
    \ar[uu, shift left=6pt, crossing over, start anchor={[yshift=3pt]}]
    \ar[from=uu, crossing over]
    \ar[uu, shift right=6pt, crossing over, start anchor={[yshift=3pt]}]
    \ar[from=uu, shift right=12pt, crossing over, end anchor={[yshift=7pt]}]
    \ar[ur, dashed]
    &&
    X \times X
    \ar[dd, shift left=6pt]
    \ar[from=dd]
    \ar[dd, shift right=6pt]
    \\
    U
    \ar[rr]
    \ar[dd, ->>]
    \ar[dr, -, shift left=1pt]
    \ar[dr, -, shift right=1pt]
    &&
    X
    \ar[dd, ->>]
    \ar[dr, -, shift left=1pt]
    \ar[dr, -, shift right=1pt]
    \\
    &
    U
    \ar[rr, crossing over]
    \ar[from=uu, shift left=6pt, crossing over]
    \ar[uu, crossing over]
    \ar[from=uu, shift right=6pt, crossing over]
    \ar[ur, dashed]
    &&
    X
    \ar[dd, ->>]
    \\
    U \times \mathbf{B}\mathcal{K}
    \ar[rr]
    \ar[dr]
       &&
    \HomotopyQuotient{X}{\mathcal{G}}
    \ar[dr]
    \\
    &
    U
    \ar[rr]
    \ar[from=uu, ->>, crossing over]
    \ar[ur, dashed]
    &&
    \ast
  \end{tikzcd}
\end{equation}

\vspace{-2mm}
\noindent
Here in all of the upper horizontal squares

\vspace{-2mm}
\begin{itemize}

\vspace{-.2cm}
\item
the left morphism is $(-1)$-connected $\twoheadrightarrow$,
since every $\infty$-group $\mathcal{K} \,\in\,\Groups(\InfinityGroupoids)$
is inhabited (Ntn. \ref{InhabitedObject})
and $\LocallyConstant$, being  lex left adjoint,
preserves this property;

\vspace{-.2cm}
\item
the right morphism is $(-1)$-truncated $\hookrightarrow$, by Lem. \ref{HigherShearMapsAreMinusOneTruncatedIfFirstShearMapIs}.
\end{itemize}
\vspace{-.2cm}

\vspace{-2mm}
\noindent
Since $n$-connected/truncated morphisms in $\infty$-presheaf categories
such as $\InfinityPresheaves(\Delta, \Topos)$ are detected objectwise,
this means that the total square diagram of simplicial objects
is a $(-1)$-connected/truncated-lifting problem,
and since $\InfinityPresheaves(\Delta, \Topos)$ is again an $\infty$-topos
this lifting problem has an essentially unique solution \eqref{MinusOneConnectedTruncatedLiftingProblem},
hence a {\it compatible} set of dashed lifts as shown above.
Therefore, forming again the $\infty$-colimit over the
vertical simplicial diagrams recovers the bottom square
(Prop. \ref{EquivalentPerspectivesOnGroupoidObjectsInAnInfinityTopos})
but now itself equipped with a lift. This is the required lift
which proves that $\HomotopyQuotient{X}{\mathcal{G}}$ is 0-truncated:
\vspace{-2mm}
\begin{equation}
  \label{ZeroTruncationOfFreeHomotopyQuotient}
    i_0
    \,
    \tau_0
    \left(
    \HomotopyQuotient{X}{\mathcal{G}}
    \right)
    \;\;
    \simeq
    \;\;
    \HomotopyQuotient{X}{\mathcal{G}}
    \,,
    \;\;\;\;\;\;\;\;\;
    \begin{tikzcd}
      \Topos_0
      \ar[rr, hook, shift right=5pt, "i_0"{below}]
      \ar[rr, phantom, "\sim"]
      &&
      \Topos
      \mathrlap{\,.}
      \ar[ll, shift right=5pt, "\mathclap{\times}"{description, pos=0}, "\tau_0"{above}]
    \end{tikzcd}
\end{equation}

\vspace{-2mm}
\noindent
With this, we may conclude as follows:
\vspace{-3mm}
$$
  \def\arraystretch{2}
  \begin{array}{lll}
    \HomotopyQuotient{X}{\mathcal{G}}
    &
    \;\simeq\;
    i_0
    \left(
      \tau_0
      \left(
        \HomotopyQuotient{X}{\mathcal{G}}
      \right)
    \right)
    &
    \proofstep{ by \eqref{ZeroTruncationOfFreeHomotopyQuotient} }
    \\
    & \;\simeq\;
    i_0
    \Big(
      \tau_0
      \big(\,
      \colimit{[n] \in \Delta^{\mathrm{op}}}
      \,
      \mathcal{G}^{\times^n}
        \times
      X
      \big)
    \Big)
    &
    \proofstep{ by \eqref{HomotopyQuotientAsHomotopyColimit} }
    \\
    &
    \;\simeq\;
    i_0
    \Big(
      \,
      \colimit{[n] \in \Delta^{\mathrm{op}}}
      \,
      \big(\tau_0(\mathcal{G})\big)^{\times^n}
        \times
      \tau_0(X)
    \Big)
    &
    \proofstep{ by \eqref{nTruncationReflection} }
    \\
    &
    \;\simeq\;
    i_0
    \Big(
      \mathrm{coeq}
      \big(
        \tau_0(\mathcal{G}) \times \tau_0(X)
        \rightrightarrows
        \tau_0(X)
      \big)
    \Big)
    &
    \proofstep{ by Ex. \ref{ColimitsOverSimplicialDiagramsOfZeroTruncatedObjectsAreCoequalizers} }
    \\
    &
    \;=\;
    i_0
    \left(
      \tau_0(X)/\tau_0(\mathcal{G})
    \right)
    \,.
  \end{array}
$$

\vspace{-7mm}
\end{proof}

\begin{definition}[Principal $\infty$-bundles]
  \label{PrincipalInfinityBundles}
%
  Given an $\infty$-topos $\Topos$,
  with $X \,\in\, \Topos$
  and
  $\mathcal{G} \,\in\, \Groups(\Topos) \xrightarrow{ \Groups(X \times(-)) }
  \Groups( \Slice{\Topos}{X})$,

  \vspace{1mm}
  \noindent {\bf (i)} we say that
  a regular $\infty$-action $\mathcal{G} \acts \, P$
  (Def. \ref{FreeInfinityAction}) of
  $\mathcal{G}$ internal to the slice topos $\SliceTopos{X}$
  (Prop. \ref{SliceInfinityTopos})
  is a
  {\it formally $\mathcal{G}$-principal $\infty$-bundle} over $X$
  and an actual {\it $\mathcal{G}$-principal $\infty$} over $X$
  if it is inhabited, $P \twoheadrightarrow \ast$ (Nota. \ref{InhabitedObject}).

\noindent {\bf (ii)}   We write
  \begin{equation}
    \label{InfinityCategoryOfPrincipalInfinityBundles}
    \PrincipalBundles{\mathcal{G}}(\Topos)_{X}
    \;\;
    \xhookrightarrow{\quad}
    \;\;
    \FormallyPrincipalBundles{\mathcal{G}}(\Topos)_{X}
    \;\;
    \xhookrightarrow{\quad}
    \;\;
    \Actions{G}(\SliceTopos{X})
  \end{equation}
  for the full sub-$\infty$-categories of
  the $\infty$-category of all $\mathcal{G}$-actions
  \eqref{InfinityCategoryOfInfinityActionsInInfinityTopos}
  on the (formally) principal ones.
\end{definition}

\begin{remark}[Local triviality]
  In contrast to the the 1-category theoretic analog
  Def. \ref{TerminologyForPrincipalBundles}
  (Rem. \ref{AssumptionOfLocalTrivializability}),
  the above
  Def. \ref{PrincipalInfinityBundles} does not include an
  explicit local triviality clause, but just formulates the
  internal notion of (formally) principal bundles
  (Ntn. \ref{InternalizationOfPrincipalBundleTheory}).
  It is instead the nature of the ambient $\infty$-toposes which
  implies the local triviality of internal principal bundles
  (Thm. \ref{OrdinaryPrincipalBundlesAmongPrincipalInfinityBundles},
  Thm. \ref{BorelClassificationOfEquivariantBundlesForResolvableSingularitiesAndEquivariantStructure} below).
\end{remark}

\begin{proposition}[Every $\infty$-action is principal over its homotopy quotient]
  \label{EveryInfinityActionIsPrincipalOverItsHomotopyQotient}
  In an $\infty$-topos $\Topos$,
  every $\infty$-action $\mathcal{G} \acts \, P \,\in\, \Actions{\mathcal{G}}(\Topos)$
  (Def. \ref{ActionObjectsInAnInfinityTopos})
  is principal (Def. \ref{PrincipalInfinityBundles})
  relative to the coprojection into its homotopy quotient
  $X \,\coloneqq\, \HomotopyQuotient{P}{\mathcal{G}}$
  \eqref{HomotopyQuotientAsHomotopyColimit}.
\end{proposition}
\begin{proof}
Since the left base change functor
$\underset{X}{\sum} \,\colon\, \Slice{\Topos}{X} \xrightarrow{\;} \Topos$
(Prop. \ref{BaseChange})

-- preserves $\infty$-colimits and hence homotopy quotients, being a left adjoint,

-- sends products to fiber products,

\noindent
the defining diagram on the left
is sent by $\sum_X$ to the following diagram on the right:
$$
  \Slice{\Topos}{X}
  \;
  \in
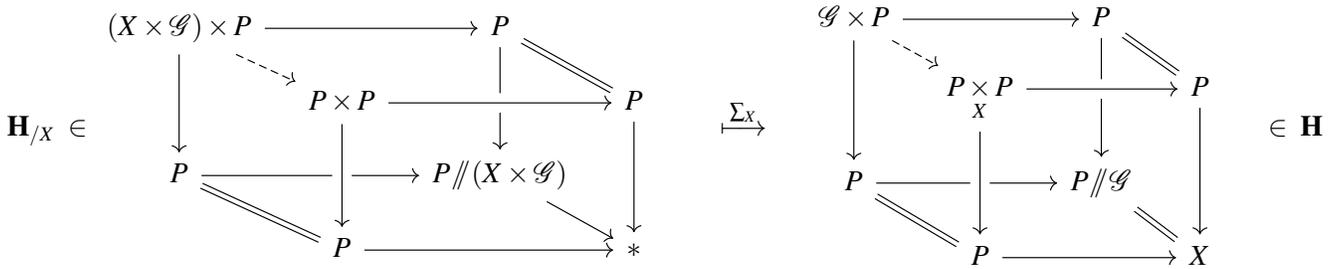
  \begin{tikzcd}[row sep=10pt, column sep=12pt]
    (X \times \mathcal{G})
    \times
    P
    \ar[dr, dashed]
    \ar[rr]
    \ar[dd]
    &&
    P
    \ar[dr, -, shift left=1pt]
    \ar[dr, -, shift right=1pt]
    \ar[dd]
    \\
    &
    P \times P
    \ar[rr, crossing over]
    &&
    P
    \ar[dd]
    \\
    P
    \ar[rr]
    \ar[dr, -, shift left=1pt]
    \ar[dr, -, shift right=1pt]
    &&
    \HomotopyQuotient{P}{(X \times \mathcal{G})}
    \ar[dr]
    \\
    &
    P
    \ar[from=uu, crossing over]
    \ar[rr]
    &&
    \ast
  \end{tikzcd}
  \;\;\;\;\;\;\;\;\;
  \overset{
    \sum_X
  }{\longmapsto}
  \;\;\;\;\;
  \begin{tikzcd}[row sep=10pt, column sep=12pt]
    \mathcal{G}
      \times
    P
    \ar[dr, dashed]
    \ar[rr]
    \ar[dd]
    &&
    P
    \ar[dr, -, shift left=1pt]
    \ar[dr, -, shift right=1pt]
    \ar[dd]
    \\
    &
    P \underset{X}{\times} P
    \ar[rr, crossing over]
    &&
    P
    \ar[dd]
    \\
    P
    \ar[rr]
    \ar[dr, -, shift left=1pt]
    \ar[dr, -, shift right=1pt]
    &&
    \HomotopyQuotient{P}{\mathcal{G}}
    \ar[dr, -, shift left=1pt]
    \ar[dr, -, shift right=1pt]
    \\
    &
    P
    \ar[from=uu, crossing over]
    \ar[rr]
    &&
    X
  \end{tikzcd}
  \;\;\;\;\;\;
  \in
  \;
  \Topos
$$
Here the dashed morphism on the right is manifestly an equivalence.
Since left base change also reflects equivalences,
so is the dashed morphism on the left,
which is the claim to be proven.
\end{proof}

\begin{theorem}
  [Delooping groupoids are moduli stacks for principal $\infty$-bundles {\cite[Thm. 3.17]{NSS12a}}]
  \label{DeloopingGroupoidsAreModuliInfinityStacksForPrincipalInfinityBundles}
  In an $\infty$-topos $\Topos$, for every $\infty$-group
  $\mathcal{G} \,\in\, \Groups(\Topos)$ (Def. \ref{GroupObjectsInAnInfinityTopos})
  there is for each $X \,\in\, \Topos$
  a natural equivalence of $\infty$-categories
  $$
    \begin{tikzcd}[column sep=huge]
      \PointsMaps{}
       { X }
       { \mathbf{B}\mathcal{G} }
      \ar[r, "\sim"{swap}]
      \ar[
        rr,
        rounded corners,
        to path={
          -- ([yshift=-10pt]\tikztostart.south)
          -- node[above]{ \scalebox{.7}{$ \mathrm{hofib} $} }
             ([yshift=-9pt]\tikztotarget.south)
          -- ([yshift=-00pt]\tikztotarget.south)
        }
      ]
      &
      \PrincipalBundles{\mathcal{G}}(\Topos)_X
      \ar[r, hook]
      &
      \Actions{\mathcal{G}}(\SliceTopos{X})
    \end{tikzcd}
$$
$$
\hspace{-3.5cm}
    (
      X \xrightarrow{\;c\;} \mathbf{B}\mathcal{G}
    )
    \;\;
    \longmapsto
    \;\;
  \bigg(
    \begin{tikzcd}[row sep=10pt]
      P
      \ar[rr]
      \ar[d]
      \ar[drr, phantom, "\mbox{\tiny \rm (pb)}"]
      &&
      \ast
      \ar[d, ->>]
      \\
      X
      \ar[rr, "c"{swap}]
      &&
      \mathbf{B}\mathcal{G}
    \end{tikzcd}
    \bigg)
  $$
  between $\mathcal{G}$-principal $\infty$-bundles over $X$
  (Def. \ref{PrincipalInfinityBundles})
  and morphisms from $X$ to the delooping
  $\mathbf{B}\mathcal{G}$ \eqref{DeloopingOfInfinityGroupAsColimit}.
\end{theorem}
\begin{proof}
  The homotopy fiber functor is
  fully faithful by Prop. \ref{EquivalentPerspectivesOnGroupoidObjectsInAnInfinityTopos},
  factors as shown
  by Prop. \ref{EveryInfinityActionIsPrincipalOverItsHomotopyQotient},
  and its factorization is essentially surjective
  by Prop. \ref{InhabitedHomotopyQuotientsOfRegularInfinityActionsAreTerminal}.
\end{proof}

\section{Cohesive homotopy theory}
\label{CohesiveHomotopyTheory}

Where general $\infty$-toposes reflect geometry
(geometric homotopy theory, \cref{ToposTheory})
in a broad sense, encompassing exotic notions of space
such as encountered, say, in arithmetic geometry,
some of their properties may be too exotic
compared to tamer notions of space
as expected in classical differential topology (\cite{Milnor64}\cite{Benedetti21});
for instance in that their shape
(Def. \ref{ShapeOfAnInfinityTopos}) may be a pro-$\infty$-groupoid
that is not represented by a plain $\infty$-groupoid.

\medskip
The notion of {\it cohesive $\infty$-toposes} (Def. \ref{CohesiveInfinityTopos} below)
narrows in on those {\it gros} $\infty$-toposes that share more of the
abstract properties of categories of spaces expected in differential topology;
for instance in that

\noindent {\bf (1)} all their slices have
genuine $\infty$-groupoidal shape
(Prop. \ref{ShapeOfSliceOfCohesiveInfinityToposIsCohesiveShape} below)
and

\noindent {\bf (2)} they contain a good supply of {\it concrete spaces},
namely of sets equipped with cohesive (e.g. topological or smooth) structure
(see Prop. \ref{DiffeologicalSpacesAreTheOneConcreteZeroTruncatedSmoothInfinityGroupoids} below).

\medskip

We first discuss the general concept of cohesive $\infty$-groupoids
and then consider standard models that encapsulate traditional
differential and equivariant topology:

\medskip

-- \cref{GeneralCohesion}: smooth cohesion,

-- \cref{GeneralSingularCohesion}: orbi-singular cohesion.

\medskip
\begin{definition}[Cohesive $\infty$-topos
{\cite[\S 3.1]{SSS09}\cite{dcct}\cite[\S 3.1.1]{SS20OrbifoldCohomology}},
following \cite{Lawvere07}]
  \label{CohesiveInfinityTopos}
  $\,$

  \noindent
  {\bf (i)}
  An $\infty$-topos $\ModalTopos{\smooth}$ over
  some base $\infty$-topos $\BaseTopos$
  is
  called {\it cohesive} if its global section geometric morphism
  $\ModalTopos{\smooth} \xrightarrow{\Points \coloneqq \GlobalSections} \BaseTopos$
  extends to an adjoint quadruple of the following form:
  \vspace{-1mm}
  \begin{equation}
    \label{CohesiveAdjointQuadruple}
    \begin{tikzcd}[column sep=large]
      \ModalTopos{\smooth}
      \ar[
        rr,
        phantom,
        shift left=9pt,
        "{\scalebox{.6}{$\bot$}}"
      ]
      \ar[
        rr,
        phantom,
        shift left=27pt,
        "{\scalebox{.6}{$\bot$}}"
      ]
      \ar[
        rr,
        phantom,
        shift left=-9pt,
        "{\scalebox{.6}{$\bot$}}"
      ]
      \ar[
        rr,
        "{\Points}"{description}
      ]
      \ar[
        rr,
        shift left=36pt,
        "{\Shape}"{description},
        "\mathclap{\times}"{description, pos=0}
      ]
      &&
      \BaseTopos
      \ar[
        ll,
        hook',
        shift right=18pt,
        "{ \Discrete }"{description}
      ]
      \ar[
        ll,
        hook',
        shift right=-18pt,
        "{ \Chaotic }"{description}
      ]
      &
      \ar[
        r,
        phantom,
        shift left=36pt,
        "{
          \mbox{
            \tiny
            \color{greenii}
            \bf
            shape
          }
        }"
      ]
      \ar[
        r,
        phantom,
        shift left=18pt,
        "{
          \mbox{
            \tiny
            \color{greenii}
            \bf
            discrete
          }
        }"
      ]
      \ar[
        r,
        phantom,
        shift left=0pt,
        "{
          \mbox{
            \tiny
            \color{greenii}
            \bf
            points
          }
        }"
      ]
      \ar[
        r,
        phantom,
        shift right=18pt,
        "{
          \mbox{
            \tiny
            \color{greenii}
            \bf
            chaotic
          }
        }"
      ]
      &[-20pt]
      {}
    \end{tikzcd}
  \end{equation}
  meaning that:

  {\bf (a)} $\Discrete \,\coloneqq\, \LocallyConstant$
  is fully faithful \eqref{FullyFaithfulInfinityFunctor}
  and has a left adjoint \eqref{AdjunctionAndHomEquivalence},
  denoted $\Shape$, which preserves
  finite products.

  {\bf (b)} $\Points \,\coloneqq\, \GlobalSections$ has a right adjoint,
  denoted $\Chaotic$.

 \noindent
 \noindent {\bf (ii)}  We denote the induced adjoint triple of endo-$\infty$-functors on $\Topos$  by
 \vspace{-2mm}
 \begin{equation}
  \label{TheModalitiesOnACohesiveInfinityTopos}
  \begin{array}{cccc}
    \raisebox{2pt}{
      \tiny
      \color{darkblue}
      \bf
      \begin{tabular}{c}
        pure shape
        \\
        aspect
      \end{tabular}
    }
    &
    \shape
    &
    \coloneqq
    &
    \Discrete \circ \Shape
    \\
    &
    \scalebox{.7}{$\bot$}
    \\
    \raisebox{2pt}{
      \tiny
      \color{darkblue}
      \bf
      \begin{tabular}{c}
        purely discrete
        \\
        aspect
      \end{tabular}
    }
    &
    \flat
    &
    \coloneqq
    &
    \Discrete \circ \Points
    \\
    &
    \scalebox{.7}{$\bot$}
    \\
    \raisebox{2pt}{
      \tiny
      \color{darkblue}
      \bf
      \begin{tabular}{c}
        purely continuous
        \\
        aspect
      \end{tabular}
    }
    &
    \sharp
    &
    \coloneqq
    &
    \Chaotic \circ \Points
  \end{array}
\end{equation}
\end{definition}
\begin{remark}[Idempotency of cohesive adjoints]
  \label{IdempotencyOfCohesiveAdjoints}
  The fact that
  the adjoints $\Discrete$ and $\Chaotic$ in \eqref{CohesiveAdjointQuadruple}
  are fully faithful is equivalent to the following (co)unit transformations
  being equivalences
  \vspace{-2mm}
  \begin{equation}
    \label{TheIdempotencyOfCohesiveAdjoints}
    \begin{tikzcd}
      \Shape \circ \Discrete
      \ar[
        r,
        "{ \epsilon^{\scalebox{.6}{$\shape$}} }"{above},
        "\sim"{below}
      ]
      &
      \mathrm{id}
      \,,
    \end{tikzcd}
    \;\;\;\;\;\;
    \begin{tikzcd}
      \mathrm{id}
      \ar[
        r,
        "{ \eta^\flat }"{above},
        "\sim"{below}
      ]
      &
      \Points
        \circ
      \Discrete
      \,,
    \end{tikzcd}
    \;\;\;\;\;\;
    \begin{tikzcd}
      \Points \circ \Chaotic
      \ar[
        r,
        "{ \epsilon^{\scalebox{.6}{$\sharp$}} }"{above},
        "\sim"{below}
      ]
      &
      \mathrm{id}
      \,,
    \end{tikzcd}
  \end{equation}

  \vspace{-2mm}
  \noindent
  and equivalent to the (co)monads \eqref{TheModalitiesOnACohesiveInfinityTopos}
  being idempotent, in that their (co)multiplications are natural equivalences:
  \vspace{-2mm}
  \begin{equation}
    \label{IdempotencyOfCohesiveModalities}
    \begin{tikzcd}
      \shape \circ \shape
      \ar[
        rr,
        " \scalebox{0.8}{$\Discrete\big( \epsilon^{\scalebox{.7}{$\shape$}}_{\Shape_{(-)}} \big)$} "{above},
        "\sim"{below}
      ]
      &&
      \shape
      \,,
    \end{tikzcd}
    \;\;\;\;\;\;\;\;
    \begin{tikzcd}
      \flat
      \ar[
        rr,
        " \scalebox{0.8}{$\Discrete\big( \eta^{\scalebox{1}{$\flat$}}_{\Points_{(-)}} \big)$} "{above},
        "\sim"{below}
      ]
      &&
      \flat \circ \flat
      \,,
    \end{tikzcd}
    \;\;\;\;\;\;\;\;
    \begin{tikzcd}
      \sharp \circ \sharp
      \ar[
        rr,
        " \scalebox{0.8}{$ \Chaotic\big( \epsilon^{\scalebox{1}{$\sharp$}}_{\Points_{(-)}} \big)$} "{above},
        "\sim"{below}
      ]
      &&
      \sharp
      \,,
    \end{tikzcd}
  \end{equation}

  \vspace{-1mm}
  \noindent
  as well as
  \begin{equation}
    \label{MixedIdempotencyOfCohesiveModalities}
    \begin{tikzcd}
      \shape \circ \flat
      \ar[
        rr,
        "\scalebox{0.8}{$\Discrete\big( \epsilon^{\scalebox{.7}{$\shape$}}_{\Points_{(-)}} \big)$} "{above},
        "\sim"{below}
      ]
      &&
      \flat
      \,,
    \end{tikzcd}
    \;\;\;\;\;\;\;\;
    \begin{tikzcd}
      \shape
      \ar[
        rr,
        " \scalebox{0.8}{$\Discrete\big( \eta^{\scalebox{1}{$\flat$}}_{\Shape_{(-)}} \big)$} "{above},
        "\sim"{below}
      ]
      &&
      \flat \circ \shape
      \,.
    \end{tikzcd}
  \end{equation}
  \vspace{-.2mm}
  \noindent
  It is this idempotency which makes these (co)monads act
  as {\it projecting out} certain qualities or {\it modes of being} of objects,
  here: {\it being pure shape}, {\it being purely discrete} and {\it being chaotic}\footnote{
  This follows the traditional terminology for ``chaotic topology'';
  see Ex. \ref{ChaoticTopolgy}.},
  while the adjointness
  $
    \shape
    \,\dashv\,
    \flat
    \,\dashv\,
    \sharp
  $ \,
  makes

  \vspace{-.2mm}
 \noindent  them project out {\it dual} such qualities, with the (co)unit
  transformations exhibiting every object $X \,\in\,\mathbf{H}$ as having
  quality in between the given opposing extremes, e.g.:
  \vspace{-3mm}
  \begin{equation}
    \label{PointsToPiecesTransform}
    \begin{tikzcd}[column sep=large]
      \overset{
        \mathclap{
        \raisebox{6pt}{
          \tiny
          \color{darkblue}
          \bf
          \begin{tabular}{c}
            purely
            \\
            discrete aspect
          \end{tabular}
        }
        }
      }{
        \flat
        \,
        X
      }
      \ar[
        r,
        "{ \epsilon^{\scalebox{.6}{$\flat$}}_X }",
        start anchor={[xshift=3pt]}
      ]
      \ar[
        rr,
        rounded corners,
        to path={
          -- ([yshift=-7.5pt]\tikztostart.south)
          --node[below]{
              \mbox{
                \tiny
                \color{greenii}
                \bf
                points-to-pieces transform
              }
            }
            ([yshift=-6pt]\tikztotarget.south)
          -- (\tikztotarget.south)
        }
      ]
      &
      X
      \ar[
        r,
        "{
           \eta^{\scalebox{.6}{$\shape$}}_X
        }",
        end anchor={[xshift=-3pt]}
      ]
      &
      \overset{
        \mathclap{
        \raisebox{6pt}{
          \tiny
          \color{darkblue}
          \bf
          \begin{tabular}{c}
            pure
            \\
            shape aspect
          \end{tabular}
        }
        }
      }{
        \shape
        \,
        X
      }.
    \end{tikzcd}
  \end{equation}
\end{remark}

The following simple class of examples of cohesive $\infty$-toposes
(Ex. \ref{ExamplesOfDiscreteCohesion})
is already interesting in itself
(as in our key Def. \ref{GEquivariantAndGloballyEquivariantHomotopyTheories} below);
and other important examples
(such as the pivotal $\SmoothInfinityGroupoids$, Ntn. \ref{SmoothInfinityGroupoids} below)
are obtained
by (co-)restricting an example in this class from $\infty$-presheaves to $\infty$-sheaves
(when possible):
\begin{example}[Discrete cohesion]
  \label{ExamplesOfDiscreteCohesion}
  Let $\InfinitySite$ be a small $\infty$-category regarded
  as an $\infty$-site with trivial Grothendieck topology. If
  $\InfinitySite$ has a terminal object,
  then its $\infty$-(pre-)sheaves
  $\mathbf{H} \,=\, \InfinitySheaves(\InfinitySite)$
  form a cohesive $\infty$-topos
  (Def. \ref{CohesiveInfinityTopos}).
    A transparent abstract way to see this is to note that the existence of a terminal object
  means equivalently that its inclusion $i$ is right adjoint to the unique functor
  $p$ to the terminal category. This implies the adjoint
  quadruple by Lem. \ref{KanExtensionOfAdjointPairOfInfinityFunctors}
  with Lem. \ref{LeftKanExtensionOfFullyFaithfulFunctorIsFullyFaithful}:
  \vspace{-1mm}
  $$
  \hspace{-2cm}
    \begin{tikzcd}
      \InfinitySite
      \ar[
        rr,
        phantom,
        "{ \scalebox{.5}{$\bot$} }"
      ]
      \ar[
        rr,
        shift left=8pt,
        "{ p }"{description}
      ]
      &&
      \ast
      \ar[
        ll,
        shift left=8pt,
        hook',
        "{ i }"{description}
      ]
    \end{tikzcd}
    \qquad
    \Rightarrow
    \qquad
    \begin{tikzcd}[column sep=large]
      \InfinityPresheaves(\InfinitySite)
      \;\;
      \ar[
        rr,
        "{ p_\ast \,\simeq\, i^\ast    }"{description}
      ]
      \ar[
        rr,
        shift left=32pt,
        "{ p_! }"{description, pos=.4},
        "\mathclap{\times}"{description, pos=0}
      ]
      &&
      \;\;
      \InfinitySheaves(\ast)
      \mathrlap{
        \;\simeq\;
        \InfinityGroupoids
      \;.}
      \ar[
        ll,
        hook',
        shift right=16pt,
        "{ p^\ast \,\simeq\, i_! }"{description}
      ]
      \ar[
        ll,
        hook',
        shift right=-16pt,
        "{ i_\ast }"{description, pos=.4}
      ]
      \ar[
        ll,
        phantom,
        shift right=8pt,
        "{\scalebox{.6}{$\bot$}}"
      ]
      \ar[
        ll,
        phantom,
        shift right=24pt,
        "{\scalebox{.6}{$\bot$}}"
      ]
      \ar[
        ll,
        phantom,
        shift right=-8pt,
        "{\scalebox{.6}{$\bot$}}"
      ]
    \end{tikzcd}
  $$

  \vspace{-1mm}
  \noindent
  Moreover:

  \noindent
  (1) the left Kan extension $i_!$ is fully-faithful since $i$ is
  (by Lem. \ref{LeftKanExtensionOfFullyFaithfulFunctorIsFullyFaithful})
  and it preserves finite products since since $p$ does
  (by  Lem. \ref{LeftKanExtensionOnRepresentablesIsOriginalFunctor},
  Lem. \ref{LeftKanExtensionPreservesBinaryProductsIfOriginalFunctorDoes})

  \noindent
  (2) a right adjoint of a right adjoint of a fully faithful functor
  is itself fully faithful.
\end{example}

\begin{remark}[The axioms on the shape modality]
\label{TheAxiomsOnTheShapeModality}
That the cohesive shape operation $\shape$ \eqref{TheModalitiesOnACohesiveInfinityTopos}
preserves finite products means that

\noindent
{\bf (a)}
it preserves binary products, in that
for all $X, Y \,\in\, \Topos$ we have natural equivalences
\vspace{-2mm}
\begin{equation}
  \label{ShapePreservesBinaryProducts}
  \shape
  (X \times Y)
  \;\simeq\;
  (\,\shape\, X)
  \times
  (\, \shape \, Y)
  \,,
  \;\;\;\;\;\;
  \mbox{equivalently}
  \;\;\;\;\;\;
  \Shape(X \times Y)
  \;\simeq\;
  \Shape(X) \times \Shape(Y)\;.
\end{equation}

\vspace{-2mm}
\noindent
{\bf (b)} it preserves the terminal object, in that
\vspace{-2mm}
\begin{equation}
  \label{ShapePreservesTheTerminalObject}
  \shape \ast
  \;\simeq\;
  \ast
  \,,
  \;\;\;\;\;
  \mbox{equivalently}:
  \;\;\;\;\;
  \mathrm{Shp}(\ast)
  \;\simeq\;
  \ast
  \,.
\end{equation}

\vspace{-2mm}
\noindent The following are some key implications of the fact/requirement that
the shape modality preserves finite products:

\noindent
{\bf (i)}
together with idempotency of modalities (Rem. \ref{IdempotencyOfCohesiveAdjoints}),
it implies the {\it projection formula}:
\vspace{-1mm}
\begin{equation}
  \label{ProjectionFormulaForShape}
  \shape
  (
    X \,\times\, \shape  \, Y
  )
  \;\simeq\;
  (
    \shape \, X
  )
  \times
  (
    \shape \, Y
  )
  \,,
  \;\;\;\;\;\;
  \mbox{equivalently}:
  \;\;\;\;\;\;
  \Shape
  (
    X
      \,\times\,
    \Discrete(S)
  )
  \;\simeq\;
  \Shape(X) \,\times\, S
  \,.
\end{equation}

\vspace{-2mm}
\noindent
{\bf (ii)}
together with the preservation \eqref{InfinityAdjointPreservesInfinityLimits}
of simplicial colimits and the effectiveness
of group objects (Prop. \ref{EquivalentPerspectivesOnGroupoidObjectsInAnInfinityTopos}),
it implies that passage to shape preserves
group objects (Def. \ref{GroupObjectsInAnInfinityTopos})
and their deloopings
\eqref{DeloopingOfInfinityGroupAsColimit}
and more generally action groupoids
(Def. \ref{ActionObjectsInAnInfinityTopos})
and their
homotopy quotients \eqref{HomotopyQuotientAsHomotopyColimit}:
\vspace{-2mm}
\begin{equation}
  \label{ShapePreservesDeloopings}
  \shape
  \,
  \mathbf{B}
  \,
  \mathcal{G}
  \;
  \simeq
  \;
  \mathbf{B}
  \,
  \shape
  \,
  \mathcal{G}
  \;
  \overset{
    \mbox{
      \tiny
      \eqref{ClassifyingSpaceIsShapeOfModuliStack}
    }
  }{
    =:
  }
  \;
  B
  \,
  \mathcal{G}
  \,,
  \;\;\;\;\;
  \shape
  (\HomotopyQuotient{X}{\mathcal{G}})
  \;\simeq\;
  \HomotopyQuotient{(\shape X)}{(\shape \mathcal{G})}
\end{equation}

\vspace{-2mm}
\noindent
The analogous statement holds for $\flat$, for the same reasons:
  \vspace{-2mm}
\begin{equation}
  \label{FlatPreservesDeloopings}
  \flat
  \,
  \mathbf{B}\mathcal{G}
  \;
  \simeq
  \;
  \mathbf{B} \flat \mathcal{G}
  \;
  \overset{
    \mbox{
      \tiny
      \eqref{ClassifyingSpaceIsShapeOfModuliStack}
    }
  }{
    =:
  }
  \;
  B \flat \mathcal{G}
  \,,
  \;\;\;\;\;
  \flat
  (\HomotopyQuotient{X}{\mathcal{G}})
  \;\simeq\;
  \HomotopyQuotient{(\flat X)}{(\flat  \mathcal{G})}
  \,,
\end{equation}

  \vspace{-2mm}
\noindent
but $\flat$ preserves also all fiber products (in particular)
and hence preserves all groupoid objects
(Def. \ref{GroupoidObjectInAnInfinityTopos},
not just the group action objects);

\noindent
{\bf (iii)}
concretely, \eqref{ShapePreservesTheTerminalObject} implies that
$\flat$ computes the underlying (geometrically discrete)
{\it $\infty$-groupoid of points}
\begin{equation}
  \label{FlatModalityComputesUnderlyingPoints}
  \def\arraystretch{1.3}
  \begin{array}{lll}
    \flat X
    \;\simeq\;
    &
    \InfinityGroupoids
    (\ast
      ,\,
      \flat X
    )
    &
    \proofstep{ by \eqref{GroupoidsAreTheirOwnGroupoidsOfPoints} }
    \\
    & \;\simeq\;
    \PointsMaps{}
      { \shape \ast }
      { X }
    &
    \proofstep{ by \eqref{CohesiveAdjointQuadruple} with \eqref{InfinityAdjointPreservesInfinityLimits} }
    \\
    & \;\simeq\;
    \PointsMaps{}
      { \ast }
      { X }
    &
    \proofstep{ by \eqref{ShapePreservesTheTerminalObject}; }
  \end{array}
\end{equation}

\noindent
{\bf (iv)}
together with the mapping-stack adjunction \eqref{InternalHomAdjunction},
it implies a canonical natural transformation
for $X, A\,\in\, \Topos$
\vspace{-2mm}
\begin{equation}
  \label{ComparisonMorphismFromShapeOfMappingStackToMappingSpaceOfShapes}
  \shape
  \,
  \Maps{}
    { X }
    { A }
  \xrightarrow{\; \scalebox{.7}{$\widetilde {\shape\mathrm{ev}}$}  \;}
  \Maps{\big}
    { \shape X }
    { \shape A }
  \,,
\end{equation}

\vspace{-2mm}
\noindent
this being the adjunct \eqref{AdjunctionAndHomEquivalence}
of the shape
of the evaluation map \eqref{EvaluationMap}:
\vspace{-2mm}
$$
  ( \shape\,  X )
  \times
  \,
  \left(
    \shape
    \,
    \Maps{}
      { X }
      { A }
  \right)
  \xrightarrow{\; \sim \;}
  \shape
  \,
  \left(
    X \times \Maps{}{X}{A}
  \right)
  \xrightarrow{\; \scalebox{.7}{$\shape \mathrm{ev}$} \;}
  \shape \,
  A
  \,.
$$
\end{remark}

Less immediate is that the shape operation also preserves certain classes of homotopy fiber products: this is the content of Prop. \ref{ShapeFunctorPreservesHomotopyFibersOverDiscreteObjects} and Prop. \ref{ShapePreservesHomotopyFibersOfDeloopingsOutOfDiscreteDomains} below.

\newpage
\noindent
{\bf Some basic properties of cohesive $\infty$-toposes.}
\begin{proposition}[Shape of slice of cohesive $\infty$-topos is cohesive shape]
  \label{ShapeOfSliceOfCohesiveInfinityToposIsCohesiveShape}
  Let $\Topos$ be a cohesive $\infty$-topos
  (Def. \ref{CohesiveInfinityTopos})
  with shape modality \eqref{TheModalitiesOnACohesiveInfinityTopos}
  \vspace{-2mm}
  $$
    \shape
      \;\colon\;
    \Topos
      \xrightarrow{\;\;\Shape \;\;}
    \InfinityGroupoids
      \xrightarrow{\;\; \Discrete \;\;}
    \Topos
    \,.
  $$

  \vspace{-1mm}
\noindent  Then for any $X \,\in\, \Topos$ the shape of the
  slice $\infty$-topos over $X$, in the sense of Def. \ref{ShapeOfAnInfinityTopos},
  is equivalently the cohesive shape of $X$, under
  the embedding \eqref{InfinityGroupoidsInsideProInfinityGroupoids}
  of $\infty$-groupoids into pro-$\infty$-groupoids:
  \vspace{-1mm}
  $$
    \Shape
    (
      \ModalTopos{/X}
    )
    \;\;
      \simeq
    \;\;
    \Shape(X)
    \;\;\;
    \in
    \;
    \InfinityGroupoids
    \xhookrightarrow{\;\;}
    \ProObjects\InfinityGroupoids
    \,.
  $$
\end{proposition}
\begin{proof}
The base change adjoint triple of the slice (Prop. \ref{BaseChange})
combines with part (1) of cohesion (Def. \ref{CohesiveInfinityTopos})
\vspace{-2mm}
\begin{equation}
  \label{GlobalSectionsGeometricMorphismsOfSliceOfCohesiveInfinityTopos}
  \begin{tikzcd}[column sep=large]
    \mathllap{
      \Gamma_X
      \;\;
      :
      \;\;\;
    }
    \mathbf{H}_{/X}
    \ar[
      rr,
      phantom,
      shift left=10+10pt,
      "{\scalebox{.7}{$\bot$}}"
    ]
    \ar[
      rr,
      phantom,
      shift left=-10+10pt,
      "{\scalebox{.7}{$\bot$}}"
    ]
    \ar[
      rr,
      shift left=20+10pt,
      "{ \sum_X }"{description}
    ]
    \ar[
      rr,
      shift left=-20+10pt,
      "{ \prod_X }"{description}
    ]
    &&
    \quad
    \mathbf{H}
    \quad
    \ar[
      ll,
      shift right=+10pt,
      "{ X \times (-) }"{description}
    ]
    \ar[
      rr,
      phantom,
      shift left=10+10pt,
      "{\scalebox{.7}{$\bot$}}"
    ]
    \ar[
      rr,
      phantom,
      shift left=-10+10pt,
      "{\scalebox{.7}{$\bot$}}"
    ]
    \ar[
      rr,
      shift left=20+10pt,
      "{ \mathrm{Shp} }"{description},
      "{\mathclap{\times}}"{description, pos=0}
    ]
    \ar[
      rr,
      shift left=-20+10pt,
      "{ \mathrm{Pnts} }"{description}
    ]
    &&
    \mathrm{Grp}_\infty
    \mathrlap{
      \;\;\;
      :
      \;\;
      \LocallyConstant_X
    }
    \ar[
      ll,
      hook',
      shift right=+10pt,
      "{ \mathrm{Disc} }"{description}
    ]
  \end{tikzcd}
\end{equation}
  and by essential uniqueness of adjoints and of
  the terminal geometric morphism, the composite
  is equivalently the global section geometric morphism
  $(\Gamma_X \,\dashv\,)$ of the slice, as shown above.

  Now consider, for $S_1, S_2 \,\in\, \InfinityGroupoids$
  the following sequence of natural equivalences:
  \vspace{-1mm}
$$
 \def\arraystretch{1.3}
  \begin{array}{lll}
       \mathrm{Grp}_\infty
    \left(
      S_1
      ,\,
      \GlobalSections_X
      \circ
      \LocallyConstant_X(S_2)
    \right)
    &
       \;\simeq\;
    \SlicePointsMaps{\big}{X}
     { \LocallyConstant_X( S_1 ) }
     { \LocallyConstant_X(S_2) }
    \big)
    &
    \mbox{\small by the composite adjunction in \eqref{GlobalSectionsGeometricMorphismsOfSliceOfCohesiveInfinityTopos} }
    \\
    & \;\simeq\;
    \SlicePointsMaps{\big}{X}
      { X \times \mathrm{Disc}(S_1) }
      { X \times \mathrm{Disc}(S_2) }
    \big)
    &
    \mbox{\small by the factorization in \eqref{GlobalSectionsGeometricMorphismsOfSliceOfCohesiveInfinityTopos} }
    \\
    &
    \;\simeq\;
    \PointsMaps{\big}
      { X \times \mathrm{Disc}(S_1) }
      { \mathrm{Disc}(S_2) }
    &
    \mbox{\small by the left adjunction in \eqref{GlobalSectionsGeometricMorphismsOfSliceOfCohesiveInfinityTopos} }
    \\
    & \;\simeq\;
    \mathrm{Grp}_\infty
    \left(
      \mathrm{Shp}
      \left(
        X \times \mathrm{Disc}(S_1)
      \right)
      ,\,
      S_2
    \right)
    &
    \mbox{\small by the right adjunction in \eqref{GlobalSectionsGeometricMorphismsOfSliceOfCohesiveInfinityTopos} }
    \\
    & \;\simeq\;
    \mathrm{Grp}_\infty
    \left(
      \mathrm{Shp}(X)
      \times
      S_1
      ,\,
      S_2
    \right)
    &
    \mbox{\small by \eqref{ProjectionFormulaForShape} }
    \\
    & \;\simeq\;
    \mathrm{Grp}_\infty
    \left(
      S_1
      ,\,
      \mathrm{Grp}_\infty
      \left(
        \mathrm{Shp}(X)
        ,\,
        S_2
      \right)
    \!\right)
    &
    \mbox{\small by \eqref{HomAdjunctionInInfinityGroupoids}\,. }
  \end{array}
$$

\vspace{-1mm}
\noindent With this, the $\infty$-Yoneda lemma (Prop. \ref{InfinityYonedaLemma})
implies the natural equivalence
that is to be proven.
\end{proof}

\begin{lemma}[Cohesive points are intrinsic points]
\label{CohesivePointsAreIntrinsicPoints}
If $\Topos$ is cohesive over $\InfinityGroupoids$,
then for all $X \in \Topos$ we have a natural equivalence

\vspace{-.5cm}
\begin{equation}
  \label{CohesivePointsArePoints}
  \mathrm{Pnts}(X)
  \;\simeq\;
  \Topos(\ast,X)\;.
\end{equation}
\end{lemma}
\begin{proof}
\vspace{-3mm}
$$
  \begin{array}{lll}
    \mathrm{Pnts}(X)
    &
    \;\simeq\;
    \InfinityGroupoids
    \left(
      \ast,
      \,
      \mathrm{Pnts}(X)
    \right)
    &
    \mbox{\small by \eqref{GroupoidsAreTheirOwnGroupoidsOfPoints}}
    \\
    & \;\simeq\;
    \Topos
    \left(
      \Discrete(\ast),
      X
    \right)
    &
    \mbox{\small by cohesion}
    \\
    & \;\simeq\;
    \Topos(\ast,X)
    &
    \mbox{\small by cohesion}
    \,.
  \end{array}
$$

\vspace{-5mm}
\end{proof}

\begin{lemma}[Mapping space consists of cohesive points of mapping stack]
\label{MappingSpaceConsistsOfCohesivePointsOfMappingStack}
If $\Topos$ is cohesive over $\BaseTopos$,
then for all $X,A \in \Topos$ we have a natural equivalence
between their hom-space \eqref{HomSpace}
and the cohesive points of their
mapping stack \eqref{InternalHomAdjunction}:
\vspace{-1mm}
\begin{equation}
  \label{PointsMapsSpaceIsPointsOfMappingStack}
  \PointsMaps{}
    { X }
    { A }
  \;\simeq\;
  \Points
  \,
  \Maps{}
    { X }
    { A }
  \,.
\end{equation}
More generally, for $B \,\in\, \Topos $
and $E_1, E_2 \,\in\, \Topos_{/B}$, the cohesive points
of the slice mapping stack (Def. \ref{SliceMappingStack})
yield the mapping space of the slice:
\vspace{-1mm}
\begin{equation}
  \label{SlicePointsMapsSpaceIsPointsOfSliceMappingStack}
  \SlicePointsMaps{}{B}
    { E_1 }
    { E_2 }
  \;\simeq\;
  \Points
  \,
  \SliceMaps{}{B}
    { E_1 }
    { E_2 }
  \,.
\end{equation}
\end{lemma}
\begin{proof}
For the first case:
\vspace{-1mm}
$$
  \begin{array}{lll}
    \PointsMaps{}
      {X}{A}
    &
    \;\simeq\;
    \PointsMaps{\big}
      {\ast}
      { \Maps{}{X}{A} }
    &
    \mbox{\small by Ex. \ref{MappingSpaceConsistsOfPointsOfMappingStack}}
    \\
    & \;\simeq\;
    \Points
    \,
    \Maps{}
      {X}{A}
    &
    \mbox{\small by Ex. \ref{CohesivePointsAreIntrinsicPoints}}
    \,.
  \end{array}
$$

  \vspace{-2mm}
\noindent More generally, for the second case:
\vspace{-2mm}
$$
  \begin{array}{lll}
    \SlicePointsMaps{}{B}
      {E_1}
      {E_2}
    &
    \;\simeq\;
    \
    \Maps{\big}
      { \ast }
      {
        \SliceMaps{}{B}
          { E_1 }
          { E_2 }
      }
    &
    \mbox{\small by Lem. \ref{PlotsOfSliceMappingStackAreSliceHoms}}
    \\
    & \;\simeq\;
    \Points
    \,
    \SliceMaps{}{B}
      { E_1 }
      { E_2 }
    &
    \mbox{\small by Ex. \ref{CohesivePointsAreIntrinsicPoints}}.
  \end{array}
$$

\vspace{-6mm}
\end{proof}

\begin{proposition}[Shape preserves fiber products over discrete objects
{\cite[Thm. 3.8.19]{dcct}\cite[Lem. 3.5]{SS20OrbifoldCohomology}}]
  \label{ShapeFunctorPreservesHomotopyFibersOverDiscreteObjects}
  $\,$

  \noindent
  For $\Topos$ a cohesive $\infty$-topos
  (Def. \ref{CohesiveInfinityTopos})
  over $\InfinityGroupoids$, for
  $B \in \InfinityGroupoids \xhookrightarrow{\Discrete} \Topos$
  a discrete object and $X,Y \in \Slice{\Topos}{B}$
  two cohesive objects fibered over $B$, there is a natural equivalence
  \vspace{-1mm}
  \begin{equation}
    \label{ShapePreservingFiberProductsOverDiscreteObjects}
    \shape\,
    \big(
      X \underset{B}{\times} Y
    \big)
    \;\simeq\;
    (\,\shape \; X)
      \underset{B}{\times}
    (\,\shape \; Y)
  \end{equation}
  between the shape \eqref{TheModalitiesOnACohesiveInfinityTopos}
  of their homotopy fiber product over $S$
  and the homotopy fiber product of their shapes.
\end{proposition}
\begin{proof}
The $\Shape \,\dashv\, \Discrete$-adjunction passes to slices
(\cite[Prop. 5.2.5.1]{Lurie09HTT}) and, by a generalization of the
Grothendieck construction, slices of $\infty$-toposes over
inverse images of $\infty$-groupoids $B$ are equivalent to the $\infty$-categories
of functors from $B$ to the them (\cite[Lem. 3.10]{BP19}):
$$
  \begin{tikzcd}
    \Functors
    \big(
      B
      ,\,
      \Topos
    \big)
    \ar[
      from=rrrrrr,
      rounded corners,
      to path={
           ([yshift=0pt]\tikztostart.south)
        -- ([yshift=-12pt]\tikztostart.south)
        -- node[yshift=-6.8pt]{ \scalebox{.8}{$\Functors(B,\Discrete)$} }
           ([yshift=-12pt]\tikztotarget.south)
        -- ([yshift=-00pt]\tikztotarget.south)
      }
    ]
    \ar[
      rrrrrr,
      rounded corners,
      to path={
           ([yshift=0pt]\tikztostart.north)
        -- ([yshift=+12pt]\tikztostart.north)
        -- node[yshift=+6.8pt]{ \scalebox{.8}{$\Functors(B,\Shape)$} }
           ([yshift=+12pt]\tikztotarget.north)
        -- ([yshift=-00pt]\tikztotarget.north)
      }
    ]
    &[-24pt]
    \simeq
    &[-24pt]
    \SliceTopos{\Discrete(B)}
    \ar[
      rr,
      shift left=7pt,
      "\Slice{\Shape}{B}"
    ]
    \ar[
      from=rr,
      shift left=7pt,
      "\Slice{\Discrete}{B}"
    ]
    \ar[
      rr,
      "{\scalebox{.8}{$\bot$}}",
      phantom
    ]
    &&
    \Slice{\InfinityGroupoids}{B}
    &[-24pt]\simeq&[-24pt]
    \Functors
    \big(
      B
      ,\,
      \InfinityGroupoids
    \big)
    \,.
  \end{tikzcd}
$$
Under these equivalences, the right adjoint is readily seen to
be given by objectwise application, as shown at the bottom,
from which it follows by uniqueness of $\infty$-adjoints that the
left adjoint is similarly given by objectwise application of the
shape operation, as shown at the top. But since
limits of $\infty$-presheaves (here: over $B$) are computed objectwise,
and since both $\Shape$ and $\Discrete$ preserve finite products,
it follows that so do their slicings over $S$.
Now the claim follows since the fiber product in question is
equivalently the plain product in  the slice $\SliceTopos{S}$.
\end{proof}

\begin{lemma}[Reverse pasting law for $\infty$-groupoids]
\label{ReversePastingLawForInfinityPullbacks}
Given a homotopy pasting diagram in
$\InfinityGroupoids$ of the form
$$
  \begin{tikzcd}
    X
    \ar[r]
    \ar[d]
    \ar[
      dr,
      phantom,
      "{\mbox{\tiny{\rm(pb)}}}"
    ]
    &
    Y
    \ar[r]
    \ar[d]
    \ar[
      dr,
      phantom
    ]
    &
    Z
    \ar[d]
    \\
    A
    \ar[
      r,
      ->>,
      "{p}"{swap}
    ]
    &
    B
    \ar[r]
    &
    C
    \mathrlap{\,,}
  \end{tikzcd}
$$
 such that

  - the total rectangle and the left square are  homotopy Cartesian \eqref{HomotopyCartesianSquare},

  -  the bottom left morphism is an effective epimorphism (Def. \ref{EffectiveEpimorphisms}),

  \noindent
  then also the right square is homotopy Cartesian.
\end{lemma}
\begin{proof}
  Using the presentation of the situation in the classical model structure on simplicial sets (Ntn. \ref{SimplicialSetsAndInfinityGroupoids}), we may consider the analogous pasting diagram in $\SimplicialSets$, such that the two bottom morphisms, say, are Kan fibrations. Since $\SimplicialSets$ is a presheaf category and hence regular (by Exp. \ref{ToposesAreRegular}) we are thus reduced to showing that the 1-category theoretic reverse pasting law from Lem. \ref{ReversePastingLawInRegularCategories} applies, which is the case if the bottom left Kan fibration (representing) $p$ is a regular epimorphism. Since $\SimplicialSets$ is in fact a topos, this is equivalent to showing that this Kan fibration is an epimorphism (again by Exp. \ref{ToposesAreRegular}), hence a degreewise surjection of simplicial sets. But the assumption that $p$ is an effective epimorphism of $\InfinityGroupoids$ means that it is a surjection on connected components (by Prop. \ref{DetectingEffectiveEpimorphismsOnZeroTruncation}), and every Kan fibration which is a surjection on connected components is already surjective in degree 0 (seen by lifting against the acyclic cofibrations $\Simplex{0} \hookrightarrow \Simplex{1}$) and hence in each degree $n$ (seen by lifting against $\Simplex{0} \hookrightarrow \Simplex{n}$).
\end{proof}

\begin{proposition}[Shape preserves homotopy fibers of deloopings out of discrete domains]
\label{ShapePreservesHomotopyFibersOfDeloopingsOutOfDiscreteDomains}
  For $\Topos$ a cohesive $\infty$-topos
  (Def. \ref{CohesiveInfinityTopos})
  over $\InfinityGroupoids$,
  the shape operation preserves homotopy fibers of morphisms
  out of cohesively discrete $\infty$-groupoids,
  $X \,\simeq\, \flat X$,  into deloopings of cohesive $\infty$-groups $\mathcal{G}$:
  $$
    \begin{tikzcd}[column sep=14pt]
      \mathrm{fib}(f)
      \ar[r]
      &
      X
      \ar[
        r,
        "{
          f
        }"
      ]
      &
      \mathbf{B}\mathcal{G}
    \end{tikzcd}
    \hspace{1cm}
    \vdash
    \hspace{1cm}
    \begin{tikzcd}[column sep=14pt]
      \shape \mathrm{fib}(f)
      \,\simeq\,
      \mathrm{fib}
      \big(
        \shape f
      \big)
      \ar[r]
      &
      X
      \ar[
        r,
        "{
          \scalebox{.8}{$
            \shape f
          $}
        }"
      ]
      &
      B\mathcal{G}
    \end{tikzcd}
  $$
\end{proposition}
\begin{proof}
  Since every $X \,\in\, \InfinityGroupoids \xhookrightarrow{\Discrete} \SmoothInfinityGroupoids$ is a coproduct of its connected components, and
  since  coproducts are preserved by homotopy pullback (Prop.
  \ref{ColimitsAreUniversalInAnInfinityTopos}) and by $\shape$ \eqref{InfinityAdjointPreservesInfinityLimits} we may assume without restriction of generality that $X$ is connected, hence (due to Prop. \ref{LoopingAndDeloopingEquivalence}) that $X \,\simeq\, B \mathcal{K}$ for some $\infty$-group $\mathcal{K}$, so that the point inclusion $\begin{tikzcd}[column sep=8pt]\ast \ar[r, ->>, shorten=-2pt] & X\end{tikzcd}$ is an effective epimorphism of (smooth) $\infty$-groupoids (by Prop. \ref{DetectingEffectiveEpimorphismsOnZeroTruncation}).
  Consider then the pasting diagram shown on the left below, where we are using the pasting law \eqref{HomotopyPastingLaw}
  and the looping/delooping equivalence \eqref{LoopingAndDelooping}
  in order to identify $\mathcal{G}$ as the fiber appearing in the top left:
  $$
    \begin{tikzcd}
      \mathcal{G}
      \ar[r]
      \ar[d]
      \ar[
        dr,
        phantom,
        "{\mbox{\tiny\rm(pb)}}"
      ]
      &
      \mathrm{fib}
      (
        f
      )
      \ar[d]
      \ar[r]
      \ar[
        dr,
        phantom,
        "{\mbox{\tiny\rm(pb)}}"
      ]
      &
      \ast
      \ar[d]
      \\
      \ast
      \ar[
        r,
        ->>
      ]
      &
      B\mathcal{K}
      \ar[
        r,
        "{f}"{swap}
      ]
      &
      \mathbf{B}\mathcal{G}
    \end{tikzcd}
    \hspace{1cm}
    \vdash
    \hspace{1cm}
    \begin{tikzcd}
      \shape \mathcal{G}
      \ar[r]
      \ar[d]
      \ar[
        dr,
        phantom,
        "{\mbox{\tiny\rm(pb)}}"
      ]
      &
      \shape
      \mathrm{fib}
      (
        f
      )
      \ar[d]
      \ar[r]
      \ar[
        dr,
        phantom,
        "{\mbox{\tiny\rm(pb)}}"
      ]
      &
      \ast
      \ar[d]
      \\
      \ast
      \ar[
        r,
        ->>
      ]
      &
      B\mathcal{K}
      \ar[
        r,
        "{
          \scalebox{.8}{$
            \raisebox{-1pt}{$\shape$} f
          $}
        }"{swap}
      ]
      &
      B\mathcal{G}
      \mathrlap{\,.}
    \end{tikzcd}
  $$
  Under applying the shape modality to this diagram, shown on the right, the total rectangle remains a homotopy pullback square with the bottom left morphism remaining an effective epimorphism (both by Rem. \ref{TheAxiomsOnTheShapeModality}), and the left square remains a pullback by Prop. \ref{ShapeFunctorPreservesHomotopyFibersOverDiscreteObjects}. Therefore the
  reverse pasting law for homotopy pullbacks (Prop. \ref{ReversePastingLawForInfinityPullbacks})
  implies that also the right square remains a homotopy pullback, which is equivalently the claim to be shown.
\end{proof}

\noindent
{\bf Cohesive charts.}

\begin{notation}[Cohesive charts {\cite[Def. 3.9]{SS20OrbifoldCohomology}}]
  \label{CohesiveCharts}
  Given a cohesive $\infty$-topos $\Topos$
  we say that an $\infty$-site $\Charts$ for $\Topos$
  is an {\it $\infty$-category of cohesive charts}
  if, under the $\infty$-Yoneda embedding
  \eqref{InfinityYonedaEmbedding},
  all its objects have contractible shape:
  \vspace{-2mm}
  \begin{equation}
    \label{CohesiveChartsAreGeometricallyContractible}
    \forall \, \TopologicalPatch \,\in\, \Charts \xhookrightarrow{ \;\;y\;\; } \Topos
    \;\;
    :
    \;\;
    \shape\, \TopologicalPatch
    \;\simeq\;
    \ast
    \,;
  \end{equation}

  \vspace{-2mm}
  \noindent
  equivalently:
  \begin{equation}
    \label{CohesiveChartsHaveContractibleShape}
    \forall \, \TopologicalPatch \,\in\, \Charts \xhookrightarrow{\; \;y\; \;} \Topos
    \;\;
    :
    \;\;
    \mathrm{Shp}(\TopologicalPatch)
    \;\simeq\;
    \ast
    \,.
  \end{equation}
\end{notation}

\begin{example}[Flat modality over a site of cohesive charts]
  \label{FlatModalityOverSiteOfCohesiveCharts}
  If a cohesive $\infty$-topos $\Topos$
  over $\InfinityGroupoids$
  has cohesive $\Charts$
  (Ntn. \ref{CohesiveCharts}) then for all $S \in \InfinityGroupoids$
  the corresponding geometrically discrete cohesive object is
  constant as an $\infty$-presheaf on the charts, in that for all
  $U \in \Charts$ there is a natural equivalence:
  \vspace{-2mm}
  \begin{equation}
    \label{DiscreteObjectsAreConstantPresheavesOnCharts}
    \left(
      \mathrm{Disc}(S)
    \right)(\TopologicalPatch)
    \;\simeq\;
    S
    \;\;\;
    \in
    \;
    \InfinityGroupoids
    \;.
  \end{equation}

    \vspace{-2mm}
\noindent
  Moreover, for $A \,\in\,$there is a natural equivalence
    \vspace{-2mm}
  $$
    (\flat A)(\TopologicalPatch)
    \;\simeq\;
    A(\ast)
    \,.
  $$
\end{example}
\begin{proof}
  The first statement is the composite of the
  following sequence of natural equivalences:
  \vspace{-2mm}
  $$
  \def\arraystretch{1.2}
    \begin{array}{lll}
    \left(
      \mathrm{Disc}(S)
    \right)(\TopologicalPatch)
    &
    \;\simeq\;
    \Topos
    \left(
      \TopologicalPatch,
      \,
      \mathrm{Disc}(S)
    \right)
    &
    \proofstep{ by \eqref{YonedaEquivalence} }
    \\
    & \;\simeq\;
    \Topos
    \left(
      \mathrm{Shp}(\ast)
      \,
      S
    \right)
    &
    \proofstep{ by \eqref{TheModalitiesOnACohesiveInfinityTopos} }
    \\
    & \;\simeq\;
    \Topos
    (
      \ast,
      \,
      S
    )
    &
    \proofstep{ by \eqref{CohesiveChartsHaveContractibleShape} }
    \\
    & \;\simeq\;
    S
    &
    \proofstep{ by \eqref{GroupoidsAreTheirOwnGroupoidsOfPoints} }
    .
    \end{array}
  $$
  From this the second statement follows by
  Lem. \ref{CohesivePointsAreIntrinsicPoints}.
\end{proof}
\begin{lemma}[Mapping stack into discrete object is discrete]
  \label{MappingStackIntoDiscreteObjectIsDiscrete}
  Let $\Topos$ be a cohesive $\infty$-topos with
  a category $\Charts$
  of cohesive charts (Ntn. \ref{CohesiveCharts}).
  Then:

  \noindent
  {\bf (i)}  the hom-equivalence of the $\shape \,\dashv\, \flat$-adjunction
  also holds internally,

  \noindent
 {\bf  (ii)}
  a mapping stack \eqref{InternalHomAdjunction} into
  a discrete object is
  itself discrete:
  \vspace{-1mm}
  \begin{equation}
    \label{EquivalenceExhibitingMappingStackIntoDiscreteObjectAsDiscrete}
    \flat A \,\simeq\, A
    \;\;\;\;\;\;\;\;\;\;
    \Rightarrow
    \;\;\;\;\;\;\;\;\;\;
    \underset{X \in \Topos}{\forall}
    \;\;\;\;
    \Maps{}
      { X }
      { A }
    \;\simeq\;
    \Maps{}
      { \shape \, X }
      { A }
    \;\;
    \simeq
    \;\;
    \flat
    \Maps{}
     { X }
     { A }
     \;\;
     \overset{
       \mathclap{
         \mbox{\tiny \eqref{MixedIdempotencyOfCohesiveModalities}}
       }
     }{\simeq}
     \;\;
     \shape
     \,
     \Maps{}
       { X }
       { A }
    \,.
  \end{equation}
\end{lemma}
\begin{proof}
  By the $\infty$-Yoneda lemma (Prop. \ref{InfinityYonedaLemma}),
  we may check the statement on $\TopologicalPatch \in \Charts$:
    \vspace{-2mm}
  $$
  \hspace{-1cm}
  \def\arraystretch{1.3}
    \begin{array}{lll}
      \Maps{}
        { X }
        { A }
      (\TopologicalPatch)
      &
      \;\simeq\;
      \Topos
      (
        \TopologicalPatch \times X, \, A
      )
      &
      \mbox{\small by \eqref{ValuesOfMappingStackAsHomSpaces}}
      \\
      & \;\simeq\;
      \Topos
      (
        \TopologicalPatch \times X, \, \flat A
      )
      &
      \proofstep{ by assumption}
      \\
      & \;\simeq\;
      \Topos
      \left(
        \shape \,
        (\TopologicalPatch \times X),
        \,
        A
      \right)
      &
      \proofstep{ by \eqref{TheModalitiesOnACohesiveInfinityTopos} }
      \\
      & \;\simeq\;
      \Topos
      \left(
        (\shape \, \TopologicalPatch )
        \times
        (\shape \, X),
        \,
        A
      \right)
      &
      \proofstep{ by \eqref{ShapePreservesBinaryProducts} }
      \\
      & \;\simeq\;
      \Topos
      \left(
        \shape \, X,
        \,
        A
      \right)
      &
      \proofstep{ by \eqref{CohesiveChartsAreGeometricallyContractible} }
      \\
      & \;\simeq\;
      \Topos
      (
        X,
        \,
        \flat A
      )
      &
      \proofstep{ by \eqref{TheModalitiesOnACohesiveInfinityTopos} }
      \\
      & \;\simeq\;
      \Topos
      (
        X,
        \,
        A
      )
      &
      \proofstep{ by assumption }
      \\
      & \;\simeq\;
      \Topos
      (
        \ast \times X,
        \,
        A
      )
      \\
      & \;\simeq\;
      \mathrm{Pts}
      \,
      \Maps{}{X}{A}
      &
      \proofstep{ by \eqref{ValuesOfMappingStackAsHomSpaces} }
      \\
      & \;\simeq\;
      \left(
      \Discrete
      \,
      \Points
      \,
      \Maps{}{X}{A}
      \right)(\TopologicalPatch)
      &
      \proofstep{ by \eqref{DiscreteObjectsAreConstantPresheavesOnCharts} }
      \\
      & \;=\;
      \left(
        \flat
        \Maps{}
          {X}{A}
      \right)(\TopologicalPatch)
      &
      \proofstep{ by \eqref{TheModalitiesOnACohesiveInfinityTopos} }
      .
    \end{array}
  $$

   \vspace{-7mm}
\end{proof}

\subsection{Smooth Cohesion}
\label{GeneralCohesion}

We turn attention now to the archetypical example of those
cohesive $\infty$-toposes (Def. \ref{CohesiveInfinityTopos})
given by generalized spaces that are built form (are $\infty$-colimits of)
{\it Cartesian spaces} $\mathbb{R}^n$ glued by smooth functions between them,
i.e. built from the elementary {\it charts} of differential topology.
The {\it concrete} 0-truncated objects among these
{\it smooth $\infty$-groupoids}
are
known as {\it diffeological spaces} (Ntn. \ref{CartesianSpacesAndDiffeologicalSpaces} below)
and we begin by recalling and developing
general topology as seen inside the category of diffeological spaces.

\medskip

\noindent
{\bf Diffeological spaces.}

\begin{definition}[Good open covers]
  \label{GoodOpenCovers}
  For $\TopologicalSpace \,\in\, \TopologicalSpaces$
  a topological manifold,

\noindent {\bf (i)}  an open cover $\{ \TopologicalPatch_i \xhookrightarrow{\;} \TopologicalSpace \}_{i \in I}$
  is {\it good} if all finite intersections of patches are
  (empty or) homeomorphic to an open ball:
  \vspace{-2mm}
  \begin{equation}
    \label{ConditionOnGoodOpenCover}
    \underset{
      n \in \mathbb{N}
    }{\forall}
    \;
    \underset{i_0, \cdots, i_n \,\in I\,}{\forall}
    \;\;\;\;
      \TopologicalPatch_{i_0}
        \cap
      \cdots
        \cap
      \TopologicalPatch_{i_n}
      \;\;
      \simeq
      \;\;
      \left\{
      \begin{array}{ll}
        \mathbb{R}^n & \mbox{or}
        \\
        \varnothing\;.
      \end{array}
      \right.
       \end{equation}

      \vspace{-2mm}
\noindent {\bf (ii)}
If $\TopologicalSpace$ is a smooth manifold, we say a good
  cover is {\it differentiably good}
  (\cite[Def. 6.3.9]{FStS12})
  if the identifications in
  \eqref{ConditionOnGoodOpenCover} can be taken to be
  diffeomorphisms instead of just homeomorphisms.
\end{definition}

\begin{notation}[Cartesian spaces and Diffeological spaces]
  \label{CartesianSpacesAndDiffeologicalSpaces}
 {\bf (i)}  We write
 \vspace{-2mm}
  \begin{equation}
    \label{CartesianSpaces}
    \CartesianSpaces \,
    \xhookrightarrow{\quad}
    \SmoothManifolds
  \end{equation}

  \vspace{-2mm}
  \noindent
  for the category whose objects are the $\mathbb{R}^n$ for $n \,\in\, \mathbb{N}$
  and whose morphisms are the smooth functions between these.

  \noindent
  {\bf (ii)}
  We regard this as a site with respect to the coverage (Grothendieck pre-topology)
  whose covers of $\mathbb{R}^n$ are the
  {\it differentiably good} open covers (\cite[Def. 6.3.9]{FStS12})
  hence the open covers
  $\{U_i \xhookrightarrow{\;} \mathbb{R}^n\}_{i \in I}$,
  such that all non-empty finite intersections
  $U_{i_0} \cap \cdots \cap U_{i_k}$, $k \,\in\, \mathbb{N}$,
  of patches are {\it diffeomorphic} to
  an open ball, and hence to $\mathbb{R}^n$.

    \noindent {\bf (iii)}  We say that a sheaf on this site is a {\it smooth set}
  and is a {\it diffeological space}
  (\cite{Souriau80}\cite{Souriau84}\cite{IglesiasZemmour85}, see \cite{BaezHoffnung08} \cite{IglesiasZemmour13})
  if it is a {\it concrete sheaf} (\cite{Dubuc79}):
  \vspace{-2mm}
  \begin{equation}
    \label{DiffeologicalSpaces}
    \DiffeologicalSpaces
    \xhookrightarrow{ \;\; i_{\sharp_1} \;\; }
    \SmoothSets
    \;\coloneqq\;
    \Sheaves(\CartesianSpaces)\;.
  \end{equation}

  \vspace{-2mm}
\noindent  This means that a functor
\vspace{-2mm}
  \begin{equation}
    \label{DiffeologicalSpacesAsFunctorOfPlots}
    \begin{tikzcd}[column sep=2pt, row sep=-4pt]
      X \,\colon\,
      \CartesianSpaces^{\mathrm{op}}
      \ar[rr]
      &&
      \Sets
      \\
  \scalebox{0.7}{$    \mathbb{R}^n $}
      &\quad \longmapsto&
    \phantom{AAAAA}    {
         \scalebox{0.7}{$  X(\mathbb{R}^n) $}
      } \;\;\;
       \scalebox{0.95}{
          \tiny
          \color{darkblue}
          \bf
              set of plots
}
    \end{tikzcd}
  \end{equation}

\vspace{-2mm}
\noindent  encodes a diffeological space if and only if it is a sheaf and
  for all $n \,\in\, \mathbb{N}$ the following
  natural morphism is a monomorphism,
  hence, for each $n$, an injection of the set of plots into the set of all functions
  of underlying point sets:
  \vspace{-2mm}
  $$
    \begin{tikzcd}[column sep=-4pt, row sep=-6pt]
      X(\mathbb{R}^n)
      &\simeq&
      \SmoothSets
      (
        \mathbb{R}^n
        ,\,
        X
      )
      \ar[
        rr,
        hook
      ]
      &[10pt]
      &
      \Sets
      \big(
        \SmoothSets(\ast ,\, \mathbb{R}^n)
        ,\,
        \SmoothSets(\ast ,\, X)
      \big)
      &=:&
      (\sharp X)(\mathbb{R}^n)\;.
      \\
      &&
       \scalebox{0.8}{$
         (\mathbb{R}^n \xrightarrow{\phi} X)
      $}
      &\longmapsto&
       \scalebox{0.8}{$
       \big(
         (\ast \xrightarrow{r} \mathbb{R}^n)
         \,\mapsto\,
         (\ast \xrightarrow{r} \mathbb{R}^n \xrightarrow{\phi} X)
       \big)
      $}
    \end{tikzcd}
  $$

  \vspace{-3mm}
  \noindent {\bf (iv)}  The category of diffeological spaces
  contains the category of smooth manifolds
  as a full subcategory, via the extended Yoneda embedding:
  \vspace{-2mm}
  \begin{equation}
    \label{ExtendedYonedaEmbeddingOfSmoothManifoldsIntoDiffeologicalSpaces}
    \CartesianSpaces
      \xhookrightarrow{\quad}
    \SmoothManifolds
      \xhookrightarrow{\;\;\; \YonedaEmbedding \;\;\;}
    \DiffeologicalSpaces
    \,.
  \end{equation}
\end{notation}

\begin{remark}[Diffeological mapping spaces]
  On general grounds,
  the category of diffeological spaces (Ntn. \ref{CartesianSpacesAndDiffeologicalSpaces})
  is Cartesian closed, in that for any $X \,\in\, \DiffeologicalSpaces$
  the cartesian product operation with $X$, given by
  \vspace{-2mm}
  \begin{equation}
    \label{InternalHomAdjunctionInDiffeologicalSpaces}
    X \times Y
      \;\colon\;
    \mathbb{R}^n
    \;\longmapsto\;
    X(\mathbb{R}^n) \times Y(\mathbb{R}^n)
    \,,
  \end{equation}

  \vspace{-2mm}
\noindent  has a right adjoint
  \vspace{-2mm}
  $$
    \begin{tikzcd}[column sep=large]
      \DiffeologicalSpaces \;
      \ar[
        rr,
        phantom,
        "{ \scalebox{.5}{$\bot$} }"
      ]
      \ar[
        rr,
        shift right=5pt,
        "{
          \Maps{}
            {X}{-}
        }"{below}
      ]
      &&
      \;
      \DiffeologicalSpaces
      \ar[
        ll,
        shift right=5pt,
        "{
          X \times (-)
        }"{above}
      ]
    \end{tikzcd}
  $$

  \vspace{-2mm}
\noindent  given by
  \begin{equation}
    \label{FormulaForDiffeologicalMappingSpace}
    \Maps{}
      { X }
      { Y }
    \;\colon\;
    \mathbb{R}^n
    \;\longmapsto\;
    \DiffeologicalSpaces
    \left(
      X \times \mathbb{R}^n
      ,\,
      Y
    \right)
    \,.
  \end{equation}
\end{remark}

\begin{example}[Continuous diffeology and D-topology]
  \label{ContinuousDiffeologyAndDTopology}
  For any
  $X \,\in\, \DiffeologicalSpaces$
  \eqref{DiffeologicalSpaces}
  the underlying set becomes a topological space
  $\DTopology(X) \,\in\, \TopologicalSpaces$
  via the {\it D-topology}
  (\cite[Def. 1.2.3]{IglesiasZemmour85}\cite[\S 2.8]{IglesiasZemmour13}\cite[\S 3]{ChristensenSinnamonWu14}),
  which is
  the final (i.e. finest) topology
  such that all plots \eqref{DiffeologicalSpacesAsFunctorOfPlots}
  become continuous functions:
  \vspace{-2mm}
  $$
    \big(
      \mathbb{R}^n \xrightarrow{\phi} X
    \big)
    \;\in\;
    \DiffeologicalSpaces(\mathbb{R}^n,\, X)
    \;\simeq\;
    X(\mathbb{R}^n)
    \;\;\;\;\;\;
    \vdash
    \;\;\;\;\;\;
    \begin{tikzcd}[column sep=20pt]
      \mathbb{R}^n
      \ar[
        rr,
        "{\phi}"{above},
        "{
          \mbox{
            \tiny
            \rm
            continuous
          }
        }"{below}
      ]
      &&
      \DTopology(X) \;.
    \end{tikzcd}
  $$

  \vspace{-2mm}
\noindent
  Conversely,
  for any $X \,\in\, \TopologicalSpaces$
  the underlying set becomes a diffeological space
  \eqref{DiffeologicalSpaces}
  via the
  {\it continuous diffeology}
  $\ContinuousDiffeology(X) \,\in\, \DiffeologicalSpaces$
  whose plots
  \eqref{DiffeologicalSpacesAsFunctorOfPlots}
  are the continuous functions
  \vspace{-2mm}
  $$
    \ContinuousDiffeology(X)
    \;:\;
    \mathbb{R}^n
    \;\longmapsto\;
    \TopologicalSpaces(\mathbb{R}^n ,\, X)
    \,.
  $$
\end{example}

\begin{notation}[D-topological spaces]
  \label{DeltaGeneratedTopologicalSpaces}
  We write
  \vspace{-2mm}
  $$
    \begin{tikzcd}
      \TopologicalSpaces
      \ar[
        rr,
        phantom,
        "{\scalebox{.5}{$\bot$}}"
      ]
      \ar[
        rr,
        shift right=5pt,
        "{\ContinuousDiffeology}"{below}
      ]
      &&
      \DTopologicalSpaces
      \ar[
        ll,
        hook',
        shift right=5pt
      ]
    \end{tikzcd}
  $$

  \vspace{-2mm}
  \noindent  for the coreflective  subcategory on the
  {\it Delta-generated spaces}
  \cite{Dugger03}
  also called {\it numerically generated spaces}
  \cite{SYH10}
  (see also \cite[\S 3.2]{ChristensenSinnamonWu14}),
  hence the topological spaces whose open subsets $U \,\subset X\,$ are precisely
  those whose pre-images $\phi^{-1}(U)$ are open under all continuous
  functions of the form $\Delta^n \xrightarrow{\phi} X$,
  or, equivalently, under all those of the form $\mathbb{R}^n \xrightarrow{\phi} X$,
  for all $n \in \mathbb{N}$.
\end{notation}
The following statement
justifies abbreviating ``Delta-generated topological spaces''
to ``D-topological spaces'':

\begin{proposition}[Adjunction between topological and diffeological spaces fixes Delta-generated spaces]
  \label{CategoryOfDeltaGeneratedTopologicalSpaces}
  $\,$

  \noindent
  The constructions
  in Ex. \ref{ContinuousDiffeologyAndDTopology}
  constitute a pair of adjoint functors
  between tological and diffeological spaces (Ntn. \ref{CartesianSpacesAndDiffeologicalSpaces})
  which factors
  through compactly generated spaces (Ntn. \ref{CompactlyGeneratedTopologicalSpaces})
  and further through
  {\it D-topological spaces} (Ntn. \ref{DeltaGeneratedTopologicalSpaces})
  as follows:
  \vspace{-3mm}
  \begin{equation}
    \label{CategoryOfDeltaGeneratedSpaces}
    \hspace{-4mm}
    \begin{tikzcd}[column sep=35pt]
    \overset{
      \mathclap{
      \raisebox{4pt}{
        \tiny
        \color{darkblue}
        \bf
        \begin{tabular}{c}
          topological spaces
        \end{tabular}
      }
      }
    }{
    \categorybox{\TopologicalSpaces}
    }
    \;\;
    \ar[
      r,
      phantom,
      "{ \scalebox{.7}{$\bot$} }"
    ]
    \ar[
      r,
      shift right=5.5pt,
      "{ k }"{below}
    ]
    &
    \;\;
    \overset{
      \mathclap{
      \raisebox{4pt}{
        \tiny
        \color{darkblue}
        \bf
        \begin{tabular}{c}
          topological k-spaces
        \end{tabular}
      }
      }
    }{
    \categorybox{\kTopologicalSpaces}
    }
    \;\;
    \ar[
      l,
      hook',
      shift right=5.5pt
    ]
    \ar[
      r,
      phantom,
      "{ \scalebox{.7}{$\bot$} }"
    ]
    \ar[
      r,
      shift right=5.5pt,
      "{
        \underset{
          \mathclap{
          \raisebox{-3pt}{
            \tiny
            \color{greenii}
            \bf
            continuous diffeology
          }
          }
        }{
          \ContinuousDiffeology
        }
      }"{below}
    ]
    &
    \;\;
    \overset{
      \mathclap{
      \raisebox{4pt}{
        \tiny
        \color{orangeii}
        \bf
        \begin{tabular}{c}
          Delta-generated
          \\
          topological spaces
        \end{tabular}
      }
      }
    }{
     \categorybox{\DTopologicalSpaces}
    }
    \;\;
    \ar[
      l,
      hook',
      shift right=5.5pt
    ]
    \ar[
      r,
      phantom,
      "{ \scalebox{.7}{$\bot$} }"
    ]
    \ar[
      r,
      hook,
      shift right=5.5pt,
    ]
    &
    \;\;
    \overset{
      \mathclap{
      \raisebox{3pt}{
        \tiny
        \color{darkblue}
        \bf
        \begin{tabular}{c}
          diffeological spaces
        \end{tabular}
      }
      }
    }{
      \categorybox{\DiffeologicalSpaces}
    }
    \ar[
      l,
      shift right=5.5pt,
      "{
        \overset{
          \raisebox{3pt}{
            \tiny
            \color{greenii}
            \bf
            D-topology
          }
        }{
          \DTopology
        }
      }"{above}
    ]
    \ar[
      r,
      hook,
      "{
        \mbox{
          \tiny
          \eqref{InjectionOfDiffeologicalSpacesAreTheOneConcreteZeroTruncatedSmoothInfinityGroupoids}
        }
      }"{below}
    ]
    &[-16pt]
    \categorybox{\SmoothInfinityGroupoids}
    \end{tikzcd}
  \end{equation}
\end{proposition}
\begin{proof}
  This is essentially the observation of
  \cite[\S 3]{SYH10}, see also \cite[\S 3.2]{ChristensenSinnamonWu14}:
  One checks that
  $\DTopology \,\dashv\, \ContinuousDiffeology$
  is an adjunction
  (\cite[Lem. 3.1]{SYH10})
  whose associated comonad is idempotent (\cite[Lem. 3.3]{SYH10})
  \vspace{-2mm}
  $$
    \begin{tikzcd}
      \DTopology
        \circ
      \ContinuousDiffeology
        \circ
      \DTopology
        \circ
      \ContinuousDiffeology
      \ar[
        rrr,
        "{
         \scalebox{0.7}{$ \DTopology( \epsilon_{\ContinuousDiffeology(-)} ) $}
        }"{above},
        "{ \sim }"{below}
      ]
      &&&
      \DTopology
        \circ
      \ContinuousDiffeology
    \end{tikzcd}
  $$

  \vspace{-3mm}
  \noindent
  and whose fixed objects
  $\begin{tikzcd}[column sep=small]
    \DTopology \circ \ContinuousDiffeology(X)
      \ar[
        r,
        "\epsilon_{X}"{above},
        "\sim"{below}
      ]
    &
    X
    \end{tikzcd}
  $
  are precisely the Delta-generated spaces
  (\cite[Prop. 3.2]{SYH10}).
  This implies the claimed factorization
  relative to all topological spaces
  by the factorization theorem for idempotent adjunctions
  (e.g. \cite[Thm. 3.8.8]{Grandis21}).
  Moreover, Delta-generated spaces
  are also coreflective among
  k-spaces (by \cite[Prop. 1.5]{Vogt71}, see also
  \cite[p. 7]{Gaucher07}), whereby the full factorization
  \eqref{CategoryOfDeltaGeneratedSpaces}
  is implied by uniqueness of adjoints.
\end{proof}

\begin{example}[Chaotic topology]
  \label{ChaoticTopolgy}
  For $S \,\in\, \Sets \xhookrightarrow{\;} \InfinityGroupoids
   \xhookrightarrow{ \Discrete } \SmoothInfinityGroupoids$
  its {\it sharp} or {\it chaotic} aspect, in the notation of
  Def. \ref{CohesiveInfinityTopos}, is the
  D-topological space obtained by equipping $S$ with the
  {\it chaotic topology}, i.e., the {\it coarsest} or {\it initial} topology
  $$
    \sharp
    \,
    \Discrete(S)
    \;\;
    \simeq
    \;\;
    \ContinuousDiffeology
    \left(
      S, \mathrm{Op}_S \coloneqq \{\varnothing, S\}
    \right)
    \,.
  $$
\end{example}

\begin{definition}[Smooth extended simplices]
  \label{SmoothExtendedSimplicies}
   Write
   \vspace{-2mm}
   $$
     \SmoothSimplex{\bullet}
     \;:\;
     \Delta
     \longrightarrow
     \CartesianSpaces
     \xhookrightarrow{\quad}
     \SmoothManifolds
     \xhookrightarrow{\;\;\YonedaEmbedding \;\;}
     \DiffeologicalSpaces
   $$

   \vspace{-2mm}
   \noindent
   for the cosimplicial diffeological space (Ntn. \ref{CartesianSpacesAndDiffeologicalSpaces})
   given by the
   extended smooth simplices:
   \vspace{-2mm}
   $$
     \SmoothSimplex{n}
     \;\coloneqq\;
     \big\{
       (x_0, \cdots, x_n)
       \,\in\,
       \mathbb{R}^{n+1}
       \,
       \big\vert
       \,
       \sum_{k = 0}^n x_k
       \,=\, 1
       \big\}
   $$

   \vspace{-2mm}
   \noindent
   co-degeneracy and co-face maps
   given by addition of consecutive variables and by insertion of zeros, respectively.
\end{definition}
\begin{notation}[Diffeological singular complex]
  \label{DiffeologicalSingularComplex}
  We write
  \vspace{-2mm}
  \begin{equation}
    \label{SmoothSingularSimplicialSet}
    \begin{tikzcd}[row sep=-3pt, column sep=2]
      \SingularSimplicialComplex
      &[-15pt]\;\colon\;&[-15pt]
      \DiffeologicalSpaces
      \ar[
        rr
      ]
      &&
      \SimplicialSets
      \\
      &&
     \scalebox{0.7}{$ X  $}
      &\longmapsto&
 \scalebox{0.7}{$     \DiffeologicalSpaces
      \left(
        \SmoothSimplex{\bullet}
        ,\,
        X
      \right)
      $}
    \end{tikzcd}
  \end{equation}

\vspace{-2mm}
\noindent
for the nerve operation induced by Def. \ref{SmoothExtendedSimplicies}.
\end{notation}

\begin{proposition}[Singular simplicial complex of continuous-diffeology represents weak homotopy type]
  \label{DiffeologicalShapeOfContinuousDiffeologyIsWeakHomotopyType}
  For $\TopologicalSpace \,\in\, \TopologicalSpaces$
  there is a natural transformation
  to the ordinary singular simplicial complex of a topological space
  from the diffeological singular simplial complex \eqref{SmoothSingularSimplicialSet}
  of its continuous-diffeological space \eqref{CategoryOfDeltaGeneratedSpaces}
  which takes values in weak homotopy equivalences:
  \vspace{-2mm}
  $$
    \begin{tikzcd}
      \SingularSimplicialComplex
      \left(
        \ContinuousDiffeology(\TopologicalSpace)
      \right)
      \ar[
        rr,
        "{ \in \WeakHomotopyEquivalences }"{below}
      ]
      &&
      \SingularSimplicialComplex
      (\TopologicalSpace)
      \;\;\;
      \in
      \;
      \SimplicialSets
      \xrightarrow{\;\;}
      \InfinityGroupoids
      \,.
    \end{tikzcd}
  $$
\end{proposition}
\begin{proof}
  The existence of these weak homotopy equivalences
  is the content of \cite[Prop. 4.14]{ChristensenWu14}.
  Inspection of the proof given there shows that these
  are indeed natural transformations.
\end{proof}

\begin{proposition}[Diffeological mapping spaces have correct underlying homotopy type]
  \label{DiffeologicalMappingSpacesHaveCorrectUnderlyingHomotopyType}
  The ordinary singular simplicial complex
  of the mapping space between a pair of k-topological spaces
  \eqref{MappingSpace}
  is naturally weakly homotopy equivalent
  to the diffeological singular simplicial complex
  \eqref{SmoothSingularSimplicialSet}
  of the diffeological mapping space
  \eqref{FormulaForDiffeologicalMappingSpace}
  of their continuous-diffeological incarnations
  \eqref{CategoryOfDeltaGeneratedSpaces}:
  \vspace{-3mm}
  $$
    X,\, Y
    \,\in\,
    \kTopologicalSpaces
    \;\;\;\;\;
    \vdash
    \;\;\;\;\;
    \begin{tikzcd}
      \SingularSimplicialComplex
      \,
      \ContinuousDiffeology
      \left(
        \Maps{}
          { X }
          { Y }
      \right)
      \ar[
        rr,
        "{ \in \WeakHomotopyEquivalences }"{below}
      ]
      &&
      \SingularSimplicialComplex\
      \,
      \Maps{\big}
        { \ContinuousDiffeology(X) }
        { \ContinuousDiffeology(Y) }
      \,.
    \end{tikzcd}
  $$
\end{proposition}
\begin{proof}
  We deduce this as a corollary of results proven in \cite{SYH10}.
  This gives the existence of a weak homotopy equivalence
  $\phi \,\in\, \WeakHomotopyEquivalences$,
  whose image under $\ContinuousDiffeology$ is of this form:
  \vspace{-2mm}
  \begin{equation}
    \label{ComparisonMorphismForIdentifyingHomotopyTypeOfDiffeologicalMappingSpace}
    \begin{tikzcd}
      \ContinuousDiffeology
      \left(
        \Maps{}
          {X}{Y}
      \right)
      \ar[
        rr,
        "{ \scalebox{0.7}{$ \ContinuousDiffeology(\phi)$} }"{above},
        "{\scalebox{0.7}{$ \in \ContinuousDiffeology(\WeakHomotopyEquivalences)$} }"{below}
      ]
      &&
      \ContinuousDiffeology
      \left(
        \mathbf{smap}(X,A)
      \right)
      \,\simeq\,
      \Maps{\big}
        { \ContinuousDiffeology(X) }
        { \ContinuousDiffeology(Y) }
      \,,
    \end{tikzcd}
  \end{equation}

  \vspace{-2mm}
\noindent  where $\mathbf{smap}$ is some topologization of the set of maps
  (defined on \cite[p. 6]{SYH10}) of which we only need to know that:

  \vspace{-3mm}
  \begin{itemize}
  \setlength\itemsep{-3pt}

  \item
  its continuous-diffeological incarnation is
  the diffeological mapping space
  (by \cite[Prop. 4.7]{SYH10})
  as shown on the right of \eqref{ComparisonMorphismForIdentifyingHomotopyTypeOfDiffeologicalMappingSpace},

  \item
  it is weakly homotopy equivalent
  to the topological mapping space with its compact-open topology
  (by \cite[Prop. 5.4]{SYH10}) hence to the k-fication of that
  (by \cite[Prop. 1.2 (h)]{Vogt71}),
  as shown on the left of \eqref{ComparisonMorphismForIdentifyingHomotopyTypeOfDiffeologicalMappingSpace}.

  \end{itemize}

  \vspace{-3mm}
\noindent
  Now consider the naturality square of $\phi$
  \eqref{ComparisonMorphismForIdentifyingHomotopyTypeOfDiffeologicalMappingSpace}
  under the
  natural weak homotopy equivalence from Prop. \ref{DiffeologicalShapeOfContinuousDiffeologyIsWeakHomotopyType}:
  \vspace{-2mm}
  $$
  \hspace{-1mm}
    \begin{tikzcd}
      \SingularSimplicialComplex
      \,
      \ContinuousDiffeology
      \left(
        \Maps{}
          { X }{ Y }
      \right)
      \ar[
        rr,
        "{
          \scalebox{0.7}{$\SingularSimplicialComplex
          \,
          \ContinuousDiffeology
          (\phi)
          $}
        }"
      ]
      \ar[
        d,
        "{
          \in \WeakHomotopyEquivalences
        }"
      ]
      &&
      \SingularSimplicialComplex
      \,
      \ContinuousDiffeology
      \left(
        \mathbf{smap}(X,\,Y)
      \right)
      \ar[
        d,
        "{
          \in \WeakHomotopyEquivalences
        }"
      ]
      &[-33pt]
      \simeq
        \;
        \SingularSimplicialComplex
        \,
        \Maps{\big}
          {\! \ContinuousDiffeology(X) }
          { \ContinuousDiffeology(Y)\! }.
      \\
      \SingularSimplicialComplex
      \left(
        \Maps{}
          { X }{ Y }
      \right)
      \ar[
        rr,
        "{\scalebox{0.7}{$ \SingularSimplicialComplex(\phi)$} }"{above},
        "{ \scalebox{0.7}{$\in \WeakHomotopyEquivalences$} }"{below}
      ]
      &&
      \SingularSimplicialComplex
      \left(
        \mathbf{smap}(X\,Y)
      \right)
    \end{tikzcd}
  $$

  \vspace{-2mm}
\noindent  With the bottom morphism and the vertical morphisms
  being weak homotopy equivalences, it follows that also the top
  morphism is a weak homotopy equivalence (by the 2-out-of-3 property),
  which is the claim to be shown.
\end{proof}
For maps out of discrete spaces, we have stronger statements,
such as the following:
\begin{lemma}[Diffeological mapping groupoid out of discrete into continuous-diffeological groupoid]
  \label{DiffeologicalMappingGroupoidOutOfDiscreteIntoContinuousDiffeologicalGroupoid}
 $\,$

\noindent Let $(S_1 \rightrightarrows S_0) \,\in\, \Groupoids(\Sets)
  \xhookrightarrow{\;\;} \Groupoids(\DTopologicalSpaces)$
  be a topologically discrete groupoid
  and
  $\Gamma \,\in\, \Groups(\kTopologicalSpaces)$
  be any topological group.

 \noindent {\bf (i)}   We have an isomorphism of diffeological groupoids:
 \vspace{-1mm}
  $$
    \Maps{\big}
      {
        (S_1 \rightrightarrows S_0)
      }
      {
        \DeloopingGroupoid
        { \ContinuousDiffeology(\Gamma) }
      \!\!}
    \;\;
      \simeq
    \;\;
    \ContinuousDiffeology
    \,
    \Maps{\big}
      {\! (S_1 \rightrightarrows S_0) }
      { \DeloopingGroupoid{\Gamma} \! }
    \;\;\;
    \in
    \;
    \Groupoids(\DiffeologicalSpaces)
    \,,
  $$
  where on the left the mapping groupoid is formed in
  $\DiffeologicalSpaces$, while on the right it is formed in $\kTopologicalSpaces$.

\noindent {\bf (ii)}  In particular, for $G \,\in\,\Groups(\Sets)$ we have an isomorphism
\vspace{-1mm}
  $$
    \Maps{\big}
      { \mathbf{E}G }
      { \mathbf{B}\ContinuousDiffeology(\Gamma) }
    \;\;
    \simeq
    \;\;
    \ContinuousDiffeology
    \,
    \Maps{}
      { \mathbf{E}G }
      { \mathbf{B}\Gamma }\;.
  $$
\end{lemma}
\begin{proof}
  The mapping groupoids are equalizers
  of maps between product spaces, e.g.
  \vspace{-2mm}
  $$
    \begin{tikzcd}
      \Maps{\big}
        {\! (S_1 \rightrightarrows S_0) }
        { \DeloopingGroupoid{ \Gamma } \!}
      \ar[rr, hook, "{\mathrm{eq}}"]
      &&[-10pt]
      \underset{ f \in S_1 }{\prod} \Gamma
      \ar[
        rrr,
        shift left=3pt,
        "F \,\mapsto\, F(f \circ f')"
      ]
      \ar[
        rrr,
        shift right=3pt,
        "F \,\mapsto\, F(f) \cdot F(f')"{swap}
      ]
      &&&
      \underset{
\scalebox{0.5}{$        { (f,f') \,\in \,}
        { S_1 \underset{S_0}{\times} S_1 }
     $} }{\prod}
      \,
      \Gamma\;.
    \end{tikzcd}
  $$

  \vspace{-3mm}
\noindent   Since
  $\ContinuousDiffeology$ is a right adjoint, it preserves all these limits.
\end{proof}

\medskip

\noindent
{\bf Smooth $\infty$-groupoids.}
Much like diffeological spaces (Ntn. \ref{CartesianSpacesAndDiffeologicalSpaces})
are ``geometric sets probeable by smooth manifolds''
so smooth $\infty$-groupoids
(Ntn. \ref{SmoothInfinityGroupoids} below)
are ``geometric $\infty$-groupoids probeable by smooth manifolds''
({\cite[Def. 3.1]{SSS09}\cite[\S 4.4]{dcct}\cite[Ex. 3.18]{SS20OrbifoldCohomology}\cite[Def. A.57]{FSS20CharacterMap}}). In fact, diffeological spaces are equivalently the
concrete 0-truncated smooth $\infty$-groupoids
(recalled as Prop. \ref{DiffeologicalSpacesAreTheOneConcreteZeroTruncatedSmoothInfinityGroupoids} below). The generalization from concrete sheaves on $\CartesianSpaces$ to
all sheaves embeds diffeological spaces into a cohesive topos, and
discarding the truncation condition makes this a cohesive $\infty$-topos
(Prop. \ref{SmoothCohesion} below).

\begin{notation}[Smooth $\infty$-groupoids]
  \label{SmoothInfinityGroupoids}
  We write
  \vspace{-2mm}
  $$
    \SmoothInfinityGroupoids
    \;\coloneqq\;
    \InfinitySheaves(\CartesianSpaces)
  $$

  \vspace{-2mm}
  \noindent
  for the hypercomplete $\infty$-topos over
  $\CartesianSpaces$ (Ntn. \ref{CartesianSpacesAndDiffeologicalSpaces}),
  presented by the local projective moel structure
  (Ntn. \ref{ModelCategoriesOfSimplicialPresheaves})
  on simplicial presheaves over $\CartesianSpaces$:
  \vspace{-2mm}
\begin{equation}
  \label{InfinityCategoryOfSmoothInfinityGroupoids}
  \begin{aligned}
    \SmoothInfinityGroupoids
    \; &\coloneqq\;
    \Localization{ \LocalWeakEquivalences }
    \big(
      \SimplicialPresheaves(\CartesianSpaces)_{\projloc}
    \big)
    \\
    \;\;
    &
    \simeq
    \;\;
    \Localization{ \LocalWeakEquivalences }
    \big(
      \SimplicialPresheaves(\SmoothManifolds)_{ \projloc }
    \big)
    \;\;
    \in
    \;
    \HomotopyTwoCategory
    (
      \InfinityToposes
    )\;.
  \end{aligned}
\end{equation}
\end{notation}
\begin{remark}[Smooth $\infty$-groupoids presented over smooth manifolds]
  \label{SmoothInfinityGroupoidsPresentedOverSmoothManifolds}
  Since every smooth manifold has, essentially by definition, an open cover by objects of $\CartesianSpaces$,
  the canonical inclusion
  $$
    \CartesianSpaces
    \xhookrightarrow[\mathrm{dense}]{}
    \SmoothManifolds
  $$
  is that of a {\it dense subsite} with respect to (differentiably good) open covers,
  and as such induces an equivalence sheaf toposes over these sites (by the {\it comparison lemma}, e.g. \cite[\S C.2.2]{Johnstone02b}), and hence of hypercomplete $\infty$-toposes:
  \begin{equation}
    \label{SmoothInfinityGroupoidsEquivalentlyStacksOverSmoothManifolds}
    \SmoothInfinityGroupoids
    \,\simeq\,
    \InfinitySheaves(\CartesianSpaces)
    \;\simeq\;
    \InfinitySheaves(\SmoothManifolds)
    \,.
  \end{equation}
  Conceptually this means that the site of $\SmoothManifolds$  has {\it redundancy} when regarded as a generating class for  smooth homotopy types, due to the fact that smooth manifolds are themselves already built from simpler smooth building blocks (charts).
  Indeed, for many purposes the presentation of $\SmoothInfinityGroupoids$ over $\CartesianSpaces$ is the more convenient (cf. for instance
  Prop. \ref{RecognitionOfFibrantObjectsOverCartSp}
  and
  Ex. \ref{CechActionGroupoidOfEquivariantGoodOpenCoverIsLocalCofibrantResolution}  below).
  But for some purposes we will invoke the presentation over $\SmoothManifolds$, such as in the proof of the smooth Oka principle (Thm. \ref{SmoothOkaPrinciple} below).
\end{remark}

\medskip

In addition to the general recognition of projectively co-fibrant simplicial presheaves (from Prop.
\ref{DuggerCofibrancyRecognition}), we have the following convenient recognition principle for fibrant simplicial presheaves presenting smooth $\infty$-groupoids (Ntn. \ref{SmoothInfinityGroupoids}), based on results in \cite{Pavlov22}:

\begin{proposition}[Recognition of fibrant simplicial presheaves over $\CartesianSpaces$]
\label{RecognitionOfFibrantObjectsOverCartSp}
 An object
$\mathcal{X} \,\in\,\SimplicialPresheaves(\CartesianSpaces)_{ \projloc }$ \eqref{InfinityCategoryOfSmoothInfinityGroupoids} is fibrant precisely if it satisfies the following conditions for all $\mathbb{R}^n \,\in\, \CartesianSpaces$:

\vspace{-.2cm}
\begin{enumerate}

\vspace{-.1cm}
\item the value $\mathcal{X}(\mathbb{R}^n) \,\in\, \SimplicialSets_{\mathrm{fib}}$ is a Kan complex, i.e. $\mathcal{X}$ is globally projectively fibrant;

\vspace{-.2cm}
\item
with
\begin{equation}
  \label{StripZigZagCoverOfCartesianSpace}
 \Big(
    U_{2i}
    :=
    (4i , 4i+3) \times \mathbb{R}^{n-1}
    \xhookrightarrow{\phantom{--}\iota_{2i}\phantom{--}}
    \mathbb{R}^n
  \Big)
  _{i \in \Integers}
  ,
  \hspace{.6cm}
  \Big(
    U_{2i+1}
    \,:=\,
    (4i+2 , 4i+5) \times \mathbb{R}^{n-1}
    \xhookrightarrow{\phantom{--}\iota_{2i+1}\phantom{--}}
    \mathbb{R}^n
  \Big)
  _{i \in \Integers}
\end{equation}
a pair of disjoint
collections of
parallel rectangular strips, which jointly cover $\mathbb{R}^n$,

\hspace{-.5cm}
\begin{tikzpicture}[xscale=1.37]

  \clip
    (-.9,2.6)
    rectangle
    (12,-2.3);

 \draw[fill=gray, fill opacity=.3, draw opacity=0]
   (-.45,-1.5)
   rectangle
   (-.35,+1.5);

 \draw[fill=gray, fill opacity=.3, draw opacity=0]
   (-.3,-1.5)
   rectangle
   (-.2,+1.5);

 \draw[fill=gray, fill opacity=.3, draw opacity=0]
   (-.15,-1.5)
   rectangle
   (-.05,+1.5);

 \draw[fill=gray, fill opacity=.3, draw opacity=0]
   (0,-1.5)
   rectangle
   (3,+1.5);

 \draw[fill=gray, fill opacity=.3, draw opacity=0]
   (2,-1.5)
   rectangle
   (5,+1.5);

 \draw[fill=gray, fill opacity=.3, draw opacity=0]
   (4,-1.5)
   rectangle
   (7,+1.5);

 \draw[fill=gray, fill opacity=.3, draw opacity=0]
   (6,-1.5)
   rectangle
   (9,+1.5);

 \draw[fill=gray, fill opacity=.3, draw opacity=0]
   (8,-1.5)
   rectangle
   (11,+1.5);

 \draw[fill=gray, fill opacity=.3, draw opacity=0]
   (11.05,-1.5)
   rectangle
   (11.15,+1.5);

 \draw[fill=gray, fill opacity=.3, draw opacity=0]
   (11.2,-1.5)
   rectangle
   (11.3,+1.5);

 \draw[fill=gray, fill opacity=.3, draw opacity=0]
   (11.35,-1.5)
   rectangle
   (11.45,+1.5);

 \draw[fill=white, draw opacity=0]
   (-.45, 1.35)
   rectangle
   (11.45, 1.3);

 \draw[fill=white, draw opacity=0]
   (-.45, 1.15)
   rectangle
   (11.45, 1.1);

 \draw[fill=white, draw opacity=0]
   (-.45, .95)
   rectangle
   (11.45, .9);

 \draw[fill=white, draw opacity=0]
   (-.45, -1.35)
   rectangle
   (11.45, -1.3);

 \draw[fill=white, draw opacity=0]
   (-.45, -1.15)
   rectangle
   (11.45, -1.1);

 \draw[fill=white, draw opacity=0]
   (-.45, -.95)
   rectangle
   (11.45, -.9);

 \draw[-Latex]
   (-.7, 1.65)
   to
   (11.8, 1.65);

 \draw[Latex-]
   (-.6, 2.4)
   to
   (-.6, -1.8);

 \draw
 (-.25, 2.15)
 node
 {\scalebox{.85}{$\mathbb{R}^{n-1}$}};

 \draw
   (4, 1.8)
     node
     {\scalebox{.85}{$
       \overset{
       4i
       }{\scalebox{.6}{$\vert$}}
     $} };

 \draw
   (5, 1.8)
     node
     {\scalebox{.85}{$
       \overset{
       4i+1
       }{\scalebox{.6}{$\vert$}}
     $} };

 \draw
   (6, 1.8)
     node
     {\scalebox{.85}{$
       \overset{
       4i+2
       }{\scalebox{.6}{$\vert$}}
     $} };

 \draw
   (7, 1.8)
     node
     {\scalebox{.85}{$
       \overset{
       4i+3
       }{\scalebox{.6}{$\vert$}}
     $} };

 \draw
   (8, 1.8)
     node
     {\scalebox{.85}{$
       \overset{
       4i+4
       }{\scalebox{.6}{$\vert$}}
     $} };

 \draw
   (9, 1.8)
     node
     {\scalebox{.85}{$
       \overset{
       4i+5
       }{\scalebox{.6}{$\vert$}}
     $} };

 \draw
   (11.5, 1.87)
     node
     {\scalebox{.85}{$
       \mathbb{R}
     $} };

 \draw
  (1.5,-1.48)
   node
   {
     \scalebox{1}{$
      \underset{U_{2(i-1)}}{
      \underbrace{\phantom{------------}}
      }
     $}
   };

 \draw
   (5.5,-1.45)
   node
   {
     \scalebox{1}{$
      \underset{
        U_{2i}
      }{
      \underbrace{\phantom{------------}}
      }
     $}
   };

 \draw
   (9.5,-1.48)
   node
   {
     \scalebox{1}{$
      \underset{U_{2(i+1)}}{
      \underbrace{\phantom{------------}}
      }
     $}
   };

 \draw
   (3.5,-1.7)
   node
   {
     \scalebox{1}{$
      \underset{U_{2(i-1)+1}}{
      \underbrace{\phantom{------------}}
      }
     $}
   };

 \draw
   (7.5,-1.7)
   node
   {
     \scalebox{1}{$
      \underset{U_{2i + 1}}{
      \underbrace{\phantom{------------}}
      }
     $}
   };

  \draw (4.5,0) node {
    \scalebox{.7}{$
      U_{2i-1}
        \cap
      U_{2i}
    $}
  };

  \draw (6.5,0) node {
    \scalebox{.7}{$
      U_{2i} \cap U_{2i+1}
    $}
  };

  \draw (5.5,-3) node {$\Gamma$};

\end{tikzpicture}

the resulting square
\begin{equation}
  \label{BrownGerstenPropertyOverCartesianSpaces}
  \begin{tikzcd}
    \mathcal{X}(\mathbb{R}^n)
    \ar[
      r,
      "{
        \big(
          \mathcal{X}(\iota_{2i})
        \big)_{i \in \Integers}
      }"
    ]
    \ar[
      d,
      "{
        \big(
          \mathcal{X}(\iota_{2i+1})
        \big)_{i \in \Integers}
      }"{swap}
    ]
    \ar[
      dr,
      phantom,
      "{\mbox{\tiny\rm(hpb)}}"
    ]
    &
    \underset{i \in \mathbb{Z}}{\prod}
    \mathcal{X}(U_{2i})
    \ar[
      d,
      shorten <=-4pt,
      "{
      }"
    ]
    \\
    \underset{i \in \mathbb{Z}}{\prod}
    \mathcal{X}(U_{2i+1})
    \ar[
      r,
      "{
      }"{swap, yshift=-1pt}
    ]
    &
    \underset{
      \mathclap{
        j \in \Integers
      }
    }{\prod}
    \mathcal{X}(U_{j} \cap U_{j+1})
  \end{tikzcd}
\end{equation}
is a homotopy pullback square.
\end{enumerate}
\end{proposition}
\begin{proof}
  On general grounds the projectively locally fibrant simplicial are those which are objectwise Kan complexes and which satisfy homotopy descent. That over the site of smooth manifolds homotopy descent is equivalent to the above {\it Brown-Gersten type condition} \eqref{BrownGerstenPropertyOverCartesianSpaces}
  follows from \cite[Thm. 4.11]{Pavlov22}.
\end{proof}

This has the following remarkable consequences (see also \cite[Prop. 4.13]{Pavlov22}):
\begin{lemma}[Local fibrancy of delooping groupoids over Cartesian spaces]
\label{LocalFibrancyOfDeloopingGroupoidsOverCartesianSpaces}
   For $\Gamma \,\in\, \Groups\big( \Sheaves(\CartesianSpaces)\big)$ any sheaf of groups over the site $\CartesianSpaces$ \eqref{CartesianSpaces},  the simplicial  nerve
  of its delooping groupoid is locally projectively fibrant:
  $$
    N\DeloopingGroupoid{\Gamma}
    \;\in\;
    \Big(
    \SimplicialPresheaves(\CartesianSpaces)_{\projloc}
    \Big)_{\mathrm{fib}}
    \,.
  $$
\end{lemma}
We thank Dmitri Pavlov for discussion of the following proof.
\begin{proof}
  In the given case,
  the homotopy pullback in \eqref{RecognitionOfFibrantObjectsOverCartSp} is readily seen to be the groupoid whose objects are tuples of {\it transition functions}
  \begin{equation}
    \label{TransitionFunctionOnZigZagCoverByCartesianStrips}
    \Big(
      \gamma_{j, j+1}
      \,\in\,
      \Gamma
      \big(
        (2j+2, 2j+3) \times \mathbb{R}^{n-1}
      \big)
    \Big)_{j \in \mathbb{Z}}
  \end{equation}
  and whose morphisms are {\it gauge transformations}
  between these, of the following form, with composition given by the group operations in the groups $\Gamma(2j, 2j+3)$:
  \begin{equation}
    \label{DescentGroupoidForStripZigZagCover}
    \Big(
    h_{j}
    \,\in\,
    \Gamma
    \big(
      (2j, 2j+3) \times \mathbb{R}^{n-1}
    \big)
    \Big)_{j \in \Integers}
    \;\;:\;\;
    (\gamma_{j, j+1})_{j \in \mathbb{Z}}
    \xrightarrow{\phantom{----}}
    \Big(
      {h_{j}^{-1}}_{\big\vert (2j+2, 2j+3)}
        \cdot
      \gamma_{j, j+1}
        \cdot
      {h_{j+1}}_{\big\vert (2j+2,2j+3)}
    \Big)_{j \in \Integers}
    \,.
  \end{equation}

\hspace{-.8cm}
\begin{tikzpicture}[xscale=1.37]

  \clip
    (-.9,2.6)
    rectangle
    (12,-4);

 \draw[fill=gray, fill opacity=.3, draw opacity=0]
   (-.45,-1.5)
   rectangle
   (-.35,+1.5);

 \draw[fill=gray, fill opacity=.3, draw opacity=0]
   (-.3,-1.5)
   rectangle
   (-.2,+1.5);

 \draw[fill=gray, fill opacity=.3, draw opacity=0]
   (-.15,-1.5)
   rectangle
   (-.05,+1.5);

 \draw[fill=gray, fill opacity=.3, draw opacity=0]
   (0,-1.5)
   rectangle
   (3,+1.5);

 \draw[fill=gray, fill opacity=.3, draw opacity=0]
   (2,-1.5)
   rectangle
   (5,+1.5);

 \draw[fill=gray, fill opacity=.3, draw opacity=0]
   (4,-1.5)
   rectangle
   (7,+1.5);

 \draw[fill=gray, fill opacity=.3, draw opacity=0]
   (6,-1.5)
   rectangle
   (9,+1.5);

 \draw[fill=gray, fill opacity=.3, draw opacity=0]
   (8,-1.5)
   rectangle
   (11,+1.5);

 \draw[fill=gray, fill opacity=.3, draw opacity=0]
   (11.05,-1.5)
   rectangle
   (11.15,+1.5);

 \draw[fill=gray, fill opacity=.3, draw opacity=0]
   (11.2,-1.5)
   rectangle
   (11.3,+1.5);

 \draw[fill=gray, fill opacity=.3, draw opacity=0]
   (11.35,-1.5)
   rectangle
   (11.45,+1.5);

 \draw[fill=white, draw opacity=0]
   (-.45, 1.35)
   rectangle
   (11.45, 1.3);

 \draw[fill=white, draw opacity=0]
   (-.45, 1.15)
   rectangle
   (11.45, 1.1);

 \draw[fill=white, draw opacity=0]
   (-.45, .95)
   rectangle
   (11.45, .9);

 \draw[fill=white, draw opacity=0]
   (-.45, -1.35)
   rectangle
   (11.45, -1.3);

 \draw[fill=white, draw opacity=0]
   (-.45, -1.15)
   rectangle
   (11.45, -1.1);

 \draw[fill=white, draw opacity=0]
   (-.45, -.95)
   rectangle
   (11.45, -.9);

 \draw[-Latex]
   (-.7, 1.65)
   to
   (11.8, 1.65);

 \draw[Latex-]
   (-.6, 2.4)
   to
   (-.6, -1.8);

 \draw
 (-.25, 2.15)
 node
 {\scalebox{.85}{$\mathbb{R}^{n-1}$}};

 \draw
   (4, 1.8)
     node
     {\scalebox{.85}{$
       \overset{
       4i
       }{\scalebox{.6}{$\vert$}}
     $} };

 \draw
   (5, 1.8)
     node
     {\scalebox{.85}{$
       \overset{
       4i+1
       }{\scalebox{.6}{$\vert$}}
     $} };

 \draw
   (6, 1.8)
     node
     {\scalebox{.85}{$
       \overset{
       4i+2
       }{\scalebox{.6}{$\vert$}}
     $} };

 \draw
   (7, 1.8)
     node
     {\scalebox{.85}{$
       \overset{
       4i+3
       }{\scalebox{.6}{$\vert$}}
     $} };

 \draw
   (8, 1.8)
     node
     {\scalebox{.85}{$
       \overset{
       4i+4
       }{\scalebox{.6}{$\vert$}}
     $} };

 \draw
   (9, 1.8)
     node
     {\scalebox{.85}{$
       \overset{
       4i+5
       }{\scalebox{.6}{$\vert$}}
     $} };

 \draw
   (11.5, 1.87)
     node
     {\scalebox{.85}{$
       \mathbb{R}
     $} };

 \draw
  (1.5,-1.49)
   node
   {
     \scalebox{1}{$
      \underset{U_{2(i-1)}}{
      \underbrace{\phantom{------------}}
      }
     $}
   };

 \draw
   (5.5,-1.47)
   node
   {
     \scalebox{1}{$
      \underset{
        U_{2i}
      }{
      \underbrace{\phantom{------------}}
      }
     $}
   };

 \draw
   (9.5,-1.49)
   node
   {
     \scalebox{1}{$
      \underset{U_{2(i+1)}}{
      \underbrace{\phantom{------------}}
      }
     $}
   };

 \draw
   (3.5,-1.7)
   node
   {
     \scalebox{1}{$
      \underset{U_{2(i-1)+1}}{
      \underbrace{\phantom{------------}}
      }
     $}
   };

 \draw
   (7.5,-1.7)
   node
   {
     \scalebox{1}{$
      \underset{U_{2i + 1}}{
      \underbrace{\phantom{------------}}
      }
     $}
   };

  \draw (4.5,0) node {
    \scalebox{.7}{$
      U_{2i-1}
        \cap
      U_{2i}
    $}
  };

  \draw (6.5,0) node {
    \scalebox{.7}{$
      U_{2i} \cap U_{2i+1}
    $}
  };

  \draw (6.2,-3.8) node {$\Gamma$};

  \draw[white, line width=3pt, bend right=20]
    (4.5,-.2)
    to
    (5.9,-3.6);
  \draw[->, bend right=20]
    (4.5,-.2)
    to
    node[near end]
     {\scalebox{.8}{\fcolorbox{white}{white}{$
       \gamma_{2i-1, 2i}
     $}}}
    (5.9,-3.6);

  \draw[white, line width=3pt, bend left=10]
    (6.5,-.2)
    to
    (6.3,-3.4);
  \draw[->, bend left=10]
    (6.5,-.2)
    to
    node[near end]
     {\scalebox{.8}{\fcolorbox{white}{white}{$
       \gamma_{2i, 2i+1}
     $}}}
    (6.3,-3.4);

  \draw[->, bend right=10]
    (3.6,-2.2)
    to
    node
     {\scalebox{.8}{\fcolorbox{white}{white}{$
       h_{2i-1}
     $}}}
    (5.76,-3.8);

  \draw[->, bend left=5]
    (5.6,-2)
    to
    node
     {\scalebox{.8}{\fcolorbox{white}{white}{$
       h_{2i}
     $}}}
    (6.1,-3.4);

  \draw[->, bend left=10]
    (7.4,-2.2)
    to
    node
     {\scalebox{.8}{\fcolorbox{white}{white}{$
       h_{2i+1}
     $}}}
    (6.5,-3.6);

\end{tikzpicture}

\vspace{.1cm}

  We need to show that this groupoid
  \eqref{DescentGroupoidForStripZigZagCover}
  is equivalent to $\Gamma(\mathbb{R}^n) \rightrightarrows \ast$, which means to show that (i) it is
  connected and (ii) that the automorphisms of the element with trivial transition functions $(\NeutralElement)_{j \in \Integers}$
  are in bijection to $\Gamma(\mathbb{R}^n)$.

  The second statement (ii) is clear: For a gauge transformation
  \eqref{DescentGroupoidForStripZigZagCover}
  to preserve the trivial transition functions $(\gamma_{j, j+1} = \NeutralElement)_{j \in \Integers}$
  its components
  $(h_{j})_{j \in \Integers}$
  must agree on intersections,
  hence they glue by the assumed sheaf property of $\Gamma$.

  For the first statement (i), the strategy is to use gauge transformations ``on the left'' to gauge away transition functions ``on the right of their domain'', and vice versa. The subtlety is to do this gauging ``on the left'' without destroying the trivial gauge already achieved ``on the right''.

  \medskip

  We close by spelling out in detail one way to do this. So let $(\gamma_{j, j+1})_{j \in \Integers}$ be any tuple of transition functions \eqref{TransitionFunctionOnZigZagCoverByCartesianStrips}. We will produce a morphism \eqref{DescentGroupoidForStripZigZagCover} of the form
  $(\gamma_{j, j+1})_{j \in \Integers} \xrightarrow{\;\;} (\NeutralElement)_{j \in \Integers}$.

  \medskip

  To this end, choose diffeomorphisms
  (by the usual methods for bump functions, e.g. as in \cite[Prop. 2.25]{Lee12}) of this form:
  $$
    \begin{tikzcd}[
      row sep=12pt,
      column sep=-5pt
    ]
      \big(
        2j, 2j+2
      \big)
      \ar[dr, hook]
      \ar[rr, "{\sim}"{swap}]
      &[-00pt]&[-63pt]
      \big(
        2j+2, 2j+2\frac{1}{4}
      \big)
      \ar[dr, hook]
      &
      \big(
        2j + 2,
        2j + 2\tfrac{1}{2}
      \big)
      \ar[rr, Rightarrow, -]
      \ar[d, hook]
      &&[-63pt]
      \big(
        2j + 2,
        2j + 2\tfrac{1}{2}
      \big)
      \ar[d, hook]
      \\
      &
      \mathllap{
        U_j
        \,=\;
      }
      \big(2j, 2j+3\big)
        \ar[
          rr,
          "{\ell_j}",
          "{\sim}"{swap}
        ]
        &&
      \big(2j+2, 2j+3\big)
      \ar[
        from=rr,
        "r_{j+1}"{swap},
        "\sim"
      ]
      &&
      \big(
        2j+2, 2j+5
      \big)
      \mathrlap{
        \;=\,
        U_{j+1}
      }
      \\
      &
      \big(2j+2\tfrac{1}{2}, 2j+3\big)
      \ar[u, hook]
      \ar[rr, Rightarrow,-]
      &&
      \big(2j+2\tfrac{1}{2}, 2j+3\big)
      \ar[u, hook]
      &
      \big(
        2j+2\tfrac{3}{4},
        2j+3
      \big)
      \ar[
        ul,
        hook
      ]
      \ar[
        from=rr,
        "{\sim}"
      ]
      &&
      \big(
        2j+3,
        2j+5
      \big)
      \ar[
        ul,
      ]
    \end{tikzcd}
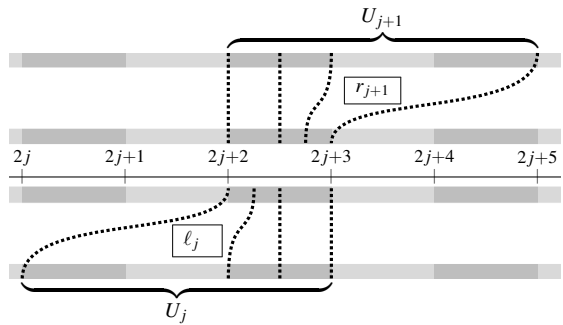
$$
\begin{center}
\begin{tikzpicture}[xscale=1.37]

  \clip
    (+3.88,2.6)
    rectangle
    (9.22,-2.3);

 \draw[fill=gray, fill opacity=.3, draw opacity=0]
   (-.45,-1.5)
   rectangle
   (-.35,+1.5);

 \draw[fill=gray, fill opacity=.3, draw opacity=0]
   (-.3,-1.5)
   rectangle
   (-.2,+1.5);

 \draw[fill=gray, fill opacity=.3, draw opacity=0]
   (-.15,-1.5)
   rectangle
   (-.05,+1.5);

 \draw[fill=gray, fill opacity=.3, draw opacity=0]
   (0,-1.5)
   rectangle
   (3,+1.5);

 \draw[fill=gray, fill opacity=.3, draw opacity=0]
   (2,-1.5)
   rectangle
   (5,+1.5);

 \draw[fill=gray, fill opacity=.3, draw opacity=0]
   (4,-1.5)
   rectangle
   (7,+1.5);

 \draw[fill=gray, fill opacity=.3, draw opacity=0]
   (6,-1.5)
   rectangle
   (9,+1.5);

 \draw[fill=gray, fill opacity=.3, draw opacity=0]
   (8,-1.5)
   rectangle
   (11,+1.5);

 \draw[fill=gray, fill opacity=.3, draw opacity=0]
   (11.05,-1.5)
   rectangle
   (11.15,+1.5);

 \draw[fill=gray, fill opacity=.3, draw opacity=0]
   (11.2,-1.5)
   rectangle
   (11.3,+1.5);

 \draw[fill=gray, fill opacity=.3, draw opacity=0]
   (11.35,-1.5)
   rectangle
   (11.45,+1.5);

 \draw[draw=white, fill=white]
   (0,1.3)
   rectangle
   (10,.5);

 \draw[draw=white, fill=white]
   (0,.28)
   rectangle
   (10,-.27);

 \draw[draw=white, fill=white]
   (0,-1.3)
   rectangle
   (10,-.5);

\begin{scope}[yshift=-1.8cm]

 \draw[-Latex]
   (-.7, 1.65)
   to
   (11.8, 1.65);

 \draw
   (4, 1.8)
     node
     {\scalebox{.85}{$
       \overset{
       2j
       }{\scalebox{.6}{$\vert$}}
     $} };

 \draw
   (5, 1.8)
     node
     {\scalebox{.85}{$
       \overset{
       2j+1
       }{\scalebox{.6}{$\vert$}}
     $} };

 \draw
   (6, 1.8)
     node
     {\scalebox{.85}{$
       \overset{
       2j+2
       }{\scalebox{.6}{$\vert$}}
     $} };

 \draw
   (7, 1.8)
     node
     {\scalebox{.85}{$
       \overset{
       2j+3
       }{\scalebox{.6}{$\vert$}}
     $} };

 \draw
   (8, 1.8)
     node
     {\scalebox{.85}{$
       \overset{
       2j + 4
       }{\scalebox{.6}{$\vert$}}
     $} };

 \draw
   (9, 1.8)
     node
     {\scalebox{.85}{$
       \overset{
       2j + 5
       }{\scalebox{.6}{$\vert$}}
     $} };

\end{scope}

 \draw
   (5.5,-1.7)
   node
   {
     \scalebox{1}{$
      \underset{
        U_{j}
      }{
      \underbrace{\phantom{------------}}
      }
     $}
   };

 \draw
   (7.5,+1.7)
   node
   {
     \scalebox{1}{$
      \overset{U_{j + 1}}{
      \overbrace{\phantom{------------}}
      }
     $}
   };

  \draw
    [line width=1.2pt, densely dotted]
    (4,-1.5)
    .. controls (4,-.6) and (6,-1) ..
    (6,-.3);

  \draw
    [line width=1.2pt, densely dotted]
    (6,-1.5)
    .. controls (6,-.6) and (6.25, -1) ..
    (6.25,-.3);

  \draw
    [line width=1.2pt, densely dotted]
    (6.5,-1.5)
    to
    (6.5,-.3);

  \draw
    [line width=1.2pt, densely dotted]
    (7,-1.5)
    to
    (7,-.3);

\begin{scope}[scale=-1, xshift=-13cm]

  \draw
    [line width=1.2pt, densely dotted]
    (4,-1.5)
    .. controls (4,-.6) and (6,-1) ..
    (6,-.3);

  \draw
    [line width=1.2pt, densely dotted]
    (6,-1.5)
    .. controls (6,-.6) and (6.25, -1) ..
    (6.25,-.3);

  \draw
    [line width=1.2pt, densely dotted]
    (6.5,-1.5)
    to
    (6.5,-.3);

  \draw
    [line width=1.2pt, densely dotted]
    (7,-1.5)
    to
    (7,-.3);

\end{scope}

  \draw
    (5.7,-1) node{
      \scalebox{.7}{
      \fcolorbox{black}{white}{
        $\ell_{j\phantom{+}}$
      }
      }
    };

  \draw
    (7.4,+1) node{
      \scalebox{.7}{
      \fcolorbox{black}{white}{
        $r_{j+1}$
      }
      }
    };
\end{tikzpicture}
\end{center}

Now, to start with, consider the morphism
$$
  \left(
  h_j
  \,:=\,
  \left\{
    \def\arraystretch{1.1}
    \begin{array}{cll}
      \ell_0^\ast(\gamma_{01})
      &\vert& j = 0
      \\
      \NeutralElement &\vert& \mbox{otherwise}
    \end{array}\right.
  \right)_{\!\!\!\! j \in \Integers}
  \;\;\;
  :\;\;\;
  (\gamma_{j, j+1})_{j \in \Integers}
  \xrightarrow{\phantom{---}}
  \big(
    \gamma^{\,\prime}_{j, j+1}
  \big)_{j \in \Integers}
$$
to new transition functions which satisfy
$$
  {\gamma^{\,\prime}_{01}}
    _{\big\vert \big(2\tfrac{1}{2}, 3\big)}
  \,=\,
  \NeutralElement
  \,;
$$
and from there the further morphism
$$
  \left(
  h_j
  \,:=\,
  \left\{
    \def\arraystretch{1.1}
    \begin{array}{cll}
      r_1^\ast(\gamma^{\,\prime}_{01})^{-1}
      &\vert& j = 1
      \\
      \NeutralElement &\vert& \mbox{otherwise}
    \end{array}\right.
  \right)_{\!\!\!\! j \in \Integers}
  \;\;\;
  :
  \;\;\;
  (\gamma^{\,\prime}_{j, j+1})_{j \in \Integers}
  \xrightarrow{\phantom{---}}
  \big(
    \gamma^{\,\prime\prime}_{j, j+1}
  \big)_{j \in \Integers}
$$
to new transition functions which satisfy
$$
  {\gamma^{\,\prime\prime}_{01}}
  \,=\,
  \NeutralElement
  \,.
$$

Next, consider the following morphism, whose components to the right are all trivial, while those to the left are defined inductively:
$$
  \left(
  h_j
  \,:=\,
  \left\{
    \def\arraystretch{1.1}
    \begin{array}{cll}
      \NeutralElement
      &\vert& j \geq 0
      \\
      \ell_{j}^\ast
      \Big(
        \gamma^{\,\prime\prime}_{j, j+1}
        \cdot
        \big(
          h_{j+1}
        \big)_{\!\big\vert (2j+2, 2j+3)}
      \Big)
      &\vert&
      j < 0
    \end{array}\right.
  \right)_{\!\!\!\! j \in \Integers}
  \;\;\;
  :
  \;\;\;
  (\gamma^{\,\prime\prime}_{j, j+1})_{j \in \Integers}
  \xrightarrow{\phantom{---}}
  \big(
    \gamma^{\,\prime\prime\prime}_{j, j+1}
  \big)_{j \in \Integers}
  \,.
$$
This yields transition functions with the property
$$
  \underset{
    \scalebox{.7}{$
      \begin{array}{c}
        j \in \Integers
        \\
        j \leq 0
      \end{array}
    $}
  }{\forall}
  \;\;\;\;\;
  {\gamma^{\,\prime\prime\prime}_{j, j+1}}
    _{\big\vert \big(2j + 2\tfrac{1}{2}, 2j + 3  \big)}
  \;=\;
  \NeutralElement
  \,;
$$
and from there we have a morphism
$$
  \left(
  h_j
  \,:=\,
  \left\{
    \def\arraystretch{1.1}
    \begin{array}{cll}
      \NeutralElement
      &\vert& j \geq 0
      \\
      r_{j+1}^\ast
      \big(
        \gamma^{\,\prime\prime\prime}_{j, j+1}
      \big)^{-1}
      &\vert&
      j < 0
    \end{array}\right.
  \right)_{\!\!\!\! j \in \Integers}
  \;\;\;
  :
  \;\;\;
  (\gamma^{\,\prime\prime\prime}_{j, j+1})_{j \in \Integers}
  \xrightarrow{\phantom{---}}
  \big(
    \gamma^{\,\prime\prime\prime\prime}_{j, j+1}
  \big)_{j \in \Integers}
$$
to transition functions with the property
$$
  \underset{
    \scalebox{.7}{$
      \begin{array}{c}
        j \in \NaturalNumbers
        \\
        j \leq 0
      \end{array}
    $}
  }{\forall}
  \;\;\;\;\;
  {\gamma^{\,\prime\prime\prime\prime}_{j, j+1}}
  \;=\;
  \NeutralElement
  \,.
$$

Finally, applying the analogous inductive construction also
mirror-symmetrically for $j \to + \infty$
serves to trivialize also the remaining transition functions to the right, giving a morphism
$(\gamma^{\,\prime\prime\prime\prime}_{j, j+1})_{j \in \Integers}
\xrightarrow{\;\;}
(\NeutralElement)_{j \in \Integers}
$. The composite morphism
$(\gamma_{j, j+1})_{j \in \Integers}  \xrightarrow{\;}
(\gamma^{\,\prime}_{j, j+1})_{j \in \Integers}  \xrightarrow{\;}
(\gamma^{\,\prime\prime}_{j, j+1})_{j \in \Integers}  \xrightarrow{\;}
(\gamma^{\,\prime\prime\prime}_{j, j+1})_{j \in \Integers}  \xrightarrow{\;}
(\gamma^{\,\prime\prime\prime\prime}_{j, j+1})_{j \in \Integers}  \xrightarrow{\;}
(\NeutralElement)_{j \in \Integers}$ is the desired one, exhibiting the connectedness of the descent groupoid.
\end{proof}
More generally:
\begin{proposition}[Local fibrancy of simplicial deloopings of locally fibrant presheaves of simplicial groups]
\label{LocalFibrancyOfSimplicialDeloopingsOfLocallyFibrantSimplicialGroups}
 Given a presheaf of simplicial groups whose underlying simplicial presheaf is locally projectively fibrant
 (Ntn. \ref{ModelCategoriesOfSimplicialPresheaves})
  then its (objectwise) simplical delooping
  \eqref{StandardSimplicialDeloopingAsQuotient}
  is itself locally projectively fibrant:
  $$
    \mathcal{G} \,\in\,
    \Groups
    \bigg(
      \Big(
      \SimplicialPresheaves(\CartesianSpaces)_{\projloc}
      \Big)_{\mathrm{fib}}
    \bigg)
    \;\;\;\;\;\;\;\;\;\;
    \Rightarrow
    \;\;\;\;\;\;\;\;\;\;
    \overline{W}\mathcal{G}
    \;\;\;
    \in
    \;
    \Big(
      \SimplicialPresheaves(\CartesianSpaces)_{\projloc}
    \Big)_{\mathrm{fib}}
    \,.
  $$
\end{proposition}
We again thank Dmitri Pavlov for discussion of the following proof.
\begin{proof}
  Since $\overline{W}(-)$ takes objectwise values in Kan complexes (Prop. \ref{BasicPropertiesOfStandardSimplicialPrincipalComplex}) it is immediate that $\overline{W}\mathcal{G}$ is globally projectively fibrant; what we have to show is that it also satisfies homotopy descent. For this, by Prop. \ref{RecognitionOfFibrantObjectsOverCartSp}, it is sufficient that for a good open cover
  $\big\{ U_j \xhookrightarrow{\;\iota_j\;} \mathbb{R}^n \big\}_{j \in \Integers}$
  of the special form \eqref{StripZigZagCoverOfCartesianSpace} the descent morphism

  \vspace{-.2cm}
  \begin{equation}
    \label{DescentMorphismForDeloopingSimplicialPresheafOnCartesianZigZagCover}
    \overline{W}\mathcal{G}
    \big(
      \mathbb{R}^n
    \big)
    \,\simeq\,
    \SimplicialPresheaves(\CartesianSpaces)
    \big(
      \mathbb{R}^n
      ,\,
      \overline{W}\mathcal{G}
    \big)
    \xrightarrow
      [\in \WeakHomotopyEquivalences]
      {\;
        \SimplicialPresheaves(\CartesianSpaces)
        \big(
          \iota_{\mathcal{U}}
          ,\,
          \overline{W}\mathcal{G}
        \big)
      \;}
    \SimplicialPresheaves(\CartesianSpaces)
    \big(
      \mathcal{U}_\bullet
      ,\,
      \overline{W}\mathcal{G}
    \big)
  \end{equation}
  \vspace{-.2cm}

  \noindent
  is a weak homotopy equivalence, where

  \vspace{-.2cm}
  $$
    \mathcal{U}_\bullet
    \;:=\;
    \Big(
      \coprod_j U_j
    \Big)^{\times^\bullet_{\mathbb{R}^n}}
    \xrightarrow{\;\;\;
      \iota_{\mathcal{U}}
    \;\;\;}
    \mathbb{R}^n
    \;\;\;
    \in
    \;
    \SimplicialPresheaves(\CartesianSpaces)
  $$
  \vspace{-.2cm}

  \noindent
  denotes the {\v C}ech nerve of the cover. Here the special property of the cover \eqref{StripZigZagCoverOfCartesianSpace} implies that it has no non-trivial intersections beyond double intersections, hence that its {\v C}ech nerve is fixed
  by the $k$-skeleton functor
  \eqref{CoSkeletaAdjunction}
  for all $k \geq 1$ :

  \vspace{-.2cm}
  \begin{equation}
    \label{CechNerveOfZigZagCoverIsOneSkeletal}
    \mathrm{sk}_{k \geq 1}
    \big(
      \mathcal{U}_\bullet
    \big)
    \;\simeq\;
    \mathcal{U}_\bullet
    \;\;\;
    \in
    \;
    \SimplicialPresheaves(\CartesianSpaces)
    \,.
  \end{equation}
  \vspace{-.3cm}

  \noindent
  Notice here, by the special properties of Cartesian products of simplicial simplices (e.g. \cite[\S B.6]{Friedman20}), that the skeletal degree is {\it additive} in a Cartesian product of simplicial sets, so that with $\mathrm{sk}_1(\Simplex{1}) \,\simeq\, \Simplex{1}$ the above property \eqref{CechNerveOfZigZagCoverIsOneSkeletal}
  gives:
  \begin{equation}
    \label{SkeletalDegreeOfProductOfCechNerveOfZigZagCoverwithOneSimplex}
    \mathrm{sk}_2
    \big(
      \mathcal{U}^{\times^\bullet_{\mathbb{R}^n}}
      \times
      \Simplex{1}
    \big)
    \;\simeq\;
      \mathcal{U}^{\times^\bullet_{\mathbb{R}^n}}
      \times
      \Simplex{1}
    \;\;\;
    \in
    \;
    \SimplicialPresheaves(\CartesianSpaces)
    \,.
  \end{equation}

  Now to see that \eqref{DescentMorphismForDeloopingSimplicialPresheafOnCartesianZigZagCover} is a weak equivalence, hence an isomorphism on all $\pi_n$, we will show that:

  \noindent
  (i) it is an isomorphism on $\pi_0$, in that, with the left hand side manifestly being connected, also the right hand side is connected;

  \noindent
  (ii) it is an isomorphism on $\pi_{\bullet \geq 1}$, in that its looping \eqref{LoopingInIntroduction} is a weak homotopy equivalence.

  \medskip

  \noindent
  Here (i) holds by the following sequence of equivalences:

  \vspace{-.7cm}
  $$
    \def\arraystretch{1.8}
    \begin{array}{ll}
      \tau_0
      \,
      \SimplicialPresheaves(\CartesianSpaces)
      \big(
        \mathcal{U}^{\times^\bullet_{\mathbb{R}^n}}
        ,\,
        \overline{W}\mathcal{G}
      \big)
      \\
      \;\simeq\;
      \SimplicialPresheaves(\CartesianSpaces)
      \big(
        \mathcal{U}^{\times^\bullet_{\mathbb{R}^n}}
        ,\,
        \overline{W}\mathcal{G}
      \big)_0
      \Big/
      \SimplicialPresheaves(\CartesianSpaces)
      \big(
        \mathcal{U}^{\times^\bullet_{\mathbb{R}^n}}
        \times
        \Simplex{1}
        ,\,
        \overline{W}\mathcal{G}
      \big)_0
      \\
      \;\simeq\;
      \SimplicialPresheaves(\CartesianSpaces)
      \Big(
        \mathrm{sk}_2
        \big(
          \mathcal{U}^{\times^\bullet_{\mathbb{R}^n}}
        \big)
        ,\,
        \overline{W}\mathcal{G}
      \Big)_0
      \Big/
      \SimplicialPresheaves(\CartesianSpaces)
      \Big(
        \mathrm{sk}_2
        \big(
        \mathcal{U}^{\times^\bullet_{\mathbb{R}^n}}
        \times
        \Simplex{1}
        \big)
        ,\,
        \overline{W}\mathcal{G}
      \Big)_0
      &
      \proofstep{
        by
        \eqref{CechNerveOfZigZagCoverIsOneSkeletal}
        \&
        \eqref{SkeletalDegreeOfProductOfCechNerveOfZigZagCoverwithOneSimplex}
      }
      \\
      \;\simeq\;
      \SimplicialPresheaves(\CartesianSpaces)
      \Big(
          \mathcal{U}^{\times^\bullet_{\mathbb{R}^n}}
        ,\,
        \mathrm{cosk}_2
        \big(
          \overline{W}\mathcal{G}
        \big)
      \Big)_0
      \Big/
      \SimplicialPresheaves(\CartesianSpaces)
      \Big(
        \mathcal{U}^{\times^\bullet_{\mathbb{R}^n}}
        \times
        \Simplex{1}
        ,\,
        \mathrm{cosk}_2
        \big(
        \overline{W}\mathcal{G}
        \big)
      \Big)_0
      &
      \proofstep{
        by
        \eqref{CoSkeletaAdjunction}
      }
      \\
      \;\simeq\;
      \tau_0
      \,
      \SimplicialPresheaves(\CartesianSpaces)
      \Big(
          \mathcal{U}^{\times^\bullet_{\mathbb{R}^n}}
        ,\,
        \mathrm{cosk}_2
        \big(
          \overline{W}\mathcal{G}
        \big)
      \Big)
      \\
      \;\simeq\;
      \tau_0
      \,
      \SimplicialPresheaves(\CartesianSpaces)
      \Big(
          \mathcal{U}^{\times^\bullet_{\mathbb{R}^n}}
        ,\,
          \overline{W}
          \big(
            \mathrm{cosk}_1
            (\mathcal{G})
          \big)
      \Big)
      &
      \proofstep{
        by
        Ex. \ref{CoprojectionsOutOfBorelConstructionAreKanFibrations}
        \&
        Lem. \ref{TruncationModeledBySimplicialCoskeleta}
      }
      \\
      \;\simeq\;
      \ast
      &
      \proofstep{
        by Prop. \ref{LocalFibrancyOfDeloopingGroupoidsOverCartesianSpaces}
      }
    \end{array}
  $$
  \vspace{-.5cm}

  \noindent
  The last step uses that $\mathrm{cosk}_{1}(\mathcal{G})$ is equivalent to a presheaf of plain groups to which the previous
  Prop. \ref{LocalFibrancyOfSimplicialDeloopingsOfLocallyFibrantSimplicialGroups} applies.

  Finally, to see that (ii) holds:
  Due to the natural weak homotopy equivalence  $\Omega \overline{W}\mathcal{G} \;\simeq\; \mathcal{G}$
  (Ex. \ref{CoprojectionsOutOfBorelConstructionAreKanFibrations})
  and using that the simplicial hom-functor preserves looping of Kan complexes (like all limits of simplicial sets) in its second argument,
  the image under looping
  \eqref{LoopingInIntroduction}
  of the descent morphism \eqref{DescentMorphismForDeloopingSimplicialPresheafOnCartesianZigZagCover}
  is
  (up to globally projective weak equivalence)
  of this form

  \vspace{-.3cm}
  $$
    \mathcal{G}
    \big(
      \mathbb{R}^n
    \big)
    \,\simeq\,
    \SimplicialPresheaves(\CartesianSpaces)
    \big(
      \mathbb{R}^n
      ,\,
      \mathcal{G}
    \big)
    \xrightarrow
      [\in \WeakHomotopyEquivalences]
      {\;
        \SimplicialPresheaves(\CartesianSpaces)
        \big(
          \iota_{\mathcal{U}}
          ,\,
          \mathcal{G}
        \big)
      \;}
    \SimplicialPresheaves(\CartesianSpaces)
    \big(
      \mathcal{U}_\bullet
      ,\,
      \mathcal{G}
    \big)
    \,.
  $$
  \vspace{-.2cm}

  \noindent
  But this is just the descent morphism for $\mathcal{G}$ itself, which is a weak equivalence by the assumption that $\mathcal{G}$ is locally projectively fibrant.
\end{proof}

\medskip

\noindent
{\bf Topological transformation groups seen in $\SmoothInfinityGroupoids$.} Before discussing the cohesion of smoooth $\infty$-groupoids (Prop.
\ref{SmoothCohesion} below), we showcase some aspects of topological transformation groups seen under the embedding of D-topological spaces into smooth $\infty$-groupoids:

\begin{example}[Stabilizer groups as the loop objects of homotopy quotient stacks]
  \label{StabilizerGroupsAsLoopObjects}
  Given a D-topological group action

  - $\Gamma \,\in\, \Groups(\DTopologicalSpaces)
    \hookrightarrow{\;}
    \Groups(\SmoothInfinityGroupoids)
  $

  -
  $
    \Gamma \acts \TopologicalSpace
    \,\in\,
    \Actions{\Gamma}
    \big(
      \DTopologicalSpaces
    \big)
    \xhookrightarrow{\;}
    \Actions{\Gamma}
    \big(
      \SmoothInfinityGroupoids
    \big)
  $

  \noindent
  seen in $\SmoothInfinityGroupoids$
  under \eqref{CategoryOfDeltaGeneratedSpaces},
  then the  simplicial nerve of their
  topological action groupoid (Exp. \ref{TopologicalActionGroupoid})

  $
    N
    \big(
      \ActionGroupoid{\TopologicalSpace}{\Gamma}
    \big)
    \;\in\;
    \simplicial\DTopologicalSpaces
    \xhookrightarrow{\;}
    \simplicial\DiffeologicalSpaces
    \xhookrightarrow{\;}
    \SimplicialPresheaves(\CartesianSpaces)
  $

  \noindent
  represents the
  $\infty$-topos theoretic homotopy quotient
  (according to
  Prop. \ref{HomotopyQuotientsAndPrincipaInfinityBundles})
  $$
    \SimplicialLocalization{\ProjectiveWeakEquivalences}
    \Big(
    N
    \big(
      \ActionGroupoid{\TopologicalSpace}{\Gamma}
    \big)
    \Big)
    \;\simeq\;
    \HomotopyQuotient
      { \TopologicalSpace }
      { \Gamma }
    \;\;\;
    \in
    \;
    \SmoothInfinityGroupoids
    \,.
  $$

  Moreover, for any $x \in \TopologicalSpace$
  its topological stabilizer group is canonically identified with the looping \eqref{LoopingInIntroduction} of the homotopy quotient at that point:
  \begin{equation}
    \label{TopologicalStabilizerGroupAsLooping}
    \begin{tikzcd}
      \mathrm{Stab}_\Gamma(x)
      \ar[rr]
      \ar[d]
      \ar[
        drr,
        phantom,
        "{\mbox{\tiny\rm(pb)}}"{pos=.37}
      ]
      &&
      \ast
      \ar[
        d,
        "{\vdash x}"{}
      ]
      \\
      \ast
      \ar[
        rr,
        "{\vdash x}"{swap}
       ]
      &&
      \HomotopyQuotient{\TopologicalSpace}{\Gamma}
    \end{tikzcd}
    \;\;\;\;\;
    \in
    \;
    \SmoothInfinityGroupoids
    \,.
  \end{equation}
\end{example}
\begin{proof}
  For the first statement, it is immediate to see that we have a pullback diagram of simplicial presheaves as shown on the left below,
  which hence gives the required homotopy pullback of $\infty$-stacks shown on the right, due to Lem. \ref{ComputingHomotopyPullbacksOfInfinityStacks}:
  $$
    \begin{tikzcd}
      N
      (
        \TopologicalSpace
        \rightrightarrows
        \TopologicalSpace
      )
      \ar[rr]
      \ar[d]
      \ar[
        drr,
        phantom,
        "\mbox{\tiny\rm(pb)}"{pos=.4}
      ]
      &&
      N
      \ActionGroupoid{\TopologicalSpace}{\Gamma}
      \ar[
        d,
        "{
          \in \ProjectiveFibrations
        }"
      ]
      \\
      N(\ast \rightrightarrows \ast)
      \ar[
        rr
      ]
      &&
      N\DeloopingGroupoid{\Gamma}
    \end{tikzcd}
    \hspace{1cm}
    \overset{
      \mathrlap{\SimplicialLocalization{\ProjectiveWeakEquivalences}}
    }{\longmapsto}
    \hspace{1cm}
    \begin{tikzcd}
      \TopologicalSpace
      \ar[rr]
      \ar[d]
      \ar[
        drr,
        phantom,
        "\mbox{\tiny\rm(pb)}"{pos=.53}
      ]
      &&
      \HomotopyQuotient{\TopologicalSpace}{\Gamma}
      \ar[
        d
      ]
      \\
      \ast
      \ar[
        rr
      ]
      &&
      \mathbf{B}\Gamma
    \end{tikzcd}
  $$
Similarly, for computing the looping,
the same Lemma \ref{ComputingHomotopyPullbacksOfInfinityStacks}
(together with direct inspection of the right morphism, or else appeal to the factorization lemma, \ref{FactorizationLemmaForProjectiveSimplicialPresheaves}) gives:
$$
  \begin{tikzcd}[column sep=6pt]
    &&
    N(\ast \rightrightarrows \ast)
    \ar[
      d,
      "{\vdash x}"{swap},
      "{\in \ProjectiveWeakEquivalences}"
    ]
    \\
    N
    \big(
      \mathrm{Stab}_\Gamma(x)
      \rightrightarrows
      \mathrm{Stab}_\Gamma(x)
    \big)
    \ar[d]
    \ar[rr]
    \ar[
      drr,
      phantom,
      "{
        \mbox{\tiny\rm(pb)}
      }"{pos=.4}
    ]
    &&
    N
    \big(
      (\{x\} \times \Gamma)
      \times
      \Gamma
      \rightrightarrows
      (\{x\} \times \Gamma)
    \big)
    \ar[
      d,
      "{
        \in \ProjectiveFibrations
      }"
    ]
    \\
    N(\ast \rightrightarrows \ast)
    \ar[
      rr,
      "{\vdash x}"{swap}
    ]
    &&
    \ActionGroupoid{\TopologicalSpace}{\Gamma}
  \end{tikzcd}
  \hspace{.7cm}
  \raisebox{-20pt}{$
    \overset{
      \mathrlap{\SimplicialLocalization{\ProjectiveWeakEquivalences}}
    }{\longmapsto}
  $}
  \hspace{.7cm}
  \begin{tikzcd}
    &&
    \phantom{N}
    \\
    \mathrm{Stab}_\Gamma(x)
    \ar[rr]
    \ar[d]
    \ar[
      drr,
      phantom,
      "{
        \mbox{\tiny\rm(pb)}
      }"
    ]
    &&
    \ast
    \ar[
      d,
      "{\vdash x}"
    ]
    \\
    \ast
    \ar[
      rr,
      "{\vdash x}"{swap}
    ]
    &&
    \HomotopyQuotient{\TopologicalSpace}{\Gamma}
  \end{tikzcd}
$$
\end{proof}

\begin{example}[Projective representations of finite groups seen via simplicial presheaves]
\label{ProjectiveRepresentationsAndTheirCentralExtensions}
The topological action groupoid (Ex. \ref{TopologicalActionGroupoid})
of the circle group $\CircleGroup$ acting by
multiplication on the unitary group $\UH$ \eqref{TheGroupUH}
carries an evident structure of a group object (by the degreewise group operations in
$\UH$ and in $\CircleGroup \times \UH$),
hence of a topological strict 2-group (Ntn. \ref{StrictTwoGroups}),
so that its simplicial nerve (Ntn. \ref{NerveOfTopologicalGroupoids})
is a simplicial topological group:
$$
  \ActionGroupoid
    { \UH }{ \CircleGroup }
  \;\;
  \in
  \;\;
  \Groups
  \big(
    \Groupoids(\kTopologicalSpaces)
  \big)
  \,,
  \;\;\;\;\;
  \SimplicialNerve
  \big(
    \ActionGroupoid
      { \UH }{ \CircleGroup }
  \big)
  \;\;\;
  \in
  \;
  \Groups
  \big(
    \SimplicialTopologicalSpaces
  \big)
  \,.
$$
Under the inclusion from Prop. \ref{CategoryOfDeltaGeneratedTopologicalSpaces}
$$
  \kTopologicalSpaces \xrightarrow{\ContinuousDiffeology} \DTopologicalSpaces
\xhookrightarrow{\;} \DiffeologicalSpaces \xhookrightarrow{\;} \SmoothInfinityGroupoids
$$
and using Lemma \ref{StabilizerGroupsAsLoopObjects},
this becomes a smooth 2-group, represented by a presheaf
of simplicial groups
on $\CartesianSpaces$:
\begin{equation}
  \label{CircleActionGroupoidOnUHAsSmoothTwoGroup}
  \ContinuousDiffeology
  \Big(
    \SimplicialNerve
    \big(
      \ActionGroupoid
        { \UH }{ \CircleGroup }
    \big)
  \Big)
  \;:\;
  \RealNumbers^n
  \;\mapsto\;
  \SimplicialNerve
  \big(
    \ActionGroupoid
      { \Maps{}{\RealNumbers^n}{\UH} }
      { \Maps{}{\RealNumbers^n}{\CircleGroup} }
  \big)
  \,.
\end{equation}
Since the $\CircleGroup$-quotient projection
$\UH \to \PUH$ is locally trivial \eqref{PUHFiberSequence}, so that
every map from an $\RealNumbers^n$ to $\PUH$ lifts to a map to $\UH$,
it follows that the canonical map from
the simplicial classifying space
\eqref{StandardSimplicialDeloopingAsQuotient} of
\eqref{CircleActionGroupoidOnUHAsSmoothTwoGroup}
to that of the group object represented by $\PUH$  \eqref{TheGroupPUH}
is a weak equivalence in the projective model structure over $\CartesianSpaces$
(Ntn. \ref{ModelCategoriesOfSimplicialPresheaves}
which presents $\SmoothInfinityGroupoids$ (Ntn. \ref{SmoothInfinityGroupoids})
of models for the delooping $\mathbf{B} \PUH$
(Lemma \ref{InfinityGroupsPresentedByPresheavesOfSimplicialGroups}):
$$
  \begin{tikzcd}
    \ContinuousDiffeology
    \Big(
      \overline{W}
      \SimplicialNerve
      \big(
        \ActionGroupoid
          { \UH }{ \CircleGroup }
      \big)
    \Big)
    \ar[
      rr,
      "{ \in \ProjectiveWeakEquivalences }"{swap}
    ]
    &&
    \ContinuousDiffeology
    \big(
      \overline{W}
      (\PUH)
    \big)
    \;\;
    \simeq
    \;
    \mathbf{B} \PUH
  \end{tikzcd}
$$
Inspection of the formula \eqref{FaceMapsOfWG}
for $\overline{W}$  shows that
this simplicial presheaf presentation of $\mathbf{B} \PUH \,\in\, \SmoothInfinityGroupoids$
looks as shown in the middle right of the following diagram, where the top
right shows an analogous presentation of $\mathbf{B} \UH$ which is such as
to make the delooped quotient $\mathbf{B}\UH \xrightarrow{\;} \mathbf{B} \PUH$
be presented by a fibration in $\SimplicialPresheaves(\CartesianSpaces)_{\mathrm{proj}}$:
\begin{equation}
  \label{DecomposingProjectiveRepresentations}
  \hspace{-1cm}
  \begin{array}{c}
  \mathllap{
    \mathbf{B}
    \widehat{G}
    \simeq
    \;
  }
  \left\{
  \!\!\!
  \adjustbox{raise=6pt}{
  \begin{tikzcd}[row sep=30pt, column sep=large]
    &
    \bullet
    \ar[
      dr,
      "{ ( g_2, \, c_2 ) }"{sloped}
    ]
    \\
    \bullet
    \ar[
      ur,
      "{
        ( g_1, \, c_1 )
      }"{sloped}
    ]
    \ar[
      rr,
      "{
          \big(
            g_1 \cdot g_2
            ,\,
            \tau(g_1, g_2) \cdot c_1 \cdot c_2
          \big)
      }"{swap},
      "{\ }"{name=t, pos=.45}
    ]
    &&
    \bullet
    \ar[
      from=t,
      to=t,
      start anchor={[yshift=+26pt, xshift=4pt]},
      end anchor={[xshift=-4pt]},
      -,
      shift left=1pt
    ]
    \ar[
      from=t,
      to=t,
      start anchor={[yshift=+26pt, xshift=4pt]},
      end anchor={[xshift=-4pt]},
      -,
      shift right=1pt
    ]
  \end{tikzcd}
  }
  \!\!\!
  \right\}
  \\
  \begin{tikzcd}[row sep=40pt, column sep=large]
    {}
    \ar[
      d,
      line width=1pt,
      "{
        \mathrlap{
        \def\arraystretch{.9}
        \begin{array}{c}
          (g_i,c_i)
          \\
          \mapsdown
          \\
          g_i
        \end{array}
        }
      }"{swap}
    ]
    \\
    {}
  \end{tikzcd}
  \\
  \mathllap{
    \mathbf{B}
    G
    \simeq
    \;
  }
  \left\{
  \!\!\!
  \adjustbox{raise=4pt}{
  \begin{tikzcd}[row sep=30pt]
    &
    \bullet
    \ar[
      dr,
      "{ g_2 }"
    ]
    \\
    \bullet
    \ar[
      ur,
      "{
        g_1
      }"
    ]
    \ar[
      rr,
      "{ g_1 \cdot g_2 }"{swap},
      "{\ }"{name=t, pos=.45}
    ]
    &&
    \bullet
    \ar[
      from=t,
      to=t,
      start anchor={[yshift=+26pt, xshift=4pt]},
      end anchor={[xshift=-4pt]},
      -,
      shift left=1pt
    ]
    \ar[
      from=t,
      to=t,
      start anchor={[yshift=+26pt, xshift=4pt]},
      end anchor={[xshift=-4pt]},
      -,
      shift right=1pt
    ]
  \end{tikzcd}
  }
  \!\!\!
  \right\}
  \\
  \begin{tikzcd}[row sep=40pt, white, column sep=large]
    {}
    \ar[
      d,
      line width=1pt,
      "{
        \mathrlap{
          \in \ProjectiveFibrations
        }
      }"
    ]
    \\
    {}
  \end{tikzcd}
  \\
  \adjustbox{raise=4pt}{
  \begin{tikzcd}[row sep=30pt, white, column sep=large]
    &
    |[alias=top]|
    \bullet
    \\
    |[alias=left]|
    \bullet
    \ar[
      rr,
      "{
        \bullet
      }"{swap},
      "{\ }"{name=t, pos=.45}
    ]
    &&
    |[alias=right]|
        \bullet
    \ar[
      from=t,
      to=t,
      start anchor={[yshift=+26pt, xshift=4pt]},
      end anchor={[xshift=-4pt]},
      Rightarrow,
      "{
        c
        =
        \tau(g_1, g_2)
      }"{description}
    ]
    \ar[
      from=left,
      to=top,
      "{ \bullet }"
    ]
    \ar[
      from=top,
      to =right,
      "{ \bullet }"
    ]
  \end{tikzcd}
  }
  \end{array}
  \hspace{-.8cm}
  \begin{array}{c}
  \adjustbox{raise=6pt}{
  \begin{tikzcd}[row sep=30pt, white, column sep=large]
    &
    \bullet
    \ar[
      dr,
      "{ ( g_2, \, c_2 ) }"{sloped}
    ]
    \\
    \phantom{\bullet}
    \ar[
      ur,
      "{
        ( g_1, \, c_1 )
      }"{sloped}
    ]
    \ar[
      rr,
      black,
      shift left=17pt,
      line width=1pt,
      "{ \widehat{\rho} }",
      "{
        \mbox{
          \tiny
          \color{greenii}
          \bf
          \def\arraystretch{.9}
          \begin{tabular}{c}
            representation of
            \\
            central extension
          \end{tabular}
        }
      }"{swap}
    ]
    \ar[
      rr,
      "{
          \big(
            g_1 \cdot g_2
            ,\,
            c(g_1, g_2) \cdot c_1 \cdot c_2
          \big)
      }"{swap},
      "{\ }"{name=t, pos=.45}
    ]
    &&
    \bullet
  \end{tikzcd}
  }
  \\
  \begin{tikzcd}[row sep=40pt]
    {}
    \ar[
      d,
      phantom,
      "{
        \mbox{
          \tiny
          \rm
          (pb)
        }
      }"
    ]
    \\
    {}
  \end{tikzcd}
  \\
  \adjustbox{raise=4pt}{
  \begin{tikzcd}[row sep=30pt, white]
    &
    \bullet
    \ar[
      dr,
      "{ g_2 }"
    ]
    \\
    \phantom{\bullet}
    \ar[
      ur,
      "{
        g_1
      }"
    ]
    \ar[
      rr,
      black,
      shift left=17pt,
      line width=1pt,
      "{ \rho }",
      "{
        \mbox{
          \tiny
          \color{greenii}
          \bf
          \begin{tabular}{c}
            projective
            \\
            representation
          \end{tabular}
        }
      }"{swap}
    ]
    \ar[
      rr,
      "{ g_1 \cdot g_2 }"{swap},
      "{\ }"{name=t, pos=.45}
    ]
    &&
    \bullet
  \end{tikzcd}
  }
  \\
  \begin{tikzcd}[row sep=40pt, column sep=large]
    {}
    \ar[
      d,
      line width=1pt,
      start anchor={[xshift=-45pt]},
      end anchor={[xshift=+10pt, yshift=-6pt]},
      "{ \tau }",
      "{
        \mbox{
          \tiny
          \color{greenii}
          \bf
          \def\arraystretch{.9}
          \begin{tabular}{c}
            twisting
            2-cocycle
          \end{tabular}
        }
      }"{swap, sloped}
    ]
    \\
    {}
  \end{tikzcd}
  \\
  \adjustbox{raise=4pt}{
  \begin{tikzcd}[row sep=30pt, white]
    &
    |[alias=top]|
    \bullet
    \\
    |[alias=left]|
    \bullet
    \ar[
      rr,
      "{
        \bullet
      }"{swap},
      "{\ }"{name=t, pos=.45}
    ]
    &&
    |[alias=right]|
        \bullet
    \ar[
      from=t,
      to=t,
      start anchor={[yshift=+26pt, xshift=4pt]},
      end anchor={[xshift=-4pt]},
      Rightarrow,
      "{
        c
        =
        \tau(g_1, g_2)
      }"{description}
    ]
    \ar[
      from=left,
      to=top,
      "{ \bullet }"
    ]
    \ar[
      from=top,
      to =right,
      "{ \bullet }"
    ]
  \end{tikzcd}
  }
  \end{array}
  \hspace{-.8cm}
  \begin{array}{c}
  \left\{
  \!\!\!
  \adjustbox{raise=6pt}{
  \begin{tikzcd}[row sep=30pt]
    &
    \bullet
    \ar[
      dr,
      "{ ( \UnitaryOperator_2, c_2 ) }"{sloped}
    ]
    \\
    \bullet
    \ar[
      ur,
      "{
        ( \UnitaryOperator_1, c_1)
      }"{sloped}
    ]
    \ar[
      rr,
      "{
          \big(
            \UnitaryOperator_1 \cdot \UnitaryOperator_2
            ,\,
            c \cdot c_1 \cdot c_2
          \big)
      }"{swap},
      "{\ }"{name=t, pos=.45}
    ]
    &&
    \bullet
    \ar[
      from=t,
      to=t,
      start anchor={[yshift=+26pt, xshift=4pt]},
      end anchor={[xshift=-4pt]},
      Rightarrow,
      "{ c }"
    ]
  \end{tikzcd}
  }
  \!\!\!
  \right\}
  \mathrlap{
    \; \simeq \,  \mathbf{B} \UH
  }
  \\
  \begin{tikzcd}[row sep=40pt]
    {}
    \ar[
      d,
      line width=1pt,
      "{
        \mathllap{
        \def\arraystretch{.9}
        \begin{array}{c}
          (U_i,c_i)
          \\
          \mapsdown
          \\
          c_i \cdot U_i
        \end{array}
        }
      }"{swap},
      "{
        \mathrlap{
          \in \ProjectiveFibrations
        }
      }"
    ]
    \\
    {}
  \end{tikzcd}
  \\
  \left\{
  \!\!\!
  \adjustbox{raise=4pt}{
  \begin{tikzcd}[row sep=30pt, column sep=large]
    &
    |[alias=top]|
    \bullet
    \\
    |[alias=left]|
    \bullet
    \ar[
      rr,
      "{
        c \cdot \UnitaryOperator_1 \cdot \UnitaryOperator_2
        =
        \tau(g_1, g_2) \cdot \rho(g_1 \cdot g_2)
      }"{swap},
      "{\ }"{name=t, pos=.45}
    ]
    &&
    |[alias=right]|
    \bullet
    \ar[
      from=t,
      to=t,
      start anchor={[yshift=+26pt, xshift=4pt]},
      end anchor={[xshift=-4pt]},
      Rightarrow,
      "{
        c
        =
        \tau(g_1, g_2)
      }"{description}
    ]
    \ar[
      from=top,
      to=right,
      "{
        \UnitaryOperator_2
        =
        \rho(g_2)
      }"{sloped}
    ]
    \ar[
      from=left,
      to=top,
      "{
        \UnitaryOperator_1
        =
        \rho(g_1)
      }"{sloped}
    ]
  \end{tikzcd}
  }
  \!\!\!
  \right\}
  \mathrlap{
    \; \simeq \,
    \mathbf{B} \PUH
  }
  \\
  \begin{tikzcd}[row sep=40pt]
    {}
    \ar[
      d,
      line width=1pt,
      "{
        \mathrlap{
          \in \ProjectiveFibrations
        }
      }"
    ]
    \\
    {}
  \end{tikzcd}
  \\
  \left\{
  \!\!\!
  \adjustbox{raise=4pt}{
  \begin{tikzcd}[row sep=30pt]
    &
    |[alias=top]|
    \bullet
    \\
    |[alias=left]|
    \bullet
    \ar[
      rr,
      "{
        \bullet
      }"{swap},
      "{\ }"{name=t, pos=.45}
    ]
    &&
    |[alias=right]|
        \bullet
    \ar[
      from=t,
      to=t,
      start anchor={[yshift=+26pt, xshift=4pt]},
      end anchor={[xshift=-4pt]},
      Rightarrow,
      "{
        c
        =
        \tau(g_1, g_2)
      }"{description}
    ]
    \ar[
      from=left,
      to=top,
      "{ \bullet }"
    ]
    \ar[
      from=top,
      to =right,
      "{ \bullet }"
    ]
  \end{tikzcd}
  }
  \!\!\!
  \right\}
  \mathrlap{
    \; \simeq \,  \mathbf{B}^2 \CircleGroup
  }
  \end{array}
\end{equation}
The bottom left of this diagram shows $\overline{W}G \,\simeq\, \mathbf{B}G$
(still by Lemma \ref{InfinityGroupsPresentedByPresheavesOfSimplicialGroups})
for a finite group $G$, which, since $G$ is discrete, is a projectively cofibrant
simplicial presheaf (by Prop. \ref{DuggerCofibrancyRecognition}).
It follows (by Lem. \ref{HomInfinityGroupoidFromCofibrantDomainAndFibrantCodomain})
that the elements of the hom-groupoid
$\SmoothInfinityGroupoids\big( \mathbf{B}G,\, \mathbf{B}\PUH \big)$
are all representable by morphisms of simplicial presheaves denoted ``$\rho$''
in this diagram, which,
again by Prop. \ref{CategoryOfDeltaGeneratedTopologicalSpaces},
are equivalently morphisms of
simplicial D-topological spaces
(in fact of topological 2-functors between D-topological 2-groupoids).

The upshot is that
the presentation for $\mathbf{B}\PUH$ chosen on the right
makes transparently manifest how such {\it projective unitary representations}
$G \xrightarrow{\rho} \PUH$
(e. g. \cite{Tappe77}\cite{Costache09}\cite[\S 5]{EspinozaUribe14},
review in \cite{Mendonca17})
encode:

\noindent
{\bf (i)} a $\CircleGroup$-valued 2-cocycle $\tau$,
by postcomposition with $\mathbf{B}\PUH \to \mathbf{B}^2 \CircleGroup$
(this construction is
the stacky refinement of Ex. \ref{HomotopyFiberSequenceOfTheProjectiveUnitaryGroup});

\noindent
{\bf (ii)}
a genuine unitary representation $\widehat{\rho}$ of the $\CircleGroup$-extension
$\widehat{G}$ of $G$ which is classified by $\tau$, this being simply the pullback
of $\rho$ along the top right projective fibration
(compare \cite[\S 5]{EspinozaUribe14}).

Abstracting away from the presentations by simplicial presheaves, the bottom right
fibration in
implies (by Lem. \ref{ComputingHomotopyPullbacksOfInfinityStacks})
that we have a homotopy fiber sequence
$\mathbf{B} \UH \xrightarrow{\;} \mathbf{B}\PUH \xrightarrow{\;} \mathbf{B}^2 \CircleGroup$
in $\SmoothInfinityGroupoids$, which, by Prop. \ref{GroupsActionsAndFiberBundles}
gives an identification
$$
  \mathbf{B} \PUH
    \;\simeq\;
  \HomotopyQuotient
    { \mathbf{B} \UH }
    { \mathbf{B} \CircleGroup }
  \;\;\;
  \in
  \;
  \SmoothInfinityGroupoids
  \,.
$$
It follows that for $\tau \,\in\, H^2_{\mathrm{Grp}}(G;\, \CircleGroup)$,
the isomorphism classes of $[\tau]$-projective unitary representations of $G$
are in natural bijection with
the connected components of the slice hom-groupoid \eqref{HomSpaceInSliceAsFiberProduct}
\begin{equation}
  \label{IsomorphismClassesOfProjectiveRepresentationsAsConnectedComponentsOfSliceHom}
  \IsomorphismClasses{
    \Representations^{[\tau]}(G)
  }
  \;\;
  \simeq
  \;\;
  \Truncation{0}
  \,
  \SlicePointsMaps{\big}{\mathbf{B}^2 \CircleGroup}
    { (\mathbf{B}G, \tau) }
    { \HomotopyQuotient{ \mathbf{B}\UH }{ \mathbf{B}\CircleGroup } }
  \,.
\end{equation}
\end{example}

\medskip

We use the occasion of Ex. \ref{ProjectiveRepresentationsAndTheirCentralExtensions} to record
that key facts of standard character theory generalize to
projective representations and to establish a fact
(Prop. \ref{TensoringProjectiveRepresentationsWithThePlainRegularRepresentation}
below) which we need below in Lem. \ref{StableProjectiveIsotorpyRepresentations}:

\begin{lemma}[Characters of projective representations of finite groups {\cite[Prop. 2.2 (1)]{Cheng15}}]
\label{CharactersOfProjectiveRepresentations}
$\,$

\noindent
  For $G$ a finite group and $[\tau] \,\in\, H^2_{\mathrm{Grp}}(G,\CircleGroup)$
  a 2-cocycle, the {\it characters}
  $$
    \rchi_{\rho}
    \;\;\;
    \in
    \;
    \Maps{}{G}{\CircleGroup}
  $$
  of $[\tau]$-projective
  unitary representations
  (Ex. \ref{ProjectiveRepresentationsAndTheirCentralExtensions})
  $$
    \begin{tikzcd}
      \mathbf{B}G
      \ar[r, "{ \rho }"]
      \ar[dr, "{ \tau }"{swap}]
      &
      \HomotopyQuotient{ \mathbf{B} \UnitaryGroup_d }{ \mathbf{B}\CircleGroup }
      \ar[r]
      &
      \HomotopyQuotient{ \mathbf{B} \UH }{ \mathbf{B}\CircleGroup }
      \mathrlap{\; \simeq\; \mathbf{B}\PUH}
      \ar[dl]
      \\
      &
      \mathbf{B}^2 \CircleGroup
    \end{tikzcd}
  $$
  of any finite dimension $d \,\in\, \NaturalNumbers$,
  which are still given by the trace of the unitary operators
  $\rho(g) \,\in\, \UnitaryGroup_d$
  \eqref{DecomposingProjectiveRepresentations}
  \begin{equation}
    \label{TraceFormulaForProjectiveCharacters}
    \Character_{\rho}(g)
    \;:=\;
    \mathrm{Tr}\big( \rho(g) \big)
    \,,
  \end{equation}
  still detect isomorphy of projective representations \eqref{IsomorphismClassesOfProjectiveRepresentationsAsConnectedComponentsOfSliceHom}:
  \begin{equation}
    \label{CharactersReflectProjectiveRepresentations}
    \Character_{\rho_1}
    \;=\;
    \Character_{\rho_2}
    \;\;\;\;\;\;\;\;
    \Leftrightarrow
    \;\;\;\;\;\;\;\;
    [\rho_1]
    \,=\,
    [\rho_2]
    \;\;\;
    \in
    \;
    \IsomorphismClasses{
      \Representations^{[\tau]}(G)
    }
    \,.
  \end{equation}
\end{lemma}
This means that projective character theory looks much like
plain character theory; for instance:
\begin{example}[Ordinary characters act on projective characters]
  The tensoring (see also \eqref{TensoringOfPUHPrimeOverUHPrime} below)
  of an ordinary unitary $G$-representation
  $\mathbf{B}G \xrightarrow{\rho} \mathbf{B}\UnitaryGroup_d$
  with a $[\tau]$-projective representation
  $\mathbf{B}G \xrightarrow{\rho} \mathbf{B}\UnitaryGroup_{d'}$
  (Ex. \ref{ProjectiveRepresentationsAndTheirCentralExtensions})
  is again a $[\tau]$-projective representation
  \begin{equation}
    \label{TensoringOfProjectiveRepresentationsWithPlainRepresentations}
    \begin{tikzcd}[row sep=0pt]
      \Representations(G)
      \times
      \Representations^{[\tau]}(G)
      \ar[rr, "{ \otimes }"]
      &&
      \Representations^{[\tau]}(G)
      \\
    \scalebox{0.7}{$  (\rho, \, \rho') $}
      &\longmapsto&
    \scalebox{0.7}{$        \rho \otimes \rho' $}
      \,,
    \end{tikzcd}
  \end{equation}
  and under this operation the characters \eqref{TraceFormulaForProjectiveCharacters}
  multiply, as usual:
  \begin{equation}
    \label{CharactersMultiplyUnderTensoringOfProjectiveWithPlainRepresentations}
    \Character_{ \rho \otimes \rho' }
    \;=\;
    \Character_{\rho} \cdot \Character_{\rho'}
    \,.
  \end{equation}
\end{example}
\begin{example}[Projective regular representation]
  \label{ProjectiveRegularRepresentation}
  For $G$ a finite group and $[\tau] \,\in\, H^2_{\mathrm{Grp}}(G;\, \CircleGroup)$,

  \noindent
  {\bf (i)}
  the {\it regular} $[\tau]$-projective representation (Ex. \ref{ProjectiveRepresentationsAndTheirCentralExtensions}) is the
  complex linear span of the underlying set of $G$
  $$
    \ComplexNumbers[G]^{[\tau]}
    \,\in\,
    \Representations^{[\tau]}(G)
  $$
  with action defined on the canonical basis elements
  $\{v_h \vert g \in G\} \,\subset \, \ComplexNumbers[G]$ by the formula
  \begin{equation}
    \label{ActionOnProjectiveRegularRepresentation}
    \rho(g)(v_k)
    \;:=\;
    \tau(g,k)\cdot v_{g \cdot k}
    \,.
  \end{equation}
  That this definition satisfies the action property is equivalently the
  cocycle condition on $\tau$:
  $$
    \begin{array}{lll}
      \rho(g_1)
      \big(
        \rho(g_2)
        (v_{g_3})
      \big)
      &
      \;=\;
      \rho(g_1)
      \big(
        \tau(g_2, g_3) \cdot v_{g_2 \cdot g_3}
      \big)
      &
      \proofstep{
        by \eqref{ActionOnProjectiveRegularRepresentation}
      }
      \\
      & \;=\;
      \tau(g_1, g_2 \cdot g_3) \cdot \tau(g_2, g_3) \cdot v_{g_1 \cdot g_2 \cdot g_3}
      &
      \proofstep{
        by \eqref{ActionOnProjectiveRegularRepresentation}
      }
      \\
      & \;=\;
      \tau(g_1, g_2)
      \cdot
      \tau( g_1 \cdot g_2, g_3 )
      \cdot
      v_{g_1 \cdot g_2 \cdot g_3}
      &
      \proofstep{
        by cocycle property
      }
      \\
      & \;=\;
      \tau(g_1, g_2)
      \cdot
      \rho(g_1 \cdot g_2)(v_{g_3})
      &
      \proofstep{
        by \eqref{ActionOnProjectiveRegularRepresentation}
      }
      .
    \end{array}
  $$

  \noindent
  {\bf (ii)}
  The character \eqref{CharactersOfProjectiveRepresentations} of
  the regular projective representation evidently is, independently of $\tau$,
  the same as that of the ordinary regular representation:
  \begin{equation}
    \label{CharacterOfTheProjectiveRegularRepresentation}
    \Character_{\ComplexNumbers[G]^{\tau}}
    \;:\;
    g
    \;\longmapsto\;
    \left\{
    \!\!\!
    \begin{array}{ccl}
      \mathrm{ord}(G) &\vert& g = \NeutralElement
      \\
      0 &\vert& \mbox{otherwise.}
    \end{array}
    \right.
  \end{equation}
\end{example}

\begin{proposition}[Regular projective representation decomposes into all projective irreps
{\cite[Prop. 2.3]{Cheng15}}]
  \label{DecompositionOfRegularProjectiveRepresentation}
  For $G$ a finite group and $[\tau] \,\in\, H^2_{\mathrm{Grp}}(G;\, \CircleGroup)$,
  the $[\tau]$-projective regular representation
  (Ex. \ref{ProjectiveRegularRepresentation})
  decomposes as the direct sum of all
  irreducible
  $[\tau]$-projective representations, each appearing with multiplicity
  equal to the complex dimension of its representation space:
  $$
    \ComplexNumbers[G]^{[\tau]}
    \;\simeq\quad
    \;
    \underset{
      \mathclap{
      \scalebox{.7}{$
        \def\arraystretch{.9}
        \begin{array}{c}
          {[\mu] \in}
          \\
          \IsomorphismClasses{
            \Representations^{[\tau]}(G)_{\mathrm{irr}}
          }
        \end{array}
      $}
      }
    }{\bigoplus}
    \quad
    \mathrm{dim}_{{}_{\ComplexNumbers}}(\mu)
    \cdot
    \mu
    \;\;\;
    \in
    \;
    \Representations^{[\tau]}(G)\;.
  $$
\end{proposition}

We deduce from this the following fact
(which will be needed  in Lem. \ref{StableProjectiveIsotorpyRepresentations} below):
\begin{proposition}[Tensoring projective representations with the plain regular representation]
  \label{TensoringProjectiveRepresentationsWithThePlainRegularRepresentation}
  For $G$ a finite group, and $[\tau] \,\in\, H^2_{\mathrm{Grp}}(G;\, \CircleGroup)$,
  the tensoring \eqref{TensoringOfProjectiveRepresentationsWithPlainRepresentations}
  of any $[\tau]$-projective representation $\rho$
  (Ex. \ref{ProjectiveRepresentationsAndTheirCentralExtensions})
  with the plain regular representation is
  the direct sum of $\mathrm{dim}_{\mathbb{C}}(\rho)$ copies of the
  $[\tau]$-projective regular representation (Ex. \ref{ProjectiveRegularRepresentation}),
  and hence a direct sum of all irreducible $[\tau]$-projective representations:
  $$
    \rho
    \,\in\,
    \Representations^{[\tau]}(G)
    \;\;\;\;\;
    \vdash
    \;\;\;\;\;
    \left\{
    \def\arraystretch{1.3}
    \begin{array}{lll}
      \mathbb{C}[G] \otimes \rho
      &
      \;\simeq\;
      \mathrm{dim}_{\mathbb{C}}(\rho)
      \cdot
      \mathbb{C}[g]^{\tau}
      \\
      &
      \;\simeq\qquad
    \underset{
      \mathclap{
      \scalebox{.7}{$
        \def\arraystretch{.9}
        \begin{array}{c}
          {[\mu] \in}
          \\
          \IsomorphismClasses{
            \Representations^{[\tau]}(G)_{\mathrm{irr}}
          }
        \end{array}
      $}
      }
    }{\bigoplus}
    \quad
    \mathrm{dim}_{{}_{\ComplexNumbers}}(\rho)
    \cdot
    \mathrm{dim}_{{}_{\ComplexNumbers}}(\mu)
    \cdot
    \mu
    &
    \proofstep{
      by Prop. \ref{DecompositionOfRegularProjectiveRepresentation}.
    }
    \end{array}
    \right.
  $$
\end{proposition}
\begin{proof}
  The character
  \eqref{TraceFormulaForProjectiveCharacters}
  of the tensor product representation is
  $$
    \def\arraystretch{1.3}
    \begin{array}{lll}
      \chi
        \big(
          \ComplexNumbers[G] \otimes \rho
        \big)
      &
      \;=\;
      \chi
      \big(
        \ComplexNumbers[G]
      \big)
      \cdot
      \chi
      \big(
        \rho
      \big)
      &
      \proofstep{
        by \eqref{CharactersMultiplyUnderTensoringOfProjectiveWithPlainRepresentations}
      }
      \\
      & \;=\;
      \big( \lvert G \rvert, 0, 0, \cdots, \big)
      \cdot
      \big( \mathrm{dim}_{\ComplexNumbers}(\rho), \mathrm{tr}(\rho(g_1)), \cdots \big)
      &
      \proofstep{
        by \eqref{CharacterOfTheProjectiveRegularRepresentation}
      }
      \\
      & \;=\;
      \big( \lvert G \rvert \cdot \mathrm{dim}_{\ComplexNumbers}(\rho) , 0, 0, \cdots\big)
      \\
      & \;=\;
      \mathrm{dim}_{\ComplexNumbers}(\rho)
      \cdot
      \big( \lvert G \rvert  , 0, 0, \cdots\big)
      \\
      & \;=\;
      \mathrm{dim}_{\ComplexNumbers}(\rho)
      \cdot
      \chi
      \big(
        \ComplexNumbers[G]^{\tau}
      \big)
      &
      \proofstep{
        by \eqref{CharacterOfTheProjectiveRegularRepresentation}
      }
      \,.
    \end{array}
  $$
  Hence the claim follows by \eqref{CharactersReflectProjectiveRepresentations}.
\end{proof}

\noindent
{\bf Smooth cohesion.} Smooth $\infty$-groupoids constitute the main example of cohesive $\infty$-toposes that we consider here (see also \cite[\S 3.1]{SS20OrbifoldCohomology} for this case and further variants):

\begin{proposition}[Smooth cohesion
  {\cite[\S 3.1]{SSS09}\cite[Prop. 4.4.8]{dcct}}
  ]
  \label{SmoothCohesion}
  The $\infty$-topos
  $\SmoothInfinityGroupoids$
  of smooth $\infty$-groupoids (Ntn. \ref{SmoothInfinityGroupoids})
  is cohesive
  in the sense of Def. \ref{CohesiveInfinityTopos}:
  The adjoint quadruple
  arises just as in Ex. \ref{ExamplesOfDiscreteCohesion} from the
  inclusion $\ast \,\simeq\,\{\mathbb{R}^0\} \xhookrightarrow{\;} \CartesianSpaces$
  of the terminal object into the site
  \eqref{CartesianSpaces}
  of smooth Cartesian spaces, and is then seen to factor through the
  full inclusion
  $\SmoothInfinityGroupoids
  \,\simeq\, \InfinitySheaves(\CartesianSpaces) \xhookrightarrow{\;}
  \InfinityPresheaves$.
\end{proposition}
Properties of smooth cohesion are
derived in \cite[\S 4.4]{dcct}
and key facts are summarized in
\cite[Ex. 3.18]{SS20OrbifoldCohomology}\cite[Ex. A.58]{FSS20CharacterMap}
and above in \cref{OverviewAndSummary}. We proceed to highlight some of these
in more detail.

To start with, contact with the differential topology discussed above, is made by observing that diffeological spaces, and with them the D-topological spaces (Prop. \ref{CategoryOfDeltaGeneratedTopologicalSpaces}), are most naturally included among smooth $\infty$-groupoids:

\begin{proposition}[Diffeological spaces are the concrete 0-truncated smooth $\infty$-groupoids
{\cite[Prop. 4.4.15]{dcct}}]
  \label{DiffeologicalSpacesAreTheOneConcreteZeroTruncatedSmoothInfinityGroupoids}
  The canonical inclusion of
  diffeological spaces (Ntn. \ref{CartesianSpacesAndDiffeologicalSpaces})
  into smooth $\infty$-groupoids (Ntn. \ref{SmoothInfinityGroupoids})
  is equivalently the inclusion of the
  concrete 0-truncated smooth $\infty$-groupoids (cf. \cite[Ex. 3.18(ii)]{SS20OrbifoldCohomology}):
  \vspace{-2mm}
  \begin{equation}
    \label{InjectionOfDiffeologicalSpacesAreTheOneConcreteZeroTruncatedSmoothInfinityGroupoids}
    \DiffeologicalSpaces
    \;\simeq\;
    \Smooth\Groupoids_{0,\sharp_1}
    \xhookrightarrow{\quad}
    \SmoothInfinityGroupoids \;.
  \end{equation}
\end{proposition}

\medskip

\noindent
{\bf Cohesive shape of smooth manifolds.}

\begin{proposition}[Quillen functor for shape modality on smooth $\infty$-groupoids]
  \label{QuillenFunctorForShapeModalityOnSmoothInfinityGroupoids}
  Consider the shape-adjunction on smooth $\infty$-groupoids (Def. \ref{SmoothInfinityGroupoids})
    \vspace{-5mm}
  $$
    \begin{tikzcd}
      \SmoothInfinityGroupoids
      \ar[
        rr,
        shift left=5pt,
        "{\Shape}"
      ]
      \ar[
        rr,
        phantom,
        "\scalebox{.6}{$\bot$}"
      ]
      &&
      \InfinityGroupoids \;.
      \ar[
        ll,
        shift left=5pt,
        "{\Discrete}"
      ]
    \end{tikzcd}
  $$

  \vspace{-2mm}
  \noindent
 {\bf (i)} This adjunction is, equivalently, the left derived functor of
  the colimit operation on simplicial presheaves over Cartesian spaces,
  regarded as functors $\CartesianSpaces^{\mathrm{op}} \xrightarrow{\;} S\SimplicialSets$,
  in that the following is a Quillen adjunction:
    \vspace{-2mm}
  \begin{equation}
    \label{ColimitQuillenAdjunctionOnLocalSimplicialPresheaves}
    \begin{tikzcd}
      \SimplicialPresheaves(\CartesianSpaces)_{ {\mathrm{proj}} \atop {\mathrm{loc}} }
      \ar[
        rr,
        shift left=5pt,
        "{\underset{\longrightarrow}{\mathrm{lim}}}"{above}
      ]
      \ar[
        rr,
        phantom,
        "{\scalebox{.7}{$\bot_{\mathrlap{\mathrm{Qu}}}$}}"
      ]
      &&
    \;   \SimplicialSets_{\mathrm{Qu}} \;.
      \ar[
        ll,
        shift left=5pt,
        "{\mathrm{const}}"{below}
      ]
    \end{tikzcd}
  \end{equation}

    \vspace{-2mm}
\noindent
{\bf (ii)}  Moreover, on a simplicial presheaf satisfying Dugger's cofibrancy condition
  \vspace{-2mm}
  (Prop. \ref{DuggerCofibrancyRecognition})
  \begin{equation}
    \label{DuggerCofibrantResolution}
    \begin{tikzcd}
      \varnothing
      \ar[
        r,
        "\in \mathrm{Cof}"
      ]
      &
      \underset{
        i_\bullet \in I_\bullet
      }{\coprod}
      \mathbb{R}^{n_{i_\bullet}}
      \ar[
        r,
        "\in \mathrm{W}"
      ]
      &
      X
    \end{tikzcd}
    \;\;
    \in
    \SimplicialPresheaves(\CartesianSpaces)_{ {\mathrm{proh}} \atop {\mathrm{loc}} }\;,
  \end{equation}

  \vspace{-3mm}
  \noindent
  the shape is given by the simplicial set obtained by contracting all copies of
  Cartesian spaces to the point:
    \vspace{-1mm}
  $$
    \shape \, X
    \;\simeq\;
      \underset{
        i_\bullet \in I_\bullet
      }{\coprod}
      \ast
      \;\;\;
      \in
      \;
      \SimplicialSets_{\mathrm{Qu}}
      \,.
  $$
\end{proposition}
\begin{proof}
  First, observe that the colimit over a representable functor is the point
  (e.g. \cite[Lem. 2.40]{SS20OrbifoldCohomology})
  \vspace{-2mm}
  \begin{equation}
    \label{ColimitOfRepresentableIsPoint}
    \underset{\longrightarrow}{\mathrm{lim}}
    \,
    \mathbb{R}^n
    \;\coloneqq\;
    \underset{\longrightarrow}{\mathrm{lim}}
    \,
    y(\mathbb{R}^n)
    \;\simeq\;
    \ast
    \;\;\;
    \in
    \;
    \Sets
    \xhookrightarrow{\;\;}
    \SimplicialSets
    \,,
  \end{equation}

  \vspace{-2mm}
  \noindent
  so that the colimit of a simplicial presheaf of the form \eqref{DuggerCofibrantResolution}
  is the simplicial set obtained by replacing all copies of Cartesian spaces
  by a point:
    \vspace{-2mm}
  \begin{equation}
    \label{ColimitOfDuggerCofibrantSimplicialPresheaf}
    \begin{array}{lll}
      \underset{\longrightarrow}{\mathrm{lim}}
      \big(\,
        \underset{ i_\bullet \in I_\bullet }{\coprod}
        \mathbb{R}^{n_{i_\bullet}}
      \big)
      & \;\simeq\;
        \big(
          \underset{\longrightarrow}{\mathrm{lim}}
          \,
          \mathbb{R}^{n_{i_\bullet}}
        \big)
      &
      \mbox{\small
        since colimits commute with coproducts
      }
      \\
      & \;\simeq\;
      \underset{ i_\bullet \in I_\bullet }{\coprod}
      \ast
      &
      \mbox{\small  \; by \eqref{ColimitOfRepresentableIsPoint}}
      \,.
    \end{array}
  \end{equation}

  \vspace{-1mm}
  \noindent
  Next, it is clear that \eqref{ColimitQuillenAdjunctionOnLocalSimplicialPresheaves}
  is a simplicial Quillen adjunction for the {\it global} projective model structure.
  To show that it is also Quillen for the local model structure it is hence sufficient,
  by \cite[Cor. A.3.7.2]{Lurie09HTT}, to see that the right adjoint preserves
  fibrant objects. By adjunction, this is equivalent to the statement that
  for $\{ U_i \hookrightarrow X \}$ a differentiably good open cover,
  with $U \,\coloneqq\, \underset{i}{\coprod}\, U_i$,
  we
  have a simplicial weak homotopy equivalence
  \vspace{-3mm}
  \begin{equation}
    \label{WeakHomotopyEquivalenceFromColimitOfCechNerveOfGoddOpenCoverOfCartesianSpace}
    \underset{\longrightarrow}{\mathrm{lim}}\; y(U^{\times^\bullet_X})
    \xrightarrow{\;\in \mathrm{W}\;}
    \ast
    \,.
  \end{equation}

  \vspace{-2mm}
 \noindent But, by \eqref{ColimitOfDuggerCofibrantSimplicialPresheaf},
  the left hand side of \eqref{WeakHomotopyEquivalenceFromColimitOfCechNerveOfGoddOpenCoverOfCartesianSpace} is
  the simplicial set obtained by contracting summands of the
  {\v C}ech nerve of the good cover to the point.
  Therefore, since any Cartesian space is contractible,
  the {\it nerve theorem}
  (\cite[Thm. 2]{McCord67}, review in \cite[Prop. 4G.3]{Hatcher02})
  implies \eqref{WeakHomotopyEquivalenceFromColimitOfCechNerveOfGoddOpenCoverOfCartesianSpace}.
    With this, the last statement follows from
  the fact that
  left derived functors may be computed on any cofibrant resolution:
    \vspace{-2mm}
  $$
    \begin{array}{lll}
      \shape \, X
      &
      \;\simeq\;
      \big(\mathbb{L} \underset{\longrightarrow}{\mathrm{lim}}\big)(X)
      &
      \mbox{\small by \eqref{ColimitQuillenAdjunctionOnLocalSimplicialPresheaves}}
      \\
      &
      \;\simeq\;
      \underset{\longrightarrow}{\mathrm{lim}}
      \big(\,
        \underset{i_\bullet \in I_\bullet}{\coprod}
        \mathbb{R}^n
      \big)
      &
      \mbox{\small by Prop. \ref{DuggerCofibrancyRecognition}}
      \\
      & \;\simeq\;
        \underset{i_\bullet \in I_\bullet}{\coprod}
        \ast
      &
      \mbox{\small by \eqref{ColimitOfDuggerCofibrantSimplicialPresheaf}}
      \,.
    \end{array}
  $$

  \vspace{-8mm}
\end{proof}

\begin{example}[Good open covers are projectively cofibrant resolutions of smooth manifolds]
  \label{GoodOpenCoversAreProjectivelyCofibrantResolutionsOfSmoothManiolds}
  $\,$

  \noindent
  Every
  $\SmoothManifold \in \SmoothManifolds \xrightarrow{\; y \;} \SimplicialPresheaves(\CartesianSpaces)$
  admits a
  {\it differentiably good open cover}
  (Def. \ref{GoodOpenCovers}),
  namely an open cover
  such that all non-empty finite intersections of patches are
  {\it diffeomorphic} to an open ball, and hence to $\mathbb{R}^{\mathrm{dim}(X)}$:
  \vspace{-2mm}
  \begin{equation}
    \label{GoodOpenCover}
    \big\{
      \TopologicalPatch_i
        \,\simeq\,
    \mathbb{R}^{\mathrm{dim}(\SmoothManifold)}
      \xhookrightarrow{\;}
    \SmoothManifold
    \big\}_{i \in I}
    \,,
    \;\;\;
    \mbox{s.t.}
    \quad
    \underset
      {
        { k \in \mathbb{N} }
        \atop
        { i_0, i_1, \cdots, i_k \in I }
      }
      {\forall}
      \;\;\;
      \TopologicalPatch_{i_0}
        \cap
      \TopologicalPatch_{i_1}
        \cap
      \cdots
        \cap
      \TopologicalPatch_{i_k}
      \;
      \underset{
        \mathclap{
          \mathrm{diff}
        }
      }{
        \simeq
      }
      \;
      \mathbb{R}^{\mathrm{dim}(\SmoothManifold)}
      \;\;\;
      \mbox{\small if non-empty}
      \,.
  \end{equation}

  \vspace{-1mm}
  \noindent
  By Dugger's recognition (Prop. \ref{DuggerCofibrancyRecognition}),
  this means that the corresponding {\v C}ech nerve is projectively cofibrant
  (hence also locally projectively cofibrant);
  moreover, its canonical morphism to $\SmoothManifold$ is
  clearly a stalkwise weak equivalence, so that it provides a
  cofibrant resolution of $X$ in the local model structure (Def. \ref{SmoothInfinityGroupoids}):
  \vspace{-3mm}
  $$
    \begin{tikzcd}[column sep=large]
      \varnothing
      \ar[
        r,
        "\in \LocalProjectiveCofibrations"{swap}
      ]
      &
      \widehat{\SmoothManifold}^{\times^\bullet_{\SmoothManifold}}
      \ar[
        r,
        "\in \LocalWeakEquivalences"{swap}
      ]
      &
      X
    \end{tikzcd}
    \;\in\;
    \SimplicialPresheaves(\CartesianSpaces)_{ {\mathrm{proj}} \atop {\mathrm{loc}} }
    \,,
    \;\;\;\;
    \mbox{for}
    \;\;\;
    \widehat{\SmoothManifold}
      \;\coloneqq
    \underset{i}{\sqcup}
    \,
    \TopologicalPatch_i
    \;.
  $$
\end{example}
In conclusion:
\vspace{-1mm}
\begin{proposition}[Cohesive shape of smooth manifolds is their homotopy type]
  \label{CohesiveShapeOfSmoothManifoldsIsTheirHomotopyType}
  For $X \in \SmoothManifolds \,\xhookrightarrow{\;y\;}\, \SmoothInfinityGroupoids$
  (Def. \ref{SmoothInfinityGroupoids}), their cohesive shape \eqref{TheModalitiesOnACohesiveInfinityTopos}
  is their standard underlying homotopy type.
\end{proposition}
\begin{proof}
  By Ex. \ref{GoodOpenCoversAreProjectivelyCofibrantResolutionsOfSmoothManiolds},
  the smooth manifold admits a differentiably good open cover whose
  {\v C}ech nerve constitutes a cofibrant resolution of the manifold
  in the local projective model structure on simplicial presheaves.
  By Prop. \ref{QuillenFunctorForShapeModalityOnSmoothInfinityGroupoids}, this
  implies that the shape of the manifold is the homotopy type of the
  simplicial set obtained by replacing each direct summand in this
  {\v C}ech nerve by a point. That this simplicial set represents the
  standard homotopy type of the manifold is the classical {\it nerve theorem}
  (\cite[Thm. 2]{McCord67}, review in \cite[Prop. 4G.3]{Hatcher02}).
\end{proof}

\begin{example}[Standard cofibrant resolution of the smooth circle]
  \label{StandardCofibrantResolutionOfTheCircle}
  Considering the smooth circle as the quotient of the real numbers by the
  integers
    \vspace{-2mm}
  $$
    \begin{tikzcd}[row sep=5pt]
      \mathbb{Z}
      \ar[
        r,
        hook,
        "i"
      ]
      &
      \mathbb{R}
      \ar[
        r,
        ->>,
        "p"
      ]
      &
      S^1
    \end{tikzcd}
    \in
    \;\;
    \SmoothManifolds
    \xhookrightarrow{\;y\;}
    \SimplicialPresheaves(\CartesianSpaces)_{ {\mathrm{proj}} \atop {\mathrm{loc}} }
    \,,
  $$

  \vspace{-2mm}
  \noindent
  Dugger's recognition (Prop. \ref{DuggerCofibrancyRecognition})
  shows that
  the {\v C}ech nerve of $p$ consistitutes a cofibrant resolution of the circle
  in the projective local model structure (Def. \ref{SmoothInfinityGroupoids})
  \vspace{-2mm}
  $$
  \begin{tikzcd}[row sep=-2pt]
    \varnothing
    \ar[
      r,
      "\in \mathrm{Cof}"{above}
    ]
    &
    \mathbb{R} \times \mathbb{Z}^{\times^\bullet}
    \ar[
      rr,
      "\sim"{above, yshift=-1pt}
    ]
    &&
    \mathbb{R}^{\times^\bullet_{S^1}}
    \ar[
      r,
      "\in \mathrm{W}"
    ]
    &
    S^1
    \\
    &
   \!\!\!   \scalebox{0.7}{$  (r,\, \vec n) $}
    &\longmapsto  \;\;\;\;&
    \hspace{-3cm}
    \mathrlap{
     \scalebox{0.7}{$  \left(
      r, (r + n_1), (r + n_1 + n_2), \cdots
    \right)
    $}.
    }
  \end{tikzcd}
  $$
\end{example}

\medskip

\noindent
{\bf Equivariant {\v C}ech cocycles from projectively cofibrant resolution of action groupoids.}
In immediate generalization of Example \ref{GoodOpenCoversAreProjectivelyCofibrantResolutionsOfSmoothManiolds}
we have the following further consequence of Dugger's cofibrancy condition
(Prop. \ref{DuggerCofibrancyRecognition}):
\begin{example}[{\v C}ech-action groupoid of equivariant good open cover is local cofibrant resolution]
  \label{CechActionGroupoidOfEquivariantGoodOpenCoverIsLocalCofibrantResolution}
  For $G \,\in\, \Groups(\Sets) \xhookrightarrow{ \; }
  \Groups\big(\SimplicialPresheaves(\CartesianSpaces)\big)$.
  let
  $G \acts \,  \SmoothManifold \,\in\,
    \Actions{G}(\SmoothManifolds)$,
  with properly equivariant good open cover
  $
    G \acts \, \widehat{\SmoothManifold}
    \coloneqq
    G \acts \,
    (
      \underset{i \in I}{\sqcup}
      \TopologicalPatch_i
    )
    \twoheadrightarrow
    G \acts \, \SmoothManifold
  $ (Def. \ref{ProperEquivariantOpenCover},
  which exists by Prop. \ref{SmoothGManifoldsAdmitProperlyEquivariantGoodOpenCovers}).
  Then the diagonal of the {\v C}ech nerve of its action groupoid
  is a cofibrant replacement of the
  simplicial nerve of the action groupoid of $\SmoothManifold$
  in $\SimplicialPresheaves(\CartesianSpaces)_{\projloc}$:
  \begin{equation}
    \label{CechActionGroupoidComponents}
    \begin{tikzcd}
      \varnothing
      \ar[
        rr,
        "\in \LocalProjectiveCofibrations"{swap}
      ]
      &&
      {}
    \end{tikzcd}
    \SimplicialNerve
    \big(\!\!
    \begin{tikzcd}
      \widehat{\SmoothManifold}
        \times_{{}_{\SmoothManifold}}
      \widehat{\SmoothManifold}
        \times
      G^{\mathrm{op}}
      \ar[
        r,
        shift left=3pt,
        "{
          \mathrm{pr}_1
        }"
      ]
      \ar[
        r,
        shift right=3pt,
        "{
          \mathrm{pr}_3 \cdot \mathrm{pr}_2
        }"{swap}
      ]
      &
      \widehat{\SmoothManifold}
    \end{tikzcd}
    \!\!\big)
    \begin{tikzcd}
      {}
      \ar[
        rr,
        "\in \LocalWeakEquivalences"{swap}
      ]
      &&
      {}
    \end{tikzcd}
    \SimplicialNerve
    \big(
     \SmoothManifold
      \times
      G^{\mathrm{op}}
      \rightrightarrows
      \SmoothManifold
    \big)
    \\
  \end{equation}
  $$
  \hspace{-2.8cm}
  \left\{
  \!\!\!\!
  \adjustbox{raise=4pt}{
  \begin{tikzcd}[column sep=50pt, row sep=30pt]
    (x ,\, i)
    \ar[
          d,
      "{
         \scalebox{0.7}{$
         \left(
          (x,\,i,\,j)
          ,\,
          \NeutralElement
        \right)
        $}
      }"{description}
    ]
    \ar[
      r,
      "{
   \scalebox{0.7}{$
   \left(
          (x,\,i,\,i)
          ,\,
          g
        \right)
        $}
      }"
    ]
    \ar[
      dr,
      "{
      \scalebox{0.7}{$  \left(
          (x,\,i,\,j)
          ,\,
          g
        \right)
        $}
      }"{sloped, description}
    ]
    &
    (g\cdot x ,\, g \cdot i)
    \ar[
      d,
      "{
         \scalebox{0.7}{$
         \left(
          (g\cdot x, \, g \cdot i, \, g\cdot j)
          ,\,
          \NeutralElement
        \right)
        $}
      }"{description}
    ]
    \\
    (x,j)
    \ar[
      r,
      "{
      \scalebox{0.7}{$  \left(
          (x,\,j,\,j)
          ,\,
          g
        \right)
        $}
      }"{swap}
    ]
    &
    (g\cdot x ,\, g \cdot j)
  \end{tikzcd}
  }
  \!\!\!\!\!
  \right\}
$$
\end{example}
\begin{remark}[Equivariant {\v C}ech cocycles]
  \label{EquivariantCechCocycles}
  Example \ref{CechActionGroupoidOfEquivariantGoodOpenCoverIsLocalCofibrantResolution}
  implies that, under the conditions stated there,
  we have for $\mathcal{G} \,\in\, \Groups\big(\SimplicialPresheaves(\CartesianSpaces)\big)$
  an equivalence of $\infty$-groupoids
  \begin{equation}
    \label{IdentifyingEquivariantBundlesWithEquivariantCechCocycles}
    \SmoothInfinityGroupoids
    \big(
      \HomotopyQuotient
        { \SmoothManifold }
        { G }
      ,\,
      \mathbf{B} \mathcal{G}
    \big)
    \;\;
      \simeq
    \;\;
    \Localization{\WeakHomotopyEquivalences}
    \SimplicialPresheaves(\CartesianSpaces)
    \Big(
      \big(\!\!
      \begin{tikzcd}
        \widehat{\SmoothManifold}
          \times_{{}_{\SmoothManifold}}
        \widehat{\SmoothManifold}
          \times
        G^{\mathrm{op}}
        \ar[
          r,
          shift left=3pt,
          "{
           \mathrm{pr}_1
          }"
        ]
        \ar[
          r,
          shift right=3pt,
          "{
            \mathrm{pr}_3 \cdot \mathrm{pr}_2
          }"{swap}
        ]
        &
       \widehat{\SmoothManifold}
      \end{tikzcd}
    \!\!\big)
    ,\,
    \overline{W}\mathcal{G}
    \Big)
    \,,
  \end{equation}
  where in the second argument we used Lemma \ref{InfinityGroupsPresentedByPresheavesOfSimplicialGroups}.
  The elements on the right may be thought of as
  {\it $G$-equivariant {\v C}ech $\mathcal{G}$-cocycles on $\SmoothManifold$}.
  These will serve to represent equivariant principal $\infty$-bundles
  in \cref{AsInfinityBundlesInternalToSliceOverBG}
  and \cref{EquivariantLocalTrivializationIsImplies} below.
\begin{equation}
  \begin{tikzcd}[column sep={between origins, 30pt}]
    &&&
    \mathclap{
      \HomotopyQuotient
        { \TopologicalSpace }
        { G }
    }
    \ar[
      rrrrrrr,
      start anchor={[xshift=10pt]},
      end anchor={[xshift=-5pt]},
      "{
        \mbox{
          \tiny
          \color{greenii}
          \bf
          modulating map
        }
      }"
    ]
    &&&&&&&
    \mathclap{
      \mathbf{B}\Gamma
    }
    &&
    \mathrlap{
      \in \, \SmoothInfinityGroupoids
    }
    \\
    &&&
    \mathclap{
      \SimplicialNerve
      \big(
        \widehat{\TopologicalSpace}
          \times_{\TopologicalSpace}
        \widehat{\TopologicalSpace}
        \times G
        \rightrightarrows
        \widehat{\TopologicalSpace}
      \big)
    }
    \ar[
      rrrrrrr,
      start anchor={[xshift=47pt]},
      end anchor={[xshift=-25pt]},
      "{
        \mbox{
          \tiny
          \color{greenii}
          \bf
          equivariant {\v C}ech cocycle
        }
      }"
    ]
    &&&&&&&
    \mathclap{
      \SimplicialNerve
      \big(
        \Gamma \rightrightarrows \ast
      \big)
    }
    &&
    \mathrlap{
      \in \,
      \SimplicialPresheaves(\CartesianSpaces)
    }
    \\
    &[-10pt]
    \phantom{g\cdot}(x,j)
    \ar[dr]
    \ar[rrr]
    &
    &
    &[-10pt]
    g\cdot(x,j)
    \ar[dr]
    &&&
    &[-10pt]
    \phantom{g\cdot}\phantom{(x,j)}\mathllap{\bullet\;\;\;\;}
    \ar[
      dr,
      "{ \gamma_{jk}(x) }"{sloped}
    ]
    \ar[
      rrr,
      start anchor={[xshift=-6pt]},
      end anchor={[xshift=+8pt]},
      "{ \rho_g(x,j) }"
    ]
    &
    &
    &[-10pt]
    \phantom{g\cdot(x,j)}\mathllap{\bullet\;\;\;\;\;\;}
    \ar[
      dr
    ]
    \\
    &&
    (x,k)
    \ar[rrr]
    &&&
    g\cdot(x,k)
    &
    \hspace{10pt}
    \longmapsto
    \hspace{-8pt}
    &
    &&
    \phantom{(x,k)}\mathllap{\bullet\;\;\;\;}
    \ar[
      rrr,
      start anchor={[xshift=-6pt]},
      end anchor={[xshift=+4pt]}
    ]
    &&&
    \phantom{g\cdot(x,k)}\mathllap{\bullet\;\;\;\;\;\;\;\;}
    &&
    \\
    \phantom{g\cdot}(x,i)
    \ar[uur, "{ \big((x,i,j), \NeutralElement\big) }"{sloped}]
    \ar[urr]
    \ar[rrr, "{ \big( (x,i,i), g \big) }"{swap}]
    &&&
    g\cdot(x,i)
    \ar[uur, crossing over]
    \ar[urr]
    &&&&
    \phantom{g\cdot} \phantom{(x,i)}\mathllap{\bullet\;\;\;\;\;}
    \ar[
      uur,
      "{ \gamma_{i j}(x) }"{sloped}
    ]
    \ar[
      urr,
      "{ \gamma_{ik}(x) }"{sloped, swap}
    ]
    \ar[
      rrr,
      "{ \rho_g(x,i) }"{swap},
      start anchor={[xshift=-8pt]},
      end anchor={[xshift=+8pt]}
    ]
    &&&
    \phantom{g\cdot(x,i)}\mathllap{\bullet\;\;\;\;\;\;}
    \ar[uur, crossing over]
    \ar[
      urr
    ]
    &&&&
  \end{tikzcd}
\end{equation}
\end{remark}

\medskip

\noindent
{\bf Smooth path $\infty$-groupoids.}
It is a familiar fact that the singular simplicial complex\footnote{
  This is the standard singular simplicial set often denoted
  $\mathrm{Sing}(-)$ or similar. Even though this notation is
  classical and standard,
  it clashes badly with any appropriate notation for ``singularities''
  in the sense of orbifolds (as well as of schemes, stacks, etc.)
  that we need below
  (in Def. \ref{Singularities}, Def. \ref{GEquivariantAndGloballyEquivariantHomotopyTheories}
  and
  Prop. \ref{GSingularCohesiveInfinityTopos}).
  Moreover, in the generalization of this
  classical construction from topological spaces to
  diffeological spaces (Ntn. \ref{DiffeologicalSingularComplex} below)
  and further to smooth $\infty$-groupoids
  (Def. \ref{SmoothPathInfinityGroupoid} below) it becomes
  increasingly manifest that the actual {\it conceptual}
  nature of the singular simplicial complex is as a model
  for the {\it path $\infty$-groupoid} of a space.
  For these reasons we choose the non-standard but more
  suggestive notation $\SingularSimplicialComplex(-)$
  instead of $\mathrm{Sing}(-)$.
}
$\SingularSimplicialComplex(\TopologicalSpace)$ of a topological space
$\TopologicalSpace$, being a Kan complex and as such representing an
$\infty$-groupoid \eqref{InfinityGroupoids} may be thought of as the
{\it $\infty$-groupoid of paths} in $\TopologicalSpace$, whose objects
are the points of $\TopologicalSpace$ and whose
$n+1$-morphisms are the continuous paths between $n$-dimensional paths
in $\TopologicalSpace$. Here we discuss how this phenomenon is subsumed
by the concept of $\infty$-groupoids of smooth paths not just in
D-topological spaces, but more generally in smooth $\infty$-groupoids
(Def. \ref{SmoothPathInfinityGroupoid} below).

\medskip
Below, Prop. \ref{SmoothShapeGivenBySmoothPathInfinityGroupoid}
\cite{BBP19} (following \cite{Pavlov14}, see also \cite[\S 3]{Bunk20} and \cite[Thm. B]{Clough21} -- a  precursor argument is in
\cite[Lem. 7.5]{BunkeNikolausVoelkl13})
says that
these $\infty$-groupoids of smooth paths in smooth $\infty$-groupoids
still reflect their smooth shape.
This a cornerstone fact for analysis of the cohesion of smooth $\infty$-groupoids.
In particular, it implies the {\it smooth Oka principle} (Thm. \ref{SmoothOkaPrinciple} below)
which in turn implies the
fundamental classification theorems for smooth principal $\infty$-bundles
(Prop. \ref{ClassificationOfPrincipalInfinityBundles} below)
and for suitable equivariant principal $\infty$-bundles (Thm. \ref{ClassificationOfGammaEquivariantPrincipalBundlesForGammaWithTruncatedClassifyingShape} below)
over smooth manifolds, as we show below in \cref{AsInfinityBundlesInternalToSliceOverBG}.

\medskip
Conceptually, this result (Prop. \ref{SmoothShapeGivenBySmoothPathInfinityGroupoid})
says that morphisms
out of the shape of
a smooth $\infty$-groupoid $X$ encode
{\it parallel transport of flat $\infty$-connections}
(\cite[p. 22]{SSS09}, see also \cite{SchreiberShulman14})
or equivalently
{\it local systems of constant $\infty$-coefficients} (\cite[\S 2.2]{FSS20CharacterMap})
\vspace{-2mm}
$$
  \overset{
    \mathclap{
    \raisebox{4pt}{
      \tiny
      \color{darkblue}
      \bf
      \begin{tabular}{c}
        smooth $\infty$-groupoid of
        \\
        smooth higher paths in $X$
      \end{tabular}
    }
    }
  }{
    \SmoothSingularSimplicialComplex(X)
  }
  \qquad \simeq\;\;\;
  \underset{
    \mathclap{
    \raisebox{-4pt}{
      \tiny
      \color{darkblue}
      \bf
      \begin{tabular}{c}
        shape of $X$
      \end{tabular}
    }
    }
  }{
    \shape \, X
  }
  \quad
  \xrightarrow{
    \mbox{
      \tiny
      \bf
      \color{greenii}
      \begin{tabular}{c}
        parallel transport of
        \\
        flat $\infty$-connection
      \end{tabular}
    }
  }
  \mathbf{B}\mathcal{G}
  \;\;\;\;\;\;\;\;\;\;\;
  \xleftrightarrow{\;\; \mbox{\tiny \eqref{AdjunctionAndHomEquivalence}} \;\;}
  \;\;\;\;\;\;\;\;\;\;\;
  X
  \xrightarrow{
    \mbox{
      \tiny
      \bf
      \color{greenii}
      \begin{tabular}{c}
        local system of
        \\
        constant $\infty$-coefficients
      \end{tabular}
    }
  }
  \flat
  \mathbf{B}\mathcal{G}
  \,.
$$

\vspace{-2mm}
\noindent Conversely, this means that in any cohesive $\infty$-topos the
shape modality encodes, in particular, flat $\infty$-connections on
principal $\infty$-bundles. This was the key motivation for
axiomatic cohesion on $\infty$-toposes in \cite{SSS09}\cite{dcct}.\footnote{
Notice that this relation
of cohesive $\infty$-toposes to higher gauge theory is invisible in
the cohesive 1-toposes considered in \cite{Lawvere07}:
The shape modality of a cohesive $(n,1)$-topos
reflects parallel transport for flat $(n-1)$-connections.
Here 1-connections are the traditional Cartan-Ehresmann connections,
while 0-connections are just functions,
and flat 0-connections are just locally constant functions.
Hence it takes at least a cohesive $(2,1)$-topos
(such as $\SmoothGroupoids$)
for its cohesion to axiomatically reflect ordinary connections/gauge fields.}

\begin{definition}[Smooth path $\infty$-groupoid]
\label{SmoothPathInfinityGroupoid}
For
$
  A
    \,\in\,
  \SmoothInfinityGroupoids
    \,\coloneqq \,
  \InfinitySheaves
  (\CartesianSpaces)
$,
we write
\vspace{-2mm}
\begin{equation}
  \label{SmoothSingularComplexOfInfinityStack}
  \hspace{-1cm}
  \SingularSimplicialComplex(A)
  \;\coloneqq\;
  {
    \colimit{ [n] \in \Delta^{\mathrm{op}} }
    A(\SmoothSimplex{n})
  }
  \;\;\;
    \in
  \;
  \InfinityGroupoids
    \hspace{2cm}
    \mathclap{
    \raisebox{0pt}{
      \tiny
      \color{darkblue}
      \bf
      \begin{tabular}{c}
        $\infty$-groupoid of
        smooth paths \\
        and their
        gauge transformations.
      \end{tabular}
    }
    }
  \end{equation}

\vspace{-2mm}
\noindent
for the $\infty$-groupoid whose $n$-morphisms are $k$-gauge transformations of
$(n-k)$-dimensional smooth paths in $X$,
in the sense of plots from the smooth extended $n$-simplices
from Def. \ref{SmoothExtendedSimplicies}.
\end{definition}

\begin{example}[Smooth shape of diffeological space given by diffeological singular simplicial complex]
  \label{SmoothShapeOfDiffeologicalSpaceGivenByDiffeologicalSingularSimplicialComplex}
  For $X \,\in\, \DiffeologicalSpaces \xhookrightarrow{\;} \SmoothInfinityGroupoids\,$,
  the smooth path $\infty$-groupoid \eqref{PathInfinityGroupoid}
  coincides with the diffeological singular simplicial complex
  \eqref{SmoothSingularSimplicialSet} up to weak homotopy equivalence:
  \vspace{-2mm}
  \begin{equation}
    \label{SingularSimplicialComplexOnDiffeologicalSpaceAsInfinityGroupoid}
    \SingularSimplicialComplex(X)
    \,\in\,
    \SimplicialSets
    \xrightarrow{\;\; \SimplicialLocalization{\WeakHomotopyEquivalences} \;\;}
    \InfinityGroupoids\;.
  \end{equation}
\end{example}

\begin{proposition}[Smooth cohesive shape given by smooth path $\infty$-groupoid
 {\cite{Pavlov14}\cite{BBP19}}]
  \label{SmoothShapeGivenBySmoothPathInfinityGroupoid}
  $\,$

  \noindent
  {\bf (i)} For $A \in \SmoothInfinityGroupoids$,
  the {\it smooth path $\infty$-groupoid} of $A$
  (Def. \ref{SmoothPathInfinityGroupoid}), namely the presheaf
  \vspace{-1mm}
  \begin{equation}
    \label{PathInfinityGroupoid}
    \SingularSimplicialComplex
    (A)
    \;:\;
    U
    \;\longmapsto\;
    \SingularSimplicialComplex
    \left(
      [U,\, A]
    \right)
    \;=\;
    \colimit{ [n] \in \Delta^{\mathrm{op}} }
    \left(
    A
    (
      \SmoothSimplex{n}
      \times
      U
    )
    \right)
    \;\;
    \in
    \;
    \InfinityGroupoids
      \end{equation}

  \vspace{-2mm}
  \noindent
  of smoothly parameterized paths (Def. \ref{SmoothPathInfinityGroupoid})
  in $A$, is in fact an $\infty$-sheaf
  \vspace{-1mm}
  $$
    \SmoothSingularSimplicialComplex(A)
    \;\in\;
    \InfinitySheaves(\SmoothManifolds)
    \;\simeq\;
    \SmoothInfinityGroupoids\;.
  $$

  \vspace{-1mm}
  \noindent
  {\bf (ii)} As such, this is equivalent to the cohesive shape \eqref{TheModalitiesOnACohesiveInfinityTopos}
  of $A$:
  \vspace{-1mm}
  $$
    \shape
    \,
    A
    \;\simeq\;
    \SmoothSingularSimplicialComplex(A)
    \,.
  $$

  \vspace{-1mm}
  \noindent
{\bf (iii)} In particular, this means,
  by evaluation on $\mathbb{R}^0 \,\in\, \CartesianSpaces$,
  that
  \vspace{-1mm}
  \begin{equation}
    \label{SmoothShapeIsSmoothSingularSmimplicialComplex}
    \Shape(A)
    \;\simeq\;
    \SingularSimplicialComplex(A)
    \;\;\;\;\;\;\;
    \mbox{and}
    \;\;\;\;\;\;\;
    \SmoothSingularSimplicialComplex(A)
    \;\simeq\;
    \Discrete
    \,
    \SingularSimplicialComplex(A)
    \,.
  \end{equation}
\end{proposition}

\begin{example}[Connected components of smooth shape are concordance classes]
  \label{ConnectedComponentsOfShapeAreConcordanceClasses}
  For $A \,\in\, \SmoothInfinityGroupoids$, there is a natural isomorphism
  \vspace{-3mm}
  \begin{equation}
    \label{ZeroTrancationOfShapeEquivalentConcordanceClassesOfZeroTruncation}
    \def\arraystretch{1.6}
    \begin{array}{lll}
      \Truncation{0}
      \,
      \shape
      \,
      A
      &
      \;\simeq\;
      \Truncation{0}
      \,
      \colimit{ [n] \in \Delta^{\mathrm{op}} }
      \,
      A
      &
      \proofstep{ by Prop. \ref{SmoothShapeGivenBySmoothPathInfinityGroupoid} }
      \\[-2pt]
      & \;\simeq\;
      \colimit{ [n] \in \Delta^{\mathrm{op}} }
      \,
      \Truncation{0}
      \,
      A
      &
      \proofstep{ by \eqref{InfinityAdjointPreservesInfinityLimits} }
      \\[10pt]
      & \;\simeq\;
      \big(
        \Truncation{0} A(\ast)
      \big)
      \big/
      \left(
        \Truncation{0}
        \,
        A(\SmoothSimplex{1})
      \right)
      &
      \proofstep{ by Lem. \ref{ColimitsOverSimplicialDiagramsOfZeroTruncatedObjectsAreCoequalizers} }
         \end{array}
  \end{equation}
  between the connected components of the shape of $A$
  and its {\it concordance classes}
  \vspace{-2mm}
  \begin{equation}
    \label{ProjectionFromIsomorphismClassesToConcordanceClasses}
    \begin{tikzcd}
    \Truncation{0}
    \,
    A(\ast)
    \;=:\;
    \IsomorphismClasses{
      A
    }
    \ar[r, ->>]
    &
    \ConcordanceClasses{
      A
    }
    \;\coloneqq\;
    \big(
      \Truncation{0}
      \,
      A(\ast)
    \big)
    \big/
    \big(
      \Truncation{0}
      \,
      A(\SmoothSimplex{1})
    \big)
    \,.
    \end{tikzcd}
  \end{equation}

  \vspace{-2mm}
  \noindent
  Notice that concordance classes are a coarsening, in general, of
  isomorphism classes on the point by isomorphism classes on the interval.
  In particular, if two smooth $\infty$-groupoids are equivalent
  after 0-truncation, then they already have the same concordance classes:
  \vspace{-2mm}
  \begin{equation}
    \label{SmoothInfinityGroupoidsWithCoinciding0TruncationHaveTheSameConcordanceClasses}
    \Truncation{0}
    \,
    A
    \;
    \simeq
    \;
    \Truncation{0}
    \,
    A'
    \;\;\;
    \in
    \;
    \SmoothInfinityGroupoids
    \;\;\;\;\;\;\;\;\;\;\;
    \Rightarrow
    \;\;\;\;\;\;\;\;\;\;\;
    \ConcordanceClasses{A}
    \;\simeq\;
    \ConcordanceClasses{A'}
    \;\;\;
    \in
    \;
    \Sets
    \,.
  \end{equation}
\end{example}

\begin{proposition}[Smooth cohesive shape of topological spaces is their weak homotopy type]
  \label{SmoothShapeOfTopologicalSpacesIsTheirWeakHomotopyType}
  The weak homotopy type of any topological space
  is naturally identified with the cohesive shape
  \eqref{TheModalitiesOnACohesiveInfinityTopos}
  seen in $\SmoothInfinityGroupoids$ (Ntn. \ref{SmoothInfinityGroupoids})
  of its continuous-diffeological incarnation \eqref{CategoryOfDeltaGeneratedSpaces}:
  \vspace{-2mm}
  $$
    \begin{tikzcd}[row sep=small]
      \kTopologicalSpaces
      \ar[
        rr,
        "{
          \ContinuousDiffeology
        }"{above},
        "{ }"{name=s, below, pos=.6}
      ]
      \ar[
        dr,
        "{ \SingularSimplicialComplex }"{below, xshift=-6pt},
        "{ }"{name=t, above, pos=.8}
      ]
      &&
      \SmoothInfinityGroupoids
      \ar[
        dl,
        "{ \Shape }"{below, xshift=4pt},
      ]
      \\
      &
      \InfinityGroupoids
      \ar[
        from=s,
        to=t,
        Rightarrow,
        "\sim"{sloped, below}
      ]
    \end{tikzcd}
  $$
\end{proposition}
\begin{proof}
  For $X \,\in\, \TopologicalSpaces$, we need to exhibit a natural
  equivalence
  \vspace{-2mm}
  $$
    \begin{tikzcd}
      \Shape \circ \ContinuousDiffeology(X)
      \ar[
        r,
        "\sim"{below}
      ]
      &
      \SingularSimplicialComplex(X)\;.
    \end{tikzcd}
  $$

  \vspace{-3mm}
  \noindent But we have:
  \vspace{-.5mm}
  $$
    \begin{array}{lll}
      \Shape \circ \ContinuousDiffeology(X)
      &
      \;\simeq\;
      \SingularSimplicialComplex \circ \ContinuousDiffeology(X)
      &
      \mbox{\small
        by Prop. \ref{SmoothShapeGivenBySmoothPathInfinityGroupoid}
      }
      \\
      & \;\simeq\;
      \SingularSimplicialComplex(X)
      &
      \mbox{\small by Prop. \ref{DiffeologicalShapeOfContinuousDiffeologyIsWeakHomotopyType}
      with Ex. \ref{SmoothShapeOfDiffeologicalSpaceGivenByDiffeologicalSingularSimplicialComplex}}.
    \end{array}
  $$

  \vspace{-6mm}
\end{proof}

More generally:
\begin{proposition}[Smooth shape of good simplicial spaces is weak homotopy type of their realization]
  \label{SmoothShapeOfGoodSimplicialSpacesIsWeakHomotopyTypeOfTheirRealization}
  The weak homotopy type
  of
  the topological realization
  $\TopologicalRealization{}{\TopologicalSpace_\bullet}$
  (Ntn. \ref{TopologicalRealizationFunctors})
  of
  any good simplicial topological space $\TopologicalSpace_\bullet$
  (Def. \ref{GoodSimplicialTopologicalSpace})
  is naturally identified with the cohesive shape
  \eqref{TheModalitiesOnACohesiveInfinityTopos}
  seen in $\SmoothInfinityGroupoids$ (Ntn. \ref{SmoothInfinityGroupoids})
  of its degreewise
  continuous-diffeological incarnation \eqref{CategoryOfDeltaGeneratedSpaces}:
  \vspace{-2.5mm}
  $$
    \begin{tikzcd}[row sep=small, column sep=large]
      \SimplicialTopologicalSpaces_\good
      \ar[
        rr,
        "{
          \ContinuousDiffeology\left((-)_\bullet\right)
        }"{above},
        "{ }"{name=s, below, pos=.6}
      ]
      \ar[
        dr,
        "{
          \SingularSimplicialComplex
          \,
          \TopologicalRealization{}{-}
        }"{swap},
        "{ }"{name=t, pos=.8}
      ]
      &&
      \SmoothInfinityGroupoids \;.
      \ar[
        dl,
        "{ \Shape }"{below, xshift=4pt},
      ]
      \\
      &
      \InfinityGroupoids
      \ar[
        from=s,
        to=t,
        Rightarrow,
        "\sim"{sloped, below}
      ]
    \end{tikzcd}
  $$
\end{proposition}
\begin{proof}
We have the following sequence of natural equivalences in
$\InfinityGroupoids$:
\vspace{-2mm}
$$
  \begin{array}{lll}
    \Shape
    \left(
      \ContinuousDiffeology(\TopologicalSpace)_\bullet
    \right)
    &
    \;\simeq\;
    \Shape
    \,
    \colimit{[n]\in\Delta^\op}
    \,
    \ContinuousDiffeology(\TopologicalSpace_n)
    &
    \proofstep{ by Lem. \ref{SimplicialPresheavesModelHomotopyColimitOfTheirComponentSheaves} }
    \\
    & \;\simeq\;
    \colimit{[n]\in\Delta^\op}
    \,
    \Shape
    \,
    \ContinuousDiffeology(\TopologicalSpace)_n
    &
    \proofstep{ by \eqref{InfinityAdjointPreservesInfinityLimits} }
    \\
    & \;\simeq\;
    \colimit{[n]\in\Delta^\op}
    \,
    \SingularSimplicialComplex (\TopologicalSpace_n)
    &
    \proofstep{ by Prop. \ref{SmoothShapeOfTopologicalSpacesIsTheirWeakHomotopyType} }
    \\
    & \;\simeq\;
    \SingularSimplicialComplex
    \TopologicalRealization{}{
      \TopologicalSpace_\bullet
    }
    &
    \proofstep{ by Rem. \ref{TopologicalRealizationOfGoodSimplicialSpacesModelsTheirHomotopyColimit}}
    \,.
  \end{array}
$$

\vspace{-7mm}
\end{proof}

\begin{proposition}[Shape of smooth $\infty$-topos over continuous-diffeology is weak homotopy type]
  \label{SmoothShapeOfSmoothInfinityToposOverContinuousDiffeologyIfWeakHomotopyType}
  For $\TopologicalSpace \,\in\, \TopologicalSpaces$, its weak homotopy type coincides
  with the shape (in the sense of Def. \ref{ShapeOfAnInfinityTopos})
  of the slice $\infty$-topos (Ntn. \ref{SliceInfinityTopos})
  over its continuous-diffeological incarnation
  $\ContinuousDiffeology(\TopologicalSpace) \,\in\, \SmoothInfinityGroupoids$
  \eqref{CategoryOfDeltaGeneratedSpaces}
  in that
  \vspace{-2mm}
  $$
    \Shape
    \big(
      (\SmoothInfinityGroupoids)_{/\ContinuousDiffeology(\TopologicalSpace)}
    \big)
    \;\;
      \simeq
    \;\;
    \SingularSimplicialComplex(\TopologicalSpace)
    \;\;
    \in
    \;\;
    \InfinityGroupoids
    \xhookrightarrow{\;}
    \ProObjects(\InfinityGroupoids)
    \,.
  $$
\end{proposition}
\begin{proof}
$\,$
\vspace{-3mm}
  $$
    \begin{array}{lll}
      \Shape
      \big(
        (\SmoothInfinityGroupoids)_{/\ContinuousDiffeology(\TopologicalSpace)}
      \big)
      &
      \;\simeq\;
      \Shape
      \left(
        \ContinuousDiffeology(\TopologicalSpace)
      \right)
      &
      \proofstep{ by Prop. \ref{ShapeOfSliceOfCohesiveInfinityToposIsCohesiveShape} }
      \\
      &
      \;\simeq\;
      \SingularSimplicialComplex(\TopologicalSpace)
      &
      \proofstep{ by Prop. \ref{SmoothShapeOfTopologicalSpacesIsTheirWeakHomotopyType} }.
    \end{array}
  $$

  \vspace{-6mm}
\end{proof}

\begin{theorem}[Smooth Oka principle for maps out of smooth manifolds]
\label{SmoothOkaPrinciple}
$\,$

\noindent Let $A \in \SmoothInfinityGroupoids$ and
$X \in \SmoothManifolds \xhookrightarrow{} \SmoothInfinityGroupoids$\,.

\noindent {\bf (i)} The shape of the mapping stack \eqref{InternalHomAdjunction}
from $X$ to $A$ is
naturally equivalent the mapping stack from $X$ to the shape of $A$:
\vspace{-2mm}
$$
  \shape
  \;
  \Maps{}
    { X }
    { A }
  \;\;\simeq\;\;
  \Maps{}
    { \shape \, X }
    { \shape \, A }
  \;\;\;\in\;
  \SmoothInfinityGroupoids \;.
$$

  \vspace{-1mm}
\noindent
{\bf (ii)} This means, equivalently, that
\vspace{-2mm}
$$
  \Shape
  \left(
    \Maps{}
      { X }
      { A }
  \right)
  \;\simeq\;
  \InfinityGroupoids
  \left(
    \Shape(X),
    \,
    \Shape(A)
  \right)
  .
$$
\end{theorem}
The following proof is a direct consequence of Prop. \ref{SmoothShapeGivenBySmoothPathInfinityGroupoid} (due to \cite{Pavlov14}\cite{BBP19}). We later learned that an alternative proof appears as \cite[Thm. B/Cor. 6.4.8]{Clough21}.
\begin{proof}
Recalling that also $\SmoothManifolds$ is a site of definition for
$\SmoothInfinityGroupoids$ (by Rem. \ref{SmoothInfinityGroupoidsPresentedOverSmoothManifolds}), observe the following sequence of
natural equivalences of $\infty$-groupoids, for $U \in \SmoothManifolds$:
\vspace{-1mm}
$$
  \begin{array}{llll}
        \Maps{}
      { X }
      { \shape \,  A }
      (U)
      &
        \;\simeq\;
    \InfinitySheaves(\SmoothManifolds)
    \left(
      X \times U
      ,\,
      \shape \, A
    \right)
    &&
    \mbox{\small by Prop. \ref{MappingStacks}}
    \\
    & \;\simeq\;
    \InfinityPresheaves(\SmoothManifolds)
    \Big(\!
      X \times U,
      \,
      \big(
      U'
      \,\mapsto\,
      \underset{\underset{[n] \in \Delta^{\mathrm{op}}}{\longrightarrow}}{\mathrm{lim}}
      A\left(U' \times \SmoothSimplex{n}\right)
     \!\! \big)
    \!\!\Big)
    &&
    \mbox{\small by Prop. \ref{SmoothShapeGivenBySmoothPathInfinityGroupoid}}
    \\
    & \;\simeq\;
    \underset{\underset{[n] \in \Delta^{\mathrm{op}}}{\longrightarrow}}{\mathrm{lim}}
    \InfinityPresheaves(\SmoothManifolds)
    \Big(\!
      X \times U,
      \,
      \big(
      U'
      \,\mapsto\,
      A\left(U' \times \SmoothSimplex{n}\right)
      \!\!\big)
    \!\!\Big)
    &&
    \mbox{\small by \cite[Cor. 5.1.2.3]{Lurie09HTT}}
    \\
    & \;\simeq\;
    \underset{\underset{[n] \in \Delta^{\mathrm{op}}}{\longrightarrow}}{\mathrm{lim}}
    \InfinityPresheaves(\SmoothManifolds)
    \big(
      X \times U
      ,\,
      \Maps{}
        { \SmoothSimplex{n} }
        { A }
    \big)
    &&
    \mbox{\small by Prop. \ref{MappingStacks}}
    \\
    & \;\simeq\;
    \underset{\underset{[n] \in \Delta^{\mathrm{op}}}{\longrightarrow}}{\mathrm{lim}}
    \InfinityPresheaves(\SmoothManifolds)
    \big(
      \SmoothSimplex{n} \times X \times U,
      \,
      A
    \big)
    &&
    \mbox{\small by Prop. \ref{MappingStacks}}
    \\
    & \;\simeq\;
    \underset{\underset{[n] \in \Delta^{\mathrm{op}}}{\longrightarrow}}{\mathrm{lim}}
    \big(
      \Maps{}
        { X }
        {A}
      \left( \SmoothSimplex{n} \times U  \right)
    \!\!\big)
    &&
    \mbox{\small by Prop. \ref{MappingStacks}}
    \\[1pt]
    & \;\simeq\;
    \left(
      \shape \;
      \Maps{}
        { X }
        { A }
    \right)(U)
    &&
    \mbox{\small by Prop. \ref{SmoothShapeGivenBySmoothPathInfinityGroupoid}}
    \,.
  \end{array}
$$

\vspace{-1mm}
\noindent By the $\infty$-Yoneda lemma (Prop. \ref{InfinityYonedaLemma}), the composite
of these natural equivalences yields the first claim.
From this the second claim follows by cohesion
(the definition $\raisebox{1pt}{\textesh} \;\simeq\; \Discrete \circ \Shape$,
the adjunction $\Shape \dashv \Discrete$ and the fully faithfulness of $\Discrete$).
\end{proof}

\medskip

\noindent
{\bf Resolvable orbi-singularities.}
The smooth Oka principle (Thm. \ref{SmoothOkaPrinciple})
applies, a priori, in the generality of domains that are smooth manifolds
$\SmoothManifold$
with maps $\SmoothManifold \to A$ into any coefficient objects $A \,\in\, \SmoothInfinityGroupoids$.
However, for particular coefficients it may happen that maps into them do not distinguish
a given domain $\infty$-stacks from its approximation by some smooth manifold,
in which case the smooth Oka principle will generalize to that particular situation as well.

\medskip
A first class of examples of this generalization occurs when
the domain is the moduli stack $\mathbf{B}G$ of a finite group $G$,
hence the orbifold $\ast \sslash G$ consisting of a single point which is
a $G$-singularity.
We may observe that
a canonical choice of approximation in this case is given by the
``blow-up'' of this point (e.g. \cite[\S 29.2]{Tu20},
Ntn. \ref{ResolvableOrbiSingularities} below)
to a
{\it smooth spherical space form} (Def. \ref{SmoothSphericalSpaceForm} below)
namely to the quotient manifold $\SmoothSphere{\,n+2}/G$ of a smooth $n$-sphere
(possibly exotic) by a smooth and free action of the group $G$:
For any given coefficients $A$, the distinction between $\mathbf{B}G$ and
$\SmoothSphere{\,n+2}/G$
is carried by morphisms $\SmoothSphere{\,n+2} \to A$, and hence disappears whenever $A$,
or at least its shape, is suitably {\it $n$-truncated}
(see Ntn. \ref{CohesiveGroupsWithTruncatedClassifyingShape} below).

\medskip
Therefore, if $G$-singularities are resolvable in this sense,
and the coefficients are suitably truncated,
we obtain an {\it orbi-smooth Oka principle},
which is Thm. \ref{OrbiSmoothOkaPrinciple} below.
Here we first recall relevant facts about
free actions of finite groups on spheres and the resulting
smooth spherical space forms.

\begin{definition}[Smooth spherical space forms\footnote{Historically and by default,
the term ``spherical space form'' refers to Riemannian geometry, where the $G$-action
on a round $\RiemannianSphere{n+2}$
is required to be (free and) by isometries
(e.g. \cite{Wolf11}\cite{Allcock18}).
In generality, if the action is just required to be (free and) continuous,
one speaks of ``topological spherical space forms'' \cite{MadsenThomasWall76},
see review in \cite{Hambleton14}.
The terminology ``smooth spherical space form'' for the case of concern here,
where the action is (free and) by diffeomorphisms seems to
originate with \cite{Madsen78} following the construction of interesting
examples in \cite{MadsenThomasWall76}, see Prop. \ref{ExistenceOfSmoothSphericaSpaceForms}.}]
  \label{SmoothSphericalSpaceForm}
  If $G \,\in\, \Groups(\FiniteSets)$ and
  $ G \acts \;  \SmoothSphere{\,n+2} \,\in\, \Actions{G}(\SmoothManifolds)$
  is a free action on a smooth $n+2$-sphere (possibly exotic),
  then the quotient smooth manifold
  \vspace{-1mm}
  $$
    \SmoothSphere{\,n+2} \!\sslash\! G
    \;\;\simeq\;\;
    \SmoothSphere{\,n+2}/G
    \;\;\;
    \in
    \;
    \SmoothManifolds
    \xhookrightarrow{\;}
    \SmoothInfinityGroupoids
    \,,
  $$

  \vspace{-1mm}
  \noindent
  is a {\it smooth spherical space form}.
\end{definition}

\newpage

\begin{proposition}[Existence of smooth spherical space forms]
\label{ExistenceOfSmoothSphericaSpaceForms}
For $G$ a finite group, the following conditions are equivalent:

\vspace{-3mm}
\begin{enumerate}[{\bf (i)}]
 \setlength\itemsep{-3pt}

\item
For each prime number $p$, if $H \subset G$ is a subgroup of order
$p^2$ or $2 p$, then $H$ is a cyclic group.

\item
$G$ has a continuous free action on a topological $d$-sphere
$\TopologicalSphere{d}$ for some $d \in \mathbb{N}$.

\item
For each $n \in \mathbb{N}$, there exists $d \geq n+2$
and a smooth manifold structure $\SmoothSphere{\,d}$ on the $d$-sphere
(possibly exotic) such that $G$ has a smooth free action on $\SmoothSphere{\,d}$,
hence such that $\SmoothSphere{\,d}/G$ is a smooth spherical space form
(Def. \ref{SmoothSphericalSpaceForm}).
\end{enumerate}
\vspace{-.2cm}

\end{proposition}

\vspace{-5mm}
\begin{proof}
  The equivalence of the first two statements is
  \cite[Thm. 0.5-0.5]{MadsenThomasWall76}.
  The implication of the third from the second statement is
  \cite[Thm. 5]{MadsenThomasWall83}, which asserts, in more detail,
  that $G$ has a smooth action (at least) on smooth spheres of dimension
  $2^{ 2 k \cdot e(G) - 1 }$, for all $k \in \mathbb{N}$,
  where $e(G) \,\in\, \mathbb{N}_+$ is the
  {\it Artin-Lam exponent} of $G$ \cite{Lam68}\cite{Yamauchi70},
  the only feature of which that matters here is that it is positive.
  Both of these statements are reviewed in \cite[Thm. 6.1]{Hambleton14}.
  Finally, the implication of the second from the third statement is trivial.
\end{proof}

\begin{notation}[Resolvable orbi-singularities]
  \label{ResolvableOrbiSingularities}
  {\bf (i)} If a finite group $G$
  satisfies the equivalent conditions of
  Prop. \ref{ExistenceOfSmoothSphericaSpaceForms},
  we say that there {\it exist resolutions of $G$-singularities}
  or that {\it $G$-singularities are resolvable}.

 \noindent {\bf (ii)} We write
 \vspace{-2mm}
 \begin{equation}
    \label{FiniteGroupsWithResolvableSingularities}
    \begin{tikzcd}[row sep=small]
      \overset{
        \mathclap{
        \raisebox{5pt}{
          \tiny
          \color{darkblue}
          \bf
          \def\arraystretch{.9}
          \begin{tabular}{c}
            finite groups that have
            \\
            resolvable singularities
          \end{tabular}
        }
        }
      }{
        \Groups(\FiniteSets)_{\Resolvable}
      }
      \ar[r, hook]
      \ar[d]
      &
      \overset{
        \mathclap{
        \raisebox{5pt}{
          \tiny
          \color{darkblue}
          \bf
          \def\arraystretch{.9}
          \begin{tabular}{c}
            finite groups with covers that
            \\
            have resolvable singularities
          \end{tabular}
        }
        }
      }{
        \Groups(\FiniteSets)_{\CoverResolvable}
      }
      \ar[r, hook]
      \ar[d]
      &
      \Groups(\FiniteSets)
      \ar[d]
      \\
      \ResolvableSingularities
      \ar[r, hook]
      &
      \CoverResolvableSingularities
      \ar[r, hook]
      &
      \Singularities
    \end{tikzcd}
  \end{equation}

  \vspace{-1mm}
  \noindent
  for the full subcategories on finite groups $G$ and on
  orbi-singularities $\orbisingularG$ (Ntn. \ref{Singularities}) on the resolvable ones,
  and, respectively, for those which admit a cover, hence a
  surjective homomorphism
      \vspace{-2mm}
  $$
  \begin{tikzcd} \widehat{G} \ar[r,->>, "p", pos=.35] &[-11pt] G \end{tikzcd}\! ,
  $$

  \vspace{-2mm}
  \noindent
  such that $\widehat{G}$ has resolvable singularities.

  \noindent
  {\bf (iii)}
  Given $G \in \Groups(\FiniteSets)_{\resolvable}$
  we say that the canonical morphism
  \vspace{-2mm}
  $$
    \begin{tikzcd}
      \SmoothSphere{\,d}/G
      \ar[rrr, "{ (\SmoothSphere{\,d} \to \ast)\,\sslash\, {G} }"]
      &&&
      \ast \sslash {G}
      \,=\,
      \mathbf{B}{G}
      \,,
    \end{tikzcd}
    $$

    \vspace{-6mm}
    $$
    \hspace{-1.5cm}
    \begin{tikzcd}
      \SmoothSphere{\,d}/{G}
      \ar[rrr, "{ \orbisingular (\SmoothSphere{\,d} \to \ast)\,\sslash\, {G} }"]
      &&&
      \orbisingularAny{G}
    \end{tikzcd}
 $$

  \vspace{-2mm}
  \noindent
  is a {\it blow-up} of the $\widehat{G}$-orbi-singularity.

  \noindent
  {\bf (iv)}
  Notice that the homotopy fiber of this morphism is the smooth sphere:
  \vspace{-2mm}
  $$
    \begin{tikzcd}[column sep=large]
      \SmoothSphere{\,d}
      \ar[rr]
      \ar[d]
      \ar[drr, phantom, "\mbox{\tiny\rm(pb)}"]
      &&
      \ast
      \ar[d]
      \\
      \SmoothSphere{\,d}/{G}
      \ar[rr, "{ (\SmoothSphere{\,d} \to \ast)\,\sslash\, {G} }"]
      &&
      \ast \sslash {G}
      \,.
    \end{tikzcd}
  $$

\end{notation}

\begin{example}[ADE-groups have (cover-)resolvable singularities]
  \label{ADEGroupsHaveSphericalSpaceForms}
  For $G \,\subset\, \SUTwo \,\simeq\, \SpOne$
  a finite subgroup,
  its left quaternion multiplication
  action on
  the unit sphere
  $\SmoothSphere{4k+3} \,\simeq\, S(\mathbb{H}^{k+1})$
  in $k+1$-Quaternion space
  is evidently smooth and free, for each $k \in \mathbb{N}$.

  \noindent {\bf (i)}
  For example, for $k = 0$ this is just the
  restriction of left multiplication action of
  $\SpOne \,\simeq\, \SmoothSphere{3}$ on itself; while
  for $k = 1$ this is the restriction of the $\SpOne$-action
  which exhibits the 7-sphere as a $\SpOne$-principal bundle over
  the 4-sphere via the quaternionic Hopf fibration $t_{\Quaternions}$
  (e.g. \cite{GluckWarnerZiller86}):
  \begin{equation}
    \label{QuaternionicHopfFibration}
    \begin{tikzcd}
      \SmoothSphere{3}
      \;\simeq\;
      \SpOne
      \ar[r, hook]
      &
      \SmoothSphere{7}
      \ar[out=180-66, in=66, looseness=3.5, "\scalebox{.77}{$\mathclap{
        G
      }$}"{description},shift right=1]
      \ar[
        d,
        "{
          t_{\Quaternions}
          \mathrlap{
            \mbox{
              \tiny
              \color{darkblue}
              \bf
              \begin{tabular}{c}
                quaternionic
                \\
                Hopf fibration
              \end{tabular}
            }
          }
        }"
      ]
      \\
      &
      \SmoothSphere{4}
    \end{tikzcd}
  \end{equation}

 \noindent {\bf (ii)}
 Indeed, by the famous {\it ADE-classification} of finite subgroups
 of $\SUTwo \simeq \SpOne$
 (for pointers to modern proofs see \cite[Rem. A.9]{HSS18})
 these must be, up to isomorphism, among the following list:

{\small
  \begin{center}
  \def\arraystretch{1.1}
  \begin{tabular}{cccc}
    \hline
    \begin{tabular}{c}
     \small \bf Label
    \end{tabular}
    &
    \begin{tabular}{c}
     \small  \bf
     Finite
     subgroup
     \\
     \bf of $\mathrm{SU}(2)$
    \end{tabular}
    &
    \begin{tabular}{c}
     \small \bf
     Order
    \end{tabular}
    &
     \begin{tabular}{c}
     \small  \bf Name of
      \\
     \small  \bf group
    \end{tabular}
    \\
    \hline
    \hline
    \rowcolor{lightgray}
    $\mathbb{A}_{\mathrlap{r}}$ &
    $\phantom{2}\mathbb{Z}_{\mathrlap{r+1}}$
    &
    $n$
    &
    Cyclic
    \\
    $\mathbb{D}_{\mathrlap{r+4}}$
    &
    $2\mathrm{D}_{\mathrlap{r+2}}$
    &
    $4 (r + 2)$
    &
    Binary dihedral
    \\
    \rowcolor{lightgray}
    $\mathbb{E}_{\mathrlap{6}}$
    & $2\mathrm{T}$
    &
    $24$
    &
    Binary tetrahedral
    \\
    $\mathbb{E}_{\mathrlap{7}}$
    &
    $2\mathrm{O}$
    &
    $ 48 $
    &
    Binary octahedral
    \\
    \rowcolor{lightgray}
    $\mathbb{E}_{\mathrlap{8}}$
    &
    $2\mathrm{I}$
    &
    $ 120 $
    &
    Binary icosahedral
    \\
    \hline
  \end{tabular}
  \end{center}
  }
Incidentally, for all these finite $G \,\subset\, \SpOne$
their integral group cohomology is
(e.g. \cite[p. 12]{EpaGanter16}):
\begin{equation}
  \label{IntegralGroupCohomologyOfADESubgroups}
  H^n_{\mathrm{Grp}}(G;\, \Integers)
  \;\simeq\;
  \left\{
  \begin{array}{lcl}
    \mathbb{Z} &\vert& n = 0
    \\
    G^{\mathrm{ab}} \,=\, G/[G,G] & \vert &  n = 2 \;\mathrm{mod}\; 4
    \\
    \CyclicGroup{\mathrm{ord}(G)}
    &\vert& \mbox{$n$ a positive multiple of 4}
    \\
    0 &\vert& \mbox{otherwise}.
  \end{array}
  \right.
\end{equation}
and this means
(by the exact sequence
$\RealNumbers \to \Integers \to \CircleGroup$
and using that $H^{\geq 1}_{\mathrm{Grp}}(G;\, \RealNumbers) \,=\, 0$ for all finite
groups) in particular that
the abelianization of $G$ is identified with its group of
multiplicative characters:
\begin{equation}
  \label{FirstUOneGroupCohomologyOfADESubgroup}
  \Homs{}{G}{\CircleGroup}
  \;\simeq\;
  H^1_{\mathrm{Grp}}(G;\, \CircleGroup)
  \;\simeq\;
  H^2_{\mathrm{Grp}}(G;\, \Integers)
  \;\simeq\;
  G^{\mathrm{ab}}
  \,.
\end{equation}

By the above ADE-classification,
the subgroups of $G \,\subset\, \SpOne$ which are not cyclic can only be of order

\vspace{1mm}
  $\bullet$ $2^2 \cdot (r + 2)$ (for $2D_{r+2}$, $r \in \mathbb{N}$)

\vspace{1mm}
  $\bullet$ $24 = 2^3 \cdot 3$ (for $2T$)

\vspace{1mm}
  $\bullet$ $48 = 2^4 \cdot 3$ (for $2 O$)

\vspace{1mm}
  $\bullet$ $120 = 2^3 \cdot 3 \cdot 5$ (for $2 I$),

\vspace{1mm}
  \noindent
  and none of these orders are of the form $p^2$ or $2p$, in accord with
  Prop. \ref{ExistenceOfSmoothSphericaSpaceForms}.

  \noindent {\bf (ii)}  This is in contrast to the plain dihedral groups
  $D_{r + 2}$ whose order $2(r + 2)$ can violate
  the condition in Prop. \ref{ExistenceOfSmoothSphericaSpaceForms},
  showing that not all finite subgroups of $\SOThree$ act freely
  (and continuously) on any sphere.

  However, every subgroup $G \,\subset\,\SOThree$ is
  (double-)covered by a finite subgroup
  $\widehat{G} \,\in\,\SUTwo$:
  \vspace{-2mm}
  $$
    \begin{tikzcd}[column sep=large]
      \widehat{G}
      \ar[r, hook]
      \ar[d, ->>]
      \ar[dr, phantom, "\mbox{\tiny\rm(pb)}"]
      &
      \SUTwo
      \ar[d, ->>]
      \\
      G
      \ar[r, hook]
      &
      \SOThree
      \,.
    \end{tikzcd}
  $$

  \vspace{-2mm}
  \noindent {\bf (iii)}
  In conclusion, the finite subgroups of $\SUTwo$ have resolvable singularities
  in the sense of Ntn. \ref{ResolvableOrbiSingularities}, and the finite subgroups
  of $\SOThree$ have covers by groups that have resolvable singularities
  \eqref{FiniteGroupsWithResolvableSingularities}:
  \vspace{-2mm}
  \begin{align*}
    & \mathrm{FinSub}\left(
      \SUTwo
    \right)
    \,\subset\,
    \Groups(\FiniteSets)_\Resolvable
    ,\,
    \\
    &    \mathrm{FinSub}\left(
      \SOThree
    \right)
    \,\subset\,
    \Groups(\FiniteSets)_\CoverResolvable
    \,.
 \end{align*}

\vspace{-1mm}
  \noindent
  These ADE equivariance groups
  already subsume key examples of general and recent interest,
  e.g. all those considered in:
  \cite{GonzalezVerdier83}\cite{HillHopkinsRavenel09}\cite[\S 8]{BrunerGreenlees10}\cite{SS19TadpoleCancellation}\cite{KrizLu20}\cite{Yu21}.
\end{example}

\begin{remark}[Further cover-resolvable singularities]
  A necessary condition for a finite to group have cover-resolvable singularities
  in the sense of Ntn. \ref{ResolvableOrbiSingularities} is
  (see \cite{Sawin21})
  that it be a {\it $\mathcal{P}'$-group} in the sense of \cite{Wall13},
  where the classification of groups with this
  necessary property is discussed.
\end{remark}

\subsection{Singular cohesion}
\label{GeneralSingularCohesion}

Here we review and develop basics of
{\it cohesive global homotopy theory} or {\it singular cohesion}
from \cite[\S 3.2]{SS20OrbifoldCohomology},
combining the cohesive homotopy theory from \cite{dcct} with the
cohesive perspective on global homotopy theory from
\cite{Rezk14}.

\medskip

\noindent
{\bf The orbit category.}
The following Ntn. \ref{Singularities}
is known under a variety of different names and notations
in the literature on global homotopy theory
(see \cite[Rem. 3.47]{SS20OrbifoldCohomology}).
Since none of the symbols in use seem particularly illuminating and
no common convention has been established,
we follow the notation in \cite{SS20OrbifoldCohomology},
which has the advantage of
transparently expressing the conceptual role of this definition in
cohesive global homotopy theory (Def. \ref{GEquivariantAndGloballyEquivariantHomotopyTheories}
below):
\begin{notation}[2-Category of orbi-singularities {\cite[Def. 3.64]{SS20OrbifoldCohomology}}]
  \label{Singularities}
  $\,$

  \vspace{0mm}
  \noindent
  {\bf (i)}  We write
  \vspace{-2mm}
  \begin{equation}
    \label{CategoryOfSingularities}
    \begin{tikzcd}[row sep=-5pt, column sep=small]
      \Singularities
      \ar[r,phantom,"\coloneqq"]
      &
      \Groupoids_{1, \geq 1}^{\mathrm{fin}}
      \ar[rr, hook]
      &&
      \InfinityGroupoids
      \\
      \scalebox{0.8}{$
        \orbisingularG
      $}
      &&\longmapsto&
     \scalebox{0.7}{$
       B G
     $}
    \end{tikzcd}
  \end{equation}

  \vspace{-2mm}
  \noindent for the
  $\infty$-site (Ntn. \ref{InfinitySite})
  of finite connected groupoids,
  namely the full sub-category
  (via Exp. \ref{InfinitySitesFromFullSubInfinityCategoriesOfPresentableInfinityCategories})
  of $\InfinityGroupoids$
  (Ntn. \ref{SimplicialSetsAndInfinityGroupoids})
  on the delooping groupoids $B G$ \eqref{DeloopingGroupoid}
  of finite groups $G \in \Groups(\FiniteSets)$.

  \noindent {\bf (ii)}
  When regarded as objects in $\Singularities$ we denote these
  by the symbol ``$\orbisingularG$''.

  \noindent {\bf (iii)}
  This means that for $\orbisingularK, \orbisingularG \,\in\, \Singularities$
  corresponding to $K, G \,\in\, \Groups(\FiniteSets)$ we have
  \vspace{-2mm}
  \begin{equation}
    \label{HomGroupoidOfSingularitiesInTermsOfGroupHomomorphisms}
    \Singularities
    \big(
      \orbisingularK
      ,\,
      \orbisingularG
    \big)
    \;\;
    \simeq
    \;\;
    \HomotopyQuotient
    {
      \Groups(K,\,G)
    }
    {
      G
    }
    \,,
  \end{equation}

 \vspace{-2mm}
 \noindent
  where on the right we have the conjugation-action groupoid
  of $G$ acting by conjugation on group homomorphisms $K \to G$.
\end{notation}

\begin{notation}[Category of $G$-orbits {\cite[\S I.3]{Bredon67}}]
  \label{GOrbitCategory}
  For $G \in \Groups(\Sets)$, the {\it category of $G$-orbits}
  or {\it $G$-orbit category} is the full subcategory of
  $G$-actions (on sets) on the connected transitive $G$-sets
  (the types of $G$-orbits),
  hence on the quotient sets $G/H$ for any subgroup $H \subset G$
  (the stabilizer subgroup of any one point in the orbit):
  \vspace{-3mm}
  $$
    \begin{tikzcd}
      G/H
      \;\in\;
      \OrbitCategory{G}
      \ar[r, hook]
      &
      \Actions{G}(\Sets)
      \,.
    \end{tikzcd}
  $$
\end{notation}
(We consider the generalization of the orbit category, from finite groups to topological groups, below in Def. \ref{ProperTopologicalOrbitInfinityCategory}.)

\begin{lemma}[0-Truncated slice of singularities]
  \label{ZeroTruncatedSliceOfSingularities}
  For $G \,\in\, \Groups$,
  the full sub-(2,1)-category of the
  slice of $\Singularities$ (Def. \ref{Singularities})
  over $\orbisingularG$ \eqref{CategoryOfSingularities}
  on the 0-truncated objects (those whose homotopy fiber is a set)
  is

  \noindent {\bf (a)} reflective and

  \noindent {\bf (b)} given by the subgroup inclusions
  $H \xhookrightarrow{ i_H } G$:
  \vspace{-4mm}
  \begin{equation}
    \label{0TruncatedObjectsInSliceOverAnOrbiSingularity}
    \begin{tikzcd}
      \Singularities_{/\scalebox{.7}{$\orbisingularG$}}
      \ar[
        rr,
        shift left=6pt,
        "{ \tau_0 }"{above}
      ]
      &&
      \big(
        \Singularities_{/\scalebox{.7}{$\orbisingularG$}}
      \big)_{\mathrm{\leq 0}}
      \ar[
        ll,
        hook',
        shift left=6pt
      ]
      \ar[
        ll,
        phantom,
        "{ \scalebox{.7}{$\bot$} }"
      ]
    \end{tikzcd}
    \;\simeq\;
    \Biggg\{\!\!
    \begin{tikzcd}[column sep=8pt, row sep=small]
      \orbisingularH
      \ar[
        rr,
        dashed,
        "{} "{name=s, below}
      ]
      \ar[
        dr,
        "{ \orbisingularAny{\,\,i_H} }"{left, yshift=-4pt},
        "{  }"{name=t, above}
      ]
      &&
      \orbisingularK
      \ar[
        dl,
        "{ \orbisingularAny{\,\,i_K} }"{right, yshift=-4pt}
      ]
      \ar[
        from=s,
        to=t,
        Rightarrow,
        "{ \sim }"{sloped, below, pos=.3}
      ]
      \\
      &
      \orbisingularG
    \end{tikzcd}
    \!\!\Biggg\}
  \end{equation}

  \vspace{-2mm}
  \noindent
  for the full sub-(2,1)-category of the
  slice of $\Singularities$ (Def. \ref{Singularities})
  over $\orbisingularG$ \eqref{CategoryOfSingularities}
  on the 0-truncated objects,
  hence on the morphisms
  $\orbisingularH \xrightarrow{\scalebox{.6}{$\orbisingulari$}} \orbisingularG$
  corresponding to subgroup inclusions $H \xhookrightarrow{i} G$.
\end{lemma}

\begin{lemma}[$G$-orbits are the 0-truncated objects in slice over $G$-orbi-singularity]
  \label{GOrbitsAre0TruncatedObjectsOverGOrbiSingularity}
  For $G \,\in\, \Groups$
  the singularities in the 0-truncated slice
  over $\orbisingularG$ \eqref{0TruncatedObjectsInSliceOverAnOrbiSingularity}
  are equivalently the $G$-orbits (Ntn. \ref{GOrbitCategory})
  in that we have an equivalence in $\InfinityCategories$
  \vspace{-2mm}
  $$
    \begin{tikzcd}[row sep=-4pt, column sep=small]
      \big(
        (
          \Groupoids^{\mathrm{fin}}_{1, \geq 1}
        )_{/\scalebox{.7}{$B G$}}
      \big)_{\leq 0}
      \ar[r, phantom, "\mbox{$\simeq$}"]
      &
      \big(
        \Singularities_{/\scalebox{.7}{$\orbisingularG$}}
      \big)_{\leq 0}
      \ar[rr,"\sim"]
      &&
      \OrbitCategory{G}
      \\
      \scalebox{0.8}{$
        (BH \xrightarrow{B i_H} BG)
      $}
      \ar[r,phantom,"\mapsto"]
      &
      \scalebox{0.8}{$
        \big(
          \orbisingularH
          \xrightarrow{
            \scalebox{.7}{$\orbisingulari$}
          }
          \orbisingularG
        \big)
      $}
      &\longmapsto&
      \scalebox{0.8}{$
        G/H
      $}
      \,,
    \end{tikzcd}
  $$
  where $i_H \,:\, H \xhookrightarrow{\;} G$ are subgroup inclusions.
\end{lemma}
\begin{proof}
  By Prop. \ref{GSetsInTheHomotopyTheoryOverBG},
  we have an equivalence
  \vspace{-4mm}
  $$
    \begin{tikzcd}
      \left(
        (
          \Groupoids_{\infty}
        )_{/ B G}
      \right)_{\leq 0}
      \ar[rr, "\mathrm{fib}"{above}, "\sim"{below}]
      &&
      \Actions{G}(\Sets)
      \,.
    \end{tikzcd}
  $$

  \vspace{-3mm}
  \noindent
  This implies the statement by the observation that the homotopy fiber
  of $B i_H \,\colon\, B H \xrightarrow{\;} B G$ is $G/H$.
\end{proof}
In summary:
\begin{proposition}[$G$-Orbit category is reflective subcategory of slice over $G$-singularity]
  \label{GOrbitCategoryIsReflectiveSubcategoryOfSliceOverGSingularity}
  For $G \,\in\, \Groups(\FiniteSets)$, the $G$-orbit category
  (Ntn. \ref{GOrbitCategory})
  is a reflective subcategory of the slice of $\Singularities$ (Ntn. \ref{Singularities})
  over $\scalebox{.8}{$\orbisingularG$}$:
  \vspace{-2mm}
  \begin{equation}
    \label{FinteProductPreservingReflectionOfSingularitiesOntoOrbitCategory}
    \begin{tikzcd}[row sep=-2pt]
      \Singularities_{/\scalebox{.7}{$\orbisingularG$}}
      \ar[
        rr,
        shift left=5pt,
        "{ \tau_0 }"{above},
        "\mathclap{\times}"{description, pos=0}
      ]
      &&
      \OrbitCategory{G}
      \ar[
        ll,
        hook',
        shift left=5pt,
        "{  }"{below}
      ]
      \ar[
        ll,
        phantom,
        "{ \scalebox{.7}{$\bot$} }"
      ]
      \\
      \scalebox{0.8}{$
      \begin{array}{c}
        \orbisingularK
        \\
        \big\downarrow{}^{
          \mathrlap{ \scalebox{0.9}{$\orbisingularAny{\;\phi_{\phantom{a}}}$} }
        }
        \\
        \orbisingularG
      \end{array}
      $}
      &\longmapsto&
    \scalebox{0.8}{$
      G/K
      \mathrlap{
        \;=\;
        G/\mathrm{im}(\phi)
        \,.
      }
    $}
    \end{tikzcd}
  \end{equation}
\end{proposition}
\begin{proof}
  This follows by Lem.  \ref{ZeroTruncatedSliceOfSingularities}
  combined with
  Lem. \ref{GOrbitsAre0TruncatedObjectsOverGOrbiSingularity}.
\end{proof}

\medskip

\noindent
{\bf Equivariant homotopy theories.}

\begin{definition}[$G$-Equivariant and globally equivariant homotopy theories]
  \label{GEquivariantAndGloballyEquivariantHomotopyTheories}
  For $\ModalTopos{\smooth}$ an $\infty$-topos and $G \,\in\, \Groups(\Sets)$,

  \noindent {\bf (i)}
  we write
  \vspace{-2mm}
  \begin{equation}
    \label{EquivariantHomotopyToposes}
    \GloballyEquivariant{(\ModalTopos{\smooth})}
    \;\coloneqq\;
    \InfinitySheaves
    \left(
      \Singularities
      ,\,
      \ModalTopos{\smooth}
    \right)
    \;\;\;\;\;\;\;\;
    \mbox{and}
    \;\;\;\;\;\;\;\;
    \GEquivariant{(\ModalTopos{\smooth})}
    \;\coloneqq\;
    \InfinitySheaves
    \left(
      \OrbitCategory{G}
      ,\,
      \ModalTopos{\smooth}
    \right)
  \end{equation}

  \vspace{-1mm}
  \noindent
  for the $\infty$-toposes of $\infty$-(pre-)sheaves
  on
  the 2-site of orbi-singularities (Def. \ref{Singularities})
  and the site of $G$-orbits (Ntn. \ref{GOrbitCategory}),
  respectively.

  \noindent {\bf (ii)}
  Since their  $\infty$-sites
  have trivial Grothendieck topology and a terminal
  object, these $\infty$-toposes \eqref{EquivariantHomotopyToposes}
  are cohesive over $\ModalTopos{\smooth}$
  (by Ex. \ref{ExamplesOfDiscreteCohesion}),
  witnessed by adjoint quadruples which we denote like this:
  \vspace{-1mm}
  \begin{equation}
    \label{AdjointQuadrupleOnAGlobalHomotopyTheory}
    \begin{tikzcd}[column sep=large]
      \GloballyEquivariant{\ModalTopos{\smooth}}
      \ar[
        rr,
        "{ \Smooth }"{description}
      ]
      \ar[
        rr,
        shift left=36pt,
        "{ \Conical }"{description},
        "\mathclap{\times}"{description, pos=0}
      ]
      &&
      \ModalTopos{\smooth}
      \ar[
        ll,
        hook',
        shift right=18pt,
        "{ \Space }"{description}
      ]
      \ar[
        ll,
        hook',
        shift right=-18pt,
        "{ \Orbisingular }"{description}
      ]
      \ar[
        ll,
        phantom,
        shift right=9pt,
        "{\scalebox{.6}{$\bot$}}"
      ]
      \ar[
        ll,
        phantom,
        shift right=27pt,
        "{\scalebox{.6}{$\bot$}}"
      ]
      \ar[
        ll,
        phantom,
        shift right=-9pt,
        "{\scalebox{.6}{$\bot$}}"
      ]
      \ar[
        r,
        phantom,
        shift left=36pt,
        "{
          \mbox{
            \tiny
            \color{greenii}
            \bf
            conical
          }
        }"
      ]
      \ar[
        r,
        phantom,
        shift left=18pt,
        "{
          \mbox{
            \tiny
            \color{greenii}
            \bf
            spatial
          }
        }"
      ]
      \ar[
        r,
        phantom,
        shift left=0pt,
        "{
          \mbox{
            \tiny
            \color{greenii}
            \bf
            smooth
          }
        }"
      ]
      \ar[
        r,
        phantom,
        shift right=18pt,
        "{
          \mbox{
            \tiny
            \color{greenii}
            \bf
            singular
          }
        }"
      ]
      &
      {}
    \end{tikzcd}
  \end{equation}

\vspace{-1mm}
\noindent {\bf (iii)}
We denote the resulting modalities by:
\vspace{-1mm}
\begin{equation}
  \label{SingularModalities}
  \def\arraystretch{1}
  \begin{array}{cccc}
    \raisebox{3pt}{
      \tiny
      \color{darkblue}
      \bf
      \begin{tabular}{c}
        purely conical
        \\
        aspect
      \end{tabular}
    }
    &
    \conical
    &
    \coloneqq
    &
    \Space
      \circ
    \Conical
    \\
    &
    \scalebox{.7}{$\bot$}
    \\
    \raisebox{3pt}{
      \tiny
      \color{darkblue}
      \bf
      \begin{tabular}{c}
        purely smooth
        \\
        aspect
      \end{tabular}
    }
    &
    \smooth
    &
    \coloneqq
    &
    \Space
      \circ
    \Smooth
    \\
    &
    \scalebox{.7}{$\bot$}
    \\
    \raisebox{3pt}{
      \tiny
      \color{darkblue}
      \bf
      \begin{tabular}{c}
        purely orbi-singular
        \\
        aspect
      \end{tabular}
    }
    &
    \orbisingular
    &
    \coloneqq
    &
    \Orbisingular
      \circ
    \Smooth
    \,.
  \end{array}
\end{equation}

\end{definition}
\begin{example}[Classical equivariant homotopy theory]
  \label{ClassicalEquivariantHomotopyTheory}
  In the base case $\ModalTopos{\smooth} \,=\, \InfinityGroupoids$,
  the $\infty$-toposes from Def. \ref{GEquivariantAndGloballyEquivariantHomotopyTheories}
  are:

  \noindent
  (i) the classical equivariant homotopy theory
  $$
    \GEquivariant{\InfinityGroupoids}
    \;\simeq\;
    \InfinityPresheaves
    \big(
      \Orbits(G)
    \big)
  $$ from
  Elmendorf's theorem in its homotopically enhanced version
  \eqref{ElmendorfTheorem} due to \cite{DwyerKan84}
  (see also \cite{CordierPorter96}, review in \cite[Thm. 1.3.8]{Blu17}),

  \noindent
  (ii) the global equivariant homotopy theory $\GloballyEquivariant{\InfinityGroupoids}$
  in its unstable form first highlighted in \cite{Rezk14}.
\end{example}

More general instances of Def. \ref{GEquivariantAndGloballyEquivariantHomotopyTheories} enhance this classical equivariant homotopy theory by cohesive geometric structure:

\medskip
\noindent
{\bf Singular-cohesive $\infty$-toposes.}

\begin{notation}[Smooth charts]
  \label{SmoothCharts}
  Given a cohesive $\infty$-topos $\ModalTopos{\smooth}$,
  we say that a {\it category of smooth charts} is
  an $\infty$-category $\Charts$ of charts for $\ModalTopos{\smooth}$
  according to Ntn. \ref{CohesiveCharts}, subject to the
  additional condition that all $U \,\in\, \Charts$ are 0-truncated,
  hence that $\Charts$ is in fact a 1-category.
\end{notation}

\begin{definition}[Singular-cohesive $\infty$-toposes]
\label{SingularCohesiveInfinityTopos}
We say that an $\infty$-topos $\Topos$ is {\it singular-cohesive} over
$\InfinityGroupoids$ if it is equivalently the
global homotopy theory (Def. \ref{GEquivariantAndGloballyEquivariantHomotopyTheories})
of a cohesive $\infty$-topos $\ModalTopos{\smooth}$ (Def. \ref{CohesiveInfinityTopos})
which admits a site $\Charts$ of smooth charts (Def. \ref{SmoothCharts}):
\vspace{-2mm}
\begin{equation}
  \label{SingularCohesiveToposAsSheavesOnSingulartiesTimesCharts}
  \Topos
  \;\simeq\;
  \InfinitySheaves
  \left(
    \Singularities
    ,\,
    \Topos_{\scalebox{.6}{$\smooth$}}
  \right)
  \;\simeq\;
  \InfinitySheaves
  \left(
    \Singularities
      \times
    \Charts
    ,\,
    \InfinityGroupoids
  \right)
  \,.
\end{equation}

This implies\footnote{
  It would be desireable to state a converse implication, hence
  to {\it characterize} singular-cohesive $\infty$-toposes by
  specification of the system of adjoint functors they carry,
  instead of resorting to a presentation over a site.
  While desireable, this is not necessary for the present purpose,
  and we leave this question for elsewhere.
}
that it carries a system of adjoint $\infty$-functors of the following form:
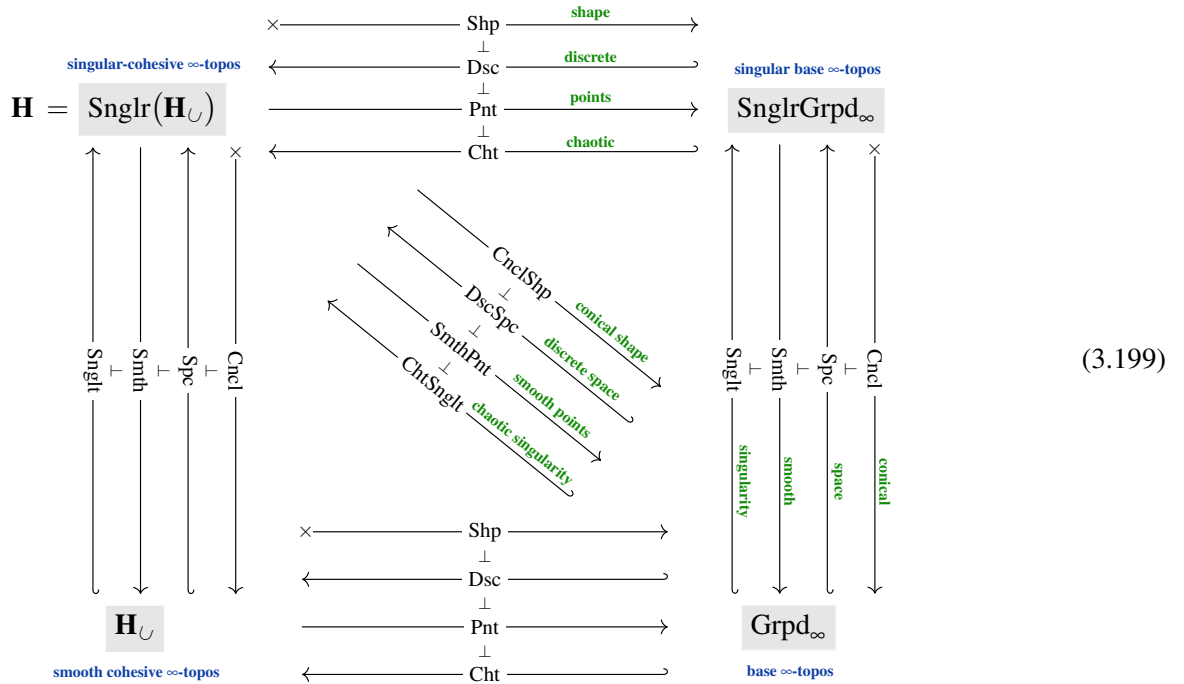
\vspace{-1mm}
\begin{equation}
  \label{SingularCohesiveAdjoints}
  \hspace{-5mm}
    \begin{tikzcd}[column sep=large]
      \mathllap{
      \Topos
        \;=\;\,
      }
      \overset{
        \mathclap{
        \raisebox{3pt}{
          \tiny
          \color{darkblue}
          \bf
          \begin{tabular}{c}
           singular-cohesive
           $\infty$-topos
          \end{tabular}
        }
        }
      }
      {
        \categorybox{
          \GloballyEquivariant{\big(\ModalTopos{\smooth}\big)}
        }
      }
      \quad
      \ar[
        rr,
        phantom,
        shift left=8pt,
        "{\scalebox{.6}{$\bot$}}"
      ]
      \ar[
        rr,
        phantom,
        shift left=24pt,
        "{\scalebox{.6}{$\bot$}}"
      ]
      \ar[
        rr,
        phantom,
        shift left=-8pt,
        "{\scalebox{.6}{$\bot$}}"
      ]
      \ar[
        rr,
        "{\Points}"{description},
        "{ \mbox{\tiny\color{greenii}\bf points} }"{above, yshift=-1pt, pos=.75}
      ]
      \ar[
        rr,
        shift left=32pt,
        "{\Shape}"{description},
        "\mathclap{\times}"{description, pos=0.01},
        "{ \mbox{\tiny\color{greenii}\bf shape} }"{above, yshift=-1pt, pos=.75}
      ]
      \ar[
        dd,
        phantom,
        shift left=9pt-10pt,
        "{\scalebox{.6}{\rotatebox{-90}{$\bot$}}}"
      ]
      \ar[
        dd,
        phantom,
        shift left=27pt-10pt,
        "{\scalebox{.6}{\rotatebox{-90}{$\bot$}}}"
      ]
      \ar[
        dd,
        phantom,
        shift left=-9pt-10pt,
        "{\scalebox{.6}{\rotatebox{-90}{$\bot$}}}"
      ]
      \ar[
        dd,
        shift left=-10pt,
        "{\scalebox{1}{$\Smooth$}}"{description, sloped}
      ]
      \ar[
        dd,
        shift left=36pt-10pt,
        "{\scalebox{1}{$\Conical$}}"{description, sloped},
        "\mathclap{\times}"{description, pos=0.01, sloped}
      ]
      \ar[
        ddrr,
        shift left=36pt,
        shorten <= 70pt,
        shorten >= 80pt,
        "{ \Conical\Shape }"{description, sloped, pos=.45},
        "{ \mbox{\tiny\color{greenii}\bf conical shape} }"{above, yshift=-1pt, sloped, pos=.6}
      ]
      \ar[
        ddrr,
        shorten <= 70pt,
        shorten >= 80pt,
        "{ \Smooth\Points }"{description, sloped, pos=.45},
        "{ \mbox{\tiny\color{greenii}\bf smooth points} }"{above, yshift=-1pt, sloped, pos=.6}
      ]
      &[60pt]
      &
      \overset{
        \raisebox{3pt}{
          \tiny
          \color{darkblue}
          \bf
          \begin{tabular}{c}
           singular base
           $\infty$-topos
          \end{tabular}
        }
      }
      {
        \categorybox{
          \GloballyEquivariant{ \InfinityGroupoids }
        }
      }
      \ar[
        ll,
        hook',
        shift right=16pt,
        "{ \Discrete }"{description},
        "{ \mbox{\tiny\color{greenii}\bf discrete} }"{above, yshift=-1pt, pos=.25}
      ]
      \ar[
        ll,
        hook',
        shift right=-16pt,
        "{ \Chaotic }"{description},
        "{ \mbox{\tiny\color{greenii}\bf chaotic} }"{above, yshift=-1pt, pos=.25}
      ]
      \ar[
        dd,
        phantom,
        shift left=9pt-10pt,
        "{\scalebox{.6}{\rotatebox{-90}{$\bot$}}}"
      ]
      \ar[
        dd,
        phantom,
        shift left=27pt-10pt,
        "{\scalebox{.6}{\rotatebox{-90}{$\bot$}}}"
      ]
      \ar[
        dd,
        phantom,
        shift left=-9pt-10pt,
        "{\scalebox{.6}{\rotatebox{-90}{$\bot$}}}"
      ]
      \ar[
        dd,
        shift left=-10pt,
        "{\scalebox{1}{$\Smooth$}}"{description, sloped},
        "{ \mbox{\tiny\color{greenii}\bf smooth} }"{above, yshift=-1pt, sloped, pos=.75}
      ]
      \ar[
        dd,
        shift left=36pt-10pt,
        "{\scalebox{1}{$\Conical$}}"{description, sloped},
        "\mathclap{\times}"{description, pos=0.01, sloped},
        "{ \mbox{\tiny\color{greenii}\bf conical} }"{above, yshift=-1pt, sloped, pos=.75}
      ]
      \\[130pt]
      \\
      \phantom{
        \InfinitySheaves
        \big)
      }
      \underset{
        \mathclap{
        \raisebox{-3pt}{
          \tiny
          \color{darkblue}
          \bf
          smooth cohesive $\infty$-topos
        }
        }
      }{
        \categorybox{
          \ModalTopos{\smooth}
        }
      }
      \phantom{
      \big(
        \Singularities
        ,\,
      }
      \quad
      \ar[
        rr,
        phantom,
        shift left=9pt,
        "{\scalebox{.6}{$\bot$}}"
      ]
      \ar[
        rr,
        phantom,
        shift left=27pt,
        "{\scalebox{.6}{$\bot$}}"
      ]
      \ar[
        rr,
        phantom,
        shift left=-9pt,
        "{\scalebox{.6}{$\bot$}}"
      ]
      \ar[
        rr,
        "{\Points}"{description}
      ]
      \ar[
        rr,
        shift left=36pt,
        "{\Shape}"{description},
        "\mathclap{\times}"{description, pos=0.01}
      ]
      \ar[
        uu,
        hook',
        shift right=18pt-10pt,
        "{\scalebox{1}{$ \Space$} }"{description, sloped, rotate=180}
      ]
      \ar[
        uu,
        hook',
        shift right=-18pt-10pt,
        "{\scalebox{1}{$ \Orbisingular$} }"{description, sloped, rotate=180}
      ]
      &&
      \phantom{
        \InfinitySheaves
       \big)
      }
      \underset{
        \mathclap{
        \raisebox{-3pt}{
          \tiny
          \color{darkblue}
          \bf
          base $\infty$-topos
        }
        }
      }{
        \categorybox{
          \InfinityGroupoids
        }
      }
      \phantom{
      \big(
        \Singularities
        ,\,
      }
      \ar[
        ll,
        hook',
        shift right=18pt,
        "{  \Discrete }"{description}
      ]
      \ar[
        ll,
        hook',
        shift right=-18pt,
        "{ \Chaotic }"{description}
      ]
      \ar[
        uu,
        hook',
        shift right=18pt-10pt,
        "{ \scalebox{1}{$\Space$} }"{description, sloped, rotate=180},
        "{ \mbox{\tiny\color{greenii}\bf space} }"{above, yshift=+1pt, sloped, rotate=180, pos=.25}
      ]
      \ar[
        uu,
        hook',
        shift right=-18pt-10pt,
        "{\scalebox{1}{$ \Orbisingular $}}"{description, sloped, rotate=180},
        "{ \mbox{\tiny\color{greenii}\bf singularity} }"{above, yshift=+1pt, sloped, rotate=180, pos=.25}
      ]
      \ar[
        uull,
        shift right=18pt,
        hook',
        shorten <= 80pt,
        shorten >= 70pt,
        "{ \Discrete\Space }"{description, sloped, pos=.55},
        "{\mbox{\tiny\color{greenii}\bf discrete space}}"{above, yshift=-1pt, sloped, pos=.4}
      ]
      \ar[
        uull,
        shift left=18pt,
        hook',
        shorten <= 80pt,
        shorten >= 70pt,
        "{ \Chaotic\Singularity }"{description, sloped, pos=.55},
        "{\mbox{\tiny\color{greenii}\bf chaotic singularity}}"{above, yshift=-1pt, sloped, pos=.4}
      ]
      \ar[
        uull,
        phantom,
        shift right=27pt,
        "\scalebox{.6}{$\bot$}"{sloped, pos=.55}
      ]
      \ar[
        uull,
        phantom,
        shift right=9pt,
        "\scalebox{.6}{$\bot$}"{sloped, pos=.55}
      ]
      \ar[
        uull,
        phantom,
        shift left=9pt,
        "\scalebox{.6}{$\bot$}"{sloped, pos=.55}
      ]
    \end{tikzcd}
\end{equation}

\end{definition}

We are going to focus on the special case relevant to the differential topology of orbifolds:
\begin{notation}[Singular-smooth $\infty$-groupoids {\cite[Ex. 3.56]{SS20OrbifoldCohomology}}]
  \label{SingularSmoothInfinityGroupoids}
  We write
  $$
    \SingularSmoothInfinityGroupoids
    \;\coloneqq\;
    \Sheaves\left(\Singularities, \, \SmoothInfinityGroupoids\right)
    \;\simeq\;
    \InfinitySheaves
    \left(
      \Singularities \times \CartesianSpaces
    \right)
  $$
  for the singular-cohesive $\infty$-topos
  (Def. \ref{SingularCohesiveInfinityTopos})
  over that of smooth $\infty$-groupoids (Ntn. \ref{SmoothInfinityGroupoids}).
\end{notation}

\begin{example}[Singular-cohesive aspects of orbi-singularities
{\cite[Lem. 3.61, Prop. 3.62]{SS20OrbifoldCohomology}}]
  \label{OrbiSingularityIsOrbiSingularizationOfHomotopyQuotient}
  For $\Topos$ a singular-cohesive $\infty$-topos (Def. \ref{SingularCohesiveInfinityTopos})
  and $G \,\in\, \Groups(\FiniteSets) \xhookrightarrow{\Groups(\Discrete)} \Groups(\Topos)$,
  the canonical $G$-orbisingularity (Ntn. \ref{Singularities})
  \vspace{-2mm}
  $$
    \orbisingularG
    \,\in\,
    \Singularities
    \xhookrightarrow{\;\;
      \overset{
        \mbox{\tiny \rm \eqref{InfinityYonedaEmbedding}}
      }{
        \YonedaEmbedding
      }
    \;\;}
    \Presheaves(\Singularities,\, \InfinityGroupoids)
    \,=\,
    \GloballyEquivariant\InfinityGroupoids
    \xhookrightarrow{ \;\;
      \overset{
        \mbox{\tiny \eqref{SingularCohesiveAdjoints} }
      }{
        \Discrete
      }
   \;\; }
    \Topos
  $$

  \vspace{0mm}
 \noindent has the following modal aspects \eqref{SingularModalities}:

  \vspace{-.2cm}
  \begin{enumerate}[{\bf (1)}]
    \setlength\itemsep{-3pt}
   \item
  its {\it purely conical aspect} is the point:
  \vspace{-1mm}
  \begin{equation}
    \label{ConicalAspectOfOrbiSingularityIsThePoint}
    \conical
    \,
    \orbisingularG
    \;\simeq\;
    \ast\;.
  \end{equation}

  \item
  Its {\it purely smooth aspect} is the delooping
  \eqref{DeloopingOfInfinityGroupAsColimit}:
  \vspace{-1mm}
  \begin{equation}
    \label{SmoothAspectOfOrbisingularity}
    \smooth
    \,
    \orbisingularG
    \;\simeq\;
    \mathbf{B}G\;.
  \end{equation}

  \item
  It is  the {\it orbi-singularization }
  \eqref{SingularModalities} of the delooping of $G$
  \eqref{DeloopingOfInfinityGroupAsColimit}:
    \vspace{-1mm}
  \begin{equation}
    \label{OrbiSingularityIsOrbiSingularizationOfDelooping}
    \orbisingularG
    \;\simeq\;
    \orbisingular \mathbf{B}G\;.
  \end{equation}

  \end{enumerate}
  \vspace{-.2cm}
\end{example}
\begin{proof}
  The first statement \eqref{ConicalAspectOfOrbiSingularityIsThePoint}
  follows via Lem. \ref{LeftKanExtensionOnRepresentablesIsOriginalFunctor},
  as $\Conical \,\simeq\, p_!$ is the left Kan extension of
  (see Ex. \ref{ExamplesOfDiscreteCohesion})
  \vspace{-2mm}
  $$
    \begin{tikzcd}
      \Singularities
      \ar[r, shift left=5pt, "p"]
      \ar[from=r, hook', shift left=5pt]
      \ar[r, phantom, "\scalebox{.7}{$\bot$}"]
      &
      \ast
      \,.
    \end{tikzcd}
  $$

    \vspace{-2mm}
\noindent
  This implies, for
  $\TopologicalPatch \times \orbisingularK \,\in\, \Charts \times \Singularities$,
  the following natural equivalences
  \vspace{-2mm}
$$
  \def\arraystretch{1.6}
  \begin{array}{lll}
  \PointsMaps{\big}
    { \TopologicalPatch \times \orbisingularK }
    { \smooth \orbisingularG }
  &
  \;\simeq\;
  \PointsMaps{\big}
    { \conical (\TopologicalPatch \times \orbisingularG) }
    { \orbisingularG }
  &
  \proofstep{ by \eqref{SingularModalities} }
  \\
  &
  \;\simeq\;
  \PointsMaps{\big}
    {
      \underset{
        \mathclap{
        \scalebox{.7}{$
        \def\arraystretch{.9}
        \begin{array}{l}
          { \simeq \, \conical \, \smooth \, U }
          \\
          { \simeq \, \TopologicalPatch }
        \end{array}
        $}
        }
      }{
        \underbrace{
          (\conical \TopologicalPatch)
        }
      }
        \times
      \underset
        {\simeq \, \orbisingularE}
        {
          \underbrace{
            (\conical \orbisingularG)
          }
        }
    }
    { \orbisingularG }
  &
  \proofstep{
    by \eqref{AdjointQuadrupleOnAGlobalHomotopyTheory}
    \& \eqref{ConicalAspectOfOrbiSingularityIsThePoint}
  }
  \\
  &
  \;\simeq\;
  \PointsMaps{\big}
    { \TopologicalPatch \times \orbisingularE }
    { \orbisingularG }
  \\
  &
  \;\simeq\;
  \PointsMaps{\big}
    { \TopologicalPatch \times \orbisingularE }
    { \flat \orbisingularG }
  \\
  &
  \;\simeq\;
  \PointsMaps{\big}
    { \shape( \TopologicalPatch \times \orbisingularE) }
    { \orbisingularG }
  \\
  &
  \;\simeq\;
  \PointsMaps{\big}
    { \orbisingularE }
    { \orbisingularG }
  \\
  &
  \;\simeq\;
  \Singularities \big(
    { \orbisingularE },
    { \orbisingularG }
    \big)
  &
  \proofstep{ by \eqref{InfinityYonedaEmbedding}}
  \\
  &
  \;\simeq\;
  \Groupoids( B 1 ,\, BG )
  \\
  & \;\simeq\;
  B G
  \\
  & \;\simeq\;
  \PointsMaps{\big}
    { \TopologicalPatch \times \orbisingularK }
    { \mathbf{B}G }\;.
 \end{array}
$$

\vspace{-2mm}
\noindent
This implies the statement \eqref{DeloopingOfInfinityGroupAsColimit}
by the $\infty$-Yoneda lemma (Prop. \ref{InfinityYonedaLemma}).

Similarly, for the last statement
\eqref{OrbiSingularityIsOrbiSingularizationOfDelooping}:
  \vspace{-2mm}
$$
\hspace{-2cm}
  \def\arraystretch{1.6}
  \begin{array}{lll}
    \PointsMaps{\big}
      { \TopologicalPatch \times \orbisingularK }
      { \orbisingular \mathbf{B}G }
    &
    \;\simeq\;
    \PointsMaps{\big}
      { \TopologicalPatch \times ( \orbisingular \mathbf{B}K ) }
      { \orbisingular \mathbf{B}G }
    \\
    & \;\simeq\;
    \PointsMaps{\big}
      { \orbisingular ( \TopologicalPatch \times \mathbf{B}K ) }
      { \orbisingular \mathbf{B}G }
    \\
    & \;\simeq\;
    \PointsMaps{}
      { \TopologicalPatch \times \mathbf{B}K }
      { \mathbf{B}G }
    \\
    & \;\simeq\;
    \PointsMaps{\big}
      { \Discrete(B K) }
      { \Discrete(B G) }
    \\
    & \;\simeq\;
    \InfinityGroupoids
    ( B K ,\, B G )
    \\
    & \;\simeq\;
    \Singularities
    \big(
      \orbisingularK
      ,\,
      \orbisingularG
    \big)
    \\
    & \;\simeq\;
    \PointsMaps{\big}
      { \orbisingularK }
      { \orbisingularG }
    \\
    & \;\simeq\;
    \PointsMaps{\big}
      { \TopologicalPatch \times \orbisingularK }
      { \orbisingularG }
    \,.
  \end{array}
$$

\vspace{-7mm}
\end{proof}

\begin{example}[Cohesive loci of orbi-singularities]
  \label{CohesiveLociOfOrbiSingularities}
  Given a singular-cohesive $\infty$-topos $\Topos$,
  and  $G \,\in\, \Groups(\FiniteSets)$,
then
\begin{equation}
  \label{ValueOnOrbisingularityIsSmoothAspectOfMappingStack}
  \mathcal{X}(\orbisingularG)
  \;\;
  \overset{
    \mbox{
      \tiny\rm
      \eqref{InfinityYonedaLemma}
    }
  }{
    \simeq
  }
  \;\;
  \PointsMaps{}
    { \orbisingularG }
    { \mathcal{X} }
  \;\;
  \overset{
    \mbox{
      \tiny
      \rm
      \eqref{HomIsGlobalPointsOfMappingStacks}
    }
  }{
    \simeq
  }
  \;\;
  \Smooth
  \left(
    \Maps{}
      { \orbisingularG }
      { \mathcal{X} }
  \right).
\end{equation}
is the cohesive mapping space of $G$-orbi-singularities into $\mathcal{X}$.

In particular, if $\mathcal{X} \,\in\, \Topos$ is smooth
$\smooth \mathcal{X} \,\simeq\, \mathcal{X}$,
then it has no orbi-singularities and hence these mapping spaces
are equivalently those out of the point:
\begin{equation}
  \label{CohesiveLocusOfOrbiSingularitiesInASmoothObject}
  \smooth(\mathcal{X}) \simeq \mathcal{X}
  \;\;\;\;\;
  \Rightarrow
  \;\;\;\;\;
  \smooth \Maps{}{ \orbisingularG }{ \mathcal{X} }
  \;\simeq\;
  \mathcal{X}
  \,.
\end{equation}
More generally, given any other smooth object $U$,
the cohesive $U$-parameterized loci of orbi-singularities in a smooth object
form the smooth mapping space out of $U$:
\begin{equation}
  \label{ParametrizedCohesiveLocusOfOrbiSingularitiesInASmoothObject}
  \begin{array}{c}
    \smooth(\mathcal{X}) \simeq \mathcal{X}
    \\
    \smooth(U) \simeq U
  \end{array}
  \;\;\;\;\;
  \Rightarrow
  \;\;\;\;\;
  \smooth \Maps{\big}{ U \times \orbisingularG }{ \mathcal{X} }
  \;\simeq\;
  \smooth \Maps{}{ U }{ \mathcal{X} }
  \,.
\end{equation}
\end{example}
\begin{proof}
The equivalence
\eqref{ParametrizedCohesiveLocusOfOrbiSingularitiesInASmoothObject}
is, under the $\infty$-Yoneda lemma,
the composite of  the following sequence of
natural equivalences:
$$
  \def\arraystretch{1.6}
  \begin{array}{lll}
  \PointsMaps{\big}
    { U \times \orbisingularG }
    { \mathcal{X} }
  &
  \;\simeq\;
  \PointsMaps{\big}
    { U \times \orbisingularG }
    { \smooth \mathcal{X} }
  &
  \proofstep{ by assumption }
  \\
  &
  \;\simeq\;
  \PointsMaps{\big}
    { \conical (U \times \orbisingularG) }
    { \mathcal{X} }
  &
  \proofstep{ by \eqref{SingularModalities} }
  \\
  &
  \;\simeq\;
  \PointsMaps{\big}
    {
      \underset{
        \mathclap{
        \scalebox{.7}{$
        \def\arraystretch{.9}
        \begin{array}{l}
          { \simeq \, \conical \, \smooth \, U }
          \\
          { \simeq \, U }
        \end{array}
        $}
        }
      }{
        \underbrace{
          (\conical U)
        }
      }
        \times
      \underset
        {\simeq \, \ast}
        {
          \underbrace{
            (\conical \orbisingularG)
          }
        }
    }
    { \mathcal{X}
  }
  &
  \proofstep{
    by \eqref{AdjointQuadrupleOnAGlobalHomotopyTheory}
    \& \eqref{ConicalAspectOfOrbiSingularityIsThePoint}
  }
  \\
  &
  \;\simeq\;
  \PointsMaps{}
    { U }
    { \mathcal{X} }\;.
  \end{array}
$$
The statement \eqref{CohesiveLocusOfOrbiSingularitiesInASmoothObject}
is the special case of
\eqref{ParametrizedCohesiveLocusOfOrbiSingularitiesInASmoothObject}
for $U \,=\, \ast$.
\end{proof}

\begin{example}[Orbi-singularization does not preserve deloopings]
  \label{OrbiSingularizationDoesNotPreserveDeloopings}
  Beware that the rightmost adjoint modality $\orbisingular$
  \eqref{SingularModalities}
  has no reason to preserve $\infty$-colimits, hence
  no reason to preserve homotopy quotients
  \eqref{HomotopyQuotientAsHomotopyColimit}
  such as deloopings $\mathbf{B}(-)$ \eqref{DeloopingOfInfinityGroupAsColimit};
  and indeed generically it does not:

  \vspace{1mm}
  \noindent {\bf (i)}
  For $G \,\in\, \Groups(\FiniteSets)$ and
  $\Gamma \,\in\, \Groups(\Sets) \xhookrightarrow{ \Groups(\Discrete) } \Topos$,
  the value of the
  orbi-singularization of
  the delooping $\mathbf{B} \Gamma$ on the
  $G$-orbi-singularity $\orbisingularG$ is the $\Gamma$-conjugation action groupoid
  of group homomorphisms $G \to \Gamma$:
  \vspace{-3mm}
  \begin{equation}
    \label{ValueOfOrbiSingularizedDeloopingOnOrbiSingularity}
    (
      \orbisingular \mathbf{B}\Gamma
    )(\orbisingularG)
    \;
    \underset{
      \mbox{
        \tiny
        \eqref{DerivingValueOfOrbisingularizedDeloopingOnOrbiSingularity}
      }
    }{
      \simeq
    }
    \;
    \Groupoids(\mathbf{B}G,\, \mathbf{B}\Gamma)
    \;\simeq\;
    \ActionGroupoid{ \Groups(G,\,\Gamma) }{\Gamma}
    \;\;\;\;\;\;
    \in
    \;
    \Groupoids
    \xhookrightarrow{\;}
    \InfinityGroupoids
    \xhookrightarrow{\; \Discrete \;}
    \Topos
    \,,
  \end{equation}
  because:
  \begin{equation}
    \label{DerivingValueOfOrbisingularizedDeloopingOnOrbiSingularity}
    \def\arraystretch{1.4}
    \begin{array}{lll}
      (
        \orbisingular \mathbf{B}\Gamma
      )(\orbisingularG)
      &
      \;\simeq\;
      \PointsMaps{}
        { \orbisingularG }
        { \orbisingular \mathbf{B}\Gamma }
      &
      \proofstep{ by Prop. \ref{InfinityYonedaLemma} }
      \\
      &
      \;\simeq\;
      \PointsMaps{}
        { \smooth \, \orbisingularG }
        { \mathbf{B}\mathcal{G} }
      &
      \proofstep{ by \eqref{SingularModalities} }
      \\
      &
      \;\simeq\;
      \PointsMaps{}
        { \mathbf{B}G }
        { \mathbf{B}\mathcal{G} }
      &
      \proofstep{ by Ex. \ref{OrbiSingularityIsOrbiSingularizationOfHomotopyQuotient} }
      \\
      &
      \;\simeq\;
      \PointsMaps{}
        { \flat \mathbf{B}G }
        { \flat \mathbf{B}\Gamma }
      &
      \proofstep{ by \eqref{FlatPreservesDeloopings} }
      \\
      & \;\simeq\;
      \InfinityGroupoids( \mathbf{B}G ,\, \mathbf{B}\Gamma)
      &
      \proofstep{ by Def. \ref{CohesiveInfinityTopos}. }
    \end{array}
  \end{equation}
  So as soon as there is more than one conjugacy class of group homomorphisms
  $G \to \Gamma$, \eqref{ValueOfOrbiSingularizedDeloopingOnOrbiSingularity} says that
  $\orbisingular \mathbf{B}\Gamma$
  has non-trivial $\pi_0$,
  hence is not connected, hence cannot be a delooping,
  by Prop. \ref{LoopingAndDeloopingEquivalence}.

  \vspace{1mm}
  \noindent {\bf (ii)}
  This failure of the rightmost singular modality $\orbisingular$
  to preserve deloopings is ultimately
  (we discuss this in \cref{EquivariantModuliStacks} below)
  the source of the rich set of connected
  components of equivariant classifying spaces
  that we had seen in
  Prop. \ref{ConnectedComponentsOfHFixedEquivariantClassifyingSpaceWhenCrossedConjugationQuotientIsDiscrete},
  Rem. \ref{PseudoPrincipalNatureOfFixedLociInUniversalEquivariantPrincipalBundles}.

  \vspace{1mm}
  \noindent {\bf (iii)}
  This state of affairs is what makes twisted equivariant cohomology rich and subtle:  While plain twisted cohomology is about sections of higher fiber bundles, twisted equivariant cohomology is about sections of the image under $\orbisingular$ of such higher fiber bundles.
\end{example}

\begin{remark}[Singularity-wise evaluation of cohesive modalities]
\label{SingularitiesWiseEvaluationOfCohesiveModalities}
Def. \ref{SingularCohesiveInfinityTopos}
means that under the identification
$$
  \Topos
  \;=\;
  \GloballyEquivariant\ModalTopos{\smooth}
  \;=\;
  \Sheaves( \Singularities ,\, \ModalTopos{\smooth} )
$$
the cohesive modalities on $\Topos$ act
functorially over each
$\orbisingularG \,\in\, \Singularities$ as the
cohesive modalities on $\ModalTopos{\smooth}$,
in that for $\mathcal{X} \,\in\, \Topos$
regarded as
$\mathcal{X} \,\colon\, \orbisingularG \,\mapsto\, \mathcal{X}(\orbisingularG)$,
we have
the following natural equivalences:
\begin{align}
  \label{SingularityWiseApplicationOfCohesiveModalities}
&  \left(\shape \mathcal{X}\right)(\orbisingularG)
  \;\simeq\;
  \shape
  \left(
    \mathcal{X}(\orbisingularG)
  \right)
  \,,
  \\
&  \left(\flat \mathcal{X}\right)(\orbisingularG)
  \;\simeq\;
  \flat
  \left(
    \mathcal{X}(\orbisingularG)
  \right)
  \,,
  \nonumber
  \\
&  \left(\sharp \mathcal{X}\right)(\orbisingularG)
  \;\simeq\;
  \sharp
  \left(
    \mathcal{X}(\orbisingularG)
  \right)
  \,.
  \nonumber
\end{align}
\end{remark}

\vspace{1mm}
\begin{proposition}[Some singular modalities commute with some cohesive modalities]
  \label{SomeSingularModalitiesCommuteWithSomeCohesiveModalities}
  In a singular-cohesive $\infty$-topos $\Topos$ (Def. \ref{SingularCohesiveInfinityTopos})
  such that $\ModalTopos{\smooth}$ admits cohesive charts
  (Ntn.  \ref{CohesiveCharts})

  -- all cohesive modalities \eqref{TheModalitiesOnACohesiveInfinityTopos}
    commute with $\smooth$ \eqref{SingularModalities},

  -- all singular modalities \eqref{SingularModalities}
    commute with $\flat$ \eqref{TheModalitiesOnACohesiveInfinityTopos},

\noindent
in that there are natural equivalences of the following form:
\vspace{-1mm}
\begin{equation}
\label{CommutativityOfSomeCohesiveModalitiesWithSomeSingularModalities}
\mbox{
\begin{tabular}{c|c|c}
\cline{2-2}
&
$
  \phantom{\mathclap{\vert^{\vert^{\vert^{\vert^{\vert^{\vert}}}}}}}
  \arraycolsep=1pt
  \begin{array}{cccc}
    &
    \flat &\circ& \conical
    \\
    \simeq
    \;\,
    &
    \conical &\circ& \flat
  \end{array}
$
&
\\
\hline
\multicolumn{1}{|c|}{
$
  \phantom{\mathclap{\vert^{\vert^{\vert^{\vert^{\vert^{\vert}}}}}}}
  \arraycolsep=1pt
  \begin{array}{cccc}
    &
    \shape &\circ& \smooth
    \\
    \simeq
    \;\,
    &
    \smooth &\circ& \shape
  \end{array}
$
}
&
\multicolumn{1}{c}{
$
  \arraycolsep=1pt
  \begin{array}{cccc}
    &
    \flat &\circ& \smooth
    \\
    \simeq
    \;\,
    &
    \smooth &\circ& \flat
  \end{array}
$
}
&
\multicolumn{1}{|c|}{
$
  \arraycolsep=1pt
  \begin{array}{cccc}
    &
    \sharp &\circ& \smooth
    \\
    \simeq
    \;\,
    &
    \smooth &\circ& \sharp
  \end{array}
$
}
\\
\hline
&
$
  \phantom{\mathclap{\vert^{\vert^{\vert^{\vert^{\vert^{\vert}}}}}}}
  \arraycolsep=1pt
  \begin{array}{cccc}
    &
    \flat &\circ& \orbisingular
    \\
    \simeq
    \;\,
    &
    \orbisingular &\circ& \flat
  \end{array}
$
\\
\cline{2-2}
\end{tabular}
}
\end{equation}
\end{proposition}
\begin{proof}
  For $\Charts$ denoting any site (of charts) for $\ModalTopos{\smooth}$
  we have, by Def. \ref{GEquivariantAndGloballyEquivariantHomotopyTheories},
  equivalences of $\infty$-categories of the following form:
  \vspace{-3mm}
  \begin{equation}
    \label{DecomposingSingularCohesiveToposOverProductSiteForProofOfCommutingModalities}
    \begin{tikzcd}[row sep=6pt,column sep=-2pt]
    &&
    \Topos
    \ar[d, phantom, "\mbox{$\simeq$}"{sloped}]
    \\
    \Sheaves
    \big(
      \Singularities
      ,\,
      \ModalTopos{\smooth}
    \big)
    &\simeq&
    \Sheaves
    \big(
      \Singularities \times \Charts
      ,\,
      \InfinityGroupoids
    \big)
    &\simeq&
    \Sheaves
    \big(
      \Charts
      ,\,
      \GloballyEquivariant\InfinityGroupoids
    \big)
    \,.
    \end{tikzcd}
  \end{equation}
  Here we may assume without restriction,
  by \eqref{ShapePreservesTheTerminalObject},
  that $\Charts$
  contains the terminal object
  $$
    \ast_{\scalebox{.7}{$\smooth$}}
    \;\in\;
    \Charts
    \xhookrightarrow{\;\; \YonedaEmbedding \;\;}
    \ModalTopos{\smooth}
    \,.
  $$
  But this means that,
  under the equivalence on the right of
  \eqref{DecomposingSingularCohesiveToposOverProductSiteForProofOfCommutingModalities},
  $\flat$ is the operation of evaluating on
  $\ast_{\scalebox{.7}{$\smooth$}}$ and then extending as a constant sheaf,
  by Ex. \ref{FlatModalityOverSiteOfCohesiveCharts}, in that there
  are natural equivalences:
  \begin{equation}
    \label{FlatModalityByConstantExtensionOverChartsOfEvaluationOnPoint}
    \underset{
      \scalebox{.7}{$
        U \,\in\, \Charts
      $}
    }{\forall}
    \;\;\;
    (\flat \mathcal{X})(U) \,\simeq\, \mathcal{X}(\ast_{\scalebox{.7}{$\smooth$}})
    \;\;\;
    \in
    \;
    \GloballyEquivariant\InfinityGroupoids
    \,.
  \end{equation}
  The analogous statement holds for $\smooth$ under the equivalence
  on the left of
  \eqref{DecomposingSingularCohesiveToposOverProductSiteForProofOfCommutingModalities},
  by the definition \eqref{AdjointQuadrupleOnAGlobalHomotopyTheory}
  via Ex. \ref{ExamplesOfDiscreteCohesion}:
  \begin{equation}
    \label{SmoothModalityByConstantExtensionOverSingularitiesOfEvaluationOnTrivialSingularity}
    \underset{
      \scalebox{.7}{$
        \orbisingularG \,\in\, \Singularities
      $}
    }{\forall}
    \;\;\;
    (\smooth \mathcal{X})(\orbisingularG)
    \;\;
    \simeq
    \;
    \mathcal{X}(\orbisingularE)
    \;\;\;
    \in
    \;
    \ModalTopos{\smooth}
    \,.
  \end{equation}
  Moreover, under the equivalence on the left
  of \eqref{DecomposingSingularCohesiveToposOverProductSiteForProofOfCommutingModalities}
  the cohesive modalities act objectwise over singularities
  (by Rem. \ref{SingularitiesWiseEvaluationOfCohesiveModalities}),
  e.g. for $\shape$
  we have natural equivalences
  \begin{equation}
    \label{CohesiveModalitiesActingObjectwiseOverSingularities}
    \left(\shape \mathcal{X}\right)(\orbisingularG)
    \;\simeq\;
    \shape
    \left(
      \mathcal{X}(\orbisingularG)
    \right)
    \;\;\;
    \in
    \;
    \ModalTopos{\smooth}
    \,;
  \end{equation}
  while
  under the equivalence on the right
  the singular modalities act objectwise over sharts, e.g. for $\orbisingular$
  we have natural equivalences
  \begin{equation}
    \label{SingularModalitiesActingObjectwiseOverCharts}
    \left(\orbisingular \mathcal{X}\right)(U)
    \;\simeq\;
    \orbisingular
    \left(
      \mathcal{X}(U)
    \right)
    \;\;\;
    \in
    \;
    \GloballyEquivariant\InfinityGroupoids
    \,.
  \end{equation}
  Using this, the claim follows with the $\infty$-Yoneda lemma
  (Prop. \ref{InfinityYonedaLemma}). For example, for the cases of
  $(\shape, \smooth)$ and
  $(\flat, \orbisingular)$ we have the following sequences of
  natural equivalences in
  $\orbisingularG \,\in\, \Singularities$ and in
  $U \,\in\, \Charts$, respectively:
  $$
     \def\arraystretch{1.4}
    \begin{array}{lll}
      \left(
        \smooth
        \,
        \shape
        \,
        \mathcal{X}
      \right)(\orbisingularG)
      &
      \;\simeq\;
      (\shape \, \mathcal{X})(\orbisingularE)
      &
      \proofstep{ by \eqref{SmoothModalityByConstantExtensionOverSingularitiesOfEvaluationOnTrivialSingularity} }
      \\
      &
      \;\simeq\;
      \shape\,
      \left(
        \mathcal{X}(\orbisingularE)
      \right)
      &
      \proofstep{ by \eqref{CohesiveModalitiesActingObjectwiseOverSingularities} }
      \\
      & \;\simeq\;
      \shape\,
      \left(
        (\smooth \, \mathcal{X})(\orbisingularG)
      \right)
      &
      \proofstep{ by \eqref{SmoothModalityByConstantExtensionOverSingularitiesOfEvaluationOnTrivialSingularity} }
      \\
      & \;\simeq\;
      \left(
      \shape
      \,
      \smooth
      \,
      \mathcal{X}
      \right)
      (\orbisingularG)
      &
      \proofstep{ by \eqref{CohesiveModalitiesActingObjectwiseOverSingularities}, }
    \end{array}
    $$
    and analogously:
    $$
    \def\arraystretch{1.3}
    \begin{array}{lll}
      \left(
      \flat
      \,
      \orbisingular
      \,
      \mathcal{X}
      \right)(U)
      &
      \;\simeq\;
      (\orbisingular \, \mathcal{X})(\ast_{\scalebox{.7}{$\smooth$}})
      &
      \proofstep{ by \eqref{FlatModalityByConstantExtensionOverChartsOfEvaluationOnPoint} }
      \\
      & \;\simeq\;
      \orbisingular
      \left(
        \mathcal{X}(\ast_{\scalebox{.7}{$\smooth$}})
      \right)
      &
      \proofstep{ by \eqref{SingularModalitiesActingObjectwiseOverCharts} }
      \\
      & \;\simeq\;
      \orbisingular
      \left(
        (\flat \mathcal{X})(U)
      \right)
      &
      \proofstep{ by \eqref{FlatModalityByConstantExtensionOverChartsOfEvaluationOnPoint} }
      \\
      & \;\simeq\;
     \left(\orbisingular \, \flat \, \mathcal{X}\right)(U)
      &
      \proofstep{ by \eqref{SingularModalitiesActingObjectwiseOverCharts}. }
    \end{array}
  $$
  The proof of the remaining cases is obtained by substitution.
\end{proof}

\begin{remark}[Noncommuting composite modalities]
  The mixed pairs of singular/cohesive modalities
  remaining besides those in Prop. \ref{SomeSingularModalitiesCommuteWithSomeCohesiveModalities}
  do {\it not} commute, in general
\vspace{-1mm}
\begin{center}
\begin{tabular}{ccc}
\cline{1-1}
\cline{3-3}
\multicolumn{1}{|c|}{
$
  \phantom{\mathclap{\vert^{\vert^{\vert^{\vert^{\vert^{\vert}}}}}}}
  \arraycolsep=1pt
  \begin{array}{cccc}
    &
    \shape &\circ& \conical
    \\
    \not\simeq
    \;\,
    &
    \conical &\circ& \shape
  \end{array}
$
}
&
\phantom{
$
  \phantom{\mathclap{\vert^{\vert^{\vert^{\vert^{\vert^{\vert}}}}}}}
  \arraycolsep=1pt
  \begin{array}{cccc}
    &
    \flat &\circ& \conical
    \\
    \simeq
    \;\,
    &
    \conical &\circ& \flat
  \end{array}
$}
&
\multicolumn{1}{|c|}{
$
  \phantom{\mathclap{\vert^{\vert^{\vert^{\vert^{\vert^{\vert}}}}}}}
  \arraycolsep=1pt
  \begin{array}{cccc}
    &
    \sharp &\circ& \conical
    \\
    \not\simeq
    \;\,
    &
    \conical &\circ& \sharp
  \end{array}
$
}
\\
\cline{1-1}
\cline{3-3}
\phantom{
$
  \phantom{\mathclap{\vert^{\vert^{\vert^{\vert^{\vert^{\vert}}}}}}}
  \arraycolsep=1pt
  \begin{array}{cccc}
    &
    \sharp &\circ& \conical
    \\
    \not\simeq
    \;\,
    &
    \conical &\circ& \sharp
  \end{array}
$
}
\\
\cline{1-1}
\cline{3-3}
\multicolumn{1}{|c|}{
$
  \phantom{\mathclap{\vert^{\vert^{\vert^{\vert^{\vert^{\vert}}}}}}}
  \arraycolsep=1pt
  \begin{array}{cccc}
    &
    \shape &\circ& \orbisingular
    \\
    \not\simeq
    \;\,
    &
    \orbisingular &\circ& \shape
  \end{array}
$
}
&
\phantom{
$
  \phantom{\mathclap{\vert^{\vert^{\vert^{\vert^{\vert^{\vert}}}}}}}
  \arraycolsep=1pt
  \begin{array}{cccc}
    &
    \flat &\circ& \conical
    \\
    \simeq
    \;\,
    &
    \conical &\circ& \flat
  \end{array}
$}
&
\multicolumn{1}{|c|}{
$
  \phantom{\mathclap{\vert^{\vert^{\vert^{\vert^{\vert^{\vert}}}}}}}
  \arraycolsep=1pt
  \begin{array}{cccc}
    &
    \sharp &\circ& \orbisingular
    \\
    \not\simeq
    \;\,
    &
    \orbisingular &\circ& \sharp
  \end{array}
$
}
\\
\cline{1-1}
\cline{3-3}
\end{tabular}
\end{center}
The failure of these equivalences reflects the key information
contained in singular cohesion.
For example:

- $\orbisingular \circ \shape \, (\HomotopyQuotient{\TopologicalSpace}{G})$
extracts homotopy fixed points of a topological $G$-space
$G \acts \, \TopologicalSpace$,

-
$\shape \circ \orbisingular (\HomotopyQuotient{\TopologicalSpace}{G})$ extracts
the geometric fixed points (Def. \ref{ShapeOfGemetricFixedLoci}).

\noindent
The fact that the latter are a finer invariant
than the former is the hallmark of proper-equivariant homotopy theory,
here brought out as an aspect of cohesive homotopy theory.
\end{remark}

\begin{example}[Smooth charts are smooth and orbisingular]
\label{SmoothChartsAreSmoothAndOrbisingular}
In a singular-cohesive $\infty$-topos $\Topos$
(Def. \ref{SingularCohesiveInfinityTopos})
all smooth charts $U \,\in\, \Charts$ (Def. \ref{SmoothCharts})
are both smooth as well as orbisingular
(\cite[Lem. 3.65]{SS20OrbifoldCohomology},
because they are 0-truncated)
\vspace{-3mm}
\begin{equation}
  \label{SmoothChartsBeingSmoothAndOrbisingular}
  U
    \,\in\,
  \Charts
  \xhookrightarrow{\;\;\;y\;\;\;}
  \Topos
  {\phantom{AAA}}
  \vdash
  {\phantom{AAA}}
  \smooth(U) \,\simeq\, U
  \,,
  \;\;\;\;
  \mbox{and}
  \;\;\;\;
  \orbisingular(U) \,\simeq\, U
  \,.
\end{equation}
\end{example}
\begin{remark}
  Ex. \ref{SmoothChartsAreSmoothAndOrbisingular}
  shows that being $\orbisingular$-modal
  is to be read as ``all singularities {\it that appear} are orbi-singularities'',
  subsuming the trivial case where no singularities are present in the first place.
\end{remark}

\begin{example}[Morphisms between orbi-singularizations of smooth objects]
  \label{MorphsimsBetweenOrbisingularizationsOfSmoothObjects}
  For
  $
    X, A
    \,\in\,
    \Topos_{\scalebox{.7}{\smooth}}
    \xhookrightarrow{\scalebox{0.7}{$\Space$}}
    \Topos
  $,
  we have that the hom-space \eqref{HomSpace}
  between their orbi-singularization is naturally equivalent to
  their hom-space as smooth objects:
  \vspace{-1mm}
  $$
    \begin{aligned}
    \Topos
    \left(
      \orbisingular X,
      \,
      \orbisingular A
    \right)
    &
    \;=\;
    \Topos
    \left(
      \Orbisingular \circ \Smooth(X),
      \,
      \Orbisingular \circ  \Smooth(A)
    \right)
    \\
    &
    \;\simeq\;
    \Topos
    \left(
      \Orbisingular (X),
      \,
      \Orbisingular (A)
    \right)
    \\
    &
    \;\simeq\;
    \Topos_{\scalebox{.7}{\smooth}}
    (
      X
      ,\,
      A
    )
    \end{aligned}
  $$

  \vspace{-3mm}
  \noindent  by the fact that
  $\Orbisingular : \Topos_{\scalebox{.7}{$\smooth$}}
    \xhookrightarrow{\;}
    \Topos
  $
  is fully faithful.
  More generally, for
  $G \in \Groups(\Topos_{\scalebox{.7}{$\smooth$}})$
  and
  $
    G \acts  \, X,\,  G \acts \, A
    \,\in\,
    \Actions{G}(\Topos_{\scalebox{.7}{\smooth}})
    \xhookrightarrow{\scalebox{.7}{$\Space$}}
    \Actions{G}(\Topos)
  $
  the slice hom-spaces
  \eqref{HomSpaceInSliceAsFiberProduct}
  of (the orbisingularizations of) the homotopy quotients
  are equivalent as follows:
  \vspace{-2mm}
  $$
    \Topos_{/\scalebox{.7}{$\orbisingular \mathbf{B}G$}}
    \left(
      \orbisingular(X \!\sslash\! G)
      ,\,
      \orbisingular(A \!\sslash\! G)
    \right)
    \;\;
    \simeq
    \;\;
    \Topos_{/\scalebox{.7}{$\mathbf{B}G$}}
    \left(
      X \!\sslash\! G,
      \,
      A \!\sslash\! G
    \right)
    \,.
  $$
\end{example}
\begin{remark}
  This means that effects of proper equivariance appear only
  for coefficients that are neither smooth nor just orbisingularizations
  of smooth objects. The canonical example of interesting
  coefficients for proper equivariance is the shape of orbi-singularizations
  of smooth objects.
\end{remark}

\medskip

\noindent
{\bf $G$-Singular cohesive $\infty$-toposes.}
Finally we connect the general singular-cohesive homotopy theory to traditional proper $G$-equivariant homotopy theory by observing, with \cite{Rezk14}, that the latter serves as a base of cohesion for the slice of the former over the generic $G$-orbi-singularity.

\begin{proposition}[$G$-Singular-cohesive $\infty$-topos]
  \label{GSingularCohesiveInfinityTopos}
  For $G \,\in\, \Groups(\Sets)$ a discrete group,
  the slice $\infty$-topos
  (Prop. \ref{SliceInfinityTopos})
  of any singular-cohesive $\infty$-topos
  $\Topos$  (Def. \ref{SingularCohesiveInfinityTopos})
  over the orbi-singularity
  $\scalebox{.8}{$\orbisingularG$}$ \eqref{CategoryOfSingularities}
  is cohesive over the
  proper $G$-equivariant homotopy theory of $\ModalTopos{\smooth}$
  (Def. \ref{GEquivariantAndGloballyEquivariantHomotopyTheories})
  in that
  there exists an adjoint quadruple of this form:
  \vspace{-1mm}
  \begin{equation}
    \label{GFixedLocFunctor}
    \begin{tikzcd}[column sep=40pt]
      \ModalTopos{/\orbisingularG}
      \;\;
      =
      \;\;
      \GloballyEquivariant(\ModalTopos{\smooth})_{\scalebox{.7}{$/\orbisingularG$}}
      \ar[
        rr,
        "{ G\Orbi\Smooth }"{description}
      ]
      \ar[
        rr,
        shift left=36pt,
        "{ \GOrbi\Conical }"{description},
        "\mathclap{\times}"{description, pos=0.01}
      ]
      &&
      \GEquivariant(\ModalTopos{\smooth})
      \,,
      \ar[
        ll,
        hook',
        shift right=18pt,
        "{ \GOrbi\Space }"{description}
      ]
      \ar[
        ll,
        hook',
        shift right=-18pt
        ,
        "{ \GOrbi\Singularity }"{description}
      ]
      \ar[
        ll,
        phantom,
        shift right=9pt,
        "{\scalebox{.6}{$\bot$}}"
      ]
      \ar[
        ll,
        phantom,
        shift right=27pt,
        "{\scalebox{.6}{$\bot$}}"
      ]
      \ar[
        ll,
        phantom,
        shift right=-9pt,
        "{\scalebox{.6}{$\bot$}}"
      ]
      &
      \ar[
        r,
        phantom,
        shift left=36pt,
        "{ \mbox{\tiny\color{greenii}\bf $G$-orbi-conical} }"
      ]
      \ar[
        r,
        phantom,
        shift left=16pt,
        "{ \mbox{\tiny\color{greenii}\bf $G$-orbi-space} }"
      ]
      \ar[
        r,
        phantom,
        shift left=0pt,
        "{ \mbox{\tiny\color{greenii}\bf $G$-orbi-smooth} }"
      ]
      \ar[
        r,
        phantom,
        shift right=16pt,
        "{ \mbox{\tiny\color{greenii}\bf $G$-orbi-singularity} }"
      ]
      &
      {}
    \end{tikzcd}
  \end{equation}
  where
  \begin{equation}
    \label{TheGFixedLociFunctor}
    G\Orbi\Smooth(\mathcal{X})
    \;\;\colon\;\;
    G/H
    \;\mapsto\;
    \Smooth
    \left(
      \SliceMaps{\big}{\orbisingularG}
        { \orbisingularH }
        { \mathcal{X} }
     \right)
  \end{equation}
  assigns the smooth aspect \eqref{AdjointQuadrupleOnAGlobalHomotopyTheory}
  of the slice mapping stack (Def. \ref{SliceMappingStack})
  (the geometric fixed loci, Def. \ref{ShapeOfGemetricFixedLoci}).
\end{proposition}
\noindent
Essentially this was observed in \cite{Rezk14}, there for the base case $\Topos \,=\, \InfinityGroupoids$ and for compact (hence finite if discrete) Lie groups $G$.
For the discrete equivariance groups considered here
the statement, for any $\Topos$, is a formal consequence
of the reflection
from Prop. \ref{FinteProductPreservingReflectionOfSingularitiesOntoOrbitCategory}
(as in \cite[\S 7.1, 7.2]{Rezk14}), as spelled out in the following proof.
\vspace{-2mm}
\begin{proof}
  With the equivalence from Prop. \ref{SystemsOfLocalSectionsOfBundlesOfInfinityPresheaves},
  the adjoint quadruple follows
  as the Kan extension from Lem. \ref{KanExtensionOfAdjointPairOfInfinityFunctors}
  of
  the adjunction of sites in
  Prop. \ref{GOrbitCategoryIsReflectiveSubcategoryOfSliceOverGSingularity}:
  \vspace{-3mm}
$$
    \phantom{AAAAAA}
    \begin{tikzcd}[row sep=0pt, column sep=20pt]
      \Singularities_{/\scalebox{.7}{$\orbisingularG$}}
      \ar[
        rr,
        shift left=5pt,
        "{ \tau }"{above}
      ]
      &&
      \OrbitCategory{G}
      \ar[
        ll,
        hook',
        shift left=5pt,
        "{ i }"{below}
      ]
      \ar[
        ll,
        phantom,  shift left=0pt,
        "{ \scalebox{.7}{$\bot$} }"
      ]
    \end{tikzcd}
   $$
   \begin{equation}
       \Rightarrow
    \;\;
    \begin{tikzcd}[column sep=40pt]
      \GloballyEquivariant(\ModalTopos{\smooth})_{\scalebox{.6}{$/\orbisingularG$}}
      \ar[
        r,
        "{
          \sim
        }"
      ]
      &[-10pt]
      \InfinityPresheaves
      \left(
        \Singularities_{/\scalebox{.7}{$\orbisingularG$}}
        ,\,
        \ModalTopos{\smooth}
      \right)
      \ar[
        rr,
        "{ \tau_\ast \,\simeq\, i^\ast    }"{description}
      ]
      \ar[
        rr,
        shift left=32pt,
        "{ \tau_! }"{description, pos=.39},
        "\mathclap{\times}"{description, pos=0}
      ]
      &&[-10pt]
      \InfinitySheaves
      \left(
        G \Orbits
        ,\,
        \ModalTopos{\smooth}
      \right)
      \;=\;
      \GEquivariant{\ModalTopos{\smooth}}\;.
      \ar[
        ll,
        hook',
        shift right=16pt,
        "{ \tau^\ast \,\simeq\, i_! }"{description}
      ]
      \ar[
        ll,
        hook',
        shift right=-15pt,
        "{ i_\ast }"{description, pos=.39}
      ]
      \ar[
        ll,
        phantom,
        shift right=8pt,
        "{\scalebox{.6}{$\bot$}}"
      ]
      \ar[
        ll,
        phantom,
        shift right=22pt,
        "{\scalebox{.6}{$\bot$}}"
      ]
      \ar[
        ll,
        phantom,
        shift right=-8pt,
        "{\scalebox{.6}{$\bot$}}"
      ]
    \end{tikzcd}
      \label{ConstructingTheRelativeSingularAdjoints}
  \end{equation}
  Moreover, that the left adjoint $\tau_!$ preserves finite products
  follows
  (by  Lem. \ref{LeftKanExtensionOnRepresentablesIsOriginalFunctor},
  Lem. \ref{LeftKanExtensionPreservesBinaryProductsIfOriginalFunctorDoes})
  since $\tau$ preserves finite products after extension to
  free cocompletions.

  The second statement \eqref{TheGFixedLociFunctor} follows
  with the explicit realization of the
  equivalence in Prop. \ref{SystemsOfLocalSectionsOfBundlesOfInfinityPresheaves}
  and using the definition of $i^\ast$:
  \vspace{-2mm}
  $$
  \hspace{-5mm}
    \def\arraystretch{2.4}
    \begin{array}{lll}
      G\Orbi\Smooth(\mathcal{X})
      &
      \;\simeq\;
      \bigg(
        G/H
        \,\mapsto\,
        \Presheaves
        \left(
          \Singularities_{\scalebox{.7}{$\orbisingularG$}}
          ,\,
          \ModalTopos{\smooth}
        \right)_{\scalebox{.7}{$\orbisingularG$}}
        \left(
          \orbisingularH
          ,\,
          \mathcal{X}
        \right)
      \!\!
      \bigg)
      &
      \mbox{\small by Prop. \ref{SystemsOfLocalSectionsOfBundlesOfInfinityPresheaves} }
      \\
      &
      \;\simeq\;
      \Bigg(
        G/H
        \,\mapsto\,
        \Presheaves
        \left(
          \Singularities_{\scalebox{.7}{$\orbisingularG$}}
          ,\,
          \ModalTopos{\smooth}
        \right)
        \left(
          \orbisingularH
          ,\,
          \mathcal{X}
        \right)
        \underset{
          \scalebox{.7}{$
          \InfinityPresheaves
          \left(
            \Singularities_{\scalebox{.7}{$\orbisingularG$}}
            ,\,
            \ModalTopos{\smooth}
          \right)
          (
            \orbisingularH
            ,\,
            \orbisingularG
          )
          $}
        }{\times}
        \big\{
          \orbisingularAny{i_H}
        \big\}
      \!\Bigg)
      &
      \mbox{\small by Prop. \ref{HomSpaceInSliceToposAsFiberProduct} }
      \\
      & \;\simeq\;
      \Big(
        G/H
        \,\mapsto\,
        \mathcal{X}(\orbisingularH)
        \underset{
          \scalebox{.7}{$
            \orbisingularG(\orbisingularH)
          $}
        }{\times}
        \big\{
          \orbisingularAny{i_H}
        \big\}
      \Big)
      &
      \mbox{\small by Prop. \ref{InfinityYonedaLemma} }
      \\
      & \;\simeq\;
      \Bigg(
        G/H
        \,\mapsto\,
        \Smooth
        \left(
        \Maps{}
          { \orbisingularH }
          { \mathcal{X} }
        \right)
        \underset{
          \scalebox{.7}{$
            \Smooth
            \left(
              \Maps{}
                { \orbisingularH }
                { \orbisingularG }
            \right)
          $}
        }{\times}
        \left\{
          \orbisingularAny{i_H}
        \right\}
     \! \Bigg)
      &
      \mbox{\small
        by \eqref{ValueOnOrbisingularityIsSmoothAspectOfMappingStack}
      }
      \\[2pt]
      & \;\simeq\;
      \Bigg(
        G/H
        \,\mapsto\,
        \Smooth
        \Big(
        \Maps{}
          { \orbisingularH }
          { \mathcal{X} }
        \underset{
          \scalebox{.7}{$
            \Maps{}
              { \orbisingularH }
              { \orbisingularG }
          $}
        }{\times}
        \big\{
          \orbisingularAny{i_H}
        \big\}
        \Big)
      \!\Bigg)
      &
      \mbox{\small by \eqref{InfinityAdjointPreservesInfinityLimits} }
      \\
      & \;\simeq\;
      \Big(
        G/H
        \,\mapsto\,
        \Smooth
        \big(
          \SliceMaps{}{\orbisingularG}
            { \orbisingularH }
            { \mathcal{X} }
        \big)
      \!\Big)
      &
      \mbox{\small by Def. \ref{SliceMappingStack}. }
    \end{array}
  $$

    \vspace{-8mm}
\end{proof}

\begin{definition}[$G$-orbi-singular modalities]
\label{ModalitieswithrespecttoGOrbiSingularities}
  Given a singular-cohesive $\infty$-topos (Def. \ref{SingularCohesiveInfinityTopos})
  and $G \,\in\, \Groups(\Sets)$, we write, in generalization of
  \eqref{SingularModalities}:
  \vspace{-2mm}
\begin{equation}
\label{GSingularModalities}
  \begin{array}{clcc}
    \raisebox{0pt}{
      \tiny
      \color{darkblue}
      \bf
      \begin{tabular}{c}
        singular aspect relative
        \\
        to $G$-orbi-singularities
      \end{tabular}
    }
    &
    \conicalrelativeG
    &
    \coloneqq
    &
    G\Orbi\Space
    \circ
    G\Orbi\Conical
    \\
    &
    \;\scalebox{.7}{$\bot$}
    \\
    \raisebox{0pt}{
      \tiny
      \color{darkblue}
      \bf
      \begin{tabular}{c}
        smooth aspect relative
        \\
        to $G$-orbi-singularities
      \end{tabular}
    }
    &
    \smoothrelativeG
    &
    \coloneqq
    &
    G\Orbi\Space
     \circ
    G\Orbi\Smooth
    \\
    &
    \;\scalebox{.7}{$\bot$}
    \\
    \raisebox{2pt}{
      \tiny
      \color{darkblue}
      \bf
      \begin{tabular}{c}
        orbisingular aspect relative
        \\
        to $G$-orbi-singularities
      \end{tabular}
    }
    &
    \orbisingularrelativeG
    &
    \coloneqq
    &
    G\Singularity
    \circ
    G\Orbi\Smooth
  \,.
  \end{array}
\end{equation}
\end{definition}

\begin{example}[Idempotency of $G$-orbi-singular modalities]
  \label{IdempotencyOfGOrbiSingularModalities}
  By the general idempotency of cohesive modalities (Rem. \ref{IdempotencyOfCohesiveAdjoints})
  we have, for instance, that $G$-orbi-spaces
  \eqref{GSingularCohesiveInfinityTopos}
  are $G$-orbi-smooth (Def. \ref{ModalitieswithrespecttoGOrbiSingularities})
  and
  hence coincide with their $G$-orbi-conical aspect \eqref{GSingularModalities}:
  $$
    G\Orbi\Space
    \;\simeq\;
    \smoothrelativeG
    \circ
    G\Orbi\Space
    \,,
    \;\;\;\;\;
    \conicalrelativeG
    \circ
    G\Orbi\Space
    \;\simeq\;
    G\Orbi\Space
    \,.
  $$
\end{example}

\begin{example}[Orbi-singularization of homotopy quotient]
  \label{OrbiSingularizationOfHomotopyQuotient}
  For $\Topos$ a singular-cohesive $\infty$-topos,
  given $G \,\acts\, X \,\in\, \Actions{G}(\ModalTopos{\smooth})
  \xhookrightarrow{\;} \Actions{G}(\Topos)$,
  we may naturally regard the orbi-singularization
  \eqref{SingularModalities}
  of its homotopy quotient (Ntn. \ref{HomotopyQuotientOfSimplicialGroupActions})
  as an object in the slice of $\Topos$ over the $G$-orbisingularity
  \eqref{GFixedLocFunctor}:
  $$
    \orbisingular
    \Bigg(
    \scalebox{0.9}{$
    \begin{array}{c}
      \HomotopyQuotient{ X }{ G }
      \\
      \downarrow{}
      \mathrlap{
          \scalebox{.7}{$
            \HomotopyQuotient{(X \to \ast)}{G}
          $}
      }
      \\
      \HomotopyQuotient{ \ast }{ G }
    \end{array}
    $}
    \;\;\;\;\;\;\;\;\;
    \Bigg)
    \;\;\;\;
    =
    \;\;\;\;
    \Bigg(\!\!
     \scalebox{0.9}{$
    \begin{array}{c}
      \orbisingular (\HomotopyQuotient{X}{G})
      \\
      \downarrow
      \\
      \orbisingularG
    \end{array}
    $}
    \!\!
    \Bigg)
    \;\;\;
    \in
    \;
    \SliceTopos{\orbisingularG}
    \,.
  $$
  If $X$ here is 0-truncated, then this construction lands in
  $\smoothrelativeG$-modal objects
  (by Prop. \ref{OrbiSpaceIncarnationOfGSpaceIsOrbisingularizationOfHomotopyQuotient}).
\end{example}

\begin{lemma}[Conical quotient of orbi-singularity relative to a $G$-orbi-singularity]
  \label{ConicalQuotientOfOrbiSingularityRelativeToGOrbiSingularity}
  Let
  $K, G \,\in\, \Groups(\FiniteSets)$
  and
  $K \xrightarrow{\phi} G$ any group homomorphism, regarded as an object
  $\orbisingularK \,\in\, \Slice{\Singularities}{\orbisingularG}
    \xhookrightarrow{\YonedaEmbedding} \SliceTopos{\orbisingularG}$.
  Then the {\it $G$-relative conical aspect} \eqref{GSingularModalities}
 of $\phi$  is given by the image factorization
  $K \twoheadrightarrow \mathrm{im}(\phi) \hookrightarrow G$,
  in that:
  $$
    \conicalrelativeG
    (\orbisingularK)
    \;\simeq\;
    \orbisingularAny{\mathrm{im}(\phi)}
    \;\;\;
    \in
    \SliceTopos{\orbisingularG}
    \,.
  $$
\end{lemma}
\begin{proof}
  By \eqref{ConstructingTheRelativeSingularAdjoints}
  with Lemma \ref{LeftKanExtensionOnRepresentablesIsOriginalFunctor}
  and Prop. \ref{GOrbitCategoryIsReflectiveSubcategoryOfSliceOverGSingularity}.
\end{proof}

\medskip

\noindent
{\bf Cohesive $G$-Orbispaces.}
The following terminology (Ntn. \ref{CohesiveGOrbispace})
is the direct generalization of that
of \cite{Rezk14} (following \cite{HenriquesGepner07})
from discrete to cohesive global homotopy theory,
as discussed in \cite[\S 4.1]{SS20OrbifoldCohomology}.
\begin{notation}[Cohesive $G$-orbispace]
  \label{CohesiveGOrbispace}
  For $\Topos$ a singular-cohesive $\infty$-topos
  and $G \,\in\, \Groups(\Sets)$,
  we say that those
  objects
  $\mathcal{X} \,\in\, \SliceTopos{\orbisingularG}$,
  which are
  $\smoothrelativeG$-modal (Def. \ref{ModalitieswithrespecttoGOrbiSingularities}),
  hence in the image of the full inclusion \eqref{GFixedLocFunctor}
  \vspace{-2mm}
  $$
    \GEquivariant(\ModalTopos{\smooth})
    \xhookrightarrow{ \; G\OrbiSpaces \; }
    \GloballyEquivariant(\ModalTopos{/\orbisingularG})
  $$

  \vspace{-2mm}
 \noindent are the {\it cohesive $G$-orbi-spaces}.
\end{notation}
\begin{remark}[Interpretation of $G$-orbispaces]
  The cohesive $G$-orbispaces (Ntn. \ref{CohesiveGOrbispace})
  may be thought of as those singular-cohesive spaces
  which are {\it smooth away from $G$-singularities}
  or {\it smooth except possibly for $G$-singularities}.
\end{remark}
\begin{example}[The terminal $G$-orbispace]
  \label{TheTerminalGOrbispace}
  For $G \,\in\, \Groups(\Sets)$ with
  $\mathbf{B}G \,\in\, \SliceTopos{\mathbf{B}G}$
  and
  $\orbisingularG \,\in\, \SliceTopos{\orbisingularG}$
  regarded as the identity maps on themselves, hence
  as the respective terminal objects (by Ex. \ref{BasicStructuresInSliceInfinityTopos}),
  we have that the
  $G$-smooth aspect \eqref{GSingularModalities}
  of the shape \eqref{TheModalitiesOnACohesiveInfinityTopos}
  of the orbi-singularization \eqref{SingularModalities}
  of
  the former terminal object $\mathbf{B}G$ is the latter terminal object $\orbisingularG$,
  hence the terminal $G$-orbispace (Ntn. \ref{CohesiveGOrbispace}):
  $$
    \smoothrelativeG
    \,
    \shape
    \,
    \orbisingular
    \,
    \mathbf{B}G
    \;\;
    \simeq
    \;\;
    \orbisingularG
    \;\;\simeq\;\;
    \ast_{\scalebox{.7}{\orbisingularG}}
    \;\;\;\;\;
    \in
    \;
    \GEquivariant\ModalTopos{\smooth}
    \xhookrightarrow{ G\Space }
    \ModalTopos{\orbisingularG}
    \,.
  $$
  This follows
  immediately since $\orbisingular$ and $\smoothrelativeG$ are right adjoints
  and hence preserve all limits \eqref{InfinityAdjointPreservesInfinityLimits}
  and since $\shape$ preserves all finite products (Rem. \ref{TheAxiomsOnTheShapeModality}).
\end{example}

\noindent
{\bf Geometric fixed loci.} We expand on the construction of
geometric fixed loci from \eqref{GFixedLocFunctor}.

\begin{definition}[Geometric fixed loci {\cite[Def. 3.69]{SS20OrbifoldCohomology}}]
  \label{ShapeOfGemetricFixedLoci}
  For

  -- $G \,\in\, \Groups(\FiniteSets)$,

  --
  $
    G \acts \, X
      \,\in\,
    \Actions{G}( \ModalTopos{\smooth} )
      \xhookrightarrow{ \; }
    \Actions{G}( \Topos )
  $:

  \noindent
  {\bf (i)}
  We say,
  for $H \subset G$ a subgroup,
  that the {\it geometric $H$-fixed locus}
  of $G \acts \, X$
  is the smooth aspect \eqref{SingularModalities}
  of the stack of sections \eqref{TheSliceMappingStack}
  of the orbi-singularization \eqref{SingularModalities}
  of its homotopy quotient \eqref{HomotopyQuotientAsHomotopyColimit}:
  \vspace{-2mm}
  \begin{align}
    \label{GeometricFixedLocus}
    H\FixedLoci(G \acts \, X)
    \;\;
    &
    \coloneqq
    \;\;
    \big(
      \FixedLoci(G \acts \, X)
    \big)(G/H)
    \;\;
    \nonumber
    \\
    &\coloneqq
    \;\;
    \smooth
    \,
    \SliceMaps{\big}{\orbisingularG}
      { \orbisingularH }
      { \orbisingular (X \!\sslash\! G) }
    \quad
    \in
    \;
    \Topos_{\scalebox{.7}{$\smooth$}}
    \xhookrightarrow{\;}
    \Topos \;.
  \end{align}
  Since orbi-singularization $\orbisingular$ is fully faithful on smooth ($\smooth$-modal) objects, this is equivalently the stacky fixed locus of Ex. \ref{FixedLociOfInfinityActions}, as it should be:
  $$
    \begin{array}{lll}
       H \FixedLoci(G \acts X)
       &
       \;\simeq\;
       \smooth
       \SliceMaps{\big}{\orbisingularG}
         { \orbisingularH }
         {
           \orbisingular(\HomotopyQuotient{X}{G})
         }
      \\
      &
      \;\simeq\;
       \SliceMaps{\big}{\mathbf{B}G}
         { \mathbf{B}H }
         {
            \HomotopyQuotient{X}{G}
         }
      \\
      &
      \;\simeq\;
      X^H
      &
      \proofstep{
        by
        \eqref{FixedPointsForSubgroupAsSliceMappingStack}
      }
    \end{array}
  $$

  \vspace{-2mm}
  \noindent
  {\bf (ii)}
  We say that the {\it shape of the geometric $H$-fixed locus} in $X$ is the
  $\infty$-groupoid of sections of
  the shape of the orbi-singularization of the homotopy quotient,
  regarded as a discrete spatial object:
  \vspace{-2mm}
  \begin{equation}
    \label{ShapeOfGeometricFixedLocus}
    \shape
    \,
    H \FixedLoci
    (
      G \acts \, X
    )
    \quad
    \overset{
      \mathclap{
      \raisebox{3pt}{
        \tiny
        Prop. \ref{ShapeOfGeometricFixedLocusIsIndeedShapeOfGeometricFixedLocus}
      }
      }
    }{
      \simeq
    }
    \qquad
    \Space
    \,
    \Discrete
    \,
    \SlicePointsMaps{\big}{\orbisingularG}
      { \orbisingularH }
      { \shape \orbisingular (X \!\sslash\! G) }
    \;\;\;
    \in
    \;
    \InfinityGroupoids
    \xhookrightarrow{\;}
    \Topos
    \;.
  \end{equation}
  \end{definition}

\begin{remark}[Orbi-spaces]
$\,$

\noindent {\bf (i)}
  By functoriality of the mapping stack construction,
  as $H \subset G$ varies the fixed loci
 \eqref{GeometricFixedLocus}
  arrange into
  an object of the $G$-equivariant homotopy theory of $\ModalTopos{\smooth}$
  (Def. \ref{GEquivariantAndGloballyEquivariantHomotopyTheories}),
  which, according to Prop. \ref{GSingularCohesiveInfinityTopos},
  coincides with the smooth aspect relative to $G$-singularities
  \eqref{TheGFixedLociFunctor} of the orbi-singularization of the
  homotopy quotient:
  \vspace{-1mm}
  \begin{equation}
    \label{SystemOfGeometricFixedLociAsObjectInGEquivariantHomotopyTheory}
    \FixedLoci(G \acts  \, X)
    \;\;
    =
    \;\;
    G\Orbi\Smooth
    \big(
      \orbisingular( \HomotopyQuotient{X}{G} )
    \big)
    \;=\;
    \Smooth
    \,
    \SliceMaps{\big}{\orbisingularG}
      { \orbisingularAny{\scalebox{.7}{$(-)$}} }
      { \orbisingular (\HomotopyQuotient{X}{G}) }
    \;\;\;
    \in
    \;
    \GEquivariant\ModalTopos{\smooth}
    \,.
  \end{equation}

  \noindent
  {\bf (ii)}
  Now, as $H \subset G$ varies, the fixed loci
  \eqref{ShapeOfGeometricFixedLocus}
  arrange into
  an object of the $G$-equivariant homotopy theory of
  $\InfinityGroupoids \xhookrightarrow{\Discrete}  \ModalTopos{\smooth}$
  (Def. \ref{GEquivariantAndGloballyEquivariantHomotopyTheories}),
  which, according to Prop. \ref{GSingularCohesiveInfinityTopos},
  coincides with the smooth aspect relative to $G$-singularities
  \eqref{TheGFixedLociFunctor} of
  {\it the shape of} the orbi-singularization of the
  homotopy quotient:
  \vspace{-2mm}
  \begin{equation}
    \label{SystemOfShapesOfGeometricFixedLociAsObjectInGEquivariantHomotopyTheory}
    \shape
    \,
    \FixedLoci(G \acts  \, X)
    \;\;
    \simeq
    \;\;
    G\Orbi\Smooth
    \big(
      \shape
      \,
      \orbisingular
      ( \HomotopyQuotient{X}{G} )
    \big)
    \;\;\;
    \in
    \;
    \begin{tikzcd}[row sep=1pt]
    &
    \GEquivariant{\ModalTopos{\smooth}}
    \ar[
      dr,
      hook,
      "{
        G\OrbiSpaces
       }"{sloped}
    ]
    \\
    \GEquivariant\InfinityGroupoids
    \ar[
      ur,
      hook,
      "{\Discrete}"{sloped}
    ]
    \ar[
      dr,
      hook,
      "{
        G\OrbiSpaces
       }"{sloped, swap, pos=.4}
    ]
    &&
    \Slice{\Singular{(\ModalTopos{\smooth})}}{\orbisingularG}
    \\
    &
    \Slice{\Singular{\InfinityGroupoids}}{\orbisingularG}
    \ar[
      ur,
      hook,
      shorten <=-7pt,
      "{\Discrete}"{sloped, swap}
    ]
    \end{tikzcd}
    \,.
  \end{equation}
  \vspace{-.5cm}

  Regarded under this embedding, we have
  \begin{equation}
    \label{OrbiSpacesAreSmoothRelativeGSingularities}
    \shape
    \,
    \FixedLoci(G \acts X)
    \;\simeq\;
    \smoothrelativeG
    \,
    \shape
    \,
    \FixedLoci(G \acts X)
    \;\simeq\;
    \smoothrelativeG
    \,
    \shape
    \orbisingular
    (\HomotopyQuotient{X}{G})
    \,,
  \end{equation}
  reflecting the fact that any such {\it $G$-orbi-space}
  (see \cite[\S 4.1]{SS20OrbifoldCohomology} for more discussion and further pointers to the literature)
  is smooth (non-singular) except for orbi-singularities corresponding to subgroups of $G$.
\end{remark}

\begin{remark}[Cocycle $\infty$-groupoids on cohesive orbi-spaces]
  \label{CocycleInfinityGroupoidsForOrbiSpaces}
  In summary, the above implies that,
for
  $G \acts \,  X \,\in\, \Actions{G}(\ModalTopos{\smooth})$
  and $\mathcal{A} \,\in\, \GEquivariant\ModalTopos{\smooth}$,
  we have a natural equivalence
  \begin{equation}
    \label{IncarnationsOfCocyclesOnCohesiveOrbispace}
    \SlicePointsMaps{\big}{\orbisingularG}
      {\!\!\orbisingular(\HomotopyQuotient{X}{G}) }
      { G\Orbi\Singularity(\mathcal{A}) }
    \;\;
    \simeq
    \;\;
    \GEquivariant{\ModalTopos{\smooth}}
    \big(
      \FixedLoci(X)
      ,\,
      \mathcal{A}
    \big)
  \end{equation}
  obtained as the following composite:
  $$
    \def\arraystretch{1.2}
    \begin{array}{lll}
      \SlicePointsMaps{\big}{\orbisingularG}
        {\!\! \orbisingular(\HomotopyQuotient{X}{G}) }
        { G\Orbi\Singularity(\mathcal{A}) }
      & \;\simeq\;
      \GEquivariant{\ModalTopos{\smooth}}
      \big(
        {
          G\Orbi\Smooth
          \left(
            \orbisingular(\HomotopyQuotient{X}{G})
          \right)
        }
       \; ,\;
        { \mathcal{A} }
      \big)
      &
      \proofstep{ by \eqref{GFixedLocFunctor} }
      \\
      & \;\simeq\;
      \GEquivariant{\ModalTopos{\smooth}}
      \Big(
        {
          \Smooth
          \,
          \SliceMaps{\big}{\orbisingularG}
            { \orbisingularAny{\scalebox{.7}{$(-)$}} }
            { \orbisingular(\HomotopyQuotient{X}{G}) }
        }
        ,\;
        { \mathcal{A} }
      \Big)
      &
      \proofstep{ by \eqref{TheGFixedLociFunctor} }
      \\
      & \;=\;
      \GEquivariant{\ModalTopos{\smooth}}
      \left(
        \FixedLoci(X)
        ,\,
        \mathcal{A}
      \right)
      &
      \proofstep{ by \eqref{SystemOfGeometricFixedLociAsObjectInGEquivariantHomotopyTheory}. }
    \end{array}
  $$
 This is the basis for the definition of proper-equivariant cohomology
  in cohesive $\infty$-toposes (Def. \ref{BorelEquivariantAndProperEquivariantCohomologyInCohesiveInfinityTopos} below).
\end{remark}

\begin{remark}[Subsuming generalized geometric fixed points]
  \label{SubsumingGeneralizedGeometricFixedPoints}
  $\,$

  \noindent
  {\bf (i)}
  If $X \,\in\, \ModalTopos{\smooth, 0}$ is smooth and 0-truncated, then
  the geometric fixed loci of $\HomotopyQuotient{X}{G}$, according to
  Def. \ref{ShapeOfGemetricFixedLoci}, consist of the expected fixed loci
  of the given $G$-$\infty$-action on $X$.

  \noindent
  {\bf (ii)}
  However, if
  $X$ is not 0-truncated, then its isotropy groups, as far as they receive
  non-trivial homomorphisms from the subgroups $H \subset G$,
  do contribute to the geometric fixed loci
  in the sense of Def. \ref{ShapeOfGemetricFixedLoci}.

  \noindent
  {\bf (iii)}
  The archetypical example of this effect occurs when
  $X \,=\, \mathbf{B}\Gamma$ is a delooping \eqref{LoopingAndDeloopingEquivalence},
  hence a single point with isotropy.
  This is exactly the case of equivariant moduli stacks,
  discussed in \cref{EquivariantModuliStacks} below.
  Their generalized geometric fixed loci is the source of all the
  interesting and characteristic structure of equivariant classifying spaces
  embodied by the Murayama-Shimakawa construction
  (Thm. \ref{MurayamaShimakawaGroupoidIsEquivariantModuliStack} below).
\end{remark}

We proceed to demonstrate the equivalence
\eqref{SystemOfShapesOfGeometricFixedLociAsObjectInGEquivariantHomotopyTheory}
in Prop. \ref{ShapeOfGeometricFixedLocusIsIndeedShapeOfGeometricFixedLocus} below:

\begin{lemma}
\label{ShapeOfFixedLocusIsSmoothAspectOfInternalFixedLocus}
The shape of the geometric fixed locus is the smooth aspect
\eqref{SingularModalities}
of the
internal fixed locus
  \vspace{-1mm}
$$
  \smooth
  \circ
  \SliceMaps{\big}{\orbisingularG}
    { \orbisingularG }
    { \shape  \,\mathcal{X} }
  \;\simeq\;
  \Space
  \circ
  \Discrete
  \circ
  \SlicePointsMaps{\big}{\orbisingularG}
    { \orbisingularG }
    { \shape \, \mathcal{X} }
  \,.
$$
\end{lemma}
\begin{proof}
  For $\orbisingularK \times U \,\in\, \Singularities \times \Charts$,
  we have the following sequence of natural equivalences:
  \vspace{-2mm}
  $$
  \hspace{-2mm}
    \begin{array}{lll}
            \Big(\!\!
        \smooth
        \,
        \SliceMaps{\big}{\orbisingularG}
          { \orbisingularG }
          { \shape \, \mathcal{X} }
      \Big)
      \left(
        \orbisingularK \times U
      \right)
           &
      \!\!\! \!\simeq
      \Big(
        \SliceMaps{\big}{\orbisingularG}
          { \orbisingularG }
          { \shape \, \mathcal{X} }
      \Big)
      (U)
      &
      \mbox{\small  by \cite[Def. 3.52]{SS20OrbifoldCohomology} }
      \\
      &    \!\!\! \!\simeq
      \SlicePointsMaps{\big}{\orbisingularG}
        { \orbisingularG \times U}
        { \shape \, \mathcal{X} }
      &
      \mbox{\small by Lem. \ref{PlotsOfSliceMappingStackAreSliceHoms} }
      \\
      &   \!\!\! \!\simeq
      \SlicePointsMaps{\big}{\orbisingularG}
        { \orbisingularG }
        { \shape \, \mathcal{X} }
      &
      \mbox{\small by \cite[Def. 3.1]{SS20OrbifoldCohomology} }
      \\
      &   \!\!\! \!\simeq
      \Big(
        \Space
        \circ
        \Discrete
        \circ
        \SlicePointsMaps{\big}{\orbisingularG}
          { \orbisingularG }
          { \shape \mathcal{X} }
      \Big)
      \left(
        \orbisingularK \times U
      \right)
      &
      \mbox{\small by \cite[Def. 3.50]{SS20OrbifoldCohomology}. }
          \end{array}
  $$

    \vspace{-2mm}
\noindent
  Therefore, the claim follows by the $\infty$-Yoneda lemma
  (Lem. \ref{InfinityYonedaLemma}).
\end{proof}

The following simple version (Prop. \ref{ShapeOfMappingStackOutOfOrbiSingularityIsMappingStackIntoShape}) of an ``orbi-smooth Oka principle'' \eqref{OrbiSmoothOkaPrincipleInIntroduction} will serve as a workhorse lemma in following proofs. Notice here that by our definitions the groups $H$ and $G$ are discrete, but that there is absolutely no condition on the geometry embodied in $\mathcal{X}$.
\begin{proposition}[Orbi-smooth Oka principle for maps out of relative orbi-singularities]
  \label{ShapeOfMappingStackOutOfOrbiSingularityIsMappingStackIntoShape}
  Let $\Topos$ be a singular-cohesive $\infty$-topos (Def. \ref{SingularCohesiveInfinityTopos}).
  Then for
  \vspace{-2mm}
  \begin{itemize}
  \setlength\itemsep{-3pt}

  \item[--]
  $\orbisingularG \,\in\,\Singularities$ {\rm (Ntn. \ref{Singularities})},

  \item[--]
  $\orbisingularH \,\in\,\Singularities_{/\scalebox{.7}{$\orbisingularG$}}$\,,

  \item[--]
  $\mathcal{X} \,\in\, \Topos_{/\scalebox{.7}{$\orbisingularG$}}$\,,
  \end{itemize}

  \vspace{-2mm}
  \noindent
  we have a natural equivalence
  \vspace{-2mm}
  \begin{equation}
    \label{ShapeOfSliceMappingStackOutOfOrbisingularity}
    \shape
    \,
    \SliceMaps{\big}{\orbisingularG}
      { \orbisingularH }
      { \mathcal{X} }
    \;\simeq\;
    \SliceMaps{\big}{\orbisingularG}
      { \shape \, \orbisingularH }
      { \shape \, \mathcal{X} }
    \;\simeq\;
    \SliceMaps{\big}{\orbisingularG}
      { \orbisingularH }
      { \shape \, \mathcal{X} }
  \end{equation}

  \vspace{-2mm}
  \noindent
  exhibiting the shape modality
  \eqref{TheModalitiesOnACohesiveInfinityTopos}
  as commuting over
  the slice mapping stack construction (Def. \ref{SliceMappingStack})
  out of orbi-singularities.
\end{proposition}
\begin{proof}
For the special case when $\orbisingularG = \ast$,
i.e. $G \,=\, 1$ the trivial group,
hence the case when $\orbisingularH \in \Singularities$,
$X \in \Topos$ and $\Maps{}{\orbisingularH}{X}$ the ordinary mapping stack
\eqref{InternalHomAdjunction},
consider the following sequence of natural equivalences:
\vspace{-3mm}
\begin{equation}
  \label{ShapeOfMapingStackOutOfOrbisingularity}
  \hspace{-1cm}
  \def\arraystretch{2.3}
  \begin{array}{lll}
    \shape \;
    \Maps{\big}
      { \orbisingularH }
      { \mathcal{X} }
    & \;=\;
    \shape
    \Big(\!\!
      \big(\orbisingularK, U \big)
      \,\mapsto\,
    \Maps{\big}
      { \orbisingularH }
      { \mathcal{X} }
    \big(U \times \orbisingularK \big)
    \!\!\Big)
    \\
    &
    \;\simeq\;
    \shape
    \Big(\!\!
      \big(\orbisingularK, U \big)
      \,\mapsto\,
      \PointsMaps{\big}
      {
        U \times \orbisingularH \times \orbisingularK
      }
      {
        \mathcal{X}
      }
   \!\! \Big)
    &
    \mbox{\small by \eqref{ValuesOfMappingStackAsHomSpaces}}
    \\
    & \;\simeq\;
    \shape
    \bigg(\!\!
      \big(\orbisingularK, U\big)
      \,\mapsto\,
      \Big(
        \mathcal{X}
        \big(
          \orbisingularH \times \orbisingularK \times U
        \big)
      \Big)
    \!\!\bigg)
    &
    \mbox{\small by \eqref{YonedaEquivalence}}
    \\
    & \;\simeq\;
    \shape
    \bigg(\!\!
      \left(\orbisingularK, U\right)
      \,\mapsto\,
      \Big(
        \mathcal{X}
        \big(
          \orbisingularH \times \orbisingularK
        \big)
      \Big)(U)
    \!\!\bigg)
    &
    \mbox{\small by \eqref{SingularCohesiveToposAsSheavesOnSingulartiesTimesCharts}}
    \\
    & \;\simeq\;
    \Bigg(\!\!
      \orbisingularK
      \,\mapsto\,
      \bigg(
        \shape
        \Big(
          \mathcal{X}
          \big(
            \orbisingularH \times \orbisingularK
          \big)
        \Big)
      \bigg)
    \!\!\Bigg)
    &
    \mbox{\small by \eqref{SingularityWiseApplicationOfCohesiveModalities}}
    \\
    & \;\simeq\;
    \bigg(\!\!
      \orbisingularK
      \,\mapsto\,
      \Big(
        \big( \shape \mathcal{X}\big)
        \big(\orbisingularH \times \orbisingularK\big)
      \Big)
    \!\!\bigg)
    &
    \mbox{\small by \eqref{SingularityWiseApplicationOfCohesiveModalities}}
    \\
    & \;\simeq\;
    \bigg(\!\!
      \big( \orbisingularK, U \big)
      \,\mapsto\,
      \Big(
        \big(\shape \mathcal{X}\big)
        \big(
          \orbisingularH \times \orbisingularK \times U
        \big)
      \Big)
    \!\!\bigg)
    &
    \mbox{\small by \eqref{SingularCohesiveToposAsSheavesOnSingulartiesTimesCharts}}
    \\
    & \;\simeq\;
    \Maps{\big}
      { \orbisingularH }
      { \shape \mathcal{X} }
    &
    \mbox{\small by \eqref{ValuesOfMappingStackAsHomSpaces}}
    \,.
  \end{array}
\end{equation}

  \vspace{-2mm}
  \noindent
The composite of these equivalence yields \eqref{ShapeOfSliceMappingStackOutOfOrbisingularity}
for the special case $\orbisingularG = \ast$.

To obtain from this the statement for general $\orbisingularG$,
consider the following sequence
of natural equivalences:
\vspace{-2mm}
$$
\def\arraystretch{1.5}
  \begin{array}{lll}
    \shape \;
    \SliceMaps{\big}{\orbisingularG}
      { \orbisingularH }
      { X }
    &
    \;=\;
    \shape
    \,
    \bigg(
      \Maps{\big}
        { \orbisingularH }
        { X }
      \underset{
        \scalebox{.7}{$
          \Maps{\big}
            { \orbisingularH }
            { \orbisingularG }
        $}
      }{\times}
      \ast
    \bigg)
    &
    \mbox{\small by Def. \ref{SliceMappingStack}}
    \\
    & \;\simeq\;
    \Big(
      \shape \;
      \Maps{\big}
        { \orbisingularH }
        { X }
    \Big)
    \underset{
      \scalebox{.7}{$
        \Maps{\big}
          { \orbisingularH }
          { \orbisingularG }
      $}}{\times}
    \big(
      \shape
      \ast
    \big)
    &
    \mbox{\small by Prop. \ref{ShapeFunctorPreservesHomotopyFibersOverDiscreteObjects} with Lem. \ref{MappingStackIntoDiscreteObjectIsDiscrete}}
    \\
    &
    \;\simeq\;
    \Maps{\big}
      { \orbisingularH }
      { \shape  \, X }
    \underset{
      \scalebox{.7}{$
        \Maps{\big}
          { \orbisingularH }
          { \orbisingularG }
      $}}{\times}
    \ast
    &
    \mbox{\small by \eqref{ShapeOfMapingStackOutOfOrbisingularity}
      and \eqref{ShapePreservesTheTerminalObject} }
    \\
    & \;=\;
    \SliceMaps{\big}{\orbisingularG}
      { \orbisingularH }
      { \shape \,  X }
    &
    \mbox{\small by Def. \ref{SliceMappingStack}}
    \,.
  \end{array}
$$

\vspace{-2mm}
\noindent The composite of these equivalences yields the desired
\eqref{ShapeOfSliceMappingStackOutOfOrbisingularity}.
\end{proof}

\begin{proposition}[Shape of geometric fixed locus]
  \label{ShapeOfGeometricFixedLocusIsIndeedShapeOfGeometricFixedLocus}
  The shape of the geometric fixed locus
  in the sense of \eqref{ShapeOfGeometricFixedLocus}
  in Def. \ref{ShapeOfGemetricFixedLoci} is indeed the
  image under $\shape$ of the geometric fixed locus
  \eqref{GeometricFixedLocus}
  \vspace{-1mm}
  $$
    \shape
    \;
    \smooth
    \,
    \SliceMaps{\big}{\orbisingularG}
      { \orbisingularG }
      { \orbisingular (X \!\sslash\! G) }
    \;\;
    \simeq
    \;\;
    \Space
    \circ
    \Discrete
    \circ
    \SlicePointsMaps{\big}{\orbisingularG}
      { \orbisingularG }
      { \shape \, \orbisingular (X \!\sslash\! G) } \;.
  $$
\end{proposition}
\begin{proof}
$\,$

\vspace{-7mm}
$$
  \begin{array}{lll}
        \shape
    \,
    \smooth
    \,
    \SliceMaps{\big}{\orbisingularG}
      { \orbisingularG }
      { \orbisingular (X \!\sslash\! G) }
    &
    \;\simeq\;
    \smooth
    \,
    \shape
    \;
    \SliceMaps{\big}{\orbisingularG}
      { \orbisingularG }
      { \orbisingular (X \!\sslash\! G) }
    &
    \mbox{\small by \cite[3.67]{SS20OrbifoldCohomology} }
    \\
    & \;\simeq\;
    \smooth
    \,
    \SliceMaps{\big}{\orbisingularG}
      { \orbisingularG }
      { \shape \orbisingular (X \!\sslash\! G) }
    &
    \mbox{\small by Prop. \ref{ShapeOfMappingStackOutOfOrbiSingularityIsMappingStackIntoShape} }
    \\
    & \;\simeq\;
    \Space
    \circ
    \Discrete
    \circ
    \SlicePointsMaps{\big}{\orbisingularG}
      { \orbisingularG }
      { \shape \orbisingular (X \!\sslash\! G) }
    &
    \mbox{\small by Lem. \ref{ShapeOfFixedLocusIsSmoothAspectOfInternalFixedLocus} }.
  \end{array}
$$

\vspace{-7mm}
\end{proof}

\begin{proposition}[Orbi-space incarnation of $G$-space is orbi-singularization of homotopy quotient]
  \label{OrbiSpaceIncarnationOfGSpaceIsOrbisingularizationOfHomotopyQuotient}
  $\,$

\noindent  For $\Topos$ a singular-cohesive $\infty$-topos (Def. \ref{SingularCohesiveInfinityTopos}),
  consider $G \acts \, X \,\in\, \Actions{G}(\ModalTopos{\smooth, 0})$

\noindent {\bf (i)} Then  there are natural equivalence
  \vspace{-2mm}
  \begin{align}
    \label{GOrbiSpaceOfSystemOfFixedLociIsOrbisingularizationOfHomotopyQuotient}
    G\Orbi\Space
    \left(
      \phantom{\shape}
      \,
      \FixedLoci(G \acts \, X)
    \right)
    &
    \;\;
    \simeq
    \;\;
    \phantom{\shape}
    \,
    \orbisingular
    \,
    (\HomotopyQuotient{ X }{ G })
    \;\;\;
    \in
    \;
    \phantom{\big(}
    \SliceTopos{\orbisingularG}
    \\
    \label{GOrbiSpaceOfSystemOfShapesOfFixedLociIsOrbisingularizationOfShapeOfHomotopyQuotient}
    G\Orbi\Space
    \left(
      \shape
      \,
      \FixedLoci(G \acts \, X)
    \right)
    &
    \;\;
    \simeq
    \;\;
    \shape
    \,
    \orbisingular
    \,
    (\HomotopyQuotient{ X }{ G })
    \;\;\;
    \in
    \;
    \Slice{(\ModalTopos{\flat})}{\orbisingularG}
    \xhookrightarrow{\;}
    \SliceTopos{\orbisingularG}
  \end{align}

  \vspace{-3mm}
\noindent  between
\vspace{-2mm}
  \begin{itemize}
   \setlength\itemsep{-3pt}
   \item[{\bf (a)}] the $G$-orbispace
  \eqref{GFixedLocFunctor}
  associated with the image
  \eqref{SystemOfGeometricFixedLociAsObjectInGEquivariantHomotopyTheory}
  of $G \acts \, X$ in cohesive $G$-equivariant homotopy theory,
  and

  \item[{\bf (b)}]  the orbisingularization of its cohesive homotopy quotient
  (Ex. \ref{OrbiSingularizationOfHomotopyQuotient}).
  \end{itemize}

\vspace{-2mm}
\noindent {\bf (ii)}   In particular, the orbi-singularization of the homotopy quotient of a
  0-truncated cohesive space is
  smooth away from $G$-singularities (Def. \ref{ModalitieswithrespecttoGOrbiSingularities}):
  \vspace{-3mm}
  \begin{align*}
&    \smoothrelativeG
    \Bigg(
    \!\!\!\!\!\!
    \scalebox{.9}{
     \begin{tikzcd}[row sep=7pt]
      \orbisingular (X \!\sslash\! G)
      \ar[d]
      \\
      \orbisingularG
    \end{tikzcd}
    }
    \!\!\!\!\!
    \Bigg)
    \;\;
    \simeq
    \;\;
    \Bigg(
    \!\!\!\!\!\!
    \scalebox{.9}{
    \begin{tikzcd}[row sep=7pt]
      \orbisingular (X \!\sslash\! G)
      \ar[d]
      \\
      \orbisingularG
    \end{tikzcd}
    }
    \!\!\!\!\!
    \Bigg)
    \;\;\;
    \in
    \;
    \SliceTopos{\orbisingularG}
    \,,
\\
  &  \smoothrelativeG
    \Bigg(
    \!\!\!\!\!
    \scalebox{.9}{
     \begin{tikzcd}[row sep=7pt]
      \shape\, \orbisingular (X \!\sslash\! G)
      \ar[d]
      \\
      \orbisingularG
    \end{tikzcd}
    }
    \!\!\!\!\!
    \Bigg)
    \;\;
    \simeq
    \;\;
    \Bigg(
    \!\!\!\!
    \scalebox{.9}{
    \begin{tikzcd}[row sep=7pt]
      \shape\, \orbisingular (X \!\sslash\! G)
      \ar[d]
      \\
      \orbisingularG
    \end{tikzcd}
    }
    \!\!\!\!\!
    \Bigg)
    \;\;\;
    \in
    \;
    \Slice{(\ModalTopos{\flat})}{\orbisingularG}
    \,.
  \end{align*}
\end{proposition}
\begin{proof}
We show the proof of \eqref{GOrbiSpaceOfSystemOfShapesOfFixedLociIsOrbisingularizationOfShapeOfHomotopyQuotient};
the proof of \eqref{GOrbiSpaceOfSystemOfFixedLociIsOrbisingularizationOfHomotopyQuotient}
follows by the exact same steps, just with all occurences of $\shape$ removed.
For $\orbisingularK \,\in\, \Slice{\Singularities}{\orbisingularG}$\,,
we have the following sequence of natural equivalences in $\ModalTopos{\smooth}$:
\vspace{-3mm}
$$
\hspace{-2mm}
\def\arraystretch{1.5}
  \begin{array}{lll}
        \GloballyEquivariant\ModalTopos{\smooth}
    \big(
      {
        \orbisingularK
      }
      ,\,
      {
        G\Orbi\Space
        \left(
          \shape
          \,
          \FixedLoci
          (\TopologicalSpace)
        \right)
      }
    \big)_{\scalebox{.7}{$\orbisingularG$}}
        &
    \simeq\;
    \GEquivariant\ModalTopos{\smooth}
    \big(
      G\Orbi\Conical(
        \orbisingularK
      )
      ,\,
      \shape
      \,
      \FixedLoci
      (\TopologicalSpace)
    \big)
    &
    \proofstep{ by \eqref{GFixedLocFunctor} }
    \\
    &
    \simeq\;
    \GEquivariant\ModalTopos{\smooth}
    \big(
      \orbisingularAny{\mathrm{im}(K)}
      ,\,
      \shape
      \,
      \FixedLoci
      (\TopologicalSpace)
    \big)
    &
    \proofstep{ by Lem. \ref{ConicalQuotientOfOrbiSingularityRelativeToGOrbiSingularity} }
    \\
    &
    \simeq\;
    \shape
    \,
    \GEquivariant\ModalTopos{\smooth}
    \big(
      \orbisingularAny{\mathrm{im}(K)}
      ,\,
      \FixedLoci
      (\TopologicalSpace)
    \big)
    &
    \proofstep{ by \eqref{SingularityWiseApplicationOfCohesiveModalities} }
    \\
    &
    \simeq\;
    \shape
    \,
    \Smooth
    \,
    \SliceMaps{\big}{\orbisingularG}
      { \orbisingularAny{\mathrm{im}(K)} }
      { \orbisingular (\HomotopyQuotient{\TopologicalSpace}{G} ) }
    &
    \proofstep{ by Prop. \ref{InfinityYonedaLemma} \&  Def. \ref{ShapeOfGemetricFixedLoci} }
    \\
    &
    \simeq\;
    \shape
    \,
    \Smooth
    \,
    \SliceMaps{\big}{\orbisingularG}
      { \orbisingularK }
      { \orbisingular (\HomotopyQuotient{\TopologicalSpace}{G} ) }
    &
    \proofstep{ by Lem. \ref{OrbiSingularizationOfHomotopyGQuotientsOfZeroTruncatedObjectsHaveOnlyGSingularities} }
    \\
    &
    \simeq\;
    \Smooth
    \,
    \shape
    \,
    \SliceMaps{\big}{\orbisingularG}
      { \orbisingularK }
      { \orbisingular (\HomotopyQuotient{\TopologicalSpace}{G} ) }
    &
    \proofstep{ by Prop. \ref{SomeSingularModalitiesCommuteWithSomeCohesiveModalities} }
    \\
    &
    \simeq\;
    \Smooth
    \,
    \SliceMaps{\big}{\orbisingularG}
      { \orbisingularK }
      { \shape \, \orbisingular (\HomotopyQuotient{\TopologicalSpace}{G} ) }
    &
    \proofstep{ by Prop. \ref{ShapeOfMappingStackOutOfOrbiSingularityIsMappingStackIntoShape} }
    \\
    &
    \simeq\;
    \GloballyEquivariant\ModalTopos{\smooth}
    \big(
      \orbisingularK
      ,\,
      \shape
      \,
      \orbisingular (\HomotopyQuotient{\TopologicalSpace}{G})
    \big)_{\scalebox{.7}{$\orbisingularG$}}
    &
    \proofstep{ by \eqref{PointsMapsSpaceIsPointsOfMappingStack} }.
  \end{array}
$$

\vspace{-2mm}
\noindent Therefore, the claim follows by the $\infty$-Yoneda lemma
(Prop. \ref{InfinityYonedaLemma}).
\end{proof}

\begin{lemma}[Orbi-singularizations of homotopy $G$-quotients of 0-truncated objects have only $G$-singularities]
  \label{OrbiSingularizationOfHomotopyGQuotientsOfZeroTruncatedObjectsHaveOnlyGSingularities}
  Given a singular-cohesive $\infty$-topos $\Topos$
  let $G \,\in\, \Groups(\FiniteSets)$
  and
  $G \acts \, X \,\in\, \Actions{G}(\ModalTopos{0})$.
  Then for $\orbisingularK \,\in\, \Singularities_{/\scalebox{.7}{$\orbisingularG$}}$
  there is a natural equivalence
  \vspace{-3mm}
  $$
    \Smooth
    \,
    \SliceMaps{\big}{\orbisingularG}
      { \orbisingularK }
      { \orbisingular (\HomotopyQuotient{X}{G}) }
    \;\;
    \simeq
    \;\;
    \Smooth
    \,
    \SliceMaps{\big}{\orbisingularG}
      { \orbisingularAny{\mathrm{im}(K)} }
      { \orbisingular (\HomotopyQuotient{X}{G}) }
    \,,
  $$

  \vspace{-2mm}
\noindent
  where
    \vspace{-2mm}
  $$
    K \twoheadrightarrow \mathrm{im}(K) \hookrightarrow G
  $$

\vspace{-1mm}
\noindent
is the image factorization of the group homomorphism
  \eqref{HomGroupoidOfSingularitiesInTermsOfGroupHomomorphisms}
  that underlies $\orbisingularK$.
\end{lemma}
\begin{proof}
The assumption that $X \,\in\, \ModalTopos{\smooth}$ is 0-truncated
implies that for
$U \,\in\, \Charts$ also $\Topos(U ,\, X) \,=\, X(U) \,\in\, \Groupoids_0
\xhookrightarrow{\;} \InfinityGroupoids$ is 0-truncated, hence a set.
This implies that the homotopy qouotient of this set of $U$-plots
(see \cite[Lem. 3.12]{SS20OrbifoldCohomology})
is a
1-groupoid and as such the disjoint union of deloopings of its isotropy groups,
which are the stabilizer subgroups $\mathrm{Stab}_G(\phi) \,\subset\, G$
of plots $U \xrightarrow{\phi} X$:
\vspace{-2mm}
\begin{equation}
  \label{DecompositionOfHomotopyQuotientOfActionOnSetofPlotsOfZeroTruncatedObject}
  X(U) \sslash G
  \;\;
  \simeq
  \;\;
  \underset{
    [x] \in \pi_0\left( \HomotopyQuotient{X(U)}{G} \right)
  }{\sqcup}
  \,
  B \mathrm{Stab}_G(x)
  \;\;
  \;
  \in
  \;
  \Groupoids_1\;.
\end{equation}

\vspace{-2mm}
\noindent
With this, we first have the following sequence of natural equivalences
for $U \,\in\, \Charts$:
\vspace{-2mm}
$$
\def\arraystretch{1.6}
  \begin{array}{lll}
       \Smooth
    \,
    \SliceMaps{\big}{\orbisingularG}
      { \orbisingularK }
      { \orbisingular (\HomotopyQuotient{X}{G}) }
    (U)
       & \;\simeq\;
    \SlicePointsMaps{\big}{\orbisingularG}
      { U \times \orbisingularK }
      { \orbisingular (\HomotopyQuotient{X}{G}) }
    &
    \proofstep{ by Lem. \ref{PlotsOfSliceMappingStackAreSliceHoms} }
    \\
    & \;\simeq\;
    \SlicePointsMaps{\big}{\orbisingularG}
      { (\orbisingular U) \times (\orbisingular B K) }
      { \orbisingular (\HomotopyQuotient{X}{G}) }
    &
    \proofstep{ by \eqref{SmoothChartsAreSmoothAndOrbisingular} }
    \\
    & \;\simeq\;
    \SlicePointsMaps{\big}{\orbisingularG}
      { \orbisingular ( U \times B K )  }
      { \orbisingular (\HomotopyQuotient{X}{G}) }
    &
    \proofstep{ by \eqref{InfinityAdjointPreservesInfinityLimits} }
    \\
    &
    \;\simeq\;
    \SlicePointsMaps{}{B G}
      { U \times B K  }
      { \HomotopyQuotient{X}{G}) }
    &
    \proofstep{ by Ex. \ref{MorphsimsBetweenOrbisingularizationsOfSmoothObjects} }
    \\
    &
    \;\simeq\;
    \InfinityGroupoids
    \big(
      { B K  }
      \,,
      { \HomotopyQuotient{X(U)}{G}) }
    \big)_{B G}
    &
    \proofstep{ by \cite[Lem. 3.12]{SS20OrbifoldCohomology} }
    \\
    &
    \;\simeq\;
    \Groupoids
    \bigg(
      { B K  }
      \,,
      {
        \underset{
          [\phi] \in \pi_0\left( \HomotopyQuotient{X(U)}{G} \right)
        }{\sqcup}
        \,
        B \mathrm{Stab}_G(x)
      }
    \bigg)_{B G}
    &
    \proofstep{ by \eqref{DecompositionOfHomotopyQuotientOfActionOnSetofPlotsOfZeroTruncatedObject} }
    \\
    &
    \;\simeq\;
    \underset{
      [\phi] \in \pi_0\left( \HomotopyQuotient{X(U)}{G} \right)
    }{\sqcup}
    \Groupoids
    \big(
      { B K  }
      \,,
      {
        \,
        B \mathrm{Stab}_G(\phi)
      }
    \big)_{B G}
    &
    \proofstep{ since $B K$ is connected }
    \\
    &
    \;\simeq\;
    \underset{
      [\phi] \in \pi_0\left( \HomotopyQuotient{X(U)}{G} \right)
    }{\sqcup}
    \Groupoids
    \big(
      { B \mathrm{im}(\phi)  }
      \,,
      {
        \,
        B \mathrm{Stab}_G(x)
      }
    \big)_{B G}
    &
    \proofstep{ by Prop. \ref{GOrbitCategoryIsReflectiveSubcategoryOfSliceOverGSingularity} }
    .
  \end{array}
$$

\vspace{-2mm}
\noindent
Now, running backwards through this chain of equivalences,
but with $K$ replaced by $\mathrm{im}(K)$ throughout, implies the claim,
by the $\infty$-Yoneda lemma (Prop. \ref{InfinityYonedaLemma}).
\end{proof}

\medskip

\noindent
{\bf Base change of proper equivariant homotopy theories.}

\begin{lemma}[Base change in equivariant homotopy theory along coverings of the equivariance group]
  \label{BaseChangeInEquivariantHomotopyTheoryAlongCoveringsOfEquivaranceGroup}
 $\,$

 \noindent  {\bf (i)} A surjective homomorphism of finite groups,
  $p : \!\!\begin{tikzcd} \widehat{G} \ar[r, ->>] &[-12pt] G \end{tikzcd}\!\!$
  \eqref{ADiscreteGroupEpimorphismInAnInfinityTopos}
  induces a reflective subcategory inclusion of orbit categories
  (Ntn. \ref{GOrbitCategory})
  $$
    \begin{tikzcd}
      \widehat{G}\Orbits
      \ar[rr, shift left=7pt, "p"]
      \ar[from=rr, shift left=7pt, hook']
      \ar[rr, phantom, "\mbox{\tiny\rm(pb)}"]
      &&
      \OrbitCategory{G}
      \,.
    \end{tikzcd}
  $$

\noindent {\bf (ii)}
  Moreover, by Kan extension (Lem. \ref{InfinityKanExtension})
  this induces,
  for $\ModalTopos{\smooth}$ any $\infty$-topos,
  an adjoint quadruple between the proper equivariant homotopy theories:
  (Def. \ref{GEquivariantAndGloballyEquivariantHomotopyTheories})
  of this form:
  \begin{equation}
    \label{BaseChangeAdjointQuadrupleOfEquivariantHomotopyTheoriesAlongEquivarianceCovering}
    \begin{tikzcd}
      \widehat{G}\InfinityGroupoids
      \ar[
        rr,
        shift left=28pt,
        "{p_!}"{description}
      ]
      \ar[
        from=rr,
        shift right=14pt,
        hook',
        "{p^\ast}"{description}
      ]
      \ar[
        rr,
        shift left=0,
        "{p_\ast}"{description}
      ]
      \ar[
        from=rr,
        shift right=-14pt,
        hook',
        "{p^!}"{description}
      ]
      \ar[rr, phantom, shift left=21pt, "\scalebox{.6}{$\bot$}"]
      \ar[rr, phantom, shift left=7pt, "\scalebox{.6}{$\bot$}"]
      \ar[rr, phantom, shift left=-7pt, "\scalebox{.6}{$\bot$}"]
      &&
      G\InfinityGroupoids \;.
    \end{tikzcd}
  \end{equation}
\end{lemma}
\begin{proof}
  The second statement follows from the first by
  Lem. \ref{KanExtensionOfAdjointPairOfInfinityFunctors}
  with Lem. \ref{LeftKanExtensionOfFullyFaithfulFunctorIsFullyFaithful}.

  A transparent way to see the adjunction in the first statement
  is to use the identification from Lem. \ref{GOrbitsAre0TruncatedObjectsOverGOrbiSingularity},
  in terms of which the adjunction is the composite
  \vspace{-2mm}
  $$
    \begin{tikzcd}[column sep=16pt]
      \widehat{G}\Orbits
      \;\simeq\;
      \big(
        (\Groupoids^{\mathrm{fin}}_{1, \geq 1})
          _{/\scalebox{.7}{$B \widehat{G}$}}
      \big)_0
      \ar[r, hook]
      &
      (\Groupoids^{\mathrm{fin}}_{1, \geq 1})
        _{/\scalebox{.7}{$B \widehat{G}$}}
      \ar[
        rr,
        shift left=5pt,
        "{ (B p)_! }"
      ]
      \ar[
        from=rr,
        shift left=5pt,
        hook',
        "{ (B p)^\ast }"
      ]
      \ar[rr, phantom, "\scalebox{.7}{$\bot$}"]
      &&
      (\Groupoids^{\mathrm{fin}}_{1, \geq 1})
        _{/\scalebox{.7}{$B G$}}
      \ar[
        rr,
        shift left=5pt,
        "{ \Truncation{0} }"
      ]
      \ar[
        from=rr,
        shift left=5pt,
        hook'
      ]
      \ar[rr, phantom, "\scalebox{.7}{$\bot$}"]
      &&
      \big(
        (\Groupoids^{\mathrm{fin}}_{1, \geq 1})
          _{/\scalebox{.7}{$B G$}}
      \big)_0
      \;\simeq\;
      \OrbitCategory{G}
    \end{tikzcd}
  $$
  of the left base change along $B p$ (Prop. \ref{BaseChange})
  with the 0-truncation reflection (Prop. \ref{nTruncation}),
  observing that the right adjoint preserves 0-truncation and hence factors.

  Under the equivalence of Lem. \ref{GOrbitsAre0TruncatedObjectsOverGOrbiSingularity},
  the right adjoint is of course given by regarding a $G$-set as a $\widehat{G}$-set
  by acting through $p$. It is immediate to see that this operation is
  fully faithful when $p$ is surjective.
\end{proof}

\newpage

\chapter{Equivariant principal $\infty$-bundles}
\label{EquivariantInfinityBundles}

In this last and main chapter we introduce the notion of equivariant principal
$\infty$-bundles, show how these subsume the traditional equivariant bundles
(from \cref{EquivariantPrincipalTopologicalBundles})
and use smooth cohesive homotopy theory
(from \cref{GeneralCohesion})
to prove a new classification theorem
for the case of truncated structure groups. Finally we
use singular-cohesive homotopy theory
(from \cref{GeneralSingularCohesion})
to show that the
equivariant moduli stacks of equivariant bundles are generally
given by the Murayama-Shimakawa construction and to
re-formulate the classification theorem in terms of
proper-equivariant homotopy theory.

\medskip

-- \cref{AsInfinityBundlesInternalToSliceOverBG}:
introducing $G$-equivariant principal $\infty$-bundles as
principal bundles internal to $\mathbf{B}G$-slices of $\infty$-toposes.

-- \cref{EquivariantLocalTrivializationIsImplies}:
recovering the classical notion of equivariant bundles and proof of the classification theorem.

-- \cref{EquivariantModuliStacks}:
equivariant moduli stacks and the proper-equivariant formulation of the classification theorem.

\medskip

\begin{notation}[Singular cohesion with cohesive charts]
 \label{SingularCohesionWithCohesiveCharts}
 In all of the following, $\Topos$
 denotes a singular-cohesive $\infty$-topos
 (Def. \ref{SingularCohesiveInfinityTopos})
 hence such that
 the cohesive sub-topos
 $\ModalTopos{\smooth} \xhookrightarrow{\;} \Topos$
 has cohesive 1-charts (Def. \ref{SmoothCharts}),
 hence a 1-site $\Charts$ of definition, such that
 $\shape U \,\simeq\, \ast$ for all $U \,\in\, \Charts$.
 For all classification results we consider the case
 $\ModalTopos{\smooth} \,=\, \SmoothInfinityGroupoids$
 from \cref{GeneralCohesion}.
\end{notation}

\begin{assumption}[Proper equivariant cohesive homotopy theory]
  \label{AssumptionOnProperEquivariantCohesiveHomotopyTheory}
  In all of the following we strengthen the
  Assumption \ref{ProperEquivariantTopology}, used in
  \cref{EquivariantTopology}, and demand that

  -- {\it equivariance groups $G$ are discrete groups (not necessariliy finite)};

  \noindent
  and for all classification results discussed in the following we also
  assume that

  -- {\it domain spaces $\TopologicalSpace$ carry the structure of smooth manifolds}.

  \medskip

  \noindent
  The first of these assumptions
  \vspace{-2mm}
  \begin{equation}
    \label{AssumptionOfDiscreteEquivarianceGroup}
    G
    \,\in\,
    \Groups(\Sets)
    \xhookrightarrow{\;\Discrete \;}
    \Groups(\Topos)
  \end{equation}
  \vspace{-2mm}
  \noindent
  implies (by Rem \ref{TheAxiomsOnTheShapeModality})
  that also the delooping \eqref{LoopingAndDelooping} of the equivariance group
  is discrete and hence pure shape:
  $$
    \begin{tikzcd}
    \flat \mathbf{B}G
    \ar[
      r,
      "{
        \scalebox{.85}{$
          \epsilon_{{}_{\mathbf{B}G}}^{\scalebox{.7}{$\flat$}}
        $}
      }",
      "\sim"{swap}
    ]
    &
    \mathbf{B}G
    \ar[
      r,
      "{
        \scalebox{.85}{$
          \eta_{{}_{\mathbf{B}G}}^{\scalebox{.7}{$\shape$}}
        $}
      }",
      "\sim"{swap}
    ]
    &
    \shape \mathbf{B}G
    \underset{
      \mathclap{
      \raisebox{-3pt}{
        \tiny
        \rm
        Def .\ref{ClassifyingSpaceOfPrincipalInfinityBundles}
      }
      }
    }{
      \; := \;
      B G
    }
    \,,
    \end{tikzcd}
  $$
  and the same holds
  (by $\flat \circ \orbisingular \,\simeq\, \orbisingular \circ \flat$, Prop. \ref{SomeSingularModalitiesCommuteWithSomeCohesiveModalities})
  for its orbi-singularization (Ex. \ref{OrbiSingularityIsOrbiSingularizationOfHomotopyQuotient}):
  $$
    \begin{tikzcd}
    \flat \orbisingularG
    \ar[
      r,
      "{
        \scalebox{.85}{$
          \epsilon_{{}_{\orbisingularG}}^{\scalebox{.7}{$\flat$}}
        $}
      }",
      "\sim"{swap}
    ]
    &
    \orbisingularG
    \ar[
      r,
      "{
        \scalebox{.85}{$
          \eta_{{}_{\orbisingularG}}^{\scalebox{.7}{$\shape$}}
        $}
      }",
      "\sim"{swap}
    ]
    &
    \shape \orbisingularG
    \,.
    \end{tikzcd}
  $$
\end{assumption}

\begin{remark}[Discrete equivariance $\infty$-groups in cohesive global homotopy theory]
  \label{NeedForDiscreteEquivarianceGroups}
  Assumption \ref{AssumptionOnProperEquivariantCohesiveHomotopyTheory}
  \eqref{AssumptionOfDiscreteEquivarianceGroup} is used in our proofs
  (in Thm. \ref{OrbiSmoothOkaPrinciple} below, just as in Prop. \ref{ShapeOfMappingStackOutOfOrbiSingularityIsMappingStackIntoShape} above)
  to guarantee, via Prop. \ref{ShapeFunctorPreservesHomotopyFibersOverDiscreteObjects},
  that $\shape$ preserves homotopy fibers over $\orbisingularG$.

\noindent {\bf (i)}   At this point it is worth recalling
  (\cite[Rem. 3.64]{SS20OrbifoldCohomology}) that in
  {\it cohesive} global homotopy theory
  according to Def. \ref{SingularCohesiveInfinityTopos}
  it is not useful to promote the
  2-category $\Singularities$
  (Ntn. \ref{Singularities}) from discrete to compact Lie groups
  with shapes of mapping stacks of their delooping between them,
  as commonly done in non-cohesive global homotopy theory.
  The reason is, conceptually, that this secretly introduces a notion
  of cohesion into the site, which does not properly interplay with the cohesion
  that is seen inside the cohesive $\infty$-topos over this site. For example,
  the crucial relations \eqref{SmoothAspectOfOrbisingularity}
  and \eqref{OrbiSingularityIsOrbiSingularizationOfDelooping} would fail, in general.

\noindent {\bf (ii)}   On the other hand, in contrast to traditional equivariant homotopy theory
(and aside from the particular classification Theorems \ref{BorelClassificationOfEquivariantBundlesForResolvableSingularitiesAndEquivariantStructure}
and
\ref{ProperClassificationOfEquivariantBundlesForResolvableSingularitiesAndEquivariantStructure}
below),
cohesive global homotopy theory works just as well for {\it higher} equivariance
groups in
$
  \Groups(\InfinityGroupoids)
  \xrightarrow{\;{\Groups(\Discrete)}\;}
  \Groups(\Topos)
  \,,
$
as long as they are geometrically discrete, hence for plain equivariance $\infty$-groups.
But if $G_i$ are topological 1-groups, then their shapes $\shape G_i$ are such
geometrically discrete $\infty$-groups, and there is the canonical comparison morphism
\eqref{ComparisonMorphismFromShapeOfMappingStackToMappingSpaceOfShapes}
between their usual hom-$\infty$-groupoids:
\begin{equation}
  \label{ComparisonMorphismForNonDiscreteEquivarianceGroups}
  \shape
  \,
  \Maps{}
    { \mathbf{B}G_1 }
    { \mathbf{B}G_2 }
  \xrightarrow{\;\widetilde{\shape \mathrm{ev}}\;}
  \Maps{}
    { \mathbf{B} \shape G_1 }
    { \mathbf{B} \shape G_2 }
  \,.
\end{equation}
\noindent {\bf (iii)}  The cohesive global homotopy theory
as set up here captures all those aspects of non-discrete equivariance groups
that are still reflected on the right of \eqref{ComparisonMorphismForNonDiscreteEquivarianceGroups}.
For these, no further condition on the topologies is required, for example the
$G_i$ could even be non-compact topological groups such as loop groups.
\end{remark}

\section{As bundles internal to $G$-$\infty$-actions}
\label{AsInfinityBundlesInternalToSliceOverBG}

In higher analogy with the discussion in \eqref{PrincipalBundlesInternalToTopologicalGActions}
we introduce equivariant principal $\infty$-bundles as
principal bundles internal to $\infty$-toposes of $G$-$\infty$-actions
(Def. \ref{GEquivariantGammaPrincipalBundles} below).
We show that ordinary topological equivariant bundles
(according to \cref{EquivariantPrincipalTopologicalBundles})
faithfully embed into this higher geometric theory,
and we use this to prove their
classification theorem for truncated structure groups
(Thm. \ref{BorelClassificationOfEquivariantBundlesForResolvableSingularitiesAndEquivariantStructure}).

\medskip
In order to warm up and to establish some preliminariries, we first
present an analogous re-proof of the classical Milgram-classification
of ordinary topological principal bundles
(Thm. \ref{OrdinaryPrincipalBundlesAmongPrincipalInfinityBundles},
Thm. \ref{ClassificationOfPrincipalBundlesAmongPrincipalInfinityBundles}):

\medskip

-- \cref{SmoothPrincipalInfinityBundles}: Smooth principal $\infty$-bundles.

-- \cref{SmoothEquivariantPrincipalInfinityBundles}: Smooth equivariant principal $\infty$-bundles.

\medskip

\subsection{Smooth principal $\infty$-bundles}
\label{SmoothPrincipalInfinityBundles}

\noindent
{\bf Principal $\infty$-bundles.}
We begin by recalling and developing some facts about
plain (i.e. not equivariant) principal $\infty$-bundles
(\cite{NSS12a}\cite{NSS12b})
with focus on their incarnation in $\SmoothInfinityGroupoids$
(Ntn. \ref{SmoothInfinityGroupoids}).

\medskip
\noindent
{\bf Recovering principal bundles among principal $\infty$-bundles.}
We discuss first how ordinary principal bundles embed into the
theory of principal $\infty$-bundles in the sense of
Prop. \ref{GroupsActionsAndFiberBundles}.
The following
Prop. \ref{OrdinaryPrincipalBundlesAmongPrincipalInfinityBundles}
is the blueprint for the analogous embedding of
equivariant principal bundles into the theory of
equivariant principal $\infty$-bundles which we prove
in Thm. \ref{BorelClassificationOfEquivariantBundlesForResolvableSingularitiesAndEquivariantStructure} below.

\begin{proposition}[Universal $\Gamma$-principal bundle over the $\Gamma$-moduli stack]
  \label{DeloopingGroupoidsOfDTopologicalGroupsAreLocallyFibrant}
  For a diffeological group
  $$\Gamma \,\in\,
    \Groups\left(\DiffeologicalSpaces\right)
    \xhookrightarrow{\;}
    \Groups\left(\Presheaves(\CartesianSpaces)\right),
  $$
  such as a D-topological group or a Lie group

  \vspace{1mm}
  \noindent
 {\bf (i)}  the simplicial nerve (Ntn. \ref{NerveOfTopologicalGroupoids})
  of the delooping groupoid of $\Gamma$
  (Ex. \ref{TopologicalDeloopingGroupoid})
  \vspace{-2mm}
  $$\mathbf{B}\Gamma
     \,=\,
     \DeloopingGroupoid{\Gamma}
    \,\in\,
    \Groupoids(\kTopologicalSpaces)
     \xhookrightarrow{\;\;
       \SimplicialNerve
       \Groupoids(\ContinuousDiffeology)
     \;\;}
     \SimplicialPresheaves(\CartesianSpaces)
    $$

  \vspace{-2mm}
  \noindent
  is fibrant in
  the local projective model structure
  of simplicial presheaves over the site of Cartesian spaces
  (Ntn. \ref{SmoothInfinityGroupoids}),
  and hence a fibrant representative of the delooping of $\Gamma$
  in the $\infty$-topos according to Prop. \ref{GroupsActionsAndFiberBundles}:
  \vspace{-2mm}
  $$
    \mathbf{B}\Gamma
    \;\in\;
    \big(
      \SimplicialPresheaves(\CartesianSpaces)_{\projloc}
    \big)_{\mathrm{fib}}
    \xrightarrow{\;\;
      \SimplicialLocalization{\LocalWeakEquivalences}
    \;\;}
    \SmoothInfinityGroupoids
    \,.
  $$

 \vspace{-2mm}
 \noindent
 {\bf (ii)} Moreover,
  pullback of the universal $\Gamma$-principal groupoid
  \eqref{TheUniversalPrincipalTopologicalGroupoid}
  \vspace{-2mm}
  $$\mathbf{E}\Gamma
     \,=\,
     \ActionGroupoid{\Gamma}{\Gamma}
    \,\in\,
    \Groupoids(\kTopologicalSpaces)
     \xhookrightarrow{\;\;
       \SimplicialNerve
       \Groupoids(\ContinuousDiffeology)
    \;\; }
     \SimplicialPresheaves(\CartesianSpaces)
    $$

  \vspace{-2mm}
  \noindent
  along morphisms of simplicial presheaves
  $X \xrightarrow{\;c\;}\mathbf{B}\Gamma$
  represents their homotopy fibers in $\SmoothInfinityGroupoids$:
  \vspace{-2mm}
  $$
    \begin{tikzcd}[row sep=4pt, column sep=15pt]
     X \underset{\mathbf{B}\Gamma}{\times} \mathbf{E}\Gamma
     \ar[r]
     \ar[d]
     \ar[dr, phantom, "\mbox{\tiny\rm(pb)}"{pos=.2}]
     &
     \mathbf{E}\Gamma
     \ar[d]
     \\
     X
     \ar[r]
     &
     \mathbf{B}\Gamma
    \end{tikzcd}
    \;\;
    \in
    \;
    \SimplicialPresheaves(\CartesianSpaces)
    \;\;\;\;\;\;
    \Rightarrow
    \;\;\;\;\;\;
    \Localization{\LocalWeakEquivalences}
    \big(
      X \underset{\mathbf{B}\Gamma}{\times} \mathbf{E}\Gamma
    \big)
    \;\;\simeq\;\;
    \Localization{\LocalWeakEquivalences}
    (X)
      \underset{
        \mathbf{B}\Gamma
      }{\times}
    \ast
    \;\;\;
    \in
    \;
    \SmoothInfinityGroupoids
    \,.
  $$
\end{proposition}
\begin{proof}
  The local fibrancy is Lemma \ref{LocalFibrancyOfDeloopingGroupoidsOverCartesianSpaces}
  and the delooping property is Lemma \ref{InfinityGroupsPresentedByPresheavesOfSimplicialGroups}in view of Ex. \ref{StandardModelOfUniversalSimplicialPrincipalComplex}.
  The second statement (as well as this delooping property) follows by
  Lem. \ref{ComputingHomotopyPullbacksOfInfinityStacks}
  since
  \vspace{-2mm}
  $$
    \begin{tikzcd}[column sep=40pt]
      \ast
      \ar[r, "\in \, \ProjectiveWeakEquivalences"{swap}]
      &
      \mathbf{E}\Gamma
      \ar[r, "\in \, \ProjectiveFibrations"{swap}]
      &
      \mathbf{B}\Gamma
    \end{tikzcd}
  $$

    \vspace{-2mm}
 \noindent
  is evidently a fibration resolution of the point inclusion in the
  global projective model structure.
\end{proof}

\begin{theorem}[Topological/smooth principal bundles embed among smooth principal $\infty$-bundles]
  \label{OrdinaryPrincipalBundlesAmongPrincipalInfinityBundles}
  Let

  -- $\Gamma \,\in\, \Groups(\DHausdorffSpaces)
    \xhookrightarrow{\Groups(\ContinuousDiffeology)}
    \Groups(\DiffeologicalSpaces)
    \xhookrightarrow{\;}
    \Groups\left(\Presheaves(\CartesianSpaces)\right),
    $

  \noindent
  or

  -- $\Gamma \,\in\, \Groups(\SmoothManifolds)
    \xhookrightarrow{\phantom{\Groups(\ContinuousDiffeology)}}
    \Groups(\DiffeologicalSpaces)
    \xhookrightarrow{\;}
    \Groups\left(\Presheaves(\CartesianSpaces)\right).
    $

  \vspace{1mm}
  \noindent
  Then the groupoid of traditional
  $\Gamma$-principal fiber bundles over $\SmoothManifold$
  (Rem. \ref{AssumptionOfLocalTrivializability}),
  respectively topological or smooth,
  is naturally equivalent to the groupoid of
  $\Gamma$-principal $\infty$-bundles
  over $\SmoothManifolds$
  internal to $\SmoothInfinityGroupoids$,
  in the sense of Prop. \ref{GroupsActionsAndFiberBundles}:

  \vspace{-2mm}
  \begin{equation}
    \label{DiffeologicalPrincipalBundlesOnSmoothManifoldEquivalentlySeenInSmoothInftyGroupoids}
    \begin{tikzcd}
      \underset{
        \mathclap{
        \raisebox{-3pt}{
          \tiny
          \color{darkblue}
          \bf
          ordinary topological/smooth-principal bundles
        }
        }
      }{
        \PrincipalFiberBundles{\Gamma}(\DiffeologicalSpaces)_X
      }
      \;\;\;\simeq\;\;\;
      \underset{
        \mathclap{
        \raisebox{-3pt}{
          \tiny
          \color{darkblue}
          \bf
          among smooth principal $\infty$-bundles
        }
        }
      }{
        \PrincipalBundles{\Gamma}(\SmoothInfinityGroupoids)_X
      }
      \;\;\;
      \in
      \;
      \Groupoids
      \xhookrightarrow{\;}
      \InfinityGroupoids
      \,.
    \end{tikzcd}
  \end{equation}
\end{theorem}
\begin{proof}
  Since $\SmoothManifold$ is assumed to be a smooth manifold, it
  admits a differentiably good open cover
  \eqref{GoodOpenCover} $\big\{ U_i \,\simeq\, \mathbb{R}^{\mathrm{dim}(\SmoothManifold)} \xhookrightarrow{\;} \SmoothManifold \big\}$.

  \noindent
  Now we consider three consecutive equivalences of groupoids:

  \noindent
  {\bf (i)}
  Since the cover is good,
  every $\Gamma$-principal bundle has a trivialization over the cover,
  which gives a diffeological
  {\v C}ech 1-cocycle with coefficients in $\Gamma$:
  \vspace{-2mm}
  \begin{equation}
    \label{CechCocycleForPrincipalBundle}
    \def\arraystretch{1}
    \begin{array}{lll}
      \mathrm{P}
      \\
      \downarrow
      \\
      \TopologicalSpace
    \end{array}
    \;\;\;\;
    \rightsquigarrow
    \;\;\;\;
    \big\{
      \mathrm{U}_i \cap \mathrm{U}_j
        \,\xrightarrow[\mathrm{dfflg}]{\gamma_{i j}}
        \,
      \Gamma
    \big\}_{i,j \in I}
    \quad
    \mbox{s.t}
    \qquad
    \underset
      { i \in I}
      { \forall }
      \left(
        \gamma_{i i }
        =
        \mathrm{const}_\NeutralElement
      \right)
    \;\;
      \mbox{and}
    \;\;
    \underset{i,j,k \in I}{\forall}
      \left(
        \gamma_{i j} \cdot \gamma_{j k}
        \,=\,
        \gamma_{i k}
      \right)
      \,,
  \end{equation}

  \vspace{-3mm}
  \noindent
  and that
  morphisms of $\Gamma$-principal bundles over $\TopologicalSpace$
  bijectively correspind to
  diffeological {\v C}ech coboundaries between these cocycles:
  \vspace{-.2cm}
  \begin{equation}
    \mathrm{P} \xrightarrow{\;h\;} \mathrm{P}^\prime
    \qquad
    \xleftrightarrow{\;\;\;}
    \qquad
    \big\{
      \begin{tikzcd}
        \TopologicalPatch_i
          \xrightarrow[\mathrm{dfflg}]{ h_i }
        \Gamma
      \end{tikzcd}
      \,\big\vert\,
      \underset{i,j \in I}{\forall}
      \;
      h_i \cdot \gamma^{\; \prime}_{i j} \,=\, \gamma_{i j} \cdot h_j
    \big\}
    \,,
  \end{equation}

  \vspace{-2mm}
  \noindent
  so that this construction constitutes a
  natural equivalence
  between the groupoid of principal fiber bundles
  and the action groupoid of {\v C}ech coboundaries acting on
  {\v C}ech cocycles:
  \vspace{-1mm}
  \begin{equation}
    \label{DiffeologicalPrincipalBundlesEquivalentToCechGroupoid}
    \PrincipalFiberBundles{\Gamma}(\DiffeologicalSpaces)_{\SmoothManifold}
    \;\;
    \simeq
    \;\;
    \ActionGroupoid
      { \hat Z^1 (\widehat{\SmoothManifold};\, \Gamma) }
      { \hat B^0(\widehat{\SmoothManifold};\,\Gamma) }
    \,.
  \end{equation}

\vspace{-1mm}
  \noindent
  {\bf (ii)}
  By direct inspection (see Rem. \ref{EquivariantCechCocycles})
  and using the defining fully-faithful embedding
  $\DiffeologicalSpaces \xhookrightarrow{\;} \Presheaves(\CartesianSpaces)$,
  the nerve of this latter groupoid  is manifestly
  isomorphic to the hom-complex \eqref{SimplicialHomComplex}
  of simplicial presheaves from the
  {\v C}ech nerve of the open cover
  \vspace{-2mm}
  \begin{equation}
    \label{CechNerveOfOpenCover}
    \SimplicialNerve
    \big(
      \widehat{\SmoothManifold}
        \times_{\SmoothManifold}
      \widehat{\SmoothManifold}
      \rightrightarrows
      \widehat{\SmoothManifold}
    \big)
    \;\;
      \simeq
    \;\;
    \widehat{\SmoothManifold}^{\times_{\SmoothManifold} \bullet}
    \;\;
    \in
    \;\;
    \SimplicialPresheaves(\CartesianSpaces)
  \end{equation}

  \vspace{-1mm}
  \noindent
  to the nerve of the delooping groupoid of $\Gamma$:
  \vspace{-1mm}
  \begin{equation}
    \label{GroupoidOfCechCocyclesAsHomComplexOfSimplicialPresheaves}
    \SimplicialNerve
    \ActionGroupoid
      { \hat Z^1 (\widehat{\SmoothManifold};\, \Gamma) }
      { \hat B^0(\widehat{\SmoothManifold};\,\Gamma) }
    \;\;
    \simeq
    \;\;
    \SimplicialPresheaves(\CartesianSpaces)
    \big(
      \widehat{\SmoothManifold}^{\times_{\SmoothManifold}^\bullet},
      \,
      \SimplicialNerve
      \DeloopingGroupoid{\Gamma}
    \!\big)
    \,.
  \end{equation}

 \vspace{-2mm}
  \noindent
  {\bf (iii)}
  This has the homotopy type of the correct hom-$\infty$-groupoid
  (Def. \ref{HomInfinityGroupoid})
  \vspace{-1mm}
  \begin{equation}
    \label{SimplicialModelForHomSpaceFromManifoldToDeloopedDiffeologicalGroup}
    \SimplicialPresheaves(\CartesianSpaces)
    \big(
      \widehat{\SmoothManifold}^{\times_{\SmoothManifold}^\bullet},
      \,
      \SimplicialNerve
      \DeloopingGroupoid{\Gamma}
    \big)
    \;\simeq\;
    \SmoothInfinityGroupoids
    (
      \TopologicalSpace,
      \,
      \mathbf{B}G
    )
    \;\;\;
    \in
    \;
    \InfinityGroupoids
    \,,
  \end{equation}

  \vspace{-3mm}
  \noindent
  because:

  \vspace{-.3cm}
  \begin{itemize}
  \setlength\itemsep{-1pt}
  \item
  $\SimplicialNerve\DeloopingGroupoid{\Gamma}$,
  which represents \eqref{PresentingObjectInModelCategory}
  the delooping \eqref{GroupsActionsAndFiberBundles}
  of
  $\Gamma$ (by Prop. \ref{InfinityGroupsPresentedByPresheavesOfSimplicialGroups}),
  is fibrant in $\SimplicialPresheaves(\CartesianSpaces)_{\projloc}$,
  by Prop. \ref{DeloopingGroupoidsOfDTopologicalGroupsAreLocallyFibrant},

  \item
  $\widehat{\SmoothManifold}^{\times^\bullet_{\SmoothManifold}}$ is a
  local projective cofibrant resolution of
  $\SmoothManifold$, by Example \ref{GoodOpenCoversAreProjectivelyCofibrantResolutionsOfSmoothManiolds}.

\end{itemize}
\vspace{-2mm}

  \noindent
  Therefore, \eqref{SimplicialModelForHomSpaceFromManifoldToDeloopedDiffeologicalGroup}
  follows by Prop. \ref{HomInfinityGroupoidFromCofibrantDomainAndFibrantCodomain}.
  In conclusion, the composite of the
  equivalences
  \eqref{DiffeologicalPrincipalBundlesEquivalentToCechGroupoid},
  \eqref{GroupoidOfCechCocyclesAsHomComplexOfSimplicialPresheaves},
  and
  \eqref{SimplicialModelForHomSpaceFromManifoldToDeloopedDiffeologicalGroup}
  yields the desired equivalence
  \eqref{DiffeologicalPrincipalBundlesOnSmoothManifoldEquivalentlySeenInSmoothInftyGroupoids}.
\end{proof}

\begin{remark}[Local triviality of principal bundles internal to $\infty$-topos is implied]
  \label{LocalTrivialityOfPrincipalBundlesInternalToInfinityToposIsImplied}
  Thm. \ref{OrdinaryPrincipalBundlesAmongPrincipalInfinityBundles}
  gives an equivalence to {\it locally trivial}
  topological principal bundles
  (Rem. \ref{AssumptionOfLocalTrivializability}),
  even though the Definition \ref{PrincipalInfinityBundles}
  of smooth principal $\infty$-bundles,
  as principal bundles internal to the $\infty$-topos
  $\SmoothInfinityGroupoids$, does not
  explicitly state a local triviality clause.
  The equivariant generalization of this phenomenon is the content of
  Thm. \ref{BorelClassificationOfEquivariantBundlesForResolvableSingularitiesAndEquivariantStructure}
  below; see
  also the discussion at the beginning of \cref{EquivariantLocalTrivializationIsImplies}.
\end{remark}

\begin{remark}[Universality of the universal principal bundle over the moduli stack]
  \label{UniversalPrincipalBundleOverTheModuliStack}
  $\,$

  \vspace{0mm}
\noindent{\bf (i)}   While the above proof of Thm. \ref{OrdinaryPrincipalBundlesAmongPrincipalInfinityBundles}
  handles principal bundles entirely in terms of their {\v C}ech cocycles,
  it entails a systematic way of reconstructing the actual bundles from these
  {\v C}ech cocycles, as follows:
  By Prop. \ref{DeloopingGroupoidsOfDTopologicalGroupsAreLocallyFibrant}
  and Prop. \ref{GroupsActionsAndFiberBundles},
  the following pullback of simplicial presheaves
  \vspace{-2mm}
  $$
    \begin{tikzcd}[row sep=10pt]
      \TopologicalPrincipalBundle
      \ar[from=rr, ->>, "\in \LocalWeakEquivalences"{yshift=-1pt}]
      \ar[d]
      &&[+10pt]
      \big(
        \Gamma
          \times
        \widehat{\SmoothManifold}
          \times_{\SmoothManifold}
        \widehat{\SmoothManifold}
        \rightrightarrows
        \Gamma \times \widehat{\SmoothManifold}
      \big)
      \ar[rr]
      \ar[d]
      \ar[drr, phantom, "\mbox{\tiny\rm(pb)}"]
      &&
      \ActionGroupoid
        { \Gamma }
        { \Gamma }
      \ar[d]
      \ar[r, phantom, "="]
      &[-10pt]
      \mathbf{E}G
      \ar[d]
      \\
      \SmoothManifold
      \ar[from=rr, ->>, "\in \LocalWeakEquivalences"]
      &&
      \big(
        \widehat{\SmoothManifold}
          \times_{\SmoothManifold}
        \widehat{\SmoothManifold}
        \rightrightarrows
        \widehat{\SmoothManifold}
      \big)
      \ar[rr, "c"{swap}]
      &&
      \DeloopingGroupoid{ \Gamma }
      \ar[r, phantom, "\simeq"]
      &
      \mathbf{B}\Gamma
    \end{tikzcd}
  $$
  computes the $\Gamma$-principal $\infty$-bundle classified
  by a {\v C}ech cocycle $c$,
  according to Thm. \ref{DeloopingGroupoidsAreModuliInfinityStacksForPrincipalInfinityBundles}.
  But inspection of the pullback diagram shows that this
  bundle $\TopologicalPrincipalBundle \xrightarrow{\;} \SmoothManifold$
  is precisely the principal bundle reconstructed from a {\v C}ech 1-cocycle
  in the traditional way.

  \vspace{0mm}
\noindent {\bf (ii)}     Therefore, the {\it groupoid} $\mathbf{E}\Gamma$
  from \eqref{TheUniversalPrincipalTopologicalGroupoid}
  (whose topological realization
  $E G \,=\, \TopologicalRealization{}{\mathbf{E}\Gamma}$
  \eqref{UniversalPrincipalBundleAsTopologicalRealizationOfUniversalPrincipalGroupoid}
  is the universal $\Gamma$-principal bundle over the classifying space $B G$)
  plays the role of the universal $\Gamma$-principal bundle over the
  {\it moduli stack} $\mathbf{B}\Gamma$.
\end{remark}

\medskip

\noindent
{\bf Concordance of principal $\infty$-bundles.}
The notion of concordance of principal bundles
(Def. \ref{ConcordanceOfEquivariantPrincipalBundles})
has an evident generalization to principal $\infty$-bundles
(and, of course, yet more generally to any contravariant structure over
spaces with cylinder objects).

\begin{definition}[Concordance of smooth principal $\infty$-bundles]
  \label{ConcordanceOfSmoothPrincipalBundles}
  Let $\mathcal{G} \,\in\, \Groups(\SmoothInfinityGroupoids)$.

  \noindent
  {\bf (i)} We say that a {\it concordance} between
  $P_1, P_2 \,\in\, \PrincipalBundles{\mathcal{G}}(\SmoothInfinityGroupoids)_X$
  is a principal $\infty$-bundle on the cylinder $X \times \mathbb{R}$
  \vspace{-2mm}
  $$
    \widehat P
    \;\in\;
    \PrincipalBundles{\mathcal{G}}(\SmoothInfinityGroupoids)_{X \times \mathbb{R}}
  $$

  \vspace{-1mm}
  \noindent
  such that over the two endpoints it restricts to the two bundles,
  respectively, up to equivalence:
  \vspace{-1mm}
  $$
    {\widehat P}\vert_{ X \times \{0\}}
    \;\simeq\;
    P_0
    \;\;\;\;\;\;
    \mbox{and}
    \;\;\;\;\;\;
    {\widehat P}\vert_{ X \times \{1\}}
    \;\simeq\;
    P_1,
    \qquad
    \;
   P_1,
    P_2 \,\in\,
    \PrincipalBundles{\mathcal{G}}(\SmoothInfinityGroupoids)_X
    \,.
  $$

  \vspace{-2mm}
  \noindent {\bf (ii)}   We write
  $$
    \left(
      \PrincipalBundles{\mathcal{G}}_X
    \right)_{/\sim_{\mathrm{conc}}}
    \;\;\;
    \in
    \;
    \Sets
  $$
  for the set of concordance classes
  of $\mathcal{G}$-principal bundles
  in $\SmoothInfinityGroupoids$.
\end{definition}

The equivalence between isomorphism classes and concordance classes
of topological principal bundles (Thm.  \ref{ConcordanceClassesOfTopologicalPrincipalBundles})
still holds after including them among
smooth principal $\infty$-bundles:

\begin{proposition}[Concordant topological principal bundles are isomorphic as smooth $\infty$-bundles]
  Given

  -- $
    \SmoothManifold
      \,\in\,
    \SmoothManifolds
      \xhookrightarrow{\;\;}
    \SmoothInfinityGroupoids
  $\,,

  \vspace{1mm}
  -- $
    \Gamma
      \,\in\,
    \Groups(\DHausdorffSpaces)
      \xhookrightarrow{\; \scalebox{0.6}{$\Groups(\ContinuousDiffeology)$}\;}
    \Groups(\SmoothInfinityGroupoids)
  $,

  \vspace{1mm}
  \noindent
  the canonical comparison morphism
  \eqref{CanonicalProjectionFromPrincipalInfinityBundlesToTheirConcordances}
  from isomorphism classes
  to concordance classes (Def. \ref{ConcordanceOfSmoothPrincipalBundles})
  of $\Gamma$-principal $\infty$-bundles over $\SmoothManifold$
  is a bijection
  \vspace{-2mm}
  $$
    \begin{tikzcd}
      \IsomorphismClasses
      {
        \PrincipalBundles{\Gamma}(\SmoothInfinityGroupoids)_{\TopologicalSpace}
      }
      \ar[r, "\sim"]
      &
      \ConcordanceClasses
      {
        \PrincipalBundles{\Gamma}(\SmoothInfinityGroupoids)_{\TopologicalSpace}
      }
      \mathrlap{\,.}
    \end{tikzcd}
  $$
\end{proposition}
\begin{proof}
This is the following composite of natural bijections:
\vspace{-.1cm}
$$
  \def\arraystretch{2}
  \begin{array}{lll}
    &
    \IsomorphismClasses
    {
      \PrincipalBundles{\Gamma}(\SmoothInfinityGroupoids)_{\TopologicalSpace}
    }
    \\
    & \;\simeq\;
    \IsomorphismClasses
    {
      \PrincipalFiberBundles{\Gamma}(\DiffeologicalSpaces)_{\TopologicalSpace}
    }
    &
    \proofstep{ by Prop. \ref{OrdinaryPrincipalBundlesAmongPrincipalInfinityBundles} }
    \\
    & \;\simeq\;
    \IsomorphismClasses
    {
      \PrincipalFiberBundles{\Gamma}(\DTopologicalSpaces)_{\TopologicalSpace}
    }
    &
    \proofstep{ by Prop. \ref{CategoryOfDeltaGeneratedTopologicalSpaces} }
    \\
    & \;\simeq\;
    \ConcordanceClasses
    {
      \PrincipalFiberBundles{\Gamma}(\DTopologicalSpaces)_{\TopologicalSpace}
    }
    &
    \proofstep{ by Thm. \ref{ConcordanceClassesOfTopologicalPrincipalBundles} }
    \\
    & \;\simeq\;
    \ConcordanceClasses
    {
      \PrincipalBundles{\Gamma}(\SmoothInfinityGroupoids)_{\TopologicalSpace}
    }
    &
    \proofstep{
      by Props.
      \ref{CategoryOfDeltaGeneratedTopologicalSpaces},
      \ref{OrdinaryPrincipalBundlesAmongPrincipalInfinityBundles}
      .
    }
  \end{array}
$$

\vspace{-8mm}
\end{proof}

\begin{proposition}[Shape of mapping stack into moduli stack is $\infty$-groupoid of concordances]
  \label{ShapeOfMappingStackIntoModuliStackISConcordances}
  Let $\Topos \,=\, \SmoothInfinityGroupoids$.
  Then, for any $\mathcal{G} \,\in\, \Groups(\Topos)$
  and $X \,\in\, \SmoothManifolds \hookrightarrow{} \Topos$, the 0-truncation of the
  shape of the mapping stack from $X$ to the
  $\mathcal{G}$-moduli stack is in natural bijection to the set
  of concordance classes (Def. \ref{ConcordanceOfSmoothPrincipalBundles})
  of $\mathcal{G}$-principal $\infty$-bundles on $X$:
  $$
    \tau_0
    \,
    \Shape
    \,
    \Maps{}
      { X }
      { \mathbf{B}\mathcal{G} }
    \;\simeq\;
    (\PrincipalBundles{\mathcal{G}}(\SmoothInfinityGroupoids)_X)_{/\sim_{\mathrm{conc}}}
    \,.
  $$
\end{proposition}
\begin{proof}
In the following sequence of natural bijections, the key step
is the ``smooth Oka principle'' (Thm. \ref{SmoothOkaPrinciple})
which allows to take the shape-operation
out of the mapping stack construction:
\vspace{-3mm}
$$
  \def\arraystretch{2}
  \begin{array}{lll}
    \tau_0
    \,
    \Shape
    \,
    \Maps{}
      { X }
      { \mathbf{B}\mathcal{G} }
    & \;\simeq\;
    \tau_0
    \,
    \SingularSimplicialComplex
    \,
    \Maps{}
      { X }
      { \mathbf{B}\mathcal{G} }
    &
    \mbox{\small by Prop. \ref{SmoothShapeGivenBySmoothPathInfinityGroupoid} }
    \\
    & \;\simeq\;
    \tau_0
    \,
    \underset{\longrightarrow}{\mathrm{lim}}
    \big(
      \Maps{}
        { X }
        { \mathbf{B}\mathcal{G} }
      (\SmoothSimplex{\bullet})
    \big)
    &
    \mbox{\small by Def. \ref{SmoothPathInfinityGroupoid}}
    \\
    & \;\simeq\;
    \tau_0
    \,
    \underset{\longrightarrow}{\mathrm{lim}}
    \Big(
    \Topos
    \big(
      X
      \times
      \SmoothSimplex{\bullet}
      ,\,
      \mathbf{B}\mathcal{G}
    \big)
   \!\! \Big)
    &
    \mbox{\small by Prop. \ref{MappingStacks}}
    \\
    & \;\simeq\;
    \tau_0
    \,
    \underset{\longrightarrow}{\mathrm{lim}}
    \left(
      \PrincipalBundles{\mathcal{G}}_{X \times \SmoothSimplex{\bullet} }
    \right)
    &
    \mbox{\small by Prop. \ref{GroupsActionsAndFiberBundles}}
    \\
    & \;\simeq\;
    \underset{\longrightarrow}{\mathrm{lim}}
    \,
    \tau_0
    \big(
    \PrincipalBundles{\mathcal{G}}_{X \times \SmoothSimplex{\bullet} }
    \big)
    &
    \mbox{\small by Prop. \ref{nTruncation}}
    \\
    & \;\simeq\;
    \big(
      \tau_0
      \left(
        \PrincipalBundles{\mathcal{G}}_{ X }
      \right)
    \big)
    \big/
    \left(
      \tau_0
      \left(
      \PrincipalBundles{\mathcal{G}}_{X \times \SmoothSimplex{1} }
      \right)
    \!\!
    \right)
    &
    \mbox{\small by Lem. \ref{ColimitsOverSimplicialDiagramsOfZeroTruncatedObjectsAreCoequalizers}}
    \\
    & \;\simeq\;
    \left(
      \PrincipalBundles{\mathcal{G}}_X
    \right)_{\sim_{\mathrm{conc}}}
    &
    \mbox{\small by Def. \ref{ConcordanceOfSmoothPrincipalBundles}}.
  \end{array}
$$

\vspace{-8mm}
\end{proof}

\begin{remark}[The $\infty$-groupoid of concordances]
\label{TheInfinityGroupoidOfConcordances}
Without applying 0-truncation, the argument in Prop. \ref{ShapeOfMappingStackIntoModuliStackISConcordances}
shows that the full shape of the mapping stack into a moduli stack is
the $\infty$-groupoid of concordances and higher order concordances
between principal $\infty$-bundles.
Here the points-to-pieces transformation $\flat \to \shape$
\eqref{PointsToPiecesTransform}
on these
mapping stacks regards higher morphisms of principal $\infty$-bundles
inside higher concordances:
\vspace{-3mm}
\begin{equation}
  \label{CanonicalProjectionFromPrincipalInfinityBundlesToTheirConcordances}
  \hspace{-7mm}
    \begin{tikzcd}[column sep=10pt]
    \PrincipalBundles{\mathcal{G}}_X
    \;\simeq
    \PointsMaps{}
      { X }
      { \mathbf{B}\mathcal{G} }
    \;\simeq\;
    \overset{
      \mathclap{
      \raisebox{6pt}{
        \tiny
        \color{darkblue}
        \bf
        \begin{tabular}{c}
          $\infty$-groupoid of {\color{orangeii}morphisms} of
          \\
          principal $\mathcal{G}$-bundles on $X$
        \end{tabular}
      }
      }
    }{
      \flat
      \,
      \Maps{}
        { X }
        { \mathbf{B}\mathcal{G} }
    }
    \quad
    \ar[
      rr,
      "\eta^{\scalebox{.5}{\textesh}}
         \,\circ\,
      \epsilon^{\flat}"
    ]
    &
    {\phantom{AA}}
    &
    \quad
    \overset{
      \mathclap{
      \raisebox{6pt}{
        \tiny
        \color{darkblue}
        \bf
        \begin{tabular}{c}
          $\infty$-groupoid of {\color{orangeii}concordances} of
          \\
          principal $\mathcal{G}$-bundles on $X$
        \end{tabular}
      }
      }
    }{
      \shape
      \,
      \Maps{}
        { X }
        { \mathbf{B}\mathcal{G} }
    }
    \quad
    \ar[
      rr,
      "\tau_0"
    ]
    &&
    (\PrincipalBundles{\mathcal{G}}_X)_{/\sim_{\mathrm{conc}}}
    \,.
    \end{tikzcd}
\end{equation}
\end{remark}

\medskip

\noindent
{\bf Classification of principal $\infty$-bundles.}
We demonstrate the existence of classifying spaces
for concordance classes of principal bundles over smooth manifolds
in the full generality of principal $\infty$-bundles internal to
$\SmoothInfinityGroupoids$.

\begin{definition}[Classifying shapes for smooth principal $\infty$-bundles]
  \label{ClassifyingSpaceOfPrincipalInfinityBundles}
  Let $\Topos$ a cohesive $\infty$-topos (Def. \ref{CohesiveInfinityTopos}).
  Then for $\mathcal{G} \,\in\, \Groups(\Topos)$ we write
  (recalling \eqref{ShapePreservesDeloopings})
  \vspace{-2mm}
  \begin{equation}
    \label{ClassifyingSpaceIsShapeOfModuliStack}
    B \mathcal{G}
    \;\coloneqq\;
    \shape
    (
      \mathbf{B}\mathcal{G}
    )
    \;\;\;
    \in
    \;
    \Topos_{\flat}
    \xhookrightarrow{\;}
    \Topos
  \end{equation}

  \vspace{-2mm}
  \noindent
  for the shape \eqref{TheModalitiesOnACohesiveInfinityTopos}
  of the universal moduli stack
  (according to Prop. \ref{GroupsActionsAndFiberBundles})
  of $\mathcal{G}$-principal bundles.
\end{definition}
\begin{example}[Classifying shapes for discrete structure $\infty$-groups]
  \label{ClassifyingShapesForDiscreteStructureInfinityGroups}
  If
  $G \,\in\, \Groups(\InfinityGroupoids) \xhookrightarrow{\Groups(\Discrete)}
  \Groups(\Topos)$
  so that $\flat \,G \;\simeq\; G \;\simeq\; \shape\, G$
  then the object $B G$ from Def. \ref{ClassifyingSpaceOfPrincipalInfinityBundles}
  coincides with the delooping $\mathbf{B}G$:
  \begin{equation}
    \label{ClassifyingShapeEquivalenceForDiscreteStructureInfinityGroups}
    \flat \mathbf{B}G
    \;\simeq\;
    \mathbf{B}G
    \simeq
    \shape \, \mathbf{B}G
    \;=\;
    B G
    \;\simeq\;
    \flat
    \,
    B G
    \,,
  \end{equation}
  by
  \eqref{ShapePreservesDeloopings}
  and
  \eqref{FlatPreservesDeloopings}.
  In particular this means that for all $\Gamma \,\in\, \Groups(\Topos)$
  $$
    B \Gamma \,\simeq\, B \,\shape\, \Gamma
    \,.
  $$
  For instance if $G \,\in\, \Groups(\Sets) \xhookrightarrow{\;} \Groups(\Topos)$
  and $n \in \mathbb{N}$, then
  $$
    B G \,\simeq\, K(G,1)
    \;\;\;
    \in
    \;
    \Groupoids
    \xhookrightarrow{\Discrete}
    \Topos
  $$
  is the homotopy type of a traditional
  Eilenberg-MacLane space, regarded among cohesive homotopy types.
\end{example}

\begin{proposition}[Milgram classifying spaces models classifying shape]
  \label{MilgramClassifyingSpaceModelClassifyingShape}
$\,$

\noindent  If $\Gamma \,\in\, \Groups(\DTopologicalSpaces)
    \xhookrightarrow{\Groups(\ContinuousDiffeology)} \Groups(\SmoothInfinityGroupoids)$
  is well-pointed (Ntn. \ref{WellPointedTopologicalGroup}),
  then the homotopy type of its Milgram classifying space $B \Gamma$
  \eqref{QuotientCoprojectionOfUniversalPrincipalBundle}
  coincides with the classifying shape
  $B \left(\ContinuousDiffeology(\Gamma)\right)$ \eqref{ClassifyingSpaceIsShapeOfModuliStack}:
 \vspace{-2mm}
  $$
    \Shape
    \left(
      \ContinuousDiffeology( B \Gamma )
    \right)
    \;\simeq\;
    B
    \left(
      \ContinuousDiffeology(\Gamma)
    \right)
    \;\;\;
    \in
    \;
    \InfinityGroupoids
    \xhookrightarrow{\;\Discrete\;}
    \SmoothInfinityGroupoids
    \,.
  $$
\end{proposition}

\begin{proof}
$\,$

\vspace{-7mm}
$$
  \def\arraystretch{1.75}
  \begin{array}{lll}
    B
    \left(
      \ContinuousDiffeology(\Gamma)
    \right)
    &
    \;=\;
    \Shape
    \left(
      \mathbf{B}
      \left(
        \ContinuousDiffeology(\Gamma)
      \right)
    \right)
    &
    \proofstep{ by Def. \ref{ClassifyingSpaceOfPrincipalInfinityBundles} }
    \\
    &
    \;\simeq\;
    \Shape
    \Big(\;
      \colimit{ [n] \in \Delta^{\mathrm{op}} }
      \ContinuousDiffeology(\Gamma)^{\times_n}
    \Big)
    &
    \proofstep{ by \eqref{DeloopingOfInfinityGroupAsColimit} }
    \\
    &
    \;\simeq\;
    \colimit{ [n] \in \Delta^{\mathrm{op}} }
    \,
    \Shape
    \left(
      \ContinuousDiffeology(\Gamma)^{\times_n}
    \right)
    &
    \proofstep{ by \eqref{InfinityAdjointPreservesInfinityLimits} }
    \\
    &
    \;\simeq\;
    \colimit{ [n] \in \Delta^{\mathrm{op}} }
    \,
    \left(
      \Shape
      \left(
        \ContinuousDiffeology(\Gamma)
      \right)
    \right)^{\times_n}
    &
    \proofstep{ by Rem. \ref{TheAxiomsOnTheShapeModality} }
    \\
    &
    \;\simeq\;
    \colimit{ [n] \in \Delta^{\mathrm{op}} }
    \,
    \left(
      \SingularSimplicialComplex
      \left(
        \ContinuousDiffeology(\Gamma)
      \right)
    \right)^{\times_n}
    &
    \proofstep{ by Prop. \ref{SmoothShapeGivenBySmoothPathInfinityGroupoid} }
    \\
    &
    \;\simeq\;
    \colimit{ [n] \in \Delta^{\mathrm{op}} }
    \,
    \left(
      \SingularSimplicialComplex(\Gamma)
    \right)^{\times_n}
    &
    \proofstep{ by Prop. \ref{DiffeologicalShapeOfContinuousDiffeologyIsWeakHomotopyType} }
    \\
    & \;\simeq\;
    \SingularSimplicialComplex
    \,
    \TopologicalRealization{}
    {
      \Gamma^{\times_\bullet}
    }
    &
    \proofstep{ by \eqref{GoodImpliesThatRealizationIsHomotopyColimit} }
    \\
    & \;\simeq\;
    \SingularSimplicialComplex
    (
      B \Gamma
    )
    &
    \proofstep{ by \eqref{QuotientCoprojectionOfUniversalPrincipalBundle} }
    \\
    & \;\simeq\;
    \SingularSimplicialComplex
    \,
    \ContinuousDiffeology
    (
      B \Gamma
    )
    &
    \proofstep{  by \eqref{DiffeologicalShapeOfContinuousDiffeologyIsWeakHomotopyType} }
    \\
    & \;\simeq\;
    \Shape
    \,
    \ContinuousDiffeology
    (
      B \Gamma
    )
    &
    \proofstep{ by Prop. \ref{SmoothShapeGivenBySmoothPathInfinityGroupoid}. }
  \end{array}
$$

\vspace{-7mm}
\end{proof}

\begin{theorem}[Concordance classification of smooth principal $\infty$-bundles]
  \label{ClassificationOfPrincipalInfinityBundles}
  $\,$

  \noindent
  For any $\mathcal{G} \,\in\, \Groups(\SmoothInfinityGroupoids)$,
  the space $B \mathcal{G}$ (Def. \ref{ClassifyingSpaceOfPrincipalInfinityBundles})
  the concordance classes (Def. \ref{ConcordanceOfSmoothPrincipalBundles})
  of
  $\mathcal{G}$-principal $\infty$-bundles
  over any $\SmoothManifold \,\in\, \SmoothManifolds \xhookrightarrow{\;}
  \SmoothInfinityGroupoids$ in that
  there is a natural bijection
  \vspace{-2mm}
  \begin{align*}
    \ConcordanceClasses{
      \PrincipalBundles{\mathcal{G}}(\SmoothInfinityGroupoids)_{\TopologicalSpace}
    }
    \;\;
    &=
    \;\;
    \Truncation{0}
    \,
    \shape
    \,
    \Maps{}
      { \SmoothManifold }
      { \mathbf{B}\mathcal{G} }
    \\
    \;\;
    & \simeq\;\;
    \Truncation{0}
    \,
    \Maps{\big}
      { \shape \TopologicalSpace }
      { B \mathcal{G} }
      \\
    \;\;
    & =
    \;\;
    H^1(X ;\, \shape \mathcal{G} )
    \;\;\;
    \in
    \;
    \Sets
    \,.
  \end{align*}
\end{theorem}
\begin{proof}
$\,$

\vspace{-8mm}
$$
  \def\arraystretch{1.6}
  \begin{array}{lll}
    \tau_0
    \,
    \PointsMaps{}
      { X }
      { B \mathcal{G} }
    )
    & \;\simeq\;
    \tau_0
    \,
    \Points
    \,
    \Maps{}
      { X }
      { B \mathcal{G} }
    &
    \proofstep{ by Lem. \ref{MappingSpaceConsistsOfCohesivePointsOfMappingStack}}
    \\
    &
    \;\simeq\;
    \tau_0
    \,
    \Points
    \,
    \Maps{}
      { X }
      { \shape \mathbf{B}\mathcal{G} }
    &
    \proofstep{ by Def. \ref{ClassifyingSpaceOfPrincipalInfinityBundles}}
    \\
    & \;\simeq\;
    \tau_0
    \,
    \Points
    \,
    \shape
    \,
    \Maps{}
      { X }
      { \mathbf{B}\mathcal{G} }
    &
    \proofstep{ by Thm. \ref{SmoothOkaPrinciple}}
    \\
    & \;\simeq\;
    \tau_0
    \,
    \Shape
    \,
    \Maps{}
      { X }
      { \mathbf{B}\mathcal{G} }
    &
    \proofstep{ by \eqref{IdempotencyOfCohesiveAdjoints} }
    \\
    & \;\simeq\;
    (\PrincipalBundles{\mathcal{G}}_X)_{/\sim_{\mathrm{conc}}}
    &
    \proofstep{by Prop. \ref{ShapeOfMappingStackIntoModuliStackISConcordances}. }
  \end{array}
$$

\vspace{-8mm}
\end{proof}

\begin{theorem}[Isomorphism classification of topological principal bundles]
  \label{ClassificationOfPrincipalBundlesAmongPrincipalInfinityBundles}
  In $\Topos \,=\, \SmoothInfinityGroupoids$\,, consider

  -- $\Gamma
      \,\in\,
    \Groups(\DHausdorffSpaces)
      \xhookrightarrow{ \Groups(\ContinuousDiffeology) \;}
    \Groups(\Topos)
    $
    being well-pointed (Ntn. \ref{WellPointedTopologicalGroup});

  -- $X \in \SmoothManifolds \xhookrightarrow{\ContinuousDiffeology\;} \Topos$.

  \vspace{1mm}
  \noindent
  Then there is a natural bijection
  \vspace{-2mm}
  $$
    \big(
      \PrincipalFiberBundles{\Gamma}(\DTopologicalSpaces)_X
    \big)_{/\sim_{\mathrm{iso}}}
    \;\;\simeq\;\;
    \tau_0
    \,
    \Maps{}
      { \shape \, X }
      { B \Gamma  }
  $$

  \vspace{-2mm}
  \noindent
  between isomorphism classes of $\Gamma$-principal fiber bundles over $\TopologicalSpace$
  (Ntn. \ref{TerminologyForPrincipalBundles})
  and homotopy classes of maps from (the shape of) $X$
  to the classifying shape of $\Gamma$
  (Def. \ref{ClassifyingSpaceOfPrincipalInfinityBundles}, Prop. \ref{MilgramClassifyingSpaceModelClassifyingShape}).
\end{theorem}
\begin{proof}
$\,$

\vspace{-7mm}
$$
 \def\arraystretch{1.8}
  \begin{array}{lll}
    \big(
      \PrincipalFiberBundles{\Gamma}(\DTopologicalSpaces)_X
    \big)_{/\sim_{\mathrm{iso}}}
    &
    \;\simeq\;
    \big(
      \PrincipalFiberBundles{\Gamma}(\DTopologicalSpaces)_X
    \big)_{/\sim_{\mathrm{conc}}}
    &
    \proofstep{ by Thm. \ref{ConcordanceClassesOfTopologicalPrincipalBundles} }
    \\
    &
    \;\simeq\;
    \big(
      \PrincipalBundles{\Gamma}(\SmoothInfinityGroupoids)_X
    \big)_{/\sim_{\mathrm{conc}}}
    &
    \proofstep{ by Prop. \ref{OrdinaryPrincipalBundlesAmongPrincipalInfinityBundles} }
    \\
    & \;\simeq\;
    \tau_0
    \,
    \Maps{}
      { \shape \TopologicalSpace }
      { B \Gamma }
    &
    \proofstep{ by Prop.  \ref{ClassificationOfPrincipalInfinityBundles}. }
  \end{array}
$$

\vspace{-7mm}
\end{proof}

\begin{remark}[Recovering the classification theorem for topological principal bundles via cohesion]
\label{RecoveringTheClassificationTheoremForTopologicalPrincipalBundlesViaCohesion}
$\,$

\vspace{-1.5mm}
\item {\bf (i)} The statement of
Thm. \ref{ClassificationOfPrincipalBundlesAmongPrincipalInfinityBundles}
reproduces the time-honored classification result for
topological principal bundles for Milgram-
(\cite{Milgram67}, review in
\cite[Thm. 3.5.1]{RudolphSchmidt17})
for the special case when the domain is a
smooth manifold.

\item {\bf (ii)}
While the statement of
Thm. \ref{ClassificationOfPrincipalBundlesAmongPrincipalInfinityBundles}
is classical,
the new proof via cohesive homotopy theory
has the advantage that it generalizes
to equivariant principal bundles,
where it provides a new classification result for
the case of truncated structure groups
(Thm. \ref{ClassificationOfGammaEquivariantPrincipalBundlesForGammaWithTruncatedClassifyingShape} below).
\end{remark}
This is what we turn to next.

\medskip

\subsection{Smooth equivariant principal $\infty$-bundles}
\label{SmoothEquivariantPrincipalInfinityBundles}

The general theory of principal $\infty$-bundles
in any $\infty$-topos,
discussed in \cref{SmoothPrincipalInfinityBundles}
{\it immediately} generalizes to the equivariant case,
by internalizing it into
the corresponding $\infty$-category of $G$-$\infty$-actions,
in direct analogy with the discussion in \cref{PrincipalBundlesInternalToTopologicalGActions},
and using that $G$-$\infty$-actions in any $\infty$-topos form
themselves again an $\infty$-topos.
Moreover, the classifying theory of smooth principal $\infty$-bundles
(from Prop.  \ref{ClassificationOfPrincipalInfinityBundles})
generalizes to a general classification result for
smooth equivariant principal $\infty$-bundles,
by the ``singular-smooth Oka principle''.

\medskip

\noindent
{\bf Equivariant $\infty$-groups.}
We generalize
(Def. \ref{GGroups} below)
the notion of {$G$-equivariant topological groups}
(Def. \ref{EquivariantTopologicalGroup})
to $\infty$-group stacks in $\infty$-toposes, and
(in higher analogy with Lem. \ref{EquivariantTopologicalGroupsAreSemidirectProductsWithG})
identify them with semidirect product $\infty$-groups
(Prop. \ref{EquivariantGroupsEquivalentToSplitExtensions}).
We currently fail to show that all $\infty$-groups internal to $G$-$\infty$-actions
are of this form (Rem. \ref{AreAllGroupsInGActionsEquivalentToSemidirectProductGroups}).

\medskip

Recall from \cite[Def. 2.101]{SS20OrbifoldCohomology}
that for $\mathcal{G} \,\in\, \Groups(\Topos)$
(Def. \ref{GGroups})
we have the
group-automorphism group
$\mathrm{Aut}_{\mathrm{Grp}}(\Topos)$ equipped with canonical
actions on $\mathcal{G}$ and on $\mathbf{B}\mathcal{G}$
preserving their basepoints:
\vspace{-2mm}
\begin{equation}
  \label{GroupAutomorphismActions}
  \hspace{-8mm}
  \arraycolsep=10pt
  \begin{array}{ccc}
    \mathrm{Aut}_{\mathrm{Grp}}(\mathcal{G})
    \;\coloneqq\;
    \mathrm{Aut}_\ast(\mathbf{B}\mathcal{G})
    &
  \begin{tikzcd}[column sep=small, row sep=large]
    \ast
    \ar[rr]
    \ar[
      d,
      "\mathrm{pt}_{\mathbf{B}G}"{left}
    ]
    \ar[
      drr,
      phantom,
      "\mbox{\tiny\rm(pb)}"
    ]
    &&
    \mathbf{B}\mathrm{Aut}_{\mathrm{Grp}}(\mathcal{G})
    \ar[
      d,
      "{
        (\mathrm{pt}_{\mathbf{B}\mathcal{G}})
        \!\sslash\!
        \mathrm{Aut}_{\mathrm{Grp}}(\mathcal{G})
      }"{description}
    ]
    \ar[
      dd,
      rounded corners,
      to path={
           -- ([xshift=+12pt]\tikztostart.east)
           --node[right]{
               \scalebox{.7}{$
                 \mathrm{id}
               $}
             } ([xshift=+12pt]\tikztotarget.east)
           -- (\tikztotarget.east)}
    ]
    &
    \\
    \mathbf{B}\mathcal{G}
    \ar[d]
    \ar[rr]
    \ar[
      drr,
      phantom,
      "\mbox{\tiny\rm(pb)}"
    ]
    &&
    (\mathbf{B}\mathcal{G}) \!\sslash\! \mathrm{Aut}_{\mathrm{Grp}}(\mathcal{G})
    \ar[d]
    &
    \\
    \ast
    \ar[rr]
    &&
    \mathbf{B} \mathrm{Aut}_{\mathrm{Grp}}(\mathcal{G})
    &
    {}
  \end{tikzcd}
  &
  \begin{tikzcd}[column sep=small, row sep=large]
    \ast
    \ar[rr]
    \ar[
      d,
      "e"{left}
    ]
    \ar[
      drr,
      phantom,
      "\mbox{\tiny\rm(pb)}"
    ]
    &&
    \mathbf{B}\mathrm{Aut}_{\mathrm{Grp}}(\mathcal{G})
    \ar[
      d,
      "{
        e
        \!\sslash\!
        \mathrm{Aut}_{\mathrm{Grp}}(\mathcal{G})
      }"{description}
    ]
    \ar[
      dd,
      rounded corners,
      to path={
           -- ([xshift=+12pt]\tikztostart.east)
           --node[right]{
               \scalebox{.7}{$
                 \mathrm{id}
               $}
             } ([xshift=+12pt]\tikztotarget.east)
           -- (\tikztotarget.east)}
    ]
    &
    \\
    \Gamma
    \ar[d]
    \ar[rr]
    \ar[
      drr,
      phantom,
      "\mbox{\tiny\rm(pb)}"
    ]
    &&
    \Gamma \!\sslash\! \mathrm{Aut}_{\mathrm{Grp}}(\mathcal{G})
    \ar[d]
    &
    \\
    \ast
    \ar[rr]
    &&
    \mathbf{B} \mathrm{Aut}_{\mathrm{Grp}}(\mathcal{G})
    &
    {}
  \end{tikzcd}
  \\
  \mbox{
    \tiny
    \color{darkblue}
    \bf
    group-automorphism group
  }
  &
  \mbox{
    \tiny
    \color{darkblue}
    \bf
    pointed action on delooping
  }
  &
  \mbox{
    \tiny
    \color{darkblue}
    \bf
    canonical action on group
  }
  \end{array}
\end{equation}

\begin{proposition}[Quotient groups by group automorphisms {\cite[Prop. 2.102]{SS20OrbifoldCohomology}}]
\label{QuotientGroupsByGroupAutomorphisms}
$\,$
Consider $G, \Gamma \in \Groups(\Topos)$ \eqref{LoopingAndDelooping}
and an action of $G$ on $\Gamma$ \eqref{EquivalenceBetweenActionsAndPrincipalBundlesAndSlices}
by group automorphisms,
i.e., a pasting diagram in $\Topos$ of this form:
\vspace{-1mm}
\begin{equation}
  \label{PullbackDiagramClassifyingGroupAutomorphismAction}
  \begin{tikzcd}[column sep=large]
    \Gamma
    \ar[
      rr,
      "\mbox{
        \tiny
        \color{greenii}
        \bf
        $G$-action on $\Gamma$ \eqref{EquivalenceBetweenActionsAndPrincipalBundlesAndSlices}
      }"{yshift=4pt, above}
    ]
    \ar[
      drr,
      phantom,
      "\mbox{\tiny\rm(pb)}"
    ]
    \ar[d]
    &&
    \Gamma \!\sslash\! G
    \ar[
      d,
      "\rho"
    ]
    \ar[
      drr,
      phantom,
      "\mbox{\tiny\rm(pb)}"
    ]
    \ar[
      rr,
      "\mbox{
        \tiny
        \color{greenii}
        \bf
        by group automorphisms
        \eqref{GroupAutomorphismActions}
      }"{yshift=4pt, above}
    ]
    &{\phantom{AAA}}&
    \Gamma \!\sslash\! \mathrm{Aut}_{\mathrm{Grp}}(\Gamma)
    \ar[
      d,
      "\rho_{\mathrm{Aut}}"
    ]
    \\
    \ast
    \ar[
      rr,
      "\mathrm{pt}_{\mathbf{B}G}"{below}
    ]
    &&
    \mathbf{B}G
    \ar[
      rr,
      "\vdash \rho"{below}
    ]
    &&
    \mathbf{B}
    \mathrm{Aut}_{\mathrm{Grp}}(\Gamma)
  \end{tikzcd}
\end{equation}

\vspace{-2mm}
\noindent
Then there is group structure on the homotopy quotient in the slice
\vspace{-1mm}
\begin{equation}
  \label{AutomorphismQuotientGroupAsGroupInSlice}
  \Gamma \!\sslash\! G
  \;\in\;
  \Groups(\Topos_{/\mathbf{B}G})
\end{equation}

\vspace{-2mm}
\noindent such that
\begin{equation}
  \label{DeloopingOfHomotopyQUotientByGroupAutomorphism}
  \mathbf{B}
  \left(
    \Gamma \!\sslash\! G
  \right)
  \;\simeq\;
  (\mathbf{B}\Gamma) \!\sslash\! G
  \;\;\;
  \in
  \;
  \Topos_{/\mathbf{B}G}
  \,.
\end{equation}
\end{proposition}
\begin{definition}[Equivariant $\infty$-groups]
\label{GGroups}
We write
\vspace{-2mm}
\begin{equation}
  \label{CategoryOfGGroups}
  \EquivariantGroups{G}(\Topos)
  \xhookrightarrow{\quad}
  \Groups(\Topos_{/\mathbf{B}G})
  \;\simeq\;
  \Groups
  \left(
    \Actions{G}(\Topos)
  \right)
\end{equation}

\vspace{-1mm}
\noindent for the $\infty$-category of group objects in the slice over $\mathbf{B}G$ that
arise via Prop. \ref{QuotientGroupsByGroupAutomorphisms}
from group objects equipped with $G$-actions by group automorphisms.
\end{definition}

\begin{example}[Equivariant $\infty$-group for trivial $G$-action]
  \label{EquivariantInfinityGroupForTrivialGAction}
  A $G$-action on $\Gamma$ is trivial precisely if
  $\Gamma \!\sslash\! G \,\simeq\, \Gamma \times \mathbf{B}G$,
  which by \eqref{DeloopingOfHomotopyQUotientByGroupAutomorphism}
  is the case
  (using Ex. \ref{HomotopyPullbackPreserbesProductProjections}) precisely if
  \vspace{-2mm}
  \begin{equation}
    \label{HomotopyQuotientForTrivialActionOnBGamma}
    \mathbf{B}(\Gamma \!\sslash\! G)
    \;\;
    \simeq
    \;\;
    (\mathbf{B}\Gamma) \times (\mathbf{B}G)
    \;\;\;
    \in
    \Topos_{/\mathbf{B}G}
    \,.
  \end{equation}
\end{example}

\begin{definition}[Split $\infty$-group extensions]
  \label{SplitGroupExtensions}
  Given $G \in \Groups(\Topos)$,

\noindent {\bf (i)}   a {\it split group extension} of $G$ is a diagram of the form
  \vspace{-2mm}
  \begin{equation}
    \label{SplitGroupExtension}
    \begin{tikzcd}[row sep=small]
      &&
      G
      \ar[
        d,
        "i\;"{left}
      ]
      \ar[
        drr,-,
        shift left=1pt
      ]
      \ar[
        drr,-,
        shift right=1pt
      ]
      \\
      \Gamma
      \ar[
        rr,
        "\mathrm{fib}(p)"{above}
      ]
      &&
      \Gamma \rtimes G
      \ar[
        rr,
        "p\;\;\;\;\;\;\;\;"{above}
      ]
      &&
      G
    \end{tikzcd}
    \;\;\;\;
    \in
    \;
    \EquivariantGroups{G}(\Topos)
    \,,
  \end{equation}

  \vspace{-1mm}
\noindent
  where the horizontal sequence is a fiber sequence.

 \noindent {\bf (ii)}  We write
  \vspace{-2mm}
  $$
    \SplitExtensions{G}
    \xhookrightarrow{\quad}
    \Groups
    (\Topos)^{G/}_{/G}
  $$

\vspace{-2mm}
\noindent
  for the full subcategory on the split extensions.
\end{definition}

\begin{proposition}[Equivariant groups are equivalent to split extensions of equivariance group]
 \label{EquivariantGroupsEquivalentToSplitExtensions}
  The operation of delooping of
  {\it split extensions}
  $\Gamma \xrightarrow{i} \Gamma \rtimes G \xrightarrow{p}$
  (Def. \ref{SplitGroupExtensions}).
  exhibits
  their equivalence with $G$-equivariant groups (Def. \ref{GGroups}):
  \vspace{-2mm}
  $$
    \begin{tikzcd}
      \SplitExtensions{G}(\Topos)
      \ar[
        r,
        "\sim"{yshift=-2pt}
      ]
      &
      \EquivariantGroups{G}(\Topos)
      \,,
    \end{tikzcd}
  $$

  \vspace{-2mm}
  \noindent
  such that \eqref{DeloopingOfHomotopyQUotientByGroupAutomorphism}
  is further identified with
  \vspace{-2mm}
  \begin{equation}
    \label{DeloopingOfSplitExtensionIsHomotopyQuotientOfBGamma}
    (\mathbf{B}\Gamma) \!\sslash\! G
    \;\simeq\;
    \mathbf{B}(\Gamma \rtimes G)\;.
  \end{equation}
\end{proposition}
\begin{proof}
First, to see that we have a functor, consider
for any split extension
the following pasting diagram:
\vspace{-2mm}
\begin{equation}
  \label{PastingDiagramFromSplitExtensionToAutomorphismAction}
  \begin{tikzcd}[column sep=large]
    \ast
    \ar[
      rr,
      "\mathrm{pt}_{\mathbf{B}G}"
    ]
    \ar[
      d,
      "\mathrm{pt}_{\mathbf{B}\Gamma}"{left}
    ]
    \ar[
      drr,
      phantom,
      "\mbox{\tiny\rm(pb)}"
    ]
    &&
    \mathbf{B}G
    \ar[
      d,
      "\mathbf{B}i"{right}
    ]
    \ar[
      dd,
      rounded corners,
      to path={
           -- ([xshift=+12pt]\tikztostart.east)
           --node[right]{
               \scalebox{.7}{$
                 \mathrm{id}
               $}
             } ([xshift=+12pt]\tikztotarget.east)
           -- (\tikztotarget.east)}
    ]
    \\
    \mathbf{B}\Gamma
    \ar[
      rr,
      "{\mathbf{B} \mathrm{fib}(p)}"{description}
    ]
    \ar[
      d
    ]
    \ar[
      drr,
      phantom,
      "\mbox{\tiny\rm(pb)}"
    ]
    &&
    \mathbf{B}(\Gamma \rtimes G)
    \ar[
      d,
      "{\mathbf{B}p}"{right}
    ]
    \\
    \ast
    \ar[
      rr,
      "\mathrm{pt}_{\mathbf{B}G}"{below}
    ]
    &&
    \mathbf{B}G
  \end{tikzcd}
\end{equation}

\vspace{-2mm}
\noindent
Here the bottom Cartesian square is the delooping of the given fiber sequence
\eqref{SplitGroupExtension}
(using \eqref{LoopingAndDelooping} in Prop. \ref{GroupsActionsAndFiberBundles}),
from which the top Cartesian square follows
by the pasting law \eqref{HomotopyPastingLaw}
and using that the pullback of an equivalence is an equivalence.
Hence
the bottom square exhibits a $G$-action on $\mathbf{B}\Gamma$
(by \eqref{EquivalenceBetweenActionsAndPrincipalBundlesAndSlices} in Prop. \ref{GroupsActionsAndFiberBundles})
with homotopy quotient as claimed in
\eqref{DeloopingOfSplitExtensionIsHomotopyQuotientOfBGamma};
and the top square exhibits $\mathrm{pt}_{\mathbf{B}\Gamma}$
as a homotopy fixed point of this action (\cite[Def. 2.97]{SS20OrbifoldCohomology}),
and hence a $G$-action on $\Gamma$ by group automorphisms
(\cite[Def. 2.101]{SS20OrbifoldCohomology}).

Similarly, $n$-morphism of split extensions are $n$-morphisms
of the diagrams \eqref{PastingDiagramFromSplitExtensionToAutomorphismAction}
which fix the top and bottom copies of $\mathbf{B}G$ and hence the
top and bottom copies of the point. These are equivalently
the pointed $n$-morphisms between $G$ actions on the $\mathbf{B}\Gamma$
and hence are equivalently $n$-morphisms between the $G$-actions
by group automorphisms on $\Gamma$.
This means that we have a fully faithful functor.

It just remains to see that this is essentially surjective.
Hence for a given $G$-action on $\Gamma$ by group automorphisms
consider the following pasting diagram:
\vspace{-2mm}
$$
  \begin{tikzcd}[row sep=23pt, column sep=large]
    \ast
    \ar[rr]
    \ar[d]
    \ar[
      drr,
      phantom,
      "\mbox{\tiny\rm(pb)}"
    ]
    &&
    \mathbf{B}G
    \ar[
      rr,
      "{
        \vdash \, \rho
      }"
    ]
    \ar[
      d
    ]
    \ar[
      drr,
      phantom,
      "{
        \mbox{\tiny\rm(pb)}
      }"
    ]
    &&
    \mathbf{B}\mathrm{Aut}_{\mathrm{Grp}}(\Gamma)
    \ar[
      d,
      "{
        (\mathrm{pt}_{\mathbf{B}\Gamma})
        \!\sslash\!
        \mathrm{Aut}_{\mathrm{Grp}}(\Gamma)
      }"{description}
    ]
    \ar[
      dd,
      rounded corners,
      to path={
           -- ([xshift=+11pt]\tikztostart.east)
           --node[right]{
               \scalebox{.7}{$
                 \mathrm{id}
               $}
             } ([xshift=+11pt]\tikztotarget.east)
           -- (\tikztotarget.east)}
    ]
    \\
    \mathbf{B}\Gamma
    \ar[rr]
    \ar[d]
    \ar[
      drr,
      phantom,
      "\mbox{\tiny\rm(pb)}"
    ]
    &&
    (\mathbf{B}\Gamma) \!\sslash\! G
    \ar[rr]
    \ar[
      d,
      "\rho"
    ]
    \ar[
      drr,
      phantom,
      "\mbox{\tiny\rm(pb)}"
    ]
    &&
    (\mathbf{B}\Gamma) \!\sslash\! \mathrm{Aut}_{\mathrm{Grp}}(\Gamma)
    \ar[
      d,
      "\rho_{\mathrm{Aut}}"{description}
    ]
    \\
    \ast
    \ar[
      rr,
      "\mathrm{pt}_{\mathbf{B}G}"{below}
    ]
    &&
    \mathbf{B}G
    \ar[
      rr,
      "\vdash \, \rho"{below}
    ]
    &&
    \mathbf{B}
    \mathrm{Aut}_{\mathrm{Grp}}(\Gamma)
  \end{tikzcd}
$$

\vspace{-2mm}
\noindent
Here the bottom square is that which defines the given
automorphism action \eqref{PullbackDiagramClassifyingGroupAutomorphismAction}
and all the other squares follow with the pasting law
\eqref{HomotopyPastingLaw}
from \eqref{GroupAutomorphismActions}.
Now the two squares on the left exhibit the delooping of a
split group extension since
$(\mathbf{B}\Gamma) \!\sslash\! G$ is connected
(by Lemma \ref{HomotopyQuotientPreservesConnectedness}), and hence
of the form \eqref{DeloopingOfSplitExtensionIsHomotopyQuotientOfBGamma}.
\end{proof}

\begin{remark}
[Are all groups in $G$-Actions equivalent to semidirect product groups?]
\label{AreAllGroupsInGActionsEquivalentToSemidirectProductGroups}
One might expect that any group object in the slice
$\Gamma \!\sslash\! G \,\in\, \Groups( \Topos_{/\mathbf{B}G})$
is an equivariant group in the sense of Def. \ref{GGroups},
hence that the full inclusion \eqref{CategoryOfGGroups} is in fact
essentially surjective and hence an equivalence.

\noindent {\bf (i)} We currently do not have a proof that this
is the case, but the following pasting diagram shows something close,
at least:
\vspace{-2mm}
$$
  \begin{tikzcd}[column sep=large, row sep=12pt]
    &
    G
    \ar[r]
    \ar[
      d,
      "{(\NeutralElement,\mathrm{id})}"{left}
    ]
    \ar[
      dr,
      phantom,
      "\mbox{\tiny\rm(d)}"
    ]
    &
    \ast
    \ar[d]
    \\
    G
    \ar[r]
    \ar[d]
    \ar[
      dr,
      phantom,
      "\mbox{\tiny\rm(d)}"
    ]
    &
    \Gamma \times G
    \ar[
      r,
      "\mathrm{pr}_1"{description}
    ]
    \ar[
      d,
      "\rho"{left}
    ]
    \ar[
      dr,
      phantom,
      "{\mbox{\tiny (c)}}"
    ]
    &
    \Gamma
    \ar[d]
    \ar[rr]
    \ar[
      drr,
      phantom,
      "{\mbox{\tiny (b)}}"
    ]
    &&
    \ast
    \ar[d]
    \\
    \ast
    \ar[r]
    &
    \Gamma
    \ar[d]
    \ar[r]
    \ar[
      dr,
      phantom,
      "{\mbox{\tiny (b)}}"
    ]
    &
    \Gamma \!\sslash\! G
    \ar[rr]
    \ar[d]
    \ar[
      drr,
      phantom,
      "\mbox{\tiny\rm(a)}"
    ]
    &&
    \mathbf{B}G
    \mathrlap{
      \;
      \simeq
      \Gamma \!\sslash\! (\Gamma \rtimes G)
    }
    \ar[d]
    \\
    &
    \ast
    \ar[r]
    &
    \mathbf{B}G
    \ar[rr]
    &&
    \underset{\mathbf{B}G}{\sum} \mathbf{B} (\Gamma \!\sslash\! G)
    \mathrlap{
      \;
      \simeq
      \mathbf{B}(\Gamma \rtimes G)\;.
    }
  \end{tikzcd}
$$

\vspace{-3mm}
\noindent
Here:

\noindent
(a) is the homotopy pullback that exhibits $\Gamma \!\sslash\! G$
as a group object in the slice over $\mathbf{B}G$;

\noindent
(b) is the homotopy pullback that exhibits $\Gamma \!\sslash\! G$
as the homotopy quotient of a $G$-action on $\Gamma$;

\noindent
(c) is the homotopy fiber product which exhibits
the shear map equivalence of
$\Gamma$ as a principal $G$-bundle over $\Gamma \!\sslash\! G$;

\noindent
(d) is a homotopy pullback implied from this by the pasting law \eqref{HomotopyPastingLaw}.

\noindent  We see from this that $\Gamma \times G$ carries a group structure,
to be denoted $\Gamma \rtimes G$, whose delooping is
the bottom right object:
\vspace{-2mm}
\begin{equation}
  \underset{\mathbf{B}G}{\sum}
  \mathbf{B}(\Gamma \!\sslash\! G)
  \;\simeq\;
  \mathbf{B}(\Gamma \rtimes G)\;.
\end{equation}

\vspace{-3mm}
\noindent Moreover,
$G \xrightarrow{(\NeutralElement,\mathrm{id})} \Gamma \times G$ is
exhibited as a homomorphism of group objects.
Also, we see that
$G \xrightarrow{(\NeutralElement,\mathrm{id})} \Gamma \times G \xrightarrow{\rho} \Gamma$
are $G$-equivariant maps for $G$ acting by right multiplication on itself.
This implies that $\rho$ is the given $G$-action.

\noindent {\bf (ii)} In conclusion, this shows almost all the structure required of
a split extension according to Prop. \ref{EquivariantGroupsEquivalentToSplitExtensions},
except that we have yet to conclude group structure on $\Gamma$ itself,
such that the morphism $\Gamma \xrightarrow{\;} \Gamma \!\sslash\! G \,\simeq\,
\Gamma \rtimes G$ is in fact a {\it homo}morphism.
\end{remark}

\begin{lemma}[Homotopy quotients by equivariant actions of equivariant groups]
\label{HomotopyQuotientsByEquivariantActionsOfEquivariantGroups}
Given

- $\HomotopyQuotient{\Gamma}{G}
     \,\in \,
       \EquivariantGroups{G}(\Topos)
       \hookrightarrow
       \Groups\big(
         \SliceTopos{\mathbf{B}G}
       \big)
       $
  (Def. \ref{GGroups})

- $(\HomotopyQuotient{\Gamma}{G})
    \acts
    (\HomotopyQuotient{A}{G})
    \,\in\,
    \Actions{(\HomotopyQuotient{\Gamma}{G})}
    \big(
      \SliceTopos{\mathbf{B}G}
    \big)
  $

\noindent
then the underlying object of the homotopy quotient of this equivariant action formed in $\SliceTopos{\mathbf{B}G} \,\simeq\, \Actions{G}{\Topos}$
is equivalent to the homotopy quotient of $A$ by $\Gamma \rtimes G$ in $\Topos$, as exhibited by the following diagram
\begin{equation}
  \begin{tikzcd}[column sep=-7pt]
    A
    \ar[r]
    \ar[d]
    \ar[
      dr,
      phantom,
      "{\mbox{\tiny\rm(pb)}}"
    ]
    &[+30pt]
    \HomotopyQuotient{A}{G}
    \ar[r]
    \ar[d]
    \ar[
      rd,
      phantom,
      "{\mbox{\tiny\rm(pb)}}"{xshift=-4pt}
    ]
    &[20pt]
    \HomotopyQuotient
      { A }
      { (\Gamma \rtimes G) }
    \ar[d]
    \ar[rr, shorten=-2pt, "\sim"]
    \ar[
      drr,
      phantom,
      "{\mbox{\tiny\rm(pb)}}"
    ]
    &&
    \HomotopyQuotient
      {
        \big(
        \HomotopyQuotient
          {A}{G}
        \big)
      }
      {
        \big(
        \HomotopyQuotient
          {\Gamma}{G}
        \big)
      }
    \ar[d]
    \\
    \ast
    \ar[r]
    &
    \mathbf{B}G
    \ar[r]
    &
    \mathbf{B}
    (\Gamma \rtimes G)
    \ar[dr]
    \ar[
      rr,
      shorten=-2pt,
      "\sim",
      "{
        \mbox{\tiny\rm Prop. \ref{EquivariantGroupsEquivalentToSplitExtensions}}
      }"{swap, pos=.4}
    ]
    &&
    \HomotopyQuotient
      { \mathbf{B}\Gamma }
      { G }
    \mathrlap{\,.}
    \ar[dl]
    \\[-10pt]
    &
    &
    &
    \mathbf{B}G
  \end{tikzcd}
\end{equation}
\end{lemma}
\begin{proof}
This uses
Prop. \ref{EquivariantGroupsEquivalentToSplitExtensions}
for the identification in the bottom right and then
Prop. \ref{HomotopyQuotientsAndPrincipaInfinityBundles} with Prop. \ref{SliceInfinityTopos},
and the pasting law \eqref{HomotopyPastingLaw} to identify the
pullbacks.
\end{proof}

\medskip

\noindent
{\bf Equivariant principal $\infty$-bundles.}
In view of Prop. \ref{GroupsActionsAndFiberBundles}
and in direct analogy with Def. \ref{EquivariantPrincipalBundle},
we consider the following evident definition of equivariant
principal $\infty$-bundles in any $\infty$-topos:

\begin{definition}[Equivariant principal $\infty$-bundles]
  \label{GEquivariantGammaPrincipalBundles}
  For

  -- $G \in \Groups(\Topos)$,

  -- $\Gamma \!\sslash\! G \,\in\, \EquivariantGroups{G}(\Topos)$
  (Def. \ref{GGroups}),

  \noindent
  the {\it $G$-equivariant $\Gamma$-principal bundles}
  on any
  $
    G \acts \, X \,\in\,
    \Actions{G}(\Topos)
    \simeq
    \Topos_{/\mathbf{B}G}
  $
  are the $\mathcal{G}$-principal bundles
  \eqref{EquivalenceBetweenActionsAndPrincipalBundlesAndSlices}
  for $\mathcal{G} \coloneqq \Gamma \!\sslash\! G$ from
  \eqref{AutomorphismQuotientGroupAsGroupInSlice},
  hence the objects of the following equivalent
  $\infty$-groupoids:
  \vspace{-3mm}
  \begin{equation}
  \label{InfinityGroupoidOfEquivariantPrincipalBundles}
  \def\arraystretch{1.5}
  \begin{array}{lll}
    \EquivariantPrincipalBundles{G}{\Gamma}(\Topos)_X
    &
    \;:=\;
    \PrincipalBundles{(\Gamma \!\sslash\! G)}
    (
      \Topos_{/\mathbf{B}G}
    )
    _{(X \!\sslash\! G)}
    &
    \mbox{\small internalization into slice}
    \\
    & \;\simeq\;
    \SlicePointsMaps{\big}{\mathbf{B}G}
      { X \!\sslash\! G }
      { \mathbf{B}(\Gamma\sslash G) }
    &
    \mbox{\small by Prop. \ref{GroupsActionsAndFiberBundles}}
    \\
    & \;\simeq\;
    \Actions{G}
    (
      \Topos
    )
    \left(
      G \acts \, X,
      \,
      G \acts \, (\mathbf{B}\Gamma)
    \right)
    &
    \mbox{\small by Prop. \ref{GroupsActionsAndFiberBundles}}
    \,.
  \end{array}
  \end{equation}
\end{definition}

\begin{remark}[Equivariant principal $\infty$-bundles and their underlying principal $\infty$-bundles]
  Let $\Gamma \in G\Groups(\Topos)$.
  The base change
  \eqref{BaseChangeAdjointTriple}
  along $\ast \xrightarrow{\mathrm{pt}_{\mathbf{B}G}} \mathbf{B}G$
  exhibits
  plain $\Gamma$-principal bundles \eqref{EquivalenceBetweenActionsAndPrincipalBundlesAndSlices}
  underlying
  (Rem. \ref{UnderlyingObjectsOfGAction})
  $G$-equivariant $\Gamma$-principal bundles on
  any $G \acts \, X$ (Def. \ref{GEquivariantGammaPrincipalBundles}):
  \vspace{-2mm}
  \begin{equation}
    \label{UnderlyingBundleOfEquivariantPrincipalBundle}
    \begin{tikzcd}[row sep=4pt]
      G\, \mathrm{Equiv}\PrincipalBundles{\Gamma}(\Topos)_X
      \ar[
        rr,
        "\underlying"{above}
      ]
      \ar[d,phantom,"\simeq"{sloped}]
      &&
      \PrincipalBundles{\Gamma}(\Topos)_X
      \ar[d,phantom,"\simeq"{sloped}]
      \\
      \Topos_{/\mathbf{B}G}
      \left(
        X \!\sslash\! G,
        \,
        (\mathbf{B}\Gamma) \!\sslash\! G
      \right)
      \ar[
        rr,
        "{
          (\mathrm{pt}_{\mathbf{B}G})^\ast
          =
          \mathrm{fib}
        }"
      ]
      &&
      \Topos
      (
        X,
        \,
        \mathbf{B}\Gamma
      )
      \\[-8pt]
      \scalebox{0.9}{$
      (\vdash P \!\sslash\! G)
      \;\simeq\;
      (\vdash P) \!\sslash\! G
      $}
        &\longmapsto&
    \scalebox{0.9}{$
      \vdash P
    $}
    \end{tikzcd}
  \end{equation}
    This situation is shown in the following diagram on the left,
    all whose squares are Cartesian:
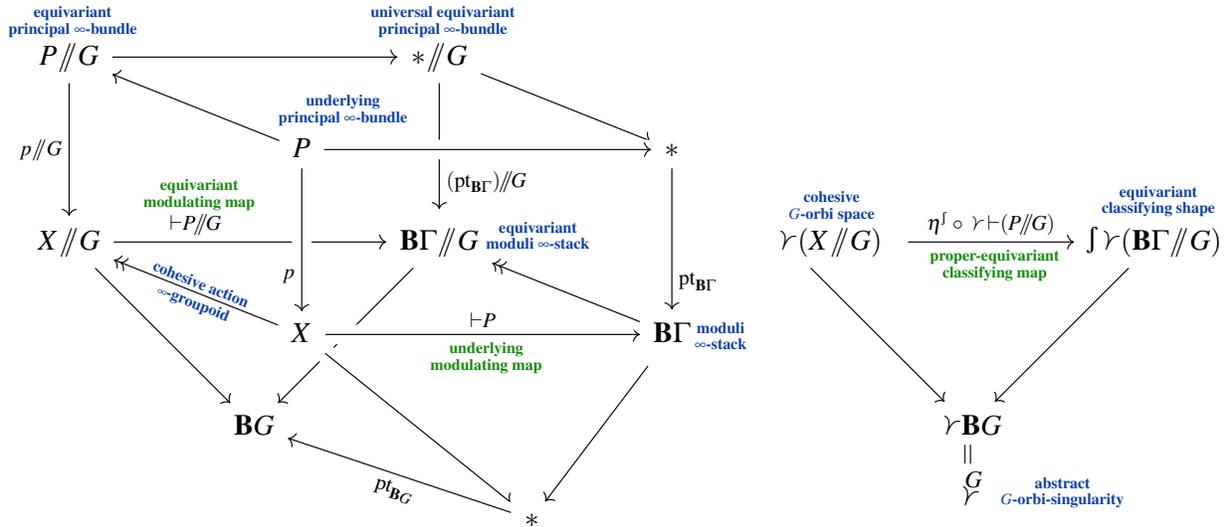
\vspace{-3mm}
$$
  \begin{tikzcd}[column sep={between origins, 35pt}, row sep={between origins, 35pt}]
    \overset{
      \mathclap{
      \raisebox{4pt}{
        \tiny
        \color{darkblue}
        \bf
        \def\arraystretch{.9}
        \begin{tabular}{c}
          equivariant
          \\
          principal $\infty$-bundle
        \end{tabular}
      }
      }
    }{
      P \!\sslash\! G
    }
    \ar[
      dd,
      "p \sslash G"{left}
    ]
    \ar[rrrr]
    \ar[from=drrr   ]
    &&&[-17pt] &[17pt]
    \overset{
      \mathclap{
      \raisebox{4pt}{
        \tiny
        \color{darkblue}
        \bf
        \def\arraystretch{.9}
        \begin{tabular}{c}
          universal equivariant
          \\
          principal $\infty$-bundle
        \end{tabular}
      }
      }
    }{
      \HomotopyQuotient
        {\ast}{G}
    }
    \ar[
      dd,
      "(\mathrm{pt}_{\mathbf{B}\Gamma}) \!\sslash\! G "{right, pos=.78}
    ]
    \ar[drrr]
    &&[-17pt]
    \\
    &&&
    \overset{
      \mathrlap{
      \;\;
      \raisebox{3pt}{
        \colorbox{white}{
        \tiny
        \color{darkblue}
        \bf
        \def\arraystretch{.9}
        \begin{tabular}{c}
          underlying
          \\
          principal
          $\infty$-bundle
        \end{tabular}
      }
      }
      }
    }{
      P
    }
    \ar[rrrr, crossing over]
    &&&&
    \ast
    \ar[
      dd,
      "\mathrm{pt}_{\mathbf{B}\Gamma}"{right, pos=.8}
    ]
    \\
    X \!\sslash\! G
    \ar[
      rrrr,
      "{
        \overset{
          \raisebox{3pt}{
            \tiny
            \color{greenii}
            \bf
            \def\arraystretch{.9}
            \begin{tabular}{c}
              equivariant
              \\
              modulating map
            \end{tabular}
          }
        }{
          \vdash \HomotopyQuotient{P}{G}
        }
      }"{above, pos=.3}
    ]
    \ar[ddrr]
    &&&&
    \HomotopyQuotient
      { \mathbf{B}\Gamma }
      { G }
    \mathrlap{
      \!\!\!\!\!\!
      \raisebox{3pt}{
        \tiny
        \color{darkblue}
        \bf
        \def\arraystretch{.9}
        \begin{tabular}{c}
          equivariant
          \\
          moduli $\infty$-stack
        \end{tabular}
      }
    }
    \ar[ddll]
    &&&
    &
    \qquad \quad
    \overset{
      \raisebox{4pt}{
        \tiny
        \color{darkblue}
        \bf
        \def\arraystretch{.9}
        \begin{tabular}{c}
          cohesive
          \\
          $G$-orbi space
        \end{tabular}
      }
    }{
      \orbisingular(
       \HomotopyQuotient
         {X}{G}
      )
    }
    \ar[ddrr]
    \ar[
      rrrr,
      "{
        \eta^{\scalebox{.5}{\textesh}}
        \;\circ\;
        \orbisingular
        \,
        \vdash(\HomotopyQuotient{P}{G})
      }"{above},
      "{
        \mbox{
          \tiny
          \color{greenii}
          \bf
          \def\arraystretch{.9}
          \begin{tabular}{c}
            proper-equivariant
            \\
            classifying map
          \end{tabular}
        }
      }"{below}
    ]
    &&{\phantom{a}}&&
    \overset{
      \mathclap{
      \raisebox{6pt}{
        \tiny
        \color{darkblue}
        \bf
        \def\arraystretch{.9}
        \begin{tabular}{c}
          equivariant
          \\
          classifying shape
        \end{tabular}
      }
      }
    }{
      \shape
      \orbisingular
      (
        \HomotopyQuotient
          { \mathbf{B}\Gamma }
          { G }
      )
    }
    \ar[ddll]
    \\
    &&&
    X
    \ar[
      rrrr,
      crossing over,
      "{ \vdash P }"{above},
      "{
        \mbox{
          \tiny
          \color{greenii}
          \bf
          \def\arraystretch{.9}
          \begin{tabular}{c}
            underlying
            \\
            modulating map
          \end{tabular}
        }
      }"{below}
    ]
    \ar[
      ulll,
      ->>,
      "{
        \mbox{
          \tiny
          \color{darkblue}
          \bf
          cohesive action
        }
      }"{sloped, above, yshift=-1pt },
      "{
        \mbox{
          \tiny
          \color{darkblue}
          \bf
          $\infty$-groupoid
        }
      }"{sloped, below, yshift=+1pt}
    ]
    \ar[ddrr, crossing over]
    \ar[
      from=uu,
      crossing over,
      "p"{left, near end}
    ]
    &&&&
    \mathbf{B}\Gamma
    \mathrlap{
      \!\!\!\!\!\!\!
      \mbox{
        \tiny
        \color{darkblue}
        \bf
        \def\arraystretch{.9}
        \begin{tabular}{c}
          moduli
          \\
          $\infty$-stack
        \end{tabular}
      }
    }
    \ar[
      ulll,
      ->>,
      crossing over
    ]
    \ar[ddll]
    \\
    &&
    \mathbf{B}G
    &&&&&&&&
    \orbisingular \mathbf{B}G
    \ar[
      d,
      phantom,
      "{
        \overset{
          \rotatebox{90}{$=$}
        }{
          \orbisingularG
          \mathrlap{
            \!\!\!
            \mbox{
              \tiny
              \color{darkblue}
              \bf
              \def\arraystretch{.9}
              \begin{tabular}{c}
                abstract
                \\
                $G$-orbi-singularity
              \end{tabular}
            }
          }
        }
      }"{pos=.37}
    ]
    \\
    &&  & &&
    \ast
    \ar[
      ulll,
      crossing over,
      "\mathrm{pt}_{\mathbf{B}G}"{sloped, below}
    ]
    & {} & {} & {} & {} & {}
  \end{tikzcd}
$$
The diagram on the right
shows the corresponding classifying map
to the proper-equivariant classifying shape, discussed in
\cref{EquivariantModuliStacks}
below
(Def. \ref{EquivariantClassifyingStack}).
\end{remark}

\begin{example}[Semidirect product-group homotopy quotient of underlying equivariant principal $\infty$-bundles]
\label{SemidirectProductHomotopyQuotientOfEquivariantInfinityBundle}
Let $G \acts \Gamma \,\in\,
       \EquivariantGroups{G}(\Topos)
       $
      (Def. \ref{GGroups})
with semidirect product $\infty$-group
$\Gamma \rtimes G$ \eqref{DeloopingOfSplitExtensionIsHomotopyQuotientOfBGamma}
then the total space $P$ of the underlying $\Gamma$-principal bundle $P$
\eqref{UnderlyingBundleOfEquivariantPrincipalBundle}
of a $G$-equivariant $\Gamma$-principal bundle satisfies
\begin{equation}
  \label{SlicedEquivalenceForSemidirectProductHomotopyQuotientOfEquivariantInfinityBundle}
  \begin{tikzcd}[column sep=10pt, row sep=14pt]
  \HomotopyQuotient
    { P }{ (\Gamma \rtimes G) }
  \ar[rr, "{\sim}"{yshift=-1pt}]
  \ar[dr]
  &&
  \HomotopyQuotient
    { X }
    { G }
  \ar[
    dl,
    "{
      \vdash
      \,
      \HomotopyQuotient{P}{G}
    }"{swap, sloped, pos=.35}
  ]
  \\
  &
  \mathbf{B}(\Gamma \rtimes G)
  \ar[
    d,
    "{
      \mathbf{B}\mathrm{pr}_2
    }"
  ]
  \\
  &
  \mathbf{B}G
  \mathrlap{\,.}
  \end{tikzcd}
\end{equation}
This follows via Prop. \ref{GroupsActionsAndFiberBundles} by applying the pasting law \eqref{HomotopyPastingLaw} to this diagram:
$$
  \begin{tikzcd}
    P
    \ar[rr]
    \ar[d]
    \ar[
      drr,
      phantom,
      "{\mbox{\tiny\rm(pb)}}"
    ]
    &&
    \ast
    \ar[d]
    \\
    \HomotopyQuotient
      { P }{ G }
    \ar[rr]
    \ar[d]
    \ar[
      drr,
      phantom,
      "{\mbox{\tiny\rm(pb)}}"
    ]
    &&
    \mathbf{B}G
    \ar[d]
    \\
    \HomotopyQuotient
      { X }{ G }
    \ar[
      rr,
      "{
        \vdash
        \,
        \HomotopyQuotient{P}{G}
      }"{description}
    ]
    \ar[dr]
    &&
    \mathbf{B}(\Gamma \rtimes G)
    \ar[dl]
    \\[-14pt]
    &
    \mathbf{B}G
    \mathrlap{\,.}
  \end{tikzcd}
$$
In particular, \eqref{SlicedEquivalenceForSemidirectProductHomotopyQuotientOfEquivariantInfinityBundle}
implies a natural identification of the fixed loci of $G$-equivariant $\Gamma$-principal bundles for
(sub-)group  $\eta \xrightarrow{\;} \Gamma \rtimes G$ with the local sections of its modulating  morphism:
\begin{equation}
  \label{FixedLociOfEquivariantPrincipalBundlesViaSliceHom}
  \def\arraystretch{1.4}
  \begin{array}{lll}
    P^{\eta}
    &
    \;\simeq\;
    \SliceMaps{\big}{\mathbf{B}(\Gamma \rtimes G)}
      {
        \mathbf{B}\eta
      }
      {
        \HomotopyQuotient
          { P }
          { (\Gamma \rtimes G) }
      }
    &
    \proofstep{
      \eqref{FixedPointsForSubgroupAsSliceMappingStack}
    }
    \\
    &
    \;\simeq\;
    \SliceMaps{\big}{\mathbf{B}(\Gamma \rtimes G)}
      {
        \mathbf{B}\eta
      }
      {
        \HomotopyQuotient
          { X }
          { G }
      }
    &
    \proofstep{
      by \eqref{SemidirectProductHomotopyQuotientOfEquivariantInfinityBundle}.
    }
  \end{array}
\end{equation}
\end{example}

\begin{proposition}[Equivariant principal $\infty$-bundles for trivial action on structure $\infty$-groups]
  \label{EquivariantBundlesForTrivialActionOnStructureGroup}
  When the $G$-action on
  the structure $\infty$-group $\Gamma$ is trivial,
  so that
  $
    (\mathbf{B}\Gamma) \!\sslash\! G
    \simeq
    (\mathbf{B}\Gamma) \times (\mathbf{B}G)
  $
  (Ex. \ref{EquivariantInfinityGroupForTrivialGAction}),
  then $G$-equivariant $\Gamma$-principal $\infty$-bundles
 (Def. \ref{GEquivariantGammaPrincipalBundles})
  are equivalently:

\noindent {\bf (i)}
$\Gamma$-principal bundles on the $G$-quotient stack:
\vspace{-3mm}
  \begin{equation}
    \label{EquivariantPrincipalBundlesAsPrincipalBundlesOnQuotientStack}
    \begin{tikzcd}[row sep=-5pt, column sep=small]
      \EquivariantPrincipalBundles{G}{\Gamma}(\Topos)_{X}
      \ar[
        rr,
        "\sim"{above, yshift=-2pt}
      ]
      &&
      \PrincipalBundles{\Gamma}(\Topos)_{X \!\sslash\! G}
      \\
      \scalebox{0.85}{$
       \big(
        G \acts \, X
        \xrightarrow{ c }
        \mathbf{B}\Gamma
      \big)
      $}
      &\longmapsto&
      \scalebox{0.85}{$
      \big(
        X \!\sslash\! G
        \xrightarrow{ c \sslash G }
        \mathbf{B}\Gamma
      \big)
      $}
    \end{tikzcd}
  \end{equation}

\vspace{-3mm}
\noindent
{\bf (ii)} $G$-Actions on $\Gamma$-principal bundles
\vspace{-2mm}
  \begin{equation}
    \label{GEquivariantPrincBundlesEquivalentToGActionsOnPrincBundles}
    \begin{tikzcd}[column sep=20pt, row sep=small]
      \EquivariantPrincipalBundles{G}{\Gamma}(\Topos)_X
      \ar[
        rr,
        "\sim"{yshift=-1pt}
      ]
      \ar[
        dr,
        "{
          \def\arraystretch{.6}
          \begin{array}{c}
           \color{greenii} \bf
            \underlying
            \\
            \mbox{
              \tiny
              \eqref{UnderlyingBundleOfEquivariantPrincipalBundle}
            }
          \end{array}
        }"{swap, pos=.7, yshift=2pt}
      ]
      &&
      \Actions{G}
      \left(
        \PrincipalBundles{\Gamma}(\Topos)
      \right)_{\scalebox{.7}{$G \acts \, X$}}
      \ar[
        dl,
        "{
          \def\arraystretch{.6}
          \begin{array}{c}
       \color{greenii} \bf
            \underlying
            \\
            \mbox{\tiny\rm\eqref{UnderlyingObjectsOfGAction}}
          \end{array}
        }"{pos=.7, yshift=2pt}
      ]
      \\
      &
      \PrincipalBundles{\Gamma}(\Topos)_X
      \mathrlap{\,,}
    \end{tikzcd}
  \end{equation}

  \vspace{-2mm}
  \noindent
  where on the right we consider $G$ base-changed
  \eqref{BaseChangeAdjointTripleToAbsoluteContext}
  to a group object in $\Gamma$-principal bundles
  \eqref{EquivalenceBetweenActionsAndPrincipalBundlesAndSlices}:
  \vspace{-1mm}
  $$
    G
      \;\in\;
    \Groups(\Topos)
    \xrightarrow{ \; \Groups( (-) \times \mathbf{B}\mathcal{G} ) \;}
    \Groups
    (
      \Topos_{/\mathbf{B}\mathcal{G}}
    )
    \xrightarrow{ \;\; \simeq \;\;}
    \Groups
    \left(
      \PrincipalBundles{\mathcal{G}}(\Topos)
    \right)
    \,.
  $$
\end{proposition}
\begin{proof}
The first equivalence in the composite of the
following sequence of natural equivalences:

  \vspace{-2mm}
  $$
    \def\arraystretch{1.2}
    \begin{array}{lll}
      \EquivariantPrincipalBundles{G}{\Gamma}(\Topos)_X
      &
      \;\simeq\;
      \Topos_{/\mathbf{B}G}
      \left(
        X \!\sslash\! G,
        \,
        \mathbf{B}(\Gamma \!\sslash\! G)
      \right)
      &
      \mbox{\small by \eqref{InfinityGroupoidOfEquivariantPrincipalBundles} }
      \\
      & \;\simeq\;
      \Topos_{/\mathbf{B}G}
      \left(
        X \!\sslash\! G,
        \,
        (\mathbf{B}\Gamma )\!\sslash\! G
      \right)
      &
      \mbox{\small by \eqref{DeloopingOfHomotopyQUotientByGroupAutomorphism}}
      \\
      & \;\simeq\;
      \Topos_{/\mathbf{B}G}
      \left(
        X \!\sslash\! G,
        \,
        (\mathbf{B}\Gamma)
        \times
        (\mathbf{B}G)
      \right)
      &
      \mbox{\small by \eqref{HomotopyQuotientForTrivialActionOnBGamma}}
      \\
      & \;\simeq\;
      \Topos
      \left(
        X \!\sslash\! G,
        \,
        \mathbf{B}\Gamma
      \right)
      &
      \mbox{\small by \eqref{BaseChangeAdjointTripleToAbsoluteContext}}
      \\
      & \;\simeq\;
      \PrincipalBundles{\Gamma}(\Topos)_{X \sslash G}
      &
      \mbox{\small by \eqref{EquivalenceBetweenActionsAndPrincipalBundlesAndSlices}}
      \,.
    \end{array}
  $$

This already implies also the second equivalence, since now
both sides are identified with maps
of the form $X \!\sslash\! G \xrightarrow{ c \sslash G } \mathbf{B}\Gamma$.
It just remains to see that the underlying objects coincide,
in that the diagram in \eqref{GEquivariantPrincBundlesEquivalentToGActionsOnPrincBundles}
commutes.
To find the underlying object of a $G$-action in
$\Gamma$-principal bundles, we need
to compute (by \eqref{UnderlyingObjectOfActionViaBaseChange})
the base change
\vspace{-2mm}
$$
  (
    \Topos_{/\mathbf{B}\Gamma}
  )_{/(\mathbf{B}\Gamma) \times (\mathbf{B}G)}
  \xrightarrow{
    \;\;
    \mathrm{pt}^\ast_{(\mathbf{B}\Gamma) \times (\mathbf{B}G)}
    \;\;
  }
  \Topos_{/\mathbf{B}\Gamma}
  \,.
$$

\vspace{0mm}
\noindent
Since fiber products in a slice are computed in the underlying
$\infty$-topos (e.g. \cite[Prop. 2.53]{SS20OrbifoldCohomology})
we are thus reduced to checking the pullback square
as shown on the
top of the following diagram on the left:
\vspace{-2mm}
$$
  \begin{tikzcd}[column sep=large]
    X
    \ar[rr]
    \ar[
      d,
      "c"{left}
    ]
    \ar[
      drr,
      phantom,
      "\mbox{\tiny\rm(pb)}"
    ]
    &&
    X \!\sslash\! G
    \ar[
      d,
      "{
        (c \sslash G, \rho_G)
      }"{right}
    ]
    \\
    \mathbf{B}\Gamma
    \ar[
      rr,
      "{
        (\mathrm{id}_{\mathbf{B}\Gamma}, \mathrm{pt}_{\mathrm{B}G})
      }"{description}
    ]
    \ar[
      dr,
      "\mathrm{id}\;\;\;"{left, near end}
    ]
    &&
    (\mathbf{B}\Gamma)
      \times
    (\mathbf{B}G)
    \ar[
      dl,
      "\;\;\;\;\;\mathrm{pr}_1"{right, near end}
    ]
    \\
    &
    \mathbf{B}\Gamma
  \end{tikzcd}
  {\phantom{AAAAAAA}}
  \begin{tikzcd}[column sep=large]
    X
    \ar[rr]
    \ar[
      d,
      "c"{left}
    ]
    \ar[
      drr,
      phantom,
      "\mbox{\tiny\rm(pb)}"
    ]
    &&
    X \!\sslash\! G
    \ar[
      d,
      "{
        (c \sslash G, \rho_G)
      }"{right}
    ]
    \\
    \mathbf{B}\Gamma
    \ar[
      rr,
      "{
        (\mathrm{id}_{\mathbf{B}\Gamma}, \mathrm{pt}_{\mathrm{B}G})
      }"{description}
    ]
    \ar[
      d
    ]
    \ar[
      drr,
      phantom,
      "\mbox{\tiny\rm(pb)}"
    ]
    &&
    (\mathbf{B}\Gamma)
      \times
    (\mathbf{B}G)
    \ar[
      d,
      "\mathrm{pr}_2"{right}
    ]
    \\
    \ast
    \ar[rr]
    &&
    \mathbf{B}G
  \end{tikzcd}
$$

\vspace{-2mm}
\noindent By forming the pasting composite shown on the right, this
follows by the pasting law \eqref{HomotopyPastingLaw}.
\end{proof}

\begin{example}[Equivariant bundle gerbes]
  \label{EquivariantBundleGerbes}
  Let $\Topos \,=\, \SmoothInfinityGroupoids$ (Ntn. \ref{SmoothInfinityGroupoids}).

  \vspace{1mm}
  \noindent
  {\bf (i)}
  If we regard the circle Lie group
  $\CircleGroup \,\in\,\Groups(\SmoothManifolds) \xhookrightarrow{\;}
  \Groups(\Topos)\,$ as equipped with trivial equivariance action,
  and consider its delooping 2-group
  $\Gamma \,=\,\mathbf{B}\CircleGroup \,\in\, \Groups(\Topos)$ (Prop. \ref{LoopingAndDeloopingEquivalence}), then
  the equivariant $\mathbf{B}\CircleGroup$-principal $\infty$-bundles
  -- in the sense of Def. \ref{GEquivariantGammaPrincipalBundles}, Prop. \ref{EquivariantBundlesForTrivialActionOnStructureGroup}
  --   are equivalently ``equivariant bundle gerbes''
  \cite{NikolausSchweigert10}\cite{MurrayRobertsStevensonVozzo15}
  (with several equivalent precursors, e.g. \cite[\S 5]{GomiTerashima09},
  and
  generalizing earlier discussion in \cite{MathaiStevenson03}\cite{Meinrenken02}).
  This follows from an analysis over action-{\v C}ech groupoids
  (Ex. \ref{CechActionGroupoidOfEquivariantGoodOpenCoverIsLocalCofibrantResolution} below)
  analogous to that in the proof of Prop. \ref{EquivariantPrincipalBundlesFromCechCocycles}
  below, see around \cite[Obs. 1.2.71]{dcct} for more.

  \vspace{1mm}
  \noindent
  {\bf (ii)} If instead we consider
  as structure 2-group the {\it discrete} underlying 2-group
  $\Gamma \,=\, \flat \mathbf{B}\CircleGroup$ \eqref{FlatModalityComputesUnderlyingPoints},
  then equivariant $\flat \mathbf{B}\CircleGroup$-principal $\infty$-bundles
  are equivalent to {\it flat} equivariant bundle gerbes,
  classified by ``discrete torsion''
  \cite{Sharpe03}\cite[\S 5]{GomiTerashima09}.

  \vspace{1mm}
  \noindent
  {\bf (iii)}
  If we regard $\ZTwo \acts \, \CircleGroup$
  as a $\ZTwo$-equivariant group with respect to the complex conjugation action,
  then $G = \ZTwo$-equivariant $\ZTwo \acts \, \mathbf{B}\CircleGroup$-principal
  $\infty$-bundles are equivalently ``Jandl gerbes''
  \cite{SSW07}\cite{GSW11} or ``real bundle gerbes'' \cite{HMSV16}.

  \noindent
  {\bf (iv)}
  Generally, if the equivariance group is given as an extension of $\ZTwo$
  \vspace{-2mm}
  $$
    \begin{tikzcd}[column sep=small]
      1
      \ar[r]
      &
      N
      \ar[r, hook]
      &
      G
      \ar[r, ->>, "p"]
      &
      \ZTwo
      \ar[r]
      &
      1
    \end{tikzcd}
  $$

  \vspace{-2mm}
  \noindent
  and we let $G$ act on $\mathbf{B}\CircleGroup$ through $p$
  (this is the situation of ``orbi-orientifolds'', see \cite[(67)]{SS19TadpoleCancellation} ),
  then $G$-equivariant $\ZTwo \acts \, \mathbf{B}\CircleGroup$-principal
  $\infty$-bundles are the equivariant generalization of Jandl gerbes.

  \noindent
  {\bf (iv)}
  Notice that the cohesive shape
  of this structure 2-group coincides with that of the
  infinite projective unitary 1-group (Ex. \ref{ProjectiveUnitarGroupOnAHilbertSpace}):  \vspace{-2mm}
  $$
    \begin{tikzcd}
      \mathbf{B}\CircleGroup
      \ar[r, "\scalebox{.9}{$ \eta^{\scalebox{.7}{$\shape$}}_{{}_{\mathbf{B}\CircleGroup(1)}} $}"]
      &
      B^2 \mathbb{Z}
      &
      \PUH \,.
      \ar[l, "\scalebox{.9}{$ \eta^{\scalebox{.7}{$\shape$}}_{{}_{\PUH}} $}"{swap}]
    \end{tikzcd}
  $$
  It is in this second guise that equivariant bundle gerbes (Jandl gerbes) model the
  general 3-twists of equivariant K-theory ($\mathrm{KR}$-theory); see
  Ex. \ref{EquivariantBundlesServingAsGeoemtricTwistsOfEquivariantKTheory} and
  Rem. \ref{GeometricTwistsOfEquivariantKTheoryAsBFieldsOnOrbiOrientifolds} below.
\end{example}

\begin{example}[Higher Jandl gerbes]
  \label{HigherJandlGerbes}
  Proceeding as in Ex. \ref{EquivariantBundleGerbes} for any $p \,\in\, \mathbb{N}$,

  \noindent
  {\bf (i)}
  the $\ZTwo$-equivariant $(p+1)$-groups which are the higher deloopings
  $\ZTwo \acts \, \Gamma \,=\, \ZTwo \acts \, \mathbf{B}^p \CircleGroup$
  of the circle group
  are the structure groups of  $\ZTwo$-equivariant $\mathbf{B}^p \CircleGroup(1)$-principal
  bundles which are, equivalently, ``higher Jandl gerbes''
  according to \cite[\S 4.4]{FSS15CField}; see also
  \cite[Ex. 2.12]{FSS20CharacterMap} for further pointers on higher bundle gerbes;

  \noindent
  {\bf (ii)} for geometrically discrete equivariant structure 2-group
  $\ZTwo \acts \,  \Gamma \,=\, \ZTwo \acts \, \flat \mathbf{B}^2 \CircleGroup$
  \eqref{FlatModalityComputesUnderlyingPoints}, these equivariant 2-bundles
  are classified by
  higher ``discrete torsion'', in the sense of \cite{Sharpe03HigherDiscreteTorsion}.
\end{example}

\medskip

\noindent {\bf Equivariant bundles under base change of the equivariance group.}
We discuss the behaviour equivariant bundles under base change of the equivariance group

-- along an inclusion: Rem. \ref{GloballyEquivariantNatureOfEquivariantPrincipalBundles};

-- along a surjection: Prop. \ref{EquivariantBundlesUnderExtensionsPfTheEquivarianceGroup}.

\begin{remark}[Globally equivariant nature of equivariant principal bundles]
  \label{GloballyEquivariantNatureOfEquivariantPrincipalBundles}
  Given $G \,\in\, \Groups(\Sets)$ and a subgroup inclusion
  $K \xhookrightarrow{\;\;i\;\;} G$,
  the induced left base change adjunction
  $(B i)_! \dashv (B i)^\ast$ (Prop. \ref{RestrictedAndLeftInducedInfinityActions})
  says, for $G$-equivariant structure groups $G \acts \, \Gamma$
  (Def. \ref{GGroups}),
  that

\noindent {\bf (i)}   K-equivariant $\Gamma$-principal bundles
  (Def. \ref{GEquivariantGammaPrincipalBundles})
  on some $K \acts \, X$
  are equivalently $G$-equivariant $\Gamma$-principal bundles on
  $X \times_K G$:
  \vspace{-2mm}
  $$
  \hspace{-5mm}
  \begin{array}{ccc}
    \EquivariantPrincipalBundles{K}{\Gamma}(\Topos)_{X}
    &\simeq&
    \qquad
    \EquivariantPrincipalBundles{G}{\Gamma}(\Topos)_{X \times_K G}
    \\
       \begin{tikzcd}[row sep=small, column sep=1pt]
      \mathllap{
        (B i)_!
        \left(
          \HomotopyQuotient
            { X }
            { K }
        \right)
        \;\simeq\;\;
      }
      \HomotopyQuotient
        { (X \times_K G) }
        { G }
      \ar[rr]
      \ar[dr]
      &&
      \HomotopyQuotient
        { \mathbf{B}\Gamma }
        { G }
      \ar[dl]
      \\
      &
      B G
    \end{tikzcd}
    \phantom{AA}
    &
    \xleftrightarrow{\;\;\;}
    &
    \begin{tikzcd}[row sep=small, column sep=1pt]
      \HomotopyQuotient
        { X }
        { K }
      \ar[rr]
      \ar[dr]
      &&
      \HomotopyQuotient
        { \mathbf{B}\Gamma }
        { K }
      \ar[dl]
      \mathrlap{
        \;\;\simeq\;
        (B i)^\ast
        \left(
          \HomotopyQuotient
            { \mathbf{B}\Gamma }
            { G }
        \right)\,.
      }
      \\
      &
      B K
          \end{tikzcd}
    \end{array}
  $$

  \vspace{-2mm}
\noindent {\bf (ii)}
 In particular, when the equivariance group action on the structure
  group $\Gamma$ is trivial, so that (by Ex. \ref{EquivariantInfinityGroupForTrivialGAction})
    \vspace{-2mm}
  \begin{align*}
    (B i)^\ast
    \HomotopyQuotient{ \mathbf{B}\Gamma }{ G }
    \;\; &\simeq\;\;
    (B i)^\ast
    \left(
      \mathbf{B}\Gamma
      \,\times\,
      \mathbf{B}G
    \right)
    \\
    \;\;&\simeq\;\;
    \mathbf{B}\Gamma
    \,\times\,
    \mathbf{B}K
    \\
    \;\;&\simeq\;\;
    \HomotopyQuotient
      { \mathbf{B}\Gamma }
      { K }
  \end{align*}

    \vspace{-2mm}
\noindent
  holds for {\it every} pair of groups and monomorphism $K \xhookrightarrow{\;i\;} G$,
  then
    \vspace{-2mm}
  $$
  \hspace{-5mm}
  \begin{array}{ccc}
    \EquivariantPrincipalBundles{K}{\Gamma}(\Topos)_{X}
    &\simeq&
\qquad    \EquivariantPrincipalBundles{G}{\Gamma}(\Topos)_{X \times_K G}
    \\
       \begin{tikzcd}[row sep=small, column sep=1pt]
      \mathllap{
        (B i)_!
        \left(
          \HomotopyQuotient
            { X }
            { K }
        \right)
        \;\simeq\;\;
      }
      \HomotopyQuotient
        { (X \times_K G) }
        { G }
      \ar[rr]
      \ar[dr]
      &&
      { \mathbf{B}\Gamma }
      \times
      \mathbf{B}G
      \ar[dl, "\mathrm{pr}_2"]
      \\
      &
      B G
    \end{tikzcd}
    \phantom{AA}
    &
    \xleftrightarrow{\;\;\;}
    &
    \begin{tikzcd}[row sep=small, column sep=1pt]
      \HomotopyQuotient
        { X }
        { K }
      \ar[rr]
      \ar[dr]
      &&
      { \mathbf{B}\Gamma }
      \times \mathbf{B}K
      \ar[dl, "\mathrm{pr}_2"]
      \mathrlap{
        \;\;\simeq\;
        (B i)^\ast
        \left(
          { \mathbf{B}\Gamma }
          \times
          \mathbf{B}G
        \right)\,.
      }
      \\
      &
      B K
          \end{tikzcd}
    \\
        \begin{tikzcd}
      \HomotopyQuotient
        { ( X \times_K G) }
        { G }
      \ar[r]
      &
      \mathbf{B}\Gamma
    \end{tikzcd}
    &
    \xleftrightarrow{\;\;\;}
    &
    \begin{tikzcd}
      \HomotopyQuotient
        { X }
        { K }
      \ar[r]
      &
      \mathbf{B}\Gamma
    \end{tikzcd}
    \end{array}
  $$

    \vspace{-2mm}
\noindent
  reflects the fact (Prop. \ref{EquivariantBundlesForTrivialActionOnStructureGroup})
  that equivariant $\Gamma$-principal bundles
  only depend on the homotopy quotient stack of their domain,
  and that $\HomotopyQuotient{ (X \times_K G) }{G} \,\simeq\, \HomotopyQuotient{ X }{K}$
  \eqref{NatureOfLeftInducedInfinityAction}.

\noindent {\bf (iii)}   This state of affairs is
  (as previously highlighted in \cite[\S 1.3]{Rezk14}\cite[p. xi]{Schwede18})
  the
  avatar
  of the incarnation of
  equivariant structures
  in {\it globally equivariant} homotopy theory,
  which we discuss as such in
  Prop. \ref{GloballyProperEquivariantClassifyingShape} below.
\end{remark}

\begin{proposition}[Equivariant bundles under extensions of the equivariance group]
\label{EquivariantBundlesUnderExtensionsPfTheEquivarianceGroup}
For $\!p : \!\!
\begin{tikzcd}\widehat{G} \ar[r,->>,"p"] & G\end{tikzcd}\!\!$ a surjection of equivariance groups
\eqref{ADiscreteGroupEpimorphismInAnInfinityTopos},
there is for

-- 1-truncated $G \acts \, X \,\in\, \Actions{G}(\ModalTopos{1})$ and

-- 0-truncated $G \acts \, \Gamma \,\, \Actions{G}\left(\Groups(\ModalTopos{0})\right)$,

\noindent
a natural full inclusion
\vspace{-2mm}
$$
  \begin{tikzcd}
    \EquivariantPrincipalBundles{G}{\Gamma}(\Topos)_X
    \ar[r, hook, "{ (B p)^\ast }"]
    &
    \EquivariantPrincipalBundles{\widehat{G}}{\Gamma}(\Topos)_X
  \end{tikzcd}
$$

\vspace{-2mm}
\noindent of the 1-groupoid $G$-equivariant
into that of $\widehat{G}$-equivariant
$\Gamma$-principal bundles on $X$, where $\widehat{G}$
acts on $\Gamma$ and on $X$ through $p$.
\end{proposition}
\begin{proof}
  Under the given assumption, both $G \acts \, X$ as well as
  $G \acts \, \mathbf{B}\Gamma$ are 1-truncated. Therefore, the statement
  follows (recalling Def. \ref{GEquivariantGammaPrincipalBundles})
  as an immediate consequence of
  \eqref{BaseChangeOfTwoActionsAlongSurjectiveHomomorphism}
  in Prop. \ref{PullbackOfActionsAlongSurjectiveGroupHomomorphismsIsFullyFaithful}.
\end{proof}

\medskip

\noindent
{\bf Conditions for classification of equivariant principal $\infty$-bundles.}
We proceed to generalize the cohesive classification result for
principal $\infty$-bundles
(Thm. \ref{ClassificationOfPrincipalBundlesAmongPrincipalInfinityBundles})
to equivariant principal $\infty$-bundles.
The crucial step in the plain proof was the smooth Oka principle
(Thm. \ref{SmoothOkaPrinciple}) and our proof strategy here is to
reduce the equivariant case to this plain case by ``blowing-up''
all orbifold singularities to smooth spherical space forms (Def. \ref{SmoothSphericalSpaceForm})
after stabilizing the dimension.
This requires assuming that
\vspace{-2mm}
\begin{itemize}
\setlength\itemsep{-2pt}
\item[{\bf (1)}]
$G$-orbisingularities are resolvable
(via Ntn. \ref{ResolvableOrbiSingularities})
and

\item[{\bf (2)}] the equivariant bundle do not detect the difference between
  an orbi-singularity and its resolution, which is
  guaranteed by

  \vspace{-.3cm}
  \begin{itemize}
    \setlength\itemsep{-2pt}
    \item[{\bf (i)}]
     a truncation condition
     on the structure group, Ntn. \ref{CohesiveGroupsWithTruncatedClassifyingShape} below,

    \item[{\bf (ii)}]
     a ``stability'' condition on the bundles, Ntn. \ref{StableEquivariantBundles} below.
  \end{itemize}

\end{itemize}

\vspace{-2mm}
\noindent The first condition is satisfied, for instance, by all finite subgroups of $\SUTwo$,
but not by all finite subgroups of $\SOThree$ (by Ex. \ref{ADEGroupsHaveSphericalSpaceForms}).
The second is mainly a truncation condition which, while rarely satisfied by
compact Lie groups,
 is the default assumption for twists in
twisted equivariant Whitehead-generalized cohomology theory, where it just means that the
twist has a bounded degree. The archetypical example of this are the
degree-3 twists of equivariant K-theory, whose classification result
is recovered by the following result (discussed in Ex. \ref{EquivariantBundlesServingAsGeoemtricTwistsOfEquivariantKTheory}).

\begin{notation}[Cohesive groups with truncated classifying shape]
  \label{CohesiveGroupsWithTruncatedClassifyingShape}
  We say that a 0-truncated smooth group
  $\Gamma \,\in\, \Groups(\SmoothInfinityGroupoids)_{\leq 0}$,
  $\Truncation{0} \, \Gamma \,\simeq\, \Gamma$,
  has {\it truncated classifying shape}
  if

  \noindent
  {\bf (i)} its shape is truncated:
    \begin{equation}
      \label{TruncationOfShapeOfGroupWithTruncatedShape}
      \underset{ n \in \mathbb{N} }{\exists}
      \;\;\;
      \Truncation{n}
      \,
      \shape
      \,
      \Gamma
      \;\simeq\;
      \shape
      \,
      \Gamma
      ;
    \end{equation}

  \noindent
  {\bf (ii)}
   isomorphism classes of $\Gamma$-principal bundles
   over smooth manifolds are already concordance classes \eqref{ProjectionFromIsomorphismClassesToConcordanceClasses}:
   \vspace{-2mm}
   $$
     \underset{
       { \SmoothManifold \,\in\,}
       \atop
       { \SmoothManifolds }
     }{\forall}
     \;\;
     \begin{tikzcd}
       \IsomorphismClasses{
         \PrincipalBundles{\Gamma}(\SmoothInfinityGroupoids)_\SmoothManifold
       }
       \ar[r, ->>,  "\sim"{above}]
       &
       \ConcordanceClasses{
         \PrincipalBundles{\Gamma}(\SmoothInfinityGroupoids)_\SmoothManifold
       }.
     \end{tikzcd}
   $$

\noindent
This second condition
is equivalent,
by Prop. \ref{ClassificationOfPrincipalInfinityBundles},
to
$B \Gamma$ \eqref{ClassifyingSpaceIsShapeOfModuliStack}
being a classifying space for
isomorphism classes of $\Gamma$-principal bundles over smooth manifolds:
\vspace{-2mm}
\begin{equation}
  \label{ForSmoothGroupsWithTruncatedClassifyingSpapeTheClassifyingSpaceIsClassifying}
   \underset{
     { \SmoothManifold \,\in\,}
     \atop
     { \SmoothManifolds }
   }{\forall}
   \;\;
  \begin{tikzcd}
     \IsomorphismClasses{
       \PrincipalBundles{\Gamma}(\SmoothInfinityGroupoids)_\SmoothManifold
     }
     \ar[r, "\sim"]
     &
     \Truncation{0}
     \,
     \Maps{}
       { \shape \, \SmoothManifold }
       { B \Gamma }
     \;\simeq\;
     H^1(X;\, \shape \, \Gamma)
     \,.
  \end{tikzcd}
\end{equation}

\vspace{-1mm}
\noindent
Notice that, by possibly increasing any truncation degree $n$ in \eqref{TruncationOfShapeOfGroupWithTruncatedShape}
a little,
Prop. \ref{ExistenceOfSmoothSphericaSpaceForms} guarantees that
we may always assume, just for notational definiteness,
  that
  $G$ has a smooth free action on some smooth $n+2$-sphere
  $\SmoothSphere{\,n+2}$ (possibly exotic):
  \vspace{-2mm}
\begin{equation}
  \label{AdjustedTruncationOfShapeOfGroupWithTruncatedShape}
  \underset{n \in \mathbb{N}}{\exists}
  \;\;\;\;\;\;
  G \acts \; \SmoothSphere{\,n+2}
  \,\in\,
  \Actions{G}(\SmoothManifolds)_{\mathrm{free}}
  \;\;\;
  \mbox{and}
  \;\;\;
  \Truncation{n}
  \,
  \shape
  \,
  \Gamma
  \;\simeq\;
  \shape
  \,
  \Gamma
  \,.
\end{equation}

\vspace{-1mm}
\noindent
We shall say for short that $\Gamma$ has {\it $n$-truncated classifying shape}
to indicate such a compatible choice of truncation degree.
\end{notation}

\begin{lemma}[Principal bundles on blow-ups for structure group with truncated classifying shape]
  \label{PrincipalBundlesOnBlowUps}
  Let $G \,\in\, \Groups(\FiniteSets)_{\resolvable}$ (Ntn. \ref{ResolvableOrbiSingularities})
  and
  $\Gamma \,\in\, \Groups(\SmoothInfinityGroupoids)_{\leq 0}$
  with $n$-truncated classifying shape (Ntn. \ref{CohesiveGroupsWithTruncatedClassifyingShape}).
  Then, for $G \acts \, \SmoothManifold \,\in\, \Actions{G}(\SmoothManifolds)$,
  pullback along the projection
  $$
    \begin{tikzcd}
      \HomotopyQuotient
        { \SmoothManifold }
        { G }
      &&
      \HomotopyQuotient
        { ( \SmoothManifold \times \SmoothSphere{\,n+2} )  }
        { G }
      \;\simeq\;
      \Quotient
        { ( \SmoothManifold \times \SmoothSphere{\,n+2} ) }
        { G }
      \ar[
        ll,
        "{
          p_{{}_{\SmoothManifold}}
          \,:=\,
          \HomotopyQuotient
            { \mathrm{pr}_1 }
            { G }
        }"{swap}
      ]
    \end{tikzcd}
  $$
  from the product with a compatible $G$-sphere
  \eqref{AdjustedTruncationOfShapeOfGroupWithTruncatedShape}
  is a natural surjection
  \vspace{-1mm}
  \begin{equation}
    \label{PullbackOfEquivariantBundleToProductWithFreeGSphereIsSurjection}
    \begin{tikzcd}
      \IsomorphismClasses
      {
        \EquivariantPrincipalBundles{G}{\Gamma}(\SmoothInfinityGroupoids)
          _{
            \SmoothManifold
          }
      }
      \ar[
        r,
        ->>,
        "{
          p_{{}_{\SmoothManifold}}^\ast
        }"
      ]
      &
      \IsomorphismClasses
      {
        \PrincipalBundles{\Gamma}(\SmoothInfinityGroupoids)
        _{
          (\SmoothManifold \times \SmoothSphere{\,n+2})/G
        }
      }
    \end{tikzcd}
  \end{equation}

  \vspace{-1mm}
  \noindent
  from
  isomorphism classes of $G$-equivariant $\Gamma$-principal bundles on
  $\SmoothManifold$ \eqref{EquivariantPrincipalBundlesAsPrincipalBundlesOnQuotientStack}
  onto the set of isomorphism classes of plain $\Gamma$-principal bundles
  on $\Quotient{ (\SmoothManifold \times \SmoothSphere{\,n+2}) }{G}$.
\end{lemma}
\begin{proof}
In terms of cocycles we need to show that
$
  \begin{tikzcd}
    \Truncation{0}
    \,
    \PointsMaps{}
      {
        \HomotopyQuotient
          { \SmoothManifold }
          { G }
      }
      { \mathbf{B}\Gamma }
    \ar[
      r,
      ->>,
      "{
        p_{{}_{\SmoothManifold}}^\ast
      }"
    ]
    &
    \Truncation{0}
    \,
    \PointsMaps{\big}
      {
        ( \SmoothManifold \times \SmoothSphere{\,n+2} )/G
      }
      { \mathbf{B}\Gamma }
    \,.
  \end{tikzcd}
$
To see this, consider the homotopy colimit diagram that exhibits the $G$-quotients:

\vspace{-.1cm}
\begin{equation}
  \label{DescentToHomotopyQuotientOfManifoldAndQuotientOfItsProductWithSphere}
  \begin{tikzcd}[row sep=40pt, column sep=30pt]
    (\SmoothManifold \times \SmoothSphere{\,n+2}) \times G \times G
    \ar[d, shift left=10pt, "\partial_0"{description}]
    \ar[d, shift left=0pt, "\partial_1"{description}]
    \ar[d, shift right=10pt, "\partial_2"{description}]
    \ar[rr, "p"]
    &&
    \SmoothManifold \times G \times G
    \ar[d, shift left=10pt, "\partial_0"{description}]
    \ar[d, shift left=0pt, "\partial_1"{description}]
    \ar[d, shift right=10pt, "\partial_2"{description}]
    \\
    (X \times \SmoothSphere{\,n+2}) \times G
    \ar[d, shift left=5pt, "\partial_0"{description}]
    \ar[d, shift right=5pt, "\partial_1"{description}]
    \ar[rr, "p"]
    &&
    \SmoothManifold \times G
    \ar[d, shift left=5pt, "\partial_0"{description}]
    \ar[d, shift right=5pt, "\partial_1"{description}]
    \ar[
      dddr,
      start anchor={[xshift=4pt]},
      end anchor={[xshift=4pt]},
      bend right=4pt,
      "\partial_1^\ast c'_0"{description},
      "{ }"{name=PartialOneAstCPrimeOneAbove, above}
    ]
    \ar[
      dddr,
      start anchor={[xshift=5pt]},
      end anchor={[xshift=5pt]},
      bend left=40,
      "\partial_0^\ast c'_0"{description},
      "{ }"{name=PartialZeroAstCPrimeOneBelow, below}
    ]
    \ar[
      from=PartialOneAstCPrimeOneAbove,
      to=PartialZeroAstCPrimeOneBelow,
      Rightarrow,
      shorten=7pt,
      "\gamma'_1"{above},
      "\sim"{sloped, below}
    ]
    \\
    \SmoothManifold \times \SmoothSphere{\,d}
    \ar[d, "q"]
    \ar[rr, "p"]
    &&
    \SmoothManifold
    \ar[d, "q"]
    \ar[
      ddr,
      "c'_0"{description, pos=.6},
      "{ }"{name=CPrimeOneBelow, below, pos=.6}
    ]
    \\
    ( \SmoothManifold \times \SmoothSphere{\,n+2} )/G
    \ar[rr, "p"]
    \ar[
      drrr,
      "c"{description, pos=.7},
      "{ }"{name=c, above, pos=.7},
      bend right=10pt,
      crossing over
    ]
    &&
    \HomotopyQuotient
      { \SmoothManifold }
      { G }
    \ar[
      dr,
      "c'"{description, pos=.4},
      "{ }"{name=cprime, above, pos=.4},
      "{ }"{name=cprimebelow, below, pos=.4}
    ]
    \ar[
      from=c,
      to=cprimebelow,
      shorten=2pt,
      Rightarrow,
      "{\sim}"{sloped, below}
    ]
    \\[-10pt]
    &&
    &[+10pt]
    \mathbf{B}\Gamma
    \ar[
      from=CPrimeOneBelow,
      to=cprime,
      Rightarrow,
      shorten=1pt,
      "\sim"{sloped, below}
    ]
  \end{tikzcd}
\end{equation}
\vspace{-.3cm}

\noindent
Hence given any cocycle
$
  c
  \,\colon
  (\SmoothManifold \times \SmoothSphere{\,n+2} )/G
  \xrightarrow{\;}
  \mathbf{B}\Gamma
$
we need to show that it is isomorphic to one in the image of $p^\ast$.
But by the universal property of the
middle vertical colimit in \eqref{DescentToHomotopyQuotientOfManifoldAndQuotientOfItsProductWithSphere},
such cocycles $c$ on the quotient
are equivalent to cocycles
$
  c_0
    \,\colon\,
  \SmoothManifold
    \times
  \SmoothSphere{\,n+2}
    \xrightarrow{\;}
  \mathbf{B}\mathcal{G}
$
on the covering space
which are equipped with {\it descent data} in the form of
an isomorphism
$
  \begin{tikzcd}
    \partial_1^\ast c_0
    \ar[r, "\sim"{below}, "\gamma_1"{above}]
    &
    \partial_0^\ast c_0
  \end{tikzcd}
$
such that\footnote{
  We are displaying the standard form of the descent condition.
  If one unwinds the universal property of the homotopy colimit in
  \eqref{DescentToHomotopyQuotientOfManifoldAndQuotientOfItsProductWithSphere}
  one finds the diagram \eqref{DescentConditionForOneCocyclesAlongGQuotientCoprojection}
  with another cocycle
  $c_1 \,:\, \SmoothManifold \times \SmoothSphere{\,n+2} \times G \xrightarrow{\;} \mathbf{B}\Gamma$
  at its center and equipped with isomorphisms to all three corners.
  These three isomorphisms (which are not shown)
  clearly carry the same information
  as the three that are shown
  in \eqref{DescentConditionForOneCocyclesAlongGQuotientCoprojection},
  up to gauge equivalence of descent data.
}

\vspace{-.1cm}
\begin{equation}
  \label{DescentConditionForOneCocyclesAlongGQuotientCoprojection}
  \begin{tikzcd}[column sep=10pt, row sep=25pt]
    &[-20pt]
    {}
    &[-20pt]
    \partial_2^\ast \partial_0^\ast c_0
    \ar[r,-,shift left=1pt]
    \ar[r,-,shift right=1pt]
    &[-0pt]
    \partial_0^\ast \partial_1^\ast c_0
    \ar[drr, "\sim"{sloped, below}, "\partial_0^\ast \gamma_1"]
    \\
    \partial_2^\ast \partial_1^\ast c_0
    \ar[urr, "\sim"{sloped, below}, "\partial_2^\ast \gamma_1"]
    \ar[dr,-,shift left=1pt]
    \ar[dr,-,shift right=1pt]
    &
    &&&
    &[-20pt]
    \partial_0^\ast \partial_0^\ast c_0
    \ar[dl, -, shift left=1pt]
    \ar[dl, -, shift right=1pt]
    \\[-20pt]
    &
    \partial_1^\ast \partial_1^\ast c_0
    \ar[rrr, "\sim"{sloped, below}, "\partial_1^\ast \gamma_1"]
    &&&
    \partial_1^\ast \partial_0^\ast c_0
  \end{tikzcd}
\end{equation}
Now, by the assumption that $\Gamma$ has truncated classifying shape,
the cocycles on $\SmoothManifold \times \SmoothSphere{\,n+2}$
must be trivial along the sphere, in that we
we have the following natural bijections:
\vspace{-1mm}
$$
  \def\arraystretch{1.5}
  \begin{array}{lll}
    \IsomorphismClasses
    {
      \PrincipalBundles{\Gamma}(\SmoothInfinityGroupoids)
        _{\SmoothManifold \times \SmoothSphere{\,n+2}}
    }
    &
    \;\simeq\;
    \Truncation{0}
    \,
    \Maps{\big}
      {
        \Shape
        (
          \SmoothManifold
          \times
          \SmoothSphere{\,n+2}
        )
      }
      { B \Gamma }
    &
    \proofstep{ by \eqref{ForSmoothGroupsWithTruncatedClassifyingSpapeTheClassifyingSpaceIsClassifying} }
    \\
    &
    \;\simeq\;
    \Truncation{0}
    \,
    \Maps{\big}
      {
        \ShapeOfSphere{n+2}
        \,\times\,
        \Shape \,\SmoothManifold
      }
      { B \Gamma }
    &
    \proofstep{ by \eqref{ShapePreservesBinaryProducts} }
    \\
    &
    \;\simeq\;
    \Truncation{0}
    \,
    \Maps{\big}
      {
        \Shape \, \SmoothManifold
      }
      {
        \Maps{}
          { \ShapeOfSphere{n+2} }
          { B \Gamma }
      }
    &
    \proofstep{ by Lem. \ref{InternalHomAdjointness} }
    \\
    &
    \;\simeq\;
    \Truncation{0}
    \,
    \Maps{}
      {
        \Shape \SmoothManifold
      }
      {
        B \Gamma
      }
    &
    \proofstep{ by Lem. \ref{TruncationAndnFoldFreeLoopSpaces} using \eqref{AdjustedTruncationOfShapeOfGroupWithTruncatedShape} }
    \\
    &
    \;\simeq\;
    \IsomorphismClasses
    {
      \PrincipalBundles{\Gamma}(\SmoothInfinityGroupoids)_{\SmoothManifold}
    }
    &
    \proofstep{ by \eqref{ForSmoothGroupsWithTruncatedClassifyingSpapeTheClassifyingSpaceIsClassifying}. }
  \end{array}
$$

\vspace{-1mm}
\noindent
This implies that the required isomorphism to a cocycle in the image of $p^\ast$
exists:
  \vspace{-2mm}
$$
  c'_0 \,\colon\, X  \xrightarrow{\;} \mathbf{B}\mathcal{G}
  \,,
  \;\;\;\;\;
  \begin{tikzcd}
    c_0
    \ar[r, "\sim"{below}, "t_1"{above}]
    &
    p^\ast c'_0
    \,.
  \end{tikzcd}
$$

\vspace{-2mm}
\noindent
Via this isomorphism, the above equivariant structure $\gamma$ on $c_0$
is transported to an equivariant structure $\gamma'$ on $c'_0$:
\vspace{-2mm}
\begin{equation}
  \label{TransportedDescentDataToCocycle}
  \begin{tikzcd}
    i^\ast \partial_0^\ast c_1
    \ar[d, "\sim"{sloped, below}, "i^\ast \gamma"{right}]
    \ar[rr, "\sim"{below}, "i^\ast \partial_0^\ast t_1"{above}]
    &&
    i^\ast \partial_0^\ast p^\ast c'_1
    \ar[r, -, shift right=1pt]
    \ar[r, -, shift left=1pt]
    &
    i^\ast p^\ast \partial_0^\ast  c'_1
    \ar[r, -, shift right=1pt]
    \ar[r, -, shift left=1pt]
    &
    \partial_0^\ast  c'_1
    \ar[d, dashed, "\gamma'_1"]
    \\
    i^\ast \partial_1^\ast c_1
    \ar[rr, "\sim"{below}, "i^\ast \partial_1^\ast t_1"{above}]
    &&
    i^\ast \partial_1^\ast p^\ast c'_1
    \ar[r, -, shift right=1pt]
    \ar[r, -, shift left=1pt]
    &
    i^\ast p^\ast \partial_1^\ast  c'_1
    \ar[r, -, shift right=1pt]
    \ar[r, -, shift left=1pt]
    &
    \partial_1^\ast  c'_1
    \,.
  \end{tikzcd}
\end{equation}

  \vspace{-1mm}
\noindent
By its conjugation nature, \eqref{TransportedDescentDataToCocycle}
clearly inherits the descent condition \eqref{DescentConditionForOneCocyclesAlongGQuotientCoprojection},
hence constitutes a co-cone under the right vertical homotopy colimiting diagram
in \eqref{DescentToHomotopyQuotientOfManifoldAndQuotientOfItsProductWithSphere}, as shown there.
Thus, by the universal property of homotopy colimits in the bottom row
of \eqref{DescentToHomotopyQuotientOfManifoldAndQuotientOfItsProductWithSphere} ,
this gives the required isomorphism
  \vspace{-2mm}
$$
  \begin{tikzcd}
  c
  \ar[r, "\sim"{below}, "t"{above}]
  &
  p^\ast c'
  \,.
  \end{tikzcd}
$$

  \vspace{-2mm}
\noindent
As this holds for every given $c$, $p^\ast$ is surjective.
\end{proof}

The point of the pullback \eqref{PullbackOfEquivariantBundleToProductWithFreeGSphereIsSurjection}
is that it gives a handle on the classification of $G$-equivariant $\Gamma$-principal bundles
on $\SmoothManifold$
through that of ordinary $\Gamma$-principal bundles on
$\Quotient{ (\SmoothManifold \times \SmoothSphere{\,n+2}) }{G}$:

\begin{lemma}[Classifying maps on blowups for structure group of truncated shape]
  \label{ConcordanceInfinityGroupoidOfPrincipalBundlesOnBlowupForStructureGroupOfTruncatedShape}
  Let $G \,\in\, \Groups(\FiniteSets)_{\resolvable}$ (Ntn. \ref{ResolvableOrbiSingularities})
  and
  $\Gamma \,\in\, \Groups(\SmoothInfinityGroupoids)_{\leq 0}$
  with $n$-truncated classifying shape (Ntn. \ref{CohesiveGroupsWithTruncatedClassifyingShape}).
  Then
  \vspace{-1mm}
  $$
    \Maps{\big}
      {
        \shape
        \,
        \Quotient{ (\SmoothManifold \times \SmoothSphere{\,n+2}) } {G}
      }
      {
        B \Gamma
      }
       \simeq
    \;\;
    \Maps{\big}
      {
        \shape
        \,
        \HomotopyQuotient{ \SmoothManifold }{G}
      }
      {
        B \Gamma
      }\;.
  $$
\end{lemma}
\begin{proof}
We compute as follows:
\vspace{-2mm}
$$
  \def\arraystretch{1.6}
  \begin{array}{lll}
        \Maps{\big}
      {
        \shape
        \,
        \Quotient{ (\SmoothManifold \times \SmoothSphere{\,n+2}) }{ G }
      }
      {
        B \Gamma
      }
        &
    \simeq\;
    \Maps{\big}
      {
        \shape
        \,
        \HomotopyQuotient{ (\SmoothManifold \times \SmoothSphere{\,n+2}) }{ G }
      }
      {
        B \Gamma
      }
    &
    \proofstep{ by Prop. \ref{HomotopyQuotientOfFreeActionsIsOrdinaryQuotientOfZeroTruncation} }
    \\
    &
    \simeq\;
    \Maps{\bigg}
      {
        \shape
        \,
        \colimit{ [n] \in \Delta^{\mathrm{op}} }
        \,
        (\SmoothManifold \times \SmoothSphere{\,n+2}) \times G^{\times_n}
      }
      {
        B \Gamma
      }
    &
    \proofstep{ by \eqref{HomotopyQuotientAsHomotopyColimit} }
    \\
    &
    \simeq\;
    \Maps{\bigg}
      {\,
        \colimit{ [n] \in \Delta^{\mathrm{op}} }
        \,
        \shape
        \,
        (\SmoothManifold \times \SmoothSphere{\,n+2}) \times G^{\times_n}
      }
      {
        B \Gamma
      }
    &
    \proofstep{ by \eqref{InfinityAdjointPreservesInfinityLimits} }
    \\
    &
    \;\simeq\;
    \limit{ [k] \in \Delta^{\mathrm{op}} }
    \,
    \Maps{\big}
      {
        \shape
        \,
        (\SmoothManifold \times \SmoothSphere{\,n+2}) \times G^{\times_k}
      }
      {
        B \Gamma
      }
    &
    \proofstep{ by Prop. \ref{MappingStackConstructionPreservesLimits} }
    \\
    &
    \simeq\;
    \limit{ [k] \in \Delta^{\mathrm{op}} }
    \,
    \Maps{\big}
      {
        (
         \shape
         \,
         \SmoothManifold
          \times
          G^{\times_k}
        )
        \times
        (
          \shape \SmoothSphere{\,n+2}
        )
      }
      {
        B \Gamma
      }
    &
    \proofstep{ by \eqref{ShapePreservesBinaryProducts} }
    \\
    &
    \simeq\;
    \limit{ [k] \in \Delta^{\mathrm{op}} }
    \,
    \Maps{\big}
    {
       \shape \,
       \SmoothManifold
        \times
        G^{\times_k}
    }
    {
      \Maps{}
        {
          \ShapeOfSphere{n+2}
        }
        {
          B \Gamma
        }
    }
    &
    \proofstep{ by Lem. \ref{InternalHomAdjointness} }
    \\
    &
    \simeq\;
    \limit{ [k] \in \Delta^{\mathrm{op}} }
    \,
    \Maps{\big}
    {
       \shape \,
       \SmoothManifold
        \times
        G^{\times_k}
    }
    {
      B \Gamma
    }
    &
    \proofstep{ by Lem. \ref{TruncationAndnFoldFreeLoopSpaces} using \eqref{AdjustedTruncationOfShapeOfGroupWithTruncatedShape} }
    \\
    &
    \simeq\;
    \Maps{\bigg}
    {\;
       \colimit{ [k] \in \Delta^{\mathrm{op}} }
       \,
       \shape \,
       \SmoothManifold
        \times
        G^{\times_k}
    }
    {
      B \Gamma
    }
    &
    \proofstep{ by Prop. \ref{MappingStackConstructionPreservesLimits} }
    \\
    &
    \simeq\;
    \Maps{}
    {
      \HomotopyQuotient
      {\,
        \shape \,
        \SmoothManifold
      }
      { G }
    }
    {
      B \Gamma
    }
   &
   \proofstep{ by \eqref{HomotopyQuotientAsHomotopyColimit}. }
  \end{array}
$$

\vspace{-7mm}
\end{proof}

While the pullback
operation \eqref{PullbackOfEquivariantBundleToProductWithFreeGSphereIsSurjection}
is always a surjection, it is often far from being injective.
The following terminology captures the idea of choosing a natural
sub-class of equivariant bundles that makes this pullback a bijection:

\begin{notation}[Blowup-stable equivariant principal bundles]
  \label{StableEquivariantBundles}
  Consider $\Gamma \,\in\, \Groups(\SmoothInfinityGroupoids)$
  of truncated classifying shape (Ntn. \ref{CohesiveGroupsWithTruncatedClassifyingShape}).

  \noindent
  {\bf(i)}
  We say that a natural system of sub-groupoids
  of $G$-equivariant $\Gamma$-principal bundles on smooth $G$-manifolds
  $\SmoothManifold$ with resolvable singularities
  (Ntn. \ref{ResolvableOrbiSingularities})
  \vspace{-2mm}
  \begin{equation}
    \label{InclusionOfStableEquivariantPrincipalBundles}
    \begin{tikzcd}
      \PointsMaps{}
        { \HomotopyQuotient{\SmoothManifold}{G} }
        { \mathbf{B}\Gamma }
      ^{\stable}
      \ar[r, hook]
      &
      \PointsMaps{}
        { \HomotopyQuotient{\SmoothManifold}{G} }
        { \mathbf{B}\Gamma }
    \end{tikzcd}
  \end{equation}

  \vspace{-2mm}
  \noindent
  is a system of {\it blowup-stable}
  equivariant bundles if restriction along this inclusion makes
  the surjections
  \eqref{PullbackOfEquivariantBundleToProductWithFreeGSphereIsSurjection}
  into bijections:
   \vspace{-2mm}
   \begin{equation}
    \label{ConditionOfCompatibleRestrictionsOnIsotropyGroups}
    \underset{
      \scalebox{.7}{$
        \begin{array}{c}
          G \in \Groups(\FiniteSets)_{\resolvable}
          \\
          G \acts \SmoothManifold \,\in\, \Actions{G}(\SmoothManifolds)
        \end{array}
      $}
    }{\forall}
    \;\;\;\;
    \begin{tikzcd}
    \Truncation{0}
    \,
    \PointsMaps{}
      { \HomotopyQuotient{\TopologicalSpace}{G} }
      { \mathbf{B}\Gamma }
    ^{\stable}
    \ar[
      r,
      hook
    ]
    \ar[
      rr,
      rounded corners,
      to path={
           ([yshift=+3pt]\tikztostart.south)
        -- ([yshift=-8pt]\tikztostart.south)
           node[swap, yshift=-6pt, xshift=+115pt]{\scalebox{.7}{$\sim$}}
        -- ([yshift=-7pt]\tikztotarget.south)
        -- ([yshift=-0pt]\tikztotarget.south)
      }
    ]
    &
    \Truncation{0}
    \,
    \PointsMaps{}
      { \HomotopyQuotient{\TopologicalSpace}{G} }
      { \mathbf{B}\Gamma }
    \ar[
      r,
      ->>,
      "{
        p_{{}_{\SmoothManifold}}^\ast
      }"
    ]
    &
    \Truncation{0}
    \,
    \PointsMaps{\big}
      {
        \Quotient
          { ( \TopologicalSpace \times \SmoothSphere{\,n+2} ) }
          { G }
      }
      { \mathbf{B}\Gamma }
    \mathrlap{\,.}
    \end{tikzcd}
   \end{equation}

  \noindent
  {\bf (ii)}
  By naturality of the inclusion
  \eqref{InclusionOfStableEquivariantPrincipalBundles} in
  the domain manifold, we have this inclusion also over products
  of the form
  $U \times \HomotopyQuotient{\SmoothManifold}{G}$ with
  $U \,\in\, \CartesianSpaces$ (Ntn. \ref{CartesianSpacesAndDiffeologicalSpaces})
  (regarded as equipped with the trivial $G$-action), so
  that a choice of stable equivariant bundles
  \eqref{ConditionOfCompatibleRestrictionsOnIsotropyGroups}
  induces a
  monomorphism of their moduli stacks:
   \vspace{-2mm}
  $$
    \begin{tikzcd}
      \Maps{}
        { \HomotopyQuotient{\SmoothManifold}{G} }
        { \mathbf{B}\Gamma }
      ^{\stable}
      \ar[r, hook]
      &
      \Maps{}
        { \HomotopyQuotient{\SmoothManifold}{G} }
        { \mathbf{B}\Gamma }
      \,,
    \end{tikzcd}
  $$

   \vspace{-2mm}
  \noindent  where, in view of \eqref{ValuesOfMappingStackAsHomSpaces}, we write
  \vspace{-2mm}
  \begin{equation}
    \label{ModuliStackOfStableEquivariantBundles}
    \Maps{}
      { \HomotopyQuotient{\SmoothManifold}{G} }
      { \mathbf{B}\Gamma }
    ^{\stable}
    \;\;\;
    \in
    \;
    \SmoothInfinityGroupoids
    \,,
 \qquad
    \Maps{}
      { \HomotopyQuotient{\SmoothManifold}{G} }
      { \mathbf{B}\Gamma }
    ^{\stable}
    (U)
    \;:=\;
    \PointsMaps{}
      { U \times \HomotopyQuotient{\SmoothManifold}{G} }
      { \mathbf{B}\Gamma }
    ^{\stable}
    \,.
  \end{equation}

  \noindent
  {\bf (iii)}
  If $\Gamma$ above is equipped with a $G$-action,
  and there is a notion of stable equivariant bundles
  for structure group $\Gamma \rtimes G$-bundles
  then we
  have also the corresponding subgroupoid of stable objects
  $$
  \def\arraystretch{1.3}
  \begin{array}{ll}
    \SlicePointsMaps{\big}{\mathbf{B}G}
      { \HomotopyQuotient{\SmoothManifold}{G} }
      { \HomotopyQuotient{\mathbf{B}\Gamma}{G} }
    ^{\stable}
     \\
    \;:=\;
    \PointsMaps{}
      { \HomotopyQuotient{\SmoothManifold}{G} }
      { \mathbf{B} (\Gamma \rtimes G) }
    ^{\stable}
    \quad
    \underset{
      \mathclap{\phantom{\vert^{\vert}}}
      \mathclap{
      \PointsMaps{}
        { \HomotopyQuotient{\SmoothManifold}{G} }
        { \mathbf{B}G }
      }
    }{\times}
    \quad \ast
    \end{array}
    \qquad
    \xhookrightarrow{\qquad}
    \quad
      \begin{array}{ll}
    \PointsMaps{}
      { \HomotopyQuotient{\SmoothManifold}{G} }
      { \mathbf{B} (\Gamma \rtimes G) }
      \quad
    \underset{
      \mathclap{\phantom{\vert^{\vert}}}
      \mathclap{
      \PointsMaps{}
        { \HomotopyQuotient{\SmoothManifold}{G} }
        { \mathbf{B}G }
      }
    }{\times}
    \quad
    \ast
    \;\;
     \\
    \overset{
      \mbox{
        \tiny
        \rm
        \eqref{HomSpaceInSliceAsFiberProduct}
      }
    }{
      \simeq
    }
    \SlicePointsMaps{\big}{\mathbf{B}G}
      { \HomotopyQuotient{\SmoothManifold}{G} }
      { \HomotopyQuotient{\mathbf{B}\Gamma}{G} }\,.
      \end{array}
      $$
  among all $G$-equivariant $G \acts \, \Gamma$-principal bundles
  according to \eqref{InfinityGroupoidOfEquivariantPrincipalBundles}.

  \noindent
  {\bf (iv)}
  Finally, in the case (iii)
  we obtain also the corresponding monomorphism of
  moduli stacks of equivariant bundles formed as slice mapping stacks
  (Def. \ref{SliceMappingStack}):
  $$
  \def\arraystretch{1.3}
  \begin{array}{ll}
    \SliceMaps{}{\mathbf{B}G}
      { \HomotopyQuotient{\SmoothManifold}{G} }
      { \HomotopyQuotient{\mathbf{B}\Gamma}{G} }
      \\
    \,:=\,
    \Maps{}
      { \HomotopyQuotient{\SmoothManifold}{G} }
      { \HomotopyQuotient{\mathbf{B}\Gamma}{G} }
    ^{\stable}
    \quad
    \underset
      {
        \mathclap{\phantom{\vert^{\vert}}}
        \mathclap{
        \Maps{}
          { \HomotopyQuotient{\SmoothManifold}{G} }
          { \mathbf{B}G }
        }
      }
      {\times}
    \quad \ast
    \end{array}
    \quad
    \xhookrightarrow{\qquad}
    \quad
    \begin{array}{ll}
    \Maps{}
      { \HomotopyQuotient{\SmoothManifold}{G} }
      { \HomotopyQuotient{\mathbf{B}\Gamma}{G} }
  \quad  \underset
      {
        \mathclap{\phantom{\vert^{\vert}}}
        \mathclap{
        \Maps{}
          { \HomotopyQuotient{\SmoothManifold}{G} }
          { \mathbf{B}G }
        }
      }
      {\times}
   \quad
   \ast
    \;
    \\
    \overset
      {
        \mbox{
          \tiny
          \eqref{TheSliceMappingStack}
        }
      }
      {=}
    \;
    \SliceMaps{}{\mathbf{B}G}
      { \HomotopyQuotient{\SmoothManifold}{G} }
      { \HomotopyQuotient{\mathbf{B}\Gamma}{G} }
    \,.
    \end{array}
  $$
\end{notation}

\begin{remark}
While
blowup-stability in Ntn. \ref{StableEquivariantBundles} may look like a strong
assumption to make, the condition can in fact be solved
nicely explicitly by ``averaging coboundaries of equivariant {\v C}ech cocycles
over blowup spheres'', at least in the following classes of examples:

\noindent {\bf (i)}  equivariant bundles with truncated compact Lie structure groups
 (Thm. \ref{ForCompactTruncatedStructureLieGroupResolutionIsInjection} below).

\noindent {\bf (ii)}  ADE-equivariant stable projective bundles are blowup-stable
(Thm. \ref{BlowupStabilityOfADEEquivariantPUHPrimePrincipalBundles},
Thm. \ref{BlowupStabilityForADEEquivariantPUHSemidirectProductZTwoPrincipalBundles} below).
\end{remark}

\medskip

\begin{theorem}[Equivariant bundles with truncated compact Lie structure are blowup-stable]
  \label{ForCompactTruncatedStructureLieGroupResolutionIsInjection}
  Let $\Gamma \,=\, \TorusGroup{r} \rtimes K$
  be the semidirect product of a connected compact abelian Lie group
  (an $r$-torus $\TorusGroup{r} \,\simeq\, (\CircleGroup)^{\times^r}$)
  with a discrete group $K$.
  If the equivariance group $G \,\in\, \Groups(\Sets)$ has a free
  action on some $\SmoothSphere{\,n+2}$ with $n \geq 1$,
  then for all $G \acts \, \SmoothManifold \,\in\, \Actions{G}(\SmoothManifolds)$
  the pullbacks
  \eqref{PullbackOfEquivariantBundleToProductWithFreeGSphereIsSurjection}
  are injections:
  $$
    \begin{tikzcd}
      \Truncation{0}
      \,
      \PointsMaps{\big}
        {
          \HomotopyQuotient
            { \SmoothManifold }
            { G }
        }
        {
          \mathbf{B}
          (
            \TorusGroup{r} \rtimes K
          )
        }
      \ar[
        rr,
        hook,
        "{p^\ast}"
      ]
      &&
      \Truncation{0}
      \,
      \PointsMaps{\big}
        {
          \Quotient
            { (\SmoothManifold \times \SmoothSphere{n + 2}) }
            { G }
        }
        {
          \mathbf{B}
          (
            \TorusGroup{r} \rtimes K
          )
        }
        \,.
    \end{tikzcd}
  $$
  With Lem. \ref{PrincipalBundlesOnBlowUps}
  this means that the full class of
  equivariant bundles with truncated compact Lie structure
  is blowup-stable in the sense of Ntn. \ref{StableEquivariantBundles}.
\end{theorem}
\begin{proof}
  We need to show the implication
  \vspace{-2mm}
  $$
    \exists
    \;\;\;
    \begin{tikzcd}[column sep=large]
      \Quotient
        { \big( \SmoothManifold \times \SmoothSphere{\,n+2} \big) }
        { G }
      \ar[
        rr,
        bend left=20pt,
        "{
          p_{{}_{\SmoothManifold}}^\ast (\vdash P)
        }",
        "{\ }"{swap, name=s}
      ]
      \ar[
        rr,
        bend right=20pt,
        "{
          p_{{}_{\SmoothManifold}}^\ast (\vdash P')
        }"{swap},
        "{\ }"{name=t}
      ]
      \ar[
        from=s,
        to=t,
        Rightarrow,
        "{f}",
        "{\sim}"{swap}
      ]
      &&
      \mathbf{B}
      (
        \TorusGroup{r} \rtimes K
      )
    \end{tikzcd}
    \hspace{.9cm}
    \Rightarrow
    \hspace{.9cm}
    \exists
    \;\;\;
    \begin{tikzcd}
      \HomotopyQuotient
        { \SmoothManifold  }
        { G }
      \ar[
        rr,
        bend left=26pt,
        "{
          \vdash P
        }",
        "{\ }"{swap, name=s}
      ]
      \ar[
        rr,
        bend right=26pt,
        "{
          \vdash P'
        }"{swap},
        "{\ }"{name=t}
      ]
      \ar[
        from=s,
        to=t,
        Rightarrow,
        "{f}",
        "{\sim}"{swap}
      ]
      &&
      \quad
      \mathbf{B}
      (
        \TorusGroup{r} \rtimes K
      )\;.
    \end{tikzcd}
  $$

  \vspace{-3mm}
  \noindent
  Choosing, by Prop. \ref{SmoothGManifoldsAdmitProperlyEquivariantGoodOpenCovers},
  any proper equivariant good open cover
  $\widehat{\SmoothManifold}$ of $\SmoothManifold$,
  we may equivalently express,
  by Ex. \ref{CechActionGroupoidOfEquivariantGoodOpenCoverIsLocalCofibrantResolution},
  the isomorphism $f$
  as a continuous natural transformation between continuous functors
  on the product of the $G$-action {\v C}ech groupoid \eqref{CechActionGroupoidComponents}
  with $\SmoothSphere{\,n+2}$:
  \vspace{-3mm}
  \begin{equation}
    \label{TransformationOfFunctorsOnActionCechGroupoid}
    \hspace{-3cm}
    \begin{tikzcd}[column sep=45pt]
      \SimplicialNerve
      \left(
        \widehat{\SmoothManifold}
          \times_{{}_{\SmoothManifold}}
        \widehat{\SmoothManifold}
          \times
        \SmoothSphere{\,n+2}
          \times
        G^\op
        \rightrightarrows
        \widehat X
          \times
        \SmoothSphere{\,n+2}
      \right)
      \qquad \qquad
      \ar[
        rr,
        bend left=12,
        "{
          p_{{}_{\SmoothManifold}}^\ast ( \vdash P )
        }",
        "{\ }"{swap, name=s}
      ]
      \ar[
        rr,
        bend right=12,
        "{
          p_{{}_{\SmoothManifold}}^\ast ( \vdash P' )
        }"{swap},
        "{\ }"{name=t}
      ]
      \ar[
        from=s,
        to=t,
        Rightarrow,
        "{f}"
      ]
      &&
      \SimplicialNerve
      \big(
        \DeloopingGroupoid
          { \TorusGroup{r} \rtimes K }
      \big)
    \end{tikzcd}
  \end{equation}
  \begin{equation}
    \label{ComponentsOfEquivariantPartOfTransformationBetweenBundlesPullbackToSphere}
    \begin{tikzcd}[column sep=-2pt]
      &
      g_1
      \cdot
      (
        \widehat{x}
        ,\,
        p
      )
      \ar[drr]
      \\
      &&&
      g_2
      \cdot
      (
        \widehat{x}
        ,\,
        p
      )
      \\
      (
        \widehat{x},
        p
      )
      \ar[urrr]
      \ar[uur]
    \end{tikzcd}
    \begin{tikzcd}
     {}
     \ar[
       rr,
       |->,
       shift right=10pt
     ]
     &&
     {}
    \end{tikzcd}
    \begin{tikzcd}[column sep=-2pt]
      &
      \mathrlap{\;\;\;\;\;\;\;\;\;\bullet}
      \phantom{
      g_2
      \cdot
      (
        \widehat{x}
        ,\,
        p
      )
      }
      \ar[drr]
      \ar[
        rrrr,
        start anchor={[xshift=-10pt]}, end anchor={[xshift=+11pt]},
        "{
          (g \cdot f)(\widehat{x}, \, p)
          \,=\,
          f(g_1 \cdot \widehat{x},\, g_1 \cdot p)
        }"
      ]
      &&&&
      \mathrlap{\;\;\;\;\;\;\;\;\;\bullet}
      \phantom{
      g_2
      \cdot
      (
        \widehat{x}
        ,\,
        p
      )
      }
      \ar[drr]
      \\
      &&&
      \mathrlap{\;\;\;\;\;\bullet}
      \phantom{
      g_2
      \cdot
      (
        \widehat{x}
        ,\,
        p
      )
      }
      \ar[rrrr, start anchor={[xshift=-18pt]}, end anchor={[xshift=+7pt]}]
      &&&&
      \mathrlap{\;\;\;\;\;\bullet}
      \phantom{
      g_2
      \cdot
      (
        \widehat{x}
        ,\,
        p
      )
      }
      \\
      \mathrlap{\;\;\;\;\; \bullet}
      \phantom{
      (
        \widehat{x},
        p
      )
      }
      \ar[
        urrr
      ]
      \ar[
        uur,
        "{
          \left(
            \rho(\widehat{x}, g_1)
            ,\,
            \alpha(\widehat{x}, g_1)
          \right)
        }"{sloped}
      ]
      \ar[
        rrrr,
        start anchor={[xshift=-2pt]}, end anchor={[xshift=+3pt]},
        "{
          f(\widehat{x},p)
        }"{swap}
      ]
      &&&&
      \mathrlap{\;\;\;\;\; \bullet}
      \phantom{
      (
        \widehat{x},
        p
      )
      }
      \ar[
        uur,
        crossing over,
        "{
          \colorbox{white}{\scalebox{.9}{$
          \left(
            \rho'(\widehat{x}, g_1)
            ,\,
            \alpha'(\widehat{x}, g_1)
          \right)
          $}}
        }"{yshift=-2pt, sloped}
      ]
      \ar[urrr]
    \end{tikzcd}
  \end{equation}
  \begin{equation}
    \label{ComponentsOfTransitionPartOfTransformationBetweenBundlesPullbackToSphere}
    \begin{tikzcd}[row sep=small, column sep=-2pt]
      &
      \left(
        (x,j)
        ,\,
        p
      \right)
      \ar[drr]
      \\
      &&&
      \left(
        (x,k)
        ,\,
        p
      \right)
      \\
      \left(
        (x,i),
        p
      \right)
      \ar[urrr]
      \ar[uur]
    \end{tikzcd}
    \begin{tikzcd}
     {}
     \ar[
       r,
       |->,
       shift right=10pt
     ]
     &
     {}
    \end{tikzcd}
    \!\!
    \begin{tikzcd}[row sep=small, column sep=-4pt]
      &
      \mathrlap{\;\;\;\;\;\;\;\;\;\bullet}
      \phantom{
      g_2
      \cdot
      (
        \widehat{x}
        ,\,
        p
      )
      }
      \ar[drr]
      \ar[
        rrrr,
        start anchor={[xshift=-10pt]}, end anchor={[xshift=+11pt]},
        "{
          f_j(x,\, p)
        }"
      ]
      &&&&
      \mathrlap{\;\;\;\;\;\;\;\;\;\bullet}
      \phantom{
      g_2
      \cdot
      (
        \widehat{x}
        ,\,
        p
      )
      }
      \ar[drr]
      \\
      &&&
      \mathrlap{\;\;\;\;\;\bullet}
      \phantom{
      g_2
      \cdot
      (
        \widehat{x}
        ,\,
        p
      )
      }
      \ar[rrrr, start anchor={[xshift=-18pt]}, end anchor={[xshift=+7pt]}]
      &&&&
      \mathrlap{\;\;\;\;\;\bullet}
      \phantom{
      g_2
      \cdot
      (
        \widehat{x}
        ,\,
        p
      )
      }
      \\
      \mathrlap{\;\;\;\;\; \bullet}
      \phantom{
      (
        \widehat{x},
        p
      )
      }
      \ar[
        urrr
      ]
      \ar[
        uur,
        "{
          g_{i j}(x)
        }"{sloped}
      ]
      \ar[
        rrrr,
        start anchor={[xshift=-2pt]}, end anchor={[xshift=+3pt]},
        "{
          f_i(x,\, p )
        }"{swap}
      ]
      &&&&
      \mathrlap{\;\;\;\;\; \bullet}
      \phantom{
      (
        \widehat{x},
        p
      )
      }
      \ar[
        uur,
        crossing over,
        "{
          \colorbox{white}{\scalebox{.9}{$
            g'_{i j}(x)
          $}}
        }"{yshift=-2pt, sloped}
      ]
      \ar[urrr]
    \end{tikzcd}
  \end{equation}
  We need to show that $f$ may compatibly be replaced by
  an equivariant functions which
  is constant on $\SmoothSphere{\,n+2}$. Our strategy is to replace the function
  on the $n+2$-sphere by its average
  $\overline{f}$ as seen by integration against the unit volume form
  on the sphere. Since this is a continuous operation, it turns
  {\v C}ech cocycles into {\v C}ech cocycles.

  The key point to achieve this is that $\TorusGroup{r} \rtimes K$ is 1-truncated
  while the dimension of $\SmoothSphere{\,n+2}$ is larger than 1, by assumption.
  This implies that
  every map from the $(n+2)$-sphere lifts to a map through the
  universal covering space:
  \vspace{-2mm}
  $$
    \begin{tikzcd}[column sep=large]
      &&
      \mathbb{R}^r \rtimes K
      \ar[d]
      \ar[r]
      \ar[
        dr,
        phantom,
        "{\mbox{\tiny\rm(pb)}}"
      ]
      &
      \ast
      \ar[d]
      \\
      \SmoothSphere{\,n+2}
      \ar[
        urr,
        dashed,
        "{
          (\vec f \left(\widehat{x}, -), \phi(\widehat{x}\,)\right)
        }"
      ]
      \ar[
        rr,
        "\;\;\;\;{
          f(\widehat{x}, -)
        }"
      ]
      &&
      \TorusGroup{r} \rtimes K
      \ar[r]
      &
      \mathbf{B}
      (
        \mathbb{Z}^r
        \rtimes
        K
      )
      \,.
    \end{tikzcd}
  $$

  \vspace{-1mm}
  \noindent
  Therefore, it makes sense to assign
  \vspace{-2mm}
  \begin{equation}
    \label{AveragingOfGaugeTransformationDataOverSphere}
    f
    \;\longmapsto\;\;
    \overline{f}(\widehat{x}\,)
    \;:=\;
    \bigg(\;
      \underset
        {\SmoothSphere{\,n+2}}
        { \int }
      \vec f(\widehat{x}, p)
      \,
      \mathrm{dvol}(p)
      ,\,
      \phi(\widehat{x}\,)
    \bigg)
    \;\;\;
    \in
    \;
    \mathbb{R}^r \rtimes K
    \twoheadrightarrow
    \TorusGroup{r} \rtimes K
    \,,
  \end{equation}

  \vspace{-1mm}
  \noindent
  which is evidently well-defined
  as a map to $\mathbb{T}^n \rtimes K$,
  in that it is independent of the choice of lift.
  Here $\phi(\widehat{x}\,) \in K$ is the $K$-component of the original function,
  which is already constant
  along the sphere, since $K$ is assumed to be discrete and $\SmoothSphere{\,n+2}$
  is connected.

  It remains to check that
  the averaged component function $\overline{f_i} : U_i  \xrightarrow{\;} \Gamma \rtimes K$
  is again a natural transformation of the form \eqref{TransformationOfFunctorsOnActionCechGroupoid}.
  The key point here is that the automorphism
  group of the $r$-torus acts linearly on component vectors
  by matrix multiplication with
  invertible integer-valued matrices, which implies that it is
  compatible with the integration operation:
    \vspace{-1mm}
  \begin{equation}
    \label{AutomorphismGroupOfTorus}
    \mathrm{Aut}_{\mathrm{Grp}}
    (
      \TorusGroup{r}
    )
    \;\;
    \simeq
    \;\;
    \mathrm{GL}(r,\mathbb{Z})
    \;\subset\;
    \mathrm{Aut}_{\mathbb{R}}
    (\mathbb{R}^r)
    \,.
  \end{equation}
  Using this we find:
    \vspace{-4mm}
  $$
    \def\arraystretch{2.3}
    \begin{array}{lll}
      &
      \overline{f}( g \cdot \widehat{x}\, )
      \\
      &
      \;=\;
      \bigg(
       \displaystyle \int_{\SmoothSphere{\,n+2}}
        \vec f( g \cdot \widehat{x}, g \cdot p )
        \,
        \mathrm{dvol}(p)
        ,\,
        \phi(g \cdot \widehat{x}\,)
      \bigg)
      \;=\;
      \bigg(
        \int_{\SmoothSphere{\,n+2}}
        \big(
          g \cdot \vec f( \widehat{x}, p )
        \big)
        \,
        \mathrm{dvol}(p)
        ,\,
        \phi(g \cdot \widehat{x}\,)
      \bigg)
      &
      \proofstep{
        by \eqref{AveragingOfGaugeTransformationDataOverSphere}
      }
      \\
      &
      \;=\;
      \Bigg(
                           {\displaystyle \int_{\SmoothSphere{\,n+2}} }
        \bigg(
          \alpha(\widehat{x}, \, g)^{-1}
          \Big(
            \vec f(\widehat{x}, \, p )
            +
            \phi(\widehat{x}\,)
            \big(
              \vec \rho'(\widehat{x}, \, g)
            \big)
            -
            \vec \rho(\widehat{x}, g)
          \Big)
        \!\!\bigg)
        \,
        \mathrm{dvol}(p)
        ,\,
        \alpha(\widehat{x},g)^{-1}
        \cdot
        \phi(\widehat{x}\,)
        \cdot
        \alpha'(\widehat{x}, \, g)
      \Bigg)
      &
      \proofstep{
        by \eqref{ComponentsOfEquivariantPartOfTransformationBetweenBundlesPullbackToSphere}
      }
      \\
      &
      \;=\;
      \bigg(
        \alpha(\widehat{x}, \, g)^{-1}
        \Big(\;
            {\displaystyle \int_ { \SmoothSphere{\,n+2} } }
            \,
            f(\widehat{x}, \, p)
            \,
            \mathrm{dvol}(p)
            \,
            +
            \phi(\widehat{x}\,)
            \big(
              \vec \rho'(\widehat{x}, \, g)
            \big)
            -
            \vec\rho(\widehat{x}, \, g)
          \Big)
          ,
          \;
          \alpha(\widehat{x},\, g)^{-1}
            \cdot
          \phi(\widehat{x}\,)
            \cdot
          \alpha(\widehat{x},\, g)
      \bigg)
      &
      \proofstep{
        by \eqref{AutomorphismGroupOfTorus}
      }
      \\
      &
      \;=\;
      \bigg(
        \alpha(\widehat{x}, \, g)^{-1}
        \Big(
          \overline{f}(\widehat{x}\,)
          +
          \phi(\widehat{x}\,)
          \big(
            \vec \rho'(\widehat{x}, \, g)
          \big)
          -
          \vec\rho(\widehat{x}, \, g)
        \Big)
        ,
        \;
        \alpha(\widehat{x},\, g)^{-1}
          \cdot
        \phi(\widehat{x}\,)
          \cdot
        \alpha(\widehat{x},\, g)
      \bigg)
      &
      \proofstep{
        by \eqref{AveragingOfGaugeTransformationDataOverSphere}
      }
      \\
      &
      \;=\;
      (g \cdot \overline{f})(u)
      &
      \proofstep{
        as in \eqref{ComponentsOfEquivariantPartOfTransformationBetweenBundlesPullbackToSphere}
      .}
          \end{array}
  $$

  \vspace{-2mm}
\noindent  This demonstrates that $\overline{f}$
  still makes the equivariance naturality diagrams
  \eqref{ComponentsOfEquivariantPartOfTransformationBetweenBundlesPullbackToSphere}
  commute.
  Verbatim the same computation, just with equivariant transitions replaced by
  transition functions, shows that it also makes the
  {\v C}ech naturality squares \eqref{ComponentsOfTransitionPartOfTransformationBetweenBundlesPullbackToSphere}
  commute. Therefore $\overline{f}$ is an isomorphism between
  $p_{{}_{\SmoothManifold}}^\ast (\vdash P)$
  and
  $p_{{}_{\SmoothManifold}}^\ast (\vdash P')$ which is itself
  the pullback to $\SmoothSphere{\,n+2}$ of an isomorphism between
  $\vdash P$ and $\vdash P'$ themselves.
\end{proof}

\medskip

\noindent
{\bf Blow-up stability of ADE-equivariant projective bundles.}
We now solve the blowup-stability condition
(Ntn. \ref{StableEquivariantBundles})
for $G$-equivariant projective bundles
in the case of ADE equivariance groups
$G \,\subset\, \mathrm{Sp}(1)$ (Ex. \ref{ADEGroupsHaveSphericalSpaceForms})
and show that this recovers
the traditional stability condition  \eqref{StabilityConditionOnGRepresentation}
due to \cite[\S 6]{AtiyahSegal04}.
This is Thm. \ref{BlowupStabilityOfADEEquivariantPUHPrimePrincipalBundles} below.
We prove this by an explicit construction on the level of equivariant
{\v C}ech cocycles. This proof immediately generalizes to the
garded and semidirect product
structure group $\GradedPUH \rtimes \ZTwo$ (Thm. \ref{BlowupStabilityForADEEquivariantPUHSemidirectProductZTwoPrincipalBundles} below).
In this generality the result eventually implies the classification theorem
for stable equivariant $\ZTwo \acts \, \GradedPUH$-principal bundles
twisting equivariant KR-theory, see Ex. \ref{EquivariantBundlesServingAsGeoemtricTwistsOfEquivariantKTheory} below.

\medskip
Traditionally,
for $\mathcal{G}$ a compact Lie group,
a projective
$\mathcal{G}$-representation $\mathcal{G} \acts \mathrm{P}(\HilbertSpace)$ is
called {\it stable} if it is stable under tensoring
with the regular Hilbert space representation $L^2(\mathcal{G})$
of square-integrable functions on the group $G$
(with respect to the Haar measure),
equipped with the canonical $\mathcal{G}$-action given by pullback of such functions
along the left multiplication action of $G$ on itself
(\cite[p. 28]{AtiyahSegal04}):
\begin{equation}
  \label{StabilityConditionOnGRepresentation}
  \mbox{
    $
      \mathcal{G} \acts \, \HilbertSpace
      \,\in\,
      \Actions{\mathcal{G}}(\HilbertSpaces)
    $
    is stable
  }
  \;\;\;\;\;\;
  \Leftrightarrow
  \;\;\;\;\;\;
  \mathcal{G} \acts \,
  \big(
    \HilbertSpace
    \,\otimes\,
    L^2(\mathcal{G})
  \big)
  \;\simeq\;
  \mathcal{G} \acts \, \HilbertSpace
  \,.
\end{equation}
This means equivalently that all projective
irreps of the given projective twist appear
as direct summands with infinite multiplicity
(see Lem. \ref{StableProjectiveIsotorpyRepresentations} below).
Likewise, an equivariant $\mathrm{PU}(\mathcal{H})$-bundle $V$ is called
{\it stable} if its total space is stable against tensoring with
$L^2{\mathcal{G}}$, hence if all its isotropy actions are stable in the above sense
(see \cite{BEJU12}\cite[\S 15]{LueckUribe14}\cite[\S 5]{EspinozaUribe15}).

This condition is traditionally motivated as
ensuring that the space of Fredholm operators on
$\mathcal{G} \acts \, \HilbertSpace$
is a classifying $\mathcal{G}$-space for equivariant $\mathrm{KU}$-theory
(see Ex. \ref{EquivariantBundlesServingAsGeoemtricTwistsOfEquivariantKTheory} below).

\medskip
Since our blow-up stability condition (Ntn. \ref{StableEquivariantBundles})
is a purely bundle-theoretic notion,
it is interesting to observe how it recovers the traditional
operator-algebraic notion:
In our situation $\mathcal{G} \,=\, \SpOne \,\simeq\, \SUTwo$
is the compact Lie group that contains our finite equivariance groups $G$,
and find that
under the
identification $L^2\left(\SpOne\right) \,\simeq\, L^2(\RiemannianSphere{3})$,
the stable $\mathcal{G}$-Hilbert space of equivariant K-theory
at a $G$-orbi-singularity
is identified with that of quantum states on the blow-up $S^3 \to \ast$ of the singularity:
The proof of Thm. \ref{BlowupStabilityOfADEEquivariantPUHPrimePrincipalBundles} below
crucially proceeds by ``quantizing'' cocycle data on the blowup $S^3$ of
an orbifold singularity, identifying it with components of the ``wavefunction''
that the cocycle takes values in, and vice versa.

\medskip

\begin{notation}[Infinite tensor product of $L^2$-Hilbert spaces]
We write $\HilbertSpaces \,\in\, \Categories$ for the category of
complex countably-dimensional Hilbert spaces.

\noindent {\bf (i)} For $(\RiemannianManifold, \VolumeElement)$ a Riemannian manifold,
regarded as a measure space,
of unit volume
\vspace{-2mm}
$$
  \underset{\RiemannianManifold}{\int}
  \, 1 \, \VolumeElement
  \;=\;
  1
$$

\vspace{-2mm}
\noindent
(these are already the most general measure spaces that we need here) we write
\vspace{-2mm}
$$
  L^2(\RiemannianManifold,\, \VolumeElement)
  \;\;\;
  \in
  \;
  \HilbertSpaces
$$

\vspace{-2mm}
\noindent
for the usual Hilbert space of its square-integrable complex-valued functions modulo those
supported on subsets of vanishing measure.

\noindent {\bf (ii)}  Recalling that the tensor product of $L^2$-Hilbert spaces
reflects the Cartesian product of underlying measure spaces
$$
  L^2(\RiemannianManifold, \VolumeElement)
    \otimes
  L^2(\RiemannianManifold', \VolumeElement')
  \;\;
  \simeq
  \;\;
  L^2
  \big(
    \RiemannianManifold
      \times
    \RiemannianManifold',
    \,
    \VolumeElement \wedge \VolumeElement'
  \big),
$$

\vspace{-1mm}
\noindent
consider the (co)limiting case of an infinite number of (tensor) products:
\begin{equation}
  \label{InfiniteTensorProductOfAnL2HilbertSpaceWithItself}
  (\RiemannianManifold,\VolumeElement)^{\times^\infty}
  \;:=\;\;\;\;
  \limit{
    \mathclap{
      { S \in }
      \atop
      { \mathrm{FinSub}(\mathbb{N}) }
    }
  }
  \;
  \RiemannianManifold^S
  \;\;\;\;\;\;\;\;\;
  \longmapsto
  \;\;\;\;\;\;\;\;\;
  \colimit{
    \mathclap{
      { S \in }
      \atop
      { \mathrm{FinSub}(\mathbb{N}) }
    }
  }
  \;
  L^2(\RiemannianManifold,\VolumeElement)^{\otimes^S}
  \;\;\;\;
  =:\;
  L^2(\RiemannianManifold,\VolumeElement)^{\otimes^\infty}
  \,,
\end{equation}

\vspace{-2mm}
\noindent
where the (co)limit is over the directed system of inclusions of
finite subsets of $\NaturalNumbers$:
\vspace{-2mm}
\begin{equation}
  \label{DirectedSystemOfFiniteSubsetsOfNAndInducedSpaces}
  \begin{tikzcd}[row sep=0pt]
    S{\,'}
      \ar[rr, hook, "{ i }"]
    &&
    S
    \\
    \RiemannianManifold^{S{\,'}}
    \ar[from=rr, ->>, "{ X^i }"{description}]
    &&
    \RiemannianManifold^{S}
    \\
    L^2(\RiemannianManifold)^{\otimes^{S{\, '}}}
    \ar[rr, "{ ( \RiemannianManifold^i )^\ast }"{description}]
    &&
    L^2(\RiemannianManifold)^{\otimes^{S}}
    \,.
  \end{tikzcd}
\end{equation}

\vspace{-2mm}
\noindent
{\bf (iii)} This means that the colimiting Hilbert space
$L^2(\RiemannianManifold, \VolumeElement)^{\otimes^\infty}$
\eqref{InfiniteTensorProductOfAnL2HilbertSpaceWithItself}
is that generated by those square-integrable
functions $f$ on $(\RiemannianManifold, \VolumeElement)^{\times^\infty}$ that depend on any finite
number of variables
(hence those which factor as
$(\RiemannianManifold, \VolumeElement)^{\times^\infty}
\xrightarrow{\;} \RiemannianManifold^{S_f} \xrightarrow{\;} \ComplexNumbers$
for some finite subset $S_f \,\subset\, \NaturalNumbers$) with the inner product
on these generating elements being
\vspace{-1mm}
\begin{equation}
\label{EqnInnerProduct}
  \langle f, \, g\rangle
  \;\;=\;
  \underset{
    s \in F_f \cup F_g
  }{\prod}
  \langle f_s, \, g_s\rangle
  \;\;=\;
  \underset{
    s \in F_f \cup F_g
  }{\prod}
  \;
  \underset
    { \RiemannianManifold }
    {\int}
  \,
  \overline{f}_s(x^s)
  \cdot
  g_s(x^s)
  \;
  \VolumeElement(x^s)
  \,.
\end{equation}
\end{notation}

\begin{remark}[Interpretations and properties of infinite tensor products of $L^2$-Hilbert spaces]
$\,$

\noindent {\bf (i)}
The space \eqref{EqnInnerProduct} may equivalently be identified with:
\vspace{-2mm}
\begin{enumerate}[{\bf (a)}]
  \setlength\itemsep{-2pt}
\item  the Hilbert space of square integrable functions on the
infinite product measure space
$(\RiemannianManifold^{\times^\infty}, \VolumeElement^{\times^\infty})$
(e.g., by \cite[Ex. 6.3.11]{Parthasarathy05}),

\item
the
(``incomplete'' \cite[\S 4.1]{vonNeumann39},
or ``Guichardet'' \cite{Guichardet69}, or ``grounded'' \cite[p. 126]{BaezSegalZhou92})
infinite tensor product of Hilbert spaces
(review in \cite[Def. 2.5.1]{Weaver01}),
if we regard each $L^2(\RiemannianManifold,\VolumeElement)$
as equipped with the vacuum state given by
\begin{equation}
  \label{CanonicalVacuumStateInLTwoHilbertSpace}
  \begin{tikzcd}[row sep=3pt]
    \ast
    \ar[
      from=rr,
      "{ p_{{}_{\TopologicalSpace}} }"{swap}
    ]
    &&
    \TopologicalSpace
    \\
    \mathllap{
      \mathbb{C}
      \;\simeq\;\,
    }
    L^2(\ast)
    \ar[rr, "{ (p_{{}_{\TopologicalSpace}})^\ast }"]
    &&
    L^2(\RiemannianManifold, \, \VolumeElement)
    \\[-3pt]
    \scalebox{0.8}{$1$}
    &\longmapsto&
    \scalebox{0.8}{$
      \mathrm{const}_1
      \mathrlap{
        \,\;=:\,
        \lvert \mathrm{vac} \rangle
        \,.
      }
    $}
  \end{tikzcd}
\end{equation}
\end{enumerate}

\noindent {\bf (ii)} In this notation, the inclusions of Hilbert spaces in
    \eqref{DirectedSystemOfFiniteSubsetsOfNAndInducedSpaces} operate by
``adding more copies of the system in its vacuum state'',
where necessary. For example:
\begin{equation}
  \label{FillingUpStatesWithVacuumStates}
  \begin{tikzcd}[row sep=0pt]
    \{1,3\}
      \ar[rr, hook, "{ i }"]
    &&
    \{1,2,3\}
      \ar[rr, hook, "{ }"]
    &&
    \NaturalNumbers
    \\
    \RiemannianManifold^2
    \ar[from=rr, ->>, "{ (\mathrm{pr}_1, \mathrm{pr}_3) }"{description}]
    &&
    \RiemannianManifold^3
    \ar[from=rr, ->>, "{ (\mathrm{pr}_1, \mathrm{pr}_2, \mathrm{pr}_3) }"{description}]
    &&
    \RiemannianManifold^{\infty}
    \\
    L^2(\RiemannianManifold)^{\otimes^2}
    \ar[rr, "{ ( \RiemannianManifold^i )^\ast }"{description}]
    &&
    L^2(\RiemannianManifold)^{\otimes^3}
    \ar[rr]
    &&
    L^2(\RiemannianManifold)^{\otimes^\infty}
    \\
    \scalebox{0.8}{$
      \lvert \psi_1 \rangle
        \otimes
      \lvert \psi_3 \rangle
    $}
    &\longmapsto&
    \scalebox{0.8}{$
      \lvert \psi_1 \rangle
        \otimes
      \lvert \mathrm{vac} \rangle
        \otimes
      \lvert \psi_3 \rangle
    $}
    &\longmapsto&
    \scalebox{0.8}{$
      \lvert \psi_1 \rangle
        \otimes
      \lvert \mathrm{vac} \rangle
        \otimes
      \lvert \psi_3 \rangle
      \otimes
      \big(
        \lvert \mathrm{vac} \rangle
      \big)^{\otimes^\infty}
    $}
    \,.
  \end{tikzcd}
\end{equation}
\noindent {\bf (iii)} This construction is functorial; in particular, for
$U \,\in\, \UnitaryGroup\left( L^2(\RiemannianManifold) \right)$
a unitary operator acting on a single copy of the Hilbert space
which fixes the vaccum state \eqref{CanonicalVacuumStateInLTwoHilbertSpace},
it extends to a diagonal unitary action on the infinite Guichardet tensor product:
\vspace{-2mm}
\begin{equation}
  \label{DiagonalUnitaryActionOnInfiniteGuichardetTensorProduct}
  \left.
  \def\arraystretch{1.3}
  \begin{array}{l}
    U \,\in\, \UnitaryGroup\big( L^2(\RiemannianManifold) \big),
    \\
    U \lvert \mathrm{vac} \rangle \,=\, \lvert \mathrm{vac} \rangle
  \end{array}
  \!\!\!
  \right\}
  \;\;\;\;\;\;
  \vdash
  \;\;\;\;\;\;
  U^{\otimes^\infty}
  \;\in\;
  \UnitaryGroup
  \big(
    L^2(\RiemannianManifold)^{\otimes^\infty}
  \big)
  .
\end{equation}
For example, a diagonal unitary action on the state
\eqref{FillingUpStatesWithVacuumStates} is of this form:
$$
  U^{\otimes^\infty}
  \Big(
    \lvert \psi_1 \rangle
      \otimes
    \lvert \mathrm{vac} \rangle
      \otimes
    \lvert \psi_3 \rangle
    \otimes
    \big(
      \lvert \mathrm{vac} \rangle
    \big)^{\otimes^\infty}
  \Big)
  \;\;
  =
  \;\;
    \lvert U(\psi_1) \rangle
      \otimes
    \lvert \mathrm{vac} \rangle
      \otimes
    \lvert U(\psi_3) \rangle
    \otimes
    \big(
      \lvert \mathrm{vac} \rangle
    \big)^{\otimes^\infty}
  .
$$
\end{remark}

\bigskip
While all countably infinite-dimensional Hilbert spaces are
abstractly isomorphic, functorial families of them -- such as
the combined tensor products and representation spaces that we
encounter in a moment
-- have
a more individual structure.
We next turn attention specifically to the
following Hilbert space:
\begin{example}[Hilbert space of finite higher spin chains]
\label{HilbertSpaceOfHigherSpinChains}
The infinite Guichardet tensor product
\eqref{InfiniteTensorProductOfAnL2HilbertSpaceWithItself}
\begin{equation}
  \label{HilbertSpaceOfBosonsOnThreeSphere}
  \HilbertSpace
  \;:=\;
  L^2\big(\RiemannianSphere{3}\big)^{\otimes^{\infty}}
  \;=\;
  L^2\big(\SpOne\big)^{\otimes^{\infty}}
  \;\;\;
  \in
  \;
  \HilbertSpaces
\end{equation}
of copies of the Hilbert space of square-integrable functions on the round unit 3-sphere:
\begin{equation}
  \label{HilbertSpaceOfSingleBosonOnThreeSphere}
  \HilbertSpace_{\mathrm{sngl}}
  \;:=\;
  L^2\big( \RiemannianSphere{3} \big)
  \;\;\;
  \in
  \;
  \HilbertSpaces
\end{equation}
is generated by the (necessarily finite-dimensional) complex irreducible representations
(irreps)  of the $\SpinGroup(3)$ group. Indeed, under the isomorphism
$$
  \RiemannianSphere{3}
  \;\simeq\;
  \SpOne
  \;\simeq\;
  \SUTwo
  \;\simeq\;
  \SpinGroup(3)
  \;\;\;
  \in
  \;
  \RiemannianManifolds
$$
(for these compact Lie groups equipped with their normalized Haar measure),
the Peter-Weyl theorem implies that
under the canonical action
\eqref{GActionOnSingleBosonHilbertSpaceOverThreeSphere}
by pullback along the (inverse) left multiplication
action
of these groups on themselves we have the following
direct sum decomposition of $G$-Hilbert spaces
(e.g. \cite[Thm. 3.28]{HofmannMorris20}):
\begin{equation}
  \label{FineStructurePeterWeylDecompositionOfLTwoSPOne}
  \SpinGroup(3) \acts \; L^2(\RiemannianSphere{3})
  \;\;\simeq\qquad \quad
  \underset{
    \mathclap{
    \scalebox{.7}{$
      \begin{array}{ll}
        {
          [\mathbf{n}]
          \,\in\,
          \IsomorphismClasses{
            \Representations(
              \SpinGroup(3)
            )_{\mathrm{irr}}
          }
        }
        \\
        \mathrm{dim}_{{}_{\mathbb{C}}}(\mathbf{n}) = n \,\in\, \NaturalNumbers_+
      \end{array}
    $}
    }
  }{\bigoplus}
  \qquad
  n
    \cdot
  \mathbf{n}
  \quad
  \in
  \;
  \Representations
  \big(
    \SpinGroup(3)
  \big)
  \,.
\end{equation}
It is suggestive to observe that the $G$-subspaces
$$
  \SpinGroup(3)
  \acts \;\,
  \mathbf{2} \otimes \mathbf{2} \otimes \mathbf{2} \otimes \cdots \otimes \mathbf{2}
\;  \xhookrightarrow{\quad} \;
  \SpOne
  \acts \;
  L^2\big( \RiemannianSphere{3} \big)^{\otimes^\infty}
  \;\;\;
  \in
  \;
  \Actions{\SpinGroup(3)}
  (
    \HilbertSpaces
  )
$$
which are finite tensor products of the fundamental representation
(the defining representation of $\SUTwo$)
$$
  \mathbf{2}
  \;\simeq\;
  \mathrm{Span}
  \big(
    \lvert \downarrow\rangle
    \, ,\,
    \lvert \uparrow\rangle
  \big)
$$
are known as the Hilbert spaces of
finite
{\it Heisenberg spin chain}-models (e.g. \cite[\S 4]{Saberi18}\cite[\S 2]{Stolz20}),
and the subspaces of finite tensor products of
any given higher dimensional representation
$$
  \SpinGroup(3)
  \acts \;\,
  \mathbf{n} \otimes \mathbf{n} \otimes \mathbf{n} \otimes \cdots \otimes \mathbf{n}
 \; \xhookrightarrow{\quad} \;
  \SpinGroup(2)
  \acts \;
  L^2\big( \RiemannianSphere{3} \big)^{\otimes^\infty}
  \;\;\;
  \in
  \;
  \Actions{\SpinGroup(3)}
  (
    \HilbertSpaces
  )
$$
are those of {\it higher spin chain} models.
Since the Hilbert space \eqref{HilbertSpaceOfBosonsOnThreeSphere}
naturally includes these finite higher spin chain Hilbert spaces,
together with all finite ``mixed higher spin'' chains,
$$
  \SpinGroup(3)
  \acts \;\,
  \mathbf{n}_1 \otimes \mathbf{n_2} \otimes \mathbf{n}_3
  \otimes \cdots \otimes \mathbf{n}_{\ell}
  \;\xhookrightarrow{\quad} \;
  \SpinGroup(3)
  \acts \;
  L^2\big( \RiemannianSphere{3} \big)
$$
we will loosely refer to it as the
{\it higher spin chain Hilbert space}.\footnote{
Some authors also consider the infinite Heisenberg spin chain given by
the Guichardet tensor product
$\mathbf{2}^{\otimes^\infty}$ with respect to
regarding $\lvert \downarrow \rangle \,\in\, \mathbf{2}$ as the vacuum state.
Beware that this space is {\it not} naturally a subspace of our Hilbert spaces, whose
selected vacuum state is instead $\lvert \mathrm{vac} \rangle \,\in\, \mathbf{1}$.}
\end{example}

\begin{definition}[Projective-unitary group on higher mixed spin chains]
\label{ModifiedUnitaryAndProjectiveUnitaryGroup}
$\,$

\noindent
{\bf (i)}
We now write
  \vspace{-3mm}
\begin{equation}
  \label{ModifiedUnitaryGroup}
  \UH
  \,=\,
  \UHPrime
  \;:=\;
  \UnitaryGroup
  \big(
    L^2(\RiemannianSphere{3})^{\otimes^\infty}
  \big)
  \;\;\;
  \in
  \;
  \Groups(\kTopologicalSpaces)
  \xrightarrow{ \Groups(\ContinuousDiffeology) }
  \Groups(\SmoothInfinityGroupoids)
\end{equation}
for the unitary group \eqref{TheGroupUH}
specifically of the Hilbert space \eqref{HilbertSpaceOfBosonsOnThreeSphere}
of ``higher spin chains'', so that
the projective unitary group $\PUH$ \eqref{TheGroupPUH}
is now specifically the
quotient of \eqref{ModifiedUnitaryGroup} by the (topological) circle group
\vspace{-2mm}
$$
  \begin{tikzcd}
    1
    \ar[r]
    &
    \CircleGroup
    \ar[r, hook]
    &
    L^2(
      \RiemannianSphere{3}
    )^{\otimes^\infty}
    \ar[r, ->>]
    &
    \PUH
    \ar[r]
    &
    1
    \,.
  \end{tikzcd}
$$

\vspace{-2mm}
\noindent
Of course, this is isomorphic to any other
construction of the bundle $\CircleGroup \to \UH \to \PUH$, but in this
incarnation we have a manifest inclusion of the following
non-isomorphic but weakly-equivalent version of the circle group.

\noindent
{\bf (ii)} We write
\vspace{-2mm}
\begin{equation}
  \label{ModifiedCircleGroup}
  \PhaseGroup
  \;:=\quad
  \colimit{
    \mathclap{
      { F \in  }
      \atop
      { \mathrm{FinSub}(\NaturalNumbers) }
    }
  }
  \,
  \Maps{\big}
    { (\SmoothSphere{3})^F }
    { \CircleGroup }
  \;\;\;
  \in
  \;
  \Groups(\DTopologicalSpaces)
\end{equation}

\vspace{-1mm}
\noindent
for the colimit of D-topological groups of $\CircleGroup$-valued
continuous functions
on the 3-sphere
over the system \eqref{DirectedSystemOfFiniteSubsetsOfNAndInducedSpaces}.
(This is the space of $\CircleGroup$-valued functions on
$(\SmoothSphere{3})^{\times^\infty}$ \eqref{InfiniteTensorProductOfAnL2HilbertSpaceWithItself}
which depend on any finite number of the factor spaces.)

\noindent
{\bf (iii)}
The system of evident multiplication actions
\vspace{-2mm}
$$
  F, F' \,\in\, \mathrm{FinSub}(\NaturalNumbers)
  \;\;
  \vdash
  \;\;
  \begin{tikzcd}
    \Maps{\big}
      { (\TopologicalSphere{3})^F }
      { \CircleGroup }
    \,\times\,
    L^2\big( (\RiemannianSphere{3})^{F'}\big)
    \ar[r]
    &
    L^2\big( (\RiemannianSphere{3})^{(F \cup F')}\big)
    \xrightarrow{\quad}
    L^2(\RiemannianSphere{3})^{\otimes^\infty}
  \end{tikzcd}
$$

\vspace{-2mm}
\noindent
passes to the double colimit to yield an action
of $\PhaseGroup$ \eqref{ModifiedCircleGroup}
on $\UHPrime$ \eqref{ModifiedUnitaryGroup}
\vspace{-2mm}
$$
  \PhaseGroup
  \,\times\,
  L^2(\RiemannianSphere{3})^{\otimes^\infty}
  \xrightarrow{\qquad}
  L^2(\RiemannianSphere{3})^{\otimes^\infty}
$$

\vspace{-2mm}
\noindent
which is manifestly unitary and hence, adjointly, an inclusion
into the unitary group \eqref{ModifiedUnitaryGroup}
\vspace{-2mm}
\begin{equation}
  \label{ModifiedCircleGroupAsSubgroupOfUnitaryGroup}
  \begin{tikzcd}
    \PhaseGroup
    \ar[r, hook]
    &
    \UnitaryGroup
    \big(
      L^2(\RiemannianSphere{3})^{\otimes^\infty}
    \big)
    \;=\;
    \UHPrime
    \,.
  \end{tikzcd}
\end{equation}

\vspace{-2mm}
\noindent
{\bf (iv)}
The quotient of the unitary group $\UHPrime$ \eqref{ModifiedUnitaryGroup}
by its subgroup $\PhaseGroup$ \eqref{ModifiedCircleGroupAsSubgroupOfUnitaryGroup}
we denote $\PUHPrime$:
\vspace{-2mm}
\begin{equation}
  \label{ModifiedProjectiveUnitaryGroup}
  \begin{tikzcd}[row sep=7pt]
    \PhaseGroup
    \ar[rr, hook]
    \ar[
      d,
      phantom,
      "{
        :=
      }"{rotate=-90}
    ]
    &&
    \UHPrime
    \ar[rr, ->>]
    \ar[
      d,
      phantom,
      "{
        :=
      }"{rotate=-90}
    ]
    &&
    \PUHPrime
    \;:=\;
    \UHPrime/\PhaseGroup
    \\
    \colimit{
      \mathclap{
        { F \in }
        \atop
        { \mathrm{FinSub}(\NaturalNumbers) }
      }
    }
    \,
    \Maps{\big}
      { (\RiemannianSphere{3})^F }
      { \CircleGroup }
    \ar[rr, hook]
    &&
    \UnitaryGroup
    \Bigg( \quad
        \colimit{
         \mathclap{
           { F \in }
             \atop
            { \mathrm{FinSub}(\NaturalNumbers) }
          }
        }
        \,
        L^2\big( \RiemannianSphere{3}  \big)^{\otimes^F}
    \Bigg)
  \end{tikzcd}
\end{equation}
\end{definition}

\begin{remark}[Adjoining a spinor system to the chain]
\label{AdjoiningAnElementaryQuantumSystemToTheTower}
The notions in Def. \ref{ModifiedUnitaryAndProjectiveUnitaryGroup}
are such as to seamlessley allow
``adjoining one more copy'' of the elementary system \eqref{HilbertSpaceOfSingleBosonOnThreeSphere}
to the infinite tower \eqref{HilbertSpaceOfBosonsOnThreeSphere}:

\noindent
{\bf (i)}
For example, the tensor product of the Hilbert space
\eqref{HilbertSpaceOfBosonsOnThreeSphere}
with one elementary system
\eqref{HilbertSpaceOfSingleBosonOnThreeSphere} from the left
is canonically re-included into the former by shifting all other
copies up:
$$
  \begin{tikzcd}[row sep=4pt]
    L^2(\RiemannianSphere{3})
    \otimes
    L^2(\RiemannianSphere{3})^{\mathrlap{ \otimes^\infty }}
    \ar[
      rr,
      hook,
      start anchor={[xshift=5pt]},
      "{
        \otimes
      }"
    ]
    &&
    L^2(\RiemannianSphere{3})^{\mathrlap{ \otimes^\infty }}
    \\[+6pt]
    L^2(\RiemannianSphere{3})
    \otimes
    L^2(\RiemannianSphere{3})^{\mathrlap{ \otimes^{ \{s_1, \cdots, s_k\} } }}
    \ar[
      rr,
      hook,
      start anchor={[xshift=30pt]}
    ]
    \ar[u,hook]
    &&
    L^2(\RiemannianSphere{3})^{\mathrlap{ \otimes^{ \{1, s_1+1, \cdots, s_k + 1\} } }}
    \ar[u,hook]
    \\
\qquad \scalebox{0.7}{$
   \lvert \psi \rangle
    \otimes
    \big(
      \lvert \psi_{s_1} \rangle
      \otimes
      \cdots
      \otimes
      \lvert \psi_{s_k} \rangle
    \big)
    $}
    &\qquad \quad \longmapsto&
    \qquad  \quad  \scalebox{0.7}{$
    \big(
      \lvert \psi \rangle
      \otimes
      \lvert \psi_{s_1} \rangle
      \otimes
      \cdots
      \otimes
      \lvert \psi_{s_k} \rangle
    \big)
    $}
    \,,
  \end{tikzcd}
$$

\vspace{-1mm}
\noindent
and this tensoring induces a corresponding inclusion of unitary
and projective unitary groups from Def. \ref{ModifiedUnitaryAndProjectiveUnitaryGroup}:
\vspace{-1mm}
\begin{equation}
  \label{TensoringOfPUHPrimeOverUHPrime}
  \begin{tikzcd}[row sep=small]
    1 \times \CircleGroup
    \ar[d]
    \ar[rr]
    &&
    \CircleGroup
    \ar[d]
    \\
    \UnitaryGroup
    \big(
      L^2(\RiemannianSphere{3})
    \big)
    \times
    \UHPrime
    \ar[
      rr,
      hook,
      "{ \otimes }"
    ]
    \ar[d]
    &&
    \UHPrime
    \ar[d]
    \\
    \UnitaryGroup
    \big(
      L^2(\RiemannianSphere{3})
    \big)
    \times
    \PUH
    \ar[rr]
    &&
    \PUH
    \mathrlap{\,,}
  \end{tikzcd}
  \hspace{1cm}
  \begin{tikzcd}[row sep=small]
    1 \times \PhaseGroup
    \ar[d]
    \ar[rr]
    &&
    \PhaseGroup
    \ar[d]
    \\
    \UnitaryGroup
    \big(
      L^2(\RiemannianSphere{3})
    \big)
    \times
    \UHPrime
    \ar[
      rr,
      hook,
      "{ \otimes }"
    ]
    \ar[d]
    &&
    \UHPrime
    \ar[d]
    \\
    \UnitaryGroup
    \big(
      L^2(\RiemannianSphere{3})
    \big)
    \times
    \PUHPrime
    \ar[rr]
    &&
    \PUHPrime
    \mathrlap{\,.}
  \end{tikzcd}
\end{equation}

\noindent
{\bf (ii)}
Similarly, since $\SmoothSphere{3}$ is compact
(certainly in the traditional sense, but also as an object
of the $\infty$-topos
$\SmoothInfinityGroupoids = \InfinitySheaves(\SmoothManifolds)$,
still in the sense that all of its covers have a finite sub-cover,
see \cite[\S 3.6.4]{dcct}) and since the colimit
\eqref{ModifiedCircleGroup} is over monomorphisms, it follows
(\cite[Prop. 3.6.61]{dcct})
that
every map from $\SmoothSphere{3}$ to $\PhaseGroup$ \eqref{ModifiedCircleGroup}
factors through one
of its finite stages, so that $\Maps{}{\SmoothSphere{3}}{\PhaseGroup}$
is naturally re-included into $\PhaseGroup$:
\begin{equation}
  \label{ReIncludingMapsIntoModifiedCircleGroupIntoModifiedCircleGroup}
  \Maps{\big}
    { \SmoothSphere{3} }
    { \PhaseGroup }
  \;\xhookrightarrow{\;\;\sim \;\;}\;
  \PhaseGroup\;,
\end{equation}

\vspace{-2mm}
\noindent
namely as follows:
\vspace{-1.5cm}
$$
\renewcommand\arraystretch{1.5}
  \begin{array}{lll}
  \\[+17pt]
  \Maps{\big}
    { \SmoothSphere{3} }
    { \PhaseGroup }
  &
  \;=\;
  \Maps{\Bigg}
  { \SmoothSphere{3} }
  {\;\;
    \colimit{
      \mathclap{
        { F \in  }
        \atop
        { \mathrm{FinSub}(\NaturalNumbers) }
      }
    }
    \,
    \Maps{\Big}
      { (\SmoothSphere{3})^F }
      { \CircleGroup }
  }
  &
  \proofstep{
    by \eqref{ModifiedCircleGroup}
  }
  \\
  &
  \;\simeq\quad
  \colimit{
    \mathclap{
      { F \in  }
      \atop
      { \mathrm{FinSub}(\NaturalNumbers) }
    }
  }
  \;
  \Maps{\Big}
  { \SmoothSphere{3} }
  {
    \Maps{\big}
      { (\TopologicalSphere{3})^F }
      { \CircleGroup }
  }
  &
  \proofstep{
    by \cite[Prop. 3.6.61]{dcct}
  }
  \\
  & \;\simeq\quad
  \colimit{
    \mathclap{
      { F \in  }
      \atop
      { \mathrm{FinSub}(\NaturalNumbers) }
    }
  }
  \;
  \Maps{\Big}
    { (\TopologicalSphere{3})^{\ast \sqcup F } }
    { \CircleGroup }
  \;\;
  =
  \quad
  \colimit{
    \mathclap{
      { F \in  }
      \atop
      { \mathrm{FinSub}(\NaturalNumbers_+) }
    }
  }
  \;
  \Maps{\Big}
    { (\TopologicalSphere{3})^{F} }
    { \CircleGroup }
  &
  \proofstep{
    by \eqref{InternalHomAdjunction}
  }
  \\
  &
  \;\hookrightarrow\quad
  \colimit{
    \mathclap{
      { F \in  }
      \atop
      { \mathrm{FinSub}(\NaturalNumbers) }
    }
  }
  \;
  \Maps{\Big}
    { (\TopologicalSphere{3})^{F} }
    { \CircleGroup }
  &
  \proofstep{
    by inspection
  }
  \\
  &
  \;=\;
  \PhaseGroup
  &
  \proofstep{
    by \eqref{ModifiedCircleGroup}.
  }
  \end{array}
$$

\end{remark}

\begin{example}[$G$-Action on space of higher spin chains]
\label{GActionOnQuantumStatesOnThreeSphere}
For $G \,\subset\, \SpOne$ any subgroup
(in particular a finite subgroup, Ex. \ref{ADEGroupsHaveSphericalSpaceForms}),
its left multiplication action on $\SmoothSphere{3} = \SpOne$
induces, by isometric pullback of square-integrable functions,
a unitary representation of $G$ on $L^2\big( \RiemannianSphere{3} \big)$
\eqref{HilbertSpaceOfSingleBosonOnThreeSphere}:
 \vspace{-2mm}
\begin{equation}
  \label{GActionOnSingleBosonHilbertSpaceOverThreeSphere}
  \begin{tikzcd}[row sep=-2pt]
    G
    \ar[
      rr,
      "{}"
    ]
    &&
    \UnitaryGroup
    \left(
      L^2
      \left(
       \RiemannianSphere{3}
      \right)
    \right)
    \\
    \scalebox{0.8}{$
      g
    $}
   &\longmapsto&
   \hspace{-1.5cm}
   \scalebox{0.8}{$
   \big( g^{-1}\big)^\ast
    \mathrlap{
      \;=:\;
      \underset{\RiemannianSphere{3}}{\displaystyle \int}
      \lvert g^{-1} \cdot p\rangle \langle p \vert
    }
    $}
      \end{tikzcd}
\end{equation}

 \vspace{-4mm}
\noindent
which diagonally extends \eqref{DiagonalUnitaryActionOnInfiniteGuichardetTensorProduct}
to an action by $\UHPrime$ \eqref{ModifiedUnitaryGroup}
on the infinite tensor product Hilbert space \eqref{HilbertSpaceOfBosonsOnThreeSphere}:
 \vspace{-2mm}
\begin{equation}
  \label{DiagonalCanonicalActionOnInfiniteGuichardetTensorProduct}
  \begin{tikzcd}[row sep=-2pt]
    G
    \ar[
      rr,
      "{}"
    ]
    &&
    \UHPrime
    =
    \UnitaryGroup
    \Big(
      L^2
      \big(
        \RiemannianSphere{3}
      \big)^{\otimes^\infty}
    \Big)
    \\
    \scalebox{0.8}{$
          g
    $}
   &\longmapsto&
   \hspace{-2cm}
   \scalebox{0.8}{$
   \big( g^{-1}\big)^\ast
    \mathrlap{
      \;=:\;
      \underset{ \NaturalNumbers }{\bigotimes}
      \,
      \underset{\RiemannianSphere{3}}{\displaystyle \int}
      \lvert g^{-1} \cdot p\rangle \langle p \vert
      \,.
    }
    $}
  \end{tikzcd}
\end{equation}

\end{example}

\begin{lemma}[Shape of projective unitary group on higher spin chains]
  \label{ShapeOfProjectiveUnitaryGroupOnBosonStatesOver3Sphere}
  $\,$

  \noindent
  {\bf (i)}
  The shape of
  the modified circle group \eqref{ModifiedCircleGroup}
  is still that of a $K(\Integers, 1)$:
   \vspace{-2mm}
  $$
    \shape
    \,
    \PhaseGroup
    \;\simeq\;
    B \Integers
    \,.
  $$

\vspace{-2mm}
  \noindent
  {\bf (ii)}
  The shape of
  the modified projective unitary group
  \eqref{ModifiedProjectiveUnitaryGroup}
  is still that
  of a $K(\mathbb{Z},2)$:
   \vspace{-2mm}
  $$
    \shape
    \,
    \PUHPrime
    \;\simeq\;
    B^2 \Integers
    \,.
  $$
\end{lemma}
\begin{proof}
  Since $\pi_{\geq 2}(\ShapeOfSphere{1}) \,=\, \ast$,
  the long exact sequence of homotopy groups induced by the homotopy fiber sequence
  \vspace{-2mm}
  \begin{equation}
    \label{MapsFrom3SphereToOneSphereWeakHomotopyEquivalentToOneSphere}
    \begin{tikzcd}
      \Omega^3 \CircleGroup
      \ar[
        rr,
        "{
          \mathrm{fib}(\mathrm{ev}_\ast)
        }"
      ]
      &&
      \Maps{}
        { \ShapeOfSphere{3} }
        { \ShapeOfSphere{1} }
      \ar[
        rr,
        "{
          \mathrm{ev}_\ast
        }",
        "{\in \WeakHomotopyEquivalences}"{swap}
      ]
      &&
      \ShapeOfSphere{1}
    \end{tikzcd}
  \end{equation}

  \vspace{-2mm}
  \noindent
  shows that the
  map on the right
  (evaluating at any chosen basepoint)
  is a weak homotopy equivalence.
  Therefore:
  \vspace{-2mm}
  \begin{equation}
    \label{ShapeOfMapsFromThreeSphereToCircle}
    \def\arraystretch{1.3}
    \begin{array}{lll}
      \shape
      \,
      \Maps{}
        { \SmoothSphere{3} }
        { \CircleGroup }
      &
      \;\simeq\;
      \Maps{\big}
        { \shape \, \SmoothSphere{3} }
        { \shape \, \CircleGroup }
      &
      \proofstep{
        by Thm. \ref{SmoothOkaPrinciple}
      }
      \\
      &
      \;\simeq\;
      \Maps{\big}
        { \ShapeOfSphere{3} }
        { \ShapeOfSphere{1} }
      &
      \proofstep{
        by Prop. \ref{CohesiveShapeOfSmoothManifoldsIsTheirHomotopyType}
        and Prop. \ref{SmoothShapeGivenBySmoothPathInfinityGroupoid}
      }
      \\
      &
      \;\simeq\;
      \ShapeOfSphere{1}
      \;\simeq\;
      B \Integers
      &
      \proofstep{
        by \eqref{MapsFrom3SphereToOneSphereWeakHomotopyEquivalentToOneSphere}.
      }
    \end{array}
  \end{equation}
  By induction on the number $k \in \NaturalNumbers$ of
  3-sphere factors, it follows more generally that:
  \begin{equation}
    \label{ShapeOfMapsFromProductOfThreeSphereToCircle}
        \def\arraystretch{1.3}
    \begin{array}{lll}
      \shape
      \,
      \Maps{\Big}
        { (\SmoothSphere{3})^{\times^{k+1}} }
        { \CircleGroup }
      &
      \;\simeq\;
      \shape
      \,
      \Maps{\Big}
        { \SmoothSphere{3} }
        {
          \Maps{\big}
            { (\SmoothSphere{3})^k  }
            { \CircleGroup }
        }
      &
      \proofstep{
        by Lem. \ref{InternalHomAdjointness}
      }
      \\
      & \;\simeq\;
      \Maps{\Big}
        { \shape\, \SmoothSphere{3} }
        {
          \shape
          \,
          \Maps{\big}
            { (\SmoothSphere{3})^k  }
            { \CircleGroup }
        }
      &
      \proofstep{
        by Thm. \ref{SmoothOkaPrinciple}
      }
      \\
      & \;\simeq\;
      \Maps{\big}
        { \ShapeOfSphere{3} }
        {
          \ShapeOfSphere{1}
        }
      &
      \proofstep{
        by induction from \eqref{ShapeOfMapsFromThreeSphereToCircle}
      }
      \\
      & \;\simeq\;
      \ShapeOfSphere{1} \,\simeq\, B \Integers
      &
      \proofstep{
        by
        \eqref{MapsFrom3SphereToOneSphereWeakHomotopyEquivalentToOneSphere}.
      }
    \end{array}
  \end{equation}

  \vspace{-2mm}
\noindent  Finally, since passing to shape commutes over colimits,
  we find that the large circle group $\PhaseGroup$ still has the
  shape of the plain circle, as claimed:
   \vspace{-2mm}
  \begin{equation}
    \label{ObtainingShapeOfModifiedCircleGroup}
    \hspace{1.5cm}
    \def\arraystretch{1.5}
    \begin{array}{lll}
      \shape
      \,
      \PhaseGroup
      &
      \;=\;
      \shape
      \;
      \colimit{
        \mathclap{
          { F \in }
          \atop
          { \mathrm{FinSub}(\NaturalNumbers) }
        }
      }
      \;
      \Maps{\big}
        { (\SmoothSphere{3})^{F} }
        { \CircleGroup }
      &
      \proofstep{
        by \eqref{ModifiedCircleGroup}
      }
      \\
      & \;\simeq\;
      \colimit{
        \mathclap{
          { F \in }
          \atop
          { \mathrm{FinSub}(\NaturalNumbers) }
        }
      }
      \;
      \shape
      \,
      \Maps{\big}
        { (\SmoothSphere{3})^{F} }
        { \CircleGroup }
      &
      \proofstep{
        by \eqref{InfinityAdjointPreservesInfinityLimits}
      }
      \\
      & \;\simeq\;
      \colimit{
        \mathclap{
          { F \in }
          \atop
          { \mathrm{FinSub}(\NaturalNumbers) }
        }
      }
      \,
      \ShapeOfSphere{1}
      &
      \proofstep{
        by
        \eqref{ShapeOfMapsFromProductOfThreeSphereToCircle}
      }
      \\
      & \;\simeq\;
      \ShapeOfSphere{1} \,\simeq\, B \Integers
      &
      \proofstep{
        colimit over constant diagram
      }
      .
    \end{array}
  \end{equation}

  \vspace{-2mm}
  \noindent
  With this, the second claim follows as for the ordinary $\PUH$:
  \vspace{-.2cm}
  \begin{equation}
    \label{AComputation}
    \def\arraystretch{1.4}
    \begin{array}{lll}
      \shape \, \PUHPrime
      &
      \;=\;
      \shape \big( \Quotient{\UHPrime}{\PhaseGroup}  \big)
      &
      \proofstep{
        by
        \eqref{ModifiedProjectiveUnitaryGroup}
      }
      \\
      &
      \;\simeq\;
      \shape \big( \HomotopyQuotient{\UHPrime}{\PhaseGroup}  \big)
      &
      \proofstep{
        by Prop. \ref{HomotopyQuotientOfFreeActionsIsOrdinaryQuotientOfZeroTruncation}
      }
      \\
      &
      \;\simeq\;
      \HomotopyQuotient{( \shape\,\UHPrime) }{ (\shape\,\PhaseGroup)}
      &
      \proofstep{
        by \eqref{InfinityAdjointPreservesInfinityLimits}
      }
      \\
      & \;\simeq\;
      \HomotopyQuotient{ \ast }{ B \Integers  }
      &
      \proofstep{
        by \eqref{UHIsContractible}
        \& \eqref{ObtainingShapeOfModifiedCircleGroup}
      }
      \\
      &
      \;\simeq\;
      B^2 \Integers
      &
      \proofstep{
        by
        \eqref{DeloopingOfInfinityGroupAsColimit}
      }
      \,.
    \end{array}
  \end{equation}

  \vspace{-7mm}
\end{proof}

\begin{example}[Inclusion of standard into modified projective unitary group is weak equivalence]
\label{InclusionOfStandardIntoModifiedProjectiveUnitaryGroupIsWeakEquivalence}
Arguing as Lem. \ref{ShapeOfProjectiveUnitaryGroupOnBosonStatesOver3Sphere},
one sees that

{\bf (a)} the canonical inclusions of the ordinary projective unitary group
$\PUH$ into $\PUHPrime$;

{\bf (b)} the self-inclusion of the latter under the infinite tensoring
\eqref{TensoringOfPUHPrimeOverUHPrime}

\noindent
are both weak homotopy equivalences:
\begin{equation}
  \label{InfiniteTensoringOfPUHIsWeakHomotopyEquivalence}
  \begin{tikzcd}[column sep=large]
    \CircleGroup
    \ar[rr, hook]
    \ar[d, hook]
    &&
    \mathrm{U}
    \big(
      L^2(\RiemannianSphere{3})^{\otimes^\infty}
    \big)
    \ar[rr, ->>]
    \ar[d, -, shift left=1pt]
    \ar[d, -, shift right=1pt]
    &&
    \PUH
    \ar[
      d,
      "{
        \in \WeakHomotopyEquivalences
      }"
    ]
    \\
    \PhaseGroup
    \ar[rr, hook]
    \ar[d, hook]
    &&
    \mathrm{U}
    \big(
      L^2(\RiemannianSphere{3})^{\otimes^\infty}
    \big)
    \ar[rr, ->>]
    \ar[
      d,
      hook,
      "{
        \UnitaryGroup
        \big(
          L^2(\RiemannianSphere{3})^{\otimes^\infty}
        \big)
        \otimes (-)
      }"{swap}]
    &&
    \PUHPrime
    \ar[d, "{ \in \WeakHomotopyEquivalences }"]
    \\
    \PhaseGroup
    \ar[rr, hook]
    &&
    \mathrm{U}
    \big(
      L^2(\RiemannianSphere{3})^{\otimes^\infty}
    \big)
    \ar[rr, ->>]
    &&
    \PUHPrime
  \end{tikzcd}
\end{equation}
\end{example}

\begin{definition}[Blowup-stable ADE-equivariant projective bundles]
  \label{BlowupStableADEEquivariantPUHPrimePrincipalBundles}
  For a finite subgroup $G \,\subset\, \SpOne$ (Ex. \ref{ADEGroupsHaveSphericalSpaceForms})
  and $G \acts \, \SmoothManifold \,\in\, \Actions{G}(\SmoothManifolds)$,
  consider any equivariant good open cover
  $
    G \acts  \; \widehat{\SmoothManifold}
    \;=\;
    G \acts \;
    \big(
      \underset{i \in I}{\sqcup}
      U_i
    \big)
  $
  (by Prop. \ref{SmoothGManifoldsAdmitProperlyEquivariantGoodOpenCovers}).

  \noindent
  {\bf (i)}
  On its {\v C}ech action groupoid (Ex. \ref{CechActionGroupoidOfEquivariantGoodOpenCoverIsLocalCofibrantResolution}),
  we have the following canonical equivariant {\v C}ech cocycle
  (Rem. \ref{EquivariantCechCocycles})
  with values in $\UnitaryGroup\big( L^2(\RiemannianSphere{3}) \big)$
  \eqref{GActionOnSingleBosonHilbertSpaceOverThreeSphere}, acting
  via the canonical unitary $G$-action from Ex. \ref{GActionOnQuantumStatesOnThreeSphere}:
  \vspace{-2mm}
  \begin{equation}
    \label{TheRegularCocycle}
    \begin{tikzcd}
      c_{{}_{\mathrm{reg}}}
      :
      &[-25pt]
      \Big(\!
        \big(
        \widehat{\TopologicalSpace}
          \times_{{}_{\TopologicalSpace}}
        \widehat{\TopologicalSpace}
        \big)
        \times
        G
        \,\rightrightarrows\,
        \widehat{\TopologicalSpace}
      \Big)
      \ar[
        r,
        ->>,
        "{
          p_{{}_{\widehat{\TopologicalSpace}}}
        }"
      ]
      &
      \big(
        \DeloopingGroupoid{ G }
      \big)
      \ar[
        rrr,
        "{
          \DeloopingGroupoid
            { \left( (-)^{-1} \right)^\ast }
        }"
      ]
      &&&
      \Big(\!\!
        \DeloopingGroupoid{
          \UnitaryGroup
          \big(
            L^2(S^3)
          \big)
        }
      \!\!\Big).
    \end{tikzcd}
  \end{equation}

  \vspace{-2mm}
  \noindent
  {\bf (ii)}
  We say that a $G$-equivariant principal bundle
  $\HomotopyQuotient{\SmoothManifold}{X} \xrightarrow{[c]}
  \mathbf{B} \PUHPrime$
  (Def. \ref{GEquivariantGammaPrincipalBundles})
  with structure group $\PUH$ or $\PUHPrime$
  (Def. \ref{ModifiedUnitaryAndProjectiveUnitaryGroup})
  is {\it stable} if the operation of
  tensoring \eqref{TensoringOfPUHPrimeOverUHPrime}
  any of its representing
  equivariant {\v C}ech cocycles \eqref{IdentifyingEquivariantBundlesWithEquivariantCechCocycles}
  \vspace{-2mm}
  $$
  \begin{tikzcd}
    c
    :
    &[-20pt]
    \Big(\!
      \big(
      \widehat{\TopologicalSpace}
        \times_{{}_{\TopologicalSpace}}
      \widehat{\TopologicalSpace}
      \big)
        \times
      G
      \,\rightrightarrows\,
      \widehat{\TopologicalSpace}
    \Big)
    \ar[rr]
    &&
    \DeloopingGroupoid
      { \PUHPrime  }
  \end{tikzcd}
$$

\vspace{-2mm}
\noindent with the canonical cocycle \eqref{TheRegularCocycle}
\vspace{-2mm}
\begin{equation}
\label{ComponentsOfTensoringWithRegularCocycle}
\begin{aligned}
  c
   & \;\longmapsto\;
  c_{{}_{\mathrm{reg}}} \otimes c
\\
\scalebox{0.8}{$
 \ProjectiveUnitaryGroup
  \Big(
    L^2(\RiemannianSphere{3})^{\otimes^\infty}
  \Big)
  \;
  \ni
  \quad
  c
  \big(
    (x,i,j),g
  \big)
  $}
   & \;\longmapsto\;\;\;
 \scalebox{0.7}{$
  \underset{
    \RiemannianSphere{3}
  }{\displaystyle \int}
  \;
  \lvert g^{-1}\cdot p \rangle \langle p \rvert
  \,\otimes\,
  c\big((x,i,j),g \big)
  \quad \in\;
  \UnitaryGroup
  \Big(
    L^2(\RiemannianSphere{3})
      \otimes
    L^2(\RiemannianSphere{3})^{\otimes^\infty}
  \Big)
  $}
  \end{aligned}
\end{equation}
is stable, in that it fixes the cocycle up to isomorphism:
\vspace{-2mm}
\begin{equation}
  \label{StabilityForEquivariantPUHBundlesExpressedThroughEquivariantCechCocycles}
  \mbox{
    [c] is stable
  }
  \;\;\;\;\;
  \Leftrightarrow
  \;\;\;\;\;
  [c]
  \;=\;
  \big[
    c_{{}_{\mathrm{reg}}}
      \otimes
    c
  \big]
  \,.
\end{equation}
\end{definition}
\begin{lemma}[Stable projective isotropy representations; cf. {\cite{BEJU12}\cite[\S 15]{LueckUribe14}\cite[\S 5]{EspinozaUribe15}}]
  \label{StableProjectiveIsotorpyRepresentations}
  For $G \,\subset\, \SpOne$ a finite subgroup
  and $[\tau] \,\in\, H^2_{\mathrm{Grp}}(G;\, \CircleGroup)$,
  a $[\tau]$-projective unitary $G$-representation
  (Ex. \ref{ProjectiveRepresentationsAndTheirCentralExtensions})
  \vspace{-2mm}
  $$
    \begin{tikzcd}[row sep=1pt]
      \mathbf{B}G
      \ar[rr, "{ \rho^{\tau} }"]
      \ar[dr, "{ \tau }"{swap}]
      &&
      \HomotopyQuotient
        { \mathbf{B}\UH }
        { \mathbf{B} \CircleGroup }
      \ar[dl]
      \\
      &
      \mathbf{B}^2
      \CircleGroup
    \end{tikzcd}
  $$

  \vspace{-2mm}
  \noindent
  is stable in the sense of Def. \ref{BlowupStableADEEquivariantPUHPrimePrincipalBundles},
  in that it is stable under tensoring with $G \acts \, L^{2}(\SpOne)$,
  iff it is
  is isomorphic to the direct sum of all
  irreducible
  $[\tau]$-projective representations
  each appearing with countably infinite multiplicity.
  $$
    \mbox{
      $\rho^\tau$ stable
    }
    \qquad
    \Leftrightarrow
    \qquad
    \rho^\tau
    \;\simeq\qquad
    \underset{
      \mathclap{
      \scalebox{.7}{$
        \def\arraystretch{.9}
        \begin{array}{c}
          {[\mu] \in }
          \\
          \IsomorphismClasses{
            \Representations^{[\tau]}(G)_{\mathrm{irr}}
          }
        \end{array}
      $}
      }
    }{\bigoplus}
    \quad
    \mu^{\otimes^\infty}
    .
  $$
\end{lemma}
\begin{proof}
  In view of \eqref{FineStructurePeterWeylDecompositionOfLTwoSPOne}
  and by Schur's Lemma,
  the restriction map
  $$
    \begin{tikzcd}[row sep=1pt]
      \SpOne
      &&
      G
      \ar[ll, hook', "{ i }"{swap}]
      \\
      G \acts \;
      L^2\big(
        \SpOne
      \big)
      \ar[rr, ->>, "{ i^\ast }"]
      &&
      L^2\big(
        G
      \big)
      \mathrlap{
        \;\simeq\;
        G \acts \; \ComplexNumbers[G]
      }
    \end{tikzcd}
  $$
  exhibits the regular $G$-representation $\ComplexNumbers[G]$ as a direct summand
  of $L^2(\SpOne)$ (which, of course,  itself is the regular $\SpOne$-representation).

  But tensoring a $[\tau]$-projective representation with
  the regular representation adjoins
  a positive finite number of
  copies of all $[\tau]$-projective irreps,
  by Prop. \ref{TensoringProjectiveRepresentationsWithThePlainRegularRepresentation}.
  This operation is thus stable if and only if an infinite number of direct summands
  of each $[\tau]$-projective irrep is already present in $\rho^{[\tau]}$.
\end{proof}

\begin{example}[Blowup-stabilization of ADE-equivariant $\PUHPrime$-principal bundles]
  \label{BlowupStabilizationOfADEEquivariantPUHPrimeBundles}
  If we equip the infinite tensor Hilbert space
  \eqref{HilbertSpaceOfBosonsOnThreeSphere}
  with the diagonal canonical $G$-action
  \eqref{DiagonalCanonicalActionOnInfiniteGuichardetTensorProduct},
  then it is stable under tensoring with one copy of its
  $G$-tensor factors, regarded as equipped with its canonical $G$-action
  \eqref{GActionOnSingleBosonHilbertSpaceOverThreeSphere}:
   \vspace{-2mm}
  $$
    \begin{tikzcd}[column sep=-6pt]
      L^2(\RiemannianSphere{3})
      \ar[out=180-66, in=66, looseness=3.5, "\scalebox{.77}{$\mathclap{
        G
      }$}"{description},shift right=1]
      &\otimes&
      L^2(\RiemannianSphere{3})^{\otimes^\infty}
      \ar[out=180-66, in=66, looseness=3.5, "\scalebox{.77}{$\mathclap{
        G
      }$}"{description},shift right=1]
      &\;\;\;\simeq\;\;\;&
      L^2(\RiemannianSphere{3})^{\otimes^\infty}
      \ar[out=180-66, in=66, looseness=3.5, "\scalebox{.77}{$\mathclap{
        G
      }$}"{description},shift right=1]
    \end{tikzcd}
    \;\;
    \in
    \;\;
    \Actions{G}(\HilbertSpaces)
    \,.
  $$

   \vspace{-2mm}
\noindent
  More generally, for $G \acts \, \HilbertSpace$ any $G$-Hilbert space, its
  tensor product with $G \acts \, L^2(\RiemannianSphere{3})^{\otimes^\infty}$ is stable
  under tensoring with a single copy of the $G$-Hilbert space
  $L^2(\RiemannianSphere{3})$.
  In this manner, if we write,
  in analogy with \eqref{TheRegularCocycle},
  \vspace{-2mm}
  \begin{equation}
    \label{InfiniteGuichardetTensorProductOfRegularCocycle}
    \begin{tikzcd}
      c_{{}_{\mathrm{reg}}}^{\otimes^\infty}
      :
      &[-25pt]
      \Big(\!
        \big(
        \widehat{\TopologicalSpace}
          \times_{{}_{\TopologicalSpace}}
        \widehat{\TopologicalSpace}
        \big)
        \times
        G
        \,\rightrightarrows\,
        \widehat{\TopologicalSpace}
      \Big)
      \ar[
        r,
        ->>,
        "{
          p_{{}_{\widehat{\TopologicalSpace}}}
        }"
      ]
      &
      \big(
        \DeloopingGroupoid{ G }
      \big)
      \ar[
        rrr,
        "{\scalebox{0.7}{$ {
          \DeloopingGroupoid
            {
              \big(\left( (-)^{-1} \right)^\ast\big)^{\otimes^\infty}
            }
        }
        $}
        }
        "
      ]
      &&&
      \bigg(\!\!
        \Big(\,
          \overset{
            \UHPrime
          }{
          \overbrace{
            \UnitaryGroup
            \big(
              L^2(S^3)^{\otimes^\infty}
            \big)
          }
          }
        \,\Big)
      \!\!\bigg)
    \end{tikzcd}
  \end{equation}

   \vspace{-2mm}
\noindent
  then tensoring
  with this cocycle
  serves to stabilize
  $G$-equivariant $\PUH$- and $\PUHPrime$-principal bundles $[c]$
  according to
  Def. \ref{BlowupStableADEEquivariantPUHPrimePrincipalBundles})
  in that:
   \vspace{-2mm}
  \begin{equation}
    \label{BlowupStabilizationofPUHPrimeBundlesByTensoringWithInfiniteRegularCocycle}
    \big[
      c_{\mathrm{reg}}
      \otimes
      (
      c_{\mathrm{reg}}^{\otimes^\infty}
      \,\otimes\,
      c
      )
    \big]
    \;\;
    =
    \;\;
    \big[
      (
      c_{\mathrm{reg}}^{\otimes^\infty}
      \,\otimes\,
      c
      )
    \big]
    \;\;\;
    \in
    \IsomorphismClasses
    {
      \EquivariantPrincipalBundles{G}{\PUHPrime}(\SmoothInfinityGroupoids)
    }
    \,.
  \end{equation}
\end{example}

\medskip

\begin{theorem}[Blowup-stability of ADE-equivariant $\PUHPrime$-principal bundles]
\label{BlowupStabilityOfADEEquivariantPUHPrimePrincipalBundles}
For $G \subset \mathrm{Sp}(1)$ a finite subgroup (Ex. \ref{ADEGroupsHaveSphericalSpaceForms}):

\noindent
{\bf (i)}
The stable $G$-equivariant $\PUHPrime$-principal bundles
according to Def. \ref{BlowupStableADEEquivariantPUHPrimePrincipalBundles}
do form a blow-up stable class in the sense of
Ntn. \ref{StableEquivariantBundles} with respect to blowing up
by the Hopf action $G \acts \, \SmoothSphere{7}$ \eqref{QuaternionicHopfFibration}.

\noindent
{\bf (ii)}
Every stable ADE-equivariant $\PUHPrime$-principal bundle reduces
to a stable equivariant $\GradedPUH$-principal bundle along the
comparison map $\GradedPUH \to \GradedPUHPrime$ \eqref{InfiniteTensoringOfPUHIsWeakHomotopyEquivalence}.
\end{theorem}
\begin{proof}
With any equivariant good open cover chosen as in
Def. \ref{BlowupStableADEEquivariantPUHPrimePrincipalBundles},
we need to show
that the pullback operation on equivariant {\v C}ech cocycles
(Rem. \ref{EquivariantCechCocycles})
\vspace{-2mm}
$$
\hspace{-2mm}
  \begin{tikzcd}
  \IsomorphismClasses{
  \Functors
  \big(
    (
      \widehat{\SmoothManifold}
        \times_{{}_{\SmoothManifold}}
      \widehat{\SmoothManifold}
        \times
      G
      \rightrightarrows
      \widehat{\SmoothManifold}
    )
    ,\,
    (\PUH \rightrightarrows \ast)
  \big)
 \! }
  \;\;
  \ar[
    r,
    start anchor={[xshift=-12pt]},
    "{
      p_{{}_{\SmoothSphere{7}}}^\ast
    }"
  ]
  &
  \!\!
  \IsomorphismClasses{
  \Functors
  \Big(
    \big(
      (
        \widehat{\SmoothManifold}
          \times_{{}_{\SmoothManifold}}
        \widehat{\SmoothManifold}
          \times
        \SmoothSphere{7}
      )
        \times
      G
      \rightrightarrows
      \widehat{\SmoothManifold}
        \times
      \SmoothSphere{7}
    \big)
    ,\,
    (\PUHPrime \rightrightarrows \ast)
  \Big)
  \!}
  \end{tikzcd}
$$

\vspace{-2mm}
\noindent
restricts on stable equivariant $\PUHPrime$-principal bundles to a bijection
and on stable equivariant $\PUH$-principal bundles at least to a surjection.
(Here ``$\mathrm{Fnctr}$'' is short for the groupoid of morphisms of
simplicial presheaves, hence, in the present case, of functors
of D-topological groupoids, Ntn. \ref{TopologicalGroupoids}).

\medskip

\noindent {\bf (1)} {\it Injectivity:}
Given a coboundary $f$ between the pullback of a pair of cocycles
\vspace{-2mm}
\begin{equation}
  \label{CoboundaryBetweenCechCocyclesPulledBackToSevenSphere}
  \hspace{-2.6cm}
  \begin{tikzcd}[column sep=8pt]
    \Big(
      \widehat{\SmoothManifold}
        \times_{{}_{\SmoothManifold}}
      \widehat{\SmoothManifold}
        \times
      \SmoothSphere{7}
        \times
      G
      \rightrightarrows
      \widehat{\SmoothManifold}
      \times
      \SmoothSphere{7}
    \Big)
    \hspace{2cm}
    \ar[
      rr,
      bend left=20,
      "{
        p_{{}_{\SmoothSphere{7}}}^\ast c
      }"{},
      "{\ }"{name=s, swap}
    ]
    \ar[
      rr,
      bend right=20,
      "{
        p_{ {}_{\SmoothSphere{7}} }^\ast  c'
      }"{swap},
      "{\ }"{name=t}
    ]
    \ar[
      from=s,
      to=t,
      Rightarrow,
      "{ f }",
      "{ \sim }"{swap, left}
    ]
    &&
    \DeloopingGroupoid{
      \PUHPrime
    }
  \end{tikzcd}
\end{equation}
\begin{equation}
  \label{EquivariantPartOfCoboundaryBetweenPUHPrimeCocyclePulledBackToSevenSphere}
  \hspace{1.5cm}
  \begin{tikzcd}[row sep=4pt]
    (\widehat{x}, \, p)
    \ar[dd]
    &[+10pt]&[+30pt]
    \bullet
    \ar[
      rr,
      "{
        f(\widehat{x}, p)
      }"
    ]
    \ar[
      dd,
      "{
        \rho(\widehat{x}, g)
      }"{left}
    ]
    &[+20pt]&
    \bullet
    \ar[
      dd,
      "{
        \rho'(\widehat{x},\, g)
      }"{right}
    ]
    \\
    &
    \mapsto
    &
    \\
    g\cdot(\widehat{x}, \, p)
    &&
    \bullet
    \ar[
      rr,
      "{
        f(g\cdot \widehat{x},\, g \cdot p)
      }"
    ]
    &&
    \bullet
  \end{tikzcd}
  \;\;\;
  \in
  \;
  \PUHPrime
\end{equation}
\begin{equation}
  \label{CechPartOfCoboundaryBetweenPUHPrimeCocyclePulledBackToSevenSphere}
  \hspace{1.5cm}
  \begin{tikzcd}[row sep=4pt]
    \big(\widehat{x}, \, p\big)
    \ar[dd]
    &[+10pt]&[+30pt]
    \bullet
    \ar[
      rr,
      "{
        f(\widehat{x}, p)
      }"
    ]
    \ar[
      dd,
      "{
        c(\widehat{x},\, \widehat{x}\,')
      }"{left}
    ]
    &[+20pt]&
    \bullet
    \ar[
      dd,
      "{
        c'(\widehat{x},\, \widehat{x}\,')
      }"{right}
    ]
    \\
    &
    \;\;\;\;\mapsto
    &
    \\
    (\widehat{x}\,', \, p)
    &&
    \bullet
    \ar[
      rr,
      "{
        f(\widehat{x}\,',\, p)
      }"
    ]
    &&
    \bullet
  \end{tikzcd}
  \;\;\;
  \in
  \;
  \PUHPrime
\end{equation}
we will construct a coboundary $\overline{f}$ between the tensoring of the
original cocycles with $c_{\mathrm{reg}}$ \eqref{ComponentsOfTensoringWithRegularCocycle}:
\begin{equation}
  \label{IsomorphismBetweenPUHPrimeCocyclesBeforePullbackToSevenSphere}
  \hspace{-1.5cm}
  \begin{tikzcd}
    \Big(
      \widehat{\SmoothManifold}
        \times_{{}_{\SmoothManifold}}
      \widehat{\SmoothManifold}
        \times
      G
      \rightrightarrows
      \widehat{\SmoothManifold}
    \Big)
    \hspace{5mm}
    \ar[
      rr,
      bend left=20,
      "{
        c_{{}_{\mathrm{reg}}} \otimes c
      }"{pos=.4},
      "{\ }"{name=s, swap}
    ]
    \ar[
      rr,
      bend right=20,
      "{
        c_{{}_{\mathrm{reg}}} \otimes c'
      }"{swap, pos=.4},
      "{\ }"{name=t}
    ]
    \ar[
      from=s,
      to=t,
      Rightarrow,
      "{ \overline{f} }",
      "{ \sim }"{swap, left}
    ]
    &[+10pt]&
    \DeloopingGroupoid{
      \PUHPrime
    }
  \end{tikzcd}
\end{equation}
\begin{equation}
  \label{EquivariantPartOfCoboundaryBetweenPUHPrimeCocycleBeforePulledBackToSevenSphere}
  \hspace{1cm}
  \begin{tikzcd}[row sep=4pt]
    \widehat{x}
    \ar[dd]
    &[+10pt]&[+20pt]
    \bullet
    \ar[
      rr,
      "{
        \overline{f}(\widehat{x}\,)
      }"
    ]
    \ar[
      dd,
      "{
        \tilde{\rho}(\widehat{x}, g)
      }"{left}
    ]
    &[+20pt]&
    \bullet
    \ar[
      dd,
      "{
        \tilde{\rho}'(\widehat{x},\, g)
      }"{right}
    ]
    \\
    & \mapsto &
    \\
    g \cdot\widehat{x}
    &&
    \bullet
    \ar[
      rr,
      "{
        \overline{f}(g\cdot \widehat{x}\,)
      }"
    ]
    &&
    \bullet
  \end{tikzcd}
  \;\;\;
  \in
  \;
  \PUHPrime
\end{equation}
\begin{equation}
  \label{CechPartOfCoboundaryBetweenPUHPrimeCocycleBeforePulledBackToSevenSphere}
  \hspace{1cm}
  \begin{tikzcd}[row sep=4pt]
    \widehat{x}
    \ar[dd]
    &[+10pt]&[+20pt]
    \bullet
    \ar[
      rr,
      "{
        \overline{f}(\widehat{x}\,)
      }"
    ]
    \ar[
      dd,
      "{
        \tilde{c}(\widehat{x}, \widehat{x}\,')
      }"{left}
    ]
    &[+20pt]&
    \bullet
    \ar[
      dd,
      "{
        \tilde{c}\,'(\widehat{x},\, \widehat{x}\,')
      }"{right}
    ]
    \\
    &\mapsto&
    \\
    \widehat{x}\,'
    &&
    \bullet
    \ar[
      rr,
      "{
        \overline{f}(g\cdot \widehat{x}\,)
      }"
    ]
    &&
    \bullet
  \end{tikzcd}
  \;\;\;
  \in
  \;
  \PUHPrime
  \,.
\end{equation}
By the assumption that both $c$ and $c'$ are stable
\eqref{StabilityForEquivariantPUHBundlesExpressedThroughEquivariantCechCocycles}, this
will imply that they were already isomorphic before pullback.

\medskip

To motivate the following construction of an isomorphism
$\overline{f}$
\eqref{IsomorphismBetweenPUHPrimeCocyclesBeforePullbackToSevenSphere},
notice that {\it if} the $G$-action on $\SmoothSphere{7}$ {\it had}
a fixed point $p_0 \xhookrightarrow{\;} \SmoothSphere{7}$,
then we could simply restrict both cocycles along its inclusion
to immediately obtain the desired isomorphism. Now, instead
of a fixed point, the left $G$-action in question
has for each point $p_0 \,\in\, \SmoothSphere{4}$ its full 3-spherical
$\SpOne$-orbit $S^3_{p_0} \xhookrightarrow{\;} \SmoothSphere{7}$,
this being the
fiber of the quaternionic Hopf fibration \eqref{QuaternionicHopfFibration}
over this point:
\begin{equation}
  \label{FiberOfQuaternionicHopfFibrationOverBasePoint}
  \begin{tikzcd}[row sep=small, column sep=large]
    S^3_{p_0}
      \ar[r, "{ i_{p_0} }"]
      \ar[d]
      \ar[dr, phantom, "{\mbox{\tiny (pb)}}"]
      &
    \SmoothSphere{7}
    \ar[d]
    \\
    \ast \ar[r, "{p_0}"{swap}]
    &
    \SmoothSphere{4}
    \,.
  \end{tikzcd}
\end{equation}
Therefore the strategy is to ``evaluate'' at some $p_0 \,\in\, \SmoothSphere{4}$
and ``average'' the value of $f$ over the
remaining $S^3_{p_0} \xhookrightarrow{\;} \SmoothSphere{7}$ over this
point, by absorbing its operator values on this 3-sphere fiber
on the left of the higher spin chain:
\begin{equation}
  \label{AveragingOfPUHPrimeCoboundaryOverThreeSphere}
  \overline{f}(\widehat{x}\,)
  \;:=\;
  \underset{
    S^3_{p_0}
  }{
    \int
  }
  \vert p \rangle \langle p \vert
    \,\otimes\,
  {f}(\widehat{x},p)
  \;\;\;
  \in
  \;
  \ProjectiveUnitaryGroup
  \Big(
    L^2(\RiemannianSphere{3})
    \otimes
    L^2(\RiemannianSphere{3})^{\otimes^\infty}
  \Big)
  \,\simeq\,
  \PUHPrime
  \,.
\end{equation}

In order to check that this
does yield a coboundary as desired,
notice that we may lift the component function $f$ from a projective-unitary to a
unitary map $\widehat{f}$:
$$
  \begin{tikzcd}
    &&
    \UHPrime
    \ar[d]
    \ar[rr]
    \ar[
      drr,
      phantom,
      "{
        \mbox{\tiny(pb)}
      }"
    ]
    &&
    \ast
    \ar[d]
    \\
    \widehat{\SmoothManifold}
      \times
    \SmoothSphere{7}
    \ar[
      rr,
      "{
        f
      }"
    ]
    \ar[
      urr,
      "{
        \widehat{f}
      }",
      dashed
    ]
    &&
    \PUHPrime
    \ar[rr]
    &&
    \mathbf{B} \PhaseGroup
    \,.
  \end{tikzcd}
$$
Namely, under the universal property of the pullback on the right,
such a lift is induced by a
trivialization of the bottom composite map shown above,
and this exists by the classification theory for $\PUHPrime$-principal bundles
(Thm. \ref{ClassificationOfPrincipalBundlesAmongPrincipalInfinityBundles})
and since $\shape \PUHPrime$ is 2-truncated
(by Lem. \ref{ShapeOfProjectiveUnitaryGroupOnBosonStatesOver3Sphere}),
while $\shape (\widehat{\SmoothManifold} \times \SmoothSphere{7})$
is is a disjoint union of copies of $\ShapeOfSphere{7}$.

Using an analogously constructed lift for $\rho(-,-)$, the equivariance condition
\eqref{EquivariantPartOfCoboundaryBetweenPUHPrimeCocyclePulledBackToSevenSphere} says
that the two products of $\widehat f$ and $\widehat{\rho}$ make
two lifts to $\UHPrime$ of the same map $\PUHPrime$, which hence differ by a
map to $\PhaseGroup$. This, of course, remains the case after
restriction to $S^3_{p_0} \xhookrightarrow{\;} S^7$, where we hence have:
\begin{equation}
  \label{LiftedEquivarianceConditionOnCocycleTransformation}
  \widehat{f}\big(\widehat{x},\,(-) \big)
  \cdot
  \widehat{\rho}'(\widehat{x},\,g)
  \;=\;
  \widehat{\rho}(\widehat{x}, \, g)
  \cdot
  \widehat{f}\big(g \cdot \widehat{x},\,g \cdot (-) \big)
  \;\;\,
  \mathrm{mod}
  \;\;\,
  \Maps{}
    { S^3_{p_0} }
    { \PhaseGroup }
  \xrightarrow{
    \mbox{
      \tiny
      \eqref{ReIncludingMapsIntoModifiedCircleGroupIntoModifiedCircleGroup}
    }
  }
  \PhaseGroup
  \,.
\end{equation}
In fact, the projectivizing maps on the right exist in smooth dependence
on $\widehat{x}$.

Using this, we may verify that $\overline{f}$ \eqref{AveragingOfPUHPrimeCoboundaryOverThreeSphere}
satisfies condition
\eqref{EquivariantPartOfCoboundaryBetweenPUHPrimeCocycleBeforePulledBackToSevenSphere}:
\begin{equation}
  \label{CheckingThatAveragedPUHPrimeCocycleIsStillEquivariant}
  \def\arraystretch{1.7}
  \begin{array}{lll}
    \overline{f}(\widehat{x}\,)
    \cdot
    \widetilde{\rho}'(\widehat{x},g)
    &
    \;=\;
    \overline{f}(\widehat{x}\,)
    \cdot
    \underset{S^3_{p_0}}{\displaystyle \int}
    \vert g^{-1} \cdot p \rangle \langle p \rvert
      \,\otimes\,
    \rho'(\widehat{x},g)
    &
    \proofstep{
      by \eqref{ComponentsOfTensoringWithRegularCocycle}
    }
    \\
    & \;=\;
    \Bigg[\;\;
    \underset{S^3_{p_0}}{\displaystyle  \int}
      \vert g^{-1} \cdot p \rangle \langle p \rvert
      \,\otimes\,
      \widehat{f}(\widehat{x}, g^{-1} \cdot p)
      \cdot
      \widehat{\rho}'(\widehat{x},g)
    \Bigg]
    &
    \proofstep{
      by \eqref{AveragingOfPUHPrimeCoboundaryOverThreeSphere}
    }
    \\
    & \;=\;
    \Bigg[\;\;
    \underset{S^3_{p_0}}{\displaystyle  \int}
      \vert g^{-1} \cdot p \rangle \langle p \rvert
      \,\otimes\,
      \widehat{\rho}(\widehat{x},g)
      \cdot
      \widehat{f}(g\cdot \widehat{x}, p)
    \Bigg]
    &
    \proofstep{
      by \eqref{LiftedEquivarianceConditionOnCocycleTransformation}
    }
    \\
    & \;=\;
    \widetilde{\rho}(\widehat{x},g)
    \cdot
    \overline{f}(g\cdot \widehat{x}\,)
    &
    \proofstep{
      by \eqref{ComponentsOfTensoringWithRegularCocycle}
      \&
      \eqref{AveragingOfPUHPrimeCoboundaryOverThreeSphere}.
    }
  \end{array}
\end{equation}
Here square brackets denote $\PhaseGroup$-equivalence classes of $\UHPrime$-operators.

The remaining condition \eqref{CechPartOfCoboundaryBetweenPUHPrimeCocycleBeforePulledBackToSevenSphere}
is verified analogously:
\begin{equation}
  \label{CheckingThatAveragedPUHPrimeCocycleIsStillACechCocycle}
  \def\arraystretch{1.7}
  \begin{array}{lll}
    \overline{f}(\widehat{x}\,)
    \cdot
    \widetilde{c}\;'(\widehat{x},  \; \widehat{x}\;')
    &
    \;=\;
    \Bigg[\;\;
      \underset{S^3_{p_0}}{\displaystyle  \int}
      \lvert p \rangle \langle p \rvert
      \,\otimes\,
      \widehat{f}(\widehat{x}, p)
      \cdot
      \widehat{c}'(\widehat{x}, \;\widehat{x}\;')
    \Bigg]
    &
    \proofstep{
      by \eqref{ComponentsOfTensoringWithRegularCocycle}
      \& \eqref{AveragingOfPUHPrimeCoboundaryOverThreeSphere}
    }
    \\
    &
    \;=\;
    \Bigg[\;\;
      \underset{S^3_{p_0}}{\displaystyle  \int}
      \lvert p \rangle \langle p \rvert
      \,\otimes\,
      \widehat{c}'(\widehat{x},  \;\widehat{x}\;')
      \cdot
      \widehat{f}(\widehat{x}\;', p)
    \Bigg]
    &
    \proofstep{
      by \eqref{CechPartOfCoboundaryBetweenPUHPrimeCocyclePulledBackToSevenSphere}
    }
    \\
    &
    \;=\;
    \widetilde{c}\;'(\widehat{x}, \; \widehat{x}\;')
    \cdot
    \overline{f}(\widehat{x}\;')
    &
    \proofstep{
      by \eqref{ComponentsOfTensoringWithRegularCocycle}
      \&
      \eqref{AveragingOfPUHPrimeCoboundaryOverThreeSphere}
    }
    \,.
  \end{array}
\end{equation}

\medskip
\noindent {\bf (2)}
{\it Surjectivity:}
First, consider the further pullback of the equivariant bundles
to the blowup by just the
3-spherical fiber of the Hopf fibration over any chose basepoint
$p_0 \,\in\, \SmoothSphere{4}$, as in \eqref{FiberOfQuaternionicHopfFibrationOverBasePoint},
hence all the way
to the last item of the following chain of morphisms:
\begin{equation}
  \label{ConsecutiveMapsFromThreeSphericalBlowup}
  \begin{tikzcd}[column sep=50pt]
    \HomotopyQuotient
      { \SmoothManifold }
      { G }
    \ar[
      from=r,
      "{
        \mathrm{pr}_1
      }"{swap}
    ]
    &
    \HomotopyQuotient
      { \SmoothManifold }
      { G }
      \times
    \SmoothSphere{4}
    \ar[
      from=r,
      "{
        \Quotient
          { ( \SmoothManifold \times t_{\Quaternions} ) }
          { G }
      }"{swap}
    ]
    &
    \Quotient
      { ( \SmoothManifold \times \SmoothSphere{7} ) }
      { G }
    \ar[
      from=r,
      "{
        \Quotient
          { ( \SmoothManifold \times i_{p_0} ) }
          { G }
      }"{swap}
    ]
    &
    \Quotient
      { \big( \SmoothManifold \times (\SmoothSphere{3})_{p_0} \big) }
      { G }\;
  \end{tikzcd}
  \!\!\!
  : p_{\SmoothSphere{3}}.
\end{equation}
Observe that on this 3-spherical blowup, the regular
unitary {\v C}ech cocycle $c_{\mathrm{reg}}$
\eqref{TheRegularCocycle}
and its infinite tensor product
$c^{\otimes^\infty}_{\mathrm{reg}}$ \eqref{InfiniteGuichardetTensorProductOfRegularCocycle}
both have a canonical trivialization after identifying $\SmoothSphere{3} \,\simeq\, \SpOne$:
$$
  \hspace{-2cm}
  \begin{tikzcd}[column sep=10pt]
    \big(
    \widehat{\SmoothManifold}
      \times_{{}_{\SmoothManifold}}
    \widehat{\SmoothManifold}
    \times
    \SpOne
    \rightrightarrows
    \widehat{\SmoothManifold}
    \times
    \SpOne
    \big)
    \hspace{2cm}
    \ar[
      rr,
      bend left=20,
      "{ c_{\mathrm{reg}} }"{pos=.4},
      "{\ }"{swap, name=s}
    ]
    \ar[
      rr,
      bend right=20,
      "{ \mathrm{const} }"{swap, pos=.4},
      "{\ }"{name=t}
    ]
    \ar[
      from=s,
      to=t,
      Rightarrow,
      "\sim"{sloped}
    ]
    &&
    \HilbertSpaces
  \end{tikzcd}
$$
$$
  \hspace{1cm}
  \begin{tikzcd}[row sep=4pt]
    p
    \ar[dd]
    &&
    L^2(\SpOne)^{\otimes^\infty}
    \ar[
      rr,
      "{
        (
          p^{-1}
        )^\ast
      }"
    ]
    \ar[dd, "{ (g^{-1})^\ast }"{left}]
    &&
    L^2(\SpOne)^{\otimes^\infty}
    \ar[dd, -, shift left=1pt]
    \ar[dd, -, shift right=1pt]
    \\
    &\longmapsto&
    \\
    g \cdot p
    &&
    L^2(\SpOne)^{\otimes^{\infty}}
    \ar[
      rr,
      "{
        \left(
          (g \cdot p)^{-1}
        \right)^\ast
      }"
    ]
    &&
    L^2(\SpOne)^{\otimes^\infty}
  \end{tikzcd}
$$
By tensoring \eqref{TensoringOfPUHPrimeOverUHPrime}
with this trivialization and
using that tensoring $\PUHPrime$
with the constant $L^2(\RiemannianSphere{3})^{\otimes^{\infty}}$
is a weak homotopy equivalence
\eqref{InfiniteTensoringOfPUHIsWeakHomotopyEquivalence},
it follows that
for every isomorphism class
$[c]$
of a $G$-equivariant $\PUHPrime$-principal bundle
over $\SmoothManifold$, its
pullback to $\SmoothManifold \times S^3_{p_0}$
coincides with the pullback of its stabilization
$[c_{\mathrm{reg}}^{\otimes^\infty} \otimes c]$ \eqref{BlowupStabilizationofPUHPrimeBundlesByTensoringWithInfiniteRegularCocycle}:
$$
  p_{\SmoothSphere{3}}^\ast
  \big[
    c_{\mathrm{reg}}^{\otimes^\infty}
    \otimes
    c
  \big]
  \;\;
  =
  \;\;
  p_{\SmoothSphere{3}}^\ast
  [c]
  \;\;\;\;
  \in
  \;
  \Truncation{0}
  \PointsMaps{\big}
    { \Quotient{ (\SmoothManifold \times \SmoothSphere{3}) }{G} }
    { \PUHPrime }
  \,.
$$
But since the intermediate pullback
along $i_{p_0}$ \eqref{ConsecutiveMapsFromThreeSphericalBlowup}
is an injection on isomorphism classes of $\PUHPrime$-principal bundles
(Lem. \ref{PullbackOfPUHBundlesFromSevenToThreeSphericalBlowupIsInjection})
\begin{equation}
  \label{IntermediatePullbackIsInjectionOnIsomorphism}
  \begin{tikzcd}
    \Truncation{0}
    \PointsMaps{\big}
      { \Quotient{ (\SmoothManifold \times \SmoothSphere{7}) }{G} }
      { \PUHPrime }
    \ar[
      rrr,
      hook,
      "{
        \big(
          \Quotient{ (\SmoothManifold \times i_{p_0}) }{G}
        \big)^\ast
      }"
    ]
    &&&
    \Truncation{0}
    \PointsMaps{\big}
      { \Quotient{ (\SmoothManifold \times \SmoothSphere{3}) }{G} }
      { \PUHPrime }
    \,,
  \end{tikzcd}
\end{equation}
we find that this equality between classes of equivariant bundles and their stabilization
must hold already after pullback to the blowup by $S^7$:
$$
  p_{\SmoothSphere{7}}^\ast
  \big[
    c_{\mathrm{reg}}^{\otimes^\infty}
    \otimes
    c
  \big]
  \;\;
  =
  \;\;
  p_{\SmoothSphere{7}}^\ast
  [c]
  \;\;\;\;
  \in
  \;
  \Truncation{0}
  \PointsMaps{\big}
    { \Quotient{ (\SmoothManifold \times \SmoothSphere{3}) }{G} }
    { \PUHPrime }
  \,.
$$
Since Lem. \ref{PrincipalBundlesOnBlowUps} says that
$p_{\SmoothSphere{7}}^\ast[-]$ is surjective, we conclude that already
$p_{\SmoothSphere{7}}^\ast\big[ c_{\mathrm{reg}}^{\otimes^\infty} \otimes (-) \big]$
is surjective. Since $c_{\mathrm{reg}}^{\ast} \otimes (-)$ produces stable
bundles (Ex. \ref{BlowupStabilizationOfADEEquivariantPUHPrimeBundles})
this is the surjectivity statement to be proven, for structure group $\PUHPrime$.
Moreover, this surjectivity argument applies verbatim also to
the structure group $\PUH$.
\end{proof}

\begin{lemma}[Pullback of $\PUH$- and $\PUHPrime$-bundles from 7-spherical to 3-spherical blowup is injection]
  \label{PullbackOfPUHBundlesFromSevenToThreeSphericalBlowupIsInjection}
  Pullback of
  equivariant bundles with structure group
  $\PUH$ or $\PUHPrime$
  along $i_{p_0}$ \eqref{ConsecutiveMapsFromThreeSphericalBlowup}
  is an injection on isomorphism classes:
$$
  \begin{tikzcd}
    \Truncation{0}
    \,
    \PointsMaps{\big}
      {
        \Quotient
          { ( \SmoothManifold \times S^7  ) }
          { G }
      }
      { \mathbf{B} \PUHPrimeOptional }
    \ar[
      rr,
      "{
        (i_{p_0})^\ast
      }",
      hook
    ]
    \ar[d,-,shift left=1pt]
    \ar[d,-,shift right=1pt]
    &&
    \Truncation{0}
    \,
    \PointsMaps{\big}
      {
        \Quotient
          { ( \SmoothManifold \times S^3  ) }
          { G }
      }
      { \mathbf{B} \PUHPrimeOptional }
    \ar[d,-,shift left=1pt]
    \ar[d,-,shift right=1pt]
    \\
    H^3
    \big(
      \HomotopyQuotient
        { \shape \, {\SmoothManifold} }
        { G }
      ;\,
      \mathbb{Z}
    \big)
    \ar[
      rr,
      hook,
      "{
        (\mathrm{id}, 0)
      }"
    ]
    &&
    H^0
    \big(
      \HomotopyQuotient
        { \shape \, {\SmoothManifold} }
        { G }
      ;\,
      \mathbb{Z}
    \big)
    \,\oplus\,
    H^3
    \big(
      \HomotopyQuotient
        { \shape \, {\SmoothManifold} }
        { G }
      ;\,
      \mathbb{Z}
    \big).
  \end{tikzcd}
$$
\end{lemma}
\begin{proof}
  Since the $G$-action both on $S^7$ as well as on $S^3$ are free,
  the standard classification of principal bundles over smooth manifolds
  (Thm. \ref{ClassificationOfPrincipalBundlesAmongPrincipalInfinityBundles})
  applies, and since
  $B \PUHPrime \,\simeq\, B  \,\shape \, \PUHPrime \,\simeq\, B^3 \Integers$
  (Lem. \ref{ShapeOfProjectiveUnitaryGroupOnBosonStatesOver3Sphere})
  this says that the classification for $\PUH'$-structure
  is the same as that of $\PUH$-bundles,
  given by integral 3-cohomology of the shape of the domain manifold.

  First observe, from the long exact sequences of homotopy groups
  induced by the fiber sequences
$$
  \begin{tikzcd}[row sep=small]
    \mathllap{
      0
      \;\simeq\;
    }
    \,
    \Omega^7 B^3 \Integers
    \ar[r]
    &
    \Maps{}
      { \ShapeOfSphere{7} }
      { B^3 \mathbb{Z} }
    \,,
    \ar[d, "{ \mathrm{ev}_\ast }"]
    \\
    &
    B^3 \Integers
  \end{tikzcd}
  \hspace{2cm}
  \begin{tikzcd}[row sep=small]
    \mathllap{
      \Integers
      \;\simeq\;
    }
    \,
    \Omega^3 B^3 \Integers
    \ar[r]
    &
    \Maps{}
      { \ShapeOfSphere{3} }
      { B^3 \Integers }
    \ar[d, "{ \mathrm{ev}_\ast }"]
    \\
    &
    B^3 \Integers
    \mathrlap{\,.}
  \end{tikzcd}
$$
  that we have equivalences
  \begin{equation}
    \label{MappingSpacesFromThreeOrSevenSphereIntoBThreeZ}
    G \acts \;
    \Maps{}
      { \ShapeOfSphere{7} }
      { B^3 \Integers }
    \;\simeq\;
    B^3 \Integers
    \;=:\;
    A_7
    \,,
    \qquad
    G \acts \;
    \Maps{}
      { \ShapeOfSphere{3} }
      { B^3 \Integers }
    \;\simeq\;
    \Integers \times B^3 \Integers
    \;=:\;
    A_3
    \;\;\;
    \in
    \;
    \Actions{G}(\InfinityGroupoids)
    \,,
  \end{equation}
  which identify the induced $G$-action (via the given $G$-action on the sphere)
  as trivial: In both cases the factor $B^3 \Integers$ comes from the maps
  that factor through the point, while the factor $\Integers$ is the
  winding of $\ShapeOfSphere{3}$ over itself under any factorization through a
  generator $\ShapeOfSphere{3} \to B^3 \Integers$  of $\pi_3(B^3\Integers) \simeq \Integers$.
  Both of these are manifestly $G$-invariant. Moreover, this shows that
  the pullback map between these mapping spaces is the inclusion of the
  0-winding sector:
  \begin{equation}
    \label{InclusionOfMapsFromSeveanSphereToThoseFromThreeSphereToBThreeZ}
    \begin{tikzcd}
      \Maps{}
        { \ShapeOfSphere{7} }
        { B^3 \Integers }
      \ar[
        rr,
        "{ (i_{p_0})^\ast }"
      ]
      \ar[d, "{\sim}"{sloped}]
      &&
      \Maps{}
        { \ShapeOfSphere{3} }
        { B^3 \Integers }
      \ar[d, "{\sim}"{sloped}]
      \\
      \mathllap{
        A_7 \,=\;
      }
      B^3 \Integers
      \ar[
        rr,
        hook,
        "{ a \,\mapsto\, (0, a) }"
      ]
      &&
      \Integers \times B^3 \Integers
      \mathrlap{
        \; =\;
        A_3
      }
      \,.
    \end{tikzcd}
  \end{equation}
  Now, for $d \in \{3,7\}$, we compute much as in Lem. \ref{ConcordanceInfinityGroupoidOfPrincipalBundlesOnBlowupForStructureGroupOfTruncatedShape}:
  $$
    \def\arraystretch{1.6}
    \begin{array}{lll}
      \Maps{\Big}
        { \shape \, \Quotient{ \big( \SmoothManifold \times \SmoothSphere{\,d} \big) }{G} }
        { B^3 \Integers }
      &
      \;\simeq\;
      \Maps{\big}
        { \HomotopyQuotient{  ( \shape \, \SmoothManifold \times \ShapeOfSphere{d} ) }{G} }
        { B^3 \Integers }
      \\
      &
      \;\simeq\;
      \Maps{\Big}
        {
          \colimit{ [n] \in \Delta^\op }
          \,
          G^{\times^n}
          \times
          \,
          \shape \, \SmoothManifold \times \ShapeOfSphere{d}
        }
        { B^3 \Integers }
      \\
      &
      \;\simeq\;
      \Maps{\Big}
        {
          \colimit{ [n] \in \Delta^\op }
          \,
          \underset{ G^{n} }{\coprod}
          \,
          \shape\, \SmoothManifold \times \ShapeOfSphere{d}
        }
        { B^3 \Integers }
      &
      \proofstep{
        by discreteness of $G$
      }
      \\
      & \;\simeq\;
      \limit{ [n] \in \Delta^\op }
      \,
      \underset{ G^n }{\prod}
      \,
      \Maps{\big}
        {
          \shape \, \SmoothManifold
        }
        { B^3 \Integers }
      &
      \proofstep{
        by \eqref{HomFunctorRespectsLimits}
      }
      \\
      & \;\simeq\;
      \limit{ [n] \in \Delta^\op }
      \,
      \underset{ G^n }{\prod}
      \,
      \Maps{\big}
        {
          \shape \,\SmoothManifold
        }
        {
          \Maps{}
            { \ShapeOfSphere{d} }
            { B^3 \Integers }
        }
      &
      \\
      &
      \;\simeq\;
      \limit{ [n] \in \Delta^\op }
      \,
      \underset{ G^n }{\prod}
      \,
      \Maps{\big}
        {
          \shape \, \SmoothManifold
        }
        {
          A_d
        }
      &
      \proofstep{
        by \eqref{MappingSpacesFromThreeOrSevenSphereIntoBThreeZ}
      }
      \\
      & \;\simeq\;
      \,
      \Maps{\Big}
        {
          \colimit{ [n] \in \Delta^\op }
          \,
          \underset{ G^n }{\coprod}
          \,
          \shape \, \SmoothManifold
        }
        {
          A_d
        }
      &
      \proofstep{
        by \eqref{HomFunctorRespectsLimits}
        \&
        \eqref{MappingSpacesFromThreeOrSevenSphereIntoBThreeZ}
      }
      \\
      & \;\simeq\;
      \,
      \Maps{\big }
        {
          \HomotopyQuotient
            { \shape \, \SmoothManifold }
            { G }
        }
        {
          A_d
        }.
      &
     \end{array}
  $$
  The claim follows by appying this natural equivalence to
  \eqref{InclusionOfMapsFromSeveanSphereToThoseFromThreeSphereToBThreeZ}.
\end{proof}

With this concrete cocycle model for stability of equivariant $\PUHPrime$-principal
bundles in hand, we may now easily generalize to
projective unitary operators on
graded Hilbert spaces and/or
to their semidirect product
with the action of complex conjugation.

\begin{notation}[Grading and Complex conjugation action on $\PUHPrime$]
  \label{ComplexConjugationActionOnPUH}
 $\,$

 \noindent {\bf (i)}  In direct generalization of \eqref{ComplexConjugationActionOnProjectiveUnitaryGroup},
  the $\ZTwo$-action by complex conjugation
  on all the groups in Ex. \ref{InclusionOfStandardIntoModifiedProjectiveUnitaryGroupIsWeakEquivalence}
  is clearly compatible with all the morphisms there.
  In particular, we have
  $\ZTwo$-equivariant weak homotopy equivalences
  $$
    \begin{tikzcd}
      \PUH
      \ar[out=180-66, in=66, looseness=3.5, "\scalebox{.77}{$\;\mathclap{
        \ZTwo
      }\;$}"{description},shift right=1]
      \ar[
        rr,
        "{
          \in \WeakHomotopyEquivalences
        }"{swap}
      ]
      &&
      \PUHPrime
      \ar[out=180-66, in=66, looseness=3.5, "\scalebox{.77}{$\;\mathclap{
        \ZTwo
      }\;$}"{description},shift right=1]
      \ar[
        rrrr,
        "{
          \UnitaryGroup
          \big(
            L^2(\RiemannianSphere{3})^{\otimes^\infty}
          \big)
          \otimes
          (-)
        }",
        "{
          \in \WeakHomotopyEquivalences
        }"{swap}
      ]
      &&&&
      \PUHPrime
      \ar[out=180-66, in=66, looseness=3.5, "\scalebox{.77}{$\;\mathclap{
        \ZTwo
      }\;$}"{description},shift right=1]
    \end{tikzcd}
    \;\;\;
    \in
    \;
    \Actions{\ZTwo}
    \big(
      \Groups(\kTopologicalSpaces)
    \big)
  $$
  and the corresponding semidirect product groups
  $$
    \begin{tikzcd}
      \PhaseGroup \rtimes \ZTwo
      \ar[r, hook]
      &
      \UHPrime \rtimes \ZTwo
      \ar[r, ->>]
      &
      \PUHPrime \rtimes \ZTwo
      \;\;\;
      \in
      \;
      \Groups(\DTopologicalSpaces)
      \,.
    \end{tikzcd}
  $$
 \noindent {\bf (ii)}   Moreover, all the above constructions on the Hilbert space
  $L^2(\RiemannianSphere{3})^{\otimes^\infty}$ evidently generalize
  to constructions on its tensor product with any
  other Hilbert space, in particular to
  \begin{equation}
    \label{GradedVersionOfInfinitL2HilbertSpace}
    L^2(\RiemannianSphere{3})^{\otimes^\infty}
    \,\otimes\,
    \ComplexNumbers^2
    \;\simeq\;
    L^2(\RiemannianSphere{3})^{\otimes^\infty}_+
    \,\oplus\,
    L^2(\RiemannianSphere{3})^{\otimes^\infty}_-
    \,.
  \end{equation}
\noindent {\bf (iii)}  Under this identification we have,
  in analogous variation of $\GradedPUH$
  \eqref{TheGroupGradedPUH},
  the modified graded projective unitary group and its semidirect product
  with complex conjugation:
  \begin{equation}
    \label{TheGroupGradedPUHPrime}
    \begin{tikzcd}[row sep=2pt]
      \PhaseGroup \rtimes \ZTwo
      \ar[r, hook]
      &
      \GradedUHPrime \rtimes \ZTwo
      \ar[r, ->>]
      &
      \GradedPUHPrime \rtimes \ZTwo
      \;\;\;
      \in
      \;
      \Groups(\DTopologicalSpaces)
      \,.
    \end{tikzcd}
  \end{equation}
   \noindent {\bf (iv)} Finally, the tensoring operation \eqref{TensoringOfPUHPrimeOverUHPrime}
  evidently lifts to these semidirect products
  \begin{equation}
    \label{TensoringOfPUHSemidirectZTwoWithUH}
    \begin{tikzcd}[row sep=6pt]
      \UnitaryGroup
      \big(
        L^2(\RiemannianSphere{3})
      \big)
        \times
      \big(
        \GradedPUHPrime
        \rtimes
        \ZTwo
      \big)
      \ar[
        rr,
        "{
          \otimes
        }"
      ]
      &&
      \GradedPUHPrime \rtimes \ZTwo
      \\
      \GradedUHPrime
        \times
      \big(
        \GradedPUHPrime
        \rtimes
        \ZTwo
      \big)
      \ar[
        rr,
        "{
          \otimes
        }"
      ]
      &&
      \GradedPUHPrime \rtimes \ZTwo
    \end{tikzcd}
  \end{equation}
\end{notation}

\begin{lemma}[Shape of graded projective groups]
  \label{ShapeOfGradedProjectiveGroups}
  The shape of the graded projective groups $\GradedPUH$
  \eqref{TheGroupGradedPUH}
  and $\GradedPUHPrime$
  \eqref{TheGroupGradedPUHPrime} is:
  $$
    \shape \,\GradedPUH\,
    \;\simeq\;
    \shape \,\GradedPUHPrime\,
    \;\simeq\;
    \ZTwo
    \,\times\,
    B^2 \Integers
    \;\;\;
    \in
    \;
    \Groups(\InfinityGroupoids)
    \,;
  $$
  and under this equivalence the complex conjugation action
  \eqref{ComplexConjugationActionOnProjectiveUnitaryGroup}
  becomes the sign involution $\ZTwo \acts \Integers_{\conjugation}$:
  $$
    \begin{tikzcd}[column sep=-7pt]
    \shape
    &
    \GradedPUH
    \ar[out=180-66, in=66, looseness=3.5, "\scalebox{.77}{$\;\mathclap{
      \ZTwo
    }\;$}"{description},shift right=1]
    &
    \;\simeq\;
    &
    \ZTwo
    \,\times\,
    &
    B^2 \Integers_{\conjugation}
    \ar[out=180-66, in=66, looseness=3.5, "\scalebox{.77}{$\;\mathclap{
      \ZTwo
    }\;$}"{description},shift right=1]
    \end{tikzcd}
    \in
    \;
    \Actions{\ZTwo}
    \big(
      \Groups(\InfinityGroupoids)
    \big)
    \,.
  $$
\end{lemma}
\begin{proof}
  This is almost the same computation as in
  \eqref{AComputation}, the only difference being that
  $\shape \, \UH \,\simeq\, \ast$
  gets replaced by $\shape \, \GradedUH \,\simeq\, \ZTwo$
  (on which the circle quotient action is trivial, by definition):
  \begin{equation}
    \def\arraystretch{1.4}
    \begin{array}{lll}
      \shape \, \GradedPUHPrime
      &
      \;=\;
      \shape \big( \Quotient{\GradedUHPrime}{\PhaseGroup}  \big)
      &
      \proofstep{
        by
        \eqref{TheGroupGradedPUHPrime}
      }
      \\
      &
      \;\simeq\;
      \shape \big( \HomotopyQuotient{\GradedUHPrime}{\PhaseGroup}  \big)
      &
      \proofstep{
        by Prop. \ref{HomotopyQuotientOfFreeActionsIsOrdinaryQuotientOfZeroTruncation}
      }
      \\
      &
      \;\simeq\;
      \HomotopyQuotient{( \shape\,\GradedUHPrime) }{ (\shape\,\PhaseGroup)}
      &
      \proofstep{
        by \eqref{InfinityAdjointPreservesInfinityLimits}
      }
      \\
      &
      \;\simeq\;
      \HomotopyQuotient
        {
          \big(
            (\shape \ZTwo )
               \times
             ( \shape \UH )
               \times
             ( \shape \UH )
          \big)
        }
        { (\shape\,\PhaseGroup) }
      &
      \proofstep{
        by \eqref{ShapePreservesBinaryProducts}
      }
      \\
      & \;\simeq\;
      \HomotopyQuotient{ \ZTwo }{ B \Integers  }
      &
      \proofstep{
        by \eqref{UHIsContractible}
        \& \eqref{ObtainingShapeOfModifiedCircleGroup}
      }
      \\
      &
      \;\simeq\;
      \ZTwo
      \times
      B^2 \Integers
      &
      \proofstep{
        by
        \eqref{DeloopingOfInfinityGroupAsColimit}
      }
      \,.
    \end{array}
  \end{equation}
  The second statement follows by observing that under the
  equivalence
  $$
    \CircleGroup
      \,\simeq\,
    \HomotopyQuotient{ \mathbb{R}}{\Integers}
  $$
  the conjugation action on the left becomes sign inversion
  on real (hence on integer) numbers on the right.
\end{proof}

Now in direct generalization of Def. \ref{BlowupStableADEEquivariantPUHPrimePrincipalBundles}:
\begin{definition}[Blowup-stable ADE-equivariant $\GradedPUHPrime \rtimes \ZTwo$-principal bundles]
\label{BlowupStableADEEquivariantPUHPrimeRtimesZTwoBundles}
  We say that an ADE-equivariant principal bundle
  with structure group $\GradedPUH \rtimes \ZTwo$
  or $\GradedPUHPrime \rtimes \ZTwo$ (Ntn. \ref{ComplexConjugationActionOnPUH})
  is {\it stable} if its class is invariant under the tensoring
  \eqref{TensoringOfPUHSemidirectZTwoWithUH} with the
  regular cocycle \eqref{TheRegularCocycle}.
\end{definition}

In direct generalization of Thm. \ref{BlowupStabilityOfADEEquivariantPUHPrimePrincipalBundles},
we have:
\begin{theorem}[Blowup-stablility of ADE-equivariant $\GradedPUHPrime \rtimes \ZTwo$-principal bundles]
\label{BlowupStabilityForADEEquivariantPUHSemidirectProductZTwoPrincipalBundles}
For $G \subset \mathrm{Sp}(1)$ a finite subgroup (Ex. \ref{ADEGroupsHaveSphericalSpaceForms}),

\noindent
{\bf (i)}
The $G$-equivariant principal bundles
with structure group $\GradedPUHPrime \rtimes \ZTwo$
(Ntn. \ref{ComplexConjugationActionOnPUH})
which are stable
according to Def. \ref{BlowupStableADEEquivariantPUHPrimeRtimesZTwoBundles}
do form a blow-up stable class in the sense of
Ntn. \ref{StableEquivariantBundles} with respect to blowing up
by the Hopf action $G \acts \, \SmoothSphere{7}$ \eqref{QuaternionicHopfFibration}.

\noindent
{\bf (ii)}
Every stable ADE-equivariant $\GradedPUHPrime \rtimes \ZTwo$-principal bundle reduces
to a stable equivariant $\PUH \rtimes \ZTwo$-principal bundle along the
comparison map $\PUH \to \PUHPrime$ \eqref{InfiniteTensoringOfPUHIsWeakHomotopyEquivalence}.
\end{theorem}
\begin{proof}
  We just need to observe that the evident generalization of the proof of
  Thm. \ref{BlowupStabilityOfADEEquivariantPUHPrimePrincipalBundles}
  still goes through.
  First, since the argument there is manifestly natural in forming
  tensor products with the underlying Hilbert space, it
  evidently generalizes to the graded structure groups $\GradedPUH$
  and $\GradedPUHPrime$.
  Moreover, since  $\ZTwo$ is discrete and $\SmoothSphere{7}$ is connected,
  the equivariant cocycles and
  the coboundary \eqref{CoboundaryBetweenCechCocyclesPulledBackToSevenSphere}
  now have components which are pairs
  $$
    \big(
      \rho(\widehat{x}, p), \sigma(\widehat{x}\,)
    \big)
    \,,\;\;\;
    \big(
      f(\widehat{x}, p), \kappa(\widehat{x}\,)
    \big)
    \;\;\;
    \in
    \;
    \GradedPUHPrime \rtimes \ZTwo
  $$
  with (not only $\sigma$ but also) $\kappa$ not depending on
  $p \,\in\, \SmoothSphere{7}$.
  This implies that their conjugation actions are compatible with
  the integrals in the computations
  \eqref{CheckingThatAveragedPUHPrimeCocycleIsStillEquivariant}
  and
  \eqref{CheckingThatAveragedPUHPrimeCocycleIsStillACechCocycle},
  so that the proof given there generalizes straightforwardly.
\end{proof}

\bigskip

\noindent
{\bf The orbi-smooth Oka principle.}
We may now prove a constrained generalization of the
``smooth Oka principle'' (Thm. \ref{SmoothOkaPrinciple})
from smooth manifolds to smooth orbifolds, subject to
the conditions discussed above
(truncated structure group and blowup-stability over resolvable singularities).

\begin{lemma}[Concordance of stable equivariant bundles on blowups]
  \label{ConcordanceOfPrincipalBundlesOnBlowupsForStructureGroupWithTruncatedShape}
  $\,$

  \noindent
  Let $G \,\in\, \Groups(\FiniteSets)_{\resolvable}$ (Ntn. \ref{ResolvableOrbiSingularities})
  and
  $\Gamma \,\in\, \Groups(\SmoothInfinityGroupoids)_{\leq 0}$
  with truncated classifying shape (Ntn. \ref{CohesiveGroupsWithTruncatedClassifyingShape}),
  so that a notion of stable $G$-equivariant $\Gamma$-principal bundles
  exists (Ntn. \ref{StableEquivariantBundles})
  then
  $$
    \Truncation{0}
    \,
    \shape
    \,
    \Maps{}
      { \HomotopyQuotient{\SmoothManifold}{G} }
      { \mathbf{B}\Gamma }
    ^{\stable}
    \;\;
    \simeq
    \;\;
    \Truncation{0}
    \,
    \shape
    \,
    \Maps{\big}
      { \Quotient{(\SmoothManifold \times \SmoothSphere{\,n+2})}{G} }
      { \mathbf{B}\Gamma }
    \,.
  $$
\end{lemma}
\begin{proof}
  By stability \eqref{ConditionOfCompatibleRestrictionsOnIsotropyGroups}
  we have such a bijection already on isomorphism classes,
  both over $\SmoothManifold$ and over $\SmoothManifold \times \SmoothSimplex{1}$.
  This implies the claimed bijection on concordance classes by
  \eqref{SmoothInfinityGroupoidsWithCoinciding0TruncationHaveTheSameConcordanceClasses}
  in Ex. \ref{ConnectedComponentsOfShapeAreConcordanceClasses}.
\end{proof}

\begin{remark}[Cardinality bound on homotopy classes of maps to structure group]
  \label{CardinalityBoundOnHomotopyClassesOfMapsFromXModGToGamma}
  In addition to the
  substantial condition of blowup-stability
  for runcated structure group over resolvable orbi-singularities,
  the following proof of the orbi-smooth Oka principle (Thm. \ref{OrbiSmoothOkaPrinciple})
  requires a mild technical condition
  (in  the step \eqref{ManifoldsForDetectingFreeHomotopyEquivalences} below):
  The cardinality of the fundamental group
  $$
    \pi_1
    \,
    \Maps{\big}
      { \HomotopyQuotient{\shape  \, \SmoothManifold}{G} }
      { B \Gamma \rtimes G }
  $$
  must be countable, at every basepoint.
  At the canonical basepoint this means that
  $$
    H^0_G \big(\SmoothManifold;\, \shape \, \Gamma \rtimes G \big)
    \;=\;
    \pi_0
    \,
    \Maps{\big}
      { \HomotopyQuotient{\shape \, \SmoothManifold}{G} }
      { \shape \, \Gamma \times G }
  $$
  must be countable, and due to the assumption that $\shape \, \Gamma$
  is truncated, this essentially just means that the Borel-equivariant cohomology of
  $\SmoothManifold$ in every degree is a countable set.
\end{remark}

\medskip

\begin{theorem}[Orbi-smooth Oka principle for
stable truncated maps out of good orbifolds with resolvable singularities]
\label{OrbiSmoothOkaPrinciple}
Given
\begin{itemize}

  \vspace{-.1cm}
  \item[--] $G \,\in\, \Groups(\FiniteSets)_{\resolvable}$ (Ntn. \ref{ResolvableOrbiSingularities}),

  \vspace{-.1cm}
  \item[--] $G \acts \, \Gamma \,\in\, \Actions{G}\big( \Groups(\kTopologicalSpaces) \big)$

    such that

    \vspace{-.2cm}
    \begin{itemize}

    \vspace{-.1cm}
    \item[--]
    $\Gamma \rtimes G $ is of truncated classifying shape (Ntn. \ref{CohesiveGroupsWithTruncatedClassifyingShape}),

    \vspace{-.1cm}
    \item[--]
    there is a notion of blowup-stability of $G$-equivariant $\Gamma \rtimes G$-principal bundles
     (Ntn. \ref{StableEquivariantBundles}),

  \end{itemize}

\end{itemize}

\vspace{-2mm}
\noindent
{\bf (i)} then,
  for all

  -- $G \acts \, \SmoothManifold \,\in\, \Actions{G}(\SmoothManifolds)$
   such that
   $\pi_1
     \Maps{\big}
       { \HomotopyQuotient{\SmoothManifold}{G} }
       { B (\Gamma \rtimes G) }
   $
   is countable (Rem. \ref{CardinalityBoundOnHomotopyClassesOfMapsFromXModGToGamma});

  \noindent
  the canonical comparison morphism \eqref{ComparisonMorphismFromShapeOfMappingStackToMappingSpaceOfShapes}
  restricted to stable equivariant bundles
  is an equivalence in $\InfinityGroupoids$:
  \vspace{-2mm}
  \begin{equation}
    \label{SlicedSmoothOkaEquivalenceForGlobalQuotientOrbifolds}
    \begin{tikzcd}
      \shape
      \,
      \SliceMaps{}{\mathbf{B}G}
        { \HomotopyQuotient{\SmoothManifold}{G} }
        { \HomotopyQuotient{\mathbf{B}\Gamma}{G} }
      ^{\stable}
      \ar[
        r,
        "\sim"{below},
        "\scalebox{.7}{$ \widetilde { \shape \, \mathrm{ev} }  $} "{above}
      ]
      &
      \SliceMaps{}{\mathbf{B}G}
        { \HomotopyQuotient{ \shape \, \SmoothManifold }{ G } }
        { \HomotopyQuotient{ \shape \, \mathbf{B}\Gamma }{ G } }
      \,.
    \end{tikzcd}
  \end{equation}

  \vspace{-1mm}
  \noindent
  {\bf (ii)} When the $G$-action on $\Gamma$ is trivial, this reduces to
  \vspace{-3mm}
  \begin{equation}
    \label{UnslicedSmoothOkaEquivalenceForGlobalQuotientOrbifolds}
    \begin{tikzcd}
      \shape
      \,
      \Maps{}
        { \HomotopyQuotient{\SmoothManifold}{G} }
        { \mathbf{B}\Gamma }
      ^{\stable}
      \ar[
        r,
        "\sim"{below},
        "\scalebox{.7}{$ \widetilde { \shape \, \mathrm{ev} }  $} "{above}
      ]
      &
      \Maps{}
        { \HomotopyQuotient{ \shape \, \SmoothManifold }{ G } }
        { \shape \, \mathbf{B}\Gamma }
      \,.
    \end{tikzcd}
  \end{equation}
\end{theorem}
\begin{proof}
  First, we may use the above lemmas to obtain the bijections
  \eqref{ClassificationOfEquivariantPrincipalBundlesForTruncatedShapeByBorelCohomology}
  in the unsliced case, for any
  $\Gamma{\, '}$ of truncated classifying shape, notably for
  $\Gamma{\,'} \,\in\, \{ \Gamma ,\, \Gamma \rtimes G \}$:
    \vspace{-2mm}
  \begin{equation}
  \label{FirstStepInClassificationOfEquivariantPrincipalBundlesWithStructureOfTruncatedClassifyingShape}
  \def\arraystretch{1.5}
  \begin{array}{lll}
    \Truncation{0}
    \,
    \shape
    \,
    \Maps{}
      { \HomotopyQuotient{\SmoothManifold}{G} }
      { \mathbf{B}\Gamma{\, '}  }
    ^{\stable}
    &
    \;\simeq\;
    \Truncation{0}
    \,
    \shape
    \,
    \Maps{\big}
      { (\SmoothManifold \times \SmoothSphere{\,n+2})/G }
      { \mathbf{B}\Gamma{\, '} }
    &
    \proofstep{
      by Lem. \ref{ConcordanceOfPrincipalBundlesOnBlowupsForStructureGroupWithTruncatedShape}
    }
    \\
    &
    \;\simeq\;
    \Truncation{0}
    \,
    \Maps{\big}
      {
        \shape
        \,
        (\SmoothManifold \times \SmoothSphere{\,n+2})/G
      }
      {
        \shape
        \,
        \mathbf{B}\Gamma{\, '}
      }
    &
    \proofstep{ by Thm. \ref{SmoothOkaPrinciple} }
    \\
    &
    \;\simeq\;
    \Truncation{0}
    \,
    \Maps{}
      {
        \HomotopyQuotient{ \shape \, \SmoothManifold }{G}
      }
      {
        \shape
        \,
        \mathbf{B}\Gamma{\, '}
      }
    &
    \proofstep{ by Lem. \ref{ConcordanceInfinityGroupoidOfPrincipalBundlesOnBlowupForStructureGroupOfTruncatedShape} }
    .
  \end{array}
\end{equation}

\vspace{-1mm}
\noindent
This implies, for $\mathrm{K} \,\in\, \SmoothManifolds \xhookrightarrow{} \Actions{G}(\SmoothManifolds)$
a smooth manifold regarded as equipped with the trivial $G$-action, that:
\vspace{-3mm}
$$
  \def\arraystretch{1.5}
  \begin{array}{lll}
    \tau_0
    \,
    \Maps{\Big}
    {
      \shape
      \,
      \mathrm{K}
    }
    {
      \shape
      \,
      \Maps{}
      {
        \HomotopyQuotient{\SmoothManifold}{G}
      }
      {
        \mathbf{B}\Gamma{\, '}
      }
      ^{\stable}
    }
       & \;\simeq\;
    \tau_0
    \,
    \shape
    \,
    \Maps{\Big}
    {
      \mathrm{K}
    }
    {
      \Maps{}
      {
        \HomotopyQuotient{\SmoothManifold}{G}
      }
      {
        \mathbf{B}\Gamma{\, '}
      }
      ^{\stable}
    }
    &
    \proofstep{ by Thm. \ref{SmoothOkaPrinciple} }
    \\
    & \;\simeq\;
    \tau_0
    \,
    \shape
    \,
      \Maps{\big}
      {
        \mathrm{K}
        \times
        \HomotopyQuotient{\SmoothManifold}{G}
      }
      {
        \mathbf{B}\Gamma{\, '}
      }
      ^{\stable}
    &
    \proofstep{ by Lem. \ref{InternalHomAdjointness} }
    \\
    & \;\simeq\;
    \tau_0
    \,
      \Maps{\big}
      {
        \shape
        (
          \mathrm{K}
          \times
          \HomotopyQuotient{\SmoothManifold}{G}
        )
      }
      {
        \shape \, \mathbf{B}\Gamma{\, '}
      }
    &
    \proofstep{ by \eqref{FirstStepInClassificationOfEquivariantPrincipalBundlesWithStructureOfTruncatedClassifyingShape} }
    \\
    & \;\simeq\;
    \tau_0
    \,
      \Maps{\big}
      {
        (\shape \, \mathrm{K})
        \times
        ( \shape \, \HomotopyQuotient{\SmoothManifold}{G})
        )
      }
      {
        \shape \, \mathbf{B}\Gamma{\, '}
      }
    &
    \proofstep{ by \eqref{ShapePreservesBinaryProducts} }
    \\
    & \;\simeq\;
    \tau_0
    \,
    \Maps{\big}
    {
      \shape \, \mathrm{K}
    }
    {
      \Maps{}
      {
        \shape \, \HomotopyQuotient{\SmoothManifold}{G})
      }
      {
        \shape \, \mathbf{B}\Gamma{\, '}
      }
    }
    &
    \proofstep{ by Lem. \ref{InternalHomAdjointness} }
    ,
  \end{array}
$$

\vspace{-2mm}
\noindent hence that
\vspace{-2mm}
\begin{equation}
  \label{ProbingTheClassificationOfEquivariantBundlesWithManifolds}
  \underset{
    \mathrm{K} \in \SmoothManifolds
  }{\forall}
  \;\;\;\;\;\;
  \tau_0
  \,
  \Maps{\Big}
    { \shape \, \mathrm{K} \;\; }
    {
      \shape
      \,
      \Maps{}
        { \HomotopyQuotient{\SmoothManifold}{G} }
        { \mathbf{B}\Gamma{\, '} }
      ^{\stable}
      \xrightarrow{\;}
      \Maps{}
        { \HomotopyQuotient{ \shape\, \SmoothManifold }{G} }
        { \shape \, \mathbf{B}\Gamma{\, '} }
    }
  \;
  \text{ is a bijection. }
\end{equation}
Considering this for the cases
\vspace{-2mm}
\begin{equation}
  \label{ManifoldsForDetectingFreeHomotopyEquivalences}
  \mathrm{K}
  \;=\;
  \ast
  \,,
  \phantom{AAAAAA}
  \mathrm{K}
  \;=\;
  \SmoothSphere{k}
  ,\,
  \;
  k \in \mathbb{N}_+
  \,,
  \phantom{AAAAAA}
  \mathrm{K}
  \;=\;
  \mathbb{R}^2
    \setminus
    \Big(
      \underset{
        \scalebox{.6}{$
          \pi_1 \Maps{\big} {\HomotopyQuotient{\SmoothManifold}{G}}{ B \Gamma' }
        $}
      }{\sqcup} \ast
    \Big),
\end{equation}

\vspace{-2mm}
\noindent
(here the last item uses the assumption that
$\pi_1 \Maps{\big} {\HomotopyQuotient{\SmoothManifold}{G}}{ B (\Gamma \rtimes G) }$
is a countable set,
Rem. \ref{CardinalityBoundOnHomotopyClassesOfMapsFromXModGToGamma},
in order to embed that many points into $\RealNumbers^2$)
and observing that, by Prop. \ref{CohesiveShapeOfSmoothManifoldsIsTheirHomotopyType},
\vspace{-2mm}
$$
  \shape \, \ast
  \;\simeq\;
  \ast
  \,,
  \phantom{AAAAAA}
  \shape \, \SmoothSphere{k}
  \;\simeq\;
  \ShapeOfSphere{k}
  ,\,
  \phantom{AAAAAA}
  \shape
  \left(
    \mathbb{R}^2
      \setminus
    \Big(
      \underset{
        \scalebox{.6}{$
          H^0(\HomotopyQuotient{\SmoothManifold}{G};\, \shape\, \Gamma')
        $}
      }{\sqcup} \ast
    \Big)
  \!\!\right)
  \;\simeq\;
  \underset{     \scalebox{.6}{$
    H^1(\HomotopyQuotient{\SmoothManifold}{G};\, \shape\, \Gamma')
    $}
  }{\vee}
  \ShapeOfSphere{1}
  \,,
$$

\vspace{-.5mm}
\noindent implies the second claim \eqref{UnslicedSmoothOkaEquivalenceForGlobalQuotientOrbifolds},
by Lem. \ref{DetectingWeakHomotopyEquivalencesViaFreeHomotopySets}
applied to \eqref{ProbingTheClassificationOfEquivariantBundlesWithManifolds}.

From this follows
the first claim \eqref{SlicedSmoothOkaEquivalenceForGlobalQuotientOrbifolds} by
the fact that the shape modality preserves fiber products over discrete objects
(Prop. \ref{ShapeFunctorPreservesHomotopyFibersOverDiscreteObjects}):
\vspace{-1mm}
$$
  \def\arraystretch{1.75}
  \begin{array}{lll}
    \shape
    \,
    \SliceMaps{\big}{BG}
      { \HomotopyQuotient{\SmoothManifold}{G} }
      { \mathbf{B}(\Gamma \rtimes G) }
    ^{\stable}
    &
    \simeq\;
    \shape
    \Big(
      \Maps{}
        { \HomotopyQuotient{\SmoothManifold}{G} }
        { \mathbf{B}(\Gamma \rtimes G) }
      ^{\stable}
      \underset{
        \Maps{}
          { \HomotopyQuotient{\SmoothManifold}{G} }
          { \mathbf{B}G }
      }{\times}
      \left\{
        (X \to \ast)\sslash G
      \right\}
    \Big)
    &
    \proofstep{ by Def. \ref{SliceMappingStack} }
    \\
    &
    \simeq\;
    \shape
    \Big(
      \Maps{}
        { \HomotopyQuotient{\SmoothManifold}{G} }
        { \mathbf{B}(\Gamma \rtimes G) }
      ^{\stable}
      \underset{
        \Maps{}
          { \scalebox{.7}{$ \HomotopyQuotient{\shape \, \SmoothManifold}{G} $} }
          { B G }
      }{\times}
      \left\{
        (\shape \, X \to \ast)\sslash G
      \right\}
    \Big)
    &
    \proofstep{ by Lem. \ref{MappingStackIntoDiscreteObjectIsDiscrete} }
    \\
    &
    \simeq\;
    \shape
    \,
    \Maps{}
      { \HomotopyQuotient{\SmoothManifold}{G} }
      { \mathbf{B}(\Gamma \rtimes G) }
    \underset{
      \Maps{}
        { \scalebox{.7}{$ \HomotopyQuotient{\shape \SmoothManifold}{G}  $} }
        { B G }
    }{\times}
    \left\{
      (\shape \, X \to \ast)\sslash G
    \right\}
    &
    \proofstep{ by Prop. \ref{ShapeFunctorPreservesHomotopyFibersOverDiscreteObjects} }
    \\
    &
    \simeq\;
    \Maps{}
      { \HomotopyQuotient{\shape \, \SmoothManifold}{G} }
      { \HomotopyQuotient{\shape \, \mathbf{B}\Gamma}{G} }
    \underset{
      \Maps{}
        { \scalebox{.7}{$ \HomotopyQuotient{\shape \, \SmoothManifold}{G} $} }
        { B G }
    }{\times}
    \left\{
      (\shape \, X \to \ast)\sslash G
    \right\}
    &
    \proofstep{ by \eqref{UnslicedSmoothOkaEquivalenceForGlobalQuotientOrbifolds} }
    \\
    &
    \simeq\;
    \SliceMaps{}{BG}
      { \HomotopyQuotient{\shape \, \SmoothManifold}{G} }
      { \HomotopyQuotient{\shape \, \mathbf{B}\Gamma}{G} }
    &
    \proofstep{ by Def. \ref{SliceMappingStack}. }
  \end{array}
$$

\vspace{-6mm}
\end{proof}

From this, we immediately obtain, in equivariant generalization of
Thm. \ref{ClassificationOfPrincipalInfinityBundles}:
\begin{theorem}[Concordance classification of smooth equivariant principal bundles]
\label{ClassificationOfGammaEquivariantPrincipalBundlesForGammaWithTruncatedClassifyingShape}
Given

  \vspace{1mm}
  --
  $G \,\in\, \Groups(\FiniteSets)_{\resolvable}$ (Ntn. \ref{ResolvableOrbiSingularities}),

  \vspace{1mm}
  --
  $\Gamma \rtimes G \,\in\, \Groups(\SmoothSets)$ of truncated classifying shape (Ntn. \ref{CohesiveGroupsWithTruncatedClassifyingShape}),

  \vspace{1mm}
  --
  a notion of blowup-stability of $G$-equivariant $\Gamma \rtimes G$-principal bundles
  (Ntn. \ref{StableEquivariantBundles});

  \vspace{1mm}
  --
  $G \acts \SmoothManifold \,\in\, \Actions{G}(\SmoothManifolds)$
  such that
  $\pi_1 \Maps{\big} {\HomotopyQuotient{\SmoothManifold}{G}}{ B (\Gamma \rtimes G) }$
  is countable (Rem. \ref{CardinalityBoundOnHomotopyClassesOfMapsFromXModGToGamma});

  \vspace{1mm}
  \noindent
  there is a natural identification of
  concordance classes
  (Def. \ref{ConcordanceOfEquivariantPrincipalBundles})
  of
  stable
  $G$-equivariant $\Gamma$-principal
  bundles
  \eqref{EquivariantPrincipalBundlesAsPrincipalBundlesOnQuotientStack}
  with
  first non-abelian Borel-equivariant cohomology:
  \vspace{-2mm}
  \begin{equation}
    \label{ClassificationOfEquivariantPrincipalBundlesForTruncatedShapeByBorelCohomology}
    \begin{tikzcd}[row sep=12pt, column sep=large]
      \ConcordanceClasses{
        \EquivariantPrincipalBundles
          {G}{\Gamma}
        (\SmoothInfinityGroupoids)_{\SmoothManifold}^{\stable}
      }
      \ar[d, -, shift left=1pt]
      \ar[d, -, shift right=1pt]
      \ar[r ,  "\sim"{below}]
      &
      H^1_G
      (
        X
        ;\,
        \shape\, \Gamma
      )
      \ar[d, -, shift left=1pt]
      \ar[d, -, shift right=1pt]
      \\
      \Truncation{0}
      \,
      \shape
      \,
      \SliceMaps{\big}{BG}
        { \HomotopyQuotient{\SmoothManifold}{G} }
        { \HomotopyQuotient{\mathbf{B}\Gamma}{G} }
      ^{\stable}
      \ar[
        r,
        "\sim"{below},
        "\scalebox{.7}{$ \widetilde { \shape \, \mathrm{ev} }  $} "{above}
      ]
      &
      \Truncation{0}
      \,
      \SliceMaps{\big}{BG}
        { \HomotopyQuotient{ \shape \, \SmoothManifold }{ G } }
        { \HomotopyQuotient{\shape \, \mathbf{B}\Gamma }{ G } }
      \,.
    \end{tikzcd}
  \end{equation}
\end{theorem}
\begin{proof}
  This is the
  restriction to connected components
  of
  \eqref{SlicedSmoothOkaEquivalenceForGlobalQuotientOrbifolds}
  in Thm. \ref{OrbiSmoothOkaPrinciple}.
\end{proof}

\vspace{-1mm}
Below in Thm. \ref{BorelClassificationOfEquivariantBundlesForResolvableSingularitiesAndEquivariantStructure}
we refine this classification statement from concordance classes
to isomorphism classes of stable equivariant bundles.

\medskip

\begin{example}[Orbi-smooth Oka principle for $\PUH$-coefficients over ADE-singularities]
  \label{OrbiSmoothOkaPrincipleForPUHCoefficientsOverThePoint}
  Let $G \,\subset \SpOne\,$ be a finite subgroup (Ex. \ref{ADEGroupsHaveSphericalSpaceForms})
  and consider the structure group $\PUH$ (Ex. \ref{ProjectiveUnitarGroupOnAHilbertSpace})
  with the corresponding stability condition (Thm. \ref{BlowupStabilityOfADEEquivariantPUHPrimePrincipalBundles}).
  We describe the shape of the
  topological space of conjugations of stable morphisms
  $\mathbf{B}G \xrightarrow{\;} \mathbf{B} \PUH$,
  hence of isomorphisms of stable $G$-equivariant $\PUH$-principal bundles over the point,
  according to the orbi-smooth Oka principle
  (Thm. \ref{OrbiSmoothOkaPrinciple}).
  This may be summarized by the following diagram:
  \begin{equation}
    \label{OrbiSmoothOkaIdentificationForEquivariantPUHBundlesOverThePoint}
    \begin{tikzcd}
      \Maps{}
        { B G }
        { B^3 \Integers }
      \ar[
        d,
        "{
          \mbox{
            \tiny
            \rm
            \eqref{SpaceOfMapsFromBGToBThreeIntegersForFiniteSubgroupsOfSpOne}
          }
        }"{swap},
        "{\sim}"
      ]
      \ar[
        rr,
        "{\sim}"{swap},
        "{
          \mbox{
            \tiny
            \rm
            Ex. \ref{HomotopyFiberSequenceOfTheProjectiveUnitaryGroup}
          }
        }"
      ]
      &[-20pt]&[-20pt]
      \Maps{}
        { \shape \, \mathbf{B}G }
        { \shape \, \mathbf{B} \PUH }
      \ar[
        rr,
        "{
          \mbox{
            \tiny
            \rm
            Thm. \ref{OrbiSmoothOkaPrinciple}
          }
        }",
        "{ \sim }"{swap}
      ]
      &&
      \shape \,
      \Maps{}
        { \mathbf{B}G }
        { \mathbf{B} \PUH }
      ^{\stable}
      \\
      B
      \Homs{}
        { G }
        { \CircleGroup }
      \,\times\,
      B^3 \Integers
      &\simeq&
      B
      \Big(
          \Homs{}
            { G }
            { \CircleGroup }
          \,\times\,
          \big(\;\;\;\;
            \underset{
              \mathclap{
                \mathrm{Irr}(G)
              }
            }{\prod}
            \,
            \UH
          \big)/ \CircleGroup
      \Big)
      \ar[
        urr,
        rounded corners,
        to path={
          -- ([xshift=4pt]\tikztostart.east)
          -- node[below]{
               \scalebox{.7}{$
                 B \mathrm{Ad} ( (-)  \otimes (-)   )
               $}
               \mbox{
                 \tiny
                 \eqref{ProjectiveEndotwinersOfStableGRepresentations}
               }
              }
             ([yshift=-33.5pt]\tikztotarget.south)
          -- ([yshift=-00pt]\tikztotarget.south)
        },
        "{
          B
          \mathrm{Ad}
          \,
          \mbox{
            \tiny
            \rm
            \eqref{ProjectiveEndotwinersOfStableGRepresentations}
          }
        }"{swap}
      ]
    \end{tikzcd}
  \end{equation}

  \vspace{-.3cm}
  \noindent
  (Here $\mathrm{Irr}(G)$ denotes the set of isomorphism classes of irredible
  unitary $G$-representations.)

  \medskip

  The interesting point to highlight is what concretely the morphisms in
  $\Maps{}{ \mathbf{B}G }{ \mathbf{B}\PUH }^{\stable}$
  (on the top right)
  are like and how exactly these
  have the shape predicted by the orbi-smooth Oka principle (on the top left).
  After some re-casting of the mapping shape on the bottom left,
  this is brought out by the assignment on the bottom right, which
  we describe in \eqref{ProjectiveEndotwinersOfStableGRepresentations} below.

  \medskip

  First to see the equivalence on the left of \eqref{OrbiSmoothOkaIdentificationForEquivariantPUHBundlesOverThePoint},
  we
  use that the integral group cohomology of $G$ is concentrated in
  even degrees \eqref{IntegralGroupCohomologyOfADESubgroups} and
  that the second integral cohomology gives the group characters
  \eqref{FirstUOneGroupCohomologyOfADESubgroup}:
  The first fact gives
  $\Truncation{0} \, \Maps{}{B G}{B^3 \Integers} \,\simeq\, \ast$
  which implies (with Prop. \ref{LoopingAndDeloopingEquivalence}) that
  $\Maps{}{ B G }{B^3 \Integers}
    \,\simeq\,
  B \Omega \Maps{}{B G}{B^3 \Integers}
  \,\simeq\,
  B \Maps{}{B G}{\Omega B^3 \Integers}
  \,\simeq\,
  B \Maps{}{ B G }{ B^2 \Integers }
  ;\,
  $
  and then the second fact gives
  $\Truncation{0} \, \Maps{}{ B G }{ B^2 \Integers }
    \,\simeq\,
    \Homs{}{G}{\CircleGroup}.
  $
  Proceeding in this manner shows that each of these connected components
  has homotopy groups concentrated in degree 2 on
  $ \Maps{}{B G}{ \Omega^3 B^3 \Integers }
   \,\simeq\,
   \Maps{}{B G}{\Integers}
   \,\simeq\,
   \Integers$,
   hence that
  $\Maps{}{B G}{B^2 \Integers} \,\simeq\, \Homs{}{G}{\CircleGroup} \times B^2 \Integers$.
  In conclusion:
  \begin{equation}
    \label{SpaceOfMapsFromBGToBThreeIntegersForFiniteSubgroupsOfSpOne}
    \Maps{}{ B G }{B^3 \Integers}
    \;\simeq\;
    B
    \big(
      H^2_{\mathrm{Grp}}(G;\, \Integers)
      \times
      B^2
      H^0_{\mathrm{Grp}}(G;\, \Integers)
    \big)
    \;\simeq\;
    B
    \big(
      \Homs{}{G}{\CircleGroup}
      \,\times\,
      B^2 \Integers
    \big)
    \;\simeq\;
    B \Homs{}{G}{ \CircleGroup }
    \,\times\,
    B^3 \Integers
    \,.
  \end{equation}

  The fact that the space \eqref{SpaceOfMapsFromBGToBThreeIntegersForFiniteSubgroupsOfSpOne}
  is connected means that there is an essentially unique stable morphism
  $\mathbf{B}G \xrightarrow{\;} \mathbf{B} \PUH$. Indeed,
  by Lem. \ref{StableProjectiveIsotorpyRepresentations} these
  stable projective representations must be the direct sums of
  infinite copies of all $[\tau]$-projective $G$-representations for
  some
  $[\tau] \,\in\, H^2_{\mathrm{Grp}}(G;\, \CircleGroup)
   \,\simeq\,
   H^3_{\mathrm{Grp}}(G;\, \Integers)$;
   but since the third integral group cohomology in the present
   case is trivial
   \eqref{IntegralGroupCohomologyOfADESubgroups},
   there is only the case $[\tau] = [0]$, and hence the only stable
   projective representation is, up to isomorphism, the projectivization of
   \begin{equation}
     \label{TheUniqueStableProjectiveGRepresentationForGFiniteSubgroupOfSpOne}
     \begin{tikzcd}[row sep=-2pt, column sep=4pt]
       G
       \ar[rr]
       &&
       \UnitaryGroup
       \Big(
       \underset{
         \mu \in \mathrm{Irr}(G)
       }{\bigoplus}
       \,
       \mu
       \;
       \otimes
       \,
       \HilbertSpace
       \Big)
       \\
  \scalebox{0.7}{$     g $}
       &\longmapsto&
     \scalebox{0.7}{$      \underset{\mu}{\oplus}
       \,
       \mu(g) \,\otimes\, \mathrm{id}
       $}
     \end{tikzcd}
   \end{equation}

   \vspace{-2mm}
\noindent   for any countably infinite-dimensional Hilbert space
   $\HilbertSpace$.

   \medskip

   Therefore the question is:
   {\it How can the group characters
   $\kappa \,\in\, \Homs{}{G}{\CircleGroup}$ encode
   projective endo-intertwiners of the stable representation \eqref{TheUniqueStableProjectiveGRepresentationForGFiniteSubgroupOfSpOne}?}
   Direct inspection shows that on $\underset{\mu}{\oplus}\, \mu$
   this is accomplished by the linear map which direct summand-wise
   is the linear identification
   \begin{equation}
     \label{TowardsTheProjectiveIntertwinersOfAStableProjectiveRep}
     \begin{tikzcd}[row sep=-2pt]
       \mu \ar[rr] && \mu \otimes \big( \kappa \times_{\CircleGroup} \ComplexNumbers\big)
       \\
       \scalebox{0.7}{$    v $}
       &\longmapsto&
         \scalebox{0.7}{$   v \otimes 1 $}
     \end{tikzcd}
   \end{equation}
   of the irrep $\mu$ with its
   tensor product with $\kappa$, the latter regarded as a
   complex linear representation on $\ComplexNumbers$ with canonical basis
   vector $1 \in \ComplexNumbers$.
   This operation fails, in general, to be a genuine unitary
   homomorphism, but the failure is exactly the phases
   $\kappa(g) \,\in\, \CircleGroup$, so that this
   does yield a projective intertwiner, as shown by the following
   commuting diagram (using the 2-groupoidal notation for the projective
   group from Ex. \ref{ProjectiveRepresentationsAndTheirCentralExtensions}):
  \begin{equation}
    \label{ProjectiveEndotwinersOfStableGRepresentations}
    \begin{tikzcd}[row sep=20pt, column sep=60pt]
      \mathbf{B}G
      \ar[
        rr,
        bend left=20,
        "{\ }"{swap,name=s1}
      ]
      \ar[
        rr,
        bend right=20,
        "{\ }"{name=t1}
      ]
      &&
      \HomotopyQuotient
        { \mathbf{B} \UH }
        { \mathbf{B} \CircleGroup }
      \mathrlap{
        \;\simeq\,
        \mathbf{B} \PUH
      }
      \ar[
        from=s1,
        to=t1,
        Rightarrow,
        "{
          \big(
            \kappa
            ,\,
            [\UnitaryOperator_\mu]_\mu
          \big)
        }"
      ]
      \\
      \\
      \bullet
      \ar[
        dd,
        "{g}"
      ]
      &
      \underset{
        \mathclap{
          \raisebox{-1pt}{
          \scalebox{.65}{
            $\mu \!\in\! \mathrm{Irr}(G)$
          }
          }
        }
      }{\bigoplus}
      \;
      \mu \,\otimes\, \HilbertSpace
      \ar[
        dd,
        start anchor={[yshift=7pt]},
        "{
          \underset{\mu}{\oplus} \, \mu(g) \,\otimes\, \mathrm{id}
        }"{description}
      ]
      \ar[
        rr,
        "{
          \underset{\mu}{\bigoplus}\,
          \Big(
            \big(
              \mu
              \xrightarrow{v \,\mapsto\, v \otimes 1}
              \mu \otimes \kappa
            \big)
            \,\otimes\,
            \UnitaryOperator_\mu
          \Big)
        }",
        "{\ }"{swap, pos=.9, name=s}
      ]
      &&
      \underset{
        \mathclap{
          \raisebox{-1pt}{
          \scalebox{.65}{
            $\mu \!\in\! \mathrm{Irr}(G)$
          }
          }
        }
      }{\bigoplus}
      \;
      \mu \,\otimes\, \HilbertSpace
      \ar[
        dd,
        start anchor={[yshift=7pt]},
        "{
          \underset{\mu}{\oplus} \, \mu(g) \,\otimes\, \mathrm{id}
        }"{description}
      ]
      \\
      {}
      \ar[
        r,
        phantom,
        "{\longmapsto}"{pos=.35}
      ]
      &
      {}
      \\
      \bullet
      &
      \underset{
        \mathclap{
          \raisebox{-1pt}{
          \scalebox{.65}{
            $\mu \!\in\! \mathrm{Irr}(G)$
          }
          }
        }
      }{\bigoplus}
      \;
      \mu \,\otimes\, \HilbertSpace
      \ar[
        rr,
        "{
          \underset{\mu}{\bigoplus}\,
          \Big(
            \big(
              \mu
              \xrightarrow{v \,\mapsto\, v \otimes 1}
              \mu \otimes \kappa
            \big)
            \,\otimes\,
            \UnitaryOperator_\mu
          \Big)
        }"{swap},
        "{\ }"{pos=.1, name=t}
      ]
      &&
      \underset{
        \mathclap{
          \raisebox{-1pt}{
          \scalebox{.65}{
            $\mu \!\in\! \mathrm{Irr}(G)$
          }
          }
        }
      }{\bigoplus}
      \;
      \mu \,\otimes\, \HilbertSpace
      \ar[
        from=t,
        to=s,
        Rightarrow,
        "{ \kappa(g) }"{description}
      ]
    \end{tikzcd}
  \end{equation}
  Here the $\UnitaryOperator_\mu \,\in\, \UH$ are any unitary operators
  acting on the Hilbert space of infinite copies of
  the irrep $\mu$.
\end{example}

\begin{example}[Orbi-smooth Oka principle for $\PUH$-coefficients over A-singularities]
It is instructive to further spell out the construction of the projective
endo-intertwiners from Ex. \ref{OrbiSmoothOkaPrincipleForPUHCoefficientsOverThePoint}
for the case that $G \,=\, \CyclicGroup{n}$ is a cyclic group.
Here $\mathrm{Irr}(G) \,=\, \mathrm{Irr}(\CyclicGroup{n}) \,\simeq\, \CyclicGroup{n}$
and
the standard $n \times n$ matrix presentation of the
regular representation
$\underset{\mu}{\oplus}\, \mu$
is
$$
  1 \leq q, q' \leq n
  \;\;\;\;\;\;\;\;
  \vdash
  \;\;\;\;\;\;\;\;
  \big(
    U_{\mathrm{reg}}([k])
  \big)_{q, q'}
  \;\;
  =
  \;\;
  \exp
  \big(
    2 \pi \ImaginaryUnit
    \cdot
    k \tfrac{q-1}{n}
  \big)
  \cdot
  \delta_{q, q'}
  \,.
$$
Moreover, the tensoring with
$\kappa \,\in\, \Homs{}{\CyclicGroup{n}}{\CircleGroup} \,\simeq\, \CyclicGroup{n}$
acts on these irreps by cyclically permuting them, so that
the matrix representation for the projective intertwiner
\eqref{TowardsTheProjectiveIntertwinersOfAStableProjectiveRep} is
$$
  1 \leq q, q' \leq n
  \;\;\;\;\;\;\;\;
  \vdash
  \;\;\;\;\;\;\;\;
  \big(
    U_{\kappa}
  \big)_{q, q'}
  \;\;
  =
  \;\;
  \delta_{q, q' - \kappa}
  \,.
$$
That this gives the required diagram \eqref{ProjectiveEndotwinersOfStableGRepresentations}
\begin{equation}
  \label{ProjectiveSelfIntertwinerOFZeroProjectiveRegularRepresentation}
  \begin{tikzcd}[row sep=30pt, column sep=40pt]
    \bullet
    \ar[d, "{ [k] }"]
    &[+20pt]
    \bullet
    \ar[
      d,
      "{ U_{\mathrm{reg}}\big([k]\big) }"{swap}
    ]
    \ar[
      rr,
      "{
        U_{([\kappa], [\vec c\,])}
      }",
      "{\ }"{swap, pos=.9, name=s}
    ]
    &&
    \bullet
    \ar[
      d,
      "{ U_{\mathrm{reg}}\big([k]\big) }"
    ]
    \\
    \bullet
    &
    \bullet
    \ar[
      rr,
      "{
        U_{([\kappa], [\vec c\,])}
      }"{swap},
      "{\ }"{pos=.1, name=t}
    ]
    &&
    \bullet
    \ar[
      from=t,
      to=s,
      Rightarrow,
      "{
        \theta_{[\kappa]}([k])
      }"{description}
    ]
  \end{tikzcd}
\end{equation}
follows from a straightforward computation:
$$
  \begin{aligned}
    \big(
      U_{\mathrm{reg}}([k])
      \cdot
      U_{[\kappa], [\vec c\,]}
    \big)_{q, q'}
    &
    \;=\;
    \exp
    \big(
      2 \pi \ImaginaryUnit
      \cdot
      k
      \tfrac{q-1}{n}
    \big)
      \cdot
    \delta_{q, q' - \kappa}
    \\
    \big(
      U_{[\kappa], [\vec c\,]}
      \cdot
      U_{\mathrm{reg}}([1])
    \big)_{q, q'}
    &
    \;=\;
    \exp
    \big(
      2 \pi \ImaginaryUnit
      \cdot
      k
      \tfrac{\kappa + q - 1}{n}
    \big)
      \cdot
    \delta_{q, q' - \kappa}
  \end{aligned}
$$
whence
$$
  \theta_{[\kappa]}([k])
  \;=\;
  \exp
  \big(
    2 \pi \ImaginaryUnit
    \cdot
    k
    \tfrac{\kappa}{n}
  \big)
  \,.
$$
  Curiously,
  the projective self-intertwining
  relation
  \eqref{ProjectiveSelfIntertwinerOFZeroProjectiveRegularRepresentation}
  of the $[0]$-projective regular $\CyclicGroup{n}$-representation
  happens to coincide, in the case $\kappa = 1$,
  with the defining equation of the
  finite matrix approximation of the ``non-commutative torus''
  as in \cite[(2.29)]{LandiLizziSzabo01}.
\end{example}

\section{Local local triviality is implied}
\label{EquivariantLocalTrivializationIsImplies}

We have seen in \cref{PrincipalBundlesInternalToTopologicalGActions}
that the definition
of principal bundles, in Cartan's original sense
and as formalized by Grothendieck's internal notion of pseudo-torsors,
namely without requiring local triviality
(Rem. \ref{AssumptionOfLocalTrivializability}),
has excellent abstract mathematical properties
that make their correct generalization to equivariant principal bundles
(or any other generalization) a purely formal matter.
While a form of local triviality of bundles may also be
formulated internally as soon as the ambient category is regular,
we have seen in \cref{NotionsOfEquivariantLocalTrivialization} that,
internal to 1-categories of group actions, the resulting notion
is overly restrictive, and that some ingenuity and labor
is required to establish a notion of equivariant local triviality
that seems more appropriate.

\medskip
However, in \cref{AsInfinityBundlesInternalToSliceOverBG},
we observed that the analogous internalization
procedure that yields equivariant bundles internal to
the category of topological spaces
applies to yield equivariant $\infty$-bundles internal to
any $\infty$-topos.
Here we prove (Thm. \ref{ProperClassificationOfEquivariantBundlesForResolvableSingularitiesAndEquivariantStructure} below)
that the latter,
when seen externally, are {\it automatically} locally trivial
in the proper equivariant sense of \cref{NotionsOfEquivariantLocalTrivialization}.
We may think of this as saying that Cartan's original definition
of principal bundles is the {\it correct} one, when regarded,
with Grothedieck, as a
formal theory waiting to be internalized (interpreted)
in an ambient category with finite limits;
and that the constraint of local
triviality is not part of the formal theory of principal bundles,
but is to be provided by the semantic context in which the theory is interpreted.

\medskip
Thm. \ref{ProperClassificationOfEquivariantBundlesForResolvableSingularitiesAndEquivariantStructure}
below shows that contexts which do provide this feature in a useful way are the
Grothendieck(-Simpson-To{\"e}n-Vezzosi-Lurie-Rezk) $\infty$-toposes.
In fact, their underlying $(2,1)$-toposes are sufficient for capturing
equivariant local triviality with respect to
equivariance 1-groups, as $(n+1,1)$-toposes will be sufficient for
capturing equivariance $n$-groups, but traditional 1-toposes
are sufficient only for equivariance 0-groups, which includes only the trivial group.
Specifically,
the model $\Topos \,=\,\SmoothInfinityGroupoids$ is such that
interpreting equivariant principal bundles inside it, in Cartan-Grothendieck's sense,
and then restricting the resulting notion along the full inclusion
$\DTopologicalSpaces \xhookrightarrow{\;} \SmoothInfinityGroupoids$,
recovers the traditional notion of equivariant principal bundles,
{\it including} now the proper equivariant local triviality constraint.

\medskip
In this sense, the traditional theory of equivariant bundles
is completed by embedding it into the more general theory of
equivariant $\infty$-bundles in singular-cohesive $\infty$-toposes.

 \vspace{1mm}
\noindent
\begin{notation}[Equivariance groups]
Throughout this section,

\vspace{0mm}
\begin{enumerate}[{\bf (i)}]

\vspace{-.2cm}
\item
$
 G \,\in\, \Groups(\Sets) \xhookrightarrow{\Groups(\Discrete)}
    \Groups \left(\SimplicialPresheaves(\CartesianSpaces)\right)
$
\\
denotes a discrete equivariance group regarded as a constant simplicial preaheaf
of groups,

\vspace{-.2cm}
\item
$
  G \acts \, \Gamma
  \,\in\,
  \Groups\left( \Actions{G}(\DHausdorffSpaces) \right)
  \xhookrightarrow{\quad}
  \Groups\left(\!
    \Actions{G}
    \left(
      \SimplicialPresheaves(\CartesianSpaces)
    \right)
 \! \right)
$
\\
denotes a $G$-equivariant
D-topological Hausdorff group $\Gamma$,
regarded as a group in $G$-actions on simplicially constant simplicial
presheaves.

\end{enumerate}
 \end{notation}

\medskip

\noindent
{\bf The universal equivariant principal bundle over the moduli stack.}
First, we generalize the construction in
Rem. \ref{UniversalPrincipalBundleOverTheModuliStack}
of universal principal bundles over the moduli stack
to the equivariant context.

\begin{definition}[Crossed adjoint action through a crossed homomorphism]
  \label{CrossedConjugationActionThroughACrossedHomomorphism}
  Given a crossed homomorphism
  $\phi \,\colon\, G \xhookrightarrow{\;} \Gamma$,
  hence a homomorphic lift
  $\widehat{(-)} \,;\, G \xrightarrow{\;} \Gamma \rtimes G$,
  we denote by
  $G \acts \, \Gamma_{\widehat{\mathrm{adj}}}
    \,\in\, \Actions{G}(\kTopologicalSpaces)$
the adjoint action which on the right is
  twisted by $\widehat{(-)}$:
   \vspace{-2mm}
  $$
    \begin{tikzcd}[row sep=-3pt]
      G \times \Gamma
      \ar[rr]
      &&
    \Gamma
      \\
    \scalebox{0.7}{$    (g,\, \gamma) $}
      &\longmapsto&
    \scalebox{0.7}{$    \alpha(g)(\gamma \cdot \phi(\gamma)^{-1})
      \mathrlap{
        \;=\;
        (\NeutralElement,\, g)
        \cdot
        (\gamma,\, \NeutralElement)
        \cdot
        \widehat{g}^{\;-1}
      }
      $}
    \end{tikzcd}
  $$
\end{definition}

\begin{lemma}[The universal $G$-equivariant $\Gamma$-principal bundle]
\label{TheUniversalGEquivariantGammaPrincipalBundle}
The action groupoid of the crossed adjoint action
$\Gamma_{\widehat{\mathrm{adj}}}$
(Def. \ref{CrossedConjugationActionThroughACrossedHomomorphism})
provides a fibrant replacement
$\SimplicialPresheaves(\CartesianSpaces)_{\proj}$
of the delooping of $G \xrightarrow{\, g\, \mapsto (\NeutralElement, g)\, } \Gamma \rtimes G$:
\vspace{-2mm}
$$
  \begin{tikzcd}[column sep=90pt, row sep=12pt]
    N\DeloopingGroupoid{ G }
    \ar[
      rr,
      "{
      \scalebox{0.8}{$  \SimplicialNerve
        \left(
          g
            \,\mapsto\,
          \left(
            \NeutralElement,
            (\NeutralElement, g)
          \right)
        \right)
        $}
      }",
      "\in \WeakEquivalences"{swap}
    ]
    \ar[
      dr,
      "{
\scalebox{0.8}{$        \SimplicialNerve
        \left(
          g
          \,\mapsto\,
          (\NeutralElement,g )
        \right)
        $}
      }"{sloped}
    ]
    &&
    N
    \big(
      \Gamma_{\widehat{\mathrm{adj}}} \times (\Gamma \rtimes G)
      \rightrightarrows
      \Gamma
    \big).
    \ar[
      dl,
      "{
        \in
        \ProjectiveFibrations
      }",
      "{
     \scalebox{0.8}{$   \SimplicialNerve
        \left(
          (\gamma{\; '}, g)
          \,\mapsfrom\,
          \left(
            \gamma,
            (\gamma{\; '}, g)
          \right)
        \right)
        $}
      }"{sloped, pos=.6}
    ]
    \\
    &
    N\DeloopingGroupoid{( \Gamma \rtimes G )}
  \end{tikzcd}
$$
\end{lemma}
\begin{proof}
  It is clear that the diagram commutes and that the right morphism
  is a projective fibration.
  To see that the top morphism is a weak equivalence it is sufficient to
  observe that we have in fact a homotopy equivalence between the
  representing groupoids:
  \vspace{-2mm}
$$
  \begin{tikzcd}[row sep=-3pt]
  \scalebox{0.7}{$    g $}
    &\longmapsfrom&
  \scalebox{0.7}{$    \left(
      \gamma,
      (\gamma{\; '}, \, g)
    \right)
    $}
    \\
    \left(
      \ast \times \, G
      \rightrightarrows
      \ast
    \right)
    \ar[
      from=rr,
      shift right=3pt
    ]
    \ar[
      rr,
      shift right=3pt
    ]
    &&
    \left(
      \Gamma
      \times
      (\Gamma \rtimes G)
      \rightrightarrows
      \Gamma
    \right)
    \\
    \scalebox{0.7}{$  g  $}
    &\longmapsto&
      \scalebox{0.7}{$ \left(\NeutralElement, (\NeutralElement, g) \right) $}
    \mathrlap{\,.}
  \end{tikzcd}
$$

\vspace{-2mm}
\noindent In one direction, the composite of these two functors equals the identity,
in the other we have the following natural transformation to the identity:
\vspace{-2mm}
$$
  \begin{tikzcd}[column sep=40pt]
    \gamma
    \ar[
      d,
      "{
        (\gamma{\; '}, \, g)
      }"{left}
    ]
    \ar[
      rr,
      "{
        (\gamma, \NeutralElement)^{-1}
      }"
    ]
    &&
    \NeutralElement
    \ar[
      d,
      "{
        (\NeutralElement, g)
      }"
    ]
    \\
    \alpha(g^{-1})(\gamma \cdot \gamma{\; '})
    \ar[
      rr,
      "{
    \scalebox{0.7}{$    \left(
          \alpha(g^{-1})(\gamma \cdot \gamma{\; '})
          ,\,
          \NeutralElement
        \right)^{-1}
        $}
      }"{swap}
    ]
    &&
    \NeutralElement
  \end{tikzcd}
$$

\vspace{-7mm}
\end{proof}

\begin{proposition}[Equivariant principal bundles from {\v C}ech cocycles]
  \label{EquivariantPrincipalBundlesFromCechCocycles}
  For

  -- $G \,\in\, \Groups(\Sets) \xhookrightarrow{\;} \Groups(\SmoothInfinityGroupoids)$,

  -- $G \acts \, \Gamma \,\in\, \Groups\left( \Actions{G}(\DTopologicalSpaces) \right)
    \xhookrightarrow{\quad} \Groups\left( \Actions{G}(\SmoothInfinityGroupoids) \right) $,

  -- $G \acts \,  \SmoothManifold \,\in\,
    \Actions{G}(\SmoothManifolds)$,

  \noindent
  the operation of forming the homotopy fiber of modulating morphisms
  lands in
  $G$-action groupoids of bundles in $\DTopologicalSpaces \xhookrightarrow{\;} \SmoothInfinityGroupoids$
  and constitutes a natural fully faithful embedding
  \vspace{-2mm}
  \begin{equation}
    \label{EmbeddingOfGroupoidOfEquivariantCechCocyclesIntoGroupoidOfEquivariantPrincipalBundles}
    \begin{tikzcd}
      \SlicePointsMaps{\big}{B G}
        { \HomotopyQuotient{ \SmoothManifold }{ G } }
        { \mathbf{B}(\Gamma \rtimes G) }
      \ar[r, phantom, "\simeq"]
      \ar[
        rr,
        rounded corners,
        to path={
          -- ([yshift=-10pt]\tikztostart.south)
          -- node[above]{ \scalebox{.7}{$ \mathrm{hofib} $} }
             ([yshift=-11pt]\tikztotarget.south)
          -- ([yshift=-00pt]\tikztotarget.south)
        }
      ]
      &
      \EquivariantPrincipalBundles{G}{\Gamma}(\SmoothInfinityGroupoids)_{\SmoothManifold}
      \ar[r, hook]
      &
      \EquivariantPrincipalBundles{G}{\Gamma}(\DTopologicalSpaces)_{\SmoothManifold}
    \end{tikzcd}
  \end{equation}
  of the groupoids of $G$-equivariant $\Gamma$-principal bundles
  internal to $\SmoothInfinityGroupoids$ and internal to $\DTopologicalSpaces$, respectively.
\end{proposition}

\begin{proof}
Abbreviating $\Topos \coloneqq \SmoothInfinityGroupoids$,
we want to compute the operation that sends a cocycle $c$ to its homotopy fiber
in $\SliceTopos{\mathbf{B}G}$, which is equivalently
given by the following homotopy pullback in $\Topos$:
\vspace{-2mm}
\begin{equation}
  \label{HomotopyPullbackExtractingEquivariantPrincipalBundleFromCocycle}
  \begin{tikzcd}[row sep=small, column sep=large]
    \HomotopyQuotient
      { \TopologicalPrincipalBundle }
      { G }
    \ar[d]
    \ar[rr]
    \ar[drr, phantom, "\mbox{\tiny\rm (pb)}"]
    &&
    \mathbf{B}G
    \ar[d]
    \\
    \HomotopyQuotient
      { \SmoothManifold }
      { G }
    \ar[dr]
    \ar[rr, "c"]
    &&
    \mathbf{B}(\Gamma \rtimes G)
    \ar[dl]
    \\[-10pt]
    &
    \mathbf{B}G
  \end{tikzcd}
\end{equation}

\vspace{-2mm}
\noindent
We will compute this in model category $\SimplicialPresheaves(\CartesianSpaces)_{\projloc}$
by using the following representatives:

\begin{enumerate}[{\bf (1)}]

\vspace{-.2cm}
\item
By Prop. \ref{SmoothGManifoldsAdmitProperlyEquivariantGoodOpenCovers}),
we may find a properly equivariant good open cover
  $
    G \acts \; \widehat{\SmoothManifold}
    \coloneqq
    G \acts \;
    (
      \underset{i \in I}{\sqcup}
      \TopologicalPatch_i
    )
    \twoheadrightarrow
    G \acts \, \SmoothManifold
  $ (Def. \ref{ProperEquivariantOpenCover})
and,
by Ex. \ref{CechActionGroupoidOfEquivariantGoodOpenCoverIsLocalCofibrantResolution},
this yields a local projective cofibrant representative
$
  N
  (
    \widehat{\SmoothManifold} \times_{\SmoothManifold} \widehat{\SmoothManifold} \times G^{\mathrm{op}}
      \rightrightarrows
    \widehat{\SmoothManifold}
  )
$ of
$\HomotopyQuotient{\SmoothManifold}{G}$.

\vspace{-.2cm}
\item
By Prop. \ref{DeloopingGroupoidsOfDTopologicalGroupsAreLocallyFibrant},
we have a local projective fibrant representative for $\mathbf{B}(\Gamma \rtimes G)$
given by $N(\Gamma \rtimes G \rightrightarrows \ast)
\,\in\, \SimplicialPresheaves(\CartesianSpaces)_{\projloc}$\,.

\vspace{-.2cm}
\item
We have the evident projective fibration
\vspace{-2mm}
$$
  \begin{tikzcd}[row sep=4pt]
    \SimplicialNerve
    \DeloopingGroupoid{ \Gamma \rtimes G }
    \ar[dd , "\in \ProjectiveFibrations"]
    &&
    \SimplicialNerve
    \DeloopingGroupoid{ \Maps{}{\mathbb{R}^n}{\Gamma \rtimes G} }
    \ar[dd, "\in \KanFibrations"]
    \\
    &
    \quad
    : \quad \mathbb{R}^n \;\;\longmapsto
    \\
    \SimplicialNerve
    \DeloopingGroupoid{ G }
    &&
    \SimplicialNerve
    \DeloopingGroupoid{ \Maps{}{\mathbb{R}^n}{\Gamma \rtimes G} }.
  \end{tikzcd}
$$

\vspace{-.2cm}
\item
By the simplicial model enrichment of $\SimplicialPresheaves(\CartesianSpaces)_{\proj}$,
this implies the following Kan fibration model:
\vspace{-2mm}
$$
  \begin{tikzcd}[row sep=small]
  \PointsMaps{\big}
    { \HomotopyQuotient{ \SmoothManifold }{ G } }
    { \mathbf{B}(\Gamma \rtimes G) }
  \ar[r, phantom, "\simeq"]
  \ar[d]
  &
  \SimplicialPresheaves(\CartesianSpaces)
  \Big(
    \SimplicialNerve
    (
    \widehat{\SmoothManifold} \times_{\SmoothManifold} \widehat{\SmoothManifold} \times G^{\mathrm{op}}
      \rightrightarrows
    \widehat{\SmoothManifold}
    )
    ,\,
    \SimplicialNerve
    ( \Gamma \rtimes G \rightrightarrows \ast )
  \Big)
  \ar[
    d,
    " \in \KanFibrations "
  ]
  \\
  \PointsMaps{\big}
    { \HomotopyQuotient{ \SmoothManifold }{ G } }
    { \mathbf{B}G }
  \ar[r, phantom, "\simeq"]
  &
  \SimplicialPresheaves(\CartesianSpaces)
  \Big(
    \SimplicialNerve
    (
    \widehat{\SmoothManifold} \times_{\SmoothManifold} \widehat{\SmoothManifold} \times G^{\mathrm{op}}
      \rightrightarrows
    \widehat{\SmoothManifold}
    )
    ,\,
    \SimplicialNerve
    ( G \rightrightarrows \ast )
  \Big).
  \end{tikzcd}
$$

\vspace{-3mm}
\item
Hence we have a presentation of the equivariant bundles by
equivariant {\v C}ech cocycles (Rem. \ref{EquivariantCechCocycles}):

\vspace{-3mm}
\begin{equation}
  \label{EquivariantCocycleRegardedInSliceMappingSpace}
  \hspace{-9mm}
  \begin{aligned}
    \SlicePointsMaps{\big}{\mathbf{B}G}
    { \HomotopyQuotient{\SmoothManifold}{G} }
    { \mathbf{B}(\Gamma \rtimes G) }
  &
  \;\simeq\;
  \mathrm{fib}
  \Bigg(\!\!
  {\small  \begin{tikzcd}[row sep=7pt]
      \mathbf{H}\left(
        { \HomotopyQuotient{\SmoothManifold}{G} },
        { \mathbf{B}(\Gamma \rtimes G) }
        \right)
      \ar[d, "\,\mathbf{B}\mathrm{pr}_2"]
      \\
      \PointsMaps{}
        { \HomotopyQuotient{\SmoothManifold}{G} }
        { \mathbf{B}G }
    \end{tikzcd}
    }
  \!\!\Bigg)
  \;\;\;\;\;\;
  \proofstep{ by Prop. \ref{HomSpaceInSliceToposAsFiberProduct} }
  \\
  & \;\simeq\;
  \mathrm{fib}
  \Bigggg(\!\!
    \begin{tikzcd}[row sep=5pt]
      \SimplicialPresheaves(\CartesianSpaces)
      \Big(
        \SimplicialNerve
        \big(
          \widehat{\SmoothManifold}
          \times_{\SmoothManifold}
          \widehat{\SmoothManifold}
          \times
          G^\op
          \rightrightarrows
          \widehat{\SmoothManifold}
        \big)
        ,\,
        \SimplicialNerve
        \DeloopingGroupoid{ \Gamma \rtimes G }
      \!\Big)
      \ar[d, shorten=-4pt]
      \\
      \SimplicialPresheaves(\CartesianSpaces)
      \Big(
          \SimplicialNerve
          \big(
            \widehat{\SmoothManifold}
            \times_{\SmoothManifold}
            \widehat{\SmoothManifold}
            \times
            G^\op
            \rightrightarrows
            \widehat{\SmoothManifold}
          \big)
          ,\,
          \SimplicialNerve
          \DeloopingGroupoid{ G }
        \big)
      \!\Big)
    \end{tikzcd}
  \!\!\Bigggg)
  \\
  &
  \;\simeq\;
      \SimplicialPresheaves(\CartesianSpaces)
      \Big(
        \SimplicialNerve
        \big(
          \widehat{\SmoothManifold}
          \times_{\SmoothManifold}
          \widehat{\SmoothManifold}
          \times
          G^\op
          \rightrightarrows
          \widehat{\SmoothManifold}
        \big)
        ,\,
        \SimplicialNerve
        \DeloopingGroupoid{ \Gamma \rtimes G }
      \Big)_{
    \SimplicialNerve
    \DeloopingGroupoid{G}
  }
  \end{aligned}
\end{equation}

\end{enumerate}

\vspace{-2mm}
\noindent
{\bf (6)} Therefore, and using
Lem. \ref{ComputingHomotopyPullbacksOfInfinityStacks},
the homotopy pullback
\eqref{HomotopyPullbackExtractingEquivariantPrincipalBundleFromCocycle}
is equivalently computed by the following pullback of
simplicial presheaves
(in equivariant generalization of Rem. \ref{UniversalPrincipalBundleOverTheModuliStack}):

\vspace{-4mm}
\begin{equation}
  \label{EquivariantPrincipalBundleExtractedFromEquivariantCechCocycle}
  \hspace{-5mm}
  \begin{tikzcd}[row sep=small, column sep=-12pt]
    \TopologicalPrincipalBundle
    \ar[out=180-66, in=66, looseness=3.5, "\scalebox{.77}{$\mathclap{
      G
    }$}"{description},shift right=1]
    \ar[d]
    &[50pt]
    \big(
      \TopologicalPrincipalBundle
        \times
      G^{\mathrm{op}}
      \rightrightarrows
      \TopologicalPrincipalBundle
    \big)
    \ar[from=r, ->>, "{ \in \LocalWeakEquivalences }"  ]
    \ar[d]
    &[50pt]
    \big(
      \Gamma
        \times
      \widehat{\SmoothManifold}
        \times_{\SmoothManifold}
      \widehat{\SmoothManifold}
        \times
      G^{\mathrm{op}}
      \rightrightarrows
      \Gamma
        \times
      \widehat{\SmoothManifold}
    \big)
    \ar[d]
    \ar[rr]
    \ar[drr, phantom, "\mbox{\tiny\rm(pb)}"]
    &&
    \left(
      \Gamma_{\widehat{\mathrm{adj}}}
      \times
      (\Gamma \rtimes G)
      \rightrightarrows
      \Gamma
    \right)
    \ar[d, "\in \ProjectiveFibrations"]
    \\
    \SmoothManifold
      \ar[shift left=-4pt, out=-180+66, in=-66, looseness=3.5, "\scalebox{.77}{$\mathclap{
        G
      }$}"{description},shift right=1]
    &
    \big(
      \SmoothManifold
        \times
      G^{\mathrm{op}}
      \rightrightarrows
      \SmoothManifold
    \big)
    \ar[from=r, ->>, "{\in \LocalWeakEquivalences}"  ]
    &
    \big(
      \widehat{\SmoothManifold}
        \times_{\SmoothManifold}
      \widehat{\SmoothManifold}
        \times
      G^{\mathrm{op}}
      \rightrightarrows
      \widehat{\SmoothManifold}
    \big)
    \ar[dr]
    \ar[
      rr
    ]
    &&
    \DeloopingGroupoid
      { (\Gamma \rtimes G) }
    \ar[dl]
    \\
    &
    &
    &
    \SimplicialNerve
    \DeloopingGroupoid
      { G }
  \end{tikzcd}
\end{equation}

\vspace{-2mm}
\noindent
The equivariant principal bundles obtained this way
are by construction all locally trivializable over $\widehat {\SmoothManifold}$,
which means that the morphisms between them are in natural bijection
to the morphsims between their local trivialization data, hence between
the {\v C}ech cocycles on the right.
\end{proof}

\begin{lemma}
[Equivariant principal bundles constructed from {\v C}ech cocycles are proper equivariant fibrations]
  \label{EquivariantPrincipalBundlesConstructedFromCechCocyclesAreProperEquivariantFibrations}
  The $G$-equivariant $\Gamma$-principal $\infty$-bundles
  obtained as in
  \eqref{EquivariantPrincipalBundleExtractedFromEquivariantCechCocycle}
  are proper equivariant Serre fibrations
  (as in Prop. \ref{ProperEquivariantModelCategoryOfGSpaces}).
\end{lemma}
\begin{proof}
  By Ex. \ref{CechActionGroupoidOfEquivariantGoodOpenCoverIsLocalCofibrantResolution},
  we may compute cocycles on {\v C}ech-action groupoids of equivariant good open covers.
  The underlying bundle
  is evidently locally trivial while only the $G$-action may vary.
  Hence fixed-point-wise we have locally trivial bundles, hence Serre fibrations.
\end{proof}

\begin{lemma}
[Concordant equivariant principal bundles constructed from {\v C}ech cocycles are isomorphic]
  \label{ConcordantEquivariantPrincipalBundlesConstructedFromCechCocyclesAreIsomorphic}
  Equivariant principal bundles in the image of
  \eqref{EmbeddingOfGroupoidOfEquivariantCechCocyclesIntoGroupoidOfEquivariantPrincipalBundles}
  are isomorphic as soon as they are concordant.
\end{lemma}
\begin{proof}
  With Lem. \ref{EquivariantPrincipalBundlesConstructedFromCechCocyclesAreProperEquivariantFibrations},
  this follows by the same proof as Thm. \ref{ConcordanceClassesOfTopologicalPrincipalBundles}.
\end{proof}

In equivariant generalization of
Thm. \ref{OrdinaryPrincipalBundlesAmongPrincipalInfinityBundles}
and Thm. \ref{ClassificationOfPrincipalBundlesAmongPrincipalInfinityBundles},
we now have:
\begin{theorem}[Isomorphism classification of stable equivariant bundles with truncated structure over good orbifolds with resolvable singularities]
  \label{BorelClassificationOfEquivariantBundlesForResolvableSingularitiesAndEquivariantStructure}
  For

  -- $G \,\in\, \Groups(\FiniteSets)_{\resolvable}$
  (Ntn. \ref{ResolvableOrbiSingularities}),

  \vspace{0.5mm}
  -- $\Gamma \,\in\, \Actions{G}\big(\Groups(\DTopologicalSpaces)\big)$
  well-pointed (Ntn. \ref{WellPointedTopologicalGroup}),
  of truncated classifying shape (Ntn. \ref{CohesiveGroupsWithTruncatedClassifyingShape}),
  with countably many connected components, and
  with a notion of blowup-stability of $G$-equivariant
  $\Gamma \rtimes G$-principal bundles (Ntn. \ref{StableEquivariantBundles})

  \noindent
  then:

  \noindent
  {\bf (i)}
  the construction
  \eqref{EmbeddingOfGroupoidOfEquivariantCechCocyclesIntoGroupoidOfEquivariantPrincipalBundles}
  from Prop. \ref{EquivariantPrincipalBundlesFromCechCocycles}
  on stable equivariant bundles
  is an equivalence
  onto the
  full subgroupoid of
  equivariantly locally trivial bundles (Def. \ref{TerminologyForPrincipalBundles}):
  \vspace{-2mm}
  \begin{equation}
    \label{EquivalenceOfGroupoidOfEquivariantCechCocyclesIntoGroupoidOfEquivariantPrincipalBundles}
    \begin{tikzcd}
      \SlicePointsMaps{\big}{\mathbf{B} G}
        { \HomotopyQuotient{ \SmoothManifold }{ G } }
        { \mathbf{B}(\Gamma \rtimes G) }
      ^{\stable}
      \ar[r, phantom, "\simeq"]
      \ar[
        rr,
        rounded corners,
        to path={
          -- ([yshift=-10pt]\tikztostart.south)
          -- node[above]{ \scalebox{.7}{$ \mathrm{hofib} $} }
             node[below]{ \scalebox{.7}{$\sim$} }
             ([yshift=-11pt]\tikztotarget.south)
          -- ([yshift=-00pt]\tikztotarget.south)
        }
      ]
      &[-20pt]
      \EquivariantPrincipalBundles
        {G}{\Gamma}
      (\SmoothInfinityGroupoids)_{\SmoothManifold}^{\stable}
      \ar[r, "\sim"]
      &[-7pt]
      \EquivariantPrincipalFiberBundles
        {G}{\Gamma}
      (\DTopologicalSpaces)_{\SmoothManifold}^{\stable}
      \,.
    \end{tikzcd}
  \end{equation}

  \vspace{-2mm}
  \noindent
  {\bf (ii)}
  Isomorphism classes of these equivariant topological bundles are classified
  by
  $\shape \, \mathbf{B}(\Gamma \rtimes G)
    \simeq B \, \shape \, \Gamma \rtimes G
    \,\in\,
    (\InfinityGroupoids)_{/\mathbf{B}G}$
  (see Ex. \ref{ClassifyingShapesForDiscreteStructureInfinityGroups}),
  hence coincide with
  equivalence classes of
  of $G$-equivariant $\shape\, \Gamma$-principal
  $\infty$-bundles (Def. \ref{GEquivariantGammaPrincipalBundles})
  in that:
  \vspace{-2mm}
  \begin{align}
    \label{IsomorphismClassificationOfEquivariantPrincipalBundlesAsCechCocycles}
    \IsomorphismClasses{
      \EquivariantPrincipalBundles
        {G}{\Gamma}
      (\SmoothInfinityGroupoids)_{\SmoothManifold}^{\stable}
    }
    \;\;
    & \simeq
    \;\;
    \Truncation{0}
    \,
    \SliceMaps{\big}{\mathbf{B}G}
      { \HomotopyQuotient{ \shape \, \SmoothManifold }{ G } }
      { \HomotopyQuotient{ B\Gamma }{ G } }
    ^{\stable}
    \nonumber
    \\
    \;\;
    & \simeq
    \;\;
    \IsomorphismClasses{
      \EquivariantPrincipalBundles
        {G}{\shape\, \Gamma}
      (\SmoothInfinityGroupoids)_{\SmoothManifold}^{\stable}
    }
    \,.
 \end{align}
\end{theorem}
\begin{proof}
The second statement is the following composite bijection:
\vspace{-2mm}
$$
  \def\arraystretch{1.5}
  \begin{array}{lll}
    \IsomorphismClasses{
      \EquivariantPrincipalBundles
        {G}{\Gamma}
      (\SmoothInfinityGroupoids)
      ^{\stable}
    }
    &
    \;\simeq\;
    \ConcordanceClasses{
      \EquivariantPrincipalBundles
        {G}{\Gamma}
      (\SmoothInfinityGroupoids)
      ^{\stable}
    }
    &
    \proofstep{
      by Lem. \ref{ConcordantEquivariantPrincipalBundlesConstructedFromCechCocyclesAreIsomorphic}
    }
    \\
    &
    \;\simeq\;
    \Truncation{0}
    \,
    \shape
    \,
    \SliceMaps{\big}{\mathbf{B}G}
      { \HomotopyQuotient{ \SmoothManifold }{ G } }
      { \HomotopyQuotient{ \mathbf{B}\Gamma }{ G } }
    ^{\stable}
    &
    \proofstep{ by Ex. \ref{ConnectedComponentsOfShapeAreConcordanceClasses} }
    \\
    &
    \;\simeq\;
    \Truncation{0}
    \,
    \SliceMaps{\big}{\mathbf{B}G}
      { \HomotopyQuotient{ \shape \, \SmoothManifold }{ G } }
      { \HomotopyQuotient{ \shape \, \mathbf{B}\Gamma }{ G } }
    ^{\stable}
    &
    \proofstep{ by Thm. \ref{OrbiSmoothOkaPrinciple}. }
  \end{array}
$$

\vspace{-2mm}
\noindent
Using this we now obtain the first statement:
Consider an equivariant principal bundle
$G \acts \, \TopologicalPrincipalBundle
  \xrightarrow{\;}
  G \acts \, \SmoothManifold$
obtained as in
\eqref{EquivariantPrincipalBundleExtractedFromEquivariantCechCocycle}
relative to a choice of proper equivariant regular good open cover
(Def. \ref{ProperEquivariantOpenCover})
$G \acts \; ( \underset{i \in I}{\sqcup} \TopologicalPatch_i )
\twoheadrightarrow  G \acts \, \SmoothManifold$.
By Thm. \ref{EquivalentNotionsOfEquivariantLocalTriviality}
it is sufficient to exhibit,
for every $x \,\in\, \SmoothManifold$ with stabilizer subgroup
$G_x \,:=\, \mathrm{Stab}_G(x) \,\subset\, G$,
a $G_x$-equivariant neighborhood
$G_x \acts \, \TopologicalPatch_x$ of $x$ over which $G_x \acts  \, \TopologicalPrincipalBundle$
restricts to a Bierstone local model
\eqref{EquivariantPrincipalHatGxTwistedProductBundleIsCartesianProductProjectionOfActions}.

Now, by the fact that we have an open cover, there exists $i \in I$ with $x \in \TopologicalPatch_i$;
and by Rem. \ref{StabilizerSubgroupOfPointAlsoFixedIndexOfPatchInRegularEquivariantOpenCover}
this is already a $G_x$-neighborhood $\TopologicalPatch_x \,\coloneqq\, \TopologicalPatch_i$.
Restricting the construction \eqref{EquivariantPrincipalBundleExtractedFromEquivariantCechCocycle}
to this patch yields the following diagram:
\vspace{-2mm}
$$
\hspace{-1mm}
  \begin{tikzcd}[column sep=76pt]
    \SimplicialNerve
    \ActionGroupoid
      { (\Gamma \times \TopologicalPatch_x)_{\rho_x} }
      { G_x^\op }
    \ar[d]
    \ar[rr]
    \ar[drr, phantom, "\mbox{\tiny\rm(pb)}"]
    &&
    \SimplicialNerve
    \left(
      \Gamma_{\widehat{\mathrm{adj}}}
      \times
      (\Gamma \rtimes G)
      \rightrightarrows
      \Gamma
    \right)
    \ar[d]
    \\
    \SimplicialNerve
    \ActionGroupoid
      { \TopologicalPatch_x }
      { G_x^\op }
    \ar[
      r,
      "{
        \SimplicialNerve
        \ActionGroupoid
          { \TopologicalPatch_i \to \ast }
          { (-)^{-1} }
      }"
    ]
    \ar[
      rr,
      bend right=15pt,
      "{\ }"{name=t},
      "{c_{\vert \TopologicalPatch_i}}"{swap}
    ]
    &
    \SimplicialNerve
    \DeloopingGroupoid{ G_x }
    \ar[
      r,
      dashed,
      "{
       \scalebox{0.7}{$ \left(
        \mathrm{pr}
        \times
        \widehat{(-)}^{\, -1}
        \rightrightarrows
        \,
        \mathrm{pr}
        \right)
        $}
      }",
      "{\ }"{swap, name=s}
    ]
    &
    \SimplicialNerve
    \DeloopingGroupoid
      { (\Gamma \rtimes G) }
    \,.
    \ar[
      from=s,
      to=t,
      Rightarrow,
      "\sim"{swap, sloped}
    ]
  \end{tikzcd}
$$

\vspace{-2mm}
\noindent
Here the factorization,
up to some simplicial homotopy shown at the bottom,
of the restricted cocycle $c_{\vert \TopologicalPatch_i}$
through a dashed morphism (as shown) follows by the
classification statement
\eqref{IsomorphismClassificationOfEquivariantPrincipalBundlesAsCechCocycles}
which we have already proven, using now that $\TopologicalPatch_x$ is a patch
in a {\it good} open cover and hence contractible, so that any cocycle
for $G_x$-equivariant $\Gamma$-principal bundles on $\TopologicalPatch_x$
factors through
\vspace{-1mm}
$$
  \shape \,
  (
    \HomotopyQuotient{\TopologicalPatch_x}{G_x}
  )
  \;\simeq\;
  \HomotopyQuotient{\ast}{G_x}
  \;\simeq\;
  \Localization{\LocalWeakEquivalences}
  \left(
    \SimplicialNerve
    \DeloopingGroupoid{ G_x }
  \right)
  \,.
$$

\vspace{-1mm}
\noindent
But this dashed morphism of simplicial presheaves
is manifestly a crossed homomorphism in the
guise \eqref{CrossedHomomorphismsAreSectionsOftheSemidirectProductProjection},
hence a lift $\widehat{G_x}$ as in Ntn. \ref{LiftsOfEquivarianceSubgroupsToSemidirectProductWithStructureGroup}:
\vspace{-2mm}
\noindent
$$
  \begin{tikzcd}[row sep=small, column sep=huge]
    \SimplicialNerve
    \DeloopingGroupoid{ G_x }
    \ar[
      dr,
      "{
        \SimplicialNerve
        \DeloopingGroupoid{ (-)^{-1} }
      }"{swap}
    ]
    \ar[
      rr,
      "{
        \SimplicialPresheaves
        \DeloopingGroupoid
          { \widehat{(-)} }
      }"
    ]
    &&
    \SimplicialNerve
    \DeloopingGroupoid{ (\Gamma \rtimes G) }\,.
    \ar[
      dl,
      "{
        \DeloopingGroupoid{ \mathrm{pr}_2 }
      }"
    ]
    \\
    &
    \SimplicialNerve
    \DeloopingGroupoid{ G }
  \end{tikzcd}
$$

\vspace{-2mm}
\noindent
With $\Gamma_{\widehat{\mathrm{adj}}}$ from
Def. \ref{CrossedConjugationActionThroughACrossedHomomorphism},
this identifies the above pullback bundle as
the trivial $\Gamma$-bundle over $\TopologicalPatch_x$
equipped with the following action:
\vspace{-2mm}
$$
  \begin{tikzcd}[row sep=-4pt]
    (\Gamma \times \TopologicalPatch_x)_{\rho_x}
      \times
    G_x^{\mathrm{op}}
    \ar[rr]
    &&
    \Gamma \times \TopologicalPatch_x
    \\
 \scalebox{0.7}{$   \left(
      (\gamma, u), g
    \right)
    $}
    &\longmapsto&
  \scalebox{0.7}{$  \mathrm{pr}_1
    \left(
      (\NeutralElement, g)
      \cdot
      (\gamma, \NeutralElement)
      \cdot
      \widehat{g}^{-1}
    \right)
    $}
  \end{tikzcd}
$$

\vspace{-2mm}
\noindent
This is precisely the action on local Bierstone models found in
\eqref{ComputingTheCartesianProductActionOnBierstoneLocalModelBundles}.
Hence the above pullback exhibits the claimed equivariant local trivialization.
\end{proof}

\begin{example}[$G$-Equivariant principal bundles over $G$-coset spaces from basic geometric homotopy theory]
\label{GEquivariantPrincipalBundlesOverGCosetSpacesViaGeometricHomotopy}
For $G \acts \Gamma$ as in Thm. \ref{ClassificationOfPrincipalBundlesAmongPrincipalInfinityBundles},
we find, with Thm. \ref{BorelClassificationOfEquivariantBundlesForResolvableSingularitiesAndEquivariantStructure},
the form of locally trivial $G$-equivariant $\Gamma$-principal bundles on the coset spaces $G/H$
as originally consider by tom Dieck in 1969 (Lem. \ref{SemidirectProductCosetBundles})
by a glance at the following homotopy Cartesian diagram in $\SmoothInfinityGroupoids$ (using just the basic facts from Ex. \ref{QuotientStackOfCosetSpacesIsDeloopingOfSubgroup}. Ex. \ref{DoubleHomotopyCosetStacks} and Prop. \ref{ConjugationGroupoidOfCrossedHomomorphismsIsSectionsOfDeloopedSemidirectProductProjection}):
\begin{equation}
  \label{EquivariantPrincipalBundlesOverCosetSpacesFromGeometricHomotopyTheory}
  \begin{tikzcd}
     \overset{
        \widehat{G}
      }{
      \overbrace{
        (\Gamma \rtimes G)
      }
    }
    /\widehat{H}
    \ar[d]
    &[-10pt]
    \HomotopyQuotient
      {
        \big(
          (\Gamma \rtimes G)
          /\widehat{H}
        \big)
      }
      { G }
    \ar[rr]
    \ar[d]
    \ar[
      drr,
      phantom,
      "{
        \overset{
          \mbox{
            \tiny
            Ex. \ref{DoubleHomotopyCosetStacks}
          }
        }{
          \mbox{\tiny\rm(pb)}
        }
      }"{yshift=4pt}
    ]
    &[-13pt]
    &[+13pt]
    \mathbf{B}G
    \ar[d]
    \\
    G/H
    &
    \HomotopyQuotient
      { (G/H) }
      { G }
    \ar[dr]
    \ar[
      r,
      "{
        \sim
      }",
      "{
        \mbox{
          \tiny
          Ex. \ref{QuotientStackOfCosetSpacesIsDeloopingOfSubgroup}
        }
      }"{swap, yshift=-1pt}
    ]
    &
    \mathbf{B}H
    \ar[
      r,
      "{
        \mathbf{B}(i_{\widehat{H}})
      }",
      "{
        \mbox{
          \tiny
          Prop. \ref{ConjugationGroupoidOfCrossedHomomorphismsIsSectionsOfDeloopedSemidirectProductProjection}
        }
      }"{swap}
    ]
    \ar[
      d,
      "{
        \mathbf{B}(i_{{H}})
      }"{description}
    ]
    &
    \mathbf{B}(\Gamma \rtimes G)
    \ar[
      dl,
      "{
        \mathbf{B}\mathrm{pr}_2
      }"
    ]
    \\
    &
    &
    \mathbf{B}G
    \mathrlap{\,.}
  \end{tikzcd}
\end{equation}

\end{example}

\begin{lemma}[$(\Gamma \rtimes G)$-CW-Complex structure on equivariant bundles over $G$-CW complexes]
\label{EquivariantCWComplexStructureOnEquivariantBundles}
Given  $G \acts \Gamma$ as in Thm. \ref{BorelClassificationOfEquivariantBundlesForResolvableSingularitiesAndEquivariantStructure}, with $\Gamma$ connected, and given a $G$-equivariant $\Gamma$-principal fiber bundle $\TopologicalPrincipalBundle$ over a $G$-CW complex $\TopologicalSpace$

- $G \acts \TopologicalSpace \,\in\, \GCWComplexes$
 (Ex. \ref{GCWComplexesAreCofibrantObjectsInProperEquivariantModelcategory})

- $(\Gamma \rtimes G) \acts \TopologicalPrincipalBundle
\,\in\, \EquivariantPrincipalFiberBundles{G}{\Gamma}_{\TopologicalSpace}
$

\noindent
then $\TopologicalPrincipalBundle$ is a proper $(\Gamma \rtimes G)$-CW complex, namely a $(\Gamma \rtimes G)$-equivariant cell complex with respect to cell attachments of the form
$$
  (\Gamma \rtimes G)/
    \times
  \partial \TopologicalSimplex{n}
  \xhookrightarrow{
    \;\;
  }
  (\Gamma \rtimes G)/
    \times
  \TopologicalSimplex{n}
$$

\end{lemma}
\begin{proof}
By the classification result of Thm. \ref{BorelClassificationOfEquivariantBundlesForResolvableSingularitiesAndEquivariantStructure} we have a homotopy Cartesian square of the following form
$$
  \begin{tikzcd}
    \HomotopyQuotient
      { \TopologicalPrincipalBundle }
      { G }
    \ar[rr]
    \ar[
      d,
      "{
        \HomotopyQuotient
         { p }
         { G }
       }"
    ]
    \ar[
      drr,
      phantom,
      "{\mbox{\tiny\rm(pb)}}"{pos=.5}
    ]
    &&
    \mathbf{B}G
    \ar[d]
    \\
    \HomotopyQuotient
      { \TopologicalSpace }
      { G }
    \ar[dr]
    \ar[
      rr,
      "{c}"
    ]
    &&
    \mathbf{B}(\Gamma \rtimes G)
    \ar[dl]
    \\[-10pt]
    &
    \mathbf{B}G
  \end{tikzcd}
$$
Denote the assumed $G$-CW-complex structure of $\TopologicalSpace$
as follows:
$$
  \TopologicalSpace_0
  \;=\; \varnothing
  \,,
  \;\;\;\;\;\;\;\;\;\;\;
  \TopologicalSpace_{\bullet+1}
  \;\;=\;\;
  \TopologicalSpace_\bullet
  \;\;
  \underset{
    \mathclap{
    \raisebox{-2pt}{
      \scalebox{.7}{$
        G/H_\bullet \times \partial \TopologicalSimplex{n_\bullet}
      $}
    }
    }
  }{\coprod}
  \;\;
  G/H_\bullet \times \TopologicalSimplex{n_\bullet}
  \,,
  \hspace{1cm}
  \TopologicalSpace
  \;=\;
  \colimit{k \in \NaturalNumbers}
  \,
  \TopologicalSpace_k
$$
Due to colimits in $\SmoothInfinityGroupoids$ being universal (Prop. \ref{ColimitsAreUniversalInAnInfinityTopos})
we find from this the above homotopy pullback that
$$
  \TopologicalPrincipalBundle_\bullet
  \;=\;
  \TopologicalSpace_\bullet
  \underset{
    \mathbf{B}\Gamma
  }{\times}
  \ast
  \,,
  \hspace{1cm}
  \TopologicalPrincipalBundle
  \;=\;
  \colimit{k \in \NaturalNumbers}
  \;
  \TopologicalPrincipalBundle_k
  \,.
$$

Hence it is now sufficient to show that
$$
  \TopologicalPrincipalBundle_{\bullet+1}
  \;\;\;\simeq\;\;\;
  \TopologicalPrincipalBundle_\bullet
  \;\;
  \underset{
    \mathclap{
    \raisebox{-3pt}{
      \scalebox{.7}{$
        (\Gamma \rtimes G)/\widehat{H_\bullet} \times \partial \TopologicalSimplex{n_\bullet}
      $}
    }
    }
  }{\coprod}
  \;\;
  (\Gamma \rtimes G)/\widehat{H_\bullet} \times \TopologicalSimplex{n_\bullet}
$$
namely that for each $k \in \NaturalNumbers$
we have homotopy Cartesian squares of this form:
$$
  \begin{tikzcd}
    \HomotopyQuotient
      { \TopologicalPrincipalBundle_{k+1} }
      { G }
    \ar[
      d,
      shorten >=-7pt,
      "{
        \HomotopyQuotient
          { p_{k+1} }
          { G }
      }"
    ]
    &[-25pt]
    \simeq
    &[-25pt]
    \HomotopyQuotient
    {
    \Bigg(
    \;\;\;\;
    \TopologicalPrincipalBundle_{k}
    \;\;\;\;\;
    \underset{
      \mathclap{
      \raisebox{-4pt}{
        \scalebox{.7}{$
          (\Gamma \rtimes G)/\widehat{H_{k+1}} \times \partial\TopologicalSimplex{n_{k+1}}
        $}
      }
      }
    }{\coprod}
    \;\;\;\;\;
    (\Gamma \rtimes G)/\widehat{H_{k+1}} \times \TopologicalSimplex{n}
    \Bigg)
    }{G}
    \ar[
      d,
      shorten <=-8pt
    ]
    \ar[rr]
    \ar[
      drr,
      phantom,
      "{\mbox{\tiny(pb)}}"{pos=.4}
    ]
    &&
    \mathbf{B}G
    \ar[d]
    \\
    \HomotopyQuotient
      { \TopologicalSpace_{k+1} }
      { G }
    &
    \simeq
    &
    \HomotopyQuotient
    {
    \Bigg(
    \;\;
    \TopologicalSpace_{k}
    \;\;\;\;\;
    \underset{
      \mathclap{
      \raisebox{-4pt}{
        \scalebox{.7}{$
          G/{H_{k+1}}
            \times
          \partial\TopologicalSimplex{n_{k+1}}
        $}
      }
      }
    }{\coprod}
    \;\;\;\;\;
      G/{H_{k+1}}
        \times
      \TopologicalSimplex{n_{k+1}}
    \Bigg)
    }{G}
    \ar[
      rr,
      "{
        \Big(
          c_{k+1}
          ,\,
          \mathbf{B}
          \big(i_{\widehat{H}_{k+1}}\big)
        \Big)
      }"{description}
    ]
    \ar[dr, shorten <=-16pt]
    &&
    \mathbf{B}(\Gamma \rtimes G)
    \ar[dl]
    \\[-10pt]
    &&
    &
    \mathbf{B}G
    \mathrlap{\,.}
  \end{tikzcd}
$$
But (again by universal colimits and using the classification result Thm. \ref{BorelClassificationOfEquivariantBundlesForResolvableSingularitiesAndEquivariantStructure} and the connectivity assumption to deduce that the cocycle must be trivializable on the simplex), this is the case as soon as we have homotopy Cartesian squares of this form:
$$
  \begin{tikzcd}
    \HomotopyQuotient
      {
        \big(
          (\Gamma \rtimes G)/\widehat{H_{k+1}}
        \big)
      }
      { G }
    \ar[d]
    \ar[rr]
    \ar[
      drr,
      phantom,
      "{\mbox{\tiny\rm(pb)}}"{pos=.4}
    ]
    &[-22pt]
    &[+22pt]
    \mathbf{B}G
    \ar[d]
    \\
    \HomotopyQuotient
      { ( G/H_{k+1}) }
      { G }
    \ar[
      r,
      "{\sim}"
    ]
    \ar[dr]
    &
    \mathbf{B}H_{k+1}
    \ar[
      r,
      "{
        \mathbf{B}
        \big(
          i_{\widehat{H_{k+1}}}
        \big)
      }"{description}
    ]
    \ar[d]
    &
    \mathbf{B}(\Gamma \rtimes G)
    \ar[
      dl
    ]
    \\
    &
    \mathbf{B}G
    \mathrlap{\,,}
  \end{tikzcd}
$$
and this holds by Ex. \ref{GEquivariantPrincipalBundlesOverGCosetSpacesViaGeometricHomotopy}.
\end{proof}

\medskip

\noindent
{\bf Equivariant fiber bundles.}
Specializing the general notion of associated fiber $\infty$-bundles (\cite[Prop. 4.10]{NSS12a} \cite[Prop. 2.92]{SS20OrbifoldCohomology}) to the equivariant context of Def. \ref{GEquivariantGammaPrincipalBundles}, we obtain the following notion of equivariant fiber $\infty$-bundles, which by Thm. \ref{BorelClassificationOfEquivariantBundlesForResolvableSingularitiesAndEquivariantStructure} are immediately guaranteed to satisfy the traditional equivariant local triviality conditions.

\begin{definition}[Equivariant associated $\infty$-bundles]
\label{EquivariantAssociatedInfinityBundles}
For

- $G\acts \Gamma \,\in \,
       \EquivariantGroups{G}(\Topos)
       $
      (Def. \ref{GGroups})

- $(\HomotopyQuotient{\Gamma}{G}) \acts
   \HomotopyQuotient{A}{G}
    \,\in\,
    \Actions{(\HomotopyQuotient{\Gamma}{G})}
    \big(
      \SliceTopos{\mathbf{B}G}
    \big)
  $ (Prop. \ref{HomotopyQuotientsAndPrincipaInfinityBundles})

\noindent
we say that

\noindent
(i) the projection
 (recall the vertical equivalences from Prop. \ref{EquivariantGroupsEquivalentToSplitExtensions}
 and Lem. \ref{HomotopyQuotientsByEquivariantActionsOfEquivariantGroups})
\begin{equation}
  \label{UniversalEquivariantAFiberInfinityBundle}
  \begin{tikzcd}[row sep=8pt]
    \HomotopyQuotient
      { (\HomotopyQuotient{A}{G}) }
      { (\HomotopyQuotient{\Gamma}{G}) }
    \ar[rr]
    \ar[d, phantom, "\simeq"{rotate=-90}]
    &&
    \HomotopyQuotient
      { \mathbf{B}\Gamma }
      { G }
    \ar[d, phantom, "\simeq"{rotate=-90}]
    \\
    \HomotopyQuotient
      { A }
      { (\Gamma \rtimes G) }
    \ar[dr]
    \ar[rr]
    &&
    \mathbf{B}(\Gamma \rtimes G)
    \ar[dl]
    \\
    &
    \mathbf{B}G
  \end{tikzcd}
\end{equation}
is the {\it universal $G$-equivariant $\Gamma$-associated $A$-fiber $\infty$-bundle};

\noindent
(ii) given a $G$-equivariant $\Gamma$-principal $\infty$-bundle $P$
(Def. \ref{GEquivariantGammaPrincipalBundles})
its {\it associated $G$-equivariant  $A$-fiber $\infty$-bundle} is the pullback of the universal such \eqref{UniversalEquivariantAFiberInfinityBundle} along the modulating map of $P$:
\begin{equation}
  \label{AssociatedEquivariantInfinityBundlesByPullback}
  \adjustbox{raise=1pt}{
  \begin{tikzcd}
    \overset{
      \mathclap{
      \raisebox{3pt}{
        \tiny
        \color{darkblue}
        \bf
        \def\arraystretch{.9}
        \begin{tabular}{c}
          equivariant
          \\
          principal $\infty$-bundle
        \end{tabular}
      }
      }
    }{
    \HomotopyQuotient
      { P }
      { G }
    }
    \ar[
      rr
    ]
    \ar[dd]
    \ar[
      ddrr,
      phantom,
      "{
        \mbox{
          \tiny
          \rm
          (pb)
        }
      }"
    ]
    &&
    \overset{
      \mathclap{
      \scalebox{.8}{$
        \def\arraystretch{.9}
        \begin{array}{c}
          \HomotopyQuotient
            { \Gamma }
            { (\Gamma \rtimes G) }
          \\
          \rotatebox[origin=c]{-90}{
            $\simeq$
          }
        \end{array}
      $}
      }
    }{
      \mathbf{B}G
    }
    \ar[dd]
    \\
    \\
    \HomotopyQuotient
      { X }
      { G }
    \ar[
      rr,
      "{
        \vdash
        \,
        \HomotopyQuotient{P}{G}
      }"
    ]
    \ar{dr}
    &&
    \mathbf{B}(\Gamma \rtimes G)
    \ar[dl]
    \\
    &
    \mathbf{B}G
  \end{tikzcd}
  }
  \hspace{1cm}
  \vdash
  \hspace{1cm}
  \begin{tikzcd}
  \overset{
      \mathclap{
      \raisebox{3pt}{
        \tiny
        \color{darkblue}
        \bf
        \def\arraystretch{.9}
        \begin{tabular}{c}
          equivariant associated
          \\
          $A$-fiber $\infty$-bundle
        \end{tabular}
      }
      }
  }{
    \HomotopyQuotient
      { E }
      { G }
    }
    \ar[
      rr
    ]
    \ar[dd]
    \ar[
      ddrr,
      phantom,
      "{
        \mbox{
          \tiny
          \rm
          (pb)
        }
      }"
    ]
    &&
    \overset{
      \mathclap{
      \raisebox{3pt}{
        \tiny
        \color{darkblue}
        \bf
        \def\arraystretch{.9}
        \begin{tabular}{c}
          universal $G$-equivariant
          \\
          $\Gamma$-associated
          \\
          $A$-fiber $\infty$-bundle
        \end{tabular}
      }
      }
    }{
    \HomotopyQuotient
      {A}
      { (\Gamma \rtimes G) }
    }
    \ar[dd]
    \\
    \\
    \HomotopyQuotient
      { X }
      { G }
    \ar[
      rr,
      "{
        \vdash
        \,
        \HomotopyQuotient{P}{G}
      }"
    ]
    \ar{dr}
    &&
    \mathbf{B}(\Gamma \rtimes G)
    \ar[dl]
    \\
    &
    \mathbf{B}G
  \end{tikzcd}
\end{equation}
Conversely, given an $A$-fiber $\infty$-bundle arising as on the right, we say that its {equivariant structure $\infty$-group} is $\Gamma$.
\end{definition}

It is a general fact that sections of $\Gamma$-associated $A$-fiber $\infty$-bundles are equivalently $\Gamma$-equivariant maps into $A$ out of the underlying $\Gamma$-principal $\infty$-bundle \cite[Cor. 4.18]{NSS12a}. The following lemma spells this out for equivariant $\infty$-bundles.

\begin{proposition}[$G$-Equivariant sections
of $\Gamma$-associated equivariant bundles as $(\Gamma \rtimes G)$-equivariant maps]
\label{EquivariantSectionsOfAssociatedEquivariantBundles}
Given $G \acts \Gamma \,\in\,\EquivariantGroups{G}(\Topos)$
(Def. \ref{GGroups})
with semidirect product \eqref{DeloopingOfSplitExtensionIsHomotopyQuotientOfBGamma} to be abbreviated as
$
  \widehat{G} \,:=\, \Gamma \rtimes G\,,
$
and
$P$ a $G$-equivariant $\Gamma$-principal $\infty$-bundle over $X \,\in\, \Topos$, modulated
 (according to Def. \ref{GEquivariantGammaPrincipalBundles})
 by
  $
    \Twist
    \,:\,
    \HomotopyQuotient{X}{G}
    \xrightarrow{\;\vdash (\HomotopyQuotient{P}{G})\;}
    \mathbf{B}(\Gamma \rtimes G)
  $
with $E$ an associated $A$-fiber $\infty$-bundle (according to Def. \ref{EquivariantAssociatedInfinityBundles}), the following sliced mapping stacks (Def. \ref{SliceMappingStack})
are naturally equivalent:
$$
  \def\arraystretch{4}
  \setlength{\arraycolsep}{2pt}
  \begin{array}{ccccccc}
  \overset{
    \mathclap{
    \raisebox{8pt}{
      \tiny
      \color{darkblue}
      \bf
      \def\arraystretch{.9}
      \begin{tabular}{c}
        $\infty$-stack of equivariant sections of
        \\
        associated
        equivariant $A$-fiber $\infty$-bundle
      \end{tabular}
    }
    }
  }{
  \SliceMaps{\big}{
    \HomotopyQuotient{X}{G}
  }
    {
      \HomotopyQuotient
        { \TopologicalSpace }
        { G }
    }
    {
      \HomotopyQuotient
        {
          E
        }
        {
          G
        }
    }
  }
  &\simeq&
  \overset{
    \mathclap{
    \raisebox{8pt}{
      \tiny
      \color{darkblue}
      \bf
      \def\arraystretch{.9}
      \begin{tabular}{c}
        $\infty$-stack of lifts of modulating map through
        \\
        universal
        equivariant
        $\Gamma$-associated
        $A$-fiber $\infty$-bundle
      \end{tabular}
    }
    }
  }{
  \SliceMaps{\big}{\mathbf{B}\widehat{G}}
    {
      \HomotopyQuotient
        { (\TopologicalSpace,\Twist) }
        { G }
    }
    {
      \HomotopyQuotient
        {
          \TopologicalCoefficients
        }
        {
          \widehat{G}
        }
    }
  }
  &\simeq&
  \overset{
    \mathclap{
    \raisebox{8pt}{
      \tiny
      \color{darkblue}
      \bf
      \def\arraystretch{.9}
      \begin{tabular}{c}
        $\infty$-stack of maps
        of homotopy quotients in slice
        \\
        from equivariant principal bundle to
        typical fiber
      \end{tabular}
    }
    }
  }{
  \SliceMaps{\big}{\mathbf{B}\widehat{G}}
    {
      \HomotopyQuotient
        { \TopologicalPrincipalBundle }
        { \widehat{G} }
    }
    {
      \HomotopyQuotient
        {
          \TopologicalCoefficients
        }
        {
          \widehat{G}
        }
    }
  }
  &\simeq&
    \overset{
    \mathclap{
    \raisebox{2pt}{
      \tiny
      \color{darkblue}
      \bf
      \def\arraystretch{.9}
      \begin{tabular}{c}
        $\infty$-stack of $\Gamma \rtimes G$-equivariant maps from
        \\
        equivariant principal bundle
        to typical fiber
      \end{tabular}
    }
    }
    }{
    \Maps{\Big}
     {
       \widehat{G}
       \acts
       \TopologicalPrincipalBundle
     }
     {
       \widehat{G}
       \acts
       \TopologicalCoefficients
     }
    ^{ \widehat{G} }
    }
    \\
  \left\{
  \hspace{-.1cm}
  \begin{tikzcd}[
    column sep=-4pt,
    row sep=10pt
  ]
    &&
    \HomotopyQuotient
      { E }{G}
    \ar[d]
    \\
    \HomotopyQuotient
      { X }{G}
    \ar[rr, Rightarrow, -]
    \ar[dr, shorten=-2pt]
    \ar[
      urr,
      dashed
    ]
    &&
    \HomotopyQuotient
      { X }{G}
    \ar[dl, shorten=-2pt]
    \\[-7pt]
    &
    \mathbf{B}G
  \end{tikzcd}
  \right\}
  &&
  \left\{
  \hspace{-.1cm}
  \begin{tikzcd}[
    column sep=-4pt,
    row sep=10pt
  ]
    &&
    \HomotopyQuotient
      { A }{(\Gamma \rtimes G)}
    \ar[d]
    \\
    \HomotopyQuotient
      { X }{G}
    \ar[rr, shorten=-2pt, "{\Twist}"{description, pos=.4}]
    \ar[dr, shorten=-2pt]
    \ar[
      urr,
      dashed
    ]
    &&
    \mathbf{B}(\Gamma \rtimes G)
    \ar[dl, shorten=-2pt]
    \\[-7pt]
    &
    \mathbf{B}G
  \end{tikzcd}
  \hspace{-.1cm}
  \right\}
  &&
  \left\{
  \hspace{-.1cm}
  \begin{tikzcd}[
    column sep=-9pt,
    row sep=10pt
  ]
    &&
    \HomotopyQuotient
      { A }{(\Gamma \rtimes G)}
    \ar[d]
    \\
    \HomotopyQuotient
      { P }{\widehat{G}}
    \ar[rr, shorten=-1pt]
    \ar[dr, shorten=-2pt]
    \ar[
      urr,
      dashed
    ]
    &&
    \mathbf{B}(\Gamma \rtimes G)
    \ar[dl, shorten=-2pt]
    \\[-7pt]
    &
    \mathbf{B}G
  \end{tikzcd}
  \hspace{-.1cm}
  \right\}
  &&
  \left\{
  \hspace{-.02cm}
  \adjustbox{raise=-7pt}{
  \begin{tikzcd}[column sep=16pt]
    P
      \ar[out=180-66, in=66, looseness=4, "\scalebox{.85}{$\;\mathclap{
        \Gamma \rtimes G
      }\;$}"{description},shift right=1]
    \ar[rr, dashed]
    &&
    A
      \ar[out=180-66, in=66, looseness=4, "\scalebox{.85}{$\;\mathclap{
        \Gamma \rtimes G
      }\;$}"{description},shift right=1]
  \end{tikzcd}
  }
  \hspace{-.02cm}
  \right\}
  \end{array}
$$
\end{proposition}
\begin{proof}
(i)
Under the identification  \eqref{PlotsOfSliceMappingStack}
of plots of the slice mapping stack with slice homs,
and using the $\infty$-Yoneda lemma (Prop. \ref{InfinityYonedaLemma}),
the first equivalence is the universal property of the defining pullback \eqref{AssociatedEquivariantInfinityBundlesByPullback}.
(ii) The second equivalence follows immediately from the identification
$
  \HomotopyQuotient
    { P }{ (\Gamma \rtimes G) }
  \;\simeq\;
  \HomotopyQuotient
    { X }
    { G }
  \;\;\;
  \in
  \;
  \SliceTopos{\mathbf{B}G}
$
\eqref{SlicedEquivalenceForSemidirectProductHomotopyQuotientOfEquivariantInfinityBundle}.
(iii) The third equivalence is a special case of Ex. \ref{EquivariantMappingStackAsSliceMappingStack}.
\end{proof}

\medskip

\noindent
{\bf Base change of concordances of equivariant bundles along coverings of the equivariance
group.}
We use the classification Theorem \ref{BorelClassificationOfEquivariantBundlesForResolvableSingularitiesAndEquivariantStructure}
to describe the base change of concordances of equivariant principal bundles
along coverings of the equivariance group (Lem. \ref{BaseChangeOfConcordancesOfEquivariantBundlesAlongCoveringOfEquivarianceGroup} below).

\begin{notation}[Equivariance group base change on mapping stacks]
  \label{BaseChangeOfEquivarianceGroupsOnMappingStacks}
  For

  -- $\Topos$ a cohesive $\infty$-topos with a 1-site $\Charts$ of charts (Ntn. \ref{CohesiveCharts}),

  -- $p : \widehat{G} \xrightarrow{\;} G \;\; \in \Groups(\Sets)$,

  -- $G \acts \, \mathcal{X}, \, G \acts \, \mathcal{A} \;\; \in \Actions{G}( \Topos)$,

  \noindent
  {\bf (i)}  we write
  \vspace{-2mm}
  $$
    (B p)^\ast_{ \mathcal{X}, \mathcal{A} }
    :
    \begin{tikzcd}
      \SliceMaps{\big}{\mathbf{B}G}
        { \HomotopyQuotient{ \mathcal{X} }{ G } }
        { \HomotopyQuotient{ \mathcal{A} }{ G } }
      \ar[r]
      &
      \SliceMaps{\big}{\mathbf{B}{\widehat{G}}}
        { \HomotopyQuotient{ \mathcal{X} }{ \widehat{G} } }
        { \HomotopyQuotient{ \mathcal{A} }{ \widehat{G} } }
    \end{tikzcd}
  $$
  for the component morphism of the
  base change
  of slice mapping stacks (Def. \ref{SliceMappingStack})
  along the delooping (Prop. \ref{LoopingAndDeloopingEquivalence}) of $p$.
  Namely:

  \noindent {\bf (ii)}
  Under Lem. \ref{PlotsOfSliceMappingStackAreSliceHoms}
  and the $\infty$-Yoneda lemma (Prop. \ref{InfinityYonedaLemma}), this
  is given at stage
  $U \,\in\, \Charts \xhookrightarrow{\YonedaEmbedding} \Topos$ by
  the $(U \times \mathcal{X}, \mathcal{Y})$-component of the base change
  functor (Prop. \ref{BaseChange}):
  \vspace{-2mm}
  $$
    (B p)^\ast_{ \mathcal{X}, \mathcal{A} }(U)
    \;=\;
    (B p)^\ast_{ U \times \mathcal{X}, \mathcal{A} }
    \;
    :
    \begin{tikzcd}
      \SlicePointsMaps{\big}{\mathbf{B}G}
        { \HomotopyQuotient{ U \times \mathcal{X} }{ G } }
        { \HomotopyQuotient{ \mathcal{A} }{ G } }
      \ar[r]
      &
      \SlicePointsMaps{\big}{\mathbf{B}{\widehat{G}}}
        { \HomotopyQuotient{ U \times \mathcal{X} }{ \widehat{G} } }
        { \HomotopyQuotient{ \mathcal{A} }{ \widehat{G} } }
      \,.
    \end{tikzcd}
  $$
\end{notation}

\begin{lemma}[Base change of concordances of equivariant bundles along covering of equivariance group]
\label{BaseChangeOfConcordancesOfEquivariantBundlesAlongCoveringOfEquivarianceGroup}
For

-- $\widehat G \,\in\, \Groups(\FiniteSets)_{\resolvable}$ (Ntn. \ref{ResolvableOrbiSingularities}),

-- $p : \begin{tikzcd} \widehat{G} \ar[r,->>] &[-12pt] G  \end{tikzcd}$
a surjective homomorphism \eqref{ADiscreteGroupEpimorphismInAnInfinityTopos},

-- $G \acts \, \SmoothManifold
   \,\in\, \Actions{G}(\SmoothManifolds)
     \xhookrightarrow{\;}
   \Actions{G}(\SmoothInfinityGroupoids)$,

-- $G \acts \, \Gamma \,\in\, \Actions{G}\left(\Groups(\kHausdorffSpaces)\right)
 \xhookrightarrow{  \ContinuousDiffeology } \Actions{G}\left( \Groups(\SmoothInfinityGroupoids)\right)$,
 well-pointed (Ntn. \ref{WellPointedTopologicalGroup})
 and of truncated classifying shape (Ntn. \ref{CohesiveGroupsWithTruncatedClassifyingShape})
 with a notion of stable equivariant bundles (Ntn. \ref{StableEquivariantBundles}),

\noindent
the
shape \eqref{CohesiveAdjointQuadruple}
of the base change
$
  (B p)^\ast_{
    \HomotopyQuotient
      { \SmoothManifold }
      { G }
    ,
    \mathbf{B}\Gamma
  }$
  on slice mapping stacks
(Ntn. \ref{BaseChangeOfEquivarianceGroupsOnMappingStacks}) --
hence the induced morphism
\eqref{ComparisonMorphismFromShapeOfMappingStackToMappingSpaceOfShapes}
of concordance $\infty$-groupoids
(Rem. \ref{TheInfinityGroupoidOfConcordances})
of equivariant $\Gamma$-principal bundles --
is a monomorphism (Ex. \ref{MonomorphismsOfInfinityGroupoids})
for stable equivariant bundles:
\vspace{-2mm}
$$
  \shape
  \,
  \SliceMaps{\big}{\mathbf{B}G}
    { \HomotopyQuotient{ \SmoothManifold }{ G } }
    { \HomotopyQuotient{ \mathbf{B}\Gamma }{ G }  }
  ^{\stable}
  \xhookrightarrow{\quad \shape (B p)^\ast_{\mathcal{X},\mathbf{B}\Gamma} \quad }
  \shape
  \,
  \SliceMaps{\big}{\mathbf{B} {\widehat{G}} }
    { \HomotopyQuotient{ \SmoothManifold }{ \widehat{G} } }
    { \HomotopyQuotient{ \mathbf{B}\Gamma }{ \widehat{G} }  }
  ^{\stable}
  \,.
$$
\end{lemma}
\begin{proof}
  For $\SmoothSimplex{\bullet} \,\in\, \SmoothManifolds \hookrightarrow \SmoothInfinityGroupoids$
  (Def. \ref{SmoothExtendedSimplicies}),
  Prop. \ref{PullbackOfActionsAlongSurjectiveGroupHomomorphismsIsFullyFaithful}
  gives that \eqref{BaseChangeOfTwoActionsAlongSurjectiveHomomorphism}
  \vspace{-2mm}
  $$
    (B p)^\ast_{\SmoothSimplex{\bullet} \times \SmoothManifold, \mathbf{B}\Gamma}
    \;\; : \;\;
    \SlicePointsMaps{\big}{\mathbf{B}G}
      { \SmoothSimplex{\bullet} \times \HomotopyQuotient{\SmoothManifold}{G} }
      { \HomotopyQuotient{ \mathbf{B}\Gamma }{G} }
    ^{\stable}
    \xhookrightarrow{\qquad} \;
    \SlicePointsMaps{\big}{\mathbf{B}{\widehat{G}}}
      { \SmoothSimplex{\bullet} \times \HomotopyQuotient{\SmoothManifold}{\widehat{G}} }
      { \HomotopyQuotient{ \mathbf{B}\Gamma }{\widehat{G}} }
    ^{\stable}
  $$
  is a degreewise monomorphism of simplicial $\infty$-groupoids.
  With
  Prop. \ref{SmoothShapeGivenBySmoothPathInfinityGroupoid}
  and Lem. \ref{PlotsOfSliceMappingStackAreSliceHoms},
  it follows that the
  morphism in question is equivalently the image under the simplicial $\infty$-colimit
  operation
  of a monomorphism of simplicial $\infty$-groupoids:
  \vspace{-2mm}
  $$
    \shape \, (B p)^\ast_{\SmoothManifold, \mathbf{B}\Gamma}
    \;\;
    \simeq
    \;\;
    \colimit{[n] \in \Delta^{\mathrm{op}}}
    \Big(
      \SlicePointsMaps{\big}{\mathbf{B}G}
        { \SmoothSimplex{n} \times \HomotopyQuotient{\SmoothManifold}{G} }
        { \HomotopyQuotient{ \mathbf{B}\Gamma }{G} }
      ^{\stable}
       \xhookrightarrow{\quad} \;
      \SlicePointsMaps{\big}{\mathbf{B}{\widehat{G}}}
        { \SmoothSimplex{n} \times \HomotopyQuotient{\SmoothManifold}{\widehat{G}} }
        { \HomotopyQuotient{ \mathbf{B}\Gamma }{\widehat{G}} }
      ^{\stable}
    \Big)
    \,.
  $$

  \vspace{-2mm}
  \noindent
  Now, since monomorphisms are characterized by a pullback property
  \eqref{InfinityMonomorphism}, it is sufficient to show that this
  homotopy pullback is preserved by the colimit operation.
  Prop. \ref{SufficientConditionsForTopologicalRealizationOfSimplicialSpacesToPreserveHomotopyPullbacks}
  says that a sufficient condition for this to be the case is that
  $$
    \Truncation{0}
    \,
      \SlicePointsMaps{\big}{\mathbf{B}{\widehat{G}}}
        { \SmoothSimplex{n} \times \HomotopyQuotient{\SmoothManifold}{\widehat{G}} }
        { \HomotopyQuotient{ \mathbf{B}\Gamma }{\widehat{G}} }
      ^{\stable}
    \;\simeq\;
    \Truncation{0}
    \,
      \SlicePointsMaps{\big}{\mathbf{B}{\widehat{G}}}
        { \HomotopyQuotient{\SmoothManifold}{\widehat{G}} }
        { \HomotopyQuotient{ \mathbf{B}\Gamma }{\widehat{G}} }
      ^{\stable}
  $$

  \vspace{-2mm}
  \noindent
  is independent of $n$. But, since
  \vspace{-2mm}
  $$
    \def\arraystretch{1.4}
    \begin{array}{lll}
      \shape \,\big( \SmoothSimplex{n} \times (-) \big)
      & \;\simeq\;
      \shape \, \SmoothSimplex{n} \,\times\, \shape \, (-)
      &
      \proofstep{ by \eqref{ShapePreservesBinaryProducts} }
      \\
      & \;\simeq\;
      \shape \, \mathbb{R}^n \,\times\, \shape \, (-)
      &
      \proofstep{ by Def. \ref{SmoothExtendedSimplicies} }
      \\
      & \;\simeq\;
      \shape \,  (-)
      &
      \proofstep{ by Prop. \ref{CohesiveShapeOfSmoothManifoldsIsTheirHomotopyType}, }
    \end{array}
  $$

  \vspace{-2mm}
\noindent  this follows by the classification result for $\widehat{G}$-equivariant
  $\Gamma$-principal bundles (Thm. \ref{BorelClassificationOfEquivariantBundlesForResolvableSingularitiesAndEquivariantStructure},
  using here the assumptions on $\widehat{G}$ and on $\Gamma$
  and the stability condition).
\end{proof}

\section{Equivariant moduli stacks}
\label{EquivariantModuliStacks}

The discussion in \cref{AsInfinityBundlesInternalToSliceOverBG},
culminating in Thm. \ref{BorelClassificationOfEquivariantBundlesForResolvableSingularitiesAndEquivariantStructure},
may be summarized as saying that
the stacky delooping
\vspace{-2mm}
$$
  \Bigg(
  \!\!\!\!\!
  \begin{tikzcd}[row sep=4pt]
    \mathbf{B}(\Gamma \rtimes G)
    \ar[d, start anchor={[yshift=+4pt]}]
    \\
    \mathbf{B}G
  \end{tikzcd}
  \!\!\!\!
  \Bigg)
  \;\;
  \simeq
  \;\;
  \Bigg(
  \!\!\!\!
  \begin{tikzcd}[row sep=4pt]
  \HomotopyQuotient
    { \mathbf{B}\Gamma }
    { G }
    \ar[d, start anchor={[yshift=+3pt]}]
    \\
    \HomotopyQuotient
      { \ast }
      { G }
  \end{tikzcd}
  \!\!\!
  \Bigg)
  \;\;\;
  \in
  \;
  \Slice{(\ModalTopos{\smooth})}{\mathbf{B}G}
$$

\vspace{-1mm}
\noindent
is the {\it moduli stack}
(as in Prop. \ref{DeloopingGroupoidsAreModuliInfinityStacksForPrincipalInfinityBundles})
for
(equivariantly locally trivial)
$G$-equivariant $\Gamma$-principal bundles.
Remarkably, this is always a Borel-equivariant object
(though with higher geometric structure).
But for these to serve as twists in a good twisted equivariant generalized cohomology theory (this is discussed in \cref{TwistedEquivariantCohomology} below) one needs to equivalently regard these moduli stacks as orbi-singularized objects in the slice $\SliceTopos{\orbisingularG}$ of a singular-cohesive $\infty$-topos (\cref{GeneralSingularCohesion}).

\medskip
This requires embedding the Borel-equivariant moduli
stacks into proper-equivariant homotopy theory,
which we achieve
(in Def. \ref{EquivariantClassifyingStack} below)
simply by forming their orbi-singularization:
\begin{equation}
  \label{OrbiSingularModuliStack}
  \orbisingular
  \Bigg(
  \!\!\!\!
  \begin{tikzcd}[row sep=4pt]
    \mathbf{B}(\Gamma \rtimes G)
    \ar[d, start anchor={[yshift=+4pt]}]
    \\
    \mathbf{B}G
  \end{tikzcd}
  \!\!\!\!
  \Bigg)
  \;\;
  \simeq
  \;\;
  \Bigg(
  \!\!\!\!\!
  \begin{tikzcd}[row sep=4pt]
  \orbisingular
  (
  \HomotopyQuotient
    { \mathbf{B}\Gamma }
    { G }
  )
    \ar[d, start anchor={[yshift=+3pt]}]
    \\
    \orbisingularG
  \end{tikzcd}
  \!\!\!\!
  \Bigg)
  \;\;\;
  \in
  \;
  \SliceTopos{\orbisingularG}
  \,.
\end{equation}

\vspace{-2mm}
\noindent
Notice that these proper equivariant moduli stacks
\eqref{OrbiSingularModuliStack}
still modulate equivariant principal bundles,
now regarded over cohesive orbi-spaces $\orbisingular (\HomotopyQuotient{X}{G})$
(recall Prop. \ref{OrbiSpaceIncarnationOfGSpaceIsOrbisingularizationOfHomotopyQuotient}
and \cite[\S 4]{SS20OrbifoldCohomology}),
in that:
\vspace{-2mm}
  $$
   \def\arraystretch{1.4}
    \begin{array}{lll}
      \SlicePointsMaps{\big}{\orbisingularG}
        { \orbisingular(X \!\sslash\! G) }
        {
          \orbisingular
          (
            \HomotopyQuotient{\mathbf{B}\Gamma}{G}
          )
        }
      & \;\simeq\;
      \SlicePointsMaps{\big}{\mathbf{B}G}
        { X \!\sslash\! G }
        { \HomotopyQuotient{\mathbf{B}\Gamma}{G} }
      &
      \proofstep{ by Ex. \ref{MorphsimsBetweenOrbisingularizationsOfSmoothObjects}}
      \\
      & \;\simeq\;
      \EquivariantPrincipalBundles{G}{\Gamma}
      (
        \ModalTopos{\smooth}
      )
      &
      \proofstep{ by Def. \ref{GEquivariantGammaPrincipalBundles}}
      \,.
    \end{array}
  $$

\vspace{-2mm}
The underlying equivariant homotopy type of
these orbi-singular
objects \eqref{OrbiSingularModuliStack}
is their pure shape \eqref{TheModalitiesOnACohesiveInfinityTopos}
hence the {\it equivariant shape}
of the
corresponding homotopy quotients:
\begin{equation}
  \label{EquivariantShape}
  \shape
  \,
  \orbisingular
  \Bigg(
  \!\!\!\!
  \begin{tikzcd}[row sep=4pt]
    \mathbf{B}(\Gamma \rtimes G)
    \ar[d, start anchor={[yshift=+4pt]}]
    \\
    \mathbf{B}G
  \end{tikzcd}
  \!\!\!\!
  \Bigg)
  \;\;
  \simeq
  \;\;
  \Bigg(
  \!\!\!
  \begin{tikzcd}[row sep=4pt]
  \shape
  \,
  \orbisingular
  (
  \HomotopyQuotient
    { \mathbf{B}\Gamma }
    { G }
  )
    \ar[d, start anchor={[yshift=+3pt]}]
    \\
    \orbisingularG
  \end{tikzcd}
  \!\!\!\!
  \Bigg)
  \;\;\;
  \in
  \;
  \SliceTopos{\orbisingularG}
  \,.
\end{equation}
We show (Thm. \ref{MurayamaShimakawaGroupoidIsEquivariantModuliStack} below)
that their $G$-(orbi-)spatial aspect
$\smoothrelativeG \;\, \shape\, \orbisingular (\HomotopyQuotient{\mathbf{B}\Gamma}{G})$
(Def. \ref{EquivariantClassifyingStack} below)
is given equivalently by the
Murayama-Shimakawa construction
from \cref{ConstructionOfUniversalEquivariantPrincipalBundles},
which serves to give that construction general abstract meaning.

\medskip
For example, this shows that the abstract reason why
equivariant classifying spaces are not connected
(see Rem. \ref{PseudoPrincipalNatureOfFixedLociInUniversalEquivariantPrincipalBundles})
even though
they arise as the images of moduli stacks which {\it are} connected
(according to Thm. \ref{DeloopingGroupoidsAreModuliInfinityStacksForPrincipalInfinityBundles}),
is that they arise so under orbi-singularization $\orbisingular$,
which is the (only) one of the three singular modalities \eqref{SingularModalities},
that does {\it not} preserve deloopings (Ex. \ref{OrbiSingularizationDoesNotPreserveDeloopings}).
This is the conceptual origin
(by
Prop. \ref{ConnectedComponentsOfHFixedEquivariantClassifyingSpaceWhenCrossedConjugationQuotientIsDiscrete}
with Thm. \ref{MurayamaShimakawaGroupoidIsEquivariantModuliStack}
and Thm. \ref{EquivariantHomotopyGroupsOfEquivariantClassifyingSpaces} below)
of the rich sets of connected components labeled by
non-abelian group cohomology classes
of  crossed homomorphisms
which is so characteristic for equivariant principal bundle theory
(as witnessed by the ubiquituous Ntn. \ref{LiftsOfEquivarianceSubgroupsToSemidirectProductWithStructureGroup}).

\medskip
Finally, we show
(Thm. \ref{ProperClassificationOfEquivariantBundlesForResolvableSingularitiesAndEquivariantStructure} below)
how the Borel-equivariant
classification of equivariant principal bundles from
Thm. \ref{BorelClassificationOfEquivariantBundlesForResolvableSingularitiesAndEquivariantStructure}
embeds into proper-equivariant theory
(Thm. \ref{ProperClassificationOfEquivariantBundlesForResolvableSingularitiesAndEquivariantStructure})
and we conclude with a list of examples and applications
of the classification result in this form
(Ex. \ref{TruncatedStructureGroupsAndTheirEquivariantClassificationResults},
 Ex. \ref{EquivariantBundlesServingAsGeoemtricTwistsOfEquivariantKTheory}).

\medskip

\noindent
{\bf Proper equivariant classifying spaces.}
In direct proper-equivariant generalization of
the notation
in Def. \ref{ClassifyingSpaceOfPrincipalInfinityBundles}, we set:

\begin{definition}[Proper equivariant classifying shape]
  \label{EquivariantClassifyingStack}
  For $\Topos$ a singular-cohesive $\infty$-topos (Def. \ref{SingularCohesiveInfinityTopos}),

  -- $G \in \Groups(\Topos)$,

  --
  $\Gamma \rtimes G
    \in
    \Groups
    (
      \SliceTopos{\mathbf{B}G}
    )
    $,

  \noindent
  {\bf (i)} we write
  \vspace{-1mm}
  \begin{equation}
    \label{EquivariantClassifyingSpaceAsShapeOfOrbisingularizationOfEquModuliStack}
    \EquivariantClassifyingShape{G}{\Gamma}
    \;\;\coloneqq\;\;
    \smoothrelativeG
    \;\,
    \shape
    \;
    \orbisingular
    \;
    \mathbf{B}(\Gamma \rtimes G)
    \;\;\;\;\;
    \in
    \;\;
    \GEquivariant\InfinityGroupoids
    \xhookrightarrow{\;\; \Discrete \, G\Orbi\Space \;\;}
    \SliceTopos{\orbisingularG}
  \end{equation}

  \vspace{-1mm}
  \noindent
  for the
  $G$-smooth aspect (Def. \ref{ModalitieswithrespecttoGOrbiSingularities})
  of the shape (Def. \ref{CohesiveInfinityTopos})
  of the orbi-singular aspect (Def. \ref{GEquivariantAndGloballyEquivariantHomotopyTheories})
  of the equivariant moduli stack
  (from Prop. \ref{DeloopingGroupoidsAreModuliInfinityStacksForPrincipalInfinityBundles},
    Thm. \ref{BorelClassificationOfEquivariantBundlesForResolvableSingularitiesAndEquivariantStructure}):
\vspace{-2mm}
$$
\hspace{-2mm}
  \begin{tikzcd}[column sep=15pt, row sep=0pt]
  \categorybox{\!  \Groups \!
    \left(
      \Slice{\!(\ModalTopos{\smooth})}{\mathbf{B}G}
    \!\right)
    \!}
    \ar[
      r,
      "{
        \overset
          {
            \raisebox{5pt}{
              \tiny
              \color{greenii}
              \bf
              delooping
            }
          }
          {\mathbf{B}}
      }"
    ]
    &
   \categorybox{    \Slice{(\ModalTopos{\smooth})}{\mathbf{B}G} }
    \ar[
      r,
      "{
        \overset
          {
            \raisebox{3pt}{
              \tiny
              \color{greenii}
              \bf
              orbi-singularization
            }
          }
          {\orbisingular}
      }"
    ]
    &
    \categorybox{   \Slice{\Topos}{\orbisingularG} }
    \ar[
      r,
      "{
        \overset
          {
            \raisebox{4pt}{
              \tiny
              \color{greenii}
              \bf
              shape
            }
          }
          {\scalebox{0.7}{$\Shape$}}
      }"
    ]
    &
   \categorybox{    \Slice{\GloballyEquivariant\InfinityGroupoids}{\orbisingularG} }
    \ar[
      rr,
      "{
        \overset
          {
            \raisebox{5pt}{
              \tiny
              \color{greenii}
              \bf
              $G$-orbi-spatial
            }
          }
          {\scalebox{0.7}{$G\Orbi\Smooth$}}
      }"
    ]
    &&
    \categorybox{   \GEquivariant\InfinityGroupoids }
    \ar[r, hook, "{\scalebox{0.7}{$\Discrete$}}"]
    &
   \categorybox{    \GEquivariant{\ModalTopos{\smooth}} }
    \ar[rr, hook, "{\scalebox{0.7}{$G\Orbi\Space$}}"]
    &&
    \categorybox{   \SliceTopos{\orbisingularG} }
    \\
\scalebox{0.7}{$
 \HomotopyQuotient{\Gamma}{G}
$}
    \ar[r, phantom, "\qquad \mapsto"]
    &
\scalebox{0.7}{$
    \mathbf{B}(\Gamma \rtimes G)
    $}
    \ar[r, phantom, "\quad \mapsto"]
    &
    \scalebox{0.7}{$
    \orbisingular \mathbf{B}(\Gamma \rtimes G)
    $}
    \ar[r, phantom, "\!\!\!\!\!\!\!\! \mapsto"]
    &
    \scalebox{0.7}{$
    \Shape \orbisingular \mathbf{B}(\Gamma \rtimes G)
    $}
    \ar[rr, phantom, "\longmapsto"]
    &&
    \scalebox{0.7}{$
    \EquivariantClassifyingShape{G}{\Gamma}
    $}
    \\[-3pt]
    \mathclap{
      \mbox{
        \tiny
        \color{darkblue}
        \bf
        \def\arraystretch{.9}
        \begin{tabular}{c}
          cohesive
          \\
          $G$-equivariant
          \\
          $\infty$-group
        \end{tabular}
      }
    }
    &
    \mathclap{
      \mbox{
        \tiny
        \color{darkblue}
        \bf
        \def\arraystretch{.9}
        \begin{tabular}{c}
          moduli stack of
          \\
          $G$-equivariant
          \\
          $\Gamma$-principal
          \\
          $\infty$-bundles
        \end{tabular}
      }
    }
    &
    \mathclap{
      \mbox{
        \tiny
        \color{darkblue}
        \bf
        \def\arraystretch{.9}
        \begin{tabular}{c}
          proper-equivariant
          \\
          moduli stack
        \end{tabular}
      }
    }
    &
    \mathclap{
      \mbox{
        \tiny
        \color{darkblue}
        \bf
        \def\arraystretch{.9}
        \begin{tabular}{c}
          global-equivariant
          \\
          homotopy type of
          \\
          proper-equivariant
          \\
          moduli stack
        \end{tabular}
      }
    }
    &&
    \mathclap{
      \mbox{
        \tiny
        \color{orangeii}
        \bf
        \def\arraystretch{.9}
        \begin{tabular}{c}
          $G$-orbi-space
          \\
          underlying
          \\
          proper-equivariant
          \\
          moduli stack
        \end{tabular}
      }
    }
  \end{tikzcd}
$$

\vspace{-3mm}
\noindent {\bf (ii)} For
$G \acts \, X \,\in\, \Actions{G}(\ModalTopos{\smooth})$
with
$X \,\in\, \ModalTopos{\smooth,0} \xhookrightarrow{\;} \Topos$,
we denote by
\vspace{-3mm}
\begin{equation}
\label{TheMapChi}
  \chi
  \;:\;
  \begin{tikzcd}
  \IsomorphismClasses
  {
    \EquivariantPrincipalBundles{G}{\Gamma}(\ModalTopos{\smooth})_{X}
  }
  \ar[r]
  &
  H^0_{\scalebox{.7}{$\orbisingularG$}}
  \left(
    X
    ;\,
    \EquivariantClassifyingShape{G}{\Gamma}
  \right)
  \end{tikzcd}
\end{equation}

\vspace{-3mm}
\noindent
the natural map from
equivalence classes of
$G$-equivariant $\Gamma$-principal $\infty$-bundles over $X$
to the
non-abelian
proper-equivariant cohomology of the $G$-orbi-space underlying $G \acts \, X$
with coefficients in $\EquivariantClassifyingShape{G}{\Gamma}$
\eqref{EquivariantClassifyingSpaceAsShapeOfOrbisingularizationOfEquModuliStack}.
\end{definition}

\begin{remark}[The equivariant characteristic class map]
The map  $\chi$ \eqref{TheMapChi} is the image under 0-truncation of the following composite:
\vspace{-2mm}
$$
  \begin{array}{lllll}
    &
    \EquivariantPrincipalBundles{G}{\Gamma}(\ModalTopos{\smooth})_{X}
    \\
    &
    \;\simeq\;
    \SlicePointsMaps{\big}{\mathbf{B}G}
      { \HomotopyQuotient{X}{G} }
      { \HomotopyQuotient{\mathbf{B}\Gamma}{G} }
    &
    \;\xrightarrow{ \shape \, \orbisingular }\;
    &
    \SlicePointsMaps{\big}{\orbisingularG}
      { \shape\, \orbisingular(\HomotopyQuotient{X}{G}) }
      { \shape\, \orbisingular(\HomotopyQuotient{\mathbf{B}\Gamma}{G}) }
    \\
    &&
    &
    \;\simeq\;
    \SlicePointsMaps{\Big}{\orbisingularG}
      {
        G\Orbi\Space
        \big(
          \Shape \, \FixedLoci(X)
        \big)
      }
      {
        \shape\,
        \orbisingular(\HomotopyQuotient{\mathbf{B}\Gamma}{G})
      }
    &
    \proofstep{ by Prop. \ref{OrbiSpaceIncarnationOfGSpaceIsOrbisingularizationOfHomotopyQuotient} }
    \\
    &&
    &
    \;\simeq\;
    \SlicePointsMaps{\Big}{\orbisingularG}
      {
        \conicalrelativeG
        \,
        G\Orbi\Space
        \big(
          \Shape \, \FixedLoci(X)
        \big)
      }
      {
        \shape\,
        \orbisingular(\HomotopyQuotient{\mathbf{B}\Gamma}{G})
      }
    &
    \proofstep{ by Ex. \ref{IdempotencyOfGOrbiSingularModalities} }
    \\
    &&
    &
    \;\simeq\;
    \SlicePointsMaps{\Big}{\orbisingularG}
      {
        G\Orbi\Space
        \big(
          \Shape \, \FixedLoci(X)
        \big)
      }
      {
        \smoothrelativeG
        \;
        \shape
        \,
        \orbisingular
        (\HomotopyQuotient{\mathbf{B}\Gamma}{G})
      }
    &
    \proofstep{ by \eqref{GSingularModalities} }
    \\
    &&
    &
    \;\simeq\;
    \GEquivariant\InfinityGroupoids
    \big(
      \Shape\, \FixedLoci(X)
      ,\,
      \EquivariantClassifyingShape{G}{\Gamma}
    \big)
    &
    \proofstep{ by \eqref{GFixedLocFunctor}. }
  \end{array}
$$
\end{remark}

We proceed to compute
(in Thm. \ref{MurayamaShimakawaGroupoidIsEquivariantModuliStack} below)
the geometric fixed loci of the proper-equivariant classifying spaces
$\EquivariantClassifyingShape{G}{\Gamma}$
\eqref{EquivariantClassifyingSpaceAsShapeOfOrbisingularizationOfEquModuliStack}.

\begin{lemma}[Geometric fixed loci of equivariant classifying spaces]
  \label{FixedLociOfEquivariantClassifyingStacks}
  In a singular-cohesive $\infty$-topos $\Topos$
  (Def. \ref{SingularCohesiveInfinityTopos}),
  given

  -- $G \,\in\, \Groups(\Sets)$,

  -- $G \acts \, \Gamma \,\in\, \Groups\left( \Actions{G}(\Topos) \right)$,

  \noindent
  {\bf (i)}
  the geometric $H$-fixed loci (Def. \ref{ShapeOfGemetricFixedLoci})
  of the equivariant moduli stack
  \eqref{OrbiSingularModuliStack}
  are
  \vspace{-1mm}
\begin{equation}
  \label{HFixedLociInEquivariantModuliStack}
  \smooth
  \SliceMaps{\big}{\orbisingularG}
    { \orbisingularH }
    { \orbisingular (\HomotopyQuotient{\mathbf{B}\Gamma}{G}) }
  \;\;
  \simeq
  \;\;
  \SliceMaps{\big}{\mathbf{B}G}
    { \mathbf{B}H }
    { \HomotopyQuotient{\mathbf{B}\Gamma}{G} };
\end{equation}

\vspace{-1mm}
\noindent
{\bf (ii)} their shape is, equivalently, the $H$-fixed locus in the
equivariant classifying shape \eqref{EquivariantClassifyingSpaceAsShapeOfOrbisingularizationOfEquModuliStack}:
\vspace{0mm}
\begin{equation}
  \label{HFixedLociInEquivariantClassifyingSpaceForInfinityGroups}
  \mathllap{
    \EquivariantClassifyingShape{G}{\Gamma}
    \;:\;
    G/H
    \;\mapsto\;
  }
  \smooth
  \SliceMaps{\big}{\orbisingularG}
    { \orbisingularH }
    {
      \shape
      \,
      \orbisingular
      \,
      \HomotopyQuotient{\mathbf{B}\Gamma}{G}
    }
  \;\;
  \simeq
  \;\;
  \shape
  \;
  \SliceMaps{\big}{\mathbf{B}G}
    { \mathbf{B}H }
    { \HomotopyQuotient{\mathbf{B}\Gamma}{G} }
  \,.
\end{equation}
\end{lemma}
\begin{proof}
For $U \times \orbisingularG \,\in\, \Charts \times \Singularities$,
we have the following sequence of natural equivalences:
\vspace{-2mm}
$$
\hspace{-3mm}
\def\arraystretch{1.9}
  \begin{array}{lll}
        \bigg(\!\!\!
    \smooth
    \Big(\!
      \SliceMaps{\big}{\orbisingularG}
        { \orbisingularH }
        {
          \orbisingular
          (\HomotopyQuotient{\mathbf{B}\Gamma}{G})
        }
    \Big)
    \!\!\bigg)
    & \!\!\!\!\!\!\!\!\!
    \left(U \times \orbisingularK \right)
    \\
    &
   \!\! \simeq\;
    \Big(
      \SliceMaps{\big}{\orbisingularG}
        { \orbisingularH }
        {
          \orbisingular
          (\HomotopyQuotient{\mathbf{B}\Gamma}{G})
        }
    \Big)
    (U)
    &
    \proofstep{ by Def. \ref{GEquivariantAndGloballyEquivariantHomotopyTheories}}
    \\
    &\!\! \simeq\;
    \PointsMaps{\big}
      { \orbisingularH \times U }
      {
        \orbisingular
        \HomotopyQuotient{\mathbf{B}\Gamma}{G}
      }
      \underset{
        \scalebox{.7}{$
          \PointsMaps{\big}
            { \orbisingularH \times U }
            { \orbisingularG  }
        $}
      }{\times}
      \{\ast\}
    &
    \proofstep{
      by Prop. \ref{HomSpaceInSliceToposAsFiberProduct}, Lem. \ref{PlotsOfSliceMappingStackAreSliceHoms}
    }
    \\
    & \!\!\simeq\;
    \PointsMaps{\big}
      { \orbisingularH \times \orbisingular U }
      {
        \orbisingular
        (\HomotopyQuotient{\mathbf{B}\Gamma}{G})
      }
      \underset{
        \scalebox{.7}{$
          \PointsMaps{}
            { \orbisingularH \times \orbisingular U }
            { \orbisingularG }
        $}
      }{\times}
    \{\ast\}
    &
    \proofstep{ by \eqref{SmoothChartsAreSmoothAndOrbisingular} }
    \\
    & \!\!\simeq\;
    \PointsMaps{\big}
      { \orbisingular(\mathbf{B}H \times U) }
      {
        \orbisingular
        (\HomotopyQuotient{\mathbf{B}\Gamma}{G})
      }
      \underset{
        \scalebox{.7}{$
          \PointsMaps{\big}
            { \orbisingular(\mathbf{B}H \times U) }
            { \orbisingular\mathbf{B}G }
        $}
      }{\times}
    \{\ast\}
    &
    \proofstep{ by \eqref{InfinityAdjointPreservesInfinityLimits} \& Ex. \ref{OrbiSingularityIsOrbiSingularizationOfHomotopyQuotient} }
    \\
    & \!\!\simeq\;
    \PointsMaps{\big}
      { \mathbf{B}H \times U }
      {
        \orbisingular
        (\HomotopyQuotient{\mathbf{B}\Gamma}{G})
      }
      \underset{
        \scalebox{.7}{$
          \PointsMaps{}
            { \mathbf{B}H \times U }
            { \mathbf{B}G }
        $}
      }{\times}
    \{\ast\}
    &
    \proofstep{ by Ex. \ref{MorphsimsBetweenOrbisingularizationsOfSmoothObjects} }
    \\
    &\!\! \simeq\;
    \PointsMaps{\big}
      { \mathbf{B}H \times U \times \orbisingularK }
      {
        \orbisingular
        (\HomotopyQuotient{\mathbf{B}\Gamma}{G})
      }
      \underset{
        \scalebox{.7}{$
          \PointsMaps{\big}
            { \mathbf{B}H \times U \times \orbisingularK }
            { \mathbf{B}G }
        $}
      }{\times}
    \{\ast\}
    &
    \proofstep{ by Ex. \ref{CohesiveLociOfOrbiSingularities} }
    \\
    & \!\!\simeq\;
    \Big(
      \SliceMaps{\big}{\mathbf{B}G}
        { \mathbf{B}H }
        {
          \orbisingular
          (\HomotopyQuotient{\mathbf{B}\Gamma}{G})
        }
    \Big)(U \times \orbisingularK)
    &
    \proofstep{
      by Props. \ref{HomSpaceInSliceToposAsFiberProduct}, \ref{PlotsOfSliceMappingStackAreSliceHoms}.
    }
   \end{array}
$$

\vspace{-2mm}
\noindent Since this is natural
in
$U \times \orbisingularK \,\in\, \CartesianSpaces \times \Singularities
\xhookrightarrow{y} \Topos$,
the first statement \eqref{HFixedLociInEquivariantModuliStack} follows by the
$\infty$-Yoneda lemma (Prop. \ref{InfinityYonedaLemma}).
With this, the second statement \eqref{HFixedLociInEquivariantClassifyingSpaceForInfinityGroups}
is implied as follows:
$$
\def\arraystretch{1.7}
  \begin{array}{lll}
    \smooth
    \Big(
      \SliceMaps{\big}{\orbisingular \mathbf{B}G}
      {
        \orbisingularH
      }
      {
        \shape
        \,
        \orbisingular
        (\HomotopyQuotient{\mathbf{B}\Gamma}{G})
      }
    \Big)
    & \;\simeq\;
    \smooth
    \,
    \shape
    \,
    \Big(
      \SliceMaps{\big}{\orbisingular \mathbf{B}G}
      {
        \orbisingularH
      }
      {
        \orbisingular
        (\HomotopyQuotient{\mathbf{B}\Gamma}{G})
      }
    \Big)
    &
    \proofstep{ by Prop. \ref{ShapeOfMappingStackOutOfOrbiSingularityIsMappingStackIntoShape}}
    \\
    & \;\simeq\;
    \shape
    \,
    \smooth
    \,
    \Big(
      \SliceMaps{\big}{\orbisingular \mathbf{B}G}
      {
        \orbisingularH
      }
      {
        \orbisingular \HomotopyQuotient{\mathbf{B}\Gamma}{G}
      }
    \Big)
    &
    \proofstep{ by Prop. \ref{SomeSingularModalitiesCommuteWithSomeCohesiveModalities} }
    \\
    & \;\simeq\;
    \shape
    \,
      \SliceMaps{\big}{\mathbf{B}G}
      { \mathbf{B}H }
      { \HomotopyQuotient{\mathbf{B}\Gamma}{G} }
    &
    \proofstep{ by \eqref{HFixedLociInEquivariantModuliStack}}.
  \end{array}
$$

\vspace{-6mm}
\end{proof}

\begin{example}[The case of trivial $G$-action on the structure $\infty$-group]
  \label{FixedLociInEquivariantModuliStackForTrivialGActionOnStructureGroup}
  In the special case when the action of $G$ on $\Gamma$ is trivial,
  in that
  $\HomotopyQuotient{\mathbf{B}\Gamma}{G} \,\simeq\, (\mathbf{B}\Gamma) \times (\mathbf{B}G)$
  (Ex. \ref{EquivariantInfinityGroupForTrivialGAction}),
  we have
  (by Ex. \ref{SliceMappingStackIntoProductProjection})
  $$
    \SliceMaps{\big}{\mathbf{B}G}
      { \mathbf{B}H }
      { \HomotopyQuotient{\mathbf{B}\Gamma}{G} }
    \;\;
      \simeq
    \;\;
    \Maps{}
      {\mathbf{B}H}
      {\mathbf{B}\Gamma}
    \,.
  $$
  Consequently, in the case
  the statement of Lem. \ref{FixedLociOfEquivariantClassifyingStacks}
  reduces to the following two equivalences
 \begin{align}
  \smooth
  \SliceMaps{\big}{\orbisingularG}
    { \orbisingularH }
    {
      \orbisingular
      (\HomotopyQuotient{\mathbf{B}\Gamma}{G})
    }
  \;\;
  &
  \simeq
  \;\;
  \Maps{}
    { \mathbf{B}H }
    { \mathbf{B}\Gamma }\;,
  \\
  \mathllap{
    \EquivariantClassifyingShape{G}{\Gamma}
    \;:\;
    G/H
    \;\mapsto\;
  }
  \smooth
  \SliceMaps{\big}{\orbisingularG}
    { \orbisingularH }
    {
      \shape
      \,
      \orbisingular
      (\HomotopyQuotient{\mathbf{B}\Gamma}{G})
    }
  \;\;
  & \simeq
  \;\;
  \shape
  \,
  \Maps{}
    { \mathbf{B}H }
    { \mathbf{B}\Gamma }
  \,.
\end{align}
\end{example}
In fact, in this case (Ex. \ref{FixedLociInEquivariantModuliStackForTrivialGActionOnStructureGroup})
the equivariant classifying shapes have a ``globally equivariant''
incarnation (in the terminology of \cite{Schwede18}), in proper-equivariant refinement of Rem.
\ref{GloballyEquivariantNatureOfEquivariantPrincipalBundles}
and along the lines of \cite[\S 1.3-1.4]{Rezk14}:
\begin{proposition}[Globally equivariant classifying shape]
 \label{GloballyProperEquivariantClassifyingShape}
If the action of the equivariance group $G$ on the structure group $\Gamma$ is trivial,
so that their semi-direct product is their direct product
\begin{equation}
  \label{TrivialGActionOnGammaInDiscussionOfGlobalEquivariance}
  \Gamma \rtimes G \,\simeq\, \Gamma \times G\;,
\end{equation}
then the equivariant shape \eqref{EquivariantShape}
of the delooping of the structure group
\begin{equation}
  \label{GlobalEquivariantClassifyingShape}
  \shape
  \,
  \orbisingular
  \,
  \mathbf{B}\Gamma
  \;\;\;
  \in
  \;
  \Topos
\end{equation}
is the \emph{globally equivariant} classifying shape,
in that, for all $G \,\in\, \Groups(\Sets)$
and $G \acts \, X \,\in\, \Actions{G}(\ModalTopos{\smooth,0})$,
we have a natural equivalence
$$
  \PointsMaps{\big}
    { \orbisingular(\HomotopyQuotient{X}{G}) }
    { \shape \, \orbisingular \, \mathbf{B}\Gamma }
  \;\;
  \simeq
  \;\;
  \SlicePointsMaps{\big}{\orbisingularG}
    { \orbisingular(\HomotopyQuotient{X}{G}) }
    { \EquivariantClassifyingShape{G}{\Gamma} }
$$
between the global hom-$\infty$-groupoid into
\eqref{GlobalEquivariantClassifyingShape}
and the slice-hom-$\infty$-groupoid into \eqref{EquivariantClassifyingSpaceAsShapeOfOrbisingularizationOfEquModuliStack}.
\end{proposition}
\begin{proof}
This is the composite of the following sequence of natural equivalences,
$$
  \def\arraystretch{1.4}
  \begin{array}{lll}
    \SlicePointsMaps{\big}{\orbisingularG}
      { \orbisingular( \HomotopyQuotient{X}{G} ) }
      { \EquivariantClassifyingShape{G}{\Gamma} }
    &
    \;=\;
    \SlicePointsMaps{\big}{\orbisingularG}
      { \orbisingular( \HomotopyQuotient{X}{G} ) }
      {
        \smoothrelativeG
        \,
        \shape
        \,
        \orbisingular
        \,
        \mathbf{B}(\Gamma \rtimes G)
      }
    &
    \proofstep{ by Def. \ref{EquivariantClassifyingStack} }
    \\
    & \;\simeq\;
    \SlicePointsMaps{\big}{\orbisingularG}
      { \conicalrelativeG\, \orbisingular( \HomotopyQuotient{X}{G} ) }
      {
        \shape
        \,
        \orbisingular
        \,
        \mathbf{B}(\Gamma \rtimes G)
      }
    &
    \proofstep{ by \eqref{GSingularModalities}, \eqref{AdjunctionAndHomEquivalence} }
    \\
    & \;\simeq\;
    \SlicePointsMaps{\big}{\orbisingularG}
      { \orbisingular( \HomotopyQuotient{X}{G} ) }
      {
        \shape
        \,
        \orbisingular
        \,
        \mathbf{B}(\Gamma \rtimes G)
      }
    &
    \proofstep{ by Prop. \ref{OrbiSpaceIncarnationOfGSpaceIsOrbisingularizationOfHomotopyQuotient}
    with Ex. \ref{IdempotencyOfGOrbiSingularModalities} }
    \\
    & \;\simeq\;
    \SlicePointsMaps{\big}{\orbisingularG}
      { \orbisingular( \HomotopyQuotient{X}{G} ) }
      {
        \shape
        \,
        \orbisingular
        \,
        \mathbf{B}(\Gamma \times G)
      }
    &
    \proofstep{ by \eqref{TrivialGActionOnGammaInDiscussionOfGlobalEquivariance} }
    \\
    & \;\simeq\;
    \SlicePointsMaps{\big}{\orbisingularG}
      { \orbisingular( \HomotopyQuotient{X}{G} ) }
      {
        \shape
        \,
        \orbisingular
        \,
        \big(
          (\mathbf{B}\Gamma)
          \times
          (\mathbf{B}G)
        \big)
      }
    &
   \proofstep{ by Ex. \ref{EquivariantInfinityGroupForTrivialGAction} or \eqref{DeloopingPreservesProducts} }
    \\
    & \;\simeq\;
    \SlicePointsMaps{\big}{\orbisingularG}
      { \orbisingular( \HomotopyQuotient{X}{G} ) }
      {
        (
        \shape
        \,
        \orbisingular
        \,
        \mathbf{B}\Gamma
        )
        \times
        (
        \shape
        \,
        \orbisingular
        \,
        \mathbf{B}G
        )
      }
    &
    \proofstep{ by \eqref{InfinityAdjointPreservesInfinityLimits} and \eqref{ShapePreservesBinaryProducts} }
     \\
    & \;\simeq\;
    \SlicePointsMaps{\big}{\orbisingularG}
      { \orbisingular( \HomotopyQuotient{X}{G} ) }
      {
        (
        \shape
        \,
        \orbisingular
        \,
        \mathbf{B}\Gamma
        )
        \times
        \orbisingularG
      }
    &
    \proofstep{ by \eqref{OrbiSingularityIsOrbiSingularizationOfDelooping}}
    \\
    & \;\simeq\;
    \PointsMaps{\big}
      { \orbisingular( \HomotopyQuotient{X}{G} ) }
      {
        \shape
        \,
        \orbisingular
        \,
        \mathbf{B}\Gamma
      }
    &
    \proofstep{ by Ex. \ref{BaseChangeToAbsoluteContext}. }
  \end{array}
$$

\vspace{-6mm}
\end{proof}

\begin{definition}[Stable proper equivariant classifying shape]
  \label{StableEquivariantClassifyingShape}
  Given $G \,\in\, \Groups(\FiniteSets)_{\resolvable}$ (Ntn. \ref{ResolvableOrbiSingularities}),
  and
  $G \acts \, \Gamma \,\in\, \Groups(\SmoothInfinityGroupoids)$
  of truncated classifying shape (Ntn. \ref{CohesiveGroupsWithTruncatedClassifyingShape})
  and with a notion of stable equivariant bundles (Ntn. \ref{StableEquivariantBundles}),
  we write
  $$
    \begin{tikzcd}[row sep=4pt]
      &[-20pt]&[-20pt]
      (
      \EquivariantClassifyingShape
        {G}{\Gamma}
      )
      ^{\stable}
      \ar[
        rr,
        hook,
      ]
      &&
      \EquivariantClassifyingShape
        {G}{\Gamma}
      \\
      G/H
      &\mapsto&
      \shape
      \SliceMaps{}{\mathbf{B}G}
        { \mathbf{B}H }
        { \HomotopyQuotient{ \mathbf{B}\Gamma }{G} }
      ^{\stable}
      \ar[
        rr,
        hook,
      ]
      &&
      \SliceMaps{}{\mathbf{B}G}
        { \mathbf{B}H }
        { \HomotopyQuotient{ \mathbf{B}\Gamma }{G} }
    \end{tikzcd}
  $$
  for the subobject of the proper-equivariant classifying shape
  (Def. \ref{EquivariantClassifyingStack})
  which is given on its values \eqref{HFixedLociInEquivariantClassifyingSpaceForInfinityGroups}
  by the stable component in the sense of Ntn. \ref{StableEquivariantBundles}.
\end{definition}

\begin{theorem}[Equivariant homotopy groups of equivariant classifying shapes]
\label{EquivariantHomotopyGroupsOfEquivariantClassifyingSpaces}
Given

  \vspace{1mm}
  -- $G \,\in\, \Groups(\FiniteSets)_{\resolvable}$ (Ntn. \ref{ResolvableOrbiSingularities}),

  \vspace{1mm}
  --
  $G \acts \Gamma \,\in\, \Actions{G}\big(\Groups(\kTopologicalSpaces)\big)$
  of truncated classifying shape (Ntn. \ref{CohesiveGroupsWithTruncatedClassifyingShape})
  and with a notion of
  of blowup-stable equivariant $\Gamma \rtimes G$-principal bundles
  (Ntn. \ref{StableEquivariantBundles});

  \vspace{1mm}
  --
  such that
  $\pi_1 \Maps{\big} {B G}{ B (\Gamma \rtimes G) }$
  is countable (Rem. \ref{CardinalityBoundOnHomotopyClassesOfMapsFromXModGToGamma});

\vspace{1mm}
\noindent
then the equivariant homotopy groups of the
blowup-stable component of the
equivariant classifying shape
$\EquivariantClassifyingShape{G}{\Gamma}$ \eqref{EquivariantClassifyingStack}
at stage $H \subset G$
are given by the non-abelian group cohomology
of $H$
with coefficients in $H \acts \shape \, \Gamma$, in that we have natural isomorphisms as follows:
\vspace{0mm}
$$
  \underset{
    { H \,\subset\, G }
    \atop
    { n \,\in\, \mathbb{N} }
  }{\forall}
  \;\;\;\;\;\;\;\;
  \pi_n^H
  \big(
    (\EquivariantClassifyingShape{G}{\Gamma})^{\stable}
  \big)
  \;\;
  \simeq
  \;\;
  H^{1-n}_{\mathrm{Grp}}
  (
    H
    ;\,
    \shape \, \Gamma
  )
  \;=\;
  H^1_{\mathrm{Grp}}
  (
    H
    ;\,
    \Omega^{n} \shape \, \Gamma
  )
  \,.
$$
\end{theorem}
\begin{proof}
$\,$

\vspace{-8mm}
$$
  \def\arraystretch{1.5}
  \begin{array}{lll}
    \pi_n^H
    \big(
      (\EquivariantClassifyingShape{G}{\Gamma})^{\stable}
    \big)
    &
    \;=\;
    \pi_n
    \Big(
      (\EquivariantClassifyingShape{G}{\Gamma})^{\stable}(G/H)
    \Big)
    \\
    &
    \;\simeq\;
    \pi_n
    \big(
      \shape
      \,
      \SliceMaps{}{B H}
        { \mathbf{B}H }
        { \HomotopyQuotient{\mathbf{B}\Gamma}{H} }
      ^{\stable}
    \big)
    &
    \proofstep{
      by Lem. \ref{FixedLociOfEquivariantClassifyingStacks}
    }
    \\
    &
    \;\simeq\;
    \pi_n
    \big(
      \SliceMaps{}{\mathbf{B}H}
        { \shape \, \mathbf{B}H }
        { \shape \, \mathbf{B}\Gamma }
    \big)
    &
    \proofstep{
      by Thm. \ref{OrbiSmoothOkaPrinciple}
    }
    \\
    & \;\simeq\;
    \pi_n
    \big(
      \SliceMaps{}{B H}
        { B H }
        {
          \HomotopyQuotient
            { B \shape \Gamma }
            { G }
        }
    \big)
    &
    \proofstep{
      by Def. \ref{ClassifyingSpaceOfPrincipalInfinityBundles},
      Ex. \ref{ClassifyingShapesForDiscreteStructureInfinityGroups}
    }
    \\
    & \;\simeq\;
    \left\{
    \begin{array}{ll}
      \Truncation{0}
      \,
      \SliceMaps{}{B H}
        { B H }
        { B \,\shape\, \Gamma }
      &
      \vert\; n = 0
      \\
      \Truncation{0}
      \,
      \Maps{}
        { B H }
        { B \Omega^n \,\shape\, \Gamma }
      &
      \vert\; n > 0
    \end{array}
    \right.
    &
    \proofstep{ by Lem. \ref{LoopingOfMappingStackIsMappingStackIntoLooping} }
    \\
    & \;=\;
    H^1_{\mathrm{Grp}}( H;\, \Omega^n \,\shape\, \Gamma )
    &
    \proofstep{
      e.g. \cite[Ex. 2.4]{FSS20CharacterMap}
    }
    \,.
  \end{array}
$$

\vspace{-6mm}
\end{proof}

\begin{remark}[Proper equivariant classifying shapes may depend on structure group beyond its shape]
  \label{ProperEquivariantClassifyingShapesMayDependOnStructureGroupBeyondItsShape}
  While plain classifying shapes (Def. \ref{ClassifyingSpaceOfPrincipalInfinityBundles})
  are sensitive only to the shape of the structure group, in that
  $B \Gamma \,\simeq\, B \, \shape \, \Gamma $
  (Ex. \ref{ClassifyingShapesForDiscreteStructureInfinityGroups}),
  for proper equivariant classifying shapes (Def. \ref{EquivariantClassifyingStack})
  this generally fails, in as far
  as the geometric fixed points of $\mathbf{B}\Gamma$
  differ from that of its shape (the ``homotopy fixed points''):
  \vspace{-1mm}
  $$
    \EquivariantClassifyingShape{G}{\Gamma}
    \;\;\;
    \underset{
      \mbox{\tiny i.g.}
    }{\not\simeq}
    \;\;\;
    \EquivariantClassifyingShape{G}{(\shape \, \Gamma)}
    \,.
  $$

    \vspace{-2mm}
\noindent  Concretely, Lem. \ref{FixedLociOfEquivariantClassifyingStacks}
  shows that the two are related by taking the shape operation into
  the mapping stack along the comparison morphism
  \eqref{ComparisonMorphismFromShapeOfMappingStackToMappingSpaceOfShapes}:
  \vspace{-3mm}
  \begin{equation}
    \label{ComparisonMorphismBetweenEquivriantClassifyingShapeForGammaAndForShapeOfGamma}
    \begin{tikzcd}[row sep=3pt]
      &
      \EquivariantClassifyingShape{G}{\Gamma}
      \ar[rr]
      &&
      \EquivariantClassifyingShape{G}{(\shape \, \Gamma)}
      \\
      G/H
      \ar[r, phantom, "\mapsto"]
      &
      \shape
      \,
      \SliceMaps{}{\mathbf{B}G}
        { \mathbf{B}H }
        { \HomotopyQuotient{\mathbf{B}\Gamma}{G} }
      \ar[
        rr,
        "{
          \scalebox{.7}{$
            \widetilde {\shape \mathrm{ev}}
          $}
        }"
      ]
      &&
      \SliceMaps{}{\mathbf{B}G}
        { \mathbf{B}H }
        { \HomotopyQuotient{B \Gamma}{G} }\;.
    \end{tikzcd}
  \end{equation}
\end{remark}

In general, there is no reason
for \eqref{ComparisonMorphismBetweenEquivriantClassifyingShapeForGammaAndForShapeOfGamma}
to be an equivalence. However, it is so under suitable conditions
where a smooth Oka principle applies
(such as in Prop. \ref{OrbiSmoothOkaPrinciple}):

\begin{proposition}[Stable equivariant classifying shape of truncated groups coincides with that
of their shape]
  \label{EquivariantClassifyingShapeOfTruncatedTopologicalGroupsCoincidesWithThatOftheirShape}
  Given

  -- $G \,\in\, \Groups(\FiniteSets)_{\resolvable}$ (Ntn. \ref{ResolvableOrbiSingularities}),

  \vspace{1mm}
  -- $\Gamma \,\in\, \Groups(\SmoothSets)$ of truncated classifying
  shape (Ntn. \ref{CohesiveGroupsWithTruncatedClassifyingShape}),

\vspace{1mm}
\noindent
then the stable proper $G$-equivariant classifying shape of $\Gamma$
(Def. \ref{StableEquivariantClassifyingShape})
is equivalent to the plain
proper equivariant classifying shape
(Def. \ref{EquivariantClassifyingStack})
of $\shape \, \Gamma$:
\vspace{-2mm}
$$
  \EquivariantClassifyingShape{G}{(\Gamma)}^{\stable}
  \;\;
  \simeq
  \;\;
  \EquivariantClassifyingShape{G}{(\shape \, \Gamma)}
  \;\;\;\;
  \in
  \;
  \GEquivariant\InfinityGroupoids
  \,.
$$
\end{proposition}
\begin{proof}
For $H \subset G$ a subgroup, we have the following sequence
of natural equivalences in
$\InfinityGroupoids \xhookrightarrow{\Discrete}  \ModalTopos{\smooth}$:
  \vspace{-2mm}
$$
  \def\arraystretch{1.3}
  \begin{array}{lll}
    \left(
      \EquivariantClassifyingShape{G}{\Gamma}
    \right)^{\stable}(G/H)
    &
    \;\simeq\;
    \shape
    \,
    \SliceMaps{}{\mathbf{B}G}
      { \mathbf{B}H }
      { \HomotopyQuotient{ \mathbf{B}\Gamma }{ G } }
    ^{\stable}
    &
    \proofstep{
      by Def. \ref{StableEquivariantClassifyingShape}
    }
    \\
    &
    \;\simeq\;
    \SliceMaps{}{\mathbf{B}G}
      { \shape \,\mathbf{B}H }
      { \HomotopyQuotient{ \shape \, \mathbf{B}\Gamma }{ G } }
    &
    \proofstep{ by Thm. \ref{OrbiSmoothOkaPrinciple} }
    \\
    &
    \;\simeq\;
    \SliceMaps{}{B G}
      { \mathbf{B} H }
      { \HomotopyQuotient{ \mathbf{B}\, \shape\, \Gamma }{ G } }
    &
    \proofstep{ by \eqref{ShapePreservesDeloopings} }
    \\
    &
    \;\simeq\;
    \shape
    \,
    \SliceMaps{}{B G}
      { \mathbf{B} H }
      { \HomotopyQuotient{ \mathbf{B} \, \shape \, \Gamma }{ G } }
    &
    \proofstep{ by \eqref{EquivalenceExhibitingMappingStackIntoDiscreteObjectAsDiscrete} }
    \\
    &
    \;\simeq\;
    \Big(
      \EquivariantClassifyingShape{G}{(\shape \, \Gamma)}
    \Big)(G/H)
    &
    \proofstep{ by Lem. \ref{FixedLociOfEquivariantClassifyingStacks}. }
  \end{array}
$$

\vspace{-1mm}
\noindent
By naturality in $G/H \,\in\, \OrbitCategory{G}$ (Ntn. \ref{GOrbitCategory}),
this
 establishes the
 claimed equivalence of
$\infty$-presheaves.
\end{proof}

\medskip

\begin{lemma}[Base change of equivariant classifying shapes along a covering of the equivariance group]
  \label{BaseChangeOfEquivariantClassifyingShapesAlongCoverings}
  For $G \acts \, \Gamma \,\in\, \Actions{G}\left( \Groups(\ModalTopos{0}) \right)$
  and  a surjective homomorphism of discrete groups
  $p :\!\! \begin{tikzcd} \widehat{G} \ar[r,->>] &[-12pt] G  \end{tikzcd}\!\!$
  \eqref{ADiscreteGroupEpimorphismInAnInfinityTopos}, there
  is a natural monomorphism \eqref{InfinityMonomorphism}
  \vspace{-2mm}
  $$
    \EquivariantClassifyingShape{G}{\Gamma}
    \xhookrightarrow{ \; (B p)^\ast \; }
    p_\ast
    \EquivariantClassifyingShape{\widehat{G}}{\Gamma}
    \;\;\;\;
    \in
    \;
    G\InfinityGroupoids
  $$

  \vspace{0mm}
\noindent
  from the equivariant classifying shape (Def. \ref{EquivariantClassifyingStack})
  of $G \acts \, \Gamma$ to the direct image $p_\ast$
  \eqref{BaseChangeAdjointQuadrupleOfEquivariantHomotopyTheoriesAlongEquivarianceCovering}
  of that of the cover
  $\widehat{G} \acts \, \Gamma \,\coloneqq\, (B p)^\ast ( G \acts \, \Gamma)$; see
  \eqref{LeftBaseChangeOfInfinityActionsAlongCoverOfEquivarianceGroup}.
\end{lemma}
\begin{proof}
  First, observe that for
  $G/H \,\in\, \OrbitCategory{G} \xhookrightarrow{\YonedaEmbedding}
  G\InfinityGroupoids$ there
  is the following sequence of natural equivalences
  \vspace{-2mm}
  $$
    \def\arraystretch{1.5}
    \begin{array}{lll}
      G\InfinityGroupoids
      \left(
        G/H
        ,\,
        p_\ast
        \,
        \EquivariantClassifyingShape{\widehat{G}}{\Gamma}
      \right)
      &
      \;\simeq\;
      G\InfinityGroupoids
      \left(
        p^\ast
        (G/H)
        ,\,
        \EquivariantClassifyingShape{\widehat{G}}{\Gamma}
      \right)
      &
      \proofstep{ by \eqref{BaseChangeAdjointQuadrupleOfEquivariantHomotopyTheoriesAlongEquivarianceCovering}
       with \eqref{AdjunctionAndHomEquivalence} }
      \\
      & \;\simeq\;
      G\InfinityGroupoids
      \left(
        \widehat{G}/\widehat{H}
        ,\,
        \EquivariantClassifyingShape{\widehat{G}}{\Gamma}
      \right)
      &
      \proofstep{ by
        Lem.
        \ref{BaseChangeInEquivariantHomotopyTheoryAlongCoveringsOfEquivaranceGroup}
        with
        Lem. \ref{LeftKanExtensionOnRepresentablesIsOriginalFunctor}
       }
      \\
      &
      \;\simeq\;
      \left(
        \EquivariantClassifyingShape{\widehat{G}}{\Gamma}
      \right)
      \left(\widehat{G}/\widehat{H}\right)
      &
      \proofstep{
        by
        Prop. \ref{InfinityYonedaLemma}
      }
      \\
      &
      \;\simeq\;
      \shape
      \,
      \SliceMaps{\big}{\mathbf{B}{\widehat{G}}}
        { \mathbf{B}\widehat{H} }
        { \HomotopyQuotient{ \mathbf{B}\Gamma }{ \widehat{G} } }
      &
      \proofstep{
        by Lem. \ref{FixedLociOfEquivariantClassifyingStacks};
      }
    \end{array}
  $$
  hence, by the $\infty$-Yoneda lemma (Prop. \ref{InfinityYonedaLemma})
  an equivalence
  \begin{equation}
    \label{BaseChangeDirectImageOfEquivariantClassifyingShape}
    p_\ast
    \,
    \EquivariantClassifyingShape{\widehat{G}}{\Gamma}
    \;\simeq\;
    \shape
    \,
    \SliceMaps{\big}{\mathbf{B}\widehat{G}}
      { \mathbf{B}\widehat{(-)} }
      { \HomotopyQuotient{ \mathbf{B}\Gamma }{\widehat{G}} }
    \;\simeq\;
    \shape
    \,
    \SliceMaps{\big}{\mathbf{B}\widehat{G}}
      { (B p)^\ast \mathbf{B}(-) }
      { (B p)^\ast \HomotopyQuotient{ \mathbf{B}\Gamma }{G} }
    \,.
  \end{equation}

  \vspace{0mm}
\noindent
  Under this identification, the natural monomorphism in question is that from
  Prop. \ref{PullbackOfActionsAlongSurjectiveGroupHomomorphismsIsFullyFaithful}:
  \vspace{-2mm}
  $$
    \def\arraystretch{1.7}
    \begin{array}{lll}
      \left(
        \EquivariantClassifyingShape{G}{\Gamma}
      \right)
      \left(G/(-)\right)
      &
      \;\simeq\;
      \shape
      \,
      \SliceMaps{\big}{\mathbf{B}G}
        { \mathbf{B}(-) }
        { \HomotopyQuotient{ \mathbf{B}\Gamma }{G} }
      &
      \proofstep{ by Lem. \ref{FixedLociOfEquivariantClassifyingStacks} }
      \\
      &
      \!\!\!\overset{ (B p)^\ast }{\hookrightarrow}\!\!
      \shape
      \,
      \SliceMaps{\big}{\mathbf{B}\widehat{G}}
        { \mathbf{B}\widehat{(-)} }
        { \HomotopyQuotient{ \mathbf{B}\Gamma }{\widehat{G}} }
      &
      \proofstep{ by Prop. \ref{PullbackOfActionsAlongSurjectiveGroupHomomorphismsIsFullyFaithful} }
      \\
      &
      \;\simeq\;
      \left(
      p_\ast
      \,
      \EquivariantClassifyingShape{\widehat{G}}{\Gamma}
      \right)
      \left(G/(-)\right)
      &
      \proofstep{
        by \eqref{BaseChangeDirectImageOfEquivariantClassifyingShape}.
      }
    \end{array}
  $$

  \vspace{-.7cm}

\end{proof}

\medskip

\noindent
{\bf Proper topological orbits.}
We use the above characterization
of equivariant classifying shapes
(Lem. \ref{FixedLociOfEquivariantClassifyingStacks})
to show
(Prop. \ref{SliceOfFiniteOrbitCategoryOverEquivariantClassifyingShape} below)
that the slice of the category of $G$-orbits over
the $G$-equivariant $\Gamma$-classifying shape
is fully faithfully embedding into the
$\infty$-site of proper $\Gamma$-orbits (Def. \ref{ProperTopologicalOrbitInfinityCategory} below) -- under a common assumption on the topological group $\Gamma$ (Def. \ref{AdmissibleEquivariantTwisting} and Prop. \ref{ExamplesOfAdmissibleGEquivariantTwistings} below).
Further below, this embedding under this admissibility condition implies the {\it twisted Elmendorf theorem} (Thm. \ref{TwistedEDKtheorem}).

\medskip

\begin{definition}[Admissible $G$-equivariant twisting]
  \label{AdmissibleEquivariantTwisting}
  For
  $G \,\in\, \Groups(\FiniteSets)$ a finite group we say that
  a $G$-equivariant
  topological group
  $G \acts \Gamma
    \,\in\,
    \Groups\big( \Actions{G}(\DTopologicalSpaces) \big)
  $
  -- or rather its delooping stack --
  is {\it admissible $G$-equivariant twisting}
  if for all (necessarily finite) subgroups $H \subset G$ and all $H \subset G$-crossed homomorphisms
  $\rho \,\in\, \CrossedHomomorphisms(H,H\acts \Gamma)$
  \eqref{CrossedHomomorphismsAndFirstNonAbelianGroupCohomology}
  the operations of

  (i) passing to shape $\shape$
   \eqref{TheModalitiesOnACohesiveInfinityTopos}

  (ii) forming the looping $\Omega_\rho$
   \eqref{LoopingInIntroduction} at $\rho$
   (via the identification of Prop. \ref{ConjugationGroupoidOfCrossedHomomorphismsIsSectionsOfDeloopedSemidirectProductProjection})

   \noindent
   commute on the slice mapping stack
   (Def. \ref{SliceMappingStack})
   from $\mathbf{B}H$ to $\mathbf{B}(\Gamma \rtimes G)$ over $\mathbf{B}G$, up to equivalence:
  \begin{equation}
    \label{AssumingShapeCommutesWithLoopingOnMappingStackBetweenDeloopedTopologicalGroups}
    \left.
      \mbox{
        \begin{tabular}{l}
          $\mathbf{B}
          (\Gamma \rtimes G)$
          is admissible
          \\
          as $G$-equivariant twisting
        \end{tabular}
      }
      \hspace{-.3cm}
    \right\}
    \hspace{.4cm}
    \Leftrightarrow
    \hspace{.0cm}
    \underset{
      \scalebox{.7}{$
        \begin{array}{c}
          H \subset G
          \\
          \rho
            \in
          \CrossedHomomorphisms
          (H,\, H\acts \Gamma)
        \end{array}
      $}
    }{\forall}
    \;
    \Omega_\rho
    \,
    \shape
    \,
    \SliceMaps{\big}{\mathbf{B}G}
      { \mathbf{B}H }
      { \mathbf{B}(\Gamma \rtimes G) }
    \;\;
    \simeq
    \;\;
    \shape
    \,
    \Omega_{\rho}
    \,
    \SliceMaps{\big}{\mathbf{B}G}
      {\mathbf{B}H}
      {\mathbf{B}(\Gamma \rtimes G)}
    \,.
  \end{equation}
  In the special case that the $G$-action on $\Gamma$ is trivial, this condition reduces to:
  \begin{equation}
    \label{AssumingShapeCommutesWithLoopingOnMappingStackBetweenDeloopedTopologicalGroupsForTrivialGActionOnGamma}
    \left.
      \mbox{
        \begin{tabular}{l}
          $\mathbf{B}
          (\Gamma \times G)$
          is admissible
          \\
          as $G$-equivariant twisting
        \end{tabular}
      }
      \hspace{-.3cm}
    \right\}
    \hspace{.4cm}
    \Leftrightarrow
    \hspace{.0cm}
    \underset{
      \scalebox{.7}{$
        \begin{array}{c}
          H \subset G
          \\
          \rho
            \in
          \Homomorphisms
          (H,\, \Gamma)
        \end{array}
      $}
    }{\forall}
    \;
    \Omega_\rho
    \,
    \shape
    \,
    \Maps{\big}
      { \mathbf{B}H }
      { \mathbf{B}\Gamma }
    \;\;
    \simeq
    \;\;
    \shape
    \,
    \Omega_{\rho}
    \,
    \Maps{\big}
      {\mathbf{B}H}
      {\mathbf{B}\Gamma}
    \,.
  \end{equation}
\end{definition}
\begin{remark}[Alternative formulations of admissible $G$-equivariant twists]
\label{AlternativeFormulationsOfAdmissibleTwists}
  Noticing that
  $$
    \def\arraystretch{1.6}
    \begin{array}{ll}
      \Omega_\rho
      \SliceMaps{}{\mathbf{B}G}
        {\mathbf{B}H}{\mathbf{B}\Gamma}
      \\
      \;\simeq\;
      \Omega_\rho
      \big(
      \HomotopyQuotient
        { \Homs{}{H}{\Gamma} }
        { \Gamma }
      \big)
      &
      \proofstep{
        by
        \eqref{ConjugationGroupoidOfCrossedHomomorphismsIsSlicedMappingGroupoidFromBGToBGammaRtimesG}
      }
      \\
      \;\simeq\;
      \mathrm{Stab}_\Gamma(\rho)
      &
      \proofstep{
        by
        \eqref{TopologicalStabilizerGroupAsLooping}
      }
    \end{array}
  $$
  is the topological stabilizer group of $\rho$ under the conjugation action by $\Gamma$, and recalling
  $$
    \def\arraystretch{1.4}
    \begin{array}{ll}
      B
      \,
      \mathrm{Stab}_\Gamma(\rho)
      \\
      \;\simeq\;
      \shape
      \,
      \mathbf{B}
      \,
      \mathrm{Stab}_\Gamma(\rho)
      &
      \proofstep{
        by Prop.
        \ref{MilgramClassifyingSpaceModelClassifyingShape}
      }
      \\
      \;\simeq\;
      B
      \,
      \shape
      \,
      \mathrm{Stab}_\Gamma(\rho)
      &
      \proofstep{
        by \eqref{ShapePreservesDeloopings}
      }
    \end{array}
  $$
we may,
under the looping/delooping equivalence (Prop. \ref{LoopingAndDeloopingEquivalence}),
equivalently rewrite \eqref{AssumingShapeCommutesWithLoopingOnMappingStackBetweenDeloopedTopologicalGroups} as follows, which proves the claim:

  \begin{equation}
    \label{ClassicalNotationForAssumingShapeCommutesWithLoopingOnMappingStackBetweenDeloopedTopologicalGroups}
    \left.
      \mbox{
        \begin{tabular}{l}
          $\mathbf{B}
          (\Gamma \rtimes G)$
          is admissible
          \\
          as $G$-equivariant twisting
        \end{tabular}
      }
      \hspace{-.3cm}
    \right\}
    \hspace{.4cm}
    \Leftrightarrow
    \hspace{.0cm}
    \underset{
      \scalebox{.7}{$
        \begin{array}{c}
          H \subset G
          \\
          \rho
            \in
          \CrossedHomomorphisms
          (H,\, H\acts \Gamma)
        \end{array}
      $}
    }{\forall}
    \;
    \Big(
    \shape
    \,
    \SliceMaps{\big}{\mathbf{B}G}
      { \mathbf{B}H }
      { \mathbf{B}(\Gamma \rtimes G) }
    \Big)_\rho
    \;\;
    \simeq
    \;\;
    B
    \,
    \mathrm{Stab}_\Gamma(\rho)
    \,,
  \end{equation}
  where on the right $(-)_\rho$ denotes the connected component of $\rho$.
  Therefore,  a sufficient condition for admissible $G$-equivariant twists is  that the $\Gamma$-conjugacy classes of crossed homomorphisms bijectively label these connected components:
  \begin{equation}
    \label{WeakerClassicalVersionForAssumingShapeCommutesWithLoopingOnMappingStackBetweenDeloopedTopologicalGroups}
    \left.
      \mbox{
        \begin{tabular}{l}
          $\mathbf{B}
          (\Gamma \rtimes G)$
          is admissible
          \\
          as $G$-equivariant twisting
        \end{tabular}
      }
      \hspace{-.3cm}
    \right\}
    \hspace{.9cm}
    \Leftarrow
    \hspace{.7cm}
    \underset{
      \scalebox{.7}{$
        H \subset G
      $}
    }{\forall}
    \;\;\;\;\;
    \shape
    \,
    \SliceMaps{\big}{\mathbf{B}G}
      { \mathbf{B}H }
      { \mathbf{B}(\Gamma \rtimes G) }
    \;\;
    \simeq
    \;\;
    \underset{
      \mathclap{
      \scalebox{.7}{$
        \begin{array}{c}
        [\rho]
        \in
        \\
        \CrossedHomomorphisms(H,H\acts \Gamma)/\Gamma
        \end{array}
      $}
      }
    }{\coprod}
    \;
    B
    \,
    \mathrm{Stab}_\Gamma(\rho)
    \,,
  \end{equation}
  Moreover, noticing that
  $$
    \def\arraystretch{1.5}
    \begin{array}{ll}
      \underset{[\rho]}{\coprod}
      B \mathrm{Stab}_\Gamma(\rho)
      \\
      \;\simeq\;
      \shape
      \,
      \underset{[\rho]}{\coprod}
      \,
      \mathbf{B}
      \,
      \mathrm{Stab}_\Gamma(\rho)
      &
      \proofstep{
        by
        Prop.
        \ref{MilgramClassifyingSpaceModelClassifyingShape}
      }
      \\
      \;\simeq\;
      \shape
      \,
      \underset{[\rho]}{\coprod}
      \;
      \HomotopyQuotient{
        \big(
        \Gamma/\mathrm{Stab}_\Gamma(\rho)
        \big)
      }{\Gamma}
      &
      \proofstep{
        by
        Exp. \ref{QuotientStackOfCosetSpacesIsDeloopingOfSubgroup}
      }
    \end{array}
  $$
  we find that a yet stronger (i.e. further from necessary but still) sufficient condition for admissibility is that the conjugacy classes already label the connected components of the action groupoid of crossed homomorphisms (i.e. even before passing to its shape):
  \begin{equation}
    \label{YetWeakerClassicalVersionForAssumingShapeCommutesWithLoopingOnMappingStackBetweenDeloopedTopologicalGroups}
    \left.
      \mbox{
        \begin{tabular}{l}
          $\mathbf{B}
          (\Gamma \rtimes G)$
          is admissible
          \\
          as $G$-equivariant twisting
        \end{tabular}
      }
      \hspace{-.3cm}
    \right\}
    \hspace{.9cm}
    \Leftarrow
    \hspace{.7cm}
    \underset{
      \scalebox{.7}{$
        H \subset G
      $}
    }{\forall}
    \;\;\;\;\;
    \SliceMaps{\big}{\mathbf{B}G}
      { \mathbf{B}H }
      { \mathbf{B}(\Gamma \rtimes G) }
    \;\;
    \simeq
    \;\;
    \underset{
      \mathclap{
      \scalebox{.7}{$
        \begin{array}{c}
        [\rho]
        \in
        \\
        \CrossedHomomorphisms(H,H\acts \Gamma)/\Gamma
        \end{array}
      $}
      }
    }{\coprod}
    \;
    \HomotopyQuotient{
      \big(
      \Gamma/\mathrm{Stab}_\Gamma(\rho)
      \big)
    }
    {\Gamma}
    \,.
  \end{equation}
  For trivial $G$-action on $\Gamma$ these specializations all simplify as in
  \eqref{AssumingShapeCommutesWithLoopingOnMappingStackBetweenDeloopedTopologicalGroupsForTrivialGActionOnGamma}.
\end{remark}
\begin{proposition}[Examples of admissible $G$-equivariant twistings]
\label{ExamplesOfAdmissibleGEquivariantTwistings}
The condition of Def. \ref{AdmissibleEquivariantTwisting} holds at least for:

(1) $\Gamma$ any Lie group (not necessarily compact) with trivial $G$-action,

(2) $\Gamma$ a compact Lie group with any $G$-action,

(3) $\Gamma = \PUH$ \eqref{TheGroupPUH}
 with trivial $G$-action.

\end{proposition}
\begin{proof}
  (1) For $\Gamma$ a Lie group
  with trivial $G$-action,
  the sufficient condition \eqref{YetWeakerClassicalVersionForAssumingShapeCommutesWithLoopingOnMappingStackBetweenDeloopedTopologicalGroups} follows readily  from the classical result \cite[Lem. 31.8]{ConnerFloyd64}, as noticed in \cite[p. 5]{Rezk14}. In the special case that $\Gamma$ is almost connected (i.e. the quotient space $\Gamma/\Gamma_{\NeutralElement}$ by the connected component of the neutral element is compact), the sufficient condition \eqref{WeakerClassicalVersionForAssumingShapeCommutesWithLoopingOnMappingStackBetweenDeloopedTopologicalGroups} is verified in \cite[Thm. 6.3 with \S 13]{LueckUribe14}.

  \noindent
  (2)
  For $\Gamma$ a compact Lie group with any $G$-action, the sufficient condition \eqref{WeakerClassicalVersionForAssumingShapeCommutesWithLoopingOnMappingStackBetweenDeloopedTopologicalGroups} follows from the combination of \cite[Thm. 10]{LashofMay86}\cite[Thm. 4.24]{GuillouMayMerling17} with \cite{MurayamaShimakawa95}\cite[Thm. 3.11]{GuillouMayMerling17}.

  \noindent
  (3)
  For $\Gamma = \PUH$ with trivial $G$-action,
  the sufficient condition \eqref{WeakerClassicalVersionForAssumingShapeCommutesWithLoopingOnMappingStackBetweenDeloopedTopologicalGroups} is verified in \cite[\S 15 with \S 13]{LueckUribe14}.
\end{proof}

\medskip

We now consider the generalization of the orbit category (Ntn. \ref{GOrbitCategory}) from finite groups $G$ to topological groups $\Gamma$.
The following Def. \ref{ProperTopologicalOrbitInfinityCategory} is the usual definition of the orbit category of a topological group as in \cite[\S 2.1 and Thm. 3.1]{DwyerKan84}, only that we"

\noindent

\vspace{-.1cm}
\begin{itemize}

\vspace{-.2cm}
\item[(i)]
consider that topological group to be decomposed as the semidirect product
$\Gamma \rtimes G$ of a $G$-equivariant topological group $\Gamma$ with a finite group $G$ (since this is the case of interest in the following);

\vspace{-.2cm}
\item[(ii)]
restrict attention to orbits given by finite subgroups (indicated by the adjective ``proper'');

\vspace{-.2cm}
\item[(ii)]
understand (as in \cite{Rezk14}) the $\SimplicialSets$-enrichment considered in \cite{DwyerKan84} as making an $\infty$-category (an $\infty$-site in the sense of Ntn. \ref{InfinitySite}, by appeal to our Prop. \ref{SmoothShapeOfTopologicalSpacesIsTheirWeakHomotopyType});

\vspace{-.2cm}
\item[(iv)]
immediately reformulate
\eqref{HomInfinityGroupoidOfProperTopologicalOrbitCategory}
the equivariant mapping space between the orbits as a slice mapping stack between the corresponding delooped subgroups, in direct generalization of Lemma \ref{GOrbitsAre0TruncatedObjectsOverGOrbiSingularity} (to which the following reduces when $\Gamma = 1$ is the trivial group).
\end{itemize}

\begin{definition}[Proper topological orbit $\infty$-category]
\label{ProperTopologicalOrbitInfinityCategory}
For any $G$-equivariant D-topological group
(not necessarily compact)

- $G \,\in\, \Groups(\FiniteSets)$

- $G \acts \Gamma \,\in\, \Groups\big(\Actions{G}(\DTopologicalSpaces)\big)$

\noindent
we say that the proper {\it orbit $\infty$-category} of its semidirect product group $\Gamma \rtimes G$
is the $\infty$-site (Ntn. \ref{InfinitySite})
$\Orbits(\Gamma \rtimes G)$ whose objects are finite subgroups
$$
   \eta
   \,\in\,
   \Groups(\FiniteSets)
   \,,
   \;\;\;\;\;\;
  \eta
  \xhookrightarrow{\;\;
    i_{\eta}
  \;}
  \Gamma \rtimes G
$$
and whose hom-$\infty$-groupoids are the shapes
of the equivariant mapping spaces between the corresponding coset spaces $(\Gamma \rtimes G)/\eta$, hence equivalently
the shapes of the slice mapping stacks
\eqref{TheSliceMappingStack}
between the deloopings of the subgroups, relative to $\mathbf{B}(\Gamma \rtimes G)$:
\begin{equation}
  \label{HomInfinityGroupoidOfProperTopologicalOrbitCategory}
  \def\arraystretch{1.9}
  \begin{array}{ll}
  \def\AmbientCategory{\Orbits(\Gamma\rtimes G)}
  \PointsMaps{\Big}
    { (\Gamma \rtimes G)/\eta_1 }
    { (\Gamma \rtimes G)/\eta_2 }
  \\
  \;:=\;
  \shape
  \bigg(
  \Maps{\Big}
    {
       (\Gamma\rtimes G)
       \acts
       \big(
         (\Gamma \rtimes G)/\eta_1
       \big)
    }
    {
       (\Gamma\rtimes G)
       \acts
       \big(
         (\Gamma \rtimes G)/\eta_2
       \big)
    }
  ^{(\Gamma \rtimes G)}
  \bigg)
  \\
  \;\simeq\;
  \shape
  \,
  \underset
    {{\mathbf{B}(\Gamma \rtimes G)}}
    {\prod}
  \Maps{\big}
    { \mathbf{B}(i_{\eta_1}) }
    { \mathbf{B}(i_{\eta_2}) }
  &
  \proofstep{
    by
    Exp. \ref{FixedLociOfInfinityActions}
  }
  \\
  \;\simeq\;
  \shape
  \Big(
  \SliceMaps{\big}{\mathbf{B}(\Gamma \rtimes G)}
    { \mathbf{B}\eta_1 }
    { \mathbf{B}\eta_2 }
  \Big)
  &
  \proofstep{
    by
    Lem. \ref{PlotsOfSliceMappingStackAreSliceHoms}.
  }
  \end{array}
\end{equation}

\end{definition}

\begin{proposition}[Slice of finite orbit category over equivariant classifying shape inside proper topological orbit category]
  \label{SliceOfFiniteOrbitCategoryOverEquivariantClassifyingShape}
  For
  $G \,\in\, \Groups(\FiniteSets)$
  and
  $G \acts \Gamma \,\in\, \Actions{G}(\DTopologicalSpaces)$
such that $\mathbf{B}(\Gamma \rtimes G)$ is admissible $G$-equivariant twisting (Def. \ref{AdmissibleEquivariantTwisting}),
then
the slice $\infty$-site (Def. \ref{SliceInfinitySite})
of the $G$-orbit category (Ntn. \ref{GOrbitCategory})
over the $G$-equivariant classifying shape
$\EquivariantClassifyingShape{G}{\Gamma}$
(Def. \ref{EquivariantClassifyingStack})
is equivelently the full sub-$\infty$-site
of the topological orbit $\infty$-category of $\Gamma \rtimes G$ (Def. \ref{ProperTopologicalOrbitInfinityCategory})
\begin{equation}
  \label{FullInclusionOfCrossedSubgroupsIntoTopologicalOrbit}
  \begin{tikzcd}[row sep=2pt]
    \Orbits(G)_{/\EquivariantClassifyingShape{G}{\Gamma}}
    \;\simeq\;
    \Orbits(\Gamma \rtimes G)_{\mathrm{crs}}
    \ar[
      r,
      hook,
      "{
        i_{\mathrm{crs}}
      }"
    ]
    &
    \Orbits(\Gamma \rtimes G)
  \end{tikzcd}
  \;\;\;\;\;
  \in
  \;
  \mathrm{Site}_\infty
  \,.
\end{equation}
on orbits of the form $\Gamma/\widehat{H}$ for
subgroups $\widehat{H} \,\subset\, \Gamma \rtimes G$
(Ntn. \ref{LiftsOfEquivarianceSubgroupsToSemidirectProductWithStructureGroup})
which are lifts of subgroups $H \subset G$:
\begin{equation}
  \label{LiftsOfSubgroupsToSemidirectProductInDiscussionOfProperOrbitCategory}
  \begin{tikzcd}
    \widehat{H}
    \ar[
      r,hook, "{i_{\widehat{H}}}"
    ]
    &
    \Gamma \rtimes G
    \ar[d, "\mathrm{pr}_2"]
    \\
    H
    \ar[
      r,hook, "{i_{H}}"
    ]
    \ar[u, "{\sim}"{sloped}]
    &
    G
    \mathrlap{\,.}
  \end{tikzcd}
\end{equation}
\end{proposition}
\begin{proof}
  Recall from Lem. \ref{FixedLociOfEquivariantClassifyingStacks} that, as an $\infty$presheaf over $\Orbits(G)$, the equivariant classifying shape is equivalently given by the assignment
$$
  \EquivariantClassifyingShape{G}{\Gamma}
  \;:\;
  G/H
  \;\longmapsto\;
  \shape
  \,
  \SliceMaps{\big}{\mathbf{B}G}
    { \mathbf{B}H }
    { \mathbf{B}(\Gamma \rtimes G) }
  \,.
$$
By the $\infty$-Yoneda lemma (Prop. \ref{InfinityYonedaLemma}) this implies that the objects of the slice
$\Slice{\Orbits(G)}{\EquivariantClassifyingShape{G}{\Gamma}}$
are lifts $\widehat{H}$ of subgroups $H \subset G$ to $\Gamma \rtimes G$ as in \eqref{LiftsOfSubgroupsToSemidirectProductInDiscussionOfProperOrbitCategory};
while the hom-$\infty$-groupoid in the slice
between a pair of such objects is naturally equivalent
to that between the corresponding $(\Gamma \rtimes G)$-orbits, as follows:
$$
  \def\arraystretch{1.8}
  \begin{array}{ll}
    \def\AmbientCategory{
    \Slice{\Orbits(G)}
      {\EquivariantClassifyingShape{G}{\Gamma}}
    }
    \PointsMaps{\Big}
      {
        \big(G/H_1, \mathbf{B}i_{\widehat{H}_1}\big)
      }
      {
        \big(
          G/H_2, \mathbf{B}i_{\widehat{H}_2}
        \big)
      }
    \\
    \;\simeq\;
    \Orbits(G)(G/H_1, G/H_2)
    \underset{
      \scalebox{.7}{$
      \def\AmbientCategory{\GEquivariant\InfinityGroupoids}
      \PointsMaps{\big}
        {
          \FixedLoci(G/H_1)
        }
        { \EquivariantClassifyingShape{G}{\Gamma} }
      $}
    }{\times}
    \Big\{
      \mathbf{B}\big(i_{\widehat{H_1}}\big)
    \Big\}
    &
    \proofstep{
      by
      \eqref{HomInfinityGroupoidOfSliceInfinitySite}
    }
    \\
    \;\simeq\;
    \SliceMaps{\big}{\mathbf{B}G}
      { \mathbf{B}H_1 }
      { \mathbf{B}H_2 }
    \underset{
      \scalebox{.7}{$
        \shape
        \,
        \SliceMaps{\big}{\mathbf{B}G}
          { \mathbf{B}H_1 }
          { \mathbf{B}(\Gamma \rtimes G) }
      $}
    }{\times}
    \Big\{
      \mathbf{B}\big(i_{\widehat{H_1}}\big)
    \Big\}
    &
    \proofstep{
      by Lem. \ref{GOrbitsAre0TruncatedObjectsOverGOrbiSingularity}
      \&
      Lem. \ref{FixedLociOfEquivariantClassifyingStacks}
    }
    \\
    \;\simeq\;
    \SliceMaps{\big}{\mathbf{B}G}
      { \mathbf{B}H_1 }
      { \mathbf{B}H_2 }
      ^{ \widehat{H}_1 }
    \underset{
      \scalebox{.7}{$
        B
        \Omega_{\widehat{H}_1}
        \shape
        \,
        \SliceMaps{\big}{\mathbf{B}G}
          { \mathbf{B}H_1 }
          { \mathbf{B}(\Gamma \rtimes G) }
      $}
    }{\times}
    \Big\{
      \mathbf{B}\big(i_{\widehat{H_1}}\big)
    \Big\}
    &
    \proofstep{
      by
      \eqref{SliceMappingStackAsFiberOverConnectedComponentOfCodomainProjection}
    }
    \\
    \;\simeq\;
    \SliceMaps{\big}{\mathbf{B}G}
      { \mathbf{B}H_1 }
      { \mathbf{B}H_2 }
     ^{\widehat{H}_1}
    \underset{
      \scalebox{.7}{$
        \shape
        \mathbf{B}
        \Omega_{\widehat{H}_1}
        \,
        \SliceMaps{\big}{\mathbf{B}G}
          { \mathbf{B}H_1 }
          { \mathbf{B}(\Gamma \rtimes G) }
      $}
    }{\times}
    \Big\{
      \mathbf{B}\big(i_{\widehat{H_1}}\big)
    \Big\}
    &
    \proofstep{
      by
      \eqref{AssumingShapeCommutesWithLoopingOnMappingStackBetweenDeloopedTopologicalGroups}
    }
    \\
    \;\simeq\;
    \shape
    \bigg(
    \SliceMaps{\big}{\mathbf{B}G}
      { \mathbf{B}H_1 }
      { \mathbf{B}H_2 }
     ^{\widehat{H}_1}
    \underset{
      \scalebox{.7}{$
        \mathbf{B}
        \Omega_{\widehat{H}_1}
        \,
        \SliceMaps{\big}{\mathbf{B}G}
          { \mathbf{B}H_1 }
          { \mathbf{B}(\Gamma \rtimes G) }
      $}
    }{\times}
    \Big\{
      \mathbf{B}\big(i_{\widehat{H_1}}\big)
    \Big\}
    \bigg)
    &
    \proofstep{
      by
      Lem.
      \ref{ShapePreservesHomotopyFibersOfDeloopingsOutOfDiscreteDomains}
    }
    \\
    \;\simeq\;
    \shape
    \bigg(
    \SliceMaps{\big}{\mathbf{B}G}
      { \mathbf{B}H_1 }
      { \mathbf{B}H_2 }
    \underset{
      \scalebox{.7}{$
        \SliceMaps{\big}{\mathbf{B}G}
          { \mathbf{B}H_1 }
          { \mathbf{B}(\Gamma \rtimes G) }
      $}
    }{\times}
    \Big\{
      \mathbf{B}\big(i_{\widehat{H_1}}\big)
    \Big\}
    \bigg)
    &
    \proofstep{
      by
      \eqref{SliceMappingStackAsFiberOverConnectedComponentOfCodomainProjection}
    }
    \\
    \;\simeq\;
    \shape
    \,
    \SliceMaps{\Big}{ \mathbf{B}(\Gamma \rtimes G) }
      {
        \big(
          \mathbf{B}H_1
          ,
          \mathbf{B}(i_{\widehat{H}_1})
        \big)
      }
      {
        \big(
          \mathbf{B}H_2
          ,
          \mathbf{B}(i_{\widehat{H}_2})
        \big)
      }
    &
    \proofstep{
      by
      Lem. \ref{SlicedSlicedMappingSpace}
    }
    \\
    \;\simeq\;
    \Orbits(\Gamma\rtimes G)
    \Big(
      (\Gamma\rtimes G)/\widehat{H_1}
      ,\,
      (\Gamma\rtimes G)/\widehat{H_2}
    \Big)
    &
    \proofstep{
      by
      \eqref{HomInfinityGroupoidOfProperTopologicalOrbitCategory}.
    }
  \end{array}
$$
Here in the second but last step we used the following observation (which could be stated much more generally):
\begin{equation}
  \label{SlicedSlicedMappingSpace}
    \SliceMaps{}{\mathbf{B}G}
      { \mathbf{B}H_1 }
      { \mathbf{B}H_2 }
    \underset{
      \SliceMaps{}{\mathbf{B}G}
        { \mathbf{B}H_1 }
        { \mathbf{B}(\Gamma \rtimes G) }
    }{\times}
    \big\{
      \mathbf{B}i_{\widehat{H}_2s}
    \big\}
    \;\;\;\;
    \simeq
    \;\;\;\;
    \SliceMaps{\Big}{ \mathbf{B}(\Gamma\rtimes G) }
      {
        \big(
          \mathbf{B}H_1,
          \,
          \mathbf{B}(i_{\widehat{H_1}})
        \big)
      }
      {
        \big(
          \mathbf{B}H_2,
          \,
          \mathbf{B}(i_{\widehat{H_2}})
        \big)
      }
      \,.
\end{equation}
To see this,
considering the following diagram in $\SmoothInfinityGroupoids$, which commutes by \eqref{LiftsOfSubgroupsToSemidirectProductInDiscussionOfProperOrbitCategory}:
$$
  \begin{tikzcd}[
    column sep=40pt,
    row sep=30pt
  ]
    \Maps{}
      { \mathbf{B}H_1 }
      { \mathbf{B}H_2 }
    \ar[
      r,
      "{
        \big(
          \mathbf{B}(i_{H_2})
        \big)
        \circ
        (-)
      }"
    ]
    \ar[
      d,
      "{
        \mathbf{B}
        \big(i_{\widehat{H_2}}\big)
        \circ
        (-)
      }"{swap}
    ]
    &
    \Maps{}
      { \mathbf{B}H_1 }
      { \mathbf{B}G }
    \ar[
      from=r,
      "{
        \vdash
        \,
        \mathbf{B}
        (i_{H_2})
      }"{swap}
    ]
    \ar[
      d,
      Rightarrow,
      -
    ]
    &
    \ast
    \ar[
      d,
      Rightarrow,
      -
    ]
    &[-15pt]
    \SliceMaps{}{\mathbf{B}G}
      {\mathbf{B}H_1}
      {\mathbf{B}H_2}
    \ar[d]
    \\
    \Maps{}
      { \mathbf{B}H_1 }
      { \mathbf{B}(\Gamma \rtimes G) }
    \ar[
      r,
      "{
        \big(
         \mathbf{B}
         (\mathrm{pr}_2)
         \big)
         \circ
         (-)
      }"
    ]
    &
    \Maps{}
      { \mathbf{B}H_1 }
      { \mathbf{B}G }
    \ar[
      from=r,
      "{
        \vdash
        \,
        \mathbf{B}
        (i_{H_2})
      }"{swap}
    ]
    &
    \ast
    &
    \SliceMaps{}{\mathbf{B}G}
      {\mathbf{B}H_1}
      {\mathbf{B}(\Gamma \rtimes G)}
    \\
    \ast
    \ar[r]
    \ar[
      u,
      "{
        \vdash
        \,
        \mathbf{B}
        (i_{\widehat{H_1}})
      }"
    ]
    &
    \ast
    \ar[
      from=r
    ]
    \ar[
      u,
      "{
        \vdash
        \,
        \mathbf{B}
        (i_{H_1})
      }"
    ]
    &
    \ast
    \ar[
      u,
      Rightarrow,
      -
    ]
    &
    \ast
    \ar[u]
    \\[-10pt]
    \SliceMaps{\Big}{ \mathbf{B}(\Gamma\rtimes G) }
      {
        \big(
          \mathbf{B}H_1,
          \,
          \mathbf{B}(i_{\widehat{H_1}})
        \big)
      }
      {
        \big(
          \mathbf{B}H_2,
          \,
          \mathbf{B}(i_{\widehat{H_2}})
        \big)
      }
      \ar[r]
      &
      \ast
      \ar[from=r]
      &
      \ast
  \end{tikzcd}
$$
By Def. \ref{SliceMappingStack}, forming limits horizontally over the rows or vertically over the columns yields the diagrams of slice mapping stacks shown on the right and on the bottom, respectively. The further limits over these two diagrams manifestly yield the two sides of the equivalence \eqref{SlicedSlicedMappingSpace}, which hence follows by the general fact that limits commute over each other.
\end{proof}

\begin{definition}[Systems of shapes of fixed loci in topological $\Gamma$-spaces]
\label{SystemsOfShapesOfFixedLociInTopologicalGammaSpaces}
For $\Gamma \,\in\, \Groups(\kTopologicalSpaces)$ and $\Gamma \acts \TopologicalSpace \,\in\, \Actions{\Gamma}(\kTopologicalSpaces)$, we write
$$
  \shape
  \,
  \FixedLoci(\TopologicalSpace)
  \,\in\,
  \InfinityPresheaves
  \big(
    \Orbits(\Gamma)
  \big)
$$
for the $\infty$-presheaf (Ntn. \ref{ModelCategoriesOfSimplicialPresheaves}) over the proper orbit category of $\Gamma$ (Def. \ref{ProperTopologicalOrbitInfinityCategory}) which assigns the shapes (presented by the classical singular simplicial complexes, via Prop. \ref{DiffeologicalMappingSpacesHaveCorrectUnderlyingHomotopyType}) of the fixed loci (Ntn. \ref{GActionOnTopologicalSpaces}) of $\TopologicalSpace$:
$$
  \FixedLoci
  \;:\;
  \Gamma/\eta
  \;\longmapsto\;
  \shape
  \,
  \Big(
  \Maps{\big}
    { \Gamma \acts (\Gamma/\eta) }
    { \Gamma \acts \TopologicalSpace }
  ^{\Gamma}
  \Big)
  \,.
$$
In the special case that $\Gamma = G$ is a discrete group, this reduces equivalently to
\eqref{ShapeOfGeometricFixedLocus}.
\end{definition}

\begin{example}[System of $(\Gamma \rtimes G)$-fixed loci of a $G$-equivariant $\Gamma$-principal topological bundle]
  \label{PresheafHomOutOfSystemOfCrossedFixedLociOfEquivariantPrincipalBundlerestrictsToCrossedSubgroups}
Let

- $G \,\in\, \Groups(Sets)$,

- $G \acts \Gamma \,\in\, \Groups\big( \Actions{G}(\kTopologicalSpaces) \big)$,

- $G \acts \TopologicalSpace \,\in\, \Actions{G}(\kTopologicalSpaces)$,

- $G \acts \TopologicalPrincipalBundle \,\in\, \EquivariantPrincipalFiberBundles{G}{\Gamma}_{\TopologicalSpace}$.

\noindent
Then the induced $(\Gamma \rtimes G)$-action on the total space $\TopologicalPrincipalBundle$ has non-empty fixed loci only for finite subgroups which are lifts $\widehat{H}$
\eqref{StabilizerOfPointInFiber}
of subgroups $H \subset G$, and for these the fixed locus is isomorphic to that of the base space:
\begin{equation}
  \label{CrossedGroupFixedLociOfEquivariantPrincipalBundle}
  \TopologicalPrincipalBundle
  \;\in\;
  \EquivariantPrincipalFiberBundles{G}{\Gamma}
  \;\;\;\;\;\;\;\;
  \Rightarrow
  \;\;\;\;\;\;\;\;
  \FixedLoci(\TopologicalPrincipalBundle)
  \;:\;
  \Gamma/\eta
  \;\mapsto\;
  \left\{
  \def\arraystretch{1.6}
  \begin{array}{lll}
    \TopologicalPrincipalBundle^{\widehat{H}}
    &\vert&
    \eta
    \;\simeq
    \!
    \adjustbox{raise=-5pt}{
    \begin{tikzcd}[row sep=2pt]
      \widehat{H}
      \ar[r, hook]
        &
        \Gamma \rtimes G
      \ar[d, shorten=-2pt, "{\mathrm{pr}_2}"]
      \\
      H
      \ar[r, hook]
      \ar[u, shorten=-2pt, "\sim"{sloped}]
      & G
    \end{tikzcd}
    }
    \\
    \varnothing &\vert&
    \mbox{otherwise.}
  \end{array}
  \right.
\end{equation}
(To see this, recall the proof of Prop. \ref{CharacterizationOfEquivariantBundlesOverCosetSpaces}: The $\Gamma$-action is free, so that for each $h \in H \subset G$ which fixes a point $x \in \TopologicalSpace$ in the base, there is a {\it unique} element $\gamma_x(h) \in \Gamma$ which compensates the $h$-action on the fiber $\TopologicalPrincipalBundle_x$.)

But this means that a morphism of presheaves out of $\shape \, \FixedLoci(\TopologicalPrincipalBundle)$ is uniquely fixed by its restriction to the full
sub-$\infty$-site $\Orbits(\Gamma \rtimes G)_{\mathrm{crs}}$ \eqref{FullInclusionOfCrossedSubgroupsIntoTopologicalOrbit} on orbits corresponding to these subgroups, in that for all $\mathcal{A} \,\in\, \InfinityPresheaves\big(\Orbits(\Gamma \rtimes G)\big)$ we have the following natural equivalences:
\begin{equation}
  \TopologicalPrincipalBundle
  \;\in\;
  \EquivariantPrincipalFiberBundles{G}{\Gamma}
    _{\TopologicalSpace}
  \;\;\;\;\;
  \Rightarrow
  \;\;\;\;\;
  \def\arraystretch{1.8}
  \begin{array}{ll}
  \def\AmbientCategory{
    \InfinityPresheaves
    \big(
      \Orbits(\Gamma\rtimes G)
    \big)
  }
  \PointsMaps{\Big}
    {
      \shape
      \,
      \FixedLoci(\TopologicalPrincipalBundle)
    }
    {
      \mathcal{A}
    }
  &
  \\
  \;
    \xrightarrow[\sim]{i_{\mathrm{crs}}^\ast}
  \;
  \def\AmbientCategory{
    \InfinityPresheaves
    \big(
      \Orbits(\Gamma \rtimes G)_{\mathrm{crs}}
    \big)
  }
  \PointsMaps{\Big}
    {
      \shape
      \,
      \FixedLoci(\TopologicalPrincipalBundle)
    }
    {
      \mathcal{A}
    }
    &
  \proofstep{
    by
    \eqref{CrossedGroupFixedLociOfEquivariantPrincipalBundle}\footnote{
In more abstract terms: Since the full inclusion $i_{\mathrm{crs}}$ is also a sieve, the left Kan extension of presheaves along $i_{\mathrm{crs}}$ is given by  extending by empty values.}
  }
  \\
  \;\simeq\;
  \def\AmbientCategory{
    \InfinityPresheaves
    \big(
      \Slice{\Orbits(G)}
        {\EquivariantClassifyingShape{G}{\Gamma}}
    \big)
  }
  \PointsMaps{\Big}
    {
      \shape
      \,
      \big(
        \FixedLoci(\TopologicalSpace),
        c
      \big)
    }
    {
      \mathcal{A}
    }
    &
    \proofstep{
      \hspace{-.3cm}
      \def\arraystretch{.9}
      \begin{tabular}{l}
      by
      Prop. \ref{SliceOfFiniteOrbitCategoryOverEquivariantClassifyingShape}
      \\
      and \eqref{FixedLociOfEquivariantPrincipalBundlesViaSliceHom}.
      \end{tabular}
    }
  \end{array}
\end{equation}
Here in the last line $c \,\colon\, \TopologicalSpace \xrightarrow{\;} \mathbf{B}\Gamma$ is the ($G$-equivariant) modulating morphism of the given equivariant bundle.
\end{example}

\def\AmbientCategory{\Topos}

\medskip
\noindent
{\bf Recovering the Murayama-Shimakawa construction.}
We show (Thm. \ref{MurayamaShimakawaGroupoidIsEquivariantModuliStack} below)
that the Murayama-Shimakawa construction of equivariant classifying spaces
(reviewed in \cref{ConstructionOfUniversalEquivariantPrincipalBundles})
may be understood as a  natural topological groupoid model for the
more abstractly defined
proper-equivariant moduli stacks \eqref{OrbiSingularModuliStack}
and hence for their proper-equivariant classifying shapes
(Def. \ref{EquivariantClassifyingStack}).

\begin{lemma}[Geometric fixed loci of equivariant classifying spaces by Murayama-Shimakawa construction]
  \label{MurayamaShimakawaGroupoidFromSingularCohesion}
  $\,$

  \noindent
  In $\Topos = \SmoothInfinityGroupoids$,

  -- $
    G \acts \, \Gamma
    \;\in\;
    \begin{tikzcd}
    \Groups
    \left(
      \Actions{G}(\kTopologicalSpaces)
    \right)
    \ar[
      rr,
      "{
        \scalebox{0.7}{$
          \Groups\left(\Actions{G}(\ContinuousDiffeology)\right)
        $}
      }"
    ]
    &&
    \Groups
    \left(
      \Actions{G}(\SmoothInfinityGroupoids)
    \right)\!\!\!
    \end{tikzcd}
  $,

  -- $H \subset G$ a subgroup,

  \noindent
  there is a natural equivalence
  \vspace{-2mm}
  \begin{equation}
   \label{MurayamaShimakawaGroupoidFormulaDerivedFromSingularCohesion}
    \shape
    \,
    \Maps{\big}
      { \mathbf{B}H }
      { (\mathbf{B}\Gamma) \!\sslash\! G }
    _{\mathbf{B}G}
    \;\;\simeq\;\;
    \SingularSimplicialComplex
    \Big(
    \TopologicalRealization{}
    {
      \Maps{}
        { \mathbf{E}G }
        { \mathbf{B}\Gamma }
    }
    ^{H}
    \Big)
    \;\;\;
    \in
    \;
    \InfinityGroupoids
    \xhookrightarrow{ \Discrete }
    \SmoothInfinityGroupoids
      \end{equation}

  \vspace{-1mm}
  \noindent
  between the shapes of the
  slice mapping stack (Def. \ref{SliceMappingStack})
  and that of the
  $H$-fixed locus of the base of the Murayama-Shimakawa construction
  \eqref{QuotientCoprojectionOfUniversalEquivariantPrincipalBundle}.
\end{lemma}
\begin{proof}
First, for
$\TopologicalPatch \in \CartesianSpaces$
we have the following sequence of
natural equivalences in
$\InfinityGroupoids$,
where in the last lines
we are using Lem. \ref{DiffeologicalMappingGroupoidOutOfDiscreteIntoContinuousDiffeologicalGroupoid}
to notationally conflate $\ContinuousDiffeology(\Gamma)$ with $\Gamma$:
\vspace{-2mm}
$$
\hspace{-2mm}
  \def\arraystretch{1.6}
  \begin{array}{lll}
    \SliceMaps{}{\mathbf{B}G}
      { \mathbf{B}H }
      { \HomotopyQuotient{\mathbf{B}\Gamma}{G} }
    (\TopologicalPatch)
    &
    \;\simeq\;
    (
      \SmoothInfinityGroupoids
    )
     _{/\mathbf{B}G}
    \left(
      \TopologicalPatch
        \times
      \mathbf{B}H,
      \,
      (\mathbf{B}\Gamma) \!\sslash\! G
    \right)
    &
    \proofstep{ by Lem. \ref{PlotsOfSliceMappingStackAreSliceHoms} }
    \\
    &
    \;\simeq\;
    (
      \SmoothInfinityGroupoids
    )
      _{/\mathbf{B}H}
    \left(
      \TopologicalPatch
        \times
      \mathbf{B}H,
      \,
      (\mathbf{B}\Gamma) \!\sslash\! H
    \right)
    &
    \proofstep{ by Prop. \ref{SliceMappingStackIsStackOfSectionsOfPullbackBundle}}
    \\
    &
    \;\simeq\;
    {\SimplicialPresheaves(\CartesianSpaces)_{/\overline{W}H}}
    \left(
      \TopologicalPatch
        \times
      \overline{W}H,
      \,
      \left(
        \overline{W}\Gamma
          \times
        W H
      \right)/H
    \right)
    &
    \proofstep{ by Lem. \ref{PresentationOfHomotopyFixedPointSpacesOfEquivariantClassifyingStacks}}
    \\
    &
    \;\simeq\;
    \Actions{H}
    \big(
      \SimplicialPresheaves(\CartesianSpaces)
    \big)
    \big(
      \TopologicalPatch
        \times
      W H
      ,
      \,
      \overline{W}\Gamma
    \big)
    &
    \proofstep{ by Prop. \ref{QuillenEquivalenceBetweenBorelModelStructureAndSliceOverClassifyingComplex}}
    \\
    &
    \;\simeq\;
    \Actions{H}
    \left(
      \SimplicialPresheaves(\CartesianSpaces)
    \right)
    \big(
      \TopologicalPatch
        \times
      W G
      ,
      \,
      \overline{W}\Gamma
    \big)
    &
    \proofstep{ by Ex. \ref{EquivariantEquivalenceOfSimplicialUniversalPrincipalComplexes}}
    \\
    &
    \;\simeq\;
    \Actions{H}
    \big(
      \SimplicialPresheaves(\CartesianSpaces)
    \big)
    \big(
     \TopologicalPatch
       \times
     \ActionGroupoid{G}{G}
      ,\,
      \DeloopingGroupoid{\Gamma}
    \big)
    &
    \proofstep{ by Ex. \ref{UniversalSimplicialPrincipalComplexForOrdoinaryGroupG}}
    \\
    &
    \;\simeq\;
    \Big(
      \ContinuousDiffeology
      \,
      \Maps{\big}
        { \ActionGroupoid{G}{G} }
        { \DeloopingGroupoid{\Gamma} }
      ^H
    \Big)
    (\TopologicalPatch)
    &
    \proofstep{ as in \cite[Lem. A.5]{SS20OrbifoldCohomology}. }
      \end{array}
$$
By naturality (in $\TopologicalPatch$) of the composite equivalence, the
Yoneda lemma (Prop. \ref{InfinityYonedaLemma}) implies that
we have an equivalence
\vspace{-2mm}
\begin{equation}
  \label{MurayamaShimakawaGroupoidAsMappingStack}
  \SliceMaps{}{\mathbf{B}G}
    { \mathbf{B}G }
    { \HomotopyQuotient{\mathbf{B}\Gamma}{G} }
  \;\;
  \simeq
  \;\;
  \ContinuousDiffeology
  \,
  \Maps{}
    { \mathbf{E}G }
    { \mathbf{B}\Gamma }
  ^H
  \;\;\;
  \in
  \;
  \Groupoids(\DiffeologicalSpaces)
  \xhookrightarrow{\;\;}
  \SmoothInfinityGroupoids
  \,.
\end{equation}

\vspace{-1mm}
\noindent Now the claim is implied from the following composite:
\vspace{-2mm}
$$
\hspace{-2mm}
  \def\arraystretch{1.5}
  \begin{array}{lll}
       \shape
    \,
    \SliceMaps{}{\mathbf{B}G}
      { \mathbf{B}H }
      { \HomotopyQuotient{\mathbf{B}\Gamma}{G} }
    &
    \!\!\!\! \simeq
    \shape
    \,
    \ContinuousDiffeology
    \big(
    \Maps{}
      { \mathbf{E}G }
      { \mathbf{B}\Gamma }
    ^H
    \big)
    &
    \proofstep{ by \eqref{MurayamaShimakawaGroupoidAsMappingStack} }
    \\
    & \!\!\!\! \simeq
    \shape
    \,
    \ContinuousDiffeology
    \Big(\!\!
      \ActionGroupoid
        { \CrossedHomomorphisms(H,\, H \acts \, \Gamma) }
        { \Gamma }
   \!\! \Big)
    &
    \proofstep{
      by \eqref{FixedLociOfBaseOfMurayamaShimakawaConstructionInTermsOfCrossedHomomorphisms}
    }
    \\
    & \!\!\!\! \simeq
    \SingularSimplicialComplex
    \,
    \TopologicalRealization{\big}
    {
      \ActionGroupoid
        { \CrossedHomomorphisms(H, \, H \acts \, \Gamma) }
        { \Gamma }
    }
    &
    \proofstep{
      by Prop. \ref{SmoothShapeOfGoodSimplicialSpacesIsWeakHomotopyTypeOfTheirRealization}
      \&
      Prop.
      \ref{NervesOfActionGroupoidsOfWellPointedTopologicalGroupActionsAreGood}
    }
    \\
    &
    \!\!\!\!\simeq
    \SingularSimplicialComplex
    \TopologicalRealization{\Big}
    {
      \Maps{}
        { \mathbf{E}G }
        { \mathbf{B}\Gamma }
       ^H
    }
    &
    \proofstep{
      by \eqref{FixedLociOfBaseOfMurayamaShimakawaConstructionInTermsOfCrossedHomomorphisms}
    }
    \\
    &
    \!\!\!\!\simeq
    \SingularSimplicialComplex
    \Big(
    \TopologicalRealization{\big}
    {
      \Maps{}
        { \mathbf{E}G }
        { \mathbf{B}\Gamma }
    }^H
    \Big)
    &
    \proofstep{
       by Lem. \ref{TopologicalRealizationPreservesFiniteLimits}.
    }
  \end{array}
$$

\vspace{-2mm}
\noindent
Here the passage through the action groupoid of crossed homomorphisms
just serves to produce a homotopy equivalent topological groupoid whose
nerve is manifestly good, by Prop. \ref{NervesOfActionGroupoidsOfWellPointedTopologicalGroupActionsAreGood},
so that its plain topological realization
is seen to have the correct weak homotopy type.
\end{proof}

\begin{theorem}[Murayama-Shimakawa construction gives shape of equivariant moduli stacks]
  \label{MurayamaShimakawaGroupoidIsEquivariantModuliStack}
  For $G \in \Groups(\Sets)$
  and $G \acts \, \Gamma \,\in\, \Actions{G}\left( \Groups(\DTopologicalSpaces) \right)$
  with $\Gamma$ well-pointed (Ntn. \ref{WellPointedTopologicalGroup}),
  the proper-equivariant classifying space from
  Def. \ref{EquivariantClassifyingStack}
  is equivalent, under Elmendorf's theorem \eqref{ElmendorfTheorem},
  to the Muryama-Shimakawa construction (Ntn. \ref{MurayamaShimakawaConstruction})
  \vspace{-2mm}
  $$
    \EquivariantClassifyingShape{G}{\Gamma}
    \;\coloneqq\;
    \smoothrelativeG
    \,
    \shape
    \,
    \orbisingular
    \,
    \mathbf{B}(\Gamma \rtimes G)
    \;\;\;
      \simeq
    \;\;\;
    \SingularSimplicialComplex
    \,
    \TopologicalRealization{}
    {
      \Maps{}
        { \mathbf{E}G }
        { \mathbf{B}\Gamma }
    }^{(-)}
    \;\;\;
    \;\;
    \in
    \;\;
    G\InfinityGroupoids
    \xhookrightarrow{ \;G\Discrete\; }
    G\ModalTopos{\smooth}
    \xhookrightarrow{\; G\Orbi\Spaces \;}
    \GloballyEquivariant{G}(\ModalTopos{\smooth})_{\scalebox{.7}{$/\orbisingularG$}}
    \,\simeq\,
    \ModalTopos{/\orbisingularG}
    \,.
  $$
\end{theorem}
\begin{proof}
For $H \subset G$, we have the following natural equivalences
in $\InfinityGroupoids$:
\vspace{-2mm}
$$
  \def\arraystretch{1.8}
  \begin{array}{lll}
    \EquivariantClassifyingShape{G}{\Gamma}
    \;:\;
    G/H
    \;\longmapsto\;
    &
    \big(
      \smoothrelativeG
      \;
      \shape
      \;
      \orbisingular
      (\HomotopyQuotient{\mathbf{B}\Gamma}{G})
    \big)
    (G/H)
    &
    \proofstep{ by Def. \ref{EquivariantClassifyingStack} }
    \\
    &
    \;\simeq\;
    \smooth
    \,
    \SliceMaps{\big}{\mathbf{B}G}
      { \orbisingularH }
      {
        \shape
        \,
        \orbisingular
        \,
        (\HomotopyQuotient{\mathbf{B}\Gamma}{G})
      }
    &
    \proofstep{ by \eqref{TheGFixedLociFunctor} }
    \\
    & \;\simeq\;
    \shape
    \,
    \SliceMaps{}{\mathbf{B}G}
      { \mathbf{B}H }
      { \HomotopyQuotient{\mathbf{B}\Gamma}{G} }
    &
    \proofstep{ by Lem. \ref{FixedLociOfEquivariantClassifyingStacks} }
    \\
    &
    \;\simeq\;
    \SingularSimplicialComplex
    \,
    \TopologicalRealization{}
    {
      \Maps{}
      { \mathbf{E}G }
      { \mathbf{B}\Gamma }
    }^H
    &
    \proofstep{ by Lem. \ref{MurayamaShimakawaGroupoidFromSingularCohesion}. }
  \end{array}
$$

\vspace{-7mm}
\end{proof}

\noindent
{\bf Classification in Borel-equivariant cohomology and in proper-equivariant cohomology.}
We now formulate the classification result for
equivariant principal bundles in terms of
equivariant cohomology.

\begin{definition}[Borel-equivariant and proper-equivariant cohomology in a cohesive $\infty$-topos]
  \label{BorelEquivariantAndProperEquivariantCohomologyInCohesiveInfinityTopos}
  Given a singular-cohesive $\infty$-topos $\Topos$
  (Def. \ref{SingularCohesiveInfinityTopos}),

  \vspace{-2mm}
  \begin{enumerate}[{\bf (i)}]
  \setlength\itemsep{-2pt}

  \item
  for $G \acts \, A \,\in\, \Actions{G}(\ModalTopos{\smooth})$, we have that
  \vspace{-2mm}
  \begin{equation}
    \label{BorelEquivariantCohomology}
    H_G^0(X;\, A)
    \;\coloneqq\;
    \Truncation{0}
    \,
    \SlicePointsMaps{}{\mathbf{B}G}
      { \HomotopyQuotient{X}{G} }
      { \HomotopyQuotient{A}{G} }
  \end{equation}
  is the {\it Borel-equivariant cohomology} of $X$ with coefficients in $A$;

  \item
  for $\mathcal{A} \,\in\, \GEquivariant\ModalTopos{\smooth}$, we have that
  \vspace{-3mm}
  \begin{equation}
    \label{ProperEquivariantCohomology}
    H_{\scalebox{.7}{$\orbisingularG$}}(X;\, \mathcal{A})
    \;\coloneqq\;
    \Truncation{0}
    \,
    \SlicePointsMaps{\big}{\orbisingularG}
      { \orbisingular (\HomotopyQuotient{X}{G}) }
      { G\Orbi\Singularity(\mathcal{A}) }
    \;\;
    \overset{
      \mathclap{
        \mbox{\tiny \eqref{IncarnationsOfCocyclesOnCohesiveOrbispace}}
      }
    }{
      \simeq
    }
    \;\;
    \Truncation{0}
    \,
    \GEquivariant{\SliceTopos{\smooth}}
    \left(
      \FixedLoci{X}
      ,\,
      \mathcal{A}
    \right)
    \,,
  \end{equation}

  \vspace{-3mm}
\noindent  is the {\it proper equivariant cohomology} of $X$ with coefficients in $\mathcal{A}$.

  \end{enumerate}
  \vspace{-.1cm}

\end{definition}

The following Prop. \ref{ProperEquivariantCohomologySubsumesBorelEquivariantCohomology}
generalizes the corresponding classical statement for bare global homotopy theory
to general singular-cohesive homotopy theory:
\begin{proposition}[Proper-equivariant cohomology subsumes Borel-equivariant cohomology]
  \label{ProperEquivariantCohomologySubsumesBorelEquivariantCohomology}
  For

  -- $G \acts \, X \,\in\, \Actions{G}(\ModalTopos{\smooth,0})$,

  -- $G \acts \, A \,\in\, \Actions{G}(\ModalTopos{\smooth})$,

  \noindent
  Borel-equivariant cohomology
  with coefficients in $A$
  is in natural bijection with the
  proper-equivariant cohomology
  (Def. \ref{BorelEquivariantAndProperEquivariantCohomologyInCohesiveInfinityTopos})
  with coefficients in
  $\FixedLoci(A)$ (Def. \ref{ShapeOfGemetricFixedLoci}):
  \vspace{-1mm}
  $$
    H^0_G(X;\, A)
    \;\simeq\;
    H^0_{\!\!\scalebox{.7}{$\orbisingularG$}}\left(X;\, \FixedLoci(A) \right)
    \,.
  $$
\end{proposition}
\begin{proof}
We have the following sequence of natural equivalences:
\vspace{-2mm}
$$
  \def\arraystretch{1.4}
  \begin{array}{lll}
    \SlicePointsMaps{}{\mathbf{B}G}
      { \HomotopyQuotient{X}{G} }
      { \HomotopyQuotient{A}{G} }
    &
    \;\simeq\;
    \SlicePointsMaps{\big}{\orbisingularG}
      {\! \orbisingular (\HomotopyQuotient{X}{G}) }
      { \orbisingular (\HomotopyQuotient{A}{G}) }
    &
    \proofstep{ by Ex. \ref{MorphsimsBetweenOrbisingularizationsOfSmoothObjects} }
    \\
    &
    \;\simeq\;
    \SlicePointsMaps{\big}{\orbisingularG}
      {
        G\Orbi\Space
        \left(
          \FixedLoci(X)
        \right)
      }
      { \orbisingular (\HomotopyQuotient{A}{G}) }
    &
    \proofstep{ by Prop. \ref{OrbiSpaceIncarnationOfGSpaceIsOrbisingularizationOfHomotopyQuotient} }
    \\
    &
    \;\simeq\;
    \GEquivariant\ModalTopos{\smooth}
    \big(
      {
        \FixedLoci(X)
      }
      ,\,
      {
        G\Orbi\Smooth
        \left(
          \orbisingular (\HomotopyQuotient{A}{G})
        \right)
      }
    \big)
    &
    \proofstep{ by \eqref{GFixedLocFunctor} }
    \\
    &
    \;\simeq\;
    \GEquivariant\ModalTopos{\smooth}
    \big(
      {
        \FixedLoci(X)
      }
      ,\,
      {
        \FixedLoci(A)
      }
    \big)
    &
    \proofstep{ by Def. \ref{ShapeOfGemetricFixedLoci}. }
  \end{array}
$$

\vspace{-2mm}
\noindent This implies the claim by passage to sets of connected components.
\end{proof}

\begin{remark}[Interpretation]
  Specialized to pure shape coefficients
  (which is essentially the only case considered in traditional literature
   on equivariant cohomology),
  Prop. \ref{ProperEquivariantCohomologySubsumesBorelEquivariantCohomology}
  means that Borel-equivariant cohomology
  only captures those proper-equivariant coefficients which are systems of
  homotopy fixed loci $\FixedLoci(\shape\, A)$
  but not coefficients
  $\shape \, \FixedLoci(A)$
  which are systems of shapes of geometric fixed loci.
\end{remark}

\begin{example}[First non-abelian Borel-equivariant cohomology as coarse proper-eqjuivariant cohomology]
\label{FirstNonAbelianBorelEquivariantCohomologyAsCoarseProperEquivariantCohomology}
Specialized to delooped
pure shape coefficients
$A \,\coloneqq\, \mathbf{B} \mathcal{G} \,\simeq\, B  \,\shape \, \mathcal{G}$
for any $G \acts \, \mathcal{G} \,\in\, \Actions{G}(\ModalTopos{\smooth})$,
Prop. \ref{ProperEquivariantCohomologySubsumesBorelEquivariantCohomology}
says that first non-abelian Borel-equivariant cohomology
is equivalently proper equivariant cohomology with
coefficients in the coarse equivariant classifying spaces from
Rem. \ref{ProperEquivariantClassifyingShapesMayDependOnStructureGroupBeyondItsShape}:
\vspace{-2mm}
$$
  H^1_G(X;\, \shape \, \Gamma)
  \;=\;
  H^0_G(X;\, B \Gamma)
  \;\simeq\;
  H^1_{\scalebox{.7}{$\orbisingularG$}}
  \left(
    X
    ;\,
    \EquivariantClassifyingShape{G}{(\shape\,\Gamma)}
  \right)
  \,.
$$
\end{example}

Finally, we obtain
the classification of stable equivariant principal bundles
(Thm. \ref{BorelClassificationOfEquivariantBundlesForResolvableSingularitiesAndEquivariantStructure})
in terms of proper equivariant cohomology
and generalized from resolvable to cover-resolvable singularities:

\begin{theorem}[Classification of equivariant bundles with truncated structure over orbifolds with cover-resolvable singularities]
  \label{ProperClassificationOfEquivariantBundlesForResolvableSingularitiesAndEquivariantStructure}
  For

  -- $G \,\in\, \Groups(\FiniteSets)_{\CoverResolvable}$
  (Ntn. \ref{ResolvableOrbiSingularities}),

  -- $\Gamma \,\in\, \Groups(\DHausdorffSpaces)$
  of truncated classifying shape (Ntn. \ref{CohesiveGroupsWithTruncatedClassifyingShape})
  and with a notion of stable equivariant bundles (Ntn. \ref{StableEquivariantBundles}),

  \noindent
  isomorphism classes of
  stable topological
  $G$-equivariant $\Gamma$-principal bundles over any
  $G \acts \, \SmoothManifold \,\in\, \Actions{G}(\SmoothManifolds)$
  are classified equivalently by

  \vspace{-3mm}
  \begin{enumerate}[{\bf (i)}]
  \setlength\itemsep{-3pt}

  \item
  proper-equivariant cohomology
  (Def. \ref{BorelEquivariantAndProperEquivariantCohomologyInCohesiveInfinityTopos})
  with coefficients in
  the stable proper-equivariant classifying shape
  $\EquivariantClassifyingShape{G}{\Gamma}$
  (Def. \ref{StableEquivariantClassifyingShape});

  \item
  proper-equivariant cohomology
  (Def. \ref{BorelEquivariantAndProperEquivariantCohomologyInCohesiveInfinityTopos})
  with coefficients in
  the coarse equivariant classifying shape
  $\EquivariantClassifyingShape{G}{(\shape \, \Gamma)}$
  (Rem. \ref{ProperEquivariantClassifyingShapesMayDependOnStructureGroupBeyondItsShape});

  \item
  first Borel-equivariant cohomology with coefficients in $\shape\,\Gamma$:
\end{enumerate}
\vspace{-2mm}
  $$
    \hspace{-1mm}
    \overset{
      \mathclap{
      \raisebox{6pt}{
        \tiny
        \color{darkblue}
        \bf
        \def\arraystretch{.9}
        \begin{tabular}{c}
          isomorphism classes of
          \\
          stable equivariant principal topological bundles
        \end{tabular}
      }
      }
    }{
    \IsomorphismClasses{
      \EquivariantPrincipalFiberBundles{G}{\Gamma}
      (\DTopologicalSpaces)_{\SmoothManifold}^{\stable}
    }
    }
        \simeq
        \overset{
      \mathclap{
      \raisebox{4pt}{
        \tiny
        \color{darkblue}
        \bf
        \begin{tabular}{c}
          Borel-equivariant cohomology
        \end{tabular}
      }
      }
    }{
      H^1_G
      (X;\, \shape\,\Gamma)
      \;
      =
      \;
      H^0_G
      (X;\, B \Gamma)
    }
    \;\;
    \simeq
    \;\;
    \overset{
      \raisebox{4pt}{
        \tiny
        \color{darkblue}
        \bf
        \begin{tabular}{c}
          proper-equivariant cohomology
        \end{tabular}
      }
    }{
      H^0_{\scalebox{.7}{$\orbisingularG$}}
      \big(X;\, \EquivariantClassifyingShape{G}{(\shape\, \Gamma)}\big)
      \;
      \simeq
      \;
      H^0_{\scalebox{.7}{$\orbisingularG$}}
      \big(
        X;\,
        (\EquivariantClassifyingShape{G}{\Gamma})^{\stable}
      \big)
    }
    \,.
  $$
\end{theorem}

\begin{proof}
  First, consider the case where $G$-singularities are resolvable.
  Then the first equivalence is
  Thm. \ref{BorelClassificationOfEquivariantBundlesForResolvableSingularitiesAndEquivariantStructure}
  in view of Def. \ref{BorelEquivariantAndProperEquivariantCohomologyInCohesiveInfinityTopos},
  the second equivalence follows by
  Ex. \ref{FirstNonAbelianBorelEquivariantCohomologyAsCoarseProperEquivariantCohomology}
  of Prop. \ref{ProperEquivariantCohomologySubsumesBorelEquivariantCohomology},
  and then the third by
  Prop. \ref{EquivariantClassifyingShapeOfTruncatedTopologicalGroupsCoincidesWithThatOftheirShape}.

  It remains to generalize this statement to the case that $G$-singularities may not
  be resolvable, but there exists a surjective homomorphism
  $p : \!\! \begin{tikzcd} \widehat{G} \ar[r, ->>] &[-12pt]  G \end{tikzcd}\!\!$
  \eqref{ADiscreteGroupEpimorphismInAnInfinityTopos} such that $\widehat{G}$
  has resolvable singularities.
  To this end, notice that, for $G \acts \, X \,\in\,
  \Actions{G}(\ModalTopos{0})$, we have the natural monomorphism
  \vspace{-2mm}
\begin{equation}
  \label{NaturalMonomorphismInProofOfClassificationTheorem}
  \begin{aligned}
    G\InfinityGroupoids
    \big(
      { X^{(-)} }
      ,\,
      { \EquivariantClassifyingShape{G}{\Gamma} }
    \big)
    &
    \;\xhookrightarrow{\quad} \;\;
    G\InfinityGroupoids
    \big(
      { X^{(-)} }
      ,\,
      {
        p_\ast \EquivariantClassifyingShape{\widehat{G}}{\Gamma}
      }
    \big)
    &
    \proofstep{
      by Lem. \ref{BaseChangeOfEquivariantClassifyingShapesAlongCoverings}
      with Ex. \ref{MonomorphismsOfInfinityGroupoids}
    }
    \\
    & \;\;\simeq \quad
    \widehat{G}\InfinityGroupoids
    \big(
      { p^\ast X^{(-)} }
      ,\,
      {
        \EquivariantClassifyingShape{\widehat{G}}{\Gamma}
      }
    \big)
    &
    \proofstep{
      by Lem. \ref{BaseChangeInEquivariantHomotopyTheoryAlongCoveringsOfEquivaranceGroup}
      with \eqref{AdjunctionAndHomEquivalence}.
    }
  \end{aligned}
\end{equation}

\vspace{-2mm}
\noindent By naturality of the base change $(B p)^\ast$, this monomorphism
is compatible with that from Prop. \ref{EquivariantBundlesUnderExtensionsPfTheEquivarianceGroup},
making a commuting square
\vspace{-2mm}
$$
  \begin{tikzcd}
    \IsomorphismClasses{
      \EquivariantPrincipalBundles{G}{\Gamma}_X
    }
    \ar[
      d,
      hook,
      "\mbox{
        \tiny
        Prop. \ref{EquivariantBundlesUnderExtensionsPfTheEquivarianceGroup}
      }"
    ]
    \ar[rr]
    &&
    \Truncation{0}
    G\InfinityGroupoids
    \left(
      { X^{(-)} }
      ,\,
      { \EquivariantClassifyingShape{G}{\Gamma} }
    \right)
    \ar[
      d,
      hook,
      "\mbox{
        \tiny
         Lem. \ref{BaseChangeOfEquivariantClassifyingShapesAlongCoverings}
         ,
         \eqref{NaturalMonomorphismInProofOfClassificationTheorem}
      }"
    ]
    \ar[r,-, shift left=1pt]
    \ar[r,-, shift right=1pt]
    &[-16pt]
    H^0_{\scalebox{.7}{$\orbisingularG$}}
    \left(
      X
      ,\,
      \EquivariantClassifyingShape{G}{\Gamma}
    \right)
    \\
    \IsomorphismClasses{
      \EquivariantPrincipalBundles{\widehat{G}}{\Gamma}_{X}
    }
    \ar[rr, "\sim"]
    &&
    \Truncation{0}
    \widehat{G}\InfinityGroupoids
    \left(
      { p^\ast X^{(-)}   }
      ,\,
      { \EquivariantClassifyingShape{\widehat{G}}{\Gamma} }
    \right)
    \ar[r,-, shift left=1pt]
    \ar[r,-, shift right=1pt]
    &[-16pt]
    H^0_{\scalebox{.7}{$\orbisingularAny{\widehat{G}}$}}
    \left(
      X
      ,\,
      \EquivariantClassifyingShape{\widehat{G}}{\Gamma}
    \right)
    \,.
  \end{tikzcd}
$$

\vspace{-2mm}
\noindent
Here the bottom map is a bijection by the previous argument,
which applies to $\widehat{G}$ by assumption; therefore the
top map is also a bijection, being its restriction to the
$G$-equivariant bundles among the $\widehat{G}$-equivariant bundles
along $(B p)^\ast$.
\end{proof}

\newpage

\part{Examples and applications}

We close by briefly indicating examples and applications of the
machinery of singular-cohesive homotopy theory in general and in view of the
classification result, Thm. \ref{ProperClassificationOfEquivariantBundlesForResolvableSingularitiesAndEquivariantStructure}, in particular.
Further discussion is relegated to
\cite{TwistedEquivariantChernCharacter}\cite{TwistedEquivariantDifferentialCohomology}.

\medskip

\section{Classification results}

\begin{example}[Equivariant bundles with truncated compact Lie structure group]
\label{TruncatedStructureGroupsAndTheirEquivariantClassificationResults}
Via Thm. \ref{ForCompactTruncatedStructureLieGroupResolutionIsInjection},
Thm. \ref{ProperClassificationOfEquivariantBundlesForResolvableSingularitiesAndEquivariantStructure}
subsumes the following special cases known in the literature
(here over smooth $G$-manifolds $G \acts \, \SmoothManifold$ with
finite equivariance group $G$ satisfying the conditions of Ntn. \ref{ResolvableOrbiSingularities}):

\vspace{-.1cm}
\begin{enumerate}[{\bf (i)}]

\vspace{-.2cm}
\item
a {\bf discrete group}
$\Gamma \,\in\, \Groups(\Sets) \xhookrightarrow{\;} \Groups(\kTopologicalSpaces)$.

In this case $\Gamma \,\simeq\, \shape \,\Gamma \,\simeq\, K(\Gamma,0)$
is an Eilenberg-MacLane space, hence truncated, and $B \Gamma$ is classifying,
by
Thm. \ref{ClassificationOfPrincipalBundlesAmongPrincipalInfinityBundles}.
Hence $\Gamma$ has truncated classifying shape, according to
Ntn. \ref{CohesiveGroupsWithTruncatedClassifyingShape}.

The statement of Thm. \ref{ProperClassificationOfEquivariantBundlesForResolvableSingularitiesAndEquivariantStructure}
in this case recovers \cite[Thm. 5]{May90}.

Notice that $\Gamma$-principal bundles for $\Gamma$ a discrete group
serve as twists for ordinary cohomology
(``local coefficient bundles'', e.g. \cite[Def. 3.1, Lem. 4.2]{BFGM03}),
see \cite[Ex. 2.32]{FSS20CharacterMap}.
In particular, in the case $\Gamma = \ZTwo$ these are the twists
of integral cohomology \cite{GS-Deligne}, reflecting the fact
that $\ZTwo$ is the group of units (see \cite[Ex. 2.37]{FSS20CharacterMap})
of the integral
Eilenberg-MacLane spectrum:
${\rm GL}_1(H \mathbb{Z}) \simeq \{\pm 1\} \simeq \ZTwo\,.$

Therefore, equivariant $\Gamma$-principal bundles in this case
constitute twists of equivariant ordinary cohomology.
Despite the relative simplicity of this example,
the classification of these equivariant twists is already quite rich,
given by the
non-abelian group cohomology of crossed homomorphisms
with coefficients in $\Gamma$ (Thm. \ref{EquivariantHomotopyGroupsOfEquivariantClassifyingSpaces}).

\vspace{-.2cm}
\item
a {\bf compact abelian Lie group}
$\Gamma \,\in\, \AbelianGroups(\Compact\SmoothManifolds)$.

Since compact abelian Lie groups are, up to isomorphism,
direct products
$\Gamma \,=\, \TorusGroup{r} \times $
of finite abelian groups $A$ with $r$-dimensional tori ($r \,\in\, \mathbb{N}$)
hence with $k$-fold direct products
of the circle group with itself, their shape is an Eilenberg-MacLane space
\vspace{-2mm}
$$
  \Shape
  \,
  \left(
    \TorusGroup{r} \times A
 \right)
  \;\;
  \simeq
  \;\;
  K(A, 0) \times K(\mathbb{Z}^r, 1)\;,
$$
and hence truncated.
Moreover, $B \Gamma$ is classifying,
by
Thm. \ref{ClassificationOfPrincipalBundlesAmongPrincipalInfinityBundles}.

The statement of Thm. \ref{ProperClassificationOfEquivariantBundlesForResolvableSingularitiesAndEquivariantStructure}
in this case recovers \cite[Lem. 1, Thm. 2]{LMS83}\cite[Thm. 3, Thm. 10]{May90}.

\vspace{-.2cm}
\item
an {\bf extension of a discrete group $K$ by a torus}
$\TorusGroup{r} \hookrightarrow \Gamma \twoheadrightarrow K$.

This is the evident joint generalization of the above cases:
The shape of $B \Gamma$ in this case is a homotopy 2-type,
now possibly with a non-trivial k-invariant.

The statement of
Thm. \ref{ProperClassificationOfEquivariantBundlesForResolvableSingularitiesAndEquivariantStructure}
in this case,
restricted to the point and for trivial $G$-action on $\Gamma$,
recovers \cite[Thm. 1.2]{Rezk16}.
\end{enumerate}

\noindent
Namely, and in more generality:

Let $\Gamma$ be a truncated compact Lie group,
hence a semidirect product
$\Gamma \,\simeq\, \TorusGroup{r} \rtimes K$
of a torus group $\TorusGroup{r} \,=\, (\CircleGroup)^{\times^r}$
with a finite group $K$;
and $G$ a finite equivariance group with resolvable singularities
(Ntn. \ref{ResolvableOrbiSingularities}),
assume that its action on $\Gamma$ is such that
$$
  (
    \TorusGroup{r} \rtimes K
  )
    \rtimes
  G
  \;\;
  \simeq
  \;\;
  \TorusGroup{r}
    \rtimes
  (
    K
      \rtimes
    G
  )
  \;\;\;
  \in
  \;
  \Groups(\SmoothInfinityGroupoids)
  \,.
$$
(The special case of this assumption which is considered in all the above references
is that the $G$-action on $\Gamma$ is trivial. On the opposite extreme, the
assumption is satisfied for $K = 1$ and general $G$-action on $\TorusGroup{r}$.)

With this assumption,
Thm. \ref{ForCompactTruncatedStructureLieGroupResolutionIsInjection}
says that  all $G$-equivariant bundles with
structure group $\Gamma \rtimes G$ are blowup-stable. Therefore the
orbi-smooth Oka principal (Thm. \ref{OrbiSmoothOkaPrinciple}) says,
for all $G \acts \, \SmoothManifold \,\in\, \Actions{G}(\SmoothManifolds)$,
that we have an equivalence
$$
  \shape
  \,
  \SliceMaps{\big}{\mathbf{B}G}
    {
      \HomotopyQuotient{ \SmoothManifold }{ G }
    }
    {
      \HomotopyQuotient
        { \mathbf{B}( \TorusGroup{r} \rtimes K ) }
        { G }
    }
  \;\simeq\;
  \SliceMaps{\big}{B G}
    {
      \HomotopyQuotient{ \shape \, \SmoothManifold }{ G }
    }
    {
      \HomotopyQuotient
        { B( \TorusGroup \rtimes K ) }
        { G }
    }
  \;\;\;
  \in
  \;
  \InfinityGroupoids
$$
and Thm. \ref{BorelClassificationOfEquivariantBundlesForResolvableSingularitiesAndEquivariantStructure}
says that this restricts on connected components to the classification theorem
for equivariantly locally trivial
$G$-equivariant $G \acts \, \Gamma$-principal bundles:
$$
  \IsomorphismClasses{
    \EquivariantPrincipalFiberBundles{G}{ (G \acts \, \Gamma) }_{\SmoothManifold}
  }
  \;\;
  \simeq
  \;\;
  \Truncation{0}
    \SliceMaps{\big}{BG}
      { \HomotopyQuotient{ \shape \SmoothManifold }{ G } }
      { \HomotopyQuotient{ B \Gamma }{ G } }
  \;\;
  =:
  \;\;
  H^1_G
  \big(
    \SmoothManifold
    ,\,
    \Gamma_G
  \big)
  \,.
$$

\end{example}

\medskip

\begin{example}[Classification of equivariant $\PUH$-bundles]
 \label{EquivariantBundlesServingAsGeoemtricTwistsOfEquivariantKTheory}
Equivariant bundles with structure group the infinite projective groups, from
Ex. \ref{ProjectiveUnitarGroupOnAHilbertSpace},
are meant to serve as geometric twists

\vspace{-2mm}
\begin{itemize}
\setlength\itemsep{-2pt}
\item[--] for equivariant $\mathrm{KU}$-theory
(going back to \cite{AdemRuan03}\cite[\S 6]{AtiyahSegal04}),

\item[--]  for $\mathrm{KO}$-theory (going back to \cite{DonovanKaroubi70},
the modern perspective is in \cite{AGG14},
considered via a variety of explicit models in
 \cite{Ros}\cite{MMS}\cite{HJ}\cite{GS-KO}\cite{GS-TwKO}\cite{GomiYamashita}),

\item[--]  and for their unification in $\mathrm{KR}$-theory \cite{Atiyah66}
(considered via a variety of explicit models in
\cite{Moutou11}\cite{Moutou12} \cite{Moutou13}\cite[\S 7]{FreedMoore13}\cite{HMSV16}\cite{Gomi17}).
\end{itemize}

\noindent
So consider the above classification result, in turn, for:

\medskip

\vspace{-.1cm}
\begin{enumerate}[{\bf (i)}]

\vspace{-.2cm}
\item
the {\bf projective unitary group on a Hilbert space}
$\Gamma \;=\; \PUH$  \eqref{TheGroupPUH}.

This has the truncated shape of an Eilenberg-MacLane space
$\shape \, \PUH \,\simeq\, K(\mathbb{Z},2)$ (by Ex. \ref{HomotopyFiberSequenceOfTheProjectiveUnitaryGroup}),
which is classifying,
by Thm. \ref{ClassificationOfPrincipalBundlesAmongPrincipalInfinityBundles}.

Via the characterization of blowup-stable projective bundles from
Thm. \ref{BlowupStabilityOfADEEquivariantPUHPrimePrincipalBundles}
for ADE-equivariance groups,
the statement of Thm. \ref{ProperClassificationOfEquivariantBundlesForResolvableSingularitiesAndEquivariantStructure}
in this case recovers
\cite[Prop. 6.3]{AtiyahSegal04}\cite[Thm. 3.8]{BEJU12}\cite[15.17]{LueckUribe14}
(see also \cite[Cor. 2.41]{TXLG04}),
and the statement of
Thm. \ref{EquivariantHomotopyGroupsOfEquivariantClassifyingSpaces}
in this case recovers \cite[Thm. 1.10]{BEJU12}.

\vspace{-.2cm}
\item
the
{\bf projective unitary group on a graded Hilbert space}
$\Gamma \;=\; \GradedPUH$ \eqref{TheGroupGradedPUH}.

Via the characterization of blowup-stable projective graded bundles from
Thm. \ref{BlowupStabilityForADEEquivariantPUHSemidirectProductZTwoPrincipalBundles}
for ADE-equivariance groups,
this special case of Thm. \ref{ProperClassificationOfEquivariantBundlesForResolvableSingularitiesAndEquivariantStructure}
gives the classification of
twists of equivariant K-theory in degrees 1 and 3 combined.

\vspace{-.2cm}
\item
the
{\bf projective unitary group with complex conjugation action}
$\ZTwo \acts \, \Gamma \;=\; \ZTwo \acts \, \PUH$ or $= \ZTwo \acts \, \GradedPUH$.

Via Thm. \ref{BlowupStabilityForADEEquivariantPUHSemidirectProductZTwoPrincipalBundles},
the statement of Thm. \ref{ProperClassificationOfEquivariantBundlesForResolvableSingularitiesAndEquivariantStructure}
in this case
generalizes the previous two cases to Atiyah's ``Real'' K-theory
($\mathrm{KR}$-theory) \cite{Atiyah66};
see also Ex. \ref{EquivariantBundleGerbes}, whence we
may refer to equivariant $\ZTwo \acts \PUH$-principal bundles
as {\it Real-equivariant} projective bundles.
Their classification result in this generality seems to be new.

\end{enumerate}
\vspace{-.2cm}

\noindent
Namely:

For $G \,\subset\, \SpOne$ (Ex. \ref{ADEGroupsHaveSphericalSpaceForms})
equipped with a homomorphism
$
  G \xrightarrow{ct}  \ZTwo
$
and for $G \acts \, \SmoothManifold \,\in\, \Actions{G}(\SmoothManifolds)$,
the isomorphism classes of
blowup-stable
(Def. \ref{BlowupStableADEEquivariantPUHPrimeRtimesZTwoBundles},
Thm. \ref{BlowupStabilityForADEEquivariantPUHSemidirectProductZTwoPrincipalBundles})
$G$-equivariant
$\ZTwo \acts \, \GradedPUHPrime$-principal bundles on $\SmoothManifold$
\eqref{BlowupStabilityForADEEquivariantPUHSemidirectProductZTwoPrincipalBundles}
are in natural bijection with the Borel-equivariant
integral cohomology of $\SmoothManifold$
with coefficients in degrees 1 and
conjugation-local coefficients in degree 3:
\begin{equation}
  \label{ClassificationStatementForStableEquivariantProjectiveGradedUnitaryBundles}
  \begin{tikzcd}
  \overset{
    \mathclap{
    \raisebox{5pt}{
      \tiny
      \color{darkblue}
      \bf
      \def\arraystretch{.9}
      \begin{tabular}{c}
        isomorphism classes of stable
        \\
        Real-equivariant projective graded bundles
      \end{tabular}
    }
    }
  }{
    \IsomorphismClasses{
      \EquivariantPrincipalBundles{G}{\GradedPUHPrime}_{\SmoothManifold}^{\stable}
    }
  }
  \ar[rr, "{\sim}"]
  \ar[
    d,
    "{\sim}"{sloped},
    "{
      \mbox{
        \tiny
        \begin{tabular}{c}
          Thm. \ref{ProperClassificationOfEquivariantBundlesForResolvableSingularitiesAndEquivariantStructure}
          \\
          via
          \\
          Thm. \ref{BlowupStabilityForADEEquivariantPUHSemidirectProductZTwoPrincipalBundles}
        \end{tabular}
      }
    }"{swap}
  ]
  &&
  \overset{
    \mathclap{
    \raisebox{6pt}{
      \tiny
      \color{darkblue}
      \bf
      \begin{tabular}{c}
        Borel-equivariant integral cohomology
        \\
        in degree 1 and conjugation-local degree3
      \end{tabular}
    }
    }
  }{
    H^1_G
    \big(
      \SmoothManifold
      ,\,
      \Integers
    \big)
    \,\times\,
    H^3_G
    \big(
      \SmoothManifold
      ,\,
      \Integers_{ct}
    \big)
  }
  \\
  \Truncation{0}
  \,
  \SliceMaps{\Big}{B \ZTwo}
    { \shape \, \HomotopyQuotient{ \SmoothManifold }{ G } }
    {
      \shape
      \,
      \HomotopyQuotient
        { \mathbf{B} \GradedPUHPrime }
        { \ZTwo }
    }
  \ar[
    rr,
    "{\sim}",
    "{
      \mbox{
        \tiny
        Lem. \ref{ShapeOfGradedProjectiveGroups}
      }
    }"{swap}
  ]
  &&
  \Truncation{0}
  \,
  \SliceMaps{\Big}{B \ZTwo}
    { \shape \, \HomotopyQuotient{ \SmoothManifold }{ G } }
    {
      B \ZTwo
      \,\times\,
      \HomotopyQuotient
        { B^3 \Integers }
        { \ZTwo }
    }
    \mathrlap{\,.}
  \ar[u,-, shift left=1pt]
  \ar[u,-, shift right=1pt]
  \end{tikzcd}
\end{equation}
(For trivial equivariance group $G = 1$ this subsumes
the classification of \cite[Thm. 3.6]{Parker88}\cite[Prop. 2.3]{AtiyahSegal04}.)

Moreover,
Thm. \ref{EquivariantHomotopyGroupsOfEquivariantClassifyingSpaces}
implies that
the equivariant homotopy groups of the equivariant classifying shape are
given by the integral group cohomology of the equivariance group with
conjugation-local coefficients:
\begin{equation}
  \label{EquivariantHomotopyGroupsOfClassifyingShapeOfStableEquivariantPUHBundles}
  \pi_k^H
  \big(
    \EquivariantClassifyingShape
      {G}
      { ( \GradedPUH )^{\stable} }
  \big)
  \;\;
  \simeq
  \;\;
  \Truncation{0}
  \,
  \SliceMaps{\big}{B H}
    { B H }
    { B^{3-k} \Integers_{\conjugation} }
  \;\;
  \simeq
  \;\;
  H^{3-k}_{\mathrm{Grp}}(H, \Integers_{\conjugation})
  \,.
\end{equation}
See also Ex. \ref{OrbiSmoothOkaPrincipleForPUHCoefficientsOverThePoint} for
how exactly this encodes the space of projective intertwiners of stable projective
$G$-representations.

\end{example}

\medskip

\section{Twisted equivariant cohomology}
\label{TwistedEquivariantCohomology}

We use the above singular-cohesive homotopy theory to obtain
(in Thm. \ref{TwistedEDKtheorem} below)
a twisted (parameterized/sliced) generalization of the classical Elmendorf-Dwyer-Kan theorem (reviewed as Prop. \ref{ElmendorfDwyerKanTheorem} below).
Before proving the twisted EDK theorem, we now explain how it serves as the foundation for a general notion of {\it twisted} \& {\it equivariant} generalized differential cohomology theory (subsuming the fairly well-known example of twisted equivariant K-theory. Ex. \ref{TwistedEquivariantKTheory} below) which we will lay out in \cite{TwistedEquivariantChernCharacter}\cite{TwistedEquivariantDifferentialCohomology}:

\medskip

\noindent
{\bf The notion of twisted equivariant cohomology.}
In asking for the proper conceptual nature for twisted \& equivariant generalized cohomology, the first step must be understanding the notion of ``generalized cohomology'' as such. The concpet  that traditionally goes by this ambitious name -- which here we rather refer to as {\it Whitehead-generalized cohomology} theories (i.e. those represented by bare spectra, for review and pointers see \cite[Ex. 1.0.13]{FSS20CharacterMap}) -- is an important example but not nearly general enough, as it fails to subsume even common generalizations such as sheaf cohomology or non-abelian cohomology, let alone the twisted equivariant differential non-abelian cohomology that is needed in many applications envisioned notably in theoretical physics.

\medskip
But some reflection shows
(see \cite[p. 6]{SS20OrbifoldCohomology})
that a notion of generalized cohomology with a good chance to live up to its name is: Hom-sets in the homotopy categories of $\infty$-toposes (as in \cref{ToposTheory}). Concretely, for $\Topos$ an $\infty$-topos and $\mathcal{A} \,\in\, \Topos$ any object, the generalized cohomology with coefficients in $\mathcal{A}$ is simply the homotopy-classes of maps into $\mathcal{A}$:
\begin{equation}
  \label{GeneralNonabelianCohomology}
  \overset{
    \mathclap{
    \raisebox{3pt}{
      \tiny
      \color{darkblue}
      \bf
      truly generalized cohomology
    }
    }
  }{
  H^0(-;\, \mathcal{A})
  \;:=\;
  \pi_0 \PointsMaps{}{-}{\mathcal{A}}
  }
  \,:\,
  \HomotopyCategory(\Topos)^{\mathrm{op}}
  \xrightarrow{\phantom{---}}
  \Sets
  \,.
\end{equation}
The Giraud-Rezk-Lurie axioms
which characterize $\infty$-toposes
(Prop. \ref{InfinityGiraudTheorem})
guarantee that these hom-sets have the main properties expected of (non-abelian) cohomology sets. In particular, if $\mathcal{A}$ is pointed-connected, it follows
(Prop. \ref{LoopingAndDeloopingEquivalence})
that $\mathcal{A} \,\simeq\, \mathbf{B}\mathcal{G}$ for an $\infty$-group object $\mathcal{G} \,\in\, \Groups(\Topos)$ and that
\begin{equation}
  \label{FirstNonabelianGeneralizedCohomology}
  \overset{
    \mathclap{
    \raisebox{3pt}{
      \tiny
      \color{darkblue}
      \bf
      first non-abelian generalized cohomology
    }
    }
  }{
  H^1(-;\, \mathcal{G})
  \,:=\,
  H^0(-;\, \mathbf{B}\mathcal{G})
  }
  \,:\,
  \HomotopyCategory(\Topos)^{\mathrm{op}}
  \xrightarrow{\phantom{---}}
  \Sets
\end{equation}
classifies $\mathcal{G}$-principal $\infty$-bundles in $\Topos$ (Prop. \ref{GroupsActionsAndFiberBundles}), in grand generalization of the classical statement from Chern-Weil theory, for $G$ a (compact) Lie group
(review in \cite[\S 6]{FSS20CharacterMap}).

\medskip
Remarkably, both {\it twisting} as well as  {\it equivariance} are native to this generalized notion
of generalized cohomology \eqref{GeneralNonabelianCohomology}, and in a unified manner: Indeed, by the {\it fundamental theorem of $\infty$-topos theory} (Prop. \ref{SliceInfinityTopos})
the {\it slice $\infty$-category} $\SliceTopos{\mathcal{B}}$
over any $\mathcal{B} \,\in\, \Topos$
is again an $\infty$-topos and its intrinsic cohomology \eqref{GeneralNonabelianCohomology}
on any $(\mathcal{X}, \tau) \,\in\, \SliceTopos{\mathcal{B}}$
with
{\it local coefficients} $(\mathcal{A}, p_{\mathcal{A}}) \,\in\, \SliceTopos{\mathcal{B}}$ is
the {\it $\tau$-twisted} cohomology of $\mathcal{X}$ :
$$
  H^\tau
  \big(
    \mathcal{X}
    ;\,
    p_{\mathcal{A}}
  \big)
  \;=\;
  H^0
  \big(
    (\mathcal{X},\tau)
    ,\,
    (\mathcal{A}, p_{\mathcal{A}})
  \big)
  \,,
$$
with its correct functoriality in pairs consisting of a domain object and a twist:
\begin{equation}
  \label{TwistedgeneralizedCohomology}
  \overset{
    \raisebox{3pt}{
      \tiny
      \color{darkblue}
      \bf
      twisted generalized cohomology
    }
  }{
  H^\bullet
  \big(
    -
    ;\,
    p_{\mathcal{A}}
  \big)
  \;:=\;
  H^0
  \big(
    (-,\bullet)
    ;\,
    (\mathcal{A}, p_{\mathcal{A}})
  \big)
  }
  \;:\;
  \HomotopyCategory
  \big(
    \SliceTopos{\mathcal{B}}
  \big)^{\mathrm{op}}
  \xrightarrow{\phantom{---}}
  \Sets
  \,.
\end{equation}

The traditional ``abelian'' namely Whitehead-generalized twisted cohomology theories, i.e., those represented by bundles of spectra, fit into this more general generalized notion of twisted cohomology \eqref{TwistedgeneralizedCohomology} most neatly, due to the remarkable fact that the $\infty$-category of bundles of spectra in any $\infty$-topos $\Topos$ is itself again an $\infty$-topos, the {\it tangent $\infty$-topos} $T\Topos$ (\cite[\S 35.5]{Joyal08Logoi}\cite[\S 4.1.2]{dcct}\cite[\S 7.3]{Lurie17}\cite{BraunackMayer21}),
which is the amalgamation of the the stabilization of all the slices of $\Topos$:
\vspace{-2mm}
\begin{equation}
  \label{TangentInfinityTopos}
  T \Topos
  \;\simeq\;
  \int_{\mathcal{B} \in \Topos}
  T_{\mathcal{B}}\Topos
  \,,
  \;\;\;\;\;\;\;\;\;\;\;\;
  \mbox{with}
  \;\;
  \begin{tikzcd}[column sep=large]
    T_{\mathcal{B}} \Topos
    \:=\;
    \Spectra
    \big(
      \SliceTopos{\mathcal{B}}
    \big)
    \ar[
      from=r, shift right=6pt,
      "{
        \Sigma_{\mathcal{B}}^\infty
      }"{swap}
    ]
    \ar[
      r, shift right=6pt,
      "{
        \Omega_{\mathcal{B}}^\infty
      }"{swap}
    ]
    \ar[
      r,
      phantom,
      "{\scalebox{.7}{$\bot$}}"
    ]
    &
    \SliceTopos{\mathcal{B}}^{\mathcal{B}/}
    \ar[
      from=r, shift right=6pt,
      "{
        (-)_+
      }"{swap}
    ]
    \ar[
      r, shift right=6pt,
      "{
        \mathrm{undrlg}
      }"{swap}
    ]
    \ar[
      r,
      phantom,
      "{\scalebox{.7}{$\bot$}}"
    ]
    &
    \SliceTopos{\mathcal{B}}
  \end{tikzcd}
  \,.
\end{equation}

\noindent Notice that for $\mathcal{E} \,\in\, T_\ast \Topos$ a (sheaf of) ring spectra, this subsumes local coefficients over the Picard groupoid of $\mathcal{E}$, and hence the usual notion of twisting and grading of multiplicative cohomology Whitehead-generalized theories (review and further pointers in \cite[Ex. 2.010]{FSS20CharacterMap}).
This means that twisted Whitehead-generalized cohomology theories together with their sheaf-hypercohomology-type generalizations
(as in \cite{Brown73})
all naturally are examples \cite[\S 4.1.2.1]{dcct}
of the truly generalized notion of cohomology \eqref{GeneralNonabelianCohomology}.

\medskip
In fact, if $\Topos$ is {\it cohesive} (Def. \ref{CohesiveInfinityTopos})
then $T\Topos$
\eqref{TangentInfinityTopos}
is also cohesive \cite[Prop. 4.1.9]{dcct} and knows all about {\it differential} generalized cohomology theory (\cite[\S 4.1.2.2]{dcct}\cite{BunkeNikolausVoelkl13}\cite{ADH21}), in that it provides any coefficient object $\mathcal{E} \in T_\ast \Topos$ with compatible notions of curvature forms, Chern characters, and flat coefficients.

\medskip
Last but not least, we have seen (Prop. \ref{HomotopyQuotientsAndPrincipaInfinityBundles})
that equivariance is just another instance of twisting by slicing
\eqref{TwistedgeneralizedCohomology}:
\vspace{-2mm}
\begin{equation}
  \label{SliceToposOverDeloopingAsContextForGeneralizedCohomology}
  \begin{tikzcd}[row sep=-3pt, column sep=4pt]
    \Actions{\mathcal{G}}(\Topos)
    \ar[rr, phantom,"{\simeq}"]
    &&
    \SliceTopos{\mathbf{B}\mathcal{G}}
    \\
    \scalebox{0.8}{$\mathcal{G}\acts \, \mathcal{X}$}
    &\longmapsto&
    \scalebox{0.8}{$ \HomotopyQuotient{\mathcal{X}}{\mathcal{G}}$}
    \,.
  \end{tikzcd}
\end{equation}

However -- and this is the major subtlety to be dealt with in the conceptual foundations of equivariant generalized cohomology  --
while the intrinsic generalized cohomology of $\SliceTopos{\mathbf{B}\mathcal{G}}$ itself is of some interest (for instance it subsumes all notions of group cohomology)
the traditional notion of {\it proper}- (as opposed to Borel-) equivariant cohomology is concerned not with the (slice) hom-$\infty$-groupoids of \eqref{SliceToposOverDeloopingAsContextForGeneralizedCohomology}, but with the {\it shape} of its {\it slice mapping stacks} (cf. Ex. \ref{EquivariantMappingStackAsSliceMappingStack}):
\vspace{-2mm}
\begin{equation}
  \label{TakingShapeIntoEquivariantMappingStack}
  \begin{tikzcd}
  \overset{
    \mathclap{
    \raisebox{2pt}{
      \tiny
      \color{darkblue}
      \bf
      shape of equivariant mapping stack
    }
    }
  }{
  \shape
  \;\SliceMaps{\big}{\mathbf{B}\mathcal{G}}
    {
      \HomotopyQuotient
        { \mathcal{X} }
        { \mathcal{G} }
    }
    {
      \HomotopyQuotient
        { \mathcal{A} }
        { \mathcal{G} }
    }
    }
    \ar[
      r,
      "{\mbox{\tiny \eqref{ComparisonMorphismFromShapeOfMappingStackToMappingSpaceOfShapes}}}"
    ]
    &
    \overset{
      \mathclap{
      \raisebox{2pt}{
        \tiny
        \color{darkblue}
        \bf
        equivariant mapping space of shapes
      }
      }
    }{
    \SliceMaps{\Big}{\mathbf{B}\shape\mathcal{G}}
    {
      \HomotopyQuotient
        { \shape\mathcal{X} }
        { \shape\mathcal{G} }
    \,}
    {\,
      \HomotopyQuotient
        { \shape\mathcal{A} }
        { \shape\mathcal{G} }
    }
    }
    \;\simeq\;
    \Discrete
    \,
    \overset{
    }{
    \SlicePointsMaps{\Big}{\mathbf{B}\shape\mathcal{G}}
    {
      \HomotopyQuotient
        { \shape\mathcal{X} }
        { \shape\mathcal{G} }
    \,}
    {\,
      \HomotopyQuotient
        { \shape\mathcal{A} }
        { \shape\mathcal{G} }
    }
    }
    \,.
  \end{tikzcd}
\end{equation}
\vspace{-2mm}

\noindent
Here the left hand side here reduces to a traditional definitions of proper {\it equivariant cohomology} (such as of equivariant K-theory, Ex. \ref{TwistedEquivariantKTheory} below):
\vspace{-1mm}
\begin{equation}
  \label{ApparentFormOfGeneralEquivariantHomotopyTheory}
  H^0_G(-; \mathrm{A})
  \,\coloneqq\,
  \pi_0 \, \shape \,
  \SliceMaps{\big}{\mathbf{B}G}
    {
      \HomotopyQuotient
        {-}
        {G}
    }
    {
      \HomotopyQuotient
        {\mathrm{A}}
        {G}
    }
  \;\; : \;\;
  \HomotopyCategory
  \big(
    \SliceTopos{\mathbf{B}G}
  \big)^{\mathrm{op}}
  \xrightarrow{\phantom{---}}
  \Sets
  \,,
\end{equation}
while the right-hand side captures only its coarse form of {\it Borel-equivariant cohomology}.
This brings out

\noindent {\bf (i)} that proper equivariant cohomology theory is intrinsically a concept of {\it cohesive} homotopy theory and

\noindent {\bf (ii)} that in order for its definition \eqref{ApparentFormOfGeneralEquivariantHomotopyTheory} to qualify as an instance of truly generalized cohomology \eqref{GeneralNonabelianCohomology} there must be {\it some} kind of equivariant Oka principle (along the lines indicated in \cref{OverviewAndSummary}) which allows to take the shape modality from outside to inside a mapping stack construction.

\medskip
The orbi-smooth Oka principle which we have proven above (Thm. \ref{OrbiSmoothOkaPrinciple}) concerns special cases when the comparison map \eqref{TakingShapeIntoEquivariantMappingStack} restricts to an equivalence. However, if we are content with the right-hand side being {\it any} generalized cohomology theory \eqref{GeneralNonabelianCohomology} -- not necessarily that in a slice of the given $\Topos$ -- then we can ask for a
more general equivariant Oka-like principle which
identifies the shape of an equivariant mapping stack with the hom-$\infty$-groupoid of the shape of its arguments regarded in some {\it other} $\infty$-topos.

\medskip
We highlight the perspective that this is what the classical Elmendorf-Dwyer-Kan theorem really accomplishes (recalled as Prop. \ref{ElmendorfDwyerKanTheorem} below) by identifying (under mild conditions) the cohesive definition \eqref{ApparentFormOfGeneralEquivariantHomotopyTheory}
of proper equivariant cohomology
with the intrinsic generalized cohomology \eqref{GeneralNonabelianCohomology} of the $\infty$-topos of $\infty$-presheaves over the corresponding orbit category -- hence (in the case of discrete equivariance groups, Ex. \ref{CohesiveFormulationOfEDKTheoremForDiscreteEquivarianceGroups} below) with the intrinsic cohomology of the $G$-equivariant $\infty$-topos $\GEquivariant{\Topos}$ (Def. \ref{GEquivariantAndGloballyEquivariantHomotopyTheories}).

\medskip

But this means, in view of \eqref{TwistedgeneralizedCohomology} that we have identified the correct general definition of twisted \& $G$-equivariant generalized cohomology: This is the intrinsic cohomology \eqref{GeneralNonabelianCohomology}
of slices of $\GEquivariant(\Topos)$:
\begin{equation}
  \label{GeneralDefinitionOfTwistedGEquivariantCohomology}
  \overset{
    \mathclap{
    \raisebox{3pt}{
      \tiny
      \color{darkblue}
      \bf
      twisted $G$-equivariant generalized cohomology
    }
    }
  }{
  H^\bullet_G
  (
    -
    ;\,
    \mathcal{A}
  )
  \;:=\;
  H^0
  \big(
    (-,\bullet)
    ;\,
    (\mathcal{A}, p_{\mathcal{A}})
  \big)
  }
  \;:\;
  \HomotopyCategory
  \big(
    \Slice{
      \GEquivariant{\Topos}
    }
    {\mathcal{B}}
  \big)^{\mathrm{op}}
  \xrightarrow{\phantom{---}}
  \Sets
\end{equation}
(for given moduli of equivariant twistings $\mathcal{B} \,\in\, \Slice{\GEquivariant(\Topos)}{\mathcal{B}}$).

\medskip
For this definition to be practically useful and to reproduce existing definitions (notably of twisted equivariant K-theory, Ex. \ref{TwistedEquivariantKTheory} below) we just need to verify one consistency condition: The equivariant principal bundles of \cref{InTopologicalSpaces} should serve, under suitable conditions, as twists in \eqref{GeneralDefinitionOfTwistedGEquivariantCohomology}, in that twisted equivariant cohomology sets in this abstract sets are given by concordance classes of equivariant sections
(as in Prop. \ref{EquivariantSectionsOfAssociatedEquivariantBundles})
of their associated equivariant fiber bundles. That this is indeed the case, under suitable conditions (``admissible twists''. Def. \ref{AdmissibleEquivariantTwisting}) is the statement of the {\it twisted EKD theorem}, which we prove as Thm. \ref{TwistedEDKtheorem} below, as an easy consequence of a couple of results in singular-cohesive homotopy theory that we have established above.

\medskip

\noindent
{\bf The twisted Elemendorf theorem.}
The  classical Elmendorf-Dywer-Kan theorem (recalled in slight singular-cohesive re-formulation as Prop. \ref{ElmendorfDwyerKanTheorem} below), expresses the shape of equivariant mapping stacks between  topological spaces equipped with topological group actions as the hom-$\infty$-groupoid between the shapes of their systems of fixed loci (which are $\infty$-presheaves over the category of orbits of the equivariance group).

\medskip
Below we use this classical EDK theorem and the singular-cohesive homotopy theory that we have developed above in order to prove a twisted/parameterized/sliced generalization (Thm. \ref{TwistedEDKtheorem}) which identifies the shape of spaces of equivariant {\it sections} of suitable equivariant fiber bundles with the {\it slice} hom-$\infty$-groupoid of systems of fixed loci slice over the respective proper equivariant classifying shape (Def. \ref{EquivariantClassifyingStack}).

\medskip

\begin{proposition}[Elmendorf-Dwyer-Kan theorem {\cite{DwyerKan84}}]
 \label{ElmendorfDwyerKanTheorem}
  For

  - $\Gamma \,\in\, \Groups(\kTopologicalSpaces)$ a topological group
  {\rm (no further conditions whatsoever)},

  - $\Gamma \acts \, \TopologicalPrincipalBundle \,\in\, \Gamma\CWComplexes \xhookrightarrow{\;\;} \Actions{\Gamma}(\kTopologicalSpaces)$ a  $\Gamma$-CW complex with finite stabilizer groups,

  - $\Gamma \acts \, \TopologicalCoefficients \,\in\, \Actions{\Gamma}(\kTopologicalSpaces)$ and topological $\Gamma$-space,

\noindent
the operation of passing to shapes of fixed loci
(Def. \ref{SystemsOfShapesOfFixedLociInTopologicalGammaSpaces}) constitutes a natural equivalence of $\infty$-groupoids
\begin{equation}
  \label{TheEDKEquivalenceOnDerivedHomSpaces}
  \begin{tikzcd}[row sep=3pt, column sep=10pt]
  \shape
  \Big(
  \Maps{\big}
    { \Gamma \acts \, \TopologicalPrincipalBundle }
    { \Gamma \acts \, \TopologicalCoefficients }
   ^{\Gamma}
  \Big)
  \ar[
    rr,
    "{\sim}"
  ]
  &&
  \def\AmbientCategory{
    \InfinityPresheaves
    \big(
      \Orbits(\Gamma)
    \big)
  }
  \PointsMaps{\Big}
    {
      \shape
      \,
      \FixedLoci(\TopologicalPrincipalBundle)
    }
    {
      \shape
      \,
      \FixedLoci(\TopologicalCoefficients)
    }
  \end{tikzcd}
\end{equation}
from the shape of the equivariant mapping stack to the hom-$\infty$-groupoid of $\infty$-presheaves over the category of proper $\Gamma$-orbits (Def. \ref{ProperTopologicalOrbitInfinityCategory}, hence here considered for finite stabilizer groups)
\end{proposition}
\begin{proof}
  This is a slight re-formulation of the result of \cite[Thm. 3.1]{DwyerKan84} applied to the computation of derived mapping spaces (see also \cite{CordierPorter96}, review in \cite[Thm. 1.3.8]{Blu17}):

  Via
  Prop. \ref{DiffeologicalMappingSpacesHaveCorrectUnderlyingHomotopyType}
  and
  Prop. \ref{CohesiveShapeOfSmoothManifoldsIsTheirHomotopyType}, the left hand side
  of \eqref{TheEDKEquivalenceOnDerivedHomSpaces}
  is equivalently the singular simplicial complex of the equivariant mapping space. (Similarly, the shapes appearing in the arguments on the right of \eqref{TheEDKEquivalenceOnDerivedHomSpaces} are equivalently the singular simplicial complexes of the systems of fixed loci.)

  By the assumption that $\TopologicalPrincipalBundle$ is a $\Gamma$-CW complex, it is cofibrant in the simplicial model category structure on $\Gamma$-spaces considered in \cite[\S 1.2 \& Thm. 2.2]{DwyerKan84}. Since every object is fibrant in this model structure, it follows by standard facts of model category theory that the left-hand side of \eqref{TheEDKEquivalenceOnDerivedHomSpaces} is equivalently the derived hom-complex -- i.e. the hom-$\infty$-groupoid -- in this model structure. Now the Quillen equivalence of
  \cite[Thm. 3.1]{DwyerKan84} identifies this hom-$\infty$-groupoid with that (between the singular simplicial complexes of the systems of fixed loci) in the projective model structure on simplicial presheaves (our Ntn. \ref{ModelCategoriesOfSimplicialPresheaves}, corresponding to \cite[\S 1.3]{DwyerKan84}) which
  (by
  Def. \ref{ProperTopologicalOrbitInfinityCategory}
  and
  Prop. \ref{PresentationOfInfinityToposesBySimplicialPresheaves}, Ntn. \ref{PresentationOfInfinityPresheavesBySimplicialPresheaves}) presents the right hand side of
  \eqref{TheEDKEquivalenceOnDerivedHomSpaces}.
\end{proof}

\medskip

\begin{example}[Elmendorf-Dwyer-Kan theorem for discrete equivariance groups in terms of singular-cohesive homotopy theory]
\label{CohesiveFormulationOfEDKTheoremForDiscreteEquivarianceGroups}
If $\Gamma = G \,\in\, \Groups(\Sets) \xhookrightarrow{\;} \Groups(\kTopologicalSpaces)$ is a discrete group then with the notation and equivalences of \cref{GeneralSingularCohesion},
the Elmendorf-Dwyer-Kan theorem
(Thm. \ref{ElmendorfDwyerKanTheorem})
reads equivalently as follows (under the assumption that $\TopologicalSpace$ is a $G$-CW-complex, Ex. \ref{GCWComplexesAreCofibrantObjectsInProperEquivariantModelcategory}):
\vspace{-2mm}
$$
  \begin{tikzcd}[row sep=3pt, column sep=10pt]
  \shape
  \Big(
  \SliceMaps{\big}{\mathbf{B}G}
    {
      \HomotopyQuotient
        { \TopologicalSpace }
        { G }
    }
    {
      \HomotopyQuotient
        { \TopologicalCoefficients }
        { G }
    }
  \Big)
  \ar[
    rr,
    "{\sim}"
  ]
  &&
  \smooth
  \SliceMaps{\big}{\orbisingularG}
    {
      \shape
      \orbisingular
      (
      \HomotopyQuotient
        { \TopologicalSpace }
        { G }
      )
    }
    {
      \shape
      \orbisingular
      (
      \HomotopyQuotient
        { \TopologicalCoefficients }
        { G }
      )
    }
    \mathrlap{\,.}
  \end{tikzcd}
$$

\vspace{-2mm}
\noindent This is the result of the following sequence of natural equivalences of $\infty$-groupoids:

\bigskip
$$
  \def\arraystretch{1.8}
  \begin{array}{ll}
  \shape
  \Big(
  \SliceMaps{\big}{\mathbf{B}G}
    {
      \HomotopyQuotient
        { \TopologicalSpace }
        { G }
    }
    {
      \HomotopyQuotient
        { \TopologicalCoefficients }
        { G }
    }
  \Big)
  \\
  \;\simeq\;
  \def\AmbientCategory{
    \InfinityPresheaves
    \big(
      \Orbits(G)
    \big)
  }
  \PointsMaps{\Big}
    {
      \shape
      \,
      \FixedLoci(\TopologicalPrincipalBundle)
    }
    {
      \shape
      \,
      \FixedLoci(\TopologicalCoefficients)
    }
    &
    \proofstep{
      by
      \eqref{TheEDKEquivalenceOnDerivedHomSpaces}
    }
  \\
  \;=\;
  \def\AmbientCategory{
    \GEquivariant\InfinityGroupoids
  }
  \PointsMaps{\Big}
    {
      \shape
      \,
      \FixedLoci(\TopologicalPrincipalBundle)
    }
    {
      \shape
      \,
      \FixedLoci(\TopologicalCoefficients)
    }
    &
    \proofstep{
      by
      Def. \ref{GEquivariantAndGloballyEquivariantHomotopyTheories},
      Ex.
      \ref{ClassicalEquivariantHomotopyTheory}
    }
  \\
  \;=\;
  \def\AmbientCategory{
    \GEquivariant\InfinityGroupoids
  }
  \PointsMaps{\Big}
    {
      G\Orbi\Smooth
      \big(
      \shape
      \orbisingular
      (\HomotopyQuotient{\TopologicalCoefficients}{G})
      \big)
    }
    {
      G\Orbi\Smooth
      \big(
      \shape
      \orbisingular
      (\HomotopyQuotient{\TopologicalCoefficients}{G})
      \big)
    }
    &
    \proofstep{
      by
      \eqref{SystemOfShapesOfGeometricFixedLociAsObjectInGEquivariantHomotopyTheory}
    }
  \\
  \;=\;
  \def\AmbientCategory{
    \Slice{
    \Singular\InfinityGroupoids
    }{\orbisingularG}
  }
  \PointsMaps{\Big}
    {
      \shape
      \orbisingular
      (\HomotopyQuotient{\TopologicalCoefficients}{G})
    }
    {
      \shape
      \orbisingular
      (\HomotopyQuotient{\TopologicalCoefficients}{G})
    }
    &
    \proofstep{
      by
      \eqref{OrbiSpacesAreSmoothRelativeGSingularities}
    }
  \\
  \;=\;
  \smooth
  \SliceMaps{\big}{\orbisingularG}
    {
      \shape
      \orbisingular
      (\HomotopyQuotient{\TopologicalCoefficients}{G})
    }
    {
      \shape
      \orbisingular
      (\HomotopyQuotient{\TopologicalCoefficients}{G})
    }
    &
    \proofstep{
      by
      \eqref{GlobalPointsOfSliceMappingStackIsSliceHom}
      .
    }
  \end{array}
$$
\end{example}

\begin{theorem}[Twisted EDK-theorem]
  \label{TwistedEDKtheorem}
  Let

  -- $G \,\in\, \Groups(\FiniteSets)$,

  -- $G \acts \, \Gamma \,\in\, \Groups\big(\Actions{G}(\DTopologicalSpaces)\big)$,
   such that $\mathbf{B}(\Gamma \rtimes G)$
   is admissible $G$-equivariant twisting
   (Def. \ref{AdmissibleEquivariantTwisting}),

 \smallskip
  -- $G \acts \, \TopologicalSpace \,\in\,
  \GCWComplexes
  \hookrightarrow
  \Actions{G}(\TopologicalSpaces)$
   a $G$-CW-complex (Exp. \ref{GCWComplexesAreCofibrantObjectsInProperEquivariantModelcategory}),

\vspace{1mm}
  -- $G \acts \, \TopologicalPrincipalBundle
  \,\in\, \EquivariantPrincipalFiberBundles{G}{\Gamma}_{\TopologicalSpace}$
  with modulating morphism
  $
   \HomotopyQuotient{\Twist}{G}
   \,\in\,
   \SlicePointsMaps{}{\mathbf{B}G}
    {
      \HomotopyQuotient
        { \TopologicalSpace }
        { G }
    }
    {
      \HomotopyQuotient
        { \mathbf{B}\Gamma }
        { G }
    }
  $,

\vspace{1mm}
  -- $(\Gamma \rtimes G) \acts \,
    \TopologicalCoefficients
    \,\in\, \Actions{G\acts \, \Gamma}
    \big(\Actions{G}(\TopologicalSpaces)\big)$.
\vspace{1mm}

\noindent
Then, for $\TopologicalFiberBundle$ denoting the associated equivariant $A$-fiber bundle (Def. \ref{EquivariantAssociatedInfinityBundles}),
we have a natural equivalence of
$\infty$-groupoids
$$
  \def\arraystretch{2}
  \begin{array}{ll}
  \overset{
    \mathclap{
    \raisebox{3pt}{
      \tiny
      \color{darkblue}
      \bf
      \begin{tabular}{c}
        shape of stack of
        $G$-equivariant sections
        of equivariant $A$-fiber bundle
      \end{tabular}
    }
    }
  }{
  \shape
  \bigg(
  \Big(
  \SliceMaps{\big}{G \acts \TopologicalSpace}
    {
      G \acts \, \TopologicalSpace
    }
    {
      G \acts \, \TopologicalFiberBundle
    }
   \Big)
   ^{G}
   \bigg)
  }
  \qquad
  &
  \qquad
  \overset{
    \mathclap{
    \raisebox{2pt}{
      \tiny
      Prop. \ref{EquivariantSectionsOfAssociatedEquivariantBundles}
    }
    }
  }{
    \;\simeq\;
  }
  \qquad
  \underset{
    \mathclap{
    \raisebox{0pt}{
      \tiny
      \color{darkblue}
      \bf
        shape of slice mapping stack
        of homotopy quotients
        from base to typical fiber
    }
    }
  }{
  \shape
  \,
  \SliceMaps{\big}{\mathbf{B}(\Gamma \rtimes G)}
    {
      \HomotopyQuotient
        { (\TopologicalSpace,\Twist) }
        { G }
    }
    {
      \HomotopyQuotient
        {
          \TopologicalCoefficients
        }
        {
          (\Gamma \rtimes G)
        }
    }
    }
    \\
    &
    \qquad
    \overset{
      \mathclap{
      \raisebox{5pt}{
        \tiny
        \color{greenii}
        \bf
        twisted EDK
      }
      }
    }{
    \;\;
    \simeq
    \;\;\;
    }
    \;\;
    \overset{
    \mathclap{
    \raisebox{3pt}{
      \tiny
      \color{darkblue}
      \bf
      \begin{tabular}{c}
        slice mapping space of
        shapes of orbi-singularizations
      \end{tabular}
    }
    }
    }{
    \smooth
    \,
    \SliceMaps{\Big}{
      \shape
      \orbisingular
      \mathbf{B}(\Gamma \rtimes G)
    }
     {
       \shape
       \orbisingular
       \big(
       \HomotopyQuotient
         {
           (\TopologicalSpace, \Twist)
         }
         { G }
       \big)
     }
     {
       \shape
       \orbisingular
       \big(
       \HomotopyQuotient
         { \mathrm{A} }
         { (\Gamma \rtimes G) }
       \big)
     }
     }
     \\
     &
   \;\; \qquad\simeq \quad
    \def\AmbientCategory{\GEquivariant\InfinityGroupoids}
    \underset{
    \mathclap{
    \raisebox{20pt}{
      \tiny
      \color{darkblue}
      \bf
        slice hom-$\infty$-groupoid of
        proper $G$-equivariant homotopy types
        of equivariant classifying shape
    }
    }
    }{
    \SlicePointsMaps{\Big}
      {\EquivariantClassifyingShape{G}{\Gamma}}
      {
        \shape
        \,
        \FixedLoci(\TopologicalSpace, \mathrm{tw})
      }
      {
        G\Orbi\Smooth
        \,
        \shape
        \,
        \orbisingular
        \big(
          \HomotopyQuotient
            { \TopologicalCoefficients }
            { (\Gamma \rtimes G) }
        \big)
      }
    }
  \end{array}
$$

\vspace{-5mm}
\noindent
between  the shape of the stack of equivariant sections of the fiber bundle and the slice hom-$\infty$-groupoid between the systems of $G$-fixed loci, sliced over the equivariant classifying shape $\EquivariantClassifyingShape{G}{\Gamma}$.
\end{theorem}

\begin{proof} This may be obtained as the composite of the following sequence of natural equivalences:
\vspace{-2mm}
$$
  \def\arraystretch{1.8}
  \begin{array}{ll}
    \shape
    \,
    \SliceMaps{\big}{\mathbf{B}(\Gamma \rtimes G)}
      {
        \HomotopyQuotient
          { (\TopologicalSpace, \Twist) }
          { G }
      }
      {
        \HomotopyQuotient
          { \TopologicalCoefficients }
          { \Gamma \rtimes G }
      }
    \\
    \;\simeq\;
    \shape
    \,
    \SliceMaps{\big}{\mathbf{B}(\Gamma \rtimes G)}
      {
        \HomotopyQuotient
          { \TopologicalPrincipalBundle }
          { \Gamma \rtimes G }
      }
      {
        \HomotopyQuotient
          { \TopologicalCoefficients }
          { \Gamma \rtimes G }
      }
    &
    \proofstep{
      by
      Ex. \ref{SemidirectProductHomotopyQuotientOfEquivariantInfinityBundle}
    }
    \\
    \;\simeq\;
    \shape
    \Big(
    \Maps{\big}
      {
        (\Gamma \rtimes G)
        \acts \,
        \TopologicalPrincipalBundle
      }
      {
        (\Gamma \rtimes G)
        \acts \,
        \TopologicalCoefficients
      }
    ^{(\Gamma \rtimes G)}
    \,
    \Big)
    &
    \proofstep{
      by
      Prop. \ref{HomotopyQuotientsAndPrincipaInfinityBundles}
    }
    \\
    \;\simeq\;
    \def\AmbientCategory{
      \InfinityPresheaves
       \big(
        \Orbits(\Gamma \rtimes G)
       \big)
   }
   \PointsMaps{\big}
     { \shape \FixedLoci(\TopologicalPrincipalBundle) }
     { \shape \FixedLoci(\TopologicalCoefficients) }
    &
    \proofstep{
      by
      Prop. \ref{ElmendorfDwyerKanTheorem}
    }
    \\
    \;\simeq\;
    \def\AmbientCategory{
      \InfinityPresheaves
       \big(
        \Slice{\Orbits(G)}{\EquivariantClassifyingShape{G}{\Gamma}}
       \big)
   }
   \PointsMaps{\Big}
     { \shape \, \FixedLoci(\TopologicalSpace, \Twist) }
     {
       G\Orbi\Smooth
       \,
       \shape\orbisingular
       \big(
         \HomotopyQuotient
           { \TopologicalCoefficients }
           { (\Gamma \rtimes G) }
       \big)
     }
    &
    \proofstep{
      by
      Prop. \ref{SliceOfFiniteOrbitCategoryOverEquivariantClassifyingShape},
      Ex. \ref{PresheafHomOutOfSystemOfCrossedFixedLociOfEquivariantPrincipalBundlerestrictsToCrossedSubgroups}
    }
    \\
    \;\simeq\;
    \def\AmbientCategory{
      \InfinityPresheaves
       \big(
         \Orbits(G)
       \big)
   }
   \SlicePointsMaps{\Big}{\EquivariantClassifyingShape{G}{\Gamma}}
     { \shape \, \FixedLoci(\TopologicalSpace, \Twist) }
     {
       G\Orbi\Smooth
       \,
       \shape \orbisingular
       \big(
         \HomotopyQuotient
           { \TopologicalCoefficients }
           { (\Gamma \rtimes G) }
       \big)
     }
    &
    \proofstep{
      by
      Prop.
      \ref{FundamentalTheoremOfPresheafToposTheory}
      .
    }
  \end{array}
$$
Here the key step is the identification of the slice $\infty$-site over the equivariant classifying shape $\EquivariantClassifyingShape{G}{\Gamma}$ via Prop. \ref{SliceOfFiniteOrbitCategoryOverEquivariantClassifyingShape}, which is the step needing the assumption that $\mathbf{B}(\Gamma \rtimes G)$ is an admissible $G$-equivariant twisting.
\vspace{-.6cm}
\end{proof}

\medskip

\begin{example}[Twisted equivariant K-theory]
\label{TwistedEquivariantKTheory}
The tradtional definition of twisted $G$-equivariant complex K-theory is as homotopy classes of sections of $\FredholmOperators$-fiber bundles which are associated to stable $G$-equivariant $\PUH$-principal bundles \cite[Def. 3.2]{AtiyahSegal04}\cite[Def. 4.2]{BEJU12}. By Prop. \ref{EquivariantSectionsOfAssociatedEquivariantBundles} we may equivalently reformulate this in terms of lifts of the modulating map of the underlying equivariant principal $\PUH$-bundle, and with Thm. \ref{TwistedEDKtheorem} (which does apply to $\PUH$-twists due to Prop. \ref{ExamplesOfAdmissibleGEquivariantTwistings}) we then find that this traditional definition of twisted equivariant K-theory does qualify as a twisted equivariant cohomology theory in the general abstract sense discussed above.

Concretely,
With Ex. \ref{EquivariantBundlesServingAsGeoemtricTwistsOfEquivariantKTheory},
the class $[\tau]$
of any stable equivariant bundle with structure group $\GradedPUHPrime \rtimes \ZTwo$
is represented by one with structure group $\GradedPUH \rtimes \ZTwo$
(Thm. \ref{BlowupStabilityForADEEquivariantPUHSemidirectProductZTwoPrincipalBundles} (2)),
and since the latter has a canonical conjugation action
\eqref{ConjugationActionOfProjectiveGradedUnitaryGroup}
on the space of Fredholm operators
\eqref{TheSpaceOfFredholmOperators}
on the graded Hilbert space
$\HilbertSpace \otimes \ComplexNumbers^2 \,=\,
L^2(\RiemannianSphere{3})^{\otimes^\infty}_+ \oplus L^2(\RiemannianSphere{3})^{\otimes^\infty}_+$
\eqref{GradedVersionOfInfinitL2HilbertSpace},
we may use these equivariant principal bundles to define twisted equivariant K-theory
in the general stacky form of \cite[Rem. 2.94]{SS20OrbifoldCohomology},
namely as the homotopy-classes of sections of the
stacky $\FredholmOperators$-fiber bundle associated
(as in \cite[\S 4]{NSS12a})
to the given equivariant principal bundle $\tau$:
\begin{equation}
  \label{TwistedEquivariantKRTheory}
  \def\arraystretch{1.9}
  \begin{aligned}
  \overset{
    \mathclap{
    \raisebox{5pt}{
      \tiny
      \color{darkblue}
      \bf
      \def\arraystretch{.9}
      \begin{tabular}{c}
        twisted equivariant\\
        KR-cohomology
      \end{tabular}
    }
    }
  }{
    \mathrm{KR}^\tau_G(\SmoothManifold)
  }
  &
  \;
  :=
  \;\;
  \overset{
    \raisebox{2pt}{
      \tiny
      \color{darkblue}
      \bf
      \hspace{-.4cm}
      \def\arraystretch{.9}
      \begin{tabular}{c}
        connected components of shape of
        sliced mapping stack
      \end{tabular}
    }
  }{
  \Truncation{0}
  \,
  \shape
  \,
  \SliceMaps{\Big}{ \mathbf{B}( \GradedPUH \rtimes \ZTwo ) }
    { ( \HomotopyQuotient{ \SmoothManifold }{ G }, \tau ) }
    { \HomotopyQuotient{ \FredholmOperators  }{ (\GradedPUH \rtimes \ZTwo) } }
  }
  \\
  &
  \;\;
  =
  \;\;
  \left\{
  \hspace{-8pt}
  \begin{tikzcd}[column sep=40pt]
    &&
    \overset{
    }{
    \HomotopyQuotient
      { \FredholmOperators }
      { ( \GradedPUH \rtimes \ZTwo ) }
    }
    \ar[
      d,
      "{
        \hspace{-10pt}
        \mbox{
          \tiny
          \color{darkblue}
          \bf
          \def\arraystretch{.9}
          \begin{tabular}{c}
            stacky universal
            \\
            Fredholm bundle
          \end{tabular}
        }
      }"
    ]
    \\
    \overset{
      \raisebox{5pt}{
        \tiny
        \color{darkblue}
        \bf
        \def\arraystretch{.9}
        \begin{tabular}{c}
          Real good
          \\
          orbifold
        \end{tabular}
      }
    }{
      \HomotopyQuotient
        { \SmoothManifold }
        { G }
    }
    \ar[
      rr,
      "{ \tau }",
      "{
        \mbox{
          \tiny
          \color{greenii}
          \bf
          \def\arraystretch{.9}
          \begin{tabular}{c}
            modulation of stable Real-equivariant
            \\
            projective graded bundle
          \end{tabular}
        }
      }"{swap}
    ]
    \ar[
      urr,
      dashed,
      "{
        \mbox{
          \tiny
          \color{orangeii}
          \bf
          \def\arraystretch{.9}
          \begin{tabular}{c}
            cocycle in
            \\
            twisted equivariant Real K-theory
          \end{tabular}
        }
      }"{xshift=4pt, sloped, pos=.4}
    ]
    \ar[
      dr,
      "{
        \mbox{
          \tiny
          \color{greenii}
          \bf
          Real structure
        }
      }"{sloped, swap}
    ]
    &&
    \mathbf{B} ( \GradedPUH \rtimes \ZTwo )
    \ar[dl]
    \\
    &
    \mathbf{B} \ZTwo
  \end{tikzcd}
  \hspace{-6pt}
  \right\}_{\big/\sim_{{}_{\mathrm{cncrdnc}}}}
  \end{aligned}
\end{equation}
Notice here that the sliced mapping space in \eqref{TwistedEquivariantKRTheory}
is formed in the cohesive $\infty$-topos $\SmoothInfinityGroupoids$:

\noindent {\bf (1)}  This makes the dashed section be an equivariant section of topological bundles,
while the image of this section in plain $\InfinityGroupoids$ under the shape operation
$\shape \widetilde{\mathrm{ev}}$
\eqref{ComparisonMorphismFromShapeOfMappingStackToMappingSpaceOfShapes}
sees only the coarser Borel-equivariant K-theory.
This cohesive-homotopy theoretic construction of twisted equivariant KR-theory generalizes
a similar approach in \cite[\S 3.19]{Pavlov14} to the twisted and Real case.

\noindent {\bf (2)} In particular, the bundle of Fredholm operators on the right of
\eqref{TwistedEquivariantKRTheory} is a stacky universal associated bundle
over the moduli stack
$\mathbf{B} (\GradedPUH \rtimes \ZTwo)$
(as for the stacky universal principal bundles in \cref{EquivariantLocalTrivializationIsImplies}),
instead of a plain universal bundle
just over the classifying space $B (\GradedPUH \rtimes \ZTwo)$
as considered in \cite[\S 15]{LueckUribe14}
(whose local triviality had remained open, see \cite[p. 2]{EspinozaUribe15}).

\end{example}

\medskip

By way of outlook, we briefly highlight how the
Real-equivariant projective graded bundles from Ex. \ref{EquivariantBundlesServingAsGeoemtricTwistsOfEquivariantKTheory}
neatly encode the structure that has been argued to capture:

-- quantum symmetries of gapped topological quantum systems (Rem. \ref{RealEquivariantProjectiveGradedBundlesAsQuantumSymmetries} below),

-- B-field configurations on orbi-orientifold stacetimes in string theory
 (Rem. \ref{GeometricTwistsOfEquivariantKTheoryAsBFieldsOnOrbiOrientifolds} below).

\noindent
While both of these applications are meanwhile classical and have been widely discussed
from a multitude of angles,
it seems that the nicely unifying
perspective of equivariant principal bundles
and for the full projective graded conjugation-equivariant structure group
$\GradedPUH \rtimes \ZTwo$ \eqref{ProjectiveGradedExtensionOfZTwoTimesZTwo}
has not been fully brought out before:

\begin{remark}[Real-equivariant projective graded bundles as quantum symmetries]
  \label{RealEquivariantProjectiveGradedBundlesAsQuantumSymmetries}
  Consider the data encoded in the
  geometric $G$-fixed locus \eqref{HFixedLociInEquivariantModuliStack}
  of the equivariant moduli stack \eqref{OrbiSingularModuliStack}
  of the Real-equivariant projective graded unitary bundles in
  \eqref{ClassificationStatementForStableEquivariantProjectiveGradedUnitaryBundles}.
  Postcomposition with the homomorphism
  $\GradedPUH \rtimes \ZTwo \xrightarrow{ (c, \mathrm{pr}_2) }  \ZTwo \times \ZTwo$
  \eqref{ProjectiveGradedExtensionOfZTwoTimesZTwo}
  projects this onto the groupoid of crossed homomorphisms
  (Prop. \ref{ConjugationGroupoidOfCrossedHomomorphismsIsSectionsOfDeloopedSemidirectProductProjection})
  from the equivariance group to $\ZTwo^2$:
  $$
    \begin{tikzcd}
      G\FixedLoci\big( G \acts \, \mathbf{B}\GradedPUH \big)
      \;\;\;
      \underset{
        \mathclap{
          \mbox{
            \tiny
            Lem. \ref{FixedLociOfEquivariantClassifyingStacks}
        }
        }
      }{
        \simeq
      }
      \;\;\;\;
      \SliceMaps{\big}{\mathbf{B} \ZTwo}
        { \mathbf{B}G }
        { \mathbf{B} ( \GradedPUH \rtimes \ZTwo ) }
      \ar[rr, "{ (c,\, ct)_\ast }"]
      &&
      \SliceMaps{\big}{\mathbf{B} \ZTwo}
        { \mathbf{B}G }
        { \mathbf{B} ( \ZTwo \times \ZTwo ) }
      \,.
    \end{tikzcd}
  $$
  This natural structure reflects exactly the data
  of ``extended QM symmetry classes'' according to
  \cite[Def. 3.7]{FreedMoore13} (review in \cite[pp. 88]{Thiang18}\cite[\S 2.2]{SS22AnyonicTopologicalOrder}):
  For $g \in G$ the elements
  $( c(g), ct(g) )$ in target group $\ZTwo \times \ZTwo$ (denoted ``$\mathcal{C}$''
  in \cite[(6.1)]{FreedMoore13}) are interpreted as the {\it c}harge-conjugation
  and its product with {\it t}ime-reversal symmetry, respectively, and
  their pre-image on the left is a projective representation
  (Ex. \ref{ProjectiveRepresentationsAndTheirCentralExtensions})
  of the quantum symmetry,
  including the datum of a $ct =: \phi$-{\it twisted extension}
  (\cite[Def. 1.7]{FreedMoore13})
  $G^\tau$ of $G$:
$$
  \begin{tikzcd}[column sep=huge]
    \CircleGroup
    \ar[r]
    \ar[d]
    \ar[dr, phantom, "{ \mbox{\tiny \rm (pb)} }"]
    &
    \quad
    \overset{
      \mathclap{
      \raisebox{7pt}{
        \tiny
        \color{darkblue}
        \bf
        \def\arraystretch{.9}
        \begin{tabular}{c}
          \color{black}
          \normalfont
          (delooping of)
          \\
          $ct =: \phi$-twisted extension of
          \\
          quantum symmetry group
        \end{tabular}
      }
      }
    }{
      \mathbf{B} G^\tau
    }
    \quad
    \ar[
      rr,
      "{
        \mbox{
          \tiny
          \color{greenii}
          \bf
          linear realization of
        }
      }"{yshift=-1pt},
      "{
        \mbox{
          \tiny
          \color{greenii}
          \bf
          quantum symmetry group
        }
      }"{swap, yshift=+1pt}
    ]
    \ar[d]
    \ar[drr, phantom, "{ \mbox{\tiny \rm (pb)} }"]
    &&
    \mathbf{B}
    ( \GradedUH \rtimes \ZTwo )
    \ar[d]
    \\
    \ast
    \ar[r]
    &
    \mathbf{B}G
    \ar[
      dr,
      "{
        \mathllap{
          \mbox{
            \tiny
            \color{greenii}
            \bf
            \def\arraystretch{.9}
            \begin{tabular}{c}
              (anti-)unitarity of
              \\
              charge+time reversal
            \end{tabular}
          }
          \!\!\!\!\!\!
        }
        \mathbf{B} ct
      }"{swap, xshift=3pt, yshift=3pt}
    ]
    \ar[
      rr,
      "{
        \tau
      }",
      "{
        \mbox{
          \tiny
          \color{orangeii}
          \bf
          \def\arraystretch{.9}
          \begin{tabular}{c}
            moduli of real-equivariant
            \\
            projective graded bundle
          \end{tabular}
        }
      }"{swap, yshift=1pt}
    ]
    &&
    \mathbf{B} ( \GradedPUH \rtimes \ZTwo )
    \ar[dl, "{ \mathbf{B} \mathrm{pr}_2 }"{xshift=-1pt, yshift=1pt}]
    \ar[
      dr,
      "{
        \mathbf{B} c
        \mathrlap{
          \!\!\!\!\!\!\!\!
          \mbox{
            \tiny
            \color{greenii}
            \bf
            \def\arraystretch{.9}
            \begin{tabular}{c}
              sign of
              \\
              charge conjugation
            \end{tabular}
          }
        }
      }"{xshift=-9pt, yshift=1pt, pos=.65}
    ]
    \\
    &
    &
    \mathbf{B} \ZTwo
    &&[-30pt]
    \mathbf{B}\ZTwo
  \end{tikzcd}
$$
\end{remark}
This and its relation to the following final remark is explained in more detail in \cite[\S 2 ]{SS22AnyonicTopologicalOrder}.

\begin{remark}[Equivariant $\ZTwo \acts \, \GradedPUH$-principal bundles as $B$-fields on orbi-orientifolds]
  \label{GeometricTwistsOfEquivariantKTheoryAsBFieldsOnOrbiOrientifolds}
  $\,$

  \noindent {\bf (i)}  A famous hypothesis asserts that
  the massless RR/NS-fields of type II string theory are classified by
  some form of geometrically twisted K-theory
  (see \cite[\S 1]{BSS18} for review and pointers and \cite{GS-RR} for recent constructions),
  specifically that
  geometrically twisted equivariant $\mathrm{KR}$-theory classifies the
  massless RR/NS-fields on orbi-orientifolds
  \cite{DFM1}\cite{DFM2}\cite{DMDR1}\cite{DMDR2}\cite{HMSV16},
  hence,
  in particular, that the geometric twists
  classified by our Ex. \ref{EquivariantBundlesServingAsGeoemtricTwistsOfEquivariantKTheory},
  reflect gauge-equivalence classes of the NS $B$-field in string theory on orbi-orientifolds.

  \noindent {\bf (ii)}
  In fact, this traditional conjecture has been argued to require corrections
  (\cite[\S 4.5.2, 4.6.5]{dBDHKMMS02}).
  Elsewhere we have shown that such corrections are
  plausibly captured by replacing twisted equivariant K-theory by
  twisted equivariant Cohomotopy
  (\cite{SS19TadpoleCancellation}\cite{BSS19} \cite[Rem. 5.30]{SS20OrbifoldCohomology}).
  The general nature of Thm. \ref{ProperClassificationOfEquivariantBundlesForResolvableSingularitiesAndEquivariantStructure}
  provides the proper theoretical backdrop for analyzing such variant definitions
  and their comparison. We hope to further discuss this elsewhere.
\end{remark}

\newpage


\end{document}